\font\fbi=cmb24
\font\fbh=cmb14
\font\fbg=cmb13
\font\fbf=cmb12
\font\fbe=cmb11
\font\fbd=cmb10
\font\fbc=cmb9
\font\fbb=cmb8
\font\fba=cmb7
\font\fri=cmr16
\font\frh=cmr14
\font\frg=cmr13
\font\frf=cmr12
\font\fre=cmr11
\font\frd=cmr10
\font\frc=cmr9
\font\frb=cmr8
\font\fra=cmr7
\font\fsg=cmsl13
\font\fsf=cmsl12
\font\fse=cmsl11
\font\fsd=cmsl10
\font\fsc=cmsl9
\font\fsb=cmsl8
\font\fsa=cmsl7
\font\fii=cmti16
\font\fih=cmti14
\font\fig=cmti13
\font\fif=cmti12
\font\fie=cmti11
\font\fid=cmti10
\font\fic=cmti9
\font\fib=cmti8
\font\fia=cmti7
\hsize 6.5 in
\vsize 8.9 in
{
\nopagenumbers
\voffset 0.5 in
\centerline{\fbi The Rational Cuboid Table }
\vskip 12 pt
\centerline{\fih - of -}
\vskip 12 pt
\baselineskip 16 pt
\centerline{\fri Maurice Borisovich Kraitchik}
\centerline{\frh (1882-1957)}
\vskip 32 pt
{
\noindent
\hskip 2.25 in
\vrule height 0 pt depth 0.1 pt width 2 in
\par
}

{
\midinsert
\vskip 36 pt
\noindent
\baselineskip 13 pt
\narrower\narrower\narrower\narrower
\fsg In honor of Maurice Kraitchik and his investigations
of rational body cuboids, this table is dedicated to both complete
and extend his original cuboid tables.
\par
\bigskip
\bigskip
\noindent
All 12,517 body cuboids with an odd side $< 2^{32}$ are given.
\endinsert
}

\vskip 3 in
\baselineskip 16 pt
\centerline{\fih Reproduced by}
\vskip 16 pt
\centerline{\frh Randall L. Rathbun}
\vskip 4 pt
\centerline{\frf 22-Nov-2001}

\vfill
\hrule width 3.08 in
\vglue 2 pt
\noindent
\eject
}
{
\pageno=-1
\voffset 0.0in
\parindent = 10 pt
\parskip 8 pt
\baselineskip 13 pt

\headline={\noindent
\underbar{\vbox to 10 pt{
\hbox to 469.75499 pt {\fbf Introduction to the Rational Cuboid Table -- \hfil \fsd of Maurice Kraitchik\frc}
\vfill}}}

\centerline{--- \hskip 5 pt \fbf History\fre\  \hskip 5 pt ---} \vskip 5 pt

\noindent\fre
In 1942 Maurice Kraitchik first introduced the rational cuboid problem,
providing a brief sketch with a parametric solution and nine other numerical
solutions[\footnote{\frc 1\fre}{\frc
Maurice Kraitchik, Question 14, Chapter 4, Arithmetico-Geometrical
Questions, \fic Mathematical Recreations\frc, 2nd rev. ed., Dover Publications,
NY, 1953, pp 107-108 }].

Three years later, in greater depth, he introduces these solutions, named
\fbe rational cuboids\fre, studies the properties of primitive cuboids, and
lists 50 solutions not derivable from the parametric
solutions[\footnote{\frc 2\fre}{\frc
ibid, ``On Certain Rational Cuboids'', \fic Scripta Mathematica\frc\ vol. 11,
1945 pg 326}].

Denoting the edges of a rectangular parallelopiped by \fie x, y, z\fre\ and the
diagonals of the faces by \fie X,Y,Z\fre, a simultaneous diophantine system in
(1) is created:

\medskip\noindent
(1) \hfill $x^2 + y^2~=~Z^2~,~~y^2 + z^2 ~=~X^2~,~~ z^2 + x^2 ~=~Y^2$ \hfill \ 
\medskip

His paper then goes on to find integer solutions of (1). He notes that from
any solution in rational numbers, another integer solution can be created by
multiplying by a suitable factor, but additionally, any solution of (1) in
integers can be divided by the common divisor of \fie x, y, z\fre. He correctly
deduces that it is no restriction of (1) for \fie x, y, z\fre\ to be integers,
having no common factor. He calls such a solution primitive.

\hangafter 2
\hangindent = -3 in
If (1) is satisfied in integers, the sides \fie x, y, z\fre\ of the
parallelopiped and the diagonals \fie X, Y, Z\fre\ of its faces form a
tetrahedron (Figure 1.) containing three Pythagorean triangles (\fie x,y,Z\fre),
(\fie x,z,Y\fre), and (\fie y,z,X\fre). Kraitchik develops the properties
of these triangles and their relationships in primitive cuboids further in his
paper. He finds one edge of the cuboid to be odd, which he labels \fie z\fre,
another edge to be even, labeled \fie x\fre, and the third edge to be
``doubly'' even with respect to the first two, labeled \fie y\fre.

\moveright 3.675 in
\vbox{
\vskip -1.2 in
\hsize 3 in
\input pictex
\beginpicture
\setcoordinatesystem units <1 in,1 in>
\setplotarea x from 0 to 2.7, y from 0 to 2.15
\plotsymbolspacing=0.2 pt
\setplotsymbol ({.})
\setlinear
\plot
	1.2	2
	0	0
	1.2	1
	0	0
	2.7	1
	1.2	2
	1.2	1
	2.7	1
/
\put {$x$} at 1.95 1.06
\put {$y$} at 1.13 1.4
\put {$z$} at 0.68 0.67
\put {$X$} at 0.54 1.1
\put {$Y$} at 1.55 0.44
\put {$Z$} at 2 1.6

\endpicture
\vskip -20 pt
\centerline{\quad Figure 1.}
}

\vskip -64 pt
\hangafter -5
\hangindent = -3 in
Working further upon the problem, Kraitchik devoted book 3 of his
\fie Th\' eorie des Nombres\fre\ to rational cuboids, and gave a table of 241
primitive solutions[\footnote{\frc 3\fre}{\frc
ibid, \fic Th\'eorie des Nombres\frc, Tome 3,
Analyse Diophantine et Applications aux Cuboides Rationnels, Gauthier-Villars,
Paris 1947 pp 122-131}],
numbering them, listing the odd side \fie z\fre, and giving the factored edges
in order \fie z, x, y\fre\ as mentioned above. He also supplies three primitive
Pythagorean generator pairs which create the Pythagorean triangles. These
generator pairs are found by examining the reduced ratios of the cuboid sides,
in order: \fie z/x\fre, then \fie x/y\fre, and finally \fie z/y\fre.

In his last published table on rational cuboids, Kraitchik makes an attempt to
find all cuboids with the smallest odd edge \fie z\fre\ less than 1,000,000, and
found 18 ``new'' ones. Unfortunately he duplicated 2 cuboids from his previous
table, thus giving only a total of 257 cuboids[\footnote{\frc 4\fre}{\frc
ibid, Sur les Cuboides Rationnels, in \fic Proceedings
of the International Congress of Mathematicians\frc\ 1954 vol. 2, Amsterdam,
1954 pp 33-34}]. Actually 416 rational cuboids do occur over this range, so
159 cuboids are missing from his lists.

\vfill\eject
\bigskip
\centerline{--- \hskip 5 pt \fbf The New Table\fre\  \hskip 5 pt ---}
\medskip

\noindent\fbe
The Rational Cuboid Table\fre\ is dedicated to Maurice Borisovich Kraitchik
in tribute to his effort upon the cuboid problem. This new table is to be
considered both a completion, restoring the missing 159 cuboids to his
original lists, and a revision, extending them from all odd sides
$<10^6$ to all odd sides $<2^{32}$. Since the author worked by hand
over 2 months finding cuboid solutions for all edges less than 10,000, he is
acquainted with the extensive manual labor that Kraitchik expended. It is
hoped that this present edition will be a suitable revision and continuation
of his work.

The format of Kraitchik's third table has been followed, but there are now
36 entries on a page, instead of 26. The table is computer generated, typeset,
and checked to insure complete accuracy for all 12,517 solutions.

\medskip
\centerline{--- \hskip 5 pt \fbf Explanation of an Individual Entry\fre\  \hskip 5 pt ---}

$$\vbox{
\frc
\baselineskip 11 pt
\parindent 5 pt
\hsize 65 pt
\hbox{%
\vbox{\hskip 65 pt}\vbox{\hskip 65 pt}\vbox{\hskip 65 pt}\vbox{\hfil Pythagorean\hfil}%
}
\hbox{%
\vbox{\hfil Entry \#\hfil}\vbox{\hfil Odd Side\hfil}\vbox{\hfil Edges\hfil}\vbox{\hfil Generators\hfil}%
}
\vskip 2 pt
\hrule
\hbox{\vrule height 10 pt depth 1 pt
\vbox{\hskip 65 pt}\vrule\vbox{\hskip 65 pt}\vrule\vbox{\fbd\quad c\hfil\frc 5.17\hfil}\vrule\vbox{\fbd\quad f\frc\hfil 11\ \ 6\hfil}\vrule%
}%
\hbox{\vrule height 10 pt depth 2 pt
\vbox{\fbd\quad a\frc\hfil 1\hfil}\vrule\vbox{\fbd\quad b\frc\hfil 85\hfil}\vrule\vbox{\fbd\quad d\hfil\frc 4.3.11\hfil}\vrule\vbox{\fbd\quad g\frc\hfil\ 5\ \ 6\hfil}\vrule%
}%
\hbox{\vrule height 9 pt depth 3 pt
\vbox{\hskip 65 pt}\vrule\vbox{\hskip 65 pt}\vrule\vbox{\fbd\quad e\hfil\frc 16.9.5\hfil}\vrule\vbox{\fbd\quad h\frc\hfil\ 9\ \ 8\hfil}\vrule%
}
\hrule
\hsize 260 pt
\vskip 8 pt
\hbox{%
\vbox{\hfil\fre Table 1. -- Example of a rational cuboid entry. \hfil}%
}
}$$

\noindent
In the table above, the items have the same identical meaning as in
Kraitchik's tables for rational cuboids:
\smallskip

\halign{\hskip 0.5 in \hfil # & # & \ \  # \hfill \cr
\fbe a &---& \fre Entry number of an odd side as it occurs in ascending numerical order. \cr
\fbe b &---& \fre Value of the odd side, \fie z\fre, of the rational cuboid. \cr
\fbe c &---& \fre Factors of the odd side, \fie z\fre, in ascending primes. \cr
&&\fsd (Powers of primes are written in full expression, i.e., see entry \fbd e \fsd above.)\fre \cr
\fbe d &---& \fre Factors of the first even side, \fie x\fre, in ascending primes. \cr
\fbe e &---& \fre Factors of the multiple even side, \fie y\fre, in ascending primes. \cr
\fbe f &---& \fre Primitive Pythagorean generators from \fie c/d\fre, in odd/even format. \cr
\fbe g &---& \fre Primitive Pythagorean generators from \fie d/e\fre, in odd/even format. \cr
\fbe h &---& \fre Primitive Pythagorean generators from \fie c/e\fre, in odd/even format. \cr
}
\smallskip

For all Pythagorean generator pairs \fie p,q\fre\ --- let ${p^2-q^2} \over 2pq$
be a Pythagorean ratio. Let the ratios of the body cuboid edges for a table
entry be \fie c/d, d/e,\fre\ and \fie c/e\fre. Since
$ {c \over d} \times {d \over e} = {c \over e}$, the relationship among three
Pythagorean ratios exists:
${{ p^2 - q^2 } \over 2pq} \times {{ r^2 - s^2 } \over 2rs} =
 {{ t^2 - u^2 } \over 2tu}$
(from the 3 generator pairs \fie f,g,h\fre) which is the fundamental property
of all body cuboids.

All entries in the Classical Rational Cuboid table follow the format in the
sample table entry above, exactly as in Kraitchik's[\footnote{\frc 3\fre}{\frc
Maurice Kraitchik, \fic Th\'eorie des Nombres\frc, Tome 3,
Analyse Diophantine et Applications aux Cuboides Rationnels, Gauthier-Villars,
Paris 1947 pp 122-131}].

\vfill\eject
}
{
\headline={\noindent
\underbar{\vbox to 10 pt{
\hbox to 469.75499 pt {\fbf Errata of Tables of Rational Cuboids -- \hfil \fsd of Maurice Kraitchik\frc}
\vfill}}}

\nointerlineskip
\baselineskip = 0 pt
\vbox{
\centerline{\fbf Omissions of Rational Cuboids}
\centerline{\fsd (as numbered in the new table)}
}%
\smallskip
$$
\vbox{
\nointerlineskip
\halign{\strut
    \vrule height 9.2 pt\hfil \ \ \frc #\ \vrule height 9.2 pt\hskip 2 pt
   &\vrule height 9.2 pt\hfil \ \ \fsc #\ \ \hfil
   &\vrule height 9.2 pt\hfil \ \ \fsc #\ \ \hfil
   &\vrule height 9.2 pt\hfil \ \ \fsc #\ \ \hfil
   &\vrule height 9.2 pt\hfil \ \ \fsc #\ \ \hfil
   &\vrule height 9.2 pt\hfil \ \ \fsc #\ \ \hfil
   &\vrule height 9.2 pt\hfil \ \ \fsc #\ \ \hfil
   &\vrule height 9.2 pt\hfil \ \ \fsc #\ \ \hfil
   &\vrule height 9.2 pt\hfil \ \ \fsc #\ \ \hfil
   &\vrule height 9.2 pt\hfil \ \ \fsc #\ \ \hfil
   &\vrule height 9.2 pt\hfil \ \ \fsc #\ \ \hfil\vrule height 9.2 pt\hskip 2 pt
   &\vrule height 9.2 pt\hfil \ \ \frc #\ \vrule height 9.2 pt\cr%
\noalign{\hrule}
&0&1&2&3&4&5&6&7&8&9&Omissions \cr
\noalign{\hrule}
0&n/a&$\cdot$&$\cdot$&$\cdot$&$\cdot$&$\cdot$&$\cdot$&$\cdot$&$\cdot$&$\cdot$&none \cr
\noalign{\hrule}
10&$\cdot$&$\cdot$&$\cdot$&$\cdot$&$\cdot$&$\cdot$&$\cdot$&$\cdot$&$\cdot$&$\cdot$&none \cr
\noalign{\hrule}
20&$\cdot$&$\cdot$&$\cdot$&$\cdot$&$\cdot$&$\cdot$&$\cdot$&$\cdot$&$\cdot$&$\cdot$&none \cr
\noalign{\hrule}
30&$\cdot$&$\cdot$&$\cdot$&$\cdot$&$\cdot$&$\cdot$&$\cdot$&$\cdot$&$\cdot$&$\cdot$&none \cr
\noalign{\hrule}
40&$\cdot$&$\cdot$&$\cdot$&$\cdot$&$\cdot$&$\cdot$&$\cdot$&$\cdot$&$\cdot$&$\cdot$&none \cr
\noalign{\hrule}
50&$\cdot$&$\cdot$&$\cdot$&$\cdot$&$\star$&$\cdot$&$\cdot$&$\cdot$&$\cdot$&$\cdot$&none \cr
\noalign{\hrule}
60&$\cdot$&$\cdot$&$\cdot$&$\cdot$&$\cdot$&$\cdot$&$\cdot$&$\cdot$&$\star$&$\cdot$&none \cr
\noalign{\hrule}
70&$\cdot$&$\cdot$&$\star$&$\cdot$&$\cdot$&$\cdot$&$\cdot$&$\cdot$&$\cdot$&$\cdot$&none \cr
\noalign{\hrule}
80&$\cdot$&$\cdot$&$\cdot$&$\star$&$\star$&$\star$& 86&$\star$&$\cdot$&$\cdot$&1 \cr
\noalign{\hrule}
90&$\cdot$&$\cdot$& 92&$\star$&$\star$&$\cdot$&$\cdot$&$\star$&$\star$&$\cdot$&1 \cr
\noalign{\hrule}
100&$\cdot$&$\cdot$&$\cdot$&$\star$&$\cdot$&$\cdot$&$\star$&$\cdot$&$\cdot$&$\cdot$&none \cr
\noalign{\hrule}
110&$\cdot$&111&$\star$&$\cdot$&$\star$&$\cdot$&$\star$&$\cdot$&118&$\cdot$& 2 \cr
\noalign{\hrule}
120&$\cdot$&$\cdot$&$\cdot$&123&$\cdot$&125&$\cdot$&127&$\cdot$&129&4 \cr
\noalign{\hrule}
130&$\cdot$&131&$\cdot$&$\cdot$&$\cdot$&$\cdot$&136&$\cdot$&138&$\cdot$&3 \cr
\noalign{\hrule}
140&140&141&$\cdot$&143&144&$\cdot$&$\cdot$&147&$\cdot$&149&6 \cr
\noalign{\hrule}
150&150&151&$\cdot$&153&154&155&156&$\cdot$&$\cdot$&159&7 \cr
\noalign{\hrule}
160&$\cdot$&$\cdot$&162&$\cdot$&$\cdot$&$\cdot$&$\cdot$&$\cdot$&$\cdot$&$\cdot$&1 \cr
\noalign{\hrule}
170&170&$\cdot$&$\cdot$&$\cdot$&174&$\cdot$&176&$\cdot$&178&$\cdot$&4 \cr
\noalign{\hrule}
180&$\cdot$&181&182&$\cdot$&$\cdot$&$\cdot$&$\cdot$&$\cdot$&$\cdot$&$\cdot$&2 \cr
\noalign{\hrule}
190&190&$\cdot$&$\cdot$&$\cdot$&194&$\cdot$&196&$\cdot$&$\cdot$&$\cdot$&3 \cr
\noalign{\hrule}
200&200&$\cdot$&$\cdot$&$\cdot$&204&205&$\cdot$&207&208&209&6 \cr
\noalign{\hrule}
210&210&$\cdot$&$\cdot$&$\cdot$&$\cdot$&$\cdot$&216&$\cdot$&$\cdot$&219&3 \cr
\noalign{\hrule}
220&220&$\cdot$&222&223&$\cdot$&$\cdot$&226&227&228&229&7 \cr
\noalign{\hrule}
230&230&231&232&233&234&235&236&$\cdot$&$\cdot$&239&8 \cr
\noalign{\hrule}
240&$\cdot$&$\cdot$&$\cdot$&$\cdot$&244&$\cdot$&$\cdot$&247&248&249&4 \cr
\noalign{\hrule}
250&250&$\cdot$&$\cdot$&253&254&$\cdot$&256&257&$\cdot$&$\cdot$&5 \cr
\noalign{\hrule}
260&260&$\cdot$&$\cdot$&263&$\cdot$&265&266&267&$\cdot$&269&6 \cr
\noalign{\hrule}
270&$\cdot$&271&$\cdot$&$\cdot$&$\cdot$&275&276&277&278&279&6 \cr
\noalign{\hrule}
280&$\cdot$&281&282&$\cdot$&284&285&$\cdot$&287&$\cdot$&289&6 \cr
\noalign{\hrule}
290&$\cdot$&291&292&293&$\cdot$&$\cdot$&296&$\cdot$&$\cdot$&299&5 \cr
\noalign{\hrule}
300&$\cdot$&$\cdot$&$\cdot$&303&$\cdot$&$\cdot$&306&307&308&$\cdot$&4 \cr
\noalign{\hrule}
310&310&$\cdot$&312&313&$\cdot$&$\cdot$&$\cdot$&$\cdot$&$\cdot$&$\cdot$&3 \cr
\noalign{\hrule}
320&$\cdot$&321&322&323&$\cdot$&325&$\cdot$&327&328&329&7 \cr
\noalign{\hrule}
330&330&$\cdot$&332&333&334&335&$\cdot$&$\cdot$&338&$\cdot$&6 \cr
\noalign{\hrule}
340&340&$\cdot$&342&343&344&$\cdot$&$\cdot$&347&$\cdot$&349&6 \cr
\noalign{\hrule}
350&350&351&352&353&354&355&356&$\cdot$&358&$\cdot$&8 \cr
\noalign{\hrule}
360&$\cdot$&361&$\cdot$&363&364&$\cdot$&366&367&$\cdot$&$\cdot$&5 \cr
\noalign{\hrule}
370&370&371&$\cdot$&$\cdot$&374&375&$\cdot$&377&$\cdot$&379&6 \cr
\noalign{\hrule}
380&380&$\cdot$&$\cdot$&383&384&385&386&$\cdot$&388&389&7 \cr
\noalign{\hrule}
390&$\cdot$&$\cdot$&$\cdot$&$\cdot$&394&395&396&397&398&$\cdot$&5 \cr
\noalign{\hrule}
400&400&$\cdot$&402&403&$\cdot$&405&$\cdot$&407&408&409&7 \cr
\noalign{\hrule}
410&410&411&412&413&$\cdot$&$\cdot$&416& & & &5 \cr
\noalign{\hrule}
}
}
$$
\vskip -3 pt
\vbox {\hskip 286 pt Total Omissions = 159}

\midinsert
\narrower\narrower
\item{$\cdot$}\fic Th\'eorie des Nombres\frc, Tome 3,
Analyse Diophantine et Applications aux Cuboides
Rationnels, Gauthier-Villars, Paris 1947 pp 122-131

\item{$\star$}Sur les Cuboides Rationnels, in \fic Proceedings of the
International Congress of Mathematicians\frc\ 1954 vol. 2, Amsterdam, 1954 pp 33-34

\endinsert
\vfill\eject

\midinsert
\narrower\narrower\narrower
\hskip 0.35 in All corrections are to the first table{\raise 3 pt\hbox{$\dag$}} unless noted.\par
\smallskip
\hskip 0.35 in The small letters reference the parts of an entry as follows:\frc
\endinsert
\vskip -10 pt
$$\vbox{
\parindent 0 pt
\hsize 72 pt
\hbox{%
\vbox{\hskip 72 pt}\vbox{\hskip 72 pt}\vbox{\hskip 72 pt}\vbox{\hfil\frc Pythagorean\hfil}%
}
\hbox{%
\vbox{\hfil Entry \#\hfil}\vbox{\hfil Odd Side\hfil}\vbox{\hfil Edges\hfil}\vbox{\hfil Generators\hfil}%
}
\vskip 2 pt
\hrule
\hbox{\vrule height 10 pt depth 1 pt
\vbox{\hskip 72 pt}\vrule\vbox{\hskip 72 pt}\vrule\vbox{\fbd\quad c\hfil\frc 5.17\hfil}\vrule\vbox{\fbd\quad f\frc\hfil 11\ \ 6\hfil}\vrule%
}%
\hbox{\vrule height 10 pt depth 2 pt
\vbox{\fbd\quad a\frc\hfil 1\hfil}\vrule\vbox{\fbd\quad b\frc\hfil 85\hfil}\vrule\vbox{\fbd\quad d\hfil\frc 4.3.11\hfil}\vrule\vbox{\fbd\quad g\frc\hfil\ 5\ \ 6\hfil}\vrule%
}%
\hbox{\vrule height 9 pt depth 3 pt
\vbox{\hskip 72 pt}\vrule\vbox{\hskip 72 pt}\vrule\vbox{\fbd\quad e\hfil\frc 16.9.5\hfil}\vrule\vbox{\fbd\quad h\frc\hfil\ 9\ \ 8\hfil}\vrule%
}
\hrule
}$$

\bigskip
\bigskip
\centerline{\fbf Misprints\frc}
\bigskip

\baselineskip 4 pt
\settabs \+ \hskip 1.7 in\ \ 13\ \ & \quad Entry\quad & \ d\ \  4.3.5.17.107 \quad & \quad \cr
\+ & \hfil Entry&& \cr
\+ & \hfil \# & \hfil Old\ Print & \hfil \quad Corrected \cr
$$
\frd
\vbox{
\nointerlineskip
\halign{\strut
    \vrule height 11 pt\ \ \hfil #\ \ \vrule \hskip 2 pt
   &\vrule \hfil \quad #\quad
   &\vrule \hfil \ \ #\ 
   &\ \ #\quad \hfil
   &\vrule \quad #\quad \hfil \vrule \cr%
\noalign{\hrule}
1&47&e&16.9.5.7\ 53&16.9.5.7.53 \cr
\noalign{\hrule}
2&94&h&17\ \ 21*&17\ \ 24 \cr
\noalign{\hrule}
3&127&d&4.3.5.17\ 107&4.3.5.17.107 \cr
\noalign{\hrule}
4&141&e&128\ 3.7.83&128.3.7.83 \cr
\noalign{\hrule}
5&145&d&16.27\ 125.7&16.27.125.7 \cr
\noalign{\hrule}
6&159&d&16.9.5.7\ 43&16.9.5.7.43 \cr
\noalign{\hrule}
7&165&d&16.9\ 121\ 89&16.9.121.89 \cr
\noalign{\hrule}
8&178&c&9.5\ 11.13.61&9.5.11.13.61 \cr
\noalign{\hrule}
9&186&e&32.3\ 5.11.31&32.3.5.11.31 \cr
\noalign{\hrule}
10&203&c&9.7\ 11\ 19\ 41&9.7.11.19.41 \cr
\noalign{\hrule}
11&216&c&5.29\ 41.107&5.29.41.107 \cr
\noalign{\hrule}
12&217&c&81.5.19,83&81.5.19.83 \cr
\noalign{\hrule}
13&221&f&29\ 42\#&29\ \ 42 \cr
\noalign{\hrule}
}%
}
$$
\baselineskip 10 pt
\centerline{$\ast$ - slight misprint, top of 4 not entirely printed}
\centerline{\# - entry correct, but move 42 right one space\ \ }

\bigskip
\bigskip
\centerline{\fbf Transpositions\frc}
\bigskip

\baselineskip 4 pt
\settabs \+ \hskip 1.84 in\ 13\ & Entry & \ \ c\quad 9.7 11 19 41 \quad & \quad Corrected \cr
\+ &Entry&& \cr
\+ & \hfil \# & \hfil Old\ Print & \hfil \quad Corrected \cr
$$
\frd
\vbox{
\nointerlineskip
\halign{\strut
    \vrule height 11 pt\ \ \hfil #\ \ \vrule \hskip 2 pt
   &\vrule \hfil \quad #\quad
   &\vrule \hfil \ \ #\ 
   &\ \ #\quad \hfil
   &\vrule \quad #\quad \hfil \vrule \cr%
\noalign{\hrule}
1&238&e&16.125.9.11.3&16.\fbd 9.125\frc .11.13 \cr
\noalign{\hrule}
2&$\ddag$&d&4.25.3.19.103&4.\fbd 3.25\frc .19.103 \cr
\noalign{\hrule}
}%
}
$$
\baselineskip 10 pt
\centerline{\ddag\ - found in second table entry for \fbd 59675\frc\hskip 0.5 in}

\vfill
\centerline{\fbf --- References ---\frc}
\medskip
\midinsert
\narrower\narrower
\item{$\dag$}\fic Th\'eorie des Nombres\frc, Tome 3,
Analyse Diophantine et Applications aux Cuboides
Rationnels, Gauthier-Villars, Paris 1947 pp 122-131
\smallskip
\item{$\ddag$}Sur les Cuboides Rationnels, in \fic Proceedings of the
International Congress of Mathematicians\frc\ 1954 vol. 2, Amsterdam, 1954 pp 33-34
\endinsert
\eject

\topinsert
\narrower\narrower\narrower
\hskip 0.35 in The small letters reference the parts of an entry as follows:\frc
\endinsert
$$
\vbox{
\parindent 0 pt
\hsize 72 pt
\hbox{%
\vbox{\hskip 72 pt}\vbox{\hskip 72 pt}\vbox{\hskip 72 pt}\vbox{\hfil\frc Pythagorean\hfil}%
}
\hbox{%
\vbox{\hfil Entry \#\hfil}\vbox{\hfil Odd Side\hfil}\vbox{\hfil Edges\hfil}\vbox{\hfil Generators\hfil}%
}
\vskip 2 pt
\hrule
\hbox{\vrule height 10 pt depth 1 pt
\vbox{\hskip 72 pt}\vrule\vbox{\hskip 72 pt}\vrule\vbox{\fbd\quad c\hfil\frc 5.17\hfil}\vrule\vbox{\fbd\quad f\frc\hfil 11\ \ 6\hfil}\vrule%
}%
\hbox{\vrule height 10 pt depth 2 pt
\vbox{\fbd\quad a\frc\hfil 1\hfil}\vrule\vbox{\fbd\quad b\frc\hfil 85\hfil}\vrule\vbox{\fbd\quad d\hfil\frc 4.3.11\hfil}\vrule\vbox{\fbd\quad g\frc\hfil\ 5\ \ 6\hfil}\vrule%
}%
\hbox{\vrule height 10 pt depth 2 pt
\vbox{\hskip 72 pt}\vrule\vbox{\hskip 72 pt}\vrule\vbox{\fbd\quad e\hfil\frc 16.9.5\hfil}\vrule\vbox{\fbd\quad h\frc\hfil\ 9\ \ 8\hfil}\vrule%
}
\hrule
}$$

\bigskip
\bigskip
\centerline{\fbf Corrections to Tome 3{\raise 5 pt\hbox{$\dag$\frc}}}
\bigskip

\baselineskip 4 pt
\settabs \+ \hskip 1.77 in\ 13\ & Entry & d\quad 16.3.5.11.13.67 \quad & \quad \quad 16.3.5.11.23.67 \cr
\+ &Entry&& \cr
\+ & \hfil \# & \hfil Old\ Print & \hfil \quad\ Corrected \cr

$$
\frd
\vbox{
\nointerlineskip
\halign{\strut
    \vrule height 11 pt\ \ \hfil #\ \ \vrule \hskip 2 pt
   &\vrule \hfil \quad #\quad
   &\vrule \hfil \ \ #\ 
   &\ \ #\quad \hfil
   &\vrule \quad #\quad \hfil \vrule \cr%
\noalign{\hrule}
1&4&e&16.9.11&\fbd 64\frd .9.11 \cr
\noalign{\hrule}
2&11&e&16.3.7.13&\fbd 64\frd .3.7.13 \cr
\noalign{\hrule}
3&32&e&16.3.7.13&\fbd 64\frd .3.7.13 \cr
\noalign{\hrule}
4&34&g&27\ \ 18*&27\ \ \fbd 2\frd 8 \cr
\noalign{\hrule}
5&45&e&64.9.7&\fbd 16\frd .9.7 \cr
\noalign{\hrule}
6&59&e&32.9.11.67&32.\fbd 3\frd .11.67 \cr
\noalign{\hrule}
7&72&e&64.5.13.29&64.5.13.\fbd 19\frd \cr
\noalign{\hrule}
8&86&e&\ 6.9.5.7.11&\fbd 1\frd 6.9.5.7.11 \cr
\noalign{\hrule}
9&100&e&128.9.5.11.31&128.\fbd 3\frd .5.11.31 \cr
\noalign{\hrule}
10&103&b&100521&10052\fbd 9\frd  \cr
\noalign{\hrule}
11&112&c&11.5.9.181&11.\fbd 59\frd .181 \cr
\noalign{\hrule}
12&119&g&85\ \ 174&85\ \ 1\fbd 5\frd 4 \cr
\noalign{\hrule}
13&133&g&576\ \ 358&5\fbd 67\frd \ \ 358 \cr
\noalign{\hrule}
14&134&f&13\ \ 1\ 8&13\ \ 1\fbd 0\frd 8 \cr
\noalign{\hrule}
15&136&g&7\ \ 80&7\ \ 8\fbd 8\frd  \cr
\noalign{\hrule}
16&148&e&64.9.5.37&64.9.5.\fbd 1\frd 7 \cr
\noalign{\hrule}
17&156&e&\ 28.3.11.23&\fbd 1\frd 28.3.11.23 \cr
\noalign{\hrule}
18&158&g&17\ \ 16&17\ \ \fbd 1\frd 16 \cr
\noalign{\hrule}
19&162&b&391645&3\fbd 0\frd 1645 \cr
\noalign{\hrule}
20&175&b&366363&3\fbd 7\frd 6363 \cr
\noalign{\hrule}
21&181&f&79\ \ 40&79\ \ \fbd 11\frd 0 \cr
\noalign{\hrule}
22&196&g&\ 7\ \ 38&\fbd 2\frd 7\ \ 38 \cr
\noalign{\hrule}
23&214&b&117419&\fbd 6\frd 17419 \cr
\noalign{\hrule}
24&241&d&16.3.5.11.13.67&16.3.5.11.\fbd 2\frd 3.67 \cr
\noalign{\hrule}
}%
}
$$

\centerline{\frc$\ast$ - Not found in all copies of \dag -- some have correct entry 28}

\vfill
\centerline{\fbf --- Reference ---\frc}
\medskip
\midinsert
\narrower\narrower
\item{$\dag$}\fic Th\'eorie des Nombres\frc, Tome 3,
Analyse Diophantine et Applications aux Cuboides
Rationnels, Gauthier-Villars, Paris 1947 pp 122-131
\endinsert
\eject
}
{
\vsize= 9.2 in
\hsize= 8.5 in
\hoffset= -1 in
\footline={\hss\frc --\ \folio\ --\hss}
\pageno=1 

\headline={
\hbox{\hskip 1 in}\vbox{\hbox to 467.5328 pt{\fbf The Rational Cuboid Table --- \hfil \fsd of Maurice Kraitchik\frb}\vskip 3 pt\hrule}\hbox{\hskip 1 in}}

\vglue -21 pt
{\noindent\hskip 1 in\hbox to 6.5 in{\ 1 -- 36 \hfil\fbd 85 -- 8789\frb}}
\vskip -9 pt
$$
\vbox{
\halign{\strut
    \vrule \ \ \hfil \frb #\ 
   &\vrule \hfil \ \ \fbb #\frb\ 
   &\vrule \hfil \ \ \frb #\ \hfil
   &\vrule \hfil \ \ \frb #\ 
   &\vrule \hfil \ \ \frb #\ \ \vrule \hskip 2 pt
   &\vrule \ \ \hfil \frb #\ 
   &\vrule \hfil \ \ \fbb #\frb\ 
   &\vrule \hfil \ \ \frb #\ \hfil
   &\vrule \hfil \ \ \frb #\ 
   &\vrule \hfil \ \ \frb #\ \vrule \cr%
\noalign{\hrule}
 & &5.17&11&6& & &27.7.11&145&152 \cr
1&85&4.3.11&5&6&19&2079&16.5.19.29&33&62 \cr
 & &16.9.5&9&8& & &64.3.11.31&31&32 \cr
\noalign{\hrule}
 & &9.13&11&2& & &3.7.103&53&50 \cr
2&117&4.11&5&6&20&2163&4.25.7.53&9&44 \cr
 & &16.3.5&5&8& & &32.9.5.11&55&48 \cr
\noalign{\hrule}
 & &11.17&5&6& & &27.5.17&13&4 \cr
3&187&4.3.5.17&11&6&21&2295&8.3.5.13&17&22 \cr
 & &16.9.11&9&8& & &32.11.17&11&16 \cr
\noalign{\hrule}
 & &3.5.13&17&22& & &9.25.11&13&2 \cr
4&195&4.11.17&9&8&22&2475&4.3.5.13&17&22 \cr
 & &64.9.11&33&32& & &16.11.17&17&8 \cr
\noalign{\hrule}
 & &3.7.11&3&4& & &9.25.13&17&22 \cr
5&231&8.9.11&1&10&23&2925&4.3.5.11.17&13&2 \cr
 & &32.5&5&16& & &16.13.17&17&8 \cr
\noalign{\hrule}
 & &25.11&7&18& & &9.7.73&33&40 \cr
6&275&4.9.7&5&2&24&4599&16.27.5.11&41&14 \cr
 & &16.3.5&3&8& & &64.7.41&41&32 \cr
\noalign{\hrule}
 & &3.11.13&5&6& & &9.7.73&33&40 \cr
7&429&4.9.5.13&11&2&25&4599&16.27.5.11&73&62 \cr
 & &16.5.11&5&8& & &64.31.73&31&32 \cr
\noalign{\hrule}
 & &9.5.11&47&52& & &169.29&99&70 \cr
8&495&8.13.47&17&30&26&4901&4.9.5.7.11&13&2 \cr
 & &32.3.5.17&17&16& & &16.3.7.13&21&8 \cr
\noalign{\hrule}
 & &9.7.11&1&10& & &289.19&135&154 \cr
9&693&4.5.7&3&4&27&5491&4.27.5.7.11&17&38 \cr
 & &32.3.5&5&16& & &16.9.17.19&9&8 \cr
\noalign{\hrule}
 & &9.5.19&11&8& & &27.11.19&91&118 \cr
10&855&16.3.5.11&13&2&28&5643&4.7.13.59&3&10 \cr
 & &64.13&13&32& & &16.3.5.59&59&40 \cr
\noalign{\hrule}
 & &5.11.17&27&28& & &27.11.19&13&14 \cr
11&935&8.27.7.17&13&4&29&5643&4.7.11.13.19&15&4 \cr
 & &64.3.7.13&91&96& & &32.3.5.7.13&91&80 \cr
\noalign{\hrule}
 & &5.13.17&9&8& & &25.11.23&7&18 \cr
12&1105&16.9.5.13&17&22&30&6325&4.9.7.23&1&22 \cr
 & &64.3.11.17&33&32& & &16.3.11&1&24 \cr
\noalign{\hrule}
 & &3.5.7.11&5&2& & &9.5.11.13&17&22 \cr
13&1155&4.25.11&7&18&31&6435&4.3.121.17&35&86 \cr
 & &16.9.7&3&8& & &16.5.7.43&43&56 \cr
\noalign{\hrule}
 & &3.5.7.11&5&6& & &11.13.53&45&98 \cr
14&1155&4.9.25.7&19&44&32&7579&4.9.5.49&13&8 \cr
 & &32.11.19&19&16& & &64.3.7.13&21&32 \cr
\noalign{\hrule}
 & &9.25.7&19&44& & &5.19.83&217&198 \cr
15&1575&8.11.19&5&6&33&7885&4.9.7.11.31&85&8 \cr
 & &32.3.5.19&19&16& & &64.3.5.17&51&32 \cr
\noalign{\hrule}
 & &27.5.13&17&22& & &9.5.11.17&13&4 \cr
16&1755&4.9.11.17&13&4&34&8415&8.5.11.13&27&28 \cr
 & &32.11.13&11&16& & &64.27.7.13&91&96 \cr
\noalign{\hrule}
 & &9.11.19&15&4& & &9.5.11.17&77&76 \cr
17&1881&8.27.5&13&14&35&8415&8.5.7.121.19&13&108 \cr
 & &32.5.7.13&91&80& & &64.27.7.13&91&96 \cr
\noalign{\hrule}
 & &5.11.37&9&46& & &11.17.47&29&18 \cr
18&2035&4.9.23&13&10&36&8789&4.9.17.29&11&40 \cr
 & &16.3.5.13&13&24& & &64.3.5.11&15&32 \cr
\noalign{\hrule}
}%
}
$$
\vfill\eject
\vglue -23 pt
\noindent\hskip 1 in\hbox to 6.5 in{\ 37 -- 72 \hfill\fbd 9045 -- 34965\frb}
\vskip -9 pt
$$
\vbox{
\nointerlineskip
\halign{\strut
    \vrule \ \ \hfil \frb #\ 
   &\vrule \hfil \ \ \fbb #\frb\ 
   &\vrule \hfil \ \ \frb #\ \hfil
   &\vrule \hfil \ \ \frb #\ 
   &\vrule \hfil \ \ \frb #\ \ \vrule \hskip 2 pt
   &\vrule \ \ \hfil \frb #\ 
   &\vrule \hfil \ \ \fbb #\frb\ 
   &\vrule \hfil \ \ \frb #\ \hfil
   &\vrule \hfil \ \ \frb #\ 
   &\vrule \hfil \ \ \frb #\ \vrule \cr%
\noalign{\hrule}
 & &27.5.67&13&14& & &27.5.11.13&7&20 \cr
37&9045&4.5.7.13.67&1&66&55&19305&8.25.7.11&7&18 \cr
 & &16.3.7.11&11&56& & &32.9.49&49&16 \cr
\noalign{\hrule}
 & &9.5.11.19&13&2& & &3.5.7.11.17&17&38 \cr
38&9405&4.3.13.19&11&8&56&19635&4.289.19&135&154 \cr
 & &64.11.13&13&32& & &16.27.5.7.11&9&8 \cr
\noalign{\hrule}
 & &9.5.11.19&31&26& & &3.11.13.47&17&30 \cr
39&9405&4.3.11.13.31&25&118&57&20163&4.9.5.11.17&47&52 \cr
 & &16.25.59&59&40& & &32.13.17.47&17&16 \cr
\noalign{\hrule}
 & &27.5.7.11&73&62& & &9.11.227&115&112 \cr
40&10395&4.7.31.73&33&40&58&22473&32.3.5.7.11.23&163&2 \cr
 & &64.3.5.11.31&31&32& & &128.163&163&64 \cr
\noalign{\hrule}
 & &27.5.7.11&29&26& & &9.7.13.29&25&4 \cr
41&10395&4.9.7.13.29&25&38&59&23751&8.3.25.13&17&22 \cr
 & &16.25.19.29&145&152& & &32.5.11.17&85&176 \cr
\noalign{\hrule}
 & &3.25.11.13&27&38& & &5.11.19.23&37&18 \cr
42&10725&4.81.5.19&7&88&60&24035&4.9.23.37&67&44 \cr
 & &64.7.11&7&32& & &32.3.11.67&67&48 \cr
\noalign{\hrule}
 & &3.25.11.13&7&18& & &3.25.17.19&33&52 \cr
43&10725&4.27.7.13&7&20&61&24225&8.9.5.11.13&17&82 \cr
 & &32.5.49&49&16& & &32.17.41&41&16 \cr
\noalign{\hrule}
 & &7.23.73&117&44& & &81.7.47&55&8 \cr
44&11753&8.9.11.13&25&14&62&26649&16.9.5.11&17&28 \cr
 & &32.3.25.7&25&48& & &128.7.17&17&64 \cr
\noalign{\hrule}
 & &3.25.7.23&1&22& & &9.25.7.17&121&104 \cr
45&12075&4.25.11&7&18&63&26775&16.7.121.13&15&106 \cr
 & &16.9.7&1&24& & &64.3.5.53&53&32 \cr
\noalign{\hrule}
 & &9.5.7.41&23&22& & &27.7.11.13&25&38 \cr
46&12915&4.7.11.23.41&17&270&64&27027&4.3.25.11.19&29&26 \cr
 & &16.27.5.17&17&24& & &16.5.13.19.29&145&152 \cr
\noalign{\hrule}
 & &27.5.109&259&286& & &5.7.13.61&107&198 \cr
47&14715&4.7.11.13.37&53&90&65&27755&4.9.11.107&37&70 \cr
 & &16.9.5.7.53&53&56& & &16.3.5.7.37&37&24 \cr
\noalign{\hrule}
 & &3.25.7.29&99&104& & &3.11.23.37&13&10 \cr
48&15225&16.27.5.11.13&139&4&66&28083&4.5.11.13.37&9&46 \cr
 & &128.139&139&64& & &16.9.13.23&13&24 \cr
\noalign{\hrule}
 & &9.7.11.23&65&142& & &9.11.13.23&29&40 \cr
49&15939&4.5.13.71&29&42&67&29601&16.3.5.13.29&53&92 \cr
 & &16.3.5.7.29&29&40& & &128.23.53&53&64 \cr
\noalign{\hrule}
 & &81.11.19&23&4& & &9.5.11.61&35&26 \cr
50&16929&8.3.11.23&29&40&68&30195&4.25.7.11.13&51&26 \cr
 & &128.5.29&145&64& & &16.3.169.17&169&136 \cr
\noalign{\hrule}
 & &3.7.19.43&7&50& & &9.25.11.13&131&144 \cr
51&17157&4.25.49&27&22&69&32175&32.81.131&25&106 \cr
 & &16.27.5.11&45&88& & &128.25.53&53&64 \cr
\noalign{\hrule}
 & &9.25.7.11&133&142& & &125.7.37&81&44 \cr
52&17325&4.49.19.71&11&60&70&32375&8.81.7.11&89&100 \cr
 & &32.3.5.11.19&19&16& & &64.3.25.89&89&96 \cr
\noalign{\hrule}
 & &3.25.13.19&17&22& & &9.7.17.31&75&44 \cr
53&18525&4.5.11.17.19&37&18&71&33201&8.27.25.11&73&62 \cr
 & &16.9.17.37&111&136& & &32.5.31.73&73&80 \cr
\noalign{\hrule}
 & &25.13.59&371&396& & &27.5.7.37&17&10 \cr
54&19175&8.9.7.11.53&65&12&72&34965&4.25.17.37&231&194 \cr
 & &64.27.5.13&27&32& & &16.3.7.11.97&97&88 \cr
\noalign{\hrule}
}%
}
$$
\eject
\vglue -23 pt
\noindent\hskip 1 in\hbox to 6.5 in{\ 73 -- 108 \hfill\fbd 35075 -- 69513\frb}
\vskip -9 pt
$$
\vbox{
\nointerlineskip
\halign{\strut
    \vrule \ \ \hfil \frb #\ 
   &\vrule \hfil \ \ \fbb #\frb\ 
   &\vrule \hfil \ \ \frb #\ \hfil
   &\vrule \hfil \ \ \frb #\ 
   &\vrule \hfil \ \ \frb #\ \ \vrule \hskip 2 pt
   &\vrule \ \ \hfil \frb #\ 
   &\vrule \hfil \ \ \fbb #\frb\ 
   &\vrule \hfil \ \ \frb #\ \hfil
   &\vrule \hfil \ \ \frb #\ 
   &\vrule \hfil \ \ \frb #\ \vrule \cr%
\noalign{\hrule}
 & &25.23.61&43&18& & &25.29.79&27&2 \cr
73&35075&4.9.23.43&33&10&91&57275&4.27.79&35&44 \cr
 & &16.27.5.11&27&88& & &32.3.5.7.11&33&112 \cr
\noalign{\hrule}
 & &11.169.19&189&20& & &3.25.19.41&13&12 \cr
74&35321&8.27.5.7&13&22&92&58425&8.9.13.19.41&61&308 \cr
 & &32.3.11.13&3&16& & &64.7.11.61&427&352 \cr
\noalign{\hrule}
 & &3.11.29.37&13&24& & &9.5.13.101&343&242 \cr
75&35409&16.9.13.29&19&10&93&59085&4.343.121&177&166 \cr
 & &64.5.13.19&247&160& & &16.3.11.59.83&649&664 \cr
\noalign{\hrule}
 & &3.49.13.19&67&66& & &25.7.11.31&103&114 \cr
76&36309&4.9.7.11.13.67&95&4&94&59675&4.3.25.19.103&3&22 \cr
 & &32.5.19.67&67&80& & &16.9.11.103&103&72 \cr
\noalign{\hrule}
 & &5.7.23.47&167&162& & &27.5.11.41&13&14 \cr
77&37835&4.81.23.167&187&20&95&60885&4.5.7.11.13.41&31&174 \cr
 & &32.9.5.11.17&153&176& & &16.3.7.29.31&203&248 \cr
\noalign{\hrule}
 & &81.25.19&1&26& & &3.5.7.11.53&13&8 \cr
78&38475&4.3.13.19&35&22&96&61215&16.11.13.53&45&98 \cr
 & &16.5.7.11&7&88& & &64.9.5.49&21&32 \cr
\noalign{\hrule}
 & &3.13.23.43&33&10& & &3.5.7.19.31&11&46 \cr
79&38571&4.9.5.11.13&31&86&97&61845&4.11.23.31&27&4 \cr
 & &16.31.43&31&8& & &32.27.11&11&144 \cr
\noalign{\hrule}
 & &9.5.13.67&1&66& & &25.37.67&31&36 \cr
80&39195&4.27.11&13&14&98&61975&8.9.5.31.37&11&26 \cr
 & &16.7.11.13&11&56& & &32.3.11.13.31&429&496 \cr
\noalign{\hrule}
 & &27.121.13&115&236& & &3.5.11.13.29&13&2 \cr
81&42471&8.5.23.59&41&18&99&62205&4.169.29&99&70 \cr
 & &32.9.5.41&41&80& & &16.9.5.7.11&21&8 \cr
\noalign{\hrule}
 & &13.17.193&103&90& & &9.5.19.73&29&44 \cr
82&42653&4.9.5.17.103&77&26&100&62415&8.3.11.19.29&23&34 \cr
 & &16.3.5.7.11.13&105&88& & &32.17.23.29&391&464 \cr
\noalign{\hrule}
 & &3.13.29.41&785&814& & &9.49.11.13&95&4 \cr
83&46371&4.5.11.37.157&171&14&101&63063&8.5.7.19&67&66 \cr
 & &16.9.7.11.19&209&168& & &32.3.5.11.67&67&80 \cr
\noalign{\hrule}
 & &9.7.11.67&233&236& & &3.5.11.13.31&41&52 \cr
84&46431&8.3.11.59.233&25&674&102&66495&8.5.169.41&187&18 \cr
 & &32.25.337&337&400& & &32.9.11.17&17&48 \cr
\noalign{\hrule}
 & &25.19.101&7&12& & &3.5.11.13.31&87&56 \cr
85&47975&8.3.5.7.101&33&68&103&66495&16.9.5.7.29&79&124 \cr
 & &64.9.11.17&187&288& & &128.31.79&79&64 \cr
\noalign{\hrule}
 & &9.13.19.23&25&44& & &3.25.17.53&69&16 \cr
86&51129&8.3.25.11.13&47&8&104&67575&32.9.5.23&11&34 \cr
 & &128.5.47&235&64& & &128.11.17&11&64 \cr
\noalign{\hrule}
 & &5.49.11.19&39&94& & &25.11.13.19&37&18 \cr
87&51205&4.3.7.13.47&69&22&105&67925&4.9.5.13.37&17&22 \cr
 & &16.9.11.23&23&72& & &16.3.11.17.37&111&136 \cr
\noalign{\hrule}
 & &3.5.49.71&11&60& & &9.5.29.53&41&46 \cr
88&52185&8.9.25.11&133&142&106&69165&4.3.23.41.53&35&88 \cr
 & &32.7.19.71&19&16& & &64.5.7.11.23&253&224 \cr
\noalign{\hrule}
 & &3.13.19.71&35&22& & &5.11.13.97&21&34 \cr
89&52611&4.5.7.11.71&53&18&107&69355&4.3.7.17.97&11&108 \cr
 & &16.9.11.53&159&88& & &32.81.11&81&16 \cr
\noalign{\hrule}
 & &3.121.149&69&80& & &3.17.29.47&11&40 \cr
90&54087&32.9.5.11.23&107&8&108&69513&16.5.11.47&29&18 \cr
 & &512.107&107&256& & &64.9.5.29&15&32 \cr
\noalign{\hrule}
}%
}
$$
\eject
\vglue -23 pt
\noindent\hskip 1 in\hbox to 6.5 in{\ 109 -- 144 \hfill\fbd 72611 -- 128205\frb}
\vskip -9 pt
$$
\vbox{
\nointerlineskip
\halign{\strut
    \vrule \ \ \hfil \frb #\ 
   &\vrule \hfil \ \ \fbb #\frb\ 
   &\vrule \hfil \ \ \frb #\ \hfil
   &\vrule \hfil \ \ \frb #\ 
   &\vrule \hfil \ \ \frb #\ \ \vrule \hskip 2 pt
   &\vrule \ \ \hfil \frb #\ 
   &\vrule \hfil \ \ \fbb #\frb\ 
   &\vrule \hfil \ \ \frb #\ \hfil
   &\vrule \hfil \ \ \frb #\ 
   &\vrule \hfil \ \ \frb #\ \vrule \cr%
\noalign{\hrule}
 & &7.11.23.41&17&270& & &25.7.19.31&3&22 \cr
109&72611&4.27.5.17&23&22&127&103075&4.3.7.11.31&103&114 \cr
 & &16.3.11.17.23&17&24& & &16.9.19.103&103&72 \cr
\noalign{\hrule}
 & &27.5.49.11&53&46& & &9.5.29.79&35&44 \cr
110&72765&4.3.5.7.23.53&109&374&128&103095&8.25.7.11.29&27&2 \cr
 & &16.11.17.109&109&136& & &32.27.7.11&33&112 \cr
\noalign{\hrule}
 & &9.7.19.61&215&212& & &3.7.121.41&2491&2470 \cr
111&73017&8.3.5.19.43.53&11&806&129&104181&4.5.13.19.47.53&41&54 \cr
 & &32.11.13.31&341&208& & &16.27.41.47.53&423&424 \cr
\noalign{\hrule}
 & &9.5.11.151&2491&2492& & &9.11.29.37&19&10 \cr
112&74745&8.3.5.7.47.53.89&41&664&130&106227&4.5.11.19.37&13&24 \cr
 & &128.41.53.83&3403&3392& & &64.3.5.13.19&247&160 \cr
\noalign{\hrule}
 & &9.5.7.13.19&13&22& & &7.11.23.61&25&36 \cr
113&77805&4.11.169.19&189&20&131&108031&8.9.25.7.23&17&52 \cr
 & &32.27.5.7&3&16& & &64.3.5.13.17&255&416 \cr
\noalign{\hrule}
 & &9.25.7.53&5&58& & &27.5.11.73&41&14 \cr
114&83475&4.125.29&77&48&132&108405&4.7.41.73&33&40 \cr
 & &128.3.7.11&11&64& & &64.3.5.11.41&41&32 \cr
\noalign{\hrule}
 & &9.7.17.79&31&48& & &9.11.19.59&13&46 \cr
115&84609&32.27.7.31&79&110&133&110979&4.3.13.19.23&35&22 \cr
 & &128.5.11.79&55&64& & &16.5.7.11.23&115&56 \cr
\noalign{\hrule}
 & &3.25.17.67&671&604& & &27.23.179&55&124 \cr
116&85425&8.11.61.151&45&106&134&111159&8.9.5.11.31&65&34 \cr
 & &32.9.5.11.53&159&176& & &32.25.13.17&221&400 \cr
\noalign{\hrule}
 & &9.25.7.61&143&82& & &5.7.11.13.23&27&38 \cr
117&96075&4.7.11.13.41&305&228&135&115115&4.27.7.19.23&109&52 \cr
 & &32.3.5.19.61&19&16& & &32.9.13.109&109&144 \cr
\noalign{\hrule}
 & &25.7.13.43&1969&1944& & &5.19.23.53&561&446 \cr
118&97825&16.243.11.179&59&238&136&115805&4.3.11.17.223&117&106 \cr
 & &64.9.7.17.59&531&544& & &16.27.13.17.53&221&216 \cr
\noalign{\hrule}
 & &11.289.31&315&26& & &289.401&145&144 \cr
119&98549&4.9.5.7.13&17&22&137&115889&32.9.5.29.401&17&418 \cr
 & &16.3.7.11.17&7&24& & &128.3.11.17.19&209&192 \cr
\noalign{\hrule}
 & &27.7.17.31&79&110& & &9.7.11.169&19&58 \cr
120&99603&4.5.11.17.79&31&48&138&117117&4.3.13.19.29&385&356 \cr
 & &128.3.5.11.31&55&64& & &32.5.7.11.89&89&80 \cr
\noalign{\hrule}
 & &81.5.13.19&35&22& & &11.59.181&415&234 \cr
121&100035&4.27.25.7.11&1&26&139&117469&4.9.5.13.83&61&22 \cr
 & &16.7.11.13&7&88& & &16.3.5.11.61&61&120 \cr
\noalign{\hrule}
 & &9.125.89&67&58& & &9.11.29.41&139&180 \cr
122&100125&4.29.67.89&11&78&140&117711&8.81.5.139&29&110 \cr
 & &16.3.11.13.29&143&232& & &32.25.11.29&25&16 \cr
\noalign{\hrule}
 & &9.5.7.11.29&97&106& & &9.5.43.61&23&38 \cr
123&100485&4.5.11.53.97&49&534&141&118035&4.3.19.23.43&55&74 \cr
 & &16.3.49.89&89&56& & &16.5.11.23.37&253&296 \cr
\noalign{\hrule}
 & &11.13.19.37&2655&2636& & &9.5.7.13.29&17&22 \cr
124&100529&8.9.5.59.659&241&418&142&118755&4.3.7.11.17.29&25&4 \cr
 & &32.3.5.11.19.241&241&240& & &32.25.11.17&85&176 \cr
\noalign{\hrule}
 & &9.5.37.61&73&110& & &121.23.43&555&434 \cr
125&101565&4.3.25.11.73&97&122&143&119669&4.3.5.7.31.37&15&22 \cr
 & &16.11.61.97&97&88& & &16.9.25.11.31&225&248 \cr
\noalign{\hrule}
 & &3.5.13.17.31&17&22& & &9.5.7.11.37&23&14 \cr
126&102765&4.11.289.31&315&26&144&128205&4.5.49.11.23&3&52 \cr
 & &16.9.5.7.13&7&24& & &32.3.13.23&23&208 \cr
\noalign{\hrule}
}%
}
$$
\eject
\vglue -23 pt
\noindent\hskip 1 in\hbox to 6.5 in{\ 145 -- 180 \hfill\fbd 131157 -- 200165\frb}
\vskip -9 pt
$$
\vbox{
\nointerlineskip
\halign{\strut
    \vrule \ \ \hfil \frb #\ 
   &\vrule \hfil \ \ \fbb #\frb\ 
   &\vrule \hfil \ \ \frb #\ \hfil
   &\vrule \hfil \ \ \frb #\ 
   &\vrule \hfil \ \ \frb #\ \ \vrule \hskip 2 pt
   &\vrule \ \ \hfil \frb #\ 
   &\vrule \hfil \ \ \fbb #\frb\ 
   &\vrule \hfil \ \ \frb #\ \hfil
   &\vrule \hfil \ \ \frb #\ 
   &\vrule \hfil \ \ \frb #\ \vrule \cr%
\noalign{\hrule}
 & &9.13.19.59&35&22& & &3.25.7.17.19&33&52 \cr
145&131157&4.3.5.7.11.59&13&46&163&169575&8.9.5.7.11.13&59&4 \cr
 & &16.5.7.13.23&115&56& & &64.13.59&59&416 \cr
\noalign{\hrule}
 & &3.5.11.19.43&69&26& & &9.7.11.13.19&3&10 \cr
146&134805&4.9.11.13.23&61&38&164&171171&4.27.5.11.19&91&118 \cr
 & &16.13.19.61&61&104& & &16.5.7.13.59&59&40 \cr
\noalign{\hrule}
 & &3.11.23.179&89&90& & &9.5.7.19.29&33&62 \cr
147&135861&4.27.5.11.23.89&179&800&165&173565&4.27.7.11.31&145&152 \cr
 & &256.125.179&125&128& & &64.5.19.29.31&31&32 \cr
\noalign{\hrule}
 & &27.25.7.29&143&172& & &27.5.7.11.17&73&62 \cr
148&137025&8.3.5.11.13.43&119&76&166&176715&4.7.17.31.73&75&44 \cr
 & &64.7.11.17.19&323&352& & &32.3.25.11.73&73&80 \cr
\noalign{\hrule}
 & &27.25.11.19&91&118& & &3.7.11.13.59&9&68 \cr
149&141075&4.25.7.13.59&369&44&167&177177&8.27.13.17&7&20 \cr
 & &32.9.11.41&41&16& & &64.5.7.17&17&160 \cr
\noalign{\hrule}
 & &9.25.17.37&231&194& & &7.11.17.137&107&30 \cr
150&141525&4.27.7.11.97&17&10&168&179333&4.3.5.17.107&11&96 \cr
 & &16.5.11.17.97&97&88& & &256.9.11&9&128 \cr
\noalign{\hrule}
 & &9.7.121.19&71&50& & &9.11.19.97&59&40 \cr
151&144837&4.3.25.19.71&59&154&169&182457&16.5.59.97&19&78 \cr
 & &16.5.7.11.59&59&40& & &64.3.5.13.19&65&32 \cr
\noalign{\hrule}
 & &9.5.11.13.23&31&86& & &3.25.37.67&11&26 \cr
152&148005&4.23.31.43&33&10&170&185925&4.5.11.13.67&31&36 \cr
 & &16.3.5.11.31&31&8& & &32.9.11.13.31&429&496 \cr
\noalign{\hrule}
 & &25.11.19.29&27&2& & &5.11.17.199&127&72 \cr
153&151525&4.27.11.19&59&40&171&186065&16.9.17.127&199&182 \cr
 & &64.3.5.59&177&32& & &64.3.7.13.199&91&96 \cr
\noalign{\hrule}
 & &9.101.167&55&46& & &9.5.11.13.29&53&92 \cr
154&151803&4.5.11.23.167&43&210&172&186615&8.3.11.23.53&29&40 \cr
 & &16.3.25.7.43&301&200& & &128.5.29.53&53&64 \cr
\noalign{\hrule}
 & &25.11.13.43&113&102& & &3.11.59.97&19&78 \cr
155&153725&4.3.5.13.17.113&393&172&173&188859&4.9.11.13.19&59&40 \cr
 & &32.9.43.131&131&144& & &64.5.13.59&65&32 \cr
\noalign{\hrule}
 & &3.11.43.109&33&76& & &5.113.337&451&114 \cr
156&154671&8.9.121.19&25&146&174&190405&4.3.11.19.41&125&84 \cr
 & &32.25.73&73&400& & &32.9.125.7&225&112 \cr
\noalign{\hrule}
 & &27.11.17.31&13&4& & &5.7.53.103&9&44 \cr
157&156519&8.3.11.13.31&25&118&175&191065&8.9.11.103&53&50 \cr
 & &32.25.59&59&400& & &32.3.25.11.53&55&48 \cr
\noalign{\hrule}
 & &7.11.13.157&621&478& & &5.49.11.71&221&276 \cr
158&157157&4.27.23.239&85&154&176&191345&8.3.7.13.17.23&71&48 \cr
 & &16.9.5.7.11.17&85&72& & &256.9.13.71&117&128 \cr
\noalign{\hrule}
 & &3.5.59.181&61&120& & &9.7.17.179&703&550 \cr
159&160185&16.9.25.61&143&82&177&191709&4.25.11.19.37&567&358 \cr
 & &64.11.13.41&533&352& & &16.81.7.179&9&8 \cr
\noalign{\hrule}
 & &9.49.13.29&19&68& & &9.5.7.13.47&295&316 \cr
160&166257&8.3.13.17.19&35&22&178&192465&8.3.25.59.79&77&2 \cr
 & &32.5.7.11.17&85&176& & &32.7.11.59&59&176 \cr
\noalign{\hrule}
 & &49.13.263&187&450& & &5.121.17.19&13&108 \cr
161&167531&4.9.25.11.17&91&96&179&195415&8.27.13.17&77&76 \cr
 & &256.27.5.7.13&135&128& & &64.3.7.11.13.19&91&96 \cr
\noalign{\hrule}
 & &5.7.11.19.23&603&442& & &5.49.19.43&27&22 \cr
162&168245&4.9.13.17.67&11&28&180&200165&4.27.11.19.43&7&50 \cr
 & &32.3.7.11.67&67&48& & &16.9.25.7.11&45&88 \cr
\noalign{\hrule}
}%
}
$$
\eject
\vglue -23 pt
\noindent\hskip 1 in\hbox to 6.5 in{\ 181 -- 216 \hfill\fbd 200385 -- 278355\frb}
\vskip -9 pt
$$
\vbox{
\nointerlineskip
\halign{\strut
    \vrule \ \ \hfil \frb #\ 
   &\vrule \hfil \ \ \fbb #\frb\ 
   &\vrule \hfil \ \ \frb #\ \hfil
   &\vrule \hfil \ \ \frb #\ 
   &\vrule \hfil \ \ \frb #\ \ \vrule \hskip 2 pt
   &\vrule \ \ \hfil \frb #\ 
   &\vrule \hfil \ \ \fbb #\frb\ 
   &\vrule \hfil \ \ \frb #\ \hfil
   &\vrule \hfil \ \ \frb #\ 
   &\vrule \hfil \ \ \frb #\ \vrule \cr%
\noalign{\hrule}
 & &9.5.61.73&481&176& & &7.11.13.229&873&730 \cr
181&200385&32.11.13.37&25&12&199&229229&4.9.5.73.97&109&182 \cr
 & &256.3.25.11&55&128& & &16.3.5.7.13.109&109&120 \cr
\noalign{\hrule}
 & &9.5.61.73&29&44& & &25.13.17.43&99&116 \cr
182&200385&8.3.11.29.61&13&74&200&237575&8.9.5.11.13.29&1&14 \cr
 & &32.11.13.37&481&176& & &32.3.7.11.29&203&528 \cr
\noalign{\hrule}
 & &27.17.19.23&3&20& & &9.121.13.17&35&86 \cr
183&200583&8.81.5.19&7&88&201&240669&4.3.5.7.13.43&17&22 \cr
 & &128.7.11&11&448& & &16.7.11.17.43&43&56 \cr
\noalign{\hrule}
 & &81.5.7.71&13&58& & &3.5.19.23.37&67&44 \cr
184&201285&4.9.7.13.29&55&62&202&242535&8.5.11.19.67&37&18 \cr
 & &16.5.11.29.31&319&248& & &32.9.37.67&67&48 \cr
\noalign{\hrule}
 & &27.7.11.97&191&488& & &3.25.11.13.23&151&174 \cr
185&201663&16.61.191&65&126&203&246675&4.9.11.29.151&125&26 \cr
 & &64.9.5.7.13&65&32& & &16.125.13.29&29&40 \cr
\noalign{\hrule}
 & &3.125.19.29&101&44& & &7.11.17.193&97&90 \cr
186&206625&8.25.11.101&63&38&204&252637&4.9.5.97.193&49&242 \cr
 & &32.9.7.11.19&77&48& & &16.3.5.49.121&105&88 \cr
\noalign{\hrule}
 & &81.5.11.47&17&28& & &23.43.263&153&110 \cr
187&209385&8.9.7.17.47&55&8&205&260107&4.9.5.11.17.23&43&26 \cr
 & &128.5.11.17&17&64& & &16.3.5.11.13.43&143&120 \cr
\noalign{\hrule}
 & &25.7.11.109&37&72& & &9.5.11.13.41&31&174 \cr
188&209825&16.9.5.11.37&83&28&206&263835&4.27.29.31&13&14 \cr
 & &128.3.7.83&83&192& & &16.7.13.29.31&203&248 \cr
\noalign{\hrule}
 & &17.31.401&185&216& & &3.5.49.359&71&76 \cr
189&211327&16.27.5.17.37&209&124&207&263865&8.19.71.359&495&854 \cr
 & &128.3.11.19.31&209&192& & &32.9.5.7.11.61&183&176 \cr
\noalign{\hrule}
 & &27.11.23.31&3431&3400& & &3.121.17.43&245&228 \cr
190&211761&16.25.17.47.73&33&1208&208&265353&8.9.5.49.11.19&17&116 \cr
 & &256.3.11.151&151&128& & &64.5.7.17.29&203&160 \cr
\noalign{\hrule}
 & &5.11.169.23&549&296& & &27.7.13.109&1517&1426 \cr
191&213785&16.9.37.61&73&110&209&267813&4.23.31.37.41&1045&102 \cr
 & &64.3.5.11.73&73&96& & &16.3.5.11.17.19&323&440 \cr
\noalign{\hrule}
 & &3.5.7.13.157&33&124& & &9.5.11.19.29&59&40 \cr
192&214305&8.9.5.11.31&7&38&210&272745&16.25.29.59&27&2 \cr
 & &32.7.11.19&19&176& & &64.27.59&177&32 \cr
\noalign{\hrule}
 & &49.11.13.31&125&216& & &11.13.23.83&413&666 \cr
193&217217&16.27.125.7&31&94&211&272987&4.9.7.37.59&85&26 \cr
 & &64.3.31.47&47&96& & &16.3.5.7.13.17&85&168 \cr
\noalign{\hrule}
 & &3.5.13.19.59&153&94& & &7.121.17.19&127&6 \cr
194&218595&4.27.5.17.47&91&44&212&273581&4.3.17.127&55&72 \cr
 & &32.7.11.13.17&119&176& & &64.27.5.11&27&160 \cr
\noalign{\hrule}
 & &5.29.37.41&107&78& & &9.25.23.53&11&34 \cr
195&219965&4.3.13.41.107&33&74&213&274275&4.5.11.17.53&69&16 \cr
 & &16.9.11.13.37&117&88& & &128.3.11.23&11&64 \cr
\noalign{\hrule}
 & &27.5.31.53&551&286& & &9.11.47.59&7&524 \cr
196&221805&4.11.13.19.29&81&62&214&274527&8.7.131&65&66 \cr
 & &16.81.29.31&29&24& & &32.3.5.7.11.13&65&112 \cr
\noalign{\hrule}
 & &3.25.29.103&517&208& & &3.7.13.19.53&33&20 \cr
197&224025&32.11.13.47&29&18&215&274911&8.9.5.7.11.19&17&116 \cr
 & &128.9.13.29&39&64& & &64.5.17.29&145&544 \cr
\noalign{\hrule}
 & &27.11.13.59&7&20& & &3.5.7.11.241&1273&1378 \cr
198&227799&8.5.7.11.59&9&68&216&278355&4.13.19.53.67&17&36 \cr
 & &64.9.5.17&17&160& & &32.9.13.17.67&871&816 \cr
\noalign{\hrule}
}%
}
$$
\eject
\vglue -23 pt
\noindent\hskip 1 in\hbox to 6.5 in{\ 217 -- 252 \hfill\fbd 287287 -- 361665\frb}
\vskip -9 pt
$$
\vbox{
\nointerlineskip
\halign{\strut
    \vrule \ \ \hfil \frb #\ 
   &\vrule \hfil \ \ \fbb #\frb\ 
   &\vrule \hfil \ \ \frb #\ \hfil
   &\vrule \hfil \ \ \frb #\ 
   &\vrule \hfil \ \ \frb #\ \ \vrule \hskip 2 pt
   &\vrule \ \ \hfil \frb #\ 
   &\vrule \hfil \ \ \fbb #\frb\ 
   &\vrule \hfil \ \ \frb #\ \hfil
   &\vrule \hfil \ \ \frb #\ 
   &\vrule \hfil \ \ \frb #\ \vrule \cr%
\noalign{\hrule}
 & &49.11.13.41&215&72& & &25.13.17.59&371&396 \cr
217&287287&16.9.5.7.43&47&82&235&325975&8.9.7.11.17.53&19&2 \cr
 & &64.3.41.47&47&96& & &32.3.11.19.53&627&848 \cr
\noalign{\hrule}
 & &3.5.11.29.61&167&138& & &3.5.19.31.37&13&44 \cr
218&291885&4.9.11.23.167&61&38&236&326895&8.5.11.13.37&171&236 \cr
 & &16.19.61.167&167&152& & &64.9.19.59&59&96 \cr
\noalign{\hrule}
 & &3.7.13.23.47&25&116& & &9.5.11.23.29&61&38 \cr
219&295113&8.25.23.29&351&374&237&330165&4.5.19.29.61&167&138 \cr
 & &32.27.11.13.17&153&176& & &16.3.19.23.167&167&152 \cr
\noalign{\hrule}
 & &9.13.43.59&25&34& & &7.121.17.23&89&72 \cr
220&296829&4.25.13.17.43&203&528&238&331177&16.9.121.89&105&16 \cr
 & &128.3.7.11.29&319&448& & &512.27.5.7&135&256 \cr
\noalign{\hrule}
 & &49.11.19.29&169&150& & &3.5.49.11.41&37&86 \cr
221&296989&4.3.25.49.169&207&38&239&331485&4.5.11.37.43&329&144 \cr
 & &16.27.5.19.23&135&184& & &128.9.7.47&141&64 \cr
\noalign{\hrule}
 & &3.5.49.11.37&3&52& & &13.17.19.79&101&120 \cr
222&299145&8.9.13.37&23&14&240&331721&16.3.5.79.101&11&90 \cr
 & &32.7.13.23&23&208& & &64.27.25.11&675&352 \cr
\noalign{\hrule}
 & &3.11.13.19.37&23&34& & &9.17.41.53&103&50 \cr
223&301587&4.13.17.23.37&45&436&241&332469&4.25.41.103&51&154 \cr
 & &32.9.5.109&327&80& & &16.3.5.7.11.17&35&88 \cr
\noalign{\hrule}
 & &5.23.43.61&33&10& & &7.19.41.61&559&600 \cr
224&301645&4.3.25.11.61&43&18&242&332633&16.3.25.7.13.43&583&492 \cr
 & &16.27.11.43&27&88& & &128.9.11.41.53&583&576 \cr
\noalign{\hrule}
 & &3.5.17.29.41&5&46& & &3.121.13.71&113&100 \cr
225&303195&4.25.23.29&351&374&243&335049&8.25.121.113&117&4 \cr
 & &16.27.11.13.17&117&88& & &64.9.25.13&25&96 \cr
\noalign{\hrule}
 & &3.11.13.23.31&25&118& & &25.7.19.101&33&68 \cr
226&305877&4.25.23.59&41&18&244&335825&8.3.5.11.17.19&7&12 \cr
 & &16.9.25.41&75&328& & &64.9.7.11.17&187&288 \cr
\noalign{\hrule}
 & &9.11.19.163&23&186& & &9.7.11.17.29&173&184 \cr
227&306603&4.27.23.31&295&326&245&341649&16.3.23.29.173&43&130 \cr
 & &16.5.59.163&59&40& & &64.5.13.23.43&1495&1376 \cr
\noalign{\hrule}
 & &9.13.43.61&375&418& & &9.11.23.151&125&26 \cr
228&306891&4.27.125.11.19&61&34&246&343827&4.125.13.23&151&174 \cr
 & &16.25.11.17.61&187&200& & &16.3.5.29.151&29&40 \cr
\noalign{\hrule}
 & &3.17.23.263&43&26& & &9.7.29.191&845&874 \cr
229&308499&4.13.43.263&153&110&247&348957&4.5.7.169.19.23&165&4 \cr
 & &16.9.5.11.13.17&143&120& & &32.3.25.11.19&209&400 \cr
\noalign{\hrule}
 & &3.121.13.67&441&430& & &3.11.13.19.43&21&22 \cr
230&316173&4.27.5.49.11.43&299&2&248&350493&4.9.7.121.13.19&1535&688 \cr
 & &16.5.7.13.23&161&40& & &128.5.43.307&307&320 \cr
\noalign{\hrule}
 & &3.5.11.19.101&23&78& & &3.25.7.11.61&51&26 \cr
231&316635&4.9.13.19.23&235&64&249&352275&4.9.13.17.61&35&26 \cr
 & &512.5.47&47&256& & &16.5.7.169.17&169&136 \cr
\noalign{\hrule}
 & &9.11.169.19&123&124& & &9.49.11.73&367&290 \cr
232&317889&8.27.11.13.31.41&5&346&250&354123&4.5.7.29.367&285&82 \cr
 & &32.5.41.173&865&656& & &16.3.25.19.41&475&328 \cr
\noalign{\hrule}
 & &25.11.13.89&217&228& & &3.5.11.17.127&199&182 \cr
233&318175&8.3.5.7.13.19.31&131&534&251&356235&4.5.7.11.13.199&127&72 \cr
 & &32.9.89.131&131&144& & &64.9.7.13.127&91&96 \cr
\noalign{\hrule}
 & &7.11.13.17.19&59&150& & &81.5.19.47&7&88 \cr
234&323323&4.3.25.17.59&13&72&252&361665&16.7.11.47&29&18 \cr
 & &64.27.5.13&27&160& & &64.9.7.29&203&32 \cr
\noalign{\hrule}
}%
}
$$
\eject
\vglue -23 pt
\noindent\hskip 1 in\hbox to 6.5 in{\ 253 -- 288 \hfill\fbd 362349 -- 462825\frb}
\vskip -9 pt
$$
\vbox{
\nointerlineskip
\halign{\strut
    \vrule \ \ \hfil \frb #\ 
   &\vrule \hfil \ \ \fbb #\frb\ 
   &\vrule \hfil \ \ \frb #\ \hfil
   &\vrule \hfil \ \ \frb #\ 
   &\vrule \hfil \ \ \frb #\ \ \vrule \hskip 2 pt
   &\vrule \ \ \hfil \frb #\ 
   &\vrule \hfil \ \ \fbb #\frb\ 
   &\vrule \hfil \ \ \frb #\ \hfil
   &\vrule \hfil \ \ \frb #\ 
   &\vrule \hfil \ \ \frb #\ \vrule \cr%
\noalign{\hrule}
 & &9.13.19.163&85&86& & &25.11.17.89&619&894 \cr
253&362349&4.5.13.17.43.163&2783&12&271&416075&4.3.149.619&235&384 \cr
 & &32.3.121.23&121&368& & &1024.9.5.47&423&512 \cr
\noalign{\hrule}
 & &3.5.361.67&91&110& & &25.97.173&99&74 \cr
254&362805&4.25.7.11.13.19&79&54&272&419525&4.9.11.37.97&65&32 \cr
 & &16.27.11.13.79&869&936& & &256.3.5.13.37&481&384 \cr
\noalign{\hrule}
 & &125.29.101&63&38& & &9.5.11.23.37&5&6 \cr
255&366125&4.9.5.7.19.29&101&44&273&421245&4.27.25.23.37&763&88 \cr
 & &32.3.7.11.101&77&48& & &64.7.11.109&109&224 \cr
\noalign{\hrule}
 & &9.125.7.47&539&586& & &9.5.17.557&253&304 \cr
256&370125&4.343.11.293&25&318&274&426105&32.3.5.11.19.23&301&136 \cr
 & &16.3.25.11.53&53&88& & &512.7.17.43&301&256 \cr
\noalign{\hrule}
 & &25.11.23.59&63&52& & &11.17.43.53&313&270 \cr
257&373175&8.9.5.7.13.59&47&12&275&426173&4.27.5.17.313&539&226 \cr
 & &64.27.13.47&611&864& & &16.3.49.11.113&339&392 \cr
\noalign{\hrule}
 & &169.17.131&19&150& & &5.17.61.83&553&858 \cr
258&376363&4.3.25.17.19&33&52&276&430355&4.3.7.11.13.79&85&6 \cr
 & &32.9.5.11.13&45&176& & &16.9.5.11.17&11&72 \cr
\noalign{\hrule}
 & &27.5.7.11.37&89&100& & &9.13.29.127&197&184 \cr
259&384615&8.125.37.89&81&44&277&430911&16.3.23.29.197&55&142 \cr
 & &64.81.11.89&89&96& & &64.5.11.23.71&1633&1760 \cr
\noalign{\hrule}
 & &3.7.11.23.73&11&12& & &27.17.23.41&209&250 \cr
260&387849&8.9.7.121.73&95&752&278&432837&4.125.11.19.23&191&246 \cr
 & &256.5.19.47&893&640& & &16.3.25.41.191&191&200 \cr
\noalign{\hrule}
 & &81.11.19.23&29&40& & &5.7.11.31.37&39&38 \cr
261&389367&16.27.5.19.29&23&4&279&441595&4.3.5.13.19.31.37&1713&1232 \cr
 & &128.5.23.29&145&64& & &128.9.7.11.571&571&576 \cr
\noalign{\hrule}
 & &9.5.11.13.61&369&424& & &9.5.11.19.47&29&18 \cr
262&392535&16.81.41.53&61&20&280&442035&4.81.5.19.29&7&88 \cr
 & &128.5.53.61&53&64& & &64.7.11.29&203&32 \cr
\noalign{\hrule}
 & &9.5.7.31.41&121&166& & &9.19.43.61&55&74 \cr
263&400365&4.121.31.83&129&212&281&448533&4.3.5.11.37.61&23&38 \cr
 & &32.3.11.43.53&583&688& & &16.11.19.23.37&253&296 \cr
\noalign{\hrule}
 & &3.5.121.13.17&103&118& & &9.7.121.59&25&34 \cr
264&401115&4.121.59.103&9&112&282&449757&4.25.7.121.17&43&78 \cr
 & &128.9.7.59&413&192& & &16.3.5.13.17.43&559&680 \cr
\noalign{\hrule}
 & &3.5.7.53.73&79&26& & &27.25.23.29&671&4 \cr
265&406245&4.13.73.79&477&550&283&450225&8.11.61&31&30 \cr
 & &16.9.25.11.53&33&40& & &32.3.5.11.31&31&176 \cr
\noalign{\hrule}
 & &5.7.13.19.47&27&638& & &9.7.13.19.29&35&22 \cr
266&406315&4.27.11.29&49&38&284&451269&4.3.5.49.11.29&19&68 \cr
 & &16.9.49.19&9&56& & &32.5.11.17.19&85&176 \cr
\noalign{\hrule}
 & &7.121.13.37&207&200& & &27.5.11.307&1387&1376 \cr
267&407407&16.9.25.11.13.23&49&16&285&455895&64.3.5.19.43.73&217&2 \cr
 & &512.3.5.49.23&805&768& & &256.7.19.31&589&896 \cr
\noalign{\hrule}
 & &27.25.13.47&319&292& & &9.5.49.11.19&11&46 \cr
268&412425&8.25.11.29.73&7&18&286&460845&4.3.7.121.23&71&50 \cr
 & &32.9.7.29.73&511&464& & &16.25.23.71&355&184 \cr
\noalign{\hrule}
 & &7.13.47.97&69&22& & &81.41.139&29&110 \cr
269&414869&4.3.11.23.97&175&78&287&461619&4.5.11.29.41&139&180 \cr
 & &16.9.25.7.13&25&72& & &32.9.25.139&25&16 \cr
\noalign{\hrule}
 & &27.7.31.71&79&110& & &9.25.121.17&29&4 \cr
270&415989&4.5.11.71.79&233&162&288&462825&8.3.11.17.29&137&50 \cr
 & &16.81.11.233&233&264& & &32.25.137&137&16 \cr
\noalign{\hrule}
}%
}
$$
\eject
\vglue -23 pt
\noindent\hskip 1 in\hbox to 6.5 in{\ 289 -- 324 \hfill\fbd 463095 -- 570843\frb}
\vskip -9 pt
$$
\vbox{
\nointerlineskip
\halign{\strut
    \vrule \ \ \hfil \frb #\ 
   &\vrule \hfil \ \ \fbb #\frb\ 
   &\vrule \hfil \ \ \frb #\ \hfil
   &\vrule \hfil \ \ \frb #\ 
   &\vrule \hfil \ \ \frb #\ \ \vrule \hskip 2 pt
   &\vrule \ \ \hfil \frb #\ 
   &\vrule \hfil \ \ \fbb #\frb\ 
   &\vrule \hfil \ \ \frb #\ \hfil
   &\vrule \hfil \ \ \frb #\ 
   &\vrule \hfil \ \ \frb #\ \vrule \cr%
\noalign{\hrule}
 & &9.5.41.251&1027&1232& & &3.125.23.59&143&202 \cr
289&463095&32.7.11.13.79&85&6&307&508875&4.25.11.13.101&63&38 \cr
 & &128.3.5.11.17&187&64& & &16.9.7.11.13.19&741&616 \cr
\noalign{\hrule}
 & &81.5.19.61&37&98& & &81.25.11.23&61&38 \cr
290&469395&4.3.49.19.37&11&122&308&512325&4.9.25.19.61&143&82 \cr
 & &16.7.11.61&7&88& & &16.11.13.19.41&247&328 \cr
\noalign{\hrule}
 & &7.11.17.361&1357&1170& & &7.17.59.73&215&198 \cr
291&472549&4.9.5.13.23.59&17&22&309&512533&4.9.5.11.43.73&247&118 \cr
 & &16.3.11.17.23.59&177&184& & &16.3.11.13.19.59&247&264 \cr
\noalign{\hrule}
 & &3.5.11.169.17&87&82& & &9.11.169.31&25&118 \cr
292&474045&4.9.11.17.29.41&247&1436&310&518661&4.3.25.13.59&121&56 \cr
 & &32.13.19.359&359&304& & &64.5.7.121&55&224 \cr
\noalign{\hrule}
 & &5.7.19.23.31&27&4& & &81.7.13.71&55&62 \cr
293&474145&8.27.5.7.19&11&46&311&523341&4.9.5.11.31.71&13&58 \cr
 & &32.9.11.23&11&144& & &16.11.13.29.31&319&248 \cr
\noalign{\hrule}
 & &81.11.13.41&61&20& & &9.49.29.41&815&374 \cr
294&474903&8.5.11.13.61&369&424&312&524349&4.5.11.17.163&3&14 \cr
 & &128.9.41.53&53&64& & &16.3.5.7.163&163&40 \cr
\noalign{\hrule}
 & &9.29.31.59&715&184& & &7.13.73.79&477&550 \cr
295&477369&16.5.11.13.23&177&122&313&524797&4.9.25.7.11.53&79&26 \cr
 & &64.3.59.61&61&32& & &16.3.5.11.13.79&33&40 \cr
\noalign{\hrule}
 & &11.169.257&79&90& & &27.101.193&55&248 \cr
296&477763&4.9.5.79.257&247&10&314&526311&16.9.5.11.31&193&148 \cr
 & &16.3.25.13.19&57&200& & &128.37.193&37&64 \cr
\noalign{\hrule}
 & &625.19.41&77&702& & &9.121.17.29&137&50 \cr
297&486875&4.27.7.11.13&25&38&315&536877&4.3.25.11.137&29&4 \cr
 & &16.3.25.11.19&33&8& & &32.29.137&137&16 \cr
\noalign{\hrule}
 & &9.5.11.23.43&61&38& & &9.7.11.19.41&43&736 \cr
298&489555&4.5.19.43.61&69&26&316&539847&64.23.43&33&10 \cr
 & &16.3.13.23.61&61&104& & &256.3.5.11&5&128 \cr
\noalign{\hrule}
 & &17.19.37.41&143&180& & &25.13.23.73&281&294 \cr
299&489991&8.9.5.11.13.41&1&40&317&545675&4.3.49.73.281&115&396 \cr
 & &128.3.25.11&75&704& & &32.27.5.7.11.23&189&176 \cr
\noalign{\hrule}
 & &3.5.7.13.361&17&22& & &27.5.7.19.31&799&286 \cr
300&492765&4.7.11.17.361&1357&1170&318&556605&4.11.13.17.47&87&134 \cr
 & &16.9.5.13.23.59&177&184& & &16.3.11.29.67&319&536 \cr
\noalign{\hrule}
 & &19.841.31&715&126& & &9.5.37.337&1001&664 \cr
301&495349&4.9.5.7.11.13&29&62&319&561105&16.7.11.13.83&113&30 \cr
 & &16.3.5.29.31&5&24& & &64.3.5.7.113&113&224 \cr
\noalign{\hrule}
 & &81.25.13.19&7&88& & &11.841.61&85&756 \cr
302&500175&16.5.7.11.13&27&38&320&564311&8.27.5.7.17&29&22 \cr
 & &64.27.7.19&7&32& & &32.9.5.11.29&45&16 \cr
\noalign{\hrule}
 & &27.5.7.13.41&29&94& & &9.29.41.53&35&88 \cr
303&503685&4.9.7.29.47&55&8&321&567153&16.3.5.7.11.29&41&46 \cr
 & &64.5.11.29&319&32& & &64.7.11.23.41&253&224 \cr
\noalign{\hrule}
 & &3.25.11.13.47&7&18& & &9.13.23.211&255&44 \cr
304&504075&4.27.7.13.47&319&292&322&567801&8.27.5.11.17&157&140 \cr
 & &32.7.11.29.73&511&464& & &64.25.7.157&1099&800 \cr
\noalign{\hrule}
 & &3.5.7.11.19.23&109&52& & &5.13.67.131&763&108 \cr
305&504735&8.5.11.13.109&27&38&323&570505&8.27.7.109&65&44 \cr
 & &32.27.19.109&109&144& & &64.9.5.11.13&99&32 \cr
\noalign{\hrule}
 & &9.5.89.127&41&86& & &9.7.13.17.41&589&958 \cr
306&508635&4.41.43.89&23&66&324&570843&4.19.31.479&55&534 \cr
 & &16.3.11.23.41&451&184& & &16.3.5.11.89&445&88 \cr
\noalign{\hrule}
}%
}
$$
\eject
\vglue -23 pt
\noindent\hskip 1 in\hbox to 6.5 in{\ 325 -- 360 \hfill\fbd 574425 -- 689481\frb}
\vskip -9 pt
$$
\vbox{
\nointerlineskip
\halign{\strut
    \vrule \ \ \hfil \frb #\ 
   &\vrule \hfil \ \ \fbb #\frb\ 
   &\vrule \hfil \ \ \frb #\ \hfil
   &\vrule \hfil \ \ \frb #\ 
   &\vrule \hfil \ \ \frb #\ \ \vrule \hskip 2 pt
   &\vrule \ \ \hfil \frb #\ 
   &\vrule \hfil \ \ \fbb #\frb\ 
   &\vrule \hfil \ \ \frb #\ \hfil
   &\vrule \hfil \ \ \frb #\ 
   &\vrule \hfil \ \ \frb #\ \vrule \cr%
\noalign{\hrule}
 & &27.25.23.37&763&88& & &3.11.13.31.47&25&118 \cr
325&574425&16.7.11.109&5&6&343&625053&4.25.47.59&171&124 \cr
 & &64.3.5.7.109&109&224& & &32.9.5.19.31&95&48 \cr
\noalign{\hrule}
 & &25.7.11.13.23&239&216& & &27.125.11.17&161&26 \cr
326&575575&16.27.5.11.239&37&202&344&631125&4.25.7.13.23&333&242 \cr
 & &64.9.37.101&909&1184& & &16.9.121.37&37&88 \cr
\noalign{\hrule}
 & &3.125.29.53&77&48& & &3.7.17.29.61&29&22 \cr
327&576375&32.9.7.11.53&5&58&345&631533&4.11.841.61&85&756 \cr
 & &128.5.11.29&11&64& & &32.27.5.7.17&45&16 \cr
\noalign{\hrule}
 & &11.13.17.239&1647&1460& & &5.29.41.107&33&74 \cr
328&581009&8.27.5.61.73&481&176&346&636115&4.3.5.11.29.37&107&78 \cr
 & &256.3.11.13.37&111&128& & &16.9.11.13.107&117&88 \cr
\noalign{\hrule}
 & &3.5.7.71.79&1007&652& & &11.17.41.83&711&700 \cr
329&588945&8.19.53.163&55&108&347&636361&8.9.25.7.41.79&31&10 \cr
 & &64.27.5.11.19&209&288& & &32.3.125.31.79&3875&3792 \cr
\noalign{\hrule}
 & &11.31.37.47&189&152& & &81.5.19.83&89&6 \cr
330&592999&16.27.7.19.47&5&52&348&638685&4.243.89&77&166 \cr
 & &128.9.5.7.13&315&832& & &16.7.11.83&7&88 \cr
\noalign{\hrule}
 & &3.293.683&781&98& & &3.5.11.169.23&37&106 \cr
331&600357&4.49.11.71&305&234&349&641355&4.5.13.37.53&437&252 \cr
 & &16.9.5.13.61&195&488& & &32.9.7.19.23&57&112 \cr
\noalign{\hrule}
 & &5.11.13.841&427&414& & &9.5.7.13.157&703&710 \cr
332&601315&4.9.5.7.11.23.61&3&58&350&642915&4.25.13.19.37.71&1&924 \cr
 & &16.27.7.23.29&161&216& & &32.3.7.11.19&19&176 \cr
\noalign{\hrule}
 & &3.25.11.17.43&201&14& & &81.25.11.29&13&68 \cr
333&603075&4.9.5.7.67&163&172&351&645975&8.5.13.17.29&47&18 \cr
 & &32.7.43.163&163&112& & &32.9.17.47&47&272 \cr
\noalign{\hrule}
 & &9.25.37.73&97&122& & &9.5.7.121.17&37&26 \cr
334&607725&4.3.37.61.97&73&110&352&647955&4.5.11.13.17.37&29&114 \cr
 & &16.5.11.73.97&97&88& & &16.3.19.29.37&551&296 \cr
\noalign{\hrule}
 & &27.11.23.89&179&800& & &49.11.23.53&65&12 \cr
335&607959&64.25.179&89&90&353&657041&8.3.5.7.13.23&57&34 \cr
 & &256.9.125.89&125&128& & &32.9.5.17.19&855&272 \cr
\noalign{\hrule}
 & &9.25.11.13.19&17&82& & &27.25.23.43&71&44 \cr
336&611325&4.5.17.19.41&33&52&354&667575&8.5.11.43.71&59&414 \cr
 & &32.3.11.13.41&41&16& & &32.9.23.59&59&16 \cr
\noalign{\hrule}
 & &5.7.13.19.71&53&18& & &9.5.11.19.71&59&154 \cr
337&613795&4.9.13.19.53&35&22&355&667755&4.3.7.121.59&71&50 \cr
 & &16.3.5.7.11.53&159&88& & &16.25.59.71&59&40 \cr
\noalign{\hrule}
 & &81.5.49.31&143&298& & &25.11.31.79&47&822 \cr
338&615195&4.9.11.13.149&25&124&356&673475&4.3.47.137&45&92 \cr
 & &32.25.13.31&65&16& & &32.27.5.23&27&368 \cr
\noalign{\hrule}
 & &3.5.7.11.13.41&305&228& & &25.23.29.41&351&374 \cr
339&615615&8.9.25.19.61&143&82&357&683675&4.27.11.13.17.41&5&46 \cr
 & &32.11.13.19.41&19&16& & &16.9.5.11.13.23&117&88 \cr
\noalign{\hrule}
 & &3.5.11.37.101&35&2& & &3.25.13.19.37&41&16 \cr
340&616605&4.25.7.101&63&38&358&685425&32.13.37.41&261&220 \cr
 & &16.9.49.19&49&456& & &256.9.5.11.29&319&384 \cr
\noalign{\hrule}
 & &11.1369.41&459&910& & &9.5.11.19.73&23&34 \cr
341&617419&4.27.5.7.13.17&37&82&359&686565&4.3.5.17.23.73&29&44 \cr
 & &16.3.13.37.41&13&24& & &32.11.17.23.29&391&464 \cr
\noalign{\hrule}
 & &9.5.7.37.53&143&122& & &9.13.71.83&835&88 \cr
342&617715&4.3.11.13.37.61&1007&214&360&689481&16.5.11.167&89&78 \cr
 & &16.19.53.107&107&152& & &64.3.5.13.89&89&160 \cr
\noalign{\hrule}
}%
}
$$
\eject
\vglue -23 pt
\noindent\hskip 1 in\hbox to 6.5 in{\ 361 -- 396 \hfill\fbd 699567 -- 873015\frb}
\vskip -9 pt
$$
\vbox{
\nointerlineskip
\halign{\strut
    \vrule \ \ \hfil \frb #\ 
   &\vrule \hfil \ \ \fbb #\frb\ 
   &\vrule \hfil \ \ \frb #\ \hfil
   &\vrule \hfil \ \ \frb #\ 
   &\vrule \hfil \ \ \frb #\ \ \vrule \hskip 2 pt
   &\vrule \ \ \hfil \frb #\ 
   &\vrule \hfil \ \ \fbb #\frb\ 
   &\vrule \hfil \ \ \frb #\ \hfil
   &\vrule \hfil \ \ \frb #\ 
   &\vrule \hfil \ \ \frb #\ \vrule \cr%
\noalign{\hrule}
 & &3.11.17.29.43&245&228& & &17.151.313&165&148 \cr
361&699567&8.9.5.49.19.29&187&1118&379&803471&8.3.5.11.37.151&57&94 \cr
 & &32.11.13.17.43&13&16& & &32.9.5.11.19.47&2585&2736 \cr
\noalign{\hrule}
 & &3.11.13.23.71&29&42& & &27.121.13.19&25&146 \cr
362&700557&4.9.7.11.23.29&65&142&380&806949&4.3.25.13.73&77&142 \cr
 & &16.5.13.29.71&29&40& & &16.5.7.11.71&71&280 \cr
\noalign{\hrule}
 & &3.25.13.17.43&1&14& & &3.25.121.89&23&98 \cr
363&712725&4.5.7.17.43&99&116&381&807675&4.49.23.89&125&36 \cr
 & &32.9.7.11.29&203&528& & &32.9.125.7&105&16 \cr
\noalign{\hrule}
 & &5.11.13.19.53&81&62& & &3.5.7.11.19.37&159&26 \cr
364&720005&4.81.5.31.53&551&286&382&811965&4.9.11.13.53&23&76 \cr
 & &16.3.11.13.19.29&29&24& & &32.13.19.23&23&208 \cr
\noalign{\hrule}
 & &3.11.13.23.73&25&14& & &5.11.13.31.37&171&236 \cr
365&720291&4.25.7.23.73&117&44&383&820105&8.9.19.31.59&13&44 \cr
 & &32.9.25.11.13&25&48& & &64.3.11.13.59&59&96 \cr
\noalign{\hrule}
 & &3.125.17.113&1421&1404& & &3.7.19.31.67&649&1118 \cr
366&720375&8.81.5.49.13.29&19&1034&384&828723&4.11.13.43.59&1197&1340 \cr
 & &32.7.11.19.47&893&1232& & &32.9.5.7.19.67&15&16 \cr
\noalign{\hrule}
 & &3.121.43.47&37&84& & &3.5.11.31.163&43&12 \cr
367&733623&8.9.7.37.43&145&188&385&833745&8.9.43.163&275&112 \cr
 & &64.5.7.29.47&203&160& & &256.25.7.11&35&128 \cr
\noalign{\hrule}
 & &81.5.23.79&163&242& & &9.169.19.29&385&356 \cr
368&735885&4.121.23.163&45&208&386&838071&8.3.5.7.11.13.89&19&58 \cr
 & &128.9.5.11.13&143&64& & &32.5.19.29.89&89&80 \cr
\noalign{\hrule}
 & &5.13.47.241&297&908& & &9.25.7.13.41&25&38 \cr
369&736255&8.27.11.227&97&130&387&839475&4.625.19.41&77&702 \cr
 & &32.9.5.13.97&97&144& & &16.27.7.11.13&33&8 \cr
\noalign{\hrule}
 & &9.7.11.23.47&37&26& & &25.11.43.71&63&8 \cr
370&749133&4.13.23.37.47&731&120&388&839575&16.9.5.7.43&29&14 \cr
 & &64.3.5.17.43&215&544& & &64.3.49.29&147&928 \cr
\noalign{\hrule}
 & &3.5.7.121.59&43&78& & &81.5.7.13.23&71&44 \cr
371&749595&4.9.13.43.59&25&34&389&847665&8.3.7.11.13.71&61&152 \cr
 & &16.25.13.17.43&559&680& & &128.11.19.61&1159&704 \cr
\noalign{\hrule}
 & &5.7.13.17.97&11&108& & &27.11.47.61&299&970 \cr
372&750295&8.27.5.11.13&21&34&390&851499&4.5.13.23.97&101&198 \cr
 & &32.81.7.17&81&16& & &16.9.5.11.101&101&40 \cr
\noalign{\hrule}
 & &9.11.13.19.31&25&118& & &5.121.17.83&147&268 \cr
373&758043&4.3.25.19.59&31&26&391&853655&8.3.49.17.67&93&26 \cr
 & &16.5.13.31.59&59&40& & &32.9.7.13.31&403&1008 \cr
\noalign{\hrule}
 & &5.11.101.139&75&64& & &9.5.7.11.13.19&53&64 \cr
374&772145&128.3.125.101&113&12&392&855855&128.5.7.19.53&199&66 \cr
 & &1024.9.113&1017&512& & &512.3.11.199&199&256 \cr
\noalign{\hrule}
 & &9.49.41.43&165&122& & &81.5.11.193&29&164 \cr
375&777483&4.27.5.7.11.61&241&430&393&859815&8.3.11.29.41&37&4 \cr
 & &16.25.43.241&241&200& & &64.29.37&37&928 \cr
\noalign{\hrule}
 & &3.7.17.37.59&1375&1634& & &27.13.23.107&1055&1406 \cr
376&779331&4.125.11.19.43&7&18&394&863811&4.5.19.37.211&13&198 \cr
 & &16.9.5.7.19.43&285&344& & &16.9.11.13.19&19&88 \cr
\noalign{\hrule}
 & &9.125.17.41&411&286& & &9.5.7.41.67&23&44 \cr
377&784125&4.27.11.13.137&55&82&395&865305&8.3.5.11.23.41&119&4 \cr
 & &16.5.121.13.41&121&104& & &64.7.11.17&17&352 \cr
\noalign{\hrule}
 & &81.71.139&29&110& & &3.5.121.13.37&49&16 \cr
378&799389&4.5.11.29.71&37&108&396&873015&32.49.11.37&207&200 \cr
 & &32.27.11.37&37&176& & &512.9.25.7.23&805&768 \cr
\noalign{\hrule}
}%
}
$$
\eject
\vglue -23 pt
\noindent\hskip 1 in\hbox to 6.5 in{\ 397 -- 432 \hfill\fbd 876645 -- 1089789\frb}
\vskip -9 pt
$$
\vbox{
\nointerlineskip
\halign{\strut
    \vrule \ \ \hfil \frb #\ 
   &\vrule \hfil \ \ \fbb #\frb\ 
   &\vrule \hfil \ \ \frb #\ \hfil
   &\vrule \hfil \ \ \frb #\ 
   &\vrule \hfil \ \ \frb #\ \ \vrule \hskip 2 pt
   &\vrule \ \ \hfil \frb #\ 
   &\vrule \hfil \ \ \fbb #\frb\ 
   &\vrule \hfil \ \ \frb #\ \hfil
   &\vrule \hfil \ \ \frb #\ 
   &\vrule \hfil \ \ \frb #\ \vrule \cr%
\noalign{\hrule}
 & &9.5.7.121.23&37&26& & &23.71.599&335&264 \cr
397&876645&4.5.11.13.23.37&219&34&415&978167&16.3.5.11.23.67&311&426 \cr
 & &16.3.13.17.73&221&584& & &64.9.71.311&311&288 \cr
\noalign{\hrule}
 & &3.7.13.41.79&1639&1600& & &5.41.47.103&169&66 \cr
398&884247&128.25.7.11.149&13&162&416&992405&4.3.11.169.41&141&310 \cr
 & &512.81.11.13&297&256& & &16.9.5.31.47&31&72 \cr
\noalign{\hrule}
 & &27.11.23.131&323&70& & &3.17.103.193&77&26 \cr
399&894861&4.9.5.7.17.19&77&94&417&1013829&4.7.11.13.193&103&90 \cr
 & &16.5.49.11.47&245&376& & &16.9.5.7.11.103&105&88 \cr
\noalign{\hrule}
 & &27.23.31.47&715&742& & &9.49.121.19&71&50 \cr
400&904797&4.5.7.11.13.23.53&17&282&418&1013859&4.3.25.7.19.71&11&46 \cr
 & &16.3.7.11.17.47&119&88& & &16.5.11.23.71&355&184 \cr
\noalign{\hrule}
 & &5.121.19.79&65&144& & &3.125.11.13.19&61&186 \cr
401&908105&32.9.25.11.13&7&18&419&1018875&4.9.11.31.61&475&196 \cr
 & &128.81.7.13&1053&448& & &32.25.49.19&49&16 \cr
\noalign{\hrule}
 & &9.5.343.59&143&388& & &7.11.17.19.41&265&1044 \cr
402&910665&8.7.11.13.97&75&68&420&1019711&8.9.5.29.53&41&12 \cr
 & &64.3.25.17.97&485&544& & &64.27.5.41&27&160 \cr
\noalign{\hrule}
 & &9.7.13.19.59&155&22& & &25.11.37.101&63&38 \cr
403&918099&4.3.5.11.13.31&25&118&421&1027675&4.9.7.11.19.37&35&2 \cr
 & &16.125.59&125&8& & &16.3.5.49.19&49&456 \cr
\noalign{\hrule}
 & &27.13.37.71&1199&718& & &27.137.281&55&82 \cr
404&922077&4.11.109.359&125&234&422&1039419&4.5.11.41.281&243&38 \cr
 & &16.9.125.11.13&125&88& & &16.243.11.19&209&72 \cr
\noalign{\hrule}
 & &9.25.11.13.29&47&18& & &5.7.121.13.19&41&54 \cr
405&933075&4.81.5.11.47&13&68&423&1046045&4.27.7.121.41&2491&2470 \cr
 & &32.13.17.47&47&272& & &16.9.5.13.19.47.53&423&424 \cr
\noalign{\hrule}
 & &27.11.47.67&115&182& & &11.13.17.433&999&1432 \cr
406&935253&4.5.7.13.23.47&37&198&424&1052623&16.27.37.179&145&34 \cr
 & &16.9.11.13.37&37&104& & &64.9.5.17.29&145&288 \cr
\noalign{\hrule}
 & &3.5.11.53.107&133&26& & &5.49.11.17.23&71&48 \cr
407&935715&4.5.7.11.13.19&9&86&425&1053745&32.3.5.7.11.71&221&276 \cr
 & &16.9.13.43&559&24& & &256.9.13.17.23&117&128 \cr
\noalign{\hrule}
 & &9.29.59.61&77&106& & &5.7.13.23.101&103&402 \cr
408&939339&4.3.7.11.53.59&203&380&426&1056965&4.3.7.67.103&55&48 \cr
 & &32.5.49.19.29&245&304& & &128.9.5.11.67&737&576 \cr
\noalign{\hrule}
 & &23.67.617&275&342& & &5.7.11.47.59&207&442 \cr
409&950797&4.9.25.11.19.23&67&48&427&1067605&4.9.7.13.17.23&1&22 \cr
 & &128.27.5.11.67&297&320& & &16.3.11.13.17&17&312 \cr
\noalign{\hrule}
 & &3.17.97.193&49&242& & &5.169.31.41&187&18 \cr
410&954771&4.49.121.17&97&90&428&1073995&4.9.11.17.31&41&52 \cr
 & &16.9.5.7.11.97&105&88& & &32.3.13.17.41&17&48 \cr
\noalign{\hrule}
 & &3.11.19.29.53&169&150& & &11.23.31.137&351&362 \cr
411&963699&4.9.25.169.53&893&628&429&1074491&4.27.13.137.181&1705&76 \cr
 & &32.5.19.47.157&785&752& & &32.3.5.11.19.31&57&80 \cr
\noalign{\hrule}
 & &3.11.13.37.61&1007&214& & &9.5.7.11.311&67&32 \cr
412&968253&4.19.53.107&27&80&430&1077615&64.67.311&189&122 \cr
 & &128.27.5.19&95&576& & &256.27.7.61&183&128 \cr
\noalign{\hrule}
 & &3.13.149.167&77&90& & &25.13.47.71&11561&11514 \cr
413&970437&4.27.5.7.11.149&169&20&431&1084525&4.3.11.19.101.1051&1081&30 \cr
 & &32.25.11.169&325&176& & &16.9.5.19.23.47&171&184 \cr
\noalign{\hrule}
 & &121.71.113&117&4& & &3.47.59.131&65&66 \cr
414&970783&8.9.13.71&113&100&432&1089789&4.9.5.11.13.47.59&7&524 \cr
 & &64.3.25.113&25&96& & &32.5.7.13.131&65&112 \cr
\noalign{\hrule}
}%
}
$$
\eject
\vglue -23 pt
\noindent\hskip 1 in\hbox to 6.5 in{\ 433 -- 468 \hfill\fbd 1090635 -- 1306305\frb}
\vskip -9 pt
$$
\vbox{
\nointerlineskip
\halign{\strut
    \vrule \ \ \hfil \frb #\ 
   &\vrule \hfil \ \ \fbb #\frb\ 
   &\vrule \hfil \ \ \frb #\ \hfil
   &\vrule \hfil \ \ \frb #\ 
   &\vrule \hfil \ \ \frb #\ \ \vrule \hskip 2 pt
   &\vrule \ \ \hfil \frb #\ 
   &\vrule \hfil \ \ \fbb #\frb\ 
   &\vrule \hfil \ \ \frb #\ \hfil
   &\vrule \hfil \ \ \frb #\ 
   &\vrule \hfil \ \ \frb #\ \vrule \cr%
\noalign{\hrule}
 & &3.5.7.13.17.47&69&22& & &9.11.79.151&89&10 \cr
433&1090635&4.9.5.11.17.23&611&424&451&1180971&4.5.89.151&73&78 \cr
 & &64.13.47.53&53&32& & &16.3.13.73.89&949&712 \cr
\noalign{\hrule}
 & &27.7.11.17.31&1073&236& & &9.25.11.13.37&1393&1382 \cr
434&1095633&8.29.37.59&15&44&452&1190475&4.3.7.13.199.691&257&2330 \cr
 & &64.3.5.11.37&185&32& & &16.5.7.233.257&1799&1864 \cr
\noalign{\hrule}
 & &81.5.11.13.19&17&28& & &25.7.13.17.31&23&198 \cr
435&1100385&8.9.7.13.17.19&79&92&453&1198925&4.9.11.23.31&205&136 \cr
 & &64.7.17.23.79&2737&2528& & &64.3.5.17.41&123&32 \cr
\noalign{\hrule}
 & &81.5.11.13.19&79&92& & &5.7.11.53.59&65&12 \cr
436&1100385&8.9.5.11.23.79&17&28&454&1203895&8.3.25.13.59&371&396 \cr
 & &64.7.17.23.79&2737&2528& & &64.27.7.11.53&27&32 \cr
\noalign{\hrule}
 & &9.25.7.19.37&103&68& & &25.7.11.17.37&177&452 \cr
437&1107225&8.5.17.37.103&57&572&455&1210825&8.3.7.59.113&263&150 \cr
 & &64.3.11.13.19&143&32& & &32.9.25.263&263&144 \cr
\noalign{\hrule}
 & &49.13.37.47&551&1188& & &9.5.49.19.29&187&1118 \cr
438&1107743&8.27.11.19.29&5&14&456&1214955&4.11.13.17.43&245&228 \cr
 & &32.3.5.7.11.29&435&176& & &32.3.5.49.13.19&13&16 \cr
\noalign{\hrule}
 & &3.25.59.251&273&22& & &9.23.59.101&425&484 \cr
439&1110675&4.9.5.7.11.13&251&134&457&1233513&8.25.121.17.23&3&118 \cr
 & &16.67.251&67&8& & &32.3.5.17.59&85&16 \cr
\noalign{\hrule}
 & &25.49.11.83&81&164& & &27.5.13.19.37&493&506 \cr
440&1118425&8.81.5.11.41&17&28&458&1233765&4.5.11.17.19.23.29&191&246 \cr
 & &64.9.7.17.41&369&544& & &16.3.17.29.41.191&5539&5576 \cr
\noalign{\hrule}
 & &3.5.13.53.109&187&78& & &3.5.7.11.29.37&57&202 \cr
441&1126515&4.9.11.169.17&89&80&459&1239315&4.9.11.19.101&559&350 \cr
 & &128.5.11.17.89&979&1088& & &16.25.7.13.43&215&104 \cr
\noalign{\hrule}
 & &5.11.67.311&189&122& & &3.5.7.11.23.47&493&494 \cr
442&1146035&4.27.5.7.11.61&67&32&460&1248555&4.5.11.13.17.19.23.29&801&134 \cr
 & &256.3.61.67&183&128& & &16.9.13.19.67.89&5963&5928 \cr
\noalign{\hrule}
 & &5.49.31.151&47&198& & &27.5.17.19.29&1529&1036 \cr
443&1146845&4.9.11.31.47&23&70&461&1264545&8.7.11.37.139&199&60 \cr
 & &16.3.5.7.11.23&33&184& & &64.3.5.11.199&199&352 \cr
\noalign{\hrule}
 & &11.23.47.97&175&78& & &5.11.13.29.61&3&58 \cr
444&1153427&4.3.25.7.13.47&69&22&462&1264835&4.3.13.841&427&414 \cr
 & &16.9.25.11.23&25&72& & &16.27.7.23.61&161&216 \cr
\noalign{\hrule}
 & &3.5.7.11.19.53&149&434& & &9.5.7.13.313&803&1388 \cr
445&1163085&4.49.31.149&99&50&463&1281735&8.11.73.347&137&210 \cr
 & &16.9.25.11.31&93&40& & &32.3.5.7.11.137&137&176 \cr
\noalign{\hrule}
 & &9.7.11.23.73&1175&1102& & &169.29.263&99&70 \cr
446&1163547&4.25.7.19.29.47&3&32&464&1288963&4.9.5.7.11.263&247&16 \cr
 & &256.3.5.19.47&893&640& & &128.3.5.13.19&57&320 \cr
\noalign{\hrule}
 & &5.11.17.29.43&245&228& & &3.5.13.29.229&187&42 \cr
447&1165945&8.3.25.49.19.29&73&102&465&1294995&4.9.7.11.13.17&41&50 \cr
 & &32.9.7.17.19.73&1197&1168& & &16.25.11.17.41&451&680 \cr
\noalign{\hrule}
 & &5.7.19.41.43&2619&2834& & &25.17.43.71&1141&66 \cr
448&1172395&4.27.13.97.109&55&42&466&1297525&4.3.7.11.163&85&78 \cr
 & &16.81.5.7.11.109&891&872& & &16.9.5.11.13.17&143&72 \cr
\noalign{\hrule}
 & &27.5.11.13.61&817&830& & &17.23.47.71&825&808 \cr
449&1177605&4.25.11.19.43.83&1&474&467&1304767&16.3.25.11.47.101&51&4 \cr
 & &16.3.79.83&83&632& & &128.9.5.17.101&505&576 \cr
\noalign{\hrule}
 & &3.13.19.37.43&47&10& & &9.5.7.11.13.29&79&124 \cr
450&1178931&4.5.13.43.47&413&198&468&1306305&8.11.13.31.79&87&56 \cr
 & &16.9.7.11.59&649&168& & &128.3.7.29.79&79&64 \cr
\noalign{\hrule}
}%
}
$$
\eject
\vglue -23 pt
\noindent\hskip 1 in\hbox to 6.5 in{\ 469 -- 504 \hfill\fbd 1307859 -- 1517549\frb}
\vskip -9 pt
$$
\vbox{
\nointerlineskip
\halign{\strut
    \vrule \ \ \hfil \frb #\ 
   &\vrule \hfil \ \ \fbb #\frb\ 
   &\vrule \hfil \ \ \frb #\ \hfil
   &\vrule \hfil \ \ \frb #\ 
   &\vrule \hfil \ \ \frb #\ \ \vrule \hskip 2 pt
   &\vrule \ \ \hfil \frb #\ 
   &\vrule \hfil \ \ \fbb #\frb\ 
   &\vrule \hfil \ \ \frb #\ \hfil
   &\vrule \hfil \ \ \frb #\ 
   &\vrule \hfil \ \ \frb #\ \vrule \cr%
\noalign{\hrule}
 & &3.343.31.41&233&110& & &3.7.17.41.97&2585&2294 \cr
469&1307859&4.5.11.31.233&287&54&487&1419789&4.5.11.31.37.47&289&1746 \cr
 & &16.27.5.7.41&9&40& & &16.9.289.97&17&24 \cr
\noalign{\hrule}
 & &9.5.7.13.17.19&61&194& & &3.7.121.13.43&71&50 \cr
470&1322685&4.3.13.61.97&539&722&488&1420419&4.25.13.43.71&999&76 \cr
 & &16.49.11.361&133&88& & &32.27.19.37&703&144 \cr
\noalign{\hrule}
 & &25.11.47.103&29&18& & &3.5.11.41.211&37&4 \cr
471&1331275&4.9.25.29.103&517&208&489&1427415&8.5.37.211&13&198 \cr
 & &128.3.11.13.47&39&64& & &32.9.11.13&3&208 \cr
\noalign{\hrule}
 & &5.19.23.617&67&48& & &27.5.13.19.43&61&34 \cr
472&1348145&32.3.67.617&275&342&490&1433835&4.13.17.43.61&375&418 \cr
 & &128.27.25.11.19&297&320& & &16.3.125.11.17.19&187&200 \cr
\noalign{\hrule}
 & &3.7.11.13.449&1577&1566& & &9.25.7.11.83&17&28 \cr
473&1348347&4.81.13.19.29.83&1309&230&491&1437975&8.5.49.17.83&81&164 \cr
 & &16.5.7.11.17.23.29&667&680& & &64.81.17.41&369&544 \cr
\noalign{\hrule}
 & &9.121.17.73&1065&992& & &3.5.7.11.29.43&135&338 \cr
474&1351449&64.27.5.31.71&1033&1168&492&1440285&4.81.25.169&203&122 \cr
 & &2048.73.1033&1033&1024& & &16.7.13.29.61&61&104 \cr
\noalign{\hrule}
 & &25.7.11.19.37&611&1314& & &9.11.13.19.59&307&320 \cr
475&1353275&4.9.13.47.73&77&64&493&1442727&128.3.5.59.307&65&242 \cr
 & &512.3.7.11.73&219&256& & &512.25.121.13&275&256 \cr
\noalign{\hrule}
 & &9.5.7.11.17.23&611&424& & &27.11.43.113&1843&1886 \cr
476&1354815&16.7.13.47.53&69&22&494&1443123&4.9.19.23.41.97&3955&22 \cr
 & &64.3.11.23.53&53&32& & &16.5.7.11.113&35&8 \cr
\noalign{\hrule}
 & &27.25.49.41&481&194& & &9.23.29.241&89&118 \cr
477&1356075&4.7.13.37.97&99&580&495&1446723&4.59.89.241&165&76 \cr
 & &32.9.5.11.29&29&176& & &32.3.5.11.19.59&1121&880 \cr
\noalign{\hrule}
 & &3.5.7.13.19.53&199&66& & &3.5.13.17.19.23&11&28 \cr
478&1374555&4.9.11.13.199&53&64&496&1448655&8.5.7.11.19.23&603&442 \cr
 & &512.53.199&199&256& & &32.9.13.17.67&67&48 \cr
\noalign{\hrule}
 & &3.49.47.199&173&26& & &7.23.47.193&153&176 \cr
479&1374891&4.13.47.173&63&110&497&1460431&32.9.11.17.193&3&190 \cr
 & &16.9.5.7.11.13&65&264& & &128.27.5.19&135&1216 \cr
\noalign{\hrule}
 & &125.11.17.59&71&54& & &9.25.17.389&187&202 \cr
480&1379125&4.27.11.59.71&125&656&498&1487925&4.3.5.11.289.101&2723&1612 \cr
 & &128.3.125.41&123&64& & &32.7.13.31.389&217&208 \cr
\noalign{\hrule}
 & &27.17.23.131&1573&1964& & &9.49.17.199&159&40 \cr
481&1382967&8.121.13.491&185&306&499&1491903&16.27.5.7.53&121&68 \cr
 & &32.9.5.13.17.37&185&208& & &128.5.121.17&121&320 \cr
\noalign{\hrule}
 & &9.25.23.269&247&22& & &243.11.13.43&505&548 \cr
482&1392075&4.11.13.19.23&135&112&500&1494207&8.3.5.11.101.137&151&14 \cr
 & &128.27.5.7.11&77&192& & &32.7.101.151&1057&1616 \cr
\noalign{\hrule}
 & &27.11.13.361&2495&2198& & &3.13.137.281&209&72 \cr
483&1393821&4.5.7.157.499&171&328&501&1501383&16.27.11.13.19&137&110 \cr
 & &64.9.5.7.19.41&205&224& & &64.5.121.137&121&160 \cr
\noalign{\hrule}
 & &25.11.13.17.23&2313&1988& & &9.125.13.103&539&436 \cr
484&1397825&8.9.7.71.257&37&34&502&1506375&8.3.5.49.11.109&107&2 \cr
 & &32.3.7.17.37.257&1799&1776& & &32.7.11.107&107&1232 \cr
\noalign{\hrule}
 & &9.7.13.29.59&2047&1670& & &81.11.13.131&25&106 \cr
485&1401309&4.5.23.89.167&429&406&503&1517373&4.25.11.13.53&131&144 \cr
 & &16.3.7.11.13.29.89&89&88& & &128.9.53.131&53&64 \cr
\noalign{\hrule}
 & &5.13.19.31.37&1713&1232& & &11.19.53.137&3735&3526 \cr
486&1416545&32.3.7.11.571&39&38&504&1517549&4.9.5.41.43.83&85&44 \cr
 & &128.9.13.19.571&571&576& & &32.3.25.11.17.83&1275&1328 \cr
\noalign{\hrule}
}%
}
$$
\eject
\vglue -23 pt
\noindent\hskip 1 in\hbox to 6.5 in{\ 505 -- 540 \hfill\fbd 1530765 -- 1743525\frb}
\vskip -9 pt
$$
\vbox{
\nointerlineskip
\halign{\strut
    \vrule \ \ \hfil \frb #\ 
   &\vrule \hfil \ \ \fbb #\frb\ 
   &\vrule \hfil \ \ \frb #\ \hfil
   &\vrule \hfil \ \ \frb #\ 
   &\vrule \hfil \ \ \frb #\ \ \vrule \hskip 2 pt
   &\vrule \ \ \hfil \frb #\ 
   &\vrule \hfil \ \ \fbb #\frb\ 
   &\vrule \hfil \ \ \frb #\ \hfil
   &\vrule \hfil \ \ \frb #\ 
   &\vrule \hfil \ \ \frb #\ \vrule \cr%
\noalign{\hrule}
 & &27.5.17.23.29&737&742& & &11.1681.89&351&1330 \cr
505&1530765&4.9.7.11.23.53.67&395&188&523&1645699&4.27.5.7.13.19&41&22 \cr
 & &32.5.7.47.67.79&5293&5264& & &16.3.5.11.13.41&65&24 \cr
\noalign{\hrule}
 & &3.13.23.29.59&875&836& & &3.7.11.13.19.29&29&62 \cr
506&1534767&8.125.7.11.19.23&143&18&524&1654653&4.19.841.31&715&126 \cr
 & &32.9.121.13.19&363&304& & &16.9.5.7.11.13&5&24 \cr
\noalign{\hrule}
 & &9.11.13.17.71&3625&3404& & &3.5.173.641&323&318 \cr
507&1553409&8.125.23.29.37&27&2&525&1663395&4.9.17.19.53.173&275&1282 \cr
 & &32.27.5.23.37&555&368& & &16.25.11.17.641&85&88 \cr
\noalign{\hrule}
 & &9.19.53.173&275&1282& & &3.11.29.37.47&273&244 \cr
508&1567899&4.25.11.641&323&318&526&1664223&8.9.7.13.37.61&47&380 \cr
 & &16.3.5.11.17.19.53&85&88& & &64.5.13.19.47&247&160 \cr
\noalign{\hrule}
 & &81.7.17.163&611&530& & &9.5.7.11.13.37&137&248 \cr
509&1571157&4.5.13.17.47.53&55&744&527&1666665&16.3.13.31.137&407&4 \cr
 & &64.3.25.11.31&341&800& & &128.11.37&1&64 \cr
\noalign{\hrule}
 & &27.25.17.137&55&82& & &9.25.17.19.23&1067&1118 \cr
510&1572075&4.125.11.17.41&411&286&528&1671525&4.3.5.11.13.43.97&253&38 \cr
 & &16.3.121.13.137&121&104& & &16.121.13.19.23&121&104 \cr
\noalign{\hrule}
 & &7921.199&3861&4060& & &3.25.11.2029&1027&1002 \cr
511&1576279&8.27.5.7.11.13.29&89&2&529&1673925&4.9.11.13.79.167&727&560 \cr
 & &32.9.5.11.89&55&144& & &128.5.7.79.727&5089&5056 \cr
\noalign{\hrule}
 & &5.49.11.19.31&9&86& & &11.29.59.89&861&850 \cr
512&1587355&4.9.7.31.43&347&304&530&1675069&4.3.25.7.17.41.89&87&2 \cr
 & &128.3.19.347&347&192& & &16.9.5.7.29.41&369&280 \cr
\noalign{\hrule}
 & &5.11.127.229&87&142& & &9.19.41.241&4931&4950 \cr
513&1599565&4.3.29.71.127&99&28&531&1689651&4.81.25.11.4931&4693&238 \cr
 & &32.27.7.11.29&203&432& & &16.5.7.13.17.361&1105&1064 \cr
\noalign{\hrule}
 & &5.37.41.211&13&198& & &27.11.29.197&577&380 \cr
514&1600435&4.9.11.13.41&37&4&532&1696761&8.9.5.19.577&203&374 \cr
 & &32.3.13.37&3&208& & &32.5.7.11.17.29&85&112 \cr
\noalign{\hrule}
 & &9.5.7.13.17.23&97&488& & &81.5.13.17.19&35&22 \cr
515&1601145&16.7.61.97&165&262&533&1700595&4.27.25.7.11.17&5699&5776 \cr
 & &64.3.5.11.131&131&352& & &128.361.41.139&2641&2624 \cr
\noalign{\hrule}
 & &5.11.19.29.53&147&118& & &5.7.11.17.263&91&96 \cr
516&1606165&4.3.49.11.19.59&291&830&534&1721335&64.3.49.13.263&187&450 \cr
 & &16.9.5.83.97&873&664& & &256.27.25.11.17&135&128 \cr
\noalign{\hrule}
 & &5.11.23.31.41&14589&14644& & &9.7.151.181&589&770 \cr
517&1607815&8.9.7.523.1621&1595&26&535&1721853&4.5.49.11.19.31&195&146 \cr
 & &32.3.5.7.11.13.29&377&336& & &16.3.25.13.19.73&1825&1976 \cr
\noalign{\hrule}
 & &9.25.13.19.29&101&146& & &27.5.7.31.59&77&78 \cr
518&1611675&4.5.29.73.101&123&22&536&1728405&4.81.49.11.13.59&1891&1000 \cr
 & &16.3.11.41.73&803&328& & &64.125.13.31.61&793&800 \cr
\noalign{\hrule}
 & &81.13.29.53&1295&242& & &3.11.13.37.109&53&90 \cr
519&1618461&4.5.7.121.37&27&28&537&1730157&4.27.5.53.109&259&286 \cr
 & &32.27.49.11.37&539&592& & &16.7.11.13.37.53&53&56 \cr
\noalign{\hrule}
 & &3.5.7.11.23.61&17&52& & &9.121.37.43&145&188 \cr
520&1620465&8.11.13.17.61&25&36&538&1732599&8.5.121.29.47&37&84 \cr
 & &64.9.25.13.17&255&416& & &64.3.5.7.29.37&203&160 \cr
\noalign{\hrule}
 & &11.13.17.23.29&315&62& & &7.121.23.89&125&36 \cr
521&1621477&4.9.5.7.17.31&27&58&539&1733809&8.9.125.121&23&98 \cr
 & &16.243.7.29&243&56& & &32.3.5.49.23&105&16 \cr
\noalign{\hrule}
 & &3.5.313.349&403&1342& & &243.25.7.41&1363&338 \cr
522&1638555&4.11.13.31.61&171&232&540&1743525&4.169.29.47&99&70 \cr
 & &64.9.11.19.29&957&608& & &16.9.5.7.11.47&47&88 \cr
\noalign{\hrule}
}%
}
$$
\eject
\vglue -23 pt
\noindent\hskip 1 in\hbox to 6.5 in{\ 541 -- 576 \hfill\fbd 1749425 -- 1987557\frb}
\vskip -9 pt
$$
\vbox{
\nointerlineskip
\halign{\strut
    \vrule \ \ \hfil \frb #\ 
   &\vrule \hfil \ \ \fbb #\frb\ 
   &\vrule \hfil \ \ \frb #\ \hfil
   &\vrule \hfil \ \ \frb #\ 
   &\vrule \hfil \ \ \frb #\ \ \vrule \hskip 2 pt
   &\vrule \ \ \hfil \frb #\ 
   &\vrule \hfil \ \ \fbb #\frb\ 
   &\vrule \hfil \ \ \frb #\ \hfil
   &\vrule \hfil \ \ \frb #\ 
   &\vrule \hfil \ \ \frb #\ \vrule \cr%
\noalign{\hrule}
 & &25.19.29.127&327&308& & &27.13.19.281&137&110 \cr
541&1749425&8.3.5.7.11.29.109&1107&92&559&1873989&4.5.11.137.281&209&72 \cr
 & &64.81.23.41&1863&1312& & &64.9.5.121.19&121&160 \cr
\noalign{\hrule}
 & &9.25.7.11.101&1591&884& & &9.7.121.13.19&1535&688 \cr
542&1749825&8.13.17.37.43&125&606&560&1882881&32.5.43.307&21&22 \cr
 & &32.3.125.101&5&16& & &128.3.5.7.11.307&307&320 \cr
\noalign{\hrule}
 & &7.11.13.17.103&207&1340& & &3.5.11.41.281&119&86 \cr
543&1752751&8.9.5.23.67&139&206&561&1900965&4.7.17.43.281&225&506 \cr
 & &32.3.103.139&139&48& & &16.9.25.7.11.23&161&120 \cr
\noalign{\hrule}
 & &9.25.11.23.31&109&98& & &9.125.19.89&371&104 \cr
544&1764675&4.25.49.31.109&769&6&562&1902375&16.3.5.7.13.53&283&88 \cr
 & &16.3.7.769&769&56& & &256.11.283&283&1408 \cr
\noalign{\hrule}
 & &27.7.11.23.37&169&238& & &9.7.13.17.137&355&1426 \cr
545&1769229&4.9.49.169.17&331&110&563&1907451&4.5.23.31.71&89&66 \cr
 & &16.5.11.13.331&331&520& & &16.3.11.71.89&979&568 \cr
\noalign{\hrule}
 & &27.7.11.23.37&505&494& & &9.5.7.11.19.29&289&376 \cr
546&1769229&4.5.7.13.19.23.101&737&576&564&1909215&16.3.11.289.47&383&416 \cr
 & &512.9.5.11.19.67&1273&1280& & &1024.13.17.383&6511&6656 \cr
\noalign{\hrule}
 & &27.7.83.113&61&52& & &3.25.11.23.101&63&38 \cr
547&1772631&8.3.7.13.61.83&3685&3868&565&1916475&4.27.7.11.19.23&383&130 \cr
 & &64.5.11.67.967&10637&10720& & &16.5.7.13.383&383&728 \cr
\noalign{\hrule}
 & &3.25.11.13.167&79&246& & &3.5.23.67.83&107&308 \cr
548&1791075&4.9.11.41.79&29&70&566&1918545&8.7.11.23.107&27&134 \cr
 & &16.5.7.29.79&203&632& & &32.27.11.67&99&16 \cr
\noalign{\hrule}
 & &3.37.67.241&87&154& & &3.13.31.37.43&295&264 \cr
549&1792317&4.9.7.11.29.37&65&268&567&1923519&16.9.5.11.37.59&469&62 \cr
 & &32.5.11.13.67&65&176& & &64.5.7.31.67&335&224 \cr
\noalign{\hrule}
 & &9.7.11.23.113&635&382& & &3.7.11.61.137&145&282 \cr
550&1801107&4.5.7.127.191&113&78&568&1930467&4.9.5.11.29.47&23&122 \cr
 & &16.3.13.113.127&127&104& & &16.23.47.61&47&184 \cr
\noalign{\hrule}
 & &9.5.49.19.43&377&358& & &7.37.73.103&1711&990 \cr
551&1801485&4.3.13.29.43.179&25&154&569&1947421&4.9.5.11.29.59&103&74 \cr
 & &16.25.7.11.13.29&377&440& & &16.3.5.11.37.103&33&40 \cr
\noalign{\hrule}
 & &11.2809.59&1729&1080& & &3.5.7.11.19.89&181&104 \cr
552&1823041&16.27.5.7.13.19&53&118&570&1953105&16.13.89.181&135&46 \cr
 & &64.3.7.53.59&21&32& & &64.27.5.13.23&207&416 \cr
\noalign{\hrule}
 & &9.49.41.101&2291&1850& & &9.31.43.163&275&112 \cr
553&1826181&4.25.29.37.79&451&1524&571&1955511&32.25.7.11.31&43&12 \cr
 & &32.3.11.41.127&127&176& & &256.3.5.7.43&35&128 \cr
\noalign{\hrule}
 & &9.5.23.29.61&31&30& & &125.13.17.71&209&1416 \cr
554&1830915&4.27.25.23.29.31&671&4&572&1961375&16.3.11.19.59&75&134 \cr
 & &32.11.31.61&31&176& & &64.9.25.67&67&288 \cr
\noalign{\hrule}
 & &19.179.541&1971&1430& & &3.11.17.31.113&65&48 \cr
555&1839941&4.27.5.11.13.73&181&38&573&1965183&32.9.5.11.13.31&449&46 \cr
 & &16.9.5.19.181&181&360& & &128.23.449&449&1472 \cr
\noalign{\hrule}
 & &3.25.11.23.97&61&36& & &3.5.11.79.151&73&78 \cr
556&1840575&8.27.11.23.61&25&646&574&1968285&4.9.11.13.73.79&89&10 \cr
 & &32.25.17.19&19&272& & &16.5.13.73.89&949&712 \cr
\noalign{\hrule}
 & &3.5.49.17.149&649&394& & &5.11.13.47.59&171&124 \cr
557&1861755&4.7.11.59.197&9&68&575&1982695&8.9.11.13.19.31&25&118 \cr
 & &32.9.17.197&197&48& & &32.3.25.19.59&95&48 \cr
\noalign{\hrule}
 & &5.13.19.37.41&261&220& & &3.11.13.41.113&1411&58 \cr
558&1873495&8.9.25.11.19.29&41&16&576&1987557&4.17.29.83&205&288 \cr
 & &256.3.11.29.41&319&384& & &256.9.5.41&15&128 \cr
\noalign{\hrule}
}%
}
$$
\eject
\vglue -23 pt
\noindent\hskip 1 in\hbox to 6.5 in{\ 577 -- 612 \hfill\fbd 1992681 -- 2329509\frb}
\vskip -9 pt
$$
\vbox{
\nointerlineskip
\halign{\strut
    \vrule \ \ \hfil \frb #\ 
   &\vrule \hfil \ \ \fbb #\frb\ 
   &\vrule \hfil \ \ \frb #\ \hfil
   &\vrule \hfil \ \ \frb #\ 
   &\vrule \hfil \ \ \frb #\ \ \vrule \hskip 2 pt
   &\vrule \ \ \hfil \frb #\ 
   &\vrule \hfil \ \ \fbb #\frb\ 
   &\vrule \hfil \ \ \frb #\ \hfil
   &\vrule \hfil \ \ \frb #\ 
   &\vrule \hfil \ \ \frb #\ \vrule \cr%
\noalign{\hrule}
 & &81.73.337&497&160& & &17.29.43.101&1925&2418 \cr
577&1992681&64.9.5.7.71&209&146&595&2141099&4.3.25.7.11.13.31&29&36 \cr
 & &256.11.19.73&209&128& & &32.27.5.11.29.31&837&880 \cr
\noalign{\hrule}
 & &81.5.7.19.37&11&122& & &9.5.29.31.53&169&1474 \cr
578&1993005&4.27.5.11.61&37&98&596&2144115&4.11.169.67&7&6 \cr
 & &16.49.11.37&7&88& & &16.3.7.11.13.67&1001&536 \cr
\noalign{\hrule}
 & &9.25.289.31&377&88& & &3.5.11.31.421&205&216 \cr
579&2015775&16.3.5.11.13.29&257&62&597&2153415&16.81.25.31.41&377&1648 \cr
 & &64.31.257&257&32& & &512.13.29.103&2987&3328 \cr
\noalign{\hrule}
 & &27.11.13.17.31&25&118& & &5.7.11.29.193&513&502 \cr
580&2034747&4.9.25.17.59&13&4&598&2154845&4.27.19.193.251&29&222 \cr
 & &32.25.13.59&59&400& & &16.81.19.29.37&703&648 \cr
\noalign{\hrule}
 & &243.5.19.89&77&166& & &3.5.169.23.37&275&206 \cr
581&2054565&4.5.7.11.19.83&89&6&599&2157285&4.125.11.13.103&1357&18 \cr
 & &16.3.7.11.89&7&88& & &16.9.23.59&177&8 \cr
\noalign{\hrule}
 & &89.101.229&10241&10140& & &25.11.29.271&147&172 \cr
582&2058481&8.3.5.49.11.169.19&47&96&600&2161225&8.3.49.43.271&15&286 \cr
 & &512.9.5.13.19.47&11609&11520& & &32.9.5.7.11.13&63&208 \cr
\noalign{\hrule}
 & &9.25.49.11.17&83&358& & &9.5.7.11.17.37&37&82 \cr
583&2061675&4.17.83.179&795&616&601&2179485&4.11.1369.41&459&910 \cr
 & &64.3.5.7.11.53&53&32& & &16.27.5.7.13.17&13&24 \cr
\noalign{\hrule}
 & &9.5.17.37.73&209&124& & &9.7.47.739&401&338 \cr
584&2066265&8.11.19.31.73&27&46&602&2188179&4.169.47.401&105&506 \cr
 & &32.27.11.23.31&713&528& & &16.3.5.7.11.13.23&299&440 \cr
\noalign{\hrule}
 & &3.5.121.17.67&3029&2666& & &3.125.67.89&11&78 \cr
585&2067285&4.13.31.43.233&783&550&603&2236125&4.9.125.11.13&67&58 \cr
 & &16.27.25.11.13.29&377&360& & &16.11.13.29.67&143&232 \cr
\noalign{\hrule}
 & &3.5.7.13.37.41&11&544& & &7.11.23.31.41&1935&1222 \cr
586&2070705&64.7.11.17&97&90&604&2250941&4.9.5.13.43.47&41&88 \cr
 & &256.9.5.97&97&384& & &64.3.5.11.13.41&65&96 \cr
\noalign{\hrule}
 & &25.11.71.107&553&624& & &9.11.13.17.103&43&60 \cr
587&2089175&32.3.25.7.13.79&601&426&605&2253537&8.27.5.11.13.43&61&412 \cr
 & &128.9.71.601&601&576& & &64.5.61.103&61&160 \cr
\noalign{\hrule}
 & &49.19.31.73&319&270& & &25.19.47.101&67&168 \cr
588&2106853&4.27.5.11.29.73&271&532&606&2254825&16.3.5.7.19.67&187&282 \cr
 & &32.3.5.7.19.271&271&240& & &64.9.11.17.47&187&288 \cr
\noalign{\hrule}
 & &7.121.47.53&397&450& & &3.25.67.449&1089&1156 \cr
589&2109877&4.9.25.47.397&551&154&607&2256225&8.27.5.121.289&161&26 \cr
 & &16.3.5.7.11.19.29&285&232& & &32.7.11.13.17.23&2737&2288 \cr
\noalign{\hrule}
 & &27.5.49.11.29&19&184& & &5.17.67.397&423&1562 \cr
590&2110185&16.9.7.19.23&715&734&608&2260915&4.9.11.47.71&257&260 \cr
 & &64.5.11.13.367&367&416& & &32.3.5.13.71.257&3341&3408 \cr
\noalign{\hrule}
 & &9.49.11.19.23&965&916& & &9.25.13.19.41&61&308 \cr
591&2119887&8.5.23.193.229&39&154&609&2278575&8.25.7.11.61&13&12 \cr
 & &32.3.7.11.13.229&229&208& & &64.3.7.11.13.61&427&352 \cr
\noalign{\hrule}
 & &3.5.31.47.97&401&304& & &3.5.7.19.31.37&17&572 \cr
592&2119935&32.19.31.401&495&94&610&2288265&8.7.11.13.17&37&54 \cr
 & &128.9.5.11.47&33&64& & &32.27.11.37&11&144 \cr
\noalign{\hrule}
 & &5.11.13.29.103&19&84& & &9.5.7.23.317&1075&1144 \cr
593&2135705&8.3.7.11.19.29&3&206&611&2296665&16.3.125.11.13.43&49&424 \cr
 & &32.9.103&1&144& & &256.49.13.53&689&896 \cr
\noalign{\hrule}
 & &3.5.49.41.71&11&60& & &3.49.13.23.53&59&10 \cr
594&2139585&8.9.25.11.41&133&92&612&2329509&4.5.13.53.59&315&374 \cr
 & &64.7.11.19.23&253&608& & &16.9.25.7.11.17&425&264 \cr
\noalign{\hrule}
}%
}
$$
\eject
\vglue -23 pt
\noindent\hskip 1 in\hbox to 6.5 in{\ 613 -- 648 \hfill\fbd 2331901 -- 2637915\frb}
\vskip -9 pt
$$
\vbox{
\nointerlineskip
\halign{\strut
    \vrule \ \ \hfil \frb #\ 
   &\vrule \hfil \ \ \fbb #\frb\ 
   &\vrule \hfil \ \ \frb #\ \hfil
   &\vrule \hfil \ \ \frb #\ 
   &\vrule \hfil \ \ \frb #\ \ \vrule \hskip 2 pt
   &\vrule \ \ \hfil \frb #\ 
   &\vrule \hfil \ \ \fbb #\frb\ 
   &\vrule \hfil \ \ \frb #\ \hfil
   &\vrule \hfil \ \ \frb #\ 
   &\vrule \hfil \ \ \frb #\ \vrule \cr%
\noalign{\hrule}
 & &11.13.23.709&4735&4482& & &7.13.19.31.47&11&600 \cr
613&2331901&4.27.5.83.947&349&598&631&2519153&16.3.25.7.11&93&82 \cr
 & &16.9.5.13.23.349&349&360& & &64.9.31.41&41&288 \cr
\noalign{\hrule}
 & &9.11.67.353&83&16& & &243.11.13.73&181&38 \cr
614&2341449&32.83.353&135&218&632&2536677&4.81.19.181&1679&1760 \cr
 & &128.27.5.109&545&192& & &256.5.11.23.73&115&128 \cr
\noalign{\hrule}
 & &29.1369.59&171&1540& & &3.11.23.31.109&19&50 \cr
615&2342359&8.9.5.7.11.19&37&58&633&2564661&4.25.11.19.109&837&362 \cr
 & &32.3.11.29.37&11&48& & &16.27.31.181&181&72 \cr
\noalign{\hrule}
 & &27.5.11.37.43&119&526& & &27.5.7.11.13.19&47&86 \cr
616&2362635&4.9.7.17.263&23&40&634&2567565&4.9.5.11.43.47&221&166 \cr
 & &64.5.23.263&263&736& & &16.13.17.47.83&799&664 \cr
\noalign{\hrule}
 & &9.7.11.19.181&95&86& & &7.289.31.41&407&120 \cr
617&2383227&4.5.7.11.361.43&219&142&635&2571233&16.3.5.11.17.37&1&186 \cr
 & &16.3.5.43.71.73&3053&2920& & &64.9.31&1&288 \cr
\noalign{\hrule}
 & &9.5.11.37.131&289&104& & &9.5.7.13.17.37&2603&3058 \cr
618&2399265&16.3.11.13.289&127&94&636&2575755&4.11.19.137.139&1389&1252 \cr
 & &64.17.47.127&2159&1504& & &32.3.11.313.463&5093&5008 \cr
\noalign{\hrule}
 & &25.7.11.29.43&73&102& & &81.49.11.59&1891&1000 \cr
619&2400475&4.3.11.17.43.73&245&228&637&2575881&16.125.31.61&77&78 \cr
 & &32.9.5.49.19.73&1197&1168& & &64.3.25.7.11.13.61&793&800 \cr
\noalign{\hrule}
 & &25.11.67.131&117&1558& & &11.17.37.373&1737&2366 \cr
620&2413675&4.9.13.19.41&29&10&638&2580787&4.9.7.169.193&205&374 \cr
 & &16.3.5.29.41&41&696& & &16.3.5.7.11.17.41&205&168 \cr
\noalign{\hrule}
 & &3.11.17.19.229&105&124& & &81.5.7.11.83&1007&1898 \cr
621&2440911&8.9.5.7.11.17.31&109&10&639&2588355&4.13.19.53.73&27&26 \cr
 & &32.25.31.109&775&1744& & &16.27.169.19.73&1387&1352 \cr
\noalign{\hrule}
 & &3.25.17.19.101&73&22& & &25.49.13.163&447&1672 \cr
622&2446725&4.5.11.73.101&289&216&640&2595775&16.3.11.19.149&65&84 \cr
 & &64.27.11.289&187&288& & &128.9.5.7.11.13&99&64 \cr
\noalign{\hrule}
 & &5.49.97.103&209&306& & &7.19.29.677&6417&6446 \cr
623&2447795&4.9.49.11.17.19&103&730&641&2611189&4.9.7.11.23.31.293&1885&7198 \cr
 & &16.3.5.73.103&73&24& & &16.3.5.13.29.59.61&2379&2360 \cr
\noalign{\hrule}
 & &9.5.11.19.263&1967&1978& & &25.7.11.23.59&47&12 \cr
624&2473515&4.3.7.19.23.43.281&13&830&642&2612225&8.3.5.11.23.47&63&52 \cr
 & &16.5.7.13.23.83&1079&1288& & &64.27.7.13.47&611&864 \cr
\noalign{\hrule}
 & &9.5.17.41.79&63&22& & &27.7.11.13.97&395&298 \cr
625&2477835&4.81.7.11.79&1&80&643&2621619&4.3.5.13.79.149&737&290 \cr
 & &128.5.7.11&7&704& & &16.25.11.29.67&725&536 \cr
\noalign{\hrule}
 & &9.25.11.17.59&551&452& & &3.7.29.59.73&103&74 \cr
626&2482425&8.25.19.29.113&69&44&644&2622963&4.7.37.73.103&1711&990 \cr
 & &64.3.11.19.23.29&667&608& & &16.9.5.11.29.59&33&40 \cr
\noalign{\hrule}
 & &49.17.31.97&11&108& & &3.7.17.53.139&69&70 \cr
627&2504831&8.27.7.11.31&5&26&645&2630019&4.9.5.49.17.23.53&461&3058 \cr
 & &32.9.5.11.13&65&1584& & &16.5.11.139.461&461&440 \cr
\noalign{\hrule}
 & &27.13.67.107&25&92& & &27.7.13.29.37&1819&820 \cr
628&2516319&8.3.25.23.107&19&88&646&2636361&8.5.17.41.107&33&74 \cr
 & &128.25.11.19&209&1600& & &32.3.5.11.17.37&85&176 \cr
\noalign{\hrule}
 & &3.25.11.43.71&59&414& & &27.19.53.97&517&490 \cr
629&2518725&4.27.5.23.59&71&44&647&2637333&4.5.49.11.47.97&81&598 \cr
 & &32.11.59.71&59&16& & &16.81.5.7.13.23&455&552 \cr
\noalign{\hrule}
 & &3.25.11.43.71&29&14& & &3.5.49.37.97&121&170 \cr
630&2518725&4.5.7.11.29.71&63&8&648&2637915&4.25.121.17.37&1827&1198 \cr
 & &64.9.49.29&147&928& & &16.9.7.29.599&599&696 \cr
\noalign{\hrule}
}%
}
$$
\eject
\vglue -23 pt
\noindent\hskip 1 in\hbox to 6.5 in{\ 649 -- 684 \hfill\fbd 2648685 -- 2866149\frb}
\vskip -9 pt
$$
\vbox{
\nointerlineskip
\halign{\strut
    \vrule \ \ \hfil \frb #\ 
   &\vrule \hfil \ \ \fbb #\frb\ 
   &\vrule \hfil \ \ \frb #\ \hfil
   &\vrule \hfil \ \ \frb #\ 
   &\vrule \hfil \ \ \frb #\ \ \vrule \hskip 2 pt
   &\vrule \ \ \hfil \frb #\ 
   &\vrule \hfil \ \ \fbb #\frb\ 
   &\vrule \hfil \ \ \frb #\ \hfil
   &\vrule \hfil \ \ \frb #\ 
   &\vrule \hfil \ \ \frb #\ \vrule \cr%
\noalign{\hrule}
 & &3.5.13.289.47&297&314& & &13.29.73.101&3875&3498 \cr
649&2648685&4.81.5.11.17.157&893&3562&667&2779621&4.3.125.11.31.53&303&38 \cr
 & &16.13.19.47.137&137&152& & &16.9.25.19.101&225&152 \cr
\noalign{\hrule}
 & &5.121.41.107&279&172& & &27.5.73.283&817&598 \cr
650&2654135&8.9.5.11.31.43&61&104&668&2788965&4.9.13.19.23.43&145&154 \cr
 & &128.3.13.31.61&1891&2496& & &16.5.7.11.19.29.43&3857&3784 \cr
\noalign{\hrule}
 & &3.121.71.103&233&130& & &3.5.121.23.67&403&202 \cr
651&2654619&4.5.13.71.233&121&1044&669&2796915&4.13.23.31.101&99&200 \cr
 & &32.9.121.29&87&16& & &64.9.25.11.31&93&160 \cr
\noalign{\hrule}
 & &25.7.17.19.47&659&234& & &5.11.17.41.73&447&488 \cr
652&2656675&4.9.7.13.659&319&340&670&2798455&16.3.61.73.149&67&6 \cr
 & &32.3.5.11.13.17.29&429&464& & &64.9.67.149&1341&2144 \cr
\noalign{\hrule}
 & &9.5.19.53.59&221&44& & &3.11.13.61.107&37&70 \cr
653&2673585&8.3.11.13.17.19&181&142&671&2800083&4.5.7.13.37.61&107&198 \cr
 & &32.11.71.181&1991&1136& & &16.9.11.37.107&37&24 \cr
\noalign{\hrule}
 & &9.5.19.53.59&1991&1136& & &3.5.7.121.13.17&15&106 \cr
654&2673585&32.11.71.181&481&300&672&2807805&4.9.25.17.53&121&104 \cr
 & &256.3.25.13.37&481&640& & &64.121.13.53&53&32 \cr
\noalign{\hrule}
 & &3.5.11.13.29.43&833&844& & &3.5.7.121.13.17&29&114 \cr
655&2674815&8.5.49.17.29.211&1971&494&673&2807805&4.9.7.11.19.29&37&26 \cr
 & &32.27.7.13.19.73&1197&1168& & &16.13.19.29.37&551&296 \cr
\noalign{\hrule}
 & &11.13.97.193&431&1692& & &9.5.11.13.19.23&47&8 \cr
656&2677103&8.9.47.431&145&286&674&2812095&16.3.19.23.47&25&44 \cr
 & &32.3.5.11.13.29&145&48& & &128.25.11.47&235&64 \cr
\noalign{\hrule}
 & &9.7.17.41.61&715&1786& & &7.11.23.37.43&15&22 \cr
657&2678571&4.5.11.13.19.47&17&30&675&2817661&4.3.5.121.23.43&555&434 \cr
 & &16.3.25.11.17.19&275&152& & &16.9.25.7.31.37&225&248 \cr
\noalign{\hrule}
 & &3.7.11.13.29.31&37&180& & &9.25.11.17.67&163&172 \cr
658&2699697&8.27.5.29.37&427&572&676&2819025&8.5.11.17.43.163&201&14 \cr
 & &64.7.11.13.61&61&32& & &32.3.7.67.163&163&112 \cr
\noalign{\hrule}
 & &9.11.19.1453&13853&13754& & &27.11.13.17.43&61&412 \cr
659&2733093&4.7.13.529.1979&725&1254&677&2822391&8.17.61.103&43&60 \cr
 & &16.3.25.7.11.13.19.29&725&728& & &64.3.5.43.61&61&160 \cr
\noalign{\hrule}
 & &9.5.49.17.73&451&206& & &81.25.23.61&143&82 \cr
660&2736405&4.11.17.41.103&73&114&678&2841075&4.9.11.13.23.41&61&38 \cr
 & &16.3.19.73.103&103&152& & &16.13.19.41.61&247&328 \cr
\noalign{\hrule}
 & &9.11.17.23.71&3319&3710& & &9.5.23.41.67&119&4 \cr
661&2748339&4.5.7.53.3319&1633&1686&679&2843145&8.3.7.17.67&23&44 \cr
 & &16.3.5.7.23.71.281&281&280& & &64.11.17.23&17&352 \cr
\noalign{\hrule}
 & &5.13.17.19.131&33&52& & &5.49.13.19.47&69&22 \cr
662&2750345&8.3.11.169.131&19&150&680&2844205&4.3.5.7.11.19.23&39&94 \cr
 & &32.9.25.11.19&45&176& & &16.9.13.23.47&23&72 \cr
\noalign{\hrule}
 & &5.7.11.13.19.29&3&206& & &9.5.13.31.157&7&38 \cr
663&2757755&4.3.5.13.103&19&84&681&2847195&4.7.13.19.157&33&124 \cr
 & &32.9.7.19&1&144& & &32.3.11.19.31&19&176 \cr
\noalign{\hrule}
 & &11.17.29.509&95&414& & &25.11.13.17.47&301&216 \cr
664&2760307&4.9.5.17.19.23&11&334&682&2856425&16.27.5.7.13.43&421&34 \cr
 & &16.3.11.167&501&8& & &64.3.17.421&421&96 \cr
\noalign{\hrule}
 & &3.5.11.97.173&65&32& & &3.5.7.11.37.67&91&646 \cr
665&2768865&64.25.13.173&99&74&683&2863245&4.49.13.17.19&99&148 \cr
 & &256.9.11.13.37&481&384& & &32.9.11.17.37&51&16 \cr
\noalign{\hrule}
 & &25.7.29.547&2277&1552& & &9.11.13.17.131&1757&470 \cr
666&2776025&32.9.11.23.97&29&40&684&2866149&4.5.7.47.251&149&102 \cr
 & &512.3.5.29.97&291&256& & &16.3.5.7.17.149&149&280 \cr
\noalign{\hrule}
}%
}
$$
\eject
\vglue -23 pt
\noindent\hskip 1 in\hbox to 6.5 in{\ 685 -- 720 \hfill\fbd 2873045 -- 3226635\frb}
\vskip -9 pt
$$
\vbox{
\nointerlineskip
\halign{\strut
    \vrule \ \ \hfil \frb #\ 
   &\vrule \hfil \ \ \fbb #\frb\ 
   &\vrule \hfil \ \ \frb #\ \hfil
   &\vrule \hfil \ \ \frb #\ 
   &\vrule \hfil \ \ \frb #\ \ \vrule \hskip 2 pt
   &\vrule \ \ \hfil \frb #\ 
   &\vrule \hfil \ \ \fbb #\frb\ 
   &\vrule \hfil \ \ \frb #\ \hfil
   &\vrule \hfil \ \ \frb #\ 
   &\vrule \hfil \ \ \frb #\ \vrule \cr%
\noalign{\hrule}
 & &5.7.23.43.83&57&358& & &25.7.11.19.83&1769&306 \cr
685&2873045&4.3.19.23.179&129&308&703&3035725&4.9.17.29.61&83&100 \cr
 & &32.9.7.11.43&99&16& & &32.3.25.29.83&29&48 \cr
\noalign{\hrule}
 & &5.11.13.37.109&17&126& & &9.7.73.661&1969&2630 \cr
686&2883595&4.9.5.7.17.37&33&152&704&3039939&4.5.11.179.263&221&42 \cr
 & &64.27.11.19&513&32& & &16.3.5.7.11.13.17&221&440 \cr
\noalign{\hrule}
 & &11.13.17.29.41&1725&2422& & &9.13.17.29.53&475&214 \cr
687&2890459&4.3.25.7.23.173&29&6&705&3057093&4.25.17.19.107&429&106 \cr
 & &16.9.5.29.173&173&360& & &16.3.5.11.13.53&11&40 \cr
\noalign{\hrule}
 & &3.5.11.13.29.47&49&38& & &5.11.17.29.113&611&1854 \cr
688&2923635&4.5.49.13.19.47&27&638&706&3063995&4.9.13.47.103&17&30 \cr
 & &16.27.7.11.29&9&56& & &16.27.5.17.103&103&216 \cr
\noalign{\hrule}
 & &27.11.59.167&349&182& & &9.25.31.443&1981&2006 \cr
689&2926341&4.3.7.11.13.349&59&290&707&3089925&4.7.17.31.59.283&33&250 \cr
 & &16.5.13.29.59&65&232& & &16.3.125.11.17.59&649&680 \cr
\noalign{\hrule}
 & &5.11.17.31.101&3151&2646& & &9.19.59.307&65&242 \cr
690&2927485&4.27.49.23.137&35&172&708&3097323&4.3.5.121.13.19&307&320 \cr
 & &32.3.5.343.43&1029&688& & &512.25.11.307&275&256 \cr
\noalign{\hrule}
 & &81.7.121.43&2885&2318& & &9.5.11.23.277&2041&1006 \cr
691&2950101&4.5.19.61.577&319&258&709&3153645&4.13.157.503&173&330 \cr
 & &16.3.5.11.19.29.43&145&152& & &16.3.5.11.13.173&173&104 \cr
\noalign{\hrule}
 & &7.11.83.463&235&228& & &7.11.17.19.127&55&72 \cr
692&2959033&8.3.5.11.19.47.83&903&10&710&3158617&16.9.5.7.121.19&127&6 \cr
 & &32.9.25.7.43&225&688& & &64.27.5.127&27&160 \cr
\noalign{\hrule}
 & &3.5.17.29.401&17&418& & &3.5.49.11.17.23&779&394 \cr
693&2965395&4.11.289.19&145&144&711&3161235&4.7.19.41.197&153&134 \cr
 & &128.9.5.11.19.29&209&192& & &16.9.17.67.197&591&536 \cr
\noalign{\hrule}
 & &3.25.7.13.19.23&37&2& & &3.25.11.23.167&113&388 \cr
694&2982525&4.5.19.23.37&473&378&712&3168825&8.23.97.113&45&68 \cr
 & &16.27.7.11.43&387&88& & &64.9.5.17.97&291&544 \cr
\noalign{\hrule}
 & &9.49.67.101&689&220& & &9.7.11.17.269&107&124 \cr
695&2984247&8.5.7.11.13.53&201&254&713&3169089&8.3.31.107.269&295&26 \cr
 & &32.3.11.67.127&127&176& & &32.5.13.31.59&2015&944 \cr
\noalign{\hrule}
 & &3.49.137.149&93&44& & &9.19.23.811&371&440 \cr
696&3000711&8.9.11.31.149&959&680&714&3189663&16.3.5.7.11.19.53&491&92 \cr
 & &128.5.7.17.137&85&64& & &128.5.23.491&491&320 \cr
\noalign{\hrule}
 & &3.5.7.11.23.113&453&338& & &5.11.97.599&3537&3052 \cr
697&3001845&4.9.11.169.151&1591&70&715&3195665&8.27.7.109.131&97&34 \cr
 & &16.5.7.37.43&43&296& & &32.3.17.97.109&327&272 \cr
\noalign{\hrule}
 & &3.5.7.11.19.137&261&124& & &9.5.7.11.13.71&31&86 \cr
698&3006465&8.27.19.29.31&241&272&716&3198195&4.7.31.43.71&915&418 \cr
 & &256.17.29.241&4097&3712& & &16.3.5.11.19.61&61&152 \cr
\noalign{\hrule}
 & &81.5.17.19.23&7&88& & &9.7.11.31.149&959&680 \cr
699&3008745&16.7.11.17.23&3&20&717&3200967&16.5.49.17.137&93&44 \cr
 & &128.3.5.7.11&11&448& & &128.3.5.11.17.31&85&64 \cr
\noalign{\hrule}
 & &3.5.11.13.23.61&63&2& & &3.11.17.29.197&1075&1092 \cr
700&3009435&4.27.7.11.23&113&140&718&3204993&8.9.25.7.13.29.43&23&238 \cr
 & &32.5.49.113&113&784& & &32.5.49.13.17.23&1127&1040 \cr
\noalign{\hrule}
 & &27.13.23.373&4223&3850& & &27.7.43.397&253&650 \cr
701&3011229&4.25.7.11.41.103&483&32&719&3226419&4.9.25.11.13.23&397&98 \cr
 & &256.3.5.49.23&245&128& & &16.5.49.397&35&8 \cr
\noalign{\hrule}
 & &3.11.19.61.79&465&404& & &81.5.31.257&1079&1234 \cr
702&3021513&8.9.5.19.31.101&2567&1978&720&3226635&4.9.13.83.617&65&682 \cr
 & &32.17.23.43.151&6493&6256& & &16.5.11.169.31&169&88 \cr
\noalign{\hrule}
}%
}
$$
\eject
\vglue -23 pt
\noindent\hskip 1 in\hbox to 6.5 in{\ 721 -- 756 \hfill\fbd 3253481 -- 3542877\frb}
\vskip -9 pt
$$
\vbox{
\nointerlineskip
\halign{\strut
    \vrule \ \ \hfil \frb #\ 
   &\vrule \hfil \ \ \fbb #\frb\ 
   &\vrule \hfil \ \ \frb #\ \hfil
   &\vrule \hfil \ \ \frb #\ 
   &\vrule \hfil \ \ \frb #\ \ \vrule \hskip 2 pt
   &\vrule \ \ \hfil \frb #\ 
   &\vrule \hfil \ \ \fbb #\frb\ 
   &\vrule \hfil \ \ \frb #\ \hfil
   &\vrule \hfil \ \ \frb #\ 
   &\vrule \hfil \ \ \frb #\ \vrule \cr%
\noalign{\hrule}
 & &7.11.29.31.47&157&360& & &27.121.17.61&3055&3116 \cr
721&3253481&16.9.5.31.157&125&32&739&3387879&8.9.5.13.19.41.47&61&308 \cr
 & &1024.3.625&1875&512& & &64.5.7.11.47.61&235&224 \cr
\noalign{\hrule}
 & &25.11.17.19.37&567&358& & &9.5.17.43.103&47&4 \cr
722&3286525&4.81.7.17.179&703&550&740&3388185&8.3.5.47.103&187&328 \cr
 & &16.9.25.11.19.37&9&8& & &128.11.17.41&451&64 \cr
\noalign{\hrule}
 & &3.13.29.41.71&1829&230& & &7.19.71.359&495&854 \cr
723&3292341&4.5.23.31.59&209&504&741&3390037&4.9.5.49.11.61&71&76 \cr
 & &64.9.7.11.19&133&1056& & &32.3.11.19.61.71&183&176 \cr
\noalign{\hrule}
 & &3.5.11.13.29.53&27&28& & &9.5.11.13.17.31&1961&674 \cr
724&3296865&8.81.7.13.29.53&1295&242&742&3391245&4.37.53.337&187&150 \cr
 & &32.5.49.121.37&539&592& & &16.3.25.11.17.53&53&40 \cr
\noalign{\hrule}
 & &25.7.11.17.101&193&312& & &3.5.23.59.167&429&406 \cr
725&3305225&16.3.5.11.13.193&69&124&743&3399285&4.9.7.11.13.29.59&2047&1670 \cr
 & &128.9.13.23.31&2691&1984& & &16.5.11.23.89.167&89&88 \cr
\noalign{\hrule}
 & &9.25.61.241&143&98& & &3.7.13.29.431&209&222 \cr
726&3307725&4.5.49.11.13.61&197&258&744&3412227&4.9.7.11.19.29.37&33695&33904 \cr
 & &16.3.7.11.43.197&2167&2408& & &128.5.13.23.163.293&18745&18752 \cr
\noalign{\hrule}
 & &9.625.19.31&73&98& & &7.31.71.223&233&264 \cr
727&3313125&4.25.49.31.73&643&132&745&3435761&16.3.11.223.233&5&228 \cr
 & &32.3.7.11.643&643&1232& & &128.9.5.11.19&1881&320 \cr
\noalign{\hrule}
 & &7.169.2803&1993&810& & &5.7.29.43.79&93&122 \cr
728&3315949&4.81.5.1993&1001&992&746&3447955&4.3.7.31.61.79&253&174 \cr
 & &256.9.5.7.11.13.31&1395&1408& & &16.9.11.23.29.31&713&792 \cr
\noalign{\hrule}
 & &11.31.71.137&135&206& & &27.5.7.13.281&17&22 \cr
729&3316907&4.27.5.103.137&121&806&747&3452085&4.9.7.11.17.281&1919&610 \cr
 & &16.3.121.13.31&33&104& & &16.5.19.61.101&1159&808 \cr
\noalign{\hrule}
 & &3.25.11.37.109&83&28& & &7.17.71.409&2035&828 \cr
730&3327225&8.5.7.83.109&37&72&748&3455641&8.9.5.11.23.37&169&238 \cr
 & &128.9.37.83&83&192& & &32.3.5.7.169.17&169&240 \cr
\noalign{\hrule}
 & &3.5.13.17.19.53&183&506& & &9.25.121.127&53&328 \cr
731&3338205&4.9.5.11.23.61&53&8&749&3457575&16.3.11.41.53&703&650 \cr
 & &64.11.23.53&253&32& & &64.25.13.19.37&481&608 \cr
\noalign{\hrule}
 & &9.13.23.29.43&14399&14282& & &29.1681.71&189&1870 \cr
732&3355677&4.7.121.17.37.193&3&190&750&3461179&4.27.5.7.11.17&41&58 \cr
 & &16.3.5.7.11.19.37&1295&1672& & &16.3.5.7.29.41&15&56 \cr
\noalign{\hrule}
 & &3.5.41.53.103&51&154& & &3.25.13.53.67&303&568 \cr
733&3357285&4.9.7.11.17.53&103&50&751&3462225&16.9.5.71.101&67&572 \cr
 & &16.25.7.11.103&35&88& & &128.11.13.67&11&64 \cr
\noalign{\hrule}
 & &5.49.11.29.43&689&732& & &5.29.71.337&187&1872 \cr
734&3360665&8.3.5.11.13.53.61&1387&2058&752&3469415&32.9.11.13.17&337&326 \cr
 & &32.9.343.19.73&1197&1168& & &128.3.163.337&163&192 \cr
\noalign{\hrule}
 & &3.5.7.11.41.71&67&138& & &25.19.73.101&289&216 \cr
735&3362205&4.9.7.11.23.67&337&400&753&3502175&16.27.5.289.19&73&22 \cr
 & &128.25.23.337&1685&1472& & &64.9.11.17.73&187&288 \cr
\noalign{\hrule}
 & &3.25.7.13.17.29&1957&682& & &27.13.79.127&215&136 \cr
736&3364725&4.11.19.31.103&243&346&754&3521583&16.5.17.43.127&171&44 \cr
 & &16.243.11.173&891&1384& & &128.9.11.17.19&323&704 \cr
\noalign{\hrule}
 & &27.25.7.23.31&407&430& & &9.125.31.101&539&236 \cr
737&3368925&4.125.7.11.37.43&201&674&755&3522375&8.3.5.49.11.59&17&38 \cr
 & &16.3.37.67.337&2479&2696& & &32.7.17.19.59&1003&2128 \cr
\noalign{\hrule}
 & &5.11.19.53.61&17&78& & &9.13.107.283&1969&578 \cr
738&3378485&4.3.11.13.17.53&67&120&756&3542877&4.11.289.179&1679&1500 \cr
 & &64.9.5.13.67&871&288& & &32.3.125.23.73&1679&2000 \cr
\noalign{\hrule}
}%
}
$$
\eject
\vglue -23 pt
\noindent\hskip 1 in\hbox to 6.5 in{\ 757 -- 792 \hfill\fbd 3559215 -- 4041205\frb}
\vskip -9 pt
$$
\vbox{
\nointerlineskip
\halign{\strut
    \vrule \ \ \hfil \frb #\ 
   &\vrule \hfil \ \ \fbb #\frb\ 
   &\vrule \hfil \ \ \frb #\ \hfil
   &\vrule \hfil \ \ \frb #\ 
   &\vrule \hfil \ \ \frb #\ \ \vrule \hskip 2 pt
   &\vrule \ \ \hfil \frb #\ 
   &\vrule \hfil \ \ \fbb #\frb\ 
   &\vrule \hfil \ \ \frb #\ \hfil
   &\vrule \hfil \ \ \frb #\ 
   &\vrule \hfil \ \ \frb #\ \vrule \cr%
\noalign{\hrule}
 & &3.5.121.37.53&2929&3484& & &3.11.169.691&365&326 \cr
757&3559215&8.13.29.67.101&315&1628&775&3853707&4.5.11.13.73.163&153&10 \cr
 & &64.9.5.7.11.37&21&32& & &16.9.25.17.73&1275&584 \cr
\noalign{\hrule}
 & &3.7.11.13.29.41&2627&3236& & &5.11.61.1151&581&570 \cr
758&3570567&8.37.71.809&369&440&776&3861605&4.3.25.7.19.61.83&1551&26 \cr
 & &128.9.5.11.37.41&185&192& & &16.9.7.11.13.47&423&728 \cr
\noalign{\hrule}
 & &5.49.11.13.103&3&718& & &11.17.23.29.31&373&120 \cr
759&3608605&4.3.7.359&183&176&777&3866599&16.3.5.31.373&171&202 \cr
 & &128.9.11.61&61&576& & &64.27.5.19.101&2565&3232 \cr
\noalign{\hrule}
 & &27.13.41.251&2185&1078& & &49.11.31.233&287&54 \cr
760&3612141&4.5.49.11.19.23&47&162&778&3893197&4.27.343.41&233&110 \cr
 & &16.81.49.47&141&392& & &16.9.5.11.233&9&40 \cr
\noalign{\hrule}
 & &3.11.17.43.151&349&382& & &9.13.37.907&1711&2618 \cr
761&3642573&4.151.191.349&171&20&779&3926403&4.7.11.17.29.59&65&54 \cr
 & &32.9.5.19.349&1047&1520& & &16.27.5.13.29.59&435&472 \cr
\noalign{\hrule}
 & &9.5.11.17.19.23&191&246& & &9.25.97.181&11&86 \cr
762&3677355&4.27.17.41.191&209&250&780&3950325&4.3.11.43.181&35&508 \cr
 & &16.125.11.19.191&191&200& & &32.5.7.127&889&16 \cr
\noalign{\hrule}
 & &3.5.13.17.19.59&153&94& & &9.11.13.17.181&989&1002 \cr
763&3716115&4.27.5.289.47&1357&88&781&3960099&4.27.17.23.43.167&265&724 \cr
 & &64.11.23.59&253&32& & &32.5.53.167.181&835&848 \cr
\noalign{\hrule}
 & &9.125.7.11.43&299&174& & &9.5.19.41.113&43&52 \cr
764&3724875&4.27.7.13.23.29&95&256&782&3961215&8.13.41.43.113&825&938 \cr
 & &2048.5.19.29&551&1024& & &32.3.25.7.11.13.67&2345&2288 \cr
\noalign{\hrule}
 & &9.13.19.23.73&235&64& & &5.17.101.463&21&484 \cr
765&3732417&128.5.47.73&81&154&783&3974855&8.3.7.121.17&69&52 \cr
 & &512.81.7.11&693&256& & &64.9.7.13.23&299&2016 \cr
\noalign{\hrule}
 & &5.49.17.29.31&473&1992& & &25.19.83.101&1001&918 \cr
766&3744335&16.3.11.43.83&435&478&784&3981925&4.27.25.7.11.13.17&83&8 \cr
 & &64.9.5.29.239&239&288& & &64.9.11.17.83&187&288 \cr
\noalign{\hrule}
 & &3.343.11.331&169&162& & &27.41.59.61&73793&73766 \cr
767&3746589&4.243.49.11.169&1655&1018&785&3984093&4.7.11.109.479.677&41&52170 \cr
 & &16.5.13.331.509&509&520& & &16.3.5.37.41.47&235&296 \cr
\noalign{\hrule}
 & &23.29.43.131&77&54& & &5.7.13.67.131&321&2024 \cr
768&3757211&4.27.7.11.29.43&923&1310&786&3993535&16.3.11.23.107&65&42 \cr
 & &16.3.5.13.71.131&355&312& & &64.9.5.7.11.13&99&32 \cr
\noalign{\hrule}
 & &27.5.7.41.97&139&148& & &27.11.13.17.61&31&20 \cr
769&3758265&8.3.5.37.97.139&451&34&787&4003857&8.9.5.13.31.61&73&476 \cr
 & &32.11.17.37.41&407&272& & &64.5.7.17.73&511&160 \cr
\noalign{\hrule}
 & &3.25.11.43.107&31&504& & &9.5.7.11.13.89&367&634 \cr
770&3795825&16.27.5.7.31&71&64&788&4009005&4.3.5.317.367&1343&242 \cr
 & &2048.31.71&2201&1024& & &16.121.17.79&187&632 \cr
\noalign{\hrule}
 & &9.5.13.73.89&241&124& & &3.13.29.53.67&359&330 \cr
771&3800745&8.31.89.241&165&76&789&4016181&4.9.5.11.67.359&1247&1984 \cr
 & &64.3.5.11.19.31&341&608& & &512.5.29.31.43&1333&1280 \cr
\noalign{\hrule}
 & &3.49.19.29.47&629&2050& & &27.11.19.23.31&295&326 \cr
772&3806859&4.25.17.37.41&63&22&790&4023459&4.5.11.19.59.163&23&186 \cr
 & &16.9.5.7.11.37&111&440& & &16.3.5.23.31.59&59&40 \cr
\noalign{\hrule}
 & &3.5.11.19.23.53&21&74& & &9.5.17.19.277&1927&2782 \cr
773&3821565&4.9.7.11.23.37&65&142&791&4026195&4.13.41.47.107&3&44 \cr
 & &16.5.13.37.71&923&296& & &32.3.11.13.107&1177&208 \cr
\noalign{\hrule}
 & &5.17.107.421&143&1962& & &5.7.19.59.103&3371&2706 \cr
774&3828995&4.9.11.13.109&157&170&792&4041205&4.3.11.41.3371&1911&1460 \cr
 & &16.3.5.11.17.157&157&264& & &32.9.5.49.13.73&949&1008 \cr
\noalign{\hrule}
}%
}
$$
\eject
\vglue -23 pt
\noindent\hskip 1 in\hbox to 6.5 in{\ 793 -- 828 \hfill\fbd 4042005 -- 4367853\frb}
\vskip -9 pt
$$
\vbox{
\nointerlineskip
\halign{\strut
    \vrule \ \ \hfil \frb #\ 
   &\vrule \hfil \ \ \fbb #\frb\ 
   &\vrule \hfil \ \ \frb #\ \hfil
   &\vrule \hfil \ \ \frb #\ 
   &\vrule \hfil \ \ \frb #\ \ \vrule \hskip 2 pt
   &\vrule \ \ \hfil \frb #\ 
   &\vrule \hfil \ \ \fbb #\frb\ 
   &\vrule \hfil \ \ \frb #\ \hfil
   &\vrule \hfil \ \ \frb #\ 
   &\vrule \hfil \ \ \frb #\ \vrule \cr%
\noalign{\hrule}
 & &3.5.121.17.131&763&678& & &3.5.121.23.101&3&118 \cr
793&4042005&4.9.7.11.109.113&131&1112&811&4216245&4.9.59.101&425&484 \cr
 & &64.7.131.139&139&224& & &32.25.121.17&85&16 \cr
\noalign{\hrule}
 & &529.79.97&957&860& & &5.11.13.17.347&1941&1876 \cr
794&4053727&8.3.5.11.23.29.43&97&570&812&4217785&8.3.7.17.67.647&893&246 \cr
 & &32.9.25.19.97&225&304& & &32.9.7.19.41.47&6251&5904 \cr
\noalign{\hrule}
 & &5.49.13.19.67&99&148& & &5.11.23.47.71&51&4 \cr
795&4054505&8.9.5.11.37.67&91&646&813&4221305&8.3.17.23.71&825&808 \cr
 & &32.3.7.13.17.19&51&16& & &128.9.25.11.101&505&576 \cr
\noalign{\hrule}
 & &13.271.1153&7359&7630& & &25.7.361.67&79&54 \cr
796&4062019&4.3.5.7.11.109.223&1157&42&814&4232725&4.27.19.67.79&91&110 \cr
 & &16.9.49.13.89&441&712& & &16.9.5.7.11.13.79&869&936 \cr
\noalign{\hrule}
 & &5.11.13.73.79&269&126& & &27.49.11.293&2273&950 \cr
797&4123405&4.9.7.73.269&25&244&815&4264029&4.25.19.2273&1139&1134 \cr
 & &32.3.25.7.61&915&112& & &16.81.5.7.17.19.67&1615&1608 \cr
\noalign{\hrule}
 & &49.13.67.97&865&396& & &3.5.29.71.139&37&108 \cr
798&4139863&8.9.5.7.11.173&97&76&816&4293015&8.81.37.139&29&110 \cr
 & &64.3.5.11.19.97&285&352& & &32.5.11.29.37&37&176 \cr
\noalign{\hrule}
 & &81.17.31.97&685&964& & &9.5.13.41.179&371&166 \cr
799&4140639&8.9.5.137.241&737&496&817&4293315&4.3.7.13.53.83&469&220 \cr
 & &256.5.11.31.67&737&640& & &32.5.49.11.67&539&1072 \cr
\noalign{\hrule}
 & &27.25.11.13.43&37&92& & &11.59.61.109&35&24 \cr
800&4150575&8.9.5.13.23.37&71&136&818&4315201&16.3.5.7.61.109&159&268 \cr
 & &128.17.37.71&1207&2368& & &128.9.5.53.67&3551&2880 \cr
\noalign{\hrule}
 & &27.25.11.13.43&1207&2368& & &9.5.139.691&341&350 \cr
801&4150575&128.17.37.71&279&350&819&4322205&4.125.7.11.31.139&7&132 \cr
 & &512.9.25.7.31&217&256& & &32.3.49.121.31&1519&1936 \cr
\noalign{\hrule}
 & &27.25.11.13.43&721&764& & &27.11.47.311&7135&6824 \cr
802&4150575&8.5.7.13.103.191&2673&4012&820&4341249&16.5.853.1427&1419&2846 \cr
 & &64.243.11.17.59&531&544& & &64.3.11.43.1423&1423&1376 \cr
\noalign{\hrule}
 & &3.49.11.17.151&537&520& & &3.5.11.361.73&221&582 \cr
803&4150839&16.9.5.7.11.13.179&817&1510&821&4348245&4.9.5.13.17.97&47&38 \cr
 & &64.25.19.43.151&817&800& & &16.13.19.47.97&611&776 \cr
\noalign{\hrule}
 & &3.5.29.41.233&2623&3322& & &27.5.13.37.67&151&34 \cr
804&4155555&4.11.43.61.151&2911&3582&822&4350645&4.3.17.67.151&25&176 \cr
 & &16.9.41.71.199&597&568& & &128.25.11.17&935&64 \cr
\noalign{\hrule}
 & &5.13.17.53.71&243&22& & &7.37.53.317&129&2090 \cr
805&4158115&4.243.11.71&5&76&823&4351459&4.3.5.11.19.43&37&18 \cr
 & &32.3.5.11.19&11&912& & &16.27.37.43&27&344 \cr
\noalign{\hrule}
 & &17.29.79.107&1199&1092& & &3.5.7.19.37.59&37&58 \cr
806&4167329&8.3.7.11.13.17.109&5&114&824&4355085&4.29.1369.59&171&1540 \cr
 & &32.9.5.11.13.19&2223&880& & &32.9.5.7.11.19&11&48 \cr
\noalign{\hrule}
 & &25.49.41.83&1089&2314& & &27.5.169.191&2759&1804 \cr
807&4168675&4.9.121.13.89&105&16&825&4357665&8.11.31.41.89&191&150 \cr
 & &128.27.5.7.13&351&64& & &32.3.25.89.191&89&80 \cr
\noalign{\hrule}
 & &9.7.113.587&2585&2698& & &5.17.19.37.73&55&18 \cr
808&4178853&4.5.7.11.19.47.71&3&74&826&4362115&4.9.25.11.17.19&43&518 \cr
 & &16.3.5.19.37.47&893&1480& & &16.3.7.37.43&21&344 \cr
\noalign{\hrule}
 & &27.5.7.11.13.31&97&20& & &5.17.19.37.73&27&46 \cr
809&4189185&8.3.25.31.97&533&242&827&4362115&4.27.5.17.23.37&209&124 \cr
 & &32.121.13.41&41&176& & &32.3.11.19.23.31&713&528 \cr
\noalign{\hrule}
 & &121.17.23.89&105&16& & &9.7.19.41.89&41&22 \cr
810&4210679&32.3.5.7.17.23&89&72&828&4367853&4.11.1681.89&351&1330 \cr
 & &512.27.5.89&135&256& & &16.27.5.7.13.19&65&24 \cr
\noalign{\hrule}
}%
}
$$
\eject
\vglue -23 pt
\noindent\hskip 1 in\hbox to 6.5 in{\ 829 -- 864 \hfill\fbd 4399241 -- 4734639\frb}
\vskip -9 pt
$$
\vbox{
\nointerlineskip
\halign{\strut
    \vrule \ \ \hfil \frb #\ 
   &\vrule \hfil \ \ \fbb #\frb\ 
   &\vrule \hfil \ \ \frb #\ \hfil
   &\vrule \hfil \ \ \frb #\ 
   &\vrule \hfil \ \ \frb #\ \ \vrule \hskip 2 pt
   &\vrule \ \ \hfil \frb #\ 
   &\vrule \hfil \ \ \fbb #\frb\ 
   &\vrule \hfil \ \ \frb #\ \hfil
   &\vrule \hfil \ \ \frb #\ 
   &\vrule \hfil \ \ \frb #\ \vrule \cr%
\noalign{\hrule}
 & &7.11.19.31.97&235&444& & &9.11.113.409&7&106 \cr
829&4399241&8.3.5.31.37.47&3&34&847&4575483&4.7.53.409&19&390 \cr
 & &32.9.5.17.47&2115&272& & &16.3.5.13.19&1235&8 \cr
\noalign{\hrule}
 & &5.169.17.307&43&264& & &9.25.7.41.71&133&92 \cr
830&4410055&16.3.5.11.13.43&269&204&848&4584825&8.49.19.23.71&11&60 \cr
 & &128.9.17.269&269&576& & &64.3.5.11.19.23&253&608 \cr
\noalign{\hrule}
 & &27.5.7.31.151&41&176& & &9.37.61.227&73&110 \cr
831&4423545&32.11.41.151&301&150&849&4611051&4.3.5.11.73.227&61&742 \cr
 & &128.3.25.7.43&215&64& & &16.5.7.53.61&53&280 \cr
\noalign{\hrule}
 & &9.29.113.151&187&74& & &81.7.79.103&407&304 \cr
832&4453443&4.11.17.37.151&1265&1302&850&4613679&32.9.7.11.19.37&395&802 \cr
 & &16.3.5.7.121.23.31&4991&4840& & &128.5.79.401&401&320 \cr
\noalign{\hrule}
 & &3.5.11.13.31.67&505&518& & &3.7.31.47.151&23&70 \cr
833&4455165&4.25.7.37.67.101&23&2502&851&4620147&4.5.49.23.151&47&198 \cr
 & &16.9.7.23.139&417&1288& & &16.9.11.23.47&33&184 \cr
\noalign{\hrule}
 & &81.17.41.79&1&80& & &11.19.67.331&105&104 \cr
834&4460103&32.5.17.41&63&22&852&4634993&16.3.5.7.13.67.331&333&2 \cr
 & &128.9.7.11&7&704& & &64.27.7.13.37&2457&1184 \cr
\noalign{\hrule}
 & &11.13.131.239&185&54& & &3.5.7.13.19.179&31&26 \cr
835&4477187&4.27.5.11.13.37&191&524&853&4642365&4.7.169.31.179&2183&3366 \cr
 & &32.3.131.191&191&48& & &16.9.11.17.37.59&3009&3256 \cr
\noalign{\hrule}
 & &27.5.29.31.37&427&572& & &11.13.19.29.59&305&246 \cr
836&4490505&8.7.11.13.31.61&37&180&854&4648787&4.3.5.11.13.41.61&31&174 \cr
 & &64.9.5.37.61&61&32& & &16.9.29.31.61&279&488 \cr
\noalign{\hrule}
 & &9.13.131.293&505&374& & &27.5.11.31.101&193&148 \cr
837&4490811&4.3.5.11.13.17.101&41&262&855&4649535&8.3.37.101.193&55&248 \cr
 & &16.5.11.41.131&205&88& & &128.5.11.31.37&37&64 \cr
\noalign{\hrule}
 & &25.13.17.19.43&117&98& & &5.7.11.43.281&119&162 \cr
838&4513925&4.9.5.49.169.17&551&44&856&4651955&4.81.5.49.11.17&1577&1118 \cr
 & &32.3.7.11.19.29&203&528& & &16.3.13.19.43.83&741&664 \cr
\noalign{\hrule}
 & &3.5.49.11.13.43&47&82& & &27.13.37.359&493&506 \cr
839&4519515&4.7.11.13.41.47&215&72&857&4662333&4.11.17.23.29.359&375&16 \cr
 & &64.9.5.43.47&47&96& & &128.3.125.11.29&1375&1856 \cr
\noalign{\hrule}
 & &9.25.7.169.17&649&2224& & &3.11.19.43.173&1649&1638 \cr
840&4524975&32.11.59.139&735&794&858&4664253&4.27.7.13.17.43.97&1211&50 \cr
 & &128.3.5.49.397&397&448& & &16.25.49.17.173&425&392 \cr
\noalign{\hrule}
 & &5.11.13.17.373&93&280& & &7.11.13.59.79&4495&3942 \cr
841&4533815&16.3.25.7.13.31&289&114&859&4665661&4.27.5.29.31.73&21&52 \cr
 & &64.9.289.19&171&544& & &32.81.5.7.13.29&405&464 \cr
\noalign{\hrule}
 & &25.121.19.79&7&18& & &3.7.11.17.29.41&3&14 \cr
842&4540525&4.9.7.11.19.79&65&144&860&4669203&4.9.49.29.41&815&374 \cr
 & &128.81.5.7.13&1053&448& & &16.5.11.17.163&163&40 \cr
\noalign{\hrule}
 & &3.7.11.17.19.61&83&100& & &27.7.149.167&169&20 \cr
843&4551393&8.25.7.11.19.83&1769&306&861&4702887&8.5.169.167&77&90 \cr
 & &32.9.17.29.61&29&48& & &32.9.25.7.11.13&325&176 \cr
\noalign{\hrule}
 & &5.7.13.17.19.31&37&54& & &25.7.11.31.79&4223&4302 \cr
844&4555915&4.27.5.19.31.37&17&572&862&4714325&4.9.7.41.103.239&991&682 \cr
 & &32.9.11.13.17&11&144& & &16.3.11.31.41.991&991&984 \cr
\noalign{\hrule}
 & &3.5.13.67.349&7&342& & &27.25.47.149&3839&3164 \cr
845&4559685&4.27.7.13.19&137&110&863&4727025&8.7.11.113.349&447&796 \cr
 & &16.5.7.11.137&959&88& & &64.3.7.149.199&199&224 \cr
\noalign{\hrule}
 & &3.5.7.11.17.233&817&118& & &27.7.13.41.47&55&8 \cr
846&4574955&4.7.19.43.59&121&180&864&4734639&16.3.5.11.13.41&29&94 \cr
 & &32.9.5.121.19&57&176& & &64.11.29.47&319&32 \cr
\noalign{\hrule}
}%
}
$$
\eject
\vglue -23 pt
\noindent\hskip 1 in\hbox to 6.5 in{\ 865 -- 900 \hfill\fbd 4751285 -- 5352399\frb}
\vskip -9 pt
$$
\vbox{
\nointerlineskip
\halign{\strut
    \vrule \ \ \hfil \frb #\ 
   &\vrule \hfil \ \ \fbb #\frb\ 
   &\vrule \hfil \ \ \frb #\ \hfil
   &\vrule \hfil \ \ \frb #\ 
   &\vrule \hfil \ \ \frb #\ \ \vrule \hskip 2 pt
   &\vrule \ \ \hfil \frb #\ 
   &\vrule \hfil \ \ \fbb #\frb\ 
   &\vrule \hfil \ \ \frb #\ \hfil
   &\vrule \hfil \ \ \frb #\ 
   &\vrule \hfil \ \ \frb #\ \vrule \cr%
\noalign{\hrule}
 & &5.49.11.41.43&377&162& & &3.25.7.73.131&291&364 \cr
865&4751285&4.81.13.29.41&125&658&883&5020575&8.9.5.49.13.97&463&22 \cr
 & &16.3.125.7.47&75&376& & &32.11.13.463&463&2288 \cr
\noalign{\hrule}
 & &3.25.121.17.31&2303&722& & &3.7.11.19.31.37&41&52 \cr
866&4782525&4.49.361.47&1287&1240&884&5034183&8.7.13.19.37.41&495&208 \cr
 & &64.9.5.7.11.13.31&91&96& & &256.9.5.11.169&507&640 \cr
\noalign{\hrule}
 & &3.7.169.19.71&2449&1100& & &11.29.97.163&911&882 \cr
867&4787601&8.25.11.31.79&13&18&885&5043709&4.9.49.97.911&795&116 \cr
 & &32.9.5.11.13.79&395&528& & &32.27.5.7.29.53&945&848 \cr
\noalign{\hrule}
 & &5.11.13.37.181&99&86& & &3.5.7.11.17.257&35&222 \cr
868&4788355&4.9.121.43.181&91&272&886&5046195&4.9.25.49.37&683&242 \cr
 & &128.3.7.13.17.43&903&1088& & &16.121.683&683&88 \cr
\noalign{\hrule}
 & &13.23.61.263&3421&2628& & &27.19.59.167&2815&1694 \cr
869&4796857&8.9.11.73.311&265&46&887&5054589&4.5.7.121.563&221&342 \cr
 & &32.3.5.11.23.53&265&528& & &16.9.5.7.13.17.19&221&280 \cr
\noalign{\hrule}
 & &5.13.17.43.101&27&532& & &3.5.7.11.23.191&113&78 \cr
870&4799015&8.27.7.17.19&77&94&888&5073915&4.9.11.13.23.113&635&382 \cr
 & &32.3.49.11.47&539&2256& & &16.5.13.127.191&127&104 \cr
\noalign{\hrule}
 & &9.7.11.169.41&535&986& & &3.5.7.11.53.83&27&26 \cr
871&4801797&4.5.7.17.29.107&13&132&889&5080845&4.81.5.7.11.13.83&1007&1898 \cr
 & &32.3.11.13.107&107&16& & &16.169.19.53.73&1387&1352 \cr
\noalign{\hrule}
 & &25.41.53.89&2871&1846& & &9.19.167.179&781&2392 \cr
872&4834925&4.9.11.13.29.71&295&82&890&5111703&16.11.13.23.71&185&114 \cr
 & &16.3.5.11.41.59&177&88& & &64.3.5.11.19.37&185&352 \cr
\noalign{\hrule}
 & &9.125.11.17.23&97&28& & &81.17.47.79&35&44 \cr
873&4838625&8.3.7.11.17.97&509&800&891&5112801&8.9.5.7.11.17.47&305&118 \cr
 & &512.25.509&509&256& & &32.25.7.59.61&3599&2800 \cr
\noalign{\hrule}
 & &7.11.13.37.131&135&124& & &9.7.127.643&893&250 \cr
874&4851847&8.27.5.13.31.131&7&124&892&5144643&4.125.7.19.47&87&88 \cr
 & &64.3.5.7.961&961&480& & &64.3.5.11.19.29.47&6815&6688 \cr
\noalign{\hrule}
 & &9.11.17.41.71&41&58& & &9.7.11.13.577&4097&3404 \cr
875&4899213&4.29.1681.71&189&1870&893&5198193&8.17.23.37.241&305&546 \cr
 & &16.27.5.7.11.17&15&56& & &32.3.5.7.13.17.61&305&272 \cr
\noalign{\hrule}
 & &5.11.17.59.89&87&2& & &25.13.19.23.37&473&378 \cr
876&4909685&4.3.11.29.59&861&850&894&5254925&4.27.5.7.11.13.43&37&2 \cr
 & &16.9.25.7.17.41&369&280& & &16.9.11.37.43&387&88 \cr
\noalign{\hrule}
 & &11.43.101.103&2177&2166& & &3.121.17.857&247&610 \cr
877&4920619&4.3.7.361.103.311&2067&110&895&5288547&4.5.13.17.19.61&235&558 \cr
 & &16.9.5.11.13.19.53&2223&2120& & &16.9.25.31.47&1457&600 \cr
\noalign{\hrule}
 & &81.31.37.53&65&28& & &9.7.11.13.19.31&25&118 \cr
878&4924071&8.27.5.7.13.53&121&68&896&5306301&4.3.25.7.19.59&155&22 \cr
 & &64.5.121.13.17&1573&2720& & &16.125.11.31&125&8 \cr
\noalign{\hrule}
 & &5.29.67.509&301&2244& & &9.125.53.89&3311&3364 \cr
879&4944935&8.3.7.11.17.43&145&156&897&5306625&8.3.5.7.11.841.43&79&2444 \cr
 & &64.9.5.13.17.29&221&288& & &64.7.13.47.79&3713&2912 \cr
\noalign{\hrule}
 & &9.5.11.17.19.31&101&70& & &3.25.121.19.31&559&464 \cr
880&4956435&4.25.7.11.17.101&1917&608&898&5345175&32.5.11.13.29.43&981&266 \cr
 & &256.27.19.71&213&128& & &128.9.7.19.109&327&448 \cr
\noalign{\hrule}
 & &49.13.73.107&3375&3268& & &5.11.31.43.73&509&294 \cr
881&4975607&8.27.125.7.19.43&671&146&899&5351995&4.3.49.31.509&363&146 \cr
 & &32.9.5.11.61.73&671&720& & &16.9.7.121.73&99&56 \cr
\noalign{\hrule}
 & &5.7.121.1181&629&552& & &81.169.17.23&109&112 \cr
882&5001535&16.3.5.11.17.23.37&893&42&900&5352399&32.27.7.13.23.109&55&244 \cr
 & &64.9.7.19.47&893&288& & &256.5.11.61.109&6649&7040 \cr
\noalign{\hrule}
}%
}
$$
\eject
\vglue -23 pt
\noindent\hskip 1 in\hbox to 6.5 in{\ 901 -- 936 \hfill\fbd 5359805 -- 5897045\frb}
\vskip -9 pt
$$
\vbox{
\nointerlineskip
\halign{\strut
    \vrule \ \ \hfil \frb #\ 
   &\vrule \hfil \ \ \fbb #\frb\ 
   &\vrule \hfil \ \ \frb #\ \hfil
   &\vrule \hfil \ \ \frb #\ 
   &\vrule \hfil \ \ \frb #\ \ \vrule \hskip 2 pt
   &\vrule \ \ \hfil \frb #\ 
   &\vrule \hfil \ \ \fbb #\frb\ 
   &\vrule \hfil \ \ \frb #\ \hfil
   &\vrule \hfil \ \ \frb #\ 
   &\vrule \hfil \ \ \frb #\ \vrule \cr%
\noalign{\hrule}
 & &5.11.19.23.223&117&106& & &3.11.13.83.157&905&822 \cr
901&5359805&4.9.5.13.19.23.53&561&446&919&5590299&4.9.5.13.137.181&1705&76 \cr
 & &16.27.11.13.17.223&221&216& & &32.25.11.19.31&775&304 \cr
\noalign{\hrule}
 & &25.11.17.31.37&2911&1764& & &9.25.7.11.17.19&59&4 \cr
902&5362225&8.9.49.41.71&187&310&920&5595975&8.5.17.19.59&33&52 \cr
 & &32.3.5.7.11.17.31&21&16& & &64.3.11.13.59&59&416 \cr
\noalign{\hrule}
 & &9.5.13.61.151&391&2354& & &27.5.343.121&25&52 \cr
903&5388435&4.11.17.23.107&135&118&921&5602905&8.125.49.11.13&207&332 \cr
 & &16.27.5.59.107&321&472& & &64.9.13.23.83&1079&736 \cr
\noalign{\hrule}
 & &3.5.7.11.31.151&153&2& & &3.5.7.17.43.73&247&118 \cr
904&5406555&4.27.7.11.17&25&52&922&5603115&4.7.13.17.19.59&215&198 \cr
 & &32.25.13.17&17&1040& & &16.9.5.11.13.19.43&247&264 \cr
\noalign{\hrule}
 & &9.7.11.13.601&3715&2896& & &3.5.49.11.17.41&73&114 \cr
905&5414409&32.5.181.743&281&462&923&5635245&4.9.5.49.19.73&451&206 \cr
 & &128.3.5.7.11.281&281&320& & &16.11.19.41.103&103&152 \cr
\noalign{\hrule}
 & &27.11.13.23.61&113&140& & &9.11.23.37.67&237&170 \cr
906&5416983&8.5.7.13.61.113&63&2&924&5644683&4.27.5.17.23.79&71&44 \cr
 & &32.9.49.113&113&784& & &32.11.17.71.79&1343&1136 \cr
\noalign{\hrule}
 & &9.7.11.41.191&29&70& & &3.7.43.61.103&95&34 \cr
907&5426883&4.5.49.29.191&615&806&925&5673549&4.5.7.17.19.103&317&198 \cr
 & &16.3.25.13.31.41&325&248& & &16.9.11.19.317&951&1672 \cr
\noalign{\hrule}
 & &9.7.13.29.229&41&50& & &9.7.19.47.101&55&8 \cr
908&5438979&4.25.29.41.229&187&42&926&5682159&16.5.11.19.101&155&54 \cr
 & &16.3.5.7.11.17.41&451&680& & &64.27.25.31&93&800 \cr
\noalign{\hrule}
 & &5.31.41.857&351&506& & &9.5.17.43.173&133&82 \cr
909&5446235&4.27.11.13.23.41&133&10&927&5690835&4.3.7.19.41.173&649&130 \cr
 & &16.9.5.7.19.23&207&1064& & &16.5.7.11.13.59&413&1144 \cr
\noalign{\hrule}
 & &27.11.13.17.83&259&820& & &9.25.19.31.43&4543&5132 \cr
910&5447871&8.9.5.7.37.41&31&32&928&5698575&8.7.11.59.1283&435&848 \cr
 & &512.5.31.37.41&7585&7936& & &256.3.5.11.29.53&1537&1408 \cr
\noalign{\hrule}
 & &9.11.17.41.79&4189&3478& & &3.7.19.103.139&2449&3422 \cr
911&5451237&4.37.47.59.71&65&6&929&5712483&4.29.31.59.79&1881&2780 \cr
 & &16.3.5.13.37.47&611&1480& & &32.9.5.11.19.139&55&48 \cr
\noalign{\hrule}
 & &3.7.11.19.29.43&1975&1766& & &9.5.11.13.29.31&71&72 \cr
912&5473083&4.25.7.79.883&279&274&930&5785065&16.81.5.29.31.71&5123&628 \cr
 & &16.9.5.31.137.883&21235&21192& & &128.47.109.157&7379&6976 \cr
\noalign{\hrule}
 & &27.25.23.353&817&242& & &9.25.149.173&2641&1084 \cr
913&5480325&4.9.121.19.43&65&56&931&5799825&8.19.139.271&205&66 \cr
 & &64.5.7.13.19.43&1729&1376& & &32.3.5.11.19.41&209&656 \cr
\noalign{\hrule}
 & &25.7.23.29.47&351&374& & &9.5.29.61.73&13&74 \cr
914&5486075&4.27.7.11.13.17.47&25&116&932&5811165&4.3.5.13.37.73&29&44 \cr
 & &32.9.25.11.17.29&153&176& & &32.11.13.29.37&481&176 \cr
\noalign{\hrule}
 & &11.31.37.439&235&204& & &9.7.19.43.113&377&440 \cr
915&5538863&8.3.5.11.17.37.47&1&186&933&5816223&16.5.11.13.29.113&129&16 \cr
 & &32.9.31.47&423&16& & &512.3.11.13.43&143&256 \cr
\noalign{\hrule}
 & &9.5.23.53.101&53&62& & &3.5.121.169.19&7&202 \cr
916&5540355&4.31.2809.101&161&2970&934&5827965&4.7.11.13.101&551&450 \cr
 & &16.27.5.7.11.23&77&24& & &16.9.25.19.29&145&24 \cr
\noalign{\hrule}
 & &27.5.7.11.13.41&43&412& & &9.49.97.137&73&24 \cr
917&5540535&8.3.11.43.103&391&82&935&5860449&16.27.73.137&55&82 \cr
 & &32.17.23.41&17&368& & &64.5.11.41.73&2993&1760 \cr
\noalign{\hrule}
 & &3.25.121.13.47&19&344& & &5.7.11.289.53&431&1014 \cr
918&5544825&16.19.43.47&45&2&936&5897045&4.3.7.169.431&261&170 \cr
 & &64.9.5.19&3&608& & &16.27.5.13.17.29&351&232 \cr
\noalign{\hrule}
}%
}
$$
\eject
\vglue -23 pt
\noindent\hskip 1 in\hbox to 6.5 in{\ 937 -- 972 \hfill\fbd 5901225 -- 6551919\frb}
\vskip -9 pt
$$
\vbox{
\nointerlineskip
\halign{\strut
    \vrule \ \ \hfil \frb #\ 
   &\vrule \hfil \ \ \fbb #\frb\ 
   &\vrule \hfil \ \ \frb #\ \hfil
   &\vrule \hfil \ \ \frb #\ 
   &\vrule \hfil \ \ \frb #\ \ \vrule \hskip 2 pt
   &\vrule \ \ \hfil \frb #\ 
   &\vrule \hfil \ \ \fbb #\frb\ 
   &\vrule \hfil \ \ \frb #\ \hfil
   &\vrule \hfil \ \ \frb #\ 
   &\vrule \hfil \ \ \frb #\ \vrule \cr%
\noalign{\hrule}
 & &3.25.11.23.311&477&788& & &9.7.121.19.43&17&116 \cr
937&5901225&8.27.5.53.197&31&166&955&6227991&8.11.17.29.43&245&228 \cr
 & &32.31.53.83&2573&848& & &64.3.5.49.19.29&203&160 \cr
\noalign{\hrule}
 & &3.13.17.59.151&35&186& & &5.11.19.47.127&13&222 \cr
938&5906667&4.9.5.7.31.59&169&110&956&6237605&4.3.13.37.127&177&304 \cr
 & &16.25.7.11.169&275&728& & &128.9.19.59&59&576 \cr
\noalign{\hrule}
 & &81.49.19.79&775&726& & &9.11.13.37.131&7&124 \cr
939&5957469&4.243.25.121.31&437&778&957&6238089&8.7.11.31.37&135&124 \cr
 & &16.5.11.19.23.389&1945&2024& & &64.27.5.961&961&480 \cr
\noalign{\hrule}
 & &9.5.13.31.331&781&874& & &7.19.31.37.41&495&208 \cr
940&6002685&4.3.11.13.19.23.71&35&178&958&6254591&32.9.5.11.13.31&41&52 \cr
 & &16.5.7.19.23.89&1691&1288& & &256.3.5.169.41&507&640 \cr
\noalign{\hrule}
 & &27.25.7.19.67&319&346& & &9.7.11.13.17.41&575&1108 \cr
941&6014925&4.5.11.29.67.173&4351&666&959&6279273&8.25.7.23.277&541&264 \cr
 & &16.9.19.37.229&229&296& & &128.3.5.11.541&541&320 \cr
\noalign{\hrule}
 & &11.13.23.31.59&41&18& & &9.125.71.79&527&598 \cr
942&6015581&4.9.11.13.31.41&25&118&960&6310125&4.13.17.23.31.79&189&110 \cr
 & &16.3.25.41.59&75&328& & &16.27.5.7.11.17.31&1023&952 \cr
\noalign{\hrule}
 & &67.89.1013&473&540& & &3.121.13.23.59&41&18 \cr
943&6040519&8.27.5.11.43.89&109&20&961&6403683&4.27.121.13.41&115&236 \cr
 & &64.9.25.11.109&2475&3488& & &32.5.23.41.59&41&80 \cr
\noalign{\hrule}
 & &3.5.43.83.113&51&164& & &25.7.23.37.43&1083&508 \cr
944&6049455&8.9.17.41.83&143&226&962&6403775&8.3.7.361.127&625&264 \cr
 & &32.11.13.17.113&143&272& & &128.9.625.11&275&576 \cr
\noalign{\hrule}
 & &27.13.41.421&6985&7406& & &9.11.17.37.103&2747&1064 \cr
945&6058611&4.5.7.11.529.127&29&6&963&6413913&16.7.19.41.67&11&30 \cr
 & &16.3.11.23.29.127&2921&2552& & &64.3.5.7.11.67&469&160 \cr
\noalign{\hrule}
 & &27.11.113.181&893&350& & &3.13.37.61.73&25&12 \cr
946&6074541&4.9.25.7.19.47&451&404&964&6425679&8.9.25.61.73&481&176 \cr
 & &32.5.7.11.41.101&1435&1616& & &256.5.11.13.37&55&128 \cr
\noalign{\hrule}
 & &5.7.31.71.79&233&162& & &9.5.11.13.17.59&871&812 \cr
947&6085765&4.81.7.31.233&79&110&965&6454305&8.5.7.169.29.67&759&424 \cr
 & &16.3.5.11.79.233&233&264& & &128.3.11.23.29.53&1537&1472 \cr
\noalign{\hrule}
 & &3.7.31.47.199&425&226& & &27.5.17.29.97&1133&2782 \cr
948&6088803&4.25.17.47.113&117&682&966&6455835&4.11.13.103.107&57&46 \cr
 & &16.9.5.11.13.31&65&264& & &16.3.13.19.23.107&2461&1976 \cr
\noalign{\hrule}
 & &49.11.17.23.29&79&240& & &9.5.11.13.19.53&53&118 \cr
949&6111721&32.3.5.7.17.79&257&138&967&6480045&4.11.2809.59&1729&1080 \cr
 & &128.9.23.257&257&576& & &64.27.5.7.13.19&21&32 \cr
\noalign{\hrule}
 & &3.5.11.131.283&469&186& & &81.49.11.149&25&124 \cr
950&6117045&4.9.7.11.31.67&131&472&968&6505191&8.9.25.49.31&143&298 \cr
 & &64.7.59.131&59&224& & &32.5.11.13.149&65&16 \cr
\noalign{\hrule}
 & &25.13.113.167&219&106& & &7.41.61.373&55&2556 \cr
951&6133075&4.3.53.73.167&57&110&969&6530111&8.9.5.11.71&41&30 \cr
 & &16.9.5.11.19.73&803&1368& & &32.27.25.41&27&400 \cr
\noalign{\hrule}
 & &31.127.1567&847&720& & &3.5.7.19.29.113&129&16 \cr
952&6169279&32.9.5.7.121.31&799&1016&970&6537615&32.9.7.19.43&377&440 \cr
 & &512.3.17.47.127&799&768& & &512.5.11.13.29&143&256 \cr
\noalign{\hrule}
 & &9.121.13.19.23&1299&1000& & &7.13.23.53.59&315&374 \cr
953&6186609&16.27.125.433&121&554&971&6544811&4.9.5.49.11.17.23&59&10 \cr
 & &64.5.121.277&277&160& & &16.3.25.11.17.59&425&264 \cr
\noalign{\hrule}
 & &27.125.7.263&219&44& & &9.11.289.229&21&208 \cr
954&6213375&8.81.5.11.73&263&628&972&6551919&32.27.7.13.17&205&16 \cr
 & &64.157.263&157&32& & &1024.5.41&41&2560 \cr
\noalign{\hrule}
}%
}
$$
\eject
\vglue -23 pt
\noindent\hskip 1 in\hbox to 6.5 in{\ 973 -- 1008 \hfill\fbd 6554275 -- 7160825\frb}
\vskip -9 pt
$$
\vbox{
\nointerlineskip
\halign{\strut
    \vrule \ \ \hfil \frb #\ 
   &\vrule \hfil \ \ \fbb #\frb\ 
   &\vrule \hfil \ \ \frb #\ \hfil
   &\vrule \hfil \ \ \frb #\ 
   &\vrule \hfil \ \ \frb #\ \ \vrule \hskip 2 pt
   &\vrule \ \ \hfil \frb #\ 
   &\vrule \hfil \ \ \fbb #\frb\ 
   &\vrule \hfil \ \ \frb #\ \hfil
   &\vrule \hfil \ \ \frb #\ 
   &\vrule \hfil \ \ \frb #\ \vrule \cr%
\noalign{\hrule}
 & &25.7.13.43.67&719&786& & &7.11.19.47.101&155&54 \cr
973&6554275&4.3.5.13.131.719&2649&946&991&6944861&4.27.5.7.31.47&55&8 \cr
 & &16.9.11.43.883&883&792& & &64.3.25.11.31&93&800 \cr
\noalign{\hrule}
 & &3.5.7.11.13.19.23&39&94& & &27.7.11.13.257&2059&1802 \cr
974&6561555&4.9.169.23.47&4949&4780&992&6945939&4.7.17.29.53.71&495&2 \cr
 & &32.5.49.101.239&1673&1616& & &16.9.5.11.53&265&8 \cr
\noalign{\hrule}
 & &9.5.41.43.83&143&226& & &3.11.169.29.43&149&20 \cr
975&6584805&4.5.11.13.43.113&51&164&993&6954519&8.5.11.29.149&85&234 \cr
 & &32.3.11.13.17.41&143&272& & &32.9.25.13.17&425&48 \cr
\noalign{\hrule}
 & &5.13.19.53.101&6293&6798& & &3.49.121.17.23&3551&2378 \cr
976&6610955&4.3.7.11.29.31.103&53&24&994&6954717&4.29.41.53.67&605&2142 \cr
 & &64.9.31.53.103&927&992& & &16.9.5.7.121.17&15&8 \cr
\noalign{\hrule}
 & &3.19.23.61.83&985&924& & &9.5.11.13.23.47&587&2702 \cr
977&6637593&8.9.5.7.11.19.197&13&184&995&6956235&4.7.193.587&197&390 \cr
 & &128.5.7.11.13.23&455&704& & &16.3.5.7.13.197&197&56 \cr
\noalign{\hrule}
 & &81.7.13.17.53&787&760& & &125.19.29.101&277&2652 \cr
978&6641271&16.3.5.19.53.787&473&314&996&6956375&8.3.13.17.277&145&132 \cr
 & &64.5.11.19.43.157&14915&15136& & &64.9.5.11.17.29&187&288 \cr
\noalign{\hrule}
 & &3.169.23.571&435&136& & &3.5.7.11.19.317&351&34 \cr
979&6658431&16.9.5.13.17.29&539&46&997&6956565&4.81.13.17.19&121&202 \cr
 & &64.49.11.23&539&32& & &16.121.13.101&143&808 \cr
\noalign{\hrule}
 & &81.5.11.19.79&91&118& & &11.31.67.307&11577&11270 \cr
980&6686955&4.3.5.7.13.59.79&29&266&998&7014029&4.3.5.49.17.23.227&599&990 \cr
 & &16.49.13.19.29&637&232& & &16.27.25.7.11.599&4725&4792 \cr
\noalign{\hrule}
 & &3.5.7.23.47.59&1&22& & &121.103.563&221&342 \cr
981&6696795&4.5.11.47.59&207&442&999&7016669&4.9.13.17.19.103&265&44 \cr
 & &16.9.13.17.23&17&312& & &32.3.5.11.19.53&285&848 \cr
\noalign{\hrule}
 & &3.11.239.859&5083&4366& & &9.11.17.59.71&125&656 \cr
982&6774933&4.13.17.23.37.59&45&436&1000&7050087&32.125.17.41&71&54 \cr
 & &32.9.5.59.109&885&1744& & &128.27.41.71&123&64 \cr
\noalign{\hrule}
 & &5.11.169.17.43&9&178& & &5.7.13.37.419&7979&7524 \cr
983&6794645&4.9.5.43.89&241&26&1001&7053865&8.9.11.19.79.101&425&444 \cr
 & &16.3.13.241&723&8& & &64.27.25.17.37.101&2727&2720 \cr
\noalign{\hrule}
 & &9.5.7.11.37.53&23&76& & &243.7.11.13.29&1525&1148 \cr
984&6794865&8.5.7.19.23.37&159&26&1002&7054047&8.25.49.41.61&1863&638 \cr
 & &32.3.13.23.53&23&208& & &32.81.11.23.29&23&16 \cr
\noalign{\hrule}
 & &27.7.13.47.59&41&50& & &3.7.23.97.151&317&166 \cr
985&6813261&4.3.25.41.47.59&1079&1694&1003&7074501&4.83.97.317&207&110 \cr
 & &16.5.7.121.13.83&605&664& & &16.9.5.11.23.83&415&264 \cr
\noalign{\hrule}
 & &9.5.19.23.347&131&154& & &27.11.17.23.61&79&470 \cr
986&6823755&4.3.7.11.131.347&23&370&1004&7083747&4.3.5.11.47.79&43&122 \cr
 & &16.5.7.11.23.37&259&88& & &16.43.47.61&47&344 \cr
\noalign{\hrule}
 & &3.25.7.13.17.59&19&2& & &27.5.19.47.59&91&44 \cr
987&6845475&4.25.13.19.59&371&396&1005&7112745&8.7.11.13.19.59&153&94 \cr
 & &32.9.7.11.19.53&627&848& & &32.9.7.11.17.47&119&176 \cr
\noalign{\hrule}
 & &25.11.13.17.113&393&172& & &5.7.23.83.107&27&134 \cr
988&6867575&8.3.5.11.43.131&113&102&1006&7149205&4.27.5.67.83&107&308 \cr
 & &32.9.17.113.131&131&144& & &32.9.7.11.107&99&16 \cr
\noalign{\hrule}
 & &3.47.137.359&515&562& & &13.17.29.1117&9683&9306 \cr
989&6934803&4.5.103.137.281&209&72&1007&7158853&4.9.11.23.47.421&565&986 \cr
 & &64.9.5.11.19.103&3399&3040& & &16.3.5.17.23.29.113&565&552 \cr
\noalign{\hrule}
 & &3.11.13.19.23.37&3&250& & &25.7.17.29.83&199&216 \cr
990&6936501&4.9.125.37&97&88&1008&7160825&16.27.5.7.29.199&209&806 \cr
 & &64.25.11.97&97&800& & &64.9.11.13.19.31&5301&4576 \cr
\noalign{\hrule}
}%
}
$$
\eject
\vglue -23 pt
\noindent\hskip 1 in\hbox to 6.5 in{\ 1009 -- 1044 \hfill\fbd 7161165 -- 7654535\frb}
\vskip -9 pt
$$
\vbox{
\nointerlineskip
\halign{\strut
    \vrule \ \ \hfil \frb #\ 
   &\vrule \hfil \ \ \fbb #\frb\ 
   &\vrule \hfil \ \ \frb #\ \hfil
   &\vrule \hfil \ \ \frb #\ 
   &\vrule \hfil \ \ \frb #\ \ \vrule \hskip 2 pt
   &\vrule \ \ \hfil \frb #\ 
   &\vrule \hfil \ \ \fbb #\frb\ 
   &\vrule \hfil \ \ \frb #\ \hfil
   &\vrule \hfil \ \ \frb #\ 
   &\vrule \hfil \ \ \frb #\ \vrule \cr%
\noalign{\hrule}
 & &9.5.11.17.23.37&7&62& & &3.5.7.11.13.17.29&167&218 \cr
1009&7161165&4.3.7.17.31.37&319&208&1027&7402395&4.13.29.109.167&495&2666 \cr
 & &128.7.11.13.29&203&832& & &16.9.5.11.31.43&129&248 \cr
\noalign{\hrule}
 & &5.49.23.31.41&3&38& & &3.25.19.59.89&693&782 \cr
1010&7162085&4.3.7.19.23.31&327&110&1028&7482675&4.27.7.11.17.19.23&13&310 \cr
 & &16.9.5.11.109&1199&72& & &16.5.7.13.23.31&403&1288 \cr
\noalign{\hrule}
 & &3.5.49.11.887&713&174& & &3.5.7.121.19.31&205&422 \cr
1011&7171395&4.9.5.23.29.31&77&68&1029&7483245&4.25.11.41.211&2863&2412 \cr
 & &32.7.11.17.23.31&391&496& & &32.9.7.67.409&1227&1072 \cr
\noalign{\hrule}
 & &25.49.11.13.41&1037&972& & &3.5.343.31.47&85&132 \cr
1012&7182175&8.243.5.11.17.61&2147&3038&1030&7496265&8.9.25.49.11.17&83&358 \cr
 & &32.3.49.19.31.113&1767&1808& & &32.17.83.179&3043&1328 \cr
\noalign{\hrule}
 & &5.11.13.19.529&2681&36& & &7.11.19.23.223&115&108 \cr
1013&7186465&8.9.7.383&193&190&1031&7503727&8.27.5.11.19.529&49&578 \cr
 & &32.3.5.7.19.193&193&336& & &32.9.5.49.289&1445&1008 \cr
\noalign{\hrule}
 & &25.7.11.53.71&9&44& & &3.7.11.19.29.59&861&850 \cr
1014&7243775&8.9.5.121.71&117&238&1032&7509579&4.9.25.49.17.19.41&11&844 \cr
 & &32.81.7.13.17&1053&272& & &32.5.11.41.211&1055&656 \cr
\noalign{\hrule}
 & &5.7.43.61.79&253&174& & &27.5.19.29.101&253&298 \cr
1015&7252595&4.3.5.11.23.29.43&93&122&1033&7512885&4.3.11.23.101.149&85&16 \cr
 & &16.9.11.23.31.61&713&792& & &128.5.11.17.149&1639&1088 \cr
\noalign{\hrule}
 & &25.311.937&3419&4356& & &3.49.17.31.97&5&26 \cr
1016&7285175&8.9.121.13.263&313&50&1034&7514493&4.5.7.13.17.97&11&108 \cr
 & &32.3.25.13.313&313&624& & &32.27.5.11.13&65&1584 \cr
\noalign{\hrule}
 & &3.5.43.89.127&23&66& & &125.17.53.67&713&2838 \cr
1017&7290435&4.9.5.11.23.127&41&86&1035&7545875&4.3.11.23.31.43&21&10 \cr
 & &16.11.23.41.43&451&184& & &16.9.5.7.23.43&387&1288 \cr
\noalign{\hrule}
 & &3.7.121.169.17&43&230& & &27.5.17.37.89&721&1166 \cr
1018&7300293&4.5.11.13.23.43&153&406&1036&7557435&4.9.7.11.53.103&37&26 \cr
 & &16.9.5.7.17.29&145&24& & &16.13.37.53.103&689&824 \cr
\noalign{\hrule}
 & &27.11.13.31.61&73&476& & &243.11.19.149&25&124 \cr
1019&7301151&8.3.7.11.17.73&31&20&1037&7567263&8.27.25.19.31&43&632 \cr
 & &64.5.7.31.73&511&160& & &128.43.79&79&2752 \cr
\noalign{\hrule}
 & &3.5.49.17.587&29&22& & &3.11.17.103.131&575&558 \cr
1020&7334565&4.5.7.11.29.587&351&2584&1038&7569573&4.27.25.23.31.131&17&638 \cr
 & &64.27.13.17.19&171&416& & &16.5.11.17.29.31&145&248 \cr
\noalign{\hrule}
 & &5.7.11.17.19.59&13&72& & &27.121.13.179&745&866 \cr
1021&7336945&16.9.7.11.13.19&59&150&1039&7602309&4.3.5.13.149.433&319&1618 \cr
 & &64.27.25.59&27&160& & &16.5.11.29.809&809&1160 \cr
\noalign{\hrule}
 & &3.11.13.17.19.53&305&358& & &3.5.7.23.47.67&37&198 \cr
1022&7344051&4.5.11.19.61.179&405&1564&1040&7604835&4.27.11.37.67&115&182 \cr
 & &32.81.25.17.23&575&432& & &16.5.7.13.23.37&37&104 \cr
\noalign{\hrule}
 & &17.1681.257&477&220& & &9.7.11.23.479&899&1378 \cr
1023&7344289&8.9.5.11.41.53&35&88&1041&7634781&4.7.13.29.31.53&297&80 \cr
 & &128.3.25.7.121&3025&1344& & &128.27.5.11.53&159&320 \cr
\noalign{\hrule}
 & &81.5.11.17.97&217&268& & &5.13.19.23.269&135&112 \cr
1024&7346295&8.27.7.11.31.67&2561&3298&1042&7640945&32.27.25.7.269&247&22 \cr
 & &32.13.17.97.197&197&208& & &128.3.7.11.13.19&77&192 \cr
\noalign{\hrule}
 & &49.19.53.149&99&50& & &3.25.7.13.19.59&369&44 \cr
1025&7352107&4.9.25.11.19.53&149&434&1043&7650825&8.27.11.19.41&91&118 \cr
 & &16.3.5.7.31.149&93&40& & &32.7.13.41.59&41&16 \cr
\noalign{\hrule}
 & &7.11.19.47.107&793&2826& & &5.49.157.199&201&44 \cr
1026&7357427&4.9.13.61.157&205&266&1044&7654535&8.3.11.67.199&133&66 \cr
 & &16.3.5.7.13.19.41&205&312& & &32.9.7.121.19&1089&304 \cr
\noalign{\hrule}
}%
}
$$
\eject
\vglue -23 pt
\noindent\hskip 1 in\hbox to 6.5 in{\ 1045 -- 1080 \hfill\fbd 7657155 -- 8430525\frb}
\vskip -9 pt
$$
\vbox{
\nointerlineskip
\halign{\strut
    \vrule \ \ \hfil \frb #\ 
   &\vrule \hfil \ \ \fbb #\frb\ 
   &\vrule \hfil \ \ \frb #\ \hfil
   &\vrule \hfil \ \ \frb #\ 
   &\vrule \hfil \ \ \frb #\ \ \vrule \hskip 2 pt
   &\vrule \ \ \hfil \frb #\ 
   &\vrule \hfil \ \ \fbb #\frb\ 
   &\vrule \hfil \ \ \frb #\ \hfil
   &\vrule \hfil \ \ \frb #\ 
   &\vrule \hfil \ \ \frb #\ \vrule \cr%
\noalign{\hrule}
 & &9.5.11.31.499&7751&7718& & &27.5.49.23.53&109&374 \cr
1045&7657155&4.3.5.17.23.227.337&341&4&1063&8063685&4.9.7.11.17.109&53&46 \cr
 & &32.11.17.31.227&227&272& & &16.17.23.53.109&109&136 \cr
\noalign{\hrule}
 & &3.23.31.37.97&1131&1100& & &5.17.151.631&369&386 \cr
1046&7676871&8.9.25.11.13.29.37&673&992&1064&8098885&4.9.41.193.631&605&26 \cr
 & &512.5.13.31.673&3365&3328& & &16.3.5.121.13.41&1599&968 \cr
\noalign{\hrule}
 & &27.5.23.37.67&71&44& & &9.5.7.29.887&77&68 \cr
1047&7697295&8.11.37.67.71&237&170&1065&8102745&8.49.11.17.887&713&174 \cr
 & &32.3.5.17.71.79&1343&1136& & &32.3.17.23.29.31&391&496 \cr
\noalign{\hrule}
 & &27.43.61.109&253&296& & &3.13.43.47.103&125&434 \cr
1048&7719489&16.3.11.23.37.109&1&110&1066&8118357&4.125.7.31.47&11&36 \cr
 & &64.5.121.23&2783&160& & &32.9.5.7.11.31&341&1680 \cr
\noalign{\hrule}
 & &7.169.31.211&3355&3186& & &9.7.11.53.223&2735&2512 \cr
1049&7738003&4.27.5.7.11.59.61&589&650&1067&8190567&32.5.7.157.547&3021&2474 \cr
 & &16.9.125.11.13.19.31&1375&1368& & &128.3.19.53.1237&1237&1216 \cr
\noalign{\hrule}
 & &9.19.131.347&23&370& & &9.11.17.31.157&2869&200 \cr
1050&7773147&4.3.5.19.23.37&131&154&1068&8191161&16.25.19.151&85&66 \cr
 & &16.7.11.37.131&259&88& & &64.3.125.11.17&125&32 \cr
\noalign{\hrule}
 & &7.121.67.137&295&174& & &9.11.17.31.157&3275&1592 \cr
1051&7774613&4.3.5.29.59.137&787&924&1069&8191161&16.25.131.199&165&34 \cr
 & &32.9.5.7.11.787&787&720& & &64.3.125.11.17&125&32 \cr
\noalign{\hrule}
 & &5.23.59.1153&2853&2912& & &13.19.73.457&2277&3664 \cr
1052&7823105&64.9.7.13.23.317&9&308&1070&8240167&32.9.11.23.229&265&494 \cr
 & &512.81.49.11&3969&2816& & &128.3.5.13.19.53&265&192 \cr
\noalign{\hrule}
 & &11.169.41.103&141&310& & &5.11.31.47.103&87&428 \cr
1053&7850557&4.3.5.31.47.103&169&66&1071&8253905&8.3.29.47.107&9&38 \cr
 & &16.9.11.169.31&31&72& & &32.27.19.107&2889&304 \cr
\noalign{\hrule}
 & &9.125.49.11.13&31&94& & &9.5.31.61.97&589&284 \cr
1054&7882875&4.7.11.13.31.47&125&216&1072&8254215&8.19.961.71&1155&194 \cr
 & &64.27.125.47&47&96& & &32.3.5.7.11.97&77&16 \cr
\noalign{\hrule}
 & &5.41.137.281&243&38& & &9.5.49.19.199&289&308 \cr
1055&7891885&4.243.19.137&55&82&1073&8337105&8.3.5.343.11.289&577&1138 \cr
 & &16.9.5.11.19.41&209&72& & &32.17.569.577&9673&9232 \cr
\noalign{\hrule}
 & &3.11.13.19.971&557&414& & &81.121.23.37&227&106 \cr
1056&7914621&4.27.19.23.557&535&22&1074&8340651&4.9.23.53.227&77&130 \cr
 & &16.5.11.23.107&115&856& & &16.5.7.11.13.227&1135&728 \cr
\noalign{\hrule}
 & &5.11.17.67.127&261&74& & &9.5.7.11.19.127&53&46 \cr
1057&7955915&4.9.29.37.127&209&172&1075&8361045&4.5.19.23.53.127&5423&6642 \cr
 & &32.3.11.19.29.43&1247&912& & &16.81.11.17.29.41&1189&1224 \cr
\noalign{\hrule}
 & &9.7.11.83.139&3379&1850& & &13.17.109.349&3195&1342 \cr
1058&7995141&4.25.31.37.109&147&38&1076&8407061&4.9.5.11.61.71&229&442 \cr
 & &16.3.5.49.19.31&665&248& & &16.3.5.13.17.229&229&120 \cr
\noalign{\hrule}
 & &5.7.11.59.353&1207&558& & &3.25.7.11.31.47&603&572 \cr
1059&8018395&4.9.7.17.31.71&143&214&1077&8414175&8.27.7.121.13.67&2825&5282 \cr
 & &16.3.11.13.31.107&1209&856& & &32.25.19.113.139&2147&2224 \cr
\noalign{\hrule}
 & &7.13.157.563&7425&6862& & &17.41.43.281&225&506 \cr
1060&8043581&4.27.25.11.47.73&49&754&1078&8421851&4.9.25.11.23.41&119&86 \cr
 & &16.9.5.49.13.29&203&360& & &16.3.5.7.17.23.43&161&120 \cr
\noalign{\hrule}
 & &3.5.1331.13.31&53&68& & &9.5.17.101.109&131&22 \cr
1061&8045895&8.11.13.17.31.53&493&90&1079&8421885&4.5.11.101.131&187&318 \cr
 & &32.9.5.289.29&867&464& & &16.3.121.17.53&121&424 \cr
\noalign{\hrule}
 & &9.19.529.89&179&350& & &9.25.89.421&323&98 \cr
1062&8050851&4.25.7.89.179&759&136&1080&8430525&4.49.17.19.89&473&150 \cr
 & &64.3.5.11.17.23&85&352& & &16.3.25.7.11.43&77&344 \cr
\noalign{\hrule}
}%
}
$$
\eject
\vglue -23 pt
\noindent\hskip 1 in\hbox to 6.5 in{\ 1081 -- 1116 \hfill\fbd 8430569 -- 8954495\frb}
\vskip -9 pt
$$
\vbox{
\nointerlineskip
\halign{\strut
    \vrule \ \ \hfil \frb #\ 
   &\vrule \hfil \ \ \fbb #\frb\ 
   &\vrule \hfil \ \ \frb #\ \hfil
   &\vrule \hfil \ \ \frb #\ 
   &\vrule \hfil \ \ \frb #\ \ \vrule \hskip 2 pt
   &\vrule \ \ \hfil \frb #\ 
   &\vrule \hfil \ \ \fbb #\frb\ 
   &\vrule \hfil \ \ \frb #\ \hfil
   &\vrule \hfil \ \ \frb #\ 
   &\vrule \hfil \ \ \frb #\ \vrule \cr%
\noalign{\hrule}
 & &7.59.137.149&4875&3916& & &125.7.11.907&261&646 \cr
1081&8430569&8.3.125.11.13.89&5&6&1099&8729875&4.9.25.17.19.29&9&484 \cr
 & &32.9.625.13.89&10413&10000& & &32.81.121&891&16 \cr
\noalign{\hrule}
 & &25.13.841.31&427&414& & &3.11.13.47.433&475&42 \cr
1082&8473075&4.9.25.7.23.31.61&319&394&1100&8730579&4.9.25.7.13.19&473&382 \cr
 & &16.3.7.11.29.61.197&4697&4728& & &16.5.11.43.191&955&344 \cr
\noalign{\hrule}
 & &9.7.31.43.101&173&44& & &5.11.43.47.79&221&174 \cr
1083&8481879&8.3.11.101.173&209&310&1101&8781245&4.3.11.13.17.29.43&425&48 \cr
 & &32.5.121.19.31&605&304& & &128.9.25.289&1445&576 \cr
\noalign{\hrule}
 & &27.5.7.11.19.43&221&166& & &9.169.53.109&89&80 \cr
1084&8492715&4.3.7.13.17.19.83&47&86&1102&8786817&32.5.53.89.109&187&78 \cr
 & &16.17.43.47.83&799&664& & &128.3.11.13.17.89&979&1088 \cr
\noalign{\hrule}
 & &3.25.17.59.113&69&44& & &31.2809.101&161&2970 \cr
1085&8500425&8.9.11.17.23.59&551&452&1103&8794979&4.27.5.7.11.23&53&62 \cr
 & &64.19.23.29.113&667&608& & &16.3.7.11.31.53&77&24 \cr
\noalign{\hrule}
 & &3.7.47.89.97&85&182& & &3.5.11.19.29.97&17&302 \cr
1086&8520771&4.5.49.13.17.47&99&146&1104&8818755&4.17.97.151&1235&1332 \cr
 & &16.9.11.13.17.73&2431&1752& & &32.9.5.13.19.37&111&208 \cr
\noalign{\hrule}
 & &27.5.17.47.79&437&832& & &25.7.31.1627&901&726 \cr
1087&8521335&128.13.17.19.23&319&72&1105&8826475&4.3.121.17.31.53&621&280 \cr
 & &2048.9.11.29&319&1024& & &64.81.5.7.11.23&891&736 \cr
\noalign{\hrule}
 & &9.25.17.23.97&83&308& & &121.19.23.167&135&302 \cr
1088&8533575&8.7.11.83.97&43&54&1106&8830459&4.27.5.121.151&529&76 \cr
 & &32.27.7.43.83&903&1328& & &32.9.19.529&23&144 \cr
\noalign{\hrule}
 & &49.11.13.23.53&57&34& & &3.13.31.71.103&95&308 \cr
1089&8541533&4.3.7.11.17.19.53&65&12&1107&8841417&8.5.7.11.19.103&899&234 \cr
 & &32.9.5.13.17.19&855&272& & &32.9.13.29.31&87&16 \cr
\noalign{\hrule}
 & &7.11.13.83.103&239&342& & &27.5.7.11.23.37&1411&1004 \cr
1090&8557549&4.9.11.13.19.239&4067&3820&1108&8846145&8.9.17.83.251&49&202 \cr
 & &32.3.5.49.83.191&573&560& & &32.49.83.101&707&1328 \cr
\noalign{\hrule}
 & &3.7.43.53.179&275&96& & &3.5.7.19.23.193&39&154 \cr
1091&8566761&64.9.25.11.43&91&134&1109&8855805&4.9.49.11.13.19&965&916 \cr
 & &256.7.11.13.67&871&1408& & &32.5.13.193.229&229&208 \cr
\noalign{\hrule}
 & &3.7.11.13.19.151&457&600& & &5.11.23.47.149&81&34 \cr
1092&8615607&16.9.25.19.457&1057&1228&1110&8858795&4.81.11.17.149&25&124 \cr
 & &128.5.7.151.307&307&320& & &32.9.25.17.31&279&1360 \cr
\noalign{\hrule}
 & &9.5.13.361.41&1645&1604& & &3.7.13.19.29.59&3431&1702 \cr
1093&8658585&8.25.7.13.47.401&363&38&1111&8874957&4.23.37.47.73&55&18 \cr
 & &32.3.7.121.19.47&847&752& & &16.9.5.11.23.47&759&1880 \cr
\noalign{\hrule}
 & &27.5.11.19.307&23&34& & &5.13.17.83.97&297&782 \cr
1094&8662005&4.9.5.17.23.307&2527&992&1112&8896355&4.27.11.289.23&95&194 \cr
 & &256.7.361.31&589&896& & &16.3.5.19.23.97&69&152 \cr
\noalign{\hrule}
 & &5.11.13.17.23.31&27&58& & &5.7.11.19.23.53&153&100 \cr
1095&8666515&4.27.11.13.23.29&315&62&1113&8916985&8.9.125.7.17.19&599&276 \cr
 & &16.243.5.7.31&243&56& & &64.27.23.599&599&864 \cr
\noalign{\hrule}
 & &3.343.11.13.59&75&68& & &11.71.73.157&469&312 \cr
1096&8681673&8.9.25.49.17.59&143&388&1114&8951041&16.3.7.13.67.73&691&180 \cr
 & &64.5.11.13.17.97&485&544& & &128.27.5.691&3455&1728 \cr
\noalign{\hrule}
 & &3.11.13.79.257&247&10& & &3.11.19.109.131&775&666 \cr
1097&8709987&4.5.11.169.19&79&90&1115&8952933&4.27.25.19.31.37&131&644 \cr
 & &16.9.25.19.79&57&200& & &32.7.23.37.131&259&368 \cr
\noalign{\hrule}
 & &5.7.17.107.137&11&96& & &5.11.17.61.157&57&728 \cr
1098&8722105&64.3.7.11.137&107&30&1116&8954495&16.3.7.13.17.19&9&8 \cr
 & &256.9.5.107&9&128& & &256.27.7.13.19&1729&3456 \cr
\noalign{\hrule}
}%
}
$$
\eject
\vglue -23 pt
\noindent\hskip 1 in\hbox to 6.5 in{\ 1117 -- 1152 \hfill\fbd 8954495 -- 9729187\frb}
\vskip -9 pt
$$
\vbox{
\nointerlineskip
\halign{\strut
    \vrule \ \ \hfil \frb #\ 
   &\vrule \hfil \ \ \fbb #\frb\ 
   &\vrule \hfil \ \ \frb #\ \hfil
   &\vrule \hfil \ \ \frb #\ 
   &\vrule \hfil \ \ \frb #\ \ \vrule \hskip 2 pt
   &\vrule \ \ \hfil \frb #\ 
   &\vrule \hfil \ \ \fbb #\frb\ 
   &\vrule \hfil \ \ \frb #\ \hfil
   &\vrule \hfil \ \ \frb #\ 
   &\vrule \hfil \ \ \frb #\ \vrule \cr%
\noalign{\hrule}
 & &5.11.17.61.157&9&8& & &25.169.31.71&3213&1012 \cr
1117&8954495&16.9.5.11.61.157&57&728&1135&9299225&8.27.7.11.17.23&65&142 \cr
 & &256.27.7.13.19&1729&3456& & &32.3.5.13.17.71&17&48 \cr
\noalign{\hrule}
 & &3.49.11.29.191&615&806& & &3.625.121.41&193&182 \cr
1118&8956563&4.9.5.11.13.31.41&29&70&1136&9301875&4.5.7.11.13.41.193&271&2394 \cr
 & &16.25.7.13.29.31&325&248& & &16.9.49.19.271&2793&2168 \cr
\noalign{\hrule}
 & &7.13.137.719&473&486& & &9.5.17.73.167&371&286 \cr
1119&8963773&4.243.11.43.719&221&940&1137&9326115&4.7.11.13.53.167&293&876 \cr
 & &32.9.5.11.13.17.47&2585&2448& & &32.3.13.73.293&293&208 \cr
\noalign{\hrule}
 & &7.11.23.61.83&3303&1394& & &5.7.13.73.281&115&396 \cr
1120&8966573&4.9.17.41.367&245&122&1138&9333415&8.9.25.11.13.23&281&294 \cr
 & &16.3.5.49.17.61&85&168& & &32.27.49.11.281&189&176 \cr
\noalign{\hrule}
 & &25.49.73.101&3599&3774& & &9.5.7.11.37.73&1649&2014 \cr
1121&9031925&4.3.7.17.37.59.61&101&528&1139&9358965&4.7.17.19.53.97&111&790 \cr
 & &128.9.11.59.101&649&576& & &16.3.5.19.37.79&79&152 \cr
\noalign{\hrule}
 & &9.7.11.13.19.53&17&116& & &9.5.43.47.103&187&328 \cr
1122&9072063&8.13.17.29.53&33&20&1140&9367335&16.3.11.17.41.43&47&4 \cr
 & &64.3.5.11.17.29&145&544& & &128.11.41.47&451&64 \cr
\noalign{\hrule}
 & &3.11.13.17.29.43&483&10& & &81.11.13.811&365&446 \cr
1123&9094371&4.9.5.7.13.23&461&344&1141&9393813&4.5.11.13.73.223&183&40 \cr
 & &64.43.461&461&32& & &64.3.25.61.73&1525&2336 \cr
\noalign{\hrule}
 & &25.19.41.467&231&236& & &9.5.13.19.23.37&97&88 \cr
1124&9094825&8.3.5.7.11.19.41.59&169&36&1142&9458865&16.11.13.19.23.97&3&250 \cr
 & &64.27.11.169.59&9971&9504& & &64.3.125.97&97&800 \cr
\noalign{\hrule}
 & &5.17.43.47.53&111&154& & &27.5.11.13.17.29&3589&1834 \cr
1125&9104605&4.3.7.11.17.37.47&71&258&1143&9517365&4.7.37.97.131&81&178 \cr
 & &16.9.37.43.71&333&568& & &16.81.89.131&393&712 \cr
\noalign{\hrule}
 & &9.49.17.23.53&461&3058& & &27.125.11.257&115&142 \cr
1126&9138843&4.11.139.461&69&70&1144&9541125&4.625.11.23.71&703&78 \cr
 & &16.3.5.7.11.23.461&461&440& & &16.3.13.19.23.37&851&1976 \cr
\noalign{\hrule}
 & &3.11.17.43.379&8165&8132& & &9.13.23.53.67&1045&174 \cr
1127&9142617&8.5.17.19.23.71.107&235&1584&1145&9555741&4.27.5.11.19.29&53&82 \cr
 & &256.9.25.11.23.47&3243&3200& & &16.11.19.41.53&451&152 \cr
\noalign{\hrule}
 & &9.17.19.23.137&295&142& & &11.23.103.367&57&310 \cr
1128&9159957&4.5.59.71.137&33&104&1146&9563653&4.3.5.19.31.103&273&242 \cr
 & &64.3.5.11.13.59&715&1888& & &16.9.7.121.13.19&1001&1368 \cr
\noalign{\hrule}
 & &9.7.19.43.179&25&154& & &49.361.541&451&90 \cr
1129&9213309&4.3.25.49.11.19&377&358&1147&9569749&4.9.5.49.11.41&299&152 \cr
 & &16.5.11.13.29.179&377&440& & &64.3.5.13.19.23&195&736 \cr
\noalign{\hrule}
 & &9.11.23.31.131&185&208& & &5.11.19.53.173&39&134 \cr
1130&9246897&32.3.5.11.13.31.37&115&226&1148&9581605&4.3.11.13.53.67&325&258 \cr
 & &128.25.13.23.113&1469&1600& & &16.9.25.169.43&1935&1352 \cr
\noalign{\hrule}
 & &9.7.11.169.79&799&722& & &3.5.7.11.43.193&69&124 \cr
1131&9252243&4.17.361.47.79&351&1150&1149&9585345&8.9.7.23.31.43&1391&58 \cr
 & &16.27.25.13.19.23&475&552& & &32.13.29.107&3103&208 \cr
\noalign{\hrule}
 & &19.3721.131&3105&616& & &9.37.83.349&299&50 \cr
1132&9261569&16.27.5.7.11.23&61&38&1150&9646011&4.3.25.13.23.37&913&62 \cr
 & &64.3.5.7.19.61&35&96& & &16.11.31.83&341&8 \cr
\noalign{\hrule}
 & &3.5.13.29.31.53&7&6& & &3.11.13.71.317&303&620 \cr
1133&9291165&4.9.5.7.29.31.53&169&1474&1151&9655503&8.9.5.11.31.101&127&28 \cr
 & &16.7.11.169.67&1001&536& & &64.7.101.127&889&3232 \cr
\noalign{\hrule}
 & &125.7.13.19.43&18689&18936& & &13.37.113.179&151&330 \cr
1134&9293375&16.9.11.263.1699&597&2296&1152&9729187&4.3.5.11.113.151&207&358 \cr
 & &256.27.7.41.199&5373&5248& & &16.27.11.23.179&297&184 \cr
\noalign{\hrule}
}%
}
$$
\eject
\vglue -23 pt
\noindent\hskip 1 in\hbox to 6.5 in{\ 1153 -- 1188 \hfill\fbd 9731799 -- 10324465\frb}
\vskip -9 pt
$$
\vbox{
\nointerlineskip
\halign{\strut
    \vrule \ \ \hfil \frb #\ 
   &\vrule \hfil \ \ \fbb #\frb\ 
   &\vrule \hfil \ \ \frb #\ \hfil
   &\vrule \hfil \ \ \frb #\ 
   &\vrule \hfil \ \ \frb #\ \ \vrule \hskip 2 pt
   &\vrule \ \ \hfil \frb #\ 
   &\vrule \hfil \ \ \fbb #\frb\ 
   &\vrule \hfil \ \ \frb #\ \hfil
   &\vrule \hfil \ \ \frb #\ 
   &\vrule \hfil \ \ \frb #\ \vrule \cr%
\noalign{\hrule}
 & &27.7.11.31.151&25&52& & &3.7.11.19.29.79&2599&4100 \cr
1153&9731799&8.25.13.31.151&153&2&1171&10055199&8.25.23.41.113&189&754 \cr
 & &32.9.5.13.17&17&1040& & &32.27.5.7.13.29&65&144 \cr
\noalign{\hrule}
 & &7.11.19.59.113&599&522& & &81.5.7.11.17.19&289&278 \cr
1154&9753821&4.9.29.113.599&2537&740&1172&10072755&4.5.4913.19.139&13603&10962 \cr
 & &32.3.5.37.43.59&555&688& & &16.27.7.29.61.223&1769&1784 \cr
\noalign{\hrule}
 & &11.23.29.31.43&507&826& & &9.5.11.19.29.37&127&82 \cr
1155&9780221&4.3.7.169.23.59&165&4&1173&10091565&4.29.37.41.127&1083&2600 \cr
 & &32.9.5.11.59&45&944& & &64.3.25.13.361&247&160 \cr
\noalign{\hrule}
 & &9.5.43.61.83&3179&556& & &3.25.11.13.23.41&31&174 \cr
1156&9796905&8.11.289.139&61&78&1174&10113675&4.9.5.23.29.31&31&176 \cr
 & &32.3.11.13.17.61&143&272& & &128.11.961&961&64 \cr
\noalign{\hrule}
 & &5.11.19.83.113&207&208& & &3.7.11.17.29.89&317&306 \cr
1157&9801055&32.9.11.13.19.23.113&251&992&1175&10135587&4.27.289.29.317&89&8470 \cr
 & &2048.3.23.31.251&23343&23552& & &16.5.7.121.89&11&40 \cr
\noalign{\hrule}
 & &3.5.11.19.53.59&339&244& & &3.5.7.13.17.19.23&9&124 \cr
1158&9803145&8.9.59.61.113&209&322&1176&10140585&8.27.13.17.31&529&308 \cr
 & &32.7.11.19.23.61&427&368& & &64.7.11.529&253&32 \cr
\noalign{\hrule}
 & &9.125.11.13.61&5&8& & &3.125.11.23.107&63&52 \cr
1159&9813375&16.3.625.11.61&1273&602&1177&10151625&8.27.25.7.13.107&37&712 \cr
 & &64.7.19.43.67&2881&4256& & &128.13.37.89&3293&832 \cr
\noalign{\hrule}
 & &27.5.7.13.17.47&1919&3674& & &5.7.11.13.19.107&313&222 \cr
1160&9815715&4.11.19.101.167&21&188&1178&10175165&4.3.11.19.37.313&261&52 \cr
 & &32.3.7.47.101&101&16& & &32.27.13.29.37&999&464 \cr
\noalign{\hrule}
 & &13.29.131.199&165&34& & &25.13.23.29.47&99&476 \cr
1161&9828013&4.3.5.11.13.17.29&393&712&1179&10188425&8.9.7.11.17.47&445&116 \cr
 & &64.9.89.131&89&288& & &64.3.5.29.89&267&32 \cr
\noalign{\hrule}
 & &3.5.11.19.43.73&509&294& & &3.5.23.109.271&1389&1118 \cr
1162&9840765&4.9.49.19.509&1825&2756&1180&10190955&4.9.5.13.43.463&1199&736 \cr
 & &32.25.13.53.73&265&208& & &256.11.13.23.109&143&128 \cr
\noalign{\hrule}
 & &3.25.7.11.17.101&59&766& & &31.41.71.113&2301&2332 \cr
1163&9915675&4.17.59.383&693&310&1181&10197233&8.3.11.13.53.59.71&1695&5884 \cr
 & &16.9.5.7.11.31&93&8& & &64.9.5.113.1471&1471&1440 \cr
\noalign{\hrule}
 & &27.5.17.61.71&493&422& & &25.11.17.37.59&7&18 \cr
1164&9939645&4.9.289.29.211&275&14&1182&10205525&4.9.7.17.37.59&1375&1634 \cr
 & &16.25.7.11.211&1477&440& & &16.3.125.11.19.43&285&344 \cr
\noalign{\hrule}
 & &3.5.7.13.37.197&209&246& & &11.17.83.659&5145&6058 \cr
1165&9949485&4.9.11.19.41.197&13&184&1183&10228339&4.3.5.343.13.233&901&264 \cr
 & &64.11.13.23.41&943&352& & &64.9.7.11.17.53&371&288 \cr
\noalign{\hrule}
 & &17.29.89.227&1177&1404& & &11.37.41.613&2613&4130 \cr
1166&9960079&8.27.11.13.17.107&23&130&1184&10229131&4.3.5.7.13.59.67&205&264 \cr
 & &32.3.5.11.169.23&2535&4048& & &64.9.25.11.13.41&325&288 \cr
\noalign{\hrule}
 & &5.17.233.503&209&294& & &9.5.11.17.23.53&4181&4234 \cr
1167&9961915&4.3.49.11.19.233&155&78&1185&10257885&4.23.29.37.73.113&51&2650 \cr
 & &16.9.5.7.13.19.31&1729&2232& & &16.3.25.17.29.53&29&40 \cr
\noalign{\hrule}
 & &9.47.137.173&847&710& & &7.11.31.59.73&377&36 \cr
1168&10025523&4.5.7.121.47.71&37&84&1186&10280809&8.9.13.29.73&95&124 \cr
 & &32.3.5.49.37.71&3479&2960& & &64.3.5.13.19.31&195&608 \cr
\noalign{\hrule}
 & &5.7.29.41.241&4433&5448& & &27.7.19.47.61&377&1270 \cr
1169&10029215&16.3.11.13.31.227&1305&1646&1187&10295397&4.5.7.13.29.127&517&372 \cr
 & &64.27.5.29.823&823&864& & &32.3.11.13.31.47&341&208 \cr
\noalign{\hrule}
 & &5.11.169.23.47&681&164& & &5.71.127.229&99&28 \cr
1170&10047895&8.3.23.41.227&585&358&1188&10324465&8.9.5.7.11.229&87&142 \cr
 & &32.27.5.13.179&179&432& & &32.27.7.29.71&203&432 \cr
\noalign{\hrule}
}%
}
$$
\eject
\vglue -23 pt
\noindent\hskip 1 in\hbox to 6.5 in{\ 1189 -- 1224 \hfill\fbd 10345335 -- 10942425\frb}
\vskip -9 pt
$$
\vbox{
\nointerlineskip
\halign{\strut
    \vrule \ \ \hfil \frb #\ 
   &\vrule \hfil \ \ \fbb #\frb\ 
   &\vrule \hfil \ \ \frb #\ \hfil
   &\vrule \hfil \ \ \frb #\ 
   &\vrule \hfil \ \ \frb #\ \ \vrule \hskip 2 pt
   &\vrule \ \ \hfil \frb #\ 
   &\vrule \hfil \ \ \fbb #\frb\ 
   &\vrule \hfil \ \ \frb #\ \hfil
   &\vrule \hfil \ \ \frb #\ 
   &\vrule \hfil \ \ \frb #\ \vrule \cr%
\noalign{\hrule}
 & &3.5.7.11.169.53&271&236& & &3.49.13.43.131&1501&3608 \cr
1189&10345335&8.11.53.59.271&351&298&1207&10764663&16.11.19.41.79&455&414 \cr
 & &32.27.13.149.271&2439&2384& & &64.9.5.7.13.19.23&437&480 \cr
\noalign{\hrule}
 & &27.13.17.37.47&253&206& & &3.5.7.37.47.59&17&312 \cr
1190&10376613&4.11.13.23.37.103&1461&3830&1208&10773105&16.9.13.17.37&2173&2156 \cr
 & &16.3.5.383.487&2435&3064& & &128.49.11.41.53&3157&3392 \cr
\noalign{\hrule}
 & &9.11.17.37.167&1475&1364& & &3.5.49.11.31.43&347&304 \cr
1191&10399257&8.3.25.121.31.59&329&34&1209&10777305&32.5.7.11.19.347&9&86 \cr
 & &32.5.7.17.31.47&1085&752& & &128.9.43.347&347&192 \cr
\noalign{\hrule}
 & &25.11.17.23.97&43&54& & &5.31.79.881&363&518 \cr
1192&10429925&4.27.25.17.23.43&83&308&1210&10787845&4.3.7.121.37.79&21&100 \cr
 & &32.3.7.11.43.83&903&1328& & &32.9.25.49.37&1665&784 \cr
\noalign{\hrule}
 & &27.121.31.103&37&3230& & &27.5.41.1951&8759&8800 \cr
1193&10431531&4.5.17.19.37&143&180&1211&10798785&64.3.125.11.19.461&1379&4 \cr
 & &32.9.25.11.13&325&16& & &512.7.19.197&3743&1792 \cr
\noalign{\hrule}
 & &27.5.11.73.97&14587&14222& & &729.7.29.73&1493&3610 \cr
1194&10515285&4.13.29.503.547&525&22&1212&10803051&4.5.361.1493&737&756 \cr
 & &16.3.25.7.11.13.29&377&280& & &32.27.5.7.11.19.67&1045&1072 \cr
\noalign{\hrule}
 & &11.13.263.281&3255&164& & &3.5.7.121.23.37&219&34 \cr
1195&10568129&8.3.5.7.31.41&221&66&1213&10811955&4.9.7.11.17.73&37&26 \cr
 & &32.9.11.13.17&9&272& & &16.13.17.37.73&221&584 \cr
\noalign{\hrule}
 & &27.5.29.37.73&1337&1364& & &9.25.13.47.79&77&2 \cr
1196&10574415&8.5.7.11.29.31.191&927&28&1214&10860525&4.3.7.11.13.47&295&316 \cr
 & &64.9.49.11.103&1133&1568& & &32.5.11.59.79&59&176 \cr
\noalign{\hrule}
 & &7.11.169.19.43&243&230& & &243.13.19.181&1679&1760 \cr
1197&10631621&4.243.5.7.13.19.23&407&160&1215&10863801&64.3.5.11.13.23.73&181&38 \cr
 & &256.3.25.11.23.37&2775&2944& & &256.5.19.23.181&115&128 \cr
\noalign{\hrule}
 & &9.25.11.31.139&703&548& & &3.11.127.2593&1233&1360 \cr
1198&10664775&8.5.11.19.37.137&139&546&1216&10867263&32.27.5.11.17.137&1313&1016 \cr
 & &32.3.7.13.19.139&133&208& & &512.5.13.101.127&1313&1280 \cr
\noalign{\hrule}
 & &27.5.7.11.13.79&1349&796& & &9.289.47.89&6391&7192 \cr
1199&10675665&8.9.19.71.199&419&220&1217&10879983&16.7.11.29.31.83&2403&170 \cr
 & &64.5.11.19.419&419&608& & &64.27.5.17.89&15&32 \cr
\noalign{\hrule}
 & &3.5.7.11.13.23.31&205&136& & &27.7.239.241&215&26 \cr
1200&10705695&16.25.7.13.17.41&23&198&1218&10886211&4.5.13.43.239&1661&1446 \cr
 & &64.9.11.23.41&123&32& & &16.3.11.151.241&151&88 \cr
\noalign{\hrule}
 & &9.7.11.13.29.41&227&150& & &27.5.49.17.97&311&1144 \cr
1201&10711701&4.27.25.41.227&47&88&1219&10908135&16.9.11.13.311&97&214 \cr
 & &64.5.11.47.227&1135&1504& & &64.11.97.107&107&352 \cr
\noalign{\hrule}
 & &81.121.1093&5447&4354& & &81.7.11.17.103&377&190 \cr
1202&10712493&4.7.13.311.419&365&54&1220&10920987&4.5.13.19.29.103&1429&1326 \cr
 & &16.27.5.7.13.73&365&728& & &16.3.169.17.1429&1429&1352 \cr
\noalign{\hrule}
 & &9.49.101.241&575&334& & &11.29.97.353&535&3348 \cr
1203&10734381&4.25.49.23.167&2651&1524&1221&10922879&8.27.5.31.107&145&176 \cr
 & &32.3.11.127.241&127&176& & &256.9.25.11.29&225&128 \cr
\noalign{\hrule}
 & &3.5.121.61.97&93&578& & &5.7.11.13.37.59&1395&1454 \cr
1204&10739355&4.9.11.289.31&95&194&1222&10925915&4.9.25.13.31.727&851&124 \cr
 & &16.5.19.31.97&31&152& & &32.3.23.961.37&961&1104 \cr
\noalign{\hrule}
 & &9.5.23.97.107&277&208& & &5.13.31.61.89&901&990 \cr
1205&10742265&32.3.13.107.277&299&22&1223&10939435&4.9.25.11.13.17.53&427&2 \cr
 & &128.11.169.23&169&704& & &16.3.7.53.61&21&424 \cr
\noalign{\hrule}
 & &5.13.37.41.109&261&220& & &27.25.13.29.43&119&76 \cr
1206&10747945&8.9.25.11.29.109&17&92&1224&10942425&8.9.5.7.17.19.29&143&172 \cr
 & &64.3.11.17.23.29&4301&2784& & &64.11.13.17.19.43&323&352 \cr
\noalign{\hrule}
}%
}
$$
\eject
\vglue -23 pt
\noindent\hskip 1 in\hbox to 6.5 in{\ 1225 -- 1260 \hfill\fbd 10966839 -- 11603163\frb}
\vskip -9 pt
$$
\vbox{
\nointerlineskip
\halign{\strut
    \vrule \ \ \hfil \frb #\ 
   &\vrule \hfil \ \ \fbb #\frb\ 
   &\vrule \hfil \ \ \frb #\ \hfil
   &\vrule \hfil \ \ \frb #\ 
   &\vrule \hfil \ \ \frb #\ \ \vrule \hskip 2 pt
   &\vrule \ \ \hfil \frb #\ 
   &\vrule \hfil \ \ \fbb #\frb\ 
   &\vrule \hfil \ \ \frb #\ \hfil
   &\vrule \hfil \ \ \frb #\ 
   &\vrule \hfil \ \ \frb #\ \vrule \cr%
\noalign{\hrule}
 & &3.13.31.47.193&825&632& & &27.7.19.29.109&31&140 \cr
1225&10966839&16.9.25.11.13.79&193&518&1243&11351151&8.3.5.49.29.31&59&88 \cr
 & &64.7.11.37.193&259&352& & &128.5.11.31.59&1705&3776 \cr
\noalign{\hrule}
 & &3.7.13.19.29.73&275&236& & &3.19.41.43.113&825&938 \cr
1226&10980879&8.25.11.19.29.59&333&1378&1244&11355483&4.9.25.7.11.19.67&43&52 \cr
 & &32.9.5.13.37.53&795&592& & &32.5.7.11.13.43.67&2345&2288 \cr
\noalign{\hrule}
 & &121.29.31.101&1525&1404& & &9.7.23.47.167&187&20 \cr
1227&10986679&8.27.25.13.31.61&629&164&1245&11373201&8.5.7.11.17.47&167&162 \cr
 & &64.9.5.17.37.41&5661&6560& & &32.81.11.17.167&153&176 \cr
\noalign{\hrule}
 & &5.23.149.643&3321&106& & &25.121.53.71&117&238 \cr
1228&11017805&4.81.41.53&1073&1100&1246&11383075&4.9.5.7.13.17.53&9&44 \cr
 & &32.3.25.11.29.37&2035&1392& & &32.81.11.13.17&1053&272 \cr
\noalign{\hrule}
 & &9.25.11.61.73&247&28& & &3.7.11.149.331&687&356 \cr
1229&11021175&8.3.7.13.19.61&185&242&1247&11392689&8.9.11.89.229&859&1660 \cr
 & &32.5.121.13.37&481&176& & &64.5.83.859&4295&2656 \cr
\noalign{\hrule}
 & &3.47.79.991&1463&1510& & &5.11.37.41.137&139&546 \cr
1230&11038749&4.5.7.11.19.79.151&963&94&1248&11430595&4.3.7.13.41.139&213&74 \cr
 & &16.9.5.19.47.107&535&456& & &16.9.13.37.71&117&568 \cr
\noalign{\hrule}
 & &27.13.23.1369&329&292& & &27.13.289.113&275&388 \cr
1231&11051937&8.7.13.37.47.73&1045&1656&1249&11462607&8.9.25.11.17.97&31&904 \cr
 & &128.9.5.7.11.19.23&665&704& & &128.5.31.113&31&320 \cr
\noalign{\hrule}
 & &9.23.41.1307&757&550& & &9.7.11.13.19.67&505&232 \cr
1232&11092509&4.25.11.41.757&481&276&1250&11468457&16.3.5.19.29.101&43&14 \cr
 & &32.3.5.11.13.23.37&715&592& & &64.5.7.43.101&505&1376 \cr
\noalign{\hrule}
 & &5.7.11.19.37.41&169&36& & &5.13.31.41.139&85&54 \cr
1233&11096855&8.9.11.169.37&223&184&1251&11483485&4.27.25.13.17.41&139&836 \cr
 & &128.3.13.23.223&2899&4416& & &32.9.11.19.139&99&304 \cr
\noalign{\hrule}
 & &11.13.31.43.59&1197&1340& & &11.17.19.53.61&67&120 \cr
1234&11246521&8.9.5.7.19.31.67&649&1118&1252&11486849&16.3.5.19.61.67&17&78 \cr
 & &32.3.5.11.13.43.59&15&16& & &64.9.13.17.67&871&288 \cr
\noalign{\hrule}
 & &3.7.11.19.31.83&85&8& & &9.5.13.17.19.61&53&8 \cr
1235&11292897&16.5.17.19.83&217&198&1253&11526255&16.13.17.19.53&183&506 \cr
 & &64.9.7.11.17.31&51&32& & &64.3.11.23.61&253&32 \cr
\noalign{\hrule}
 & &81.11.31.409&241&650& & &3.25.121.31.41&4553&4522 \cr
1236&11296989&4.25.13.31.241&43&198&1254&11534325&4.7.17.19.29.41.157&45&1144 \cr
 & &16.9.5.11.13.43&65&344& & &64.9.5.11.13.17.19&663&608 \cr
\noalign{\hrule}
 & &3.5.13.149.389&2603&2454& & &29.41.89.109&649&540 \cr
1237&11302395&4.9.5.19.137.409&1639&406&1255&11534489&8.27.5.11.59.89&5&94 \cr
 & &16.7.11.19.29.149&551&616& & &32.3.25.47.59&3525&944 \cr
\noalign{\hrule}
 & &3.5.7.37.41.71&187&310& & &29.41.89.109&649&540 \cr
1238&11309235&4.25.11.17.31.37&2911&1764&1256&11534489&8.27.5.11.59.89&133&44 \cr
 & &32.9.49.41.71&21&16& & &64.9.5.7.121.19&5445&4256 \cr
\noalign{\hrule}
 & &9.121.13.17.47&655&434& & &81.7.19.29.37&127&386 \cr
1239&11311443&4.5.7.31.47.131&89&42&1257&11559429&4.3.29.127.193&275&304 \cr
 & &16.3.5.49.31.89&2759&1960& & &128.25.11.19.127&1397&1600 \cr
\noalign{\hrule}
 & &7.47.173.199&63&110& & &71.167.977&405&572 \cr
1240&11326483&4.9.5.49.11.199&173&26&1258&11584289&8.81.5.11.13.71&391&248 \cr
 & &16.3.5.11.13.173&65&264& & &128.9.5.17.23.31&6417&5440 \cr
\noalign{\hrule}
 & &27.25.103.163&339&176& & &27.5.7.13.23.41&19&46 \cr
1241&11332575&32.81.5.11.113&419&824&1259&11584755&4.7.19.529.41&121&408 \cr
 & &512.103.419&419&256& & &64.3.121.17.19&2057&608 \cr
\noalign{\hrule}
 & &9.5.17.37.401&209&124& & &3.11.13.17.37.43&255&218 \cr
1242&11350305&8.11.19.31.401&185&216&1260&11603163&4.9.5.13.289.109&203&86 \cr
 & &128.27.5.11.19.37&209&192& & &16.5.7.29.43.109&1015&872 \cr
\noalign{\hrule}
}%
}
$$
\eject
\vglue -23 pt
\noindent\hskip 1 in\hbox to 6.5 in{\ 1261 -- 1296 \hfill\fbd 11604285 -- 12427349\frb}
\vskip -9 pt
$$
\vbox{
\nointerlineskip
\halign{\strut
    \vrule \ \ \hfil \frb #\ 
   &\vrule \hfil \ \ \fbb #\frb\ 
   &\vrule \hfil \ \ \frb #\ \hfil
   &\vrule \hfil \ \ \frb #\ 
   &\vrule \hfil \ \ \frb #\ \ \vrule \hskip 2 pt
   &\vrule \ \ \hfil \frb #\ 
   &\vrule \hfil \ \ \fbb #\frb\ 
   &\vrule \hfil \ \ \frb #\ \hfil
   &\vrule \hfil \ \ \frb #\ 
   &\vrule \hfil \ \ \frb #\ \vrule \cr%
\noalign{\hrule}
 & &9.5.7.11.17.197&1157&222& & &5.7.11.139.223&117&106 \cr
1261&11604285&4.27.13.37.89&493&506&1279&11933845&4.9.5.7.13.53.139&211&484 \cr
 & &16.11.17.23.29.89&667&712& & &32.3.121.53.211&2321&2544 \cr
\noalign{\hrule}
 & &9.5.11.31.757&133&890& & &9.7.29.47.139&2353&4180 \cr
1262&11616165&4.3.25.7.19.89&757&668&1280&11935791&8.5.11.13.19.181&557&348 \cr
 & &32.7.167.757&167&112& & &64.3.13.29.557&557&416 \cr
\noalign{\hrule}
 & &7.17.23.31.137&743&216& & &11.23.149.317&4465&2826 \cr
1263&11624039&16.27.23.743&3355&3332&1281&11949949&4.9.5.19.47.157&149&8 \cr
 & &128.3.5.49.11.17.61&2135&2112& & &64.3.5.19.149&285&32 \cr
\noalign{\hrule}
 & &9.5.13.19.1049&895&154& & &3.11.13.17.31.53&137&84 \cr
1264&11659635&4.3.25.7.11.179&1247&722&1282&11982399&8.9.7.11.31.137&223&1730 \cr
 & &16.361.29.43&817&232& & &32.5.173.223&865&3568 \cr
\noalign{\hrule}
 & &27.13.139.239&323&3430& & &9.49.31.877&4055&3838 \cr
1265&11660571&4.5.343.17.19&341&324&1283&11989467&4.5.7.19.101.811&759&52 \cr
 & &32.81.49.11.31&1023&784& & &32.3.5.11.13.19.23&2717&1840 \cr
\noalign{\hrule}
 & &9.125.11.23.41&1&124& & &5.49.11.61.73&261&250 \cr
1266&11669625&8.3.11.23.31&1075&1064&1284&12000835&4.9.625.7.29.61&99&526 \cr
 & &128.25.7.19.43&817&448& & &16.81.11.29.263&2349&2104 \cr
\noalign{\hrule}
 & &9.5.11.17.19.73&47&38& & &3.125.17.31.61&473&442 \cr
1267&11671605&4.11.361.47.73&221&582&1285&12055125&4.25.11.13.289.43&307&18 \cr
 & &16.3.13.17.47.97&611&776& & &16.9.11.43.307&1419&2456 \cr
\noalign{\hrule}
 & &9.13.17.61.97&1025&964& & &9.5.49.13.421&1919&1870 \cr
1268&11768913&8.25.41.97.241&363&122&1286&12067965&4.25.11.13.17.19.101&47&2478 \cr
 & &32.3.5.121.41.61&605&656& & &16.3.7.19.47.59&893&472 \cr
\noalign{\hrule}
 & &9.11.139.857&769&760& & &27.13.17.19.107&1925&2248 \cr
1269&11793177&16.5.19.769.857&813&44&1287&12130911&16.9.25.7.11.281&1&10 \cr
 & &128.3.5.11.19.271&1355&1216& & &64.125.7.281&1967&4000 \cr
\noalign{\hrule}
 & &3.11.13.59.467&4009&2062& & &27.5.11.13.17.37&137&322 \cr
1270&11820237&4.19.211.1031&621&410&1288&12142845&4.7.11.13.23.137&551&408 \cr
 & &16.27.5.19.23.41&2185&2952& & &64.3.17.19.23.29&551&736 \cr
\noalign{\hrule}
 & &3.121.29.1123&75&46& & &3.5.121.37.181&91&272 \cr
1271&11821821&4.9.25.23.1123&1079&44&1289&12155055&32.5.7.13.17.37&99&86 \cr
 & &32.5.11.13.83&83&1040& & &128.9.7.11.17.43&903&1088 \cr
\noalign{\hrule}
 & &9.7.23.41.199&503&440& & &9.5.11.13.31.61&475&196 \cr
1272&11822391&16.5.11.199.503&351&152&1290&12168585&8.125.49.13.19&61&186 \cr
 & &256.27.5.11.13.19&2145&2432& & &32.3.49.31.61&49&16 \cr
\noalign{\hrule}
 & &5.7.29.43.271&15&286& & &5.11.13.19.29.31&67&78 \cr
1273&11827795&4.3.25.11.13.29&147&172&1291&12212915&4.3.169.19.31.67&379&210 \cr
 & &32.9.49.13.43&63&208& & &16.9.5.7.67.379&3411&3752 \cr
\noalign{\hrule}
 & &9.5.7.529.71&87&442& & &3.5.7.11.23.463&473&10 \cr
1274&11831085&4.27.7.13.17.29&275&184&1292&12299595&4.25.121.43&477&598 \cr
 & &64.25.11.23.29&145&352& & &16.9.13.23.53&39&424 \cr
\noalign{\hrule}
 & &25.7.31.37.59&1377&452& & &5.47.179.293&301&594 \cr
1275&11842775&8.81.7.17.113&29&22&1293&12325045&4.27.7.11.43.47&29&358 \cr
 & &32.27.11.29.113&3277&4752& & &16.3.11.29.179&29&264 \cr
\noalign{\hrule}
 & &3.25.7.13.37.47&209&246& & &27.13.23.29.53&187&164 \cr
1276&11868675&4.9.5.11.19.41.47&169&686&1294&12408201&8.11.17.29.41.53&45&538 \cr
 & &16.343.169.41&533&392& & &32.9.5.41.269&1345&656 \cr
\noalign{\hrule}
 & &243.343.11.13&1655&1018& & &9.7.11.19.23.41&59&40 \cr
1277&11918907&4.5.7.331.509&169&162&1295&12416481&16.5.7.23.41.59&561&1504 \cr
 & &16.81.5.169.509&509&520& & &1024.3.11.17.47&799&512 \cr
\noalign{\hrule}
 & &9.11.13.73.127&607&680& & &11.19.97.613&4293&2450 \cr
1278&11931777&16.5.17.127.607&621&14&1296&12427349&4.81.25.49.53&227&38 \cr
 & &64.27.7.17.23&1173&224& & &16.3.5.7.19.227&227&840 \cr
\noalign{\hrule}
}%
}
$$
\eject
\vglue -23 pt
\noindent\hskip 1 in\hbox to 6.5 in{\ 1297 -- 1332 \hfill\fbd 12439791 -- 13091881\frb}
\vskip -9 pt
$$
\vbox{
\nointerlineskip
\halign{\strut
    \vrule \ \ \hfil \frb #\ 
   &\vrule \hfil \ \ \fbb #\frb\ 
   &\vrule \hfil \ \ \frb #\ \hfil
   &\vrule \hfil \ \ \frb #\ 
   &\vrule \hfil \ \ \frb #\ \ \vrule \hskip 2 pt
   &\vrule \ \ \hfil \frb #\ 
   &\vrule \hfil \ \ \fbb #\frb\ 
   &\vrule \hfil \ \ \frb #\ \hfil
   &\vrule \hfil \ \ \frb #\ 
   &\vrule \hfil \ \ \frb #\ \vrule \cr%
\noalign{\hrule}
 & &27.7.13.61.83&89&28& & &9.13.23.37.127&235&616 \cr
1297&12439791&8.3.49.83.89&2167&1900&1315&12645009&16.3.5.7.11.13.47&25&116 \cr
 & &64.25.11.19.197&5225&6304& & &128.125.11.29&3625&704 \cr
\noalign{\hrule}
 & &27.13.529.67&107&94& & &9.25.11.53.97&931&136 \cr
1298&12440493&4.9.529.47.107&4895&134&1316&12723975&16.3.5.49.17.19&43&62 \cr
 & &16.5.11.67.89&89&440& & &64.7.17.31.43&1333&3808 \cr
\noalign{\hrule}
 & &3.101.181.227&141&40& & &27.5.7.97.139&451&34 \cr
1299&12449361&16.9.5.47.227&1045&998&1317&12741435&4.9.7.11.17.41&139&148 \cr
 & &64.25.11.19.499&9481&8800& & &32.11.17.37.139&407&272 \cr
\noalign{\hrule}
 & &81.25.11.13.43&203&122& & &9.7.11.53.347&41&118 \cr
1300&12451725&4.7.11.29.43.61&135&338&1318&12744963&4.3.41.59.347&689&1730 \cr
 & &16.27.5.169.61&61&104& & &16.5.13.53.173&173&520 \cr
\noalign{\hrule}
 & &9.5.31.79.113&11303&10738& & &3.7.11.13.31.137&115&102 \cr
1301&12453165&4.7.13.59.89.127&1023&134&1319&12753741&4.9.5.11.17.23.137&2303&26 \cr
 & &16.3.11.31.59.67&649&536& & &16.5.49.13.47&35&376 \cr
\noalign{\hrule}
 & &9.5.13.61.349&3553&412& & &3.5.7.13.17.19.29&109&486 \cr
1302&12454065&8.11.17.19.103&61&42&1320&12785955&4.729.19.109&355&374 \cr
 & &32.3.7.11.17.61&77&272& & &16.5.11.17.71.109&781&872 \cr
\noalign{\hrule}
 & &3.5.49.11.23.67&97&104& & &9.7.121.23.73&95&752 \cr
1303&12458985&16.5.7.11.13.23.97&1387&3618&1321&12799017&32.5.19.23.47&11&12 \cr
 & &64.27.19.67.73&657&608& & &256.3.5.11.19.47&893&640 \cr
\noalign{\hrule}
 & &9.5.7.11.59.61&299&3544& & &3.11.23.101.167&43&210 \cr
1304&12470535&16.13.23.443&233&210&1322&12802053&4.9.5.7.43.101&55&46 \cr
 & &64.3.5.7.13.233&233&416& & &16.25.7.11.23.43&301&200 \cr
\noalign{\hrule}
 & &3.11.13.17.29.59&49&270& & &9.11.31.37.113&455&568 \cr
1305&12478323&4.81.5.49.59&247&166&1323&12831489&16.3.5.7.13.37.71&155&118 \cr
 & &16.5.7.13.19.83&415&1064& & &64.25.31.59.71&1775&1888 \cr
\noalign{\hrule}
 & &13.151.6361&3187&3174& & &9.13.19.53.109&667&340 \cr
1306&12486643&4.3.529.151.3187&143&3330&1324&12842271&8.3.5.13.17.23.29&1643&1672 \cr
 & &16.27.5.11.13.23.37&3105&3256& & &128.11.19.23.31.53&713&704 \cr
\noalign{\hrule}
 & &9.11.13.17.571&4553&2870& & &3.11.361.23.47&477&40 \cr
1307&12492909&4.5.7.29.41.157&17&24&1325&12877953&16.27.5.19.53&283&230 \cr
 & &64.3.5.17.29.157&785&928& & &64.25.23.283&283&800 \cr
\noalign{\hrule}
 & &5.7.11.17.23.83&97&90& & &9.11.19.41.167&443&3616 \cr
1308&12494405&4.9.25.23.83.97&2167&258&1326&12879207&64.113.443&165&278 \cr
 & &16.27.11.43.197&1161&1576& & &256.3.5.11.139&139&640 \cr
\noalign{\hrule}
 & &25.13.137.281&4403&2622& & &25.121.17.251&1387&1638 \cr
1309&12511525&4.3.7.17.19.23.37&325&66&1327&12907675&4.9.7.13.17.19.73&177&44 \cr
 & &16.9.25.11.13.19&209&72& & &32.27.11.59.73&1593&1168 \cr
\noalign{\hrule}
 & &3.25.7.11.41.53&23&182& & &3.5.13.131.509&33&98 \cr
1310&12549075&4.5.49.11.13.23&333&382&1328&13002405&4.9.49.11.509&475&34 \cr
 & &16.9.23.37.191&2553&1528& & &16.25.11.17.19&1615&88 \cr
\noalign{\hrule}
 & &9.11.23.37.149&169&238& & &9.5.7.19.41.53&1419&754 \cr
1311&12553101&4.3.7.169.17.149&55&94&1329&13005405&4.27.11.13.29.43&7&20 \cr
 & &16.5.7.11.13.17.47&1547&1880& & &32.5.7.11.29.43&319&688 \cr
\noalign{\hrule}
 & &3.7.11.13.53.79&75&68& & &5.7.169.31.71&531&314 \cr
1312&12573561&8.9.25.17.53.79&41&436&1330&13018915&4.9.59.71.157&1859&2330 \cr
 & &64.5.17.41.109&4469&2720& & &16.3.5.11.169.233&233&264 \cr
\noalign{\hrule}
 & &11.19.23.37.71&91&162& & &7.19.29.31.109&15&14 \cr
1313&12627989&4.81.7.13.19.37&365&1364&1331&13032803&4.3.5.49.19.31.109&4017&638 \cr
 & &32.3.5.11.31.73&465&1168& & &16.9.11.13.29.103&927&1144 \cr
\noalign{\hrule}
 & &9.13.23.37.127&1251&400& & &11.61.109.179&5459&5460 \cr
1314&12645009&32.81.25.139&29&110&1332&13091881&8.3.5.7.11.13.53.103.109&1077&6854 \cr
 & &128.125.11.29&3625&704& & &32.9.5.13.23.149.359&107341&107280 \cr
\noalign{\hrule}
}%
}
$$
\eject
\vglue -23 pt
\noindent\hskip 1 in\hbox to 6.5 in{\ 1333 -- 1368 \hfill\fbd 13097175 -- 13771615\frb}
\vskip -9 pt
$$
\vbox{
\nointerlineskip
\halign{\strut
    \vrule \ \ \hfil \frb #\ 
   &\vrule \hfil \ \ \fbb #\frb\ 
   &\vrule \hfil \ \ \frb #\ \hfil
   &\vrule \hfil \ \ \frb #\ 
   &\vrule \hfil \ \ \frb #\ \ \vrule \hskip 2 pt
   &\vrule \ \ \hfil \frb #\ 
   &\vrule \hfil \ \ \fbb #\frb\ 
   &\vrule \hfil \ \ \frb #\ \hfil
   &\vrule \hfil \ \ \frb #\ 
   &\vrule \hfil \ \ \frb #\ \vrule \cr%
\noalign{\hrule}
 & &3.25.7.13.19.101&83&8& & &5.11.13.19.23.43&45&254 \cr
1333&13097175&16.19.83.101&1001&918&1351&13435565&4.9.25.43.127&347&728 \cr
 & &64.27.7.11.13.17&187&288& & &64.3.7.13.347&1041&224 \cr
\noalign{\hrule}
 & &25.11.19.23.109&837&362& & &5.11.13.19.23.43&351&466 \cr
1334&13099075&4.27.23.31.181&19&50&1352&13435565&4.27.11.169.233&119&1978 \cr
 & &16.9.25.19.181&181&72& & &16.3.7.17.23.43&119&24 \cr
\noalign{\hrule}
 & &9.5.121.29.83&103&158& & &9.125.17.19.37&711&1414 \cr
1335&13106115&4.11.79.83.103&405&508&1353&13444875&4.81.7.79.101&37&44 \cr
 & &32.81.5.79.127&1143&1264& & &32.11.37.79.101&1111&1264 \cr
\noalign{\hrule}
 & &81.5.11.13.227&1379&664& & &5.17.29.43.127&429&302 \cr
1336&13146705&16.9.7.83.197&275&472&1354&13461365&4.3.5.11.13.29.151&2451&1696 \cr
 & &256.25.7.11.59&413&640& & &256.9.19.43.53&1007&1152 \cr
\noalign{\hrule}
 & &3.7.11.19.31.97&169&510& & &3.29.43.59.61&95&34 \cr
1337&13197723&4.9.5.169.17.19&83&88&1355&13463859&4.5.17.19.29.59&99&394 \cr
 & &64.11.169.17.83&2873&2656& & &16.9.11.19.197&2167&456 \cr
\noalign{\hrule}
 & &5.49.11.169.29&207&38& & &81.5.11.13.233&19&214 \cr
1338&13208195&4.9.11.19.23.29&169&150&1356&13494195&4.27.11.19.107&1165&868 \cr
 & &16.27.25.169.23&135&184& & &32.5.7.31.233&31&112 \cr
\noalign{\hrule}
 & &3.5.23.109.353&1077&1430& & &5.11.13.79.239&269&126 \cr
1339&13274565&4.9.25.11.13.359&1417&1058&1357&13499915&4.9.7.239.269&583&1300 \cr
 & &16.169.529.109&169&184& & &32.3.25.11.13.53&159&80 \cr
\noalign{\hrule}
 & &49.11.71.347&169&3648& & &3.5.7.19.67.101&4697&4898 \cr
1340&13279343&128.3.169.19&85&84&1358&13500165&4.49.11.31.61.79&2881&990 \cr
 & &1024.9.5.7.17.19&2907&2560& & &16.9.5.121.43.67&363&344 \cr
\noalign{\hrule}
 & &9.25.7.11.13.59&251&134& & &3.25.13.17.19.43&131&116 \cr
1341&13288275&4.5.59.67.251&273&22&1359&13541775&8.5.17.29.43.131&693&38 \cr
 & &16.3.7.11.13.67&67&8& & &32.9.7.11.19.29&203&528 \cr
\noalign{\hrule}
 & &3.11.17.19.29.43&641&90& & &9.5.17.37.479&401&364 \cr
1342&13291773&4.27.5.11.641&469&172&1360&13558095&8.7.13.401.479&39&440 \cr
 & &32.5.7.43.67&67&560& & &128.3.5.7.11.169&1183&704 \cr
\noalign{\hrule}
 & &3.125.7.13.17.23&333&242& & &3.343.47.281&1573&730 \cr
1343&13342875&4.27.5.121.17.37&161&26&1361&13590003&4.5.7.121.13.73&47&558 \cr
 & &16.7.11.13.23.37&37&88& & &16.9.13.31.47&31&312 \cr
\noalign{\hrule}
 & &49.11.17.31.47&71&258& & &9.23.251.263&1703&4070 \cr
1344&13350491&4.3.7.31.43.71&915&418&1362&13664691&4.5.11.13.37.131&269&138 \cr
 & &16.9.5.11.19.61&1159&360& & &16.3.5.13.23.269&269&520 \cr
\noalign{\hrule}
 & &3.7.11.17.41.83&31&10& & &3.11.13.19.23.73&225&212 \cr
1345&13363581&4.5.11.17.31.83&711&700&1363&13685529&8.27.25.11.53.73&331&34 \cr
 & &32.9.125.7.31.79&3875&3792& & &32.5.17.53.331&5627&4240 \cr
\noalign{\hrule}
 & &5.7.19.101.199&4897&4698& & &3.13.17.127.163&355&134 \cr
1346&13365835&4.81.7.29.59.83&101&682&1364&13724763&4.5.67.71.127&99&28 \cr
 & &16.3.11.31.59.101&649&744& & &32.9.5.7.11.67&2211&560 \cr
\noalign{\hrule}
 & &5.11.17.41.349&647&1098& & &9.25.7.11.13.61&287&262 \cr
1347&13378915&4.9.17.61.647&349&298&1365&13738725&4.49.11.13.41.131&3&536 \cr
 & &16.3.61.149.349&447&488& & &64.3.67.131&131&2144 \cr
\noalign{\hrule}
 & &5.7.11.19.31.59&129&284& & &9.11.13.289.37&203&86 \cr
1348&13379135&8.3.11.19.43.71&83&126&1366&13761891&4.7.11.29.37.43&255&218 \cr
 & &32.27.7.71.83&1917&1328& & &16.3.5.7.17.29.109&1015&872 \cr
\noalign{\hrule}
 & &9.19.59.1327&1981&2000& & &9.5.7.11.29.137&23&122 \cr
1349&13388103&32.3.125.7.59.283&1327&88&1367&13766445&4.7.23.61.137&145&282 \cr
 & &512.25.11.1327&275&256& & &16.3.5.23.29.47&47&184 \cr
\noalign{\hrule}
 & &25.13.19.41.53&341&666& & &5.121.13.17.103&9&112 \cr
1350&13418275&4.9.11.31.37.41&43&80&1368&13771615&32.9.5.7.13.17&103&118 \cr
 & &128.3.5.11.31.43&1333&2112& & &128.3.7.59.103&413&192 \cr
\noalign{\hrule}
}%
}
$$
\eject
\vglue -23 pt
\noindent\hskip 1 in\hbox to 6.5 in{\ 1369 -- 1404 \hfill\fbd 13776425 -- 14835483\frb}
\vskip -9 pt
$$
\vbox{
\nointerlineskip
\halign{\strut
    \vrule \ \ \hfil \frb #\ 
   &\vrule \hfil \ \ \fbb #\frb\ 
   &\vrule \hfil \ \ \frb #\ \hfil
   &\vrule \hfil \ \ \frb #\ 
   &\vrule \hfil \ \ \frb #\ \ \vrule \hskip 2 pt
   &\vrule \ \ \hfil \frb #\ 
   &\vrule \hfil \ \ \fbb #\frb\ 
   &\vrule \hfil \ \ \frb #\ \hfil
   &\vrule \hfil \ \ \frb #\ 
   &\vrule \hfil \ \ \frb #\ \vrule \cr%
\noalign{\hrule}
 & &25.13.19.23.97&111&136& & &25.49.13.29.31&727&792 \cr
1369&13776425&16.3.17.23.37.97&801&1430&1387&14316575&16.9.5.11.29.727&581&146 \cr
 & &64.27.5.11.13.89&979&864& & &64.3.7.11.73.83&2739&2336 \cr
\noalign{\hrule}
 & &7.13.47.53.61&18821&18450& & &9.11.361.401&131&230 \cr
1370&13827541&4.9.25.11.29.41.59&13&42&1388&14331339&4.5.23.131.401&7733&7332 \cr
 & &16.27.5.7.13.41.59&1593&1640& & &32.3.11.13.19.37.47&611&592 \cr
\noalign{\hrule}
 & &27.17.97.311&385&74& & &9.169.17.557&725&946 \cr
1371&13846653&4.5.7.11.37.97&281&204&1389&14402349&4.3.25.11.13.29.43&1513&2228 \cr
 & &32.3.17.37.281&281&592& & &32.5.17.89.557&89&80 \cr
\noalign{\hrule}
 & &9.5.11.67.419&2405&2204& & &27.7.23.31.107&325&388 \cr
1372&13896135&8.3.25.13.19.29.37&419&506&1390&14418999&8.3.25.13.97.107&841&550 \cr
 & &32.11.13.19.23.419&299&304& & &32.625.11.841&9251&10000 \cr
\noalign{\hrule}
 & &9.5.169.31.59&121&56& & &49.11.17.19.83&1437&26 \cr
1373&13909545&16.3.7.121.13.31&25&118&1391&14450051&4.3.7.13.479&285&194 \cr
 & &64.25.7.11.59&55&224& & &16.9.5.19.97&873&40 \cr
\noalign{\hrule}
 & &27.49.53.199&121&68& & &7.121.169.101&551&450 \cr
1374&13953681&8.7.121.17.199&159&40&1392&14457443&4.9.25.11.13.19.29&7&202 \cr
 & &128.3.5.121.53&121&320& & &16.3.5.7.29.101&145&24 \cr
\noalign{\hrule}
 & &3.49.11.13.23.29&171&148& & &3.25.11.13.19.71&233&508 \cr
1375&14021007&8.27.49.13.19.37&575&62&1393&14468025&8.71.127.233&99&28 \cr
 & &32.25.23.31.37&925&496& & &64.9.7.11.233&699&224 \cr
\noalign{\hrule}
 & &11.563.2267&1963&4230& & &9.5.11.23.31.41&1075&1064 \cr
1376&14039531&4.9.5.13.47.151&1139&1126&1394&14470335&16.3.125.7.19.41.43&1&124 \cr
 & &16.3.17.47.67.563&1139&1128& & &128.7.19.31.43&817&448 \cr
\noalign{\hrule}
 & &3.97.211.229&57&154& & &27.13.19.41.53&7&20 \cr
1377&14060829&4.9.7.11.19.229&1055&826&1395&14491737&8.5.7.19.41.53&1419&754 \cr
 & &16.5.49.59.211&295&392& & &32.3.11.13.29.43&319&688 \cr
\noalign{\hrule}
 & &5.7.11.13.29.97&1271&1368& & &3.13.29.41.313&423&110 \cr
1378&14079065&16.9.5.11.19.31.41&113&1158&1396&14514123&4.27.5.11.29.47&133&650 \cr
 & &64.27.113.193&5211&3616& & &16.125.7.13.19&125&1064 \cr
\noalign{\hrule}
 & &3.25.31.73.83&8701&9476& & &3.7.29.113.211&795&682 \cr
1379&14087175&8.7.11.23.103.113&13&90&1397&14520387&4.9.5.11.29.31.53&161&422 \cr
 & &32.9.5.13.23.113&1469&1104& & &16.5.7.23.31.211&155&184 \cr
\noalign{\hrule}
 & &27.5.11.13.17.43&7&466& & &3.5.7.11.23.547&29&40 \cr
1380&14111955&4.5.7.13.233&227&228&1398&14531055&16.25.7.29.547&2277&1552 \cr
 & &32.3.19.227.233&4313&3728& & &512.9.11.23.97&291&256 \cr
\noalign{\hrule}
 & &9.5.31.67.151&307&28& & &27.5.7.13.29.41&47&88 \cr
1381&14113215&8.7.151.307&375&682&1399&14606865&16.7.11.13.29.47&227&150 \cr
 & &32.3.125.11.31&25&176& & &64.3.25.47.227&1135&1504 \cr
\noalign{\hrule}
 & &9.5.13.19.31.41&721&550& & &3.5.49.11.13.139&4013&5542 \cr
1382&14127165&4.125.7.11.13.103&1379&246&1400&14609595&4.17.163.4013&1925&2088 \cr
 & &16.3.49.41.197&197&392& & &64.9.25.7.11.17.29&493&480 \cr
\noalign{\hrule}
 & &11.23.29.41.47&273&244& & &19.41.89.211&9779&9000 \cr
1383&14138399&8.3.7.13.23.41.61&435&968&1401&14628841&16.9.125.7.11.127&279&356 \cr
 & &128.9.5.7.121.29&495&448& & &128.81.25.31.89&2025&1984 \cr
\noalign{\hrule}
 & &3.5.11.13.17.389&67&322& & &27.31.89.197&1675&1084 \cr
1384&14184885&4.7.11.13.23.67&163&306&1402&14675121&8.9.25.67.271&979&376 \cr
 & &16.9.17.23.163&163&552& & &128.5.11.47.89&517&320 \cr
\noalign{\hrule}
 & &27.25.7.31.97&533&242& & &25.7.11.13.19.31&131&534 \cr
1385&14208075&4.9.7.121.13.41&97&20&1403&14739725&4.3.5.11.89.131&217&228 \cr
 & &32.5.11.41.97&41&176& & &32.9.7.19.31.131&131&144 \cr
\noalign{\hrule}
 & &3.5.11.361.239&403&642& & &9.13.23.37.149&55&94 \cr
1386&14236035&4.9.13.19.31.107&187&776&1404&14835483&4.3.5.11.23.37.47&169&238 \cr
 & &64.11.13.17.97&1261&544& & &16.5.7.169.17.47&1547&1880 \cr
\noalign{\hrule}
}%
}
$$
\eject
\vglue -23 pt
\noindent\hskip 1 in\hbox to 6.5 in{\ 1405 -- 1440 \hfill\fbd 14835645 -- 15406435\frb}
\vskip -9 pt
$$
\vbox{
\nointerlineskip
\halign{\strut
    \vrule \ \ \hfil \frb #\ 
   &\vrule \hfil \ \ \fbb #\frb\ 
   &\vrule \hfil \ \ \frb #\ \hfil
   &\vrule \hfil \ \ \frb #\ 
   &\vrule \hfil \ \ \frb #\ \ \vrule \hskip 2 pt
   &\vrule \ \ \hfil \frb #\ 
   &\vrule \hfil \ \ \fbb #\frb\ 
   &\vrule \hfil \ \ \frb #\ \hfil
   &\vrule \hfil \ \ \frb #\ 
   &\vrule \hfil \ \ \frb #\ \vrule \cr%
\noalign{\hrule}
 & &9.5.11.17.41.43&167&38& & &3.5.13.19.61.67&1337&1276 \cr
1405&14835645&4.3.11.17.19.167&157&344&1423&15142335&8.5.7.11.19.29.191&487&468 \cr
 & &64.19.43.157&157&608& & &64.9.7.11.13.29.487&10227&10208 \cr
\noalign{\hrule}
 & &27.7.17.31.149&703&3916& & &3.11.31.113.131&3913&410 \cr
1406&14840847&8.11.19.37.89&149&60&1424&15143469&4.5.7.13.41.43&297&262 \cr
 & &64.3.5.37.149&185&32& & &16.27.11.41.131&41&72 \cr
\noalign{\hrule}
 & &3.5.121.13.17.37&5&116& & &9.13.23.43.131&1045&658 \cr
1407&14841255&8.25.13.17.29&473&252&1425&15158403&4.5.7.11.19.23.47&1591&2028 \cr
 & &64.9.7.11.43&43&672& & &32.3.5.169.37.43&185&208 \cr
\noalign{\hrule}
 & &9.25.11.17.353&26999&27010& & &3.13.19.97.211&193&440 \cr
1408&14852475&4.125.49.19.29.37.73&3601&24&1426&15166047&16.5.11.97.193&437&630 \cr
 & &64.3.13.19.37.277&9139&8864& & &64.9.25.7.19.23&525&736 \cr
\noalign{\hrule}
 & &5.11.29.67.139&333&362& & &9.7.11.43.509&1445&4154 \cr
1409&14854235&4.9.11.37.67.181&2235&244&1427&15167691&4.5.289.31.67&25&42 \cr
 & &32.27.5.61.149&1647&2384& & &16.3.125.7.17.31&527&1000 \cr
\noalign{\hrule}
 & &9.11.31.37.131&115&226& & &49.11.47.599&569&30 \cr
1410&14875443&4.3.5.23.113.131&185&208&1428&15174467&4.3.5.47.569&261&308 \cr
 & &128.25.13.37.113&1469&1600& & &32.27.5.7.11.29&135&464 \cr
\noalign{\hrule}
 & &3.11.13.17.23.89&105&194& & &5.11.17.37.439&1&186 \cr
1411&14928771&4.9.5.7.11.17.97&3211&2144&1429&15187205&4.3.31.439&235&204 \cr
 & &256.169.19.67&1273&1664& & &32.9.5.17.47&423&16 \cr
\noalign{\hrule}
 & &9.5.11.19.37.43&609&1426& & &27.19.59.503&5575&3982 \cr
1412&14963355&4.27.7.23.29.31&25&2&1430&15224301&4.25.11.181.223&63&118 \cr
 & &16.25.7.29.31&145&1736& & &16.9.5.7.59.223&223&280 \cr
\noalign{\hrule}
 & &9.5.11.13.17.137&29&166& & &9.25.11.47.131&697&744 \cr
1413&14987115&4.3.11.17.29.83&923&1484&1431&15238575&16.27.25.17.31.41&931&94 \cr
 & &32.7.13.53.71&371&1136& & &64.49.17.19.47&833&608 \cr
\noalign{\hrule}
 & &13.17.29.2339&1161&1178& & &3.11.29.37.431&375&56 \cr
1414&14990651&4.27.13.19.29.31.43&1265&9356&1432&15261279&16.9.125.7.37&431&494 \cr
 & &32.3.5.11.23.2339&253&240& & &64.5.13.19.431&247&160 \cr
\noalign{\hrule}
 & &11.17.19.41.103&405&728& & &125.49.11.227&141&134 \cr
1415&15004319&16.81.5.7.13.41&43&412&1433&15294125&4.3.5.7.47.67.227&627&962 \cr
 & &128.9.43.103&387&64& & &16.9.11.13.19.37.47&5499&5624 \cr
\noalign{\hrule}
 & &3.5.17.19.29.107&429&106& & &13.29.97.419&5445&6706 \cr
1416&15034035&4.9.11.13.29.53&475&214&1434&15322411&4.9.5.7.121.479&421&58 \cr
 & &16.25.11.19.107&11&40& & &16.3.5.7.29.421&421&840 \cr
\noalign{\hrule}
 & &49.121.43.59&85&36& & &3.5.41.61.409&143&266 \cr
1417&15041873&8.9.5.17.43.59&83&212&1435&15343635&4.5.7.11.13.19.61&587&648 \cr
 & &64.3.17.53.83&4233&1696& & &64.81.7.11.587&6457&6048 \cr
\noalign{\hrule}
 & &3.13.31.37.337&187&150& & &27.5.41.47.59&1079&1694 \cr
1418&15075021&4.9.25.11.13.17.31&1961&674&1436&15348555&4.9.7.121.13.83&41&50 \cr
 & &16.5.37.53.337&53&40& & &16.25.121.41.83&605&664 \cr
\noalign{\hrule}
 & &3.7.11.361.181&219&142& & &3.49.11.13.17.43&67&80 \cr
1419&15093771&4.9.71.73.181&95&86&1437&15366351&32.5.11.17.43.67&333&806 \cr
 & &16.5.19.43.71.73&3053&2920& & &128.9.5.13.31.37&1147&960 \cr
\noalign{\hrule}
 & &25.49.11.19.59&1377&1318& & &9.11.17.41.223&217&234 \cr
1420&15105475&4.81.5.17.19.659&3773&2158&1438&15387669&4.81.7.13.31.223&415&638 \cr
 & &16.9.343.11.13.83&747&728& & &16.5.7.11.29.31.83&4495&4648 \cr
\noalign{\hrule}
 & &81.5.11.43.79&439&34& & &11.23.83.733&3077&4986 \cr
1421&15133635&4.17.79.439&259&180&1439&15392267&4.9.17.181.277&325&506 \cr
 & &32.9.5.7.17.37&119&592& & &16.3.25.11.13.17.23&325&408 \cr
\noalign{\hrule}
 & &27.121.41.113&497&610& & &5.11.19.23.641&47&162 \cr
1422&15136011&4.5.7.121.61.71&601&246&1440&15406435&4.81.47.641&391&250 \cr
 & &16.3.41.61.601&601&488& & &16.27.125.17.23&425&216 \cr
\noalign{\hrule}
}%
}
$$
\eject
\vglue -23 pt
\noindent\hskip 1 in\hbox to 6.5 in{\ 1441 -- 1476 \hfill\fbd 15418377 -- 16150855\frb}
\vskip -9 pt
$$
\vbox{
\nointerlineskip
\halign{\strut
    \vrule \ \ \hfil \frb #\ 
   &\vrule \hfil \ \ \fbb #\frb\ 
   &\vrule \hfil \ \ \frb #\ \hfil
   &\vrule \hfil \ \ \frb #\ 
   &\vrule \hfil \ \ \frb #\ \ \vrule \hskip 2 pt
   &\vrule \ \ \hfil \frb #\ 
   &\vrule \hfil \ \ \fbb #\frb\ 
   &\vrule \hfil \ \ \frb #\ \hfil
   &\vrule \hfil \ \ \frb #\ 
   &\vrule \hfil \ \ \frb #\ \vrule \cr%
\noalign{\hrule}
 & &27.169.31.109&19&20& & &9.5.17.107.193&767&4048 \cr
1441&15418377&8.9.5.13.19.31.109&11&1406&1459&15798015&32.11.13.23.59&257&510 \cr
 & &32.11.361.37&3971&592& & &128.3.5.17.257&257&64 \cr
\noalign{\hrule}
 & &3.11.47.61.163&177&340& & &13.31.197.199&4365&1804 \cr
1442&15421593&8.9.5.17.59.61&235&296&1460&15798809&8.9.5.11.41.97&69&28 \cr
 & &128.25.17.37.47&925&1088& & &64.27.5.7.11.23&4347&1760 \cr
\noalign{\hrule}
 & &27.17.151.223&341&118& & &9.7.23.67.163&275&208 \cr
1443&15455907&4.11.31.59.151&4439&4470&1461&15824529&32.3.25.11.13.163&119&44 \cr
 & &16.3.5.11.23.149.193&17135&16984& & &256.7.121.13.17&1573&2176 \cr
\noalign{\hrule}
 & &5.11.13.53.409&10845&10832& & &5.121.13.43.47&45&2 \cr
1444&15499055&32.9.25.11.241.677&17&258&1462&15895165&4.9.25.121.13&19&344 \cr
 & &128.27.17.43.677&29111&29376& & &64.3.19.43&3&608 \cr
\noalign{\hrule}
 & &27.25.7.17.193&713&638& & &9.5.53.59.113&209&322 \cr
1445&15502725&4.9.11.17.23.29.31&1589&1930&1463&15900795&4.5.7.11.19.23.53&339&244 \cr
 & &16.5.7.29.193.227&227&232& & &32.3.7.23.61.113&427&368 \cr
\noalign{\hrule}
 & &9.17.61.1663&13861&14410& & &3.7.13.19.37.83&231&250 \cr
1446&15520779&4.5.11.83.131.167&39&874&1464&15929277&4.9.125.49.11.83&13&428 \cr
 & &16.3.13.19.23.131&3013&1976& & &32.25.11.13.107&1177&400 \cr
\noalign{\hrule}
 & &5.49.19.47.71&3303&2948& & &9.11.23.47.149&25&124 \cr
1447&15533735&8.9.7.11.67.367&299&68&1465&15945831&8.25.23.31.47&81&34 \cr
 & &64.3.13.17.23.67&5083&6432& & &32.81.5.17.31&279&1360 \cr
\noalign{\hrule}
 & &3.125.343.121&207&332& & &25.19.59.569&845&276 \cr
1448&15563625&8.27.7.11.23.83&25&52&1466&15946225&8.3.125.169.23&1881&2006 \cr
 & &64.25.13.23.83&1079&736& & &32.27.11.17.19.59&297&272 \cr
\noalign{\hrule}
 & &9.25.7.11.29.31&587&312& & &9.5.11.13.37.67&13667&13602 \cr
1449&15575175&16.27.7.13.587&935&1522&1467&15952365&4.27.79.173.2267&67&2200 \cr
 & &64.5.11.17.761&761&544& & &64.25.11.67.173&173&160 \cr
\noalign{\hrule}
 & &81.11.83.211&445&446& & &9.25.17.43.97&253&38 \cr
1450&15604083&4.5.83.89.211.223&135&18644&1468&15954075&4.3.5.11.17.19.23&1067&1118 \cr
 & &32.27.25.59.79&1475&1264& & &16.121.13.43.97&121&104 \cr
\noalign{\hrule}
 & &27.5.7.13.31.41&1513&242& & &3.7.121.61.103&25&696 \cr
1451&15614235&4.7.121.17.89&333&1180&1469&15965103&16.9.25.11.29&103&158 \cr
 & &32.9.5.37.59&59&592& & &64.5.79.103&395&32 \cr
\noalign{\hrule}
 & &3.7.31.103.233&11&710& & &9.25.7.11.13.71&1539&2036 \cr
1452&15623349&4.5.11.31.71&135&206&1470&15990975&8.729.19.509&619&110 \cr
 & &16.27.25.103&225&8& & &32.5.11.19.619&619&304 \cr
\noalign{\hrule}
 & &5.7.11.13.53.59&351&298& & &3.11.19.107.239&367&260 \cr
1453&15650635&4.27.5.7.169.149&271&236&1471&16034271&8.5.13.239.367&781&414 \cr
 & &32.9.59.149.271&2439&2384& & &32.9.11.13.23.71&897&1136 \cr
\noalign{\hrule}
 & &9.5.11.103.307&1027&106& & &9.7.11.19.23.53&65&142 \cr
1454&15652395&4.3.5.13.53.79&1231&836&1472&16050573&4.5.13.19.53.71&21&74 \cr
 & &32.11.19.1231&1231&304& & &16.3.7.13.37.71&923&296 \cr
\noalign{\hrule}
 & &27.19.127.241&575&4004& & &31.1849.281&909&940 \cr
1455&15701391&8.25.7.11.13.23&63&52&1473&16106639&8.9.5.47.101.281&11&292 \cr
 & &64.9.5.49.169&845&1568& & &64.3.5.11.47.73&3431&5280 \cr
\noalign{\hrule}
 & &9.7.389.641&2071&1430& & &81.17.19.617&461&1078 \cr
1456&15708987&4.5.7.11.13.19.109&645&554&1474&16142571&4.49.11.17.461&171&290 \cr
 & &16.3.25.19.43.277&6925&6536& & &16.9.5.7.11.19.29&385&232 \cr
\noalign{\hrule}
 & &27.5.7.11.17.89&37&26& & &3.5.17.29.37.59&429&574 \cr
1457&15727635&4.3.5.13.17.37.89&721&1166&1475&16143285&4.9.7.11.13.37.41&55&314 \cr
 & &16.7.11.13.53.103&689&824& & &16.5.121.13.157&2041&968 \cr
\noalign{\hrule}
 & &5.49.19.31.109&4017&638& & &5.7.19.149.163&65&84 \cr
1458&15729245&4.3.11.13.29.103&15&14&1476&16150855&8.3.25.49.13.163&447&1672 \cr
 & &16.9.5.7.11.13.103&927&1144& & &128.9.11.19.149&99&64 \cr
\noalign{\hrule}
}%
}
$$
\eject
\vglue -23 pt
\noindent\hskip 1 in\hbox to 6.5 in{\ 1477 -- 1512 \hfill\fbd 16175475 -- 17235603\frb}
\vskip -9 pt
$$
\vbox{
\nointerlineskip
\halign{\strut
    \vrule \ \ \hfil \frb #\ 
   &\vrule \hfil \ \ \fbb #\frb\ 
   &\vrule \hfil \ \ \frb #\ \hfil
   &\vrule \hfil \ \ \frb #\ 
   &\vrule \hfil \ \ \frb #\ \ \vrule \hskip 2 pt
   &\vrule \ \ \hfil \frb #\ 
   &\vrule \hfil \ \ \fbb #\frb\ 
   &\vrule \hfil \ \ \frb #\ \hfil
   &\vrule \hfil \ \ \frb #\ 
   &\vrule \hfil \ \ \frb #\ \vrule \cr%
\noalign{\hrule}
 & &9.25.29.37.67&419&506& & &7.23.107.967&11495&10746 \cr
1477&16175475&4.3.11.23.67.419&2405&2204&1495&16658509&4.27.5.121.19.199&413&214 \cr
 & &32.5.13.19.23.29.37&299&304& & &16.9.5.7.11.59.107&495&472 \cr
\noalign{\hrule}
 & &3.7.19.23.41.43&33&10& & &243.17.37.109&113&130 \cr
1478&16179051&4.9.5.7.11.19.41&43&736&1496&16660323&4.5.13.37.109.113&2299&1734 \cr
 & &256.5.23.43&5&128& & &16.3.121.13.289.19&2057&1976 \cr
\noalign{\hrule}
 & &9.13.109.1277&8791&7810& & &243.25.13.211&649&406 \cr
1479&16285581&4.5.11.59.71.149&697&48&1497&16663725&4.5.7.11.13.29.59&161&216 \cr
 & &128.3.17.41.71&1207&2624& & &64.27.49.23.59&1357&1568 \cr
\noalign{\hrule}
 & &9.7.29.37.241&65&268& & &9.41.199.227&4753&4554 \cr
1480&16291359&8.5.13.67.241&87&154&1498&16668837&4.81.49.11.23.97&10127&2270 \cr
 & &32.3.5.7.11.13.29&65&176& & &16.5.13.19.41.227&95&104 \cr
\noalign{\hrule}
 & &5.7.11.17.47.53&71&258& & &25.7.17.71.79&209&288 \cr
1481&16303595&4.3.5.43.53.71&111&154&1499&16686775&64.9.25.11.17.19&43&518 \cr
 & &16.9.7.11.37.71&333&568& & &256.3.7.37.43&1591&384 \cr
\noalign{\hrule}
 & &9.11.13.19.23.29&655&104& & &5.49.17.19.211&3531&3854 \cr
1482&16310151&16.3.5.169.131&319&188&1500&16697485&4.3.7.11.41.47.107&633&116 \cr
 & &128.5.11.29.47&235&64& & &32.9.29.41.211&369&464 \cr
\noalign{\hrule}
 & &9.5.37.9817&14707&14744& & &27.49.19.23.29&715&734 \cr
1483&16345305&16.3.5.7.11.19.97.191&9817&6188&1501&16766379&4.3.5.7.11.13.29.367&19&184 \cr
 & &128.49.13.17.9817&833&832& & &64.13.19.23.367&367&416 \cr
\noalign{\hrule}
 & &25.13.59.853&4807&6282& & &9.25.29.31.83&77&68 \cr
1484&16356275&4.9.11.19.23.349&397&650&1502&16788825&8.5.7.11.17.31.83&999&86 \cr
 & &16.3.25.13.19.397&397&456& & &32.27.17.37.43&1591&816 \cr
\noalign{\hrule}
 & &27.5.11.73.151&73&62& & &3.5.7.11.13.19.59&85&124 \cr
1485&16369155&4.31.5329.151&5005&324&1503&16831815&8.25.7.17.31.59&589&414 \cr
 & &32.81.5.7.11.13&39&112& & &32.9.19.23.961&961&1104 \cr
\noalign{\hrule}
 & &27.23.61.433&221&212& & &11.361.31.137&351&10 \cr
1486&16402473&8.3.13.17.23.53.61&2165&946&1504&16864837&4.27.5.13.137&913&868 \cr
 & &32.5.11.13.43.433&715&688& & &32.3.7.11.31.83&249&112 \cr
\noalign{\hrule}
 & &9.5.7.11.47.101&943&1178& & &3.13.19.23.991&649&662 \cr
1487&16448355&4.3.11.19.23.31.41&65&188&1505&16889613&4.11.59.331.991&9269&10260 \cr
 & &32.5.13.19.31.47&403&304& & &32.27.5.11.13.19.23.31&495&496 \cr
\noalign{\hrule}
 & &5.23.31.41.113&57&98& & &9.5.11.19.23.79&139&70 \cr
1488&16516645&4.3.49.19.23.113&275&162&1506&17088885&4.3.25.7.79.139&779&1196 \cr
 & &16.243.25.49.11&2673&1960& & &32.7.13.19.23.41&287&208 \cr
\noalign{\hrule}
 & &9.121.17.19.47&1475&582& & &9.11.13.53.251&1725&1036 \cr
1489&16532109&4.27.25.59.97&19&116&1507&17120961&8.27.25.7.23.37&31&4 \cr
 & &32.5.19.29.59&145&944& & &64.5.23.31.37&4255&992 \cr
\noalign{\hrule}
 & &81.7.19.29.53&23&110& & &125.23.59.101&63&38 \cr
1490&16558101&4.27.5.11.23.53&3451&3704&1508&17132125&4.9.5.7.19.23.59&143&202 \cr
 & &64.7.17.29.463&463&544& & &16.3.7.11.13.19.101&741&616 \cr
\noalign{\hrule}
 & &5.11.13.19.23.53&441&142& & &9.5.23.73.227&5947&4268 \cr
1491&16560115&4.9.5.49.19.71&391&106&1509&17150985&8.11.19.97.313&377&690 \cr
 & &16.3.7.17.23.53&119&24& & &32.3.5.13.19.23.29&377&304 \cr
\noalign{\hrule}
 & &5.7.13.59.617&1159&1926& & &3.49.11.13.19.43&3&46 \cr
1492&16563365&4.9.7.19.61.107&803&1230&1510&17174157&4.9.11.13.19.23&1255&1462 \cr
 & &16.27.5.11.41.73&2993&2376& & &16.5.17.43.251&251&680 \cr
\noalign{\hrule}
 & &3.7.13.17.43.83&339&220& & &9.11.29.53.113&161&422 \cr
1493&16563729&8.9.5.11.83.113&289&206&1511&17194419&4.7.23.113.211&795&682 \cr
 & &32.289.103.113&1921&1648& & &16.3.5.11.23.31.53&155&184 \cr
\noalign{\hrule}
 & &41.3721.109&4095&374& & &9.49.121.17.19&101&222 \cr
1494&16629149&4.9.5.7.11.13.17&61&82&1512&17235603&4.27.49.37.101&1207&2530 \cr
 & &16.3.5.17.41.61&85&24& & &16.5.11.17.23.71&355&184 \cr
\noalign{\hrule}
}%
}
$$
\eject
\vglue -23 pt
\noindent\hskip 1 in\hbox to 6.5 in{\ 1513 -- 1548 \hfill\fbd 17242045 -- 18299325\frb}
\vskip -9 pt
$$
\vbox{
\nointerlineskip
\halign{\strut
    \vrule \ \ \hfil \frb #\ 
   &\vrule \hfil \ \ \fbb #\frb\ 
   &\vrule \hfil \ \ \frb #\ \hfil
   &\vrule \hfil \ \ \frb #\ 
   &\vrule \hfil \ \ \frb #\ \ \vrule \hskip 2 pt
   &\vrule \ \ \hfil \frb #\ 
   &\vrule \hfil \ \ \fbb #\frb\ 
   &\vrule \hfil \ \ \frb #\ \hfil
   &\vrule \hfil \ \ \frb #\ 
   &\vrule \hfil \ \ \frb #\ \vrule \cr%
\noalign{\hrule}
 & &5.31.173.643&1521&1694& & &11.17.361.263&4221&250 \cr
1513&17242045&4.9.7.121.169.31&25&118&1531&17754341&4.9.125.7.67&221&154 \cr
 & &16.3.25.7.11.13.59&3003&2360& & &16.3.49.11.13.17&49&312 \cr
\noalign{\hrule}
 & &41.5329.79&1045&4284& & &27.13.17.19.157&1375&1532 \cr
1514&17260631&8.9.5.7.11.17.19&73&158&1532&17799561&8.3.125.11.13.383&23&406 \cr
 & &32.3.19.73.79&57&16& & &32.125.7.23.29&4669&2000 \cr
\noalign{\hrule}
 & &289.529.113&3267&3380& & &5.23.37.53.79&7917&6698 \cr
1515&17275553&8.27.5.121.169.23&71&94&1533&17815685&4.3.7.13.17.29.197&53&66 \cr
 & &32.9.11.169.47.71&30033&29744& & &16.9.11.29.53.197&2167&2088 \cr
\noalign{\hrule}
 & &27.5.11.13.29.31&59&202& & &25.49.13.19.59&583&642 \cr
1516&17355195&4.3.5.31.59.101&203&262&1534&17851925&4.3.11.13.19.53.107&945&1088 \cr
 & &16.7.29.101.131&707&1048& & &512.81.5.7.17.53&4293&4352 \cr
\noalign{\hrule}
 & &9.5.7.13.31.137&407&4& & &7.23.29.43.89&45&44 \cr
1517&17391465&8.3.5.7.11.37&137&248&1535&17868263&8.9.5.7.11.23.29.43&13&1002 \cr
 & &128.31.137&1&64& & &32.27.11.13.167&2171&4752 \cr
\noalign{\hrule}
 & &3.5.7.11.13.19.61&963&1172& & &9.11.37.67.73&265&338 \cr
1518&17402385&8.27.13.107.293&29&322&1536&17915733&4.5.11.169.37.53&219&626 \cr
 & &32.7.23.29.107&667&1712& & &16.3.53.73.313&313&424 \cr
\noalign{\hrule}
 & &5.7.43.71.163&85&78& & &27.125.47.113&3107&2768 \cr
1519&17417365&4.3.25.13.17.43.71&1141&66&1537&17924625&32.9.13.173.239&775&2332 \cr
 & &16.9.7.11.13.163&143&72& & &256.25.11.31.53&1643&1408 \cr
\noalign{\hrule}
 & &5.7.71.79.89&187&258& & &9.5.7.19.31.97&83&88 \cr
1520&17472035&4.3.7.11.17.43.79&585&284&1538&17996895&16.7.11.31.83.97&169&510 \cr
 & &32.27.5.13.17.71&351&272& & &64.3.5.169.17.83&2873&2656 \cr
\noalign{\hrule}
 & &81.125.7.13.19&397&478& & &27.5.11.17.23.31&1589&1930 \cr
1521&17506125&4.13.19.239.397&4851&310&1539&17999685&4.3.25.7.193.227&713&638 \cr
 & &16.9.5.49.11.31&77&248& & &16.11.23.29.31.227&227&232 \cr
\noalign{\hrule}
 & &3.25.11.13.23.71&37&34& & &27.25.13.29.71&971&946 \cr
1522&17513925&4.25.11.13.17.23.37&2313&1988&1540&18067725&4.11.13.29.43.971&425&48 \cr
 & &32.9.7.37.71.257&1799&1776& & &128.3.25.17.971&971&1088 \cr
\noalign{\hrule}
 & &9.5.13.79.379&521&506& & &9.5.7.11.23.227&163&2 \cr
1523&17515485&4.3.11.23.379.521&12245&262&1541&18090765&4.3.163.227&115&112 \cr
 & &16.5.31.79.131&131&248& & &128.5.7.23.163&163&64 \cr
\noalign{\hrule}
 & &5.19.239.773&267&506& & &5.7.19.23.29.41&401&378 \cr
1524&17550965&4.3.5.11.19.23.89&501&1546&1542&18185755&4.27.5.49.29.401&5357&1748 \cr
 & &16.9.167.773&167&72& & &32.3.11.19.23.487&487&528 \cr
\noalign{\hrule}
 & &27.5.13.17.19.31&529&308& & &9.11.23.61.131&61&38 \cr
1525&17572815&8.5.7.11.19.529&9&124&1543&18195507&4.19.3721.131&3105&616 \cr
 & &64.9.11.23.31&253&32& & &64.27.5.7.11.23&35&96 \cr
\noalign{\hrule}
 & &27.49.97.137&55&82& & &5.11.169.37.53&437&252 \cr
1526&17581347&4.5.49.11.41.97&73&24&1544&18227495&8.9.7.11.13.19.23&37&106 \cr
 & &64.3.5.11.41.73&2993&1760& & &32.3.7.19.37.53&57&112 \cr
\noalign{\hrule}
 & &5.31.73.1559&3861&3934& & &9.13.19.29.283&149&700 \cr
1527&17640085&4.27.7.11.13.31.281&1&280&1545&18244161&8.3.25.7.13.149&17&22 \cr
 & &64.3.5.49.11.13&1911&352& & &32.5.7.11.17.149&5215&2992 \cr
\noalign{\hrule}
 & &27.11.19.31.101&299&812& & &7.17.23.59.113&619&738 \cr
1528&17668233&8.7.13.23.29.31&535&132&1546&18247579&4.9.41.113.619&3245&1388 \cr
 & &64.3.5.7.11.107&749&160& & &32.3.5.11.59.347&1041&880 \cr
\noalign{\hrule}
 & &3.11.29.59.313&481&168& & &3.5.7.13.43.311&119&76 \cr
1529&17672919&16.9.7.13.29.37&65&268&1547&18254145&8.49.17.19.311&261&572 \cr
 & &128.5.169.67&845&4288& & &64.9.11.13.19.29&627&928 \cr
\noalign{\hrule}
 & &9.13.19.73.109&8525&9506& & &3.25.11.41.541&37&578 \cr
1530&17688411&4.25.49.11.31.97&513&1192&1548&18299325&4.5.11.289.37&287&342 \cr
 & &64.27.5.7.19.149&745&672& & &16.9.7.17.19.41&133&408 \cr
\noalign{\hrule}
}%
}
$$
\eject
\vglue -23 pt
\noindent\hskip 1 in\hbox to 6.5 in{\ 1549 -- 1584 \hfill\fbd 18343121 -- 19452675\frb}
\vskip -9 pt
$$
\vbox{
\nointerlineskip
\halign{\strut
    \vrule \ \ \hfil \frb #\ 
   &\vrule \hfil \ \ \fbb #\frb\ 
   &\vrule \hfil \ \ \frb #\ \hfil
   &\vrule \hfil \ \ \frb #\ 
   &\vrule \hfil \ \ \frb #\ \ \vrule \hskip 2 pt
   &\vrule \ \ \hfil \frb #\ 
   &\vrule \hfil \ \ \fbb #\frb\ 
   &\vrule \hfil \ \ \frb #\ \hfil
   &\vrule \hfil \ \ \frb #\ 
   &\vrule \hfil \ \ \frb #\ \vrule \cr%
\noalign{\hrule}
 & &23.601.1327&675&652& & &3.25.23.103.107&16819&16244 \cr
1549&18343121&8.27.25.163.601&1309&494&1567&19011225&8.121.31.131.139&5&126 \cr
 & &32.9.5.7.11.13.17.19&15561&14960& & &32.9.5.7.31.139&973&1488 \cr
\noalign{\hrule}
 & &5.11.19.73.241&449&354& & &5.7.13.19.31.71&13&18 \cr
1550&18384685&4.3.59.241.449&345&104&1568&19027645&4.9.7.169.19.71&2449&1100 \cr
 & &64.9.5.13.23.59&2691&1888& & &32.3.25.11.31.79&395&528 \cr
\noalign{\hrule}
 & &3.11.17.71.463&835&372& & &27.5.11.13.23.43&397&98 \cr
1551&18441753&8.9.5.11.31.167&133&34&1569&19092645&4.3.49.43.397&253&650 \cr
 & &32.5.7.17.19.31&1085&304& & &16.25.7.11.13.23&35&8 \cr
\noalign{\hrule}
 & &3.11.29.97.199&305&14& & &3.13.17.19.37.41&1&40 \cr
1552&18472971&4.5.7.61.199&117&82&1570&19109649&16.5.17.19.37&143&180 \cr
 & &16.9.13.41.61&2501&312& & &128.9.25.11.13&75&704 \cr
\noalign{\hrule}
 & &7.11.293.823&185&108& & &3.5.7.11.13.19.67&1593&1124 \cr
1553&18567703&8.27.5.37.823&245&578&1571&19114095&8.81.5.59.281&343&62 \cr
 & &32.3.25.49.289&2023&1200& & &32.343.31.59&1829&784 \cr
\noalign{\hrule}
 & &9.25.7.11.29.37&431&494& & &9.5.11.13.19.157&4709&1726 \cr
1554&18589725&4.11.13.19.29.431&375&56&1572&19195605&4.17.277.863&293&570 \cr
 & &64.3.125.7.13.19&247&160& & &16.3.5.17.19.293&293&136 \cr
\noalign{\hrule}
 & &5.7.11.193.251&29&222& & &9.121.289.61&95&194 \cr
1555&18650555&4.3.5.7.11.29.37&513&502&1573&19197981&4.5.11.19.61.97&93&578 \cr
 & &16.81.19.37.251&703&648& & &16.3.289.19.31&31&152 \cr
\noalign{\hrule}
 & &9.5.121.23.149&107&8& & &9.11.19.59.173&271&260 \cr
1556&18660015&16.11.107.149&69&80&1574&19199367&8.5.13.19.173.271&3405&118 \cr
 & &512.3.5.23.107&107&256& & &32.3.25.59.227&227&400 \cr
\noalign{\hrule}
 & &27.7.121.19.43&255&46& & &9.7.13.53.443&253&190 \cr
1557&18683973&4.81.5.11.17.23&167&86&1575&19229301&4.5.11.13.19.23.53&17&282 \cr
 & &16.5.17.43.167&167&680& & &16.3.11.17.19.47&517&2584 \cr
\noalign{\hrule}
 & &31.47.101.127&1419&4550& & &3.7.11.169.17.29&95&282 \cr
1558&18688939&4.3.25.7.11.13.43&101&114&1576&19246227&4.9.5.7.13.19.47&473&382 \cr
 & &16.9.5.7.11.19.101&665&792& & &16.11.43.47.191&2021&1528 \cr
\noalign{\hrule}
 & &3.5.11.19.43.139&1557&1084& & &27.7.11.13.23.31&9125&10126 \cr
1559&18737895&8.27.5.173.271&203&68&1577&19270251&4.125.61.73.83&2961&7414 \cr
 & &64.7.17.29.173&3451&5536& & &16.9.7.11.47.337&337&376 \cr
\noalign{\hrule}
 & &9.23.197.461&2173&1976& & &125.11.101.139&113&12 \cr
1560&18799119&16.13.19.23.41.53&999&220&1578&19303625&8.3.11.113.139&75&64 \cr
 & &128.27.5.11.13.37&2035&2496& & &1024.9.25.113&1017&512 \cr
\noalign{\hrule}
 & &9.11.53.59.61&203&380& & &59.331.991&9269&10260 \cr
1561&18883953&8.3.5.7.19.29.61&77&106&1579&19353239&8.27.5.13.19.23.31&649&662 \cr
 & &32.5.49.11.19.53&245&304& & &32.9.5.11.31.59.331&495&496 \cr
\noalign{\hrule}
 & &3.13.19.71.359&191&550& & &5.11.31.41.277&159&118 \cr
1562&18887349&4.25.11.71.191&983&792&1580&19363685&4.3.5.11.31.53.59&3783&4432 \cr
 & &64.9.121.983&2949&3872& & &128.9.13.97.277&873&832 \cr
\noalign{\hrule}
 & &3.11.19.47.641&391&250& & &5.11.37.89.107&169&276 \cr
1563&18889629&4.125.11.17.19.23&47&162&1581&19379305&8.3.11.169.23.37&353&54 \cr
 & &16.81.25.17.47&425&216& & &32.81.13.353&4589&1296 \cr
\noalign{\hrule}
 & &25.7.13.53.157&9&166& & &3.25.7.17.41.53&9&44 \cr
1564&18930275&4.9.13.53.83&469&220&1582&19394025&8.27.5.11.17.41&511&596 \cr
 & &32.3.5.7.11.67&737&48& & &64.7.11.73.149&1639&2336 \cr
\noalign{\hrule}
 & &27.5.49.47.61&247&58& & &3.5.11.19.41.151&1209&1660 \cr
1565&18965205&4.7.13.19.29.47&69&22&1583&19408785&8.9.25.13.31.83&23&302 \cr
 & &16.3.11.19.23.29&437&2552& & &32.23.83.151&83&368 \cr
\noalign{\hrule}
 & &9.11.31.41.151&301&150& & &3.25.11.17.19.73&43&518 \cr
1566&19000179&4.27.25.7.31.43&41&176&1584&19452675&4.7.37.43.73&55&18 \cr
 & &128.5.11.41.43&215&64& & &16.9.5.7.11.43&21&344 \cr
\noalign{\hrule}
}%
}
$$
\eject
\vglue -23 pt
\noindent\hskip 1 in\hbox to 6.5 in{\ 1585 -- 1620 \hfill\fbd 19531655 -- 20452047\frb}
\vskip -9 pt
$$
\vbox{
\nointerlineskip
\halign{\strut
    \vrule \ \ \hfil \frb #\ 
   &\vrule \hfil \ \ \fbb #\frb\ 
   &\vrule \hfil \ \ \frb #\ \hfil
   &\vrule \hfil \ \ \frb #\ 
   &\vrule \hfil \ \ \frb #\ \ \vrule \hskip 2 pt
   &\vrule \ \ \hfil \frb #\ 
   &\vrule \hfil \ \ \fbb #\frb\ 
   &\vrule \hfil \ \ \frb #\ \hfil
   &\vrule \hfil \ \ \frb #\ 
   &\vrule \hfil \ \ \frb #\ \vrule \cr%
\noalign{\hrule}
 & &5.11.13.59.463&93&556& & &3.7.11.13.61.109&61&82 \cr
1585&19531655&8.3.5.13.31.139&37&102&1603&19966947&4.41.3721.109&4095&374 \cr
 & &32.9.17.31.37&629&4464& & &16.9.5.7.11.13.17&85&24 \cr
\noalign{\hrule}
 & &3.11.31.61.313&187&126& & &27.125.61.97&143&82 \cr
1586&19532139&4.27.7.121.17.31&1525&1742&1604&19969875&4.3.5.11.13.41.97&1891&2086 \cr
 & &16.25.13.17.61.67&1675&1768& & &16.7.11.31.61.149&1639&1736 \cr
\noalign{\hrule}
 & &5.7.17.107.307&737&4482& & &3.5.7.11.13.31.43&307&338 \cr
1587&19545155&4.27.11.67.83&4075&4142&1605&20014995&4.7.11.2197.307&6001&9378 \cr
 & &16.3.25.19.109.163&6213&6520& & &16.9.17.353.521&8857&8472 \cr
\noalign{\hrule}
 & &11.13.41.47.71&137745&137816& & &3.5.11.29.53.79&489&94 \cr
1588&19564831&16.9.5.7.23.107.3061&3361&5822&1606&20034795&4.9.29.47.163&3961&3700 \cr
 & &64.3.5.7.41.71.3361&3361&3360& & &32.25.17.37.233&3145&3728 \cr
\noalign{\hrule}
 & &9.17.97.1319&165&1484& & &3.7.67.109.131&65&44 \cr
1589&19575279&8.27.5.7.11.53&97&92&1607&20090553&8.5.11.13.67.131&763&108 \cr
 & &64.11.23.53.97&583&736& & &64.27.7.11.109&99&32 \cr
\noalign{\hrule}
 & &3.5.13.31.41.79&407&802& & &3.25.13.23.29.31&319&394 \cr
1590&19579755&4.11.37.41.401&959&558&1608&20160075&4.11.13.841.197&427&414 \cr
 & &16.9.7.11.31.137&411&616& & &16.9.7.11.23.61.197&4697&4728 \cr
\noalign{\hrule}
 & &27.25.7.11.13.29&139&4& & &9.5.17.23.31.37&319&208 \cr
1591&19594575&8.5.7.29.139&99&104&1609&20181465&32.3.5.11.13.23.29&7&62 \cr
 & &128.9.11.13.139&139&64& & &128.7.13.29.31&203&832 \cr
\noalign{\hrule}
 & &9.13.23.67.109&407&1010& & &121.169.23.43&153&406 \cr
1592&19652373&4.5.11.23.37.101&1995&1742&1610&20224061&4.9.7.11.13.17.29&43&230 \cr
 & &16.3.25.7.13.19.67&175&152& & &16.3.5.23.29.43&145&24 \cr
\noalign{\hrule}
 & &9.5.49.169.53&4601&4356& & &9.31.1369.53&731&638 \cr
1593&19750185&8.81.121.43.107&8215&4732&1611&20243403&4.3.11.17.29.43.53&259&1160 \cr
 & &64.5.7.169.31.53&31&32& & &64.5.7.841.37&841&1120 \cr
\noalign{\hrule}
 & &3.11.13.29.37.43&6137&7580& & &3.73.193.479&349&130 \cr
1594&19793631&8.5.17.361.379&9&370&1612&20245893&4.5.13.193.349&657&308 \cr
 & &32.9.25.17.37&425&48& & &32.9.7.11.13.73&143&336 \cr
\noalign{\hrule}
 & &7.11.361.23.31&5185&2658& & &3.5.7.11.13.19.71&391&106 \cr
1595&19819261&4.3.5.17.61.443&69&374&1613&20255235&4.11.13.17.23.53&441&142 \cr
 & &16.9.11.289.23&289&72& & &16.9.49.17.71&119&24 \cr
\noalign{\hrule}
 & &27.5.19.59.131&213&82& & &7.17.23.31.239&3509&3900 \cr
1596&19824885&4.81.19.41.71&715&634&1614&20278433&8.3.25.7.121.13.29&239&36 \cr
 & &16.5.11.13.41.317&4121&3608& & &64.27.11.13.239&297&416 \cr
\noalign{\hrule}
 & &27.25.7.13.17.19&113&22& & &5.17.19.29.433&2881&5346 \cr
1597&19840275&4.5.11.17.19.113&141&46&1615&20279555&4.243.11.43.67&433&304 \cr
 & &16.3.23.47.113&2599&376& & &128.81.19.433&81&64 \cr
\noalign{\hrule}
 & &9.7.11.23.29.43&325&578& & &25.169.61.79&297&4522 \cr
1598&19875933&4.3.25.13.289.29&1&86&1616&20360275&4.27.7.11.17.19&65&122 \cr
 & &16.5.13.17.43&13&680& & &16.9.5.7.13.61&7&72 \cr
\noalign{\hrule}
 & &9.7.11.23.29.43&923&1310& & &5.11.23.71.227&4563&658 \cr
1599&19875933&4.5.13.23.71.131&77&54&1617&20388005&4.27.7.169.47&277&230 \cr
 & &16.27.5.7.11.13.71&355&312& & &16.9.5.7.23.277&277&504 \cr
\noalign{\hrule}
 & &5.7.11.13.41.97&27&38& & &29.31.37.613&8965&10038 \cr
1600&19904885&4.27.7.19.41.97&155&524&1618&20390219&4.3.5.7.11.163.239&1885&744 \cr
 & &32.3.5.19.31.131&2489&1488& & &64.9.25.13.29.31&325&288 \cr
\noalign{\hrule}
 & &27.7.19.23.241&143&580& & &3.5.37.83.443&429&14 \cr
1601&19904913&8.9.5.7.11.13.29&337&482&1619&20406795&4.9.7.11.13.37&125&134 \cr
 & &32.11.241.337&337&176& & &16.125.11.13.67&871&2200 \cr
\noalign{\hrule}
 & &81.125.11.179&97&178& & &3.7.11.29.43.71&1425&1628 \cr
1602&19936125&4.5.89.97.179&133&312&1620&20452047&8.9.25.121.19.37&3013&5312 \cr
 & &64.3.7.13.19.97&1729&3104& & &1024.23.83.131&10873&11776 \cr
\noalign{\hrule}
}%
}
$$
\eject
\vglue -23 pt
\noindent\hskip 1 in\hbox to 6.5 in{\ 1621 -- 1656 \hfill\fbd 20563881 -- 21845635\frb}
\vskip -9 pt
$$
\vbox{
\nointerlineskip
\halign{\strut
    \vrule \ \ \hfil \frb #\ 
   &\vrule \hfil \ \ \fbb #\frb\ 
   &\vrule \hfil \ \ \frb #\ \hfil
   &\vrule \hfil \ \ \frb #\ 
   &\vrule \hfil \ \ \frb #\ \ \vrule \hskip 2 pt
   &\vrule \ \ \hfil \frb #\ 
   &\vrule \hfil \ \ \fbb #\frb\ 
   &\vrule \hfil \ \ \frb #\ \hfil
   &\vrule \hfil \ \ \frb #\ 
   &\vrule \hfil \ \ \frb #\ \vrule \cr%
\noalign{\hrule}
 & &3.13.31.73.233&3401&3388& & &9.11.169.31.41&191&150 \cr
1621&20563881&8.7.121.19.179.233&169&1800&1639&21265101&4.27.25.169.191&2759&1804 \cr
 & &128.9.25.11.169.19&3575&3648& & &32.5.11.31.41.89&89&80 \cr
\noalign{\hrule}
 & &9.5.49.47.199&579&814& & &27.11.23.53.59&221&428 \cr
1622&20623365&4.27.7.11.37.193&1175&176&1640&21360537&8.3.13.17.53.107&505&184 \cr
 & &128.25.121.47&121&320& & &128.5.17.23.101&505&1088 \cr
\noalign{\hrule}
 & &5.7.1331.443&1773&442& & &5.7.11.19.23.127&16507&17142 \cr
1623&20637155&4.9.7.13.17.197&109&88&1641&21367115&4.3.17.971.2857&943&1914 \cr
 & &64.3.11.13.17.109&1853&1248& & &16.9.11.17.23.29.41&1189&1224 \cr
\noalign{\hrule}
 & &9.5.49.11.23.37&793&58& & &11.31.131.479&8145&6704 \cr
1624&20641005&4.3.11.13.29.61&901&230&1642&21397409&32.9.5.181.419&481&62 \cr
 & &16.5.17.23.53&17&424& & &128.3.5.13.31.37&481&960 \cr
\noalign{\hrule}
 & &9.11.23.31.293&1835&4904& & &7.13.23.37.277&1819&1782 \cr
1625&20681991&16.5.367.613&1349&1716&1643&21451157&4.81.7.11.17.23.107&185&1634 \cr
 & &128.3.11.13.19.71&923&1216& & &16.9.5.11.19.37.43&1881&1720 \cr
\noalign{\hrule}
 & &81.5.11.59.79&29&266& & &3.7.13.17.41.113&87&200 \cr
1626&20764755&4.27.7.11.19.29&91&118&1644&21501753&16.9.25.13.17.29&539&46 \cr
 & &16.49.13.29.59&637&232& & &64.5.49.11.23&161&1760 \cr
\noalign{\hrule}
 & &3.5.11.13.89.109&21&1178& & &3.3125.121.19&1591&1534 \cr
1627&20808645&4.9.5.7.19.31&109&46&1645&21553125&4.121.13.37.43.59&9&1582 \cr
 & &16.19.23.109&437&8& & &16.9.7.59.113&413&2712 \cr
\noalign{\hrule}
 & &9.5.17.163.167&1827&1012& & &11.67.83.353&135&218 \cr
1628&20824065&8.81.7.11.23.29&167&86&1646&21593363&4.27.5.11.67.109&83&16 \cr
 & &32.7.29.43.167&301&464& & &128.3.5.83.109&545&192 \cr
\noalign{\hrule}
 & &3.25.49.13.19.23&1497&1562& & &23.67.107.131&65&42 \cr
1629&20877675&4.9.5.7.11.71.499&989&3484&1647&21600197&4.3.5.7.13.67.131&321&2024 \cr
 & &32.11.13.23.43.67&737&688& & &64.9.11.23.107&99&32 \cr
\noalign{\hrule}
 & &27.7.17.67.97&209&260& & &3.11.31.37.571&115&456 \cr
1630&20881287&8.9.5.11.13.19.97&1&98&1648&21612921&16.9.5.19.23.37&455&248 \cr
 & &32.5.49.13.19&1729&80& & &256.25.7.13.31&325&896 \cr
\noalign{\hrule}
 & &81.7.11.17.197&1007&1160& & &3.5.101.109.131&187&318 \cr
1631&20887713&16.9.5.7.19.29.53&1367&170&1649&21632685&4.9.11.17.53.109&131&22 \cr
 & &64.25.17.1367&1367&800& & &16.121.53.131&121&424 \cr
\noalign{\hrule}
 & &27.5.37.47.89&865&874& & &3.49.13.17.23.29&4895&1514 \cr
1632&20894085&4.3.25.19.23.89.173&4699&374&1650&21668829&4.5.11.89.757&351&406 \cr
 & &16.11.17.23.37.127&2159&2024& & &16.27.7.13.29.89&89&72 \cr
\noalign{\hrule}
 & &9.25.7.11.17.71&299&86& & &17.19.229.293&837&4730 \cr
1633&20911275&4.3.5.13.17.23.43&107&22&1651&21672331&4.27.5.11.31.43&493&152 \cr
 & &16.11.13.23.107&1391&184& & &64.9.17.19.29&261&32 \cr
\noalign{\hrule}
 & &9.25.13.23.311&5137&2338& & &9.125.23.841&157&418 \cr
1634&20922525&4.7.11.167.467&317&150&1652&21760875&4.5.11.19.29.157&233&552 \cr
 & &16.3.25.7.11.317&317&616& & &64.3.19.23.233&233&608 \cr
\noalign{\hrule}
 & &3.5.29.127.379&4543&6448& & &7.13.23.101.103&55&48 \cr
1635&20937855&32.7.11.13.31.59&377&36&1653&21773479&32.3.5.11.13.23.101&103&402 \cr
 & &256.9.169.29&507&128& & &128.9.11.67.103&737&576 \cr
\noalign{\hrule}
 & &3.25.13.17.31.41&139&836& & &3.11.17.71.547&129&58 \cr
1636&21066825&8.11.19.31.139&85&54&1654&21787557&4.9.29.43.547&5885&5338 \cr
 & &32.27.5.11.17.19&99&304& & &16.5.11.17.107.157&785&856 \cr
\noalign{\hrule}
 & &9.7.53.59.107&2629&3042& & &81.5.11.59.83&3391&3332 \cr
1637&21079107&4.81.11.169.239&361&530&1655&21816135&8.5.49.11.17.3391&1577&4968 \cr
 & &16.5.361.53.239&1805&1912& & &128.27.7.19.23.83&437&448 \cr
\noalign{\hrule}
 & &27.5.7.17.1319&97&92& & &5.7.59.71.149&4257&6322 \cr
1638&21189735&8.17.23.97.1319&165&1484&1656&21845635&4.9.11.29.43.109&65&22 \cr
 & &64.3.5.7.11.23.53&583&736& & &16.3.5.121.13.109&1573&2616 \cr
\noalign{\hrule}
}%
}
$$
\eject
\vglue -23 pt
\noindent\hskip 1 in\hbox to 6.5 in{\ 1657 -- 1692 \hfill\fbd 21846129 -- 22752405\frb}
\vskip -9 pt
$$
\vbox{
\nointerlineskip
\halign{\strut
    \vrule \ \ \hfil \frb #\ 
   &\vrule \hfil \ \ \fbb #\frb\ 
   &\vrule \hfil \ \ \frb #\ \hfil
   &\vrule \hfil \ \ \frb #\ 
   &\vrule \hfil \ \ \frb #\ \ \vrule \hskip 2 pt
   &\vrule \ \ \hfil \frb #\ 
   &\vrule \hfil \ \ \fbb #\frb\ 
   &\vrule \hfil \ \ \frb #\ \hfil
   &\vrule \hfil \ \ \frb #\ 
   &\vrule \hfil \ \ \frb #\ \vrule \cr%
\noalign{\hrule}
 & &3.599.12157&299&300& & &81.13.17.29.43&41&418 \cr
1657&21846129&8.9.25.13.23.12157&649&12806&1675&22322547&4.3.11.19.41.43&169&40 \cr
 & &32.5.11.19.59.337&18535&17936& & &64.5.169.41&65&1312 \cr
\noalign{\hrule}
 & &121.13.17.19.43&1195&378& & &5.13.17.89.227&1433&1518 \cr
1658&21847397&4.27.5.7.17.239&299&418&1676&22324315&4.3.11.23.89.1433&227&1206 \cr
 & &16.9.5.11.13.19.23&115&72& & &16.27.23.67.227&621&536 \cr
\noalign{\hrule}
 & &9.25.13.31.241&283&42& & &25.7.11.13.19.47&93&82 \cr
1659&21852675&4.27.7.31.283&79&110&1677&22347325&4.3.13.19.31.41.47&11&600 \cr
 & &16.5.11.79.283&869&2264& & &64.9.25.11.41&41&288 \cr
\noalign{\hrule}
 & &3.25.11.101.263&49&214& & &625.7.19.269&303&572 \cr
1660&21914475&4.5.49.101.107&621&86&1678&22360625&8.3.5.11.13.19.101&1173&62 \cr
 & &16.27.7.23.43&2709&184& & &32.9.17.23.31&4743&368 \cr
\noalign{\hrule}
 & &9.5.121.29.139&59&62& & &9.11.23.43.229&265&494 \cr
1661&21948795&4.3.5.29.31.59.139&187&1898&1679&22421619&4.3.5.13.19.43.53&335&1342 \cr
 & &16.11.13.17.31.73&2263&1768& & &16.25.11.61.67&1525&536 \cr
\noalign{\hrule}
 & &11.13.29.47.113&17&30& & &27.25.149.223&4873&1148 \cr
1662&22024717&4.3.5.11.17.29.113&611&1854&1680&22428225&8.7.11.41.443&447&4 \cr
 & &16.27.13.47.103&103&216& & &64.3.7.149&1&224 \cr
\noalign{\hrule}
 & &81.5.13.53.79&163&242& & &3.49.11.17.19.43&1495&956 \cr
1663&22044555&4.121.13.53.163&2147&4266&1681&22458513&8.5.13.17.23.239&77&162 \cr
 & &16.27.19.79.113&113&152& & &32.81.7.11.13.23&299&432 \cr
\noalign{\hrule}
 & &27.7.11.13.19.43&833&586& & &27.5.7.13.31.59&817&1012 \cr
1664&22081059&4.9.343.17.293&95&248&1682&22469265&8.9.7.11.19.23.43&47&124 \cr
 & &64.5.19.31.293&1465&992& & &64.23.31.43.47&1081&1376 \cr
\noalign{\hrule}
 & &17.19.67.1021&10269&9130& & &25.7.19.67.101&187&282 \cr
1665&22095461&4.9.5.7.11.83.163&1021&772&1683&22500275&4.3.5.11.17.47.101&67&168 \cr
 & &32.3.5.7.193.1021&579&560& & &64.9.7.11.17.67&187&288 \cr
\noalign{\hrule}
 & &9.5.13.23.31.53&17&282& & &9.3125.11.73&16859&17516 \cr
1666&22106565&4.27.17.31.47&715&742&1684&22584375&8.23.29.151.733&33&700 \cr
 & &16.5.7.11.13.17.53&119&88& & &64.3.25.7.11.151&151&224 \cr
\noalign{\hrule}
 & &27.5.11.13.31.37&4601&4606& & &9.7.11.13.23.109&555&446 \cr
1667&22142835&4.49.13.37.43.47.107&23&3936&1685&22585563&4.27.5.23.37.223&247&868 \cr
 & &256.3.7.23.41.47&6601&6016& & &32.7.13.19.31.37&703&496 \cr
\noalign{\hrule}
 & &3.49.13.67.173&97&76& & &25.7.31.43.97&429&646 \cr
1668&22150401&8.7.13.19.67.97&865&396&1686&22627675&4.3.11.13.17.19.97&657&410 \cr
 & &64.9.5.11.19.173&285&352& & &16.27.5.17.41.73&2993&3672 \cr
\noalign{\hrule}
 & &27.5.13.19.23.29&71&74& & &121.23.47.173&1305&1478 \cr
1669&22241115&4.9.13.19.23.37.71&671&2020&1687&22628573&4.9.5.29.47.739&299&440 \cr
 & &32.5.11.37.61.101&6161&6512& & &64.3.25.11.13.23.29&975&928 \cr
\noalign{\hrule}
 & &11.13.23.67.101&99&200& & &25.19.23.31.67&349&924 \cr
1670&22256663&16.9.25.121.67&403&202&1688&22691225&8.3.7.11.31.349&159&190 \cr
 & &64.3.5.13.31.101&93&160& & &32.9.5.7.11.19.53&693&848 \cr
\noalign{\hrule}
 & &27.7.169.17.41&5071&1858& & &5.7.13.23.41.53&183&22 \cr
1671&22262877&4.11.461.929&4879&5340&1689&22740445&4.3.11.13.53.61&375&314 \cr
 & &32.3.5.7.17.41.89&89&80& & &16.9.125.11.157&2475&1256 \cr
\noalign{\hrule}
 & &5.7.11.13.61.73&589&204& & &9.5.7.17.31.137&241&286 \cr
1672&22287265&8.3.17.19.31.73&957&430&1690&22742685&4.7.11.13.137.241&1087&2046 \cr
 & &32.9.5.11.29.43&261&688& & &16.3.121.31.1087&1087&968 \cr
\noalign{\hrule}
 & &27.11.13.23.251&5033&1744& & &3.13.1681.347&4777&266 \cr
1673&22289553&32.7.109.719&305&414&1691&22748973&4.7.17.19.281&1271&990 \cr
 & &128.9.5.7.23.61&305&448& & &16.9.5.11.31.41&93&440 \cr
\noalign{\hrule}
 & &5.71.227.277&31617&31262& & &9.5.13.19.23.89&965&1082 \cr
1674&22322045&4.27.49.11.29.1171&25&4&1692&22752405&4.25.19.193.541&33&508 \cr
 & &32.9.25.7.11.1171&8197&7920& & &32.3.11.127.193&1397&3088 \cr
\noalign{\hrule}
}%
}
$$
\eject
\vglue -23 pt
\noindent\hskip 1 in\hbox to 6.5 in{\ 1693 -- 1728 \hfill\fbd 22855959 -- 23647519\frb}
\vskip -9 pt
$$
\vbox{
\nointerlineskip
\halign{\strut
    \vrule \ \ \hfil \frb #\ 
   &\vrule \hfil \ \ \fbb #\frb\ 
   &\vrule \hfil \ \ \frb #\ \hfil
   &\vrule \hfil \ \ \frb #\ 
   &\vrule \hfil \ \ \frb #\ \ \vrule \hskip 2 pt
   &\vrule \ \ \hfil \frb #\ 
   &\vrule \hfil \ \ \fbb #\frb\ 
   &\vrule \hfil \ \ \frb #\ \hfil
   &\vrule \hfil \ \ \frb #\ 
   &\vrule \hfil \ \ \frb #\ \vrule \cr%
\noalign{\hrule}
 & &27.7.31.47.83&793&1780& & &11.13.37.53.83&245&162 \cr
1693&22855959&8.9.5.13.61.89&11&2&1711&23275109&4.81.5.49.13.53&397&1034 \cr
 & &32.5.11.61.89&979&4880& & &16.3.5.11.47.397&1191&1880 \cr
\noalign{\hrule}
 & &3.11.37.97.193&437&630& & &7.11.43.79.89&585&284 \cr
1694&22858341&4.27.5.7.19.23.37&97&902&1712&23279641&8.9.5.13.71.89&187&258 \cr
 & &16.11.19.41.97&41&152& & &32.27.11.13.17.43&351&272 \cr
\noalign{\hrule}
 & &81.125.7.17.19&37&44& & &3.49.19.31.269&869&62 \cr
1695&22892625&8.125.11.17.19.37&711&1414&1713&23290827&4.11.961.79&441&520 \cr
 & &32.9.7.11.79.101&1111&1264& & &64.9.5.49.11.13&429&160 \cr
\noalign{\hrule}
 & &3.7.11.13.29.263&247&16& & &81.19.109.139&29&110 \cr
1696&22903881&32.169.19.29&99&70&1714&23317389&4.5.11.19.29.109&377&168 \cr
 & &128.9.5.7.11.19&57&320& & &64.3.7.13.841&841&2912 \cr
\noalign{\hrule}
 & &19.31.97.401&495&94& & &5.49.11.13.23.29&103&12 \cr
1697&22910333&4.9.5.11.47.97&401&304&1715&23368345&8.3.7.11.29.103&201&520 \cr
 & &128.3.11.19.401&33&64& & &128.9.5.13.67&603&64 \cr
\noalign{\hrule}
 & &81.7.11.13.283&431&418& & &3.7.13.23.37.101&737&576 \cr
1698&22945923&4.27.7.121.19.431&32545&29528&1716&23464623&128.27.11.37.67&505&494 \cr
 & &64.5.23.283.3691&3691&3680& & &512.5.13.19.67.101&1273&1280 \cr
\noalign{\hrule}
 & &3.5.49.11.17.167&667&166& & &25.13.17.31.137&1823&1602 \cr
1699&22953315&4.5.11.23.29.83&81&334&1717&23464675&4.9.31.89.1823&1355&4114 \cr
 & &16.81.29.167&29&216& & &16.3.5.121.17.271&813&968 \cr
\noalign{\hrule}
 & &5.49.17.37.149&1829&2574& & &3.7.11.17.43.139&57&244 \cr
1700&22961645&4.9.7.11.13.31.59&437&2384&1718&23471679&8.9.19.61.139&595&656 \cr
 & &128.3.19.23.149&437&192& & &256.5.7.17.19.41&779&640 \cr
\noalign{\hrule}
 & &9.5.121.41.103&2671&2774& & &9.11.13.17.29.37&305&188 \cr
1701&22994235&4.19.41.73.2671&1725&946&1719&23476167&8.5.11.37.47.61&2451&416 \cr
 & &16.3.25.11.23.43.73&1679&1720& & &512.3.13.19.43&817&256 \cr
\noalign{\hrule}
 & &11.13.19.37.229&1269&1250& & &3.5.7.17.59.223&137&86 \cr
1702&23021141&4.27.625.13.37.47&437&5062&1720&23485245&4.5.7.43.59.137&99&314 \cr
 & &16.3.5.19.23.2531&2531&2760& & &16.9.11.137.157&1727&3288 \cr
\noalign{\hrule}
 & &1331.13.31.43&595&738& & &13.23.31.43.59&15713&14946 \cr
1703&23064899&4.9.5.7.121.17.41&731&116&1721&23515453&4.3.19.47.53.827&33&860 \cr
 & &32.3.289.29.43&867&464& & &32.9.5.11.43.53&477&880 \cr
\noalign{\hrule}
 & &81.11.19.29.47&3071&2990& & &81.11.13.19.107&1165&868 \cr
1704&23074227&4.5.13.23.37.47.83&1&1080&1722&23548239&8.3.5.7.13.31.233&19&214 \cr
 & &64.27.25.37&925&32& & &32.7.19.31.107&31&112 \cr
\noalign{\hrule}
 & &9.25.7.11.31.43&1717&1748& & &3.7.13.361.239&2095&1012 \cr
1705&23094225&8.5.17.19.23.43.101&1609&2046&1723&23554167&8.5.7.11.23.419&171&248 \cr
 & &32.3.11.31.101.1609&1609&1616& & &128.9.5.19.23.31&713&960 \cr
\noalign{\hrule}
 & &9.5.11.13.37.97&73&112& & &9.7.11.41.829&1387&1100 \cr
1706&23095215&32.3.7.11.73.97&233&1300&1724&23554377&8.3.25.121.19.73&203&82 \cr
 & &256.25.13.233&233&640& & &32.5.7.29.41.73&365&464 \cr
\noalign{\hrule}
 & &25.11.13.29.223&67&78& & &7.13.17.97.157&961&1080 \cr
1707&23119525&4.3.5.169.67.223&473&642&1725&23559263&16.27.5.961.97&917&44 \cr
 & &16.9.11.43.67.107&4601&4824& & &128.3.5.7.11.131&1965&704 \cr
\noalign{\hrule}
 & &3.13.59.89.113&77&190& & &7.11.19.89.181&135&46 \cr
1708&23141157&4.5.7.11.13.19.59&279&488&1726&23567467&4.27.5.7.11.19.23&181&104 \cr
 & &64.9.5.7.31.61&1891&3360& & &64.9.13.23.181&207&416 \cr
\noalign{\hrule}
 & &5.11.19.67.331&333&2& & &27.7.29.31.139&4945&914 \cr
1709&23174965&4.9.11.19.37&105&104&1727&23617629&4.5.23.43.457&21&22 \cr
 & &64.27.5.7.13.37&2457&1184& & &16.3.5.7.11.23.457&2285&2024 \cr
\noalign{\hrule}
 & &5.121.19.43.47&227&246& & &7.23.191.769&495&4888 \cr
1710&23231395&4.3.5.11.41.47.227&361&156&1728&23647519&16.9.5.11.13.47&289&322 \cr
 & &32.9.13.361.227&2951&2736& & &64.3.5.7.289.23&289&480 \cr
\noalign{\hrule}
}%
}
$$
\eject
\vglue -23 pt
\noindent\hskip 1 in\hbox to 6.5 in{\ 1729 -- 1764 \hfill\fbd 23667309 -- 24951379\frb}
\vskip -9 pt
$$
\vbox{
\nointerlineskip
\halign{\strut
    \vrule \ \ \hfil \frb #\ 
   &\vrule \hfil \ \ \fbb #\frb\ 
   &\vrule \hfil \ \ \frb #\ \hfil
   &\vrule \hfil \ \ \frb #\ 
   &\vrule \hfil \ \ \frb #\ \ \vrule \hskip 2 pt
   &\vrule \ \ \hfil \frb #\ 
   &\vrule \hfil \ \ \fbb #\frb\ 
   &\vrule \hfil \ \ \frb #\ \hfil
   &\vrule \hfil \ \ \frb #\ 
   &\vrule \hfil \ \ \frb #\ \vrule \cr%
\noalign{\hrule}
 & &81.37.53.149&7015&5054& & &9.529.47.109&479&502 \cr
1729&23667309&4.5.7.361.23.61&583&222&1747&24390603&4.23.47.251.479&11407&390 \cr
 & &16.3.11.37.53.61&61&88& & &16.3.5.11.13.17.61&3355&1768 \cr
\noalign{\hrule}
 & &3.961.43.191&2095&1904& & &3.5.7.19.37.331&807&1510 \cr
1730&23678079&32.5.7.17.31.419&473&54&1748&24432765&4.9.25.151.269&247&22 \cr
 & &128.27.5.7.11.43&385&576& & &16.11.13.19.151&151&1144 \cr
\noalign{\hrule}
 & &5.49.11.59.149&9&68& & &9.23.83.1423&385&362 \cr
1731&23691745&8.9.5.7.17.149&649&394&1749&24448563&4.5.7.11.181.1423&1345&78 \cr
 & &32.3.11.59.197&197&48& & &16.3.25.11.13.269&3575&2152 \cr
\noalign{\hrule}
 & &5.7.31.43.509&363&146& & &9.11.29.83.103&5&314 \cr
1732&23747395&4.3.5.121.43.73&509&294&1750&24544179&4.3.5.83.157&203&46 \cr
 & &16.9.49.11.509&99&56& & &16.5.7.23.29&7&920 \cr
\noalign{\hrule}
 & &3.7.121.17.19.29&927&1130& & &9.25.79.1381&803&578 \cr
1733&23801547&4.27.5.19.103.113&317&2464&1751&24547275&4.11.289.73.79&1073&270 \cr
 & &256.5.7.11.317&317&640& & &16.27.5.17.29.37&493&888 \cr
\noalign{\hrule}
 & &11.289.59.127&755&1404& & &3.11.13.23.47.53&285&232 \cr
1734&23820247&8.27.5.13.17.151&33&118&1752&24578697&16.9.5.13.19.23.29&235&64 \cr
 & &32.81.11.13.59&81&208& & &2048.25.29.47&725&1024 \cr
\noalign{\hrule}
 & &3.5.7.17.361.37&657&638& & &3.5.11.29.53.97&49&534 \cr
1735&23842245&4.27.11.17.19.29.73&2575&3068&1753&24599685&4.9.49.29.89&97&106 \cr
 & &32.25.13.59.73.103&21535&21424& & &16.7.53.89.97&89&56 \cr
\noalign{\hrule}
 & &5.11.23.113.167&45&68& & &125.19.101.103&7397&5478 \cr
1736&23871815&8.9.25.11.17.167&113&388&1754&24707125&4.3.11.13.83.569&255&824 \cr
 & &64.3.17.97.113&291&544& & &64.9.5.11.17.103&187&288 \cr
\noalign{\hrule}
 & &121.17.29.401&1955&1554& & &11.29.71.1091&581&510 \cr
1737&23920853&4.3.5.7.289.23.37&3971&2676&1755&24710059&4.3.5.7.11.17.29.83&1269&2182 \cr
 & &32.9.11.361.223&3249&3568& & &16.81.5.47.1091&405&376 \cr
\noalign{\hrule}
 & &5.49.11.83.107&311&228& & &7.13.439.619&19665&20284 \cr
1738&23934295&8.3.5.19.107.311&101&6&1756&24728431&8.9.5.11.19.23.461&2173&2634 \cr
 & &32.9.101.311&909&4976& & &32.27.5.41.53.439&2173&2160 \cr
\noalign{\hrule}
 & &9.31.97.887&1947&1060& & &5.121.41.997&493&504 \cr
1739&24004881&8.27.5.11.53.59&97&38&1757&24730585&16.9.5.7.11.17.29.41&1843&248 \cr
 & &32.11.19.53.97&583&304& & &256.3.7.19.31.97&12901&11904 \cr
\noalign{\hrule}
 & &9.7.17.29.773&565&208& & &3.5.7.11.13.17.97&509&800 \cr
1740&24008607&32.3.5.13.29.113&517&952&1758&24759735&64.125.13.509&1067&558 \cr
 & &512.7.11.17.47&517&256& & &256.9.11.31.97&93&128 \cr
\noalign{\hrule}
 & &3.5.7.17.103.131&317&198& & &81.5.17.59.61&37&4978 \cr
1741&24085005&4.27.11.131.317&25&3512&1759&24779115&4.19.37.131&351&352 \cr
 & &64.25.439&2195&32& & &256.27.11.13.131&1703&1408 \cr
\noalign{\hrule}
 & &5.343.13.23.47&6169&7884& & &81.25.13.23.41&3553&3922 \cr
1742&24100895&8.27.31.73.199&833&1430&1760&24824475&4.9.11.17.19.37.53&5&328 \cr
 & &32.9.5.49.11.13.17&153&176& & &64.5.11.41.53&53&352 \cr
\noalign{\hrule}
 & &3.7.31.137.271&13081&12122& & &3.5.13.29.53.83&157&92 \cr
1743&24169677&4.11.19.29.103.127&21235&19278&1761&24876345&8.23.29.53.157&1037&2574 \cr
 & &16.81.5.7.17.31.137&135&136& & &32.9.11.13.17.61&561&976 \cr
\noalign{\hrule}
 & &9.5.47.73.157&77&80& & &27.25.17.41.53&511&596 \cr
1744&24240015&32.3.25.7.11.47.73&697&478&1762&24935175&8.5.7.53.73.149&9&44 \cr
 & &128.7.11.17.41.239&28441&28864& & &64.9.11.73.149&1639&2336 \cr
\noalign{\hrule}
 & &3.5.49.19.37.47&4121&6424& & &31.5329.151&5005&324 \cr
1745&24285135&16.11.13.73.317&35&108&1763&24945049&8.81.5.7.11.13&73&62 \cr
 & &128.27.5.7.317&317&576& & &32.3.7.13.31.73&39&112 \cr
\noalign{\hrule}
 & &5.7.13.197.271&603&1958& & &61.449.911&425&486 \cr
1746&24291085&4.9.7.11.67.89&95&172&1764&24951379&4.243.25.17.449&427&22 \cr
 & &32.3.5.19.43.67&2451&1072& & &16.3.5.7.11.17.61&231&680 \cr
\noalign{\hrule}
}%
}
$$
\eject
\vglue -23 pt
\noindent\hskip 1 in\hbox to 6.5 in{\ 1765 -- 1800 \hfill\fbd 24960851 -- 26190045\frb}
\vskip -9 pt
$$
\vbox{
\nointerlineskip
\halign{\strut
    \vrule \ \ \hfil \frb #\ 
   &\vrule \hfil \ \ \fbb #\frb\ 
   &\vrule \hfil \ \ \frb #\ \hfil
   &\vrule \hfil \ \ \frb #\ 
   &\vrule \hfil \ \ \frb #\ \ \vrule \hskip 2 pt
   &\vrule \ \ \hfil \frb #\ 
   &\vrule \hfil \ \ \fbb #\frb\ 
   &\vrule \hfil \ \ \frb #\ \hfil
   &\vrule \hfil \ \ \frb #\ 
   &\vrule \hfil \ \ \frb #\ \vrule \cr%
\noalign{\hrule}
 & &19.29.89.509&22375&22926& & &9.11.19.107.127&377&250 \cr
1765&24960851&4.3.125.179.3821&1463&2358&1783&25560909&4.3.125.13.29.107&323&2 \cr
 & &16.27.25.7.11.19.131&7425&7336& & &16.5.17.19.29&17&1160 \cr
\noalign{\hrule}
 & &3.25.11.157.193&669&1454& & &9.11.19.107.127&7303&5270 \cr
1766&24998325&4.9.5.223.727&341&386&1784&25560909&4.5.17.31.67.109&17927&17382 \cr
 & &16.11.31.193.223&223&248& & &16.3.7.13.197.2897&20279&20488 \cr
\noalign{\hrule}
 & &5.49.19.41.131&1903&2682& & &5.7.11.13.19.269&2097&862 \cr
1767&25002005&4.9.7.11.149.173&781&262&1785&25580555&4.9.7.233.431&2231&2662 \cr
 & &16.3.121.71.131&363&568& & &16.3.1331.23.97&2783&2328 \cr
\noalign{\hrule}
 & &9.5.13.61.701&5929&3184& & &13.19.173.599&5537&2250 \cr
1768&25015185&32.49.121.199&61&138&1786&25595869&4.9.125.49.113&209&356 \cr
 & &128.3.7.11.23.61&253&448& & &32.3.25.11.19.89&825&1424 \cr
\noalign{\hrule}
 & &25.7.11.13.17.59&589&414& & &5.7.11.19.31.113&851&2094 \cr
1769&25100075&4.9.11.13.19.23.31&85&124&1787&25624445&4.3.7.23.37.349&99&62 \cr
 & &32.3.5.17.23.961&961&1104& & &16.27.11.31.349&349&216 \cr
\noalign{\hrule}
 & &5.7.11.37.41.43&329&144& & &49.37.67.211&7975&6162 \cr
1770&25113935&32.9.49.41.47&37&86&1788&25630381&4.3.25.11.13.29.79&2271&1876 \cr
 & &128.3.37.43.47&141&64& & &32.9.5.7.67.757&757&720 \cr
\noalign{\hrule}
 & &3.13.41.79.199&1719&1520& & &9.25.157.727&341&386 \cr
1771&25137879&32.27.5.13.19.191&1639&1990&1789&25681275&4.5.11.31.157.193&669&1454 \cr
 & &128.25.11.149.199&1639&1600& & &16.3.31.223.727&223&248 \cr
\noalign{\hrule}
 & &11.13.19.37.251&1371&1390& & &3.5.7.19.37.349&657&638 \cr
1772&25232779&4.3.5.13.37.139.457&399&5542&1790&25761435&4.27.11.29.73.349&323&26 \cr
 & &16.9.5.7.17.19.163&2771&2520& & &16.13.17.19.29.73&2117&1768 \cr
\noalign{\hrule}
 & &25.13.17.23.199&33&358& & &9.43.61.1093&455&638 \cr
1773&25287925&4.3.11.179.199&885&1084&1791&25802451&4.3.5.7.11.13.29.43&1093&1702 \cr
 & &32.9.5.59.271&2439&944& & &16.11.23.37.1093&253&296 \cr
\noalign{\hrule}
 & &5.11.31.37.401&263&78& & &7.11.13.17.37.41&97&90 \cr
1774&25297085&4.3.13.263.401&2475&2738&1792&25814789&4.9.5.13.37.41.97&11&544 \cr
 & &16.27.25.11.1369&185&216& & &256.3.11.17.97&97&384 \cr
\noalign{\hrule}
 & &3.25.7.11.13.337&199&186& & &27.25.49.11.71&689&1228 \cr
1775&25300275&4.9.5.31.199.337&53&1738&1793&25831575&8.25.13.53.307&9&316 \cr
 & &16.11.31.53.79&1643&632& & &64.9.53.79&79&1696 \cr
\noalign{\hrule}
 & &81.5.7.169.53&341&712& & &3.7.13.31.43.71&915&418 \cr
1776&25393095&16.5.11.13.31.89&1431&274&1794&25837539&4.9.5.11.13.19.61&31&86 \cr
 & &64.27.53.137&137&32& & &16.19.31.43.61&61&152 \cr
\noalign{\hrule}
 & &3.11.29.101.263&1333&1596& & &27.125.13.19.31&583&1172 \cr
1777&25420791&8.9.7.11.19.31.43&59&40&1795&25842375&8.25.11.53.293&9&284 \cr
 & &128.5.7.31.43.59&12803&13760& & &64.9.53.71&71&1696 \cr
\noalign{\hrule}
 & &3.25.49.11.17.37&11023&9802& & &3.5.23.31.41.59&5411&1846 \cr
1778&25427325&4.169.29.73.151&1715&2664&1796&25871205&4.7.13.71.773&341&432 \cr
 & &64.9.5.343.13.37&91&96& & &128.27.11.31.71&781&576 \cr
\noalign{\hrule}
 & &169.151.997&583&414& & &27.41.97.241&1189&1430 \cr
1779&25442443&4.9.11.23.53.151&65&518&1797&25878339&4.5.11.13.29.1681&10767&11086 \cr
 & &16.3.5.7.13.23.37&483&1480& & &16.3.5.23.37.97.241&185&184 \cr
\noalign{\hrule}
 & &9.125.1331.17&5947&6032& & &5.7.121.17.361&1095&962 \cr
1780&25455375&32.25.13.19.29.313&319&6&1798&25990195&4.3.25.13.19.37.73&231&1156 \cr
 & &128.3.11.19.841&841&1216& & &32.9.7.11.13.289&153&208 \cr
\noalign{\hrule}
 & &9.5.11.13.17.233&1159&938& & &3.5.13.31.61.71&199&204 \cr
1781&25489035&4.5.7.11.19.61.67&699&1436&1799&26180895&8.9.17.61.71.199&11501&638 \cr
 & &32.3.19.233.359&359&304& & &32.7.11.29.31.53&1537&1232 \cr
\noalign{\hrule}
 & &3.5.7.29.83.101&549&466& & &9.5.7.29.47.61&341&646 \cr
1782&25526235&4.27.61.101.233&65&6226&1800&26190045&4.3.11.17.19.29.31&137&50 \cr
 & &16.5.11.13.283&143&2264& & &16.25.19.31.137&2603&1240 \cr
\noalign{\hrule}
}%
}
$$
\eject
\vglue -23 pt
\noindent\hskip 1 in\hbox to 6.5 in{\ 1801 -- 1836 \hfill\fbd 26276373 -- 27461511\frb}
\vskip -9 pt
$$
\vbox{
\nointerlineskip
\halign{\strut
    \vrule \ \ \hfil \frb #\ 
   &\vrule \hfil \ \ \fbb #\frb\ 
   &\vrule \hfil \ \ \frb #\ \hfil
   &\vrule \hfil \ \ \frb #\ 
   &\vrule \hfil \ \ \frb #\ \ \vrule \hskip 2 pt
   &\vrule \ \ \hfil \frb #\ 
   &\vrule \hfil \ \ \fbb #\frb\ 
   &\vrule \hfil \ \ \frb #\ \hfil
   &\vrule \hfil \ \ \frb #\ 
   &\vrule \hfil \ \ \frb #\ \vrule \cr%
\noalign{\hrule}
 & &27.17.19.23.131&77&94& & &9.53.71.797&425&372 \cr
1801&26276373&4.3.7.11.23.47.131&323&70&1819&26991999&8.27.25.17.31.71&671&536 \cr
 & &16.5.49.17.19.47&245&376& & &128.5.11.31.61.67&20801&21440 \cr
\noalign{\hrule}
 & &9.5.11.13.61.67&893&22& & &3.25.43.83.101&539&536 \cr
1802&26299845&4.3.121.19.47&1079&1220&1820&27035175&16.49.11.67.83.101&6579&188 \cr
 & &32.5.13.61.83&83&16& & &128.9.7.17.43.47&987&1088 \cr
\noalign{\hrule}
 & &3.169.23.37.61&73&110& & &9.7.13.19.37.47&5&14 \cr
1803&26318877&4.5.11.169.23.73&549&296&1821&27060579&4.5.49.13.37.47&551&1188 \cr
 & &64.9.37.61.73&73&96& & &32.27.5.11.19.29&435&176 \cr
\noalign{\hrule}
 & &3.25.49.13.19.29&4799&5176& & &3.5.11.13.73.173&181&38 \cr
1804&26324025&16.7.647.4799&135&4664&1822&27089205&4.5.19.173.181&1287&2152 \cr
 & &256.27.5.11.53&477&1408& & &64.9.11.13.269&269&96 \cr
\noalign{\hrule}
 & &9.5.37.97.163&455&418& & &27.5.169.29.41&99&70 \cr
1805&26325315&4.25.7.11.13.19.163&97&2022&1823&27127035&4.243.25.7.11.41&1363&338 \cr
 & &16.3.19.97.337&337&152& & &16.11.169.29.47&47&88 \cr
\noalign{\hrule}
 & &5.7.11.37.1849&943&906& & &9.25.11.79.139&779&1196 \cr
1806&26339005&4.3.5.7.11.23.41.151&1161&104&1824&27177975&8.3.11.13.19.23.41&139&70 \cr
 & &64.81.13.41.43&1053&1312& & &32.5.7.13.41.139&287&208 \cr
\noalign{\hrule}
 & &49.11.17.43.67&333&806& & &49.11.109.463&215&324 \cr
1807&26398603&4.9.49.13.31.37&67&80&1825&27201713&8.81.5.43.463&1199&736 \cr
 & &128.3.5.31.37.67&1147&960& & &512.9.11.23.109&207&256 \cr
\noalign{\hrule}
 & &3.5.7.17.83.179&795&616& & &3.31.227.1289&629&660 \cr
1808&26519745&16.9.25.49.11.53&83&358&1826&27212079&8.9.5.11.17.37.227&53&280 \cr
 & &64.53.83.179&53&32& & &128.25.7.11.17.53&9911&11200 \cr
\noalign{\hrule}
 & &243.5.13.23.73&133&110& & &49.19.23.31.41&327&110 \cr
1809&26519805&4.25.7.11.13.19.73&709&5934&1827&27215923&4.3.5.7.11.41.109&3&38 \cr
 & &16.3.23.43.709&709&344& & &16.9.11.19.109&1199&72 \cr
\noalign{\hrule}
 & &11.13.29.43.149&85&234& & &81.5.17.59.67&5221&206 \cr
1810&26569829&4.9.5.169.17.43&149&20&1828&27216405&4.23.103.227&1071&1298 \cr
 & &32.3.25.17.149&425&48& & &16.9.7.11.17.59&11&56 \cr
\noalign{\hrule}
 & &25.49.109.199&2869&2106& & &3.13.37.113.167&253&86 \cr
1811&26571475&4.81.7.13.19.151&583&1640&1829&27230853&4.11.13.23.37.43&2151&3140 \cr
 & &64.9.5.11.41.53&2173&3168& & &32.9.5.157.239&2355&3824 \cr
\noalign{\hrule}
 & &9.5.7.19.61.73&185&242& & &3.25.13.17.31.53&1021&1046 \cr
1812&26651205&4.3.25.121.37.73&247&28&1830&27232725&4.17.31.523.1021&16785&572 \cr
 & &32.7.11.13.19.37&481&176& & &32.9.5.11.13.373&373&528 \cr
\noalign{\hrule}
 & &243.11.169.59&3649&6322& & &7.11.41.89.97&565&114 \cr
1813&26652483&4.29.41.89.109&649&540&1831&27254381&4.3.5.19.89.113&789&902 \cr
 & &32.27.5.11.59.89&89&80& & &16.9.5.11.41.263&263&360 \cr
\noalign{\hrule}
 & &9.7.11.361.107&65&296& & &7.19.29.37.191&603&470 \cr
1814&26768511&16.3.5.13.37.107&973&418&1832&27257419&4.9.5.47.67.191&407&598 \cr
 & &64.7.11.19.139&139&32& & &16.3.11.13.23.37.47&1833&2024 \cr
\noalign{\hrule}
 & &27.5.11.13.19.73&77&142& & &25.169.29.223&473&642 \cr
1815&26776035&4.9.7.121.19.71&25&146&1833&27323075&4.3.5.11.29.43.107&67&78 \cr
 & &16.25.7.71.73&71&280& & &16.9.13.43.67.107&4601&4824 \cr
\noalign{\hrule}
 & &27.17.19.37.83&239&220& & &3.25.11.89.373&337&782 \cr
1816&26782191&8.5.11.37.83.239&221&2850&1834&27387525&4.5.11.17.23.337&749&936 \cr
 & &32.3.125.13.17.19&125&208& & &64.9.7.13.23.107&4173&5152 \cr
\noalign{\hrule}
 & &3.5.41.193.227&35&158& & &9.11.13.83.257&79&178 \cr
1817&26943765&4.25.7.79.227&193&1782&1835&27452997&4.13.79.83.89&537&620 \cr
 & &16.81.11.193&11&216& & &32.3.5.31.79.179&5549&6320 \cr
\noalign{\hrule}
 & &27.5.121.13.127&523&874& & &81.49.11.17.37&173&184 \cr
1818&26969085&4.5.11.19.23.523&47&162&1836&27461511&16.27.7.23.37.173&331&520 \cr
 & &16.81.47.523&523&1128& & &256.5.13.173.331&21515&22144 \cr
\noalign{\hrule}
}%
}
$$
\eject
\vglue -23 pt
\noindent\hskip 1 in\hbox to 6.5 in{\ 1837 -- 1872 \hfill\fbd 27462435 -- 29171303\frb}
\vskip -9 pt
$$
\vbox{
\nointerlineskip
\halign{\strut
    \vrule \ \ \hfil \frb #\ 
   &\vrule \hfil \ \ \fbb #\frb\ 
   &\vrule \hfil \ \ \frb #\ \hfil
   &\vrule \hfil \ \ \frb #\ 
   &\vrule \hfil \ \ \frb #\ \ \vrule \hskip 2 pt
   &\vrule \ \ \hfil \frb #\ 
   &\vrule \hfil \ \ \fbb #\frb\ 
   &\vrule \hfil \ \ \frb #\ \hfil
   &\vrule \hfil \ \ \frb #\ 
   &\vrule \hfil \ \ \frb #\ \vrule \cr%
\noalign{\hrule}
 & &3.5.7.11.13.31.59&203&262& & &3.49.11.13.17.79&183&370 \cr
1837&27462435&4.49.11.13.29.131&81&458&1855&28231203&4.9.5.7.13.37.61&79&506 \cr
 & &16.81.131.229&3537&1832& & &16.11.23.37.79&37&184 \cr
\noalign{\hrule}
 & &9.11.23.47.257&3451&6278& & &54289.521&27405&26884 \cr
1838&27503883&4.7.17.29.43.73&771&470&1856&28284569&8.27.5.7.11.13.29.47&233&2 \cr
 & &16.3.5.29.47.257&29&40& & &32.9.13.29.233&377&144 \cr
\noalign{\hrule}
 & &17.19.31.41.67&1&1272& & &3.49.103.1871&96283&96430 \cr
1839&27505711&16.3.17.53&451&450&1857&28328811&4.5.11.8753.9643&445&9198 \cr
 & &64.27.25.11.41&275&864& & &16.9.25.7.11.73.89&6497&6600 \cr
\noalign{\hrule}
 & &9.121.13.29.67&139&238& & &9.5.7.11.13.17.37&303&292 \cr
1840&27507051&4.7.11.17.67.139&195&1334&1858&28333305&8.27.13.37.73.101&157&1156 \cr
 & &16.3.5.7.13.23.29&161&40& & &64.289.73.157&2669&2336 \cr
\noalign{\hrule}
 & &9.49.11.13.19.23&1255&1462& & &49.11.23.29.79&2025&208 \cr
1841&27558531&4.5.49.17.43.251&3&46&1859&28401527&32.81.25.7.13&67&158 \cr
 & &16.3.5.17.23.251&251&680& & &128.9.67.79&603&64 \cr
\noalign{\hrule}
 & &3.5.7.13.17.29.41&13&132& & &9.7.37.73.167&55&18 \cr
1842&27590745&8.9.11.169.41&535&986&1860&28417221&4.81.5.7.11.167&1643&2812 \cr
 & &32.5.17.29.107&107&16& & &32.19.31.37.53&1007&496 \cr
\noalign{\hrule}
 & &243.25.7.11.59&8131&6206& & &1331.13.31.53&493&90 \cr
1843&27598725&4.29.47.107.173&77&30&1861&28428829&4.9.5.121.17.29&53&68 \cr
 & &16.3.5.7.11.29.173&173&232& & &32.3.289.29.53&867&464 \cr
\noalign{\hrule}
 & &7.11.29.83.149&5067&7300& & &3.5.13.211.691&2387&5842 \cr
1844&27615511&8.9.25.73.563&29&44&1862&28431195&4.7.11.23.31.127&109&108 \cr
 & &64.3.5.11.29.563&563&480& & &32.27.11.23.109.127&22563&22352 \cr
\noalign{\hrule}
 & &81.25.11.17.73&161&26& & &9.5.13.19.31.83&1189&1384 \cr
1845&27643275&4.3.5.7.13.23.73&77&142&1863&28598895&16.3.19.29.41.173&649&130 \cr
 & &16.49.11.23.71&1633&392& & &64.5.11.13.29.59&649&928 \cr
\noalign{\hrule}
 & &5.289.127.151&33&118& & &3.25.11.13.17.157&8033&5992 \cr
1846&27710765&4.3.11.17.59.127&755&1404&1864&28625025&16.7.29.107.277&1413&1690 \cr
 & &32.81.5.13.151&81&208& & &64.9.5.7.169.157&91&96 \cr
\noalign{\hrule}
 & &5.11.169.29.103&981&878& & &9.7.17.73.367&1905&664 \cr
1847&27764165&4.9.5.29.109.439&483&2678&1865&28693161&16.27.5.83.127&803&1438 \cr
 & &16.27.7.13.23.103&161&216& & &64.11.73.719&719&352 \cr
\noalign{\hrule}
 & &3.13.31.83.277&145&132& & &3.13.19.47.827&33&860 \cr
1848&27796119&8.9.5.11.29.31.83&823&76&1866&28801929&8.9.5.11.13.43&49&94 \cr
 & &64.5.11.19.823&4115&6688& & &32.49.43.47&43&784 \cr
\noalign{\hrule}
 & &27.5.19.61.179&517&338& & &3.7.11.29.59.73&95&124 \cr
1849&28007235&4.3.11.169.47.61&589&82&1867&28852593&8.5.7.11.19.31.59&377&36 \cr
 & &16.19.31.41.47&1271&376& & &64.9.5.13.19.29&195&608 \cr
\noalign{\hrule}
 & &3.121.113.683&963&280& & &9.5.13.361.137&913&868 \cr
1850&28015977&16.27.5.7.11.107&523&226&1868&28932345&8.7.11.361.31.83&351&10 \cr
 & &64.5.113.523&523&160& & &32.27.5.7.13.83&249&112 \cr
\noalign{\hrule}
 & &3.5.7.11.83.293&603&310& & &9.25.11.13.17.53&139&86 \cr
1851&28088445&4.27.25.7.31.67&103&572&1869&28989675&4.11.13.17.43.139&373&186 \cr
 & &32.11.13.31.103&403&1648& & &16.3.31.139.373&4309&2984 \cr
\noalign{\hrule}
 & &9.13.31.61.127&1771&120& & &9.5.11.961.61&755&206 \cr
1852&28098369&16.27.5.7.11.23&113&140&1870&29017395&4.25.11.103.151&213&62 \cr
 & &128.25.49.113&5537&1600& & &16.3.31.71.103&103&568 \cr
\noalign{\hrule}
 & &9.49.23.47.59&5&418& & &3.7.11.191.661&235&426 \cr
1853&28126539&4.5.7.11.19.23&181&204&1871&29163981&4.9.5.7.11.47.71&2015&1322 \cr
 & &32.3.17.19.181&3077&304& & &16.25.13.31.661&325&248 \cr
\noalign{\hrule}
 & &27.5.13.17.23.41&1073&682& & &7.17.29.79.107&2223&880 \cr
1854&28134405&4.11.29.31.37.41&391&798&1872&29171303&32.9.5.7.11.13.19&69&316 \cr
 & &16.3.7.17.19.23.31&217&152& & &256.27.23.79&621&128 \cr
\noalign{\hrule}
}%
}
$$
\eject
\vglue -23 pt
\noindent\hskip 1 in\hbox to 6.5 in{\ 1873 -- 1908 \hfill\fbd 29174145 -- 30514133\frb}
\vskip -9 pt
$$
\vbox{
\nointerlineskip
\halign{\strut
    \vrule \ \ \hfil \frb #\ 
   &\vrule \hfil \ \ \fbb #\frb\ 
   &\vrule \hfil \ \ \frb #\ \hfil
   &\vrule \hfil \ \ \frb #\ 
   &\vrule \hfil \ \ \frb #\ \ \vrule \hskip 2 pt
   &\vrule \ \ \hfil \frb #\ 
   &\vrule \hfil \ \ \fbb #\frb\ 
   &\vrule \hfil \ \ \frb #\ \hfil
   &\vrule \hfil \ \ \frb #\ 
   &\vrule \hfil \ \ \frb #\ \vrule \cr%
\noalign{\hrule}
 & &3.5.7.11.13.29.67&43&188& & &9.11.13.19.23.53&235&64 \cr
1873&29174145&8.13.43.47.67&1135&1746&1891&29808207&128.5.11.47.53&285&232 \cr
 & &32.9.5.97.227&681&1552& & &2048.3.25.19.29&725&1024 \cr
\noalign{\hrule}
 & &27.125.17.509&299&1826& & &9.5.7.17.19.293&389&682 \cr
1874&29203875&4.9.11.13.23.83&185&68&1892&29811285&4.5.11.19.31.389&667&3612 \cr
 & &32.5.17.37.83&83&592& & &32.3.7.23.29.43&1247&368 \cr
\noalign{\hrule}
 & &5.11.37.113.127&1081&954& & &3.5.49.41.991&639&352 \cr
1875&29204285&4.9.23.47.53.113&91&22&1893&29863785&64.27.5.7.11.71&863&82 \cr
 & &16.3.7.11.13.47.53&2067&2632& & &256.41.863&863&128 \cr
\noalign{\hrule}
 & &25.7.47.53.67&9&44& & &27.5.7.31.1021&377&5482 \cr
1876&29206975&8.9.5.11.47.67&943&608&1894&29910195&4.13.29.2741&1377&1364 \cr
 & &512.3.19.23.41&2337&5888& & &32.81.11.17.29.31&493&528 \cr
\noalign{\hrule}
 & &3.5.11.29.31.197&243&98& & &81.7.121.19.23&167&86 \cr
1877&29221995&4.729.49.197&325&1054&1895&29981259&4.7.11.19.43.167&255&46 \cr
 & &16.25.7.13.17.31&455&136& & &16.3.5.17.23.167&167&680 \cr
\noalign{\hrule}
 & &9.11.23.71.181&8143&9776& & &81.5.7.71.149&61&88 \cr
1878&29261727&32.13.17.47.479&545&66&1896&29991465&16.3.5.7.11.61.71&209&706 \cr
 & &128.3.5.11.17.109&545&1088& & &64.121.19.353&6707&3872 \cr
\noalign{\hrule}
 & &9.11.23.79.163&45&208& & &121.29.83.103&405&508 \cr
1879&29320929&32.81.5.13.79&163&242&1897&29998441&8.81.5.11.29.127&103&158 \cr
 & &128.121.13.163&143&64& & &32.9.79.103.127&1143&1264 \cr
\noalign{\hrule}
 & &5.121.89.547&763&216& & &27.125.7.31.41&2761&6136 \cr
1880&29453215&16.27.5.7.11.109&29&356&1898&30027375&16.11.13.59.251&199&450 \cr
 & &128.9.29.89&29&576& & &64.9.25.13.199&199&416 \cr
\noalign{\hrule}
 & &9.19.23.59.127&5713&1780& & &81.13.17.23.73&7135&8998 \cr
1881&29469969&8.5.29.89.197&143&54&1899&30055779&4.5.11.409.1427&509&918 \cr
 & &32.27.5.11.13.29&1885&528& & &16.27.5.11.17.509&509&440 \cr
\noalign{\hrule}
 & &27.5.13.17.23.43&341&556& & &13.61.83.457&187&270 \cr
1882&29506815&8.9.11.17.31.139&7&146&1900&30079283&4.27.5.11.13.17.61&125&1162 \cr
 & &32.7.11.31.73&511&5456& & &16.3.625.7.83&625&168 \cr
\noalign{\hrule}
 & &3.125.11.13.19.29&241&716& & &25.49.79.311&121&432 \cr
1883&29547375&8.5.13.179.241&153&88&1901&30097025&32.27.25.7.121&457&632 \cr
 & &128.9.11.17.179&537&1088& & &512.3.79.457&457&768 \cr
\noalign{\hrule}
 & &27.11.169.19.31&5&346& & &7.11.59.61.109&159&268 \cr
1884&29563677&4.5.13.19.173&123&124&1902&30206407&8.3.11.53.59.67&35&24 \cr
 & &32.3.5.31.41.173&865&656& & &128.9.5.7.53.67&3551&2880 \cr
\noalign{\hrule}
 & &9.5.11.13.43.107&911&266& & &3.7.23.31.43.47&41&88 \cr
1885&29607435&4.3.7.13.19.911&47&86&1903&30260433&16.7.11.23.31.41&1935&1222 \cr
 & &16.43.47.911&911&376& & &64.9.5.13.43.47&65&96 \cr
\noalign{\hrule}
 & &9.23.37.53.73&895&784& & &9.25.11.13.23.41&501&524 \cr
1886&29632671&32.3.5.49.53.179&3053&3212&1904&30341025&8.27.11.13.131.167&703&1000 \cr
 & &256.7.11.43.71.73&5467&5504& & &128.125.19.37.167&6179&6080 \cr
\noalign{\hrule}
 & &5.11.29.43.433&57&376& & &25.11.31.43.83&639&694 \cr
1887&29697305&16.3.5.19.43.47&69&26&1905&30425725&4.9.5.71.83.347&5549&344 \cr
 & &64.9.13.23.47&5499&736& & &64.3.31.43.179&179&96 \cr
\noalign{\hrule}
 & &9.5.49.13.17.61&79&506& & &27.5.7.11.29.101&745&38 \cr
1888&29725605&4.7.11.17.23.79&183&370&1906&30446955&4.25.11.19.149&1057&582 \cr
 & &16.3.5.23.37.61&37&184& & &16.3.7.97.151&97&1208 \cr
\noalign{\hrule}
 & &17.37.151.313&57&94& & &9.125.11.23.107&37&62 \cr
1889&29728427&4.3.17.19.47.313&165&148&1907&30454875&4.5.23.31.37.107&89&624 \cr
 & &32.9.5.11.19.37.47&2585&2736& & &128.3.13.37.89&3293&832 \cr
\noalign{\hrule}
 & &9.7.17.37.751&3509&3250& & &2197.17.19.43&1507&690 \cr
1890&29759877&4.125.121.13.17.29&1759&1266&1908&30514133&4.3.5.11.17.23.137&949&558 \cr
 & &16.3.5.13.211.1759&13715&14072& & &16.27.5.13.31.73&1971&1240 \cr
\noalign{\hrule}
}%
}
$$
\eject
\vglue -23 pt
\noindent\hskip 1 in\hbox to 6.5 in{\ 1909 -- 1944 \hfill\fbd 30533503 -- 31943795\frb}
\vskip -9 pt
$$
\vbox{
\nointerlineskip
\halign{\strut
    \vrule \ \ \hfil \frb #\ 
   &\vrule \hfil \ \ \fbb #\frb\ 
   &\vrule \hfil \ \ \frb #\ \hfil
   &\vrule \hfil \ \ \frb #\ 
   &\vrule \hfil \ \ \frb #\ \ \vrule \hskip 2 pt
   &\vrule \ \ \hfil \frb #\ 
   &\vrule \hfil \ \ \fbb #\frb\ 
   &\vrule \hfil \ \ \frb #\ \hfil
   &\vrule \hfil \ \ \frb #\ 
   &\vrule \hfil \ \ \frb #\ \vrule \cr%
\noalign{\hrule}
 & &7.121.13.47.59&2173&600& & &3.5.7.13.17.19.71&367&296 \cr
1909&30533503&16.3.25.7.41.53&183&188&1927&31303545&16.5.7.19.37.367&473&3042 \cr
 & &128.9.5.41.47.61&2745&2624& & &64.9.11.169.43&559&1056 \cr
\noalign{\hrule}
 & &5.11.37.43.349&621&970& & &3.25.7.169.353&33&58 \cr
1910&30539245&4.27.25.11.23.97&329&1396&1928&31319925&4.9.11.13.29.353&245&3638 \cr
 & &32.9.7.47.349&329&144& & &16.5.49.17.107&107&952 \cr
\noalign{\hrule}
 & &9.25.7.11.41.43&71&544& & &3.5.193.10831&5319&5512 \cr
1911&30543975&64.3.5.7.17.71&551&656&1929&31355745&16.81.5.13.53.197&5423&5018 \cr
 & &2048.19.29.41&551&1024& & &64.11.169.17.29.193&5423&5408 \cr
\noalign{\hrule}
 & &5.11.29.71.271&1455&1526& & &9.5.13.17.29.109&1643&1672 \cr
1912&30689395&4.3.25.7.29.97.109&2769&44&1930&31436145&16.3.11.19.31.53.109&667&340 \cr
 & &32.9.7.11.13.71&63&208& & &128.5.11.17.23.29.31&713&704 \cr
\noalign{\hrule}
 & &9.107.127.251&55&52& & &3.289.41.887&877&10 \cr
1913&30697551&8.3.5.11.13.127.251&1075&322&1931&31530189&4.5.41.877&459&418 \cr
 & &32.125.7.13.23.43&12857&14000& & &16.27.5.11.17.19&95&792 \cr
\noalign{\hrule}
 & &9.5.49.11.31.41&109&232& & &25.11.131.877&6461&3186 \cr
1914&30828105&16.3.5.49.29.109&79&166&1932&31593925&4.27.7.13.59.71&373&550 \cr
 & &64.79.83.109&6557&3488& & &16.9.25.7.11.373&373&504 \cr
\noalign{\hrule}
 & &9.11.13.289.83&95&194& & &25.13.31.43.73&1771&1368 \cr
1915&30871269&4.5.13.19.83.97&297&782&1933&31625425&16.9.25.7.11.19.23&353&122 \cr
 & &16.27.11.17.19.23&69&152& & &64.3.23.61.353&8119&5856 \cr
\noalign{\hrule}
 & &9.7.11.13.23.149&29&178& & &9.5.7.11.13.19.37&1009&454 \cr
1916&30873843&4.7.11.13.29.89&535&444&1934&31666635&4.3.13.227.1009&2989&38 \cr
 & &32.3.5.29.37.107&3959&2320& & &16.49.19.61&7&488 \cr
\noalign{\hrule}
 & &3.49.11.169.113&591&592& & &3.25.11.19.43.47&821&4 \cr
1917&30879849&32.9.7.11.37.113.197&23951&1690&1935&31679175&8.47.821&387&434 \cr
 & &128.5.169.43.557&2785&2752& & &32.9.7.31.43&21&496 \cr
\noalign{\hrule}
 & &13.17.29.61.79&7845&9614& & &3.25.7.11.289.19&461&406 \cr
1918&30884971&4.3.5.11.19.23.523&157&366&1936&31710525&4.5.49.19.29.461&153&398 \cr
 & &16.9.5.23.61.157&1035&1256& & &16.9.17.199.461&1383&1592 \cr
\noalign{\hrule}
 & &3.5.73.89.317&871&5626& & &25.11.31.3721&46683&46342 \cr
1919&30893235&4.13.29.67.97&341&1602&1937&31721525&4.27.7.13.17.19.29.47&305&682 \cr
 & &16.9.11.31.89&341&24& & &16.9.5.11.17.19.31.61&153&152 \cr
\noalign{\hrule}
 & &5.7.11.13.23.269&81&172& & &9.17.19.67.163&325&814 \cr
1920&30965935&8.81.5.43.269&833&1102&1938&31747347&4.3.25.11.13.19.37&335&368 \cr
 & &32.9.49.17.19.29&3451&2736& & &128.125.13.23.67&1625&1472 \cr
\noalign{\hrule}
 & &3.25.49.11.13.59&183&92& & &11.169.19.29.31&379&210 \cr
1921&31005975&8.9.7.23.59.61&1243&2600&1939&31753579&4.3.5.7.11.29.379&67&78 \cr
 & &128.25.11.13.113&113&64& & &16.9.7.13.67.379&3411&3752 \cr
\noalign{\hrule}
 & &25.13.17.43.131&693&38& & &9.7.31.43.379&173&44 \cr
1922&31122325&4.9.5.7.11.13.19&131&116&1940&31828041&8.3.11.173.379&1825&2344 \cr
 & &32.3.7.11.29.131&203&528& & &128.25.73.293&7325&4672 \cr
\noalign{\hrule}
 & &25.13.17.43.131&693&38& & &13.359.6823&5745&1078 \cr
1923&31122325&4.9.5.7.11.13.19&2227&2162&1941&31842941&4.3.5.49.11.383&1313&1368 \cr
 & &16.3.17.23.47.131&141&184& & &64.27.7.13.19.101&3591&3232 \cr
\noalign{\hrule}
 & &5.7.13.17.37.109&33&152& & &9.23.149.1033&14905&15938 \cr
1924&31195255&16.3.11.13.19.109&17&126&1942&31860819&4.5.11.13.271.613&5133&1610 \cr
 & &64.27.7.17.19&513&32& & &16.3.25.7.23.29.59&1475&1624 \cr
\noalign{\hrule}
 & &5.7.121.53.139&123&262& & &27.7.13.41.317&1133&1450 \cr
1925&31199245&4.3.11.41.53.131&1807&366&1943&31933629&4.3.25.11.13.29.103&317&8 \cr
 & &16.9.13.61.139&117&488& & &64.11.29.317&319&32 \cr
\noalign{\hrule}
 & &9.49.11.47.137&173&250& & &5.13.31.83.191&1449&1034 \cr
1926&31235589&4.125.7.137.173&47&912&1944&31943795&4.9.7.11.23.31.47&3071&1300 \cr
 & &128.3.25.19.47&475&64& & &32.3.25.13.37.83&111&80 \cr
\noalign{\hrule}
}%
}
$$
\eject
\vglue -23 pt
\noindent\hskip 1 in\hbox to 6.5 in{\ 1945 -- 1980 \hfill\fbd 31954027 -- 33342435\frb}
\vskip -9 pt
$$
\vbox{
\nointerlineskip
\halign{\strut
    \vrule \ \ \hfil \frb #\ 
   &\vrule \hfil \ \ \fbb #\frb\ 
   &\vrule \hfil \ \ \frb #\ \hfil
   &\vrule \hfil \ \ \frb #\ 
   &\vrule \hfil \ \ \frb #\ \ \vrule \hskip 2 pt
   &\vrule \ \ \hfil \frb #\ 
   &\vrule \hfil \ \ \fbb #\frb\ 
   &\vrule \hfil \ \ \frb #\ \hfil
   &\vrule \hfil \ \ \frb #\ 
   &\vrule \hfil \ \ \frb #\ \vrule \cr%
\noalign{\hrule}
 & &49.29.113.199&611&810& & &5.343.17.19.59&333&10 \cr
1945&31954027&4.81.5.13.47.113&203&814&1963&32682755&4.9.25.37.59&551&374 \cr
 & &16.9.5.7.11.29.37&495&296& & &16.3.11.17.19.29&319&24 \cr
\noalign{\hrule}
 & &3.7.13.41.47.61&17&30& & &11.13.19.23.523&135&388 \cr
1946&32090331&4.9.5.7.17.41.61&715&1786&1964&32682793&8.27.5.13.19.97&101&146 \cr
 & &16.25.11.13.19.47&275&152& & &32.3.73.97.101&7373&4656 \cr
\noalign{\hrule}
 & &9.5.7.11.13.23.31&65&142& & &9.7.13.23.37.47&18821&19672 \cr
1947&32117085&4.25.169.31.71&3213&1012&1965&32757543&16.11.29.59.2459&905&1554 \cr
 & &32.27.7.11.17.23&17&48& & &64.3.5.7.29.37.181&905&928 \cr
\noalign{\hrule}
 & &17.31.47.1297&11753&10296& & &3.13.31.37.733&165&568 \cr
1948&32125393&16.9.7.11.13.23.73&625&1054&1966&32789289&16.9.5.11.37.71&31&68 \cr
 & &64.3.625.7.17.31&625&672& & &128.5.17.31.71&355&1088 \cr
\noalign{\hrule}
 & &25.121.13.19.43&381&854& & &27.11.17.67.97&1&98 \cr
1949&32128525&4.3.5.7.11.61.127&129&256&1967&32813451&4.3.49.17.67&209&260 \cr
 & &2048.9.43.61&549&1024& & &32.5.7.11.13.19&1729&80 \cr
\noalign{\hrule}
 & &3.49.13.19.887&125&762& & &5.11.13.19.41.59&31&174 \cr
1950&32206083&4.9.125.19.127&553&572&1968&32862115&4.3.19.29.31.59&305&246 \cr
 & &32.7.11.13.79.127&1397&1264& & &16.9.5.31.41.61&279&488 \cr
\noalign{\hrule}
 & &3.25.13.19.37.47&77&64& & &7.13.17.89.239&9&230 \cr
1951&32214975&128.25.7.11.19.37&611&1314&1969&32906237&4.9.5.7.23.89&923&946 \cr
 & &512.9.13.47.73&219&256& & &16.3.5.11.13.43.71&2343&1720 \cr
\noalign{\hrule}
 & &81.25.19.839&657&182& & &9.5.49.11.23.59&199&494 \cr
1952&32280525&4.729.7.13.73&839&110&1970&32914035&4.7.13.19.23.199&43&204 \cr
 & &16.5.7.11.839&7&88& & &32.3.17.43.199&3383&688 \cr
\noalign{\hrule}
 & &27.7.17.19.529&65&464& & &9.5.7.11.13.17.43&227&228 \cr
1953&32293863&32.9.5.13.17.29&539&46&1971&32927895&8.27.11.17.19.43.227&7&466 \cr
 & &128.49.11.23&11&448& & &32.7.19.227.233&4313&3728 \cr
\noalign{\hrule}
 & &27.17.97.727&1673&946& & &3.11.31.103.313&6511&3068 \cr
1954&32368221&4.7.11.17.43.239&485&246&1972&32980497&8.13.17.59.383&693&310 \cr
 & &16.3.5.7.11.41.97&287&440& & &32.9.5.7.11.13.31&105&208 \cr
\noalign{\hrule}
 & &11.17.19.23.397&971&5778& & &5.7.31.113.269&411&380 \cr
1955&32442443&4.27.107.971&325&646&1973&32980745&8.3.25.19.137.269&2061&4664 \cr
 & &16.9.25.13.17.19&325&72& & &128.27.11.53.229&15741&14656 \cr
\noalign{\hrule}
 & &5.19.23.83.179&129&308& & &27.7.169.17.61&577&460 \cr
1956&32462545&8.3.5.7.11.43.83&57&358&1974&33122817&8.3.5.7.13.23.577&43&22 \cr
 & &32.9.11.19.179&99&16& & &32.11.23.43.577&10879&9232 \cr
\noalign{\hrule}
 & &3.11.71.83.167&89&78& & &9.5.7.11.61.157&493&920 \cr
1957&32476323&4.9.13.71.83.89&835&88&1975&33184305&16.25.11.17.23.29&471&196 \cr
 & &64.5.11.89.167&89&160& & &128.3.49.17.157&119&64 \cr
\noalign{\hrule}
 & &7.97.151.317&207&110& & &5.61.89.1223&459&764 \cr
1958&32501693&4.9.5.7.11.23.151&317&166&1976&33198335&8.27.17.89.191&1223&2024 \cr
 & &16.3.5.11.83.317&415&264& & &128.3.11.23.1223&253&192 \cr
\noalign{\hrule}
 & &9.5.13.19.29.101&343&242& & &7.11.53.79.103&65&12 \cr
1959&32555835&4.343.121.19.29&237&314&1977&33207097&8.3.5.13.79.103&249&146 \cr
 & &16.3.49.11.79.157&7693&6952& & &32.9.13.73.83&6059&1872 \cr
\noalign{\hrule}
 & &5.7.23.71.571&611&2244& & &3.49.11.19.23.47&117&1010 \cr
1960&32635505&8.3.7.11.13.17.47&71&258&1978&33211563&4.27.5.11.13.101&1435&1292 \cr
 & &32.9.13.43.71&117&688& & &32.25.7.17.19.41&425&656 \cr
\noalign{\hrule}
 & &3.25.13.19.41.43&231&16& & &289.31.47.79&1865&1848 \cr
1961&32659575&32.9.5.7.11.41&107&98&1979&33264767&16.3.5.7.11.17.31.373&2623&12 \cr
 & &128.343.11.107&3773&6848& & &128.9.11.43.61&2623&6336 \cr
\noalign{\hrule}
 & &81.25.13.17.73&77&142& & &81.5.7.19.619&259&254 \cr
1962&32669325&4.27.5.7.11.17.71&161&26&1980&33342435&4.3.49.37.127.619&1835&22 \cr
 & &16.49.13.23.71&1633&392& & &16.5.11.127.367&1397&2936 \cr
\noalign{\hrule}
}%
}
$$
\eject
\vglue -23 pt
\noindent\hskip 1 in\hbox to 6.5 in{\ 1981 -- 2016 \hfill\fbd 33367515 -- 34768767\frb}
\vskip -9 pt
$$
\vbox{
\nointerlineskip
\halign{\strut
    \vrule \ \ \hfil \frb #\ 
   &\vrule \hfil \ \ \fbb #\frb\ 
   &\vrule \hfil \ \ \frb #\ \hfil
   &\vrule \hfil \ \ \frb #\ 
   &\vrule \hfil \ \ \frb #\ \ \vrule \hskip 2 pt
   &\vrule \ \ \hfil \frb #\ 
   &\vrule \hfil \ \ \fbb #\frb\ 
   &\vrule \hfil \ \ \frb #\ \hfil
   &\vrule \hfil \ \ \frb #\ 
   &\vrule \hfil \ \ \frb #\ \vrule \cr%
\noalign{\hrule}
 & &3.5.17.19.71.97&749&458& & &3.11.13.19.53.79&35&22 \cr
1981&33367515&4.5.7.19.107.229&639&506&1999&34128237&4.5.7.121.53.79&3483&2930 \cr
 & &16.9.11.23.71.107&759&856& & &16.81.25.43.293&7911&8600 \cr
\noalign{\hrule}
 & &11.19.59.2711&15471&14350& & &27.17.23.53.61&2165&946 \cr
1982&33429341&4.81.25.7.41.191&1031&76&2000&34130781&4.9.5.11.43.433&221&212 \cr
 & &32.3.5.7.19.1031&1031&1680& & &32.5.11.13.17.43.53&715&688 \cr
\noalign{\hrule}
 & &37.73.79.157&1383&1540& & &3.23.37.59.227&697&660 \cr
1983&33500503&8.3.5.7.11.73.461&413&48&2001&34192329&8.9.5.11.17.41.227&71&298 \cr
 & &256.9.49.11.59&4851&7552& & &32.5.11.17.71.149&13277&11920 \cr
\noalign{\hrule}
 & &13.17.19.79.101&11&90& & &9.5.7.19.59.97&3239&2366 \cr
1984&33503821&4.9.5.11.13.17.19&101&120&2002&34252155&4.49.169.41.79&291&242 \cr
 & &64.27.25.11.101&675&352& & &16.3.121.13.79.97&1027&968 \cr
\noalign{\hrule}
 & &5.7.17.23.31.79&4279&4806& & &27.5.7.19.23.83&4573&4972 \cr
1985&33514565&4.27.7.11.89.389&79&310&2003&34276095&8.9.11.17.113.269&133&20 \cr
 & &16.9.5.31.79.89&89&72& & &64.5.7.11.19.269&269&352 \cr
\noalign{\hrule}
 & &31.61.113.157&8901&8840& & &11.529.71.83&1273&360 \cr
1986&33548231&16.9.5.13.17.23.31.43&791&5874&2004&34291367&16.9.5.19.23.67&721&284 \cr
 & &64.27.7.11.89.113&2403&2464& & &128.3.7.71.103&721&192 \cr
\noalign{\hrule}
 & &3.25.11.23.29.61&471&196& & &25.13.23.43.107&30483&31042 \cr
1987&33566775&8.9.49.61.157&493&920&2005&34392475&4.27.11.17.83.1129&1237&2150 \cr
 & &128.5.7.17.23.29&119&64& & &16.9.25.17.43.1237&1237&1224 \cr
\noalign{\hrule}
 & &5.7.13.17.43.101&29&36& & &3.7.11.137.1087&6065&1544 \cr
1988&33593105&8.9.17.29.43.101&1925&2418&2006&34400289&16.5.193.1213&1089&124 \cr
 & &32.27.25.7.11.13.31&837&880& & &128.9.121.31&1023&64 \cr
\noalign{\hrule}
 & &3.13.23.157.239&85&154& & &125.7.23.29.59&143&18 \cr
1989&33658131&4.5.7.11.13.17.157&621&478&2007&34433875&4.9.11.13.29.59&875&836 \cr
 & &16.27.5.17.23.239&85&72& & &32.3.125.7.121.19&363&304 \cr
\noalign{\hrule}
 & &27.7.17.47.223&1261&2530& & &81.25.7.11.13.17&5699&5776 \cr
1990&33675453&4.5.7.11.13.23.97&489&190&2008&34459425&32.3.13.361.41.139&35&22 \cr
 & &16.3.25.11.19.163&1793&3800& & &128.5.7.11.19.41.139&2641&2624 \cr
\noalign{\hrule}
 & &5.343.13.17.89&3737&3828& & &81.7.31.37.53&121&68 \cr
1991&33732335&8.3.49.11.29.37.101&75&26&2009&34468497&8.3.121.17.31.37&65&28 \cr
 & &32.9.25.11.13.29.37&2871&2960& & &64.5.7.121.13.17&1573&2720 \cr
\noalign{\hrule}
 & &9.11.397.859&2185&2182& & &27.7.19.29.331&5885&404 \cr
1992&33761277&4.3.5.19.23.859.1091&16343&22&2010&34470009&8.5.11.101.107&609&502 \cr
 & &16.11.23.59.277&1357&2216& & &32.3.5.7.29.251&251&80 \cr
\noalign{\hrule}
 & &27.5.11.529.43&413&116& & &27.19.23.29.101&85&16 \cr
1993&33779295&8.5.7.29.43.59&409&1656&2011&34559271&32.9.5.17.19.29&253&298 \cr
 & &128.9.23.409&409&64& & &128.11.17.23.149&1639&1088 \cr
\noalign{\hrule}
 & &5.7.11.19.31.149&5967&1348& & &3.13.59.83.181&61&22 \cr
1994&33787985&8.27.13.17.337&227&110&2012&34567923&4.11.59.61.181&415&234 \cr
 & &32.3.5.11.17.227&227&816& & &16.9.5.13.61.83&61&120 \cr
\noalign{\hrule}
 & &3.11.19.31.37.47&5&52& & &27.5.13.17.19.61&3487&3182 \cr
1995&33800943&8.5.11.13.31.37&189&152&2013&34578765&4.11.17.37.43.317&551&78 \cr
 & &128.27.5.7.13.19&315&832& & &16.3.13.19.29.317&317&232 \cr
\noalign{\hrule}
 & &9.7.11.19.31.83&295&46& & &3.5.11.43.67.73&29&14 \cr
1996&33878691&4.3.5.7.19.23.59&3071&3484&2014&34701645&4.7.11.29.67.73&6399&7202 \cr
 & &32.13.37.67.83&871&592& & &16.81.13.79.277&7479&8216 \cr
\noalign{\hrule}
 & &3.25.7.17.37.103&57&572& & &27.5.361.23.31&539&544 \cr
1997&34013175&8.9.5.7.11.13.19&103&68&2015&34748055&64.9.49.11.17.23.31&379&1898 \cr
 & &64.11.13.17.103&143&32& & &256.13.17.73.379&27667&28288 \cr
\noalign{\hrule}
 & &5.13.71.73.101&209&714& & &3.11.29.47.773&1319&1000 \cr
1998&34026395&4.3.7.11.17.19.73&725&516&2016&34768767&16.125.47.1319&1247&72 \cr
 & &32.9.25.7.29.43&1935&3248& & &256.9.5.29.43&645&128 \cr
\noalign{\hrule}
}%
}
$$
\eject
\vglue -23 pt
\noindent\hskip 1 in\hbox to 6.5 in{\ 2017 -- 2052 \hfill\fbd 34810875 -- 36504657\frb}
\vskip -9 pt
$$
\vbox{
\nointerlineskip
\halign{\strut
    \vrule \ \ \hfil \frb #\ 
   &\vrule \hfil \ \ \fbb #\frb\ 
   &\vrule \hfil \ \ \frb #\ \hfil
   &\vrule \hfil \ \ \frb #\ 
   &\vrule \hfil \ \ \frb #\ \ \vrule \hskip 2 pt
   &\vrule \ \ \hfil \frb #\ 
   &\vrule \hfil \ \ \fbb #\frb\ 
   &\vrule \hfil \ \ \frb #\ \hfil
   &\vrule \hfil \ \ \frb #\ 
   &\vrule \hfil \ \ \frb #\ \vrule \cr%
\noalign{\hrule}
 & &9.125.11.29.97&1343&1082& & &9.25.7.11.29.71&1167&892 \cr
2017&34810875&4.5.11.17.79.541&291&104&2035&35672175&8.27.7.223.389&17&206 \cr
 & &64.3.13.97.541&541&416& & &32.17.103.389&1751&6224 \cr
\noalign{\hrule}
 & &3.5.7.19.101.173&39&134& & &3.11.13.19.41.107&245&206 \cr
2018&34858635&4.9.7.13.67.101&689&220&2036&35758437&4.5.49.19.103.107&1377&656 \cr
 & &32.5.11.169.53&583&2704& & &128.81.5.7.17.41&945&1088 \cr
\noalign{\hrule}
 & &9.11.13.37.733&31&68& & &3.5.7.17.29.691&2665&2172 \cr
2019&34904727&8.13.17.31.733&165&568&2037&35769615&8.9.25.13.41.181&11101&10076 \cr
 & &128.3.5.11.17.71&355&1088& & &64.11.17.229.653&7183&7328 \cr
\noalign{\hrule}
 & &27.5.49.11.13.37&9617&7582& & &9.7.11.13.41.97&155&524 \cr
2020&34999965&4.17.59.163.223&111&52&2038&35828793&8.5.11.13.31.131&27&38 \cr
 & &32.3.13.17.37.223&223&272& & &32.27.19.31.131&2489&1488 \cr
\noalign{\hrule}
 & &3.5.103.139.163&157&260& & &7.19.29.71.131&65&66 \cr
2021&35005065&8.25.13.157.163&1017&3058&2039&35873957&4.3.5.7.11.13.19.29.71&14017&822 \cr
 & &32.9.11.113.139&339&176& & &16.9.107.131.137&963&1096 \cr
\noalign{\hrule}
 & &27.343.19.199&5945&572& & &5.11.17.83.463&6603&8918 \cr
2022&35015841&8.5.11.13.29.41&209&168&2040&35931115&4.3.343.13.31.71&61&30 \cr
 & &128.3.5.7.121.19&121&320& & &16.9.5.49.61.71&3479&4392 \cr
\noalign{\hrule}
 & &9.7.11.19.2663&1363&1300& & &5.13.59.83.113&32113&31548 \cr
2023&35063721&8.25.11.13.19.29.47&81&128&2041&35968465&8.3.11.17.239.1889&10509&10270 \cr
 & &2048.81.25.13.29&9425&9216& & &32.9.5.13.17.31.79.113&2449&2448 \cr
\noalign{\hrule}
 & &5.19.37.67.149&93&242& & &3.7.11.13.67.179&411&590 \cr
2024&35090245&4.3.121.19.31.37&2533&1944&2042&36014979&4.9.5.59.67.137&5549&2534 \cr
 & &64.729.17.149&729&544& & &16.7.31.179.181&181&248 \cr
\noalign{\hrule}
 & &7.121.19.37.59&433&414& & &3.125.139.691&7&132 \cr
2025&35131019&4.9.23.37.59.433&1705&8254&2043&36018375&8.9.7.11.691&341&350 \cr
 & &16.3.5.11.31.4127&4127&3720& & &32.25.49.121.31&1519&1936 \cr
\noalign{\hrule}
 & &9.13.467.643&5929&142& & &9.11.43.61.139&595&656 \cr
2026&35132877&4.49.121.71&61&60&2044&36095103&32.5.7.11.17.41.43&57&244 \cr
 & &32.3.5.49.61.71&3479&4880& & &256.3.5.19.41.61&779&640 \cr
\noalign{\hrule}
 & &3.7.11.17.23.389&79&310& & &125.11.13.43.47&531&544 \cr
2027&35134869&4.5.17.23.31.79&4279&4806&2045&36125375&64.9.5.11.17.47.59&211&24 \cr
 & &16.27.11.89.389&89&72& & &1024.27.59.211&12449&13824 \cr
\noalign{\hrule}
 & &13.23.41.47.61&435&968& & &3.13.19.59.827&385&442 \cr
2028&35146553&16.3.5.121.29.47&273&244&2046&36155613&4.5.7.11.169.17.59&1431&428 \cr
 & &128.9.5.7.11.13.61&495&448& & &32.27.5.7.53.107&4815&5936 \cr
\noalign{\hrule}
 & &81.5.11.41.193&37&4& & &27.49.11.41.61&241&430 \cr
2029&35252415&8.27.5.37.193&29&164&2047&36397053&4.5.7.41.43.241&165&122 \cr
 & &64.29.37.41&37&928& & &16.3.25.11.61.241&241&200 \cr
\noalign{\hrule}
 & &27.5.11.13.31.59&203&262& & &7.13.61.79.83&85&6 \cr
2030&35308845&4.9.7.11.13.29.131&59&202&2048&36397907&4.3.5.17.61.83&553&858 \cr
 & &16.7.59.101.131&707&1048& & &16.9.7.11.13.79&11&72 \cr
\noalign{\hrule}
 & &3.5.7.11.23.31.43&793&540& & &729.7.37.193&235&494 \cr
2031&35411145&8.81.25.7.13.61&1333&3608&2049&36440523&4.5.13.19.47.193&27&220 \cr
 & &128.11.31.41.43&41&64& & &32.27.25.11.47&1175&176 \cr
\noalign{\hrule}
 & &9.11.2197.163&8365&8366& & &5.11.17.43.907&6527&8892 \cr
2032&35452989&4.5.7.13.47.89.163.239&272297&402&2050&36465935&8.9.13.19.61.107&1507&2666 \cr
 & &16.3.23.67.11839&11839&12328& & &32.3.11.31.43.137&411&496 \cr
\noalign{\hrule}
 & &17.19.31.53.67&451&450& & &5.11.31.73.293&331&1134 \cr
2033&35556163&4.9.25.11.19.31.41.67&1&1272&2051&36468245&4.81.7.31.331&253&584 \cr
 & &64.27.25.11.53&275&864& & &64.3.7.11.23.73&69&224 \cr
\noalign{\hrule}
 & &3.125.7.107.127&1441&806& & &9.49.23.59.61&1243&2600 \cr
2034&35671125&4.25.11.13.31.131&1467&1808&2052&36504657&16.25.7.11.13.113&183&92 \cr
 & &128.9.13.113.163&6357&7232& & &128.3.23.61.113&113&64 \cr
\noalign{\hrule}
}%
}
$$
\eject
\vglue -23 pt
\noindent\hskip 1 in\hbox to 6.5 in{\ 2053 -- 2088 \hfill\fbd 36509165 -- 37857105\frb}
\vskip -9 pt
$$
\vbox{
\nointerlineskip
\halign{\strut
    \vrule \ \ \hfil \frb #\ 
   &\vrule \hfil \ \ \fbb #\frb\ 
   &\vrule \hfil \ \ \frb #\ \hfil
   &\vrule \hfil \ \ \frb #\ 
   &\vrule \hfil \ \ \frb #\ \ \vrule \hskip 2 pt
   &\vrule \ \ \hfil \frb #\ 
   &\vrule \hfil \ \ \fbb #\frb\ 
   &\vrule \hfil \ \ \frb #\ \hfil
   &\vrule \hfil \ \ \frb #\ 
   &\vrule \hfil \ \ \frb #\ \vrule \cr%
\noalign{\hrule}
 & &5.49.11.19.23.31&107&138& & &9.49.11.13.19.31&391&50 \cr
2053&36509165&4.3.11.19.529.107&369&160&2071&37144107&4.25.13.17.19.23&251&186 \cr
 & &256.27.5.41.107&4387&3456& & &16.3.5.17.31.251&251&680 \cr
\noalign{\hrule}
 & &3.5.11.17.29.449&2223&2716& & &25.11.29.59.79&687&38 \cr
2054&36523905&8.27.5.7.13.19.97&29&484&2072&37171475&4.3.19.79.229&865&636 \cr
 & &64.121.29.97&97&352& & &32.9.5.53.173&477&2768 \cr
\noalign{\hrule}
 & &3.5.13.37.61.83&6083&688& & &3.11.31.37.983&493&530 \cr
2055&36529545&32.7.11.43.79&61&18&2073&37207533&4.5.17.29.53.983&277&1260 \cr
 & &128.9.7.11.61&33&448& & &32.9.25.7.17.277&5817&6800 \cr
\noalign{\hrule}
 & &7.11.13.17.19.113&2319&398& & &121.13.41.577&1075&498 \cr
2056&36535499&4.3.7.199.773&1083&310&2074&37212461&4.3.25.41.43.83&187&228 \cr
 & &16.9.5.361.31&1395&152& & &32.9.5.11.17.19.43&2907&3440 \cr
\noalign{\hrule}
 & &729.7.13.19.29&355&374& & &3.5.13.263.727&677&638 \cr
2057&36552789&4.5.7.11.13.17.29.71&109&486&2075&37284195&4.11.29.677.727&25&702 \cr
 & &16.243.11.71.109&781&872& & &16.27.25.11.13.29&495&232 \cr
\noalign{\hrule}
 & &49.13.19.3023&1635&1388& & &3.11.61.97.191&65&126 \cr
2058&36587369&8.3.5.49.109.347&149&198&2076&37294851&4.27.5.7.11.13.97&191&488 \cr
 & &32.27.5.11.109.149&20115&19184& & &64.5.13.61.191&65&32 \cr
\noalign{\hrule}
 & &11.47.193.367&355&162& & &9.5.11.23.29.113&6109&106 \cr
2059&36619627&4.81.5.71.367&2717&2788&2077&37308645&4.41.53.149&1161&1012 \cr
 & &32.27.11.13.17.19.41&8721&8528& & &32.27.11.23.43&129&16 \cr
\noalign{\hrule}
 & &9.5.7.11.19.557&3379&3936& & &3.5.19.173.757&119&638 \cr
2060&36670095&64.27.31.41.109&2107&836&2078&37323885&4.5.7.11.17.19.29&83&468 \cr
 & &512.49.11.19.43&301&256& & &32.9.13.17.83&249&3536 \cr
\noalign{\hrule}
 & &9.25.7.11.29.73&6431&6344& & &9.11.169.23.97&9481&7250 \cr
2061&36677025&16.3.11.13.59.61.109&265&1682&2079&37326861&4.125.19.29.499&213&338 \cr
 & &64.5.841.53.61&1769&1696& & &16.3.169.71.499&499&568 \cr
\noalign{\hrule}
 & &3.25.29.101.167&279&446& & &3.11.17.19.31.113&273&254 \cr
2062&36685725&4.27.31.101.223&1445&4576&2080&37338477&4.9.7.11.13.113.127&1207&190 \cr
 & &256.5.11.13.289&3179&1664& & &16.5.7.13.17.19.71&455&568 \cr
\noalign{\hrule}
 & &9.11.17.139.157&95&62& & &27.125.11.19.53&59&1066 \cr
2063&36728109&4.3.5.17.19.31.139&55&472&2081&37384875&4.3.11.13.41.59&475&1124 \cr
 & &64.25.11.19.59&475&1888& & &32.25.19.281&281&16 \cr
\noalign{\hrule}
 & &25.7.11.361.53&8551&6024& & &81.7.169.17.23&55&244 \cr
2064&36831025&16.3.17.251.503&2385&1882&2082&37466793&8.3.5.11.13.17.61&109&112 \cr
 & &64.27.5.53.941&941&864& & &256.5.7.11.61.109&6649&7040 \cr
\noalign{\hrule}
 & &3.5.17.23.61.103&1333&418& & &9.7.11.13.37.113&155&118 \cr
2065&36849795&4.11.19.23.31.43&793&540&2083&37666629&4.3.5.11.31.59.113&455&568 \cr
 & &32.27.5.13.19.61&247&144& & &64.25.7.13.59.71&1775&1888 \cr
\noalign{\hrule}
 & &25.17.29.41.73&1859&1134& & &3.121.17.29.211&1407&914 \cr
2066&36888725&4.81.7.11.169.17&205&16&2084&37760349&4.9.7.11.67.457&2425&1688 \cr
 & &128.3.5.11.13.41&429&64& & &64.25.7.97.211&679&800 \cr
\noalign{\hrule}
 & &9.11.113.3299&1819&1480& & &9.25.7.11.37.59&353&178 \cr
2067&36905913&16.3.5.11.17.37.107&1699&2260&2085&37820475&4.11.37.89.353&295&3588 \cr
 & &128.25.113.1699&1699&1600& & &32.3.5.13.23.59&23&208 \cr
\noalign{\hrule}
 & &27.11.13.17.563&67&496& & &9.19.841.263&335&506 \cr
2068&36953631&32.9.17.31.67&565&38&2086&37822293&4.5.11.23.67.263&289&1026 \cr
 & &128.5.19.113&113&6080& & &16.27.289.19.23&289&552 \cr
\noalign{\hrule}
 & &81.23.83.239&8987&10850& & &17.23.29.47.71&1335&2002 \cr
2069&36956331&4.25.7.11.19.31.43&69&26&2087&37838243&4.3.5.7.11.13.17.89&87&2 \cr
 & &16.3.5.7.11.13.23.31&1001&1240& & &16.9.7.11.13.29&63&1144 \cr
\noalign{\hrule}
 & &49.11.23.29.103&201&520& & &27.5.11.13.37.53&213&268 \cr
2070&37029839&16.3.5.7.13.23.67&103&12&2088&37857105&8.81.53.67.71&4595&832 \cr
 & &128.9.67.103&603&64& & &1024.5.13.919&919&512 \cr
\noalign{\hrule}
}%
}
$$
\eject
\vglue -23 pt
\noindent\hskip 1 in\hbox to 6.5 in{\ 2089 -- 2124 \hfill\fbd 37905219 -- 39821925\frb}
\vskip -9 pt
$$
\vbox{
\nointerlineskip
\halign{\strut
    \vrule \ \ \hfil \frb #\ 
   &\vrule \hfil \ \ \fbb #\frb\ 
   &\vrule \hfil \ \ \frb #\ \hfil
   &\vrule \hfil \ \ \frb #\ 
   &\vrule \hfil \ \ \frb #\ \ \vrule \hskip 2 pt
   &\vrule \ \ \hfil \frb #\ 
   &\vrule \hfil \ \ \fbb #\frb\ 
   &\vrule \hfil \ \ \frb #\ \hfil
   &\vrule \hfil \ \ \frb #\ 
   &\vrule \hfil \ \ \frb #\ \vrule \cr%
\noalign{\hrule}
 & &27.11.23.31.179&65&34& & &27.11.13.19.529&413&116 \cr
2089&37905219&4.3.5.13.17.23.179&55&124&2107&38806911&8.7.13.19.29.59&759&1880 \cr
 & &32.25.11.13.17.31&221&400& & &128.3.5.11.23.47&235&64 \cr
\noalign{\hrule}
 & &9.5.11.13.43.137&85&1592& & &9.5.7.11.17.661&3661&3610 \cr
2090&37908585&16.3.25.17.199&137&62&2108&38936205&4.3.25.49.361.523&1397&172 \cr
 & &64.17.31.137&527&32& & &32.11.361.43.127&5461&5776 \cr
\noalign{\hrule}
 & &25.289.29.181&6237&988& & &9.23.31.59.103&11575&9746 \cr
2091&37924025&8.81.7.11.13.19&85&58&2109&38996109&4.25.11.443.463&1339&876 \cr
 & &32.3.5.7.17.19.29&133&48& & &32.3.5.11.13.73.103&949&880 \cr
\noalign{\hrule}
 & &27.11.289.443&95&194& & &5.13.19.23.1373&537&836 \cr
2092&38024019&4.3.5.19.97.443&713&616&2110&39000065&8.3.5.11.361.179&1887&82 \cr
 & &64.5.7.11.19.23.31&4123&3680& & &32.9.17.37.41&6273&592 \cr
\noalign{\hrule}
 & &27.7.11.13.17.83&31&32& & &5.121.19.41.83&9657&7358 \cr
2093&38135097&64.3.11.13.17.31.83&259&820&2111&39117485&4.9.13.29.37.283&41&70 \cr
 & &512.5.7.31.37.41&7585&7936& & &16.3.5.7.13.41.283&849&728 \cr
\noalign{\hrule}
 & &25.7.11.83.239&6183&208& & &3.11.17.31.37.61&1129&92 \cr
2094&38186225&32.27.13.229&337&350&2112&39251487&8.23.31.1129&549&580 \cr
 & &128.9.25.7.337&337&576& & &64.9.5.23.29.61&667&480 \cr
\noalign{\hrule}
 & &3.11.13.29.37.83&963&1444& & &5.7.11.13.47.167&11611&9774 \cr
2095&38206311&8.27.11.361.107&845&332&2113&39284245&4.27.17.181.683&613&70 \cr
 & &64.5.169.19.83&247&160& & &16.9.5.7.17.613&613&1224 \cr
\noalign{\hrule}
 & &7.11.17.19.29.53&41&12& & &9.5.7.11.17.23.29&1&86 \cr
2096&38226727&8.3.7.11.17.19.41&265&1044&2114&39289635&4.3.7.11.23.43&325&578 \cr
 & &64.27.5.29.53&27&160& & &16.25.13.289&13&680 \cr
\noalign{\hrule}
 & &9.5.11.13.43.139&4453&4582& & &11.13.17.19.23.37&45&436 \cr
2097&38461995&4.3.11.29.61.73.79&1205&3614&2115&39306839&8.9.5.11.19.109&23&34 \cr
 & &16.5.13.29.139.241&241&232& & &32.3.5.17.23.109&327&80 \cr
\noalign{\hrule}
 & &9.5.11.23.31.109&211&6206& & &81.13.37.1009&8057&5060 \cr
2098&38469915&4.29.107.211&39&68&2116&39311649&8.5.7.11.23.1151&449&702 \cr
 & &32.3.13.17.211&3587&208& & &32.27.5.7.13.449&449&560 \cr
\noalign{\hrule}
 & &49.17.113.409&985&936& & &9.7.11.17.47.71&6631&2158 \cr
2099&38498761&16.9.5.13.197.409&121&1106&2117&39313197&4.13.19.83.349&365&714 \cr
 & &64.3.7.121.13.79&3081&3872& & &16.3.5.7.17.19.73&365&152 \cr
\noalign{\hrule}
 & &27.25.13.23.191&1649&6424& & &11.17.23.89.103&7137&9506 \cr
2100&38548575&16.11.17.73.97&353&450&2118&39427267&4.9.49.13.61.97&1025&1012 \cr
 & &64.9.25.17.353&353&544& & &32.3.25.7.11.23.41.61&4575&4592 \cr
\noalign{\hrule}
 & &3.11.17.97.709&635&74& & &7.11.31.71.233&135&206 \cr
2101&38581653&4.5.37.97.127&179&306&2119&39488141&4.27.5.7.103.233&11&710 \cr
 & &16.9.17.37.179&179&888& & &16.9.25.11.71&225&8 \cr
\noalign{\hrule}
 & &81.5.19.47.107&221&202& & &41.47.103.199&2321&2520 \cr
2102&38698155&4.9.5.13.17.101.107&187&722&2120&39497719&16.9.5.7.11.41.211&2441&1990 \cr
 & &16.11.13.289.361&2717&2312& & &64.3.25.199.2441&2441&2400 \cr
\noalign{\hrule}
 & &17.19.37.41.79&7595&4356& & &27.5.17.59.293&817&1820 \cr
2103&38709289&8.9.5.49.121.31&41&8&2121&39673665&8.3.25.7.13.19.43&2363&88 \cr
 & &128.3.5.11.31.41&1023&320& & &128.11.17.139&1529&64 \cr
\noalign{\hrule}
 & &9.5.13.97.683&2257&6622& & &25.11.37.47.83&321&86 \cr
2104&38756835&4.7.11.37.43.61&145&156&2122&39692675&4.3.5.43.83.107&1517&2052 \cr
 & &32.3.5.13.29.37.61&1073&976& & &32.81.19.37.41&779&1296 \cr
\noalign{\hrule}
 & &3.7.17.961.113&20651&19690& & &3.107.337.367&559&452 \cr
2105&38767701&4.5.11.107.179.193&1071&106&2123&39700959&8.13.43.113.367&4815&44 \cr
 & &16.9.7.17.53.179&537&424& & &64.9.5.11.107&15&352 \cr
\noalign{\hrule}
 & &7.47.191.617&4797&4180& & &3.25.11.13.47.79&161&114 \cr
2106&38771663&8.9.5.7.11.13.19.41&47&86&2124&39821925&4.9.7.13.19.23.79&203&2020 \cr
 & &32.3.5.11.41.43.47&2255&2064& & &32.5.49.29.101&2929&784 \cr
\noalign{\hrule}
}%
}
$$
\eject
\vglue -23 pt
\noindent\hskip 1 in\hbox to 6.5 in{\ 2125 -- 2160 \hfill\fbd 39831935 -- 41484645\frb}
\vskip -9 pt
$$
\vbox{
\nointerlineskip
\halign{\strut
    \vrule \ \ \hfil \frb #\ 
   &\vrule \hfil \ \ \fbb #\frb\ 
   &\vrule \hfil \ \ \frb #\ \hfil
   &\vrule \hfil \ \ \frb #\ 
   &\vrule \hfil \ \ \frb #\ \ \vrule \hskip 2 pt
   &\vrule \ \ \hfil \frb #\ 
   &\vrule \hfil \ \ \fbb #\frb\ 
   &\vrule \hfil \ \ \frb #\ \hfil
   &\vrule \hfil \ \ \frb #\ 
   &\vrule \hfil \ \ \frb #\ \vrule \cr%
\noalign{\hrule}
 & &5.11.13.17.29.113&2679&2744& & &3.25.7.19.31.131&41&434 \cr
2125&39831935&16.3.343.19.47.113&2225&78&2143&40508475&4.49.961.41&1485&524 \cr
 & &64.9.25.7.13.89&801&1120& & &32.27.5.11.131&11&144 \cr
\noalign{\hrule}
 & &5.13.31.53.373&33&20& & &9.11.17.101.239&43&196 \cr
2126&39834535&8.3.25.11.31.373&201&574&2144&40625937&8.49.11.43.101&2133&2210 \cr
 & &32.9.7.11.41.67&4059&7504& & &32.27.5.7.13.17.79&1365&1264 \cr
\noalign{\hrule}
 & &3.23.961.601&545&1258& & &9.7.29.37.601&689&88 \cr
2127&39851709&4.5.17.31.37.109&531&616&2145&40626999&16.3.11.13.29.53&785&752 \cr
 & &64.9.7.11.59.109&6431&7392& & &512.5.13.47.157&10205&12032 \cr
\noalign{\hrule}
 & &29.31.101.439&65&36& & &27.73.107.193&4195&3616 \cr
2128&39860761&8.9.5.13.31.439&385&824&2146&40703121&64.9.5.113.839&2123&2962 \cr
 & &128.3.25.7.11.103&7931&4800& & &256.11.193.1481&1481&1408 \cr
\noalign{\hrule}
 & &9.7.29.131.167&185&316& & &3.25.11.19.23.113&91&1334 \cr
2129&39969279&8.3.5.7.29.37.79&323&286&2147&40739325&4.7.13.529.29&279&250 \cr
 & &32.5.11.13.17.19.79&16511&17680& & &16.9.125.7.13.31&403&840 \cr
\noalign{\hrule}
 & &3.37.41.59.149&3025&3084& & &25.121.13.17.61&8253&5228 \cr
2130&40007841&8.9.25.121.37.257&23&14&2148&40780025&8.9.7.131.1307&457&850 \cr
 & &32.25.7.121.23.257&44975&44528& & &32.3.25.7.17.457&457&336 \cr
\noalign{\hrule}
 & &81.37.43.311&755&836& & &3.5.11.13.23.827&425&402 \cr
2131&40078881&8.5.11.19.151.311&3145&276&2149&40800045&4.9.125.11.13.17.67&29&1654 \cr
 & &64.3.25.17.23.37&425&736& & &16.29.67.827&67&232 \cr
\noalign{\hrule}
 & &9.25.23.61.127&589&814& & &9.125.7.29.179&6533&1342 \cr
2132&40090725&4.11.19.31.37.127&5589&890&2150&40879125&4.11.47.61.139&39&100 \cr
 & &16.243.5.23.89&89&216& & &32.3.25.11.13.47&143&752 \cr
\noalign{\hrule}
 & &81.11.23.37.53&77&130& & &5.121.19.3557&3969&7526 \cr
2133&40186773&4.9.5.7.121.13.37&227&106&2151&40887715&4.81.49.53.71&2257&340 \cr
 & &16.5.7.13.53.227&1135&728& & &32.3.5.17.37.61&629&2928 \cr
\noalign{\hrule}
 & &27.7.11.289.67&1267&1334& & &7.11.13.17.29.83&2779&1368 \cr
2134&40255677&4.3.49.11.23.29.181&1273&2890&2152&40959919&16.9.49.19.397&1595&1198 \cr
 & &16.5.289.19.29.67&145&152& & &64.3.5.11.29.599&599&480 \cr
\noalign{\hrule}
 & &3.5.23.841.139&363&478& & &5.13.37.113.151&207&358 \cr
2135&40330155&4.9.121.139.239&745&506&2153&41036515&4.9.13.23.37.179&151&330 \cr
 & &16.5.1331.23.149&1331&1192& & &16.27.5.11.23.151&297&184 \cr
\noalign{\hrule}
 & &3.7.11.17.43.239&299&418& & &81.11.47.983&39553&40070 \cr
2136&40357779&4.121.13.19.23.43&1195&378&2154&41165091&4.5.37.1069.4007&1469&2538 \cr
 & &16.27.5.7.23.239&115&72& & &16.27.5.13.37.47.113&1469&1480 \cr
\noalign{\hrule}
 & &27.11.23.61.97&25&646& & &49.41.103.199&1415&2808 \cr
2137&40419027&4.25.17.19.97&61&36&2155&41178473&16.27.5.7.13.283&2189&1906 \cr
 & &32.9.17.19.61&19&272& & &64.3.11.199.953&953&1056 \cr
\noalign{\hrule}
 & &625.71.911&1343&432& & &5.31.37.43.167&243&6422 \cr
2138&40425625&32.27.25.17.79&143&568&2156&41183035&4.243.169.19&341&172 \cr
 & &512.3.11.13.71&143&768& & &32.9.11.31.43&9&176 \cr
\noalign{\hrule}
 & &27.343.11.397&223&174& & &9.7.11.13.23.199&19&180 \cr
2139&40442787&4.81.7.11.29.223&115&6352&2157&41234193&8.81.5.11.13.19&31&112 \cr
 & &128.5.23.397&23&320& & &256.5.7.19.31&155&2432 \cr
\noalign{\hrule}
 & &3.17.53.71.211&343&290& & &9.5.7.11.17.19.37&5459&5200 \cr
2140&40493643&4.5.343.17.29.71&495&2&2158&41410215&32.3.125.13.53.103&6517&6358 \cr
 & &16.9.25.49.11&3675&88& & &128.343.11.13.289.19&833&832 \cr
\noalign{\hrule}
 & &81.7.11.43.151&377&76& & &5.7.11.19.53.107&9&86 \cr
2141&40496841&8.27.11.13.19.29&335&16&2159&41483365&4.9.43.53.107&133&26 \cr
 & &256.5.19.67&6365&128& & &16.3.7.13.19.43&559&24 \cr
\noalign{\hrule}
 & &9.5.11.19.59.73&1513&1732& & &9.5.29.83.383&9821&7414 \cr
2142&40507335&8.3.17.19.89.433&295&28&2160&41484645&4.7.11.23.61.337&533&870 \cr
 & &64.5.7.59.433&433&224& & &16.3.5.7.11.13.29.41&533&616 \cr
\noalign{\hrule}
}%
}
$$
\eject
\vglue -23 pt
\noindent\hskip 1 in\hbox to 6.5 in{\ 2161 -- 2196 \hfill\fbd 41525825 -- 42762875\frb}
\vskip -9 pt
$$
\vbox{
\nointerlineskip
\halign{\strut
    \vrule \ \ \hfil \frb #\ 
   &\vrule \hfil \ \ \fbb #\frb\ 
   &\vrule \hfil \ \ \frb #\ \hfil
   &\vrule \hfil \ \ \frb #\ 
   &\vrule \hfil \ \ \frb #\ \ \vrule \hskip 2 pt
   &\vrule \ \ \hfil \frb #\ 
   &\vrule \hfil \ \ \fbb #\frb\ 
   &\vrule \hfil \ \ \frb #\ \hfil
   &\vrule \hfil \ \ \frb #\ 
   &\vrule \hfil \ \ \frb #\ \vrule \cr%
\noalign{\hrule}
 & &25.11.29.41.127&7479&3796& & &9.49.11.13.23.29&425&334 \cr
2161&41525825&8.27.13.73.277&95&22&2179&42063021&4.3.25.7.17.29.167&149&352 \cr
 & &32.3.5.11.19.277&277&912& & &256.25.11.17.149&2533&3200 \cr
\noalign{\hrule}
 & &27.5.17.23.787&203&2158& & &9.5.7.13.43.239&1661&1446 \cr
2162&41541795&4.9.7.13.29.83&85&176&2180&42084315&4.27.7.11.151.241&215&26 \cr
 & &128.5.11.17.83&83&704& & &16.5.11.13.43.151&151&88 \cr
\noalign{\hrule}
 & &3.5.19.139.1049&477&572& & &3.5.49.11.41.127&59&186 \cr
2163&41556135&8.27.11.13.53.139&589&940&2181&42098595&4.9.11.31.41.59&325&2744 \cr
 & &64.5.19.31.47.53&1643&1504& & &64.25.343.13&35&416 \cr
\noalign{\hrule}
 & &841.73.677&759&82& & &3.5.11.169.17.89&49&138 \cr
2164&41563061&4.3.11.23.41.73&435&508&2182&42190005&4.9.5.49.169.23&73&242 \cr
 & &32.9.5.11.29.127&1397&720& & &16.7.121.23.73&511&2024 \cr
\noalign{\hrule}
 & &3.7.121.169.97&207&1054& & &3.5.83.107.317&109&426 \cr
2165&41654613&4.27.13.17.23.31&3773&4300&2183&42229155&4.9.71.83.109&319&428 \cr
 & &32.25.343.11.43&1075&784& & &32.11.29.71.107&781&464 \cr
\noalign{\hrule}
 & &9.5.23.67.601&1201&1804& & &13.19.41.43.97&55&42 \cr
2166&41676345&8.11.23.41.1201&129&1072&2184&42239717&4.3.5.7.11.19.41.43&2619&2834 \cr
 & &256.3.11.43.67&473&128& & &16.81.11.13.97.109&891&872 \cr
\noalign{\hrule}
 & &121.29.11881&7695&4186& & &9.17.31.59.151&169&110 \cr
2167&41690429&4.81.5.7.13.19.23&109&242&2185&42255387&4.5.11.169.17.151&35&186 \cr
 & &16.3.5.121.23.109&69&40& & &16.3.25.7.11.13.31&275&728 \cr
\noalign{\hrule}
 & &5.71.257.457&871&414& & &9.125.49.13.59&8041&8066 \cr
2168&41694395&4.9.13.23.67.71&2827&1930&2186&42280875&4.3.5.7.11.17.37.43.109&2717&428 \cr
 & &16.3.5.11.193.257&193&264& & &32.121.13.19.43.107&12947&13072 \cr
\noalign{\hrule}
 & &9.5.13.23.29.107&6479&476& & &3.5.7.13.139.223&211&484 \cr
2169&41750865&8.7.11.17.19.31&161&162&2187&42310905&8.121.211.223&117&106 \cr
 & &32.81.49.11.23.31&1519&1584& & &32.9.11.13.53.211&2321&2544 \cr
\noalign{\hrule}
 & &121.17.23.883&587&1470& & &9.25.7.11.31.79&9809&7334 \cr
2170&41775613&4.3.5.49.23.587&213&374&2188&42428925&4.17.19.193.577&105&88 \cr
 & &16.9.5.7.11.17.71&315&568& & &64.3.5.7.11.19.577&577&608 \cr
\noalign{\hrule}
 & &9.11.19.97.229&1055&826& & &3.7.17.271.439&209&230 \cr
2171&41782653&4.5.7.59.97.211&57&154&2189&42471933&4.5.11.17.19.23.271&297&26 \cr
 & &16.3.5.49.11.19.59&295&392& & &16.27.5.121.13.23&2783&4680 \cr
\noalign{\hrule}
 & &5.139.157.383&11869&9954& & &13.31.263.401&2475&2738 \cr
2172&41791045&4.9.7.11.13.79.83&157&74&2190&42501589&4.9.25.11.31.1369&263&78 \cr
 & &16.3.13.37.79.157&1027&888& & &16.27.5.13.37.263&185&216 \cr
\noalign{\hrule}
 & &27.7.23.59.163&14725&16082& & &25.11.13.73.163&67&882 \cr
2173&41805099&4.25.11.17.19.31.43&199&1134&2191&42538925&4.9.5.49.11.67&3&52 \cr
 & &16.81.5.7.19.199&597&760& & &32.27.13.67&1809&16 \cr
\noalign{\hrule}
 & &11.13.23.29.439&1921&6750& & &9.11.17.19.31.43&505&226 \cr
2174&41872259&4.27.125.17.113&551&466&2192&42625341&4.5.11.19.101.113&1581&338 \cr
 & &16.3.25.19.29.233&1425&1864& & &16.3.5.169.17.31&169&40 \cr
\noalign{\hrule}
 & &3.19.31.151.157&85&66& & &9.25.7.11.23.107&1111&1136 \cr
2175&41890269&4.9.5.11.17.31.157&2869&200&2193&42636825&32.3.121.23.71.101&3661&8560 \cr
 & &64.125.19.151&125&32& & &1024.5.7.107.523&523&512 \cr
\noalign{\hrule}
 & &27.5.23.103.131&17&638& & &3.11.13.23.61.71&119&180 \cr
2176&41895765&4.11.17.29.103&575&558&2194&42733977&8.27.5.7.11.17.71&161&620 \cr
 & &16.9.25.23.29.31&145&248& & &64.25.49.23.31&775&1568 \cr
\noalign{\hrule}
 & &25.13.23.71.79&189&110& & &9.5.49.11.41.43&527&978 \cr
2177&41927275&4.27.125.7.11.71&527&598&2195&42761565&4.27.7.17.31.163&13&176 \cr
 & &16.3.7.11.13.17.23.31&1023&952& & &128.11.13.17.31&221&1984 \cr
\noalign{\hrule}
 & &3.25.13.17.43.59&203&528& & &125.157.2179&8723&10902 \cr
2178&42050775&32.9.7.11.29.59&25&34&2196&42762875&4.3.11.13.23.61.79&157&96 \cr
 & &128.25.7.11.17.29&319&448& & &256.9.13.79.157&1027&1152 \cr
\noalign{\hrule}
}%
}
$$
\eject
\vglue -23 pt
\noindent\hskip 1 in\hbox to 6.5 in{\ 2197 -- 2232 \hfill\fbd 42865515 -- 44127291\frb}
\vskip -9 pt
$$
\vbox{
\nointerlineskip
\halign{\strut
    \vrule \ \ \hfil \frb #\ 
   &\vrule \hfil \ \ \fbb #\frb\ 
   &\vrule \hfil \ \ \frb #\ \hfil
   &\vrule \hfil \ \ \frb #\ 
   &\vrule \hfil \ \ \frb #\ \ \vrule \hskip 2 pt
   &\vrule \ \ \hfil \frb #\ 
   &\vrule \hfil \ \ \fbb #\frb\ 
   &\vrule \hfil \ \ \frb #\ \hfil
   &\vrule \hfil \ \ \frb #\ 
   &\vrule \hfil \ \ \frb #\ \vrule \cr%
\noalign{\hrule}
 & &9.5.7.11.89.139&299&6554& & &9.5.7.19.37.197&13&184 \cr
2197&42865515&4.13.23.29.113&401&1068&2215&43624665&16.5.7.13.23.37&209&246 \cr
 & &32.3.89.401&401&16& & &64.3.11.19.23.41&943&352 \cr
\noalign{\hrule}
 & &3.25.7.127.643&87&88& & &49.17.23.43.53&1199&1080 \cr
2198&42872025&16.9.11.29.127.643&893&250&2216&43663361&16.27.5.7.11.23.109&1703&68 \cr
 & &64.125.11.19.29.47&6815&6688& & &128.9.13.17.131&1703&576 \cr
\noalign{\hrule}
 & &11.13.41.71.103&725&6588& & &49.11.53.1531&1275&2806 \cr
2199&42876119&8.27.25.29.61&143&82&2217&43736077&4.3.25.7.17.23.61&93&212 \cr
 & &32.3.11.13.29.41&87&16& & &32.9.5.23.31.53&1035&496 \cr
\noalign{\hrule}
 & &3.13.37.113.263&635&154& & &9.5.49.43.463&1199&736 \cr
2200&42884517&4.5.7.11.113.127&95&18&2218&43899345&64.49.11.23.109&215&324 \cr
 & &16.9.25.19.127&381&3800& & &512.81.5.23.43&207&256 \cr
\noalign{\hrule}
 & &11.169.17.23.59&3573&7460& & &5.7.11.13.31.283&207&1208 \cr
2201&42885271&8.9.5.373.397&409&782&2219&43908865&16.9.23.31.151&5225&5194 \cr
 & &32.3.5.17.23.409&409&240& & &64.3.25.49.11.19.53&1855&1824 \cr
\noalign{\hrule}
 & &3.23.29.89.241&165&76& & &5.11.13.239.257&571&714 \cr
2202&42919449&8.9.5.11.19.23.29&89&118&2220&43917445&4.3.7.17.239.571&33&4030 \cr
 & &32.5.11.19.59.89&1121&880& & &16.9.5.11.13.31&9&248 \cr
\noalign{\hrule}
 & &9.5.7.17.29.277&143&134& & &7.11.17.37.907&65&54 \cr
2203&43016715&4.5.7.11.13.17.29.67&1077&62&2221&43928731&4.27.5.13.37.907&1711&2618 \cr
 & &16.3.11.13.31.359&4667&2728& & &16.3.5.7.11.17.29.59&435&472 \cr
\noalign{\hrule}
 & &3125.7.11.179&17277&17098& & &243.5.11.19.173&1073&830 \cr
2204&43071875&4.3.7.13.83.103.443&179&900&2222&43930755&4.25.19.29.37.83&1773&1298 \cr
 & &32.27.25.179.443&443&432& & &16.9.11.29.59.197&1711&1576 \cr
\noalign{\hrule}
 & &5.7.11.19.43.137&119&2484& & &9.5.13.19.37.107&973&418 \cr
2205&43092665&8.27.49.17.23&209&232&2223&44004285&4.3.7.11.361.139&65&296 \cr
 & &128.3.11.17.19.29&493&192& & &64.5.13.37.139&139&32 \cr
\noalign{\hrule}
 & &9.11.13.19.41.43&85&44& & &27.49.17.19.103&1243&920 \cr
2206&43110639&8.3.5.121.13.17.19&305&58&2224&44014887&16.9.5.7.11.23.113&103&1346 \cr
 & &32.25.17.29.61&1769&6800& & &64.5.103.673&673&160 \cr
\noalign{\hrule}
 & &9.7.11.19.29.113&655&542& & &9.5.11.19.31.151&299&290 \cr
2207&43148259&4.5.11.29.131.271&409&3390&2225&44024805&4.25.11.13.23.29.151&3977&402 \cr
 & &16.3.25.113.409&409&200& & &16.3.23.41.67.97&6499&7544 \cr
\noalign{\hrule}
 & &3.625.7.11.13.23&73&502& & &11.13.37.53.157&125&282 \cr
2208&43168125&4.25.7.73.251&3393&2882&2226&44026411&4.3.125.13.47.53&1083&1408 \cr
 & &16.9.11.13.29.131&393&232& & &1024.9.5.11.361&3249&2560 \cr
\noalign{\hrule}
 & &27.5.11.43.677&3737&3710& & &3.7.121.13.31.43&1037&2714 \cr
2209&43229835&4.25.7.37.43.53.101&677&1602&2227&44032989&4.7.17.23.59.61&715&288 \cr
 & &16.9.7.89.101.677&707&712& & &256.9.5.11.13.23&345&128 \cr
\noalign{\hrule}
 & &27.11.13.103.109&209&100& & &9.11.31.113.127&1207&190 \cr
2210&43347447&8.9.25.121.13.19&857&1442&2228&44043219&4.5.17.19.31.71&273&254 \cr
 & &32.5.7.103.857&857&560& & &16.3.5.7.13.71.127&455&568 \cr
\noalign{\hrule}
 & &9.25.17.83.137&1007&2418& & &9.5.11.13.41.167&29&70 \cr
2211&43494075&4.27.13.19.31.53&137&110&2229&44060445&4.25.7.13.29.167&79&246 \cr
 & &16.5.11.31.53.137&341&424& & &16.3.7.29.41.79&203&632 \cr
\noalign{\hrule}
 & &3.5.19.43.53.67&509&308& & &7.169.23.1621&1133&2754 \cr
2212&43517505&8.5.7.11.53.509&387&122&2230&44105789&4.81.7.11.17.103&535&598 \cr
 & &32.9.7.11.43.61&671&336& & &16.9.5.13.17.23.107&765&856 \cr
\noalign{\hrule}
 & &11.17.31.73.103&9667&9594& & &3.125.11.289.37&377&252 \cr
2213&43587643&4.9.7.13.31.41.1381&55&1326&2231&44108625&8.27.7.11.13.17.29&1867&4100 \cr
 & &16.27.5.7.11.169.17&1183&1080& & &64.25.41.1867&1867&1312 \cr
\noalign{\hrule}
 & &3.5.11.3721.71&5617&5546& & &3.13.17.19.31.113&283&686 \cr
2214&43591515&4.5.11.41.47.59.137&117&568&2232&44127291&4.343.113.283&85&198 \cr
 & &64.9.13.47.59.71&1833&1888& & &16.9.5.343.11.17&1029&440 \cr
\noalign{\hrule}
}%
}
$$
\eject
\vglue -23 pt
\noindent\hskip 1 in\hbox to 6.5 in{\ 2233 -- 2268 \hfill\fbd 44146025 -- 46069425\frb}
\vskip -9 pt
$$
\vbox{
\nointerlineskip
\halign{\strut
    \vrule \ \ \hfil \frb #\ 
   &\vrule \hfil \ \ \fbb #\frb\ 
   &\vrule \hfil \ \ \frb #\ \hfil
   &\vrule \hfil \ \ \frb #\ 
   &\vrule \hfil \ \ \frb #\ \ \vrule \hskip 2 pt
   &\vrule \ \ \hfil \frb #\ 
   &\vrule \hfil \ \ \fbb #\frb\ 
   &\vrule \hfil \ \ \frb #\ \hfil
   &\vrule \hfil \ \ \frb #\ 
   &\vrule \hfil \ \ \frb #\ \vrule \cr%
\noalign{\hrule}
 & &25.7.11.17.19.71&27&44& & &13.19.31.61.97&12925&11034 \cr
2233&44146025&8.27.25.7.121.19&1937&362&2251&45306469&4.9.25.11.47.613&661&1178 \cr
 & &32.3.13.149.181&7059&2384& & &16.3.25.19.31.661&661&600 \cr
\noalign{\hrule}
 & &5.11.13.127.487&421&294& & &9.7.11.29.37.61&47&380 \cr
2234&44222035&4.3.49.421.487&33&454&2252&45358929&8.5.11.19.29.47&273&244 \cr
 & &16.9.49.11.227&2043&392& & &64.3.5.7.13.19.61&247&160 \cr
\noalign{\hrule}
 & &3.5.19.29.53.101&6809&1456& & &9.5.7.11.13.19.53&601&634 \cr
2235&44242545&32.7.11.13.619&191&810&2253&45360315&4.3.7.53.317.601&14711&2090 \cr
 & &128.81.5.191&191&1728& & &16.5.11.19.47.313&313&376 \cr
\noalign{\hrule}
 & &3.5.11.13.23.29.31&177&122& & &9.125.13.29.107&2219&2596 \cr
2236&44352165&4.9.29.31.59.61&715&184&2254&45381375&8.25.7.11.59.317&7003&9222 \cr
 & &64.5.11.13.23.61&61&32& & &32.3.29.47.53.149&2491&2384 \cr
\noalign{\hrule}
 & &9.5.7.17.361.23&1781&746& & &3.29.37.103.137&1425&1562 \cr
2237&44462565&4.13.17.137.373&297&76&2255&45423309&4.9.25.11.19.37.71&3451&824 \cr
 & &32.27.11.19.137&137&528& & &64.7.11.17.29.103&187&224 \cr
\noalign{\hrule}
 & &9.11.29.37.419&7787&4016& & &9.7.11.29.31.73&271&532 \cr
2238&44509113&32.13.251.599&425&174&2256&45479511&8.49.19.31.271&319&270 \cr
 & &128.3.25.13.17.29&425&832& & &32.27.5.11.29.271&271&240 \cr
\noalign{\hrule}
 & &3.53.193.1451&805&646& & &43.53.83.241&2221&2178 \cr
2239&44526837&4.5.7.17.19.23.193&99&292&2257&45586837&4.9.121.241.2221&215&2436 \cr
 & &32.9.5.7.11.19.73&4389&5840& & &32.27.5.7.11.29.43&2079&2320 \cr
\noalign{\hrule}
 & &49.19.31.1543&1237&306& & &9.5.31.97.337&2167&2198 \cr
2240&44532523&4.9.17.31.1237&665&572&2258&45601155&4.7.11.157.197.337&2543&1164 \cr
 & &32.3.5.7.11.13.17.19&663&880& & &32.3.97.157.2543&2543&2512 \cr
\noalign{\hrule}
 & &25.7.227.1123&3399&2276& & &9.7.13.29.41.47&1717&3080 \cr
2241&44611175&8.3.7.11.103.569&323&246&2259&45768177&16.5.49.11.17.101&3&52 \cr
 & &32.9.17.19.41.103&13243&14832& & &128.3.13.17.101&101&1088 \cr
\noalign{\hrule}
 & &9.625.7.17.67&31&94& & &3.5.121.151.167&529&76 \cr
2242&44848125&4.5.17.31.47.67&959&3036&2260&45768855&8.19.529.167&135&302 \cr
 & &32.3.7.11.23.137&1507&368& & &32.27.5.23.151&23&144 \cr
\noalign{\hrule}
 & &3.11.29.31.37.41&207&244& & &25.49.11.41.83&469&444 \cr
2243&45004839&8.27.23.29.31.61&85&752&2261&45855425&8.3.343.37.41.67&7719&4972 \cr
 & &256.5.17.47.61&3995&7808& & &64.9.11.31.83.113&1017&992 \cr
\noalign{\hrule}
 & &3.25.11.169.17.19&717&1142& & &27.5.11.17.23.79&43&122 \cr
2244&45034275&4.9.19.239.571&299&4840&2262&45870165&4.9.17.23.43.61&79&470 \cr
 & &64.5.121.13.23&253&32& & &16.5.43.47.79&47&344 \cr
\noalign{\hrule}
 & &3.7.11.23.61.139&25&36& & &5.49.11.29.587&351&2584 \cr
2245&45048927&8.27.25.7.23.139&3961&764&2263&45876985&16.27.7.13.17.19&29&22 \cr
 & &64.17.191.233&3247&7456& & &64.9.11.13.19.29&171&416 \cr
\noalign{\hrule}
 & &9.5.71.103.137&121&806& & &27.125.7.29.67&5203&5672 \cr
2246&45084645&4.121.13.31.71&135&206&2264&45903375&16.9.121.43.709&899&190 \cr
 & &16.27.5.11.13.103&33&104& & &64.5.19.29.31.43&817&992 \cr
\noalign{\hrule}
 & &3.5.7.23.53.353&187&166& & &3.5.49.11.13.19.23&251&186 \cr
2247&45182235&4.5.11.17.23.53.83&3177&1222&2265&45930885&4.9.49.11.31.251&391&50 \cr
 & &16.9.11.13.47.353&429&376& & &16.25.17.23.251&251&680 \cr
\noalign{\hrule}
 & &81.25.11.19.107&221&2254& & &3.5.49.11.13.19.23&8627&5338 \cr
2248&45285075&4.9.49.13.17.23&605&214&2266&45930885&4.17.157.8627&4235&4392 \cr
 & &16.5.7.121.107&7&88& & &64.9.5.7.121.17.61&1037&1056 \cr
\noalign{\hrule}
 & &81.25.11.19.107&51&158& & &17.41.149.443&711&6820 \cr
2249&45285075&4.243.25.17.79&3053&1078&2267&46006879&8.9.5.11.31.79&41&38 \cr
 & &16.49.11.43.71&2107&568& & &32.3.5.11.19.31.41&1045&1488 \cr
\noalign{\hrule}
 & &121.23.41.397&2085&7046& & &27.25.131.521&34103&34148 \cr
2250&45298891&4.3.5.13.139.271&205&66&2268&46069425&8.3.5.67.509.8537&451&8086 \cr
 & &16.9.25.11.13.41&325&72& & &32.11.13.41.67.311&44473&43952 \cr
\noalign{\hrule}
}%
}
$$
\eject
\vglue -23 pt
\noindent\hskip 1 in\hbox to 6.5 in{\ 2269 -- 2304 \hfill\fbd 46135947 -- 48841533\frb}
\vskip -9 pt
$$
\vbox{
\nointerlineskip
\halign{\strut
    \vrule \ \ \hfil \frb #\ 
   &\vrule \hfil \ \ \fbb #\frb\ 
   &\vrule \hfil \ \ \frb #\ \hfil
   &\vrule \hfil \ \ \frb #\ 
   &\vrule \hfil \ \ \frb #\ \ \vrule \hskip 2 pt
   &\vrule \ \ \hfil \frb #\ 
   &\vrule \hfil \ \ \fbb #\frb\ 
   &\vrule \hfil \ \ \frb #\ \hfil
   &\vrule \hfil \ \ \frb #\ 
   &\vrule \hfil \ \ \frb #\ \vrule \cr%
\noalign{\hrule}
 & &3.11.13.41.43.61&95&34& & &27.7.53.67.71&2585&1178 \cr
2269&46135947&4.5.11.13.17.19.41&173&360&2287&47650869&4.9.5.11.19.31.47&1207&1378 \cr
 & &64.9.25.19.173&4325&1824& & &16.13.17.31.53.71&221&248 \cr
\noalign{\hrule}
 & &13.29.31.59.67&231&172& & &25.11.29.43.139&111&584 \cr
2270&46198711&8.3.7.11.29.43.67&1455&1426&2288&47666575&16.3.5.29.37.73&109&36 \cr
 & &32.9.5.7.11.23.31.97&11155&11088& & &128.27.37.109&999&6976 \cr
\noalign{\hrule}
 & &11.13.19.29.587&3503&7650& & &9.17.19.529.31&265&58 \cr
2271&46251491&4.9.25.17.31.113&1007&914&2289&47671893&4.5.23.29.31.53&383&330 \cr
 & &16.3.25.19.53.457&3975&3656& & &16.3.25.11.29.383&4213&5800 \cr
\noalign{\hrule}
 & &11.17.19.29.449&351&142& & &81.5.23.47.109&7117&1712 \cr
2272&46263613&4.27.13.71.449&95&544&2290&47720745&32.11.107.647&377&270 \cr
 & &256.3.5.13.17.19&65&384& & &128.27.5.11.13.29&377&704 \cr
\noalign{\hrule}
 & &27.5.37.59.157&3667&2882& & &3.5.7.11.961.43&2413&2392 \cr
2273&46268685&4.9.11.19.131.193&1147&976&2291&47728065&16.11.13.19.23.43.127&315&18476 \cr
 & &128.31.37.61.131&4061&3904& & &128.9.5.7.31.149&149&192 \cr
\noalign{\hrule}
 & &243.7.11.19.131&785&916& & &49.13.137.547&45&592 \cr
2274&46571679&8.5.11.19.157.229&2751&232&2292&47736143&32.9.5.37.137&187&224 \cr
 & &128.3.5.7.29.131&145&64& & &2048.3.5.7.11.17&2805&1024 \cr
\noalign{\hrule}
 & &3.11.43.107.307&1635&1742& & &9.25.7.121.251&1643&1118 \cr
2275&46612731&4.9.5.13.43.67.109&1787&1228&2293&47834325&4.3.11.13.31.43.53&5789&1790 \cr
 & &32.109.307.1787&1787&1744& & &16.5.7.179.827&827&1432 \cr
\noalign{\hrule}
 & &121.37.53.197&6851&438& & &25.11.31.41.137&613&162 \cr
2276&46744357&4.3.13.17.31.73&605&636&2294&47884925&4.81.137.613&923&310 \cr
 & &32.9.5.121.13.53&117&80& & &16.9.5.13.31.71&117&568 \cr
\noalign{\hrule}
 & &9.5.11.17.67.83&1469&2216& & &3.7.47.67.727&2665&484 \cr
2277&46795815&16.13.17.113.277&249&28&2295&48075783&8.5.7.121.13.41&603&398 \cr
 & &128.3.7.83.113&113&448& & &32.9.11.67.199&199&528 \cr
\noalign{\hrule}
 & &9.25.11.13.31.47&439&956& & &5.49.11.19.23.41&153&134 \cr
2278&46878975&8.5.13.239.439&2651&456&2296&48286315&4.9.5.7.11.17.23.67&779&394 \cr
 & &128.3.11.19.241&241&1216& & &16.3.19.41.67.197&591&536 \cr
\noalign{\hrule}
 & &27.5.13.17.19.83&137&110& & &5.49.13.59.257&83&174 \cr
2279&47049795&4.25.11.17.83.137&1007&2418&2297&48294155&4.3.5.7.29.59.83&143&438 \cr
 & &16.3.11.13.19.31.53&341&424& & &16.9.11.13.29.73&803&2088 \cr
\noalign{\hrule}
 & &3.25.7.11.13.17.37&239&424& & &9.29.31.43.139&21&22 \cr
2280&47222175&16.5.7.11.53.239&621&2294&2298&48359907&4.27.7.11.29.31.139&4945&914 \cr
 & &64.27.23.31.37&713&288& & &16.5.11.23.43.457&2285&2024 \cr
\noalign{\hrule}
 & &9.121.13.47.71&239&850& & &9.7.13.23.31.83&2365&456 \cr
2281&47241909&4.25.17.71.239&297&58&2299&48467601&16.27.5.11.19.43&1519&1046 \cr
 & &16.27.5.11.17.29&85&696& & &64.49.31.523&523&224 \cr
\noalign{\hrule}
 & &9.7.43.101.173&209&310& & &9.11.17.29.997&295&266 \cr
2282&47334357&4.3.5.7.11.19.31.43&173&44&2300&48660579&4.3.5.7.19.59.997&121&1118 \cr
 & &32.5.121.19.173&605&304& & &16.5.121.13.19.43&2717&1720 \cr
\noalign{\hrule}
 & &9.125.7.11.547&169&106& & &25.11.23.43.179&377&198 \cr
2283&47383875&4.5.169.53.547&149&696&2301&48683525&4.9.121.13.29.43&203&160 \cr
 & &64.3.29.53.149&4321&1696& & &256.3.5.7.13.841&5887&4992 \cr
\noalign{\hrule}
 & &243.5.11.53.67&247&490& & &27.7.13.83.239&15403&13730 \cr
2284&47459115&4.25.49.13.19.53&257&432&2302&48739509&4.5.73.211.1373&581&792 \cr
 & &128.27.7.19.257&1799&1216& & &64.9.5.7.11.73.83&365&352 \cr
\noalign{\hrule}
 & &9.11.13.19.29.67&43&14& & &9.5.11.19.29.179&1183&1502 \cr
2285&47512179&4.3.7.11.13.43.67&505&232&2303&48821355&4.3.7.169.19.751&1991&262 \cr
 & &64.5.29.43.101&505&1376& & &16.11.13.131.181&2353&1048 \cr
\noalign{\hrule}
 & &23.29.191.373&7605&3212& & &9.13.19.127.173&2563&2390 \cr
2286&47519081&8.9.5.11.169.73&181&38&2304&48841533&4.3.5.11.19.233.239&6157&1730 \cr
 & &32.3.5.13.19.181&7059&1520& & &16.25.47.131.173&1175&1048 \cr
\noalign{\hrule}
}%
}
$$
\eject
\vglue -23 pt
\noindent\hskip 1 in\hbox to 6.5 in{\ 2305 -- 2340 \hfill\fbd 48922797 -- 51131223\frb}
\vskip -9 pt
$$
\vbox{
\nointerlineskip
\halign{\strut
    \vrule \ \ \hfil \frb #\ 
   &\vrule \hfil \ \ \fbb #\frb\ 
   &\vrule \hfil \ \ \frb #\ \hfil
   &\vrule \hfil \ \ \frb #\ 
   &\vrule \hfil \ \ \frb #\ \ \vrule \hskip 2 pt
   &\vrule \ \ \hfil \frb #\ 
   &\vrule \hfil \ \ \fbb #\frb\ 
   &\vrule \hfil \ \ \frb #\ \hfil
   &\vrule \hfil \ \ \frb #\ 
   &\vrule \hfil \ \ \frb #\ \vrule \cr%
\noalign{\hrule}
 & &3.7.11.29.67.109&475&2686& & &27.7.11.71.337&5525&5596 \cr
2305&48922797&4.25.7.17.19.79&43&36&2323&49744233&8.9.25.7.13.17.1399&2747&1348 \cr
 & &32.9.25.17.19.43&3225&5168& & &64.5.17.41.67.337&2747&2720 \cr
\noalign{\hrule}
 & &27.25.7.13.797&3989&4786& & &25.13.17.29.311&3699&4076 \cr
2306&48955725&4.7.2393.3989&12765&15158&2324&49829975&8.27.17.137.1019&3421&5750 \cr
 & &16.3.5.11.13.23.37.53&1961&2024& & &32.3.125.11.23.311&253&240 \cr
\noalign{\hrule}
 & &5.13.19.151.263&32103&32858& & &9.49.11.43.239&2125&2132 \cr
2307&49045555&4.27.7.29.41.2347&83&286&2325&49853727&8.125.7.13.17.41.239&3041&66 \cr
 & &16.3.11.13.83.2347&7041&7304& & &32.3.5.11.41.3041&3041&3280 \cr
\noalign{\hrule}
 & &9.625.11.13.61&1273&602& & &3.5.71.173.271&133&138 \cr
2308&49066875&4.3.7.13.19.43.67&5&8&2326&49930395&4.9.7.19.23.71.173&925&286 \cr
 & &64.5.7.19.43.67&2881&4256& & &16.25.11.13.19.23.37&9361&9880 \cr
\noalign{\hrule}
 & &3.5.17.199.967&5561&8944& & &49.11.19.67.73&2675&7566 \cr
2309&49070415&32.13.43.67.83&901&1980&2327&50088731&4.3.25.13.97.107&453&938 \cr
 & &256.9.5.11.17.53&583&384& & &16.9.5.7.67.151&151&360 \cr
\noalign{\hrule}
 & &5.7.11.29.53.83&355&558& & &5.121.41.43.47&361&156 \cr
2310&49114835&4.9.25.31.53.71&709&616&2328&50130905&8.3.11.13.361.43&227&246 \cr
 & &64.3.7.11.71.709&2127&2272& & &32.9.13.19.41.227&2951&2736 \cr
\noalign{\hrule}
 & &9.5.13.83.1013&157&92& & &9.19.23.53.241&143&580 \cr
2311&49186215&8.3.23.157.1013&3325&286&2329&50236209&8.3.5.11.13.29.53&109&268 \cr
 & &32.25.7.11.13.19&133&880& & &64.5.11.67.109&7303&1760 \cr
\noalign{\hrule}
 & &5.31.43.83.89&711&622& & &121.19.47.467&37&84 \cr
2312&49234355&4.9.5.79.83.311&5611&946&2330&50460751&8.3.7.19.37.467&585&118 \cr
 & &16.3.11.31.43.181&181&264& & &32.27.5.7.13.59&1593&7280 \cr
\noalign{\hrule}
 & &9.11.17.19.23.67&1255&286& & &3.11.23.101.659&9245&5912 \cr
2313&49276557&4.3.5.121.13.251&437&316&2331&50518281&16.5.1849.739&923&2772 \cr
 & &32.5.13.19.23.79&395&208& & &128.9.7.11.13.71&1491&832 \cr
\noalign{\hrule}
 & &3.11.13.19.73.83&387&1190& & &11.17.23.61.193&1305&1976 \cr
2314&49386909&4.27.5.7.13.17.43&185&374&2332&50635673&16.9.5.13.19.23.29&193&106 \cr
 & &16.25.11.289.37&925&2312& & &64.3.5.19.53.193&795&608 \cr
\noalign{\hrule}
 & &27.5.11.29.31.37&301&598& & &3.49.11.17.19.97&103&730 \cr
2315&49395555&4.5.7.13.23.37.43&71&114&2333&50662227&4.5.73.97.103&209&306 \cr
 & &16.3.7.13.19.23.71&5681&3976& & &16.9.11.17.19.73&73&24 \cr
\noalign{\hrule}
 & &27.11.43.53.73&16549&17980& & &81.11.19.41.73&221&148 \cr
2316&49410999&8.5.13.19.29.31.67&43&24&2334&50668497&8.9.11.13.17.19.37&995&292 \cr
 & &128.3.5.13.29.31.43&1885&1984& & &64.5.17.73.199&995&544 \cr
\noalign{\hrule}
 & &27.7.11.13.31.59&149&500& & &7.17.23.97.191&6185&9432 \cr
2317&49432383&8.125.7.31.149&171&46&2335&50708399&16.9.5.131.1237&291&946 \cr
 & &32.9.19.23.149&437&2384& & &64.27.11.43.97&473&864 \cr
\noalign{\hrule}
 & &5.49.13.103.151&59&696& & &5.11.23.67.599&311&426 \cr
2318&49536305&16.3.29.59.103&81&22&2336&50768245&4.3.71.311.599&335&264 \cr
 & &64.243.11.29&7047&352& & &64.9.5.11.67.311&311&288 \cr
\noalign{\hrule}
 & &9.5.7.23.41.167&2883&3718& & &3.5.11.13.19.29.43&41&206 \cr
2319&49606515&4.27.11.169.961&305&656&2337&50821485&4.29.41.43.103&2375&612 \cr
 & &128.5.11.13.41.61&671&832& & &32.9.125.17.19&425&48 \cr
\noalign{\hrule}
 & &27.5.7.11.17.281&241&1726& & &9.125.11.23.179&551&574 \cr
2320&49656915&4.17.241.863&1617&2480&2338&50947875&4.7.11.19.29.41.179&3711&1742 \cr
 & &128.3.5.49.11.31&217&64& & &16.3.13.29.67.1237&16081&15544 \cr
\noalign{\hrule}
 & &81.7.11.13.613&25&38& & &5.11.169.23.239&1737&892 \cr
2321&49702653&4.9.25.11.19.613&3281&2236&2339&51094615&8.9.23.193.223&2275&2854 \cr
 & &32.5.13.17.43.193&3655&3088& & &32.3.25.7.13.1427&1427&1680 \cr
\noalign{\hrule}
 & &27.5.121.17.179&10829&10830& & &27.11.13.17.19.41&61&308 \cr
2322&49707405&4.81.25.49.13.289.361&1&2024&2340&51131223&8.3.7.121.17.61&3055&3116 \cr
 & &64.7.11.13.361.23&4693&5152& & &64.5.7.13.19.41.47&235&224 \cr
\noalign{\hrule}
}%
}
$$
\eject
\vglue -23 pt
\noindent\hskip 1 in\hbox to 6.5 in{\ 2341 -- 2376 \hfill\fbd 51141519 -- 53798625\frb}
\vskip -9 pt
$$
\vbox{
\nointerlineskip
\halign{\strut
    \vrule \ \ \hfil \frb #\ 
   &\vrule \hfil \ \ \fbb #\frb\ 
   &\vrule \hfil \ \ \frb #\ \hfil
   &\vrule \hfil \ \ \frb #\ 
   &\vrule \hfil \ \ \frb #\ \ \vrule \hskip 2 pt
   &\vrule \ \ \hfil \frb #\ 
   &\vrule \hfil \ \ \fbb #\frb\ 
   &\vrule \hfil \ \ \frb #\ \hfil
   &\vrule \hfil \ \ \frb #\ 
   &\vrule \hfil \ \ \frb #\ \vrule \cr%
\noalign{\hrule}
 & &9.11.13.79.503&563&464& & &343.17.89.101&75&26 \cr
2341&51141519&32.29.503.563&533&30&2359&52414859&4.3.25.7.13.17.89&3737&3828 \cr
 & &128.3.5.13.29.41&1189&320& & &32.9.5.11.29.37.101&2871&2960 \cr
\noalign{\hrule}
 & &3.5.49.13.53.101&201&254& & &9.5.343.43.79&241&4 \cr
2342&51147915&4.9.7.67.101.127&689&220&2360&52432695&8.3.7.43.241&2509&2552 \cr
 & &32.5.11.13.53.127&127&176& & &128.11.13.29.193&5597&9152 \cr
\noalign{\hrule}
 & &9.125.7.67.97&507&172& & &5.17.61.67.151&45&106 \cr
2343&51179625&8.27.25.169.43&743&418&2361&52456645&4.9.25.17.53.67&671&604 \cr
 & &32.11.13.19.743&8173&3952& & &32.3.11.53.61.151&159&176 \cr
\noalign{\hrule}
 & &9.7.11.17.19.229&109&10& & &125.37.59.193&2431&2394 \cr
2344&51259131&4.5.19.109.229&105&124&2362&52664875&4.9.5.7.11.13.17.19.59&37&258 \cr
 & &32.3.25.7.31.109&775&1744& & &16.27.7.11.19.37.43&3591&3784 \cr
\noalign{\hrule}
 & &7.23.29.101.109&9711&6550& & &81.11.13.23.199&31&112 \cr
2345&51401021&4.9.25.13.83.131&11&404&2363&53015391&32.7.23.31.199&19&180 \cr
 & &32.3.5.11.13.101&715&48& & &256.9.5.19.31&155&2432 \cr
\noalign{\hrule}
 & &27.25.7.11.23.43&8381&8174& & &5.7.11.23.53.113&453&338 \cr
2346&51403275&4.3.5.289.29.61.67&3139&1196&2364&53032595&4.3.11.169.53.151&1175&486 \cr
 & &32.13.23.43.61.73&949&976& & &16.729.25.13.47&3645&4888 \cr
\noalign{\hrule}
 & &9.49.11.13.19.43&373&100& & &25.11.41.53.89&2133&2584 \cr
2347&51522471&8.3.25.7.19.373&227&892&2365&53184175&16.27.25.17.19.79&41&1384 \cr
 & &64.5.223.227&1115&7264& & &256.9.41.173&173&1152 \cr
\noalign{\hrule}
 & &3.5.11.59.67.79&12259&14206& & &49.19.211.271&571&360 \cr
2348&51527355&4.13.23.41.7103&4023&3080&2366&53235511&16.9.5.271.571&121&692 \cr
 & &64.27.5.7.11.13.149&1937&2016& & &128.3.5.121.173&2595&7744 \cr
\noalign{\hrule}
 & &81.7.11.17.487&659&650& & &9.5.7.131.1291&65&66 \cr
2349&51636123&4.9.25.13.487.659&3839&544&2367&53273115&4.27.25.7.11.13.1291&131&8906 \cr
 & &256.5.11.13.17.349&1745&1664& & &16.11.61.73.131&803&488 \cr
\noalign{\hrule}
 & &9.5.7.41.4003&307&308& & &9.25.37.43.149&959&7366 \cr
2350&51698745&8.3.49.11.307.4003&44581&548&2368&53338275&4.7.29.127.137&2431&1542 \cr
 & &64.109.137.409&14933&13088& & &16.3.11.13.17.257&4369&1144 \cr
\noalign{\hrule}
 & &5.121.17.47.107&7657&2628& & &5.49.37.71.83&7479&1586 \cr
2351&51723265&8.9.13.19.31.73&35&22&2369&53420045&4.27.13.61.277&385&446 \cr
 & &32.3.5.7.11.31.73&651&1168& & &16.9.5.7.11.13.223&1287&1784 \cr
\noalign{\hrule}
 & &3.5.11.13.361.67&1271&534& & &5.7.11.19.71.103&899&234 \cr
2352&51881115&4.9.13.31.41.89&667&134&2370&53494595&4.9.13.29.31.71&95&308 \cr
 & &16.23.29.31.67&899&184& & &32.3.5.7.11.19.29&87&16 \cr
\noalign{\hrule}
 & &37.683.2053&11609&13662& & &3.11.73.97.229&109&182 \cr
2353&51881363&4.27.11.13.19.23.47&685&74&2371&53511117&4.7.11.13.109.229&873&730 \cr
 & &16.9.5.19.37.137&685&1368& & &16.9.5.73.97.109&109&120 \cr
\noalign{\hrule}
 & &9.11.29.67.271&185&86& & &9.5.11.13.53.157&451&334 \cr
2354&52128747&4.5.29.37.43.67&503&570&2372&53545635&4.121.41.53.167&217&6630 \cr
 & &16.3.25.19.43.503&9557&8600& & &16.3.5.7.13.17.31&527&56 \cr
\noalign{\hrule}
 & &81.5.11.13.17.53&97&124& & &5.343.11.17.167&639&296 \cr
2355&52181415&8.3.5.11.31.53.97&931&136&2373&53557735&16.9.37.71.167&215&286 \cr
 & &128.49.17.19.31&931&1984& & &64.3.5.11.13.37.43&1591&1248 \cr
\noalign{\hrule}
 & &7.11.13.17.37.83&100757&101016& & &9.5.7.11.13.29.41&113&92 \cr
2356&52259207&16.9.19.23.61.5303&913&4390&2374&53558505&8.3.11.13.23.29.113&119&548 \cr
 & &64.3.5.11.23.83.439&2195&2208& & &64.7.17.113.137&2329&3616 \cr
\noalign{\hrule}
 & &17.19.29.37.151&11477&8910& & &5.13.19.157.277&1159&882 \cr
2357&52333429&4.81.5.11.23.499&629&130&2375&53708915&4.9.5.49.361.61&277&638 \cr
 & &16.27.25.13.17.37&351&200& & &16.3.49.11.29.277&539&696 \cr
\noalign{\hrule}
 & &3.5.7.11.113.401&1387&986& & &9.125.17.29.97&637&1012 \cr
2358&52336515&4.5.11.17.19.29.73&1143&98&2376&53798625&8.3.49.11.13.23.29&347&550 \cr
 & &16.9.49.29.127&381&1624& & &32.25.7.121.347&2429&1936 \cr
\noalign{\hrule}
}%
}
$$
\eject
\vglue -23 pt
\noindent\hskip 1 in\hbox to 6.5 in{\ 2377 -- 2412 \hfill\fbd 53836025 -- 55368885\frb}
\vskip -9 pt
$$
\vbox{
\nointerlineskip
\halign{\strut
    \vrule \ \ \hfil \frb #\ 
   &\vrule \hfil \ \ \fbb #\frb\ 
   &\vrule \hfil \ \ \frb #\ \hfil
   &\vrule \hfil \ \ \frb #\ 
   &\vrule \hfil \ \ \frb #\ \ \vrule \hskip 2 pt
   &\vrule \ \ \hfil \frb #\ 
   &\vrule \hfil \ \ \fbb #\frb\ 
   &\vrule \hfil \ \ \frb #\ \hfil
   &\vrule \hfil \ \ \frb #\ 
   &\vrule \hfil \ \ \frb #\ \vrule \cr%
\noalign{\hrule}
 & &25.17.19.59.113&121&444& & &9.5.47.149.173&755&802 \cr
2377&53836025&8.3.5.121.37.59&351&56&2395&54518355&4.25.149.151.401&6237&16262 \cr
 & &128.81.7.11.13&7371&704& & &16.81.7.11.47.173&77&72 \cr
\noalign{\hrule}
 & &7.17.43.53.199&10701&12980& & &125.13.19.29.61&803&822 \cr
2378&53968999&8.9.5.11.29.41.59&221&98&2396&54617875&4.3.11.29.61.73.137&131&1638 \cr
 & &32.3.5.49.13.17.59&1239&1040& & &16.27.7.13.73.131&3537&4088 \cr
\noalign{\hrule}
 & &9.5.11.13.37.227&527&1970& & &9.5.13.223.419&229&190 \cr
2379&54047565&4.3.25.17.31.197&259&166&2397&54660645&4.3.25.19.223.229&3493&9218 \cr
 & &16.7.37.83.197&1379&664& & &16.7.11.419.499&499&616 \cr
\noalign{\hrule}
 & &3.25.47.67.229&2147&1002& & &3.5.11.13.17.19.79&65&122 \cr
2380&54084075&4.9.5.19.113.167&469&1034&2398&54733965&4.25.169.61.79&297&4522 \cr
 & &16.7.11.19.47.67&77&152& & &16.27.7.11.17.19&7&72 \cr
\noalign{\hrule}
 & &27.5.49.13.17.37&10637&12932& & &81.5.17.73.109&3497&3388 \cr
2381&54090855&8.11.53.61.967&819&148&2399&54783945&8.7.121.13.73.269&9&82 \cr
 & &64.9.7.13.37.53&53&32& & &32.9.121.41.269&4961&4304 \cr
\noalign{\hrule}
 & &3.25.7.11.17.19.29&83&468& & &7.11.13.19.43.67&5&138 \cr
2382&54094425&8.27.5.13.17.83&43&178&2400&54793739&4.3.5.23.43.67&461&528 \cr
 & &32.43.83.89&7387&688& & &128.9.5.11.461&2305&576 \cr
\noalign{\hrule}
 & &9.5.11.17.59.109&133&428& & &27.5.13.17.19.97&4393&4328 \cr
2383&54116865&8.3.7.19.107.109&83&26&2401&54985905&16.23.97.191.541&845&1386 \cr
 & &32.7.13.83.107&7553&1712& & &64.9.5.7.11.169.191&2483&2464 \cr
\noalign{\hrule}
 & &3.25.43.97.173&287&578& & &25.7.13.361.67&211&36 \cr
2384&54118725&4.5.7.289.41.43&873&572&2402&55025425&8.9.19.67.211&2015&1804 \cr
 & &32.9.11.13.41.97&451&624& & &64.3.5.11.13.31.41&1353&992 \cr
\noalign{\hrule}
 & &5.19.47.53.229&561&332& & &9.13.53.83.107&469&220 \cr
2385&54191705&8.3.5.11.17.53.83&2709&1796&2403&55071081&8.3.5.7.11.67.107&689&488 \cr
 & &64.27.7.43.449&12123&9632& & &128.5.7.13.53.61&427&320 \cr
\noalign{\hrule}
 & &3.49.17.529.41&377&1210& & &3.25.19.29.31.43&429&904 \cr
2386&54200811&4.5.121.13.29.41&37&414&2404&55086225&16.9.11.13.29.113&1075&2318 \cr
 & &16.9.5.11.23.37&55&888& & &64.25.19.43.61&61&32 \cr
\noalign{\hrule}
 & &25.29.37.43.47&22361&21114& & &25.7.11.23.29.43&221&1026 \cr
2387&54213325&4.27.17.23.59.379&1247&110&2405&55210925&4.27.5.11.13.17.19&23&232 \cr
 & &16.9.5.11.17.29.43&99&136& & &64.9.13.23.29&13&288 \cr
\noalign{\hrule}
 & &3.25.13.19.29.101&123&22& & &125.7.17.47.79&1763&7638 \cr
2388&54259725&4.9.5.11.13.19.41&101&146&2406&55230875&4.3.19.41.43.67&55&12 \cr
 & &16.11.41.73.101&803&328& & &32.9.5.11.19.41&369&3344 \cr
\noalign{\hrule}
 & &3.7.13.19.37.283&639&1342& & &27.25.7.11.1063&59&1004 \cr
2389&54313077&4.27.11.13.61.71&215&566&2407&55249425&8.5.11.59.251&153&98 \cr
 & &16.5.43.61.283&215&488& & &32.9.49.17.59&119&944 \cr
\noalign{\hrule}
 & &13.361.67.173&723&550& & &9.7.13.19.53.67&1067&1000 \cr
2390&54396563&4.3.25.11.13.19.241&1943&708&2408&55257111&16.3.125.7.11.19.97&43&62 \cr
 & &32.9.5.29.59.67&1305&944& & &64.25.11.31.43.97&33325&34144 \cr
\noalign{\hrule}
 & &3.5.17.37.73.79&14473&14362& & &9.7.13.19.53.67&33325&34144 \cr
2391&54411645&4.17.41.43.167.353&423&6424&2409&55257111&64.25.11.31.43.97&21&76 \cr
 & &64.9.11.43.47.73&1419&1504& & &512.3.5.7.19.31.43&1333&1280 \cr
\noalign{\hrule}
 & &3.5.49.11.53.127&3139&3084& & &9.5.7.13.23.587&289&298 \cr
2392&54420135&8.9.43.53.73.257&17&2296&2410&55286595&4.5.7.13.289.23.149&5149&66 \cr
 & &128.7.17.41.73&2993&1088& & &16.3.11.17.19.271&4607&1672 \cr
\noalign{\hrule}
 & &3.13.37.67.563&717&154& & &5.7.19.29.47.61&895&468 \cr
2393&54431403&4.9.7.11.37.239&1003&670&2411&55290095&8.9.25.13.19.179&4939&5264 \cr
 & &16.5.11.17.59.67&1003&440& & &256.3.7.11.47.449&1347&1408 \cr
\noalign{\hrule}
 & &7.23.31.67.163&27&4& & &3.5.11.13.83.311&551&694 \cr
2394&54506711&8.27.7.67.163&671&470&2412&55368885&4.19.29.311.347&329&18 \cr
 & &32.9.5.11.47.61&2867&7920& & &16.9.7.19.29.47&987&4408 \cr
\noalign{\hrule}
}%
}
$$
\eject
\vglue -23 pt
\noindent\hskip 1 in\hbox to 6.5 in{\ 2413 -- 2448 \hfill\fbd 55371225 -- 57470523\frb}
\vskip -9 pt
$$
\vbox{
\nointerlineskip
\halign{\strut
    \vrule \ \ \hfil \frb #\ 
   &\vrule \hfil \ \ \fbb #\frb\ 
   &\vrule \hfil \ \ \frb #\ \hfil
   &\vrule \hfil \ \ \frb #\ 
   &\vrule \hfil \ \ \frb #\ \ \vrule \hskip 2 pt
   &\vrule \ \ \hfil \frb #\ 
   &\vrule \hfil \ \ \fbb #\frb\ 
   &\vrule \hfil \ \ \frb #\ \hfil
   &\vrule \hfil \ \ \frb #\ 
   &\vrule \hfil \ \ \frb #\ \vrule \cr%
\noalign{\hrule}
 & &3.25.49.13.19.61&1081&444& & &3.13.17.23.29.127&3395&3014 \cr
2413&55371225&8.9.19.23.37.47&385&52&2431&56162067&4.5.7.11.23.97.137&2159&5310 \cr
 & &64.5.7.11.13.47&47&352& & &16.9.25.17.59.127&177&200 \cr
\noalign{\hrule}
 & &121.37.89.139&2813&7956& & &27.7.11.17.37.43&23&40 \cr
2414&55384967&8.9.13.17.29.97&259&1390&2432&56230713&16.3.5.11.23.37.43&119&526 \cr
 & &32.3.5.7.37.139&105&16& & &64.7.17.23.263&263&736 \cr
\noalign{\hrule}
 & &9.7.11.127.631&7471&530& & &3.5.7.11.23.29.73&3&32 \cr
2415&55534941&4.5.31.53.241&147&94&2433&56238105&64.9.11.23.73&1175&1102 \cr
 & &16.3.5.49.31.47&1645&248& & &256.25.19.29.47&893&640 \cr
\noalign{\hrule}
 & &11.13.17.73.313&9&8& & &81.11.17.47.79&305&118 \cr
2416&55545919&16.9.11.13.73.313&3239&830&2434&56240811&4.9.5.59.61.79&35&44 \cr
 & &64.3.5.41.79.83&10209&12640& & &32.25.7.11.59.61&3599&2800 \cr
\noalign{\hrule}
 & &3.961.37.521&741&220& & &125.11.13.23.137&603&1178 \cr
2417&55575591&8.9.5.11.13.19.37&271&62&2435&56324125&4.9.5.11.19.31.67&481&146 \cr
 & &32.5.13.31.271&1355&208& & &16.3.13.31.37.73&2263&888 \cr
\noalign{\hrule}
 & &3.7.47.199.283&123&76& & &49.121.13.17.43&635&96 \cr
2418&55584879&8.9.7.19.41.283&1397&3980&2436&56343287&64.3.5.11.13.127&421&294 \cr
 & &64.5.11.127.199&635&352& & &256.9.49.421&421&1152 \cr
\noalign{\hrule}
 & &3.25.11.13.71.73&451&1374& & &25.7.11.19.23.67&1327&858 \cr
2419&55587675&4.9.121.41.229&245&124&2437&56362075&4.3.5.121.13.1327&361&966 \cr
 & &32.5.49.31.229&1519&3664& & &16.9.7.13.361.23&171&104 \cr
\noalign{\hrule}
 & &9.19.29.103.109&905&1166& & &9.25.23.67.163&119&44 \cr
2420&55674693&4.5.11.53.103.181&39&142&2438&56516175&8.3.7.11.17.23.67&275&208 \cr
 & &16.3.5.11.13.53.71&3763&5720& & &256.25.121.13.17&1573&2176 \cr
\noalign{\hrule}
 & &25.13.19.71.127&99&28& & &625.29.31.101&1503&1628 \cr
2421&55679975&8.9.25.7.11.13.19&233&508&2439&56749375&8.9.5.11.29.37.167&697&5482 \cr
 & &64.3.7.127.233&699&224& & &32.3.17.41.2741&8223&11152 \cr
\noalign{\hrule}
 & &5.49.17.43.311&261&572& & &9.11.13.19.23.101&235&64 \cr
2422&55698545&8.9.5.11.13.29.43&119&76&2440&56804319&128.5.11.47.101&23&78 \cr
 & &64.3.7.11.17.19.29&627&928& & &512.3.13.23.47&47&256 \cr
\noalign{\hrule}
 & &25.7.19.103.163&10057&6732& & &9.7.17.29.31.59&15&44 \cr
2423&55823425&8.9.11.17.89.113&19&70&2441&56806911&8.27.5.7.11.17.31&1073&236 \cr
 & &32.3.5.7.11.19.113&339&176& & &64.5.29.37.59&185&32 \cr
\noalign{\hrule}
 & &9.5.11.13.19.457&1057&1228& & &9.5.11.13.37.239&191&524 \cr
2424&55875105&8.7.11.13.151.307&457&600&2442&56904705&8.131.191.239&185&54 \cr
 & &128.3.25.307.457&307&320& & &32.27.5.37.191&191&48 \cr
\noalign{\hrule}
 & &11.13.529.739&1275&2014& & &81.125.7.11.73&263&628 \cr
2425&55903133&4.3.25.17.19.23.53&169&222&2443&56912625&8.25.7.157.263&219&44 \cr
 & &16.9.25.169.19.37&4329&3800& & &64.3.11.73.157&157&32 \cr
\noalign{\hrule}
 & &5.31.277.1303&667&636& & &3.5.19.29.61.113&715&1054 \cr
2426&55944305&8.3.5.23.29.53.277&211&66&2444&56970645&4.25.11.13.17.19.31&4041&3616 \cr
 & &32.9.11.23.53.211&20889&19504& & &256.9.11.113.449&1347&1408 \cr
\noalign{\hrule}
 & &9.5.37.151.223&16929&16744& & &11.19.29.9409&7735&1674 \cr
2427&56065545&16.729.7.11.13.19.23&3791&4228&2445&57027949&4.27.5.7.13.17.31&97&58 \cr
 & &128.49.13.17.151.223&833&832& & &16.9.7.17.29.97&153&56 \cr
\noalign{\hrule}
 & &3.13.37.59.659&241&418& & &9.49.31.59.71&4741&8930 \cr
2428&56105283&4.11.13.19.37.241&2655&2636&2446&57267819&4.5.11.19.47.431&221&210 \cr
 & &32.9.5.59.241.659&241&240& & &16.3.25.7.13.17.19.47&6175&6392 \cr
\noalign{\hrule}
 & &81.11.61.1033&3211&8152& & &121.43.73.151&7663&1170 \cr
2429&56144583&16.169.19.1019&425&594&2447&57352669&4.9.5.13.79.97&151&86 \cr
 & &64.27.25.11.17.19&425&608& & &16.3.43.97.151&97&24 \cr
\noalign{\hrule}
 & &3.25.343.37.59&551&374& & &3.121.17.67.139&195&1334 \cr
2430&56157675&4.343.11.17.19.29&333&10&2448&57470523&4.9.5.11.13.23.29&139&238 \cr
 & &16.9.5.11.29.37&319&24& & &16.5.7.17.23.139&161&40 \cr
\noalign{\hrule}
}%
}
$$
\eject
\vglue -23 pt
\noindent\hskip 1 in\hbox to 6.5 in{\ 2449 -- 2484 \hfill\fbd 57547035 -- 59705217\frb}
\vskip -9 pt
$$
\vbox{
\nointerlineskip
\halign{\strut
    \vrule \ \ \hfil \frb #\ 
   &\vrule \hfil \ \ \fbb #\frb\ 
   &\vrule \hfil \ \ \frb #\ \hfil
   &\vrule \hfil \ \ \frb #\ 
   &\vrule \hfil \ \ \frb #\ \ \vrule \hskip 2 pt
   &\vrule \ \ \hfil \frb #\ 
   &\vrule \hfil \ \ \fbb #\frb\ 
   &\vrule \hfil \ \ \frb #\ \hfil
   &\vrule \hfil \ \ \frb #\ 
   &\vrule \hfil \ \ \frb #\ \vrule \cr%
\noalign{\hrule}
 & &9.5.7.169.23.47&319&526& & &9.11.67.83.107&689&488 \cr
2449&57547035&4.7.11.29.47.263&33&296&2467&58907673&16.3.13.53.61.83&469&220 \cr
 & &64.3.121.29.37&4477&928& & &128.5.7.11.61.67&427&320 \cr
\noalign{\hrule}
 & &9.5.11.17.41.167&157&344& & &11.19.521.541&255&266 \cr
2450&57617505&16.3.5.41.43.157&167&38&2468&58908949&4.3.5.7.17.361.541&451&90 \cr
 & &64.19.157.167&157&608& & &16.27.25.7.11.17.41&4879&5400 \cr
\noalign{\hrule}
 & &3.49.13.19.37.43&33&670& & &3.5.7.19.109.271&1221&676 \cr
2451&57767619&4.9.5.11.43.67&1463&1418&2469&58930305&8.9.11.169.19.37&271&62 \cr
 & &16.7.121.19.709&709&968& & &32.169.31.271&169&496 \cr
\noalign{\hrule}
 & &27.13.23.67.107&19&88& & &27.11.13.17.29.31&265&112 \cr
2452&57875337&16.9.11.13.19.67&25&92&2470&59007663&32.3.5.7.11.31.53&1399&244 \cr
 & &128.25.11.19.23&209&1600& & &256.61.1399&1399&7808 \cr
\noalign{\hrule}
 & &9.11.59.61.163&235&296& & &9.7.23.83.491&8891&4472 \cr
2453&58077063&16.5.11.37.47.163&177&340&2471&59051097&16.13.17.43.523&253&270 \cr
 & &128.3.25.17.37.59&925&1088& & &64.27.5.11.13.23.43&1677&1760 \cr
\noalign{\hrule}
 & &9.5.11.19.37.167&6509&10024& & &25.7.29.103.113&351&2926 \cr
2454&58113495&16.7.23.179.283&231&52&2472&59067925&4.27.49.11.13.19&215&226 \cr
 & &128.3.49.11.13.23&637&1472& & &16.3.5.13.19.43.113&559&456 \cr
\noalign{\hrule}
 & &3.25.121.13.17.29&577&698& & &27.17.23.29.193&2959&1480 \cr
2455&58161675&4.13.29.349.577&8811&1310&2473&59087529&16.9.5.11.37.269&193&2228 \cr
 & &16.9.5.11.89.131&393&712& & &128.193.557&557&64 \cr
\noalign{\hrule}
 & &9.17.31.71.173&517&690& & &9.7.11.17.29.173&43&130 \cr
2456&58258269&4.27.5.11.23.31.47&3287&4556&2474&59105277&4.3.5.7.11.13.17.43&173&184 \cr
 & &32.5.17.19.67.173&335&304& & &64.5.13.23.43.173&1495&1376 \cr
\noalign{\hrule}
 & &27.19.197.577&203&374& & &25.11.17.47.269&177&92 \cr
2457&58312197&4.3.7.11.17.29.197&577&380&2475&59106025&8.3.5.11.23.47.59&2019&754 \cr
 & &32.5.7.17.19.577&85&112& & &32.9.13.29.673&8749&4176 \cr
\noalign{\hrule}
 & &125.11.59.719&297&422& & &11.19.23.31.397&8957&4590 \cr
2458&58328875&4.27.121.59.211&13031&12500&2476&59159749&4.27.5.169.17.53&31&190 \cr
 & &32.3.3125.83.157&3925&3984& & &16.9.25.13.19.31&325&72 \cr
\noalign{\hrule}
 & &9.25.49.11.13.37&779&446& & &9.5.7.13.31.467&19&26 \cr
2459&58333275&4.11.13.19.41.223&3213&5930&2477&59283315&4.169.19.31.467&5633&8844 \cr
 & &16.27.5.7.17.593&593&408& & &32.3.11.43.67.131&8777&7568 \cr
\noalign{\hrule}
 & &9.5.7.19.41.239&407&372& & &5.7.11.23.37.181&397&408 \cr
2460&58647015&8.27.11.31.37.239&9025&182&2478&59301935&16.3.17.37.181.397&6723&26 \cr
 & &32.25.7.13.361&65&304& & &64.243.13.83&1079&7776 \cr
\noalign{\hrule}
 & &11.13.17.23.1049&453&596& & &27.5.11.13.17.181&85&58 \cr
2461&58652737&8.3.17.23.149.151&3003&470&2479&59401485&4.25.289.29.181&6237&988 \cr
 & &32.9.5.7.11.13.47&329&720& & &32.81.7.11.13.19&133&48 \cr
\noalign{\hrule}
 & &3.5.49.17.37.127&1267&638& & &9.5.23.137.419&35&172 \cr
2462&58714005&4.343.11.29.181&5969&3978&2480&59412105&8.25.7.43.419&297&122 \cr
 & &16.9.13.17.47.127&141&104& & &32.27.11.43.61&2623&528 \cr
\noalign{\hrule}
 & &3.5.7.11.17.41.73&73&158& & &3.25.7.19.59.101&193&92 \cr
2463&58767555&4.41.5329.79&1045&4284&2481&59441025&8.5.7.23.59.193&2717&4068 \cr
 & &32.9.5.7.11.17.19&57&16& & &64.9.11.13.19.113&1469&1056 \cr
\noalign{\hrule}
 & &81.7.97.1069&187&7670& & &11.43.71.1777&1125&652 \cr
2464&58793931&4.5.11.13.17.59&389&378&2482&59676991&8.9.125.71.163&969&806 \cr
 & &16.27.5.7.17.389&389&680& & &32.27.5.13.17.19.31&15903&17680 \cr
\noalign{\hrule}
 & &9.5.13.17.61.97&539&722& & &27.11.169.29.41&125&658 \cr
2465&58844565&4.3.5.49.11.17.361&61&194&2483&59679477&4.125.7.11.13.47&421&96 \cr
 & &16.7.11.19.61.97&133&88& & &256.3.5.7.421&421&4480 \cr
\noalign{\hrule}
 & &9.13.17.101.293&41&262& & &9.11.13.23.2017&23839&22552 \cr
2466&58860477&4.3.41.131.293&505&374&2484&59705217&16.31.769.2819&1025&1794 \cr
 & &16.5.11.17.41.101&205&88& & &64.3.25.13.23.31.41&1025&992 \cr
\noalign{\hrule}
}%
}
$$
\eject
\vglue -23 pt
\noindent\hskip 1 in\hbox to 6.5 in{\ 2485 -- 2520 \hfill\fbd 59781085 -- 61490715\frb}
\vskip -9 pt
$$
\vbox{
\nointerlineskip
\halign{\strut
    \vrule \ \ \hfil \frb #\ 
   &\vrule \hfil \ \ \fbb #\frb\ 
   &\vrule \hfil \ \ \frb #\ \hfil
   &\vrule \hfil \ \ \frb #\ 
   &\vrule \hfil \ \ \frb #\ \ \vrule \hskip 2 pt
   &\vrule \ \ \hfil \frb #\ 
   &\vrule \hfil \ \ \fbb #\frb\ 
   &\vrule \hfil \ \ \frb #\ \hfil
   &\vrule \hfil \ \ \frb #\ 
   &\vrule \hfil \ \ \frb #\ \vrule \cr%
\noalign{\hrule}
 & &5.7.13.37.53.67&10561&7194& & &27.49.11.23.181&1273&2890 \cr
2485&59781085&4.3.11.59.109.179&25&84&2503&60584139&4.9.5.289.19.67&1267&1334 \cr
 & &32.9.25.7.11.179&895&1584& & &16.5.7.19.23.29.181&145&152 \cr
\noalign{\hrule}
 & &3.5.13.29.71.149&209&238& & &9.11.37.73.227&61&742 \cr
2486&59824245&4.5.7.11.13.17.19.71&943&162&2504&60699573&4.3.7.37.53.61&73&110 \cr
 & &16.81.7.19.23.41&4347&6232& & &16.5.7.11.53.73&53&280 \cr
\noalign{\hrule}
 & &27.17.101.1291&2795&1078& & &27.49.17.37.73&101&2600 \cr
2487&59849469&4.9.5.49.11.13.43&101&200&2505&60748191&16.9.25.13.101&343&242 \cr
 & &64.125.7.13.101&875&416& & &64.5.343.121&121&1120 \cr
\noalign{\hrule}
 & &3.5.13.37.43.193&547&418& & &3.5.23.353.499&77&422 \cr
2488&59877285&4.11.13.19.37.547&1561&7578&2506&60770715&4.7.11.211.353&2047&1836 \cr
 & &16.9.7.223.421&2947&5352& & &32.27.7.17.23.89&1071&1424 \cr
\noalign{\hrule}
 & &125.13.19.29.67&1813&942& & &11.23.59.61.67&1431&2168 \cr
2489&59990125&4.3.25.49.37.157&87&1012&2507&61006649&16.27.23.53.271&1829&610 \cr
 & &32.9.7.11.23.29&1449&176& & &64.3.5.31.59.61&93&160 \cr
\noalign{\hrule}
 & &9.7.13.61.1201&635&7772& & &5.17.29.53.467&217&684 \cr
2490&60000759&8.5.29.67.127&61&66&2508&61011215&8.9.5.7.19.29.31&121&34 \cr
 & &32.3.11.29.61.67&737&464& & &32.3.7.121.17.19&399&1936 \cr
\noalign{\hrule}
 & &27.25.7.11.13.89&701&1702& & &81.17.23.41.47&137&560 \cr
2491&60135075&4.25.23.37.701&63&638&2509&61030017&32.9.5.7.23.137&583&376 \cr
 & &16.9.7.11.29.37&29&296& & &512.5.11.47.53&583&1280 \cr
\noalign{\hrule}
 & &3.7.11.19.71.193&4625&4818& & &3.13.73.89.241&165&76 \cr
2492&60142467&4.9.125.121.37.73&893&772&2510&61065303&8.9.5.11.13.19.73&241&124 \cr
 & &32.25.19.47.73.193&1175&1168& & &64.11.19.31.241&341&608 \cr
\noalign{\hrule}
 & &9.7.13.23.31.103&3165&1826& & &9.11.179.3449&5263&5084 \cr
2493&60146541&4.27.5.11.83.211&109&26&2511&61119729&8.3.11.19.31.41.277&20275&13796 \cr
 & &16.11.13.109.211&2321&872& & &64.25.811.3449&811&800 \cr
\noalign{\hrule}
 & &81.5.7.13.23.71&61&152& & &7.121.17.31.137&743&216 \cr
2494&60184215&16.27.5.19.23.61&71&44&2512&61152553&16.27.121.743&553&190 \cr
 & &128.11.19.61.71&1159&704& & &64.9.5.7.19.79&711&3040 \cr
\noalign{\hrule}
 & &5.7.17.23.53.83&3177&1222& & &81.7.13.43.193&187&380 \cr
2495&60200315&4.9.7.13.47.353&187&166&2513&61171929&8.5.11.13.17.19.43&113&360 \cr
 & &16.3.11.13.17.47.83&429&376& & &128.9.25.17.113&2825&1088 \cr
\noalign{\hrule}
 & &9.7.11.233.373&443&676& & &3.5.11.13.361.79&833&668 \cr
2496&60227937&8.3.7.11.169.443&1165&164&2514&61173255&8.49.13.17.19.167&207&40 \cr
 & &64.5.13.41.233&533&160& & &128.9.5.49.17.23&2499&1472 \cr
\noalign{\hrule}
 & &27.19.41.47.61&2311&3470& & &3.169.23.29.181&539&358 \cr
2497&60301611&4.9.5.347.2311&2717&406&2515&61208589&4.49.11.13.29.179&1265&1062 \cr
 & &16.5.7.11.13.19.29&1885&616& & &16.9.5.7.121.23.59&2541&2360 \cr
\noalign{\hrule}
 & &27.7.13.41.599&115&484& & &3.25.7.11.13.19.43&2227&2162 \cr
2498&60341463&8.3.5.7.121.13.23&599&116&2516&61336275&4.5.17.23.43.47.131&693&38 \cr
 & &64.11.29.599&319&32& & &16.9.7.11.19.23.47&141&184 \cr
\noalign{\hrule}
 & &125.13.73.509&7871&1254& & &9.19.31.71.163&209&280 \cr
2499&60380125&4.3.11.17.19.463&325&138&2517&61348473&16.3.5.7.11.361.31&413&52 \cr
 & &16.9.25.13.19.23&437&72& & &128.49.11.13.59&8437&3136 \cr
\noalign{\hrule}
 & &9.17.37.59.181&2847&3850& & &169.43.79.107&9063&9020 \cr
2500&60453819&4.27.25.7.11.13.73&181&116&2518&61427951&8.9.5.11.19.41.53.79&20683&20604 \cr
 & &32.5.7.29.73.181&1015&1168& & &64.27.5.11.13.17.37.43.101&54945&54944 \cr
\noalign{\hrule}
 & &3.5.11.13.89.317&203&114& & &9.49.11.19.23.29&289&338 \cr
2501&60516885&4.9.5.7.11.13.19.29&2219&1174&2519&61476723&4.3.169.289.23.29&665&2 \cr
 & &16.49.317.587&587&392& & &16.5.7.13.17.19&13&680 \cr
\noalign{\hrule}
 & &9.7.11.23.29.131&387&1054& & &3.5.11.13.109.263&667&122 \cr
2502&60552261&4.81.7.17.31.43&305&262&2520&61490715&4.11.13.23.29.61&2043&2104 \cr
 & &16.5.17.31.61.131&1037&1240& & &64.9.23.227.263&681&736 \cr
\noalign{\hrule}
}%
}
$$
\eject
\vglue -23 pt
\noindent\hskip 1 in\hbox to 6.5 in{\ 2521 -- 2556 \hfill\fbd 61516455 -- 62984375\frb}
\vskip -9 pt
$$
\vbox{
\nointerlineskip
\halign{\strut
    \vrule \ \ \hfil \frb #\ 
   &\vrule \hfil \ \ \fbb #\frb\ 
   &\vrule \hfil \ \ \frb #\ \hfil
   &\vrule \hfil \ \ \frb #\ 
   &\vrule \hfil \ \ \frb #\ \ \vrule \hskip 2 pt
   &\vrule \ \ \hfil \frb #\ 
   &\vrule \hfil \ \ \fbb #\frb\ 
   &\vrule \hfil \ \ \frb #\ \hfil
   &\vrule \hfil \ \ \frb #\ 
   &\vrule \hfil \ \ \frb #\ \vrule \cr%
\noalign{\hrule}
 & &3.5.7.11.13.17.241&817&388& & &9.343.113.179&20185&18574 \cr
2521&61516455&8.7.17.19.43.97&11&108&2539&62440749&4.5.11.37.251.367&113&13692 \cr
 & &64.27.11.19.43&171&1376& & &32.3.7.113.163&163&16 \cr
\noalign{\hrule}
 & &5.7.19.37.41.61&2613&2308& & &3.25.7.11.29.373&3&8 \cr
2522&61537105&8.3.13.41.67.577&555&22&2540&62468175&16.9.5.7.29.373&19&1846 \cr
 & &32.9.5.11.37.67&603&176& & &64.13.19.71&71&7904 \cr
\noalign{\hrule}
 & &9.7.13.19.37.107&1199&7940& & &7.11.47.61.283&123&2990 \cr
2523&61605999&8.5.11.109.397&401&798&2541&62474797&4.3.5.7.13.23.41&517&426 \cr
 & &32.3.5.7.19.401&401&80& & &16.9.5.11.47.71&355&72 \cr
\noalign{\hrule}
 & &5.17.23.139.227&2277&1582& & &3.5.7.23.41.631&235&396 \cr
2524&61686115&4.9.7.11.529.113&737&850&2542&62478465&8.27.25.11.41.47&1141&34 \cr
 & &16.3.25.7.121.17.67&2541&2680& & &32.7.11.17.163&163&2992 \cr
\noalign{\hrule}
 & &9.25.13.47.449&14663&14522& & &9.25.7.97.409&371&856 \cr
2525&61726275&4.3.5.11.31.43.53.137&73&338&2543&62484975&16.3.5.49.53.107&341&194 \cr
 & &16.11.169.31.43.73&10439&10664& & &64.11.31.53.97&583&992 \cr
\noalign{\hrule}
 & &3.5.7.19.47.659&319&340& & &9.25.193.1439&2167&728 \cr
2526&61791135&8.25.11.17.19.29.47&659&234&2544&62488575&16.3.5.7.11.13.197&607&772 \cr
 & &32.9.11.13.29.659&429&464& & &128.13.193.607&607&832 \cr
\noalign{\hrule}
 & &13.19.23.73.149&1837&1590& & &3.5.49.11.59.131&259&390 \cr
2527&61792237&4.3.5.11.53.73.167&4389&4462&2545&62488965&4.9.25.343.13.37&59&284 \cr
 & &16.9.5.7.121.19.23.97&5445&5432& & &32.13.37.59.71&481&1136 \cr
\noalign{\hrule}
 & &3.13.19.29.43.67&253&124& & &3.13.79.103.197&10925&13486 \cr
2528&61909809&8.11.19.23.31.67&3285&4558&2546&62516571&4.25.11.19.23.613&1751&1314 \cr
 & &32.9.5.43.53.73&1095&848& & &16.9.5.11.17.73.103&1241&1320 \cr
\noalign{\hrule}
 & &3.25.7.13.29.313&319&6& & &27.7.11.17.29.61&32651&32680 \cr
2529&61950525&4.9.7.11.841&767&74&2547&62521767&16.3.5.11.19.43.103.317&391&82 \cr
 & &16.13.37.59&2183&8& & &64.5.17.19.23.41.317&30115&30176 \cr
\noalign{\hrule}
 & &3.11.41.163.281&3325&3358& & &3.169.17.53.137&1341&440 \cr
2530&61971459&4.25.7.19.23.73.281&1757&3582&2548&62582559&16.27.5.11.13.149&197&548 \cr
 & &16.9.49.23.199.251&36897&36616& & &128.11.137.197&197&704 \cr
\noalign{\hrule}
 & &3.125.7.13.23.79&339&214& & &25.7.11.13.41.61&10541&6966 \cr
2531&62005125&4.9.13.23.107.113&11809&10340&2549&62587525&4.81.43.83.127&85&42 \cr
 & &32.5.49.11.47.241&3619&3856& & &16.243.5.7.17.83&1411&1944 \cr
\noalign{\hrule}
 & &27.529.43.101&12737&10010& & &11.5329.1069&133&936 \cr
2532&62031069&4.5.7.11.13.47.271&69&22&2550&62663711&16.9.7.13.19.73&643&890 \cr
 & &16.3.5.121.23.271&1355&968& & &64.3.5.89.643&3215&8544 \cr
\noalign{\hrule}
 & &9.25.11.13.1931&37&62& & &9.11.13.29.1681&10767&11086 \cr
2533&62129925&4.13.31.37.1931&13125&11978&2551&62739963&4.27.23.37.97.241&1189&1430 \cr
 & &16.3.625.7.53.113&2825&2968& & &16.5.11.13.23.29.37.41&185&184 \cr
\noalign{\hrule}
 & &9.25.37.59.127&67&8& & &19.23.31.41.113&275&162 \cr
2534&62379225&16.3.37.67.127&209&172&2552&62763251&4.81.25.11.31.41&57&98 \cr
 & &128.11.19.43.67&8987&4288& & &16.243.5.49.11.19&2673&1960 \cr
\noalign{\hrule}
 & &9.25.37.59.127&209&172& & &3.11.17.23.31.157&8651&3470 \cr
2535&62379225&8.3.25.11.19.43.59&67&8&2553&62798901&4.5.41.211.347&279&68 \cr
 & &128.11.19.43.67&8987&4288& & &32.9.5.17.31.41&205&48 \cr
\noalign{\hrule}
 & &27.125.11.1681&503&1178& & &9.125.11.13.17.23&29&1654 \cr
2536&62407125&4.5.11.19.31.503&123&218&2554&62902125&4.23.29.827&425&402 \cr
 & &16.3.41.109.503&503&872& & &16.3.25.17.29.67&67&232 \cr
\noalign{\hrule}
 & &9.31.467.479&5083&9394& & &5.7.19.31.43.71&4039&2706 \cr
2537&62410347&4.7.11.13.17.23.61&485&186&2555&62937595&4.3.49.11.41.577&19&558 \cr
 & &16.3.5.7.17.31.97&485&952& & &16.27.19.31.41&27&328 \cr
\noalign{\hrule}
 & &9.5.11.13.31.313&1691&1378& & &15625.29.139&5797&9828 \cr
2538&62438805&4.5.169.19.53.89&2781&1936&2556&62984375&8.27.7.11.13.17.31&125&278 \cr
 & &128.27.121.19.103&1957&2112& & &32.3.125.7.11.139&77&48 \cr
\noalign{\hrule}
}%
}
$$
\eject
\vglue -23 pt
\noindent\hskip 1 in\hbox to 6.5 in{\ 2557 -- 2592 \hfill\fbd 62987873 -- 65103555\frb}
\vskip -9 pt
$$
\vbox{
\nointerlineskip
\halign{\strut
    \vrule \ \ \hfil \frb #\ 
   &\vrule \hfil \ \ \fbb #\frb\ 
   &\vrule \hfil \ \ \frb #\ \hfil
   &\vrule \hfil \ \ \frb #\ 
   &\vrule \hfil \ \ \frb #\ \ \vrule \hskip 2 pt
   &\vrule \ \ \hfil \frb #\ 
   &\vrule \hfil \ \ \fbb #\frb\ 
   &\vrule \hfil \ \ \frb #\ \hfil
   &\vrule \hfil \ \ \frb #\ 
   &\vrule \hfil \ \ \frb #\ \vrule \cr%
\noalign{\hrule}
 & &13.17.257.1109&7337&7080& & &25.7.11.29.31.37&111&34 \cr
2557&62987873&16.3.5.11.17.23.29.59&7217&7710&2575&64031275&4.3.5.17.31.1369&633&2002 \cr
 & &64.9.25.7.257.1031&7217&7200& & &16.9.7.11.13.211&211&936 \cr
\noalign{\hrule}
 & &3.7.13.23.79.127&41&340& & &13.29.43.59.67&1455&1426 \cr
2558&62997207&8.5.7.17.41.79&319&234&2576&64082083&4.3.5.13.23.31.59.97&231&172 \cr
 & &32.9.11.13.29.41&1353&464& & &32.9.5.7.11.23.43.97&11155&11088 \cr
\noalign{\hrule}
 & &5.11.61.89.211&1205&1116& & &53.89.107.127&27&80 \cr
2559&63003545&8.9.25.31.61.241&7&68&2577&64099313&32.27.5.89.127&451&184 \cr
 & &64.3.7.17.31.241&11067&7712& & &512.9.11.23.41&8487&2816 \cr
\noalign{\hrule}
 & &9.5.17.41.2011&8947&9152& & &3.5.11.41.53.179&639&1534 \cr
2560&63075015&128.11.13.17.23.389&4679&378&2578&64179555&4.27.11.13.59.71&4025&164 \cr
 & &512.27.7.4679&4679&5376& & &32.25.7.23.41&23&560 \cr
\noalign{\hrule}
 & &25.289.59.149&783&8008& & &3.5.49.19.43.107&367&382 \cr
2561&63514975&16.27.7.11.13.29&85&118&2579&64252965&4.7.19.43.191.367&4811&20592 \cr
 & &64.9.5.13.17.59&117&32& & &128.9.11.13.17.283&11037&11968 \cr
\noalign{\hrule}
 & &25.7.11.17.29.67&117&202& & &5.17.29.89.293&22671&21206 \cr
2562&63584675&4.9.5.7.13.67.101&451&956&2580&64279805&4.9.11.23.229.461&703&680 \cr
 & &32.3.11.13.41.239&3107&1968& & &64.3.5.11.17.19.37.229&13053&13024 \cr
\noalign{\hrule}
 & &5.19.23.37.787&7293&10808& & &27.11.53.61.67&115&182 \cr
2563&63625015&16.3.7.11.13.17.193&3&190&2581&64333467&4.5.7.13.23.53.61&8987&9252 \cr
 & &64.9.5.7.13.19&819&32& & &32.9.7.11.19.43.257&4883&4816 \cr
\noalign{\hrule}
 & &81.5.7.11.13.157&29507&28408& & &3.49.11.13.37.83&2033&2034 \cr
2564&63648585&16.19.53.67.1553&1413&140&2582&64555491&4.27.11.13.19.37.107.113&1411&17570 \cr
 & &128.9.5.7.53.157&53&64& & &16.5.7.17.83.107.251&4267&4280 \cr
\noalign{\hrule}
 & &25.7.19.41.467&169&36& & &9.25.7.17.19.127&671&1742 \cr
2565&63663775&8.9.5.169.467&231&236&2583&64608075&4.25.11.13.61.67&1131&394 \cr
 & &64.27.7.11.169.59&9971&9504& & &16.3.169.29.197&4901&1576 \cr
\noalign{\hrule}
 & &27.7.31.83.131&3151&910& & &3.5.7.31.103.193&481&484 \cr
2566&63704907&4.5.49.13.23.137&93&44&2584&64706145&8.7.121.13.31.37.103&5085&36424 \cr
 & &32.3.5.11.13.23.31&299&880& & &128.9.5.29.113.157&9831&10048 \cr
\noalign{\hrule}
 & &25.7.11.79.419&7293&7372& & &81.5.7.11.31.67&2561&3298 \cr
2567&63719425&8.3.5.121.13.17.19.97&27243&31442&2585&64771245&4.3.5.13.17.97.197&217&268 \cr
 & &32.81.79.199.1009&16119&16144& & &32.7.13.31.67.197&197&208 \cr
\noalign{\hrule}
 & &5.11.23.29.37.47&167&684& & &3.25.7.169.17.43&551&44 \cr
2568&63795215&8.9.5.19.29.167&1127&1628&2586&64857975&8.5.11.19.29.43&117&98 \cr
 & &64.3.49.11.23.37&49&96& & &32.9.49.11.13.29&203&528 \cr
\noalign{\hrule}
 & &3.5.71.139.431&611&682& & &3.5.19.29.47.167&1127&1628 \cr
2569&63803085&4.5.11.13.31.47.139&3447&862&2587&64871985&8.49.11.23.37.47&167&684 \cr
 & &16.9.13.383.431&383&312& & &64.9.49.19.167&49&96 \cr
\noalign{\hrule}
 & &9.13.31.73.241&1073&1190& & &25.7.19.109.179&177&2 \cr
2570&63809811&4.5.7.17.29.37.241&1139&66&2588&64874075&4.3.19.59.109&615&506 \cr
 & &16.3.7.11.289.67&5159&2312& & &16.9.5.11.23.41&207&3608 \cr
\noalign{\hrule}
 & &3.47.167.2711&1285&1426& & &27.11.13.17.23.43&265&724 \cr
2571&63835917&4.5.23.31.167.257&367&5544&2589&64914993&8.5.11.13.53.181&989&1002 \cr
 & &64.9.5.7.11.367&2569&5280& & &32.3.5.23.43.53.167&835&848 \cr
\noalign{\hrule}
 & &9.7.11.13.41.173&893&1010& & &3.5.49.11.29.277&139&106 \cr
2572&63900837&4.5.7.19.41.47.101&1677&242&2590&64946805&4.29.53.139.277&333&7700 \cr
 & &16.3.121.13.43.47&517&344& & &32.9.25.7.11.37&37&240 \cr
\noalign{\hrule}
 & &5.7.19.23.37.113&99&62& & &5.13.53.113.167&57&110 \cr
2573&63948395&4.9.5.11.19.31.113&851&2094&2591&65010595&4.3.25.11.13.19.113&219&106 \cr
 & &16.27.23.37.349&349&216& & &16.9.11.19.53.73&803&1368 \cr
\noalign{\hrule}
 & &9.11.13.17.37.79&437&590& & &3.5.11.229.1723&18031&19754 \cr
2574&63952317&4.5.11.19.23.37.59&101&750&2592&65103555&4.7.13.17.19.73.83&1485&244 \cr
 & &16.3.625.19.101&1919&5000& & &32.27.5.11.61.83&747&976 \cr
\noalign{\hrule}
}%
}
$$
\eject
\vglue -23 pt
\noindent\hskip 1 in\hbox to 6.5 in{\ 2593 -- 2628 \hfill\fbd 65117533 -- 67222155\frb}
\vskip -9 pt
$$
\vbox{
\nointerlineskip
\halign{\strut
    \vrule \ \ \hfil \frb #\ 
   &\vrule \hfil \ \ \fbb #\frb\ 
   &\vrule \hfil \ \ \frb #\ \hfil
   &\vrule \hfil \ \ \frb #\ 
   &\vrule \hfil \ \ \frb #\ \ \vrule \hskip 2 pt
   &\vrule \ \ \hfil \frb #\ 
   &\vrule \hfil \ \ \fbb #\frb\ 
   &\vrule \hfil \ \ \frb #\ \hfil
   &\vrule \hfil \ \ \frb #\ 
   &\vrule \hfil \ \ \frb #\ \vrule \cr%
\noalign{\hrule}
 & &13.59.73.1163&965&198& & &7.151.167.373&715&342 \cr
2593&65117533&4.9.5.11.73.193&107&472&2611&65841587&4.9.5.11.13.19.167&2611&106 \cr
 & &64.3.11.59.107&107&1056& & &16.3.7.53.373&159&8 \cr
\noalign{\hrule}
 & &5.13.17.19.29.107&161&162& & &3.5.7.11.23.37.67&31&438 \cr
2594&65147485&4.81.5.7.13.23.29.107&6479&476&2612&65854635&4.9.5.23.31.73&79&286 \cr
 & &32.9.49.11.17.19.31&1519&1584& & &16.11.13.31.79&2449&104 \cr
\noalign{\hrule}
 & &9.121.19.47.67&1079&1220& & &27.7.121.43.67&299&2 \cr
2595&65155959&8.3.5.13.61.67.83&893&22&2613&65885589&4.11.13.23.67&441&430 \cr
 & &32.11.19.47.83&83&16& & &16.9.5.49.23.43&161&40 \cr
\noalign{\hrule}
 & &27.25.13.17.19.23&177&398& & &3.5.7.17.19.29.67&309&242 \cr
2596&65189475&4.81.19.59.199&869&670&2614&65896845&4.9.5.7.121.17.103&36461&35534 \cr
 & &16.5.11.59.67.79&4661&5896& & &16.361.101.109.163&16463&16568 \cr
\noalign{\hrule}
 & &25.11.13.17.29.37&473&252& & &9.7.11.23.41.101&65&188 \cr
2597&65211575&8.9.7.121.37.43&5&116&2615&66003399&8.3.5.7.13.47.101&943&1178 \cr
 & &64.3.5.7.29.43&43&672& & &32.13.19.23.31.41&403&304 \cr
\noalign{\hrule}
 & &3.7.67.151.307&375&682& & &7.13.23.101.313&10165&6096 \cr
2598&65224299&4.9.125.11.31.67&307&28&2616&66166009&32.3.5.19.107.127&939&1474 \cr
 & &32.25.7.11.307&25&176& & &128.9.11.67.313&737&576 \cr
\noalign{\hrule}
 & &7.11.13.19.47.73&639&310& & &3.5.7.13.139.349&53&402 \cr
2599&65254189&4.9.5.11.19.31.71&493&1274&2617&66217515&4.9.53.67.139&169&308 \cr
 & &16.3.5.49.13.17.29&435&952& & &32.7.11.169.67&871&176 \cr
\noalign{\hrule}
 & &27.5.13.19.37.53&173&308& & &3.121.31.71.83&8323&1720 \cr
2600&65389545&8.7.11.19.53.173&1935&1352&2618&66313929&16.5.7.29.41.43&117&88 \cr
 & &128.9.5.7.169.43&559&448& & &256.9.7.11.13.43&1677&896 \cr
\noalign{\hrule}
 & &9.25.11.13.19.107&683&708& & &3.43.47.97.113&33&80 \cr
2601&65411775&8.27.11.19.59.683&3055&15386&2619&66456543&32.9.5.11.43.97&1501&434 \cr
 & &32.5.49.13.47.157&2303&2512& & &128.7.19.31.79&2449&8512 \cr
\noalign{\hrule}
 & &3.25.7.11.17.23.29&67&186& & &5.7.11.13.53.251&3103&342 \cr
2602&65482725&4.9.25.29.31.67&4301&2224&2620&66581515&4.9.7.19.29.107&251&358 \cr
 & &128.11.17.23.139&139&64& & &16.3.19.179.251&537&152 \cr
\noalign{\hrule}
 & &3.5.7.11.17.47.71&533&674& & &11.89.149.457&4117&9144 \cr
2603&65521995&4.5.7.11.13.41.337&1343&342&2621&66663047&16.9.23.127.179&1371&1550 \cr
 & &16.9.17.19.41.79&1501&984& & &64.27.25.31.457&775&864 \cr
\noalign{\hrule}
 & &81.11.251.293&523&230& & &3.11.19.31.47.73&481&1070 \cr
2604&65526813&4.27.5.11.23.523&481&1004&2622&66688347&4.5.13.37.73.107&423&58 \cr
 & &32.13.23.37.251&299&592& & &16.9.29.47.107&107&696 \cr
\noalign{\hrule}
 & &27.5.11.19.23.101&383&130& & &9.13.23.137.181&1705&76 \cr
2605&65543445&4.25.13.101.383&63&38&2623&66728727&8.5.11.19.23.31&351&362 \cr
 & &16.9.7.13.19.383&383&728& & &32.27.5.13.19.181&57&80 \cr
\noalign{\hrule}
 & &243.17.59.269&433&374& & &5.11.23.47.1123&731&1854 \cr
2606&65563101&4.81.11.289.433&301&590&2624&66767965&4.9.17.23.43.103&2019&2410 \cr
 & &16.5.7.43.59.433&2165&2408& & &16.27.5.241.673&6507&5384 \cr
\noalign{\hrule}
 & &59.10201.109&1885&8316& & &5.11.13.41.43.53&333&118 \cr
2607&65602631&8.27.5.7.11.13.29&101&218&2625&66808885&4.9.13.37.53.59&161&320 \cr
 & &32.3.5.7.101.109&21&80& & &512.3.5.7.23.59&4071&1792 \cr
\noalign{\hrule}
 & &7.121.13.59.101&8165&1026& & &3.5.11.13.529.59&1177&410 \cr
2608&65614549&4.27.5.19.23.71&101&538&2626&66947595&4.25.121.41.107&1143&3818 \cr
 & &16.3.5.101.269&807&40& & &16.9.23.83.127&249&1016 \cr
\noalign{\hrule}
 & &729.5.13.19.73&839&110& & &3.5.11.47.89.97&2625&1558 \cr
2609&65722995&4.25.11.19.839&657&182&2627&66948915&4.9.625.7.19.41&979&1604 \cr
 & &16.9.7.11.13.73&7&88& & &32.11.19.89.401&401&304 \cr
\noalign{\hrule}
 & &3.7.17.19.31.313&19&12& & &3.5.7.1331.13.37&323&158 \cr
2610&65815449&8.9.17.361.313&4477&1660&2628&67222155&4.7.121.17.19.79&3551&5850 \cr
 & &64.5.121.37.83&10043&5920& & &16.9.25.13.53.67&795&536 \cr
\noalign{\hrule}
}%
}
$$
\eject
\vglue -23 pt
\noindent\hskip 1 in\hbox to 6.5 in{\ 2629 -- 2664 \hfill\fbd 67297545 -- 70079961\frb}
\vskip -9 pt
$$
\vbox{
\nointerlineskip
\halign{\strut
    \vrule \ \ \hfil \frb #\ 
   &\vrule \hfil \ \ \fbb #\frb\ 
   &\vrule \hfil \ \ \frb #\ \hfil
   &\vrule \hfil \ \ \frb #\ 
   &\vrule \hfil \ \ \frb #\ \ \vrule \hskip 2 pt
   &\vrule \ \ \hfil \frb #\ 
   &\vrule \hfil \ \ \fbb #\frb\ 
   &\vrule \hfil \ \ \frb #\ \hfil
   &\vrule \hfil \ \ \frb #\ 
   &\vrule \hfil \ \ \frb #\ \vrule \cr%
\noalign{\hrule}
 & &9.5.7.29.53.139&629&484& & &27.25.17.43.139&329&746 \cr
2629&67297545&8.3.121.17.37.139&1363&1000&2647&68586075&4.9.7.17.47.373&263&110 \cr
 & &128.125.29.37.47&1739&1600& & &16.5.7.11.47.263&1841&4136 \cr
\noalign{\hrule}
 & &81.13.17.53.71&5&76& & &9.5.241.6337&3191&3146 \cr
2630&67361463&8.5.13.17.19.53&243&22&2648&68724765&4.121.13.241.3191&2921&270 \cr
 & &32.243.11.19&11&912& & &16.27.5.11.13.23.127&3289&3048 \cr
\noalign{\hrule}
 & &9.5.7.121.29.61&103&158& & &3.7.23.181.787&10951&7150 \cr
2631&67425435&4.7.11.61.79.103&25&696&2649&68801901&4.25.11.13.47.233&189&422 \cr
 & &64.3.25.29.79&395&32& & &16.27.25.7.11.211&2475&1688 \cr
\noalign{\hrule}
 & &3.13.17.19.53.101&2515&2838& & &7.11.13.89.773&5977&6750 \cr
2632&67431741&4.9.5.11.13.43.503&251&5282&2650&68865797&4.27.125.7.43.139&383&3092 \cr
 & &16.5.19.139.251&695&2008& & &32.3.5.383.773&383&240 \cr
\noalign{\hrule}
 & &9.11.13.23.2281&15941&13660& & &5.11.169.17.19.23&36935&36918 \cr
2633&67519881&8.5.19.683.839&761&78&2651&69052555&4.9.25.7.11.83.89.293&24397&78 \cr
 & &32.3.5.13.19.761&761&1520& & &16.27.7.13.31.787&5859&6296 \cr
\noalign{\hrule}
 & &3.7.11.13.113.199&703&690& & &3.7.11.13.109.211&109&122 \cr
2634&67528461&4.9.5.11.19.23.37.113&65&2212&2652&69065997&4.61.11881.211&495&12376 \cr
 & &32.25.7.13.37.79&925&1264& & &64.9.5.7.11.13.17&51&160 \cr
\noalign{\hrule}
 & &9.125.19.29.109&923&1148& & &3.5.11.169.37.67&219&626 \cr
2635&67566375&8.5.7.13.29.41.71&513&2398&2653&69126915&4.9.67.73.313&265&338 \cr
 & &32.27.7.11.19.109&77&48& & &16.5.169.53.313&313&424 \cr
\noalign{\hrule}
 & &25.13.43.47.103&11&36& & &3.5.13.37.43.223&1509&1286 \cr
2636&67652975&8.9.11.13.43.103&125&434&2654&69184635&4.9.37.503.643&85&418 \cr
 & &32.3.125.7.11.31&341&1680& & &16.5.11.17.19.643&3553&5144 \cr
\noalign{\hrule}
 & &9.5.7.31.53.131&47&418& & &5.29.383.1249&2431&8676 \cr
2637&67798395&4.3.11.19.47.131&3475&3992&2655&69363215&8.9.11.13.17.241&433&290 \cr
 & &64.25.139.499&2495&4448& & &32.3.5.17.29.433&433&816 \cr
\noalign{\hrule}
 & &81.5.7.71.337&209&146& & &11.169.107.349&763&414 \cr
2638&67833045&4.9.11.19.73.337&497&160&2656&69420637&4.9.7.169.23.109&755&428 \cr
 & &256.5.7.11.19.71&209&128& & &32.3.5.23.107.151&755&1104 \cr
\noalign{\hrule}
 & &27.11.17.89.151&1267&1300& & &5.49.11.19.23.59&8983&3348 \cr
2639&67853511&8.9.25.7.13.89.181&19&604&2657&69485185&8.27.13.31.691&259&950 \cr
 & &64.5.19.151.181&905&608& & &32.9.25.7.19.37&333&80 \cr
\noalign{\hrule}
 & &5.7.13.23.43.151&1089&874& & &9.25.7.13.43.79&1573&1652 \cr
2640&67949245&4.9.7.121.19.529&3001&702&2658&69553575&8.3.49.121.169.59&11&158 \cr
 & &16.243.13.3001&3001&1944& & &32.1331.59.79&1331&944 \cr
\noalign{\hrule}
 & &27.7.13.17.23.71&275&184& & &3.25.11.37.43.53&677&1602 \cr
2641&68208777&16.25.11.529.71&87&442&2659&69566475&4.27.11.89.677&3737&3710 \cr
 & &64.3.5.11.13.17.29&145&352& & &16.5.7.37.53.89.101&707&712 \cr
\noalign{\hrule}
 & &9.5.49.19.23.71&89&44& & &243.5.23.47.53&6557&6322 \cr
2642&68414535&8.7.11.23.71.89&775&1272&2660&69610995&4.23.29.79.83.109&2025&4532 \cr
 & &128.3.25.11.31.53&2915&1984& & &32.81.25.11.29.103&1595&1648 \cr
\noalign{\hrule}
 & &25.17.233.691&2961&8786& & &5.7.11.17.29.367&1077&758 \cr
2643&68426275&4.9.7.23.47.191&3497&4070&2661&69658435&4.3.7.17.359.379&249&130 \cr
 & &16.3.5.11.13.37.269&5291&6456& & &16.9.5.13.83.359&4667&5976 \cr
\noalign{\hrule}
 & &9.121.37.1699&30887&31976& & &3.89.311.839&553&286 \cr
2644&68457807&16.7.67.461.571&55&516&2662&69668043&4.7.11.13.79.311&575&1602 \cr
 & &128.3.5.7.11.43.67&2345&2752& & &16.9.25.11.23.89&575&264 \cr
\noalign{\hrule}
 & &9.11.13.17.31.101&353&50& & &3.5.11.13.17.19.101&261&244 \cr
2645&68503149&4.3.25.11.17.353&9139&8864&2663&69976335&8.27.11.13.19.29.61&7675&13736 \cr
 & &256.13.19.37.277&5263&4736& & &128.25.17.101.307&307&320 \cr
\noalign{\hrule}
 & &5.7.13.31.43.113&297&262& & &3.7.19.37.47.101&9581&4660 \cr
2646&68536195&4.27.11.31.113.131&3913&410&2664&70079961&8.5.11.13.67.233&1717&846 \cr
 & &16.9.5.7.13.41.43&41&72& & &32.9.5.17.47.101&85&48 \cr
\noalign{\hrule}
}%
}
$$
\eject
\vglue -23 pt
\noindent\hskip 1 in\hbox to 6.5 in{\ 2665 -- 2700 \hfill\fbd 70080153 -- 72557485\frb}
\vskip -9 pt
$$
\vbox{
\nointerlineskip
\halign{\strut
    \vrule \ \ \hfil \frb #\ 
   &\vrule \hfil \ \ \fbb #\frb\ 
   &\vrule \hfil \ \ \frb #\ \hfil
   &\vrule \hfil \ \ \frb #\ 
   &\vrule \hfil \ \ \frb #\ \ \vrule \hskip 2 pt
   &\vrule \ \ \hfil \frb #\ 
   &\vrule \hfil \ \ \fbb #\frb\ 
   &\vrule \hfil \ \ \frb #\ \hfil
   &\vrule \hfil \ \ \frb #\ 
   &\vrule \hfil \ \ \frb #\ \vrule \cr%
\noalign{\hrule}
 & &3.11.13.29.43.131&83&476& & &3.5.11.37.89.131&269&138 \cr
2665&70080153&8.7.11.17.29.83&131&450&2683&71178195&4.9.5.23.89.269&1073&272 \cr
 & &32.9.25.17.131&425&48& & &128.17.23.29.37&493&1472 \cr
\noalign{\hrule}
 & &5.11.13.17.29.199&393&712& & &27.37.43.1657&23165&21574 \cr
2666&70146505&16.3.89.131.199&165&34&2684&71179749&4.5.7.23.41.67.113&4257&3314 \cr
 & &64.9.5.11.17.89&89&288& & &16.9.5.7.11.43.1657&55&56 \cr
\noalign{\hrule}
 & &27.49.29.31.59&1507&86& & &7.53.401.479&425&54 \cr
2667&70173243&4.11.31.43.137&735&598&2685&71261309&4.27.25.17.401&583&182 \cr
 & &16.3.5.49.11.13.23&715&184& & &16.3.5.7.11.13.53&39&440 \cr
\noalign{\hrule}
 & &3.5.11.19.61.367&497&130& & &9.25.49.11.19.31&227&172 \cr
2668&70183245&4.25.7.13.61.71&3303&1028&2686&71430975&8.3.5.7.31.43.227&893&242 \cr
 & &32.9.257.367&257&48& & &32.121.19.43.47&473&752 \cr
\noalign{\hrule}
 & &9.11.31.89.257&1553&1274& & &5.79.113.1601&1083&518 \cr
2669&70197237&4.49.13.89.1553&465&1088&2687&71460635&4.3.7.361.37.79&141&220 \cr
 & &512.3.5.7.13.17.31&1547&1280& & &32.9.5.7.11.37.47&3663&5264 \cr
\noalign{\hrule}
 & &3.5.11.13.71.461&2997&2074& & &9.11.13.19.37.79&1655&952 \cr
2670&70207995&4.243.5.17.37.61&293&922&2688&71476119&16.3.5.7.13.17.331&2923&4028 \cr
 & &16.61.293.461&293&488& & &128.19.37.53.79&53&64 \cr
\noalign{\hrule}
 & &3.7.11.13.41.571&17&24& & &9.5.7.11.23.29.31&247&218 \cr
2671&70303233&16.9.11.13.17.571&4553&2870&2689&71645805&4.3.7.11.13.19.23.109&841&358 \cr
 & &64.5.7.29.41.157&785&928& & &16.13.19.841.179&3401&3016 \cr
\noalign{\hrule}
 & &7.17.43.59.233&121&180& & &27.289.29.317&89&8470 \cr
2672&70343399&8.9.5.121.17.233&817&118&2690&71732979&4.5.7.121.89&317&306 \cr
 & &32.3.11.19.43.59&57&176& & &16.9.5.11.17.317&11&40 \cr
\noalign{\hrule}
 & &9.5.7.19.61.193&137&442& & &13.41.157.859&2365&8802 \cr
2673&70461405&4.3.7.13.17.19.137&965&1364&2691&71881979&4.27.5.11.43.163&157&58 \cr
 & &32.5.11.13.31.193&341&208& & &16.3.29.157.163&163&696 \cr
\noalign{\hrule}
 & &27.49.11.29.167&4141&4042& & &121.19.23.29.47&8515&12024 \cr
2674&70480179&4.3.29.41.43.47.101&5357&3430&2692&72071351&16.9.5.13.131.167&601&1102 \cr
 & &16.5.343.11.43.487&2435&2408& & &64.3.5.19.29.601&601&480 \cr
\noalign{\hrule}
 & &27.5.37.103.137&121&806& & &5.11.169.19.409&1107&938 \cr
2675&70484445&4.3.121.13.31.37&157&250&2693&72231445&4.27.7.11.19.41.67&4225&4016 \cr
 & &16.125.11.13.157&3925&1144& & &128.9.25.7.169.251&2259&2240 \cr
\noalign{\hrule}
 & &3.5.343.47.293&25&318& & &9.25.11.289.101&2723&1612 \cr
2676&70851795&4.9.125.47.53&539&586&2694&72242775&8.3.5.7.13.31.389&187&202 \cr
 & &16.49.11.53.293&53&88& & &32.7.11.13.17.31.101&217&208 \cr
\noalign{\hrule}
 & &81.5.11.17.937&1157&220& & &3.5.13.17.19.31.37&439&1142 \cr
2677&70963695&8.25.121.13.89&1557&1468&2695&72243795&4.5.13.439.571&187&252 \cr
 & &64.9.13.173.367&4771&5536& & &32.9.7.11.17.571&1713&1232 \cr
\noalign{\hrule}
 & &3.5.7.11.13.29.163&19&184& & &43.89.127.149&9361&3900 \cr
2678&70975905&16.13.19.23.163&841&1278&2696&72418321&8.3.25.11.13.23.37&51&14 \cr
 & &64.9.841.71&213&928& & &32.9.5.7.11.17.23&8855&2448 \cr
\noalign{\hrule}
 & &27.59.109.409&695&286& & &9.25.7.11.37.113&359&884 \cr
2679&71017533&4.3.5.11.13.59.139&259&436&2697&72435825&8.3.13.17.37.359&853&224 \cr
 & &32.7.11.13.37.109&1001&592& & &512.7.13.853&853&3328 \cr
\noalign{\hrule}
 & &7.19.487.1097&177&310& & &9.25.13.17.31.47&487&722 \cr
2680&71053787&4.3.5.31.59.1097&1463&366&2698&72449325&4.3.5.17.361.487&6721&584 \cr
 & &16.9.5.7.11.19.61&305&792& & &64.11.13.47.73&73&352 \cr
\noalign{\hrule}
 & &3.125.7.11.23.107&585&592& & &9.5.49.11.29.103&215&226 \cr
2681&71061375&32.27.625.13.23.37&623&2&2699&72449685&4.25.29.43.103.113&351&2926 \cr
 & &128.7.13.37.89&3293&832& & &16.27.7.11.13.19.43&559&456 \cr
\noalign{\hrule}
 & &27.25.121.13.67&37&238& & &5.49.11.13.19.109&439&978 \cr
2682&71138925&4.9.7.11.13.17.37&103&40&2700&72557485&4.3.5.19.163.439&3763&4578 \cr
 & &64.5.17.37.103&3811&544& & &16.9.7.53.71.109&477&568 \cr
\noalign{\hrule}
}%
}
$$
\eject
\vglue -23 pt
\noindent\hskip 1 in\hbox to 6.5 in{\ 2701 -- 2736 \hfill\fbd 72662931 -- 74999925\frb}
\vskip -9 pt
$$
\vbox{
\nointerlineskip
\halign{\strut
    \vrule \ \ \hfil \frb #\ 
   &\vrule \hfil \ \ \fbb #\frb\ 
   &\vrule \hfil \ \ \frb #\ \hfil
   &\vrule \hfil \ \ \frb #\ 
   &\vrule \hfil \ \ \frb #\ \ \vrule \hskip 2 pt
   &\vrule \ \ \hfil \frb #\ 
   &\vrule \hfil \ \ \fbb #\frb\ 
   &\vrule \hfil \ \ \frb #\ \hfil
   &\vrule \hfil \ \ \frb #\ 
   &\vrule \hfil \ \ \frb #\ \vrule \cr%
\noalign{\hrule}
 & &9.11.37.83.239&221&2850& & &11.13.29.107.167&5229&10072 \cr
2701&72662931&4.27.25.13.17.19&239&220&2719&74102743&16.9.7.83.1259&505&754 \cr
 & &32.125.11.13.239&125&208& & &64.3.5.7.13.29.101&505&672 \cr
\noalign{\hrule}
 & &9.25.7.11.13.17.19&141&46& & &5.19.61.67.191&487&468 \cr
2702&72747675&4.27.5.7.13.23.47&113&22&2720&74158615&8.9.13.61.67.487&1337&1276 \cr
 & &16.11.23.47.113&2599&376& & &64.3.7.11.29.191.487&10227&10208 \cr
\noalign{\hrule}
 & &243.25.7.29.59&401&814& & &19.23.277.613&9009&2638 \cr
2703&72760275&4.5.11.29.37.401&2133&2278&2721&74203037&4.9.7.11.13.1319&445&874 \cr
 & &16.27.17.37.67.79&5293&5032& & &16.3.5.7.19.23.89&445&168 \cr
\noalign{\hrule}
 & &9.5.19.71.1201&629&572& & &19.67.173.337&5805&5786 \cr
2704&72906705&8.3.5.11.13.17.37.71&367&2402&2722&74217173&4.27.5.11.43.263.337&11591&44954 \cr
 & &32.17.367.1201&367&272& & &16.3.7.169.19.67.173&169&168 \cr
\noalign{\hrule}
 & &27.7.11.61.577&319&258& & &3.5.343.11.13.101&177&166 \cr
2705&73174563&4.81.7.121.29.43&2885&2318&2723&74309235&4.9.5.13.59.83.101&343&242 \cr
 & &16.5.19.29.61.577&145&152& & &16.343.121.59.83&649&664 \cr
\noalign{\hrule}
 & &27.17.31.37.139&4543&9686& & &3.5.11.17.41.647&349&298 \cr
2706&73179747&4.7.11.29.59.167&185&18&2724&74408235&4.5.11.41.149.349&647&1098 \cr
 & &16.9.5.11.37.59&55&472& & &16.9.61.149.647&447&488 \cr
\noalign{\hrule}
 & &27.5.7.13.59.101&31&382& & &9.49.17.19.523&1201&368 \cr
2707&73206315&4.5.31.101.191&427&528&2725&74497689&32.3.19.23.1201&55&1256 \cr
 & &128.3.7.11.31.61&671&1984& & &512.5.11.157&8635&256 \cr
\noalign{\hrule}
 & &5.7.67.157.199&133&66& & &27.11.17.29.509&385&124 \cr
2708&73264835&4.3.5.49.11.19.157&201&44&2726&74528289&8.3.5.7.121.17.31&241&122 \cr
 & &32.9.121.19.67&1089&304& & &32.5.31.61.241&9455&3856 \cr
\noalign{\hrule}
 & &7.13.23.101.347&605&8586& & &7.13.79.97.107&1037&8700 \cr
2709&73353371&4.81.5.121.53&13&2&2727&74614631&8.3.25.17.29.61&11&6 \cr
 & &16.27.11.13.53&53&2376& & &32.9.5.11.29.61&1595&8784 \cr
\noalign{\hrule}
 & &25.13.37.41.149&1431&506& & &5.29.31.37.449&349&798 \cr
2710&73460725&4.27.11.23.41.53&353&230&2728&74675435&4.3.5.7.19.29.349&333&682 \cr
 & &16.9.5.529.353&4761&2824& & &16.27.11.19.31.37&297&152 \cr
\noalign{\hrule}
 & &9.7.853.1367&5725&6578& & &19.23.29.71.83&16717&14310 \cr
2711&73461213&4.25.7.11.13.23.229&853&918&2729&74681989&4.27.5.53.73.229&217&3652 \cr
 & &16.27.5.17.229.853&687&680& & &32.9.7.11.31.83&341&1008 \cr
\noalign{\hrule}
 & &11.17.19.23.29.31&8909&2430& & &9.23.31.61.191&1049&842 \cr
2712&73465381&4.243.5.59.151&79&374&2730&74764467&4.191.421.1049&429&620 \cr
 & &16.81.11.17.79&81&632& & &32.3.5.11.13.31.421&2105&2288 \cr
\noalign{\hrule}
 & &3.25.11.29.37.83&5909&5894& & &11.17.19.37.569&8865&1946 \cr
2713&73473675&4.5.7.19.83.311.421&11&1566&2731&74801309&4.9.5.7.139.197&1037&342 \cr
 & &16.27.7.11.29.421&421&504& & &16.81.17.19.61&81&488 \cr
\noalign{\hrule}
 & &9.7.11.13.41.199&425&8584& & &3.13.31.103.601&1103&700 \cr
2714&73504431&16.25.17.29.37&339&154&2732&74840727&8.25.7.103.1103&2573&5148 \cr
 & &64.3.5.7.11.113&565&32& & &64.9.11.13.31.83&249&352 \cr
\noalign{\hrule}
 & &3.5.7.11.23.47.59&153&176& & &9.13.67.73.131&575&374 \cr
2715&73664745&32.27.5.121.17.59&4141&874&2733&74964357&4.3.25.11.17.23.131&169&224 \cr
 & &128.19.23.41.101&1919&2624& & &256.5.7.169.17.23&5083&4480 \cr
\noalign{\hrule}
 & &3.125.7.13.17.127&23&198& & &5.7.19.137.823&891&68 \cr
2716&73675875&4.27.5.11.23.127&169&466&2734&74979415&8.81.5.11.17.19&259&254 \cr
 & &16.169.23.233&233&2392& & &32.3.7.11.17.37.127&6919&6096 \cr
\noalign{\hrule}
 & &7.361.31.941&793&1734& & &9.5.11.127.1193&9269&3304 \cr
2717&73715117&4.3.13.289.31.61&133&660&2735&74997945&16.7.13.23.31.59&231&172 \cr
 & &32.9.5.7.11.17.19&765&176& & &128.3.49.11.23.43&2107&1472 \cr
\noalign{\hrule}
 & &3.11.19.529.223&49&578& & &81.25.7.11.13.37&361&46 \cr
2718&73965309&4.49.289.223&115&108&2736&74999925&4.9.5.13.361.23&473&112 \cr
 & &32.27.5.7.289.23&1445&1008& & &128.7.11.23.43&989&64 \cr
\noalign{\hrule}
}%
}
$$
\eject
\vglue -23 pt
\noindent\hskip 1 in\hbox to 6.5 in{\ 2737 -- 2772 \hfill\fbd 75040719 -- 78235255\frb}
\vskip -9 pt
$$
\vbox{
\nointerlineskip
\halign{\strut
    \vrule \ \ \hfil \frb #\ 
   &\vrule \hfil \ \ \fbb #\frb\ 
   &\vrule \hfil \ \ \frb #\ \hfil
   &\vrule \hfil \ \ \frb #\ 
   &\vrule \hfil \ \ \frb #\ \ \vrule \hskip 2 pt
   &\vrule \ \ \hfil \frb #\ 
   &\vrule \hfil \ \ \fbb #\frb\ 
   &\vrule \hfil \ \ \frb #\ \hfil
   &\vrule \hfil \ \ \frb #\ 
   &\vrule \hfil \ \ \frb #\ \vrule \cr%
\noalign{\hrule}
 & &3.13.29.43.1543&707&836& & &81.25.7.11.17.29&47&272 \cr
2737&75040719&8.7.11.13.19.29.101&2279&360&2755&76871025&32.9.7.289.47&1625&1336 \cr
 & &128.9.5.11.43.53&583&960& & &512.125.13.167&2171&1280 \cr
\noalign{\hrule}
 & &79.9409.101&715&8694& & &9.5.11.13.17.19.37&31&734 \cr
2738&75074411&4.27.5.7.11.13.23&97&202&2756&76904685&4.11.13.31.367&2033&2400 \cr
 & &16.9.11.97.101&9&88& & &256.3.25.19.107&535&128 \cr
\noalign{\hrule}
 & &3.5.11.13.17.29.71&49&38& & &27.5.169.31.109&11&1406 \cr
2739&75081435&4.5.49.13.17.19.71&801&304&2757&77091885&4.3.11.13.19.37&19&20 \cr
 & &128.9.7.361.89&7581&5696& & &32.5.11.361.37&3971&592 \cr
\noalign{\hrule}
 & &3.5.11.13.17.29.71&337&326& & &27.5.29.53.373&529&902 \cr
2740&75081435&4.5.29.71.163.337&187&1872&2758&77395635&4.5.11.529.29.41&217&102 \cr
 & &128.9.11.13.17.163&163&192& & &16.3.7.17.23.31.41&6601&4216 \cr
\noalign{\hrule}
 & &3.5.17.37.79.101&603&2320& & &3.29.31.59.487&693&206 \cr
2741&75281865&32.27.25.29.67&1309&634&2759&77492901&4.27.7.11.59.103&1363&230 \cr
 & &128.7.11.17.317&2219&704& & &16.5.7.23.29.47&1645&184 \cr
\noalign{\hrule}
 & &27.7.13.23.31.43&319&670& & &3.5.13.17.97.241&363&122 \cr
2742&75329163&4.5.7.11.29.31.67&215&684&2760&77494755&4.9.121.13.17.61&1025&964 \cr
 & &32.9.25.11.19.43&475&176& & &32.25.121.41.241&605&656 \cr
\noalign{\hrule}
 & &5.19.67.71.167&8967&2222& & &9.25.11.17.19.97&1211&1696 \cr
2743&75469805&4.3.49.11.61.101&75&26&2761&77544225&64.5.7.11.53.173&4731&4784 \cr
 & &16.9.25.11.13.61&3965&792& & &2048.3.7.13.19.23.83&13363&13312 \cr
\noalign{\hrule}
 & &3.5.7.13.23.29.83&85&2& & &9.25.49.31.227&893&242 \cr
2744&75567765&4.25.7.13.17.23&333&242&2762&77582925&4.3.5.7.121.19.47&227&172 \cr
 & &16.9.121.17.37&363&5032& & &32.11.43.47.227&473&752 \cr
\noalign{\hrule}
 & &121.79.7921&819&8740& & &3.25.13.29.41.67&77&948 \cr
2745&75716839&8.9.5.7.13.19.23&89&158&2763&77671425&8.9.7.11.29.79&407&146 \cr
 & &32.3.5.7.79.89&105&16& & &32.121.37.73&8833&592 \cr
\noalign{\hrule}
 & &3.11.13.17.101.103&105&116& & &3.5.49.23.43.107&27&188 \cr
2746&75869079&8.9.5.7.29.101.103&815&94&2764&77779905&8.81.7.47.107&1529&2278 \cr
 & &32.25.29.47.163&7661&11600& & &32.11.17.67.139&9313&2992 \cr
\noalign{\hrule}
 & &9.37.41.67.83&143&226& & &27.125.17.23.59&2189&1186 \cr
2747&75924333&4.11.13.37.67.113&505&1974&2765&77857875&4.11.23.199.593&423&170 \cr
 & &16.3.5.7.11.47.101&2585&5656& & &16.9.5.17.47.199&199&376 \cr
\noalign{\hrule}
 & &3.5.169.17.41.43&555&2318& & &3.25.49.17.29.43&143&158 \cr
2748&75976485&4.9.25.19.37.61&143&82&2766&77906325&4.5.7.11.13.17.29.79&1179&164 \cr
 & &16.11.13.19.37.41&407&152& & &32.9.11.13.41.131&5109&7216 \cr
\noalign{\hrule}
 & &625.13.47.199&981&1606& & &9.125.7.11.17.53&451&26 \cr
2749&75993125&4.9.11.47.73.109&125&16&2767&78049125&4.5.7.121.13.41&391&456 \cr
 & &128.3.125.11.73&803&192& & &64.3.17.19.23.41&437&1312 \cr
\noalign{\hrule}
 & &3.25.11.13.41.173&1031&1564& & &27.5.7.11.73.103&8957&1438 \cr
2750&76072425&8.5.11.17.23.1031&117&1148&2768&78160005&4.169.53.719&33&20 \cr
 & &64.9.7.13.17.41&51&224& & &32.3.5.11.13.719&719&208 \cr
\noalign{\hrule}
 & &3.5.7.11.13.37.137&295&706& & &81.25.1331.29&2623&1292 \cr
2751&76111035&4.25.37.59.353&639&286&2769&78162975&8.3.5.17.19.43.61&121&164 \cr
 & &16.9.11.13.59.71&213&472& & &64.121.17.41.61&1037&1312 \cr
\noalign{\hrule}
 & &3.5.29.37.47.101&429&934& & &5.49.11.13.23.97&1387&3618 \cr
2752&76402965&4.9.11.13.37.467&4747&544&2770&78163085&4.27.7.19.67.73&97&104 \cr
 & &256.17.47.101&17&128& & &64.9.13.19.73.97&657&608 \cr
\noalign{\hrule}
 & &9.5.11.37.53.79&2071&2116& & &25.11.13.79.277&179&1206 \cr
2753&76684905&8.11.19.529.37.109&15&422&2771&78231725&4.9.5.11.67.179&117&62 \cr
 & &32.3.5.23.109.211&2507&3376& & &16.81.13.31.67&2511&536 \cr
\noalign{\hrule}
 & &9.7.23.29.31.59&863&1276& & &5.7.19.71.1657&14499&16984 \cr
2754&76856409&8.3.11.841.863&1715&874&2772&78235255&16.81.11.179.193&793&1330 \cr
 & &32.5.343.11.19.23&1045&784& & &64.27.5.7.13.19.61&793&864 \cr
\noalign{\hrule}
}%
}
$$
\eject
\vglue -23 pt
\noindent\hskip 1 in\hbox to 6.5 in{\ 2773 -- 2808 \hfill\fbd 78287391 -- 80544555\frb}
\vskip -9 pt
$$
\vbox{
\nointerlineskip
\halign{\strut
    \vrule \ \ \hfil \frb #\ 
   &\vrule \hfil \ \ \fbb #\frb\ 
   &\vrule \hfil \ \ \frb #\ \hfil
   &\vrule \hfil \ \ \frb #\ 
   &\vrule \hfil \ \ \frb #\ \ \vrule \hskip 2 pt
   &\vrule \ \ \hfil \frb #\ 
   &\vrule \hfil \ \ \fbb #\frb\ 
   &\vrule \hfil \ \ \frb #\ \hfil
   &\vrule \hfil \ \ \frb #\ 
   &\vrule \hfil \ \ \frb #\ \vrule \cr%
\noalign{\hrule}
 & &81.7.169.19.43&407&160& & &3.5.49.29.37.101&323&286 \cr
2773&78287391&64.5.11.13.37.43&243&230&2791&79654155&4.5.7.11.13.17.19.101&783&12802 \cr
 & &256.243.25.23.37&2775&2944& & &16.27.29.37.173&173&72 \cr
\noalign{\hrule}
 & &7.41.47.5807&3047&2760& & &9.11.17.19.47.53&475&1274 \cr
2774&78330623&16.3.5.11.23.47.277&2927&3444&2792&79654707&4.3.25.49.13.361&187&548 \cr
 & &128.9.5.7.41.2927&2927&2880& & &32.5.11.13.17.137&685&208 \cr
\noalign{\hrule}
 & &3.125.11.31.613&1631&2244& & &25.11.23.43.293&2149&1074 \cr
2775&78387375&8.9.7.121.17.233&475&596&2793&79688675&4.3.7.23.179.307&4157&2904 \cr
 & &64.25.19.149.233&4427&4768& & &64.9.121.4157&4157&3168 \cr
\noalign{\hrule}
 & &27.5.13.23.29.67&53&82& & &27.25.23.37.139&121&454 \cr
2776&78429195&4.13.23.41.53.67&1045&174&2794&79845075&4.3.121.139.227&95&322 \cr
 & &16.3.5.11.19.29.41&451&152& & &16.5.7.121.19.23&847&152 \cr
\noalign{\hrule}
 & &3.7.11.13.17.29.53&83&2150& & &9.5.13.311.439&719&836 \cr
2777&78465387&4.25.17.43.83&1827&1742&2795&79869465&8.11.19.439.719&8125&216 \cr
 & &16.9.5.7.13.29.67&67&120& & &128.27.625.13&375&64 \cr
\noalign{\hrule}
 & &3.97.211.1283&1819&2030& & &7.11.43.101.239&2133&2210 \cr
2778&78777483&4.5.7.17.29.97.107&1199&450&2796&79924229&4.27.5.13.17.79.239&43&196 \cr
 & &16.9.125.11.29.109&10875&9592& & &32.3.5.49.13.43.79&1365&1264 \cr
\noalign{\hrule}
 & &5.49.11.23.31.41&53&108& & &9.25.13.23.29.41&31&176 \cr
2779&78782935&8.27.7.31.41.53&293&76&2797&79989975&32.5.11.13.31.41&31&174 \cr
 & &64.3.19.53.293&3021&9376& & &128.3.29.961&961&64 \cr
\noalign{\hrule}
 & &3.5.11.31.73.211&507&296& & &3.49.11.17.41.71&1157&460 \cr
2780&78786345&16.9.5.169.31.37&19&136&2798&80020479&8.5.13.23.71.89&1039&594 \cr
 & &256.13.17.19.37&9139&2176& & &32.27.11.13.1039&1039&1872 \cr
\noalign{\hrule}
 & &27.25.11.13.19.43&527&518& & &9.25.49.169.43&11&158 \cr
2781&78860925&4.3.5.7.13.17.31.37.43&1349&16&2799&80118675&4.3.25.11.43.79&1573&1652 \cr
 & &128.17.19.37.71&1207&2368& & &32.7.1331.13.59&1331&944 \cr
\noalign{\hrule}
 & &27.5.11.79.673&1037&982& & &9.121.29.43.59&1747&790 \cr
2782&78952995&4.9.17.61.79.491&4619&200&2800&80120997&4.3.5.11.79.1747&2177&430 \cr
 & &64.25.17.31.149&2635&4768& & &16.25.7.43.311&311&1400 \cr
\noalign{\hrule}
 & &3.25.7.31.43.113&649&1724& & &3.5.11.19.37.691&2103&1412 \cr
2783&79080225&8.11.31.59.431&3285&1456&2801&80152545&8.9.11.353.701&127&226 \cr
 & &256.9.5.7.13.73&949&384& & &32.113.127.701&14351&11216 \cr
\noalign{\hrule}
 & &3.5.13.433.937&6307&5874& & &3.25.7.11.17.19.43&23&232 \cr
2784&79115595&4.9.5.7.11.17.53.89&4043&38&2802&80208975&16.5.7.23.29.43&221&1026 \cr
 & &16.13.17.19.311&323&2488& & &64.27.13.17.19&13&288 \cr
\noalign{\hrule}
 & &9.11.19.23.31.59&2975&2326& & &5.7.11.13.17.23.41&29&6 \cr
2785&79128027&4.25.7.17.23.1163&787&1950&2803&80235155&4.3.11.13.17.29.41&1725&2422 \cr
 & &16.3.625.13.787&8125&6296& & &16.9.25.7.23.173&173&360 \cr
\noalign{\hrule}
 & &27.25.7.13.1291&131&8906& & &3.5.17.37.67.127&209&172 \cr
2786&79299675&4.61.73.131&65&66&2804&80282415&8.5.11.17.19.43.67&261&74 \cr
 & &16.3.5.11.13.61.73&803&488& & &32.9.19.29.37.43&1247&912 \cr
\noalign{\hrule}
 & &9.13.53.67.191&5497&4626& & &3.11.617.3947&1957&1990 \cr
2787&79354197&4.81.23.239.257&803&1060&2805&80364867&4.5.19.103.199.617&4387&16110 \cr
 & &32.5.11.53.73.239&4015&3824& & &16.9.25.41.107.179&22017&21400 \cr
\noalign{\hrule}
 & &11.17.29.97.151&1235&1332& & &27.7.11.23.1681&247&40 \cr
2788&79430681&8.9.5.11.13.19.29.37&17&302&2806&80380377&16.3.5.11.13.19.41&245&206 \cr
 & &32.3.13.17.37.151&111&208& & &64.25.49.19.103&2575&4256 \cr
\noalign{\hrule}
 & &27.49.19.29.109&59&88& & &5.89.227.797&169&966 \cr
2789&79458057&16.9.11.19.59.109&31&140&2807&80508955&4.3.7.169.23.89&227&396 \cr
 & &128.5.7.11.31.59&1705&3776& & &32.27.11.23.227&621&176 \cr
\noalign{\hrule}
 & &5.11.169.43.199&323&522& & &9.5.7.169.17.89&73&242 \cr
2790&79537315&4.9.11.17.19.29.43&985&434&2808&80544555&4.121.17.73.89&49&138 \cr
 & &16.3.5.7.17.31.197&3349&5208& & &16.3.49.11.23.73&511&2024 \cr
\noalign{\hrule}
}%
}
$$
\eject
\vglue -23 pt
\noindent\hskip 1 in\hbox to 6.5 in{\ 2809 -- 2844 \hfill\fbd 80570175 -- 83360395\frb}
\vskip -9 pt
$$
\vbox{
\nointerlineskip
\halign{\strut
    \vrule \ \ \hfil \frb #\ 
   &\vrule \hfil \ \ \fbb #\frb\ 
   &\vrule \hfil \ \ \frb #\ \hfil
   &\vrule \hfil \ \ \frb #\ 
   &\vrule \hfil \ \ \frb #\ \ \vrule \hskip 2 pt
   &\vrule \ \ \hfil \frb #\ 
   &\vrule \hfil \ \ \fbb #\frb\ 
   &\vrule \hfil \ \ \frb #\ \hfil
   &\vrule \hfil \ \ \frb #\ 
   &\vrule \hfil \ \ \frb #\ \vrule \cr%
\noalign{\hrule}
 & &3.25.7.1849.83&837&1012& & &9.49.31.43.139&551&968 \cr
2809&80570175&8.81.11.23.31.83&301&2210&2827&81711567&16.3.121.19.29.43&7&50 \cr
 & &32.5.7.11.13.17.43&143&272& & &64.25.7.121.29&3509&800 \cr
\noalign{\hrule}
 & &5.19.43.47.421&127&108& & &9.25.7.17.43.71&107&22 \cr
2810&80829895&8.27.43.127.421&4625&836&2828&81744075&4.3.5.7.11.71.107&299&86 \cr
 & &64.3.125.11.19.37&1221&800& & &16.13.23.43.107&1391&184 \cr
\noalign{\hrule}
 & &5.7.17.199.683&639&44& & &9.25.11.19.37.47&169&686 \cr
2811&80870615&8.9.11.71.199&7&206&2829&81776475&4.5.343.169.37&209&246 \cr
 & &32.3.7.11.103&3399&16& & &16.3.49.11.13.19.41&533&392 \cr
\noalign{\hrule}
 & &27.5.7.121.709&899&190& & &9.49.19.43.227&19903&20130 \cr
2812&81070605&4.3.25.7.19.29.31&611&814&2830&81787419&4.27.5.11.13.61.1531&1589&58 \cr
 & &16.11.13.31.37.47&1739&3224& & &16.5.7.11.13.29.227&377&440 \cr
\noalign{\hrule}
 & &27.5.11.31.41.43&1885&1928& & &5.13.17.19.47.83&389&1188 \cr
2813&81159705&16.9.25.11.13.29.241&137&2788&2831&81901495&8.27.5.11.13.389&37&28 \cr
 & &128.17.29.41.137&2329&1856& & &64.3.7.11.37.389&8547&12448 \cr
\noalign{\hrule}
 & &9.25.121.19.157&161&4& & &9.5.49.13.47.61&69&22 \cr
2814&81212175&8.3.5.7.11.19.23&1003&992&2832&82182555&4.27.5.7.11.23.61&247&58 \cr
 & &512.17.23.31.59&12121&15104& & &16.11.13.19.23.29&437&2552 \cr
\noalign{\hrule}
 & &3.7.11.13.17.37.43&125&606& & &5.11.19.31.43.59&83&126 \cr
2815&81222141&4.9.125.7.11.101&1591&884&2833&82186115&4.9.5.7.31.59.83&129&284 \cr
 & &32.5.13.17.37.43&5&16& & &32.27.43.71.83&1917&1328 \cr
\noalign{\hrule}
 & &9.25.7.11.13.361&961&844& & &41.43.149.313&22437&24200 \cr
2816&81306225&8.5.7.11.961.211&1591&114&2834&82221031&16.81.25.121.277&41&14 \cr
 & &32.3.19.31.37.43&1333&592& & &64.3.5.7.11.41.277&3047&3360 \cr
\noalign{\hrule}
 & &9.11.43.97.197&2227&2030& & &9.25.11.29.31.37&923&2518 \cr
2817&81347013&4.5.7.17.29.97.131&1283&366&2835&82325925&4.3.5.13.71.1259&6293&7552 \cr
 & &16.3.5.29.61.1283&8845&10264& & &1024.7.29.31.59&413&512 \cr
\noalign{\hrule}
 & &9.25.19.79.241&131&106& & &25.19.229.757&4329&10054 \cr
2818&81391725&4.3.19.53.131.241&3397&10340&2836&82342675&4.9.11.13.37.457&493&950 \cr
 & &32.5.11.43.47.79&517&688& & &16.3.25.11.17.19.29&493&264 \cr
\noalign{\hrule}
 & &729.121.13.71&875&1798& & &13.163.167.233&20515&18396 \cr
2819&81416907&4.3.125.7.11.29.31&1349&884&2837&82452409&8.9.5.7.11.73.373&5&2 \cr
 & &32.25.13.17.19.71&425&304& & &32.3.25.11.73.373&20075&17904 \cr
\noalign{\hrule}
 & &3.5.11.169.23.127&2261&1626& & &1331.23.37.73&779&900 \cr
2820&81452085&4.9.7.11.17.19.271&21463&20000&2838&82685713&8.9.25.11.19.37.41&73&2182 \cr
 & &256.625.169.127&125&128& & &32.3.5.73.1091&1091&240 \cr
\noalign{\hrule}
 & &3.5.61.269.331&693&962& & &43.79.97.251&5445&5348 \cr
2821&81470685&4.27.7.11.13.37.61&323&1324&2839&82706759&8.9.5.7.121.79.191&1261&76 \cr
 & &32.17.19.37.331&323&592& & &64.3.121.13.19.97&1573&1824 \cr
\noalign{\hrule}
 & &27.5.7.121.23.31&1711&1556& & &81.11.19.59.83&11063&6284 \cr
2822&81527985&8.7.23.29.59.389&18315&21038&2840&82901313&8.13.23.37.1571&545&1026 \cr
 & &32.9.5.11.37.67.157&2479&2512& & &32.27.5.19.23.109&545&368 \cr
\noalign{\hrule}
 & &9.11.13.19.47.71&12203&12250& & &81.25.17.19.127&1247&1166 \cr
2823&81599661&4.125.49.71.12203&111&12314&2841&83067525&4.25.11.17.29.43.53&17591&684 \cr
 & &16.3.5.7.37.47.131&1295&1048& & &32.9.49.19.359&359&784 \cr
\noalign{\hrule}
 & &3.7.31.283.443&33&250& & &3.7.11.521.691&1825&1822 \cr
2824&81615219&4.9.125.11.443&1981&2006&2842&83162541&4.25.11.73.691.911&169&50274 \cr
 & &16.5.7.11.17.59.283&649&680& & &16.27.5.49.169.19&3211&2520 \cr
\noalign{\hrule}
 & &27.25.13.71.131&10823&14098& & &5.13.19.23.29.101&323&990 \cr
2825&81616275&4.7.19.53.79.137&819&682&2843&83198245&4.9.25.11.17.361&593&232 \cr
 & &16.9.49.11.13.31.53&2597&2728& & &64.3.17.29.593&1779&544 \cr
\noalign{\hrule}
 & &3.13.17.59.2087&1545&542& & &5.23.31.67.349&159&190 \cr
2826&81637179&4.9.5.13.103.271&2225&1298&2844&83360395&4.3.25.19.23.53.67&349&924 \cr
 & &16.125.11.59.89&1375&712& & &32.9.7.11.53.349&693&848 \cr
\noalign{\hrule}
}%
}
$$
\eject
\vglue -23 pt
\noindent\hskip 1 in\hbox to 6.5 in{\ 2845 -- 2880 \hfill\fbd 83426205 -- 86185121\frb}
\vskip -9 pt
$$
\vbox{
\nointerlineskip
\halign{\strut
    \vrule \ \ \hfil \frb #\ 
   &\vrule \hfil \ \ \fbb #\frb\ 
   &\vrule \hfil \ \ \frb #\ \hfil
   &\vrule \hfil \ \ \frb #\ 
   &\vrule \hfil \ \ \frb #\ \ \vrule \hskip 2 pt
   &\vrule \ \ \hfil \frb #\ 
   &\vrule \hfil \ \ \fbb #\frb\ 
   &\vrule \hfil \ \ \frb #\ \hfil
   &\vrule \hfil \ \ \frb #\ 
   &\vrule \hfil \ \ \frb #\ \vrule \cr%
\noalign{\hrule}
 & &3.5.83.113.593&255&338& & &3.11.47.227.241&97&130 \cr
2845&83426205&4.9.25.169.17.113&2407&418&2863&84850557&4.5.13.47.97.241&297&908 \cr
 & &16.11.13.19.29.83&551&1144& & &32.27.11.97.227&97&144 \cr
\noalign{\hrule}
 & &3.25.59.109.173&247&272& & &9.13.29.127.197&55&142 \cr
2846&83442225&32.13.17.19.59.109&207&1210&2864&84889467&4.3.5.11.13.71.127&197&184 \cr
 & &128.9.5.121.19.23&2783&3648& & &64.5.11.23.71.197&1633&1760 \cr
\noalign{\hrule}
 & &3.11.17.19.41.191&2063&1490& & &3.49.41.73.193&4031&10058 \cr
2847&83470629&4.5.41.149.2063&2023&4086&2865&84914403&4.29.47.107.139&6755&8118 \cr
 & &16.9.5.7.289.227&1589&2040& & &16.9.5.7.11.41.193&33&40 \cr
\noalign{\hrule}
 & &9.121.13.61.97&413&380& & &5.11.13.113.1051&2699&2556 \cr
2848&83766969&8.3.5.7.11.19.59.97&1343&694&2866&84915545&8.9.71.113.2699&37&8060 \cr
 & &32.5.17.19.79.347&27413&25840& & &64.3.5.13.31.37&1147&96 \cr
\noalign{\hrule}
 & &5.11.19.29.47.59&333&1378& & &9.17.37.43.349&209&838 \cr
2849&84035765&4.9.13.37.47.53&949&790&2867&84954627&4.3.11.19.43.419&8695&9322 \cr
 & &16.3.5.169.73.79&5767&4056& & &16.5.37.47.59.79&2773&3160 \cr
\noalign{\hrule}
 & &3.17.19.29.41.73&6665&6578& & &27.25.17.31.239&3367&4042 \cr
2850&84106293&4.5.11.13.23.31.43.73&551&252&2868&85018275&4.7.13.17.37.43.47&1195&396 \cr
 & &32.9.5.7.19.29.31.43&1505&1488& & &32.9.5.7.11.13.239&143&112 \cr
\noalign{\hrule}
 & &3.5.13.43.79.127&171&44& & &3.7.121.23.31.47&137&850 \cr
2851&84126705&8.27.11.13.19.79&215&136&2869&85151451&4.25.121.17.137&371&234 \cr
 & &128.5.11.17.19.43&323&704& & &16.9.5.7.13.17.53&663&2120 \cr
\noalign{\hrule}
 & &9.5.11.13.23.569&419&166& & &3.5.17.19.79.223&273&50 \cr
2852&84214845&4.83.419.569&75&494&2870&85354365&4.9.125.7.13.79&293&418 \cr
 & &16.3.25.13.19.83&95&664& & &16.7.11.13.19.293&1001&2344 \cr
\noalign{\hrule}
 & &3.49.11.31.1681&1625&1532& & &5.7.11.463.479&39&424 \cr
2853&84263487&8.125.7.13.41.383&2673&8&2871&85384145&16.3.13.53.479&213&266 \cr
 & &128.243.25.11&81&1600& & &64.9.7.13.19.71&1349&3744 \cr
\noalign{\hrule}
 & &9.11.17.361.139&3725&3364& & &729.5.23.1019&1313&2332 \cr
2854&84451257&8.3.25.11.841.149&361&86&2872&85427865&8.11.13.23.53.101&441&142 \cr
 & &32.361.841.43&841&688& & &32.9.49.71.101&3479&1616 \cr
\noalign{\hrule}
 & &25.7.41.79.149&13&162& & &3.25.13.23.37.103&2387&18 \cr
2855&84456925&4.81.13.41.79&1639&1600&2873&85461675&4.27.5.7.11.31&851&634 \cr
 & &512.27.25.11.149&297&256& & &16.23.37.317&317&8 \cr
\noalign{\hrule}
 & &9.5.121.361.43&193&1612& & &9.25.53.71.101&67&572 \cr
2856&84522735&8.3.11.13.31.193&361&218&2874&85514175&8.5.11.13.53.67&303&568 \cr
 & &32.361.31.109&109&496& & &128.3.11.71.101&11&64 \cr
\noalign{\hrule}
 & &5.7.121.13.29.53&2571&2252& & &5.7.11.31.71.101&387&394 \cr
2857&84619535&8.3.5.11.563.857&5239&954&2875&85585885&4.9.5.31.43.101.197&3043&88 \cr
 & &32.27.169.31.53&351&496& & &64.3.11.17.43.179&3043&4128 \cr
\noalign{\hrule}
 & &3.11.169.43.353&5803&5846& & &9.7.37.83.443&125&134 \cr
2858&84653283&4.7.169.37.79.829&45&124&2876&85708539&4.125.67.83.443&429&14 \cr
 & &32.9.5.7.31.37.829&40145&39792& & &16.3.25.7.11.13.67&871&2200 \cr
\noalign{\hrule}
 & &3.25.7.121.31.43&153&148& & &9.5.121.19.829&203&82 \cr
2859&84678825&8.27.5.121.17.31.37&2951&316&2877&85764195&4.3.7.29.41.829&1387&1100 \cr
 & &64.13.37.79.227&17933&15392& & &32.25.11.19.29.73&365&464 \cr
\noalign{\hrule}
 & &11.53.337.431&123&460& & &11.13.17.23.29.53&171&148 \cr
2860&84679001&8.3.5.23.41.431&195&236&2878&85938281&8.9.13.17.19.37.53&829&140 \cr
 & &64.9.25.13.23.59&13275&9568& & &64.3.5.7.37.829&12435&8288 \cr
\noalign{\hrule}
 & &27.5.17.19.29.67&77&94& & &5.11.43.73.499&651&152 \cr
2861&84724515&4.3.5.7.11.29.47.67&3097&52&2879&86149855&16.3.5.7.19.31.43&331&486 \cr
 & &32.11.13.19.163&2119&176& & &64.729.7.331&5103&10592 \cr
\noalign{\hrule}
 & &3.5.7.11.17.29.149&85&118& & &11.17.19.127.191&13689&10060 \cr
2862&84842835&4.25.289.59.149&783&8008&2880&86185121&8.81.5.169.503&1519&1016 \cr
 & &64.27.7.11.13.29&117&32& & &128.27.49.31.127&1519&1728 \cr
\noalign{\hrule}
}%
}
$$
\eject
\vglue -23 pt
\noindent\hskip 1 in\hbox to 6.5 in{\ 2881 -- 2916 \hfill\fbd 86288709 -- 89501335\frb}
\vskip -9 pt
$$
\vbox{
\nointerlineskip
\halign{\strut
    \vrule \ \ \hfil \frb #\ 
   &\vrule \hfil \ \ \fbb #\frb\ 
   &\vrule \hfil \ \ \frb #\ \hfil
   &\vrule \hfil \ \ \frb #\ 
   &\vrule \hfil \ \ \frb #\ \ \vrule \hskip 2 pt
   &\vrule \ \ \hfil \frb #\ 
   &\vrule \hfil \ \ \fbb #\frb\ 
   &\vrule \hfil \ \ \frb #\ \hfil
   &\vrule \hfil \ \ \frb #\ 
   &\vrule \hfil \ \ \frb #\ \vrule \cr%
\noalign{\hrule}
 & &3.13.19.23.61.83&93&154& & &9.19.31.61.271&2849&2300 \cr
2881&86288709&4.9.7.11.23.31.83&2489&580&2899&87630831&8.25.7.11.23.31.37&813&38 \cr
 & &32.5.7.19.29.131&917&2320& & &32.3.7.11.19.271&77&16 \cr
\noalign{\hrule}
 & &9.5.11.37.53.89&389&56& & &9.5.11.13.31.443&7751&12184 \cr
2882&86391855&16.7.11.53.389&97&486&2900&88371855&16.23.337.1523&593&930 \cr
 & &64.243.7.97&189&3104& & &64.3.5.23.31.593&593&736 \cr
\noalign{\hrule}
 & &7.11.53.59.359&15385&12258& & &7.11.19.193.313&635&828 \cr
2883&86439661&4.27.5.17.181.227&99&82&2901&88378367&8.9.5.23.127.313&991&1930 \cr
 & &16.243.5.11.41.227&9307&9720& & &32.3.25.193.991&991&1200 \cr
\noalign{\hrule}
 & &3.7.11.37.67.151&779&628& & &5.13.41.79.421&539&1566 \cr
2884&86469999&8.11.19.37.41.157&3015&3422&2902&88635235&4.27.49.11.29.41&79&530 \cr
 & &32.9.5.19.29.59.67&2755&2832& & &16.9.5.7.53.79&477&56 \cr
\noalign{\hrule}
 & &9.5.11.179.977&571&406& & &5.7.11.19.61.199&711&284 \cr
2885&86567085&4.3.7.29.179.571&167&370&2903&88796785&8.9.11.19.71.79&223&14 \cr
 & &16.5.37.167.571&6179&4568& & &32.3.7.71.223&213&3568 \cr
\noalign{\hrule}
 & &11.29.269.1009&1649&9450& & &9.5.7.13.529.41&121&408 \cr
2886&86583299&4.27.25.7.17.97&377&496&2904&88816455&16.27.5.121.13.17&19&46 \cr
 & &128.3.25.13.29.31&975&1984& & &64.121.17.19.23&2057&608 \cr
\noalign{\hrule}
 & &3.5.7.67.97.127&411&478& & &5.49.169.19.113&183&748 \cr
2887&86664165&4.9.5.97.137.239&8509&14674&2905&88896535&8.3.11.169.17.61&5061&5248 \cr
 & &16.11.23.29.67.127&253&232& & &2048.9.7.41.241&9881&9216 \cr
\noalign{\hrule}
 & &41.839.2521&15939&18460& & &7.11.23.61.823&31655&31716 \cr
2888&86719879&8.9.5.7.11.13.23.71&841&82&2906&88909513&8.9.5.13.23.487.881&889&8 \cr
 & &32.3.5.7.841.41&841&1680& & &128.3.5.7.127.487&7305&8128 \cr
\noalign{\hrule}
 & &3.11.13.529.383&193&190& & &9.5.7.11.67.383&943&206 \cr
2889&86918403&4.5.11.13.19.529.193&2681&36&2907&88915365&4.3.5.7.23.41.103&3439&1276 \cr
 & &32.9.7.193.383&193&336& & &32.11.19.29.181&551&2896 \cr
\noalign{\hrule}
 & &25.19.397.461&583&9342& & &13.17.19.127.167&12447&15620 \cr
2890&86933075&4.27.11.53.173&325&152&2908&89056591&8.27.5.11.71.461&1543&762 \cr
 & &64.3.25.11.13.19&429&32& & &32.81.127.1543&1543&1296 \cr
\noalign{\hrule}
 & &9.169.19.23.131&319&188& & &9.11.19.23.29.71&1391&610 \cr
2891&87072687&8.3.11.19.23.29.47&655&104&2909&89078517&4.3.5.13.19.61.107&667&724 \cr
 & &128.5.13.47.131&235&64& & &32.5.23.29.61.181&905&976 \cr
\noalign{\hrule}
 & &7.29.347.1237&9361&702& & &9.5.13.19.71.113&731&618 \cr
2892&87135517&4.27.11.13.23.37&445&406&2910&89175645&4.27.5.13.17.43.103&2201&11704 \cr
 & &16.9.5.7.11.29.89&445&792& & &64.7.11.19.31.71&341&224 \cr
\noalign{\hrule}
 & &27.5.23.61.461&329&1976& & &3.11.23.41.47.61&2363&1846 \cr
2893&87315705&16.7.13.19.23.47&153&176&2911&89218173&4.13.17.41.71.139&1525&1386 \cr
 & &512.9.11.13.17.19&4199&2816& & &16.9.25.7.11.13.17.61&975&952 \cr
\noalign{\hrule}
 & &7.13.151.6361&10051&3690& & &9.25.11.13.47.59&49&226 \cr
2894&87406501&4.9.5.19.529.41&481&462&2912&89221275&4.3.49.13.47.113&1969&2308 \cr
 & &16.27.5.7.11.13.23.37&3105&3256& & &32.7.11.179.577&4039&2864 \cr
\noalign{\hrule}
 & &9.5.11.13.289.47&25&586& & &5.7.23.29.43.89&13&1002 \cr
2895&87406605&4.3.125.17.293&2303&2678&2913&89341315&4.3.13.89.167&45&44 \cr
 & &16.49.13.47.103&103&392& & &32.27.5.11.13.167&2171&4752 \cr
\noalign{\hrule}
 & &5.49.11.17.23.83&81&334& & &9.125.7.11.1033&1987&1112 \cr
2896&87460835&4.81.49.17.167&667&166&2914&89483625&16.3.11.139.1987&3287&1300 \cr
 & &16.27.23.29.83&29&216& & &128.25.13.19.173&2249&1216 \cr
\noalign{\hrule}
 & &3.5.43.137.991&14729&14726& & &3.13.19.23.59.89&341&400 \cr
2897&87569715&4.11.13.37.103.199.991&3177&98896&2915&89492793&32.25.11.23.31.89&171&1876 \cr
 & &128.9.7.11.353.883&67991&67776& & &256.9.5.7.19.67&1407&640 \cr
\noalign{\hrule}
 & &5.17.19.193.281&5533&10872& & &5.7.11.37.61.103&1577&1272 \cr
2898&87586295&16.9.11.151.503&579&1082&2916&89501335&16.3.19.53.83.103&1221&3178 \cr
 & &64.27.193.541&541&864& & &64.9.7.11.37.227&227&288 \cr
\noalign{\hrule}
}%
}
$$
\eject
\vglue -23 pt
\noindent\hskip 1 in\hbox to 6.5 in{\ 2917 -- 2952 \hfill\fbd 89516735 -- 92581775\frb}
\vskip -9 pt
$$
\vbox{
\nointerlineskip
\halign{\strut
    \vrule \ \ \hfil \frb #\ 
   &\vrule \hfil \ \ \fbb #\frb\ 
   &\vrule \hfil \ \ \frb #\ \hfil
   &\vrule \hfil \ \ \frb #\ 
   &\vrule \hfil \ \ \frb #\ \ \vrule \hskip 2 pt
   &\vrule \ \ \hfil \frb #\ 
   &\vrule \hfil \ \ \fbb #\frb\ 
   &\vrule \hfil \ \ \frb #\ \hfil
   &\vrule \hfil \ \ \frb #\ 
   &\vrule \hfil \ \ \frb #\ \vrule \cr%
\noalign{\hrule}
 & &5.7.11.41.53.107&13897&16812& & &9.5.7.11.41.643&629&14 \cr
2917&89516735&8.9.13.467.1069&253&214&2935&91347795&4.3.49.11.17.37&325&214 \cr
 & &32.3.11.23.107.1069&1069&1104& & &16.25.13.17.107&221&4280 \cr
\noalign{\hrule}
 & &27.5.7.13.37.197&493&506& & &3.5.31.47.53.79&157&78 \cr
2918&89545365&4.5.7.11.17.23.29.197&1157&222&2936&91506885&4.9.13.31.53.157&205&484 \cr
 & &16.3.13.23.29.37.89&667&712& & &32.5.121.41.157&4961&2512 \cr
\noalign{\hrule}
 & &3.5.11.13.83.503&683&826& & &9.11.29.71.449&95&544 \cr
2919&89551605&4.5.7.59.83.683&1&414&2937&91524609&64.5.11.17.19.29&351&142 \cr
 & &16.9.23.683&69&5464& & &256.27.5.13.71&65&384 \cr
\noalign{\hrule}
 & &3.25.7.11.13.1193&5959&4766& & &11.19.31.71.199&4335&10504 \cr
2920&89564475&4.7.59.101.2383&19665&22048&2938&91541791&16.3.5.13.289.101&161&60 \cr
 & &256.9.5.13.19.23.53&3021&2944& & &128.9.25.7.17.23&3825&10304 \cr
\noalign{\hrule}
 & &27.7.11.79.547&1577&1030& & &11.29.239.1201&10071&3140 \cr
2921&89839827&4.9.5.7.19.83.103&547&650&2939&91565441&8.27.5.157.373&481&638 \cr
 & &16.125.13.83.547&1079&1000& & &32.9.5.11.13.29.37&481&720 \cr
\noalign{\hrule}
 & &5.11.23.29.31.79&4887&5002& & &9.5.7.11.13.19.107&47&86 \cr
2922&89841565&4.27.41.61.79.181&4105&18404&2940&91576485&4.3.5.11.43.47.107&911&266 \cr
 & &32.3.5.43.107.821&13803&13136& & &16.7.19.47.911&911&376 \cr
\noalign{\hrule}
 & &9.7.17.23.41.89&429&268& & &9.25.11.101.367&1187&1288 \cr
2923&89885817&8.27.11.13.67.89&415&1394&2941&91740825&16.7.23.367.1187&8375&66 \cr
 & &32.5.13.17.41.83&415&208& & &64.3.125.11.67&335&32 \cr
\noalign{\hrule}
 & &81.5.13.17.19.53&185&1192& & &27.5.7.13.31.241&79&110 \cr
2924&90131535&16.25.13.37.149&1431&506&2942&91781235&4.25.11.13.79.241&283&42 \cr
 & &64.27.11.23.53&253&32& & &16.3.7.11.79.283&869&2264 \cr
\noalign{\hrule}
 & &9.25.49.13.17.37&209&124& & &9.5.11.13.17.839&10349&3914 \cr
2925&90151425&8.5.49.11.13.19.31&1377&142&2943&91782405&4.19.79.103.131&61&42 \cr
 & &32.81.11.17.71&781&144& & &16.3.7.61.79.131&4819&7336 \cr
\noalign{\hrule}
 & &3.29.41.131.193&10855&5258& & &9.25.43.53.179&91&134 \cr
2926&90184461&4.5.11.13.167.239&1683&488&2944&91786725&4.7.13.53.67.179&275&96 \cr
 & &64.9.121.17.61&2057&5856& & &256.3.25.11.13.67&871&1408 \cr
\noalign{\hrule}
 & &3.11.89.97.317&203&114& & &3.7.13.29.41.283&1585&396 \cr
2927&90309813&4.9.7.11.19.29.97&575&1268&2945&91860951&8.27.5.11.13.317&2089&1772 \cr
 & &32.25.23.29.317&667&400& & &64.5.443.2089&10445&14176 \cr
\noalign{\hrule}
 & &9.25.11.13.29.97&421&646& & &27.5.49.13.1069&3401&2332 \cr
2928&90508275&4.13.17.19.29.421&1947&7420&2946&91928655&8.3.5.11.19.53.179&149&434 \cr
 & &32.3.5.7.11.53.59&413&848& & &32.7.31.149.179&4619&2864 \cr
\noalign{\hrule}
 & &3.7.11.13.97.311&1303&874& & &5.11.19.37.2383&1219&1164 \cr
2929&90591501&4.19.23.97.1303&1573&270&2947&92138695&8.3.19.23.37.53.97&929&78 \cr
 & &16.27.5.121.13.23&495&184& & &32.9.13.97.929&11349&14864 \cr
\noalign{\hrule}
 & &3.13.17.281.487&3999&4280& & &5.7.13.31.47.139&69&22 \cr
2930&90729561&16.9.5.13.31.43.107&487&4114&2948&92147965&4.3.5.11.23.31.139&1521&1676 \cr
 & &64.5.121.17.487&121&160& & &32.27.11.169.419&5447&4752 \cr
\noalign{\hrule}
 & &9.5.11.31.61.97&323&226& & &81.43.103.257&1793&1690 \cr
2931&90796365&4.5.11.17.19.31.113&213&128&2949&92198493&4.5.11.169.163.257&2567&774 \cr
 & &1024.3.19.71.113&8023&9728& & &16.9.5.13.17.43.151&1105&1208 \cr
\noalign{\hrule}
 & &25.11.169.37.53&207&482& & &9.19.29.43.433&539&278 \cr
2932&91137475&4.9.13.23.37.241&467&430&2950&92331621&4.49.11.139.433&703&270 \cr
 & &16.3.5.43.241.467&10363&11208& & &16.27.5.7.11.19.37&407&840 \cr
\noalign{\hrule}
 & &9.5.13.19.29.283&17&22& & &5.13.19.37.43.47&413&198 \cr
2933&91220805&4.3.11.17.19.29.283&149&700&2951&92349595&4.9.7.11.19.37.59&47&10 \cr
 & &32.25.7.11.17.149&5215&2992& & &16.3.5.7.11.47.59&649&168 \cr
\noalign{\hrule}
 & &9.5.13.43.3631&3213&418& & &25.11.13.19.29.47&119&24 \cr
2934&91337805&4.243.7.11.17.19&689&1012&2952&92581775&16.3.5.7.17.29.47&237&92 \cr
 & &32.121.13.23.53&2783&848& & &128.9.17.23.79&1817&9792 \cr
\noalign{\hrule}
}%
}
$$
\eject
\vglue -23 pt
\noindent\hskip 1 in\hbox to 6.5 in{\ 2953 -- 2988 \hfill\fbd 92721519 -- 95267469\frb}
\vskip -9 pt
$$
\vbox{
\nointerlineskip
\halign{\strut
    \vrule \ \ \hfil \frb #\ 
   &\vrule \hfil \ \ \fbb #\frb\ 
   &\vrule \hfil \ \ \frb #\ \hfil
   &\vrule \hfil \ \ \frb #\ 
   &\vrule \hfil \ \ \frb #\ \ \vrule \hskip 2 pt
   &\vrule \ \ \hfil \frb #\ 
   &\vrule \hfil \ \ \fbb #\frb\ 
   &\vrule \hfil \ \ \frb #\ \hfil
   &\vrule \hfil \ \ \frb #\ 
   &\vrule \hfil \ \ \frb #\ \vrule \cr%
\noalign{\hrule}
 & &9.11.17.37.1489&10825&14488& & &25.121.173.179&213&392 \cr
2953&92721519&16.25.433.1811&689&1122&2971&93675175&16.3.5.49.71.173&783&428 \cr
 & &64.3.25.11.13.17.53&689&800& & &128.81.7.29.107&8667&12992 \cr
\noalign{\hrule}
 & &3.7.169.17.29.53&625&506& & &81.13.79.1129&1091&38 \cr
2954&92731821&4.625.11.13.23.53&967&252&2972&93918123&4.19.79.1091&585&506 \cr
 & &32.9.125.7.967&2901&2000& & &16.9.5.11.13.19.23&209&920 \cr
\noalign{\hrule}
 & &7.11.59.137.149&5&6& & &9.11.13.19.23.167&901&602 \cr
2955&92736259&4.3.5.7.59.137.149&4875&3916&2973&93923973&4.7.11.17.19.43.53&267&740 \cr
 & &32.9.625.11.13.89&10413&10000& & &32.3.5.7.17.37.89&4403&7120 \cr
\noalign{\hrule}
 & &9.5.13.19.61.137&221&84& & &27.25.37.53.71&253&1178 \cr
2956&92888055&8.27.7.169.17.19&341&172&2974&93980925&4.11.19.23.31.71&91&162 \cr
 & &64.7.11.17.31.43&5117&10912& & &16.81.7.13.19.31&651&1976 \cr
\noalign{\hrule}
 & &5.11.29.101.577&129&448& & &7.23.641.911&135&776 \cr
2957&92951815&128.3.5.7.43.101&461&246&2975&94016111&16.27.5.7.23.97&517&356 \cr
 & &512.9.41.461&4149&10496& & &128.3.5.11.47.89&4183&10560 \cr
\noalign{\hrule}
 & &125.49.67.227&627&962& & &27.19.31.61.97&529&1118 \cr
2958&93155125&4.3.25.7.11.13.19.37&141&134&2976&94098051&4.13.529.43.97&1395&836 \cr
 & &16.9.13.19.37.47.67&5499&5624& & &32.9.5.11.19.23.31&115&176 \cr
\noalign{\hrule}
 & &49.71.73.367&65&432& & &5.121.29.41.131&333&322 \cr
2959&93205889&32.27.5.7.13.73&33&40&2977&94234195&4.9.7.11.23.29.37.41&7&950 \cr
 & &512.81.25.11.13&11583&6400& & &16.3.25.49.19.37&931&4440 \cr
\noalign{\hrule}
 & &49.71.73.367&33&40& & &25.11.13.23.31.37&413&438 \cr
2960&93205889&16.3.5.7.11.71.367&65&432&2978&94312075&4.3.7.11.13.31.59.73&63&4370 \cr
 & &512.81.25.11.13&11583&6400& & &16.27.5.49.19.23&931&216 \cr
\noalign{\hrule}
 & &7.17.529.1481&545&936& & &3.5.121.61.853&1127&1432 \cr
2961&93230431&16.9.5.7.13.23.109&1111&306&2979&94439895&16.49.121.23.179&85&36 \cr
 & &64.81.11.17.101&1111&2592& & &128.9.5.17.23.179&4117&3264 \cr
\noalign{\hrule}
 & &27.5.121.29.197&5311&11024& & &3.5.13.29.73.229&147&82 \cr
2962&93321855&32.13.47.53.113&1189&1302&2980&94534635&4.9.49.29.41.73&583&2410 \cr
 & &128.3.7.13.29.31.41&2821&2624& & &16.5.7.11.53.241&1687&4664 \cr
\noalign{\hrule}
 & &3.11.17.43.53.73&1213&1196& & &9.13.41.109.181&4419&3002 \cr
2963&93331887&8.13.23.43.53.1213&1289&990&2981&94640013&4.81.19.79.491&205&286 \cr
 & &32.9.5.11.1213.1289&19335&19408& & &16.5.11.13.19.41.79&869&760 \cr
\noalign{\hrule}
 & &3.5.7.13.19.59.61&349&444& & &5.17.29.107.359&1089&730 \cr
2964&93340065&8.9.7.37.59.349&8315&4598&2982&94688045&4.9.25.121.29.73&1391&434 \cr
 & &32.5.121.19.1663&1663&1936& & &16.3.7.11.13.31.107&651&1144 \cr
\noalign{\hrule}
 & &27.49.13.61.89&323&314& & &3.5.7.73.83.149&29&44 \cr
2965&93373371&4.3.17.19.61.89.157&2849&1690&2983&94793055&8.7.11.29.83.149&5067&7300 \cr
 & &16.5.7.11.169.37.157&5809&5720& & &64.9.25.73.563&563&480 \cr
\noalign{\hrule}
 & &3.25.7.11.19.23.37&3111&3364& & &81.5.7.13.31.83&457&3278 \cr
2966&93375975&8.9.17.19.841.61&335&506&2984&94827915&4.9.11.149.457&899&442 \cr
 & &32.5.11.17.23.61.67&1037&1072& & &16.11.13.17.29.31&493&88 \cr
\noalign{\hrule}
 & &5.7.11.19.53.241&17&36& & &9.11.19.29.37.47&325&748 \cr
2967&93434495&8.9.5.7.11.17.241&1273&1378&2985&94860711&8.25.121.13.17.19&597&1702 \cr
 & &32.3.13.17.19.53.67&871&816& & &32.3.5.23.37.199&995&368 \cr
\noalign{\hrule}
 & &5.83.113.1993&22451&24444& & &9.25.7.29.31.67&4301&2224 \cr
2968&93461735&8.9.7.11.13.97.157&549&452&2986&94866975&32.7.11.17.23.139&67&186 \cr
 & &64.81.61.113.157&4941&5024& & &128.3.31.67.139&139&64 \cr
\noalign{\hrule}
 & &25.11.17.37.541&287&342& & &3.49.19.59.577&495&82 \cr
2969&93579475&4.9.5.7.19.41.541&37&578&2987&95082099&4.27.5.7.11.19.41&295&484 \cr
 & &16.3.7.289.19.37&133&408& & &32.25.1331.59&1331&400 \cr
\noalign{\hrule}
 & &25.11.13.17.23.67&1413&542& & &3.11.89.163.199&1139&1050 \cr
2970&93654275&4.9.5.11.157.271&1541&186&2988&95267469&4.9.25.7.17.67.163&319&7654 \cr
 & &16.27.23.31.67&27&248& & &16.5.11.29.43.89&145&344 \cr
\noalign{\hrule}
}%
}
$$
\eject
\vglue -23 pt
\noindent\hskip 1 in\hbox to 6.5 in{\ 2989 -- 3024 \hfill\fbd 95285619 -- 97978725\frb}
\vskip -9 pt
$$
\vbox{
\nointerlineskip
\halign{\strut
    \vrule \ \ \hfil \frb #\ 
   &\vrule \hfil \ \ \fbb #\frb\ 
   &\vrule \hfil \ \ \frb #\ \hfil
   &\vrule \hfil \ \ \frb #\ 
   &\vrule \hfil \ \ \frb #\ \ \vrule \hskip 2 pt
   &\vrule \ \ \hfil \frb #\ 
   &\vrule \hfil \ \ \fbb #\frb\ 
   &\vrule \hfil \ \ \frb #\ \hfil
   &\vrule \hfil \ \ \frb #\ 
   &\vrule \hfil \ \ \frb #\ \vrule \cr%
\noalign{\hrule}
 & &27.11.13.23.29.37&1159&9830& & &27.7.11.19.31.79&15265&11674 \cr
2989&95285619&4.5.19.61.983&461&522&3007&96737949&4.5.13.43.71.449&237&686 \cr
 & &16.9.5.19.29.461&461&760& & &16.3.5.343.43.79&245&344 \cr
\noalign{\hrule}
 & &9.5.11.289.23.29&589&2590& & &27.5.53.83.163&941&526 \cr
2990&95417685&4.3.25.7.19.31.37&667&1442&3008&96799995&4.3.53.263.941&391&550 \cr
 & &16.49.23.29.103&103&392& & &16.25.11.17.23.263&6049&7480 \cr
\noalign{\hrule}
 & &9.17.419.1489&133&286& & &9.25.13.157.211&1617&3658 \cr
2991&95455323&4.7.11.13.19.1489&1245&244&3009&96896475&4.27.49.11.31.59&157&184 \cr
 & &32.3.5.19.61.83&5795&1328& & &64.49.23.59.157&1357&1568 \cr
\noalign{\hrule}
 & &27.25.11.13.23.43&71&136& & &9.121.19.43.109&25&146 \cr
2992&95463225&16.3.5.11.17.43.71&37&92&3010&96978717&4.25.43.73.109&33&76 \cr
 & &128.17.23.37.71&1207&2368& & &32.3.25.11.19.73&73&400 \cr
\noalign{\hrule}
 & &3.13.79.89.349&55&1102& & &49.11.31.37.157&27675&25862 \cr
2993&95698941&4.5.11.19.29.79&387&482&3011&97062581&4.27.25.41.67.193&2821&74 \cr
 & &16.9.29.43.241&1247&5784& & &16.9.5.7.13.31.37&117&40 \cr
\noalign{\hrule}
 & &11.31.131.2143&2025&2036& & &23.29.41.53.67&605&2142 \cr
2994&95729953&8.81.25.509.2143&24497&29078&3012&97109197&4.9.5.7.121.17.23&67&186 \cr
 & &32.9.7.11.17.31.67.131&1071&1072& & &16.27.5.11.31.67&297&1240 \cr
\noalign{\hrule}
 & &9.11.23.137.307&2135&628& & &9.25.7.17.19.191&12949&16274 \cr
2995&95768343&8.5.7.23.61.157&561&538&3013&97166475&4.23.79.103.563&627&1190 \cr
 & &32.3.5.11.17.61.269&4573&4880& & &16.3.5.7.11.17.19.103&103&88 \cr
\noalign{\hrule}
 & &11.361.59.409&309&100& & &3.11.23.181.709&481&228 \cr
2996&95824201&8.3.25.19.59.103&303&818&3014&97401711&8.9.13.19.37.181&95&86 \cr
 & &32.9.5.101.409&505&144& & &32.5.13.361.37.43&23465&25456 \cr
\noalign{\hrule}
 & &3.7.13.17.19.1087&341&322& & &5.13.19.23.47.73&81&154 \cr
2997&95850573&4.49.11.23.31.1087&1107&20&3015&97457555&4.81.7.11.13.19.23&235&64 \cr
 & &32.27.5.11.31.41&1705&5904& & &512.9.5.7.11.47&693&256 \cr
\noalign{\hrule}
 & &9.5.47.61.743&319&2548& & &25.7.47.71.167&911&864 \cr
2998&95858145&8.3.5.49.11.13.29&109&94&3016&97523825&64.27.7.167.911&2209&1298 \cr
 & &32.7.11.13.47.109&1001&1744& & &256.9.11.2209.59&5841&6016 \cr
\noalign{\hrule}
 & &3.17.19.167.593&1123&1716& & &9.121.157.571&9497&9500 \cr
2999&95961039&8.9.11.13.19.1123&5177&4930&3017&97625583&8.3.125.19.571.9497&49933&21442 \cr
 & &32.5.11.17.29.31.167&899&880& & &32.13.19.23.71.151.167&403351&403472 \cr
\noalign{\hrule}
 & &7.11.19.29.31.73&11349&11938& & &81.13.31.41.73&9917&10450 \cr
3000&96012301&4.9.7.13.47.97.127&12617&1060&3018&97700499&4.9.25.11.19.47.211&449&26 \cr
 & &32.3.5.11.31.37.53&795&592& & &16.11.13.211.449&2321&3592 \cr
\noalign{\hrule}
 & &5.13.17.19.23.199&957&658& & &49.13.131.1171&3247&4950 \cr
3001&96094115&4.3.7.11.29.47.199&4695&1076&3019&97716437&4.9.25.7.11.17.191&299&656 \cr
 & &32.9.5.269.313&2421&5008& & &128.3.5.11.13.23.41&3795&2624 \cr
\noalign{\hrule}
 & &5.11.31.157.359&1247&9882& & &3.7.11.59.71.101&597&184 \cr
3002&96098915&4.81.29.43.61&95&34&3020&97733559&16.9.23.101.199&355&554 \cr
 & &16.27.5.17.19.29&459&4408& & &64.5.23.71.277&1385&736 \cr
\noalign{\hrule}
 & &27.7.47.83.131&4199&1958& & &3.5.49.11.107.113&221&956 \cr
3003&96584859&4.7.11.13.17.19.89&153&1310&3021&97755735&8.13.17.113.239&63&176 \cr
 & &16.9.5.289.131&289&40& & &256.9.7.11.13.17&663&128 \cr
\noalign{\hrule}
 & &3.23.521.2687&1309&1378& & &27.125.7.41.101&1637&1738 \cr
3004&96594963&4.7.11.13.17.53.521&2115&1532&3022&97831125&4.7.11.41.79.1637&2333&20340 \cr
 & &32.9.5.13.17.47.383&31161&30640& & &32.9.5.113.2333&2333&1808 \cr
\noalign{\hrule}
 & &5.13.23.37.1747&1023&724& & &9.25.11.13.17.179&3031&2494 \cr
3005&96635305&8.3.5.11.31.37.181&927&1108&3023&97908525&4.3.7.11.29.43.433&2225&806 \cr
 & &64.27.31.103.277&28531&26784& & &16.25.13.29.31.89&899&712 \cr
\noalign{\hrule}
 & &5.13.19.29.37.73&539&534& & &9.25.13.19.41.43&107&98 \cr
3006&96736315&4.3.49.11.13.19.73.89&115&1272&3024&97978725&4.5.49.13.19.43.107&231&16 \cr
 & &64.9.5.49.11.23.53&13409&14112& & &128.3.343.11.107&3773&6848 \cr
\noalign{\hrule}
}%
}
$$
\eject
\vglue -23 pt
\noindent\hskip 1 in\hbox to 6.5 in{\ 3025 -- 3060 \hfill\fbd 98080125 -- 100370925\frb}
\vskip -9 pt
$$
\vbox{
\nointerlineskip
\halign{\strut
    \vrule \ \ \hfil \frb #\ 
   &\vrule \hfil \ \ \fbb #\frb\ 
   &\vrule \hfil \ \ \frb #\ \hfil
   &\vrule \hfil \ \ \frb #\ 
   &\vrule \hfil \ \ \frb #\ \ \vrule \hskip 2 pt
   &\vrule \ \ \hfil \frb #\ 
   &\vrule \hfil \ \ \fbb #\frb\ 
   &\vrule \hfil \ \ \frb #\ \hfil
   &\vrule \hfil \ \ \frb #\ 
   &\vrule \hfil \ \ \frb #\ \vrule \cr%
\noalign{\hrule}
 & &3.125.11.13.31.59&27727&25898& & &9.7.11.13.43.257&329&230 \cr
3025&98080125&4.7.17.23.233.563&1485&1252&3043&99558459&4.5.49.23.47.257&251&6 \cr
 & &32.27.5.11.313.563&5067&5008& & &16.3.23.47.251&5773&376 \cr
\noalign{\hrule}
 & &5.31.443.1429&10439&3294& & &27.5.7.121.13.67&103&40 \cr
3026&98122285&4.27.11.13.61.73&493&310&3044&99594495&16.3.25.11.67.103&37&238 \cr
 & &16.9.5.13.17.29.31&493&936& & &64.7.17.37.103&3811&544 \cr
\noalign{\hrule}
 & &3.47.137.5081&7645&7598& & &9.49.17.97.137&8103&8200 \cr
3027&98149677&4.5.11.29.131.137.139&823&684&3045&99627633&16.27.25.7.37.41.73&187&2 \cr
 & &32.9.5.19.29.131.823&119335&119472& & &64.5.11.17.41.73&2993&1760 \cr
\noalign{\hrule}
 & &3.11.19.233.673&13&686& & &27.25.11.89.151&511&244 \cr
3028&98319243&4.343.11.13.19&279&260&3046&99784575&8.9.5.7.11.61.73&493&178 \cr
 & &32.9.5.7.169.31&5915&1488& & &32.17.29.73.89&493&1168 \cr
\noalign{\hrule}
 & &125.7.13.17.509&1067&558& & &5.11.23.157.503&711&554 \cr
3029&98427875&4.9.7.11.17.31.97&509&800&3047&99898315&4.9.79.277.503&133&370 \cr
 & &256.3.25.31.509&93&128& & &16.3.5.7.19.37.277&4921&6648 \cr
\noalign{\hrule}
 & &7.11.19.31.41.53&3213&3266& & &9.5.7.19.59.283&1327&88 \cr
3030&98552069&4.27.49.17.23.41.71&325&3154&3048&99931545&16.3.11.19.1327&1981&2000 \cr
 & &16.9.25.13.17.19.83&5525&5976& & &512.125.7.11.283&275&256 \cr
\noalign{\hrule}
 & &5.7.11.13.17.19.61&1021&594& & &23.37.41.47.61&26919&24992 \cr
3031&98613515&4.27.121.13.1021&745&2318&3049&100032497&64.27.11.71.997&1105&1886 \cr
 & &16.9.5.19.61.149&149&72& & &256.9.5.13.17.23.41&1105&1152 \cr
\noalign{\hrule}
 & &3.5.121.17.23.139&2627&570& & &5.7.13.31.41.173&407&1842 \cr
3032&98643435&4.9.25.19.37.71&143&782&3050&100046765&4.3.11.31.37.307&169&138 \cr
 & &16.11.13.17.19.23&13&152& & &16.9.11.169.23.37&3663&2392 \cr
\noalign{\hrule}
 & &3.11.17.53.3319&1633&1686& & &9.29.31.89.139&10231&2140 \cr
3033&98683827&4.9.11.17.23.71.281&3319&3710&3051&100093761&8.5.13.107.787&2021&1914 \cr
 & &16.5.7.53.281.3319&281&280& & &32.3.11.13.29.43.47&2021&2288 \cr
\noalign{\hrule}
 & &9.5.11.13.67.229&5363&4218& & &9.5.7.121.37.71&893&772 \cr
3034&98732205&4.27.19.31.37.173&67&770&3052&100128105&8.7.19.47.71.193&4625&4818 \cr
 & &16.5.7.11.67.173&173&56& & &32.3.125.11.37.47.73&1175&1168 \cr
\noalign{\hrule}
 & &3.5.37.73.2441&1111&1330& & &9.17.29.107.211&845&634 \cr
3035&98897115&4.25.7.11.19.37.101&611&1314&3053&100174149&4.3.5.169.107.317&319&2 \cr
 & &16.9.13.47.73.101&1313&1128& & &16.5.11.169.29&1859&40 \cr
\noalign{\hrule}
 & &3.361.29.47.67&309&242& & &9.5.11.13.37.421&571&1834 \cr
3036&98900643&4.9.121.19.47.103&5963&11650&3054&100237995&4.3.7.11.131.571&181&50 \cr
 & &16.25.67.89.233&2225&1864& & &16.25.181.571&905&4568 \cr
\noalign{\hrule}
 & &27.17.23.83.113&133&20& & &11.17.37.43.337&11615&5886 \cr
3037&99014103&8.3.5.7.19.23.83&4573&4972&3055&100263229&4.27.5.23.101.109&16571&16456 \cr
 & &64.11.17.113.269&269&352& & &64.9.121.17.73.227&7227&7264 \cr
\noalign{\hrule}
 & &3.5.17.43.83.109&93&638& & &17.23.29.37.239&11869&6372 \cr
3038&99200355&4.9.11.29.31.83&1417&1156&3056&100270777&8.27.11.13.59.83&2507&2390 \cr
 & &32.11.13.289.109&143&272& & &32.3.5.11.23.109.239&545&528 \cr
\noalign{\hrule}
 & &27.19.41.53.89&3335&314& & &3.25.7.13.61.241&197&258 \cr
3039&99212661&4.9.5.23.29.157&209&52&3057&100334325&4.9.5.43.197.241&143&98 \cr
 & &32.5.11.13.19.23&115&2288& & &16.49.11.13.43.197&2167&2408 \cr
\noalign{\hrule}
 & &343.19.97.157&1663&1320& & &3.25.13.97.1061&7139&6654 \cr
3040&99247393&16.3.5.11.97.1663&1365&298&3058&100344075&4.9.5.121.59.1109&7427&4772 \cr
 & &64.9.25.7.13.149&2925&4768& & &32.7.11.1061.1193&1193&1232 \cr
\noalign{\hrule}
 & &5.49.47.89.97&99&146& & &3.5.7.11.17.19.269&97&90 \cr
3041&99408995&4.9.11.73.89.97&85&182&3059&100354485&4.27.25.19.97.269&19459&6634 \cr
 & &16.3.5.7.11.13.17.73&2431&1752& & &16.11.29.31.61.107&6527&7192 \cr
\noalign{\hrule}
 & &27.5.7.139.757&11371&15124& & &9.25.61.71.103&143&82 \cr
3042&99435735&8.19.83.137.199&8327&8190&3060&100370925&4.11.13.41.71.103&725&6588 \cr
 & &32.9.5.7.11.13.19.757&209&208& & &32.27.25.29.61&87&16 \cr
\noalign{\hrule}
}%
}
$$
\eject
\vglue -23 pt
\noindent\hskip 1 in\hbox to 6.5 in{\ 3061 -- 3096 \hfill\fbd 100856847 -- 103604501\frb}
\vskip -9 pt
$$
\vbox{
\nointerlineskip
\halign{\strut
    \vrule \ \ \hfil \frb #\ 
   &\vrule \hfil \ \ \fbb #\frb\ 
   &\vrule \hfil \ \ \frb #\ \hfil
   &\vrule \hfil \ \ \frb #\ 
   &\vrule \hfil \ \ \frb #\ \ \vrule \hskip 2 pt
   &\vrule \ \ \hfil \frb #\ 
   &\vrule \hfil \ \ \fbb #\frb\ 
   &\vrule \hfil \ \ \frb #\ \hfil
   &\vrule \hfil \ \ \frb #\ 
   &\vrule \hfil \ \ \frb #\ \vrule \cr%
\noalign{\hrule}
 & &3.49.13.89.593&12793&16264& & &9.7.11.13.17.23.29&665&2 \cr
3061&100856847&16.11.19.107.1163&7&1170&3079&102153051&4.3.5.49.11.19&289&338 \cr
 & &64.9.5.7.13.19&15&608& & &16.5.169.289&13&680 \cr
\noalign{\hrule}
 & &9.5.19.41.43.67&3509&3844& & &3.5.11.19.67.487&3315&2042 \cr
3062&100993455&8.121.29.961.41&1075&114&3080&102291915&4.9.25.13.17.1021&4757&4432 \cr
 & &32.3.25.121.19.43&121&80& & &128.17.67.71.277&4709&4544 \cr
\noalign{\hrule}
 & &729.97.1429&1079&350& & &25.7.11.17.31.101&1917&608 \cr
3063&101048877&4.25.7.13.83.97&297&782&3081&102461975&64.27.19.31.71&101&70 \cr
 & &16.27.5.7.11.17.23&595&2024& & &256.3.5.7.71.101&213&128 \cr
\noalign{\hrule}
 & &9.5.11.13.19.827&49&94& & &27.25.47.53.61&127&2618 \cr
3064&101113155&4.49.19.47.827&33&860&3082&102566925&4.3.5.7.11.17.127&53&32 \cr
 & &32.3.5.49.11.43&43&784& & &256.11.53.127&1397&128 \cr
\noalign{\hrule}
 & &5.7.107.113.239&63&176& & &7.11.73.101.181&1189&78 \cr
3065&101141215&32.9.5.49.11.107&221&956&3083&102757501&4.3.13.29.41.73&303&230 \cr
 & &256.3.13.17.239&663&128& & &16.9.5.23.29.101&1035&232 \cr
\noalign{\hrule}
 & &5.7.11.13.37.547&11313&15142& & &5.7.19.23.53.127&5423&6642 \cr
3066&101296195&4.27.67.113.419&245&358&3084&102950645&4.81.7.11.17.29.41&53&46 \cr
 & &16.3.5.49.179.419&3759&3352& & &16.9.17.23.29.41.53&1189&1224 \cr
\noalign{\hrule}
 & &27.5.13.197.293&77&274& & &25.7.29.103.197&11869&6156 \cr
3067&101300355&4.5.7.11.137.293&3621&3914&3085&102976825&8.81.11.13.19.83&985&1732 \cr
 & &16.3.7.17.19.71.103&12257&10792& & &64.9.5.197.433&433&288 \cr
\noalign{\hrule}
 & &27.5.31.53.457&551&286& & &81.5.11.19.23.53&3451&3704 \cr
3068&101364885&4.11.13.19.29.457&6095&1068&3086&103182255&16.3.7.17.19.29.463&23&110 \cr
 & &32.3.5.23.53.89&89&368& & &64.5.11.17.23.463&463&544 \cr
\noalign{\hrule}
 & &9.7.13.19.61.107&405&344& & &25.13.29.47.233&1397&1632 \cr
3069&101566647&16.729.5.13.19.43&253&982&3087&103213175&64.3.5.11.17.29.127&1877&282 \cr
 & &64.11.23.43.491&11293&15136& & &256.9.47.1877&1877&1152 \cr
\noalign{\hrule}
 & &9.5.49.11.59.71&2743&3098& & &3.5.11.13.17.19.149&49&100 \cr
3070&101604195&4.49.13.211.1549&213&424&3088&103232415&8.125.49.11.13.19&6739&5364 \cr
 & &64.3.53.71.1549&1549&1696& & &64.9.23.149.293&879&736 \cr
\noalign{\hrule}
 & &5.11.13.23.37.167&1717&2124& & &13.19.43.71.137&5275&5346 \cr
3071&101613655&8.9.5.13.17.59.101&7&58&3089&103310467&4.243.25.11.137.211&659&26 \cr
 & &32.3.7.29.59.101&5133&11312& & &16.81.5.11.13.659&4455&5272 \cr
\noalign{\hrule}
 & &9.5.7.11.13.37.61&1007&214& & &81.7.13.37.379&3629&1298 \cr
3072&101666565&4.3.5.7.19.53.107&143&122&3090&103363533&4.9.11.19.59.191&7025&5306 \cr
 & &16.11.13.19.61.107&107&152& & &16.25.7.281.379&281&200 \cr
\noalign{\hrule}
 & &9.13.529.31.53&869&340& & &9.11.47.71.313&15517&15470 \cr
3073&101690199&8.3.5.11.17.53.79&437&358&3091&103403619&4.5.7.13.17.59.71.263&87&15604 \cr
 & &32.11.17.19.23.179&3553&2864& & &32.3.5.7.29.47.83&1015&1328 \cr
\noalign{\hrule}
 & &3.25.11.31.41.97&303&148& & &9.25.29.83.191&11189&4664 \cr
3074&101711775&8.9.5.37.97.101&697&212&3092&103440825&16.11.53.67.167&285&452 \cr
 & &64.17.37.41.53&629&1696& & &128.3.5.19.53.113&2147&3392 \cr
\noalign{\hrule}
 & &25.43.173.547&261&286& & &9.17.31.139.157&55&472 \cr
3075&101728325&4.9.11.13.29.43.173&5251&10960&3093&103506489&16.3.5.11.59.157&95&62 \cr
 & &128.3.5.59.89.137&12193&11328& & &64.25.19.31.59&475&1888 \cr
\noalign{\hrule}
 & &361.521.541&451&90& & &9.5.31.47.1579&867&712 \cr
3076&101751821&4.9.5.11.41.521&255&266&3094&103527135&16.27.289.47.89&253&206 \cr
 & &16.27.25.7.17.19.41&4879&5400& & &64.11.17.23.89.103&34799&36256 \cr
\noalign{\hrule}
 & &25.7.37.79.199&1159&234& & &125.289.47.61&711&2156 \cr
3077&101793475&4.9.13.19.61.79&605&422&3095&103570375&8.9.25.49.11.79&689&136 \cr
 & &16.3.5.121.19.211&2299&5064& & &128.3.7.13.17.53&371&2496 \cr
\noalign{\hrule}
 & &5.7.11.37.71.101&221&276& & &7.11.13.29.43.83&131&450 \cr
3078&102150895&8.3.13.17.23.37.101&971&342&3096&103604501&4.9.25.13.43.131&83&476 \cr
 & &32.27.19.23.971&11799&15536& & &32.3.25.7.17.83&425&48 \cr
\noalign{\hrule}
}%
}
$$
\eject
\vglue -23 pt
\noindent\hskip 1 in\hbox to 6.5 in{\ 3097 -- 3132 \hfill\fbd 103822425 -- 107027487\frb}
\vskip -9 pt
$$
\vbox{
\nointerlineskip
\halign{\strut
    \vrule \ \ \hfil \frb #\ 
   &\vrule \hfil \ \ \fbb #\frb\ 
   &\vrule \hfil \ \ \frb #\ \hfil
   &\vrule \hfil \ \ \frb #\ 
   &\vrule \hfil \ \ \frb #\ \ \vrule \hskip 2 pt
   &\vrule \ \ \hfil \frb #\ 
   &\vrule \hfil \ \ \fbb #\frb\ 
   &\vrule \hfil \ \ \frb #\ \hfil
   &\vrule \hfil \ \ \frb #\ 
   &\vrule \hfil \ \ \frb #\ \vrule \cr%
\noalign{\hrule}
 & &27.25.49.43.73&187&488& & &5.101.109.1913&307&198 \cr
3097&103822425&16.7.11.17.61.73&2265&2188&3115&105301085&4.9.11.307.1913&1417&496 \cr
 & &128.3.5.17.151.547&9299&9664& & &128.3.11.13.31.109&1209&704 \cr
\noalign{\hrule}
 & &3.11.23.41.47.71&455&488& & &9.7.31.199.271&3305&5096 \cr
3098&103844103&16.5.7.13.47.61.71&1107&4444&3116&105323337&16.5.343.13.661&527&1188 \cr
 & &128.27.5.11.41.101&505&576& & &128.27.11.13.17.31&561&832 \cr
\noalign{\hrule}
 & &9.17.29.41.571&309&880& & &3.5.29.37.79.83&1669&1254 \cr
3099&103874607&32.27.5.11.17.103&943&808&3117&105534915&4.9.11.19.29.1669&235&26 \cr
 & &512.11.23.41.101&2323&2816& & &16.5.13.47.1669&1669&4888 \cr
\noalign{\hrule}
 & &9.5.11.13.19.23.37&191&246& & &3.5.17.29.109.131&29419&27566 \cr
3100&104047515&4.27.13.37.41.191&493&506&3118&105593205&4.7.11.13.31.73.179&161&18 \cr
 & &16.11.17.23.29.41.191&5539&5576& & &16.9.49.23.31.73&6789&9016 \cr
\noalign{\hrule}
 & &9.7.13.29.41.107&33&74& & &27.5.7.23.43.113&203&418 \cr
3101&104195637&4.27.7.11.13.29.37&1819&820&3119&105610365&4.49.11.19.29.113&1087&156 \cr
 & &32.5.11.17.41.107&85&176& & &32.3.13.29.1087&1087&6032 \cr
\noalign{\hrule}
 & &27.49.11.13.19.29&575&62& & &9.11.13.103.797&2407&1610 \cr
3102&104243139&4.25.11.23.29.31&171&148&3120&105651117&4.3.5.7.11.23.29.83&1703&1036 \cr
 & &32.9.25.19.31.37&925&496& & &32.5.49.13.37.131&4847&3920 \cr
\noalign{\hrule}
 & &9.5.11.169.29.43&1513&2228& & &9.25.7.11.17.359&1417&1058 \cr
3103&104317785&8.3.13.17.89.557&725&946&3121&105734475&4.7.13.17.529.109&2077&660 \cr
 & &32.25.11.29.43.89&89&80& & &32.3.5.11.23.31.67&713&1072 \cr
\noalign{\hrule}
 & &9.5.17.19.43.167&193&194& & &3.7.17.31.73.131&117&100 \cr
3104&104375835&4.5.17.19.97.167.193&5709&2536&3122&105833721&8.27.25.13.73.131&1837&134 \cr
 & &64.3.11.173.193.317&61181&60896& & &32.25.11.67.167&11189&4400 \cr
\noalign{\hrule}
 & &7.11.13.29.59.61&129&190& & &9.37.41.43.181&275&94 \cr
3105&104475371&4.3.5.7.13.19.43.59&83&330&3123&106261299&4.25.11.37.43.47&1383&208 \cr
 & &16.9.25.11.43.83&3569&1800& & &128.3.11.13.461&5993&704 \cr
\noalign{\hrule}
 & &9.5.7.11.97.311&281&204& & &9.5.7.17.31.641&437&158 \cr
3106&104528655&8.27.17.281.311&385&74&3124&106409205&4.19.23.79.641&539&102 \cr
 & &32.5.7.11.37.281&281&592& & &16.3.49.11.17.79&553&88 \cr
\noalign{\hrule}
 & &3.25.17.23.43.83&1767&3178& & &9.5.7.11.13.17.139&95&44 \cr
3107&104660925&4.9.5.7.19.31.227&451&2494&3125&106441335&8.3.25.7.121.13.19&1649&1376 \cr
 & &16.7.11.29.41.43&451&1624& & &512.17.19.43.97&4171&4864 \cr
\noalign{\hrule}
 & &3.13.841.31.103&1695&1292& & &9.25.131.3613&9229&8836 \cr
3108&104727207&8.9.5.17.19.29.113&77&68&3126&106493175&8.3.5.11.2209.839&67&772 \cr
 & &64.7.11.289.19.113&38437&39776& & &64.11.47.67.193&12931&16544 \cr
\noalign{\hrule}
 & &9.7.11.13.29.401&535&158& & &23.37.41.43.71&1335&1292 \cr
3109&104765661&4.5.79.107.401&161&240&3127&106522223&8.3.5.17.19.23.41.89&2971&924 \cr
 & &128.3.25.7.23.107&2675&1472& & &64.9.7.11.17.2971&32681&34272 \cr
\noalign{\hrule}
 & &5.49.11.17.29.79&257&138& & &27.11.19.43.439&2215&10556 \cr
3110&104962165&4.3.7.11.23.29.257&79&240&3128&106522911&8.5.7.13.29.443&207&236 \cr
 & &128.9.5.79.257&257&576& & &64.9.5.7.13.23.59&3835&5152 \cr
\noalign{\hrule}
 & &27.5.7.277.401&10757&3278& & &3.7.79.239.269&583&1300 \cr
3111&104967765&4.11.31.149.347&401&4218&3129&106658769&8.25.11.13.53.79&269&126 \cr
 & &16.3.19.37.401&19&296& & &32.9.5.7.53.269&159&80 \cr
\noalign{\hrule}
 & &3.11.169.47.401&129&272& & &27.5.31.71.359&319&40 \cr
3112&105109719&32.9.13.17.43.47&2005&16&3130&106671465&16.3.25.11.29.71&31&244 \cr
 & &1024.5.401&5&512& & &128.29.31.61&61&1856 \cr
\noalign{\hrule}
 & &3.11.19.271.619&2569&2580& & &5.121.13.41.331&2049&2254 \cr
3113&105178623&8.9.5.7.43.367.619&3199&104&3131&106735915&4.3.49.121.23.683&765&82 \cr
 & &128.49.13.43.457&27391&29248& & &16.27.5.7.17.23.41&621&952 \cr
\noalign{\hrule}
 & &9.11.17.19.37.89&715&86& & &81.7.23.29.283&209&412 \cr
3114&105300261&4.5.121.13.19.43&293&312&3132&107027487&8.3.11.19.103.283&455&172 \cr
 & &64.3.169.43.293&7267&9376& & &64.5.7.13.43.103&4429&2080 \cr
\noalign{\hrule}
}%
}
$$
\eject
\vglue -23 pt
\noindent\hskip 1 in\hbox to 6.5 in{\ 3133 -- 3168 \hfill\fbd 107039445 -- 110562771\frb}
\vskip -9 pt
$$
\vbox{
\nointerlineskip
\halign{\strut
    \vrule \ \ \hfil \frb #\ 
   &\vrule \hfil \ \ \fbb #\frb\ 
   &\vrule \hfil \ \ \frb #\ \hfil
   &\vrule \hfil \ \ \frb #\ 
   &\vrule \hfil \ \ \frb #\ \ \vrule \hskip 2 pt
   &\vrule \ \ \hfil \frb #\ 
   &\vrule \hfil \ \ \fbb #\frb\ 
   &\vrule \hfil \ \ \frb #\ \hfil
   &\vrule \hfil \ \ \frb #\ 
   &\vrule \hfil \ \ \frb #\ \vrule \cr%
\noalign{\hrule}
 & &3.5.19.47.61.131&2499&3658& & &9.11.17.19.43.79&10537&1550 \cr
3133&107039445&4.9.5.49.17.31.59&517&3172&3151&108625869&4.25.31.41.257&663&622 \cr
 & &32.7.11.13.47.61&77&208& & &16.3.5.13.17.31.311&2015&2488 \cr
\noalign{\hrule}
 & &5.49.13.19.29.61&1219&1770& & &9.11.17.19.43.79&1217&1690 \cr
3134&107051035&4.3.25.13.23.53.59&671&96&3152&108625869&4.5.169.79.1217&95&1122 \cr
 & &256.9.11.53.61&477&1408& & &16.3.25.11.13.17.19&25&104 \cr
\noalign{\hrule}
 & &81.5.7.23.31.53&13243&16808& & &5.11.13.19.53.151&1779&1090 \cr
3135&107131815&16.11.17.19.41.191&75&116&3153&108720755&4.3.25.11.109.593&159&434 \cr
 & &128.3.25.11.17.19.29&6061&5440& & &16.9.7.31.53.109&981&1736 \cr
\noalign{\hrule}
 & &11.83.239.491&211&702& & &29.37.41.2473&51233&50160 \cr
3136&107139637&4.27.13.211.239&1069&830&3154&108794689&32.3.5.7.11.13.19.563&7419&1226 \cr
 & &16.3.5.13.83.1069&1069&1560& & &128.9.613.2473&613&576 \cr
\noalign{\hrule}
 & &5.11.17.19.37.163&23&186& & &3.23.31.151.337&395&58 \cr
3137&107140715&4.3.5.17.23.31.37&171&356&3155&108847293&4.5.23.29.31.79&317&396 \cr
 & &32.27.19.23.89&2403&368& & &32.9.5.11.29.317&3487&6960 \cr
\noalign{\hrule}
 & &11.17.47.73.167&797&2634& & &9.25.7.257.269&247&22 \cr
3138&107146699&4.3.17.439.797&3333&4130&3156&108884475&4.7.11.13.19.257&83&174 \cr
 & &16.9.5.7.11.59.101&4545&3304& & &16.3.11.19.29.83&6061&664 \cr
\noalign{\hrule}
 & &125.11.13.17.353&2643&3358& & &7.13.47.71.359&297&626 \cr
3139&107267875&4.3.25.23.73.881&153&728&3157&109016453&4.27.11.313.359&1505&2444 \cr
 & &64.27.7.13.17.73&511&864& & &32.9.5.7.13.43.47&215&144 \cr
\noalign{\hrule}
 & &25.1369.3137&93&92& & &9.5.107.139.163&679&572 \cr
3140&107363825&8.3.5.23.31.37.3137&559&3696&3158&109093455&8.5.7.11.13.97.163&69&8896 \cr
 & &256.9.7.11.13.31.43&39897&38528& & &1024.3.23.139&23&512 \cr
\noalign{\hrule}
 & &5.11.29.31.41.53&10647&2432& & &7.13.17.19.47.79&243&3956 \cr
3141&107443985&256.9.7.169.19&35&22&3159&109136209&8.243.7.23.43&191&110 \cr
 & &1024.3.5.49.11.13&637&1536& & &32.3.5.11.23.191&10505&1104 \cr
\noalign{\hrule}
 & &9.5.7.47.53.137&11033&3772& & &27.11.193.1907&805&1102 \cr
3142&107499105&8.11.17.23.41.59&147&106&3160&109311147&4.5.7.19.23.29.193&237&430 \cr
 & &32.3.49.17.53.59&413&272& & &16.3.25.7.19.43.79&10507&8600 \cr
\noalign{\hrule}
 & &9.11.19.23.47.53&305&2186& & &81.5.19.23.619&737&118 \cr
3143&107768133&4.5.23.61.1093&141&164&3161&109553715&4.9.11.23.59.67&133&74 \cr
 & &32.3.41.47.1093&1093&656& & &16.7.11.19.37.67&2479&616 \cr
\noalign{\hrule}
 & &9.7.11.29.31.173&425&598& & &25.23.107.1783&341&234 \cr
3144&107780211&4.3.25.7.13.17.23.29&347&550&3162&109699075&4.9.11.13.31.1783&1403&380 \cr
 & &16.625.11.17.347&5899&5000& & &32.3.5.13.19.23.61&1159&624 \cr
\noalign{\hrule}
 & &243.5.7.19.23.29&859&1004& & &9.7.23.191.397&1451&10582 \cr
3145&107783865&8.3.7.19.251.859&4147&10160&3163&109873323&4.11.13.37.1451&485&966 \cr
 & &256.5.11.13.29.127&1651&1408& & &16.3.5.7.11.23.97&485&88 \cr
\noalign{\hrule}
 & &9.5.41.71.827&2387&1748& & &3.5.11.37.67.269&101&168 \cr
3146&108332865&8.7.11.19.23.31.41&213&500&3164&110030415&16.9.5.7.11.37.101&377&118 \cr
 & &64.3.125.11.19.71&475&352& & &64.13.29.59.101&5959&12064 \cr
\noalign{\hrule}
 & &19.23.29.83.103&455&1122& & &3.11.13.17.79.191&1919&2110 \cr
3147&108341477&4.3.5.7.11.13.17.103&19&84&3165&110044077&4.5.11.13.19.101.211&1817&504 \cr
 & &32.9.49.11.17.19&539&2448& & &64.9.5.7.19.23.79&1311&1120 \cr
\noalign{\hrule}
 & &3.5.19.23.61.271&11543&4988& & &13.43.47.59.71&21231&18458 \cr
3148&108360705&8.7.17.29.43.97&295&198&3166&110057597&4.9.7.11.337.839&925&86 \cr
 & &32.9.5.7.11.43.59&2537&3696& & &16.3.25.7.11.37.43&925&1848 \cr
\noalign{\hrule}
 & &27.5.121.29.229&1519&1748& & &3.25.7.13.19.23.37&2941&3016 \cr
3149&108480735&8.5.49.19.23.29.31&1&666&3167&110353425&16.169.17.19.29.173&207&3080 \cr
 & &32.9.7.31.37&37&3472& & &256.9.5.7.11.23.29&319&384 \cr
\noalign{\hrule}
 & &25.7.121.23.223&6213&6052& & &3.11.23.31.37.127&2071&850 \cr
3150&108606575&8.3.5.11.17.19.89.109&1449&404&3168&110562771&4.25.17.19.31.109&549&226 \cr
 & &64.27.7.23.89.101&2727&2848& & &16.9.61.109.113&6649&2712 \cr
\noalign{\hrule}
}%
}
$$
\eject
\vglue -23 pt
\noindent\hskip 1 in\hbox to 6.5 in{\ 3169 -- 3204 \hfill\fbd 110755385 -- 113733927\frb}
\vskip -9 pt
$$
\vbox{
\nointerlineskip
\halign{\strut
    \vrule \ \ \hfil \frb #\ 
   &\vrule \hfil \ \ \fbb #\frb\ 
   &\vrule \hfil \ \ \frb #\ \hfil
   &\vrule \hfil \ \ \frb #\ 
   &\vrule \hfil \ \ \frb #\ \ \vrule \hskip 2 pt
   &\vrule \ \ \hfil \frb #\ 
   &\vrule \hfil \ \ \fbb #\frb\ 
   &\vrule \hfil \ \ \frb #\ \hfil
   &\vrule \hfil \ \ \frb #\ 
   &\vrule \hfil \ \ \frb #\ \vrule \cr%
\noalign{\hrule}
 & &5.13.71.103.233&121&1044& & &7.11.19.193.397&2117&2250 \cr
3169&110755385&8.9.121.29.103&233&130&3187&112096523&4.9.125.29.73.193&763&5588 \cr
 & &32.3.5.13.29.233&87&16& & &32.3.5.7.11.109.127&1635&2032 \cr
\noalign{\hrule}
 & &3.121.37.73.113&167&240& & &27.37.107.1049&691&358 \cr
3170&110792319&32.9.5.11.113.167&7&106&3188&112130757&4.3.107.179.691&185&506 \cr
 & &128.5.7.53.167&1855&10688& & &16.5.11.23.37.179&1969&920 \cr
\noalign{\hrule}
 & &9.11.23.181.269&145&398& & &7.47.569.599&261&308 \cr
3171&110864853&4.3.5.29.199.269&3289&2482&3189&112133399&8.9.49.11.29.599&569&30 \cr
 & &16.5.11.13.17.23.73&949&680& & &32.27.5.29.569&135&464 \cr
\noalign{\hrule}
 & &3.5.7.19.23.41.59&561&1504& & &9.5.169.113.131&989&976 \cr
3172&110995815&64.9.11.17.19.47&59&40&3190&112576815&32.3.13.23.43.61.113&2489&110 \cr
 & &1024.5.17.47.59&799&512& & &128.5.11.19.43.131&817&704 \cr
\noalign{\hrule}
 & &25.7.11.529.109&477&2248& & &9.5.7.13.107.257&41&50 \cr
3173&110997425&16.9.23.53.281&469&750&3191&112608405&4.125.41.107.257&1419&11956 \cr
 & &64.27.125.7.67&335&864& & &32.3.49.11.43.61&2623&1232 \cr
\noalign{\hrule}
 & &3.7.17.31.79.127&105083&105610& & &3.25.11.17.29.277&101&176 \cr
3174&111035211&4.5.11.41.59.179.233&553&612&3192&112662825&32.121.17.29.101&2493&436 \cr
 & &32.9.7.11.17.41.79.179&1969&1968& & &256.9.109.277&109&384 \cr
\noalign{\hrule}
 & &9.25.7.11.13.17.29&457&138& & &27.25.11.17.19.47&221&202 \cr
3175&111035925&4.27.5.13.23.457&649&1106&3193&112718925&4.3.25.11.13.289.101&2077&1102 \cr
 & &16.7.11.23.59.79&1357&632& & &16.19.29.31.67.101&6767&7192 \cr
\noalign{\hrule}
 & &5.7.11.17.361.47&1287&1240& & &7.11.43.67.509&145&156 \cr
3176&111049015&16.9.25.121.13.17.31&2303&722&3194&112915033&8.3.5.13.29.67.509&301&2244 \cr
 & &64.3.49.13.361.47&91&96& & &64.9.7.11.13.17.43&221&288 \cr
\noalign{\hrule}
 & &5.11.19.97.1097&1071&26& & &3.19.23.43.2003&1023&980 \cr
3177&111197405&4.9.7.13.17.97&2519&2428&3195&112915119&8.9.5.49.11.19.23.31&5&166 \cr
 & &32.3.11.229.607&1821&3664& & &32.25.7.11.31.83&6391&12400 \cr
\noalign{\hrule}
 & &9.25.343.11.131&59&284& & &9.5.13.17.31.367&4343&2508 \cr
3178&111209175&8.11.59.71.131&259&390&3196&113144265&8.27.11.19.43.101&1259&1468 \cr
 & &32.3.5.7.13.37.71&481&1136& & &64.43.367.1259&1259&1376 \cr
\noalign{\hrule}
 & &5.7.53.173.347&87&260& & &3.19.37.223.241&455&214 \cr
3179&111357505&8.3.25.7.13.29.53&1111&1164&3197&113343987&4.5.7.13.19.37.107&2453&1062 \cr
 & &64.9.11.29.97.101&30943&29088& & &16.9.7.11.59.223&231&472 \cr
\noalign{\hrule}
 & &125.17.19.31.89&2783&1092& & &9.5.11.29.53.149&1409&1462 \cr
3180&111394625&8.3.7.121.13.17.23&267&124&3198&113361435&4.5.17.43.149.1409&319&6726 \cr
 & &64.9.7.11.31.89&99&224& & &16.3.11.17.19.29.59&323&472 \cr
\noalign{\hrule}
 & &27.17.23.61.173&3575&634& & &25.7.13.19.43.61&583&492 \cr
3181&111408021&4.9.25.11.13.317&61&56&3199&113379175&8.3.11.19.41.53.61&559&600 \cr
 & &64.5.7.11.61.317&2219&1760& & &128.9.25.11.13.43.53&583&576 \cr
\noalign{\hrule}
 & &25.19.103.2281&29491&27534& & &3.5.7.13.29.47.61&305&682 \cr
3182&111597925&4.3.7.11.13.353.383&3725&1254&3200&113490195&4.25.11.31.3721&46683&46342 \cr
 & &16.9.25.121.19.149&1089&1192& & &16.27.7.13.17.19.29.47&153&152 \cr
\noalign{\hrule}
 & &3.11.13.37.79.89&11&100& & &3.7.121.13.19.181&153&1420 \cr
3183&111603063&8.25.121.13.79&999&2026&3201&113600487&8.27.5.17.19.71&175&148 \cr
 & &32.27.37.1013&1013&144& & &64.125.7.37.71&4625&2272 \cr
\noalign{\hrule}
 & &3.7.17.23.31.439&1015&302& & &9.7.11.19.89.97&575&1268 \cr
3184&111743499&4.5.49.17.29.151&5707&7128&3202&113670711&8.25.23.89.317&203&114 \cr
 & &64.81.11.13.439&297&416& & &32.3.25.7.19.23.29&667&400 \cr
\noalign{\hrule}
 & &9.7.13.19.23.313&125&148& & &25.49.13.37.193&109&84 \cr
3185&112023639&8.3.125.19.37.313&869&556&3203&113720425&8.3.343.13.37.109&213&4246 \cr
 & &64.5.11.37.79.139&25715&27808& & &32.9.11.71.193&781&144 \cr
\noalign{\hrule}
 & &9.5.1327.1877&28919&30796& & &9.289.73.599&321&920 \cr
3186&112085055&8.121.239.7699&42225&42464&3204&113733927&16.27.5.17.23.107&1001&1460 \cr
 & &512.3.25.11.563.1327&2815&2816& & &128.25.7.11.13.73&1001&1600 \cr
\noalign{\hrule}
}%
}
$$
\eject
\vglue -23 pt
\noindent\hskip 1 in\hbox to 6.5 in{\ 3205 -- 3240 \hfill\fbd 113947119 -- 117195111\frb}
\vskip -9 pt
$$
\vbox{
\nointerlineskip
\halign{\strut
    \vrule \ \ \hfil \frb #\ 
   &\vrule \hfil \ \ \fbb #\frb\ 
   &\vrule \hfil \ \ \frb #\ \hfil
   &\vrule \hfil \ \ \frb #\ 
   &\vrule \hfil \ \ \frb #\ \ \vrule \hskip 2 pt
   &\vrule \ \ \hfil \frb #\ 
   &\vrule \hfil \ \ \fbb #\frb\ 
   &\vrule \hfil \ \ \frb #\ \hfil
   &\vrule \hfil \ \ \frb #\ 
   &\vrule \hfil \ \ \frb #\ \vrule \cr%
\noalign{\hrule}
 & &9.11.13.29.43.71&425&48& & &81.13.17.41.157&4213&2224 \cr
3205&113947119&32.27.25.17.71&971&946&3223&115228737&32.9.11.139.383&119&20 \cr
 & &128.11.17.43.971&971&1088& & &256.5.7.17.383&2681&640 \cr
\noalign{\hrule}
 & &5.7.11.17.23.757&351&406& & &5.49.31.43.353&285&68 \cr
3206&113954995&4.27.49.13.17.23.29&4895&1514&3224&115284505&8.3.25.7.17.19.43&363&62 \cr
 & &16.9.5.11.89.757&89&72& & &32.9.121.19.31&171&1936 \cr
\noalign{\hrule}
 & &11.29.271.1319&819&500& & &9.5.7.13.47.599&7733&13232 \cr
3207&114026231&8.9.125.7.13.271&233&38&3225&115286535&32.11.19.37.827&617&210 \cr
 & &32.3.25.7.19.233&3325&11184& & &128.3.5.7.19.617&617&1216 \cr
\noalign{\hrule}
 & &3.7.11.13.191.199&545&1938& & &9.5.53.139.349&169&308 \cr
3208&114141027&4.9.5.11.17.19.109&523&676&3226&115698735&8.5.7.11.169.349&53&402 \cr
 & &32.5.169.19.523&2615&3952& & &32.3.11.13.53.67&871&176 \cr
\noalign{\hrule}
 & &9.5.11.13.41.433&167&266& & &3.5.7.23.191.251&779&26 \cr
3209&114240555&4.5.7.13.19.41.167&1917&748&3227&115777515&4.13.19.41.191&5159&4968 \cr
 & &32.27.11.17.19.71&1349&816& & &64.27.7.11.23.67&603&352 \cr
\noalign{\hrule}
 & &7.83.239.823&121&702& & &9.5.121.89.239&2767&2678 \cr
3210&114280957&4.27.121.13.239&425&664&3228&115820595&4.13.103.239.2767&2937&170 \cr
 & &64.3.25.13.17.83&425&1248& & &16.3.5.11.17.89.103&103&136 \cr
\noalign{\hrule}
 & &9.7.17.19.41.137&667&530& & &11.31.59.73.79&21&52 \cr
3211&114300333&4.5.17.23.29.41.53&15&682&3229&116026273&8.3.7.11.13.59.79&4495&3942 \cr
 & &16.3.25.11.31.53&583&6200& & &32.81.5.29.31.73&405&464 \cr
\noalign{\hrule}
 & &27.5.23.137.269&1177&2522& & &3.7.11.13.29.31.43&803&530 \cr
3212&114428565&4.11.13.23.97.107&569&822&3230&116086971&4.5.121.29.53.73&3689&180 \cr
 & &16.3.97.137.569&569&776& & &32.9.25.7.17.31&425&48 \cr
\noalign{\hrule}
 & &9.11.13.23.53.73&629&590& & &5.7.17.19.43.239&77&162 \cr
3213&114526269&4.3.5.11.17.37.59.73&2951&1064&3231&116181485&4.81.49.11.19.43&1495&956 \cr
 & &64.7.13.19.59.227&7847&7264& & &32.27.5.13.23.239&299&432 \cr
\noalign{\hrule}
 & &7.13.19.23.43.67&461&528& & &5.23.31.127.257&999&286 \cr
3214&114568727&32.3.7.11.13.19.461&5&138&3232&116358035&4.27.11.13.37.127&775&368 \cr
 & &128.9.5.23.461&2305&576& & &128.3.25.13.23.31&65&192 \cr
\noalign{\hrule}
 & &27.5.11.19.31.131&1139&302& & &31.47.109.733&1323&2056 \cr
3215&114581115&4.5.17.19.67.151&259&1014&3233&116409929&16.27.49.47.257&221&550 \cr
 & &16.3.7.169.17.37&6253&952& & &64.9.25.7.11.13.17&13923&8800 \cr
\noalign{\hrule}
 & &81.13.17.37.173&12151&5750& & &841.97.1427&293&1134 \cr
3216&114584301&4.125.23.29.419&153&572&3234&116410379&4.81.7.97.293&583&290 \cr
 & &32.9.5.11.13.17.23&253&80& & &16.9.5.7.11.29.53&693&2120 \cr
\noalign{\hrule}
 & &5.7.19.23.59.127&3901&11394& & &7.13.19.31.41.53&963&1210 \cr
3217&114605435&4.27.47.83.211&1001&2900&3235&116470627&4.9.5.7.121.31.107&169&48 \cr
 & &32.3.25.7.11.13.29&1885&528& & &128.27.5.169.107&2889&4160 \cr
\noalign{\hrule}
 & &3.11.13.397.673&265&408& & &5.7.11.23.59.223&169&54 \cr
3218&114620649&16.9.5.17.53.397&437&40&3236&116505235&4.27.7.11.169.59&437&850 \cr
 & &256.25.17.19.23&10925&2176& & &16.3.25.13.17.19.23&1105&456 \cr
\noalign{\hrule}
 & &25.13.31.59.193&1177&1332& & &7.17.73.89.151&573&484 \cr
3219&114724025&8.9.5.11.37.59.107&1537&2422&3237&116744593&8.3.121.17.73.191&2793&6040 \cr
 & &32.3.7.11.29.53.173&32277&30448& & &128.9.5.49.19.151&665&576 \cr
\noalign{\hrule}
 & &81.7.31.61.107&851&1040& & &3.25.13.257.467&201&266 \cr
3220&114725079&32.3.5.13.23.37.107&19&88&3238&117018525&4.9.5.7.19.67.257&3949&934 \cr
 & &512.5.11.13.19.37&9139&14080& & &16.7.11.359.467&359&616 \cr
\noalign{\hrule}
 & &3.5.11.829.839&1653&2492& & &27.23.29.67.97&233&434 \cr
3221&114762615&8.9.7.11.19.29.89&629&1252&3239&117040491&4.9.7.31.97.233&455&2552 \cr
 & &64.17.29.37.313&11581&15776& & &64.5.49.11.13.29&539&2080 \cr
\noalign{\hrule}
 & &9.13.17.151.383&35&186& & &9.11.47.89.283&1763&784 \cr
3222&115029837&4.27.5.7.31.383&97&286&3240&117195111&32.49.41.43.47&2033&270 \cr
 & &16.5.11.13.31.97&3007&440& & &128.27.5.19.107&535&3648 \cr
\noalign{\hrule}
}%
}
$$
\eject
\vglue -23 pt
\noindent\hskip 1 in\hbox to 6.5 in{\ 3241 -- 3276 \hfill\fbd 117207125 -- 120125655\frb}
\vskip -9 pt
$$
\vbox{
\nointerlineskip
\halign{\strut
    \vrule \ \ \hfil \frb #\ 
   &\vrule \hfil \ \ \fbb #\frb\ 
   &\vrule \hfil \ \ \frb #\ \hfil
   &\vrule \hfil \ \ \frb #\ 
   &\vrule \hfil \ \ \frb #\ \ \vrule \hskip 2 pt
   &\vrule \ \ \hfil \frb #\ 
   &\vrule \hfil \ \ \fbb #\frb\ 
   &\vrule \hfil \ \ \frb #\ \hfil
   &\vrule \hfil \ \ \frb #\ 
   &\vrule \hfil \ \ \frb #\ \vrule \cr%
\noalign{\hrule}
 & &125.7.29.31.149&3519&2774& & &27.5.13.17.23.173&61&56 \cr
3241&117207125&4.9.25.17.19.23.73&1519&6556&3259&118713465&16.3.7.17.23.61.173&3575&634 \cr
 & &32.3.49.11.31.149&77&48& & &64.25.7.11.13.317&2219&1760 \cr
\noalign{\hrule}
 & &9.7.19.29.31.109&1495&1666& & &3.7.13.257.1693&2585&2494 \cr
3242&117295227&4.5.343.13.17.23.31&1133&4698&3260&118782573&4.5.11.29.43.47.257&31&288 \cr
 & &16.81.11.13.29.103&927&1144& & &256.9.5.31.43.47&10105&11904 \cr
\noalign{\hrule}
 & &81.7.29.37.193&275&304& & &9.5.11.17.29.487&1085&598 \cr
3243&117419463&32.27.25.7.11.19.37&127&386&3261&118845045&4.25.7.13.23.29.31&2217&7208 \cr
 & &128.25.11.127.193&1397&1600& & &64.3.17.53.739&739&1696 \cr
\noalign{\hrule}
 & &13.19.29.47.349&165&212& & &243.5.7.11.31.41&437&778 \cr
3244&117494689&8.3.5.11.19.53.349&4437&598&3262&118908405&4.7.19.23.41.389&1747&7200 \cr
 & &32.27.13.17.23.29&621&272& & &256.9.25.1747&1747&640 \cr
\noalign{\hrule}
 & &3.17.47.191.257&2255&2114& & &3.7.11.19.41.661&325&336 \cr
3245&117661539&4.5.7.11.41.151.191&9729&1898&3263&118946289&32.9.25.49.13.19.41&71&4726 \cr
 & &16.9.5.13.23.47.73&1495&1752& & &128.5.17.71.139&11815&4544 \cr
\noalign{\hrule}
 & &3.5.121.241.269&491&854& & &169.47.71.211&3519&11462 \cr
3246&117664635&4.7.61.241.491&1089&598&3264&118994083&4.9.11.17.23.521&65&456 \cr
 & &16.9.121.13.23.61&793&552& & &64.27.5.11.13.19&297&3040 \cr
\noalign{\hrule}
 & &9.5.11.29.59.139&187&1898& & &23.29.31.73.79&891&1558 \cr
3247&117725355&4.3.121.13.17.73&59&62&3265&119244259&4.81.11.19.41.73&1247&140 \cr
 & &16.13.17.31.59.73&2263&1768& & &32.3.5.7.11.29.43&385&2064 \cr
\noalign{\hrule}
 & &5.49.13.103.359&183&176& & &3.25.49.13.41.61&1863&638 \cr
3248&117771745&32.3.5.7.11.13.61.103&3&718&3266&119485275&4.243.11.13.23.29&1525&1148 \cr
 & &128.9.61.359&61&576& & &32.25.7.23.41.61&23&16 \cr
\noalign{\hrule}
 & &25.19.43.73.79&1807&1332& & &49.13.163.1151&1635&484 \cr
3249&117790975&8.9.13.37.79.139&2365&558&3267&119509481&8.3.5.49.121.109&85&36 \cr
 & &32.81.5.11.31.43&891&496& & &64.27.25.17.109&11475&3488 \cr
\noalign{\hrule}
 & &5.11.13.37.61.73&1747&954& & &9.11.37.127.257&775&368 \cr
3250&117804115&4.9.5.11.53.1747&2331&584&3268&119556657&32.25.23.31.257&999&286 \cr
 & &64.81.7.37.73&81&224& & &128.27.5.11.13.37&65&192 \cr
\noalign{\hrule}
 & &7.11.13.29.31.131&599&2040& & &243.25.7.29.97&241&484 \cr
3251&117886769&16.3.5.17.31.599&377&222&3269&119622825&8.7.121.97.241&1665&986 \cr
 & &64.9.13.17.29.37&629&288& & &32.9.5.11.17.29.37&407&272 \cr
\noalign{\hrule}
 & &9.7.11.169.19.53&913&94& & &9.121.13.43.197&4285&4186 \cr
3252&117936819&4.121.13.47.83&245&366&3270&119923947&4.5.7.11.169.23.857&3349&16362 \cr
 & &16.3.5.49.61.83&581&2440& & &16.81.5.17.101.197&765&808 \cr
\noalign{\hrule}
 & &13.23.43.67.137&2331&550& & &5.7.13.47.71.79&1277&4332 \cr
3253&118014403&4.9.25.7.11.23.37&43&302&3271&119948465&8.3.7.361.1277&705&572 \cr
 & &16.3.5.11.43.151&1661&120& & &64.9.5.11.13.19.47&209&288 \cr
\noalign{\hrule}
 & &27.25.17.41.251&3367&2908& & &3.625.7.13.19.37&947&1428 \cr
3254&118089225&8.7.13.37.41.727&395&1122&3272&119949375&8.9.5.49.17.947&1111&1094 \cr
 & &32.3.5.7.11.13.17.79&1027&1232& & &32.11.101.547.947&95647&96272 \cr
\noalign{\hrule}
 & &27.5.19.29.37.43&539&1052& & &9.37.41.59.149&23&14 \cr
3255&118346535&8.5.49.11.29.263&639&376&3273&120023523&4.7.23.41.59.149&3025&3084 \cr
 & &128.9.7.11.47.71&3337&4928& & &32.3.25.7.121.23.257&44975&44528 \cr
\noalign{\hrule}
 & &9.7.11.13.59.223&437&850& & &27.125.7.13.17.23&209&2666 \cr
3256&118531413&4.25.17.19.23.223&169&54&3274&120085875&4.11.17.19.31.43&5243&4770 \cr
 & &16.27.5.169.17.19&1105&456& & &16.9.5.49.53.107&749&424 \cr
\noalign{\hrule}
 & &29.37.47.2351&44175&42812& & &3.7.11.13.23.37.47&731&120 \cr
3257&118563281&8.3.25.7.11.19.31.139&4033&1392&3275&120110991&16.9.5.7.11.17.43&37&26 \cr
 & &256.9.11.29.37.109&1199&1152& & &64.5.13.17.37.43&215&544 \cr
\noalign{\hrule}
 & &3.11.13.23.53.227&441&142& & &9.5.13.17.47.257&1133&1180 \cr
3258&118709877&4.27.49.71.227&845&1072&3276&120125655&8.25.11.13.17.59.103&27&248 \cr
 & &128.5.49.169.67&4355&3136& & &128.27.31.59.103&5487&6592 \cr
\noalign{\hrule}
}%
}
$$
\eject
\vglue -23 pt
\noindent\hskip 1 in\hbox to 6.5 in{\ 3277 -- 3312 \hfill\fbd 120158125 -- 122735041\frb}
\vskip -9 pt
$$
\vbox{
\nointerlineskip
\halign{\strut
    \vrule \ \ \hfil \frb #\ 
   &\vrule \hfil \ \ \fbb #\frb\ 
   &\vrule \hfil \ \ \frb #\ \hfil
   &\vrule \hfil \ \ \frb #\ 
   &\vrule \hfil \ \ \frb #\ \ \vrule \hskip 2 pt
   &\vrule \ \ \hfil \frb #\ 
   &\vrule \hfil \ \ \fbb #\frb\ 
   &\vrule \hfil \ \ \frb #\ \hfil
   &\vrule \hfil \ \ \frb #\ 
   &\vrule \hfil \ \ \frb #\ \vrule \cr%
\noalign{\hrule}
 & &625.17.43.263&10967&342& & &9.11.19.53.1217&2117&1534 \cr
3277&120158125&4.9.11.19.997&497&500&3295&121326381&4.3.13.19.29.59.73&737&650 \cr
 & &32.3.125.7.11.19.71&2343&2128& & &16.25.11.169.59.67&11323&11800 \cr
\noalign{\hrule}
 & &81.5.47.59.107&77&30& & &81.7.353.607&1195&1276 \cr
3278&120167955&4.243.25.7.11.59&8131&6206&3296&121491657&8.5.11.29.239.607&127&6804 \cr
 & &16.29.47.107.173&173&232& & &64.243.5.7.127&127&480 \cr
\noalign{\hrule}
 & &3.11.13.23.73.167&589&1090& & &9.5.11.29.43.197&23&238 \cr
3279&120288597&4.5.11.13.19.31.109&11189&11592&3297&121601205&4.7.11.17.23.197&1075&1092 \cr
 & &64.9.5.7.23.67.167&469&480& & &32.3.25.49.13.23.43&1127&1040 \cr
\noalign{\hrule}
 & &3.103.263.1481&2353&2090& & &9.121.19.5881&2891&2990 \cr
3280&120356427&4.5.11.13.19.103.181&1431&526&3298&121683771&4.5.49.11.13.19.23.59&411&26 \cr
 & &16.27.11.13.53.263&689&792& & &16.3.7.169.59.137&9971&7672 \cr
\noalign{\hrule}
 & &27.25.11.13.29.43&497&62& & &25.11.23.71.271&403&378 \cr
3281&120366675&4.9.5.7.11.31.71&493&592&3299&121699325&4.27.7.13.23.31.271&4925&3476 \cr
 & &128.17.29.37.71&1207&2368& & &32.3.25.11.13.79.197&3081&3152 \cr
\noalign{\hrule}
 & &9.5.11.19.23.557&301&136& & &3.5.7.11.31.41.83&129&212 \cr
3282&120487455&16.3.7.17.43.557&253&304&3300&121844415&8.9.5.7.41.43.53&121&166 \cr
 & &512.7.11.19.23.43&301&256& & &32.121.43.53.83&583&688 \cr
\noalign{\hrule}
 & &3.5.17.31.79.193&105&88& & &27.7.13.83.599&11737&3950 \cr
3283&120527535&16.9.25.7.11.31.79&9809&7334&3301&122154669&4.25.121.79.97&89&186 \cr
 & &64.17.19.193.577&577&608& & &16.3.11.31.79.89&2759&6952 \cr
\noalign{\hrule}
 & &5.11.361.59.103&303&818& & &5.11.13.17.19.529&11289&2296 \cr
3284&120658835&4.3.11.19.101.409&309&100&3302&122169905&16.3.7.41.53.71&221&150 \cr
 & &32.9.25.101.103&505&144& & &64.9.25.13.17.41&369&160 \cr
\noalign{\hrule}
 & &9.5.49.229.239&473&244& & &9.5.11.53.59.79&149&116 \cr
3285&120681855&8.3.5.49.11.43.61&703&32&3303&122281335&8.3.29.59.79.149&4151&8812 \cr
 & &512.19.37.43&703&11008& & &64.7.593.2203&15421&18976 \cr
\noalign{\hrule}
 & &5.7.17.43.53.89&2841&986& & &27.25.7.11.13.181&769&860 \cr
3286&120684445&4.3.289.29.947&3717&4664&3304&122297175&8.3.125.11.43.769&49&424 \cr
 & &64.27.7.11.53.59&649&864& & &128.49.53.769&5383&3392 \cr
\noalign{\hrule}
 & &3.25.23.31.37.61&813&38& & &9.5.71.101.379&20737&17542 \cr
3287&120693075&4.9.19.61.271&2849&2300&3305&122301405&4.49.89.179.233&18673&23034 \cr
 & &32.25.7.11.23.37&77&16& & &16.3.11.71.263.349&2893&2792 \cr
\noalign{\hrule}
 & &9.11.13.37.43.59&469&62& & &3.25.19.23.37.101&187&288 \cr
3288&120809403&4.7.13.31.43.67&295&264&3306&122480175&64.27.11.17.23.37&695&304 \cr
 & &64.3.5.7.11.59.67&335&224& & &2048.5.11.19.139&1529&1024 \cr
\noalign{\hrule}
 & &9.11.289.41.103&19097&18910& & &11.13.19.23.37.53&485&522 \cr
3289&120824253&4.5.169.17.31.61.113&15&1906&3307&122544851&4.9.5.11.13.23.29.97&1273&5420 \cr
 & &16.3.25.169.953&4225&7624& & &32.3.25.19.67.271&5025&4336 \cr
\noalign{\hrule}
 & &243.5.83.1201&1223&22& & &3.5.11.13.19.31.97&97&58 \cr
3290&121114845&4.81.11.1223&1057&166&3308&122550285&4.11.19.29.9409&7735&1674 \cr
 & &16.7.83.151&7&1208& & &16.27.5.7.13.17.31&153&56 \cr
\noalign{\hrule}
 & &27.7.11.167.349&59&290& & &27.11.19.37.587&139&158 \cr
3291&121170357&4.9.5.29.59.167&349&182&3309&122560317&4.37.79.139.587&5369&16350 \cr
 & &16.5.7.13.29.349&65&232& & &16.3.25.7.13.59.109&9919&11800 \cr
\noalign{\hrule}
 & &9.25.47.73.157&697&478& & &27.5.7.121.29.37&8471&7624 \cr
3292&121200075&4.3.17.41.157.239&77&80&3310&122692185&16.9.43.197.953&1363&410 \cr
 & &128.5.7.11.17.41.239&28441&28864& & &64.5.29.41.43.47&1763&1504 \cr
\noalign{\hrule}
 & &9.25.841.641&533&308& & &5.11.169.43.307&269&204 \cr
3293&121293225&8.7.11.13.41.641&177&464&3311&122703295&8.3.13.17.269.307&43&264 \cr
 & &256.3.11.13.29.59&767&1408& & &128.9.11.43.269&269&576 \cr
\noalign{\hrule}
 & &7.17.37.59.467&1325&858& & &11.13.19.199.227&981&1208 \cr
3294&121315859&4.3.25.7.11.13.17.53&1073&852&3312&122735041&16.9.13.19.109.151&1909&3980 \cr
 & &32.9.29.37.53.71&3763&4176& & &128.3.5.23.83.199&1245&1472 \cr
\noalign{\hrule}
}%
}
$$
\eject
\vglue -23 pt
\noindent\hskip 1 in\hbox to 6.5 in{\ 3313 -- 3348 \hfill\fbd 122849325 -- 126196571\frb}
\vskip -9 pt
$$
\vbox{
\nointerlineskip
\halign{\strut
    \vrule \ \ \hfil \frb #\ 
   &\vrule \hfil \ \ \fbb #\frb\ 
   &\vrule \hfil \ \ \frb #\ \hfil
   &\vrule \hfil \ \ \frb #\ 
   &\vrule \hfil \ \ \frb #\ \ \vrule \hskip 2 pt
   &\vrule \ \ \hfil \frb #\ 
   &\vrule \hfil \ \ \fbb #\frb\ 
   &\vrule \hfil \ \ \frb #\ \hfil
   &\vrule \hfil \ \ \frb #\ 
   &\vrule \hfil \ \ \frb #\ \vrule \cr%
\noalign{\hrule}
 & &27.25.23.41.193&15143&10318& & &5.7.11.23.59.239&2869&1674 \cr
3313&122849325&4.7.11.19.67.797&1035&238&3331&124864355&4.27.19.23.31.151&7453&6094 \cr
 & &16.9.5.49.11.17.23&187&392& & &16.3.11.29.257.277&7453&6648 \cr
\noalign{\hrule}
 & &5.7.11.19.53.317&37&18& & &3.5.49.43.59.67&1195&1342 \cr
3314&122899315&4.9.7.37.53.317&129&2090&3332&124934565&4.25.11.61.67.239&5031&944 \cr
 & &16.27.5.11.19.43&27&344& & &128.9.11.13.43.59&143&192 \cr
\noalign{\hrule}
 & &3.121.19.71.251&767&14& & &9.11.31.83.491&385&106 \cr
3315&122911437&4.7.11.13.19.59&315&334&3333&125070957&4.5.7.121.53.83&2139&2260 \cr
 & &16.9.5.49.13.167&1911&6680& & &32.3.25.7.23.31.113&2825&2576 \cr
\noalign{\hrule}
 & &9.37.191.1933&1253&680& & &41.12769.239&1485&11284 \cr
3316&122944599&16.3.5.7.17.37.179&1133&1910&3334&125123431&8.27.5.7.11.13.31&113&82 \cr
 & &64.25.11.103.191&1133&800& & &32.9.7.11.41.113&63&176 \cr
\noalign{\hrule}
 & &3.5.7.17.361.191&1443&1804& & &25.11.13.17.29.71&27&2 \cr
3317&123077535&8.9.5.7.11.13.37.41&361&46&3335&125135725&4.27.11.13.17.71&3625&3404 \cr
 & &32.13.361.23.41&299&656& & &32.3.125.23.29.37&555&368 \cr
\noalign{\hrule}
 & &9.5.13.29.53.137&913&868& & &7.11.13.31.37.109&1349&11268 \cr
3318&123182865&8.7.11.29.31.53.83&2055&518&3336&125148023&8.9.19.71.313&505&434 \cr
 & &32.3.5.49.11.37.137&539&592& & &32.3.5.7.19.31.101&505&912 \cr
\noalign{\hrule}
 & &9.5.13.17.79.157&1403&638& & &27.5.13.23.29.107&2387&502 \cr
3319&123347835&4.11.23.29.61.79&845&924&3337&125252595&4.7.11.23.31.251&1235&522 \cr
 & &32.3.5.7.121.169.23&2093&1936& & &16.9.5.11.13.19.29&19&88 \cr
\noalign{\hrule}
 & &9.7.11.41.43.101&125&986& & &9.5.41.149.457&9145&9592 \cr
3320&123397659&4.3.125.17.29.43&2933&808&3338&125631585&16.3.25.11.31.59.109&301&26 \cr
 & &64.7.101.419&419&32& & &64.7.13.31.43.59&17329&13216 \cr
\noalign{\hrule}
 & &27.7.11.23.29.89&2125&4706& & &3.625.7.11.13.67&197&428 \cr
3321&123415677&4.125.7.13.17.181&3&178&3339&125750625&8.13.67.107.197&7295&5904 \cr
 & &16.3.5.13.17.89&13&680& & &256.9.5.41.1459&4377&5248 \cr
\noalign{\hrule}
 & &5.11.23.251.389&3081&2692& & &3.5.11.13.89.659&21373&21462 \cr
3322&123513335&8.3.5.11.13.79.673&21&694&3340&125806395&4.9.49.121.29.67.73&3293&2636 \cr
 & &32.9.7.79.347&2429&11376& & &32.29.37.67.89.659&1073&1072 \cr
\noalign{\hrule}
 & &3.7.11.17.163.193&899&410& & &3.7.121.13.37.103&69&190 \cr
3323&123539493&4.5.29.31.41.193&117&76&3341&125888763&4.9.5.13.19.23.103&1369&854 \cr
 & &32.9.5.13.19.29.31&6045&8816& & &16.7.23.1369.61&851&488 \cr
\noalign{\hrule}
 & &9.5.11.37.43.157&559&226& & &81.7.13.19.29.31&1357&1154 \cr
3324&123644565&4.11.13.1849.113&9435&10904&3342&125904051&4.13.19.23.59.577&165&602 \cr
 & &64.3.5.17.29.37.47&799&928& & &16.3.5.7.11.43.577&2885&3784 \cr
\noalign{\hrule}
 & &5.11.83.103.263&173&90& & &9.5.49.11.29.179&1883&988 \cr
3325&123661285&4.9.25.11.103.173&2911&4814&3343&125907705&8.343.13.19.269&4785&326 \cr
 & &16.3.29.41.71.83&1189&1704& & &32.3.5.11.29.163&163&16 \cr
\noalign{\hrule}
 & &9.7.11.23.43.181&325&578& & &25.7.13.23.29.83&333&242 \cr
3326&124053237&4.3.25.13.289.181&121&784&3344&125946275&4.9.121.29.37.83&85&2 \cr
 & &128.5.49.121.17&1309&320& & &16.3.5.121.17.37&363&5032 \cr
\noalign{\hrule}
 & &13.23.31.59.227&12339&5302& & &9.125.19.71.83&671&754 \cr
3327&124139717&4.27.11.241.457&565&806&3345&125962875&4.3.5.11.13.29.61.71&5377&17272 \cr
 & &16.9.5.11.13.31.113&565&792& & &64.17.19.127.283&4811&4064 \cr
\noalign{\hrule}
 & &9.5.7.11.13.31.89&109&46& & &9.7.13.137.1123&3725&4136 \cr
3328&124279155&4.11.13.23.89.109&21&1178&3346&126003969&16.3.25.11.13.47.149&3989&3014 \cr
 & &16.3.7.19.23.31&437&8& & &64.121.137.3989&3989&3872 \cr
\noalign{\hrule}
 & &7.11.61.71.373&41&30& & &5.169.17.67.131&11571&2794 \cr
3329&124390651&4.3.5.7.41.61.373&55&2556&3347&126081605&4.3.7.11.19.29.127&169&150 \cr
 & &32.27.25.11.71&27&400& & &16.9.25.7.169.127&635&504 \cr
\noalign{\hrule}
 & &27.25.31.59.101&671&166& & &89.10201.139&1085&11286 \cr
3330&124692075&4.5.11.59.61.83&909&4154&3348&126196571&4.27.5.7.11.19.31&101&178 \cr
 & &16.9.31.67.101&67&8& & &16.3.5.19.89.101&19&120 \cr
\noalign{\hrule}
}%
}
$$
\eject
\vglue -23 pt
\noindent\hskip 1 in\hbox to 6.5 in{\ 3349 -- 3384 \hfill\fbd 126220325 -- 129055455\frb}
\vskip -9 pt
$$
\vbox{
\nointerlineskip
\halign{\strut
    \vrule \ \ \hfil \frb #\ 
   &\vrule \hfil \ \ \fbb #\frb\ 
   &\vrule \hfil \ \ \frb #\ \hfil
   &\vrule \hfil \ \ \frb #\ 
   &\vrule \hfil \ \ \frb #\ \ \vrule \hskip 2 pt
   &\vrule \ \ \hfil \frb #\ 
   &\vrule \hfil \ \ \fbb #\frb\ 
   &\vrule \hfil \ \ \frb #\ \hfil
   &\vrule \hfil \ \ \frb #\ 
   &\vrule \hfil \ \ \frb #\ \vrule \cr%
\noalign{\hrule}
 & &25.49.11.17.19.29&153&398& & &3.7.19.23.53.263&2431&590 \cr
3349&126220325&4.9.5.11.289.199&461&406&3367&127918203&4.5.11.13.17.23.59&41&18 \cr
 & &16.3.7.29.199.461&1383&1592& & &16.9.5.11.13.17.41&6765&1768 \cr
\noalign{\hrule}
 & &27.13.31.41.283&361&11242& & &5.343.17.23.191&891&1846 \cr
3350&126252243&4.7.11.361.73&697&690&3368&128077915&4.81.49.11.13.71&191&3670 \cr
 & &16.3.5.11.17.19.23.41&1955&1672& & &16.3.5.191.367&367&24 \cr
\noalign{\hrule}
 & &49.361.37.193&1087&726& & &5.49.13.19.29.73&795&592 \cr
3351&126317149&4.3.121.193.1087&1605&518&3369&128110255&32.3.25.7.13.37.53&33&292 \cr
 & &16.9.5.7.11.37.107&495&856& & &256.9.11.53.73&477&1408 \cr
\noalign{\hrule}
 & &81.5.11.13.37.59&1189&994& & &5.7.11.23.43.337&813&2498 \cr
3352&126428445&4.27.7.11.29.41.71&163&944&3370&128317805&4.3.23.271.1249&3741&2492 \cr
 & &128.7.29.59.163&1141&1856& & &32.9.7.29.43.89&801&464 \cr
\noalign{\hrule}
 & &3.5.7.19.61.1039&391&524& & &25.11.23.53.383&4377&4432 \cr
3353&126441105&8.17.23.131.1039&715&324&3371&128391175&32.3.5.53.277.1459&271&6 \cr
 & &64.81.5.11.13.131&3537&4576& & &128.9.271.1459&13131&17344 \cr
\noalign{\hrule}
 & &3.5.17.43.83.139&61&78& & &11.19.53.67.173&325&258 \cr
3354&126503205&4.9.5.13.43.61.83&3179&556&3372&128393507&4.3.25.13.19.43.173&39&134 \cr
 & &32.11.13.289.139&143&272& & &16.9.5.169.43.67&1935&1352 \cr
\noalign{\hrule}
 & &3.7.11.47.89.131&1633&2550& & &9.7.121.23.733&61&38 \cr
3355&126581763&4.9.25.11.17.23.71&163&262&3373&128516157&4.7.11.19.61.733&4615&9312 \cr
 & &16.23.71.131.163&1633&1304& & &256.3.5.13.71.97&6305&9088 \cr
\noalign{\hrule}
 & &5.11.17.43.47.67&441&1580& & &9.5.11.169.29.53&893&628 \cr
3356&126605545&8.9.25.49.11.79&1081&844&3374&128577735&8.11.19.29.47.157&169&150 \cr
 & &64.3.7.23.47.211&1477&2208& & &32.3.25.169.47.157&785&752 \cr
\noalign{\hrule}
 & &3.121.41.47.181&13975&8288& & &5.11.13.17.71.149&86673&85928 \cr
3357&126609681&64.25.7.13.37.43&181&144&3375&128587745&16.3.23.167.173.467&38475&42316 \cr
 & &2048.9.7.43.181&903&1024& & &128.243.25.19.71.149&1215&1216 \cr
\noalign{\hrule}
 & &9.5.19.89.1667&1261&406& & &7.11.23.59.1231&291&940 \cr
3358&126850365&4.7.13.29.89.97&319&942&3376&128625959&8.3.5.7.23.47.97&129&32 \cr
 & &16.3.11.841.157&1727&6728& & &512.9.5.43.47&2021&11520 \cr
\noalign{\hrule}
 & &9.11.13.151.653&343&310& & &9.5.49.19.37.83&13&428 \cr
3359&126902061&4.3.5.343.13.31.151&1133&76&3377&128659545&8.13.19.37.107&231&250 \cr
 & &32.5.49.11.19.103&1957&3920& & &32.3.125.7.11.107&1177&400 \cr
\noalign{\hrule}
 & &3.5.13.19.43.797&1163&1228& & &25.7.41.79.227&193&1782 \cr
3360&126974055&8.19.43.307.1163&173&990&3378&128669275&4.81.11.41.193&35&158 \cr
 & &32.9.5.11.173.307&5709&4912& & &16.27.5.7.11.79&11&216 \cr
\noalign{\hrule}
 & &25.13.23.29.587&351&374& & &3.17.337.7487&3769&3718 \cr
3361&127246925&4.27.11.169.17.587&1143&1730&3379&128679069&4.11.169.337.3769&4075&306 \cr
 & &16.243.5.11.127.173&21971&21384& & &16.9.25.11.13.17.163&3575&3912 \cr
\noalign{\hrule}
 & &9.11.53.79.307&9503&14750& & &7.17.23.31.37.41&2343&3614 \cr
3362&127255491&4.125.13.17.43.59&921&154&3380&128712899&4.3.11.13.17.71.139&555&368 \cr
 & &16.3.5.7.11.17.307&119&40& & &128.9.5.23.37.139&695&576 \cr
\noalign{\hrule}
 & &27.5.7.11.13.23.41&599&116& & &3.5.7.11.17.79.83&3067&3016 \cr
3363&127432305&8.9.29.41.599&115&484&3381&128746695&16.5.13.29.83.3067&19&396 \cr
 & &64.5.121.23.29&319&32& & &128.9.11.19.3067&3067&3648 \cr
\noalign{\hrule}
 & &3.5.7.11.163.677&11591&2626& & &19.311.21799&10909&10890 \cr
3364&127455405&4.13.67.101.173&137&36&3382&128810291&4.9.5.121.311.10909&323&10586 \cr
 & &32.9.13.67.137&871&6576& & &16.3.5.11.17.19.67.79&11055&10744 \cr
\noalign{\hrule}
 & &9.41.277.1249&323&46& & &27.5.11.17.19.269&7199&8134 \cr
3365&127664037&4.17.19.23.1249&429&820&3383&129027195&4.9.49.23.83.313&845&1972 \cr
 & &32.3.5.11.13.19.41&209&1040& & &32.5.169.17.29.83&2407&2704 \cr
\noalign{\hrule}
 & &9.7.13.17.67.137&869&2008& & &9.5.83.109.317&319&428 \cr
3366&127799217&16.3.11.13.79.251&137&890&3384&129055455&8.5.11.29.107.317&109&426 \cr
 & &64.5.11.89.137&979&160& & &32.3.11.29.71.109&781&464 \cr
\noalign{\hrule}
}%
}
$$
\eject
\vglue -23 pt
\noindent\hskip 1 in\hbox to 6.5 in{\ 3385 -- 3420 \hfill\fbd 129058475 -- 131990859\frb}
\vskip -9 pt
$$
\vbox{
\nointerlineskip
\halign{\strut
    \vrule \ \ \hfil \frb #\ 
   &\vrule \hfil \ \ \fbb #\frb\ 
   &\vrule \hfil \ \ \frb #\ \hfil
   &\vrule \hfil \ \ \frb #\ 
   &\vrule \hfil \ \ \frb #\ \ \vrule \hskip 2 pt
   &\vrule \ \ \hfil \frb #\ 
   &\vrule \hfil \ \ \fbb #\frb\ 
   &\vrule \hfil \ \ \frb #\ \hfil
   &\vrule \hfil \ \ \frb #\ 
   &\vrule \hfil \ \ \frb #\ \vrule \cr%
\noalign{\hrule}
 & &25.7.13.17.47.71&393&218& & &3.169.29.53.167&359&330 \cr
3385&129058475&4.3.17.71.109.131&5577&3724&3403&130136253&4.9.5.11.13.167.359&2701&530 \cr
 & &32.9.49.11.169.19&1881&1456& & &16.25.11.37.53.73&2701&2200 \cr
\noalign{\hrule}
 & &27.7.11.13.17.281&1919&610& & &3.5.19.23.31.641&539&102 \cr
3386&129107979&4.3.5.13.19.61.101&17&22&3404&130254405&4.9.5.49.11.17.31&437&158 \cr
 & &16.11.17.19.61.101&1159&808& & &16.7.11.19.23.79&553&88 \cr
\noalign{\hrule}
 & &9.5.97.107.277&299&22& & &27.49.11.169.53&103&194 \cr
3387&129374235&4.3.5.11.13.23.97&277&208&3405&130351221&4.7.13.53.97.103&855&484 \cr
 & &128.11.169.277&169&704& & &32.9.5.121.19.97&1045&1552 \cr
\noalign{\hrule}
 & &9.7.23.31.43.67&64931&65920& & &9.25.13.17.43.61&143&82 \cr
3388&129411639&256.5.29.103.2239&2613&374&3406&130428675&4.11.169.17.41.43&555&2318 \cr
 & &1024.3.5.11.13.17.67&2431&2560& & &16.3.5.11.19.37.61&407&152 \cr
\noalign{\hrule}
 & &9.23.41.101.151&1391&2750& & &5.7.11.13.67.389&163&306 \cr
3389&129435237&4.125.11.13.23.107&301&876&3407&130445315&4.9.5.17.163.389&67&322 \cr
 & &32.3.5.7.13.43.73&3139&7280& & &16.3.7.23.67.163&163&552 \cr
\noalign{\hrule}
 & &263.677.727&25&702& & &81.5.29.41.271&1199&1240 \cr
3390&129443077&4.27.25.13.263&677&638&3408&130498695&16.9.25.11.29.31.109&1687&788 \cr
 & &16.9.5.11.29.677&495&232& & &128.7.109.197.241&47477&48832 \cr
\noalign{\hrule}
 & &3.13.29.239.479&2285&11606& & &9.11.13.29.31.113&1075&2318 \cr
3391&129478011&4.5.7.457.829&473&3672&3409&130742469&4.25.19.31.43.61&429&904 \cr
 & &64.27.11.17.43&6579&352& & &64.3.11.13.61.113&61&32 \cr
\noalign{\hrule}
 & &3.343.19.37.179&1019&4420& & &27.125.11.13.271&301&3824 \cr
3392&129486273&8.5.7.13.17.1019&513&506&3410&130791375&32.9.7.43.239&925&1226 \cr
 & &32.27.5.11.13.17.19.23&6435&6256& & &128.25.37.613&613&2368 \cr
\noalign{\hrule}
 & &27.11.23.61.311&3527&3626& & &3.5.11.13.17.37.97&95&386 \cr
3393&129590901&4.3.49.37.61.3527&8959&15730&3411&130872885&4.25.11.17.19.193&4171&504 \cr
 & &16.5.7.121.13.289.31&18785&19096& & &64.9.7.43.97&43&672 \cr
\noalign{\hrule}
 & &9.49.17.59.293&817&1820& & &81.5.7.11.13.17.19&139&4 \cr
3394&129600639&8.5.343.13.19.43&451&108&3412&130945815&8.3.7.17.19.139&3611&3478 \cr
 & &64.27.5.11.19.41&1045&3936& & &32.23.37.47.157&7379&13616 \cr
\noalign{\hrule}
 & &3.25.11.31.37.137&139&546& & &5.11.13.19.31.311&1633&1788 \cr
3395&129639675&4.9.5.7.13.31.139&703&548&3413&130972985&8.3.13.19.23.71.149&1643&294 \cr
 & &32.7.13.19.37.137&133&208& & &32.9.49.23.31.53&2597&3312 \cr
\noalign{\hrule}
 & &25.11.13.19.23.83&289&36& & &7.11.17.31.53.61&97&90 \cr
3396&129668825&8.9.289.19.83&443&526&3414&131191907&4.9.5.31.53.61.97&5423&494 \cr
 & &32.3.17.263.443&13413&7088& & &16.3.5.11.13.17.19.29&1131&760 \cr
\noalign{\hrule}
 & &9.25.11.19.31.89&757&668& & &9.25.7.239.349&11891&3166 \cr
3397&129741975&8.3.11.31.167.757&133&890&3415&131372325&4.11.23.47.1583&533&1050 \cr
 & &32.5.7.19.89.167&167&112& & &16.3.25.7.13.23.41&533&184 \cr
\noalign{\hrule}
 & &7.19.89.97.113&789&902& & &27.13.37.67.151&25&176 \cr
3398&129745357&4.3.7.11.41.97.263&565&114&3416&131389479&32.9.25.11.13.37&151&34 \cr
 & &16.9.5.19.113.263&263&360& & &128.5.11.17.151&935&64 \cr
\noalign{\hrule}
 & &27.5.7.17.41.197&179&26& & &9.5.7.11.17.23.97&589&104 \cr
3399&129757005&4.3.7.13.179.197&341&250&3417&131417055&16.13.17.19.23.31&329&108 \cr
 & &16.125.11.31.179&5549&2200& & &128.27.7.31.47&1457&192 \cr
\noalign{\hrule}
 & &27.11.13.19.29.61&7675&13736& & &9.5.7.11.17.23.97&509&800 \cr
3400&129772071&16.25.17.101.307&261&244&3418&131417055&64.3.125.23.509&97&28 \cr
 & &128.9.5.29.61.307&307&320& & &512.7.97.509&509&256 \cr
\noalign{\hrule}
 & &9.5.11.13.17.29.41&247&1436& & &27.5.49.11.13.139&16813&17242 \cr
3401&130070655&8.5.169.19.359&87&82&3419&131486355&4.9.17.23.37.43.233&13&220 \cr
 & &32.3.19.29.41.359&359&304& & &32.5.11.13.17.37.43&731&592 \cr
\noalign{\hrule}
 & &9.7.11.19.41.241&4209&5672& & &9.343.11.169.23&491&1030 \cr
3402&130103127&16.27.23.61.709&665&44&3420&131990859&4.5.7.23.103.491&351&454 \cr
 & &128.5.7.11.19.61&305&64& & &16.27.13.227.491&1473&1816 \cr
\noalign{\hrule}
}%
}
$$
\eject
\vglue -23 pt
\noindent\hskip 1 in\hbox to 6.5 in{\ 3421 -- 3456 \hfill\fbd 132011451 -- 135438303\frb}
\vskip -9 pt
$$
\vbox{
\nointerlineskip
\halign{\strut
    \vrule \ \ \hfil \frb #\ 
   &\vrule \hfil \ \ \fbb #\frb\ 
   &\vrule \hfil \ \ \frb #\ \hfil
   &\vrule \hfil \ \ \frb #\ 
   &\vrule \hfil \ \ \frb #\ \ \vrule \hskip 2 pt
   &\vrule \ \ \hfil \frb #\ 
   &\vrule \hfil \ \ \fbb #\frb\ 
   &\vrule \hfil \ \ \frb #\ \hfil
   &\vrule \hfil \ \ \frb #\ 
   &\vrule \hfil \ \ \frb #\ \vrule \cr%
\noalign{\hrule}
 & &243.11.13.29.131&3589&3458& & &3.125.7.11.4639&25327&25702 \cr
3421&132011451&4.7.11.169.19.37.97&1125&58&3439&133951125&4.7.19.31.43.71.181&8613&830 \cr
 & &16.9.125.19.29.37&703&1000& & &16.27.5.11.29.31.83&2407&2232 \cr
\noalign{\hrule}
 & &5.19.31.193.233&17171&12744& & &9.25.11.113.479&2303&2782 \cr
3422&132433705&16.27.7.11.59.223&193&1754&3440&133964325&4.5.49.11.13.47.107&113&498 \cr
 & &64.9.193.877&877&288& & &16.3.7.83.107.113&581&856 \cr
\noalign{\hrule}
 & &27.5.7.121.19.61&1619&460& & &3.25.961.1861&2093&232 \cr
3423&132525855&8.25.11.23.1619&139&114&3441&134131575&16.7.13.23.29.31&495&404 \cr
 & &32.3.19.139.1619&1619&2224& & &128.9.5.11.23.101&3333&1472 \cr
\noalign{\hrule}
 & &3.25.11.13.47.263&407&382& & &27.19.23.83.137&17&154 \cr
3424&132571725&4.121.13.37.47.191&83&1656&3442&134166429&4.3.7.11.17.23.83&1235&674 \cr
 & &64.9.23.83.191&4393&7968& & &16.5.7.13.19.337&455&2696 \cr
\noalign{\hrule}
 & &27.5.121.23.353&65&56& & &27.11.13.19.31.59&5&346 \cr
3425&132623865&16.3.25.7.13.23.353&817&242&3443&134173611&4.5.19.59.173&5151&5056 \cr
 & &64.7.121.13.19.43&1729&1376& & &512.3.17.79.101&7979&4352 \cr
\noalign{\hrule}
 & &9.125.7.19.887&553&572& & &5.13.61.97.349&921&824 \cr
3426&132717375&8.49.11.13.79.887&125&762&3444&134227145&16.3.13.61.103.307&5137&1146 \cr
 & &32.3.125.11.79.127&1397&1264& & &64.9.11.191.467&18909&14944 \cr
\noalign{\hrule}
 & &125.11.13.17.19.23&267&124& & &27.5.7.11.37.349&323&26 \cr
3427&132793375&8.3.125.19.31.89&2783&1092&3445&134230635&4.5.7.13.17.19.37&657&638 \cr
 & &64.9.7.121.13.23&99&224& & &16.9.11.13.17.29.73&2117&1768 \cr
\noalign{\hrule}
 & &9.11.83.103.157&203&46& & &9.11.19.23.29.107&611&2492 \cr
3428&132877107&4.3.7.11.23.29.103&5&314&3446&134245089&8.7.13.23.47.89&3509&4590 \cr
 & &16.5.7.23.157&7&920& & &32.27.5.121.17.29&165&272 \cr
\noalign{\hrule}
 & &9.5.7.11.19.43.47&25&448& & &9.49.11.19.31.47&965&916 \cr
3429&133052535&128.125.49.19&403&528&3447&134290233&8.5.31.47.193.229&8085&986 \cr
 & &4096.3.11.13.31&403&2048& & &32.3.25.49.11.17.29&493&400 \cr
\noalign{\hrule}
 & &3.5.19.421.1109&3445&4554& & &81.125.11.17.71&2021&104 \cr
3430&133063365&4.27.25.11.13.23.53&3973&11552&3448&134429625&16.3.11.13.43.47&1&142 \cr
 & &256.361.29.137&2603&3712& & &64.43.71&1&1376 \cr
\noalign{\hrule}
 & &9.7.121.13.17.79&199&1858& & &3.25.49.23.37.43&793&58 \cr
3431&133089957&4.3.13.199.929&4345&3416&3449&134479275&4.5.13.29.43.61&207&352 \cr
 & &64.5.7.11.61.79&61&160& & &256.9.11.23.61&183&1408 \cr
\noalign{\hrule}
 & &3.5.13.19.103.349&61&42& & &3.5.13.53.83.157&469&220 \cr
3432&133183635&4.9.5.7.13.61.349&3553&412&3450&134675385&8.25.7.11.67.157&9&166 \cr
 & &32.7.11.17.19.103&77&272& & &32.9.11.67.83&737&48 \cr
\noalign{\hrule}
 & &5.17.19.23.37.97&801&1430& & &5.23.31.103.367&273&242 \cr
3433&133313405&4.9.25.11.13.19.89&111&136&3451&134760565&4.3.7.121.13.23.367&57&310 \cr
 & &64.27.11.17.37.89&979&864& & &16.9.5.7.11.13.19.31&1001&1368 \cr
\noalign{\hrule}
 & &31.41.61.1721&85&1806& & &9.25.29.139.149&283&134 \cr
3434&133430851&4.3.5.7.17.41.43&165&122&3452&135139275&4.3.25.29.67.283&11473&13148 \cr
 & &16.9.25.11.17.61&225&1496& & &32.7.11.19.149.173&1903&2128 \cr
\noalign{\hrule}
 & &27.11.13.17.19.107&1&10& & &9.7.11.41.67.71&337&400 \cr
3435&133440021&4.3.5.13.17.19.107&1925&2248&3453&135160641&32.25.41.71.337&67&138 \cr
 & &64.125.7.11.281&1967&4000& & &128.3.5.23.67.337&1685&1472 \cr
\noalign{\hrule}
 & &9.125.7.11.23.67&1073&698& & &17.31.101.2543&21565&21666 \cr
3436&133489125&4.3.29.37.67.349&1495&448&3454&135356261&4.3.5.19.23.31.157.227&3939&374 \cr
 & &512.5.7.13.23.37&481&256& & &16.9.11.13.17.101.157&1287&1256 \cr
\noalign{\hrule}
 & &25.11.53.89.103&2507&2952& & &5.7.11.13.17.37.43&4189&12366 \cr
3437&133609025&16.9.5.11.23.41.109&2867&1602&3455&135370235&4.27.59.71.229&25&34 \cr
 & &64.81.47.61.89&2867&2592& & &16.3.25.17.71.229&1145&1704 \cr
\noalign{\hrule}
 & &243.13.19.23.97&9185&14774& & &3.49.11.13.17.379&4427&500 \cr
3438&133906851&4.5.11.83.89.167&39&874&3456&135438303&8.125.7.19.233&1053&578 \cr
 & &16.3.13.19.23.89&89&8& & &32.81.5.13.289&459&80 \cr
\noalign{\hrule}
}%
}
$$
\eject
\vglue -23 pt
\noindent\hskip 1 in\hbox to 6.5 in{\ 3457 -- 3492 \hfill\fbd 135462391 -- 139145193\frb}
\vskip -9 pt
$$
\vbox{
\nointerlineskip
\halign{\strut
    \vrule \ \ \hfil \frb #\ 
   &\vrule \hfil \ \ \fbb #\frb\ 
   &\vrule \hfil \ \ \frb #\ \hfil
   &\vrule \hfil \ \ \frb #\ 
   &\vrule \hfil \ \ \frb #\ \ \vrule \hskip 2 pt
   &\vrule \ \ \hfil \frb #\ 
   &\vrule \hfil \ \ \fbb #\frb\ 
   &\vrule \hfil \ \ \frb #\ \hfil
   &\vrule \hfil \ \ \frb #\ 
   &\vrule \hfil \ \ \frb #\ \vrule \cr%
\noalign{\hrule}
 & &71.83.127.181&99&28& & &7.13.17.19.31.151&2545&5112 \cr
3457&135462391&8.9.7.11.83.181&215&34&3475&137588633&16.9.5.7.71.509&3473&3982 \cr
 & &32.3.5.7.11.17.43&1155&11696& & &64.3.11.23.151.181&1991&2208 \cr
\noalign{\hrule}
 & &9.5.13.29.61.131&109&22& & &5.11.13.29.61.109&2043&2104 \cr
3458&135567315&4.3.5.11.13.61.109&8251&7336&3476&137867015&16.9.5.109.227.263&667&122 \cr
 & &64.7.37.131.223&1561&1184& & &64.3.23.29.61.227&681&736 \cr
\noalign{\hrule}
 & &81.7.43.67.83&275&194& & &5.7.11.29.53.233&27&292 \cr
3459&135582741&4.25.11.43.83.97&197&4368&3477&137876585&8.27.7.73.233&605&94 \cr
 & &128.3.5.7.13.197&2561&320& & &32.9.5.121.47&47&1584 \cr
\noalign{\hrule}
 & &9.5.11.169.1621&4313&12418& & &7.13.19.173.461&297&164 \cr
3460&135604755&4.7.19.227.887&557&330&3478&137892937&8.27.11.13.41.173&145&28 \cr
 & &16.3.5.7.11.19.557&557&1064& & &64.3.5.7.11.29.41&1353&4640 \cr
\noalign{\hrule}
 & &3.49.169.53.103&855&484& & &25.7.13.19.31.103&993&2332 \cr
3461&135617937&8.27.5.7.121.13.19&103&194&3479&138017425&8.3.11.31.53.331&987&656 \cr
 & &32.5.11.19.97.103&1045&1552& & &256.9.7.11.41.47&4059&6016 \cr
\noalign{\hrule}
 & &3.5.121.13.73.79&151&86& & &9.7.11.41.43.113&815&428 \cr
3462&136072365&4.121.43.73.151&7663&1170&3480&138058767&8.5.7.41.107.163&1469&5214 \cr
 & &16.9.5.13.79.97&97&24& & &32.3.11.13.79.113&79&208 \cr
\noalign{\hrule}
 & &27.25.49.13.317&5497&3278& & &27.7.11.181.367&425&4462 \cr
3463&136302075&4.7.11.23.149.239&1395&1234&3481&138101733&4.25.7.17.23.97&323&162 \cr
 & &16.9.5.31.149.617&4619&4936& & &16.81.5.289.19&289&2280 \cr
\noalign{\hrule}
 & &3.5.13.37.18899&13057&5842& & &27.5.337.3037&21229&24266 \cr
3464&136356285&4.11.23.127.1187&593&594&3482&138168315&4.11.13.23.71.1103&1013&90 \cr
 & &16.27.121.23.127.593&122751&122936& & &16.9.5.11.23.1013&1013&2024 \cr
\noalign{\hrule}
 & &9.17.19.61.769&869&100& & &5.41.43.61.257&231&26 \cr
3465&136364463&8.3.25.11.61.79&1079&896&3483&138192755&4.3.7.11.13.43.61&375&418 \cr
 & &2048.7.11.13.83&6391&13312& & &16.9.125.7.121.19&4275&6776 \cr
\noalign{\hrule}
 & &625.49.61.73&99&526& & &5.41.47.83.173&13299&6206 \cr
3466&136373125&4.9.7.11.73.263&261&250&3484&138348965&4.3.11.13.29.31.107&173&204 \cr
 & &16.81.125.29.263&2349&2104& & &32.9.11.17.107.173&1819&1584 \cr
\noalign{\hrule}
 & &29.59.199.401&40755&39044& & &3.25.7.11.289.83&29&386 \cr
3467&136536089&8.3.5.11.13.19.43.227&9&218&3485&138524925&4.5.11.17.29.193&2331&3266 \cr
 & &32.27.5.13.43.109&14715&8944& & &16.9.7.23.37.71&2627&552 \cr
\noalign{\hrule}
 & &27.25.13.79.197&551&434& & &9.5.121.13.19.103&1369&854 \cr
3468&136565325&4.3.5.7.19.29.31.79&143&1042&3486&138526245&4.7.121.1369.61&69&190 \cr
 & &16.7.11.13.19.521&3647&1672& & &16.3.5.19.23.37.61&851&488 \cr
\noalign{\hrule}
 & &3.625.7.11.13.73&1127&698& & &9.7.11.29.61.113&1865&1978 \cr
3469&137011875&4.25.343.23.349&8307&418&3487&138528621&4.5.11.23.29.43.373&31&10848 \cr
 & &16.9.11.13.19.71&57&568& & &256.3.5.31.113&31&640 \cr
\noalign{\hrule}
 & &25.13.19.53.419&12593&4632& & &9.49.29.37.293&11815&974 \cr
3470&137128225&16.3.49.193.257&125&132&3488&138645549&4.5.17.139.487&201&286 \cr
 & &128.9.125.7.11.193&6755&6336& & &16.3.11.13.67.139&1529&6968 \cr
\noalign{\hrule}
 & &125.13.29.41.71&513&2398& & &9.5.121.13.37.53&1873&2456 \cr
3471&137180875&4.27.25.11.19.109&923&1148&3489&138809385&16.5.11.307.1873&9379&7506 \cr
 & &32.3.7.11.13.41.71&77&48& & &64.27.83.113.139&11537&10848 \cr
\noalign{\hrule}
 & &5.7.17.23.79.127&319&234& & &25.11.293.1723&13139&5814 \cr
3472&137301605&4.9.11.13.23.29.127&41&340&3490&138830725&4.9.7.17.19.1877&1685&3946 \cr
 & &32.3.5.11.17.29.41&1353&464& & &16.3.5.337.1973&5919&2696 \cr
\noalign{\hrule}
 & &9.5.31.241.409&43&198& & &9.5.49.13.37.131&2425&2684 \cr
3473&137503755&4.81.11.43.409&241&650&3491&138939255&8.3.125.7.11.61.97&47&1328 \cr
 & &16.25.13.43.241&65&344& & &256.47.83.97&3901&12416 \cr
\noalign{\hrule}
 & &3.11.23.37.59.83&85&26& & &9.11.23.53.1153&2405&1252 \cr
3474&137522451&4.5.11.13.17.23.83&413&666&3492&139145193&8.3.5.11.13.37.313&2645&1424 \cr
 & &16.9.5.7.17.37.59&85&168& & &256.25.529.89&2047&3200 \cr
\noalign{\hrule}
}%
}
$$
\eject
\vglue -23 pt
\noindent\hskip 1 in\hbox to 6.5 in{\ 3493 -- 3528 \hfill\fbd 139152285 -- 143819949\frb}
\vskip -9 pt
$$
\vbox{
\nointerlineskip
\halign{\strut
    \vrule \ \ \hfil \frb #\ 
   &\vrule \hfil \ \ \fbb #\frb\ 
   &\vrule \hfil \ \ \frb #\ \hfil
   &\vrule \hfil \ \ \frb #\ 
   &\vrule \hfil \ \ \frb #\ \ \vrule \hskip 2 pt
   &\vrule \ \ \hfil \frb #\ 
   &\vrule \hfil \ \ \fbb #\frb\ 
   &\vrule \hfil \ \ \frb #\ \hfil
   &\vrule \hfil \ \ \frb #\ 
   &\vrule \hfil \ \ \frb #\ \vrule \cr%
\noalign{\hrule}
 & &9.5.61.163.311&581&886& & &27.25.7.121.13.19&41&206 \cr
3493&139152285&4.7.83.311.443&377&66&3511&141216075&4.9.5.7.11.41.103&377&8 \cr
 & &16.3.7.11.13.29.83&2639&7304& & &64.13.29.103&29&3296 \cr
\noalign{\hrule}
 & &27.5.31.79.421&803&382& & &3.25.61.89.347&2461&2114 \cr
3494&139188915&4.9.11.31.73.191&689&1030&3512&141289725&4.7.23.89.107.151&3025&6498 \cr
 & &16.5.13.53.73.103&5459&7592& & &16.9.25.7.121.361&2541&2888 \cr
\noalign{\hrule}
 & &7.11.13.23.73.83&405&106& & &81.5.41.67.127&13&28 \cr
3495&139496357&4.81.5.11.53.83&377&536&3513&141291945&8.27.7.13.67.127&4763&1334 \cr
 & &64.27.5.13.29.67&3915&2144& & &32.11.23.29.433&12557&4048 \cr
\noalign{\hrule}
 & &5.49.17.19.41.43&899&4554& & &9.5.43.103.709&1127&418 \cr
3496&139515005&4.9.7.11.23.29.31&325&326&3514&141307245&4.3.49.11.19.23.43&1175&244 \cr
 & &16.3.25.11.13.23.29.163&43355&43032& & &32.25.23.47.61&5405&976 \cr
\noalign{\hrule}
 & &3.5.11.13.17.43.89&241&26& & &9.5.59.139.383&1397&1258 \cr
3497&139551555&4.11.169.17.241&9&178&3515&141344235&4.11.17.37.127.383&1027&3186 \cr
 & &16.9.89.241&723&8& & &16.27.13.37.59.79&1027&888 \cr
\noalign{\hrule}
 & &19.22201.331&14245&7956& & &5.11.13.17.89.131&1477&36 \cr
3498&139622089&8.9.5.7.11.13.17.37&149&38&3516&141715145&8.9.5.7.13.211&1319&1424 \cr
 & &32.3.5.7.13.19.149&91&240& & &256.3.89.1319&1319&384 \cr
\noalign{\hrule}
 & &9.25.7.79.1123&9449&8326& & &27.5.7.13.83.139&779&194 \cr
3499&139729275&4.7.11.23.181.859&925&5088&3517&141732045&4.3.19.41.83.97&187&104 \cr
 & &256.3.25.11.37.53&1961&1408& & &64.11.13.17.19.41&3553&1312 \cr
\noalign{\hrule}
 & &3.7.121.13.19.223&675&898& & &5.13.17.19.29.233&2761&3996 \cr
3500&139960821&4.81.25.7.19.449&19847&16522&3518&141863215&8.27.11.17.37.251&2337&1930 \cr
 & &16.11.89.223.751&751&712& & &32.81.5.19.41.193&3321&3088 \cr
\noalign{\hrule}
 & &3.7.13.29.89.199&89&2& & &9.13.71.109.157&803&11950 \cr
3501&140217987&4.7921.199&3861&4060&3519&142157691&4.25.11.73.239&101&174 \cr
 & &32.27.5.7.11.13.29&55&144& & &16.3.29.101.239&6931&808 \cr
\noalign{\hrule}
 & &3.5.11.37.103.223&163&60& & &3.11.17.89.2857&15847&15580 \cr
3502&140225745&8.9.25.11.37.163&3869&206&3520&142647153&8.5.13.17.19.23.41.53&693&4 \cr
 & &32.53.73.103&53&1168& & &64.9.5.7.11.19.23&2185&672 \cr
\noalign{\hrule}
 & &3.13.101.179.199&5159&12920& & &27.13.29.107.131&1333&370 \cr
3503&140311119&16.5.7.11.17.19.67&12469&12402&3521&142679043&4.3.5.29.31.37.43&1441&1256 \cr
 & &64.9.5.13.37.53.337&17861&17760& & &64.11.43.131.157&1727&1376 \cr
\noalign{\hrule}
 & &9.11.31.131.349&111&3950& & &29.71.103.673&351&322 \cr
3504&140311611&4.27.25.37.79&1085&1048&3522&142727821&4.27.7.13.23.71.103&12035&2662 \cr
 & &64.125.7.31.131&125&224& & &16.3.5.1331.29.83&3993&3320 \cr
\noalign{\hrule}
 & &27.31.359.467&11869&740& & &25.13.29.59.257&1071&814 \cr
3505&140325561&8.5.11.13.37.83&521&558&3523&142911275&4.9.5.7.11.17.37.59&311&2494 \cr
 & &32.9.5.11.31.521&521&880& & &16.3.7.29.43.311&933&2408 \cr
\noalign{\hrule}
 & &5.13.17.73.1741&17171&12426& & &27.23.47.59.83&65569&66650 \cr
3506&140437765&4.3.7.11.19.109.223&73&150&3524&142928739&4.25.7.17.19.29.31.43&141&1474 \cr
 & &16.9.25.19.73.109&981&760& & &16.3.5.7.11.29.47.67&2345&2552 \cr
\noalign{\hrule}
 & &11.13.23.31.1381&2961&12230& & &5.7.11.23.31.521&137513&136992 \cr
3507&140805379&4.9.5.7.47.1223&541&682&3525&143017105&64.3.17.1427.8089&8085&16174 \cr
 & &16.3.5.7.11.31.541&541&840& & &256.9.5.49.11.8087&8087&8064 \cr
\noalign{\hrule}
 & &9.5.11.41.53.131&1075&3248& & &13.17.19.73.467&221753&221430 \cr
3508&140908185&32.3.125.7.29.43&689&1936&3526&143148109&4.3.5.7.121.61.79.401&4615&204 \cr
 & &1024.121.13.53&143&512& & &32.9.25.7.11.13.17.71&4473&4400 \cr
\noalign{\hrule}
 & &3.289.23.73.97&1&290& & &25.7.11.29.31.83&999&86 \cr
3509&141202221&4.5.23.29.73&1001&1116&3527&143637725&4.27.5.29.37.43&77&68 \cr
 & &32.9.7.11.13.31&217&6864& & &32.3.7.11.17.37.43&1591&816 \cr
\noalign{\hrule}
 & &27.19.409.673&7007&5780& & &3.49.13.17.19.233&407&426 \cr
3510&141206841&8.9.5.49.11.13.289&109&10&3528&143819949&4.9.11.13.37.71.233&6125&4028 \cr
 & &32.25.7.13.17.109&12971&5200& & &32.125.49.19.37.53&1961&2000 \cr
\noalign{\hrule}
}%
}
$$
\eject
\vglue -23 pt
\noindent\hskip 1 in\hbox to 6.5 in{\ 3529 -- 3564 \hfill\fbd 143871651 -- 148062915\frb}
\vskip -9 pt
$$
\vbox{
\nointerlineskip
\halign{\strut
    \vrule \ \ \hfil \frb #\ 
   &\vrule \hfil \ \ \fbb #\frb\ 
   &\vrule \hfil \ \ \frb #\ \hfil
   &\vrule \hfil \ \ \frb #\ 
   &\vrule \hfil \ \ \frb #\ \ \vrule \hskip 2 pt
   &\vrule \ \ \hfil \frb #\ 
   &\vrule \hfil \ \ \fbb #\frb\ 
   &\vrule \hfil \ \ \frb #\ \hfil
   &\vrule \hfil \ \ \frb #\ 
   &\vrule \hfil \ \ \frb #\ \vrule \cr%
\noalign{\hrule}
 & &9.7.11.31.37.181&10015&15626& & &9.25.7.121.13.59&33263&30988 \cr
3529&143871651&4.5.13.601.2003&701&1302&3547&146171025&8.29.31.37.61.127&713&2970 \cr
 & &16.3.5.7.13.31.701&701&520& & &32.27.5.11.23.961&961&1104 \cr
\noalign{\hrule}
 & &25.7.11.37.43.47&923&1098& & &27.5.7.13.79.151&4351&6314 \cr
3530&143945725&4.9.11.13.37.61.71&43&2300&3548&146547765&4.49.11.19.41.229&9815&426 \cr
 & &32.3.25.13.23.43&39&368& & &16.3.5.13.71.151&71&8 \cr
\noalign{\hrule}
 & &9.7.13.19.47.197&5115&6494& & &9.25.59.61.181&143&82 \cr
3531&144079299&4.27.5.11.17.31.191&12367&1862&3549&146569275&4.11.13.41.59.181&61&120 \cr
 & &16.49.19.83.149&581&1192& & &64.3.5.11.13.41.61&533&352 \cr
\noalign{\hrule}
 & &27.125.11.169.23&1353&1522& & &27.5.7.13.17.19.37&463&482 \cr
3532&144304875&4.81.121.41.761&20501&10700&3550&146818035&4.13.17.37.241.463&1881&1252 \cr
 & &32.25.13.19.83.107&1577&1712& & &32.9.11.19.313.463&5093&5008 \cr
\noalign{\hrule}
 & &5.49.11.17.23.137&1697&4848& & &9.5.17.401.479&583&182 \cr
3533&144363065&32.3.7.101.1697&495&1202&3551&146940435&4.7.11.13.53.479&425&54 \cr
 & &128.27.5.11.601&601&1728& & &16.27.25.11.13.17&39&440 \cr
\noalign{\hrule}
 & &9.5.29.31.43.83&1417&1156& & &121.23.101.523&135&388 \cr
3534&144383895&8.5.13.289.43.109&93&638&3552&147006409&8.27.5.11.97.101&1&100 \cr
 & &32.3.11.13.17.29.31&143&272& & &64.3.125.97&12125&96 \cr
\noalign{\hrule}
 & &7.11.29.71.911&615&296& & &25.121.41.1187&919&2106 \cr
3535&144432673&16.3.5.7.37.41.71&2001&484&3553&147217675&4.81.13.41.919&275&644 \cr
 & &128.9.121.23.29&207&704& & &32.9.25.7.11.13.23&819&368 \cr
\noalign{\hrule}
 & &49.11.331.811&169&162& & &9.5.7.61.79.97&1181&4736 \cr
3536&144689699&4.81.7.11.169.811&5627&50&3554&147244545&256.7.37.1181&587&594 \cr
 & &16.27.25.17.331&675&136& & &1024.27.11.37.587&19371&18944 \cr
\noalign{\hrule}
 & &49.13.29.47.167&963&1340& & &27.11.17.163.179&923&760 \cr
3537&144994577&8.9.5.67.107.167&77&244&3555&147314673&16.3.5.13.19.71.179&17&196 \cr
 & &64.3.5.7.11.61.67&3685&5856& & &128.5.49.13.17.19&3185&1216 \cr
\noalign{\hrule}
 & &9.5.17.373.509&127&382& & &3.11.13.19.101.179&925&388 \cr
3538&145240605&4.3.127.191.373&473&100&3556&147361929&8.25.11.19.37.97&567&358 \cr
 & &32.25.11.43.127&5461&880& & &32.81.7.97.179&679&432 \cr
\noalign{\hrule}
 & &5.17.53.59.547&2931&196& & &27.7.13.19.29.109&109&242 \cr
3539&145389865&8.3.49.17.977&429&548&3557&147564963&4.121.29.11881&7695&4186 \cr
 & &64.9.7.11.13.137&12467&3168& & &16.81.5.7.13.19.23&69&40 \cr
\noalign{\hrule}
 & &5.121.19.53.239&2403&2138& & &9.17.19.43.1181&11413&11026 \cr
3540&145607165&4.27.121.89.1069&105&16&3558&147626181&4.17.37.101.113.149&3135&602 \cr
 & &128.81.5.7.1069&7483&5184& & &16.3.5.7.11.19.43.113&565&616 \cr
\noalign{\hrule}
 & &9.5.49.169.17.23&359&814& & &9.11.43.79.439&259&180 \cr
3541&145704195&4.3.7.11.13.37.359&535&1978&3559&147637017&8.81.5.7.11.37.43&439&34 \cr
 & &16.5.11.23.43.107&1177&344& & &32.7.17.37.439&119&592 \cr
\noalign{\hrule}
 & &3.17.31.37.47.53&1221&422& & &5.49.23.73.359&15921&10286 \cr
3542&145716027&4.9.11.1369.211&265&1634&3560&147676445&4.9.29.37.61.139&1679&1540 \cr
 & &16.5.11.19.43.53&215&1672& & &32.3.5.7.11.23.61.73&183&176 \cr
\noalign{\hrule}
 & &13.17.37.71.251&2945&318& & &27.25.23.31.307&49&44 \cr
3543&145722317&4.3.5.17.19.31.53&303&286&3561&147751425&8.9.5.49.11.23.307&5263&1798 \cr
 & &16.9.5.11.13.53.101&5247&4040& & &32.7.19.29.31.277&3857&4432 \cr
\noalign{\hrule}
 & &27.11.13.23.31.53&23825&20536& & &81.13.19.83.89&1309&230 \cr
3544&145903329&16.25.17.151.953&689&264&3562&147791709&4.5.7.11.17.23.89&171&82 \cr
 & &256.3.11.13.53.151&151&128& & &16.9.5.7.17.19.41&205&952 \cr
\noalign{\hrule}
 & &9.5.7.11.73.577&26737&23852& & &5.7.11.23.59.283&9&68 \cr
3545&145949265&8.67.89.26737&10387&16350&3563&147851935&8.9.5.17.23.283&1927&2318 \cr
 & &32.3.25.13.17.47.109&9265&9776& & &32.3.19.41.47.61&7503&14288 \cr
\noalign{\hrule}
 & &81.121.13.31.37&1297&1700& & &9.5.7.11.13.19.173&947&782 \cr
3546&146142711&8.25.121.17.1297&1677&380&3564&148062915&4.3.17.23.173.947&2891&50 \cr
 & &64.3.125.13.19.43&2375&1376& & &16.25.49.23.59&2065&184 \cr
\noalign{\hrule}
}%
}
$$
\eject
\vglue -23 pt
\noindent\hskip 1 in\hbox to 6.5 in{\ 3565 -- 3600 \hfill\fbd 148337397 -- 152920351\frb}
\vskip -9 pt
$$
\vbox{
\nointerlineskip
\halign{\strut
    \vrule \ \ \hfil \frb #\ 
   &\vrule \hfil \ \ \fbb #\frb\ 
   &\vrule \hfil \ \ \frb #\ \hfil
   &\vrule \hfil \ \ \frb #\ 
   &\vrule \hfil \ \ \frb #\ \ \vrule \hskip 2 pt
   &\vrule \ \ \hfil \frb #\ 
   &\vrule \hfil \ \ \fbb #\frb\ 
   &\vrule \hfil \ \ \frb #\ \hfil
   &\vrule \hfil \ \ \frb #\ 
   &\vrule \hfil \ \ \frb #\ \vrule \cr%
\noalign{\hrule}
 & &9.13.67.127.149&97&30& & &3.11.17.19.37.383&2097&2116 \cr
3565&148337397&4.27.5.13.97.149&341&2278&3583&151048689&8.27.17.529.37.233&625&4 \cr
 & &16.5.11.17.31.67&155&1496& & &64.625.23.233&14375&7456 \cr
\noalign{\hrule}
 & &3.17.29.31.41.79&627&644& & &3.25.11.23.79.101&2171&354 \cr
3566&148504911&8.9.7.11.19.23.29.79&205&2086&3584&151401525&4.9.11.13.59.167&2533&4370 \cr
 & &32.5.49.23.41.149&3427&3920& & &16.5.17.19.23.149&323&1192 \cr
\noalign{\hrule}
 & &9.5.11.67.4481&2323&2158& & &3.17.19.37.41.103&11&30 \cr
3567&148612365&4.3.13.23.67.83.101&895&646&3585&151407219&4.9.5.11.17.37.103&2747&1064 \cr
 & &16.5.13.17.19.101.179&24947&24344& & &64.5.7.19.41.67&469&160 \cr
\noalign{\hrule}
 & &27.7.11.73.983&107677&107600& & &9.5.7.121.29.137&3007&3692 \cr
3568&149186961&32.9.25.29.47.79.269&65&2356&3586&151430895&8.3.11.13.31.71.97&851&1918 \cr
 & &256.125.13.19.31.47&76375&75392& & &32.7.23.31.37.137&713&592 \cr
\noalign{\hrule}
 & &9.101.139.1181&221&82& & &27.5.7.11.13.19.59&47&124 \cr
3569&149220531&4.3.13.17.41.1181&209&1390&3587&151486335&8.3.5.13.31.47.59&817&1012 \cr
 & &16.5.11.17.19.139&1045&136& & &64.11.19.23.43.47&1081&1376 \cr
\noalign{\hrule}
 & &27.11.13.23.1681&245&206& & &7.13.29.71.809&369&440 \cr
3570&149277843&4.9.5.49.23.41.103&247&40&3588&151581521&16.9.5.7.11.13.29.41&2627&3236 \cr
 & &64.25.7.13.19.103&2575&4256& & &128.3.5.37.71.809&185&192 \cr
\noalign{\hrule}
 & &9.49.11.13.23.103&15745&15052& & &3.25.13.29.31.173&1267&748 \cr
3571&149396247&8.5.7.47.53.67.71&73&2418&3589&151638825&8.5.7.11.17.29.181&117&202 \cr
 & &32.3.13.31.71.73&2263&1136& & &32.9.7.13.101.181&2121&2896 \cr
\noalign{\hrule}
 & &3.5.7.37.137.281&137&122& & &9.5.43.181.433&539&106 \cr
3572&149560845&4.61.18769.281&17955&814&3590&151651755&4.3.49.11.53.181&4891&2900 \cr
 & &16.27.5.7.11.19.37&209&72& & &32.25.29.67.73&1943&5840 \cr
\noalign{\hrule}
 & &9.5.49.11.31.199&26413&29108& & &5.13.17.37.47.79&803&1602 \cr
3573&149629095&8.19.61.383.433&363&796&3591&151806005&4.9.11.73.79.89&163&74 \cr
 & &64.3.121.199.383&383&352& & &16.3.11.37.73.163&2409&1304 \cr
\noalign{\hrule}
 & &3.23.643.3373&20497&23870& & &27.11.289.29.61&9959&9500 \cr
3574&149649891&4.5.7.11.31.103.199&2691&502&3592&151838577&8.125.17.19.23.433&429&4 \cr
 & &16.9.5.7.13.23.251&1365&2008& & &64.3.5.11.13.19.23&2185&416 \cr
\noalign{\hrule}
 & &5.17.53.79.421&7799&14514& & &3.11.13.17.83.251&1341&1090 \cr
3575&149831795&4.3.11.41.59.709&9&50&3593&151935069&4.27.5.83.109.149&1493&748 \cr
 & &16.27.25.11.709&1485&5672& & &32.11.17.109.1493&1493&1744 \cr
\noalign{\hrule}
 & &9.25.7.169.563&12361&1714& & &5.49.11.13.43.101&17&522 \cr
3576&149856525&4.47.263.857&297&560&3594&152157005&4.9.13.17.29.43&3269&3140 \cr
 & &128.27.5.7.11.47&141&704& & &32.3.5.7.157.467&1401&2512 \cr
\noalign{\hrule}
 & &9.25.7.11.13.23.29&1963&2062& & &9.5.7.11.13.31.109&2939&3056 \cr
3577&150225075&4.169.29.151.1031&89309&84930&3595&152207055&32.7.31.191.2939&13&204 \cr
 & &16.3.5.11.19.23.149.353&2831&2824& & &256.3.13.17.2939&2939&2176 \cr
\noalign{\hrule}
 & &25.31.131.1481&18447&18578& & &9.11.37.61.683&15487&9340 \cr
3578&150358525&4.3.7.11.13.31.43.1327&3&1330&3596&152611569&8.5.17.467.911&689&222 \cr
 & &16.9.5.49.11.13.19&5733&1672& & &32.3.5.13.17.37.53&689&1360 \cr
\noalign{\hrule}
 & &11.19.173.4159&21231&24518& & &9.5.121.289.97&767&1834 \cr
3579&150376963&4.9.7.13.23.41.337&6897&6920&3597&152639685&4.5.7.11.13.59.131&97&552 \cr
 & &64.27.5.7.121.13.19.173&2079&2080& & &64.3.23.97.131&131&736 \cr
\noalign{\hrule}
 & &25.121.17.29.101&2493&436& & &3.19.23.113.1031&9625&9964 \cr
3580&150623825&8.9.25.109.277&101&176&3598&152735433&8.125.7.11.23.47.53&1197&22 \cr
 & &256.3.11.101.109&109&384& & &32.9.5.49.121.19&1815&784 \cr
\noalign{\hrule}
 & &31.71.89.769&429&340& & &9.7.29.241.347&14425&7436 \cr
3581&150638641&8.3.5.11.13.17.31.71&299&228&3599&152786529&8.25.11.169.577&641&1218 \cr
 & &64.9.5.11.169.19.23&35321&33120& & &32.3.25.7.29.641&641&400 \cr
\noalign{\hrule}
 & &9.121.31.41.109&2779&1690& & &61.11881.211&495&12376 \cr
3582&150868971&4.5.7.169.31.397&1591&1188&3600&152920351&16.9.5.7.11.13.17&109&122 \cr
 & &32.27.5.11.13.37.43&2405&2064& & &64.3.5.17.61.109&51&160 \cr
\noalign{\hrule}
}%
}
$$
\eject
\vglue -23 pt
\noindent\hskip 1 in\hbox to 6.5 in{\ 3601 -- 3636 \hfill\fbd 153295461 -- 157415427\frb}
\vskip -9 pt
$$
\vbox{
\nointerlineskip
\halign{\strut
    \vrule \ \ \hfil \frb #\ 
   &\vrule \hfil \ \ \fbb #\frb\ 
   &\vrule \hfil \ \ \frb #\ \hfil
   &\vrule \hfil \ \ \frb #\ 
   &\vrule \hfil \ \ \frb #\ \ \vrule \hskip 2 pt
   &\vrule \ \ \hfil \frb #\ 
   &\vrule \hfil \ \ \fbb #\frb\ 
   &\vrule \hfil \ \ \frb #\ \hfil
   &\vrule \hfil \ \ \frb #\ 
   &\vrule \hfil \ \ \frb #\ \vrule \cr%
\noalign{\hrule}
 & &9.11.71.113.193&98531&97750& & &3.5.11.67.101.139&217&520 \cr
3601&153295461&4.125.17.23.37.2663&2769&106&3619&155201145&16.25.7.13.31.139&1809&3616 \cr
 & &16.3.13.17.37.53.71&1961&1768& & &1024.27.67.113&1017&512 \cr
\noalign{\hrule}
 & &5.361.841.101&18245&18216& & &27.5.7.23.37.193&97&902 \cr
3602&153318505&16.9.25.11.23.29.41.89&247&2828&3620&155209635&4.11.41.97.193&437&630 \cr
 & &128.3.7.11.13.19.23.101&2093&2112& & &16.9.5.7.19.23.41&41&152 \cr
\noalign{\hrule}
 & &3.5.11.53.89.197&357&622& & &25.49.11.83.139&969&12506 \cr
3603&153326085&4.9.7.17.197.311&275&3074&3621&155461075&4.3.169.17.19.37&231&250 \cr
 & &16.25.7.11.29.53&35&232& & &16.9.125.7.11.13.17&221&360 \cr
\noalign{\hrule}
 & &3.25.11.13.17.841&2921&2734& & &3.13.47.73.1163&335&1498 \cr
3604&153335325&4.5.23.29.127.1367&9&136&3622&155619867&4.5.7.67.73.107&423&88 \cr
 & &64.9.17.23.1367&1367&2208& & &64.9.11.47.107&107&1056 \cr
\noalign{\hrule}
 & &3.31.41.131.307&6097&6490& & &49.13.41.59.101&75&26 \cr
3605&153347421&4.5.7.11.13.31.59.67&1107&722&3623&155631203&4.3.25.169.41.59&2727&4202 \cr
 & &16.27.13.361.41.67&4693&4824& & &16.81.11.101.191&891&1528 \cr
\noalign{\hrule}
 & &7.11.17.233.503&155&78& & &27.5.7.13.19.23.29&337&482 \cr
3606&153413491&4.3.5.13.17.31.503&209&294&3624&155687805&4.3.19.23.241.337&143&580 \cr
 & &16.9.49.11.13.19.31&1729&2232& & &32.5.11.13.29.337&337&176 \cr
\noalign{\hrule}
 & &169.769.1181&675&506& & &27.5.7.37.61.73&2321&2132 \cr
3607&153483941&4.27.25.11.23.769&247&522&3625&155699145&8.5.11.13.37.41.211&13&198 \cr
 & &16.243.13.19.23.29&5589&4408& & &32.9.121.169.41&4961&2704 \cr
\noalign{\hrule}
 & &49.23.31.53.83&2119&2280& & &25.7.13.17.29.139&23&198 \cr
3608&153687863&16.3.5.7.13.19.31.163&2959&138&3626&155898925&4.9.11.23.29.139&931&1070 \cr
 & &64.9.5.11.23.269&2959&1440& & &16.3.5.49.11.19.107&1463&2568 \cr
\noalign{\hrule}
 & &27.7.61.67.199&935&874& & &27.13.19.67.349&137&110 \cr
3609&153716157&4.5.7.11.17.19.23.199&159&40&3627&155941227&4.5.11.67.137.349&7&342 \cr
 & &64.3.25.11.19.23.53&13409&15200& & &16.9.7.11.19.137&959&88 \cr
\noalign{\hrule}
 & &3.121.53.61.131&14105&6724& & &5.11.53.59.907&481&426 \cr
3610&153738849&8.5.7.13.31.1681&99&304&3628&155990395&4.3.13.37.53.59.71&2877&250 \cr
 & &256.9.7.11.19.41&2337&896& & &16.9.125.7.13.137&3425&6552 \cr
\noalign{\hrule}
 & &3.11.13.19.113.167&155&1992& & &625.13.71.271&1449&2074 \cr
3611&153817521&16.9.5.13.31.83&251&334&3629&156333125&4.9.7.17.23.61.71&6721&3100 \cr
 & &64.31.167.251&251&992& & &32.3.25.11.13.31.47&1023&752 \cr
\noalign{\hrule}
 & &5.7.11.19.107.197&655&522& & &3.25.121.19.907&3417&1118 \cr
3612&154192885&4.9.25.29.131.197&1219&4494&3630&156389475&4.9.5.13.17.43.67&361&26 \cr
 & &16.27.7.23.53.107&621&424& & &16.169.17.361&2873&152 \cr
\noalign{\hrule}
 & &3.13.29.31.53.83&2055&518& & &5.11.97.149.197&1809&12644 \cr
3613&154233339&4.9.5.7.13.37.137&913&868&3631&156598255&8.27.29.67.109&97&164 \cr
 & &32.49.11.31.37.83&539&592& & &64.3.41.97.109&327&1312 \cr
\noalign{\hrule}
 & &9.25.121.53.107&481&844& & &9.125.7.11.1811&793&782 \cr
3614&154392975&8.3.13.37.107.211&2135&6094&3632&156877875&4.5.13.17.23.61.1811&1653&158 \cr
 & &32.5.7.11.61.277&1939&976& & &16.3.17.19.29.61.79&19703&18328 \cr
\noalign{\hrule}
 & &3.11.17.29.37.257&665&408& & &9.19.67.71.193&833&904 \cr
3615&154701921&16.9.5.7.11.289.19&95&194&3633&156995271&16.49.17.19.67.113&1053&220 \cr
 & &64.25.7.361.97&16975&11552& & &128.81.5.11.13.113&6215&7488 \cr
\noalign{\hrule}
 & &81.5.7.193.283&187&380& & &9.5.13.37.53.137&41053&39092 \cr
3616&154844865&8.25.11.17.19.283&207&5018&3634&157164345&8.29.61.337.673&1221&548 \cr
 & &32.9.13.23.193&299&16& & &64.3.11.37.137.337&337&352 \cr
\noalign{\hrule}
 & &9.25.7.131.751&3419&1166& & &5.7.11.13.23.1367&853&918 \cr
3617&154950075&4.3.5.11.13.53.263&427&262&3635&157362205&4.27.17.853.1367&5725&6578 \cr
 & &16.7.61.131.263&263&488& & &16.3.25.11.13.17.23.229&687&680 \cr
\noalign{\hrule}
 & &9.13.19.101.691&913&1160& & &27.13.17.23.31.37&125&1022 \cr
3618&155145393&16.3.5.11.29.83.101&95&2834&3636&157415427&4.9.125.7.17.73&473&598 \cr
 & &64.25.13.19.109&109&800& & &16.11.13.23.43.73&803&344 \cr
\noalign{\hrule}
}%
}
$$
\eject
\vglue -23 pt
\noindent\hskip 1 in\hbox to 6.5 in{\ 3637 -- 3672 \hfill\fbd 157437885 -- 161868949\frb}
\vskip -9 pt
$$
\vbox{
\nointerlineskip
\halign{\strut
    \vrule \ \ \hfil \frb #\ 
   &\vrule \hfil \ \ \fbb #\frb\ 
   &\vrule \hfil \ \ \frb #\ \hfil
   &\vrule \hfil \ \ \frb #\ 
   &\vrule \hfil \ \ \frb #\ \ \vrule \hskip 2 pt
   &\vrule \ \ \hfil \frb #\ 
   &\vrule \hfil \ \ \fbb #\frb\ 
   &\vrule \hfil \ \ \frb #\ \hfil
   &\vrule \hfil \ \ \frb #\ 
   &\vrule \hfil \ \ \frb #\ \vrule \cr%
\noalign{\hrule}
 & &3.5.11.71.89.151&1091&1846& & &9.11.19.29.37.79&251&460 \cr
3637&157437885&4.13.5041.1091&4571&9612&3655&159446727&8.5.23.29.37.251&5767&1512 \cr
 & &32.27.7.89.653&653&1008& & &128.27.7.73.79&511&192 \cr
\noalign{\hrule}
 & &27.5.37.139.227&95&322& & &243.25.7.121.31&1391&1634 \cr
3638&157607235&4.9.25.7.19.23.37&121&454&3656&159511275&4.7.13.19.31.43.107&1531&198 \cr
 & &16.7.121.19.227&847&152& & &16.9.11.107.1531&1531&856 \cr
\noalign{\hrule}
 & &9.13.19.29.31.79&1397&370& & &3.13.29.113.1249&1359&110 \cr
3639&157879683&4.3.5.11.29.37.127&437&118&3657&159625947&4.27.5.11.29.151&439&1222 \cr
 & &16.19.23.59.127&2921&472& & &16.5.13.47.439&235&3512 \cr
\noalign{\hrule}
 & &7.13.19.29.47.67&111&440& & &13.19.41.97.163&12969&2842 \cr
3640&157894009&16.3.5.11.13.37.67&207&274&3658&160117997&4.9.49.11.29.131&39&10 \cr
 & &64.27.5.11.23.137&15755&9504& & &16.27.5.11.13.131&3537&440 \cr
\noalign{\hrule}
 & &5.11.13.37.47.127&177&304& & &25.49.11.73.163&3&52 \cr
3641&157909895&32.3.5.11.19.47.59&13&222&3659&160339025&8.3.5.13.73.163&67&882 \cr
 & &128.9.13.37.59&59&576& & &32.27.49.67&1809&16 \cr
\noalign{\hrule}
 & &3.5.7.41.73.503&231&272& & &9.41.43.67.151&25&176 \cr
3642&158075295&32.9.5.49.11.17.73&403&38&3660&160526439&32.3.25.11.41.43&161&290 \cr
 & &128.11.13.17.19.31&6479&14144& & &128.125.7.23.29&2875&12992 \cr
\noalign{\hrule}
 & &25.11.13.17.19.137&447&238& & &3.125.11.17.29.79&291&104 \cr
3643&158197325&4.3.5.7.13.289.149&1691&246&3661&160656375&16.9.25.13.29.97&1343&1082 \cr
 & &16.9.7.19.41.89&623&2952& & &64.13.17.79.541&541&416 \cr
\noalign{\hrule}
 & &13.97.269.467&10011&16082& & &25.11.41.53.269&3071&3654 \cr
3644&158410603&4.3.11.17.43.47.71&1401&620&3662&160747675&4.9.7.29.37.41.83&967&106 \cr
 & &32.9.5.17.31.467&765&496& & &16.3.53.83.967&967&1992 \cr
\noalign{\hrule}
 & &9.5.13.19.53.269&1751&1694& & &5.7.17.43.61.103&317&198 \cr
3645&158466555&4.3.7.121.17.103.269&13091&628&3663&160750555&4.9.11.43.61.317&95&34 \cr
 & &32.7.13.19.53.157&157&112& & &16.3.5.11.17.19.317&951&1672 \cr
\noalign{\hrule}
 & &41.79.173.283&55&228& & &9.5.7.11.13.43.83&289&206 \cr
3646&158578201&8.3.5.11.19.41.79&265&186&3664&160765605&4.7.13.289.43.103&339&220 \cr
 & &32.9.25.19.31.53&14725&7632& & &32.3.5.11.17.103.113&1921&1648 \cr
\noalign{\hrule}
 & &3.25.7.17.23.773&591&1364& & &11.29.37.43.317&45&362 \cr
3647&158677575&8.9.5.7.11.31.197&13&328&3665&160886693&4.9.5.29.43.181&337&308 \cr
 & &128.13.41.197&533&12608& & &32.3.7.11.181.337&3801&5392 \cr
\noalign{\hrule}
 & &3.5.23.67.6869&3377&3492& & &9.11.29.47.1193&65&76 \cr
3648&158776935&8.27.11.67.97.307&751&1058&3666&160979841&8.3.5.13.19.29.1193&31&1162 \cr
 & &32.11.529.97.751&17273&17072& & &32.5.7.19.31.83&12865&2128 \cr
\noalign{\hrule}
 & &81.7.169.1657&1045&12644& & &3.25.11.401.487&3423&988 \cr
3649&158778711&8.5.11.19.29.109&221&330&3667&161111775&8.9.5.7.13.19.163&401&1066 \cr
 & &32.3.25.121.13.17&2057&400& & &32.169.41.401&169&656 \cr
\noalign{\hrule}
 & &243.19.163.211&203&40& & &3.25.13.37.41.109&17&92 \cr
3650&158792481&16.5.7.19.29.211&1903&2106&3668&161219175&8.13.17.23.37.41&261&220 \cr
 & &64.81.5.11.13.173&2249&1760& & &64.9.5.11.17.23.29&4301&2784 \cr
\noalign{\hrule}
 & &5.7.29.349.449&333&682& & &3.13.43.79.1217&95&1122 \cr
3651&159051515&4.9.11.31.37.449&349&798&3669&161231811&4.9.5.11.17.19.43&1217&1690 \cr
 & &16.27.7.11.19.349&297&152& & &16.25.169.1217&25&104 \cr
\noalign{\hrule}
 & &27.11.19.71.397&2219&2148& & &5.193.349.479&657&308 \cr
3652&159059241&8.81.7.19.179.317&859&680&3670&161320015&8.9.7.11.73.479&349&130 \cr
 & &128.5.7.17.317.859&102221&101440& & &32.3.5.7.11.13.349&143&336 \cr
\noalign{\hrule}
 & &5.49.13.23.41.53&551&2046& & &9.7.11.23.73.139&1027&224 \cr
3653&159183115&4.3.11.19.29.31.41&255&296&3671&161733033&64.49.13.23.79&2085&1786 \cr
 & &64.9.5.11.17.31.37&6919&8928& & &256.3.5.19.47.139&893&640 \cr
\noalign{\hrule}
 & &5.19.53.103.307&825&1132& & &11.31.479.991&69&410 \cr
3654&159211735&8.3.125.11.53.283&4869&1756&3672&161868949&4.3.5.23.41.991&507&484 \cr
 & &64.27.439.541&11853&17312& & &32.9.5.121.169.41&6929&7920 \cr
\noalign{\hrule}
}%
}
$$
\eject
\vglue -23 pt
\noindent\hskip 1 in\hbox to 6.5 in{\ 3673 -- 3708 \hfill\fbd 161889585 -- 165069905\frb}
\vskip -9 pt
$$
\vbox{
\nointerlineskip
\halign{\strut
    \vrule \ \ \hfil \frb #\ 
   &\vrule \hfil \ \ \fbb #\frb\ 
   &\vrule \hfil \ \ \frb #\ \hfil
   &\vrule \hfil \ \ \frb #\ 
   &\vrule \hfil \ \ \frb #\ \ \vrule \hskip 2 pt
   &\vrule \ \ \hfil \frb #\ 
   &\vrule \hfil \ \ \fbb #\frb\ 
   &\vrule \hfil \ \ \frb #\ \hfil
   &\vrule \hfil \ \ \frb #\ 
   &\vrule \hfil \ \ \frb #\ \vrule \cr%
\noalign{\hrule}
 & &3.5.11.13.71.1063&1985&1204& & &3.13.47.257.347&215&826 \cr
3673&161889585&8.25.7.13.43.397&6919&3006&3691&163465107&4.5.7.43.59.257&99&158 \cr
 & &32.9.11.17.37.167&2839&1776& & &16.9.5.7.11.43.79&8295&3784 \cr
\noalign{\hrule}
 & &9.11.17.23.47.89&1865&182& & &3.7.23.31.67.163&671&470 \cr
3674&161919747&4.5.7.13.47.373&1423&1188&3692&163520133&4.5.11.23.31.47.61&27&4 \cr
 & &32.27.11.13.1423&1423&624& & &32.27.5.11.47.61&2867&7920 \cr
\noalign{\hrule}
 & &3.25.23.109.863&31109&33616& & &81.13.17.41.223&415&638 \cr
3675&162265575&32.11.13.191.2393&45&2438&3693&163668843&4.5.11.17.29.41.83&217&234 \cr
 & &128.9.5.11.23.53&159&704& & &16.9.5.7.13.29.31.83&4495&4648 \cr
\noalign{\hrule}
 & &3.5.11.59.79.211&963&1358& & &9.11.29.127.449&12797&224 \cr
3676&162272715&4.27.7.59.97.107&1171&422&3694&163713033&64.7.67.191&635&702 \cr
 & &16.97.211.1171&1171&776& & &256.27.5.13.127&65&384 \cr
\noalign{\hrule}
 & &3.5.11.31.61.521&103&568& & &9.7.11.349.677&13775&13098 \cr
3677&162559815&16.71.103.521&225&296&3695&163737189&4.27.25.19.29.37.59&349&164 \cr
 & &256.9.25.37.103&3811&1920& & &32.5.29.41.59.349&2419&2320 \cr
\noalign{\hrule}
 & &27.7.41.139.151&2405&1348& & &27.25.37.79.83&1093&982 \cr
3678&162643761&8.5.13.37.41.337&261&220&3696&163761075&4.9.79.491.1093&601&110 \cr
 & &64.9.25.11.29.337&9773&8800& & &16.5.11.601.1093&6611&8744 \cr
\noalign{\hrule}
 & &49.19.41.4261&1665&2596& & &7.11.13.17.23.419&1425&1006 \cr
3679&162646631&8.9.5.11.37.41.59&61&62&3697&163992829&4.3.25.7.19.23.503&33&128 \cr
 & &32.3.5.11.31.37.59.61&104005&104784& & &1024.9.5.11.503&2515&4608 \cr
\noalign{\hrule}
 & &7.11.19.31.37.97&3&34& & &289.2209.257&3289&1080 \cr
3680&162771917&4.3.7.11.17.19.97&235&444&3698&164069057&16.27.5.11.13.17.23&1249&1028 \cr
 & &32.9.5.17.37.47&2115&272& & &128.3.5.257.1249&1249&960 \cr
\noalign{\hrule}
 & &5.11.17.19.89.103&1469&3426& & &5.11.19.71.2213&1781&432 \cr
3681&162851755&4.3.13.17.113.571&675&1246&3699&164193535&32.27.5.11.13.137&359&326 \cr
 & &16.81.25.7.13.89&405&728& & &128.9.13.163.359&19071&22976 \cr
\noalign{\hrule}
 & &3.11.107.193.239&1475&1154& & &25.11.61.97.101&3211&2706 \cr
3682&162874437&4.25.59.193.577&67&126&3700&164344675&4.3.5.121.169.19.41&183&2482 \cr
 & &16.9.25.7.67.577&11725&13848& & &16.9.13.17.61.73&657&1768 \cr
\noalign{\hrule}
 & &5.7.23.31.47.139&1521&1676& & &3.49.457.2447&25465&25922 \cr
3683&163031015&8.9.7.169.47.419&69&22&3701&164387013&4.5.7.11.13.463.997&4113&872 \cr
 & &32.27.11.13.23.419&5447&4752& & &64.9.11.13.109.457&1199&1248 \cr
\noalign{\hrule}
 & &11.17.23.83.457&261&652& & &9.11.17.23.31.137&2303&26 \cr
3684&163141231&8.9.29.163.457&147&310&3702&164397123&4.49.13.31.47&115&102 \cr
 & &32.27.5.49.29.31&6615&14384& & &16.3.5.7.17.23.47&35&376 \cr
\noalign{\hrule}
 & &3.5.11.47.53.397&551&154& & &3.7.11.17.163.257&2285&486 \cr
3685&163172955&4.7.121.19.29.53&397&450&3703&164505957&4.729.5.11.457&3781&4238 \cr
 & &16.9.25.19.29.397&285&232& & &16.5.13.19.163.199&1235&1592 \cr
\noalign{\hrule}
 & &3.5.13.17.53.929&473&456& & &27.11.13.17.23.109&3437&4636 \cr
3686&163220655&16.9.5.11.13.19.43.53&6503&4118&3704&164551959&8.7.17.19.61.491&1035&8294 \cr
 & &64.7.11.29.71.929&2233&2272& & &32.9.5.11.13.23.29&29&80 \cr
\noalign{\hrule}
 & &3.11.29.281.607&175&782& & &5.11.13.17.19.23.31&123&2308 \cr
3687&163232619&4.25.7.17.23.281&353&72&3705&164663785&8.3.31.41.577&273&304 \cr
 & &64.9.7.23.353&7413&736& & &256.9.7.13.19.41&369&896 \cr
\noalign{\hrule}
 & &27.5.7.11.41.383&145&3302& & &9.7.13.83.2423&3925&3344 \cr
3688&163232685&4.3.25.13.29.127&5687&5362&3706&164708271&32.3.25.11.13.19.157&449&292 \cr
 & &16.7.121.47.383&47&88& & &256.25.11.73.449&32777&35200 \cr
\noalign{\hrule}
 & &3.5.13.31.113.239&113&82& & &125.7.11.17.19.53&599&276 \cr
3689&163257315&4.41.12769.239&1485&11284&3707&164770375&8.3.11.23.53.599&153&100 \cr
 & &32.27.5.7.11.13.31&63&176& & &64.27.25.17.599&599&864 \cr
\noalign{\hrule}
 & &19.24649.349&9009&15640& & &5.7.11.169.43.59&2921&2448 \cr
3690&163447519&16.9.5.7.11.13.17.23&157&38&3708&165069905&32.9.5.13.17.23.127&649&14 \cr
 & &64.3.11.19.23.157&253&96& & &128.3.7.11.23.59&23&192 \cr
\noalign{\hrule}
}%
}
$$
\eject
\vglue -23 pt
\noindent\hskip 1 in\hbox to 6.5 in{\ 3709 -- 3744 \hfill\fbd 165363471 -- 168789005\frb}
\vskip -9 pt
$$
\vbox{
\nointerlineskip
\halign{\strut
    \vrule \ \ \hfil \frb #\ 
   &\vrule \hfil \ \ \fbb #\frb\ 
   &\vrule \hfil \ \ \frb #\ \hfil
   &\vrule \hfil \ \ \frb #\ 
   &\vrule \hfil \ \ \frb #\ \ \vrule \hskip 2 pt
   &\vrule \ \ \hfil \frb #\ 
   &\vrule \hfil \ \ \fbb #\frb\ 
   &\vrule \hfil \ \ \frb #\ \hfil
   &\vrule \hfil \ \ \frb #\ 
   &\vrule \hfil \ \ \frb #\ \vrule \cr%
\noalign{\hrule}
 & &27.7.13.17.37.107&6023&14200& & &7.13.17.23.43.109&3267&1850 \cr
3709&165363471&16.25.19.71.317&3531&3214&3727&166768147&4.27.25.121.23.37&109&224 \cr
 & &64.3.5.11.107.1607&1607&1760& & &256.3.5.7.121.109&605&384 \cr
\noalign{\hrule}
 & &5.121.13.17.1237&8343&7738& & &27.13.23.73.283&145&154 \cr
3710&165393085&4.81.17.53.73.103&11&62&3728&166780107&4.3.5.7.11.29.73.283&817&598 \cr
 & &16.27.11.31.53.103&5459&6696& & &16.7.11.13.19.23.29.43&3857&3784 \cr
\noalign{\hrule}
 & &9.5.7.23.53.431&2867&1012& & &27.7.19.23.43.47&14333&6206 \cr
3711&165497535&8.11.529.47.61&71&600&3729&166920453&4.11.29.107.1303&705&598 \cr
 & &128.3.25.47.71&3337&320& & &16.3.5.11.13.23.29.47&377&440 \cr
\noalign{\hrule}
 & &9.31.37.43.373&209&178& & &5.11.13.29.83.97&297&782 \cr
3712&165570597&4.11.19.37.89.373&405&3698&3730&166937485&4.27.121.17.23.29&985&1072 \cr
 & &16.81.5.19.1849&387&760& & &128.9.5.23.67.197&13869&12608 \cr
\noalign{\hrule}
 & &7.289.23.3559&1975&1584& & &243.59.89.131&1339&12998 \cr
3713&165596711&32.9.25.7.11.17.79&2791&184&3731&167155083&4.13.67.97.103&55&42 \cr
 & &512.3.23.2791&2791&768& & &16.3.5.7.11.67.103&7931&2680 \cr
\noalign{\hrule}
 & &27.5.7.169.17.61&43&22& & &5.7.169.23.1229&1661&432 \cr
3714&165614085&4.9.11.13.17.43.61&577&460&3732&167199305&32.27.5.11.13.151&1037&322 \cr
 & &32.5.11.23.43.577&10879&9232& & &128.3.7.17.23.61&1037&192 \cr
\noalign{\hrule}
 & &7.11.41.137.383&1115&1566& & &25.7.11.13.41.163&603&538 \cr
3715&165650947&4.27.5.29.137.223&2617&3850&3733&167242075&4.9.5.11.41.67.269&61&676 \cr
 & &16.3.125.7.11.2617&2617&3000& & &32.3.169.61.269&2379&4304 \cr
\noalign{\hrule}
 & &81.5.449.911&427&22& & &11.13.17.23.41.73&551&252 \cr
3716&165660795&4.7.11.61.911&425&486&3734&167347609&8.9.7.17.19.29.41&6665&6578 \cr
 & &16.243.25.7.11.17&231&680& & &32.3.5.7.11.13.23.31.43&1505&1488 \cr
\noalign{\hrule}
 & &3.7.11.29.53.467&1475&1794& & &125.11.13.17.19.29&123&370 \cr
3717&165806949&4.9.25.13.23.53.59&1271&1856&3735&167435125&4.3.625.11.37.41&87&538 \cr
 & &512.5.23.29.31.41&6355&5888& & &16.9.29.37.269&2421&296 \cr
\noalign{\hrule}
 & &9.5.11.19.31.569&2855&2266& & &3.25.7.17.19.23.43&1609&2046 \cr
3718&165894795&4.25.121.103.571&1227&1798&3736&167709675&4.9.5.7.11.31.1609&1717&1748 \cr
 & &16.3.29.31.103.409&2987&3272& & &32.17.19.23.101.1609&1609&1616 \cr
\noalign{\hrule}
 & &23.37.109.1789&18557&22590& & &9.5.11.17.19.1049&439&610 \cr
3719&165945851&4.9.5.7.11.241.251&949&1702&3737&167719365&4.25.11.17.61.439&11191&216 \cr
 & &16.3.5.7.13.23.37.73&949&840& & &64.27.361.31&57&992 \cr
\noalign{\hrule}
 & &5.7.11.53.79.103&249&146& & &5.49.61.103.109&28037&28098 \cr
3720&166035485&4.3.7.11.53.73.83&65&12&3738&167787515&4.9.343.529.53.223&3811&80300 \cr
 & &32.9.5.13.73.83&6059&1872& & &32.3.25.11.37.73.103&2701&2640 \cr
\noalign{\hrule}
 & &9.7.83.113.281&15503&16250& & &49.11.19.529.31&369&160 \cr
3721&166036437&4.625.7.13.37.419&33&292&3739&167942159&64.9.5.49.31.41&107&138 \cr
 & &32.3.25.11.73.419&10475&12848& & &256.27.23.41.107&4387&3456 \cr
\noalign{\hrule}
 & &5.23.37.103.379&1089&12934& & &9.25.13.17.31.109&4983&5092 \cr
3722&166102435&4.9.121.29.223&5&6&3740&168020775&8.27.11.17.19.67.151&3815&262 \cr
 & &16.27.5.11.29.223&2453&6264& & &32.5.7.67.109.131&917&1072 \cr
\noalign{\hrule}
 & &3.37.1849.811&341&470& & &3.5.7.13.257.479&2319&1034 \cr
3723&166448829&4.5.11.31.37.43.47&3249&4706&3741&168035595&4.9.11.13.47.773&257&1030 \cr
 & &16.9.11.13.361.181&14079&15928& & &16.5.47.103.257&103&376 \cr
\noalign{\hrule}
 & &27.11.53.71.149&377&404& & &9.5.7.13.17.41.59&197&418 \cr
3724&166524039&8.13.29.53.101.149&1737&200&3742&168398685&4.3.7.11.19.59.197&1703&464 \cr
 & &128.9.25.101.193&4825&6464& & &128.13.19.29.131&3799&1216 \cr
\noalign{\hrule}
 & &7.29.37.67.331&6039&3560& & &25.7.11.13.53.127&1107&218 \cr
3725&166571447&16.9.5.7.11.61.89&37&268&3743&168443275&4.27.11.13.41.109&235&92 \cr
 & &128.3.37.67.89&267&64& & &32.9.5.23.41.47&1081&5904 \cr
\noalign{\hrule}
 & &3.11.13.17.73.313&3239&830& & &5.7.11.289.37.41&1&186 \cr
3726&166637757&4.5.17.41.79.83&9&8&3744&168789005&4.3.7.17.31.41&407&120 \cr
 & &64.9.5.41.79.83&10209&12640& & &64.9.5.11.37&1&288 \cr
\noalign{\hrule}
}%
}
$$
\eject
\vglue -23 pt
\noindent\hskip 1 in\hbox to 6.5 in{\ 3745 -- 3780 \hfill\fbd 168903007 -- 173299269\frb}
\vskip -9 pt
$$
\vbox{
\nointerlineskip
\halign{\strut
    \vrule \ \ \hfil \frb #\ 
   &\vrule \hfil \ \ \fbb #\frb\ 
   &\vrule \hfil \ \ \frb #\ \hfil
   &\vrule \hfil \ \ \frb #\ 
   &\vrule \hfil \ \ \frb #\ \ \vrule \hskip 2 pt
   &\vrule \ \ \hfil \frb #\ 
   &\vrule \hfil \ \ \fbb #\frb\ 
   &\vrule \hfil \ \ \frb #\ \hfil
   &\vrule \hfil \ \ \frb #\ 
   &\vrule \hfil \ \ \frb #\ \vrule \cr%
\noalign{\hrule}
 & &7.13.17.23.47.101&99&200& & &9.19.61.83.197&13&184 \cr
3745&168903007&16.9.25.7.11.17.47&207&592&3763&170557281&16.13.23.61.83&985&924 \cr
 & &512.81.5.23.37&2997&1280& & &128.3.5.7.11.13.197&455&704 \cr
\noalign{\hrule}
 & &3.5.49.41.71.79&969&1040& & &81.17.73.1697&1177&520 \cr
3746&169027215&32.9.25.13.17.19.79&143&568&3764&170584137&16.9.5.11.13.17.107&1537&146 \cr
 & &512.11.169.19.71&3211&2816& & &64.5.29.53.73&1537&160 \cr
\noalign{\hrule}
 & &5.7.11.23.29.659&1215&556& & &5.11.17.41.61.73&67&6 \cr
3747&169227905&8.243.25.29.139&403&322&3765&170705755&4.3.5.11.17.41.67&447&488 \cr
 & &32.3.7.13.23.31.139&1807&1488& & &64.9.61.67.149&1341&2144 \cr
\noalign{\hrule}
 & &3.5.11.17.23.37.71&13091&1376& & &9.13.17.31.47.59&4715&4312 \cr
3748&169480905&64.13.19.43.53&253&306&3766&170980407&16.5.49.11.23.41.47&177&1258 \cr
 & &256.9.11.17.19.23&57&128& & &64.3.7.11.17.37.59&259&352 \cr
\noalign{\hrule}
 & &3.7.11.19.529.73&2253&1450& & &5.7.13.19.43.461&12123&17842 \cr
3749&169490013&4.9.25.19.29.751&803&52&3767&171369835&4.27.11.449.811&1253&3686 \cr
 & &32.5.11.13.29.73&65&464& & &16.9.7.19.97.179&873&1432 \cr
\noalign{\hrule}
 & &9.5.7.13.181.229&1351&1364& & &9.11.17.19.31.173&425&598 \cr
3750&169733655&8.3.49.11.31.193.229&8507&950&3768&171492651&4.3.25.13.289.19.23&173&11098 \cr
 & &32.25.19.31.47.181&1457&1520& & &16.31.173.179&179&8 \cr
\noalign{\hrule}
 & &125.7.11.31.569&148903&149472& & &5.19.53.67.509&387&122 \cr
3751&169775375&64.27.17.19.173.461&8855&17614&3769&171708605&4.9.19.43.61.67&509&308 \cr
 & &256.3.5.7.11.23.8807&8807&8832& & &32.3.7.11.61.509&671&336 \cr
\noalign{\hrule}
 & &3.7.11.43.71.241&737&950& & &3.49.37.131.241&16585&10166 \cr
3752&169963563&4.25.121.19.43.67&1687&612&3770&171714669&4.5.13.17.23.31.107&3725&5544 \cr
 & &32.9.7.17.67.241&201&272& & &64.9.125.7.11.149&4917&4000 \cr
\noalign{\hrule}
 & &5.11.13.17.71.197&517&588& & &47.1103.3313&27577&24264 \cr
3753&170011985&8.3.49.121.47.197&159&38&3771&171749233&16.9.11.23.109.337&1105&94 \cr
 & &32.9.49.19.47.53&20727&16112& & &64.3.5.13.17.23.47&1105&2208 \cr
\noalign{\hrule}
 & &3.25.11.31.61.109&2561&1036& & &3.49.11.157.677&4185&3508 \cr
3754&170048175&8.7.13.31.37.197&2813&3294&3772&171869313&8.81.5.11.31.877&7&884 \cr
 & &32.27.7.29.61.97&1827&1552& & &64.5.7.13.17.31&221&4960 \cr
\noalign{\hrule}
 & &9.7.17.19.61.137&965&1364& & &5.7.121.13.53.59&183&188 \cr
3755&170056593&8.3.5.11.31.61.193&137&442&3773&172156985&8.3.121.13.47.59.61&2173&600 \cr
 & &32.11.13.17.31.137&341&208& & &128.9.25.41.53.61&2745&2624 \cr
\noalign{\hrule}
 & &25.7.13.17.53.83&837&242& & &3.5.19.23.97.271&13301&2146 \cr
3756&170131325&4.27.5.121.31.53&1729&86&3774&172311285&4.29.37.47.283&5917&4554 \cr
 & &16.9.7.13.19.43&43&1368& & &16.9.11.23.61.97&183&88 \cr
\noalign{\hrule}
 & &11.47.191.1723&189&1912& & &3.25.7.29.47.241&387&628 \cr
3757&170141081&16.27.7.47.239&335&382&3775&172453575&8.27.5.43.47.157&91&44 \cr
 & &64.9.5.7.67.191&603&1120& & &64.7.11.13.43.157&6751&4576 \cr
\noalign{\hrule}
 & &25.11.13.181.263&989&1002& & &11.29.59.67.137&787&924 \cr
3758&170180725&4.3.25.23.43.167.263&303&6878&3776&172757959&8.3.7.121.67.787&295&174 \cr
 & &16.9.19.23.101.181&1919&1656& & &32.9.5.29.59.787&787&720 \cr
\noalign{\hrule}
 & &25.7.19.83.617&221&396& & &25.13.19.101.277&145&132 \cr
3759&170276575&8.9.11.13.17.19.83&827&86&3777&172757975&8.3.125.11.19.29.101&277&2652 \cr
 & &32.3.17.43.827&2193&13232& & &64.9.11.13.17.277&187&288 \cr
\noalign{\hrule}
 & &3.5.23.53.67.139&13673&4082& & &9.5.11.71.4919&2779&2140 \cr
3760&170288205&4.121.13.113.157&117&4&3778&172878255&8.25.7.11.107.397&5551&4374 \cr
 & &32.9.169.157&169&7536& & &32.2187.49.13.61&11907&12688 \cr
\noalign{\hrule}
 & &11.17.29.53.593&8307&8890& & &5.49.13.17.31.103&891&860 \cr
3761&170439467&4.9.5.7.13.17.71.127&551&656&3779&172884985&8.81.25.49.11.13.43&31&1256 \cr
 & &128.3.13.19.29.41.127&15621&15808& & &128.9.31.43.157&1413&2752 \cr
\noalign{\hrule}
 & &25.7.13.137.547&15543&8432& & &3.11.13.31.83.157&369&710 \cr
3762&170486225&32.9.11.17.31.157&125&32&3780&173299269&4.27.5.41.71.157&3575&664 \cr
 & &2048.3.125.11.17&2805&1024& & &64.125.11.13.83&125&32 \cr
\noalign{\hrule}
}%
}
$$
\eject
\vglue -23 pt
\noindent\hskip 1 in\hbox to 6.5 in{\ 3781 -- 3816 \hfill\fbd 173410853 -- 176955075\frb}
\vskip -9 pt
$$
\vbox{
\nointerlineskip
\halign{\strut
    \vrule \ \ \hfil \frb #\ 
   &\vrule \hfil \ \ \fbb #\frb\ 
   &\vrule \hfil \ \ \frb #\ \hfil
   &\vrule \hfil \ \ \frb #\ 
   &\vrule \hfil \ \ \frb #\ \ \vrule \hskip 2 pt
   &\vrule \ \ \hfil \frb #\ 
   &\vrule \hfil \ \ \fbb #\frb\ 
   &\vrule \hfil \ \ \frb #\ \hfil
   &\vrule \hfil \ \ \frb #\ 
   &\vrule \hfil \ \ \frb #\ \vrule \cr%
\noalign{\hrule}
 & &343.11.19.41.59&15&34& & &9.5.17.31.83.89&7189&5676 \cr
3781&173410853&4.3.5.7.11.17.41.59&817&4062&3799&175182705&8.27.7.11.13.43.79&1303&830 \cr
 & &16.9.19.43.677&677&3096& & &32.5.7.13.83.1303&1303&1456 \cr
\noalign{\hrule}
 & &7.71.113.3089&1149&1940& & &5.7.11.13.17.29.71&801&304 \cr
3782&173481329&8.3.5.71.97.383&4401&2486&3800&175190015&32.9.11.19.29.89&49&38 \cr
 & &32.81.11.113.163&1793&1296& & &128.3.49.361.89&7581&5696 \cr
\noalign{\hrule}
 & &5.11.17.31.53.113&3929&2286& & &3.5.11.13.151.541&2183&522 \cr
3783&173591165&4.9.17.127.3929&5203&1274&3801&175227195&4.27.13.29.37.59&775&8 \cr
 & &16.3.49.121.13.43&1677&4312& & &64.25.31.37&31&5920 \cr
\noalign{\hrule}
 & &3.5.7.11.17.53.167&293&876& & &5.7.23.31.79.89&13689&3454 \cr
3784&173789385&8.9.5.17.73.293&371&286&3802&175458605&4.81.11.169.157&7&20 \cr
 & &32.7.11.13.53.293&293&208& & &32.3.5.7.11.13.157&471&2288 \cr
\noalign{\hrule}
 & &5.7.41.281.431&2211&806& & &25.121.19.43.71&1687&612 \cr
3785&173794285&4.3.11.13.31.41.67&843&3590&3803&175471175&8.9.7.17.71.241&737&950 \cr
 & &16.9.5.281.359&359&72& & &32.3.25.11.17.19.67&201&272 \cr
\noalign{\hrule}
 & &3.125.11.17.37.67&1081&806& & &5.13.17.67.2371&1755&616 \cr
3786&173839875&4.5.13.23.31.47.67&369&34&3804&175536985&16.27.25.7.11.169&1829&646 \cr
 & &16.9.17.23.41.47&1927&552& & &64.3.17.19.31.59&3363&992 \cr
\noalign{\hrule}
 & &27.5.13.23.59.73&209&740& & &27.11.17.43.809&6683&7070 \cr
3787&173852055&8.3.25.11.19.23.37&1169&1606&3805&175639563&4.3.5.7.11.41.101.163&1247&106 \cr
 & &32.7.121.73.167&1169&1936& & &16.5.29.43.53.101&2929&2120 \cr
\noalign{\hrule}
 & &27.19.37.89.103&275&238& & &49.11.17.29.661&2925&4346 \cr
3788&173998827&4.25.7.11.17.89.103&13113&11362&3806&175645547&4.9.25.13.17.41.53&77&128 \cr
 & &16.9.7.13.19.23.31.47&4991&4888& & &1024.3.5.7.11.13.53&2067&2560 \cr
\noalign{\hrule}
 & &7.11.17.37.3593&26013&22420& & &25.11.23.37.751&387&364 \cr
3789&174019769&8.3.5.13.19.23.29.59&425&126&3807&175752775&8.9.25.7.11.13.37.43&73&1148 \cr
 & &32.27.125.7.17.59&1593&2000& & &64.3.49.13.41.73&10731&17056 \cr
\noalign{\hrule}
 & &9.5.23.337.499&341&4& & &9.125.17.29.317&3173&3202 \cr
3790&174048705&8.3.11.31.499&7751&7718&3808&175816125&4.3.19.167.317.1601&325&1276 \cr
 & &32.17.23.227.337&227&272& & &32.25.11.13.19.29.167&2717&2672 \cr
\noalign{\hrule}
 & &3.7.11.13.19.43.71&115&158& & &9.5.11.13.23.29.41&119&548 \cr
3791&174195021&4.5.11.19.23.71.79&1141&360&3809&175977945&8.3.5.7.17.41.137&113&92 \cr
 & &64.9.25.7.23.163&3749&2400& & &64.17.23.113.137&2329&3616 \cr
\noalign{\hrule}
 & &27.125.13.41.97&1891&2086& & &9.11.47.157.241&619&2032 \cr
3792&174490875&4.9.25.7.31.61.149&143&82&3810&176055561&32.47.127.619&2675&3294 \cr
 & &16.7.11.13.31.41.149&1639&1736& & &128.27.25.61.107&6527&4800 \cr
\noalign{\hrule}
 & &3.5.11.17.19.29.113&113&142& & &27.11.361.31.53&4963&14170 \cr
3793&174647715&4.11.19.71.12769&1035&13804&3811&176157531&4.5.7.13.109.709&627&82 \cr
 & &32.9.5.7.17.23.29&161&48& & &16.3.7.11.13.19.41&533&56 \cr
\noalign{\hrule}
 & &3.7.17.23.89.239&923&946& & &11.17.19.179.277&23&300 \cr
3794&174656181&4.11.13.17.43.71.239&9&230&3812&176168399&8.3.25.11.23.179&1921&2196 \cr
 & &16.9.5.11.23.43.71&2343&1720& & &64.27.17.61.113&1647&3616 \cr
\noalign{\hrule}
 & &9.7.121.13.41.43&3145&2144& & &27.5.7.11.361.47&731&214 \cr
3795&174711537&64.3.5.11.17.37.67&553&586&3813&176371965&4.17.361.43.107&8671&6852 \cr
 & &256.5.7.37.79.293&23147&23680& & &32.3.13.23.29.571&7423&10672 \cr
\noalign{\hrule}
 & &5.11.17.31.37.163&171&356& & &11.37249.431&20995&16254 \cr
3796&174808535&8.9.11.19.89.163&23&186&3814&176597509&4.27.5.7.13.17.19.43&193&22 \cr
 & &32.27.23.31.89&2403&368& & &16.3.7.11.13.17.193&119&312 \cr
\noalign{\hrule}
 & &3.5.7.11.13.89.131&837&2278& & &5.11.17.43.53.83&1827&1742 \cr
3797&175059885&4.81.13.17.31.67&1691&820&3815&176861795&4.9.7.11.13.29.53.67&83&2150 \cr
 & &32.5.17.19.41.89&697&304& & &16.3.25.43.67.83&67&120 \cr
\noalign{\hrule}
 & &9.7.11.157.1609&13403&2140& & &9.25.11.19.53.71&2161&1886 \cr
3798&175060809&8.5.13.107.1031&569&462&3816&176955075&4.3.23.41.53.2161&2155&4328 \cr
 & &32.3.5.7.11.13.569&569&1040& & &64.5.23.431.541&12443&13792 \cr
\noalign{\hrule}
}%
}
$$
\eject
\vglue -23 pt
\noindent\hskip 1 in\hbox to 6.5 in{\ 3817 -- 3852 \hfill\fbd 176960781 -- 181563921\frb}
\vskip -9 pt
$$
\vbox{
\nointerlineskip
\halign{\strut
    \vrule \ \ \hfil \frb #\ 
   &\vrule \hfil \ \ \fbb #\frb\ 
   &\vrule \hfil \ \ \frb #\ \hfil
   &\vrule \hfil \ \ \frb #\ 
   &\vrule \hfil \ \ \frb #\ \ \vrule \hskip 2 pt
   &\vrule \ \ \hfil \frb #\ 
   &\vrule \hfil \ \ \fbb #\frb\ 
   &\vrule \hfil \ \ \frb #\ \hfil
   &\vrule \hfil \ \ \frb #\ 
   &\vrule \hfil \ \ \frb #\ \vrule \cr%
\noalign{\hrule}
 & &81.23.43.2209&535&454& & &3.7.11.13.841.71&621&302 \cr
3817&176960781&4.5.2209.107.227&17667&6622&3835&179312133&4.81.7.23.29.151&4331&8710 \cr
 & &16.9.7.11.13.43.151&1001&1208& & &16.5.13.61.67.71&305&536 \cr
\noalign{\hrule}
 & &3.25.23.83.1237&3731&2494& & &5.29.719.1721&501&1220 \cr
3818&177107475&4.7.13.23.29.41.43&1089&100&3836&179422855&8.3.25.29.61.167&4731&5456 \cr
 & &32.9.25.7.121.13&1573&336& & &256.9.11.19.31.83&23157&26752 \cr
\noalign{\hrule}
 & &9.5.7.11.29.41.43&3439&128& & &7.13.17.19.31.197&2645&1098 \cr
3819&177155055&256.3.5.19.181&233&52&3837&179503051&4.9.5.529.31.61&187&342 \cr
 & &2048.13.233&233&13312& & &16.81.11.17.19.61&891&488 \cr
\noalign{\hrule}
 & &9.13.19.173.461&145&28& & &9.5.121.19.37.47&22379&17914 \cr
3820&177290919&8.5.7.19.29.461&297&164&3838&179908245&4.7.169.23.53.139&15&154 \cr
 & &64.27.5.11.29.41&1353&4640& & &16.3.5.49.11.23.53&1127&424 \cr
\noalign{\hrule}
 & &3.5.7.457.3697&14041&4444& & &9.25.13.19.41.79&281&44 \cr
3821&177400545&8.11.19.101.739&4115&4014&3839&180007425&8.3.11.19.41.281&1993&3346 \cr
 & &32.9.5.19.223.823&12711&13168& & &32.7.239.1993&13951&3824 \cr
\noalign{\hrule}
 & &5.11.23.239.587&413&174& & &5.7.17.179.1693&17523&11258 \cr
3822&177470645&4.3.5.7.11.23.29.59&5283&3572&3840&180312965&4.27.11.13.59.433&757&1190 \cr
 & &32.27.19.47.587&513&752& & &16.9.5.7.13.17.757&757&936 \cr
\noalign{\hrule}
 & &9.5.7.169.47.71&271&226& & &7.1331.289.67&279&568 \cr
3823&177645195&4.169.47.113.271&19283&11340&3841&180405071&16.9.11.31.67.71&133&2210 \cr
 & &32.81.5.7.11.1753&1753&1584& & &64.3.5.7.13.17.19&1235&96 \cr
\noalign{\hrule}
 & &25.13.19.47.613&1221&1834& & &7.13.19.29.59.61&83&330 \cr
3824&177907925&4.3.5.7.11.19.37.131&2091&2756&3842&180457459&4.3.5.11.29.61.83&129&190 \cr
 & &32.9.11.13.17.41.53&8109&7216& & &16.9.25.19.43.83&3569&1800 \cr
\noalign{\hrule}
 & &3.125.7.19.43.83&671&754& & &529.31.73.151&1501&3180 \cr
3825&178003875&4.5.7.11.13.29.43.61&5893&342&3843&180766177&8.3.5.19.23.53.79&209&186 \cr
 & &16.9.11.19.71.83&213&88& & &32.9.11.361.31.53&5247&5776 \cr
\noalign{\hrule}
 & &9.13.17.37.41.59&217&550& & &3.5.7.19.31.37.79&1749&1196 \cr
3826&178021467&4.25.7.11.17.31.41&277&318&3844&180772935&8.9.11.13.23.37.53&511&178 \cr
 & &16.3.5.11.31.53.277&15235&13144& & &32.7.11.23.73.89&6497&4048 \cr
\noalign{\hrule}
 & &27.5.11.17.23.307&2527&992& & &5.7.13.17.97.241&11&108 \cr
3827&178254945&64.3.7.11.361.31&23&34&3845&180821095&8.27.5.11.13.241&817&388 \cr
 & &256.7.17.19.23.31&589&896& & &64.9.19.43.97&171&1376 \cr
\noalign{\hrule}
 & &9.5.13.19.61.263&223&1012& & &3.11.23.419.569&75&494 \cr
3828&178317945&8.3.11.23.61.223&367&1036&3846&180953949&4.9.25.11.13.19.23&419&166 \cr
 & &64.7.11.37.367&13579&2464& & &16.5.19.83.419&95&664 \cr
\noalign{\hrule}
 & &13.19.41.67.263&825&46& & &9.7.19.37.61.67&1819&660 \cr
3829&178447867&4.3.25.11.23.263&143&120&3847&181009143&8.27.5.7.11.17.107&1147&1742 \cr
 & &64.9.125.121.13&1125&3872& & &32.11.13.31.37.67&341&208 \cr
\noalign{\hrule}
 & &5.11.13.31.83.97&2913&1652& & &25.7.19.107.509&267&242 \cr
3830&178450415&8.3.7.31.59.971&377&594&3848&181089475&4.3.7.121.19.89.107&5493&4030 \cr
 & &32.81.11.13.29.59&1711&1296& & &16.9.5.11.13.31.1831&23803&24552 \cr
\noalign{\hrule}
 & &27.11.17.23.29.53&395&188& & &27.23.487.599&11&610 \cr
3831&178487199&8.3.5.17.29.47.79&737&742&3849&181153773&4.5.11.61.487&249&238 \cr
 & &32.7.11.47.53.67.79&5293&5264& & &16.3.5.7.17.61.83&2905&8296 \cr
\noalign{\hrule}
 & &25.11.13.23.41.53&441&142& & &3.5.7.13.23.29.199&969&3608 \cr
3832&178674925&4.9.25.49.41.71&283&242&3850&181180545&16.9.5.11.17.19.41&19253&19082 \cr
 & &16.3.7.121.71.283&5467&6792& & &64.7.13.29.47.1481&1481&1504 \cr
\noalign{\hrule}
 & &5.7.13.17.61.379&1369&1284& & &81.7.11.17.29.59&247&166 \cr
3833&178825465&8.3.13.1369.61.107&1583&2376&3851&181415619&4.11.13.17.19.29.83&49&270 \cr
 & &128.81.11.37.1583&58571&57024& & &16.27.5.49.19.83&415&1064 \cr
\noalign{\hrule}
 & &81.7.11.23.29.43&305&262& & &9.7.11.37.73.97&233&1300 \cr
3834&178883397&4.5.11.23.29.61.131&387&1054&3852&181563921&8.3.25.13.37.233&73&112 \cr
 & &16.9.5.17.31.43.61&1037&1240& & &256.5.7.73.233&233&640 \cr
\noalign{\hrule}
}%
}
$$
\eject
\vglue -23 pt
\noindent\hskip 1 in\hbox to 6.5 in{\ 3853 -- 3888 \hfill\fbd 181695745 -- 184191651\frb}
\vskip -9 pt
$$
\vbox{
\nointerlineskip
\halign{\strut
    \vrule \ \ \hfil \frb #\ 
   &\vrule \hfil \ \ \fbb #\frb\ 
   &\vrule \hfil \ \ \frb #\ \hfil
   &\vrule \hfil \ \ \frb #\ 
   &\vrule \hfil \ \ \frb #\ \ \vrule \hskip 2 pt
   &\vrule \ \ \hfil \frb #\ 
   &\vrule \hfil \ \ \fbb #\frb\ 
   &\vrule \hfil \ \ \frb #\ \hfil
   &\vrule \hfil \ \ \frb #\ 
   &\vrule \hfil \ \ \frb #\ \vrule \cr%
\noalign{\hrule}
 & &5.7.11.289.23.71&14687&5832& & &3.5.23.37.113.127&91&22 \cr
3853&181695745&16.729.19.773&2695&1922&3871&183190515&4.5.7.11.13.37.127&1081&954 \cr
 & &64.3.5.49.11.961&961&672& & &16.9.7.13.23.47.53&2067&2632 \cr
\noalign{\hrule}
 & &3.25.19.29.53.83&671&754& & &9.17.61.67.293&1855&1562 \cr
3854&181788675&4.11.13.841.53.61&765&76&3872&183216123&4.3.5.7.11.53.61.71&1139&874 \cr
 & &32.9.5.11.17.19.61&561&976& & &16.7.17.19.23.67.71&1349&1288 \cr
\noalign{\hrule}
 & &81.5.17.29.911&449&44& & &5.49.19.961.41&1485&524 \cr
3855&181894815&8.11.449.911&231&680&3873&183411655&8.27.25.11.19.131&41&434 \cr
 & &128.3.5.7.121.17&847&64& & &32.9.7.11.31.41&11&144 \cr
\noalign{\hrule}
 & &9.125.7.121.191&2977&1258& & &27.25.7.13.29.103&409&374 \cr
3856&181999125&4.25.13.17.37.229&1701&1276&3874&183476475&4.5.11.13.17.103.409&1119&14 \cr
 & &32.243.7.11.29.37&783&592& & &16.3.7.373.409&409&2984 \cr
\noalign{\hrule}
 & &53.337.10193&5123&5070& & &3.121.13.19.23.89&89&158 \cr
3857&182057173&4.3.5.169.47.109.337&371&4752&3875&183536067&4.121.79.7921&819&8740 \cr
 & &128.81.5.7.11.13.53&5265&4928& & &32.9.5.7.13.19.23&105&16 \cr
\noalign{\hrule}
 & &5.7.121.19.31.73&151&360& & &3.5.41.43.53.131&1209&554 \cr
3858&182092295&16.9.25.11.31.151&91&184&3876&183607635&4.9.13.31.53.277&1591&902 \cr
 & &256.3.7.13.23.151&5889&2944& & &16.11.31.37.41.43&407&248 \cr
\noalign{\hrule}
 & &25.11.47.73.193&14573&5502& & &9.5.11.13.17.23.73&7687&8446 \cr
3859&182100325&4.3.7.13.19.59.131&73&60&3877&183674205&4.3.5.41.103.7687&973&22088 \cr
 & &32.9.5.59.73.131&1179&944& & &64.7.11.139.251&1757&4448 \cr
\noalign{\hrule}
 & &25.7.13.23.3481&1683&1798& & &9.5.11.13.17.23.73&241&254 \cr
3860&182143325&4.9.5.7.11.13.17.29.31&59&214&3878&183674205&4.17.23.73.127.241&29575&1032 \cr
 & &16.3.11.17.29.59.107&3531&3944& & &64.3.25.7.169.43&559&1120 \cr
\noalign{\hrule}
 & &3.5.7.19.241.379&2757&4444& & &11.13.19.53.1277&15921&8342 \cr
3861&182221305&8.9.5.11.101.919&427&482&3879&183889277&4.9.29.43.61.97&2915&3002 \cr
 & &32.7.61.241.919&919&976& & &16.3.5.11.19.43.53.79&645&632 \cr
\noalign{\hrule}
 & &5.7.17.31.41.241&741&946& & &11.169.19.41.127&965&432 \cr
3862&182255045&4.3.11.13.17.19.31.43&27&500&3880&183916447&32.27.5.13.19.193&907&328 \cr
 & &32.81.125.13.19&6175&1296& & &512.9.41.907&907&2304 \cr
\noalign{\hrule}
 & &9.5.121.19.41.43&1255&3706& & &81.5.49.13.23.31&1661&1444 \cr
3863&182391165&4.3.25.17.109.251&589&164&3881&183943305&8.3.7.11.13.361.151&899&158 \cr
 & &32.19.31.41.109&109&496& & &32.11.19.29.31.79&2291&3344 \cr
\noalign{\hrule}
 & &81.13.17.23.443&2995&2972& & &25.13.19.83.359&1161&914 \cr
3864&182393289&8.3.5.443.599.743&209&2006&3882&183996475&4.27.43.359.457&865&506 \cr
 & &32.11.17.19.59.743&12331&11888& & &16.9.5.11.23.43.173&10879&12456 \cr
\noalign{\hrule}
 & &3.7.11.13.31.37.53&177&230& & &9.5.19.31.53.131&2057&2004 \cr
3865&182555373&4.9.5.7.13.23.31.59&53&766&3883&184024215&8.27.5.121.17.19.167&2201&364 \cr
 & &16.5.53.59.383&1915&472& & &64.7.11.13.17.31.71&6461&5984 \cr
\noalign{\hrule}
 & &25.19.1681.229&18837&23188& & &27.11.23.29.929&10541&11470 \cr
3866&182850775&8.9.7.11.13.17.23.31&205&458&3884&184033971&4.9.5.31.37.83.127&715&1858 \cr
 & &32.3.5.7.31.41.229&93&112& & &16.25.11.13.37.929&325&296 \cr
\noalign{\hrule}
 & &9.5.7.43.59.229&2739&5276& & &11.169.97.1021&1141&120 \cr
3867&183006495&8.27.11.83.1319&511&808&3885&184109783&16.3.5.7.11.13.163&419&582 \cr
 & &128.7.73.83.101&7373&5312& & &64.9.5.97.419&419&1440 \cr
\noalign{\hrule}
 & &7.11.13.17.31.347&5307&5450& & &5.31.929.1279&1683&2962 \cr
3868&183051869&4.3.25.7.17.29.61.109&11&1026&3886&184169605&4.9.11.17.31.1481&12077&13100 \cr
 & &16.81.5.11.19.109&1539&4360& & &32.3.25.13.131.929&655&624 \cr
\noalign{\hrule}
 & &7.11.19.47.2663&81&128& & &5.11.23.41.53.67&639&304 \cr
3869&183110543&256.81.7.2663&1363&1300&3887&184172615&32.9.11.19.53.71&1675&668 \cr
 & &2048.9.25.13.29.47&9425&9216& & &256.3.25.67.167&835&384 \cr
\noalign{\hrule}
 & &3.13.1369.47.73&1045&1656& & &81.7.17.97.197&9595&9514 \cr
3870&183184521&16.27.5.11.19.23.37&329&292&3888&184191651&4.5.7.17.19.67.71.101&237&902 \cr
 & &128.5.7.11.19.47.73&665&704& & &16.3.11.41.71.79.101&45551&44872 \cr
\noalign{\hrule}
}%
}
$$
\eject
\vglue -23 pt
\noindent\hskip 1 in\hbox to 6.5 in{\ 3889 -- 3924 \hfill\fbd 184320045 -- 188724891\frb}
\vskip -9 pt
$$
\vbox{
\nointerlineskip
\halign{\strut
    \vrule \ \ \hfil \frb #\ 
   &\vrule \hfil \ \ \fbb #\frb\ 
   &\vrule \hfil \ \ \frb #\ \hfil
   &\vrule \hfil \ \ \frb #\ 
   &\vrule \hfil \ \ \frb #\ \ \vrule \hskip 2 pt
   &\vrule \ \ \hfil \frb #\ 
   &\vrule \hfil \ \ \fbb #\frb\ 
   &\vrule \hfil \ \ \frb #\ \hfil
   &\vrule \hfil \ \ \frb #\ 
   &\vrule \hfil \ \ \frb #\ \vrule \cr%
\noalign{\hrule}
 & &9.5.7.13.19.23.103&509&418& & &3.17.361.29.349&3069&3562 \cr
3889&184320045&4.5.11.361.23.509&15781&4074&3907&186337731&4.27.11.13.19.31.137&4945&698 \cr
 & &16.3.7.43.97.367&4171&2936& & &16.5.13.23.43.349&559&920 \cr
\noalign{\hrule}
 & &3.11.29.37.41.127&1083&2600& & &9.5.13.31.43.239&3553&446 \cr
3890&184374663&16.9.25.11.13.361&127&82&3908&186373395&4.3.5.11.17.19.223&1679&1456 \cr
 & &64.5.13.19.41.127&247&160& & &128.7.13.17.23.73&2737&4672 \cr
\noalign{\hrule}
 & &11.23.29.41.613&16235&8898& & &9.5.7.43.47.293&29&358 \cr
3891&184400821&4.3.5.17.191.1483&1219&264&3909&186528195&4.5.29.179.293&301&594 \cr
 & &64.9.11.17.23.53&901&288& & &16.27.7.11.29.43&29&264 \cr
\noalign{\hrule}
 & &3.5.7.11.23.53.131&20557&24638& & &3.5.19.29.67.337&309&242 \cr
3892&184440795&4.61.97.127.337&3825&3922&3910&186615435&4.9.5.121.103.337&4171&464 \cr
 & &16.9.25.17.37.53.337&5055&5032& & &128.11.29.43.97&1067&2752 \cr
\noalign{\hrule}
 & &81.5.37.109.113&2299&1734& & &243.125.11.13.43&56657&58282 \cr
3893&184570245&4.243.121.289.19&113&130&3911&186775875&4.7.23.53.181.1069&75&1144 \cr
 & &16.5.121.13.17.19.113&2057&1976& & &64.3.25.7.11.13.181&181&224 \cr
\noalign{\hrule}
 & &5.11.23.37.59.67&959&2994& & &9.5.47.241.367&167&2002 \cr
3894&185020165&4.3.7.23.137.499&171&3322&3912&187065405&4.7.11.13.47.167&7857&7340 \cr
 & &16.27.11.19.151&151&4104& & &32.81.5.97.367&97&144 \cr
\noalign{\hrule}
 & &3.5.19.29.73.307&291&1826& & &121.529.37.79&3147&1330 \cr
3895&185226915&4.9.11.19.83.97&899&982&3913&187098307&4.3.5.7.19.23.1049&605&444 \cr
 & &16.29.31.97.491&3007&3928& & &32.9.25.121.19.37&225&304 \cr
\noalign{\hrule}
 & &13.113.257.491&26071&29412& & &3.5.49.37.71.97&87&158 \cr
3896&185368703&8.9.19.841.31.43&385&514&3914&187291965&4.9.29.37.79.97&1925&1664 \cr
 & &32.3.5.7.11.19.29.257&3135&3248& & &1024.25.7.11.13.79&5135&5632 \cr
\noalign{\hrule}
 & &5.11.361.47.199&1541&648& & &9.25.7.11.29.373&19&1846 \cr
3897&185703815&16.81.5.19.23.67&1997&188&3915&187404525&4.5.11.13.19.71&3&8 \cr
 & &128.3.47.1997&1997&192& & &64.3.13.19.71&71&7904 \cr
\noalign{\hrule}
 & &27.5.11.19.29.227&343&208& & &13.23.43.61.239&891&98 \cr
3898&185739345&32.343.11.13.227&981&3478&3916&187442203&4.81.49.11.239&1159&920 \cr
 & &128.9.37.47.109&5123&2368& & &64.3.5.7.19.23.61&399&160 \cr
\noalign{\hrule}
 & &9.11.13.19.71.107&667&724& & &9.25.7.11.79.137&5549&5274 \cr
3899&185769441&8.3.11.23.29.71.181&1391&610&3917&187508475&4.81.7.31.179.293&137&430 \cr
 & &32.5.13.61.107.181&905&976& & &16.5.31.43.137.179&1333&1432 \cr
\noalign{\hrule}
 & &27.7.17.151.383&97&286& & &11.19.31.103.281&89&120 \cr
3900&185817429&4.11.13.17.97.151&35&186&3918&187521697&16.3.5.89.103.281&117&398 \cr
 & &16.3.5.7.11.31.97&3007&440& & &64.27.13.89.199&17711&11232 \cr
\noalign{\hrule}
 & &9.5.121.13.37.71&6749&6386& & &9.5.7.11.13.23.181&3109&2204 \cr
3901&185952195&4.3.13.17.31.103.397&1265&74&3919&187522335&8.3.13.19.29.3109&1925&1184 \cr
 & &16.5.11.17.23.31.37&391&248& & &512.25.7.11.29.37&1073&1280 \cr
\noalign{\hrule}
 & &81.13.17.19.547&77&94& & &125.7.11.101.193&493&618 \cr
3902&186045093&4.9.7.11.13.47.547&323&4600&3920&187620125&4.3.7.17.29.103.193&1551&200 \cr
 & &64.25.11.17.19.23&253&800& & &64.9.25.11.29.47&423&928 \cr
\noalign{\hrule}
 & &25.13.23.37.673&1727&13752& & &11.13.19.37.1873&20439&15148 \cr
3903&186134975&16.9.11.157.191&577&1150&3921&188290817&8.27.7.541.757&865&1406 \cr
 & &64.3.25.23.577&577&96& & &32.9.5.7.19.37.173&865&1008 \cr
\noalign{\hrule}
 & &3.5.289.29.1481&987&494& & &25.11.19.43.839&14981&5994 \cr
3904&186183915&4.9.5.7.13.17.19.47&467&638&3922&188502325&4.81.37.71.211&2945&4862 \cr
 & &16.7.11.29.47.467&3269&4136& & &16.3.5.11.13.17.19.31&527&312 \cr
\noalign{\hrule}
 & &11.13.17.19.29.139&1065&742& & &9.5.7.11.13.53.79&41&436 \cr
3905&186187859&4.3.5.7.11.29.53.71&1159&378&3923&188603415&8.7.11.13.41.109&75&68 \cr
 & &16.81.5.49.19.61&3969&2440& & &64.3.25.17.41.109&4469&2720 \cr
\noalign{\hrule}
 & &9.5.7.19.841.37&335&506& & &3.19.53.179.349&5203&14690 \cr
3906&186235245&4.25.7.11.23.37.67&3111&3364&3924&188724891&4.5.121.13.43.113&901&342 \cr
 & &32.3.17.841.61.67&1037&1072& & &16.9.5.11.17.19.53&85&264 \cr
\noalign{\hrule}
}%
}
$$
\eject
\vglue -23 pt
\noindent\hskip 1 in\hbox to 6.5 in{\ 3925 -- 3960 \hfill\fbd 188727795 -- 193615965\frb}
\vskip -9 pt
$$
\vbox{
\nointerlineskip
\halign{\strut
    \vrule \ \ \hfil \frb #\ 
   &\vrule \hfil \ \ \fbb #\frb\ 
   &\vrule \hfil \ \ \frb #\ \hfil
   &\vrule \hfil \ \ \frb #\ 
   &\vrule \hfil \ \ \frb #\ \ \vrule \hskip 2 pt
   &\vrule \ \ \hfil \frb #\ 
   &\vrule \hfil \ \ \fbb #\frb\ 
   &\vrule \hfil \ \ \frb #\ \hfil
   &\vrule \hfil \ \ \frb #\ 
   &\vrule \hfil \ \ \frb #\ \vrule \cr%
\noalign{\hrule}
 & &9.5.17.29.47.181&3601&4906& & &3.13.107.109.421&157&170 \cr
3925&188727795&4.11.13.17.223.277&905&3804&3943&191494797&4.5.17.107.157.421&143&1962 \cr
 & &32.3.5.11.181.317&317&176& & &16.9.11.13.109.157&157&264 \cr
\noalign{\hrule}
 & &9.5.7.11.529.103&1377&992& & &3.5.49.11.137.173&47&912 \cr
3926&188797455&64.729.17.23.31&721&8&3944&191622585&32.9.7.11.19.47&173&250 \cr
 & &1024.7.17.103&17&512& & &128.125.19.173&475&64 \cr
\noalign{\hrule}
 & &7.17.67.151.157&5137&4980& & &25.13.59.73.137&79&216 \cr
3927&189015911&8.3.5.7.11.17.83.467&97&90&3945&191769175&16.27.5.13.73.79&77&142 \cr
 & &32.27.25.83.97.467&201275&201744& & &64.9.7.11.71.79&7821&15904 \cr
\noalign{\hrule}
 & &27.25.11.71.359&31&244& & &9.25.7.13.17.19.29&2449&4334 \cr
3928&189255825&8.9.31.61.359&319&40&3946&191789325&4.3.5.11.31.79.197&493&98 \cr
 & &128.5.11.29.61&61&1856& & &16.49.11.17.29.31&77&248 \cr
\noalign{\hrule}
 & &7.17.19.31.37.73&23815&27504& & &9.49.11.13.17.179&817&1510 \cr
3929&189315791&32.9.5.11.191.433&949&3814&3947&191900709&4.5.7.17.19.43.151&537&520 \cr
 & &128.3.13.73.1907&1907&2496& & &64.3.25.13.19.43.179&817&800 \cr
\noalign{\hrule}
 & &125.13.23.37.137&303&178& & &5.11.17.19.79.137&93&230 \cr
3930&189453875&4.3.23.89.101.137&5895&7942&3948&192270595&4.3.25.11.23.31.79&817&1158 \cr
 & &16.27.5.11.361.131&9747&11528& & &16.9.19.23.43.193&4439&3096 \cr
\noalign{\hrule}
 & &11.19.71.12769&1035&13804& & &3.5.19.43.113.139&1573&878 \cr
3931&189479191&8.9.5.7.17.23.29&113&142&3949&192489285&4.121.13.113.439&3475&2232 \cr
 & &32.3.7.23.71.113&161&48& & &64.9.25.11.31.139&465&352 \cr
\noalign{\hrule}
 & &81.7.13.47.547&323&4600& & &25.7.149.7393&24013&27738 \cr
3932&189501039&16.9.25.17.19.23&77&94&3950&192772475&4.9.11.23.37.59.67&1205&1274 \cr
 & &64.25.7.11.23.47&253&800& & &16.3.5.49.11.13.59.241&18557&18408 \cr
\noalign{\hrule}
 & &3.11.29.37.53.101&247&160& & &27.11.19.47.727&4255&3742 \cr
3933&189544377&64.5.13.19.53.101&153&1160&3951&192815667&4.5.23.37.47.1871&1805&66 \cr
 & &1024.9.25.17.29&425&1536& & &16.3.25.11.361.23&575&152 \cr
\noalign{\hrule}
 & &27.5.7.11.13.23.61&58381&56734& & &9.11.17.19.37.163&335&368 \cr
3934&189594405&4.19.79.739.1493&1497&4&3952&192853287&32.3.5.17.23.67.163&325&814 \cr
 & &32.3.499.739&739&7984& & &128.125.11.13.23.37&1625&1472 \cr
\noalign{\hrule}
 & &3.5.17.47.83.191&1001&2246& & &5.11.13.29.71.131&621&302 \cr
3935&189998205&4.7.11.13.47.1123&1577&2700&3953&192856235&4.27.5.23.131.151&9373&10408 \cr
 & &32.27.25.11.19.83&495&304& & &64.3.7.13.103.1301&9107&9888 \cr
\noalign{\hrule}
 & &9.7.13.29.53.151&359&1178& & &81.19.29.61.71&19765&19214 \cr
3936&190079253&4.19.31.151.359&5489&5640&3954&193296861&4.9.5.13.59.67.739&1349&5302 \cr
 & &64.3.5.11.19.47.499&27445&28576& & &16.5.11.13.19.71.241&1205&1144 \cr
\noalign{\hrule}
 & &25.11.41.47.359&1387&2562& & &3.49.169.31.251&109&60 \cr
3937&190243075&4.3.7.19.41.61.73&9447&8060&3955&193303383&8.9.5.31.109.251&19943&18962 \cr
 & &32.9.5.13.31.47.67&2077&1872& & &32.49.11.19.37.499&7733&7984 \cr
\noalign{\hrule}
 & &9.17.37.151.223&91&242& & &9.11.13.97.1549&5267&14870 \cr
3938&190622853&4.7.121.13.17.223&795&2104&3956&193375611&4.5.23.229.1487&629&858 \cr
 & &64.3.5.11.53.263&2893&8480& & &16.3.5.11.13.17.23.37&629&920 \cr
\noalign{\hrule}
 & &9.49.11.19.2069&115&94& & &27.7.13.17.41.113&1195&404 \cr
3939&190697661&4.3.5.7.23.47.2069&793&1276&3957&193515777&8.9.5.17.101.239&1577&2486 \cr
 & &32.5.11.13.29.47.61&8845&9776& & &32.5.11.19.83.113&1577&880 \cr
\noalign{\hrule}
 & &3.29.101.103.211&91&120& & &9.5.11.13.17.29.61&271&766 \cr
3940&190967871&16.9.5.7.13.101.103&817&110&3958&193519755&4.13.29.271.383&7315&3792 \cr
 & &64.25.11.13.19.43&6175&15136& & &128.3.5.7.11.19.79&553&1216 \cr
\noalign{\hrule}
 & &27.25.11.289.89&3629&3596& & &9.25.11.17.43.107&1043&1150 \cr
3941&190978425&8.9.19.29.31.89.191&431&2150&3959&193587075&4.3.625.7.11.23.149&197&428 \cr
 & &32.25.19.31.43.431&13361&13072& & &32.23.107.149.197&3427&3152 \cr
\noalign{\hrule}
 & &81.17.19.67.109&9889&8750& & &9.5.79.107.509&847&338 \cr
3942&191068389&4.9.625.7.11.29.31&109&46&3960&193615965&4.3.7.121.169.107&245&76 \cr
 & &16.125.11.23.29.109&2875&2552& & &32.5.343.121.19&6517&1936 \cr
\noalign{\hrule}
}%
}
$$
\eject
\vglue -23 pt
\noindent\hskip 1 in\hbox to 6.5 in{\ 3961 -- 3996 \hfill\fbd 193995245 -- 199863345\frb}
\vskip -9 pt
$$
\vbox{
\nointerlineskip
\halign{\strut
    \vrule \ \ \hfil \frb #\ 
   &\vrule \hfil \ \ \fbb #\frb\ 
   &\vrule \hfil \ \ \frb #\ \hfil
   &\vrule \hfil \ \ \frb #\ 
   &\vrule \hfil \ \ \frb #\ \ \vrule \hskip 2 pt
   &\vrule \ \ \hfil \frb #\ 
   &\vrule \hfil \ \ \fbb #\frb\ 
   &\vrule \hfil \ \ \frb #\ \hfil
   &\vrule \hfil \ \ \frb #\ 
   &\vrule \hfil \ \ \frb #\ \vrule \cr%
\noalign{\hrule}
 & &5.17.59.101.383&693&310& & &27.5.17.23.37.101&695&304 \cr
3961&193995245&4.9.25.7.11.31.101&59&766&3979&197257545&32.25.19.101.139&187&288 \cr
 & &16.3.31.59.383&93&8& & &2048.9.11.17.139&1529&1024 \cr
\noalign{\hrule}
 & &11.37.367.1301&75335&74034& & &25.11.19.37.1021&97173&96152 \cr
3962&194329069&4.81.5.13.19.61.457&15047&6364&3980&197384825&16.27.7.17.59.61.101&739&2860 \cr
 & &32.3.5.37.41.43.367&645&656& & &128.9.5.11.13.17.739&9607&9792 \cr
\noalign{\hrule}
 & &9.11.13.29.41.127&575&614& & &9.11.73.109.251&1729&530 \cr
3963&194340861&4.3.25.11.23.127.307&1435&5626&3981&197723493&4.5.7.13.19.53.73&29&978 \cr
 & &16.125.7.29.41.97&875&776& & &16.3.5.7.29.163&163&8120 \cr
\noalign{\hrule}
 & &3.49.31.179.239&19063&19780& & &3.11.47.173.739&299&440 \cr
3964&194953017&8.5.7.11.23.43.1733&1017&716&3982&198290697&16.5.121.13.23.173&1305&1478 \cr
 & &64.9.5.11.23.113.179&3729&3680& & &64.9.25.13.29.739&975&928 \cr
\noalign{\hrule}
 & &121.23.163.431&4905&14818& & &71.18769.149&4095&14674 \cr
3965&195514099&4.9.5.31.109.239&283&44&3983&198557251&4.9.5.7.11.13.23.29&137&298 \cr
 & &32.3.5.11.31.283&4245&496& & &16.3.11.13.137.149&33&104 \cr
\noalign{\hrule}
 & &9.25.23.29.1303&869&434& & &17.79.353.419&881&462 \cr
3966&195547725&4.3.5.7.11.23.31.79&319&164&3984&198639101&4.3.7.11.353.881&795&1676 \cr
 & &32.121.29.41.79&3239&1936& & &32.9.5.11.53.419&495&848 \cr
\noalign{\hrule}
 & &9.11.13.361.421&175&186& & &5.49.11.13.53.107&333&5338 \cr
3967&195599547&4.27.25.7.13.31.421&10703&178&3985&198683485&4.9.7.17.37.157&397&380 \cr
 & &16.49.11.89.139&4361&1112& & &32.3.5.19.157.397&8949&6352 \cr
\noalign{\hrule}
 & &9.1681.67.193&505&2242& & &27.11.31.113.191&6305&4204 \cr
3968&195633099&4.5.19.41.59.101&957&1462&3986&198714681&8.9.5.13.97.1051&467&584 \cr
 & &16.3.11.17.19.29.43&8987&3944& & &128.5.73.97.467&34091&31040 \cr
\noalign{\hrule}
 & &5.7.11.13.19.29.71&14017&822& & &3.5.13.47.109.199&125&16 \cr
3969&195800605&4.3.107.131.137&65&66&3987&198798015&32.625.13.199&981&1606 \cr
 & &16.9.5.11.13.107.137&963&1096& & &128.9.11.73.109&803&192 \cr
\noalign{\hrule}
 & &3.5.7.11.223.761&715&46& & &27.5.7.43.59.83&49&34 \cr
3970&196006965&4.25.7.121.13.23&1467&1558&3988&198989595&4.9.343.17.43.59&275&2812 \cr
 & &16.9.19.23.41.163&6683&10488& & &32.25.11.17.19.37&1615&6512 \cr
\noalign{\hrule}
 & &3.5.11.13.23.29.137&103&34& & &9.5.289.61.251&147&142 \cr
3971&196007955&4.5.11.13.17.29.103&2331&1816&3989&199119555&4.27.49.61.71.251&14365&946 \cr
 & &64.9.7.17.37.227&11577&8288& & &16.5.7.11.169.17.43&1859&2408 \cr
\noalign{\hrule}
 & &3.7.11.19.59.757&375&382& & &25.13.31.53.373&201&574 \cr
3972&196025907&4.9.125.11.19.59.191&23467&4558&3990&199172675&4.3.7.13.41.53.67&33&20 \cr
 & &16.5.31.43.53.757&1643&1720& & &32.9.5.7.11.41.67&4059&7504 \cr
\noalign{\hrule}
 & &243.25.13.47.53&1853&638& & &5.7.23.31.61.131&87&218 \cr
3973&196726725&4.5.11.13.17.29.109&1223&1938&3991&199415405&4.3.7.23.29.31.109&1893&4400 \cr
 & &16.3.289.19.1223&5491&9784& & &128.9.25.11.631&6941&2880 \cr
\noalign{\hrule}
 & &3.49.13.29.53.67&275&1146& & &11.41.467.947&4365&14782 \cr
3974&196792869&4.9.25.11.53.191&83&182&3992&199454299&4.9.5.19.97.389&51&1894 \cr
 & &16.5.7.13.83.191&415&1528& & &16.27.17.947&459&8 \cr
\noalign{\hrule}
 & &5.11.169.59.359&1591&2358& & &25.49.11.13.17.67&23&198 \cr
3975&196877395&4.9.5.13.37.43.131&511&1166&3993&199524325&4.9.7.121.23.67&2041&500 \cr
 & &16.3.7.11.37.53.73&5883&4088& & &32.3.125.13.157&785&48 \cr
\noalign{\hrule}
 & &9.25.11.13.29.211&161&216& & &3.25.61.149.293&2009&1716 \cr
3976&196878825&16.243.5.7.23.211&649&406&3994&199730775&8.9.49.11.13.41.61&745&74 \cr
 & &64.49.11.23.29.59&1357&1568& & &32.5.7.37.41.149&287&592 \cr
\noalign{\hrule}
 & &27.7.17.29.2113&665&1448& & &5.7.11.19.59.463&1689&1556 \cr
3977&196883001&16.5.49.17.19.181&22011&22334&3995&199823855&8.3.389.463.563&513&50 \cr
 & &64.3.11.13.23.29.859&9449&9568& & &32.81.25.19.389&1945&1296 \cr
\noalign{\hrule}
 & &9.169.29.41.109&649&540& & &3.5.11.67.101.179&353&152 \cr
3978&197123121&8.243.5.11.169.59&3649&6322&3996&199863345&16.11.19.179.353&1161&808 \cr
 & &32.5.29.41.89.109&89&80& & &256.27.19.43.101&817&1152 \cr
\noalign{\hrule}
}%
}
$$
\eject
\vglue -23 pt
\noindent\hskip 1 in\hbox to 6.5 in{\ 3997 -- 4032 \hfill\fbd 200024445 -- 204405135\frb}
\vskip -9 pt
$$
\vbox{
\nointerlineskip
\halign{\strut
    \vrule \ \ \hfil \frb #\ 
   &\vrule \hfil \ \ \fbb #\frb\ 
   &\vrule \hfil \ \ \frb #\ \hfil
   &\vrule \hfil \ \ \frb #\ 
   &\vrule \hfil \ \ \frb #\ \ \vrule \hskip 2 pt
   &\vrule \ \ \hfil \frb #\ 
   &\vrule \hfil \ \ \fbb #\frb\ 
   &\vrule \hfil \ \ \frb #\ \hfil
   &\vrule \hfil \ \ \frb #\ 
   &\vrule \hfil \ \ \frb #\ \vrule \cr%
\noalign{\hrule}
 & &3.5.23.41.79.179&561&382& & &27.7.11.29.3359&31825&28466 \cr
3997&200024445&4.9.5.11.17.79.191&3401&154&4015&202517469&4.25.19.43.67.331&613&2268 \cr
 & &16.7.121.19.179&121&1064& & &32.81.5.7.19.613&1839&1520 \cr
\noalign{\hrule}
 & &9.7.121.13.43.47&3589&14600& & &5.13.17.23.79.101&859&858 \cr
3998&200279079&16.25.37.73.97&141&44&4016&202786285&4.3.5.11.169.23.79.859&122823&24038 \cr
 & &128.3.5.11.47.73&365&64& & &16.81.7.17.101.4549&4549&4536 \cr
\noalign{\hrule}
 & &3.11.13.19.103.239&4067&3820& & &81.25.7.41.349&443&418 \cr
3999&200653167&8.5.49.83.103.191&239&342&4017&202830075&4.27.11.19.349.443&2665&9296 \cr
 & &32.9.5.7.19.191.239&573&560& & &128.5.7.11.13.41.83&913&832 \cr
\noalign{\hrule}
 & &3.11.169.17.29.73&18977&12626& & &3.25.11.289.23.37&133&422 \cr
4000&200710653&4.7.59.107.2711&981&1730&4018&202899675&4.5.7.11.19.23.211&603&442 \cr
 & &16.9.5.59.109.173&10207&13080& & &16.9.13.17.67.211&2613&1688 \cr
\noalign{\hrule}
 & &5.11.13.17.83.199&989&3576& & &5.7.17.29.61.193&11659&18414 \cr
4001&200764135&16.3.17.23.43.149&291&440&4019&203143115&4.27.11.31.89.131&205&74 \cr
 & &256.9.5.11.23.97&2231&1152& & &16.3.5.11.37.41.89&4551&7832 \cr
\noalign{\hrule}
 & &9.13.17.29.3481&51469&49480& & &3.59.71.103.157&14689&3542 \cr
4002&200787561&16.5.11.1237.4679&1721&2958&4020&203220957&4.7.11.23.37.397&325&72 \cr
 & &64.3.5.11.17.29.1721&1721&1760& & &64.9.25.7.13.37&975&8288 \cr
\noalign{\hrule}
 & &25.11.13.53.1061&7029&6764& & &3.5.121.31.3613&1295&2318 \cr
4003&201032975&8.9.5.121.19.71.89&233&1582&4021&203285445&4.25.7.11.19.37.61&423&248 \cr
 & &32.3.7.89.113.233&26329&29904& & &64.9.19.31.37.47&2109&1504 \cr
\noalign{\hrule}
 & &27.5.17.79.1109&427&682& & &5.11.17.23.9463&4927&4536 \cr
4004&201067245&4.9.7.11.31.61.79&1499&950&4022&203501815&16.81.5.7.11.13.379&4759&5474 \cr
 & &16.25.7.11.19.1499&10493&8360& & &64.3.49.17.23.4759&4759&4704 \cr
\noalign{\hrule}
 & &9.7.11.17.19.29.31&1865&1586& & &27.41.251.733&105&146 \cr
4005&201231261&4.5.11.13.19.61.373&7905&818&4023&203669181&4.81.5.7.73.733&391&5522 \cr
 & &16.3.25.17.31.409&409&200& & &16.5.11.17.23.251&187&920 \cr
\noalign{\hrule}
 & &81.17.37.59.67&1727&1690& & &3.5.121.13.53.163&2147&4266 \cr
4006&201401397&4.27.5.11.169.59.157&67&1660&4024&203837205&4.81.5.19.79.113&163&242 \cr
 & &32.25.169.67.83&2075&2704& & &16.121.19.113.163&113&152 \cr
\noalign{\hrule}
 & &27.25.29.41.251&25597&25858& & &13.19.61.83.163&1507&1590 \cr
4007&201446325&4.3.5.7.11.13.179.1847&6119&3116&4025&203841443&4.3.5.11.13.53.61.137&389&7968 \cr
 & &32.19.29.41.179.211&3401&3376& & &256.9.5.83.389&1945&1152 \cr
\noalign{\hrule}
 & &3.11.43.173.821&151&670& & &9.11.19.29.37.101&559&350 \cr
4008&201544827&4.5.11.43.67.151&15921&15770&4026&203849613&4.25.7.13.29.37.43&57&202 \cr
 & &16.9.25.19.29.61.83&41325&40504& & &16.3.5.13.19.43.101&215&104 \cr
\noalign{\hrule}
 & &25.11.19.47.821&387&434& & &9.11.13.23.61.113&301&370 \cr
4009&201617075&4.9.25.7.11.19.31.43&821&4&4027&204039693&4.3.5.7.13.37.43.113&3151&1708 \cr
 & &32.3.7.31.821&21&496& & &32.5.49.23.61.137&685&784 \cr
\noalign{\hrule}
 & &13.61.127.2003&605&1398& & &3.7.31.53.61.97&5423&494 \cr
4010&201724133&4.3.5.121.127.233&1281&116&4028&204154251&4.7.11.13.17.19.29&97&90 \cr
 & &32.9.7.11.29.61&1827&176& & &16.9.5.13.19.29.97&1131&760 \cr
\noalign{\hrule}
 & &9.5.29.359.431&539&1616& & &3.5.49.11.13.29.67&1957&2828 \cr
4011&201921345&32.3.49.11.29.101&559&862&4029&204219015&8.343.19.101.103&1943&8460 \cr
 & &128.11.13.43.431&473&832& & &64.9.5.29.47.67&47&96 \cr
\noalign{\hrule}
 & &81.5.11.289.157&893&3562& & &9.5.7.17.37.1031&253&5408 \cr
4012&202136715&4.13.17.19.47.137&297&314&4030&204277185&64.7.11.169.23&87&74 \cr
 & &16.27.11.19.137.157&137&152& & &256.3.11.13.29.37&377&1408 \cr
\noalign{\hrule}
 & &3.17.41.109.887&95&792& & &5.11.29.37.3463&1891&1572 \cr
4013&202164153&16.27.5.11.19.109&887&1184&4031&204368945&8.3.5.31.37.61.131&801&1456 \cr
 & &1024.5.37.887&185&512& & &256.27.7.13.31.89&31239&27776 \cr
\noalign{\hrule}
 & &5.49.11.13.53.109&6487&6498& & &3.5.11.19.113.577&1273&1612 \cr
4014&202397195&4.9.169.361.109.499&22379&32012&4032&204405135&8.11.13.361.31.67&805&10386 \cr
 & &32.3.7.19.23.53.139.151&20989&20976& & &32.9.5.7.23.577&161&48 \cr
\noalign{\hrule}
}%
}
$$
\eject
\vglue -23 pt
\noindent\hskip 1 in\hbox to 6.5 in{\ 4033 -- 4068 \hfill\fbd 204459255 -- 208918395\frb}
\vskip -9 pt
$$
\vbox{
\nointerlineskip
\halign{\strut
    \vrule \ \ \hfil \frb #\ 
   &\vrule \hfil \ \ \fbb #\frb\ 
   &\vrule \hfil \ \ \frb #\ \hfil
   &\vrule \hfil \ \ \frb #\ 
   &\vrule \hfil \ \ \frb #\ \ \vrule \hskip 2 pt
   &\vrule \ \ \hfil \frb #\ 
   &\vrule \hfil \ \ \fbb #\frb\ 
   &\vrule \hfil \ \ \frb #\ \hfil
   &\vrule \hfil \ \ \frb #\ 
   &\vrule \hfil \ \ \frb #\ \vrule \cr%
\noalign{\hrule}
 & &27.5.7.11.13.17.89&19&604& & &9.5.7.53.89.139&1087&998 \cr
4033&204459255&8.3.11.17.19.151&1267&1300&4051&206533845&4.3.7.53.499.1087&13&1100 \cr
 & &64.25.7.13.19.181&905&608& & &32.25.11.13.499&6487&880 \cr
\noalign{\hrule}
 & &9.25.7.13.17.19.31&269&320& & &3.5.7.13.17.37.241&1881&1252 \cr
4034&205016175&128.3.125.7.13.269&3379&1496&4052&206918985&8.27.5.7.11.19.313&463&482 \cr
 & &2048.11.17.31.109&1199&1024& & &32.11.241.313.463&5093&5008 \cr
\noalign{\hrule}
 & &9.7.11.169.17.103&15853&15750& & &9.25.7.11.17.19.37&59&94 \cr
4035&205071867&4.81.125.49.83.191&8041&1318&4053&207051075&4.5.11.19.37.47.59&12367&6762 \cr
 & &16.125.11.17.43.659&5375&5272& & &16.3.49.23.83.149&3427&4648 \cr
\noalign{\hrule}
 & &7.11.17.19.73.113&905&14352& & &3.7.13.59.61.211&589&650 \cr
4036&205160879&32.3.5.13.23.181&657&838&4054&207313197&4.25.169.19.31.211&3355&3186 \cr
 & &128.27.73.419&419&1728& & &16.27.125.11.19.59.61&1375&1368 \cr
\noalign{\hrule}
 & &5.7.11.17.23.29.47&801&134& & &3.343.17.71.167&215&286 \cr
4037&205179205&4.9.7.47.67.89&493&494&4055&207414501&4.5.343.11.13.17.43&639&296 \cr
 & &16.3.13.17.19.29.67.89&5963&5928& & &64.9.13.37.43.71&1591&1248 \cr
\noalign{\hrule}
 & &9.5.11.13.19.23.73&5959&5156& & &9.11.41.137.373&745&608 \cr
4038&205282935&8.23.59.101.1289&1323&34&4056&207418959&64.3.5.19.149.373&37&410 \cr
 & &32.27.49.17.101&833&4848& & &256.25.19.37.41&925&2432 \cr
\noalign{\hrule}
 & &3.5.11.13.17.43.131&747&188& & &27.5.7.13.61.277&6469&4334 \cr
4039&205407345&8.27.47.83.131&4199&1958&4057&207579645&4.9.11.197.6469&4121&2348 \cr
 & &32.11.13.17.19.89&89&304& & &32.11.13.317.587&6457&5072 \cr
\noalign{\hrule}
 & &9.121.13.73.199&7225&16058& & &11.17.67.73.227&2021&14550 \cr
4040&205658739&4.25.7.289.31.37&199&726&4058&207618059&4.3.25.43.47.97&21&22 \cr
 & &16.3.7.121.17.199&17&56& & &16.9.25.7.11.47.97&8225&6984 \cr
\noalign{\hrule}
 & &81.5.47.101.107&187&722& & &9.11.17.19.73.89&295&28 \cr
4041&205711245&4.9.11.17.361.47&221&202&4059&207754569&8.3.5.7.11.59.73&1513&1732 \cr
 & &16.11.13.289.19.101&2717&2312& & &64.7.17.89.433&433&224 \cr
\noalign{\hrule}
 & &27.25.11.53.523&1691&3016& & &243.11.289.269&301&590 \cr
4042&205813575&16.3.11.13.19.29.89&2513&5230&4060&207801693&4.3.5.7.43.59.269&433&374 \cr
 & &64.5.7.359.523&359&224& & &16.5.7.11.17.43.433&2165&2408 \cr
\noalign{\hrule}
 & &81.7.13.17.31.53&191&880& & &11.13.529.41.67&18547&11670 \cr
4043&205879401&32.9.5.11.31.191&911&1190&4061&207802309&4.3.5.17.389.1091&427&1518 \cr
 & &128.25.7.17.911&911&1600& & &16.9.7.11.17.23.61&1037&504 \cr
\noalign{\hrule}
 & &9.19.31.47.827&33&860& & &81.25.13.53.149&26359&22066 \cr
4044&206044569&8.27.5.11.31.43&565&596&4062&207888525&4.11.17.43.59.613&3393&3350 \cr
 & &64.25.11.113.149&16837&8800& & &16.9.25.13.17.29.59.67&3953&3944 \cr
\noalign{\hrule}
 & &3.25.17.19.47.181&1863&2662& & &3.5.7.13.127.1201&61&66 \cr
4045&206082075&4.243.1331.19.23&1645&3944&4063&208199355&4.9.7.11.13.61.1201&635&7772 \cr
 & &64.5.7.11.17.29.47&319&224& & &32.5.11.29.67.127&737&464 \cr
\noalign{\hrule}
 & &11.17.383.2879&667&3546& & &243.5.13.79.167&119&286 \cr
4046&206196859&4.9.17.23.29.197&3443&2270&4064&208383435&4.3.7.11.169.17.79&2117&3460 \cr
 & &16.3.5.11.227.313&4695&1816& & &32.5.7.29.73.173&5017&8176 \cr
\noalign{\hrule}
 & &27.121.17.47.79&2305&962& & &5.289.37.47.83&3523&378 \cr
4047&206216307&4.5.13.37.47.461&75&536&4065&208566965&4.27.7.13.17.271&517&296 \cr
 & &64.3.125.37.67&8375&1184& & &64.9.7.11.37.47&63&352 \cr
\noalign{\hrule}
 & &11.113.127.1307&45&1352& & &7.11.13.257.811&5233&5310 \cr
4048&206324327&16.9.5.169.113&85&254&4066&208635427&4.9.5.59.257.5233&15431&268 \cr
 & &64.3.25.17.127&1275&32& & &32.3.5.13.67.1187&3561&5360 \cr
\noalign{\hrule}
 & &9.11.13.19.23.367&1241&1450& & &49.11.13.83.359&23635&27702 \cr
4049&206407773&4.25.17.29.73.367&1&366&4067&208787579&4.729.5.19.29.163&175&338 \cr
 & &16.3.5.17.29.61&145&8296& & &16.27.125.7.169.29&3375&3016 \cr
\noalign{\hrule}
 & &27.11.19.23.37.43&25&2& & &9.5.7.19.67.521&11&46 \cr
4050&206494299&4.25.11.19.37.43&609&1426&4068&208918395&4.3.11.23.67.521&9367&7826 \cr
 & &16.3.5.7.23.29.31&145&1736& & &16.7.13.17.19.29.43&1247&1768 \cr
\noalign{\hrule}
}%
}
$$
\eject
\vglue -23 pt
\noindent\hskip 1 in\hbox to 6.5 in{\ 4069 -- 4104 \hfill\fbd 209384747 -- 213469725\frb}
\vskip -9 pt
$$
\vbox{
\nointerlineskip
\halign{\strut
    \vrule \ \ \hfil \frb #\ 
   &\vrule \hfil \ \ \fbb #\frb\ 
   &\vrule \hfil \ \ \frb #\ \hfil
   &\vrule \hfil \ \ \frb #\ 
   &\vrule \hfil \ \ \frb #\ \ \vrule \hskip 2 pt
   &\vrule \ \ \hfil \frb #\ 
   &\vrule \hfil \ \ \fbb #\frb\ 
   &\vrule \hfil \ \ \frb #\ \hfil
   &\vrule \hfil \ \ \frb #\ 
   &\vrule \hfil \ \ \frb #\ \vrule \cr%
\noalign{\hrule}
 & &11.169.163.691&153&10& & &81.7.13.23.29.43&4675&4054 \cr
4069&209384747&4.9.5.13.17.691&365&326&4087&211407651&4.3.25.11.13.17.2027&2861&19436 \cr
 & &16.3.25.17.73.163&1275&584& & &32.43.113.2861&2861&1808 \cr
\noalign{\hrule}
 & &3.5.13.17.83.761&539&124& & &7.11.53.197.263&23865&27946 \cr
4070&209385345&8.49.11.31.761&365&396&4088&211440691&4.3.5.37.43.89.157&23&66 \cr
 & &64.9.5.49.121.73&5929&7008& & &16.9.5.11.23.37.157&7659&6280 \cr
\noalign{\hrule}
 & &9.7.13.31.37.223&5621&1292& & &25.7.11.13.79.107&601&426 \cr
4071&209484639&8.49.11.17.19.73&277&1110&4089&211536325&4.3.11.71.107.601&553&624 \cr
 & &32.3.5.11.37.277&1385&176& & &128.9.7.13.79.601&601&576 \cr
\noalign{\hrule}
 & &81.11.17.101.137&8719&16900& & &49.11.13.109.277&27575&31176 \cr
4072&209589039&8.25.169.8719&4275&4444&4090&211562351&16.9.25.433.1103&1201&98 \cr
 & &64.9.625.11.19.101&625&608& & &64.3.25.49.1201&1201&2400 \cr
\noalign{\hrule}
 & &9.5.13.19.109.173&2461&826& & &3.5.13.41.103.257&231&26 \cr
4073&209595555&4.3.7.13.23.59.107&9515&9424&4091&211635645&4.9.7.11.169.103&5375&5272 \cr
 & &128.5.11.19.23.31.173&713&704& & &64.125.11.43.659&16475&15136 \cr
\noalign{\hrule}
 & &7.11.13.19.103.107&1377&656& & &25.11.23.109.307&11457&18518 \cr
4074&209608399&32.81.11.13.17.41&245&206&4092&211653475&4.9.19.47.67.197&7189&4040 \cr
 & &128.27.5.49.17.103&945&1088& & &64.3.5.7.13.79.101&7189&9696 \cr
\noalign{\hrule}
 & &3.125.13.19.31.73&319&1944& & &5.11.37.41.43.59&91&2274 \cr
4075&209610375&16.729.11.19.29&469&260&4093&211674595&4.3.7.13.41.379&333&46 \cr
 & &128.5.7.13.29.67&469&1856& & &16.27.13.23.37&27&2392 \cr
\noalign{\hrule}
 & &9.5.49.13.71.103&473&164& & &625.49.11.17.37&663&712 \cr
4076&209627145&8.3.5.11.41.43.71&349&1414&4094&211894375&16.3.5.13.289.37.89&11039&7746 \cr
 & &32.7.11.101.349&1111&5584& & &64.9.7.19.83.1291&24529&23904 \cr
\noalign{\hrule}
 & &5.7.11.37.41.359&2015&498& & &27.5.7.11.19.29.37&3127&1832 \cr
4077&209672155&4.3.25.11.13.31.83&469&444&4095&211922865&16.3.11.53.59.229&989&760 \cr
 & &32.9.7.13.31.37.67&2077&1872& & &256.5.19.23.43.59&2537&2944 \cr
\noalign{\hrule}
 & &3.5.7.41.53.919&36649&39406& & &5.7.11.31.109.163&43&120 \cr
4078&209683635&4.17.19.61.67.547&5863&27504&4096&212049145&16.3.25.31.43.109&333&442 \cr
 & &128.9.11.13.41.191&2483&2112& & &64.27.13.17.37.43&20683&14688 \cr
\noalign{\hrule}
 & &9.5.29.37.43.101&3781&562& & &27.7.13.23.53.71&43099&43450 \cr
4079&209701755&4.3.5.19.199.281&1001&404&4097&212650893&4.25.49.11.47.79.131&1349&954 \cr
 & &32.7.11.13.19.101&209&1456& & &16.9.5.11.19.53.71.131&1045&1048 \cr
\noalign{\hrule}
 & &9.7.71.173.271&925&286& & &9.11.71.79.383&3475&3554 \cr
4080&209707659&4.25.11.13.37.271&133&138&4098&212676453&4.25.139.383.1777&69&1846 \cr
 & &16.3.5.7.11.13.19.23.37&9361&9880& & &16.3.5.13.23.71.139&1495&1112 \cr
\noalign{\hrule}
 & &25.7.11.13.37.227&4167&4232& & &3.11.19.37.53.173&1935&1352 \cr
4081&210184975&16.9.5.7.11.529.463&4747&5902&4099&212711631&16.27.5.169.37.43&173&308 \cr
 & &64.3.13.23.47.101.227&3243&3232& & &128.7.11.13.43.173&559&448 \cr
\noalign{\hrule}
 & &3.5.11.17.31.41.59&20391&20732& & &243.7.17.53.139&29&110 \cr
4082&210344145&8.9.5.7.71.73.971&1891&1394&4100&213031539&4.3.5.7.11.17.29.53&1301&2780 \cr
 & &32.17.31.41.61.971&971&976& & &32.25.139.1301&1301&400 \cr
\noalign{\hrule}
 & &11.13.361.61.67&155&516& & &11.13.41.163.223&155&378 \cr
4083&210983201&8.3.5.13.31.43.67&231&1102&4101&213114187&4.27.5.7.11.31.163&571&1222 \cr
 & &32.9.5.7.11.19.29&1015&144& & &16.9.5.13.47.571&2115&4568 \cr
\noalign{\hrule}
 & &343.11.17.37.89&213&620& & &3.5.19.29.107.241&5401&15566 \cr
4084&211216313&8.3.5.7.31.71.89&1131&1628&4102&213129555&4.11.43.181.491&9&482 \cr
 & &64.9.5.11.13.29.37&585&928& & &16.9.181.241&543&8 \cr
\noalign{\hrule}
 & &3.5.11.13.19.71.73&1233&3950& & &9.25.61.103.151&8051&7502 \cr
4085&211233165&4.27.125.79.137&7099&3724&4103&213464925&4.25.121.31.83.97&867&1208 \cr
 & &32.49.19.31.229&1519&3664& & &64.3.11.289.97.151&3179&3104 \cr
\noalign{\hrule}
 & &9.5.7.17.19.31.67&649&1118& & &3.25.49.29.2003&3367&2642 \cr
4086&211324365&4.3.5.11.13.17.43.59&1147&1862&4104&213469725&4.343.13.37.1321&8415&8758 \cr
 & &16.49.19.31.37.43&259&344& & &16.9.5.11.17.29.37.151&4983&5032 \cr
\noalign{\hrule}
}%
}
$$
\eject
\vglue -23 pt
\noindent\hskip 1 in\hbox to 6.5 in{\ 4105 -- 4140 \hfill\fbd 213711147 -- 217862645\frb}
\vskip -9 pt
$$
\vbox{
\nointerlineskip
\halign{\strut
    \vrule \ \ \hfil \frb #\ 
   &\vrule \hfil \ \ \fbb #\frb\ 
   &\vrule \hfil \ \ \frb #\ \hfil
   &\vrule \hfil \ \ \frb #\ 
   &\vrule \hfil \ \ \frb #\ \ \vrule \hskip 2 pt
   &\vrule \ \ \hfil \frb #\ 
   &\vrule \hfil \ \ \fbb #\frb\ 
   &\vrule \hfil \ \ \frb #\ \hfil
   &\vrule \hfil \ \ \frb #\ 
   &\vrule \hfil \ \ \frb #\ \vrule \cr%
\noalign{\hrule}
 & &9.169.23.41.149&43&490& & &7.121.29.67.131&1125&14726 \cr
4105&213711147&4.3.5.49.13.23.43&1639&1156&4123&215589451&4.9.125.37.199&2213&2412 \cr
 & &32.7.11.289.149&289&1232& & &32.81.67.2213&2213&1296 \cr
\noalign{\hrule}
 & &5.7.11.23.101.239&153&958& & &125.13.19.29.241&153&88 \cr
4106&213750845&4.9.17.239.479&1313&2750&4124&215785375&16.9.25.11.17.19.29&241&716 \cr
 & &16.3.125.11.13.101&39&200& & &128.3.17.179.241&537&1088 \cr
\noalign{\hrule}
 & &9.5.7.13.109.479&649&170& & &27.5.11.17.43.199&169&304 \cr
4107&213804045&4.25.11.17.59.109&189&1664&4125&216021465&32.169.17.19.199&261&62 \cr
 & &1024.27.7.11.13&33&512& & &128.9.169.29.31&5239&1856 \cr
\noalign{\hrule}
 & &3.25.19.223.673&2451&3124& & &27.19.23.73.251&781&530 \cr
4108&213862575&8.9.11.361.43.71&13&374&4126&216193077&4.9.5.11.53.71.73&137&502 \cr
 & &32.121.13.17.71&2057&14768& & &16.11.53.137.251&583&1096 \cr
\noalign{\hrule}
 & &3.25.11.17.101.151&63&38& & &9.25.11.17.53.97&1411&3836 \cr
4109&213895275&4.27.7.11.17.19.151&601&1060&4127&216307575&8.7.289.83.137&435&146 \cr
 & &32.5.7.19.53.601&7049&9616& & &32.3.5.29.73.137&3973&1168 \cr
\noalign{\hrule}
 & &3.169.17.103.241&6061&6230& & &9.7.17.241.839&3013&2860 \cr
4110&213949437&4.5.7.11.19.29.89.103&19755&1154&4128&216555129&8.5.11.13.23.131.241&29&270 \cr
 & &16.9.25.439.577&10975&13848& & &32.27.25.11.29.131&7975&6288 \cr
\noalign{\hrule}
 & &3.5.169.29.41.71&673&1862& & &3.17.19.23.71.137&33&104 \cr
4111&214002165&4.49.19.71.673&585&88&4129&216785649&16.9.11.13.17.19.23&295&142 \cr
 & &64.9.5.7.11.13.19&133&1056& & &64.5.11.13.59.71&715&1888 \cr
\noalign{\hrule}
 & &9.7.11.169.31.59&283&224& & &81.11.19.23.557&103&4910 \cr
4112&214206993&64.3.49.11.31.283&3835&722&4130&216877419&4.9.5.103.491&253&262 \cr
 & &256.5.13.361.59&361&640& & &16.11.23.131.491&491&1048 \cr
\noalign{\hrule}
 & &27.5.13.17.43.167&2599&2432& & &9.25.11.23.37.103&203&2572 \cr
4113&214245135&256.3.5.17.19.23.113&11&334&4131&216941175&8.3.7.11.29.643&437&206 \cr
 & &1024.11.113.167&1243&512& & &32.19.23.29.103&29&304 \cr
\noalign{\hrule}
 & &11.43.53.83.103&22359&16900& & &3.5.7.13.23.31.223&5137&15602 \cr
4114&214314881&8.3.25.169.29.257&831&454&4132&217033635&4.11.29.269.467&8251&5292 \cr
 & &32.9.5.13.227.277&26559&22160& & &32.27.49.37.223&259&144 \cr
\noalign{\hrule}
 & &9.7.11.19.41.397&2573&1794& & &7.13.19.23.53.103&14837&1746 \cr
4115&214319259&4.27.7.13.23.31.83&55&244&4133&217088053&4.9.37.97.401&583&620 \cr
 & &32.5.11.31.61.83&5063&2480& & &32.3.5.11.31.53.97&3201&2480 \cr
\noalign{\hrule}
 & &9.5.7.37.53.347&817&1612& & &27.5.11.31.53.89&1561&1508 \cr
4116&214347105&8.3.13.19.31.37.43&55&1388&4134&217147095&8.3.5.7.13.29.89.223&61&206 \cr
 & &64.5.11.19.347&19&352& & &32.7.13.61.103.223&43981&46384 \cr
\noalign{\hrule}
 & &9.7.121.107.263&7657&5290& & &5.19.73.173.181&1287&2152 \cr
4117&214518843&4.5.7.13.19.529.31&121&282&4135&217155655&16.9.11.13.73.269&181&38 \cr
 & &16.3.5.121.19.23.47&1081&760& & &64.3.19.181.269&269&96 \cr
\noalign{\hrule}
 & &27.5.13.37.3307&493&506& & &9.17.23.37.1669&3665&1996 \cr
4118&214740045&4.5.11.17.23.29.3307&3321&14&4136&217308807&8.5.23.499.733&6105&5372 \cr
 & &16.81.7.11.17.41&451&2856& & &64.3.25.11.17.37.79&869&800 \cr
\noalign{\hrule}
 & &27.25.49.73.89&587&2990& & &9.5.17.59.61.79&121&274 \cr
4119&214888275&4.125.13.23.587&231&356&4137&217505565&4.121.59.61.137&4503&3854 \cr
 & &32.3.7.11.13.23.89&143&368& & &16.3.11.19.41.47.79&1927&1672 \cr
\noalign{\hrule}
 & &25.7.17.29.47.53&22009&16416& & &81.11.13.89.211&15181&3598 \cr
4120&214911025&64.27.13.19.1693&265&1958&4138&217517157&4.7.17.19.47.257&355&5238 \cr
 & &256.3.5.11.53.89&979&384& & &16.27.5.71.97&97&2840 \cr
\noalign{\hrule}
 & &5.11.17.317.727&205&522& & &25.7.31.137.293&1859&1566 \cr
4121&215479165&4.9.25.11.17.29.41&4121&554&4139&217764925&4.27.7.11.169.29.31&145&548 \cr
 & &16.3.13.277.317&831&104& & &32.3.5.13.841.137&841&624 \cr
\noalign{\hrule}
 & &5.7.11.29.97.199&117&82& & &5.7.11.13.19.29.79&15237&3782 \cr
4122&215517995&4.9.11.13.29.41.97&305&14&4140&217862645&4.9.31.61.1693&77&16 \cr
 & &16.3.5.7.13.41.61&2501&312& & &128.3.7.11.1693&1693&192 \cr
\noalign{\hrule}
}%
}
$$
\eject
\vglue -23 pt
\noindent\hskip 1 in\hbox to 6.5 in{\ 4141 -- 4176 \hfill\fbd 217985625 -- 223277873\frb}
\vskip -9 pt
$$
\vbox{
\nointerlineskip
\halign{\strut
    \vrule \ \ \hfil \frb #\ 
   &\vrule \hfil \ \ \fbb #\frb\ 
   &\vrule \hfil \ \ \frb #\ \hfil
   &\vrule \hfil \ \ \frb #\ 
   &\vrule \hfil \ \ \frb #\ \ \vrule \hskip 2 pt
   &\vrule \ \ \hfil \frb #\ 
   &\vrule \hfil \ \ \fbb #\frb\ 
   &\vrule \hfil \ \ \frb #\ \hfil
   &\vrule \hfil \ \ \frb #\ 
   &\vrule \hfil \ \ \frb #\ \vrule \cr%
\noalign{\hrule}
 & &9.625.11.13.271&331&956& & &9.5.7.11.17.53.71&203&698 \cr
4141&217985625&8.239.271.331&301&30&4159&221659515&4.49.29.71.349&3321&6800 \cr
 & &32.3.5.7.43.239&301&3824& & &128.81.25.17.41&205&576 \cr
\noalign{\hrule}
 & &3.5.149.239.409&903&1142& & &3.7.11.31.83.373&565&348 \cr
4142&218473485&4.9.7.43.149.571&385&956&4160&221697399&8.9.5.29.113.373&3317&40 \cr
 & &32.5.49.11.43.239&539&688& & &128.25.31.107&2675&64 \cr
\noalign{\hrule}
 & &5.11.31.41.53.59&3783&4432& & &27.25.29.47.241&91&44 \cr
4143&218592935&32.3.13.41.97.277&159&118&4161&221726025&8.5.7.11.13.29.241&387&628 \cr
 & &128.9.13.53.59.97&873&832& & &64.9.11.13.43.157&6751&4576 \cr
\noalign{\hrule}
 & &9.11.13.23.83.89&171&82& & &9.7.19.31.43.139&7&50 \cr
4144&218662587&4.81.13.19.41.83&1309&230&4162&221788539&4.3.25.49.31.139&551&968 \cr
 & &16.5.7.11.17.23.41&205&952& & &64.25.121.19.29&3509&800 \cr
\noalign{\hrule}
 & &3.125.7.121.13.53&391&456& & &27.7.121.31.313&1525&1742 \cr
4145&218843625&16.9.25.17.19.23.53&451&26&4163&221897907&4.25.13.61.67.313&187&126 \cr
 & &64.11.13.19.23.41&437&1312& & &16.9.25.7.11.13.17.67&1675&1768 \cr
\noalign{\hrule}
 & &27.5.31.59.887&97&38& & &343.11.103.571&3089&3192 \cr
4146&219013605&4.19.31.97.887&1947&1060&4164&221901449&16.3.2401.19.3089&2745&344 \cr
 & &32.3.5.11.19.53.59&583&304& & &256.27.5.19.43.61&31293&27520 \cr
\noalign{\hrule}
 & &3.5.29.37.43.317&337&308& & &3.625.7.11.29.53&27&292 \cr
4147&219390945&8.7.11.37.317.337&45&362&4165&221904375&8.81.125.7.73&221&346 \cr
 & &32.9.5.7.181.337&3801&5392& & &32.13.17.73.173&2941&15184 \cr
\noalign{\hrule}
 & &9.25.7.11.19.23.29&2689&2486& & &3.7.13.41.79.251&85&6 \cr
4148&219559725&4.121.19.113.2689&195&2494&4166&221945997&4.9.5.17.41.251&1027&1232 \cr
 & &16.3.5.13.29.43.113&559&904& & &128.7.11.13.17.79&187&64 \cr
\noalign{\hrule}
 & &25.41.2209.97&29601&25624& & &243.5.31.71.83&979&1222 \cr
4149&219629825&16.9.11.13.23.3203&1153&2050&4167&221959845&4.5.11.13.47.83.89&159&76 \cr
 & &64.3.25.11.41.1153&1153&1056& & &32.3.11.13.19.53.89&13091&15664 \cr
\noalign{\hrule}
 & &3.25.11.251.1063&153&98& & &3.125.7.19.61.73&251&176 \cr
4150&220120725&4.27.5.49.17.1063&59&1004&4168&222093375&32.5.11.19.73.251&377&4392 \cr
 & &32.7.17.59.251&119&944& & &512.9.13.29.61&1131&256 \cr
\noalign{\hrule}
 & &5.7.29.31.47.149&3289&1926& & &9.5.11.37.67.181&459&278 \cr
4151&220349395&4.9.11.13.23.31.107&25&118&4169&222106005&4.243.5.17.37.139&1451&1694 \cr
 & &16.3.25.23.59.107&2461&7080& & &16.7.121.139.1451&10703&11608 \cr
\noalign{\hrule}
 & &9.5.343.109.131&28613&29158& & &3.5.7.13.29.31.181&18533&23782 \cr
4152&220396365&4.7.13.31.61.71.239&109&1782&4170&222111435&4.11.23.43.47.431&5215&4698 \cr
 & &16.81.11.13.71.109&923&792& & &16.81.5.7.29.43.149&1161&1192 \cr
\noalign{\hrule}
 & &27.5.61.73.367&1169&802& & &7.11.19.137.1109&659&450 \cr
4153&220623885&4.5.7.61.167.401&12111&1924&4171&222277979&4.9.25.7.137.659&1109&1768 \cr
 & &32.3.11.13.37.367&481&176& & &64.3.25.13.17.1109&663&800 \cr
\noalign{\hrule}
 & &81.11.13.23.829&15325&3742& & &125.121.61.241&3933&3692 \cr
4154&220853061&4.25.613.1871&629&1242&4172&222352625&8.9.121.13.19.23.71&305&58 \cr
 & &16.27.25.17.23.37&629&200& & &32.3.5.23.29.61.71&1633&1392 \cr
\noalign{\hrule}
 & &11.13.19.127.641&1295&1422& & &9.11.13.29.59.101&101&218 \cr
4155&221182819&4.9.5.7.37.79.641&597&44&4173&222407757&4.59.10201.109&1885&8316 \cr
 & &32.27.5.11.37.199&5373&2960& & &32.27.5.7.11.13.29&21&80 \cr
\noalign{\hrule}
 & &5.49.11.19.29.149&109&1530& & &13.17.43.47.499&6975&16478 \cr
4156&221256805&4.9.25.17.19.109&1391&1334&4174&222873859&4.9.25.7.11.31.107&313&338 \cr
 & &16.3.13.17.23.29.107&4173&3128& & &16.3.11.169.107.313&10329&11128 \cr
\noalign{\hrule}
 & &25.13.43.71.223&929&14904& & &27.11.13.17.43.79&1073&1060 \cr
4157&221266175&16.81.23.929&775&154&4175&222968889&8.5.11.17.29.37.43.53&21385&186 \cr
 & &64.3.25.7.11.31&77&2976& & &32.3.25.7.13.31.47&1457&2800 \cr
\noalign{\hrule}
 & &3.49.11.169.811&5627&50& & &7.13.19.29.61.73&795&592 \cr
4158&221624403&4.25.7.17.331&169&162&4176&223277873&32.3.5.13.37.53.61&1649&1584 \cr
 & &16.81.25.169.17&675&136& & &1024.27.11.17.37.97&96903&95744 \cr
\noalign{\hrule}
}%
}
$$
\eject
\vglue -23 pt
\noindent\hskip 1 in\hbox to 6.5 in{\ 4177 -- 4212 \hfill\fbd 223290053 -- 228379785\frb}
\vskip -9 pt
$$
\vbox{
\nointerlineskip
\halign{\strut
    \vrule \ \ \hfil \frb #\ 
   &\vrule \hfil \ \ \fbb #\frb\ 
   &\vrule \hfil \ \ \frb #\ \hfil
   &\vrule \hfil \ \ \frb #\ 
   &\vrule \hfil \ \ \frb #\ \ \vrule \hskip 2 pt
   &\vrule \ \ \hfil \frb #\ 
   &\vrule \hfil \ \ \fbb #\frb\ 
   &\vrule \hfil \ \ \frb #\ \hfil
   &\vrule \hfil \ \ \frb #\ 
   &\vrule \hfil \ \ \frb #\ \vrule \cr%
\noalign{\hrule}
 & &7.17.29.89.727&3835&1254& & &9.11.17.29.41.113&395&5028 \cr
4177&223290053&4.3.5.11.13.17.19.59&1363&1068&4195&226122831&8.27.5.79.419&857&1276 \cr
 & &32.9.19.29.47.89&423&304& & &64.5.11.29.857&857&160 \cr
\noalign{\hrule}
 & &9.11.109.127.163&9751&10950& & &27.25.49.13.17.31&1881&2284 \cr
4178&223384491&4.27.25.49.73.199&109&136&4196&226596825&8.243.5.11.19.571&18011&13394 \cr
 & &64.5.17.73.109.199&6205&6368& & &32.7.31.37.83.181&3071&2896 \cr
\noalign{\hrule}
 & &9.125.13.17.29.31&125&278& & &3.5.13.37.89.353&6533&16412 \cr
4179&223513875&4.15625.29.139&5797&9828&4197&226673655&8.11.47.139.373&445&72 \cr
 & &32.27.7.11.13.17.31&77&48& & &128.9.5.89.139&417&64 \cr
\noalign{\hrule}
 & &7.11.17.29.71.83&1269&2182& & &3.11.13.23.83.277&1783&1506 \cr
4180&223704173&4.27.47.71.1091&581&510&4198&226852197&4.9.83.251.1783&1265&518 \cr
 & &16.81.5.7.17.47.83&405&376& & &16.5.7.11.23.37.251&1757&1480 \cr
\noalign{\hrule}
 & &27.13.59.101.107&15799&21758& & &3.11.17.197.2053&35549&32200 \cr
4181&223802163&4.7.11.23.37.43.61&3&40&4199&226891401&16.25.7.19.23.1871&2465&594 \cr
 & &64.3.5.7.11.23.61&8855&1952& & &64.27.125.11.17.29&1125&928 \cr
\noalign{\hrule}
 & &3.11.59.103.1117&33003&32900& & &31.37.53.3733&56881&58842 \cr
4182&224004297&8.27.25.7.11.19.47.193&515&2&4200&226932803&4.9.7.11.467.5171&1885&3286 \cr
 & &32.125.7.103.193&1351&2000& & &16.3.5.7.11.13.29.31.53&1885&1848 \cr
\noalign{\hrule}
 & &49.289.71.223&15723&4796& & &9.5.7.121.59.101&19&866 \cr
4183&224211113&8.9.11.109.1747&1037&710&4201&227127285&4.3.19.101.433&3095&2662 \cr
 & &32.3.5.11.17.61.71&305&528& & &16.5.1331.619&619&88 \cr
\noalign{\hrule}
 & &9.5.11.17.149.179&763&848& & &9.125.37.53.103&6517&6358 \cr
4184&224436465&32.7.11.53.109.149&525&674&4202&227230875&4.3.343.11.289.19.37&5459&5200 \cr
 & &128.3.25.49.53.337&16513&16960& & &128.25.49.13.17.53.103&833&832 \cr
\noalign{\hrule}
 & &27.5.7.11.17.31.41&71&226& & &3.11.19.349.1039&109099&108052 \cr
4185&224604765&4.7.17.41.71.113&747&460&4203&227357097&8.7.17.79.227.1381&20705&2772 \cr
 & &32.9.5.23.83.113&2599&1328& & &64.9.5.49.11.41.101&10045&9696 \cr
\noalign{\hrule}
 & &9.25.11.29.31.101&2627&302& & &9.7.41.107.823&37&786 \cr
4186&224727525&4.3.11.37.71.151&97&310&4204&227461563&4.27.37.41.131&535&572 \cr
 & &16.5.31.97.151&97&1208& & &32.5.11.13.107.131&1703&880 \cr
\noalign{\hrule}
 & &3.7.11.17.23.47.53&13711&13690& & &9.5.11.13.23.29.53&9331&8954 \cr
4187&224989611&4.5.17.23.1369.13711&14089&378&4205&227483685&4.3.7.1331.31.37.43&2665&1334 \cr
 & &16.27.5.7.37.73.193&13505&13896& & &16.5.7.13.23.29.37.41&287&296 \cr
\noalign{\hrule}
 & &9.5.7.13.17.53.61&127&188& & &243.13.61.1181&649&532 \cr
4188&225065295&8.13.17.47.53.127&55&744&4206&227577519&8.27.7.11.19.59.61&4747&3100 \cr
 & &128.3.5.11.31.127&3937&704& & &64.25.11.31.47.101&34441&37600 \cr
\noalign{\hrule}
 & &3.529.47.3023&3133&110& & &9.11.23.67.1493&1469&3010 \cr
4189&225482547&4.5.11.13.23.241&1337&1314&4207&227770587&4.3.5.7.11.13.43.113&67&406 \cr
 & &16.9.5.7.13.73.191&13943&10920& & &16.5.49.13.29.67&1885&392 \cr
\noalign{\hrule}
 & &25.149.151.401&6237&16262& & &3.7.17.31.59.349&4579&2750 \cr
4190&225552475&4.81.7.11.47.173&755&802&4208&227880597&4.125.11.17.19.241&2077&2502 \cr
 & &16.9.5.7.11.151.401&77&72& & &16.9.5.11.31.67.139&3685&3336 \cr
\noalign{\hrule}
 & &9.7.19.29.67.97&767&506& & &3.5.17.59.109.139&1329&1034 \cr
4191&225599787&4.7.11.13.23.59.97&255&158&4209&227946795&4.9.11.47.109.443&2927&1946 \cr
 & &16.3.5.11.13.17.23.79&14773&11960& & &16.7.47.139.2927&2927&2632 \cr
\noalign{\hrule}
 & &3.17.19.23.53.191&2171&2222& & &3.13.23.43.61.97&1395&836 \cr
4192&225611301&4.11.13.19.53.101.167&153&1160&4210&228224607&8.27.5.11.19.31.61&529&1118 \cr
 & &64.9.5.11.17.29.167&4843&5280& & &32.5.11.13.529.43&115&176 \cr
\noalign{\hrule}
 & &3.5.49.121.43.59&83&212& & &9.5.7.11.13.37.137&551&408 \cr
4193&225628095&8.49.121.53.83&85&36&4211&228333105&16.27.5.17.19.29.37&137&322 \cr
 & &64.9.5.17.53.83&4233&1696& & &64.7.19.23.29.137&551&736 \cr
\noalign{\hrule}
 & &5.11.17.41.71.83&2967&2926& & &3.5.17.29.89.347&15317&14872 \cr
4194&225908155&4.3.5.7.121.17.19.23.43&107&498&4212&228379785&16.11.169.4913.53&1333&6246 \cr
 & &16.9.7.19.43.83.107&6741&6536& & &64.9.13.31.43.347&1333&1248 \cr
\noalign{\hrule}
}%
}
$$
\eject
\vglue -23 pt
\noindent\hskip 1 in\hbox to 6.5 in{\ 4213 -- 4248 \hfill\fbd 228404071 -- 233957475\frb}
\vskip -9 pt
$$
\vbox{
\nointerlineskip
\halign{\strut
    \vrule \ \ \hfil \frb #\ 
   &\vrule \hfil \ \ \fbb #\frb\ 
   &\vrule \hfil \ \ \frb #\ \hfil
   &\vrule \hfil \ \ \frb #\ 
   &\vrule \hfil \ \ \frb #\ \ \vrule \hskip 2 pt
   &\vrule \ \ \hfil \frb #\ 
   &\vrule \hfil \ \ \fbb #\frb\ 
   &\vrule \hfil \ \ \frb #\ \hfil
   &\vrule \hfil \ \ \frb #\ 
   &\vrule \hfil \ \ \frb #\ \vrule \cr%
\noalign{\hrule}
 & &7.37.41.137.157&2761&3048& & &9.11.13.37.43.113&3151&1708 \cr
4213&228404071&16.3.11.127.137.251&5&132&4231&231380721&8.3.7.11.23.61.137&301&370 \cr
 & &128.9.5.121.251&5445&16064& & &32.5.49.37.43.137&685&784 \cr
\noalign{\hrule}
 & &27.49.13.97.137&515&164& & &3.11.23.37.73.113&51&2650 \cr
4214&228557511&8.5.7.41.103.137&429&292&4232&231656667&4.9.25.11.17.53&4181&4234 \cr
 & &64.3.5.11.13.41.73&2993&1760& & &16.5.29.37.73.113&29&40 \cr
\noalign{\hrule}
 & &5.7.13.289.37.47&6231&5942& & &11.17.43.151.191&171&20 \cr
4215&228669805&4.3.5.13.31.67.2971&3663&692&4233&231910481&8.9.5.11.17.19.43&349&382 \cr
 & &32.27.11.31.37.173&5363&4752& & &32.3.5.19.191.349&1047&1520 \cr
\noalign{\hrule}
 & &31.107.151.457&1887&1430& & &3.5.121.23.67.83&1037&872 \cr
4216&228896219&4.3.5.11.13.17.37.151&2033&372&4234&232143945&16.11.17.61.67.109&905&234 \cr
 & &32.9.17.19.31.107&153&304& & &64.9.5.13.109.181&4251&5792 \cr
\noalign{\hrule}
 & &5.7.121.13.23.181&21257&14922& & &3.11.17.19.139.157&1707&1846 \cr
4217&229193965&4.9.29.733.829&877&1610&4235&232611357&4.9.13.71.157.569&4405&2992 \cr
 & &16.3.5.7.23.29.877&877&696& & &128.5.11.17.71.881&4405&4544 \cr
\noalign{\hrule}
 & &125.13.73.1933&979&954& & &5.11.67.179.353&1161&808 \cr
4218&229302125&4.9.5.11.13.53.73.89&27451&3866&4236&232844095&16.27.5.43.67.101&353&152 \cr
 & &16.3.97.283.1933&849&776& & &256.9.19.43.353&817&1152 \cr
\noalign{\hrule}
 & &9.11.13.19.83.113&1631&3100& & &9.11.53.199.223&1295&712 \cr
4219&229344687&8.3.25.7.11.31.233&1577&1918&4237&232846119&16.5.7.37.89.199&5129&1836 \cr
 & &32.5.49.19.83.137&685&784& & &128.27.17.23.223&391&192 \cr
\noalign{\hrule}
 & &3.29.961.41.67&1075&114& & &5.49.11.13.289.23&369&226 \cr
4220&229668429&4.9.25.19.43.67&3509&3844&4238&232877645&4.9.7.17.23.41.113&591&200 \cr
 & &32.5.121.29.961&121&80& & &64.27.25.41.197&8077&4320 \cr
\noalign{\hrule}
 & &9.7.11.19.73.239&685&1196& & &3.5.7.11.13.361.43&127&146 \cr
4221&229724649&8.5.13.23.137.239&3091&2406&4239&233077845&4.5.11.19.43.73.127&27671&18684 \cr
 & &32.3.11.13.281.401&5213&4496& & &32.27.7.59.67.173&10207&9648 \cr
\noalign{\hrule}
 & &3.11.13.19.71.397&215&182& & &3.5.11.19.23.53.61&141&164 \cr
4222&229752237&4.5.7.169.19.43.71&79&3132&4240&233115465&8.9.11.19.41.47.53&305&2186 \cr
 & &32.27.5.7.29.79&1305&8848& & &32.5.41.61.1093&1093&656 \cr
\noalign{\hrule}
 & &27.25.7.17.47.61&53&32& & &5.23.67.79.383&4365&4444 \cr
4223&230291775&64.9.5.47.53.61&127&2618&4241&233130185&8.9.25.11.67.97.101&2171&4596 \cr
 & &256.7.11.17.127&1397&128& & &64.27.11.13.167.383&4509&4576 \cr
\noalign{\hrule}
 & &3.11.289.19.31.41&83&206& & &9.11.13.37.59.83&217&550 \cr
4224&230309013&4.11.19.31.83.103&1035&7514&4242&233190243&4.25.7.121.31.83&147&268 \cr
 & &16.9.5.13.289.23&69&520& & &32.3.5.343.31.67&10385&5488 \cr
\noalign{\hrule}
 & &7.121.289.23.41&20687&14040& & &9.5.49.11.59.163&111&52 \cr
4225&230830369&16.27.5.13.137.151&2023&242&4243&233260335&8.27.5.49.11.13.37&9617&7582 \cr
 & &64.9.7.121.289&9&32& & &32.17.59.163.223&223&272 \cr
\noalign{\hrule}
 & &9.5.7.11.163.409&5339&4930& & &3.13.17.29.61.199&253&240 \cr
4226&231001155&4.25.11.17.19.29.281&5001&224&4244&233396553&32.9.5.11.23.61.199&731&5308 \cr
 & &256.3.7.29.1667&1667&3712& & &256.5.17.43.1327&6635&5504 \cr
\noalign{\hrule}
 & &25.7.23.137.419&297&122& & &5.23.41.59.839&891&52 \cr
4227&231047075&4.27.11.23.61.137&35&172&4245&233397215&8.81.5.11.13.59&277&218 \cr
 & &32.3.5.7.11.43.61&2623&528& & &32.9.13.109.277&3601&15696 \cr
\noalign{\hrule}
 & &3.5.29.37.53.271&4245&5782& & &5.7.11.17.53.673&2449&2262 \cr
4228&231172485&4.9.25.49.59.283&319&94&4246&233453605&4.3.5.13.29.31.53.79&121&1764 \cr
 & &16.7.11.29.47.283&3619&2264& & &32.27.49.121.79&2133&1232 \cr
\noalign{\hrule}
 & &9.5.13.17.43.541&1265&724& & &5.11.19.41.53.103&2065&108 \cr
4229&231350535&8.25.11.23.43.181&2757&1768&4247&233890855&8.27.25.7.11.59&437&212 \cr
 & &128.3.11.13.17.919&919&704& & &64.3.7.19.23.53&483&32 \cr
\noalign{\hrule}
 & &5.7.23.31.73.127&15719&3966& & &9.25.37.157.179&731&194 \cr
4230&231357805&4.3.11.661.1429&2921&4350&4248&233957475&4.3.17.43.97.157&2395&2552 \cr
 & &16.9.25.23.29.127&145&72& & &64.5.11.29.43.479&13717&15328 \cr
\noalign{\hrule}
}%
}
$$
\eject
\vglue -23 pt
\noindent\hskip 1 in\hbox to 6.5 in{\ 4249 -- 4284 \hfill\fbd 234235287 -- 240038799\frb}
\vskip -9 pt
$$
\vbox{
\nointerlineskip
\halign{\strut
    \vrule \ \ \hfil \frb #\ 
   &\vrule \hfil \ \ \fbb #\frb\ 
   &\vrule \hfil \ \ \frb #\ \hfil
   &\vrule \hfil \ \ \frb #\ 
   &\vrule \hfil \ \ \frb #\ \ \vrule \hskip 2 pt
   &\vrule \ \ \hfil \frb #\ 
   &\vrule \hfil \ \ \fbb #\frb\ 
   &\vrule \hfil \ \ \frb #\ \hfil
   &\vrule \hfil \ \ \frb #\ 
   &\vrule \hfil \ \ \frb #\ \vrule \cr%
\noalign{\hrule}
 & &27.11.13.19.31.103&25&272& & &9.7.11.13.23.31.37&53&766 \cr
4249&234235287&32.25.17.31.103&953&798&4267&237666429&4.11.37.53.383&177&230 \cr
 & &128.3.5.7.19.953&953&2240& & &16.3.5.23.59.383&1915&472 \cr
\noalign{\hrule}
 & &3.5.11.17.29.43.67&433&304& & &9.25.11.13.83.89&257&2482 \cr
4250&234354945&32.5.17.19.29.433&2881&5346&4268&237676725&4.3.13.17.73.257&89&860 \cr
 & &128.243.11.43.67&81&64& & &32.5.17.43.89&731&16 \cr
\noalign{\hrule}
 & &9.5.11.19.97.257&185&442& & &7.13.29.251.359&473575&473826 \cr
4251&234457245&4.3.25.13.17.37.97&17&308&4269&237797651&4.3.25.19.157.503.997&102157&54372 \cr
 & &32.7.11.289.37&259&4624& & &32.9.5.11.23.37.197.251&36445&36432 \cr
\noalign{\hrule}
 & &3.5.169.17.43.127&649&14& & &13.23.59.103.131&6363&7130 \cr
4252&235341795&4.7.11.13.43.59&2921&2448&4270&238030013&4.9.5.7.529.31.101&3417&286 \cr
 & &128.9.17.23.127&23&192& & &16.27.5.11.13.17.67&5049&2680 \cr
\noalign{\hrule}
 & &27.5.11.41.53.73&1829&344& & &5.13.53.257.269&473&216 \cr
4253&235564065&16.31.43.59.73&2175&2132&4271&238163185&16.27.5.11.43.269&1157&1802 \cr
 & &128.3.25.13.29.31.41&1885&1984& & &64.9.13.17.53.89&801&544 \cr
\noalign{\hrule}
 & &3.25.11.31.61.151&213&62& & &3.25.7.23.109.181&239&566 \cr
4254&235571325&4.9.961.61.71&755&206&4272&238227675&4.5.181.239.283&333&572 \cr
 & &16.5.71.103.151&103&568& & &32.9.11.13.37.283&11037&6512 \cr
\noalign{\hrule}
 & &27.7.37.41.823&535&572& & &3.5.11.169.83.103&4489&4060 \cr
4255&235964799&8.5.7.11.13.107.823&37&786&4273&238388865&8.25.7.13.29.4489&649&1026 \cr
 & &32.3.5.11.13.37.131&1703&880& & &32.27.7.11.19.59.67&8911&8496 \cr
\noalign{\hrule}
 & &27.41.383.557&1027&644& & &27.25.7.11.43.107&71&64 \cr
4256&236157417&8.9.7.13.23.41.79&383&2200&4274&239136975&128.5.11.43.71.107&31&504 \cr
 & &128.25.11.13.383&275&832& & &2048.9.7.31.71&2201&1024 \cr
\noalign{\hrule}
 & &7.11.19.109.1481&533&666& & &9.49.11.31.37.43&2885&778 \cr
4257&236170627&4.9.13.37.41.1481&59&1540&4275&239256171&4.5.31.389.577&273&304 \cr
 & &32.3.5.7.11.37.59&295&1776& & &128.3.5.7.13.19.389&5057&6080 \cr
\noalign{\hrule}
 & &3.7.11.17.137.439&12151&13468& & &5.11.29.31.47.103&9&38 \cr
4258&236181561&8.49.13.29.37.419&1485&16988&4276&239363245&4.9.5.11.19.31.103&87&428 \cr
 & &64.27.5.11.31.137&279&160& & &32.27.19.29.107&2889&304 \cr
\noalign{\hrule}
 & &9.121.17.19.673&1027&1030& & &19.23.47.89.131&1953&1060 \cr
4259&236725731&4.3.5.13.19.79.103.673&259&1760&4277&239464201&8.9.5.7.31.53.89&55&34 \cr
 & &256.25.7.11.13.37.103&33475&33152& & &32.3.25.11.17.31.53&22525&16368 \cr
\noalign{\hrule}
 & &243.5.11.13.29.47&655&2018& & &3.5.7.11.17.29.421&1443&1022 \cr
4260&236814435&4.25.13.131.1009&667&342&4278&239723715&4.9.49.11.13.37.73&1513&2150 \cr
 & &16.9.19.23.29.131&437&1048& & &16.25.17.43.73.89&3139&3560 \cr
\noalign{\hrule}
 & &9.5.11.169.19.149&57&112& & &9.41.233.2789&1579&1210 \cr
4261&236827305&32.27.7.361.149&3275&748&4279&239789853&4.5.121.233.1579&2071&492 \cr
 & &256.25.11.17.131&655&2176& & &32.3.5.11.19.41.109&1045&1744 \cr
\noalign{\hrule}
 & &27.5.41.113.379&34937&34558& & &7.11.83.157.239&669&430 \cr
4262&237047445&4.9.49.23.31.37.467&565&286&4280&239809493&4.3.5.11.43.83.223&7077&2512 \cr
 & &16.5.49.11.13.113.467&5137&5096& & &128.9.7.157.337&337&576 \cr
\noalign{\hrule}
 & &9.5.7.169.61.73&319&136& & &243.11.23.47.83&179&1730 \cr
4263&237055455&16.3.11.13.17.29.73&181&38&4281&239829579&4.81.5.173.179&15281&15686 \cr
 & &64.17.19.29.181&3439&15776& & &16.7.11.23.31.37.59&1829&2072 \cr
\noalign{\hrule}
 & &3.19.961.61.71&1155&194& & &25.7.13.19.31.179&36963&33562 \cr
4264&237239187&4.9.5.7.11.61.97&589&284&4282&239855525&4.27.1369.97.173&1463&94 \cr
 & &32.7.11.19.31.71&77&16& & &16.3.7.11.19.47.97&1551&776 \cr
\noalign{\hrule}
 & &9.7.17.19.107.109&253&146& & &27.5.7.17.67.223&439&506 \cr
4265&237330387&4.3.11.17.23.73.109&3115&608&4283&240027165&4.11.17.23.223.439&1167&6296 \cr
 & &256.5.7.11.19.89&445&1408& & &64.3.11.389.787&8657&12448 \cr
\noalign{\hrule}
 & &7.11.19.37.41.107&2123&1836& & &3.49.11.13.19.601&1165&564 \cr
4266&237472697&8.27.121.17.19.193&2015&1652&4284&240038799&8.9.5.7.11.47.233&463&230 \cr
 & &64.9.5.7.13.17.31.59&23777&24480& & &32.25.23.47.463&10649&18800 \cr
\noalign{\hrule}
}%
}
$$
\eject
\vglue -23 pt
\noindent\hskip 1 in\hbox to 6.5 in{\ 4285 -- 4320 \hfill\fbd 240288191 -- 244002745\frb}
\vskip -9 pt
$$
\vbox{
\nointerlineskip
\halign{\strut
    \vrule \ \ \hfil \frb #\ 
   &\vrule \hfil \ \ \fbb #\frb\ 
   &\vrule \hfil \ \ \frb #\ \hfil
   &\vrule \hfil \ \ \frb #\ 
   &\vrule \hfil \ \ \frb #\ \ \vrule \hskip 2 pt
   &\vrule \ \ \hfil \frb #\ 
   &\vrule \hfil \ \ \fbb #\frb\ 
   &\vrule \hfil \ \ \frb #\ \hfil
   &\vrule \hfil \ \ \frb #\ 
   &\vrule \hfil \ \ \frb #\ \vrule \cr%
\noalign{\hrule}
 & &11.13.101.127.131&2465&15696& & &5.11.29.47.53.61&13&42 \cr
4285&240288191&32.9.5.17.29.109&127&18&4303&242361845&4.3.7.13.47.53.61&18821&18450 \cr
 & &128.81.17.127&1377&64& & &16.27.25.11.29.41.59&1593&1640 \cr
\noalign{\hrule}
 & &25.11.13.23.37.79&6723&16898& & &125.13.29.37.139&541&1266 \cr
4286&240343675&4.81.7.17.71.83&13&14&4304&242363875&4.3.5.37.211.541&13&198 \cr
 & &16.3.49.13.17.71.83&10437&11288& & &16.27.11.13.541&297&4328 \cr
\noalign{\hrule}
 & &5.13.841.53.83&765&76& & &9.5.23.421.557&5353&7458 \cr
4287&240471335&8.9.25.17.19.83&671&754&4305&242704395&4.27.11.53.101.113&205&92 \cr
 & &32.3.11.13.17.29.61&561&976& & &32.5.23.41.53.101&2173&1616 \cr
\noalign{\hrule}
 & &3.5.7.13.23.79.97&97&202& & &9.5.43.271.463&1199&736 \cr
4288&240579885&4.79.9409.101&715&8694&4306&242790255&64.11.23.109.271&1389&1118 \cr
 & &16.27.5.7.11.13.23&9&88& & &256.3.11.13.43.463&143&128 \cr
\noalign{\hrule}
 & &23.37.277.1021&15703&22074& & &27.97.137.677&5629&7660 \cr
4289&240677267&4.3.13.41.283.383&185&198&4307&242909631&8.9.5.13.383.433&1507&1940 \cr
 & &16.27.5.11.37.41.283&11603&11880& & &64.25.11.13.97.137&325&352 \cr
\noalign{\hrule}
 & &25.7.13.23.43.107&1233&158& & &27.5.7.11.97.241&1665&986 \cr
4290&240747325&4.9.7.23.79.137&125&286&4308&243003915&4.243.25.17.29.37&241&484 \cr
 & &16.3.125.11.13.79&869&120& & &32.121.17.37.241&407&272 \cr
\noalign{\hrule}
 & &3.11.17.29.59.251&2989&1278& & &9.7.23.43.47.83&205&44 \cr
4291&240927621&4.27.49.11.61.71&1075&1004&4309&243059607&8.3.5.11.41.43.47&37&178 \cr
 & &32.25.7.43.61.251&2623&2800& & &32.11.37.41.89&1517&15664 \cr
\noalign{\hrule}
 & &81.5.7.11.59.131&949&536& & &9.5.7.13.23.29.89&1931&2554 \cr
4292&241028865&16.3.13.67.73.131&2591&2518&4310&243091485&4.3.29.1277.1931&21413&15620 \cr
 & &64.67.1259.2591&84353&82912& & &32.5.49.11.19.23.71&1349&1232 \cr
\noalign{\hrule}
 & &9.169.23.61.113&2489&110& & &3.25.7.11.13.41.79&459&74 \cr
4293&241137819&4.3.5.11.13.19.131&989&976&4311&243167925&4.81.5.17.37.79&163&242 \cr
 & &128.11.19.23.43.61&817&704& & &16.121.17.37.163&1793&5032 \cr
\noalign{\hrule}
 & &3.49.11.17.31.283&3835&722& & &9.13.17.19.41.157&961&638 \cr
4294&241160997&4.5.13.17.361.59&1&18&4312&243260667&4.3.11.29.961.157&835&3718 \cr
 & &16.9.5.13.19.59&195&8968& & &16.5.121.169.167&2171&4840 \cr
\noalign{\hrule}
 & &25.11.13.19.53.67&8287&8262& & &9.11.13.421.449&283&166 \cr
4295&241201675&4.243.11.17.53.8287&6767&1520&4313&243280323&4.11.83.283.421&9429&14060 \cr
 & &128.27.5.17.19.67.101&1717&1728& & &32.3.5.7.19.37.449&703&560 \cr
\noalign{\hrule}
 & &5.11.17.841.307&2537&2682& & &3.11.19.173.2243&12423&12250 \cr
4296&241404845&4.9.11.29.43.59.149&85&234&4314&243300453&4.9.125.49.19.41.101&173&2 \cr
 & &16.81.5.13.17.43.59&4779&4472& & &16.5.7.41.101.173&1435&808 \cr
\noalign{\hrule}
 & &9.25.7.11.13.29.37&973&248& & &3.25.13.257.971&353&418 \cr
4297&241666425&16.3.49.13.31.139&1375&3182&4315&243308325&4.5.11.19.353.971&4369&486 \cr
 & &64.125.11.37.43&43&160& & &16.243.17.19.257&323&648 \cr
\noalign{\hrule}
 & &27.11.13.17.29.127&931&1228& & &3.5.11.361.61.67&1263&2422 \cr
4298&241741071&8.49.13.19.29.307&35&342&4316&243442155&4.9.7.19.173.421&251&3538 \cr
 & &32.9.5.343.361&1805&5488& & &16.7.29.61.251&203&2008 \cr
\noalign{\hrule}
 & &27.7.289.19.233&299&5192& & &5.343.19.31.241&6087&1508 \cr
4299&241807167&16.9.11.13.23.59&185&68&4317&243442535&8.3.7.13.29.2029&913&1116 \cr
 & &128.5.17.37.59&185&3776& & &64.27.11.13.31.83&3861&2656 \cr
\noalign{\hrule}
 & &5.17.59.139.347&9889&10584& & &9.7.11.17.23.29.31&13&680 \cr
4300&241888495&16.27.49.11.17.29.31&19&338&4318&243595737&16.5.13.289.31&1863&1894 \cr
 & &64.9.7.169.19.31&10647&18848& & &64.81.5.23.947&947&1440 \cr
\noalign{\hrule}
 & &9.67.229.1753&533&1220& & &9.13.23.31.37.79&939&88 \cr
4301&242066511&8.3.5.13.41.61.67&1753&748&4319&243839583&16.27.11.31.313&575&262 \cr
 & &64.11.13.17.1753&187&416& & &64.25.11.23.131&3275&352 \cr
\noalign{\hrule}
 & &11.17.31.37.1129&549&580& & &5.7.23.61.4969&87659&86256 \cr
4302&242158081&8.9.5.11.17.29.37.61&1129&92&4320&244002745&32.9.11.13.599.613&2473&4270 \cr
 & &64.3.5.23.29.1129&667&480& & &128.3.5.7.13.61.2473&2473&2496 \cr
\noalign{\hrule}
}%
}
$$
\eject
\vglue -23 pt
\noindent\hskip 1 in\hbox to 6.5 in{\ 4321 -- 4356 \hfill\fbd 244206435 -- 249954915\frb}
\vskip -9 pt
$$
\vbox{
\nointerlineskip
\halign{\strut
    \vrule \ \ \hfil \frb #\ 
   &\vrule \hfil \ \ \fbb #\frb\ 
   &\vrule \hfil \ \ \frb #\ \hfil
   &\vrule \hfil \ \ \frb #\ 
   &\vrule \hfil \ \ \frb #\ \ \vrule \hskip 2 pt
   &\vrule \ \ \hfil \frb #\ 
   &\vrule \hfil \ \ \fbb #\frb\ 
   &\vrule \hfil \ \ \frb #\ \hfil
   &\vrule \hfil \ \ \frb #\ 
   &\vrule \hfil \ \ \frb #\ \vrule \cr%
\noalign{\hrule}
 & &3.5.121.157.857&731&126& & &9.11.13.43.61.73&7&664 \cr
4321&244206435&4.27.7.17.43.157&1685&1528&4339&246433473&16.7.13.43.83&1739&1830 \cr
 & &64.5.43.191.337&8213&10784& & &64.3.5.37.47.61&235&1184 \cr
\noalign{\hrule}
 & &9.7.139.27889&11575&11638& & &3.25.11.59.61.83&909&4154 \cr
4322&244223973&4.25.11.529.167.463&239&2076&4340&246441525&4.27.5.31.67.101&671&166 \cr
 & &32.3.5.529.173.239&41347&42320& & &16.11.61.67.83&67&8 \cr
\noalign{\hrule}
 & &5.13.37.157.647&6867&17072& & &3.7.11.43.103.241&589&830 \cr
4323&244297495&32.9.7.11.97.109&697&370&4341&246566859&4.5.7.19.31.83.103&3049&144 \cr
 & &128.3.5.7.17.37.41&697&1344& & &128.9.19.3049&9147&1216 \cr
\noalign{\hrule}
 & &7.11.13.467.523&23171&17100& & &9.5.11.13.193.199&1141&146 \cr
4324&244485241&8.9.25.17.19.29.47&22769&22774&4342&247149045&4.7.73.163.193&757&594 \cr
 & &32.3.5.29.59.193.22769&1343371&1343280& & &16.27.11.73.757&757&1752 \cr
\noalign{\hrule}
 & &3.5.167.233.419&433&266& & &3.7.13.19.43.1109&3473&4290 \cr
4325&244555635&4.5.7.19.419.433&2563&468&4343&247352469&4.9.5.11.169.23.151&1591&70 \cr
 & &32.9.11.13.19.233&627&208& & &16.25.7.23.37.43&925&184 \cr
\noalign{\hrule}
 & &27.49.17.73.149&4189&6688& & &9.11.29.79.1091&1981&890 \cr
4326&244634607&64.9.11.19.59.71&115&56&4344&247448619&4.5.7.79.89.283&135&418 \cr
 & &1024.5.7.11.23.71&8165&5632& & &16.27.25.11.19.89&1691&600 \cr
\noalign{\hrule}
 & &13.17.19.101.577&4825&6138& & &3.25.7.11.19.37.61&113&146 \cr
4327&244705123&4.9.25.11.17.31.193&7501&6974&4345&247649325&4.25.19.61.73.113&311&1836 \cr
 & &16.3.121.13.317.577&951&968& & &32.27.17.73.311&5287&10512 \cr
\noalign{\hrule}
 & &9.25.361.23.131&341&134& & &23.29.373.997&6057&16874 \cr
4328&244730925&4.11.19.31.67.131&2283&206&4346&248044627&4.9.11.13.59.673&685&1334 \cr
 & &16.3.11.103.761&761&9064& & &16.3.5.13.23.29.137&685&312 \cr
\noalign{\hrule}
 & &13.19.31.113.283&85&198& & &7.13.23.29.61.67&785&756 \cr
4329&244863203&4.9.5.11.13.17.19.31&283&686&4347&248068639&8.27.5.49.13.61.157&593&44 \cr
 & &16.3.5.343.11.283&1029&440& & &64.3.5.11.157.593&19569&25120 \cr
\noalign{\hrule}
 & &3.5.7.71.107.307&4873&2724& & &27.5.7.181.1451&421&484 \cr
4330&244889295&8.9.5.11.227.443&4249&8236&4348&248186295&8.3.121.421.1451&1357&94 \cr
 & &64.7.29.71.607&607&928& & &32.121.23.47.59&7139&17296 \cr
\noalign{\hrule}
 & &9.5.7.19.151.271&347&1012& & &25.11.29.163.191&893&702 \cr
4331&244912185&8.11.23.271.347&5025&1208&4349&248285675&4.27.5.13.19.47.163&1351&116 \cr
 & &128.3.25.67.151&335&64& & &32.3.7.29.47.193&987&3088 \cr
\noalign{\hrule}
 & &5.7.11.19.107.313&261&52& & &27.49.13.17.23.37&331&110 \cr
4332&244986665&8.9.5.7.13.29.107&313&222&4350&248817933&4.3.5.11.23.37.331&169&238 \cr
 & &32.27.29.37.313&999&464& & &16.5.7.169.17.331&331&520 \cr
\noalign{\hrule}
 & &27.5.7.23.29.389&239&428& & &5.13.23.31.41.131&45369&48034 \cr
4333&245192535&8.5.107.239.389&403&792&4351&248918995&4.9.7.47.5041.73&7667&2626 \cr
 & &128.9.11.13.31.107&4433&6848& & &16.3.7.11.13.17.41.101&1717&1848 \cr
\noalign{\hrule}
 & &25.13.19.151.263&1639&324& & &3.5.23.29.139.179&781&114 \cr
4334&245227775&8.81.5.11.19.149&1357&1208&4352&248934405&4.9.11.19.71.139&389&250 \cr
 & &128.3.11.23.59.151&1947&1472& & &16.125.11.19.389&7391&2200 \cr
\noalign{\hrule}
 & &9.25.509.2143&24497&29078& & &13.29.37.107.167&2475&2368 \cr
4335&245427075&4.7.11.17.31.67.131&2025&2036&4353&249254681&128.9.25.11.13.1369&3841&266 \cr
 & &32.81.25.7.17.67.509&1071&1072& & &512.3.7.19.23.167&1311&1792 \cr
\noalign{\hrule}
 & &243.25.11.13.283&4067&388& & &9.5.11.13.17.43.53&17869&8366 \cr
4336&245849175&8.3.5.49.83.97&283&962&4354&249311205&4.47.89.107.167&837&8686 \cr
 & &32.7.13.37.283&259&16& & &16.27.31.43.101&93&808 \cr
\noalign{\hrule}
 & &81.19.23.6947&3255&3692& & &3.13.23.47.61.97&101&198 \cr
4337&245902959&8.243.5.7.13.31.71&979&1222&4355&249454803&4.27.11.47.61.101&299&970 \cr
 & &32.5.7.11.169.47.89&65065&66928& & &16.5.13.23.97.101&101&40 \cr
\noalign{\hrule}
 & &25.11.17.23.29.79&3699&1724& & &3.5.7.23.29.43.83&149&66 \cr
4338&246339775&8.27.23.137.431&319&112&4356&249954915&4.9.7.11.23.29.149&13&680 \cr
 & &256.3.7.11.29.137&959&384& & &64.5.13.17.149&149&7072 \cr
\noalign{\hrule}
}%
}
$$
\eject
\vglue -23 pt
\noindent\hskip 1 in\hbox to 6.5 in{\ 4357 -- 4392 \hfill\fbd 250128615 -- 255454353\frb}
\vskip -9 pt
$$
\vbox{
\nointerlineskip
\halign{\strut
    \vrule \ \ \hfil \frb #\ 
   &\vrule \hfil \ \ \fbb #\frb\ 
   &\vrule \hfil \ \ \frb #\ \hfil
   &\vrule \hfil \ \ \frb #\ 
   &\vrule \hfil \ \ \frb #\ \ \vrule \hskip 2 pt
   &\vrule \ \ \hfil \frb #\ 
   &\vrule \hfil \ \ \fbb #\frb\ 
   &\vrule \hfil \ \ \frb #\ \hfil
   &\vrule \hfil \ \ \frb #\ 
   &\vrule \hfil \ \ \frb #\ \vrule \cr%
\noalign{\hrule}
 & &3.5.11.31.79.619&1781&76& & &3.5.121.43.53.61&179&126 \cr
4357&250128615&8.13.19.79.137&819&682&4375&252319485&4.27.7.121.43.179&185&5018 \cr
 & &32.9.7.11.169.31&507&112& & &16.5.7.13.37.193&7141&728 \cr
\noalign{\hrule}
 & &5.11.43.53.1997&1235&762& & &25.11.23.73.547&2397&16072 \cr
4358&250313965&4.3.25.13.19.53.127&341&666&4376&252563575&16.3.49.17.41.47&547&1380 \cr
 & &16.27.11.31.37.127&4699&6696& & &128.9.5.23.547&9&64 \cr
\noalign{\hrule}
 & &3.25.11.73.4157&2491&1666& & &3.5.11.361.31.137&719&2226 \cr
4359&250355325&4.49.17.47.53.73&1501&990&4377&252972555&4.9.7.19.53.719&445&274 \cr
 & &16.9.5.7.11.17.19.79&2261&1896& & &16.5.7.53.89.137&371&712 \cr
\noalign{\hrule}
 & &25.49.11.29.641&6063&988& & &25.7.17.97.877&99&778 \cr
4360&250486775&8.3.7.13.19.43.47&155&174&4378&253080275&4.9.25.11.17.389&2921&1754 \cr
 & &32.9.5.13.29.31.43&1333&1872& & &16.3.23.127.877&127&552 \cr
\noalign{\hrule}
 & &3.37.43.47.1117&311&1428& & &9.5.49.121.13.73&3293&2636 \cr
4361&250577727&8.9.7.17.43.311&1765&1034&4379&253197945&8.5.13.37.89.659&21373&21462 \cr
 & &32.5.7.11.47.353&2471&880& & &32.3.49.11.29.37.67.73&1073&1072 \cr
\noalign{\hrule}
 & &9.5.17.157.2087&53&104& & &25.49.11.19.23.43&1107&118 \cr
4362&250659135&16.3.5.13.53.2087&1441&646&4380&253208725&4.27.11.19.41.59&2531&3310 \cr
 & &64.11.13.17.19.131&2489&4576& & &16.3.5.331.2531&7593&2648 \cr
\noalign{\hrule}
 & &5.17.19.23.43.157&7337&6552& & &9.5.7.11.13.17.331&2923&4028 \cr
4363&250765895&16.9.7.11.13.529.29&323&1910&4381&253468215&8.3.11.19.37.53.79&1655&952 \cr
 & &64.3.5.13.17.19.191&573&416& & &128.5.7.17.53.331&53&64 \cr
\noalign{\hrule}
 & &27.11.61.109.127&36379&36760& & &81.25.349.359&517&158 \cr
4364&250793631&16.9.5.7.919.5197&1537&6734&4382&253714275&4.3.11.47.79.349&5395&5744 \cr
 & &64.5.49.13.29.37.53&75313&76960& & &128.5.11.13.83.359&913&832 \cr
\noalign{\hrule}
 & &3125.49.11.149&837&802& & &25.11.13.19.37.101&567&358 \cr
4365&250971875&4.27.625.7.31.401&13&638&4383&253835725&4.81.7.13.101.179&925&388 \cr
 & &16.9.11.13.29.401&5213&2088& & &32.27.25.7.37.97&679&432 \cr
\noalign{\hrule}
 & &9.25.7.169.23.41&811&788& & &25.23.109.4051&221067&220492 \cr
4366&251003025&8.3.25.7.13.197.811&79&2354&4384&253896425&8.9.7.121.29.199.277&25415&8102 \cr
 & &32.11.79.107.197&21079&13904& & &32.3.5.7.13.17.23.4051&273&272 \cr
\noalign{\hrule}
 & &9.25.7.11.43.337&4727&2368& & &3.7.17.19.89.421&473&150 \cr
4367&251056575&128.3.5.29.37.163&337&152&4385&254152227&4.9.25.11.43.421&323&98 \cr
 & &2048.19.29.337&551&1024& & &16.49.11.17.19.43&77&344 \cr
\noalign{\hrule}
 & &5.7.17.19.97.229&639&506& & &9.5.11.17.19.37.43&4667&16618 \cr
4368&251117965&4.9.11.17.23.71.97&749&458&4386&254377035&4.7.13.359.1187&663&1850 \cr
 & &16.3.7.11.23.107.229&759&856& & &16.3.25.169.17.37&169&40 \cr
\noalign{\hrule}
 & &9.841.139.239&745&506& & &3.25.17.31.41.157&15829&16356 \cr
4369&251449749&4.5.11.23.841.149&363&478&4387&254422425&8.9.5.11.29.47.1439&67&1372 \cr
 & &16.3.1331.149.239&1331&1192& & &64.343.11.47.67&22981&16544 \cr
\noalign{\hrule}
 & &3.25.7.121.17.233&3959&934& & &5.7.13.17.167.197&1089&1082 \cr
4370&251622525&4.17.37.107.467&33885&33418&4388&254473765&4.9.5.121.17.197.541&22477&14362 \cr
 & &16.27.5.49.11.31.251&1953&2008& & &16.3.7.11.169.19.43.167&1677&1672 \cr
\noalign{\hrule}
 & &27.11.169.29.173&1273&976& & &5.11.47.83.1187&16281&3224 \cr
4371&251818281&32.13.19.29.61.67&575&1368&4389&254676785&16.243.13.31.67&107&94 \cr
 & &512.9.25.361.23&8303&6400& & &64.81.31.47.107&2511&3424 \cr
\noalign{\hrule}
 & &27.5.13.43.47.71&1903&1292& & &5.17.89.151.223&84939&83426 \cr
4372&251826705&8.3.11.17.19.43.173&169&40&4390&254736245&4.3.7.23.59.101.1231&51&110 \cr
 & &128.5.169.17.173&2249&1088& & &16.9.5.11.17.101.1231&9999&9848 \cr
\noalign{\hrule}
 & &9.49.17.19.29.61&845&2266& & &3.17.23.43.61.83&985&924 \cr
4373&251981667&4.3.5.11.169.19.103&527&812&4391&255372657&8.9.5.7.11.17.43.197&1655&10126 \cr
 & &32.7.11.13.17.29.31&341&208& & &32.25.61.83.331&331&400 \cr
\noalign{\hrule}
 & &9.5.29.157.1231&1793&562& & &9.7.121.23.31.47&3071&1300 \cr
4374&252213435&4.3.11.29.163.281&22423&23380&4392&255454353&8.3.25.11.13.37.83&47&8 \cr
 & &32.5.7.17.167.1319&19873&21104& & &128.5.37.47.83&415&2368 \cr
\noalign{\hrule}
}%
}
$$
\eject
\vglue -23 pt
\noindent\hskip 1 in\hbox to 6.5 in{\ 4393 -- 4428 \hfill\fbd 255863685 -- 261065849\frb}
\vskip -9 pt
$$
\vbox{
\nointerlineskip
\halign{\strut
    \vrule \ \ \hfil \frb #\ 
   &\vrule \hfil \ \ \fbb #\frb\ 
   &\vrule \hfil \ \ \frb #\ \hfil
   &\vrule \hfil \ \ \frb #\ 
   &\vrule \hfil \ \ \frb #\ \ \vrule \hskip 2 pt
   &\vrule \ \ \hfil \frb #\ 
   &\vrule \hfil \ \ \fbb #\frb\ 
   &\vrule \hfil \ \ \frb #\ \hfil
   &\vrule \hfil \ \ \frb #\ 
   &\vrule \hfil \ \ \frb #\ \vrule \cr%
\noalign{\hrule}
 & &3.5.7.11.17.83.157&1569&158& & &81.25.7.13.23.61&1333&3608 \cr
4393&255863685&4.9.5.7.79.523&229&166&4411&258537825&16.11.23.31.41.43&793&540 \cr
 & &16.83.229.523&523&1832& & &128.27.5.13.41.61&41&64 \cr
\noalign{\hrule}
 & &5.13.19.43.61.79&199&594& & &5.23.29.31.41.61&297&602 \cr
4394&255912995&4.27.11.19.43.199&3965&4592&4412&258565885&4.27.7.11.23.41.43&329&122 \cr
 & &128.9.5.7.13.41.61&369&448& & &16.3.49.43.47.61&2303&1032 \cr
\noalign{\hrule}
 & &11.17.41.173.193&855&1048& & &3.5.7.13.29.47.139&305&682 \cr
4395&255993463&16.9.5.17.19.41.131&107&148&4413&258608805&4.25.11.31.61.139&377&1152 \cr
 & &128.3.19.37.107.131&42051&44992& & &1024.9.13.29.61&183&512 \cr
\noalign{\hrule}
 & &9.41.827.839&20729&13178& & &9.11.23.1369.83&11635&19852 \cr
4396&256031757&4.11.19.599.1091&845&246&4414&258728679&8.5.7.13.179.709&809&1518 \cr
 & &16.3.5.11.169.19.41&845&1672& & &32.3.5.7.11.23.809&809&560 \cr
\noalign{\hrule}
 & &3.25.7.121.29.139&1411&1614& & &27.49.23.67.127&221&248 \cr
4397&256069275&4.9.17.83.139.269&1331&80&4415&258920361&16.7.13.17.23.31.127&2871&50 \cr
 & &128.5.1331.269&269&704& & &64.9.25.11.17.29&5423&800 \cr
\noalign{\hrule}
 & &5.7.11.293.2273&1139&1134& & &5.11.31.383.397&4379&16686 \cr
4398&256405765&4.81.49.11.17.67.293&2273&950&4416&259246955&4.81.29.103.151&175&278 \cr
 & &16.3.25.17.19.67.2273&1615&1608& & &16.27.25.7.29.139&3753&8120 \cr
\noalign{\hrule}
 & &9.5.11.17.131.233&817&118& & &49.13.19.29.739&2875&6732 \cr
4399&256851045&4.3.19.43.59.131&2465&5002&4417&259379393&8.9.125.7.11.17.23&817&58 \cr
 & &16.5.17.29.41.61&1769&328& & &32.3.17.19.29.43&51&688 \cr
\noalign{\hrule}
 & &27.25.13.19.23.67&107&1648& & &27.25.7.13.41.103&317&8 \cr
4400&256923225&32.5.19.103.107&759&1274&4418&259397775&16.9.7.41.317&1133&1450 \cr
 & &128.3.49.11.13.23&539&64& & &64.25.11.29.103&319&32 \cr
\noalign{\hrule}
 & &7.11.13.19.59.229&393&374& & &5.13.29.157.877&17819&7614 \cr
4401&256965709&4.3.7.121.17.131.229&21551&5700&4419&259543765&4.81.47.103.173&377&550 \cr
 & &32.9.25.19.23.937&8433&9200& & &16.9.25.11.13.29.47&517&360 \cr
\noalign{\hrule}
 & &5.7.19.23.67.251&6237&1468& & &27.11.169.31.167&717&1120 \cr
4402&257216015&8.81.49.11.367&1357&2680&4420&259849161&64.81.5.7.13.239&1363&310 \cr
 & &128.3.5.23.59.67&59&192& & &256.25.29.31.47&1363&3200 \cr
\noalign{\hrule}
 & &9.11.13.19.67.157&2597&386& & &9.49.13.137.331&19807&25540 \cr
4403&257221107&4.3.49.13.53.193&335&244&4421&259974351&8.5.29.683.1277&297&980 \cr
 & &32.5.7.53.61.67&2135&848& & &64.27.25.49.11.29&725&1056 \cr
\noalign{\hrule}
 & &9.5.19.43.47.149&2659&3404& & &3.11.23.41.61.137&3707&4650 \cr
4404&257465295&8.3.19.23.37.2659&275&2384&4422&260061483&4.9.25.121.31.337&2501&7946 \cr
 & &256.25.11.23.149&253&640& & &16.5.29.41.61.137&29&40 \cr
\noalign{\hrule}
 & &3.5.343.11.29.157&709&390& & &7.11.13.289.29.31&1863&1894 \cr
4405&257677035&4.9.25.49.13.709&383&58&4423&260070811&4.81.7.11.23.29.947&13&680 \cr
 & &16.29.383.709&709&3064& & &64.9.5.13.17.947&947&1440 \cr
\noalign{\hrule}
 & &3.11.13.17.23.29.53&45&538& & &13.19.67.79.199&109&90 \cr
4406&257814843&4.27.5.13.23.269&187&164&4424&260166829&4.9.5.13.67.79.109&3355&1938 \cr
 & &32.5.11.17.41.269&1345&656& & &16.27.25.11.17.19.61&7425&8296 \cr
\noalign{\hrule}
 & &961.127.2113&3025&912& & &9.5.121.13.29.127&2395&998 \cr
4407&257885311&32.3.25.121.19.31&1537&762&4425&260701155&4.25.11.479.499&377&102 \cr
 & &128.9.29.53.127&1537&576& & &16.3.13.17.29.499&499&136 \cr
\noalign{\hrule}
 & &25.7.13.23.4933&2869&2064& & &3.19.37.337.367&4225&16694 \cr
4408&258119225&32.3.5.13.19.43.151&2629&3864&4426&260839011&4.25.169.17.491&209&2664 \cr
 & &512.9.7.11.23.239&2629&2304& & &64.9.5.11.19.37&15&352 \cr
\noalign{\hrule}
 & &5.11.13.31.43.271&1143&872& & &3.7.47.241.1097&187&910 \cr
4409&258288745&16.9.11.43.109.127&4439&248&4427&260940099&4.5.49.11.13.17.47&387&152 \cr
 & &256.3.23.31.193&579&2944& & &64.9.13.17.19.43&9503&1824 \cr
\noalign{\hrule}
 & &9.17.47.83.433&35075&31174& & &121.31.79.881&21&100 \cr
4410&258437349&4.25.11.13.23.61.109&1301&102&4428&261065849&8.3.25.7.31.881&363&518 \cr
 & &16.3.25.13.17.1301&1301&2600& & &32.9.5.49.121.37&1665&784 \cr
\noalign{\hrule}
}%
}
$$
\eject
\vglue -23 pt
\noindent\hskip 1 in\hbox to 6.5 in{\ 4429 -- 4464 \hfill\fbd 261080475 -- 268141775\frb}
\vskip -9 pt
$$
\vbox{
\nointerlineskip
\halign{\strut
    \vrule \ \ \hfil \frb #\ 
   &\vrule \hfil \ \ \fbb #\frb\ 
   &\vrule \hfil \ \ \frb #\ \hfil
   &\vrule \hfil \ \ \frb #\ 
   &\vrule \hfil \ \ \frb #\ \ \vrule \hskip 2 pt
   &\vrule \ \ \hfil \frb #\ 
   &\vrule \hfil \ \ \fbb #\frb\ 
   &\vrule \hfil \ \ \frb #\ \hfil
   &\vrule \hfil \ \ \frb #\ 
   &\vrule \hfil \ \ \frb #\ \vrule \cr%
\noalign{\hrule}
 & &3.25.17.23.29.307&521&1014& & &5.991.53569&493&498 \cr
4429&261080475&4.9.5.169.23.521&319&526&4447&265434395&4.3.17.29.83.53569&6325&47244 \cr
 & &16.11.29.263.521&2893&4168& & &32.9.25.11.23.31.127&26289&27280 \cr
\noalign{\hrule}
 & &9.25.49.19.29.43&2323&902& & &9.5.7.11.19.37.109&271&62 \cr
4430&261215325&4.3.11.19.23.41.101&157&280&4448&265512555&4.5.7.31.109.271&1221&676 \cr
 & &64.5.7.11.101.157&1727&3232& & &32.3.11.169.31.37&169&496 \cr
\noalign{\hrule}
 & &11.13.17.41.43.61&173&360& & &3.17.79.149.443&41&38 \cr
4431&261437033&16.9.5.43.61.173&95&34&4449&265942203&4.17.19.41.149.443&711&6820 \cr
 & &64.3.25.17.19.173&4325&1824& & &32.9.5.11.19.31.79&1045&1488 \cr
\noalign{\hrule}
 & &5.7.13.37.41.379&707&1188& & &25.7.11.17.79.103&7479&6346 \cr
4432&261599065&8.27.49.11.41.101&2291&1850&4450&266283325&4.27.17.19.167.277&245&78 \cr
 & &32.3.25.11.29.37.79&2607&2320& & &16.81.5.49.13.277&3601&4536 \cr
\noalign{\hrule}
 & &27.13.41.131.139&20141&25840& & &3.11.131.229.269&409&278 \cr
4433&262045719&32.5.11.17.19.1831&393&1438&4451&266301123&4.11.139.269.409&899&630 \cr
 & &128.3.17.131.719&719&1088& & &16.9.5.7.29.31.409&18879&16360 \cr
\noalign{\hrule}
 & &9.11.19.23.73.83&455&202& & &9.5.17.29.41.293&7239&4774 \cr
4434&262130517&4.5.7.13.19.83.101&4071&2494&4452&266508405&4.27.7.11.19.31.127&85&104 \cr
 & &16.3.7.23.29.43.59&2537&1624& & &64.5.11.13.17.31.127&4433&4064 \cr
\noalign{\hrule}
 & &5.7.17.593.743&1111&1854& & &3.25.7.121.13.17.19&7919&8174 \cr
4435&262156405&4.9.7.11.17.101.103&7709&4310&4453&266741475&4.5.13.61.67.7919&6137&1782 \cr
 & &16.3.5.13.431.593&431&312& & &16.81.11.17.361.61&513&488 \cr
\noalign{\hrule}
 & &11.23.31.107.313&18387&8684& & &9.49.47.61.211&4895&5444 \cr
4436&262669913&8.81.13.167.227&745&1426&4454&266777217&8.5.11.47.89.1361&939&422 \cr
 & &32.27.5.23.31.149&745&432& & &32.3.5.89.211.313&1565&1424 \cr
\noalign{\hrule}
 & &3.11.71.313.359&1505&2444& & &3.25.23.137.1129&5291&4162 \cr
4437&263275881&8.5.7.13.43.47.71&297&626&4455&266810925&4.25.11.13.37.2081&1503&578 \cr
 & &32.27.5.11.43.313&215&144& & &16.9.11.13.289.167&11271&14696 \cr
\noalign{\hrule}
 & &27.5.7.11.13.1949&13015&12322& & &9.25.49.11.31.71&6071&4366 \cr
4438&263378115&4.3.25.19.61.101.137&3451&26&4456&266926275&4.3.5.13.37.59.467&869&574 \cr
 & &16.7.13.17.29.101&1717&232& & &16.7.11.41.79.467&3239&3736 \cr
\noalign{\hrule}
 & &27.13.569.1321&485&836& & &25.7.11.13.47.227&873&262 \cr
4439&263828799&8.5.11.19.97.569&333&236&4457&266991725&4.9.5.7.11.97.131&3859&464 \cr
 & &64.9.5.11.19.37.59&7733&9440& & &128.3.17.29.227&87&1088 \cr
\noalign{\hrule}
 & &3.7.11.17.31.41.53&285&166& & &27.11.361.47.53&283&230 \cr
4440&264534501&4.9.5.19.31.53.83&451&4850&4458&267077547&4.5.11.19.23.47.283&477&40 \cr
 & &16.125.11.41.97&125&776& & &64.9.25.53.283&283&800 \cr
\noalign{\hrule}
 & &27.41.43.67.83&1895&1508& & &9.25.7.121.23.61&139&114 \cr
4441&264709161&8.3.5.13.29.67.379&451&5378&4459&267376725&4.27.7.11.19.61.139&1619&460 \cr
 & &32.5.11.41.2689&2689&880& & &32.5.23.139.1619&1619&2224 \cr
\noalign{\hrule}
 & &3.7.11.13.227.389&28747&23690& & &3.5.19.29.179.181&4325&924 \cr
4442&265173909&4.5.17.19.23.89.103&4279&4176&4460&267777735&8.9.125.7.11.173&841&716 \cr
 & &128.9.11.17.23.29.389&1479&1472& & &64.7.11.841.179&319&224 \cr
\noalign{\hrule}
 & &9.7.11.29.67.197&229&240& & &27.5.13.331.461&55&406 \cr
4443&265260303&32.27.5.29.197.229&101&884&4461&267797205&4.25.7.11.29.331&2183&1458 \cr
 & &256.13.17.101.229&23129&28288& & &16.729.7.37.59&999&3304 \cr
\noalign{\hrule}
 & &9.11.13.361.571&103709&102422& & &5.7.23.167.1993&5063&4902 \cr
4444&265290597&4.83.137.617.757&10849&73680&4462&267928955&4.3.19.43.61.83.167&275&2898 \cr
 & &128.3.5.19.307.571&307&320& & &16.27.25.7.11.23.83&913&1080 \cr
\noalign{\hrule}
 & &3.5.11.29.31.1789&3481&1886& & &3.17.19.43.59.109&1335&518 \cr
4445&265371315&4.23.31.41.3481&43&1314&4463&267960477&4.9.5.7.37.59.89&1793&3458 \cr
 & &16.9.43.59.73&9417&472& & &16.49.11.13.19.163&2119&4312 \cr
\noalign{\hrule}
 & &11.13.19.23.31.137&111&98& & &25.11.361.37.73&231&1156 \cr
4446&265399277&4.3.49.23.31.37.137&39&5030&4464&268141775&8.3.7.121.289.19&1095&962 \cr
 & &16.9.5.7.13.503&4527&280& & &32.9.5.13.17.37.73&153&208 \cr
\noalign{\hrule}
}%
}
$$
\eject
\vglue -23 pt
\noindent\hskip 1 in\hbox to 6.5 in{\ 4465 -- 4500 \hfill\fbd 268949725 -- 275732261\frb}
\vskip -9 pt
$$
\vbox{
\nointerlineskip
\halign{\strut
    \vrule \ \ \hfil \frb #\ 
   &\vrule \hfil \ \ \fbb #\frb\ 
   &\vrule \hfil \ \ \frb #\ \hfil
   &\vrule \hfil \ \ \frb #\ 
   &\vrule \hfil \ \ \frb #\ \ \vrule \hskip 2 pt
   &\vrule \ \ \hfil \frb #\ 
   &\vrule \hfil \ \ \fbb #\frb\ 
   &\vrule \hfil \ \ \frb #\ \hfil
   &\vrule \hfil \ \ \frb #\ 
   &\vrule \hfil \ \ \frb #\ \vrule \cr%
\noalign{\hrule}
 & &25.121.67.1327&361&966& & &9.7.31.1681.83&1143&1430 \cr
4465&268949725&4.3.5.7.361.23.67&1327&858&4483&272488419&4.81.5.11.13.41.127&25079&18094 \cr
 & &16.9.11.13.19.1327&171&104& & &16.31.83.109.809&809&872 \cr
\noalign{\hrule}
 & &7.17.19.271.439&297&26& & &25.7.587.2657&281&306 \cr
4466&268988909&4.27.7.11.13.439&209&230&4484&272940325&4.9.7.17.281.2657&16465&2134 \cr
 & &16.9.5.121.13.19.23&2783&4680& & &16.3.5.11.37.89.97&10767&7832 \cr
\noalign{\hrule}
 & &9.5.7.11.23.31.109&841&358& & &13.157.181.739&9405&19012 \cr
4467&269289405&4.3.5.841.31.179&247&218&4485&273002119&8.9.5.49.11.19.97&181&104 \cr
 & &16.13.19.29.109.179&3401&3016& & &128.3.7.13.97.181&291&448 \cr
\noalign{\hrule}
 & &9.5.23.29.47.191&395&968& & &7.11.17.23.47.193&3&190 \cr
4468&269444655&16.3.25.121.23.79&923&802&4486&273100597&4.3.5.7.19.23.47&153&176 \cr
 & &64.13.71.79.401&28471&32864& & &128.27.5.11.17.19&135&1216 \cr
\noalign{\hrule}
 & &3.7.11.47.103.241&51&52& & &81.13.53.59.83&469&220 \cr
4469&269503311&8.9.7.11.13.17.47.241&155&2806&4487&273296673&8.27.5.7.11.59.67&8843&10922 \cr
 & &32.5.13.17.23.31.61&25415&30256& & &32.37.43.127.239&30353&25456 \cr
\noalign{\hrule}
 & &9.7.11.19.59.347&17&116& & &27.11.31.113.263&4355&4324 \cr
4470&269567991&8.17.29.59.347&675&328&4488&273622833&8.9.5.13.23.47.67.113&217&1252 \cr
 & &128.27.25.29.41&1025&5568& & &64.7.31.47.67.313&14711&15008 \cr
\noalign{\hrule}
 & &9.25.11.13.289.29&257&62& & &3.23.43.257.359&51&308 \cr
4471&269658675&4.3.5.289.31.257&377&88&4489&273744321&8.9.7.11.17.23.43&359&1090 \cr
 & &64.11.13.29.257&257&32& & &32.5.11.109.359&109&880 \cr
\noalign{\hrule}
 & &5.17.19.23.53.137&1353&3682& & &81.5.7.13.17.19.23&605&214 \cr
4472&269709845&4.3.7.11.23.41.263&603&340&4490&273795795&4.9.25.121.19.107&221&2254 \cr
 & &32.27.5.7.11.17.67&2079&1072& & &16.49.11.13.17.23&7&88 \cr
\noalign{\hrule}
 & &729.5.11.23.293&16441&326& & &9.7.11.13.113.269&1711&710 \cr
4473&270200205&4.41.163.401&3141&3542&4491&273846573&4.5.29.59.71.113&4867&3156 \cr
 & &16.9.7.11.23.349&349&56& & &32.3.5.31.157.263&8153&12560 \cr
\noalign{\hrule}
 & &5.31.59.127.233&553&612& & &81.25.43.47.67&2453&428 \cr
4474&270609695&8.9.7.17.31.79.127&105083&105610&4492&274199175&8.11.47.107.223&3741&1288 \cr
 & &32.3.5.11.41.59.179.233&1969&1968& & &128.3.7.23.29.43&203&1472 \cr
\noalign{\hrule}
 & &29.89.101.1039&10571&19560& & &7.11.13.23.43.277&177&100 \cr
4475&270847559&16.3.5.11.961.163&89&252&4493&274226953&8.3.25.13.23.43.59&399&958 \cr
 & &128.27.5.7.31.89&945&1984& & &32.9.25.7.19.479&4311&7600 \cr
\noalign{\hrule}
 & &9.125.7.11.31.101&17&38& & &29.97.211.463&1301&1512 \cr
4476&271222875&4.3.25.17.19.31.101&539&236&4494&274810409&16.27.7.463.1301&1345&44 \cr
 & &32.49.11.17.19.59&1003&2128& & &128.9.5.7.11.269&16947&3520 \cr
\noalign{\hrule}
 & &25.7.43.109.331&2509&2178& & &3.139.769.857&813&44 \cr
4477&271494475&4.9.25.7.121.13.193&799&1324&4495&274816761&8.9.11.139.271&769&760 \cr
 & &32.3.11.13.17.47.331&2431&2256& & &128.5.19.271.769&1355&1216 \cr
\noalign{\hrule}
 & &9.25.11.13.23.367&1799&1776& & &121.841.37.73&1411&1290 \cr
4478&271589175&32.27.7.37.257.367&9709&200&4496&274856461&4.3.5.17.841.43.83&55&786 \cr
 & &512.25.49.19.73&3577&4864& & &16.9.25.11.83.131&3275&5976 \cr
\noalign{\hrule}
 & &9.5.7.11.131.599&97&34& & &9.25.11.13.43.199&137&62 \cr
4479&271895085&4.5.11.17.97.599&3537&3052&4497&275321475&4.3.11.13.31.43.137&85&1592 \cr
 & &32.27.7.17.109.131&327&272& & &64.5.17.31.199&527&32 \cr
\noalign{\hrule}
 & &9.7.17.23.43.257&359&1090& & &9.25.49.13.17.113&3229&7054 \cr
4480&272219283&4.5.109.257.359&51&308&4498&275327325&4.7.3229.3527&13959&10730 \cr
 & &32.3.5.7.11.17.109&109&880& & &16.27.5.11.29.37.47&4089&3256 \cr
\noalign{\hrule}
 & &5.11.19.23.47.241&2249&3294& & &9.5.19.23.107.131&31031&26216 \cr
4481&272244445&4.27.13.47.61.173&4807&1940&4499&275644305&16.7.11.13.29.31.113&173&30 \cr
 & &32.9.5.11.19.23.97&97&144& & &64.3.5.31.113.173&3503&5536 \cr
\noalign{\hrule}
 & &3.5.49.17.19.31.37&6409&7556& & &49.29.61.3181&47619&44630 \cr
4482&272303535&8.13.289.29.1889&333&44&4500&275732261&4.9.5.11.13.37.4463&1621&2842 \cr
 & &64.9.11.37.1889&1889&1056& & &16.3.5.49.13.29.1621&1621&1560 \cr
\noalign{\hrule}
}%
}
$$
\eject
\vglue -23 pt
\noindent\hskip 1 in\hbox to 6.5 in{\ 4501 -- 4536 \hfill\fbd 275930655 -- 282339995\frb}
\vskip -9 pt
$$
\vbox{
\nointerlineskip
\halign{\strut
    \vrule \ \ \hfil \frb #\ 
   &\vrule \hfil \ \ \fbb #\frb\ 
   &\vrule \hfil \ \ \frb #\ \hfil
   &\vrule \hfil \ \ \frb #\ 
   &\vrule \hfil \ \ \frb #\ \ \vrule \hskip 2 pt
   &\vrule \ \ \hfil \frb #\ 
   &\vrule \hfil \ \ \fbb #\frb\ 
   &\vrule \hfil \ \ \frb #\ \hfil
   &\vrule \hfil \ \ \frb #\ 
   &\vrule \hfil \ \ \frb #\ \vrule \cr%
\noalign{\hrule}
 & &3.5.7.11.13.17.23.47&445&116& & &5.43.61.101.211&4587&4486 \cr
4501&275930655&8.25.13.23.29.89&99&476&4519&279493765&4.3.5.11.61.139.2243&9847&1368 \cr
 & &64.9.7.11.17.89&267&32& & &64.27.11.19.43.229&6183&6688 \cr
\noalign{\hrule}
 & &27.11.17.47.1163&793&370& & &27.5.11.29.43.151&2041&2338 \cr
4502&275983389&4.3.5.11.13.17.37.61&1163&58&4520&279621045&4.5.7.13.43.157.167&453&6298 \cr
 & &16.29.61.1163&29&488& & &16.3.13.47.67.151&871&376 \cr
\noalign{\hrule}
 & &5.11.83.103.587&1011&1924& & &25.11.13.17.43.107&747&1072 \cr
4503&276004465&8.3.13.37.103.337&189&292&4521&279625775&32.9.11.43.67.83&195&278 \cr
 & &64.81.7.73.337&24601&18144& & &128.27.5.13.67.139&3753&4288 \cr
\noalign{\hrule}
 & &9.5.7.13.19.53.67&43&62& & &121.29.173.461&2449&2622 \cr
4504&276285555&4.3.13.31.43.53.67&1067&1000&4522&279853277&4.3.11.19.23.29.31.79&1065&1384 \cr
 & &64.125.11.31.43.97&33325&34144& & &64.9.5.19.23.71.173&6745&6624 \cr
\noalign{\hrule}
 & &3.13.37.43.61.73&441&2698& & &3.25.13.41.47.149&321&854 \cr
4505&276304197&4.27.49.13.19.71&1595&5074&4523&279944925&4.9.7.61.107.149&1045&296 \cr
 & &16.5.11.29.43.59&295&2552& & &64.5.11.19.37.61&7733&1952 \cr
\noalign{\hrule}
 & &125.169.23.569&1881&2006& & &25.7.11.19.79.97&3639&4024 \cr
4506&276462875&4.9.11.17.19.59.569&845&276&4524&280274225&16.3.5.19.503.1213&7811&1746 \cr
 & &32.27.5.11.169.17.23&297&272& & &64.27.73.97.107&2889&2336 \cr
\noalign{\hrule}
 & &25.121.13.79.89&999&2026& & &7.11.71.51283&25647&25636 \cr
4507&276494075&4.27.37.89.1013&11&100&4525&280364161&8.3.7.13.17.29.71.83.103&1965&16378 \cr
 & &32.9.25.11.1013&1013&144& & &32.9.5.19.103.131.431&282305&281808 \cr
\noalign{\hrule}
 & &3.5.11.13.191.677&8579&18734& & &17.31.37.73.197&605&636 \cr
4508&277363515&4.17.19.23.29.373&9&382&4526&280415119&8.3.5.121.37.53.197&6851&438 \cr
 & &16.9.19.29.191&19&696& & &32.9.5.13.17.31.73&117&80 \cr
\noalign{\hrule}
 & &7.13.23.41.53.61&375&314& & &3.25.13.53.61.89&177&88 \cr
4509&277433429&4.3.125.7.23.41.157&183&22&4527&280543575&16.9.5.11.13.59.61&323&262 \cr
 & &16.9.25.11.61.157&2475&1256& & &64.11.17.19.59.131&27379&32096 \cr
\noalign{\hrule}
 & &27.7.11.13.43.239&1159&920& & &3.17.19.23.12601&6105&6496 \cr
4510&277756479&16.5.13.19.23.43.61&891&98&4528&280838487&64.9.5.7.11.19.29.37&6329&5474 \cr
 & &64.81.5.49.11.19&399&160& & &256.49.17.23.6329&6329&6272 \cr
\noalign{\hrule}
 & &27.5.7.491.599&16679&506& & &5.49.11.13.71.113&375&262 \cr
4511&277933005&4.11.13.23.1283&653&630&4529&281085805&4.3.625.11.71.131&703&78 \cr
 & &16.9.5.7.11.13.653&653&1144& & &16.9.13.19.37.131&4847&1368 \cr
\noalign{\hrule}
 & &9.5.11.17.173.191&18487&7982& & &9.5.7.121.47.157&931&884 \cr
4512&278056845&4.7.13.19.139.307&2013&3820&4530&281250585&8.3.343.13.17.19.157&1535&506 \cr
 & &32.3.5.7.11.61.191&61&112& & &32.5.11.17.19.23.307&5833&6256 \cr
\noalign{\hrule}
 & &9.149.283.733&1025&1174& & &81.5.11.233.271&157&76 \cr
4513&278175699&4.3.25.41.283.587&1043&1892&4531&281302065&8.5.11.19.157.271&387&658 \cr
 & &32.5.7.11.41.43.149&3157&3440& & &32.9.7.43.47.157&7379&4816 \cr
\noalign{\hrule}
 & &3.5.7.31.37.2311&1927&1958& & &27.11.31.127.241&4307&3164 \cr
4514&278325285&4.11.41.47.89.2311&3247&936&4532&281798649&8.3.7.11.59.73.113&155&494 \cr
 & &64.9.11.13.17.41.191&42211&43296& & &32.5.7.13.19.31.73&2555&3952 \cr
\noalign{\hrule}
 & &3.125.19.89.439&1771&454& & &27.7.1331.19.59&3583&3556 \cr
4515&278380875&4.5.7.11.19.23.227&979&1206&4533&281997639&8.49.11.19.127.3583&68265&188 \cr
 & &16.9.7.121.67.89&1407&968& & &64.9.5.37.41.47&1927&5920 \cr
\noalign{\hrule}
 & &5.7.11.13.19.29.101&53&24& & &3.17.23.37.67.97&653&1826 \cr
4516&278533255&16.3.5.13.19.53.101&6293&6798&4534&282063099&4.11.83.97.653&207&860 \cr
 & &64.9.7.11.29.31.103&927&992& & &32.9.5.23.43.83&415&2064 \cr
\noalign{\hrule}
 & &13.137.229.683&1147&1830& & &25.11.13.23.47.73&7209&6844 \cr
4517&278560867&4.3.5.31.37.61.137&231&916&4535&282113975&8.81.5.11.29.59.89&18907&24158 \cr
 & &32.9.7.11.61.229&549&1232& & &32.3.7.37.47.73.257&1799&1776 \cr
\noalign{\hrule}
 & &81.5.11.29.2161&635&1526& & &5.7.289.103.271&4887&4598 \cr
4518&279190395&4.25.7.29.109.127&73&102&4536&282339995&4.27.121.19.103.181&1073&11390 \cr
 & &16.3.17.73.109.127&13843&9928& & &16.9.5.17.29.37.67&1943&2664 \cr
\noalign{\hrule}
}%
}
$$
\eject
\vglue -23 pt
\noindent\hskip 1 in\hbox to 6.5 in{\ 4537 -- 4572 \hfill\fbd 282565283 -- 287467739\frb}
\vskip -9 pt
$$
\vbox{
\nointerlineskip
\halign{\strut
    \vrule \ \ \hfil \frb #\ 
   &\vrule \hfil \ \ \fbb #\frb\ 
   &\vrule \hfil \ \ \frb #\ \hfil
   &\vrule \hfil \ \ \frb #\ 
   &\vrule \hfil \ \ \frb #\ \ \vrule \hskip 2 pt
   &\vrule \ \ \hfil \frb #\ 
   &\vrule \hfil \ \ \fbb #\frb\ 
   &\vrule \hfil \ \ \frb #\ \hfil
   &\vrule \hfil \ \ \frb #\ 
   &\vrule \hfil \ \ \frb #\ \vrule \cr%
\noalign{\hrule}
 & &7.11.13.19.83.179&629&450& & &27.7.11.23.67.89&415&1394 \cr
4537&282565283&4.9.25.7.11.17.19.37&181&104&4555&285132771&4.5.7.17.23.41.83&429&268 \cr
 & &64.3.5.13.17.37.181&6697&8160& & &32.3.5.11.13.67.83&415&208 \cr
\noalign{\hrule}
 & &125.49.11.13.17.19&6739&5364& & &3.125.11.23.31.97&1207&1218 \cr
4538&282907625&8.9.17.23.149.293&49&100&4556&285289125&4.9.5.7.17.23.29.31.71&95161&34556 \cr
 & &64.3.25.49.23.293&879&736& & &32.11.41.53.163.211&34393&34768 \cr
\noalign{\hrule}
 & &81.11.41.61.127&973&424& & &3.5.11.19.181.503&63&118 \cr
4539&283005657&16.9.7.41.53.139&115&254&4557&285419805&4.27.7.19.59.503&5575&3982 \cr
 & &64.5.7.23.53.127&805&1696& & &16.25.7.11.181.223&223&280 \cr
\noalign{\hrule}
 & &25.7.11.13.43.263&1219&1674& & &27.37.43.61.109&1&110 \cr
4540&283007725&4.27.5.23.31.43.53&6487&4208&4558&285621093&4.9.5.11.43.61&253&296 \cr
 & &128.9.13.263.499&499&576& & &64.5.121.23.37&2783&160 \cr
\noalign{\hrule}
 & &9.5.11.13.29.37.41&391&798& & &5.11.19.167.1637&1341&296 \cr
4541&283094955&4.27.5.7.13.17.19.23&1073&682&4559&285681055&16.9.37.149.167&325&176 \cr
 & &16.7.11.19.29.31.37&217&152& & &512.3.25.11.13.37&2405&768 \cr
\noalign{\hrule}
 & &3.25.13.17.19.29.31&407&492& & &49.11.13.361.113&755&1392 \cr
4542&283117575&8.9.5.11.13.19.37.41&97&682&4560&285836551&32.3.5.11.19.29.151&1667&2712 \cr
 & &32.121.31.37.97&3589&1936& & &512.9.113.1667&1667&2304 \cr
\noalign{\hrule}
 & &19.31.557.863&8085&18668& & &3.5.19.61.109.151&48433&51302 \cr
4543&283126999&8.3.5.49.11.13.359&173&186&4561&286139715&4.7.11.17.37.113.227&151&8550 \cr
 & &32.9.5.49.11.31.173&8477&7920& & &16.9.25.17.19.151&85&24 \cr
\noalign{\hrule}
 & &5.13.23.131.1447&94031&95526& & &3.29.103.109.293&20449&11488 \cr
4544&283387715&4.27.49.19.29.61.101&1483&286&4562&286187457&64.121.169.359&2059&1890 \cr
 & &16.3.7.11.13.101.1483&16313&16968& & &256.27.5.7.11.29.71&5467&5760 \cr
\noalign{\hrule}
 & &3.7.11.17.23.43.73&771&470& & &81.5.169.47.89&32131&36314 \cr
4545&283517619&4.9.5.11.23.47.257&3451&6278&4563&286305435&4.11.23.67.127.271&169&102 \cr
 & &16.5.7.17.29.43.73&29&40& & &16.3.11.169.17.23.127&2159&2024 \cr
\noalign{\hrule}
 & &3.7.11.23.83.643&2665&4408& & &5.7.13.359.1753&17677&5112 \cr
4546&283549497&16.5.13.19.23.29.41&747&196&4564&286343785&16.9.11.71.1607&697&910 \cr
 & &128.9.5.49.13.83&195&448& & &64.3.5.7.11.13.17.41&697&1056 \cr
\noalign{\hrule}
 & &3.11.13.239.2767&2937&170& & &243.125.11.857&1399&1274 \cr
4547&283703277&4.9.5.121.17.89&2767&2678&4565&286345125&4.49.13.857.1399&3699&2300 \cr
 & &16.13.17.103.2767&103&136& & &32.27.25.7.13.23.137&2093&2192 \cr
\noalign{\hrule}
 & &81.17.31.61.109&19&1872& & &5.7.11.13.19.23.131&567&698 \cr
4548&283825863&32.729.13.19&355&374&4566&286521235&4.81.49.13.19.349&575&62 \cr
 & &128.5.11.13.17.71&781&4160& & &16.3.25.23.31.349&1745&744 \cr
\noalign{\hrule}
 & &81.5.11.17.23.163&167&86& & &3.7.191.199.359&102157&103550 \cr
4549&283930515&4.5.17.43.163.167&1827&1012&4567&286549851&4.25.11.19.37.109.251&9&28 \cr
 & &32.9.7.11.23.29.43&301&464& & &32.9.25.7.11.109.251&18825&19184 \cr
\noalign{\hrule}
 & &27.7.37.179.227&143&116& & &9.121.271.971&317&46 \cr
4550&284146569&8.11.13.29.179.227&85&2412&4568&286560549&4.3.23.317.971&451&520 \cr
 & &64.9.5.17.29.67&5695&928& & &64.5.11.13.41.317&4121&6560 \cr
\noalign{\hrule}
 & &27.11.13.17.61.71&161&620& & &49.23.43.61.97&2825&6996 \cr
4551&284273847&8.5.7.13.23.31.61&119&180&4569&286743737&8.3.25.7.11.53.113&97&468 \cr
 & &64.9.25.49.17.31&775&1568& & &64.27.5.11.13.97&1755&352 \cr
\noalign{\hrule}
 & &9.11.23.29.31.139&799&730& & &3.25.7.13.23.31.59&59&214 \cr
4552&284536197&4.3.5.17.29.31.47.73&2495&12788&4570&287107275&4.5.23.3481.107&1683&1798 \cr
 & &32.25.23.139.499&499&400& & &16.9.11.17.29.31.107&3531&3944 \cr
\noalign{\hrule}
 & &9.5.19.23.43.337&2211&1874& & &19.23.47.71.197&4411&4848 \cr
4553&284965515&4.27.11.23.67.937&779&158&4571&287278993&32.3.11.71.101.401&591&190 \cr
 & &16.11.19.41.67.79&5293&3608& & &128.9.5.19.101.197&505&576 \cr
\noalign{\hrule}
 & &25.7.11.17.31.281&16549&16890& & &13.17.19.223.307&121&102 \cr
4554&285067475&4.3.125.13.19.67.563&7847&528&4572&287467739&4.3.121.13.289.307&3567&190 \cr
 & &128.9.7.11.361.59&3249&3776& & &16.9.5.11.19.29.41&4059&1160 \cr
\noalign{\hrule}
}%
}
$$
\eject
\vglue -23 pt
\noindent\hskip 1 in\hbox to 6.5 in{\ 4573 -- 4608 \hfill\fbd 287476609 -- 295330035\frb}
\vskip -9 pt
$$
\vbox{
\nointerlineskip
\halign{\strut
    \vrule \ \ \hfil \frb #\ 
   &\vrule \hfil \ \ \fbb #\frb\ 
   &\vrule \hfil \ \ \frb #\ \hfil
   &\vrule \hfil \ \ \frb #\ 
   &\vrule \hfil \ \ \frb #\ \ \vrule \hskip 2 pt
   &\vrule \ \ \hfil \frb #\ 
   &\vrule \hfil \ \ \fbb #\frb\ 
   &\vrule \hfil \ \ \frb #\ \hfil
   &\vrule \hfil \ \ \frb #\ 
   &\vrule \hfil \ \ \frb #\ \vrule \cr%
\noalign{\hrule}
 & &7.23.31.239.241&28897&23400& & &9.5.11.31.67.283&131&472 \cr
4573&287476609&16.9.25.11.13.37.71&241&2386&4591&290956545&16.5.59.131.283&469&186 \cr
 & &64.3.5.241.1193&1193&480& & &64.3.7.31.59.67&59&224 \cr
\noalign{\hrule}
 & &27.25.11.17.43.53&277&148& & &25.7.11.19.73.109&999&926 \cr
4574&287666775&8.9.11.37.53.277&125&458&4592&291027275&4.27.19.37.109.463&17885&754 \cr
 & &32.125.229.277&1385&3664& & &16.3.5.49.13.29.73&377&168 \cr
\noalign{\hrule}
 & &3.7.13.41.47.547&17411&3922& & &9.121.127.2111&3151&1040 \cr
4575&287760837&4.23.37.53.757&405&352&4593&291957633&32.3.5.11.13.23.137&511&1270 \cr
 & &256.81.5.11.23.37&10989&14720& & &128.25.7.73.127&511&1600 \cr
\noalign{\hrule}
 & &3.5.23.37.107.211&13&198& & &9.5.11.13.289.157&437&148 \cr
4576&288195405&4.27.11.13.23.107&1055&1406&4594&291975255&8.11.19.23.37.157&2325&2482 \cr
 & &16.5.11.19.37.211&19&88& & &32.3.25.17.31.37.73&2701&2480 \cr
\noalign{\hrule}
 & &27.5.7.47.73.89&2159&2024& & &5.19.53.59.983&693&428 \cr
4577&288564255&16.7.11.17.23.73.127&43&846&4595&292014895&8.9.7.11.107.983&1675&5206 \cr
 & &64.9.17.23.43.47&989&544& & &32.3.25.19.67.137&2055&1072 \cr
\noalign{\hrule}
 & &25.169.23.2971&36279&32054& & &9.11.23.181.709&95&86 \cr
4578&288706925&4.9.11.29.31.47.139&41&52&4596&292205133&4.5.11.19.23.43.709&481&228 \cr
 & &32.3.13.29.41.47.139&19599&19024& & &32.3.5.13.361.37.43&23465&25456 \cr
\noalign{\hrule}
 & &27.7.11.13.289.37&89&100& & &27.25.13.17.37.53&453&28 \cr
4579&288999711&8.25.13.289.37.89&8239&2454&4597&292532175&8.81.7.53.151&4285&3718 \cr
 & &32.3.5.7.11.107.409&2045&1712& & &32.5.11.169.857&857&2288 \cr
\noalign{\hrule}
 & &27.5.19.137.823&259&254& & &9.19.23.163.457&5075&3608 \cr
4580&289206315&4.7.37.127.137.823&891&68&4598&292973103&16.25.7.11.23.29.41&19&96 \cr
 & &32.81.11.17.37.127&6919&6096& & &1024.3.5.19.29.41&1189&2560 \cr
\noalign{\hrule}
 & &3.5.11.13.19.47.151&1407&1462& & &11.19.37.137.277&2563&2700 \cr
4581&289238235&4.9.7.13.17.43.47.67&547&6644&4599&293459617&8.27.25.121.37.233&89&274 \cr
 & &32.11.43.151.547&547&688& & &32.9.5.89.137.233&4005&3728 \cr
\noalign{\hrule}
 & &25.7.11.13.31.373&289&114& & &25.121.37.43.61&2449&2754 \cr
4582&289364075&4.3.11.289.19.373&93&280&4600&293579275&4.81.5.17.31.37.79&2101&5246 \cr
 & &64.9.5.7.17.19.31&171&544& & &16.27.11.43.61.191&191&216 \cr
\noalign{\hrule}
 & &81.11.157.2069&1085&3154& & &19.37.67.6241&825&2098 \cr
4583&289426203&4.3.5.7.11.19.31.83&157&74&4601&293957341&4.3.25.11.79.1049&327&722 \cr
 & &16.5.19.31.37.157&589&1480& & &16.9.5.11.361.109&1881&4360 \cr
\noalign{\hrule}
 & &3.5.49.17.97.239&9139&11176& & &9.5.49.11.17.23.31&2291&2452 \cr
4584&289671585&16.7.11.13.19.37.127&549&340&4602&293994855&8.5.7.11.29.79.613&27&4318 \cr
 & &128.9.5.13.17.37.61&2257&2496& & &32.27.17.29.127&127&1392 \cr
\noalign{\hrule}
 & &3.5.23.67.83.151&14129&13676& & &25.13.19.29.31.53&1067&168 \cr
4585&289700295&8.13.23.71.199.263&36783&15554&4603&294220225&16.3.5.7.11.53.97&1989&1406 \cr
 & &32.9.7.11.61.67.101&4697&4848& & &64.27.13.17.19.37&459&1184 \cr
\noalign{\hrule}
 & &5.7.13.59.101.107&737&30& & &11.13.59.139.251&21663&13226 \cr
4586&290113915&4.3.25.11.67.107&14679&14746&4604&294358493&4.9.17.29.83.389&35&52 \cr
 & &16.27.7.73.101.233&1971&1864& & &32.3.5.7.13.83.389&8169&6640 \cr
\noalign{\hrule}
 & &9.7.13.23.73.211&255&44& & &125.17.19.23.317&4897&4878 \cr
4587&290146311&8.27.5.7.11.17.73&107&404&4605&294374125&4.9.5.59.83.271.317&79439&506 \cr
 & &64.5.17.101.107&1717&17120& & &16.3.11.19.23.37.113&1221&904 \cr
\noalign{\hrule}
 & &81.5.7.17.19.317&121&202& & &5.7.11.431.1777&2231&17316 \cr
4588&290278485&4.5.7.121.101.317&351&34&4606&294866495&8.9.13.23.37.97&673&770 \cr
 & &16.27.11.13.17.101&143&808& & &32.3.5.7.11.23.673&673&1104 \cr
\noalign{\hrule}
 & &3.7.23.59.61.167&589&650& & &11.17.19.61.1361&623&414 \cr
4589&290298939&4.25.13.19.23.31.167&7&2178&4607&294973613&4.9.7.23.89.1361&1615&254 \cr
 & &16.9.5.7.121.31&121&3720& & &16.3.5.17.19.23.127&345&1016 \cr
\noalign{\hrule}
 & &3.169.19.109.277&3245&2968& & &3.5.7.11.169.17.89&399&1114 \cr
4590&290849169&16.5.7.11.169.53.59&2493&634&4608&295330035&4.9.49.13.19.557&3841&3400 \cr
 & &64.9.5.7.277.317&951&1120& & &64.25.17.19.23.167&3173&3680 \cr
\noalign{\hrule}
}%
}
$$
\eject
\vglue -23 pt
\noindent\hskip 1 in\hbox to 6.5 in{\ 4609 -- 4644 \hfill\fbd 295500205 -- 301167801\frb}
\vskip -9 pt
$$
\vbox{
\nointerlineskip
\halign{\strut
    \vrule \ \ \hfil \frb #\ 
   &\vrule \hfil \ \ \fbb #\frb\ 
   &\vrule \hfil \ \ \frb #\ \hfil
   &\vrule \hfil \ \ \frb #\ 
   &\vrule \hfil \ \ \frb #\ \ \vrule \hskip 2 pt
   &\vrule \ \ \hfil \frb #\ 
   &\vrule \hfil \ \ \fbb #\frb\ 
   &\vrule \hfil \ \ \frb #\ \hfil
   &\vrule \hfil \ \ \frb #\ 
   &\vrule \hfil \ \ \frb #\ \vrule \cr%
\noalign{\hrule}
 & &5.7.11.13.17.23.151&171&886& & &3.5.13.19.37.41.53&43&80 \cr
4609&295500205&4.9.17.19.23.443&4379&3050&4627&297885705&32.25.13.19.43.53&341&666 \cr
 & &16.3.25.29.61.151&305&696& & &128.9.11.31.37.43&1333&2112 \cr
\noalign{\hrule}
 & &9.11.17.29.73.83&9983&9070& & &5.11.17.37.79.109&3303&3412 \cr
4610&295721613&4.5.17.67.149.907&1023&116&4628&297897545&8.9.11.37.367.853&347&754 \cr
 & &32.3.5.11.29.31.149&745&496& & &32.3.13.29.347.853&74211&72176 \cr
\noalign{\hrule}
 & &5.7.31.227.1201&21131&20904& & &25.49.23.71.149&13&162 \cr
4611&295800295&16.3.11.13.17.31.67.113&177&2254&4629&298063325&4.81.7.13.23.71&703&220 \cr
 & &64.9.49.23.59.113&18193&16992& & &32.27.5.11.19.37&7733&432 \cr
\noalign{\hrule}
 & &27.5.11.89.2239&377&1862& & &7.13.29.37.43.71&1089&17300 \cr
4612&295917435&4.49.13.19.29.89&755&666&4630&298104079&8.9.25.121.173&59&62 \cr
 & &16.9.5.13.19.37.151&2869&3848& & &32.3.25.31.59.173&16089&23600 \cr
\noalign{\hrule}
 & &7.13.841.53.73&765&76& & &5.17.19.29.6367&1369&1386 \cr
4613&296098439&8.9.5.7.17.19.73&33&40&4631&298198445&4.9.7.11.1369.6367&44873&304 \cr
 & &128.27.25.11.17.19&8075&19008& & &128.3.19.23.1951&1951&4416 \cr
\noalign{\hrule}
 & &7.13.841.53.73&33&40& & &5.11.289.19.23.43&36537&36580 \cr
4614&296098439&16.3.5.11.13.841.53&765&76&4632&298682945&8.3.25.361.31.59.641&5423&14448 \cr
 & &128.27.25.11.17.19&8075&19008& & &256.9.7.11.17.29.43.59&3717&3712 \cr
\noalign{\hrule}
 & &9.25.17.29.2671&1927&4598& & &9.5.49.23.43.137&209&232 \cr
4615&296280675&4.121.17.19.41.47&1057&870&4633&298762065&16.5.11.19.29.43.137&119&2484 \cr
 & &16.3.5.7.11.19.29.151&1463&1208& & &128.27.7.17.23.29&493&192 \cr
\noalign{\hrule}
 & &9.7.13.19.137.139&555&418& & &9.1331.13.19.101&1873&2120 \cr
4616&296328123&4.27.5.11.13.361.37&8479&4508&4634&298840113&16.3.5.53.101.1873&133&5486 \cr
 & &32.5.49.23.61.139&805&976& & &64.5.7.13.19.211&211&1120 \cr
\noalign{\hrule}
 & &9.25.7.13.41.353&1507&3082& & &9.125.11.19.31.41&123&218 \cr
4617&296334675&4.11.23.41.67.137&1247&1500&4635&298843875&4.27.25.1681.109&503&1178 \cr
 & &32.3.125.29.43.137&3973&3440& & &16.19.31.109.503&503&872 \cr
\noalign{\hrule}
 & &3.25.29.31.53.83&709&616& & &3.5.11.17.197.541&22477&14362 \cr
4618&296602575&16.7.11.29.83.709&355&558&4636&298948485&4.7.169.19.43.167&1089&1082 \cr
 & &64.9.5.31.71.709&2127&2272& & &16.9.121.13.19.43.541&1677&1672 \cr
\noalign{\hrule}
 & &9.25.53.149.167&6479&18404& & &9.625.23.2311&1343&968 \cr
4619&296729775&8.11.19.31.43.107&2067&1250&4637&298985625&16.3.5.121.17.23.79&2311&494 \cr
 & &32.3.625.11.13.53&325&176& & &64.11.13.19.2311&209&416 \cr
\noalign{\hrule}
 & &3.11.13.17.23.29.61&193&106& & &3.25.29.79.1741&25199&18326 \cr
4620&296730291&4.11.17.53.61.193&1305&1976&4638&299147325&4.49.11.17.113.223&1741&180 \cr
 & &64.9.5.13.19.29.53&795&608& & &32.9.5.7.11.1741&33&112 \cr
\noalign{\hrule}
 & &27.7.13.31.47.83&11&2& & &9.121.41.53.127&703&650 \cr
4621&297127467&4.3.7.11.31.47.83&793&1780&4639&300532419&4.3.25.11.13.19.37.127&53&328 \cr
 & &32.5.11.13.61.89&979&4880& & &64.13.19.37.41.53&481&608 \cr
\noalign{\hrule}
 & &5.7.11.19.97.419&1387&708& & &5.7.11.41.137.139&213&74 \cr
4622&297303545&8.3.11.361.59.73&221&582&4640&300593755&4.3.5.11.37.71.137&139&546 \cr
 & &32.9.13.17.59.97&1989&944& & &16.9.7.13.71.139&117&568 \cr
\noalign{\hrule}
 & &3.5.7.73.79.491&137&374& & &3.5.7.11.43.73.83&689&1054 \cr
4623&297317685&4.5.11.17.137.491&1159&1296&4641&300920235&4.11.13.17.31.43.53&657&74 \cr
 & &128.81.11.17.19.61&18117&20672& & &16.9.13.31.37.73&1209&296 \cr
\noalign{\hrule}
 & &9.5.11.19.29.1091&8023&8342& & &3.5.7.11.43.73.83&10127&7718 \cr
4624&297564795&4.3.19.43.71.97.113&1091&752&4642&300920235&4.7.13.17.19.41.227&10043&19350 \cr
 & &128.43.47.71.1091&3053&3008& & &16.9.25.121.43.83&33&40 \cr
\noalign{\hrule}
 & &3.5.7.11.439.587&3009&3448& & &9.25.11.19.37.173&6731&3444 \cr
4625&297635415&16.9.5.7.17.59.431&6187&13208&4643&301007025&8.27.7.41.53.127&2405&7612 \cr
 & &256.13.23.127.269&34163&38272& & &64.5.11.13.37.173&13&32 \cr
\noalign{\hrule}
 & &81.5.11.13.53.97&931&136& & &81.11.17.59.337&871&812 \cr
4626&297741015&16.27.49.13.17.19&97&124&4644&301167801&8.9.7.13.29.67.337&545&2488 \cr
 & &128.49.19.31.97&931&1984& & &128.5.7.13.109.311&33899&29120 \cr
\noalign{\hrule}
}%
}
$$
\eject
\vglue -23 pt
\noindent\hskip 1 in\hbox to 6.5 in{\ 4645 -- 4680 \hfill\fbd 301215537 -- 308289575\frb}
\vskip -9 pt
$$
\vbox{
\nointerlineskip
\halign{\strut
    \vrule \ \ \hfil \frb #\ 
   &\vrule \hfil \ \ \fbb #\frb\ 
   &\vrule \hfil \ \ \frb #\ \hfil
   &\vrule \hfil \ \ \frb #\ 
   &\vrule \hfil \ \ \frb #\ \ \vrule \hskip 2 pt
   &\vrule \ \ \hfil \frb #\ 
   &\vrule \hfil \ \ \fbb #\frb\ 
   &\vrule \hfil \ \ \frb #\ \hfil
   &\vrule \hfil \ \ \frb #\ 
   &\vrule \hfil \ \ \frb #\ \vrule \cr%
\noalign{\hrule}
 & &27.7.17.241.389&1025&1144& & &5.13.19.37.41.163&7983&7502 \cr
4645&301215537&16.3.25.11.13.41.389&71&1096&4663&305379685&4.9.121.31.41.887&6307&21190 \cr
 & &256.11.13.71.137&19591&9088& & &16.3.5.7.13.17.53.163&371&408 \cr
\noalign{\hrule}
 & &17.23.37.83.251&11583&20870& & &3.13.19.29.41.347&1661&2850 \cr
4646&301391011&4.81.5.11.13.2087&829&1258&4664&305724003&4.9.25.11.361.151&2665&1306 \cr
 & &16.27.5.17.37.829&829&1080& & &16.125.13.41.653&653&1000 \cr
\noalign{\hrule}
 & &9.41.43.83.229&6625&3222& & &27.5.19.23.71.73&137&502 \cr
4647&301584069&4.81.125.53.179&229&176&4665&305771085&4.3.19.23.137.251&781&530 \cr
 & &128.25.11.179.229&1969&1600& & &16.5.11.53.71.137&583&1096 \cr
\noalign{\hrule}
 & &9.7.11.31.101.139&101&178& & &9.5.7.121.71.113&601&246 \cr
4648&301599837&4.89.10201.139&1085&11286&4666&305796645&4.27.41.113.601&497&610 \cr
 & &16.27.5.7.11.19.31&19&120& & &16.5.7.61.71.601&601&488 \cr
\noalign{\hrule}
 & &3.13.19.29.101.139&12773&21560& & &9.7.11.37.79.151&731&1600 \cr
4649&301684071&16.5.49.11.53.241&65&306&4667&305871489&128.25.17.43.151&97&54 \cr
 & &64.9.25.7.11.13.17&2975&1056& & &512.27.25.17.97&7275&4352 \cr
\noalign{\hrule}
 & &243.5.7.17.2087&79&164& & &27.11.13.19.43.97&1211&50 \cr
4650&301748895&8.7.41.79.2087&767&1320&4668&305980389&4.25.7.11.19.173&1649&1638 \cr
 & &128.3.5.11.13.41.59&5863&3776& & &16.9.25.49.13.17.97&425&392 \cr
\noalign{\hrule}
 & &81.25.11.71.191&2027&17498& & &3.7.17.31.89.311&397&130 \cr
4651&302071275&4.13.673.2027&677&1350&4669&306323493&4.5.7.13.311.397&1287&890 \cr
 & &16.27.25.13.677&677&104& & &16.9.25.11.169.89&825&1352 \cr
\noalign{\hrule}
 & &3.49.31.151.439&767&550& & &27.11.67.73.211&507&296 \cr
4652&302078973&4.25.7.11.13.59.151&1863&202&4670&306504297&16.81.169.37.67&413&5840 \cr
 & &16.81.5.13.23.101&3105&10504& & &512.5.7.59.73&2065&256 \cr
\noalign{\hrule}
 & &47.1283.5011&27645&32656& & &9.5.7.23.101.419&1639&1294 \cr
4653&302168311&32.3.5.13.19.97.157&1443&1540&4671&306601155&4.3.11.101.149.647&475&172 \cr
 & &256.9.25.7.11.169.37&117117&118400& & &32.25.11.19.43.149&8987&11920 \cr
\noalign{\hrule}
 & &27.11.37.79.349&1085&1048& & &27.11.13.19.37.113&1411&17570 \cr
4654&302977719&16.5.7.11.31.131.349&111&3950&4672&306713979&4.5.7.17.83.251&333&248 \cr
 & &64.3.125.7.37.79&125&224& & &64.9.31.37.251&251&992 \cr
\noalign{\hrule}
 & &343.43.59.349&2615&17622& & &7.529.79.1049&605&444 \cr
4655&303696659&4.9.5.11.89.523&295&1274&4673&306871313&8.3.5.121.23.37.79&3147&1330 \cr
 & &16.3.25.49.13.59&325&24& & &32.9.25.7.19.1049&225&304 \cr
\noalign{\hrule}
 & &27.5.49.19.41.59&15103&23068& & &11.31.37.97.251&141&110 \cr
4656&304032015&8.11.73.79.1373&1121&252&4674&307186099&4.3.5.121.37.47.97&41&4518 \cr
 & &64.9.7.19.59.73&73&32& & &16.27.5.41.251&1107&40 \cr
\noalign{\hrule}
 & &27.49.47.67.73&8473&5324& & &3.5.13.17.23.29.139&931&1070 \cr
4657&304127271&8.7.1331.37.229&1005&598&4675&307343595&4.25.49.13.17.19.107&23&198 \cr
 & &32.3.5.121.13.23.67&1495&1936& & &16.9.7.11.19.23.107&1463&2568 \cr
\noalign{\hrule}
 & &3.5.7.43.89.757&1581&2204& & &3.11.19.43.101.113&1581&338 \cr
4658&304189095&8.9.17.19.29.31.43&539&278&4676&307705893&4.9.169.17.31.43&505&226 \cr
 & &32.49.11.17.31.139&10703&8432& & &16.5.169.101.113&169&40 \cr
\noalign{\hrule}
 & &3.23.37.43.47.59&697&660& & &3.5.7.11.169.19.83&827&86 \cr
4659&304417167&8.9.5.11.17.41.43.47&1423&692&4677&307822515&4.5.7.13.43.827&14193&14752 \cr
 & &64.11.41.173.1423&58343&60896& & &256.9.19.83.461&461&384 \cr
\noalign{\hrule}
 & &3.7.23.43.107.137&125&286& & &5.11.13.47.89.103&2349&1834 \cr
4660&304452771&4.125.11.13.43.107&1233&158&4678&308057035&4.81.7.11.13.29.131&89&1090 \cr
 & &16.9.5.11.79.137&869&120& & &16.9.5.29.89.109&261&872 \cr
\noalign{\hrule}
 & &3.25.13.19.41.401&363&38& & &5.7.121.83.877&1891&1014 \cr
4661&304569525&4.9.121.361.41&1645&1604&4679&308269885&4.3.121.169.31.61&145&24 \cr
 & &32.5.7.121.47.401&847&752& & &64.9.5.29.31.61&8091&1952 \cr
\noalign{\hrule}
 & &9.7.17.461.617&171&290& & &25.11.841.31.43&1273&432 \cr
4662&304632027&4.81.5.19.29.617&461&1078&4680&308289575&32.27.5.19.43.67&3509&3844 \cr
 & &16.5.49.11.29.461&385&232& & &256.3.121.29.961&341&384 \cr
\noalign{\hrule}
}%
}
$$
\eject
\vglue -23 pt
\noindent\hskip 1 in\hbox to 6.5 in{\ 4681 -- 4716 \hfill\fbd 308339605 -- 313137475\frb}
\vskip -9 pt
$$
\vbox{
\nointerlineskip
\halign{\strut
    \vrule \ \ \hfil \frb #\ 
   &\vrule \hfil \ \ \fbb #\frb\ 
   &\vrule \hfil \ \ \frb #\ \hfil
   &\vrule \hfil \ \ \frb #\ 
   &\vrule \hfil \ \ \frb #\ \ \vrule \hskip 2 pt
   &\vrule \ \ \hfil \frb #\ 
   &\vrule \hfil \ \ \fbb #\frb\ 
   &\vrule \hfil \ \ \frb #\ \hfil
   &\vrule \hfil \ \ \frb #\ 
   &\vrule \hfil \ \ \frb #\ \vrule \cr%
\noalign{\hrule}
 & &5.49.59.83.257&143&438& & &3.11.13.59.71.173&20387&10234 \cr
4681&308339605&4.3.7.11.13.73.257&83&174&4699&310895013&4.7.17.19.29.37.43&641&90 \cr
 & &16.9.11.29.73.83&803&2088& & &16.9.5.7.37.641&1923&10360 \cr
\noalign{\hrule}
 & &3.37.71.109.359&125&234& & &25.7.11.13.289.43&933&2378 \cr
4682&308391411&4.27.125.13.37.71&1199&718&4700&310985675&4.3.5.13.29.41.311&291&86 \cr
 & &16.125.11.109.359&125&88& & &16.9.43.97.311&2799&776 \cr
\noalign{\hrule}
 & &49.121.61.853&85&36& & &27.25.7.13.61.83&19&16 \cr
4683&308503657&8.9.5.17.61.853&1127&1432&4701&310994775&32.9.5.13.19.61.83&7511&374 \cr
 & &128.3.49.17.23.179&4117&3264& & &128.7.11.17.29.37&6919&1856 \cr
\noalign{\hrule}
 & &3.37.53.137.383&2343&2726& & &5.11.59.239.401&273&922 \cr
4684&308686893&4.9.11.29.47.53.71&955&4292&4702&310997555&4.3.7.13.401.461&1593&4400 \cr
 & &32.5.841.37.191&4205&3056& & &128.81.25.11.59&81&320 \cr
\noalign{\hrule}
 & &3.25.7.13.31.1459&1217&242& & &81.25.11.61.229&3043&1898 \cr
4685&308687925&4.7.121.31.1217&585&1802&4703&311159475&4.5.11.13.17.73.179&1537&432 \cr
 & &16.9.5.11.13.17.53&187&1272& & &128.27.29.53.73&3869&1856 \cr
\noalign{\hrule}
 & &25.11.17.29.43.53&1823&1092& & &9.121.17.67.251&437&316 \cr
4686&308975425&8.3.5.7.13.29.1823&1419&404&4704&311333121&8.3.17.19.23.67.79&1255&286 \cr
 & &64.9.11.13.43.101&909&416& & &32.5.11.13.79.251&395&208 \cr
\noalign{\hrule}
 & &121.19.43.53.59&18081&25220& & &3.5.17.19.239.269&4557&16 \cr
4687&309125839&8.9.5.49.13.41.97&23&26&4705&311489895&32.9.5.49.31&737&782 \cr
 & &32.3.5.169.23.41.97&59655&62192& & &128.11.17.23.67&253&4288 \cr
\noalign{\hrule}
 & &3.11.23.37.101.109&1995&1742& & &9.5.7.13.29.43.61&1093&1702 \cr
4688&309165747&4.9.5.7.13.19.67.109&407&1010&4706&311494365&4.3.23.37.61.1093&455&638 \cr
 & &16.25.7.11.19.37.101&175&152& & &16.5.7.11.13.23.29.37&253&296 \cr
\noalign{\hrule}
 & &3.13.73.313.347&137&210& & &11.59.449.1069&5909&5850 \cr
4689&309215517&4.9.5.7.13.137.313&803&1388&4707&311507669&4.9.25.13.19.311.449&413&6322 \cr
 & &32.11.73.137.347&137&176& & &16.3.5.7.13.29.59.109&5655&6104 \cr
\noalign{\hrule}
 & &3.7.121.13.17.19.29&755&92& & &3.5.109.139.1373&339&1034 \cr
4690&309420111&8.5.19.23.29.151&351&200&4708&312034845&4.9.11.47.109.113&131&1112 \cr
 & &128.27.125.13.23&2875&576& & &64.47.131.139&131&1504 \cr
\noalign{\hrule}
 & &9.625.7.29.271&7843&10282& & &17.19.479.2017&847&1170 \cr
4691&309448125&4.7.11.23.31.53.97&3125&4344&4709&312064189&4.9.5.7.121.13.479&769&2584 \cr
 & &64.3.3125.31.181&905&992& & &64.3.13.17.19.769&769&1248 \cr
\noalign{\hrule}
 & &81.11.41.61.139&115&254& & &49.13.31.97.163&18975&20116 \cr
4692&309746349&4.9.5.11.23.61.127&973&424&4710&312219817&8.3.25.7.11.23.47.107&3077&2328 \cr
 & &64.5.7.23.53.139&805&1696& & &128.9.5.11.17.97.181&9955&9792 \cr
\noalign{\hrule}
 & &9.7.121.13.53.59&19781&20870& & &5.49.11.17.19.359&1629&166 \cr
4693&309882573&4.5.7.131.151.2087&81&836&4711&312504115&4.9.7.17.83.181&4667&5210 \cr
 & &32.81.11.19.2087&2087&2736& & &16.3.5.13.359.521&521&312 \cr
\noalign{\hrule}
 & &25.11.19.23.29.89&291&146& & &17.29.31.97.211&17875&2592 \cr
4694&310171675&4.3.5.11.73.89.97&247&732&4712&312797161&64.81.125.11.13&31&86 \cr
 & &32.9.13.19.61.73&4453&1872& & &256.9.25.31.43&387&3200 \cr
\noalign{\hrule}
 & &5.7.59.359.419&75181&73086& & &11.13.17.23.29.193&1095&1414 \cr
4695&310619365&4.3.13.937.75181&475&462&4713&312945061&4.3.5.7.17.23.73.101&13&378 \cr
 & &16.9.25.7.11.19.75181&75181&75240& & &16.81.49.13.101&4949&648 \cr
\noalign{\hrule}
 & &3.11.13.79.89.103&4185&2846& & &7.121.113.3271&18285&4612 \cr
4696&310678797&4.81.5.11.31.1423&1157&266&4714&313070681&8.3.5.23.53.1153&565&588 \cr
 & &16.5.7.13.19.31.89&589&280& & &64.9.25.49.53.113&1575&1696 \cr
\noalign{\hrule}
 & &9.5.7.19.23.37.61&385&52& & &27.5.49.19.47.53&715&292 \cr
4697&310687335&8.25.49.11.13.61&1081&444&4715&313081335&8.3.25.49.11.13.73&1363&2312 \cr
 & &64.3.11.23.37.47&47&352& & &128.11.289.29.47&3179&1856 \cr
\noalign{\hrule}
 & &9.5.43.347.463&2371&2834& & &25.7.37.137.353&639&286 \cr
4698&310879035&4.3.13.43.109.2371&347&2024&4716&313137475&4.9.7.11.13.71.137&295&706 \cr
 & &64.11.23.109.347&1199&736& & &16.3.5.59.71.353&213&472 \cr
\noalign{\hrule}
}%
}
$$
\eject
\vglue -23 pt
\noindent\hskip 1 in\hbox to 6.5 in{\ 4717 -- 4752 \hfill\fbd 313168925 -- 317578275\frb}
\vskip -9 pt
$$
\vbox{
\nointerlineskip
\halign{\strut
    \vrule \ \ \hfil \frb #\ 
   &\vrule \hfil \ \ \fbb #\frb\ 
   &\vrule \hfil \ \ \frb #\ \hfil
   &\vrule \hfil \ \ \frb #\ 
   &\vrule \hfil \ \ \frb #\ \ \vrule \hskip 2 pt
   &\vrule \ \ \hfil \frb #\ 
   &\vrule \hfil \ \ \fbb #\frb\ 
   &\vrule \hfil \ \ \frb #\ \hfil
   &\vrule \hfil \ \ \frb #\ 
   &\vrule \hfil \ \ \frb #\ \vrule \cr%
\noalign{\hrule}
 & &25.19.37.103.173&8701&8874& & &9.5.11.13.19.29.89&2219&1174 \cr
4717&313168925&4.9.7.11.17.29.103.113&8881&950&4735&315565965&4.7.89.317.587&203&114 \cr
 & &16.3.25.17.19.83.107&1819&1992& & &16.3.49.19.29.587&587&392 \cr
\noalign{\hrule}
 & &9.149.409.571&385&956& & &27.7.11.31.59.83&1433&520 \cr
4718&313175799&8.5.7.11.239.409&903&1142&4736&315606753&16.3.5.13.59.1433&805&628 \cr
 & &32.3.49.11.43.571&539&688& & &128.25.7.13.23.157&7475&10048 \cr
\noalign{\hrule}
 & &11.23.53.61.383&35073&39286& & &81.7.121.43.107&8215&4732 \cr
4719&313274467&4.81.13.433.1511&1405&106&4737&315660807&8.5.49.169.31.53&4601&4356 \cr
 & &16.27.5.13.53.281&1405&2808& & &64.9.121.31.43.107&31&32 \cr
\noalign{\hrule}
 & &9.25.11.353.359&1417&1058& & &19.23.41.67.263&143&120 \cr
4720&313649325&4.13.529.109.353&1077&1430&4738&315715457&16.3.5.11.13.19.41.67&825&46 \cr
 & &16.3.5.11.169.23.359&169&184& & &64.9.125.121.23&1125&3872 \cr
\noalign{\hrule}
 & &9.11.13.43.53.107&485&98& & &17.19.59.73.227&5093&18486 \cr
4721&313838811&4.5.49.13.97.107&313&222&4739&315793547&4.9.11.13.79.463&415&454 \cr
 & &16.3.7.37.97.313&11581&5432& & &16.3.5.83.227.463&2315&1992 \cr
\noalign{\hrule}
 & &9.13.289.37.251&46487&49750& & &27.13.19.23.29.71&671&2020 \cr
4722&314021331&4.125.7.29.199.229&187&42&4740&315823833&8.3.5.11.29.61.101&71&74 \cr
 & &16.3.25.49.11.17.199&4975&4312& & &32.11.37.61.71.101&6161&6512 \cr
\noalign{\hrule}
 & &3.17.103.163.367&235&132& & &3.37.43.239.277&98509&100100 \cr
4723&314239713&8.9.5.11.17.47.163&1133&334&4741&315986919&8.25.7.11.13.23.4283&43&342 \cr
 & &32.5.121.103.167&835&1936& & &32.9.5.19.43.4283&4283&4560 \cr
\noalign{\hrule}
 & &5.7.11.53.73.211&7719&7684& & &9.25.7.11.17.29.37&857&738 \cr
4724&314298215&8.3.11.17.31.53.83.113&171&1072&4742&316025325&4.81.5.37.41.857&24157&7552 \cr
 & &256.27.19.31.67.83&105659&107136& & &1024.49.17.29.59&413&512 \cr
\noalign{\hrule}
 & &125.13.17.59.193&37&258& & &25.49.17.43.353&363&62 \cr
4725&314565875&4.3.25.37.43.193&2431&2394&4743&316102675&4.3.7.121.31.353&285&68 \cr
 & &16.27.7.11.13.17.19.43&3591&3784& & &32.9.5.121.17.19&171&1936 \cr
\noalign{\hrule}
 & &3.13.41.47.53.79&2865&374& & &25.7.13.43.53.61&9&44 \cr
4726&314665611&4.9.5.11.13.17.191&6023&3922&4744&316268225&8.9.5.11.13.43.61&737&178 \cr
 & &16.19.37.53.317&703&2536& & &32.3.121.67.89&17889&1936 \cr
\noalign{\hrule}
 & &3.125.11.23.31.107&89&624& & &3.25.7.121.13.383&443&404 \cr
4727&314700375&32.9.25.11.13.89&37&62&4745&316290975&8.25.101.383.443&5729&5346 \cr
 & &128.13.31.37.89&3293&832& & &32.243.11.17.101.337&27297&27472 \cr
\noalign{\hrule}
 & &5.11.13.17.19.29.47&291&226& & &9.13.773.3499&10721&20770 \cr
4728&314778035&4.3.17.19.29.97.113&797&10164&4746&316453059&4.5.31.67.71.151&2173&2508 \cr
 & &32.9.7.121.797&5579&1584& & &32.3.11.19.41.53.71&23903&21584 \cr
\noalign{\hrule}
 & &3.5.7.13.31.43.173&1843&978& & &3.49.11.23.67.127&6565&18962 \cr
4729&314781285&4.9.19.43.97.163&275&112&4747&316458219&4.5.13.19.101.499&373&126 \cr
 & &128.25.7.11.19.97&1843&3520& & &16.9.5.7.101.373&1119&4040 \cr
\noalign{\hrule}
 & &3.5.11.13.17.89.97&333&112& & &9.625.7.11.17.43&197&428 \cr
4730&314802345&32.27.7.11.37.97&2963&4030&4748&316614375&8.3.17.43.107.197&1043&1150 \cr
 & &128.5.13.31.2963&2963&1984& & &32.25.7.23.149.197&3427&3152 \cr
\noalign{\hrule}
 & &3.5.37.41.101.137&4605&464& & &27.7.31.191.283&3725&3916 \cr
4731&314860935&32.9.25.29.307&143&118&4749&316696527&8.25.7.11.31.89.149&2223&3202 \cr
 & &128.11.13.59.307&8437&19648& & &32.9.13.19.149.1601&30419&30992 \cr
\noalign{\hrule}
 & &27.5.19.43.2857&2009&848& & &9.5.11.67.73.131&169&224 \cr
4732&315112815&32.5.49.19.41.53&1419&754&4750&317156895&64.3.7.169.67.73&575&374 \cr
 & &128.3.7.11.13.29.43&1001&1856& & &256.25.7.11.13.17.23&5083&4480 \cr
\noalign{\hrule}
 & &9.25.11.347.367&1411&1064& & &3.625.11.89.173&1889&14 \cr
4733&315188775&16.7.17.19.83.367&8425&1452&4751&317563125&4.7.89.1889&989&900 \cr
 & &128.3.25.121.337&337&704& & &32.9.25.7.23.43&161&2064 \cr
\noalign{\hrule}
 & &9.5.17.29.59.241&2491&1606& & &9.25.7.17.29.409&697&4378 \cr
4734&315448515&4.3.11.29.47.53.73&67&20&4752&317578275&4.11.289.41.199&81&370 \cr
 & &32.5.11.53.67.73&4891&9328& & &16.81.5.37.199&199&2664 \cr
\noalign{\hrule}
}%
}
$$
\eject
\vglue -23 pt
\noindent\hskip 1 in\hbox to 6.5 in{\ 4753 -- 4788 \hfill\fbd 317748717 -- 324257439\frb}
\vskip -9 pt
$$
\vbox{
\nointerlineskip
\halign{\strut
    \vrule \ \ \hfil \frb #\ 
   &\vrule \hfil \ \ \fbb #\frb\ 
   &\vrule \hfil \ \ \frb #\ \hfil
   &\vrule \hfil \ \ \frb #\ 
   &\vrule \hfil \ \ \frb #\ \ \vrule \hskip 2 pt
   &\vrule \ \ \hfil \frb #\ 
   &\vrule \hfil \ \ \fbb #\frb\ 
   &\vrule \hfil \ \ \frb #\ \hfil
   &\vrule \hfil \ \ \frb #\ 
   &\vrule \hfil \ \ \frb #\ \vrule \cr%
\noalign{\hrule}
 & &27.11.13.17.47.103&301&310& & &729.7.11.29.197&325&1054 \cr
4753&317748717&4.3.5.7.11.17.31.43.103&1387&364&4771&320687829&4.25.11.13.17.29.31&243&98 \cr
 & &32.5.49.13.19.43.73&15695&14896& & &16.243.5.49.13.17&455&136 \cr
\noalign{\hrule}
 & &9.125.11.17.1511&1157&2668& & &7.13.37.47.2029&18405&4202 \cr
4754&317876625&8.5.11.13.23.29.89&2781&2114&4772&321087221&4.9.5.11.191.409&709&518 \cr
 & &32.27.7.13.103.151&9373&7248& & &16.3.5.7.11.37.709&709&1320 \cr
\noalign{\hrule}
 & &13.19.29.37.1201&7825&7788& & &9.5.11.13.19.37.71&367&2402 \cr
4755&318302231&8.3.25.11.19.29.59.313&1221&1534&4773&321190155&4.3.19.367.1201&629&572 \cr
 & &32.9.5.121.13.37.3481&17405&17424& & &32.11.13.17.37.367&367&272 \cr
\noalign{\hrule}
 & &27.11.13.17.23.211&157&140& & &3.5.121.23.43.179&203&160 \cr
4756&318536361&8.5.7.13.23.157.211&255&44&4774&321311265&64.25.7.23.29.179&377&198 \cr
 & &64.3.25.7.11.17.157&1099&800& & &256.9.7.11.13.841&5887&4992 \cr
\noalign{\hrule}
 & &27.13.19.137.349&1675&2024& & &27.7.11.19.79.103&547&650 \cr
4757&318864897&16.25.11.13.19.23.67&387&88&4775&321419637&4.3.25.11.13.79.547&1577&1030 \cr
 & &256.9.121.43.67&8107&5504& & &16.125.13.19.83.103&1079&1000 \cr
\noalign{\hrule}
 & &81.7.281.2003&7051&6970& & &61.18769.281&17955&814 \cr
4758&319131981&4.5.11.17.41.281.641&25083&1198&4776&321719429&4.27.5.7.11.19.37&137&122 \cr
 & &16.27.11.599.929&6589&7432& & &16.9.11.19.61.137&209&72 \cr
\noalign{\hrule}
 & &81.23.37.41.113&3751&430& & &9.25.31.61.757&1357&1388 \cr
4759&319357323&4.5.121.23.31.43&793&540&4777&322084575&8.5.23.59.347.757&2379&42284 \cr
 & &32.27.25.11.13.61&3575&976& & &64.3.11.13.961.61&341&416 \cr
\noalign{\hrule}
 & &47.109.127.491&4617&18460& & &27.121.47.2099&1099&1000 \cr
4760&319454911&8.243.5.13.19.71&715&634&4778&322299351&16.3.125.7.11.47.157&4427&3328 \cr
 & &32.3.25.11.169.317&23775&29744& & &8192.25.13.19.233&75725&77824 \cr
\noalign{\hrule}
 & &5.49.11.19.6241&387&482& & &3.25.11.13.17.29.61&4489&4234 \cr
4761&319570405&4.9.49.43.79.241&625&10988&4779&322532925&4.5.841.4489.73&4347&142 \cr
 & &32.3.625.41.67&8241&2000& & &16.27.7.23.71.73&5183&11592 \cr
\noalign{\hrule}
 & &3.5.49.251.1733&603&652& & &81.7.11.13.23.173&5&86 \cr
4762&319712505&8.27.67.163.1733&1771&38&4780&322621299&4.5.11.23.43.173&1521&3424 \cr
 & &32.7.11.19.23.163&1793&6992& & &256.9.169.107&107&1664 \cr
\noalign{\hrule}
 & &5.7.13.19.163.227&9723&11842& & &7.11.17.37.59.113&263&150 \cr
4763&319873645&4.3.49.31.191.463&2497&11856&4781&322902811&4.3.25.11.17.37.263&177&452 \cr
 & &128.9.11.13.19.227&99&64& & &32.9.59.113.263&263&144 \cr
\noalign{\hrule}
 & &9.5.13.241.2269&19771&8426& & &5.11.19.47.6581&441&452 \cr
4764&319894965&4.11.17.383.1163&773&390&4782&323225815&8.9.5.49.113.6581&3971&23714 \cr
 & &16.3.5.11.13.17.773&773&1496& & &32.3.11.361.71.167&3173&3408 \cr
\noalign{\hrule}
 & &5.11.13.487.919&10773&20882& & &27.5.7.23.107.139&5687&5548 \cr
4765&320000395&4.81.7.19.53.197&481&110&4783&323264655&8.9.121.19.23.47.73&1&208 \cr
 & &16.27.5.11.13.19.37&703&216& & &256.11.13.47.73&3431&18304 \cr
\noalign{\hrule}
 & &27.7.17.37.2693&4843&2150& & &3.5.23.67.71.197&1309&324 \cr
4766&320146533&4.25.17.29.43.167&1057&1782&4784&323309505&8.243.7.11.17.67&1219&482 \cr
 & &16.81.7.11.43.151&1661&1032& & &32.17.23.53.241&901&3856 \cr
\noalign{\hrule}
 & &27.29.449.911&231&680& & &3.5.19.59.71.271&1417&2772 \cr
4767&320277537&16.81.5.7.11.17.29&449&44&4785&323537415&8.27.7.11.13.19.109&193&2264 \cr
 & &128.7.121.449&847&64& & &128.11.193.283&2123&18112 \cr
\noalign{\hrule}
 & &11.17.29.31.1907&2847&18130& & &9.5.11.317.2063&4199&6262 \cr
4768&320591491&4.3.5.49.13.37.73&1037&1518&4786&323715645&4.3.5.13.17.19.31.101&2063&3982 \cr
 & &16.9.7.11.17.23.61&1403&504& & &16.11.17.181.2063&181&136 \cr
\noalign{\hrule}
 & &7.11.169.71.347&85&84& & &9.7.11.13.127.283&1615&1498 \cr
4769&320601281&8.3.5.49.11.17.71.347&169&3648&4787&323792469&4.5.49.17.19.107.127&1415&3828 \cr
 & &1024.9.5.169.17.19&2907&2560& & &32.3.25.11.17.29.283&425&464 \cr
\noalign{\hrule}
 & &3.5.11.71.101.271&691&420& & &3.11.17.19.29.1049&249&800 \cr
4770&320651265&8.9.25.7.71.691&3997&2222&4788&324257439&64.9.25.11.17.83&1049&634 \cr
 & &32.49.11.101.571&571&784& & &256.5.317.1049&317&640 \cr
\noalign{\hrule}
}%
}
$$
\eject
\vglue -23 pt
\noindent\hskip 1 in\hbox to 6.5 in{\ 4789 -- 4824 \hfill\fbd 324662975 -- 330780169\frb}
\vskip -9 pt
$$
\vbox{
\nointerlineskip
\halign{\strut
    \vrule \ \ \hfil \frb #\ 
   &\vrule \hfil \ \ \fbb #\frb\ 
   &\vrule \hfil \ \ \frb #\ \hfil
   &\vrule \hfil \ \ \frb #\ 
   &\vrule \hfil \ \ \frb #\ \ \vrule \hskip 2 pt
   &\vrule \ \ \hfil \frb #\ 
   &\vrule \hfil \ \ \fbb #\frb\ 
   &\vrule \hfil \ \ \frb #\ \hfil
   &\vrule \hfil \ \ \frb #\ 
   &\vrule \hfil \ \ \frb #\ \vrule \cr%
\noalign{\hrule}
 & &25.49.13.19.29.37&33&292& & &3.125.7.121.13.79&1177&1098 \cr
4789&324662975&8.3.7.11.19.29.73&795&592&4807&326200875&4.27.5.1331.61.107&1651&34286 \cr
 & &256.9.5.11.37.53&477&1408& & &16.7.13.31.79.127&127&248 \cr
\noalign{\hrule}
 & &13.17.71.127.163&99&28& & &11.19.41.113.337&125&84 \cr
4790&324819391&8.9.7.11.13.17.163&355&134&4808&326316089&8.3.125.7.113.337&451&114 \cr
 & &32.3.5.7.11.67.71&2211&560& & &32.9.25.7.11.19.41&225&112 \cr
\noalign{\hrule}
 & &13.23.43.131.193&26037&31670& & &9.7.13.53.73.103&33&20 \cr
4791&325063531&4.9.5.11.263.3167&1189&1978&4809&326377233&8.27.5.7.11.73.103&8957&1438 \cr
 & &16.3.5.11.23.29.41.43&1189&1320& & &32.169.53.719&719&208 \cr
\noalign{\hrule}
 & &11.13.17.173.773&11151&21200& & &81.7.17.19.1783&4495&6278 \cr
4792&325095199&32.27.25.7.53.59&17&8&4810&326540403&4.5.17.29.31.43.73&259&2376 \cr
 & &512.3.7.17.53.59&3127&5376& & &64.27.7.11.37.43&407&1376 \cr
\noalign{\hrule}
 & &27.5.49.11.41.109&299&152& & &3.13.17.23.29.739&6039&5300 \cr
4793&325186785&16.9.5.13.19.23.109&1681&826&4811&326801319&8.27.25.11.13.53.61&739&2494 \cr
 & &64.7.13.1681.59&767&1312& & &32.5.11.29.43.739&215&176 \cr
\noalign{\hrule}
 & &3.25.49.11.13.619&1541&316& & &7.169.23.41.293&27621&20882 \cr
4794&325299975&8.11.13.23.67.79&145&882&4812&326861717&4.81.11.31.53.197&145&52 \cr
 & &32.9.5.49.23.29&667&48& & &32.27.5.11.13.29.53&7155&5104 \cr
\noalign{\hrule}
 & &27.5.49.101.487&4807&8342& & &5.11.23.29.37.241&1751&900 \cr
4795&325372005&4.7.11.19.23.43.97&4383&2540&4813&327120145&8.9.125.17.29.103&4439&814 \cr
 & &32.9.5.11.127.487&127&176& & &32.3.11.23.37.193&193&48 \cr
\noalign{\hrule}
 & &9.13.43.71.911&15433&23740& & &3.5.11.169.61.193&303&368 \cr
4796&325410111&8.5.11.23.61.1187&929&258&4814&328290105&32.9.13.23.101.193&1525&212 \cr
 & &32.3.5.23.43.929&929&1840& & &256.25.23.53.61&1219&640 \cr
\noalign{\hrule}
 & &3.7.11.13.107.1013&667&346& & &27.5.7.59.71.83&2921&1976 \cr
4797&325498173&4.7.11.13.23.29.173&95&1998&4815&328564215&16.13.19.23.71.127&99&28 \cr
 & &16.27.5.19.29.37&703&10440& & &128.9.7.11.13.19.23&2717&1472 \cr
\noalign{\hrule}
 & &9.7.13.17.149.157&319&18050& & &27.7.11.13.43.283&421&138 \cr
4798&325700739&4.25.11.361.29&353&372&4816&328891563&4.81.7.11.23.421&8773&11720 \cr
 & &32.3.11.19.31.353&10943&3344& & &64.5.31.283.293&1465&992 \cr
\noalign{\hrule}
 & &13.23.67.71.229&3069&2198& & &3.5.11.17.41.47.61&10541&20996 \cr
4799&325716547&4.9.7.11.31.71.157&575&916&4817&329719335&8.29.83.127.181&27&154 \cr
 & &32.3.25.23.157.229&471&400& & &32.27.7.11.29.83&747&3248 \cr
\noalign{\hrule}
 & &9.5.29.37.43.157&693&6058& & &5.7.11.19.23.37.53&4891&4914 \cr
4800&325972035&4.81.7.11.13.233&157&76&4818&329928445&4.27.49.11.13.19.67.73&115&1502 \cr
 & &32.7.11.13.19.157&209&1456& & &16.9.5.13.23.67.751&6759&6968 \cr
\noalign{\hrule}
 & &3.5.7.11.19.83.179&181&104& & &5.11.17.29.43.283&4469&3738 \cr
4801&326036865&16.13.83.179.181&629&450&4819&329962435&4.3.5.7.11.41.89.109&377&822 \cr
 & &64.9.25.17.37.181&6697&8160& & &16.9.7.13.29.41.137&5617&6552 \cr
\noalign{\hrule}
 & &3.25.7.11.131.431&8749&5474& & &361.311.2939&1289&1650 \cr
4802&326062275&4.49.13.17.23.673&655&18&4820&329964469&4.3.25.11.311.1289&133&1422 \cr
 & &16.9.5.17.23.131&23&408& & &16.27.5.7.11.19.79&1485&4424 \cr
\noalign{\hrule}
 & &27.5.7.11.13.19.127&899&6086& & &7.13.17.29.37.199&1203&6974 \cr
4803&326080755&4.9.17.29.31.179&133&394&4821&330326269&4.3.7.11.317.401&3315&1096 \cr
 & &16.7.19.179.197&179&1576& & &64.9.5.13.17.137&1233&160 \cr
\noalign{\hrule}
 & &9.11.169.101.193&1525&212& & &5.11.13.41.59.191&159&2260 \cr
4804&326137383&8.25.11.13.53.61&303&368&4822&330350735&8.3.25.13.53.113&3157&2832 \cr
 & &256.3.5.23.53.101&1219&640& & &256.9.7.11.41.59&63&128 \cr
\noalign{\hrule}
 & &27.841.53.271&451&1988& & &3.5.11.13.29.47.113&589&1106 \cr
4805&326140641&8.3.7.11.29.41.71&271&590&4823&330370755&4.7.13.19.29.31.79&23077&23454 \cr
 & &32.5.59.71.271&295&1136& & &16.9.7.47.491.1303&10311&10424 \cr
\noalign{\hrule}
 & &5.11.19.29.47.229&17&534& & &17.41.677.701&6105&5404 \cr
4806&326172715&4.3.5.17.89.229&617&528&4824&330780169&8.3.5.7.11.37.41.193&83&1434 \cr
 & &128.9.11.17.617&5553&1088& & &32.9.5.11.83.239&19837&7920 \cr
\noalign{\hrule}
}%
}
$$
\eject
\vglue -23 pt
\noindent\hskip 1 in\hbox to 6.5 in{\ 4825 -- 4860 \hfill\fbd 331357625 -- 337368213\frb}
\vskip -9 pt
$$
\vbox{
\nointerlineskip
\halign{\strut
    \vrule \ \ \hfil \frb #\ 
   &\vrule \hfil \ \ \fbb #\frb\ 
   &\vrule \hfil \ \ \frb #\ \hfil
   &\vrule \hfil \ \ \frb #\ 
   &\vrule \hfil \ \ \frb #\ \ \vrule \hskip 2 pt
   &\vrule \ \ \hfil \frb #\ 
   &\vrule \hfil \ \ \fbb #\frb\ 
   &\vrule \hfil \ \ \frb #\ \hfil
   &\vrule \hfil \ \ \frb #\ 
   &\vrule \hfil \ \ \frb #\ \vrule \cr%
\noalign{\hrule}
 & &125.17.19.29.283&5291&2916& & &27.25.11.19.23.103&109&406 \cr
4825&331357625&8.729.11.13.17.37&4237&3940&4843&334206675&4.5.7.19.23.29.109&711&1474 \cr
 & &64.27.5.19.197.223&6021&6304& & &16.9.11.29.67.79&1943&632 \cr
\noalign{\hrule}
 & &27.7.17.19.61.89&2849&1690& & &5.49.19.29.37.67&1105&2178 \cr
4826&331424163&4.9.5.49.11.169.37&323&314&4844&334652605&4.9.25.121.13.17.19&43&518 \cr
 & &16.5.11.13.17.19.37.157&5809&5720& & &16.3.7.11.13.37.43&129&1144 \cr
\noalign{\hrule}
 & &3.19.59.241.409&203&206& & &7.13.19.97.1997&1863&134 \cr
4827&331487547&4.7.19.29.59.103.241&13959&260&4845&334922861&4.81.23.67.97&539&2080 \cr
 & &32.27.5.11.13.29.47&14355&9776& & &256.3.5.49.11.13&77&1920 \cr
\noalign{\hrule}
 & &11.169.19.41.229&3105&1246& & &27.5.7.19.47.397&20987&10268 \cr
4828&331628869&4.27.5.7.23.41.89&605&338&4846&335022345&8.17.31.151.677&2679&2002 \cr
 & &16.9.25.7.121.169&275&504& & &32.3.7.11.13.17.19.47&221&176 \cr
\noalign{\hrule}
 & &3.25.49.11.43.191&531&806& & &27.11.361.53.59&9775&11524 \cr
4829&332010525&4.27.7.13.31.43.59&13561&19420&4847&335267559&8.9.25.17.23.43.67&361&26 \cr
 & &32.5.71.191.971&971&1136& & &32.5.13.17.361.23&391&1040 \cr
\noalign{\hrule}
 & &3.5.41.53.61.167&4991&5196& & &27.289.97.443&713&616 \cr
4830&332045265&8.9.7.23.31.53.433&6215&2318&4848&335302713&16.9.7.11.289.23.31&95&194 \cr
 & &32.5.11.19.31.61.113&3503&3344& & &64.5.7.19.23.31.97&4123&3680 \cr
\noalign{\hrule}
 & &9.121.37.73.113&7&106& & &3.11.13.199.3929&1075&1114 \cr
4831&332376957&4.7.11.37.53.73&167&240&4849&335422659&4.25.43.557.3929&3357&572 \cr
 & &128.3.5.7.53.167&1855&10688& & &32.9.5.11.13.43.373&1865&2064 \cr
\noalign{\hrule}
 & &5.11.19.23.61.227&1649&4146& & &25.11.41.83.359&469&444 \cr
4832&332812645&4.3.17.23.97.691&6989&4758&4850&335961175&8.3.7.37.41.67.359&2015&498 \cr
 & &16.9.13.29.61.241&3133&2088& & &32.9.5.13.31.67.83&2077&1872 \cr
\noalign{\hrule}
 & &3.5.7.121.17.23.67&20821&10034& & &81.5.11.19.29.137&5389&5708 \cr
4833&332832885&4.29.47.173.443&63&110&4851&336294585&8.5.17.19.317.1427&1521&94 \cr
 & &16.9.5.7.11.29.443&443&696& & &32.9.169.47.317&7943&5072 \cr
\noalign{\hrule}
 & &7.11.37.179.653&13065&11096& & &11.17.47.101.379&3777&2666 \cr
4834&333011063&16.3.5.7.13.19.67.73&111&358&4852&336434131&4.3.31.43.47.1259&707&750 \cr
 & &64.9.5.37.73.179&365&288& & &16.9.125.7.101.1259&8813&9000 \cr
\noalign{\hrule}
 & &3.5.29.43.47.379&5425&5566& & &5.7.17.43.59.223&99&314 \cr
4835&333192165&4.125.7.121.23.31.43&7767&22892&4853&336621845&4.9.11.17.157.223&137&86 \cr
 & &32.9.7.59.97.863&40061&41424& & &16.3.11.43.137.157&1727&3288 \cr
\noalign{\hrule}
 & &25.49.31.67.131&3287&3132& & &9.25.13.31.47.79&193&518 \cr
4836&333306575&8.27.5.19.29.67.173&11&56&4854&336676275&4.7.31.37.47.193&825&632 \cr
 & &128.3.7.11.19.29.173&15051&13376& & &64.3.25.7.11.37.79&259&352 \cr
\noalign{\hrule}
 & &27.7.13.19.37.193&27&220& & &169.17.19.31.199&729&4510 \cr
4837&333363303&8.729.5.7.11.37&235&494&4855&336747203&4.729.5.11.17.41&299&398 \cr
 & &32.25.11.13.19.47&1175&176& & &16.81.5.13.23.199&405&184 \cr
\noalign{\hrule}
 & &27.13.23.67.617&2783&2770& & &9.125.17.79.223&293&418 \cr
4838&333729747&4.3.5.121.529.67.277&8021&86&4856&336925125&4.11.17.19.223.293&273&50 \cr
 & &16.13.43.277.617&277&344& & &16.3.25.7.11.13.293&1001&2344 \cr
\noalign{\hrule}
 & &3.5.29.31.53.467&121&34& & &9.25.7.13.101.163&1411&4946 \cr
4839&333767235&4.121.17.53.467&217&684&4857&337079925&4.3.5.17.83.2473&1859&614 \cr
 & &32.9.7.121.19.31&399&1936& & &16.11.169.17.307&5219&1144 \cr
\noalign{\hrule}
 & &81.11.13.19.37.41&995&292& & &3.343.11.13.29.79&1145&276 \cr
4840&333856809&8.9.5.41.73.199&221&148&4858&337113777&8.9.5.7.13.23.229&373&212 \cr
 & &64.5.13.17.37.199&995&544& & &64.53.229.373&19769&7328 \cr
\noalign{\hrule}
 & &9.5.7.11.211.457&30719&35746& & &27.7.121.23.641&2353&430 \cr
4841&334119555&4.13.17.61.139.293&77&216&4859&337157667&4.9.5.7.13.43.181&1807&902 \cr
 & &64.27.7.11.13.17.61&1037&1248& & &16.11.169.41.139&5699&1352 \cr
\noalign{\hrule}
 & &25.11.53.101.227&859&276& & &27.7.13.17.41.197&341&250 \cr
4842&334161025&8.3.5.23.101.859&3309&986&4860&337368213&4.9.125.11.17.31.41&179&26 \cr
 & &32.9.17.29.1103&18751&4176& & &16.25.11.13.31.179&5549&2200 \cr
\noalign{\hrule}
}%
}
$$
\eject
\vglue -23 pt
\noindent\hskip 1 in\hbox to 6.5 in{\ 4861 -- 4896 \hfill\fbd 337451625 -- 344345265\frb}
\vskip -9 pt
$$
\vbox{
\nointerlineskip
\halign{\strut
    \vrule \ \ \hfil \frb #\ 
   &\vrule \hfil \ \ \fbb #\frb\ 
   &\vrule \hfil \ \ \frb #\ \hfil
   &\vrule \hfil \ \ \frb #\ 
   &\vrule \hfil \ \ \frb #\ \ \vrule \hskip 2 pt
   &\vrule \ \ \hfil \frb #\ 
   &\vrule \hfil \ \ \fbb #\frb\ 
   &\vrule \hfil \ \ \frb #\ \hfil
   &\vrule \hfil \ \ \frb #\ 
   &\vrule \hfil \ \ \frb #\ \vrule \cr%
\noalign{\hrule}
 & &9.125.7.73.587&231&356& & &9.11.13.19.73.191&15689&24230 \cr
4861&337451625&8.27.49.11.73.89&587&2990&4879&340948179&4.5.29.541.2423&941&1482 \cr
 & &32.5.11.13.23.587&143&368& & &16.3.5.13.19.29.941&941&1160 \cr
\noalign{\hrule}
 & &9.11.107.127.251&1075&322& & &9.7.11.361.29.47&3425&7044 \cr
4862&337673061&4.3.25.7.23.43.107&55&52&4880&340985799&8.27.25.137.587&9541&8954 \cr
 & &32.125.7.11.13.23.43&12857&14000& & &32.5.7.121.29.37.47&185&176 \cr
\noalign{\hrule}
 & &7.289.199.839&6435&7828& & &3.25.7.13.23.41.53&283&242 \cr
4863&337762103&8.9.5.11.13.17.19.103&1945&2254&4881&341106675&4.121.13.23.53.283&441&142 \cr
 & &32.3.25.49.11.23.389&29175&28336& & &16.9.49.11.71.283&5467&6792 \cr
\noalign{\hrule}
 & &3.125.7.23.29.193&1739&2314& & &3.5.19.61.67.293&253&52 \cr
4864&337918875&4.5.13.29.37.47.89&4851&17756&4882&341284935&8.11.13.19.23.293&1139&4428 \cr
 & &32.9.49.11.23.193&77&48& & &64.27.17.41.67&153&1312 \cr
\noalign{\hrule}
 & &23.31.479.991&507&484& & &5.121.47.61.197&159&38 \cr
4865&338453257&8.3.121.169.31.479&69&410&4883&341703395&4.3.5.19.47.53.61&697&462 \cr
 & &32.9.5.11.169.23.41&6929&7920& & &16.9.7.11.17.41.53&6307&2952 \cr
\noalign{\hrule}
 & &81.29.31.4649&3499&1150& & &5.11.23.37.67.109&1025&174 \cr
4866&338535531&4.25.23.31.3499&1827&1672&4884&341816915&4.3.125.29.41.67&1591&3534 \cr
 & &64.9.5.7.11.19.23.29&2185&2464& & &16.9.19.31.37.43&589&3096 \cr
\noalign{\hrule}
 & &11.23.47.71.401&591&190& & &5.7.97.263.383&1089&826 \cr
4867&338548661&4.3.5.19.23.47.197&4411&4848&4885&341974955&4.9.49.121.59.97&10205&1532 \cr
 & &128.9.5.11.101.401&505&576& & &32.3.5.13.157.383&471&208 \cr
\noalign{\hrule}
 & &9.125.359.839&143&982& & &7.23.47.53.853&819&1672 \cr
4868&338851125&4.11.13.359.491&367&5034&4886&342096503&16.9.49.11.13.19.23&235&304 \cr
 & &16.3.367.839&367&8& & &512.3.5.13.361.47&4693&3840 \cr
\noalign{\hrule}
 & &9.7.11.17.107.269&295&26& & &9.7.17.29.103.107&1315&1672 \cr
4869&339092523&4.3.5.7.11.13.17.59&107&124&4887&342301239&16.3.5.11.19.107.263&497&818 \cr
 & &32.5.13.31.59.107&2015&944& & &64.7.11.19.71.409&14839&13088 \cr
\noalign{\hrule}
 & &7.53.71.79.163&55&108& & &3.5.121.41.43.107&61&104 \cr
4870&339193057&8.27.5.7.11.71.79&1007&652&4888&342383415&16.11.13.41.61.107&279&172 \cr
 & &64.9.11.19.53.163&209&288& & &128.9.13.31.43.61&1891&2496 \cr
\noalign{\hrule}
 & &9.25.11.13.59.179&21161&28064& & &27.13.19.89.577&913&1490 \cr
4871&339800175&64.7.877.3023&1073&1950&4889&342473157&4.5.11.13.19.83.149&1731&986 \cr
 & &256.3.25.7.13.29.37&1073&896& & &16.3.17.29.83.577&493&664 \cr
\noalign{\hrule}
 & &3.5.23.929.1061&11737&12666& & &27.13.19.53.971&77&18526 \cr
4872&340055805&4.9.5.121.97.2111&13793&9428&4890&343206747&4.7.11.59.157&555&544 \cr
 & &32.11.13.1061.2357&2357&2288& & &256.3.5.17.37.59&10915&2176 \cr
\noalign{\hrule}
 & &81.5.23.59.619&133&74& & &3.25.7.169.53.73&9757&9588 \cr
4873&340193115&4.9.5.7.19.37.619&737&118&4891&343277025&8.9.5.7.11.17.47.887&50297&8608 \cr
 & &16.7.11.37.59.67&2479&616& & &512.13.53.73.269&269&256 \cr
\noalign{\hrule}
 & &9.11.73.197.239&415&176& & &11.59.389.1361&18961&3990 \cr
4874&340268841&32.3.5.121.73.83&4541&4292&4892&343599421&4.3.5.7.19.67.283&2655&2722 \cr
 & &256.5.19.29.37.239&3515&3712& & &16.27.25.7.59.1361&189&200 \cr
\noalign{\hrule}
 & &81.25.49.3433&1937&1496& & &9.11.169.19.23.47&4949&4780 \cr
4875&340639425&16.9.25.11.13.17.149&2729&196&4893&343638009&8.5.49.11.19.101.239&39&94 \cr
 & &128.49.11.2729&2729&704& & &32.3.7.13.47.101.239&1673&1616 \cr
\noalign{\hrule}
 & &9.25.7.13.127.131&157&1808& & &3.5.19.23.131.401&7733&7332 \cr
4876&340642575&32.3.5.7.113.157&1397&298&4894&344340705&8.9.11.13.361.37.47&131&230 \cr
 & &128.11.127.149&1639&64& & &32.5.13.23.37.47.131&611&592 \cr
\noalign{\hrule}
 & &27.5.7.43.83.101&1271&970& & &9.11.29.31.53.73&24225&23422 \cr
4877&340643205&4.25.31.41.97.101&3283&858&4895&344344869&4.27.25.49.17.19.239&2743&1798 \cr
 & &16.3.49.11.13.31.67&5159&3224& & &16.5.7.13.17.29.31.211&3587&3640 \cr
\noalign{\hrule}
 & &3.31.97.107.353&145&176& & &9.5.11.361.41.47&783&1144 \cr
4878&340732191&32.5.11.29.97.353&535&3348&4896&344345265&16.243.5.121.13.29&179&1394 \cr
 & &256.27.25.31.107&225&128& & &64.17.29.41.179&3043&928 \cr
\noalign{\hrule}
}%
}
$$
\eject
\vglue -23 pt
\noindent\hskip 1 in\hbox to 6.5 in{\ 4897 -- 4932 \hfill\fbd 344369025 -- 351623025\frb}
\vskip -9 pt
$$
\vbox{
\nointerlineskip
\halign{\strut
    \vrule \ \ \hfil \frb #\ 
   &\vrule \hfil \ \ \fbb #\frb\ 
   &\vrule \hfil \ \ \frb #\ \hfil
   &\vrule \hfil \ \ \frb #\ 
   &\vrule \hfil \ \ \frb #\ \ \vrule \hskip 2 pt
   &\vrule \ \ \hfil \frb #\ 
   &\vrule \hfil \ \ \fbb #\frb\ 
   &\vrule \hfil \ \ \frb #\ \hfil
   &\vrule \hfil \ \ \frb #\ 
   &\vrule \hfil \ \ \frb #\ \vrule \cr%
\noalign{\hrule}
 & &9.25.7.121.13.139&1649&1376& & &3.25.7.11.157.383&137&962 \cr
4897&344369025&64.3.17.43.97.139&95&44&4915&347256525&4.13.37.137.383&45&5024 \cr
 & &512.5.11.19.43.97&4171&4864& & &256.9.5.157&1&384 \cr
\noalign{\hrule}
 & &3.11.13.19.29.31.47&2911&3810& & &3.5.7.11.31.89.109&20519&8584 \cr
4898&344404203&4.9.5.19.41.71.127&961&182&4916&347344305&16.289.29.37.71&279&350 \cr
 & &16.5.7.13.961.71&1085&568& & &64.9.25.7.17.29.31&435&544 \cr
\noalign{\hrule}
 & &3.5.49.11.13.29.113&81&458& & &27.5.49.11.17.281&1577&1118 \cr
4899&344429085&4.243.5.113.229&551&664&4917&347598405&4.13.19.43.83.281&119&162 \cr
 & &64.19.29.83.229&4351&2656& & &16.81.7.13.17.19.83&741&664 \cr
\noalign{\hrule}
 & &289.43.53.523&14225&13494& & &7.121.13.131.241&2025&1108 \cr
4900&344464013&4.3.25.13.17.173.569&129&44&4918&347628281&8.81.25.121.277&2227&5252 \cr
 & &32.9.5.11.13.43.569&7397&7920& & &64.3.13.17.101.131&303&544 \cr
\noalign{\hrule}
 & &41.151.179.311&7139&19890& & &5.11.13.29.31.541&111&430 \cr
4901&344646779&4.9.5.121.13.17.59&287&716&4919&347746685&4.3.25.13.31.37.43&261&664 \cr
 & &32.3.5.7.11.41.179&231&80& & &64.27.29.43.83&3569&864 \cr
\noalign{\hrule}
 & &9.11.19.1681.109&18565&74& & &61.29929.191&9139&20790 \cr
4902&344653749&4.5.37.47.79&4387&4308&4920&348702779&4.27.5.7.11.13.19.37&173&382 \cr
 & &32.3.41.107.359&359&1712& & &16.9.7.13.173.191&117&56 \cr
\noalign{\hrule}
 & &9.11.17.29.37.191&1585&2078& & &9.5.11.13.17.31.103&953&798 \cr
4903&344919069&4.5.191.317.1039&319&636&4921&349298235&4.27.7.11.13.19.953&25&272 \cr
 & &32.3.11.29.53.1039&1039&848& & &128.25.7.17.953&953&2240 \cr
\noalign{\hrule}
 & &9.5.11.67.101.103&5117&1784& & &13.19.31.71.643&445&198 \cr
4904&345015495&16.3.5.7.17.43.223&53&32&4922&349565021&4.9.5.11.31.71.89&313&468 \cr
 & &1024.43.53.223&11819&22016& & &32.81.13.89.313&7209&5008 \cr
\noalign{\hrule}
 & &3.5.49.29.97.167&935&1102& & &25.17.29.113.251&3427&3852 \cr
4905&345281685&4.25.7.11.17.19.841&1067&1908&4923&349573975&8.9.23.107.113.149&91&22 \cr
 & &32.9.121.19.53.97&2299&2544& & &32.3.7.11.13.107.149&24717&30992 \cr
\noalign{\hrule}
 & &3.5.169.17.71.113&517&588& & &3.49.17.19.37.199&715&118 \cr
4906&345751185&8.9.49.11.13.47.113&19&5518&4924&349602603&4.5.11.13.19.37.59&1139&1044 \cr
 & &32.11.19.31.89&18601&496& & &32.9.11.13.17.29.67&4147&3216 \cr
\noalign{\hrule}
 & &3.13.19.23.79.257&18623&890& & &9.5.7.11.13.19.409&973&1072 \cr
4907&346024029&4.5.11.89.1693&819&874&4925&350044695&32.49.13.19.67.139&955&318 \cr
 & &16.9.7.13.19.23.89&89&168& & &128.3.5.53.139.191&10123&8896 \cr
\noalign{\hrule}
 & &27.11.41.43.661&269&392& & &13.17.29.31.41.43&147&550 \cr
4908&346106871&16.9.49.11.43.269&445&2866&4926&350271077&4.3.25.49.11.29.43&689&732 \cr
 & &64.5.7.89.1433&7165&19936& & &32.9.25.11.13.53.61&11925&10736 \cr
\noalign{\hrule}
 & &7.17.101.151.191&1233&1334& & &3.5.19.59.67.311&8843&20702 \cr
4909&346639979&4.9.7.23.29.137.191&35855&30316&4927&350374155&4.11.37.239.941&351&590 \cr
 & &32.3.5.11.13.53.71.101&11289&11440& & &16.27.5.11.13.37.59&481&792 \cr
\noalign{\hrule}
 & &19.23.73.83.131&2295&718& & &5.11.19.31.79.137&387&482 \cr
4910&346859573&4.27.5.17.73.359&5405&5764&4928&350611085&4.9.31.43.137.241&6681&790 \cr
 & &32.3.25.11.23.47.131&825&752& & &16.27.5.17.79.131&459&1048 \cr
\noalign{\hrule}
 & &25.49.13.19.31.37&1377&142& & &9.5.13.19.73.433&2423&1474 \cr
4911&347053525&4.81.5.17.37.71&209&124&4929&351334035&4.5.11.19.67.2423&1971&4394 \cr
 & &32.9.11.19.31.71&781&144& & &16.27.11.2197.73&169&264 \cr
\noalign{\hrule}
 & &9.25.7.11.13.23.67&1231&494& & &5.13.17.487.653&201&286 \cr
4912&347071725&4.3.7.169.19.1231&1255&2438&4930&351402155&4.3.11.169.67.653&9253&2070 \cr
 & &16.5.19.23.53.251&1007&2008& & &16.27.5.19.23.487&437&216 \cr
\noalign{\hrule}
 & &9.5.11.13.17.19.167&21&188& & &25.11.19.137.491&169&306 \cr
4913&347110335&8.27.5.7.13.17.47&1919&3674&4931&351470075&4.9.11.169.17.491&31&460 \cr
 & &32.11.19.101.167&101&16& & &32.3.5.13.17.23.31&897&8432 \cr
\noalign{\hrule}
 & &9.7.19.29.73.137&535&2068& & &81.25.13.361.37&473&112 \cr
4914&347164713&8.3.5.11.29.47.107&2993&2036&4932&351623025&32.9.5.7.11.37.43&361&46 \cr
 & &64.5.41.73.509&2545&1312& & &128.361.23.43&989&64 \cr
\noalign{\hrule}
}%
}
$$
\eject
\vglue -23 pt
\noindent\hskip 1 in\hbox to 6.5 in{\ 4933 -- 4968 \hfill\fbd 351969885 -- 359001643\frb}
\vskip -9 pt
$$
\vbox{
\nointerlineskip
\halign{\strut
    \vrule \ \ \hfil \frb #\ 
   &\vrule \hfil \ \ \fbb #\frb\ 
   &\vrule \hfil \ \ \frb #\ \hfil
   &\vrule \hfil \ \ \frb #\ 
   &\vrule \hfil \ \ \frb #\ \ \vrule \hskip 2 pt
   &\vrule \ \ \hfil \frb #\ 
   &\vrule \hfil \ \ \fbb #\frb\ 
   &\vrule \hfil \ \ \frb #\ \hfil
   &\vrule \hfil \ \ \frb #\ 
   &\vrule \hfil \ \ \frb #\ \vrule \cr%
\noalign{\hrule}
 & &9.5.79.181.547&319&866& & &7.19.529.31.163&2169&1580 \cr
4933&351969885&4.3.11.29.181.433&10255&5492&4951&355513921&8.9.5.7.23.79.241&1&22 \cr
 & &32.5.7.293.1373&9611&4688& & &32.3.5.11.79.241&19039&2640 \cr
\noalign{\hrule}
 & &27.5.29.47.1913&1591&322& & &81.47.61.1531&1705&3236 \cr
4934&352001565&4.5.7.23.29.37.43&4307&4422&4952&355539537&8.5.11.31.47.809&287&522 \cr
 & &16.3.11.37.59.67.73&29711&31624& & &32.9.7.11.29.31.41&8323&5456 \cr
\noalign{\hrule}
 & &11.13.19.173.751&1075&828& & &81.5.17.113.457&26767&18998 \cr
4935&353000791&8.9.25.23.43.751&53&698&4953&355548285&4.7.13.23.29.59.71&363&304 \cr
 & &32.3.5.23.53.349&5235&19504& & &128.3.7.121.13.19.71&29887&31808 \cr
\noalign{\hrule}
 & &9.25.7.29.59.131&23&154& & &9.5.19.541.769&3731&3190 \cr
4936&353022075&4.3.25.49.11.23.29&1651&524&4954&355704795&4.25.7.11.13.19.29.41&2441&1416 \cr
 & &32.11.13.127.131&1651&176& & &64.3.11.13.59.2441&26851&24544 \cr
\noalign{\hrule}
 & &27.25.11.199.239&187&52& & &25.11.13.29.47.73&19359&4366 \cr
4937&353140425&8.5.121.13.17.199&531&1526&4955&355708925&4.81.37.59.239&733&1450 \cr
 & &32.9.7.13.59.109&6431&1456& & &16.27.25.29.733&733&216 \cr
\noalign{\hrule}
 & &9.5.7.121.13.23.31&47&8& & &9.5.7.11.13.53.149&1147&998 \cr
4938&353287935&16.3.7.11.23.31.47&3071&1300&4956&355720365&4.3.7.31.37.53.499&17&1130 \cr
 & &128.25.13.37.83&415&2368& & &16.5.17.113.499&8483&904 \cr
\noalign{\hrule}
 & &3.5.49.121.29.137&2885&1378& & &9.5.23.37.71.131&4873&26 \cr
4939&353338755&4.25.11.13.53.577&2511&3836&4957&356181795&4.3.5.11.13.443&851&1294 \cr
 & &32.81.7.13.31.137&351&496& & &16.23.37.647&647&8 \cr
\noalign{\hrule}
 & &3.5.13.17.19.31.181&141&2494& & &27.13.61.127.131&551&1100 \cr
4940&353408835&4.9.19.29.43.47&3395&4642&4958&356214807&8.3.25.11.19.29.131&41&434 \cr
 & &16.5.7.11.97.211&1477&8536& & &32.7.11.29.31.41&2387&19024 \cr
\noalign{\hrule}
 & &3.5.121.19.37.277&89&274& & &5.7.11.13.19.23.163&841&1278 \cr
4941&353436765&4.19.89.137.277&2563&2700&4959&356511155&4.9.5.7.11.841.71&19&184 \cr
 & &32.27.25.11.89.233&4005&3728& & &64.3.19.23.29.71&213&928 \cr
\noalign{\hrule}
 & &3.5.13.31.59.991&121&56& & &9.5.7.11.29.53.67&83&182 \cr
4942&353445105&16.7.121.31.991&1593&5344&4960&356822235&4.49.13.29.67.83&275&1146 \cr
 & &1024.27.59.167&1503&512& & &16.3.25.11.83.191&415&1528 \cr
\noalign{\hrule}
 & &7.11.19.47.53.97&81&598& & &5.11.17.43.83.107&195&278 \cr
4943&353500301&4.81.13.19.23.53&517&490&4961&357060605&4.3.25.13.17.107.139&747&1072 \cr
 & &16.3.5.49.11.13.23.47&455&552& & &128.27.67.83.139&3753&4288 \cr
\noalign{\hrule}
 & &9.11.529.43.157&2359&2402& & &27.7.11.283.607&2701&19090 \cr
4944&353556621&4.7.11.157.337.1201&19849&33060&4962&357132699&4.5.23.37.73.83&5733&2662 \cr
 & &32.3.5.7.19.23.29.863&16397&16240& & &16.9.49.1331.13&121&728 \cr
\noalign{\hrule}
 & &11.103.457.683&32661&37688& & &9.25.19.29.43.67&3509&3844 \cr
4945&353644423&16.9.7.19.191.673&1105&914&4963&357171975&8.5.121.841.961&1273&432 \cr
 & &64.3.5.7.13.17.19.457&3705&3808& & &256.27.11.19.31.67&341&384 \cr
\noalign{\hrule}
 & &5.121.13.19.23.103&2623&666& & &81.25.7.19.1327&4309&7634 \cr
4946&354011515&4.9.5.11.37.43.61&169&136&4964&357394275&4.9.11.31.139.347&455&796 \cr
 & &64.3.169.17.37.43&4773&7072& & &32.5.7.13.199.347&4511&3184 \cr
\noalign{\hrule}
 & &3.7.11.17.31.41.71&747&460& & &3.5.11.13.19.67.131&29&10 \cr
4947&354376407&8.27.5.11.23.31.83&71&226&4965&357706635&4.25.11.29.67.131&117&1558 \cr
 & &32.23.71.83.113&2599&1328& & &16.9.13.19.29.41&41&696 \cr
\noalign{\hrule}
 & &9.7.13.41.59.179&20491&27830& & &121.29.79.1291&211&1080 \cr
4948&354627819&4.5.121.23.31.661&1061&1722&4966&357879401&16.27.5.11.29.211&1109&790 \cr
 & &16.3.5.7.31.41.1061&1061&1240& & &64.3.25.79.1109&1109&2400 \cr
\noalign{\hrule}
 & &3.5.19.71.89.197&2501&2572& & &81.11.17.19.29.43&169&40 \cr
4949&354780255&8.5.41.61.197.643&171&814&4967&358877871&16.27.5.169.17.29&41&418 \cr
 & &32.9.11.19.37.41.61&6771&7216& & &64.5.11.13.19.41&65&1312 \cr
\noalign{\hrule}
 & &5.11.19.151.2251&603&1648& & &7.11.13.29.83.149&227&150 \cr
4950&355196545&32.9.67.103.151&25&176&4968&359001643&4.3.25.83.149.227&6751&5616 \cr
 & &1024.3.25.11.103&309&2560& & &128.81.5.13.43.157&12717&13760 \cr
\noalign{\hrule}
}%
}
$$
\eject
\vglue -23 pt
\noindent\hskip 1 in\hbox to 6.5 in{\ 4969 -- 5004 \hfill\fbd 359119475 -- 365099735\frb}
\vskip -9 pt
$$
\vbox{
\nointerlineskip
\halign{\strut
    \vrule \ \ \hfil \frb #\ 
   &\vrule \hfil \ \ \fbb #\frb\ 
   &\vrule \hfil \ \ \frb #\ \hfil
   &\vrule \hfil \ \ \frb #\ 
   &\vrule \hfil \ \ \frb #\ \ \vrule \hskip 2 pt
   &\vrule \ \ \hfil \frb #\ 
   &\vrule \hfil \ \ \fbb #\frb\ 
   &\vrule \hfil \ \ \frb #\ \hfil
   &\vrule \hfil \ \ \frb #\ 
   &\vrule \hfil \ \ \frb #\ \vrule \cr%
\noalign{\hrule}
 & &25.11.13.17.19.311&249&62& & &11.19.529.29.113&279&250 \cr
4969&359119475&4.3.25.13.19.31.83&511&264&4987&362308397&4.9.125.11.19.31.113&91&1334 \cr
 & &64.9.7.11.73.83&6059&2016& & &16.3.5.7.13.23.29.31&403&840 \cr
\noalign{\hrule}
 & &7.13.113.181.193&2669&19140& & &3.5.11.17.47.2749&1417&1332 \cr
4970&359216039&8.3.5.11.17.29.157&5&6&4988&362414415&8.27.11.13.37.47.109&4883&27490 \cr
 & &32.9.25.17.29.157&12325&22608& & &32.5.19.257.2749&257&304 \cr
\noalign{\hrule}
 & &13.41.43.61.257&375&418& & &3.17.37.401.479&39&440 \cr
4971&359301163&4.3.125.11.19.41.257&231&26&4989&362453073&16.9.5.11.13.17.37&401&364 \cr
 & &16.9.25.7.121.13.19&4275&6776& & &128.7.11.169.401&1183&704 \cr
\noalign{\hrule}
 & &27.5.121.361.61&455&94& & &3.13.289.53.607&21721&10450 \cr
4972&359713035&4.3.25.7.121.13.47&3743&7268&4990&362599341&4.25.7.11.19.29.107&731&864 \cr
 & &32.19.23.79.197&4531&1264& & &256.27.5.17.43.107&4815&5504 \cr
\noalign{\hrule}
 & &11.23.53.139.193&12469&22698& & &9.11.79.89.521&6381&650 \cr
4973&359723243&4.9.13.37.97.337&265&746&4991&362651949&4.81.25.13.709&2299&1246 \cr
 & &16.3.5.53.97.373&1455&2984& & &16.5.7.121.19.89&95&616 \cr
\noalign{\hrule}
 & &3.5.11.17.19.43.157&14271&15088& & &3.13.17.29.43.439&14345&23848 \cr
4974&359794545&32.9.5.23.41.67.71&343&14212&4992&362948079&16.5.11.19.151.271&783&878 \cr
 & &256.343.11.17.19&343&128& & &64.27.29.271.439&271&288 \cr
\noalign{\hrule}
 & &9.5.11.19.101.379&427&482& & &81.5.7.11.19.613&3281&2236 \cr
4975&360013995&4.7.19.61.241.379&2757&4444&4993&363211695&8.9.7.13.17.43.193&25&38 \cr
 & &32.3.11.61.101.919&919&976& & &32.25.17.19.43.193&3655&3088 \cr
\noalign{\hrule}
 & &81.5.49.11.13.127&251&134& & &11.13.17.23.67.97&819&722 \cr
4976&360405045&4.9.7.67.127.251&307&1450&4994&363378587&4.9.7.11.169.17.361&201&2060 \cr
 & &16.25.29.67.307&8903&2680& & &32.27.5.19.67.103&2781&1520 \cr
\noalign{\hrule}
 & &3.43.337.8293&535067&534730& & &3.25.19.47.61.89&697&462 \cr
4977&360521589&4.5.7.13.79.521.7639&225&7414&4995&363607275&4.9.5.7.11.17.41.89&407&38 \cr
 & &16.9.125.11.337.521&4125&4168& & &16.7.121.17.19.37&4477&952 \cr
\noalign{\hrule}
 & &3.23.1873.2791&971&902& & &5.7.11.841.1123&27939&33826 \cr
4978&360700467&4.11.41.971.2791&1881&910&4996&363610555&4.3.13.67.139.1301&253&1554 \cr
 & &16.9.5.7.121.13.19.41&33033&31160& & &16.9.7.11.23.37.67&1541&2664 \cr
\noalign{\hrule}
 & &9.5.7.11.79.1319&1829&3148& & &9.5.169.137.349&913&868 \cr
4979&361056465&8.5.11.31.59.787&459&1246&4997&363617865&8.7.11.13.31.83.349&2725&1812 \cr
 & &32.27.7.17.59.89&1513&2832& & &64.3.25.7.31.109.151&23653&24160 \cr
\noalign{\hrule}
 & &81.25.11.13.29.43&61&412& & &9.5.13.17.23.37.43&109&224 \cr
4980&361100025&8.3.25.29.61.103&5243&3718&4998&363917385&64.7.13.17.43.109&3267&1850 \cr
 & &32.49.11.169.107&637&1712& & &256.27.25.121.37&605&384 \cr
\noalign{\hrule}
 & &27.41.137.2381&763&6380& & &9.5.7.11.17.37.167&1633&698 \cr
4981&361100079&8.9.5.7.11.29.109&137&182&4999&363973995&4.23.71.167.349&91&258 \cr
 & &32.49.13.109.137&637&1744& & &16.3.7.13.23.43.71&1633&4472 \cr
\noalign{\hrule}
 & &27.5.13.227.907&13733&12826& & &7.11.19.23.31.349&1625&2214 \cr
4982&361335195&4.3.5.121.31.53.443&11&454&5000&364048531&4.27.125.7.13.23.41&1003&2728 \cr
 & &16.1331.53.227&1331&424& & &64.9.5.11.17.31.59&1003&1440 \cr
\noalign{\hrule}
 & &27.17.19.181.229&7105&3212& & &9.11.169.19.31.37&3437&226 \cr
4983&361476729&8.9.5.49.11.29.73&1027&1090&5001&364618683&4.7.31.113.491&12015&12506 \cr
 & &32.25.7.11.13.79.109&60277&57200& & &16.27.5.169.37.89&89&120 \cr
\noalign{\hrule}
 & &5.49.11.19.23.307&459&668& & &9.49.121.6841&6385&456 \cr
4984&361558505&8.27.5.17.167.307&281&26&5002&365042601&16.27.5.19.1277&1921&644 \cr
 & &32.9.13.167.281&19539&4496& & &128.7.17.23.113&391&7232 \cr
\noalign{\hrule}
 & &9.13.19.29.71.79&3175&2434& & &27.13.361.43.67&36001&36560 \cr
4985&361595403&4.3.25.29.127.1217&1463&2188&5003&365054391&32.9.5.7.37.139.457&1273&22 \cr
 & &32.7.11.19.127.547&9779&8752& & &128.11.19.67.457&457&704 \cr
\noalign{\hrule}
 & &9.11.13.17.71.233&6125&4028& & &5.49.11.13.17.613&4275&4888 \cr
4986&361944297&8.125.49.17.19.53&407&426&5004&365099735&16.9.125.169.19.47&10727&13102 \cr
 & &32.3.125.11.37.53.71&1961&2000& & &64.3.17.631.6551&19653&20192 \cr
\noalign{\hrule}
}%
}
$$
\eject
\vglue -23 pt
\noindent\hskip 1 in\hbox to 6.5 in{\ 5005 -- 5040 \hfill\fbd 365100183 -- 371561775\frb}
\vskip -9 pt
$$
\vbox{
\nointerlineskip
\halign{\strut
    \vrule \ \ \hfil \frb #\ 
   &\vrule \hfil \ \ \fbb #\frb\ 
   &\vrule \hfil \ \ \frb #\ \hfil
   &\vrule \hfil \ \ \frb #\ 
   &\vrule \hfil \ \ \frb #\ \ \vrule \hskip 2 pt
   &\vrule \ \ \hfil \frb #\ 
   &\vrule \hfil \ \ \fbb #\frb\ 
   &\vrule \hfil \ \ \frb #\ \hfil
   &\vrule \hfil \ \ \frb #\ 
   &\vrule \hfil \ \ \frb #\ \vrule \cr%
\noalign{\hrule}
 & &27.7.23.47.1787&27&20& & &25.7.11.53.3607&15201&24476 \cr
5005&365100183&8.729.5.23.1787&12851&3916&5023&368004175&8.27.29.211.563&9113&7214 \cr
 & &64.11.71.89.181&16109&24992& & &32.3.13.701.3607&701&624 \cr
\noalign{\hrule}
 & &11.23.37.103.379&5&6& & &27.121.17.19.349&659&310 \cr
5006&365425357&4.3.5.23.37.103.379&1089&12934&5024&368279109&4.9.5.121.31.659&215&874 \cr
 & &16.27.121.29.223&2453&6264& & &16.25.19.23.31.43&1333&4600 \cr
\noalign{\hrule}
 & &3.25.7.11.167.379&761&3414& & &27.11.13.17.31.181&3979&1988 \cr
5007&365517075&4.9.11.569.761&235&334&5025&368289207&8.7.23.31.71.173&1865&3498 \cr
 & &16.5.47.167.761&761&376& & &32.3.5.7.11.53.373&2611&4240 \cr
\noalign{\hrule}
 & &3.5.11.31.233.307&19277&2392& & &25.121.13.17.19.29&1541&2034 \cr
5008&365881065&16.13.23.37.521&279&242&5026&368357275&4.9.11.19.23.67.113&593&650 \cr
 & &64.9.121.13.23.31&429&736& & &16.3.25.13.23.67.593&4623&4744 \cr
\noalign{\hrule}
 & &9.5.11.13.19.41.73&145&596& & &3.17.19.23.61.271&547&490 \cr
5009&365939145&8.3.25.29.73.149&1313&5038&5027&368426397&4.5.49.23.271.547&11033&1548 \cr
 & &32.11.13.101.229&229&1616& & &32.9.7.11.17.43.59&2537&3696 \cr
\noalign{\hrule}
 & &9.5.11.13.29.37.53&3067&152& & &27.7.23.29.37.79&36355&30874 \cr
5010&365952015&16.3.13.19.3067&4477&4724&5028&368482149&4.5.11.43.359.661&567&94 \cr
 & &128.121.37.1181&1181&704& & &16.81.5.7.47.359&1795&1128 \cr
\noalign{\hrule}
 & &9.17.361.6637&4943&1694& & &5.7.11.13.17.61.71&219&716 \cr
5011&366581421&4.7.121.17.4943&2565&2378&5029&368503135&8.3.13.61.73.179&6137&6930 \cr
 & &16.27.5.7.11.19.29.41&3157&3480& & &32.27.5.7.11.17.361&361&432 \cr
\noalign{\hrule}
 & &5.13.19.37.71.113&6691&6444& & &5.13.59.157.613&89397&91438 \cr
5012&366610985&8.9.113.179.6691&77&20150&5030&369084235&4.27.7.11.43.131.349&551&628 \cr
 & &32.3.25.7.11.13.31&2387&240& & &32.3.19.29.43.157.349&19893&19952 \cr
\noalign{\hrule}
 & &9.5.11.19.127.307&403&232& & &9.125.7.173.271&1177&1448 \cr
5013&366691545&16.11.13.29.31.307&17&324&5031&369203625&16.3.11.107.173.181&329&848 \cr
 & &128.81.13.17.29&6409&576& & &512.7.47.53.181&9593&12032 \cr
\noalign{\hrule}
 & &7.13.29.163.853&70839&68200& & &27.7.11.13.79.173&527&500 \cr
5014&366923921&16.9.25.11.17.31.463&61&402&5032&369378009&8.125.7.11.17.31.173&277&2664 \cr
 & &64.27.25.17.61.67&41175&36448& & &128.9.125.37.277&10249&8000 \cr
\noalign{\hrule}
 & &27.25.7.11.23.307&5263&1798& & &3.5.23.37.53.547&137&22 \cr
5015&366995475&4.3.5.19.29.31.277&49&44&5033&370070115&4.11.37.137.547&2261&2808 \cr
 & &32.49.11.19.29.277&3857&4432& & &64.27.7.11.13.17.19&13167&7072 \cr
\noalign{\hrule}
 & &9.7.19.83.3697&12151&13728& & &9.5.11.169.19.233&119&1978 \cr
5016&367300647&64.27.11.13.29.419&2785&7394&5034&370340685&4.5.7.17.19.23.43&351&466 \cr
 & &256.5.557.3697&557&640& & &16.27.7.13.17.233&119&24 \cr
\noalign{\hrule}
 & &9.7.13.17.23.31.37&3547&856& & &5.13.23.29.83.103&19&84 \cr
5017&367302663&16.31.107.3547&115&3432&5035&370641895&8.3.7.19.23.29.83&455&1122 \cr
 & &256.3.5.11.13.23&11&640& & &32.9.5.49.11.13.17&539&2448 \cr
\noalign{\hrule}
 & &25.13.361.31.101&781&3912& & &9.5.11.23.29.1123&1079&44 \cr
5018&367344575&16.3.25.11.71.163&969&806&5036&370775295&8.121.13.29.83&75&46 \cr
 & &64.9.11.13.17.19.31&187&288& & &32.3.25.13.23.83&83&1040 \cr
\noalign{\hrule}
 & &9.19.43.107.467&427&40& & &49.11.41.107.157&8625&8174 \cr
5019&367422057&16.5.7.19.61.107&3311&3216&5037&371241101&4.3.125.49.23.61.67&533&594 \cr
 & &512.3.49.11.43.67&3283&2816& & &16.81.125.11.13.41.67&8375&8424 \cr
\noalign{\hrule}
 & &9.343.23.31.167&255&88& & &9.49.13.17.37.103&3509&7320 \cr
5020&367572177&16.27.5.11.17.23.31&233&388&5038&371423871&16.27.5.121.29.61&931&2578 \cr
 & &128.11.17.97.233&22601&11968& & &64.5.49.19.1289&1289&3040 \cr
\noalign{\hrule}
 & &5.11.13.23.79.283&38367&33232& & &5.7.11.53.131.139&111&806 \cr
5021&367660865&32.27.49.29.31.67&293&2990&5039&371554645&4.3.11.13.31.37.53&3197&1236 \cr
 & &128.9.5.13.23.293&293&576& & &32.9.23.103.139&207&1648 \cr
\noalign{\hrule}
 & &3.5.7.11.37.79.109&247&138& & &3.25.13.17.29.773&7567&20708 \cr
5022&367991085&4.9.13.19.23.37.79&2573&350&5040&371561775&8.7.23.31.47.167&27115&27828 \cr
 & &16.25.7.23.31.83&713&3320& & &64.9.5.11.17.29.773&33&32 \cr
\noalign{\hrule}
}%
}
$$
\eject
\vglue -23 pt
\noindent\hskip 1 in\hbox to 6.5 in{\ 5041 -- 5076 \hfill\fbd 371767175 -- 378203007\frb}
\vskip -9 pt
$$
\vbox{
\nointerlineskip
\halign{\strut
    \vrule \ \ \hfil \frb #\ 
   &\vrule \hfil \ \ \fbb #\frb\ 
   &\vrule \hfil \ \ \frb #\ \hfil
   &\vrule \hfil \ \ \frb #\ 
   &\vrule \hfil \ \ \frb #\ \ \vrule \hskip 2 pt
   &\vrule \ \ \hfil \frb #\ 
   &\vrule \hfil \ \ \fbb #\frb\ 
   &\vrule \hfil \ \ \frb #\ \hfil
   &\vrule \hfil \ \ \frb #\ 
   &\vrule \hfil \ \ \frb #\ \vrule \cr%
\noalign{\hrule}
 & &25.13.53.113.191&3157&2832& & &27.11.17.191.389&475&86 \cr
5041&371767175&32.3.7.11.41.59.191&159&2260&5059&375135651&4.9.25.19.43.191&847&2782 \cr
 & &256.9.5.7.53.113&63&128& & &16.5.7.121.13.107&5885&728 \cr
\noalign{\hrule}
 & &25.11.13.61.1709&513&19312& & &7.13.19.29.59.127&61&828 \cr
5042&372690175&32.27.17.19.71&2015&2032&5060&375706513&8.9.19.23.29.61&229&1540 \cr
 & &1024.9.5.13.31.127&3937&4608& & &64.3.5.7.11.229&165&7328 \cr
\noalign{\hrule}
 & &7.11.13.29.37.347&575&498& & &27.125.11.13.19.41&1579&19454 \cr
5043&372703331&4.3.25.13.23.83.347&1293&6688&5061&375964875&4.71.137.1579&825&754 \cr
 & &256.9.5.11.19.431&8189&5760& & &16.3.25.11.13.29.137&137&232 \cr
\noalign{\hrule}
 & &11.13.101.131.197&1215&1346& & &9.5.11.17.23.29.67&3571&3824 \cr
5044&372730501&4.243.5.11.101.673&3349&16&5062&376057935&32.3.67.239.3571&19&220 \cr
 & &128.81.17.197&1377&64& & &256.5.11.19.3571&3571&2432 \cr
\noalign{\hrule}
 & &3.5.11.31.269.271&1247&1712& & &25.7.289.43.173&873&572 \cr
5045&372878385&32.29.43.107.271&4275&7378&5063&376227425&8.9.5.11.13.97.173&287&578 \cr
 & &128.9.25.7.17.19.31&1615&1344& & &32.3.7.11.13.289.41&451&624 \cr
\noalign{\hrule}
 & &3.19.37.149.1187&945&242& & &9.25.7.197.1213&4223&4268 \cr
5046&373004067&4.81.5.7.121.149&481&562&5064&376363575&8.5.11.41.97.103.197&203&20088 \cr
 & &16.5.121.13.37.281&3653&4840& & &128.81.7.11.29.31&3069&1856 \cr
\noalign{\hrule}
 & &13.17.37.43.1061&1323&19360& & &3.5.49.37.83.167&2881&2964 \cr
5047&373059271&64.27.5.49.121&37&26&5065&376949895&8.9.7.13.19.37.43.67&1199&1132 \cr
 & &256.3.5.7.11.13.37&165&896& & &64.11.13.19.43.109.283&339317&339872 \cr
\noalign{\hrule}
 & &9.7.13.37.109.113&1895&478& & &5.11.13.349.1511&43821&39284 \cr
5048&373242051&4.3.5.37.239.379&12317&14212&5066&377047385&8.81.7.23.61.541&1513&110 \cr
 & &32.11.17.19.109.113&209&272& & &32.27.5.7.11.17.89&1513&3024 \cr
\noalign{\hrule}
 & &17.41.127.4217&495&4712& & &11.169.43.53.89&1629&650 \cr
5049&373284623&16.9.5.11.17.19.31&273&254&5067&377062829&4.9.25.2197.181&1551&646 \cr
 & &64.27.5.7.11.13.127&2457&1760& & &16.27.5.11.17.19.47&2565&6392 \cr
\noalign{\hrule}
 & &3.5.19.41.43.743&37937&32648& & &9.625.7.11.13.67&427&310 \cr
5050&373324065&16.7.11.53.59.643&2575&1926&5068&377251875&4.3125.49.31.61&803&2322 \cr
 & &64.9.25.53.103.107&16377&17120& & &16.27.11.43.61.73&2623&1752 \cr
\noalign{\hrule}
 & &9.25.11.13.17.683&4577&4302& & &27.7.13.19.59.137&31625&24956 \cr
5051&373583925&4.81.17.23.199.239&2963&1100&5069&377338689&8.125.11.17.23.367&1233&3068 \cr
 & &32.25.11.199.2963&2963&3184& & &64.9.25.13.59.137&25&32 \cr
\noalign{\hrule}
 & &9.5.11.43.97.181&35&508& & &27.11.37.61.563&1357&19474 \cr
5052&373700745&8.3.25.7.97.127&11&86&5070&377395227&4.7.13.23.59.107&725&666 \cr
 & &32.7.11.43.127&889&16& & &16.9.25.7.23.29.37&725&1288 \cr
\noalign{\hrule}
 & &9.5.37.347.647&20977&8138& & &11.19.47.83.463&903&10 \cr
5053&373807485&4.11.13.313.1907&797&1110&5071&377488067&4.3.5.7.43.463&235&228 \cr
 & &16.3.5.11.13.37.797&797&1144& & &32.9.25.19.43.47&225&688 \cr
\noalign{\hrule}
 & &3.5.7.61.71.823&2703&3058& & &3.11.13.17.37.1399&2413&1784 \cr
5054&374263365&4.9.11.17.53.61.139&38855&27448&5072&377507559&16.11.13.19.127.223&1295&1422 \cr
 & &64.5.19.47.73.409&29857&28576& & &64.9.5.7.37.79.223&7805&7584 \cr
\noalign{\hrule}
 & &169.17.151.863&74207&71640& & &27.29.31.103.151&8165&4972 \cr
5055&374389249&16.9.5.7.199.10601&5599&5002&5073&377517969&8.9.5.11.23.71.113&941&302 \cr
 & &64.3.5.7.11.41.61.509&217343&216480& & &32.5.23.151.941&941&1840 \cr
\noalign{\hrule}
 & &9.5.49.11.13.29.41&3397&1388& & &9.11.23.31.53.101&13895&23894 \cr
5056&374909535&8.3.13.43.79.347&161&398&5074&377852211&4.5.7.13.397.919&5797&636 \cr
 & &32.7.23.199.347&7981&3184& & &32.3.5.11.17.31.53&85&16 \cr
\noalign{\hrule}
 & &9.11.17.19.37.317&6349&5380& & &27.5.13.17.19.23.29&539&46 \cr
5057&375058233&8.3.5.7.11.269.907&217&52&5075&378098955&4.3.49.11.19.529&65&464 \cr
 & &64.49.13.31.907&28117&20384& & &128.5.7.11.13.29&11&448 \cr
\noalign{\hrule}
 & &3.5.7.19.229.821&8333&3982& & &3.7.13.23.29.31.67&215&684 \cr
5058&375077955&4.7.11.13.181.641&821&180&5076&378203007&8.27.5.13.19.23.43&319&670 \cr
 & &32.9.5.181.821&181&48& & &32.25.11.19.29.67&475&176 \cr
\noalign{\hrule}
}%
}
$$
\eject
\vglue -23 pt
\noindent\hskip 1 in\hbox to 6.5 in{\ 5077 -- 5112 \hfill\fbd 378401283 -- 383372919\frb}
\vskip -9 pt
$$
\vbox{
\nointerlineskip
\halign{\strut
    \vrule \ \ \hfil \frb #\ 
   &\vrule \hfil \ \ \fbb #\frb\ 
   &\vrule \hfil \ \ \frb #\ \hfil
   &\vrule \hfil \ \ \frb #\ 
   &\vrule \hfil \ \ \frb #\ \ \vrule \hskip 2 pt
   &\vrule \ \ \hfil \frb #\ 
   &\vrule \hfil \ \ \fbb #\frb\ 
   &\vrule \hfil \ \ \frb #\ \hfil
   &\vrule \hfil \ \ \frb #\ 
   &\vrule \hfil \ \ \frb #\ \vrule \cr%
\noalign{\hrule}
 & &9.13.289.361.31&2059&1698& & &5.11.169.17.19.127&801&814 \cr
5077&378401283&4.27.29.31.71.283&3685&4522&5095&381290195&4.9.121.13.37.89.127&4709&10 \cr
 & &16.5.7.11.17.19.67.71&3905&3752& & &16.3.5.17.89.277&277&2136 \cr
\noalign{\hrule}
 & &3.121.29.103.349&21413&8950& & &9.49.13.227.293&1343&1294 \cr
5078&378414069&4.25.49.19.23.179&377&198&5096&381307563&4.13.17.79.227.647&473&50640 \cr
 & &16.9.49.11.13.19.29&741&392& & &128.3.5.11.43.211&11605&2752 \cr
\noalign{\hrule}
 & &9.25.11.13.79.149&773&1202& & &5.17.29.43.59.61&99&394 \cr
5079&378731925&4.3.149.601.773&163&610&5097&381476005&4.9.11.43.61.197&95&34 \cr
 & &16.5.61.163.601&9943&4808& & &16.3.5.11.17.19.197&2167&456 \cr
\noalign{\hrule}
 & &3.25.13.31.83.151&2773&2622& & &25.11.13.841.127&9&136 \cr
5080&378809925&4.9.5.19.23.31.47.59&143&1972&5098&381835025&16.9.5.11.13.17.29&2921&2734 \cr
 & &32.11.13.17.19.23.29&7429&5104& & &64.3.23.127.1367&1367&2208 \cr
\noalign{\hrule}
 & &3.5.11.17.29.59.79&287&582& & &3.5.11.67.179.193&1537&1358 \cr
5081&379149045&4.9.7.17.29.41.97&6175&4526&5099&381917085&4.7.11.29.53.67.97&3043&27900 \cr
 & &16.25.7.13.19.31.73&14105&11096& & &32.9.25.17.31.179&465&272 \cr
\noalign{\hrule}
 & &17.41.157.3469&32705&26268& & &9.5.13.53.97.127&223&158 \cr
5082&379609201&8.3.5.11.31.199.211&205&6&5100&381950595&4.3.53.79.97.223&30353&34540 \cr
 & &32.9.25.11.31.41&3069&400& & &32.5.11.127.157.239&2629&2512 \cr
\noalign{\hrule}
 & &9.25.49.11.31.101&703&692& & &3.11.13.19.173.271&3405&118 \cr
5083&379712025&8.5.49.19.37.101.173&459&4196&5101&382143333&4.9.5.11.59.227&271&260 \cr
 & &64.27.17.173.1049&17833&16608& & &32.25.13.227.271&227&400 \cr
\noalign{\hrule}
 & &9.11.83.113.409&27217&19000& & &3.25.29.43.61.67&5863&2982 \cr
5084&379765089&16.125.17.19.1601&3237&4838&5102&382236675&4.9.5.7.11.13.41.71&67&32 \cr
 & &64.3.5.13.41.59.83&2419&2080& & &256.13.41.67.71&2911&1664 \cr
\noalign{\hrule}
 & &27.5.17.19.23.379&1763&1648& & &3.5.11.13.19.83.113&251&992 \cr
5085&380104785&32.3.17.19.41.43.103&13265&22&5103&382241145&64.5.31.83.251&207&208 \cr
 & &128.5.7.11.379&11&448& & &2048.9.13.23.31.251&23343&23552 \cr
\noalign{\hrule}
 & &81.5.7.11.73.167&1643&2812& & &3.5.11.373.6211&9503&9130 \cr
5086&380176335&8.19.31.37.53.73&55&18&5104&382255995&4.25.121.13.17.43.83&161&5364 \cr
 & &32.9.5.11.19.31.53&1007&496& & &32.9.7.23.83.149&13363&7152 \cr
\noalign{\hrule}
 & &3.31.131.157.199&165&34& & &9.49.13.131.509&475&34 \cr
5087&380633469&4.9.5.11.17.31.157&3275&1592&5105&382270707&4.25.13.17.19.131&33&98 \cr
 & &64.125.131.199&125&32& & &16.3.5.49.11.17.19&1615&88 \cr
\noalign{\hrule}
 & &25.49.11.13.41.53&333&382& & &27.31.41.71.157&3575&664 \cr
5088&380655275&4.9.5.37.41.53.191&23&182&5106&382531599&16.25.11.13.31.83&369&710 \cr
 & &16.3.7.13.23.37.191&2553&1528& & &64.9.125.41.71&125&32 \cr
\noalign{\hrule}
 & &3.5.121.169.17.73&3149&5684& & &5.49.59.103.257&803&2088 \cr
5089&380658135&8.49.17.29.47.67&2041&1242&5107&382638305&16.9.11.29.73.103&247&556 \cr
 & &32.27.13.23.29.157&4553&3312& & &128.3.13.19.29.139&21489&8896 \cr
\noalign{\hrule}
 & &5.31.79.137.227&2691&1556& & &3.125.7.17.23.373&4483&1858 \cr
5090&380807255&8.9.13.23.79.389&6677&8494&5108&382837875&4.23.929.4483&12925&8442 \cr
 & &32.3.11.31.137.607&607&528& & &16.9.25.7.11.47.67&737&1128 \cr
\noalign{\hrule}
 & &9.529.29.31.89&1555&1026& & &9.5.11.13.17.31.113&449&46 \cr
5091&380932371&4.243.5.19.31.311&6721&812&5109&383210685&4.17.23.113.449&65&48 \cr
 & &32.5.7.11.13.29.47&1001&3760& & &128.3.5.13.23.449&449&1472 \cr
\noalign{\hrule}
 & &3.5.17.89.103.163&19&70& & &5.7.841.47.277&1651&288 \cr
5092&381026355&4.25.7.19.103.163&10057&6732&5110&383214265&64.9.5.13.29.127&1133&752 \cr
 & &32.9.11.17.89.113&339&176& & &2048.3.11.47.103&1133&3072 \cr
\noalign{\hrule}
 & &9.25.47.137.263&21593&9232& & &9.169.59.4271&2401&1870 \cr
5093&381027825&32.11.13.151.577&4155&2192&5111&383275269&4.5.2401.11.169.17&89&2784 \cr
 & &1024.3.5.137.277&277&512& & &256.3.49.29.89&1421&11392 \cr
\noalign{\hrule}
 & &81.19.31.61.131&1519&970& & &81.41.241.479&9779&9860 \cr
5094&381242619&4.9.5.49.961.97&917&44&5112&383372919&8.5.7.11.17.29.127.241&61&180 \cr
 & &32.5.343.11.131&1715&176& & &64.9.25.11.29.61.127&44225&44704 \cr
\noalign{\hrule}
}%
}
$$
\eject
\vglue -23 pt
\noindent\hskip 1 in\hbox to 6.5 in{\ 5113 -- 5148 \hfill\fbd 383429475 -- 390938715\frb}
\vskip -9 pt
$$
\vbox{
\nointerlineskip
\halign{\strut
    \vrule \ \ \hfil \frb #\ 
   &\vrule \hfil \ \ \fbb #\frb\ 
   &\vrule \hfil \ \ \frb #\ \hfil
   &\vrule \hfil \ \ \frb #\ 
   &\vrule \hfil \ \ \frb #\ \ \vrule \hskip 2 pt
   &\vrule \ \ \hfil \frb #\ 
   &\vrule \hfil \ \ \fbb #\frb\ 
   &\vrule \hfil \ \ \frb #\ \hfil
   &\vrule \hfil \ \ \frb #\ 
   &\vrule \hfil \ \ \frb #\ \vrule \cr%
\noalign{\hrule}
 & &9.25.11.13.17.701&3857&4558& & &7.53.439.2383&12825&10442 \cr
5113&383429475&4.5.7.13.19.29.43.53&6777&3982&5131&388116827&4.27.25.7.19.23.227&31&1166 \cr
 & &16.27.11.19.181.251&4769&4344& & &16.3.5.11.23.31.53&345&2728 \cr
\noalign{\hrule}
 & &5.11.13.29.107.173&1813&3204& & &9.5.17.317.1601&325&1276 \cr
5114&383825585&8.9.5.49.11.37.89&107&338&5132&388250505&8.3.125.11.13.17.29&3173&3202 \cr
 & &32.3.7.169.37.107&273&592& & &32.11.13.19.167.1601&2717&2672 \cr
\noalign{\hrule}
 & &9.5.11.41.127.149&13837&12198& & &3.11.13.19.29.31.53&5289&4600 \cr
5115&384042285&4.27.19.101.107.137&2665&62&5133&388370697&16.9.25.19.23.41.43&583&196 \cr
 & &16.5.13.31.41.107&1391&248& & &128.25.49.11.23.53&1225&1472 \cr
\noalign{\hrule}
 & &11.17.23.179.499&22475&10998& & &81.5.11.29.31.97&259&578 \cr
5116&384169621&4.9.25.13.29.31.47&1243&3128&5134&388489365&4.3.5.7.289.37.97&583&872 \cr
 & &64.3.5.11.17.23.113&339&160& & &64.7.11.37.53.109&5777&8288 \cr
\noalign{\hrule}
 & &3.121.17.19.29.113&593&650& & &3.5.289.19.53.89&387&5104 \cr
5117&384224973&4.25.11.13.17.29.593&1541&2034&5135&388515705&32.27.5.11.29.43&119&1366 \cr
 & &16.9.23.67.113.593&4623&4744& & &128.7.17.683&4781&64 \cr
\noalign{\hrule}
 & &5.7.23.29.109.151&10263&7102& & &625.13.17.29.97&1449&1364 \cr
5118&384235355&4.3.7.11.53.67.311&261&208&5136&388545625&8.9.125.7.11.13.23.31&359&1734 \cr
 & &128.27.11.13.29.311&8397&9152& & &32.27.289.31.359&9693&8432 \cr
\noalign{\hrule}
 & &125.11.23.29.419&147&272& & &7.11.13.17.53.431&261&170 \cr
5119&384275375&32.3.49.11.17.23.29&321&988&5137&388719331&4.9.5.11.289.29.53&431&1014 \cr
 & &256.9.7.13.19.107&9737&21888& & &16.27.169.29.431&351&232 \cr
\noalign{\hrule}
 & &3.11.19.71.89.97&19631&1030& & &25.29.53.67.151&1349&594 \cr
5120&384315261&4.5.67.103.293&4183&2718&5138&388745725&4.27.5.11.19.53.71&29&1036 \cr
 & &16.9.47.89.151&141&1208& & &32.9.7.11.29.37&3663&112 \cr
\noalign{\hrule}
 & &5.11.17.23.31.577&273&304& & &5.11.169.17.23.107&2667&206 \cr
5121&384659935&32.3.5.7.11.13.17.19.23&123&2308&5139&388874915&4.3.5.7.11.103.127&249&884 \cr
 & &256.9.7.41.577&369&896& & &32.9.7.13.17.83&63&1328 \cr
\noalign{\hrule}
 & &3.25.29.47.53.71&15419&15466& & &9.25.11.31.37.137&23999&13676 \cr
5122&384672675&4.5.11.17.19.37.53.907&9891&86&5140&388919025&8.13.103.233.263&31309&29970 \cr
 & &16.9.7.17.19.43.157&20881&17544& & &32.81.5.37.131.239&2151&2096 \cr
\noalign{\hrule}
 & &9.5.11.17.19.29.83&1049&634& & &5.7.17.19.127.271&131&2028 \cr
5123&384843195&4.19.29.317.1049&249&800&5141&389084185&8.3.5.169.19.131&451&204 \cr
 & &256.3.25.83.317&317&640& & &64.9.11.13.17.41&4797&352 \cr
\noalign{\hrule}
 & &27.5.11.17.101.151&601&1060& & &5.13.37.43.53.71&253&306 \cr
5124&385011495&8.25.53.101.601&63&38&5142&389150645&4.9.5.11.17.23.37.71&13091&1376 \cr
 & &32.9.7.19.53.601&7049&9616& & &256.3.13.19.43.53&57&128 \cr
\noalign{\hrule}
 & &3.11.17.29.131.181&1495&1582& & &25.7.11.31.47.139&377&1152 \cr
5125&385754259&4.5.7.11.13.23.113.131&1737&34&5143&389856775&256.9.7.13.29.47&305&682 \cr
 & &16.9.5.17.113.193&2895&904& & &1024.3.5.11.31.61&183&512 \cr
\noalign{\hrule}
 & &9.11.13.41.71.103&4333&2980& & &9.5.49.47.53.71&741&2596 \cr
5126&385885071&8.3.5.7.13.149.619&5771&3914&5144&389978505&8.27.7.11.13.19.59&1471&2120 \cr
 & &32.7.19.29.103.199&3781&3248& & &128.5.13.53.1471&1471&832 \cr
\noalign{\hrule}
 & &3.5.11.31.163.463&5073&20& & &3.25.121.19.31.73&91&184 \cr
5127&386023935&8.9.25.19.89&163&638&5145&390197775&16.7.11.13.19.23.73&151&360 \cr
 & &32.11.29.163&1&464& & &256.9.5.13.23.151&5889&2944 \cr
\noalign{\hrule}
 & &9.5.19.41.103.107&759&1274& & &11.241.293.503&1863&1360 \cr
5128&386341155&4.27.49.11.13.23.41&401&50&5146&390701729&32.81.5.17.23.241&503&262 \cr
 & &16.25.49.23.401&2005&9016& & &128.9.23.131.503&1179&1472 \cr
\noalign{\hrule}
 & &3.25.19.43.59.107&9061&9878& & &3.5.19.29.41.1153&56551&52984 \cr
5129&386829075&4.25.11.13.17.41.449&2537&2988&5147&390711345&16.11.37.53.97.179&2779&810 \cr
 & &32.9.43.59.83.449&1347&1328& & &64.81.5.7.53.397&10719&11872 \cr
\noalign{\hrule}
 & &19.961.79.269&441&520& & &9.5.17.47.83.131&2981&754 \cr
5130&388022009&16.9.5.49.13.19.269&869&62&5148&390938715&4.11.13.29.47.271&6557&6180 \cr
 & &64.3.5.11.13.31.79&429&160& & &32.3.5.11.79.83.103&1133&1264 \cr
\noalign{\hrule}
}%
}
$$
\eject
\vglue -23 pt
\noindent\hskip 1 in\hbox to 6.5 in{\ 5149 -- 5184 \hfill\fbd 391785381 -- 398731473\frb}
\vskip -9 pt
$$
\vbox{
\nointerlineskip
\halign{\strut
    \vrule \ \ \hfil \frb #\ 
   &\vrule \hfil \ \ \fbb #\frb\ 
   &\vrule \hfil \ \ \frb #\ \hfil
   &\vrule \hfil \ \ \frb #\ 
   &\vrule \hfil \ \ \frb #\ \ \vrule \hskip 2 pt
   &\vrule \ \ \hfil \frb #\ 
   &\vrule \hfil \ \ \fbb #\frb\ 
   &\vrule \hfil \ \ \frb #\ \hfil
   &\vrule \hfil \ \ \frb #\ 
   &\vrule \hfil \ \ \frb #\ \vrule \cr%
\noalign{\hrule}
 & &9.13.23.41.53.67&361&3190& & &9.25.7.11.137.167&23579&22346 \cr
5149&391785381&4.3.5.11.13.361.29&469&1184&5167&396378675&4.7.17.19.73.11173&1243&9930 \cr
 & &256.7.19.37.67&259&2432& & &16.3.5.11.19.113.331&2147&2648 \cr
\noalign{\hrule}
 & &3.7.19.23.79.541&7865&2414& & &81.5.11.23.53.73&443&178 \cr
5150&392215803&4.5.7.121.13.17.71&81&10&5168&396437085&4.3.11.73.89.443&599&380 \cr
 & &16.81.25.121.17&11475&968& & &32.5.19.443.599&11381&7088 \cr
\noalign{\hrule}
 & &9.5.31.317.887&2059&2376& & &9.11.13.23.59.227&1375&1316 \cr
5151&392244705&16.243.11.29.31.71&887&1786&5169&396446193&8.125.7.121.47.227&5681&6 \cr
 & &64.19.47.71.887&1349&1504& & &32.3.5.7.13.19.23&665&16 \cr
\noalign{\hrule}
 & &9.13.29.37.53.59&497&1570& & &125.13.17.113.127&1421&1404 \cr
5152&392566707&4.3.5.7.59.71.157&1859&2330&5170&396446375&8.27.5.49.169.29.127&867&22 \cr
 & &16.25.7.11.169.233&5825&8008& & &32.81.7.11.289.29&6237&7888 \cr
\noalign{\hrule}
 & &11.13.89.109.283&437&720& & &9.5.31.59.61.79&16247&25702 \cr
5153&392589769&32.9.5.11.19.23.109&37&290&5171&396627795&4.7.11.71.181.211&527&1794 \cr
 & &128.3.25.19.29.37&20387&4800& & &16.3.13.17.23.31.71&1207&2392 \cr
\noalign{\hrule}
 & &3.5.7.11.13.17.23.67&139&206& & &27.5.17.29.67.89&3913&3652 \cr
5154&393347955&4.7.11.13.17.103.139&207&1340&5172&396867465&8.3.7.11.13.43.67.83&247&1160 \cr
 & &32.9.5.23.67.139&139&48& & &128.5.169.19.29.43&3211&2752 \cr
\noalign{\hrule}
 & &3.5.7.11.461.739&353&386& & &9.5.11.23.67.521&9367&7826 \cr
5155&393484245&4.5.7.193.353.461&477&1828&5173&397416195&4.3.5.7.13.17.19.29.43&11&46 \cr
 & &32.9.53.353.457&18709&21936& & &16.11.13.17.23.29.43&1247&1768 \cr
\noalign{\hrule}
 & &9.5.49.29.61.101&3301&2860& & &5.7.11.13.19.37.113&13161&13294 \cr
5156&393965145&8.25.11.13.29.3301&2013&1288&5174&397592195&4.3.289.23.41.107.113&121&11970 \cr
 & &128.3.7.121.13.23.61&1573&1472& & &16.27.5.7.121.19.23&297&184 \cr
\noalign{\hrule}
 & &9.5.121.13.19.293&10199&10492& & &3.7.19.29.163.211&1903&2106 \cr
5157&394060095&8.5.7.13.31.43.47.61&2319&304&5175&397961403&4.243.11.13.163.173&203&40 \cr
 & &256.3.7.19.47.773&5411&6016& & &64.5.7.11.13.29.173&2249&1760 \cr
\noalign{\hrule}
 & &3.11.13.29.79.401&161&240& & &27.25.7.11.13.19.31&11&236 \cr
5158&394118439&32.9.5.7.11.13.23.29&535&158&5176&397972575&8.3.7.121.31.59&2185&356 \cr
 & &128.25.23.79.107&2675&1472& & &64.5.19.23.89&2047&32 \cr
\noalign{\hrule}
 & &3.5.7.11.23.37.401&731&472& & &9.121.29.12611&4551&8060 \cr
5159&394144905&16.5.11.17.23.43.59&27&962&5177&398267991&8.27.5.13.31.37.41&493&506 \cr
 & &64.27.13.37.59&117&1888& & &32.5.11.17.23.29.31.41&6355&6256 \cr
\noalign{\hrule}
 & &9.13.23.29.31.163&235&142& & &9.5.13.43.71.223&6179&3410 \cr
5160&394331067&4.3.5.23.47.71.163&2407&1342&5178&398279115&4.3.25.11.31.37.167&43931&42094 \cr
 & &16.11.29.47.61.83&3901&5368& & &16.13.197.223.1619&1619&1576 \cr
\noalign{\hrule}
 & &9.5.11.47.71.239&2275&5612& & &9.5.7.361.31.113&803&214 \cr
5161&394783785&8.3.125.7.13.23.61&2627&248&5179&398343645&4.5.7.11.19.73.107&3093&3842 \cr
 & &128.7.31.37.71&217&2368& & &16.3.11.17.113.1031&1031&1496 \cr
\noalign{\hrule}
 & &9.5.13.19.961.37&427&428& & &5.13.113.227.239&2121&986 \cr
5162&395216055&8.7.13.961.37.61.107&1653&14146&5180&398488285&4.3.7.17.29.101.113&2497&780 \cr
 & &32.3.11.19.29.107.643&18647&18832& & &32.9.5.7.11.13.227&63&176 \cr
\noalign{\hrule}
 & &3.25.49.41.43.61&1863&638& & &25.11.23.29.41.53&20961&6386 \cr
5163&395220525&4.243.11.23.29.43&5065&5384&5181&398582525&4.9.17.31.103.137&1391&1802 \cr
 & &64.5.23.673.1013&23299&21536& & &16.3.13.289.53.107&3757&2568 \cr
\noalign{\hrule}
 & &3.47.71.127.311&43&170& & &27.23.263.2441&24521&31622 \cr
5164&395404467&4.5.17.43.47.311&855&1166&5182&398671443&4.7.31.97.113.163&12205&23166 \cr
 & &16.9.25.11.17.19.53&14025&8056& & &16.81.5.11.13.2441&143&120 \cr
\noalign{\hrule}
 & &27.25.343.29.59&10051&10186& & &9.7.11.31.67.277&545&478 \cr
5165&396139275&4.5.11.19.529.29.463&879&4214&5183&398702997&4.3.5.7.109.239.277&10117&15934 \cr
 & &16.3.49.19.23.43.293&6739&6536& & &16.5.31.67.151.257&1285&1208 \cr
\noalign{\hrule}
 & &3.7.11.19.137.659&1109&1768& & &9.49.19.23.2069&793&1276 \cr
5166&396252087&16.11.13.17.19.1109&659&450&5184&398731473&8.3.7.11.13.19.29.61&115&94 \cr
 & &64.9.25.13.17.659&663&800& & &32.5.13.23.29.47.61&8845&9776 \cr
\noalign{\hrule}
}%
}
$$
\eject
\vglue -23 pt
\noindent\hskip 1 in\hbox to 6.5 in{\ 5185 -- 5220 \hfill\fbd 399229935 -- 404556867\frb}
\vskip -9 pt
$$
\vbox{
\nointerlineskip
\halign{\strut
    \vrule \ \ \hfil \frb #\ 
   &\vrule \hfil \ \ \fbb #\frb\ 
   &\vrule \hfil \ \ \frb #\ \hfil
   &\vrule \hfil \ \ \frb #\ 
   &\vrule \hfil \ \ \frb #\ \ \vrule \hskip 2 pt
   &\vrule \ \ \hfil \frb #\ 
   &\vrule \hfil \ \ \fbb #\frb\ 
   &\vrule \hfil \ \ \frb #\ \hfil
   &\vrule \hfil \ \ \frb #\ 
   &\vrule \hfil \ \ \frb #\ \vrule \cr%
\noalign{\hrule}
 & &3.5.13.31.211.313&495121&495524& & &9.11.13.53.71.83&4891&4820 \cr
5185&399229935&8.11.19.23.73.103.1697&9&1688&5203&401967423&8.5.11.53.67.73.241&3783&232 \cr
 & &128.9.11.19.103.211&1957&2112& & &128.3.13.29.97.241&6989&6208 \cr
\noalign{\hrule}
 & &3.25.11.169.47.61&5649&2294& & &5.7.17.541.1249&565&684 \cr
5186&399731475&4.9.5.7.31.37.269&3337&5002&5204&402046855&8.9.25.19.113.541&233&308 \cr
 & &16.7.41.47.61.71&497&328& & &64.3.7.11.19.113.233&23617&22368 \cr
\noalign{\hrule}
 & &97.139.149.199&5&144& & &3.25.13.37.71.157&1&924 \cr
5187&399784433&32.9.5.97.199&143&342&5205&402128025&8.9.7.11.157&703&710 \cr
 & &128.81.11.13.19&20007&704& & &32.5.11.19.37.71&19&176 \cr
\noalign{\hrule}
 & &25.31.457.1129&7029&21196& & &3.7.23.31.107.251&1235&522 \cr
5188&399863575&8.9.7.11.71.757&793&1550&5206&402129861&4.27.5.13.19.29.107&2387&502 \cr
 & &32.3.25.7.13.31.61&793&336& & &16.7.11.19.31.251&19&88 \cr
\noalign{\hrule}
 & &9.25.29.101.607&231&376& & &81.17.37.53.149&5447&2450 \cr
5189&400028175&16.27.5.7.11.47.101&17&118&5207&402344253&4.25.49.13.17.419&2783&2664 \cr
 & &64.7.11.17.47.59&5593&20768& & &64.9.25.7.121.23.37&4025&3872 \cr
\noalign{\hrule}
 & &3.25.7.29.41.641&177&464& & &9.13.23.157.953&1775&20144 \cr
5190&400128225&32.9.25.841.59&533&308&5208&402630111&32.25.71.1259&665&594 \cr
 & &256.7.11.13.41.59&767&1408& & &128.27.125.7.11.19&4389&8000 \cr
\noalign{\hrule}
 & &25.11.37.67.587&1099&1836& & &25.7.23.71.1409&2717&1308 \cr
5191&400172575&8.27.5.7.17.37.157&559&226&5209&402656975&8.3.11.13.19.71.109&75&1274 \cr
 & &32.3.7.13.17.43.113&14577&24752& & &32.9.25.49.169&1521&112 \cr
\noalign{\hrule}
 & &5.7.121.17.67.83&93&26& & &9.5.1331.53.127&1679&2314 \cr
5192&400364195&4.3.5.121.13.31.83&147&268&5210&403153245&4.3.13.23.53.73.89&6601&7550 \cr
 & &32.9.49.13.31.67&403&1008& & &16.25.7.529.41.151&30955&29624 \cr
\noalign{\hrule}
 & &9.5.11.17.19.23.109&1407&208& & &7.11.19.239.1153&46661&42120 \cr
5193&400831695&32.27.7.13.23.67&277&344&5211&403154521&16.81.5.13.29.1609&1153&2762 \cr
 & &512.7.13.43.277&11911&23296& & &64.3.13.1153.1381&1381&1248 \cr
\noalign{\hrule}
 & &81.5.11.17.67.79&1561&6854& & &7.11.13.17.19.29.43&515&486 \cr
5194&400865355&4.9.7.23.149.223&85&76&5212&403183781&4.243.5.17.19.43.103&4669&584 \cr
 & &32.5.17.19.149.223&2831&3568& & &64.81.7.23.29.73&1863&2336 \cr
\noalign{\hrule}
 & &243.25.11.17.353&2459&2216& & &3.7.11.17.23.41.109&4553&8854 \cr
5195&401016825&16.277.353.2459&1053&1406&5213&403644549&4.7.19.29.157.233&3445&3312 \cr
 & &64.81.13.19.37.277&9139&8864& & &128.9.5.13.23.53.157&10335&10048 \cr
\noalign{\hrule}
 & &5.11.79.127.727&403&324& & &7.11.23.457.499&18983&19440 \cr
5196&401169505&8.81.5.11.13.31.127&281&1424&5214&403864153&32.243.5.23.41.463&19223&18280 \cr
 & &256.9.13.89.281&25009&14976& & &512.3.25.47.409.457&19223&19200 \cr
\noalign{\hrule}
 & &27.25.7.17.19.263&2269&2728& & &5.7.11.13.17.47.101&207&592 \cr
5197&401384025&16.25.7.11.31.2269&2097&172&5215&403898495&32.9.13.23.37.101&99&200 \cr
 & &128.9.31.43.233&7223&2752& & &512.81.25.11.37&2997&1280 \cr
\noalign{\hrule}
 & &27.5.17.29.37.163&1079&586& & &9.5.7.11.31.53.71&107&1598 \cr
5198&401393205&4.3.13.83.163.293&18389&22198&5216&404202645&4.3.17.47.53.107&1085&1406 \cr
 & &16.7.11.37.71.1009&7063&6248& & &16.5.7.17.19.31.37&629&152 \cr
\noalign{\hrule}
 & &3.5.17.19.41.43.47&17941&18672& & &9.5.23.31.43.293&5183&3718 \cr
5199&401461545&32.9.5.7.11.233.389&2021&76&5217&404238915&4.11.169.31.71.73&645&304 \cr
 & &256.7.11.19.43.47&77&128& & &128.3.5.13.19.43.71&923&1216 \cr
\noalign{\hrule}
 & &3.5.13.17.31.3907&3611&296& & &3.121.47.137.173&37&84 \cr
5200&401502855&16.23.31.37.157&495&652&5218&404362761&8.9.7.37.137.173&847&710 \cr
 & &128.9.5.11.23.163&5379&1472& & &32.5.49.121.37.71&3479&2960 \cr
\noalign{\hrule}
 & &81.5.17.23.43.59&347&642& & &27.7.37.151.383&8671&5500 \cr
5201&401746635&4.243.17.107.347&10051&15950&5219&404426169&8.9.125.11.13.23.29&31&86 \cr
 & &16.25.11.19.529.29&2755&2024& & &32.25.23.29.31.43&24725&14384 \cr
\noalign{\hrule}
 & &7.11.169.17.23.79&519&350& & &9.11.13.23.79.173&557&730 \cr
5202&401958557&4.3.25.49.17.23.173&649&9126&5220&404556867&4.5.23.73.79.557&2301&484 \cr
 & &16.81.11.169.59&81&472& & &32.3.121.13.59.73&649&1168 \cr
\noalign{\hrule}
}%
}
$$
\eject
\vglue -23 pt
\noindent\hskip 1 in\hbox to 6.5 in{\ 5221 -- 5256 \hfill\fbd 405161757 -- 410947371\frb}
\vskip -9 pt
$$
\vbox{
\nointerlineskip
\halign{\strut
    \vrule \ \ \hfil \frb #\ 
   &\vrule \hfil \ \ \fbb #\frb\ 
   &\vrule \hfil \ \ \frb #\ \hfil
   &\vrule \hfil \ \ \frb #\ 
   &\vrule \hfil \ \ \frb #\ \ \vrule \hskip 2 pt
   &\vrule \ \ \hfil \frb #\ 
   &\vrule \hfil \ \ \fbb #\frb\ 
   &\vrule \hfil \ \ \frb #\ \hfil
   &\vrule \hfil \ \ \frb #\ 
   &\vrule \hfil \ \ \frb #\ \vrule \cr%
\noalign{\hrule}
 & &81.7.11.13.19.263&10585&10718& & &9.11.13.19.59.283&1397&1150 \cr
5221&405161757&4.5.11.13.23.29.73.233&309&10&5239&408291741&4.25.121.23.59.127&3211&10704 \cr
 & &16.3.25.73.103.233&17009&20600& & &128.3.5.169.19.223&1115&832 \cr
\noalign{\hrule}
 & &9.5.7.11.19.61.101&5627&2150& & &25.29.31.37.491&3113&2622 \cr
5222&405609435&4.3.125.17.43.331&1159&1034&5240&408303325&4.3.5.11.19.23.29.283&111&3224 \cr
 & &16.11.19.47.61.331&331&376& & &64.9.13.19.31.37&117&608 \cr
\noalign{\hrule}
 & &3.17.23.277.1249&429&820& & &9.11.13.41.71.109&29155&28942 \cr
5223&405826329&8.9.5.11.13.41.277&323&46&5241&408363813&4.3.5.343.11.17.29.499&199&8284 \cr
 & &32.5.11.13.17.19.23&209&1040& & &32.7.19.29.109.199&3781&3248 \cr
\noalign{\hrule}
 & &5.11.13.529.29.37&387&686& & &7.121.41.61.193&3015&4366 \cr
5224&405846155&4.9.5.343.11.23.43&241&494&5242&408840971&4.9.5.37.41.59.67&2509&4928 \cr
 & &16.3.7.13.19.43.241&4579&7224& & &512.3.5.7.11.13.193&195&256 \cr
\noalign{\hrule}
 & &3.5.11.47.199.263&629&366& & &3.5.121.169.31.43&4541&2726 \cr
5225&405873435&4.9.11.17.37.47.61&3265&398&5243&408877755&4.19.29.31.47.239&14157&2924 \cr
 & &16.5.17.199.653&653&136& & &32.9.121.13.17.43&17&48 \cr
\noalign{\hrule}
 & &27.5.31.293.331&253&584& & &9.5.7.11.13.31.293&145&548 \cr
5226&405873855&16.5.11.23.73.293&331&1134&5244&409143735&8.25.29.137.293&1859&1566 \cr
 & &64.81.7.23.331&69&224& & &32.27.11.169.841&841&624 \cr
\noalign{\hrule}
 & &81.11.23.29.683&151&470& & &5.11.13.17.151.223&10759&22914 \cr
5227&405904851&4.3.5.47.151.683&8671&1574&5245&409295315&4.9.7.19.29.53.67&403&604 \cr
 & &16.13.23.29.787&787&104& & &32.3.7.13.29.31.151&651&464 \cr
\noalign{\hrule}
 & &3.49.19.23.71.89&775&1272& & &81.5.7.353.409&7139&21454 \cr
5228&405926241&16.9.25.7.19.31.53&89&44&5246&409308795&4.121.17.59.631&409&222 \cr
 & &128.5.11.31.53.89&2915&1984& & &16.3.11.37.59.409&649&296 \cr
\noalign{\hrule}
 & &5.13.53.179.659&88481&86154& & &7.23.67.163.233&24383&13596 \cr
5229&406375645&4.3.23.83.173.3847&3913&66&5247&409679473&8.3.11.37.103.659&33&70 \cr
 & &16.9.7.11.13.43.83&6391&3096& & &32.9.5.7.121.659&5931&9680 \cr
\noalign{\hrule}
 & &3.71.83.127.181&215&34& & &25.17.173.5573&4257&1316 \cr
5230&406387173&4.5.17.43.71.127&99&28&5248&409754825&8.9.25.7.11.43.47&1313&238 \cr
 & &32.9.5.7.11.17.43&1155&11696& & &32.3.49.13.17.101&4949&624 \cr
\noalign{\hrule}
 & &81.11.17.47.571&3217&3064& & &27.5.7.353.1229&791&438 \cr
5231&406500039&16.9.47.383.3217&403&20&5249&409975965&4.81.5.49.73.113&13519&4366 \cr
 & &128.5.13.31.3217&16085&25792& & &16.11.37.59.1229&649&296 \cr
\noalign{\hrule}
 & &27.125.49.23.107&4687&3938& & &9.11.23.101.1783&8645&10968 \cr
5232&406987875&4.9.7.11.43.109.179&2461&850&5250&410048991&16.27.5.7.13.19.457&485&28 \cr
 & &16.25.17.23.107.109&109&136& & &128.25.49.13.97&4753&20800 \cr
\noalign{\hrule}
 & &11.13.17.29.53.109&399&290& & &27.5.17.19.23.409&1421&194 \cr
5233&407272723&4.3.5.7.11.17.19.841&1075&234&5251&410192235&4.9.49.23.29.97&103&770 \cr
 & &16.27.125.13.19.43&3375&6536& & &16.5.343.11.103&343&9064 \cr
\noalign{\hrule}
 & &25.11.23.37.1741&26443&17082& & &81.25.11.113.163&419&824 \cr
5234&407437525&4.9.13.31.73.853&997&1850&5252&410283225&16.5.103.163.419&339&176 \cr
 & &16.3.25.31.37.997&997&744& & &512.3.11.113.419&419&256 \cr
\noalign{\hrule}
 & &3.5.11.13.17.53.211&3895&3684& & &27.5.11.13.89.239&1333&2528 \cr
5235&407788095&8.9.25.17.19.41.307&5291&2384&5253&410636655&64.31.43.79.89&23&66 \cr
 & &256.11.13.37.41.149&5513&5248& & &256.3.11.23.31.79&1817&3968 \cr
\noalign{\hrule}
 & &3.25.17.19.113.149&157&182& & &27.5.37.83.991&12731&23936 \cr
5236&407876325&4.7.13.17.19.149.157&81&2750&5254&410853735&256.11.17.29.439&2429&2400 \cr
 & &16.81.125.7.11.13&455&2376& & &16384.3.25.7.17.347&41293&40960 \cr
\noalign{\hrule}
 & &11.37.503.1993&20267&1656& & &3.11.19.47.73.191&18775&18966 \cr
5237&408008953&16.9.13.23.1559&745&814&5255&410886267&4.9.25.19.29.109.751&803&52 \cr
 & &64.3.5.11.13.37.149&745&1248& & &32.5.11.13.29.73.109&1885&1744 \cr
\noalign{\hrule}
 & &11.13.37.179.431&4825&20772& & &27.19.23.29.1201&16225&16202 \cr
5238&408195359&8.9.25.193.577&769&962&5256&410947371&4.25.11.19.29.59.8101&2673&5428 \cr
 & &32.3.25.13.37.769&769&1200& & &32.243.5.121.23.3481&17405&17424 \cr
\noalign{\hrule}
}%
}
$$
\eject
\vglue -23 pt
\noindent\hskip 1 in\hbox to 6.5 in{\ 5257 -- 5292 \hfill\fbd 410975675 -- 416281635\frb}
\vskip -9 pt
$$
\vbox{
\nointerlineskip
\halign{\strut
    \vrule \ \ \hfil \frb #\ 
   &\vrule \hfil \ \ \fbb #\frb\ 
   &\vrule \hfil \ \ \frb #\ \hfil
   &\vrule \hfil \ \ \frb #\ 
   &\vrule \hfil \ \ \frb #\ \ \vrule \hskip 2 pt
   &\vrule \ \ \hfil \frb #\ 
   &\vrule \hfil \ \ \fbb #\frb\ 
   &\vrule \hfil \ \ \frb #\ \hfil
   &\vrule \hfil \ \ \frb #\ 
   &\vrule \hfil \ \ \frb #\ \vrule \cr%
\noalign{\hrule}
 & &25.11.841.1777&10507&10518& & &11.19.29.961.71&4131&1930 \cr
5257&410975675&4.3.7.19.79.1753.1777&1765&12&5275&413548091&4.243.5.17.31.193&2059&1222 \cr
 & &32.9.5.7.19.79.353&22239&24016& & &16.9.5.13.29.47.71&585&376 \cr
\noalign{\hrule}
 & &27.49.11.13.41.53&277&200& & &9.7.13.31.43.379&21505&22838 \cr
5258&411107697&16.3.25.7.13.41.277&1583&1492&5276&413764533&4.5.7.11.17.19.23.601&27&3032 \cr
 & &128.277.373.1583&103321&101312& & &64.27.11.17.379&51&352 \cr
\noalign{\hrule}
 & &3.7.19.59.101.173&3585&2374& & &3.13.17.19.107.307&3305&1914 \cr
5259&411331893&4.9.5.19.239.1187&3071&7612&5277&413798853&4.9.5.11.19.29.661&2149&2810 \cr
 & &32.5.11.37.83.173&2035&1328& & &16.25.7.11.281.307&1925&2248 \cr
\noalign{\hrule}
 & &49.41.239.857&1265&408& & &5.49.17.23.29.149&10293&6842 \cr
5260&411489407&16.3.5.7.11.17.23.41&1433&2868&5278&413930195&4.3.7.11.47.73.311&279&232 \cr
 & &128.9.239.1433&1433&576& & &64.27.11.29.31.311&9207&9952 \cr
\noalign{\hrule}
 & &729.11.19.37.73&30475&22742& & &49.121.13.41.131&6021&5900 \cr
5261&411527061&4.25.23.53.83.137&11&126&5279&413980567&8.27.25.7.41.59.223&1703&142 \cr
 & &16.9.5.7.11.53.83&581&2120& & &32.3.5.13.59.71.131&1065&944 \cr
\noalign{\hrule}
 & &9.169.167.1621&1571&50& & &3.11.31.47.79.109&2527&1070 \cr
5262&411745347&4.25.167.1571&869&702&5280&414025491&4.5.7.361.79.107&1125&376 \cr
 & &16.27.25.11.13.79&275&1896& & &64.9.625.19.47&625&1824 \cr
\noalign{\hrule}
 & &27.5.13.41.59.97&1207&1412& & &9.5.7.29.61.743&109&94 \cr
5263&411798465&8.13.17.59.71.353&4389&200&5281&414025605&4.3.47.61.109.743&319&2548 \cr
 & &128.3.25.7.11.17.19&1463&5440& & &32.49.11.13.29.109&1001&1744 \cr
\noalign{\hrule}
 & &3.5.7.121.17.1907&1861&46& & &9.19.23.31.43.79&1289&1160 \cr
5264&411883395&4.7.17.23.1861&871&990&5282&414172431&16.3.5.19.23.29.1289&11&1300 \cr
 & &16.9.5.11.13.23.67&299&1608& & &128.125.11.13.29&1625&20416 \cr
\noalign{\hrule}
 & &9.25.7.11.37.643&437&206& & &7.29.53.139.277&333&7700 \cr
5265&412179075&4.3.25.19.23.37.103&203&2572&5283&414253777&8.9.25.49.11.37&139&106 \cr
 & &32.7.19.29.643&29&304& & &32.3.5.37.53.139&37&240 \cr
\noalign{\hrule}
 & &9.7.107.193.317&8899&11752& & &27.5.49.31.43.47&331&754 \cr
5266&412421121&16.7.11.13.113.809&4755&3946&5284&414436365&4.3.7.13.29.43.331&10819&2090 \cr
 & &64.3.5.13.317.1973&1973&2080& & &16.5.11.19.31.349&349&1672 \cr
\noalign{\hrule}
 & &5.7.11.13.23.3583&7313&10602& & &9.5.23.37.79.137&24559&13736 \cr
5267&412457045&4.9.7.19.31.71.103&407&1756&5285&414466785&16.17.41.101.599&1771&2370 \cr
 & &32.3.11.31.37.439&3441&7024& & &64.3.5.7.11.17.23.79&187&224 \cr
\noalign{\hrule}
 & &9.13.41.113.761&22211&8990& & &67.89.197.353&65&132 \cr
5268&412508421&4.5.7.19.29.31.167&351&484&5286&414672983&8.3.5.11.13.89.353&399&46 \cr
 & &32.27.121.13.29.31&2697&1936& & &32.9.7.11.13.19.23&3059&20592 \cr
\noalign{\hrule}
 & &9.17.41.157.419&3421&20600& & &9.25.49.121.311&457&632 \cr
5269&412656759&16.25.11.103.311&3321&4454&5287&414881775&16.7.79.311.457&121&432 \cr
 & &64.81.17.41.131&131&288& & &512.27.121.457&457&768 \cr
\noalign{\hrule}
 & &3.25.7.11.13.23.239&37&202& & &27.7.11.89.2243&11401&9158 \cr
5270&412687275&4.5.7.13.23.37.101&239&216&5288&415024533&4.9.13.19.241.877&5767&22430 \cr
 & &64.27.37.101.239&909&1184& & &16.5.73.79.2243&395&584 \cr
\noalign{\hrule}
 & &11.23.37.157.281&1289&438& & &5.7.121.263.373&14229&14492 \cr
5271&412979237&4.3.73.281.1289&785&504&5289&415449265&8.27.5.11.17.31.3623&59&26 \cr
 & &64.27.5.7.73.157&1971&1120& & &32.9.13.31.59.3623&112313&110448 \cr
\noalign{\hrule}
 & &5.11.13.41.73.193&10189&22734& & &3.25.11.13.17.43.53&373&186 \cr
5272&413019035&4.27.23.421.443&853&410&5290&415518675&4.9.25.31.53.373&139&86 \cr
 & &16.9.5.23.41.853&853&1656& & &16.31.43.139.373&4309&2984 \cr
\noalign{\hrule}
 & &27.25.49.13.961&11837&12188& & &9.5.7.23.67.857&5029&2684 \cr
5273&413205975&8.343.11.19.89.277&33&310&5291&416000655&8.11.23.47.61.107&113&1290 \cr
 & &32.3.5.121.19.31.89&1691&1936& & &32.3.5.43.47.113&5311&688 \cr
\noalign{\hrule}
 & &5.7.13.17.19.29.97&109&24& & &9.5.7.11.17.37.191&361&46 \cr
5274&413412545&16.3.13.29.97.109&2211&950&5292&416281635&4.17.361.23.191&1443&1804 \cr
 & &64.9.25.11.19.67&495&2144& & &32.3.11.13.23.37.41&299&656 \cr
\noalign{\hrule}
}%
}
$$
\eject
\vglue -23 pt
\noindent\hskip 1 in\hbox to 6.5 in{\ 5293 -- 5328 \hfill\fbd 416396241 -- 423853881\frb}
\vskip -9 pt
$$
\vbox{
\nointerlineskip
\halign{\strut
    \vrule \ \ \hfil \frb #\ 
   &\vrule \hfil \ \ \fbb #\frb\ 
   &\vrule \hfil \ \ \frb #\ \hfil
   &\vrule \hfil \ \ \frb #\ 
   &\vrule \hfil \ \ \frb #\ \ \vrule \hskip 2 pt
   &\vrule \ \ \hfil \frb #\ 
   &\vrule \hfil \ \ \fbb #\frb\ 
   &\vrule \hfil \ \ \frb #\ \hfil
   &\vrule \hfil \ \ \frb #\ 
   &\vrule \hfil \ \ \frb #\ \vrule \cr%
\noalign{\hrule}
 & &27.109.151.937&53105&49028& & &7.121.29.71.241&801&7790 \cr
5293&416396241&8.5.7.13.17.19.43.103&4983&6322&5311&420297493&4.9.5.7.19.41.89&241&374 \cr
 & &32.3.11.29.43.109.151&473&464& & &16.3.11.17.89.241&89&408 \cr
\noalign{\hrule}
 & &9.5.7.19.257.271&31831&37816& & &9.7.31.283.761&9295&8534 \cr
5294&416837295&16.29.139.163.229&14649&8008&5312&420603939&4.5.11.169.17.31.251&5327&2454 \cr
 & &256.3.7.11.13.19.257&143&128& & &16.3.5.7.11.409.761&409&440 \cr
\noalign{\hrule}
 & &11.17.23.31.53.59&387&970& & &61.449.15361&7711&7650 \cr
5295&416926037&4.9.5.17.31.43.97&839&742&5313&420722429&4.9.25.11.17.449.701&7223&488 \cr
 & &16.3.5.7.43.53.839&5873&5160& & &64.3.5.17.31.61.233&7905&7456 \cr
\noalign{\hrule}
 & &3.625.83.2687&437&188& & &3.289.19.107.239&583&1450 \cr
5296&418164375&8.19.23.47.2687&1125&1562&5314&421264029&4.25.11.29.53.239&889&306 \cr
 & &32.9.125.11.47.71&1551&1136& & &16.9.5.7.17.29.127&1015&3048 \cr
\noalign{\hrule}
 & &9.5.7.17.23.43.79&5&124& & &81.5.11.13.19.383&581&568 \cr
5297&418391505&8.3.25.23.31.79&2057&82&5315&421447455&16.27.5.7.11.19.71.83&1817&13022 \cr
 & &32.121.17.41&4961&16& & &64.7.17.23.79.383&2737&2528 \cr
\noalign{\hrule}
 & &61.109.113.557&375947&375390& & &19.37.43.73.191&15019&15210 \cr
5298&418494709&4.9.5.121.13.43.97.239&595&122&5316&421482947&4.9.5.169.23.73.653&209&3056 \cr
 & &16.3.25.7.11.13.17.61.97&31525&31416& & &128.3.11.13.19.23.191&759&832 \cr
\noalign{\hrule}
 & &27.25.7.283.313&7733&92& & &25.7.121.43.463&477&598 \cr
5299&418535775&8.7.11.19.23.37&313&390&5317&421573075&4.9.7.13.23.53.463&473&10 \cr
 & &32.3.5.13.23.313&299&16& & &16.3.5.11.13.43.53&39&424 \cr
\noalign{\hrule}
 & &3.5.11.17.19.29.271&145&178& & &27.5.17.397.463&20273&19082 \cr
5300&418845405&4.25.841.89.271&1547&22572&5318&421846245&4.9.7.11.19.29.47.97&397&496 \cr
 & &32.27.7.11.13.17.19&63&208& & &128.7.29.31.97.397&6293&6208 \cr
\noalign{\hrule}
 & &11.31.37.89.373&405&3698& & &5.343.17.41.353&58797&62282 \cr
5301&418846549&4.81.5.31.1849&209&178&5319&421960315&4.9.11.19.47.139.149&145&4 \cr
 & &16.9.5.11.19.43.89&387&760& & &32.3.5.11.19.29.139&4031&10032 \cr
\noalign{\hrule}
 & &5.7.13.19.47.1031&27&638& & &9.343.11.289.43&6533&8216 \cr
5302&418910765&4.27.11.29.1031&559&472&5320&421983639&16.13.17.47.79.139&675&3038 \cr
 & &64.9.11.13.43.59&5841&1376& & &64.27.25.49.13.31&775&1248 \cr
\noalign{\hrule}
 & &27.25.13.23.31.67&59&266& & &81.7.11.71.953&3781&2890 \cr
5303&419190525&4.3.7.19.31.59.67&1015&814&5321&422014131&4.5.289.19.71.199&1287&2494 \cr
 & &16.5.49.11.19.29.37&11803&7448& & &16.9.5.11.13.17.29.43&3655&3016 \cr
\noalign{\hrule}
 & &27.25.67.73.127&22327&21692& & &3.11.13.17.103.563&265&44 \cr
5304&419280975&8.3.5.11.17.29.83.269&127&142&5322&422913777&8.5.121.53.563&221&342 \cr
 & &32.11.17.29.71.83.127&22649&22576& & &32.9.5.13.17.19.53&285&848 \cr
\noalign{\hrule}
 & &9.5.11.17.19.43.61&1421&262& & &11.841.149.307&85&234 \cr
5305&419378355&4.5.49.29.43.131&177&478&5323&423168493&4.9.5.13.17.29.307&2537&2682 \cr
 & &16.3.7.29.59.239&11977&1912& & &16.81.13.43.59.149&4779&4472 \cr
\noalign{\hrule}
 & &5.7.169.89.797&227&396& & &9.7.13.37.61.229&5497&21296 \cr
5306&419568695&8.9.5.11.227.797&169&966&5324&423302607&32.1331.23.239&785&546 \cr
 & &32.27.7.11.169.23&621&176& & &128.3.5.7.13.23.157&785&1472 \cr
\noalign{\hrule}
 & &9.5.11.23.29.31.41&329&122& & &3.2197.17.19.199&5159&4990 \cr
5307&419639715&4.5.7.29.31.47.61&297&602&5325&423649707&4.5.7.11.13.19.67.499&61&6426 \cr
 & &16.27.49.11.43.47&2303&1032& & &16.27.49.11.17.61&2989&792 \cr
\noalign{\hrule}
 & &25.17.31.53.601&63&838& & &25.13.37.131.269&28143&34868 \cr
5308&419663275&4.9.7.419.601&1767&1166&5326&423748975&8.9.23.53.59.379&269&110 \cr
 & &16.27.11.19.31.53&513&88& & &32.3.5.11.23.59.269&759&944 \cr
\noalign{\hrule}
 & &3.7.17.23.29.41.43&269&398& & &3.11.17.73.79.131&205&598 \cr
5309&419803797&4.7.17.41.199.269&2305&2574&5327&423822597&4.5.13.17.23.41.79&1143&200 \cr
 & &16.9.5.11.13.199.461&38805&40568& & &64.9.125.13.127&1651&12000 \cr
\noalign{\hrule}
 & &5.53.59.107.251&1617&1510& & &3.11.19.43.79.199&3965&4592 \cr
5310&419909195&4.3.25.49.11.151.251&6837&562&5328&423853881&32.5.7.13.41.61.79&199&594 \cr
 & &16.9.11.43.53.281&2529&3784& & &128.27.7.11.41.199&369&448 \cr
\noalign{\hrule}
}%
}
$$
\eject
\vglue -23 pt
\noindent\hskip 1 in\hbox to 6.5 in{\ 5329 -- 5364 \hfill\fbd 424599609 -- 432763285\frb}
\vskip -9 pt
$$
\vbox{
\nointerlineskip
\halign{\strut
    \vrule \ \ \hfil \frb #\ 
   &\vrule \hfil \ \ \fbb #\frb\ 
   &\vrule \hfil \ \ \frb #\ \hfil
   &\vrule \hfil \ \ \frb #\ 
   &\vrule \hfil \ \ \frb #\ \ \vrule \hskip 2 pt
   &\vrule \ \ \hfil \frb #\ 
   &\vrule \hfil \ \ \fbb #\frb\ 
   &\vrule \hfil \ \ \frb #\ \hfil
   &\vrule \hfil \ \ \frb #\ 
   &\vrule \hfil \ \ \frb #\ \vrule \cr%
\noalign{\hrule}
 & &3.7.73.173.1601&1853&1780& & &9.25.7.31.67.131&11&56 \cr
5329&424599609&8.5.17.89.109.1601&1557&44&5347&428537025&16.5.49.11.31.131&3287&3132 \cr
 & &64.9.5.11.109.173&545&1056& & &128.27.11.19.29.173&15051&13376 \cr
\noalign{\hrule}
 & &9.5.11.73.79.149&7611&3266& & &9.5.19.31.103.157&40513&40342 \cr
5330&425345085&4.27.23.43.59.71&35&8&5348&428612355&4.11.23.29.31.127.877&28761&55948 \cr
 & &64.5.7.23.59.71&9499&2272& & &32.3.11.71.197.9587&153857&153392 \cr
\noalign{\hrule}
 & &81.25.11.29.659&403&322& & &27.11.17.19.41.109&887&1184 \cr
5331&425697525&4.7.11.13.23.31.659&1215&556&5349&428715639&64.17.37.41.887&95&792 \cr
 & &32.243.5.13.31.139&1807&1488& & &1024.9.5.11.19.37&185&512 \cr
\noalign{\hrule}
 & &27.11.17.37.43.53&125&458& & &7.13.41.137.839&891&890 \cr
5332&425746827&4.3.125.17.43.229&277&148&5350&428852333&4.81.5.7.11.41.89.839&2197&5354 \cr
 & &32.5.37.229.277&1385&3664& & &16.9.5.2197.89.2677&120465&120328 \cr
\noalign{\hrule}
 & &27.125.7.13.19.73&10703&7328& & &5.7.29.173.2447&1325&1122 \cr
5333&425982375&64.49.11.139.229&6375&4846&5351&429680965&4.3.125.11.17.53.173&783&2158 \cr
 & &256.3.125.17.2423&2423&2176& & &16.81.13.29.53.83&6723&5512 \cr
\noalign{\hrule}
 & &7.13.23.47.61.71&1107&4444& & &9.5.13.31.37.641&487&154 \cr
5334&426044801&8.27.11.23.41.101&455&488&5352&430107795&4.5.7.11.13.31.487&1419&1016 \cr
 & &128.9.5.7.13.61.101&505&576& & &64.3.7.121.43.127&15367&9632 \cr
\noalign{\hrule}
 & &9.5.31.47.67.97&4981&5404& & &9.17.47.163.367&1133&334 \cr
5335&426106935&8.7.17.97.193.293&11&108&5353&430172811&4.11.103.167.367&235&132 \cr
 & &64.27.11.193.293&6369&9376& & &32.3.5.121.47.167&835&1936 \cr
\noalign{\hrule}
 & &9.5.7.13.29.37.97&3307&88& & &9.17.19.317.467&5863&3010 \cr
5336&426211695&16.3.11.13.3307&1673&1634&5354&430349373&4.5.7.11.13.17.41.43&373&186 \cr
 & &64.7.11.19.43.239&10277&6688& & &16.3.5.7.31.41.373&8897&14920 \cr
\noalign{\hrule}
 & &7.13.19.23.71.151&10505&7032& & &3.125.19.23.37.71&39&76 \cr
5337&426342007&16.3.5.7.11.191.293&4077&6178&5355&430499625&8.9.25.13.361.71&737&2512 \cr
 & &64.81.151.3089&3089&2592& & &256.11.13.67.157&22451&8576 \cr
\noalign{\hrule}
 & &9.5.49.11.73.241&263&2432& & &9.5.13.17.19.43.53&6503&4118 \cr
5338&426718215&256.19.73.263&95&168&5356&430628445&4.7.17.29.71.929&473&456 \cr
 & &4096.3.5.7.361&361&2048& & &64.3.7.11.19.29.43.71&2233&2272 \cr
\noalign{\hrule}
 & &11.13.23.31.53.79&441&142& & &27.5.19.31.61.89&34349&40144 \cr
5339&426902333&4.9.49.31.71.79&3425&3922&5357&431686935&32.49.13.193.701&305&396 \cr
 & &16.3.25.7.37.53.137&5069&4200& & &256.9.5.7.11.61.193&1351&1408 \cr
\noalign{\hrule}
 & &27.5.11.19.109.139&377&168& & &243.5.13.109.251&1577&1686 \cr
5340&427485465&16.81.7.13.29.139&29&110&5358&432135405&4.729.5.19.83.281&1067&338 \cr
 & &64.5.7.11.13.841&841&2912& & &16.11.169.19.83.97&13871&12616 \cr
\noalign{\hrule}
 & &9.25.49.13.19.157&1319&94& & &27.25.7.13.31.227&3071&3058 \cr
5341&427538475&4.13.19.47.1319&8127&9020&5359&432247725&4.25.7.11.31.37.83.139&1989&26314 \cr
 & &32.27.5.7.11.41.43&1353&688& & &16.9.13.17.37.59.223&8251&8024 \cr
\noalign{\hrule}
 & &49.169.51659&21689&29970& & &25.11.169.71.131&1739&36 \cr
5342&427788179&4.81.5.529.37.41&3809&3850&5360&432263975&8.9.11.13.37.47&655&1084 \cr
 & &16.9.125.7.11.13.23.293&25875&25784& & &64.3.5.131.271&813&32 \cr
\noalign{\hrule}
 & &9.5.7.11.169.17.43&1121&562& & &3.17.37.107.2141&161&1980 \cr
5343&428062635&4.5.7.13.19.59.281&193&102&5361&432287169&8.27.5.7.11.23.37&493&358 \cr
 & &16.3.17.19.193.281&3667&2248& & &32.7.11.17.29.179&5191&1232 \cr
\noalign{\hrule}
 & &3.7.31.59.71.157&1859&2330& & &9.5.11.13.17.59.67&1147&1862 \cr
5344&428145123&4.5.7.11.169.31.233&531&314&5362&432438435&4.3.49.19.31.37.67&649&1118 \cr
 & &16.9.11.59.157.233&233&264& & &16.7.11.13.37.43.59&259&344 \cr
\noalign{\hrule}
 & &243.5.7.11.23.199&3139&1438& & &13.17.31.83.761&365&396 \cr
5345&428201235&4.5.11.43.73.719&509&294&5363&432729713&8.9.5.11.13.17.73.83&539&124 \cr
 & &16.3.49.509.719&3563&5752& & &64.3.49.121.31.73&5929&7008 \cr
\noalign{\hrule}
 & &9.25.11.131.1321&22997&10028& & &5.19.23.37.53.101&301&402 \cr
5346&428301225&8.13.23.29.61.109&857&1650&5364&432763285&4.3.5.7.23.43.53.67&8383&9372 \cr
 & &32.3.25.11.29.857&857&464& & &32.9.7.11.71.83.101&8217&7952 \cr
\noalign{\hrule}
}%
}
$$
\eject
\vglue -23 pt
\noindent\hskip 1 in\hbox to 6.5 in{\ 5365 -- 5400 \hfill\fbd 432777345 -- 439450165\frb}
\vskip -9 pt
$$
\vbox{
\nointerlineskip
\halign{\strut
    \vrule \ \ \hfil \frb #\ 
   &\vrule \hfil \ \ \fbb #\frb\ 
   &\vrule \hfil \ \ \frb #\ \hfil
   &\vrule \hfil \ \ \frb #\ 
   &\vrule \hfil \ \ \frb #\ \ \vrule \hskip 2 pt
   &\vrule \ \ \hfil \frb #\ 
   &\vrule \hfil \ \ \fbb #\frb\ 
   &\vrule \hfil \ \ \frb #\ \hfil
   &\vrule \hfil \ \ \frb #\ 
   &\vrule \hfil \ \ \frb #\ \vrule \cr%
\noalign{\hrule}
 & &3.5.7.11.13.19.37.41&223&184& & &9.5.7.13.37.43.67&1463&1418 \cr
5365&432777345&16.5.7.19.23.41.223&169&36&5383&436514715&4.49.11.13.19.37.709&33&670 \cr
 & &128.9.169.23.223&2899&4416& & &16.3.5.121.67.709&709&968 \cr
\noalign{\hrule}
 & &27.7.11.17.37.331&6047&6200& & &9.3125.11.17.83&17893&16482 \cr
5366&432845721&16.3.25.7.11.31.6047&557&6604&5384&436528125&4.27.29.41.67.617&83&700 \cr
 & &128.25.13.127.557&41275&35648& & &32.25.7.41.67.83&469&656 \cr
\noalign{\hrule}
 & &3.11.31.379.1117&729&388& & &3.5.49.29.103.199&209&806 \cr
5367&433079889&8.2187.97.379&19475&17288&5385&436893555&4.7.11.13.19.31.103&603&736 \cr
 & &128.25.19.41.2161&54025&49856& & &256.9.11.23.31.67&16951&11904 \cr
\noalign{\hrule}
 & &27.5.29.59.1877&1591&286& & &11.31.37.59.587&10535&11184 \cr
5368&433558845&4.3.11.13.37.43.59&1075&1108&5386&436964561&32.3.5.49.31.43.233&225&8 \cr
 & &32.25.13.1849.277&24037&22160& & &512.27.125.7.43&23625&11008 \cr
\noalign{\hrule}
 & &9.25.49.23.29.59&1651&524& & &3.11.67.163.1213&975&238 \cr
5369&433866825&8.3.13.59.127.131&23&154&5387&437156709&4.9.25.7.13.17.163&331&484 \cr
 & &32.7.11.13.23.127&1651&176& & &32.5.7.121.13.331&3641&7280 \cr
\noalign{\hrule}
 & &25.31.1369.409&10773&23452& & &9.25.11.13.107.127&323&2 \cr
5370&433938775&8.81.7.11.13.19.41&185&62&5388&437226075&4.3.11.17.19.127&377&250 \cr
 & &32.27.5.7.11.31.37&297&112& & &16.125.13.17.29&17&1160 \cr
\noalign{\hrule}
 & &25.193.293.307&3509&3816& & &23.41.47.71.139&1525&1386 \cr
5371&434013575&16.9.121.29.53.193&293&1830&5389&437403949&4.9.25.7.11.23.47.61&2363&1846 \cr
 & &64.27.5.11.61.293&671&864& & &16.3.25.7.13.17.71.139&975&952 \cr
\noalign{\hrule}
 & &13.29.31.53.701&12733&24420& & &9.5.121.17.29.163&791&298 \cr
5372&434207111&8.3.5.7.11.17.37.107&689&60&5390&437554755&4.5.7.113.149.163&11271&5566 \cr
 & &64.9.25.11.13.53&99&800& & &16.3.121.13.289.23&299&136 \cr
\noalign{\hrule}
 & &25.11.13.29.59.71&371&396& & &9.11.23.41.43.109&21535&30022 \cr
5373&434294575&8.9.7.121.29.53.71&6295&118&5391&437564259&4.5.17.59.73.883&943&60 \cr
 & &32.3.5.7.59.1259&1259&336& & &32.3.25.23.41.73&73&400 \cr
\noalign{\hrule}
 & &9.7.13.59.89.101&28225&19034& & &3.5.11.13.29.31.227&3723&772 \cr
5374&434357469&4.25.31.307.1129&231&76&5392&437736585&8.9.11.17.73.193&775&962 \cr
 & &32.3.5.7.11.19.1129&5645&3344& & &32.25.13.31.37.73&365&592 \cr
\noalign{\hrule}
 & &5.11.23.47.71.103&1161&1208& & &3.5.11.19.29.61.79&123&428 \cr
5375&434794415&16.27.5.11.43.71.151&59&414&5393&438119385&8.9.11.41.79.107&13987&15164 \cr
 & &64.243.23.59.151&8909&7776& & &64.17.71.197.223&43931&38624 \cr
\noalign{\hrule}
 & &3.11.13.17.19.43.73&5155&994& & &9.7.13.29.59.313&8023&3954 \cr
5376&434961813&4.5.7.17.71.1031&473&558&5394&438609717&4.27.71.113.659&1015&902 \cr
 & &16.9.7.11.31.43.71&497&744& & &16.5.7.11.29.41.659&3295&3608 \cr
\noalign{\hrule}
 & &121.13.43.59.109&1197&1340& & &9.11.23.37.41.127&2567&1624 \cr
5377&434986409&8.9.5.7.11.19.67.109&205&532&5395&438684543&16.3.7.17.29.37.151&115&4 \cr
 & &64.3.25.49.361.41&44403&39200& & &128.5.23.29.151&4379&320 \cr
\noalign{\hrule}
 & &17.37.41.47.359&135&494& & &3.5.49.13.19.41.59&291&242 \cr
5378&435137797&4.27.5.13.19.41.47&1073&1606&5396&439157355&4.9.5.121.19.59.97&3239&2366 \cr
 & &16.9.5.11.29.37.73&4015&2088& & &16.7.121.169.41.79&1027&968 \cr
\noalign{\hrule}
 & &9.29.43.71.547&5885&5338& & &81.53.101.1013&43703&38350 \cr
5379&435867651&4.5.11.17.71.107.157&129&58&5397&439229709&4.25.11.13.29.59.137&901&606 \cr
 & &16.3.5.29.43.107.157&785&856& & &16.3.5.13.17.29.53.101&493&520 \cr
\noalign{\hrule}
 & &25.7.11.29.73.107&1053&542& & &3.43.59.197.293&187&10 \cr
5380&436049075&4.81.5.13.107.271&4169&646&5398&439314531&4.5.11.17.43.293&367&1098 \cr
 & &16.9.11.17.19.379&6443&1368& & &16.9.11.61.367&12111&488 \cr
\noalign{\hrule}
 & &9.7.11.13.79.613&2675&4514& & &9.5.17.361.37.43&185&202 \cr
5381&436278843&4.3.25.11.37.61.107&613&58&5399&439378515&4.25.361.1369.101&35343&1118 \cr
 & &16.5.29.107.613&145&856& & &16.27.7.11.13.17.43&273&88 \cr
\noalign{\hrule}
 & &5.121.13.19.23.127&47&162& & &5.7.11.571.1999&961&1038 \cr
5382&436499635&4.81.11.13.47.127&523&874&5400&439450165&4.3.5.961.173.571&4109&1254 \cr
 & &16.3.19.23.47.523&523&1128& & &16.9.7.11.19.31.587&5283&4712 \cr
\noalign{\hrule}
}%
}
$$
\eject
\vglue -23 pt
\noindent\hskip 1 in\hbox to 6.5 in{\ 5401 -- 5436 \hfill\fbd 439521225 -- 448489965\frb}
\vskip -9 pt
$$
\vbox{
\nointerlineskip
\halign{\strut
    \vrule \ \ \hfil \frb #\ 
   &\vrule \hfil \ \ \fbb #\frb\ 
   &\vrule \hfil \ \ \frb #\ \hfil
   &\vrule \hfil \ \ \frb #\ 
   &\vrule \hfil \ \ \frb #\ \ \vrule \hskip 2 pt
   &\vrule \ \ \hfil \frb #\ 
   &\vrule \hfil \ \ \fbb #\frb\ 
   &\vrule \hfil \ \ \frb #\ \hfil
   &\vrule \hfil \ \ \frb #\ 
   &\vrule \hfil \ \ \frb #\ \vrule \cr%
\noalign{\hrule}
 & &3.25.11.13.107.383&2627&2352& & &3.7.11.13.31.67.71&575&916 \cr
5401&439521225&32.9.49.37.71.107&2299&1660&5419&442843401&8.25.13.23.67.229&3069&2198 \cr
 & &256.5.49.121.19.83&10241&10624& & &32.9.25.7.11.31.157&471&400 \cr
\noalign{\hrule}
 & &27.5.121.13.19.109&857&1442& & &5.13.19.23.67.233&1733&1296 \cr
5402&439787205&4.3.7.103.109.857&209&100&5420&443430455&32.81.5.67.1733&699&1034 \cr
 & &32.25.7.11.19.857&857&560& & &128.243.11.47.233&2673&3008 \cr
\noalign{\hrule}
 & &3.11.13.23.841.53&22137&22436& & &27.7.43.197.277&7975&496 \cr
5403&439801791&8.9.11.47.71.79.157&115&754&5421&443482263&32.25.7.11.29.31&83&258 \cr
 & &32.5.13.23.29.47.157&785&752& & &128.3.29.43.83&83&1856 \cr
\noalign{\hrule}
 & &3.11.13.43.107.223&861&2038& & &25.11.23.29.41.59&1197&1222 \cr
5404&440163867&4.9.7.41.43.1019&845&1864&5422&443705075&4.9.7.11.13.19.23.29.47&28735&32882 \cr
 & &64.5.169.41.233&3029&6560& & &16.3.5.49.41.401.821&19649&19704 \cr
\noalign{\hrule}
 & &9.5.49.11.41.443&787&542& & &9.11.13.41.47.179&865&686 \cr
5405&440543565&4.3.11.41.271.787&20605&11662&5423&443928771&4.3.5.343.13.41.173&3149&484 \cr
 & &16.5.343.13.17.317&2219&1768& & &32.49.121.47.67&539&1072 \cr
\noalign{\hrule}
 & &11.37.71.79.193&6825&22072& & &3.11.13.37.83.337&113&30 \cr
5406&440592559&16.3.25.7.13.31.89&71&84&5424&443983683&4.9.5.37.113.337&1001&664 \cr
 & &128.9.5.49.71.89&4005&3136& & &64.7.11.13.83.113&113&224 \cr
\noalign{\hrule}
 & &3.5.11.13.37.67.83&41&42& & &5.49.11.13.19.23.29&4623&1438 \cr
5407&441348765&4.9.5.7.11.13.37.41.67&25813&3502&5425&443998555&4.3.529.67.719&95&624 \cr
 & &16.7.17.83.103.311&5287&5768& & &128.9.5.13.19.67&603&64 \cr
\noalign{\hrule}
 & &3.5.13.17.19.43.163&2783&12& & &3.19.37.43.59.83&671&920 \cr
5408&441461865&8.9.121.19.23&85&86&5426&444094239&16.5.11.19.23.59.61&99&1258 \cr
 & &32.5.121.17.23.43&121&368& & &64.9.5.121.17.37&2057&480 \cr
\noalign{\hrule}
 & &81.5.19.103.557&253&262& & &13.841.179.227&25783&14850 \cr
5409&441469845&4.9.11.19.23.131.557&103&4910&5427&444240589&4.27.25.11.19.23.59&87&28 \cr
 & &16.5.103.131.491&491&1048& & &32.81.5.7.11.19.29&4455&2128 \cr
\noalign{\hrule}
 & &25.7.31.97.839&407&432& & &7.169.17.23.961&1793&1080 \cr
5410&441502775&32.27.7.11.31.37.97&5375&1786&5428&444513433&16.27.5.7.11.31.163&1927&460 \cr
 & &128.9.125.19.43.47&10105&10944& & &128.3.25.23.41.47&1927&4800 \cr
\noalign{\hrule}
 & &5.49.11.17.23.419&235639&236574& & &9.5.7.13.23.29.163&773&1228 \cr
5411&441519155&4.27.13.67.337.3517&4813&8330&5429&445212495&8.3.163.307.773&1241&1078 \cr
 & &16.9.5.49.17.67.4813&4813&4824& & &32.49.11.17.73.307&22411&20944 \cr
\noalign{\hrule}
 & &9.7.11.47.71.191&2015&1322& & &27.25.7.121.19.41&377&8 \cr
5412&441695331&4.5.13.31.191.661&235&426&5430&445373775&16.3.5.11.13.19.29&41&206 \cr
 & &16.3.25.13.31.47.71&325&248& & &64.29.41.103&29&3296 \cr
\noalign{\hrule}
 & &3.5.7.59.113.631&28067&18602& & &9.13.61.73.857&749&200 \cr
5413&441722085&4.13.17.71.127.131&27&44&5431&446497857&16.25.7.107.857&803&54 \cr
 & &32.27.11.13.127.131&16637&20592& & &64.27.25.11.73&825&32 \cr
\noalign{\hrule}
 & &27.5.7.13.17.29.73&205&16& & &9.11.17.19.61.229&6475&3956 \cr
5414&442124865&32.25.29.41.73&1859&1134&5432&446686713&8.25.7.17.23.37.43&6047&7638 \cr
 & &128.81.7.11.169&429&64& & &32.3.5.19.67.6047&6047&5360 \cr
\noalign{\hrule}
 & &5.11.13.29.83.257&1953&874& & &3.23.53.71.1721&1677&44 \cr
5415&442298285&4.9.5.7.19.23.29.31&33&62&5433&446852487&8.9.11.13.43.53&2645&2602 \cr
 & &16.27.7.11.23.961&4347&7688& & &32.5.13.529.1301&6505&4784 \cr
\noalign{\hrule}
 & &9.13.1511.2503&89645&87142& & &27.11.23.29.37.61&25&646 \cr
5416&442497861&4.5.11.17.233.17929&18867&938&5434&447109443&4.25.17.19.29.37&579&494 \cr
 & &16.3.7.11.19.67.331&14003&18536& & &16.3.5.13.361.193&12545&2888 \cr
\noalign{\hrule}
 & &3.5.7.529.31.257&2367&2624& & &7.31.83.149.167&6409&7452 \cr
5417&442527015&128.27.5.23.41.263&3839&6944&5435&448167713&8.81.13.17.23.29.31&679&220 \cr
 & &8192.7.11.31.349&3839&4096& & &64.3.5.7.11.13.23.97&11155&13728 \cr
\noalign{\hrule}
 & &3.5.49.11.19.43.67&899&564& & &3.5.7.11.19.107.191&799&806 \cr
5418&442564815&8.9.7.29.31.43.47&2147&814&5436&448489965&4.11.13.17.19.31.47.191&6057&3956 \cr
 & &32.11.19.29.37.113&1073&1808& & &32.9.13.23.43.47.673&139449&139984 \cr
\noalign{\hrule}
}%
}
$$
\eject
\vglue -23 pt
\noindent\hskip 1 in\hbox to 6.5 in{\ 5437 -- 5472 \hfill\fbd 448498053 -- 454269733\frb}
\vskip -9 pt
$$
\vbox{
\nointerlineskip
\halign{\strut
    \vrule \ \ \hfil \frb #\ 
   &\vrule \hfil \ \ \fbb #\frb\ 
   &\vrule \hfil \ \ \frb #\ \hfil
   &\vrule \hfil \ \ \frb #\ 
   &\vrule \hfil \ \ \frb #\ \ \vrule \hskip 2 pt
   &\vrule \ \ \hfil \frb #\ 
   &\vrule \hfil \ \ \fbb #\frb\ 
   &\vrule \hfil \ \ \frb #\ \hfil
   &\vrule \hfil \ \ \frb #\ 
   &\vrule \hfil \ \ \frb #\ \vrule \cr%
\noalign{\hrule}
 & &243.37.83.601&21203&1034& & &11.13.19.961.173&1075&828 \cr
5437&448498053&4.7.11.13.47.233&4975&5976&5455&451709401&8.9.25.23.961.43&919&3886 \cr
 & &64.9.25.83.199&199&800& & &32.3.5.29.67.919&26651&16080 \cr
\noalign{\hrule}
 & &3.49.19.31.71.73&1107&1156& & &3.5.49.11.29.41.47&3&52 \cr
5438&448759689&8.81.289.19.41.71&5621&130&5456&451814055&8.9.13.29.41.47&1717&3080 \cr
 & &32.5.7.11.13.41.73&451&1040& & &128.5.7.11.17.101&101&1088 \cr
\noalign{\hrule}
 & &3.25.7.11.23.31.109&769&6& & &3.7.169.23.29.191&165&4 \cr
5439&448815675&4.9.11.23.769&109&98&5457&452131953&8.9.5.11.29.191&845&874 \cr
 & &16.49.109.769&769&56& & &32.25.11.169.19.23&209&400 \cr
\noalign{\hrule}
 & &25.11.71.127.181&5727&7124& & &9.11.13.17.23.29.31&1855&5482 \cr
5440&448821175&8.3.25.13.23.83.137&2667&758&5458&452392083&4.5.7.17.53.2741&1311&1430 \cr
 & &32.9.7.13.127.379&2653&1872& & &16.3.25.11.13.19.23.53&475&424 \cr
\noalign{\hrule}
 & &5.13.29.37.41.157&171&14& & &37.109.151.743&447&296 \cr
5441&448948565&4.9.7.13.19.29.41&785&814&5459&452474369&16.3.1369.109.149&8805&7436 \cr
 & &16.3.5.7.11.19.37.157&209&168& & &128.9.5.11.169.587&58113&54080 \cr
\noalign{\hrule}
 & &11.13.19.37.41.109&63&470& & &81.5.13.31.47.59&2603&2662 \cr
5442&449264101&4.9.5.7.19.47.109&2419&2704&5460&452595195&4.1331.19.31.47.137&25&1482 \cr
 & &128.3.7.169.41.59&1239&832& & &16.3.25.121.13.361&1805&968 \cr
\noalign{\hrule}
 & &19.31.59.67.193&391&198& & &3.5.7.13.17.109.179&253&74 \cr
5443&449365181&4.9.11.17.23.59.67&133&870&5461&452752755&4.5.7.11.13.17.23.37&2151&1216 \cr
 & &16.27.5.7.19.23.29&1015&4968& & &512.9.19.23.239&13623&5888 \cr
\noalign{\hrule}
 & &25.343.13.37.109&213&4246& & &3.25.7.11.13.37.163&97&2022 \cr
5444&449578675&4.3.25.11.71.193&109&84&5462&452777325&4.9.37.97.337&455&418 \cr
 & &32.9.7.11.71.109&781&144& & &16.5.7.11.13.19.337&337&152 \cr
\noalign{\hrule}
 & &3.11.41.43.59.131&893&3430& & &27.5.59.139.409&259&436 \cr
5445&449665491&4.5.343.19.41.47&23&18&5463&452818215&8.9.7.37.109.409&695&286 \cr
 & &16.9.343.19.23.47&16121&10488& & &32.5.7.11.13.37.139&1001&592 \cr
\noalign{\hrule}
 & &9.23.41.197.269&2327&2204& & &27.5.13.289.19.47&1357&88 \cr
5446&449751591&8.3.13.19.29.179.269&935&6046&5464&452925135&16.11.13.19.23.59&153&94 \cr
 & &32.5.11.17.29.3023&33253&39440& & &64.9.11.17.23.47&253&32 \cr
\noalign{\hrule}
 & &81.25.11.17.29.41&91&96& & &7.11.13.23.103.191&1685&684 \cr
5447&450244575&64.243.5.7.13.29.41&13&1202&5465&452931479&8.9.5.19.191.337&601&410 \cr
 & &256.7.169.601&4207&21632& & &32.3.25.19.41.601&19475&28848 \cr
\noalign{\hrule}
 & &9.11.19.443.541&43819&48692& & &9.5.7.17.19.61.73&417&94 \cr
5448&450806103&8.7.29.37.47.1511&1625&114&5466&453070485&4.27.5.47.61.139&5291&9044 \cr
 & &32.3.125.7.13.19.29&1625&3248& & &32.7.11.13.17.19.37&481&176 \cr
\noalign{\hrule}
 & &7.11.17.283.1217&551&21240& & &3.7.11.13.17.83.107&841&1590 \cr
5449&450833999&16.9.5.19.29.59&1407&1348&5467&453383931&4.9.5.841.53.83&1547&2852 \cr
 & &128.27.7.67.337&9099&4288& & &32.7.13.17.23.29.31&713&464 \cr
\noalign{\hrule}
 & &3.5.7.11.13.59.509&1019&928& & &25.19.37.83.311&11&1566 \cr
5450&450915465&64.5.29.509.1019&763&1782&5468&453663475&4.27.5.11.29.37&5909&5894 \cr
 & &256.81.7.11.29.109&3161&3456& & &16.9.7.19.311.421&421&504 \cr
\noalign{\hrule}
 & &3.25.7.2197.17.23&1731&2156& & &3.11.23.31.101.191&4313&1990 \cr
5451&450989175&8.9.343.11.13.577&367&4826&5469&453897939&4.5.19.31.199.227&191&36 \cr
 & &32.11.19.127.367&6973&22352& & &32.9.19.191.199&199&912 \cr
\noalign{\hrule}
 & &9.29.41.149.283&3641&7962& & &3.5.11.13.17.59.211&47&8 \cr
5452&451229067&4.27.11.331.1327&11767&2830&5470&453952785&16.17.47.59.211&5825&6624 \cr
 & &16.5.7.1681.283&205&56& & &1024.9.25.23.233&5359&7680 \cr
\noalign{\hrule}
 & &5.13.47.71.2081&4897&5508& & &9.25.11.223.823&347&322 \cr
5453&451379305&8.81.17.59.71.83&517&1724&5471&454234275&4.3.7.11.23.347.823&9589&17570 \cr
 & &64.3.11.47.59.431&4741&5664& & &16.5.49.43.223.251&2107&2008 \cr
\noalign{\hrule}
 & &81.7.11.19.37.103&395&802& & &17.59.313.1447&1225&222 \cr
5454&451614933&4.9.5.79.103.401&407&304&5472&454269733&4.3.25.49.37.313&1003&1188 \cr
 & &128.5.11.19.37.401&401&320& & &32.81.5.7.11.17.59&567&880 \cr
\noalign{\hrule}
}%
}
$$
\eject
\vglue -23 pt
\noindent\hskip 1 in\hbox to 6.5 in{\ 5473 -- 5508 \hfill\fbd 454270685 -- 462257361\frb}
\vskip -9 pt
$$
\vbox{
\nointerlineskip
\halign{\strut
    \vrule \ \ \hfil \frb #\ 
   &\vrule \hfil \ \ \fbb #\frb\ 
   &\vrule \hfil \ \ \frb #\ \hfil
   &\vrule \hfil \ \ \frb #\ 
   &\vrule \hfil \ \ \frb #\ \ \vrule \hskip 2 pt
   &\vrule \ \ \hfil \frb #\ 
   &\vrule \hfil \ \ \fbb #\frb\ 
   &\vrule \hfil \ \ \frb #\ \hfil
   &\vrule \hfil \ \ \frb #\ 
   &\vrule \hfil \ \ \frb #\ \vrule \cr%
\noalign{\hrule}
 & &5.11.17.53.89.103&5191&474& & &17.113.179.1327&2185&858 \cr
5473&454270685&4.3.17.29.79.179&3267&1924&5491&456300893&4.3.5.11.13.19.23.113&1827&358 \cr
 & &32.81.121.13.37&1053&6512& & &16.27.7.11.29.179&297&1624 \cr
\noalign{\hrule}
 & &5.11.19.29.53.283&663&344& & &27.625.11.23.107&623&2 \cr
5474&454544695&16.3.5.13.17.43.283&1073&342&5492&456823125&4.7.11.89.107&585&592 \cr
 & &64.27.13.19.29.37&999&416& & &128.9.5.13.37.89&3293&832 \cr
\noalign{\hrule}
 & &27.5.49.13.17.311&97&214& & &7.23.37.41.1871&18381&16510 \cr
5475&454655565&4.3.5.49.17.97.107&311&1144&5493&456967427&4.3.5.7.11.13.127.557&779&222 \cr
 & &64.11.13.107.311&107&352& & &16.9.5.19.37.41.127&855&1016 \cr
\noalign{\hrule}
 & &9.49.11.13.7211&10547&11086& & &9.5.7.59.67.367&70829&67526 \cr
5476&454747293&4.3.13.23.53.199.241&202895&208438&5494&456986565&4.11.19.47.137.1777&413&2190 \cr
 & &16.5.7.11.17.31.89.1171&46903&46840& & &16.3.5.7.11.47.59.73&517&584 \cr
\noalign{\hrule}
 & &5.49.11.19.83.107&101&6& & &81.49.13.83.107&46321&41078 \cr
5477&454751605&4.3.49.11.83.101&311&228&5495&458232957&4.11.19.23.47.4211&1565&2646 \cr
 & &32.9.19.101.311&909&4976& & &16.27.5.49.11.19.313&1565&1672 \cr
\noalign{\hrule}
 & &3.5.49.11.101.557&589&5538& & &25.43.53.83.97&111&154 \cr
5478&454837845&4.9.5.13.19.31.71&101&146&5496&458705725&4.3.5.7.11.37.83.97&5777&12168 \cr
 & &16.31.71.73.101&2201&584& & &64.27.169.53.109&4563&3488 \cr
\noalign{\hrule}
 & &11.13.17.23.79.103&21375&1388& & &9.25.11.13.53.269&779&196 \cr
5479&454964081&8.9.125.19.347&31&316&5497&458718975&8.3.49.19.41.269&65&334 \cr
 & &64.3.25.31.79&31&2400& & &32.5.7.13.41.167&287&2672 \cr
\noalign{\hrule}
 & &343.11.83.1453&10025&5958& & &3.19.41.349.563&6873&7436 \cr
5480&455020027&4.9.25.7.331.401&1079&3886&5498&459190119&8.9.11.169.19.29.79&1&170 \cr
 & &16.3.5.13.29.67.83&2613&1160& & &32.5.11.17.29.79&1595&21488 \cr
\noalign{\hrule}
 & &9.25.7.11.13.43.47&421&34& & &9.5.7.31.103.457&6413&22378 \cr
5481&455179725&4.5.11.17.47.421&301&216&5499&459648315&4.121.53.67.167&8479&372 \cr
 & &64.27.7.43.421&421&96& & &32.3.31.61.139&61&2224 \cr
\noalign{\hrule}
 & &5.7.11.31.37.1031&25041&13106& & &3.5.121.227.1117&427&2924 \cr
5482&455284445&4.3.17.491.6553&897&7450&5500&460209585&8.5.7.11.17.43.61&2043&1312 \cr
 & &16.9.25.13.23.149&9685&1656& & &512.9.7.41.227&861&256 \cr
\noalign{\hrule}
 & &3.13.193.241.251&5015&4774& & &27.125.7.17.31.37&473&598 \cr
5483&455315757&4.5.7.11.17.31.59.193&3263&18&5501&460663875&4.3.11.13.23.31.37.43&125&1022 \cr
 & &16.9.7.13.31.251&217&24& & &16.125.7.11.43.73&803&344 \cr
\noalign{\hrule}
 & &7.17.31.311.397&1287&890& & &3.5.13.529.41.109&5673&1204 \cr
5484&455469763&4.9.5.11.13.17.31.89&397&130&5502&460999695&8.9.5.7.31.43.61&1&44 \cr
 & &16.3.25.11.169.397&825&1352& & &64.7.11.31.61&217&21472 \cr
\noalign{\hrule}
 & &9.25.361.71.79&557&628& & &25.23.37.109.199&173&372 \cr
5485&455591025&8.3.5.361.157.557&67&1738&5503&461476025&8.3.5.23.31.37.173&7&858 \cr
 & &32.11.67.79.157&737&2512& & &32.9.7.11.13.31&21483&208 \cr
\noalign{\hrule}
 & &27.625.49.19.29&611&814& & &3.7.73.227.1327&667&660 \cr
5486&455608125&4.9.25.7.11.13.37.47&11107&9368&5504&461784057&8.9.5.11.23.29.73.227&1861&182 \cr
 & &64.11.29.383.1171&12881&12256& & &32.5.7.11.13.29.1861&24193&25520 \cr
\noalign{\hrule}
 & &3.7.11.17.19.41.149&1399&1432& & &9.25.7.11.53.503&24713&25084 \cr
5487&455810817&16.7.17.41.179.1399&18029&39330&5505&461867175&8.25.13.1901.6271&1113&788 \cr
 & &64.9.5.121.19.23.149&345&352& & &64.3.7.53.197.6271&6271&6304 \cr
\noalign{\hrule}
 & &3.5.11.17.19.43.199&261&62& & &9.7.23.197.1619&7975&6596 \cr
5488&456045315&4.27.5.11.29.31.43&169&304&5506&462148407&8.25.11.17.23.29.97&1619&1716 \cr
 & &128.169.19.29.31&5239&1856& & &64.3.5.121.13.17.1619&2057&2080 \cr
\noalign{\hrule}
 & &5.49.19.29.31.109&927&1144& & &9.5.7.13.157.719&3421&2702 \cr
5489&456148105&16.9.5.7.11.13.29.103&2071&916&5507&462255885&4.3.5.49.11.193.311&31&1586 \cr
 & &128.3.13.19.109.229&687&832& & &16.13.31.61.193&5983&488 \cr
\noalign{\hrule}
 & &9.625.11.73.101&1457&952& & &81.29.47.53.79&12355&10064 \cr
5490&456204375&16.3.125.7.17.31.47&541&916&5508&462257361&32.9.5.7.17.37.353&403&2068 \cr
 & &128.7.17.229.541&27251&34624& & &256.11.13.17.31.47&2431&3968 \cr
\noalign{\hrule}
}%
}
$$
\eject
\vglue -23 pt
\noindent\hskip 1 in\hbox to 6.5 in{\ 5509 -- 5544 \hfill\fbd 462742401 -- 469973075\frb}
\vskip -9 pt
$$
\vbox{
\nointerlineskip
\halign{\strut
    \vrule \ \ \hfil \frb #\ 
   &\vrule \hfil \ \ \fbb #\frb\ 
   &\vrule \hfil \ \ \frb #\ \hfil
   &\vrule \hfil \ \ \frb #\ 
   &\vrule \hfil \ \ \frb #\ \ \vrule \hskip 2 pt
   &\vrule \ \ \hfil \frb #\ 
   &\vrule \hfil \ \ \fbb #\frb\ 
   &\vrule \hfil \ \ \frb #\ \hfil
   &\vrule \hfil \ \ \frb #\ 
   &\vrule \hfil \ \ \frb #\ \vrule \cr%
\noalign{\hrule}
 & &3.11.47.61.67.73&3879&1012& & &5.13.19.37.61.167&5445&4742 \cr
5509&462742401&8.27.121.23.431&4427&7210&5527&465494965&4.9.25.121.13.2371&1769&5344 \cr
 & &32.5.7.19.103.233&23999&10640& & &256.3.11.29.61.167&319&384 \cr
\noalign{\hrule}
 & &9.7.11.13.83.619&1385&472& & &9.7.29.37.71.97&1925&1664 \cr
5510&462855393&16.3.5.7.13.59.277&331&436&5528&465554313&256.25.49.11.13.71&87&158 \cr
 & &128.109.277.331&30193&21184& & &1024.3.5.11.13.29.79&5135&5632 \cr
\noalign{\hrule}
 & &3.5.361.233.367&67&300& & &3.49.121.17.23.67&2041&500 \cr
5511&463042065&8.9.125.361.67&2563&5812&5529&465966039&8.125.7.13.17.157&23&198 \cr
 & &64.11.233.1453&1453&352& & &32.9.5.11.23.157&785&48 \cr
\noalign{\hrule}
 & &9.13.67.113.523&247&770& & &5.7.13.19.31.37.47&955&2412 \cr
5512&463277061&4.5.7.11.169.19.67&759&424&5530&466043305&8.9.25.19.67.191&6281&4606 \cr
 & &64.3.121.19.23.53&6413&13984& & &32.3.49.11.47.571&1713&1232 \cr
\noalign{\hrule}
 & &9.13.17.109.2137&24541&11788& & &5.11.19.59.67.113&903&3050 \cr
5513&463303737&8.7.11.23.97.421&905&1326&5531&466790005&4.3.125.7.11.43.61&201&674 \cr
 & &32.3.5.7.11.13.17.181&905&1232& & &16.9.61.67.337&549&2696 \cr
\noalign{\hrule}
 & &9.5.49.13.19.23.37&2987&8668& & &3.625.49.5081&703&4378 \cr
5514&463484385&8.7.11.29.103.197&2183&804&5532&466816875&4.25.11.19.37.199&1379&6354 \cr
 & &64.3.11.37.59.67&649&2144& & &16.9.7.197.353&353&4728 \cr
\noalign{\hrule}
 & &3.5.7.13.361.941&2533&2172& & &27.5.13.23.43.269&2297&1738 \cr
5515&463691865&8.9.7.13.17.149.181&361&1628&5533&466901955&4.9.11.23.79.2297&331&538 \cr
 & &64.11.361.37.149&1639&1184& & &16.269.331.2297&2297&2648 \cr
\noalign{\hrule}
 & &3.11.31.43.53.199&3281&21788& & &3.11.13.31.41.857&133&10 \cr
5516&463951983&8.13.17.193.419&3465&3658&5534&467286963&4.5.7.19.31.857&351&506 \cr
 & &32.9.5.7.11.13.31.59&1239&1040& & &16.27.7.11.13.19.23&207&1064 \cr
\noalign{\hrule}
 & &3.5.29.59.101.179&65&36& & &625.11.101.673&243&868 \cr
5517&463997535&8.27.25.13.59.179&1441&3034&5535&467314375&8.243.7.31.673&1313&640 \cr
 & &32.11.13.37.41.131&16687&27248& & &2048.27.5.13.101&351&1024 \cr
\noalign{\hrule}
 & &9.11.19.37.59.113&157&860& & &3.5.19.101.109.149&1183&736 \cr
5518&464003199&8.5.11.43.59.157&405&2132&5536&467497185&64.5.7.169.23.109&1111&306 \cr
 & &64.81.25.13.41&533&7200& & &256.9.11.13.17.101&663&1408 \cr
\noalign{\hrule}
 & &9.11.17.29.37.257&3601&21842& & &3.5.7.17.37.73.97&401&110 \cr
5519&464105763&4.13.67.163.277&145&132&5537&467664645&4.25.11.17.37.401&1679&8496 \cr
 & &32.3.5.11.29.67.163&815&1072& & &128.9.23.59.73&177&1472 \cr
\noalign{\hrule}
 & &27.25.23.37.809&21791&3184& & &27.5.7.11.19.23.103&799&664 \cr
5520&464709825&32.7.11.199.283&145&138&5538&467889345&16.17.23.47.83.103&3135&766 \cr
 & &128.3.5.11.23.29.199&2189&1856& & &64.3.5.11.17.19.383&383&544 \cr
\noalign{\hrule}
 & &9.5.23.37.61.199&1577&5786& & &3.7.11.17.97.1229&178165&179474 \cr
5521&464863005&4.3.5.11.19.83.263&3023&1708&5539&468149451&4.5.13.19.2741.4723&179537&127458 \cr
 & &32.7.11.61.3023&3023&1232& & &16.9.17.59.73.97.179&4307&4296 \cr
\noalign{\hrule}
 & &11.13.17.961.199&893&1296& & &3.5.7.13.23.43.347&467&1962 \cr
5522&464902009&32.81.17.19.31.47&725&2182&5540&468444795&4.27.43.109.467&347&814 \cr
 & &128.9.25.29.1091&27275&16704& & &16.11.37.109.347&1199&296 \cr
\noalign{\hrule}
 & &9.125.49.11.13.59&5321&3116& & &9.5.11.53.61.293&1139&874 \cr
5523&465089625&8.25.17.19.41.313&177&602&5541&468898155&4.3.17.19.23.67.293&1855&1562 \cr
 & &32.3.7.43.59.313&313&688& & &16.5.7.11.19.23.53.71&1349&1288 \cr
\noalign{\hrule}
 & &5.7.41.179.1811&31603&31782& & &27.5.17.89.2297&10199&12496 \cr
5524&465182515&4.3.11.169.17.41.5297&25669&32598&5542&469173735&32.9.7.11.31.47.71&1513&874 \cr
 & &16.27.7.17.19.193.1811&3667&3672& & &128.17.19.23.47.89&1081&1216 \cr
\noalign{\hrule}
 & &13.19.31.89.683&3487&4170& & &7.121.13.43.991&6403&6480 \cr
5525&465446059&4.3.5.11.89.139.317&203&114&5543&469211743&32.81.5.11.19.43.337&9043&56 \cr
 & &16.9.5.7.11.19.29.139&11165&10008& & &512.3.5.7.9043&9043&3840 \cr
\noalign{\hrule}
 & &9.25.11.13.17.23.37&791&934& & &25.121.13.17.19.37&207&218 \cr
5526&465475725&4.3.7.17.37.113.467&10241&2302&5544&469973075&4.9.11.13.19.23.37.109&2699&590 \cr
 & &16.343.11.19.1151&6517&9208& & &16.3.5.59.109.2699&19293&21592 \cr
\noalign{\hrule}
}%
}
$$
\eject
\vglue -23 pt
\noindent\hskip 1 in\hbox to 6.5 in{\ 5545 -- 5580 \hfill\fbd 470108925 -- 476160009\frb}
\vskip -9 pt
$$
\vbox{
\nointerlineskip
\halign{\strut
    \vrule \ \ \hfil \frb #\ 
   &\vrule \hfil \ \ \fbb #\frb\ 
   &\vrule \hfil \ \ \frb #\ \hfil
   &\vrule \hfil \ \ \frb #\ 
   &\vrule \hfil \ \ \frb #\ \ \vrule \hskip 2 pt
   &\vrule \ \ \hfil \frb #\ 
   &\vrule \hfil \ \ \fbb #\frb\ 
   &\vrule \hfil \ \ \frb #\ \hfil
   &\vrule \hfil \ \ \frb #\ 
   &\vrule \hfil \ \ \frb #\ \vrule \cr%
\noalign{\hrule}
 & &9.25.11.13.19.769&9983&9242& & &9.5.11.71.97.139&34751&41006 \cr
5545&470108925&4.3.11.67.149.4621&47975&61838&5563&473860035&4.7.19.29.31.59.101&629&78 \cr
 & &16.25.49.19.101.631&4949&5048& & &16.3.13.17.31.37.59&14911&8024 \cr
\noalign{\hrule}
 & &3.1331.23.47.109&3953&1170& & &27.11.23.173.401&235&62 \cr
5546&470491197&4.27.5.11.13.59.67&301&436&5564&473886963&4.5.23.31.47.401&741&340 \cr
 & &32.7.13.43.59.109&2537&1456& & &32.3.25.13.17.19.31&13175&3952 \cr
\noalign{\hrule}
 & &27.13.31.83.521&7049&7018& & &9.13.61.127.523&353&170 \cr
5547&470527083&4.7.121.13.19.29.53.83&4871&15630&5565&474046677&4.3.5.13.17.127.353&3247&1342 \cr
 & &16.3.5.121.521.4871&4871&4840& & &16.11.289.61.191&3179&1528 \cr
\noalign{\hrule}
 & &9.5.31.43.47.167&209&256& & &17.37.691.1091&22057&3510 \cr
5548&470822265&512.3.11.19.43.167&877&2296&5566&474191149&4.27.5.7.13.23.137&629&330 \cr
 & &8192.7.41.877&35957&28672& & &16.81.25.11.17.37&891&200 \cr
\noalign{\hrule}
 & &27.49.11.13.19.131&575&62& & &9.121.13.19.41.43&305&58 \cr
5549&470891421&4.25.11.23.31.131&567&698&5567&474217029&4.3.5.29.41.43.61&85&44 \cr
 & &16.81.5.7.31.349&1745&744& & &32.25.11.17.29.61&1769&6800 \cr
\noalign{\hrule}
 & &243.5.29.43.311&6565&6808& & &7.13.17.29.97.109&2211&950 \cr
5550&471197655&16.25.13.23.29.37.101&1687&28512&5568&474336499&4.3.25.7.11.17.19.67&109&24 \cr
 & &1024.81.7.11.241&2651&3584& & &64.9.5.11.67.109&495&2144 \cr
\noalign{\hrule}
 & &25.289.37.41.43&19879&19896& & &243.125.7.23.97&1203&1028 \cr
5551&471293975&16.3.17.41.103.193.829&763&66&5569&474366375&8.729.5.257.401&1367&638 \cr
 & &64.9.7.11.103.109.193&218669&219744& & &32.11.29.257.1367&39643&45232 \cr
\noalign{\hrule}
 & &3.7.11.19.163.659&1027&950& & &3.961.1849.89&517&2366 \cr
5552&471453213&4.25.13.361.79.163&18171&5294&5570&474429363&4.7.11.169.47.89&837&320 \cr
 & &16.27.5.673.2647&23823&26920& & &512.27.5.7.13.31&455&2304 \cr
\noalign{\hrule}
 & &25.7.29.61.1523&6215&4446& & &25.13.289.31.163&671&144 \cr
5553&471482725&4.9.125.11.13.19.113&107&232&5571&474603025&32.9.5.11.13.17.61&2263&1702 \cr
 & &64.3.11.13.19.29.107&6099&4576& & &128.3.23.31.37.73&2553&4672 \cr
\noalign{\hrule}
 & &27.5.53.233.283&3287&3004& & &9.121.271.1609&1349&260 \cr
5554&471793545&8.5.19.53.173.751&2893&2142&5572&474846471&8.5.13.19.71.271&183&88 \cr
 & &32.9.7.11.17.173.263&31297&30448& & &128.3.11.13.61.71&923&3904 \cr
\noalign{\hrule}
 & &9.11.67.83.857&1495&4066& & &9.5.169.197.317&583&3538 \cr
5555&471811923&4.3.5.11.13.19.23.107&133&166&5573&474924645&4.3.11.13.29.53.61&317&476 \cr
 & &16.5.7.361.83.107&3745&2888& & &32.7.11.17.29.317&1309&464 \cr
\noalign{\hrule}
 & &3.25.7.13.257.269&83&174& & &5.49.121.17.23.41&481&366 \cr
5556&471832725&4.9.25.29.83.269&247&22&5574&475238995&4.3.7.13.17.37.41.61&477&2024 \cr
 & &16.11.13.19.29.83&6061&664& & &64.27.11.23.37.53&999&1696 \cr
\noalign{\hrule}
 & &9.11.17.31.83.109&24613&32830& & &27.25.11.169.379&407&3818 \cr
5557&472009131&4.5.49.67.151.163&85&78&5575&475578675&4.3.121.23.37.83&1267&1516 \cr
 & &16.3.25.7.13.17.67.151&13741&13400& & &32.7.37.181.379&1267&592 \cr
\noalign{\hrule}
 & &81.11.169.43.73&36155&28888& & &27.5.11.17.83.227&10823&8018 \cr
5558&472667481&16.5.7.23.157.1033&1289&2322&5576&475641045&4.9.19.79.137.211&17&154 \cr
 & &64.27.5.7.43.1289&1289&1120& & &16.7.11.17.79.211&553&1688 \cr
\noalign{\hrule}
 & &81.13.23.29.673&2351&21868& & &3.25.13.19.61.421&607&186 \cr
5559&472682223&8.7.11.71.2351&785&1566&5577&475740525&4.9.25.19.31.607&4631&10094 \cr
 & &32.27.5.7.29.157&785&112& & &16.49.11.103.421&1133&392 \cr
\noalign{\hrule}
 & &9.7.13.23.41.613&2087&2204& & &9.5.7.11.17.41.197&1975&2162 \cr
5560&473430321&8.19.23.29.41.2087&1515&572&5578&475775685&4.3.125.23.41.47.79&3349&5276 \cr
 & &64.3.5.11.13.19.29.101&9595&10208& & &32.17.79.197.1319&1319&1264 \cr
\noalign{\hrule}
 & &7.11.17.19.79.241&575&2076& & &5.49.961.43.47&26419&14904 \cr
5561&473518969&8.3.25.7.17.23.173&1801&936&5579&475834345&16.81.23.29.911&4433&3766 \cr
 & &128.27.5.13.1801&23413&8640& & &64.9.7.11.13.31.269&3497&3168 \cr
\noalign{\hrule}
 & &11.19.23.29.43.79&7&36& & &3.13.71.359.479&2575&3652 \cr
5562&473551991&8.9.7.11.19.23.79&221&1090&5580&476160009&8.25.11.71.83.103&6227&2322 \cr
 & &32.3.5.7.13.17.109&1853&21840& & &32.27.5.13.43.479&215&144 \cr
\noalign{\hrule}
}%
}
$$
\eject
\vglue -23 pt
\noindent\hskip 1 in\hbox to 6.5 in{\ 5581 -- 5616 \hfill\fbd 476692425 -- 485891721\frb}
\vskip -9 pt
$$
\vbox{
\nointerlineskip
\halign{\strut
    \vrule \ \ \hfil \frb #\ 
   &\vrule \hfil \ \ \fbb #\frb\ 
   &\vrule \hfil \ \ \frb #\ \hfil
   &\vrule \hfil \ \ \frb #\ 
   &\vrule \hfil \ \ \frb #\ \ \vrule \hskip 2 pt
   &\vrule \ \ \hfil \frb #\ 
   &\vrule \hfil \ \ \fbb #\frb\ 
   &\vrule \hfil \ \ \frb #\ \hfil
   &\vrule \hfil \ \ \frb #\ 
   &\vrule \hfil \ \ \frb #\ \vrule \cr%
\noalign{\hrule}
 & &27.25.11.19.31.109&131&644& & &17.37.43.71.251&625&582 \cr
5581&476692425&8.7.11.23.109.131&775&666&5599&482004587&4.3.625.37.97.251&1419&10706 \cr
 & &32.9.25.7.23.31.37&259&368& & &16.9.5.11.43.53.101&5247&4040 \cr
\noalign{\hrule}
 & &27.11.13.19.67.97&21895&1622& & &9.25.49.107.409&341&194 \cr
5582&476760141&4.5.29.151.811&481&330&5600&482487075&4.3.5.11.31.97.409&371&856 \cr
 & &16.3.25.11.13.29.37&1073&200& & &64.7.11.31.53.107&583&992 \cr
\noalign{\hrule}
 & &27.11.13.23.41.131&703&1000& & &25.11.13.337.401&27&28 \cr
5583&476960913&16.125.19.23.37.41&501&524&5601&483114775&8.27.5.7.13.337.401&817&5198 \cr
 & &128.3.5.19.37.131.167&6179&6080& & &32.9.7.19.23.43.113&43731&48944 \cr
\noalign{\hrule}
 & &5.7.11.13.19.29.173&1075&828& & &3.13.19.31.109.193&169807&171224 \cr
5584&477091615&8.9.125.7.23.29.43&2147&2522&5602&483240927&16.11.17.43.359.1259&60635&6498 \cr
 & &32.3.13.19.43.97.113&4859&4656& & &64.9.5.361.67.181&10317&10720 \cr
\noalign{\hrule}
 & &19.29.31.73.383&1323&940& & &121.23.29.53.113&1335&1448 \cr
5585&477567679&8.27.5.49.19.29.47&111&440&5603&483354223&16.3.5.29.53.89.181&15481&10764 \cr
 & &128.81.25.7.11.37&14175&26048& & &128.27.13.23.113.137&1781&1728 \cr
\noalign{\hrule}
 & &9.11.37.73.1787&565&1222& & &27.5.7.13.23.29.59&561&106 \cr
5586&477842013&4.5.11.13.37.47.113&283&1752&5604&483451605&4.81.11.17.53.59&2009&1118 \cr
 & &64.3.47.73.283&283&1504& & &16.49.13.17.41.43&697&2408 \cr
\noalign{\hrule}
 & &3.7.11.17.73.1667&22601&5738& & &25.11.169.101.103&19&84 \cr
5587&477880557&4.19.97.151.233&18513&16670&5605&483479425&8.3.5.7.11.13.19.101&783&328 \cr
 & &16.9.5.121.17.1667&33&40& & &128.81.19.29.41&22591&5184 \cr
\noalign{\hrule}
 & &9.5.11.17.29.37.53&229&178& & &11.17.23.107.1051&7011&4550 \cr
5588&478552635&4.3.5.29.53.89.229&949&2486&5606&483677557&4.9.25.7.13.17.19.41&253&32 \cr
 & &16.11.13.73.89.113&10057&7592& & &256.3.5.7.11.23.41&1435&384 \cr
\noalign{\hrule}
 & &243.5.49.11.17.43&487&244& & &3.7.11.89.101.233&8339&650 \cr
5589&478720935&8.5.49.11.61.487&117&422&5607&483814947&4.25.7.13.31.269&243&212 \cr
 & &32.9.13.211.487&2743&7792& & &32.243.5.53.269&14257&6480 \cr
\noalign{\hrule}
 & &27.25.7.11.61.151&493&178& & &31.58081.269&24871&33210 \cr
5590&478741725&4.3.5.17.29.89.151&511&244&5608&484337459&4.81.5.7.11.17.19.41&241&538 \cr
 & &32.7.17.29.61.73&493&1168& & &16.3.5.7.17.241.269&105&136 \cr
\noalign{\hrule}
 & &9.43.709.1747&2975&2266& & &7.11.29.41.67.79&6305&6774 \cr
5591&479347101&4.3.25.7.11.17.43.103&709&194&5609&484590029&4.3.5.13.79.97.1129&49245&50374 \cr
 & &16.5.11.17.97.709&1067&680& & &16.9.25.49.67.89.283&17829&17800 \cr
\noalign{\hrule}
 & &5.7.121.19.59.101&3007&2952& & &3.289.191.2927&1177&1750 \cr
5592&479490935&16.9.7.11.19.31.41.97&14809&6240&5610&484702419&4.125.7.11.289.107&2927&252 \cr
 & &1024.27.5.13.59.251&6777&6656& & &32.9.5.49.2927&49&240 \cr
\noalign{\hrule}
 & &7.11.83.193.389&34343&40734& & &25.11.67.83.317&1037&4524 \cr
5593&479817107&4.9.31.61.73.563&83&10&5611&484780175&8.3.25.13.17.29.61&27&2 \cr
 & &16.3.5.61.83.563&2815&1464& & &32.81.13.17.61&793&22032 \cr
\noalign{\hrule}
 & &27.5.19.23.79.103&7171&6734& & &3.13.19.23.149.191&1023&1460 \cr
5594&480042315&4.7.13.37.71.79.101&759&7220&5612&485026737&8.9.5.11.31.73.149&667&3952 \cr
 & &32.3.5.11.361.23.37&407&304& & &256.11.13.19.23.29&319&128 \cr
\noalign{\hrule}
 & &5.11.97.127.709&179&306& & &9.7.11.23.83.367&841&260 \cr
5595&480379405&4.9.11.17.179.709&635&74&5613&485517879&8.3.5.11.13.23.841&191&568 \cr
 & &16.3.5.37.127.179&179&888& & &128.5.29.71.191&5539&22720 \cr
\noalign{\hrule}
 & &5.11.23.41.73.127&801&596& & &3.49.19.31.71.79&9955&5452 \cr
5596&480840415&8.9.23.73.89.149&8509&2368&5614&485644047&8.5.7.11.29.47.181&1207&1026 \cr
 & &1024.3.37.67.127&2479&1536& & &32.27.5.17.19.47.71&799&720 \cr
\noalign{\hrule}
 & &5.121.17.101.463&69&52& & &25.49.11.13.47.59&2193&1426 \cr
5597&480957455&8.3.5.13.23.101.463&21&484&5615&485760275&4.3.25.7.17.23.31.43&293&132 \cr
 & &64.9.7.121.13.23&299&2016& & &32.9.11.31.43.293&9083&6192 \cr
\noalign{\hrule}
 & &9.25.11.19.37.277&1139&2524& & &9.7.23.97.3457&2165&1292 \cr
5598&481959225&8.5.17.19.67.631&1739&1416&5616&485891721&8.5.7.17.19.23.433&297&136 \cr
 & &128.3.37.47.59.67&3149&3776& & &128.27.5.11.289.19&9537&6080 \cr
\noalign{\hrule}
}%
}
$$
\eject
\vglue -23 pt
\noindent\hskip 1 in\hbox to 6.5 in{\ 5617 -- 5652 \hfill\fbd 485895141 -- 493403001\frb}
\vskip -9 pt
$$
\vbox{
\nointerlineskip
\halign{\strut
    \vrule \ \ \hfil \frb #\ 
   &\vrule \hfil \ \ \fbb #\frb\ 
   &\vrule \hfil \ \ \frb #\ \hfil
   &\vrule \hfil \ \ \frb #\ 
   &\vrule \hfil \ \ \frb #\ \ \vrule \hskip 2 pt
   &\vrule \ \ \hfil \frb #\ 
   &\vrule \hfil \ \ \fbb #\frb\ 
   &\vrule \hfil \ \ \frb #\ \hfil
   &\vrule \hfil \ \ \frb #\ 
   &\vrule \hfil \ \ \frb #\ \vrule \cr%
\noalign{\hrule}
 & &9.41.43.113.271&26675&15022& & &27.7.19.23.31.191&325&388 \cr
5617&485895141&4.25.7.11.29.37.97&2911&678&5635&489033153&8.3.25.13.19.97.191&511&2354 \cr
 & &16.3.25.41.71.113&71&200& & &32.5.7.11.13.73.107&7811&11440 \cr
\noalign{\hrule}
 & &5.29.89.101.373&2397&532& & &9.7.11.47.83.181&5725&5678 \cr
5618&486170065&8.3.7.17.19.47.89&4849&5742&5636&489314133&4.25.11.17.83.167.229&39&874 \cr
 & &32.27.11.13.29.373&297&208& & &16.3.5.13.17.19.23.229&28405&31144 \cr
\noalign{\hrule}
 & &5.11.29.41.43.173&25477&25650& & &729.31.59.367&473&22126 \cr
5619&486473405&4.27.125.11.19.73.349&1211&164&5637&489336147&4.11.13.23.37.43&5&6 \cr
 & &32.9.7.19.41.73.173&1197&1168& & &16.3.5.13.23.37.43&12857&1480 \cr
\noalign{\hrule}
 & &5.19.29.41.59.73&151&2268& & &13.17.529.53.79&225&304 \cr
5620&486497185&8.81.5.7.19.151&1199&1066&5638&489497983&32.9.25.13.17.19.53&557&132 \cr
 & &32.27.11.13.41.109&3861&1744& & &256.27.11.19.557&15039&26752 \cr
\noalign{\hrule}
 & &9.5.7.89.97.179&649&604& & &7.17.19.29.31.241&29049&36520 \cr
5621&486771705&8.11.59.89.97.151&8771&138&5639&489865999&16.3.5.11.23.83.421&85&168 \cr
 & &32.3.49.11.23.179&77&368& & &256.9.25.7.17.421&3789&3200 \cr
\noalign{\hrule}
 & &81.19.37.83.103&725&22& & &9.125.11.289.137&371&1136 \cr
5622&486805707&4.9.25.11.29.103&2477&2158&5640&489963375&32.25.7.17.53.71&4887&3562 \cr
 & &16.5.13.83.2477&2477&520& & &128.27.13.137.181&543&832 \cr
\noalign{\hrule}
 & &3.25.7.619.1499&737&762& & &9.5.49.11.17.29.41&163&1846 \cr
5623&487137525&4.9.7.11.67.127.619&5365&206&5641&490266315&4.5.13.29.71.163&12279&11356 \cr
 & &16.5.29.37.103.127&13081&8584& & &32.3.17.167.4093&4093&2672 \cr
\noalign{\hrule}
 & &13.17.23.37.2591&27555&16492& & &3.5.11.43.257.269&1157&1802 \cr
5624&487291961&8.3.5.7.11.19.31.167&557&612&5642&490498635&4.13.17.53.89.257&473&216 \cr
 & &64.27.17.19.31.557&17267&16416& & &64.27.11.17.43.89&801&544 \cr
\noalign{\hrule}
 & &5.19.23.83.2687&1125&1562& & &25.11.13.19.31.233&3449&3774 \cr
5625&487300885&4.9.625.11.71.83&437&188&5643&490622275&4.3.11.17.19.37.3449&3501&52 \cr
 & &32.3.11.19.23.47.71&1551&1136& & &32.27.13.37.389&999&6224 \cr
\noalign{\hrule}
 & &9.11.19.53.67.73&43&24& & &5.11.31.71.4057&60931&64836 \cr
5626&487598463&16.27.11.43.53.73&16549&17980&5644&491120135&8.9.13.43.109.1801&1993&3410 \cr
 & &128.5.13.19.29.31.67&1885&1984& & &32.3.5.11.31.43.1993&1993&2064 \cr
\noalign{\hrule}
 & &3.49.19.41.4261&61&62& & &27.5.11.31.47.227&30067&37352 \cr
5627&487939893&4.49.19.31.61.4261&1665&2596&5645&491147415&16.7.23.29.107.281&39&242 \cr
 & &32.9.5.11.31.37.59.61&104005&104784& & &64.3.121.13.23.107&2461&4576 \cr
\noalign{\hrule}
 & &81.11.19.127.227&1247&1166& & &9.11.61.97.839&9295&8456 \cr
5628&488046141&4.121.29.43.53.227&85&6498&5646&491471937&16.3.5.7.121.169.151&839&734 \cr
 & &16.9.5.17.361.43&323&1720& & &64.13.151.367.839&4771&4832 \cr
\noalign{\hrule}
 & &25.11.23.31.47.53&6931&7644& & &27.17.19.131.431&32591&23870 \cr
5629&488422825&8.3.49.13.29.47.239&155&1518&5647&492396381&4.5.7.11.13.23.31.109&2489&1944 \cr
 & &32.9.5.7.11.13.23.31&91&144& & &64.243.7.19.23.131&207&224 \cr
\noalign{\hrule}
 & &5.11.31.37.61.127&2427&2272& & &9.5.7.19.281.293&407&5746 \cr
5630&488719495&64.3.11.61.71.809&69&740&5648&492763005&4.3.5.11.169.17.37&1301&586 \cr
 & &512.9.5.23.37.71&1633&2304& & &16.13.293.1301&1301&104 \cr
\noalign{\hrule}
 & &3.11.181.223.367&3589&3770& & &5.17.19.37.73.113&1143&778 \cr
5631&488836293&4.5.13.29.37.97.367&2169&8474&5649&492918995&4.9.19.37.127.389&935&232 \cr
 & &16.9.19.37.223.241&2109&1928& & &64.3.5.11.17.29.127&1397&2784 \cr
\noalign{\hrule}
 & &9.5.7.13.19.61.103&527&812& & &3.5.343.11.31.281&57&1462 \cr
5632&488848815&8.3.49.17.29.31.61&845&2266&5650&492999045&4.9.7.11.17.19.43&925&538 \cr
 & &32.5.11.169.31.103&341&208& & &16.25.17.37.269&9953&680 \cr
\noalign{\hrule}
 & &9.5.11.13.17.41.109&931&268& & &3.11.19.29.43.631&2847&4094 \cr
5633&488886255&8.3.5.49.19.41.67&73&542&5651&493359339&4.9.13.19.23.73.89&5663&6820 \cr
 & &32.7.19.73.271&1897&22192& & &32.5.7.11.23.31.809&18607&17360 \cr
\noalign{\hrule}
 & &243.7.11.17.29.53&1301&2780& & &3.49.13.19.107.127&1415&3828 \cr
5634&488899719&8.81.5.139.1301&29&110&5652&493403001&8.9.5.11.13.29.283&1615&1498 \cr
 & &32.25.11.29.1301&1301&400& & &32.25.7.17.19.29.107&425&464 \cr
\noalign{\hrule}
}%
}
$$
\eject
\vglue -23 pt
\noindent\hskip 1 in\hbox to 6.5 in{\ 5653 -- 5688 \hfill\fbd 493420707 -- 503034675\frb}
\vskip -9 pt
$$
\vbox{
\nointerlineskip
\halign{\strut
    \vrule \ \ \hfil \frb #\ 
   &\vrule \hfil \ \ \fbb #\frb\ 
   &\vrule \hfil \ \ \frb #\ \hfil
   &\vrule \hfil \ \ \frb #\ 
   &\vrule \hfil \ \ \frb #\ \ \vrule \hskip 2 pt
   &\vrule \ \ \hfil \frb #\ 
   &\vrule \hfil \ \ \fbb #\frb\ 
   &\vrule \hfil \ \ \frb #\ \hfil
   &\vrule \hfil \ \ \frb #\ 
   &\vrule \hfil \ \ \frb #\ \vrule \cr%
\noalign{\hrule}
 & &27.13.31.137.331&7999&8030& & &27.25.13.29.37.53&253&1178 \cr
5653&493420707&4.3.5.11.19.73.331.421&13&6302&5671&499025475&4.11.13.19.23.29.31&2881&3180 \cr
 & &16.11.13.23.73.137&253&584& & &32.3.5.31.43.53.67&1333&1072 \cr
\noalign{\hrule}
 & &3.5.11.13.31.41.181&1233&782& & &27.25.11.23.37.79&13&14 \cr
5654&493459395&4.27.17.23.137.181&1279&3608&5672&499175325&4.25.7.11.13.23.37.79&6723&16898 \cr
 & &64.11.23.41.1279&1279&736& & &16.81.49.17.71.83&10437&11288 \cr
\noalign{\hrule}
 & &9.11.113.131.337&365&28& & &5.7.11.13.73.1367&645&722 \cr
5655&493872489&8.3.5.7.11.73.113&379&412&5673&499453955&4.3.25.13.361.43.73&4101&19624 \cr
 & &64.5.73.103.379&27667&16480& & &64.9.11.223.1367&223&288 \cr
\noalign{\hrule}
 & &3.121.13.17.47.131&89&42& & &9.25.61.151.241&38437&44462 \cr
5656&493933011&4.9.7.121.13.17.89&655&434&5674&499466475&4.7.11.289.19.43.47&591&302 \cr
 & &16.5.49.31.89.131&2759&1960& & &16.3.7.11.43.151.197&2167&2408 \cr
\noalign{\hrule}
 & &9.7.11.71.89.113&16523&27710& & &5.11.13.47.89.167&2449&1734 \cr
5657&494834571&4.5.13.17.31.41.163&1695&424&5675&499471115&4.3.289.31.79.167&123&44 \cr
 & &64.3.25.17.53.113&1325&544& & &32.9.11.289.31.41&8959&5904 \cr
\noalign{\hrule}
 & &3.25.7.13.59.1229&569&660& & &5.7.43.61.5441&2871&2570 \cr
5658&494887575&8.9.125.11.59.569&1127&6248&5676&499511005&4.9.25.11.29.61.257&6751&8926 \cr
 & &128.49.121.23.71&11431&7744& & &16.3.11.43.157.4463&13389&13816 \cr
\noalign{\hrule}
 & &27.5.7.11.19.23.109&1681&826& & &81.7.11.13.61.101&437&356 \cr
5659&495145035&4.3.49.11.1681.59&299&152&5677&499540041&8.7.11.19.23.89.101&1845&74 \cr
 & &64.13.19.23.41.59&767&1312& & &32.9.5.37.41.89&7585&1424 \cr
\noalign{\hrule}
 & &3.5.7.11.169.43.59&193&102& & &3.5.7.23.29.37.193&41&152 \cr
5660&495209715&4.9.11.13.17.43.193&1121&562&5678&500119935&16.5.7.19.23.29.41&3069&5254 \cr
 & &16.19.59.193.281&3667&2248& & &64.9.11.31.37.71&2343&992 \cr
\noalign{\hrule}
 & &5.7.47.211.1429&187&1242& & &3.11.13.17.19.23.157&2325&2482 \cr
5661&495998755&4.27.7.11.17.23.47&3569&50&5679&500365437&4.9.25.13.289.31.73&437&148 \cr
 & &16.3.25.43.83&43&9960& & &32.5.19.23.31.37.73&2701&2480 \cr
\noalign{\hrule}
 & &27.5.11.19.73.241&34889&41396& & &9.25.13.23.43.173&5575&3326 \cr
5662&496386495&8.79.131.139.251&21931&10950&5680&500458725&4.625.223.1663&519&1144 \cr
 & &32.3.25.7.13.73.241&91&80& & &64.3.11.13.173.223&223&352 \cr
\noalign{\hrule}
 & &27.5.7.11.163.293&61&754& & &81.625.11.29.31&1079&296 \cr
5663&496454805&4.3.13.29.61.293&425&454&5681&500630625&16.3.5.13.31.37.83&47&202 \cr
 & &16.25.13.17.61.227&13481&9080& & &64.13.37.47.101&22607&3232 \cr
\noalign{\hrule}
 & &81.5.7.13.97.139&451&34& & &25.49.169.41.59&2727&4202 \cr
5664&496915965&4.27.7.11.13.17.41&1577&880&5682&500793475&4.27.49.11.101.191&75&26 \cr
 & &128.5.121.19.83&10043&1216& & &16.81.25.11.13.191&891&1528 \cr
\noalign{\hrule}
 & &9.5.49.11.13.19.83&1139&226& & &9.11.43.191.617&3107&3680 \cr
5665&497251755&4.3.7.17.19.67.113&169&188&5683&501674679&64.3.5.13.23.43.239&3587&1910 \cr
 & &32.169.47.67.113&7571&9776& & &256.25.17.191.211&3587&3200 \cr
\noalign{\hrule}
 & &5.37.41.137.479&8759&8964& & &243.25.7.11.29.37&8611&1564 \cr
5666&497750455&8.27.19.83.137.461&17&154&5684&501922575&8.7.17.23.79.109&1053&290 \cr
 & &32.3.7.11.17.83.461&46563&51632& & &32.81.5.13.23.29&23&208 \cr
\noalign{\hrule}
 & &3.7.11.17.109.1163&345&1508& & &59.83.89.1153&53757&48860 \cr
5667&497814009&8.9.5.7.11.13.23.29&31&130&5685&502515449&8.27.5.7.11.181.349&1157&110 \cr
 & &32.25.169.29.31&22475&2704& & &32.9.25.121.13.89&3025&1872 \cr
\noalign{\hrule}
 & &3.49.11.13.19.29.43&305&2412& & &27.5.13.29.41.241&2489&644 \cr
5668&498050553&8.27.5.29.61.67&5089&3146&5686&502893495&8.3.7.19.23.29.131&4961&7712 \cr
 & &32.7.121.13.727&727&176& & &512.121.41.241&121&256 \cr
\noalign{\hrule}
 & &3.13.31.337.1223&627&596& & &27.7.23.37.53.59&4433&4100 \cr
5669&498290559&8.9.11.13.19.149.337&2875&158&5687&502943553&8.3.25.11.13.31.41.59&203&262 \cr
 & &32.125.23.79.149&18625&29072& & &32.5.7.11.13.29.41.131&41789&42640 \cr
\noalign{\hrule}
 & &9.5.7.169.17.19.29&473&382& & &3.25.11.13.17.31.89&427&2 \cr
5670&498652245&4.11.13.17.29.43.191&95&282&5688&503034675&4.7.31.61.89&901&990 \cr
 & &16.3.5.19.43.47.191&2021&1528& & &16.9.5.7.11.17.53&21&424 \cr
\noalign{\hrule}
}%
}
$$
\eject
\vglue -23 pt
\noindent\hskip 1 in\hbox to 6.5 in{\ 5689 -- 5724 \hfill\fbd 503250825 -- 513293495\frb}
\vskip -9 pt
$$
\vbox{
\nointerlineskip
\halign{\strut
    \vrule \ \ \hfil \frb #\ 
   &\vrule \hfil \ \ \fbb #\frb\ 
   &\vrule \hfil \ \ \frb #\ \hfil
   &\vrule \hfil \ \ \frb #\ 
   &\vrule \hfil \ \ \frb #\ \ \vrule \hskip 2 pt
   &\vrule \ \ \hfil \frb #\ 
   &\vrule \hfil \ \ \fbb #\frb\ 
   &\vrule \hfil \ \ \frb #\ \hfil
   &\vrule \hfil \ \ \frb #\ 
   &\vrule \hfil \ \ \frb #\ \vrule \cr%
\noalign{\hrule}
 & &3.25.49.11.59.211&1551&1340& & &243.11.13.19.769&5375&9236 \cr
5689&503250825&8.9.125.121.47.67&21733&6608&5707&507717639&8.9.125.43.2309&3683&5992 \cr
 & &256.7.59.103.211&103&128& & &128.5.7.29.107.127&25781&34240 \cr
\noalign{\hrule}
 & &27.7.433.6151&2981&3170& & &9.125.11.17.41.59&32701&23674 \cr
5690&503379387&4.5.11.271.317.433&947&1218&5708&508897125&4.7.19.53.89.617&2667&2050 \cr
 & &16.3.7.11.29.317.947&27463&27896& & &16.3.25.49.19.41.127&931&1016 \cr
\noalign{\hrule}
 & &9.7.11.13.17.19.173&6283&23300& & &3.25.7.13.17.41.107&313&222 \cr
5691&503413911&8.25.61.103.233&141&374&5709&509001675&4.9.5.17.37.41.313&3919&2354 \cr
 & &32.3.5.11.17.47.61&235&976& & &16.11.37.107.3919&3919&3256 \cr
\noalign{\hrule}
 & &3.17.43.163.1409&121&610& & &7.13.61.107.857&803&54 \cr
5692&503659731&4.5.121.61.1409&765&644&5710&509021149&4.27.11.13.61.73&749&200 \cr
 & &32.9.25.7.17.23.61&4209&2800& & &64.3.25.7.11.107&825&32 \cr
\noalign{\hrule}
 & &7.11.13.19.71.373&1545&1172& & &9.47.1091.1103&20675&30602 \cr
5693&503680177&8.3.5.7.71.103.293&1389&1096&5711&509026779&4.25.11.13.107.827&1109&282 \cr
 & &128.9.103.137.463&63431&59328& & &16.3.25.11.47.1109&1109&2200 \cr
\noalign{\hrule}
 & &27.25.11.13.17.307&1901&1476& & &25.11.13.131.1087&4221&22954 \cr
5694&503763975&8.243.13.41.1901&629&2530&5712&509069275&4.9.7.23.67.499&169&330 \cr
 & &32.5.11.17.23.37.41&943&592& & &16.27.5.11.169.67&871&216 \cr
\noalign{\hrule}
 & &3.5.121.13.23.929&679&250& & &7.11.17.19.59.347&675&328 \cr
5695&504154365&4.625.7.11.23.97&4553&2322&5713&509183983&16.27.25.7.11.19.41&17&116 \cr
 & &16.27.7.29.43.157&6751&14616& & &128.3.25.17.29.41&1025&5568 \cr
\noalign{\hrule}
 & &9.7.13.41.83.181&61&22& & &27.5.11.193.1777&1631&146 \cr
5696&504457317&4.3.7.11.41.61.181&3151&650&5714&509297085&4.7.73.193.233&559&792 \cr
 & &16.25.11.13.23.137&6325&1096& & &64.9.11.13.43.73&559&2336 \cr
\noalign{\hrule}
 & &11.29.41.173.223&117&106& & &9.5.7.11.13.17.23.29&461&344 \cr
5697&504574741&4.9.13.29.41.53.173&649&20630&5715&510765255&16.11.17.29.43.461&483&10 \cr
 & &16.3.5.11.59.2063&2063&7080& & &64.3.5.7.23.461&461&32 \cr
\noalign{\hrule}
 & &3.7.13.19.271.359&11993&8470& & &5.11.19.23.89.239&501&1546 \cr
5698&504638043&4.5.49.121.67.179&359&180&5716&511248485&4.3.167.239.773&267&506 \cr
 & &32.9.25.11.67.359&825&1072& & &16.9.11.23.89.167&167&72 \cr
\noalign{\hrule}
 & &9.13.17.19.31.431&389&820& & &7.11.13.17.151.199&11493&10100 \cr
5699&504925551&8.3.5.17.19.41.389&4741&2650&5717&511343833&8.9.25.17.101.1277&1057&2774 \cr
 & &32.125.11.53.431&1375&848& & &32.3.25.7.19.73.151&1387&1200 \cr
\noalign{\hrule}
 & &27.5.7.11.13.37.101&157&1156& & &27.5.7.13.17.31.79&47&506 \cr
5700&504999495&8.5.7.11.289.157&303&292&5718&511461405&4.5.11.13.23.31.47&6953&4938 \cr
 & &64.3.17.73.101.157&2669&2336& & &16.3.17.409.823&823&3272 \cr
\noalign{\hrule}
 & &81.7.19.107.439&163&8504& & &3.5.121.17.59.281&3429&3710 \cr
5701&506040129&16.7.163.1063&535&528&5719&511545045&4.81.25.7.17.53.127&1463&562 \cr
 & &512.3.5.11.107.163&1793&1280& & &16.49.11.19.127.281&931&1016 \cr
\noalign{\hrule}
 & &9.5.23.31.43.367&1045&56& & &3.49.13.19.59.239&925&748 \cr
5702&506333385&16.3.25.7.11.19.31&2977&3502&5720&511993209&8.25.7.11.13.17.19.37&673&4248 \cr
 & &64.13.17.103.229&23587&7072& & &128.9.17.59.673&2019&1088 \cr
\noalign{\hrule}
 & &9.5.11.23.79.563&637&74& & &13.17.89.131.199&23331&5620 \cr
5703&506370645&4.5.49.11.13.23.37&29&2064&5721&512751161&8.3.5.7.11.101.281&3&4 \cr
 & &128.3.7.29.43&43&12992& & &64.9.5.11.101.281&27819&16160 \cr
\noalign{\hrule}
 & &27.7.31.137.631&2165&1534& & &11.17.53.191.271&20925&30836 \cr
5704&506492973&4.5.7.13.31.59.433&275&492&5722&513003271&8.27.25.13.31.593&1759&1166 \cr
 & &32.3.125.11.41.433&17753&22000& & &32.3.11.31.53.1759&1759&1488 \cr
\noalign{\hrule}
 & &25.841.89.271&1547&22572& & &27.2203.8629&7619&1010 \cr
5705&507101975&8.27.7.11.13.17.19&145&178&5723&513261549&4.9.5.19.101.401&101167&101338 \cr
 & &32.9.5.7.13.29.89&63&208& & &16.11.17.23.541.2203&4301&4328 \cr
\noalign{\hrule}
 & &3.7.23.31.41.827&213&500& & &5.1331.13.17.349&8643&8660 \cr
5706&507689511&8.9.125.71.827&2387&1748&5724&513293495&8.3.25.43.67.349.433&363&712 \cr
 & &64.25.7.11.19.23.31&475&352& & &128.9.121.67.89.433&38537&38592 \cr
\noalign{\hrule}
}%
}
$$
\eject
\vglue -23 pt
\noindent\hskip 1 in\hbox to 6.5 in{\ 5725 -- 5760 \hfill\fbd 513422845 -- 522016495\frb}
\vskip -9 pt
$$
\vbox{
\nointerlineskip
\halign{\strut
    \vrule \ \ \hfil \frb #\ 
   &\vrule \hfil \ \ \fbb #\frb\ 
   &\vrule \hfil \ \ \frb #\ \hfil
   &\vrule \hfil \ \ \frb #\ 
   &\vrule \hfil \ \ \frb #\ \ \vrule \hskip 2 pt
   &\vrule \ \ \hfil \frb #\ 
   &\vrule \hfil \ \ \fbb #\frb\ 
   &\vrule \hfil \ \ \frb #\ \hfil
   &\vrule \hfil \ \ \frb #\ 
   &\vrule \hfil \ \ \frb #\ \vrule \cr%
\noalign{\hrule}
 & &5.169.19.113.283&989&1158& & &49.11.31.83.373&4085&18 \cr
5725&513422845&4.3.5.23.43.193.283&341&624&5743&517293931&4.9.5.19.31.43&99&56 \cr
 & &128.9.11.13.23.31.43&14663&13248& & &64.81.7.11.19&19&2592 \cr
\noalign{\hrule}
 & &3.5.11.17.31.61.97&213&128& & &81.343.103.181&4081&4262 \cr
5726&514512735&256.9.61.71.97&323&226&5744&517958469&4.2401.11.53.2131&135&2266 \cr
 & &1024.17.19.71.113&8023&9728& & &16.27.5.121.53.103&605&424 \cr
\noalign{\hrule}
 & &9.5.7.19.127.677&391&286& & &3.5.11.17.31.59.101&1021&494 \cr
5727&514584315&4.3.11.13.17.19.23.127&4253&700&5745&518164845&4.11.13.19.59.1021&279&488 \cr
 & &32.25.7.23.4253&4253&1840& & &64.9.31.61.1021&3063&1952 \cr
\noalign{\hrule}
 & &25.11.17.31.53.67&21&10& & &81.25.7.13.29.97&947&1672 \cr
5728&514628675&4.3.125.7.17.53.67&713&2838&5746&518365575&16.3.7.11.13.19.947&14761&16490 \cr
 & &16.9.7.11.23.31.43&387&1288& & &64.5.17.29.97.509&509&544 \cr
\noalign{\hrule}
 & &27.25.11.139.499&8531&3944& & &3.25.23.271.1109&16137&15028 \cr
5729&515005425&16.9.17.19.29.449&139&310&5747&518429775&8.27.5.11.13.289.163&2263&508 \cr
 & &64.5.17.29.31.139&493&992& & &64.11.17.31.73.127&38471&44704 \cr
\noalign{\hrule}
 & &11.13.17.23.61.151&135&118& & &5.7.11.37.59.617&18297&29212 \cr
5730&515014643&4.27.5.13.59.61.151&391&2354&5748&518560735&8.9.19.67.107.109&797&476 \cr
 & &16.3.11.17.23.59.107&321&472& & &64.3.7.17.109.797&13549&10464 \cr
\noalign{\hrule}
 & &3.49.11.13.107.229&8585&346& & &3.7.31.47.71.239&335&382 \cr
5731&515077563&4.5.7.17.101.173&2977&3078&5749&519200493&4.5.7.31.67.71.191&1305&14102 \cr
 & &16.81.13.17.19.229&323&216& & &16.9.25.11.29.641&7975&15384 \cr
\noalign{\hrule}
 & &81.97.173.379&11375&25388& & &243.5.49.11.13.61&493&178 \cr
5732&515159919&8.125.7.11.13.577&1499&2076&5750&519323805&4.27.7.13.17.29.89&12017&6050 \cr
 & &64.3.5.7.173.1499&1499&1120& & &16.25.121.61.197&197&440 \cr
\noalign{\hrule}
 & &3.5.11.13.17.71.199&3999&616& & &3.25.11.29.31.701&587&312 \cr
5733&515213985&16.9.7.121.31.43&563&284&5751&519914175&16.9.13.587.701&409&292 \cr
 & &128.43.71.563&563&2752& & &128.73.409.587&29857&37568 \cr
\noalign{\hrule}
 & &5.11.13.17.19.23.97&329&108& & &17.841.41.887&95&792 \cr
5734&515238295&8.27.5.7.11.47.97&589&104&5752&519938999&16.9.5.11.19.841&943&102 \cr
 & &128.3.13.19.31.47&1457&192& & &64.27.17.23.41&23&864 \cr
\noalign{\hrule}
 & &3.53.59.179.307&65&242& & &5.7.13.17.23.37.79&53&66 \cr
5735&515514093&4.5.121.13.53.179&2043&4370&5753&520016315&4.3.5.11.23.37.53.79&7917&6698 \cr
 & &16.9.25.19.23.227&4313&13800& & &16.9.7.11.13.17.29.197&2167&2088 \cr
\noalign{\hrule}
 & &3.5.49.11.23.47.59&181&204& & &3.125.11.13.17.571&443&118 \cr
5736&515653215&8.9.7.17.47.59.181&5&418&5754&520537875&4.5.59.443.571&1649&1206 \cr
 & &32.5.11.17.19.181&3077&304& & &16.9.17.59.67.97&3953&2328 \cr
\noalign{\hrule}
 & &3.5.7.11.29.89.173&107&338& & &125.49.167.509&3177&2668 \cr
5737&515722515&4.169.29.107.173&1813&3204&5755&520643375&8.9.25.7.23.29.353&3553&11672 \cr
 & &32.9.49.13.37.89&273&592& & &128.3.11.17.19.1459&27721&35904 \cr
\noalign{\hrule}
 & &9.7.19.29.107.139&679&572& & &9.5.13.17.23.43.53&385&1604 \cr
5738&516286449&8.49.11.13.19.29.97&19795&4164&5756&521287065&8.25.7.11.43.401&5163&4862 \cr
 & &64.3.5.37.107.347&1735&1184& & &32.3.121.13.17.1721&1721&1936 \cr
\noalign{\hrule}
 & &27.59.439.739&427&1166& & &243.5.19.97.233&1529&314 \cr
5739&516802653&4.7.11.53.61.439&2229&2600&5757&521744085&4.11.139.157.233&747&980 \cr
 & &64.3.25.13.61.743&19825&23776& & &32.9.5.49.83.139&4067&2224 \cr
\noalign{\hrule}
 & &9.5.11.17.19.53.61&899&1906& & &81.5.11.193.607&151&1972 \cr
5740&516908205&4.3.29.31.61.953&2375&484&5758&521907705&8.27.5.17.29.151&607&148 \cr
 & &32.125.121.19.29&725&176& & &64.29.37.607&37&928 \cr
\noalign{\hrule}
 & &5.13.31.43.47.127&101&114& & &5.121.29.151.197&159&38 \cr
5741&517184005&4.3.19.31.47.101.127&1419&4550&5759&521911115&4.3.5.19.29.53.151&1507&1362 \cr
 & &16.9.25.7.11.13.19.43&665&792& & &16.9.11.53.137.227&12031&9864 \cr
\noalign{\hrule}
 & &7.11.17.151.2617&4959&23278& & &5.7.11.169.71.113&75&68 \cr
5742&517273603&4.9.19.29.103.113&2617&660&5760&522016495&8.3.125.13.17.71.113&3199&4824 \cr
 & &32.27.5.11.2617&135&16& & &128.27.7.17.67.457&30619&29376 \cr
\noalign{\hrule}
}%
}
$$
\eject
\vglue -23 pt
\noindent\hskip 1 in\hbox to 6.5 in{\ 5761 -- 5796 \hfill\fbd 522321345 -- 529907939\frb}
\vskip -9 pt
$$
\vbox{
\nointerlineskip
\halign{\strut
    \vrule \ \ \hfil \frb #\ 
   &\vrule \hfil \ \ \fbb #\frb\ 
   &\vrule \hfil \ \ \frb #\ \hfil
   &\vrule \hfil \ \ \frb #\ 
   &\vrule \hfil \ \ \frb #\ \ \vrule \hskip 2 pt
   &\vrule \ \ \hfil \frb #\ 
   &\vrule \hfil \ \ \fbb #\frb\ 
   &\vrule \hfil \ \ \frb #\ \hfil
   &\vrule \hfil \ \ \frb #\ 
   &\vrule \hfil \ \ \frb #\ \vrule \cr%
\noalign{\hrule}
 & &27.5.7.13.17.41.61&377&172& & &27.47.59.79.89&7675&644 \cr
5761&522321345&8.3.7.169.17.29.43&5225&10126&5779&526418001&8.9.25.7.23.307&979&1784 \cr
 & &32.25.11.19.61.83&1577&880& & &128.5.11.89.223&2453&320 \cr
\noalign{\hrule}
 & &11.76729.619&41769&34960& & &5.121.23.157.241&999&1784 \cr
5762&522447761&32.27.5.7.13.17.19.23&277&22&5780&526501855&16.27.37.223.241&973&1196 \cr
 & &128.9.7.11.19.277&171&448& & &128.3.7.13.23.37.139&10101&8896 \cr
\noalign{\hrule}
 & &3.49.13.79.3461&4873&5510& & &5.11.43.251.887&25973&12168 \cr
5763&522503709&4.5.11.19.29.79.443&6381&2036&5781&526536505&16.9.169.19.1367&937&430 \cr
 & &32.9.29.509.709&20561&24432& & &64.3.5.19.43.937&937&1824 \cr
\noalign{\hrule}
 & &9.11.131.173.233&22865&202& & &7.37.47.181.239&15543&27716 \cr
5764&522767421&4.5.17.101.269&4427&4158&5782&526591807&8.9.11.169.41.157&141&310 \cr
 & &16.27.7.11.19.233&21&152& & &32.27.5.31.47.157&4185&2512 \cr
\noalign{\hrule}
 & &81.11.17.19.23.79&13&4& & &3.13.29.89.5233&100753&103334 \cr
5765&522919881&8.9.11.13.19.23.79&6751&1070&5783&526748547&4.7.121.53.61.1901&21771&860 \cr
 & &32.5.43.107.157&23005&2512& & &32.9.5.11.41.43.59&19393&14160 \cr
\noalign{\hrule}
 & &25.17.29.31.1369&633&2002& & &25.11.37.103.503&2509&3024 \cr
5766&523060675&4.3.5.7.11.13.29.211&111&34&5784&527156575&32.27.5.7.13.37.193&3229&824 \cr
 & &16.9.13.17.37.211&211&936& & &512.9.103.3229&3229&2304 \cr
\noalign{\hrule}
 & &9.125.7.13.19.269&3379&1496& & &13.289.841.167&2299&1458 \cr
5767&523238625&16.3.11.17.19.31.109&269&320&5785&527659379&4.729.121.19.167&1105&3404 \cr
 & &2048.5.11.109.269&1199&1024& & &32.27.5.13.17.23.37&999&1840 \cr
\noalign{\hrule}
 & &9.17.31.211.523&53&580& & &81.25.11.19.29.43&121&164 \cr
5768&523404279&8.3.5.29.53.523&341&182&5786&527761575&8.27.5.1331.29.41&2623&1292 \cr
 & &32.5.7.11.13.29.31&4147&560& & &64.17.19.41.43.61&1037&1312 \cr
\noalign{\hrule}
 & &5.7.11.67.103.197&36963&29032& & &9.11.13.17.19.31.41&299&398 \cr
5769&523406345&16.27.19.1369.191&23&14&5787&528355971&4.169.19.23.31.199&729&4510 \cr
 & &64.3.7.19.23.37.191&16169&18336& & &16.729.5.11.23.41&405&184 \cr
\noalign{\hrule}
 & &27.11.31.101.563&695&416& & &9.41.71.20173&715957&716326 \cr
5770&523537641&64.3.5.13.139.563&73&490&5788&528512427&4.121.13.61.97.27551&10085&17466 \cr
 & &256.25.49.13.73&23725&6272& & &16.3.5.13.41.71.97.2017&10085&10088 \cr
\noalign{\hrule}
 & &3.19.23.401.997&235&166& & &5.11.19.29.31.563&3015&3046 \cr
5771&524133867&4.5.19.47.83.997&1287&290&5789&528913165&4.9.25.67.563.1523&7&1682 \cr
 & &16.9.25.11.13.29.47&15275&7656& & &16.3.7.841.1523&4569&1624 \cr
\noalign{\hrule}
 & &729.11.13.47.107&301&310& & &25.7.11.17.19.23.37&603&442 \cr
5772&524258163&4.81.5.7.11.31.43.107&5989&2678&5790&529130525&4.9.5.13.289.37.67&133&422 \cr
 & &16.5.13.31.53.103.113&27295&28024& & &16.3.7.13.19.67.211&2613&1688 \cr
\noalign{\hrule}
 & &27.5.19.149.1373&67229&63206& & &5.11.47.257.797&619&666 \cr
5773&524740005&4.11.169.17.23.37.79&545&798&5791&529482965&4.9.11.37.619.797&2209&4600 \cr
 & &16.3.5.7.169.19.37.109&6253&6104& & &64.3.25.23.37.2209&4255&4512 \cr
\noalign{\hrule}
 & &81.5.47.79.349&5395&5744& & &2209.53.4523&2235&2288 \cr
5774&524813985&32.27.25.13.83.359&517&158&5792&529539271&32.3.5.11.13.2209.149&285&1924 \cr
 & &128.11.13.47.79.83&913&832& & &256.9.25.169.19.37&56277&60800 \cr
\noalign{\hrule}
 & &3.5.49.11.43.1511&6023&4512& & &81.13.17.101.293&739&638 \cr
5775&525306705&64.9.11.19.47.317&2275&5762&5793&529744293&4.11.13.29.293.739&5969&2160 \cr
 & &256.25.7.13.43.67&871&640& & &128.27.5.29.47.127&6815&8128 \cr
\noalign{\hrule}
 & &7.29.43.139.433&703&270& & &7.23.29.233.487&3135&3622 \cr
5776&525372323&4.27.5.19.29.37.43&539&278&5794&529796099&4.3.5.7.11.19.23.1811&825&986 \cr
 & &16.3.5.49.11.37.139&407&840& & &16.9.125.121.17.19.29&18513&19000 \cr
\noalign{\hrule}
 & &25.841.131.191&23023&1998& & &5.7.11.19.23.47.67&533&204 \cr
5777&526066525&4.27.7.11.13.23.37&145&262&5795&529803505&8.3.5.13.17.19.23.41&519&424 \cr
 & &16.3.5.7.23.29.131&23&168& & &128.9.13.17.53.173&38233&30528 \cr
\noalign{\hrule}
 & &7.41.61.107.281&2407&1980& & &11.47.53.83.233&3481&918 \cr
5778&526382969&8.9.5.11.29.83.281&265&16&5796&529907939&4.27.17.47.3481&429&370 \cr
 & &256.3.25.11.29.53&23925&6784& & &16.81.5.11.13.37.59&10915&8424 \cr
\noalign{\hrule}
}%
}
$$
\eject
\vglue -23 pt
\noindent\hskip 1 in\hbox to 6.5 in{\ 5797 -- 5832 \hfill\fbd 529925445 -- 538839791\frb}
\vskip -9 pt
$$
\vbox{
\nointerlineskip
\halign{\strut
    \vrule \ \ \hfil \frb #\ 
   &\vrule \hfil \ \ \fbb #\frb\ 
   &\vrule \hfil \ \ \frb #\ \hfil
   &\vrule \hfil \ \ \frb #\ 
   &\vrule \hfil \ \ \frb #\ \ \vrule \hskip 2 pt
   &\vrule \ \ \hfil \frb #\ 
   &\vrule \hfil \ \ \fbb #\frb\ 
   &\vrule \hfil \ \ \frb #\ \hfil
   &\vrule \hfil \ \ \frb #\ 
   &\vrule \hfil \ \ \frb #\ \vrule \cr%
\noalign{\hrule}
 & &9.5.49.17.67.211&3301&286& & &9.7.11.409.1889&6251&10750 \cr
5797&529925445&4.49.11.13.3301&1675&1626&5815&535412493&4.125.49.19.43.47&591&634 \cr
 & &16.3.25.11.13.67.271&1355&1144& & &16.3.5.19.47.197.317&46295&48184 \cr
\noalign{\hrule}
 & &9.25.7.11.71.431&2077&1802& & &23.37.47.59.227&53911&48690 \cr
5798&530162325&4.7.17.31.53.67.71&30615&2684&5816&535679821&4.9.5.11.169.29.541&1741&118 \cr
 & &32.3.5.11.13.61.157&2041&976& & &16.3.5.29.59.1741&1741&3480 \cr
\noalign{\hrule}
 & &9.5.11.13.17.37.131&127&94& & &3.5.7.361.67.211&2015&1804 \cr
5799&530237565&4.3.5.37.47.127.131&289&104&5817&535862985&8.25.7.11.13.19.31.41&211&36 \cr
 & &64.13.289.47.127&2159&1504& & &64.9.11.31.41.211&1353&992 \cr
\noalign{\hrule}
 & &3.25.11.41.61.257&3243&3182& & &25.11.169.83.139&277&138 \cr
5800&530274525&4.9.11.23.37.41.43.47&257&2020&5818&536182075&4.3.5.11.169.23.277&1101&1946 \cr
 & &32.5.37.47.101.257&1739&1616& & &16.9.7.23.139.367&2569&1656 \cr
\noalign{\hrule}
 & &9.11.23.29.83.97&471&442& & &7.19.73.101.547&4843&5550 \cr
5801&531631683&4.27.13.17.23.97.157&45733&1210&5819&536393123&4.3.25.29.37.73.167&3413&8778 \cr
 & &16.5.121.19.29.83&55&152& & &16.9.5.7.11.19.3413&3413&3960 \cr
\noalign{\hrule}
 & &5.7.11.19.23.29.109&393&370& & &3.5.7.11.137.3391&2449&942 \cr
5802&531822445&4.3.25.11.19.29.37.131&2743&5232&5820&536574885&4.9.5.7.31.79.157&3809&3256 \cr
 & &128.9.13.37.109.211&7807&7488& & &64.11.13.31.37.293&10841&12896 \cr
\noalign{\hrule}
 & &3.5.13.61.97.461&627&166& & &9.5.11.13.89.937&3811&874 \cr
5803&531908715&4.9.5.11.19.83.97&4453&3598&5821&536633955&4.3.13.19.23.37.103&1325&14 \cr
 & &16.7.11.61.73.257&2827&4088& & &16.25.7.37.53&53&10360 \cr
\noalign{\hrule}
 & &11.29.41.67.607&8433&9170& & &9.5.11.17.227.281&9657&5798 \cr
5804&531909851&4.9.5.7.41.131.937&899&38&5822&536767605&4.81.13.29.37.223&281&200 \cr
 & &16.3.5.19.29.31.131&1965&4712& & &64.25.29.223.281&1115&928 \cr
\noalign{\hrule}
 & &125.49.13.41.163&2871&23246& & &3.121.23.131.491&185&306 \cr
5805&532133875&4.9.11.29.59.197&325&266&5823&537016029&4.27.5.17.23.37.131&1573&1964 \cr
 & &16.3.25.7.11.13.19.29&551&264& & &32.5.121.13.37.491&185&208 \cr
\noalign{\hrule}
 & &9.11.29.151.1229&67&386& & &3.5.13.37.163.457&247&210 \cr
5806&532797309&4.3.67.193.1229&715&514&5824&537452565&4.9.25.7.169.19.163&19943&16732 \cr
 & &16.5.11.13.193.257&3341&7720& & &32.343.11.37.47.89&16121&15664 \cr
\noalign{\hrule}
 & &9.7.13.37.43.409&1265&1598& & &9.7.11.17.43.1061&37&26 \cr
5807&532938861&4.5.11.13.17.23.43.47&11767&5046&5825&537484563&4.13.17.37.43.1061&1323&19360 \cr
 & &16.3.5.7.841.1681&8405&6728& & &256.27.5.49.121&165&896 \cr
\noalign{\hrule}
 & &25.29.41.131.137&65151&69124& & &81.11.17.131.271&2107&2500 \cr
5808&533474575&8.27.11.19.127.1571&2881&1310&5826&537733647&8.27.625.49.11.43&289&586 \cr
 & &32.9.5.19.43.67.131&2881&2736& & &32.5.7.289.43.293&10255&11696 \cr
\noalign{\hrule}
 & &9.7.11.19.107.379&519&7720& & &11.13.19.37.53.101&153&1160 \cr
5809&533961351&16.27.5.173.193&2173&3038&5827&538131737&16.9.5.11.17.29.37&247&160 \cr
 & &64.49.31.41.53&8897&1696& & &1024.3.25.13.17.19&425&1536 \cr
\noalign{\hrule}
 & &27.11.31.131.443&133&164& & &5.13.19.59.83.89&371&786 \cr
5810&534309831&8.7.19.41.131.443&2795&306&5828&538253755&4.3.7.19.53.59.131&319&2808 \cr
 & &32.9.5.13.17.41.43&3655&8528& & &64.81.7.11.13.29&891&6496 \cr
\noalign{\hrule}
 & &125.7.11.19.37.79&1391&984& & &81.11.13.23.43.47&8375&8122 \cr
5811&534543625&16.3.7.13.41.79.107&657&370&5829&538412589&4.3.125.31.43.67.131&299&5932 \cr
 & &64.27.5.37.73.107&2889&2336& & &32.125.13.23.1483&1483&2000 \cr
\noalign{\hrule}
 & &5.11.37.41.43.149&333&118& & &9.11.13.17.151.163&3103&3254 \cr
5812&534568045&4.9.1369.59.149&13871&12502&5830&538507827&4.3.11.17.29.107.1627&53605&6422 \cr
 & &16.3.7.11.13.19.47.97&12901&14664& & &16.5.169.19.71.151&923&760 \cr
\noalign{\hrule}
 & &9.5.13.43.89.239&23&66& & &3.5.7.11.169.31.89&1431&274 \cr
5813&535072005&4.27.5.11.13.23.239&1333&2528&5831&538543005&4.81.7.13.53.137&341&712 \cr
 & &256.23.31.43.79&1817&3968& & &64.11.31.89.137&137&32 \cr
\noalign{\hrule}
 & &9.125.7.11.37.167&131&1706& & &7.361.23.73.127&625&264 \cr
5814&535255875&4.5.37.131.853&519&334&5832&538839791&16.3.625.11.23.73&977&702 \cr
 & &16.3.131.167.173&173&1048& & &64.81.25.13.977&24425&33696 \cr
\noalign{\hrule}
}%
}
$$
\eject
\vglue -23 pt
\noindent\hskip 1 in\hbox to 6.5 in{\ 5833 -- 5868 \hfill\fbd 539236035 -- 548565435\frb}
\vskip -9 pt
$$
\vbox{
\nointerlineskip
\halign{\strut
    \vrule \ \ \hfil \frb #\ 
   &\vrule \hfil \ \ \fbb #\frb\ 
   &\vrule \hfil \ \ \frb #\ \hfil
   &\vrule \hfil \ \ \frb #\ 
   &\vrule \hfil \ \ \frb #\ \ \vrule \hskip 2 pt
   &\vrule \ \ \hfil \frb #\ 
   &\vrule \hfil \ \ \fbb #\frb\ 
   &\vrule \hfil \ \ \frb #\ \hfil
   &\vrule \hfil \ \ \frb #\ 
   &\vrule \hfil \ \ \frb #\ \vrule \cr%
\noalign{\hrule}
 & &81.5.13.23.61.73&1295&3158& & &3.13.47.401.739&105&506 \cr
5833&539236035&4.25.7.13.37.1579&33&292&5851&543189387&4.9.5.7.11.23.739&401&338 \cr
 & &32.3.11.73.1579&1579&176& & &16.5.11.169.23.401&299&440 \cr
\noalign{\hrule}
 & &9.5.49.11.13.29.59&8455&10366& & &13.17.23.89.1201&607&594 \cr
5834&539503965&4.3.25.19.71.73.89&59059&64384&5852&543316787&4.27.11.17.23.89.607&1865&182 \cr
 & &1024.7.11.13.59.503&503&512& & &16.3.5.7.13.373.607&13055&14568 \cr
\noalign{\hrule}
 & &9.25.41.233.251&2233&8058& & &9.11.71.167.463&133&34 \cr
5835&539505675&4.27.7.11.17.29.79&233&260&5853&543489309&4.7.17.19.71.463&835&372 \cr
 & &32.5.7.11.13.79.233&1027&1232& & &32.3.5.7.19.31.167&1085&304 \cr
\noalign{\hrule}
 & &9.11.289.113.167&29725&18538& & &7.13.67.191.467&1199&2070 \cr
5836&539918181&4.25.13.23.29.31.41&969&302&5854&543834109&4.9.5.11.23.109.191&233&1868 \cr
 & &16.3.25.13.17.19.151&1963&3800& & &32.3.23.233.467&233&1104 \cr
\noalign{\hrule}
 & &3.5.7.13.19.59.353&1019&746& & &9.7.11.13.17.53.67&4553&6100 \cr
5837&540148245&4.19.59.373.1019&747&374&5855&543846303&8.3.25.11.29.61.157&161&4 \cr
 & &16.9.11.17.83.1019&17323&21912& & &64.5.7.23.29.61&8845&736 \cr
\noalign{\hrule}
 & &3.25.121.13.19.241&305&58& & &81.13.41.43.293&833&274 \cr
5838&540207525&4.125.29.61.241&3933&3692&5856&543936627&4.3.49.17.137.293&18821&21320 \cr
 & &32.9.13.19.23.29.71&1633&1392& & &64.5.11.13.29.41.59&1595&1888 \cr
\noalign{\hrule}
 & &27.169.109.1087&4183&5600& & &3.5.11.13.19.31.431&633&602 \cr
5839&540637929&64.3.25.7.13.47.89&17&22&5857&544527555&4.9.7.11.43.211.431&19711&1178 \cr
 & &256.5.7.11.17.47.89&61523&56960& & &16.7.19.23.31.857&857&1288 \cr
\noalign{\hrule}
 & &81.29.43.53.101&1921&1820& & &3.125.11.841.157&233&552 \cr
5840&540690471&8.27.5.7.13.17.53.113&4679&1628&5858&544652625&16.9.25.23.29.233&157&418 \cr
 & &64.5.11.13.37.4679&60827&65120& & &64.11.19.157.233&233&608 \cr
\noalign{\hrule}
 & &27.5.11.13.109.257&1409&1418& & &13.19.23.41.2341&28369&16110 \cr
5841&540790965&4.3.5.13.109.709.1409&1397&5648&5859&545268061&4.9.5.11.179.2579&1021&1558 \cr
 & &128.11.127.353.709&44831&45376& & &16.3.5.11.19.41.1021&1021&1320 \cr
\noalign{\hrule}
 & &9.11.13.23.101.181&14103&11780& & &27.17.151.7873&3707&4166 \cr
5842&541135881&8.81.5.19.31.1567&14927&14846&5860&545669757&4.11.151.337.2083&157435&157098 \cr
 & &32.5.11.13.23.31.59.571&9145&9136& & &16.3.5.11.23.1369.26183&602209&602360 \cr
\noalign{\hrule}
 & &27.5.11.13.289.97&31&904& & &9.5.17.43.53.313&539&226 \cr
5843&541177065&16.3.13.17.31.113&275&388&5861&545695155&4.49.11.43.53.113&313&270 \cr
 & &128.25.11.31.97&31&320& & &16.27.5.49.113.313&339&392 \cr
\noalign{\hrule}
 & &27.5.73.137.401&847&2852& & &3.25.7.13.29.31.89&451&74 \cr
5844&541404135&8.7.121.23.31.73&801&878&5862&546075075&4.11.31.37.41.89&1125&146 \cr
 & &32.9.11.31.89.439&13609&15664& & &16.9.125.37.73&2701&120 \cr
\noalign{\hrule}
 & &7.13.31.47.61.67&963&494& & &7.11.19.29.61.211&169&150 \cr
5845&541883069&4.9.169.19.61.107&7303&10780&5863&546077917&4.3.25.7.169.61.211&1501&24 \cr
 & &32.3.5.49.11.67.109&1635&1232& & &64.9.169.19.79&711&5408 \cr
\noalign{\hrule}
 & &17.41.877.887&459&418& & &3.7.17.41.67.557&869&802 \cr
5846&542195603&4.27.11.289.19.887&877&10&5864&546238203&4.7.11.17.41.79.401&4645&234 \cr
 & &16.9.5.11.19.877&95&792& & &16.9.5.13.79.929&15405&7432 \cr
\noalign{\hrule}
 & &7.11.13.41.73.181&303&230& & &9.25.17.19.73.103&421&506 \cr
5847&542274733&4.3.5.7.11.23.101.181&1189&78&5865&546443325&4.5.11.19.23.73.421&61&1326 \cr
 & &16.9.5.13.23.29.41&1035&232& & &16.3.13.17.61.421&793&3368 \cr
\noalign{\hrule}
 & &25.169.19.29.233&2761&3996& & &5.11.31.107.3001&3159&158 \cr
5848&542418175&8.27.5.11.13.37.251&151&34&5866&547487435&4.243.5.11.13.79&181&214 \cr
 & &32.3.11.17.151.251&12801&26576& & &16.81.13.107.181&2353&648 \cr
\noalign{\hrule}
 & &3.7.29.31.41.701&4237&670& & &9.25.23.29.41.89&2527&302 \cr
5849&542601339&4.5.19.31.67.223&3157&1080&5867&547623675&4.3.7.361.29.151&3649&6820 \cr
 & &64.27.25.7.11.41&99&800& & &32.5.11.31.41.89&31&176 \cr
\noalign{\hrule}
 & &11.13.17.29.43.179&369&190& & &9.5.11.17.19.47.73&593&958 \cr
5850&542630803&4.9.5.11.17.19.29.41&827&130&5868&548565435&4.3.17.19.479.593&781&188 \cr
 & &16.3.25.13.19.827&1425&6616& & &32.11.47.71.479&479&1136 \cr
\noalign{\hrule}
}%
}
$$
\eject
\vglue -23 pt
\noindent\hskip 1 in\hbox to 6.5 in{\ 5869 -- 5904 \hfill\fbd 548789931 -- 557902345\frb}
\vskip -9 pt
$$
\vbox{
\nointerlineskip
\halign{\strut
    \vrule \ \ \hfil \frb #\ 
   &\vrule \hfil \ \ \fbb #\frb\ 
   &\vrule \hfil \ \ \frb #\ \hfil
   &\vrule \hfil \ \ \frb #\ 
   &\vrule \hfil \ \ \frb #\ \ \vrule \hskip 2 pt
   &\vrule \ \ \hfil \frb #\ 
   &\vrule \hfil \ \ \fbb #\frb\ 
   &\vrule \hfil \ \ \frb #\ \hfil
   &\vrule \hfil \ \ \frb #\ 
   &\vrule \hfil \ \ \frb #\ \vrule \cr%
\noalign{\hrule}
 & &27.31.53.89.139&4921&9230& & &3.5.7.13.31.103.127&209&426 \cr
5869&548789931&4.9.5.7.13.19.37.71&11&106&5887&553522515&4.9.11.13.19.71.103&3175&5132 \cr
 & &16.7.11.37.53.71&497&3256& & &32.25.11.127.1283&1283&880 \cr
\noalign{\hrule}
 & &3.7.11.17.23.59.103&225&166& & &27.11.17.19.29.199&175&148 \cr
5870&548880717&4.27.25.7.11.83.103&3757&10148&5888&553617801&8.25.7.11.29.37.199&101&1494 \cr
 & &32.5.13.289.43.59&731&1040& & &32.9.5.37.83.101&8383&2960 \cr
\noalign{\hrule}
 & &3.7.29.53.73.233&605&94& & &9.5.7.13.17.73.109&9&82 \cr
5871&548999493&4.5.121.29.47.53&27&292&5889&553926555&4.81.5.17.41.109&3497&3388 \cr
 & &32.27.11.47.73&47&1584& & &32.7.121.13.41.269&4961&4304 \cr
\noalign{\hrule}
 & &7.11.23.31.73.137&1885&378& & &125.7.11.13.19.233&17&108 \cr
5872&549064901&4.27.5.49.13.23.29&1319&682&5890&553928375&8.27.11.17.19.233&2921&1040 \cr
 & &16.9.5.11.31.1319&1319&360& & &256.3.5.13.23.127&381&2944 \cr
\noalign{\hrule}
 & &9.25.31.227.347&5401&5356& & &27.625.11.29.103&1043&7918 \cr
5873&549413775&8.5.11.13.103.227.491&1&1134&5891&554461875&4.9.7.37.107.149&3445&3296 \cr
 & &32.81.7.13.491&819&7856& & &256.5.13.37.53.103&1961&1664 \cr
\noalign{\hrule}
 & &125.37.211.563&39091&31284& & &9.5.11.37.157.193&69&124 \cr
5874&549417625&8.9.11.13.31.79.97&2825&376&5892&554962815&8.27.23.31.37.157&6721&12530 \cr
 & &128.3.25.13.47.113&5311&2496& & &32.5.7.11.13.47.179&4277&2864 \cr
\noalign{\hrule}
 & &7.11.13.41.59.227&24463&6984& & &9.25.49.11.23.199&361&214 \cr
5875&549662113&16.9.17.97.1439&865&574&5893&555075675&4.3.11.361.107.199&1955&4144 \cr
 & &64.3.5.7.17.41.173&865&1632& & &128.5.7.17.19.23.37&703&1088 \cr
\noalign{\hrule}
 & &3.5.7.13.17.23.1031&87&74& & &25.19.71.101.163&969&806 \cr
5876&550260165&4.9.5.17.29.37.1031&253&5408&5894&555214675&4.3.13.17.361.31.101&781&3912 \cr
 & &256.11.169.23.29&377&1408& & &64.9.11.17.71.163&187&288 \cr
\noalign{\hrule}
 & &27.25.19.29.1481&27887&15062& & &27.11.23.31.43.61&4105&3764 \cr
5877&550820925&4.17.79.353.443&495&848&5895&555449103&8.9.5.23.821.941&47&988 \cr
 & &128.9.5.11.53.443&4873&3392& & &64.13.19.47.821&11609&26272 \cr
\noalign{\hrule}
 & &25.7.17.23.83.97&2167&258& & &729.11.19.41.89&18275&18256 \cr
5878&550889675&4.3.7.11.17.43.197&97&90&5896&555965289&32.9.25.7.17.43.89.163&103&698 \cr
 & &16.27.5.43.97.197&1161&1576& & &128.5.43.103.163.349&284435&283456 \cr
\noalign{\hrule}
 & &3.5.31.43.47.587&12597&12644& & &9.25.113.131.167&12351&9526 \cr
5879&551642055&8.9.5.13.17.19.29.31.109&1&154&5897&556222725&4.27.11.23.179.433&113&140 \cr
 & &32.7.11.13.19.29.109&26923&35728& & &32.5.7.113.179.433&3031&2864 \cr
\noalign{\hrule}
 & &9.11.13.193.2221&925&23506& & &81.5.13.19.67.83&517&562 \cr
5880&551676411&4.25.7.23.37.73&821&858&5898&556294635&4.9.11.19.47.67.281&5395&13432 \cr
 & &16.3.25.7.11.13.821&821&1400& & &64.5.11.13.23.73.83&803&736 \cr
\noalign{\hrule}
 & &3.29.41.271.571&197&74& & &7.13.29.433.487&3267&3064 \cr
5881&551961147&4.29.37.197.571&4635&11924&5899&556488569&16.27.121.383.433&4055&158 \cr
 & &32.9.5.11.103.271&309&880& & &64.3.5.11.79.811&26763&12640 \cr
\noalign{\hrule}
 & &3.25.121.127.479&377&102& & &9.5.13.289.37.89&8239&2454 \cr
5882&552059475&4.9.11.13.17.29.127&2395&998&5900&556731045&4.27.7.11.107.409&89&100 \cr
 & &16.5.17.479.499&499&136& & &32.25.89.107.409&2045&1712 \cr
\noalign{\hrule}
 & &31.67.73.3643&783&2860& & &3.13.23.613.1013&11343&11956 \cr
5883&552355303&8.27.5.11.13.29.73&251&1054&5901&557009193&8.9.49.13.19.61.199&1265&13802 \cr
 & &32.3.13.17.31.251&4267&624& & &32.5.7.11.23.67.103&5665&7504 \cr
\noalign{\hrule}
 & &81.31.359.613&24041&5038& & &7.13.19.23.107.131&165&1868 \cr
5884&552588237&4.11.29.229.829&845&1674&5902&557414039&8.3.5.7.11.23.467&153&314 \cr
 & &16.27.5.169.29.31&845&232& & &32.27.5.11.17.157&2669&23760 \cr
\noalign{\hrule}
 & &9.11.19.83.3541&229&20& & &3.25.121.13.29.163&851&3876 \cr
5885&552831543&8.3.5.229.3541&8509&9196&5903&557667825&8.9.13.17.19.23.37&103&220 \cr
 & &64.121.19.67.127&1397&2144& & &64.5.11.23.37.103&851&3296 \cr
\noalign{\hrule}
 & &9.49.11.13.31.283&3835&722& & &5.7.11.13.17.79.83&24707&7248 \cr
5886&553251699&4.3.5.169.361.59&283&224&5904&557902345&32.3.31.151.797&2739&1942 \cr
 & &256.5.7.361.283&361&640& & &128.9.11.83.971&971&576 \cr
\noalign{\hrule}
}%
}
$$
\eject
\vglue -23 pt
\noindent\hskip 1 in\hbox to 6.5 in{\ 5905 -- 5940 \hfill\fbd 558145203 -- 565834731\frb}
\vskip -9 pt
$$
\vbox{
\nointerlineskip
\halign{\strut
    \vrule \ \ \hfil \frb #\ 
   &\vrule \hfil \ \ \fbb #\frb\ 
   &\vrule \hfil \ \ \frb #\ \hfil
   &\vrule \hfil \ \ \frb #\ 
   &\vrule \hfil \ \ \frb #\ \ \vrule \hskip 2 pt
   &\vrule \ \ \hfil \frb #\ 
   &\vrule \hfil \ \ \fbb #\frb\ 
   &\vrule \hfil \ \ \frb #\ \hfil
   &\vrule \hfil \ \ \frb #\ 
   &\vrule \hfil \ \ \frb #\ \vrule \cr%
\noalign{\hrule}
 & &3.7.11.43.83.677&1037&994& & &7.11.101.233.311&1205&972 \cr
5905&558145203&4.49.11.17.61.71.83&873&40&5923&563544751&8.243.5.11.101.241&2537&10718 \cr
 & &64.9.5.61.71.97&17751&11360& & &32.3.23.43.59.233&1357&2064 \cr
\noalign{\hrule}
 & &11.19.29.37.47.53&973&390& & &9.7.13.59.107.109&9515&9424 \cr
5906&558624187&4.3.5.7.13.19.37.139&1827&1688&5924&563567823&32.3.5.11.19.31.109.173&2461&826 \cr
 & &64.27.49.13.29.211&10339&11232& & &128.7.11.23.31.59.107&713&704 \cr
\noalign{\hrule}
 & &9.43.47.59.521&1529&1008& & &3.7.31.71.89.137&1189&3058 \cr
5907&559111671&32.81.7.11.47.139&6385&148&5925&563572653&4.11.29.41.71.139&405&376 \cr
 & &256.5.37.1277&1277&23680& & &64.81.5.41.47.139&32665&35424 \cr
\noalign{\hrule}
 & &5.529.73.2897&2771&126& & &7.11.19.47.59.139&3473&3060 \cr
5908&559367245&4.9.7.17.73.163&1495&1276&5926&563908961&8.9.5.11.17.19.23.151&1123&188 \cr
 & &32.3.5.7.11.13.23.29&609&2288& & &64.3.47.151.1123&3369&4832 \cr
\noalign{\hrule}
 & &9.11.43.47.2797&3185&5206& & &5.11.19.29.37.503&43497&49558 \cr
5909&559620963&4.3.5.49.11.13.19.137&2797&80&5927&564006355&4.243.71.179.349&1643&11066 \cr
 & &128.25.7.2797&25&448& & &16.9.11.31.53.503&279&424 \cr
\noalign{\hrule}
 & &27.5.49.11.43.179&2461&850& & &3.5.11.73.139.337&179&516 \cr
5910&560072205&4.3.125.7.17.23.107&4687&3938&5928&564223935&8.9.11.43.73.179&16459&18070 \cr
 & &16.11.17.43.109.179&109&136& & &32.5.13.109.139.151&1963&1744 \cr
\noalign{\hrule}
 & &27.5.13.17.19.23.43&11&334& & &5.11.19.29.43.433&69&26 \cr
5911&560629485&4.9.11.13.43.167&2599&2432&5929&564248795&4.3.11.13.23.29.433&57&376 \cr
 & &1024.11.19.23.113&1243&512& & &64.9.13.19.23.47&5499&736 \cr
\noalign{\hrule}
 & &27.125.7.11.17.127&169&466& & &7.17.41.109.1061&55385&60264 \cr
5912&561070125&4.25.7.169.17.233&23&198&5930&564251471&16.243.5.11.19.31.53&7&88 \cr
 & &16.9.11.13.23.233&233&2392& & &256.3.7.121.31.53&11253&6784 \cr
\noalign{\hrule}
 & &9.5.49.19.59.227&901&220& & &3.625.7.19.31.73&643&132 \cr
5913&561099735&8.3.25.49.11.17.53&1487&262&5931&564335625&8.9.25.11.19.643&73&98 \cr
 & &32.17.131.1487&2227&23792& & &32.49.11.73.643&643&1232 \cr
\noalign{\hrule}
 & &9.5.7.13.23.67.89&767&1102& & &9.5.17.37.127.157&559&226 \cr
5914&561625155&4.3.169.19.23.29.59&623&4510&5932&564373395&4.13.17.43.113.127&429&302 \cr
 & &16.5.7.11.19.41.89&451&152& & &16.3.11.169.113.151&19097&13288 \cr
\noalign{\hrule}
 & &49.19.37.47.347&9545&3294& & &27.7.121.23.29.37&6665&24674 \cr
5915&561796123&4.27.5.7.23.61.83&1807&3256&5933&564384051&4.5.169.31.43.73&1771&1368 \cr
 & &64.3.5.11.13.37.139&2145&4448& & &64.9.5.7.11.13.19.23&247&160 \cr
\noalign{\hrule}
 & &81.5.17.79.1033&3959&2926& & &9.11.31.71.2591&149&490 \cr
5916&561864195&4.7.11.19.37.79.107&873&2050&5934&564576309&4.5.49.149.2591&923&1668 \cr
 & &16.9.25.7.19.41.97&3977&5320& & &32.3.49.13.71.139&1807&784 \cr
\noalign{\hrule}
 & &9.7.121.13.53.107&5917&7030& & &11.13.97.193.211&437&630 \cr
5917&561990429&4.3.5.13.19.37.61.97&30899&11984&5935&564868733&4.9.5.7.13.19.23.211&193&440 \cr
 & &128.7.11.2809.107&53&64& & &64.3.25.7.11.23.193&525&736 \cr
\noalign{\hrule}
 & &27.5.53.127.619&3887&3268& & &9.25.11.53.59.73&6929&6346 \cr
5918&562476015&8.169.19.23.43.127&5571&110&5936&564970725&4.169.19.41.73.167&71921&59730 \cr
 & &32.9.5.11.13.619&143&16& & &16.3.5.11.23.53.59.181&181&184 \cr
\noalign{\hrule}
 & &3.1331.13.37.293&4755&6086& & &3.7.13.53.139.281&44935&50836 \cr
5919&562745469&4.9.5.13.17.179.317&1153&836&5937&565144671&8.5.11.19.43.71.179&3653&252 \cr
 & &32.5.11.19.179.1153&17005&18448& & &64.9.7.13.43.281&43&96 \cr
\noalign{\hrule}
 & &5.13.83.151.691&49473&54868& & &27.25.13.73.883&3193&1222 \cr
5920&562919695&8.9.11.23.29.43.239&4567&3320&5938&565627725&4.5.169.31.47.103&407&438 \cr
 & &128.3.5.23.83.4567&4567&4416& & &16.3.11.37.47.73.103&3811&4136 \cr
\noalign{\hrule}
 & &27.5.17.31.41.193&2915&3068& & &9.5.11.17.23.37.79&6641&7826 \cr
5921&562970385&8.3.25.11.13.41.53.59&1537&62&5939&565732035&4.3.7.11.13.29.43.229&5&38 \cr
 & &32.11.29.31.2809&2809&5104& & &16.5.7.13.19.29.229&7163&12824 \cr
\noalign{\hrule}
 & &27.11.13.139.1049&589&940& & &3.7.11.53.113.409&19&390 \cr
5922&562976271&8.5.19.31.47.1049&477&572&5940&565834731&4.9.5.11.13.19.113&7&106 \cr
 & &64.9.11.13.31.47.53&1643&1504& & &16.5.7.13.19.53&1235&8 \cr
\noalign{\hrule}
}%
}
$$
\eject
\vglue -23 pt
\noindent\hskip 1 in\hbox to 6.5 in{\ 5941 -- 5976 \hfill\fbd 565957161 -- 574757391\frb}
\vskip -9 pt
$$
\vbox{
\nointerlineskip
\halign{\strut
    \vrule \ \ \hfil \frb #\ 
   &\vrule \hfil \ \ \fbb #\frb\ 
   &\vrule \hfil \ \ \frb #\ \hfil
   &\vrule \hfil \ \ \frb #\ 
   &\vrule \hfil \ \ \frb #\ \ \vrule \hskip 2 pt
   &\vrule \ \ \hfil \frb #\ 
   &\vrule \hfil \ \ \fbb #\frb\ 
   &\vrule \hfil \ \ \frb #\ \hfil
   &\vrule \hfil \ \ \frb #\ 
   &\vrule \hfil \ \ \frb #\ \vrule \cr%
\noalign{\hrule}
 & &9.7.11.19.53.811&491&92& & &3.5.17.37.191.317&319&636 \cr
5941&565957161&8.3.23.491.811&371&440&5959&571260945&8.9.11.17.29.37.53&1585&2078 \cr
 & &128.5.7.11.53.491&491&320& & &32.5.53.317.1039&1039&848 \cr
\noalign{\hrule}
 & &9.25.13.19.10193&25853&25112& & &81.5.23.83.739&1813&1882 \cr
5942&566475975&16.3.5.43.73.103.251&27227&27742&5960&571354155&4.27.49.37.83.941&1955&23452 \cr
 & &64.11.13.19.43.97.1433&45881&45856& & &32.5.7.11.13.17.23.41&3731&2992 \cr
\noalign{\hrule}
 & &9.5.169.23.41.79&319&526& & &11.13.19.23.41.223&813&5050 \cr
5943&566549685&4.11.29.41.79.263&59267&59346&5961&571355213&4.3.25.23.101.271&39&62 \cr
 & &16.27.7.13.29.47.97.157&59073&59032& & &16.9.25.13.31.271&6975&2168 \cr
\noalign{\hrule}
 & &9.49.121.13.19.43&851&722& & &5.169.37.101.181&57&44 \cr
5944&566747181&4.3.49.6859.23.37&6149&710&5962&571555465&8.3.5.11.13.19.37.181&207&2198 \cr
 & &16.5.11.13.23.43.71&355&184& & &32.27.7.19.23.157&25277&8208 \cr
\noalign{\hrule}
 & &9.125.7.17.19.223&2501&1736& & &9.13.41.101.1181&209&1390 \cr
5945&567228375&16.25.49.31.41.61&1863&638&5963&572190957&4.3.5.11.19.101.139&221&82 \cr
 & &64.81.11.23.29.31&6417&10208& & &16.5.11.13.17.19.41&1045&136 \cr
\noalign{\hrule}
 & &9.17.23.359.449&6149&2108& & &3.13.41.359.997&301&1298 \cr
5946&567231129&8.11.13.289.31.43&135&424&5964&572318877&4.7.11.43.59.359&1835&1476 \cr
 & &128.27.5.11.31.53&8215&2112& & &32.9.5.41.59.367&1835&2832 \cr
\noalign{\hrule}
 & &5.31.43.53.1607&111&154& & &3.19.31.227.1427&22005&22232 \cr
5947&567664715&4.3.7.11.31.37.1607&695&912&5965&572382543&16.81.5.7.19.163.397&7&88 \cr
 & &128.9.5.11.19.37.139&23769&26048& & &256.49.11.163.397&64711&68992 \cr
\noalign{\hrule}
 & &7.121.19.47.751&1525&774& & &5.11.13.23.61.571&1577&1278 \cr
5948&568034621&4.9.25.7.43.47.61&143&284&5966&572795795&4.9.11.19.61.71.83&211&460 \cr
 & &32.3.25.11.13.43.71&5325&8944& & &32.3.5.19.23.71.211&4047&3376 \cr
\noalign{\hrule}
 & &9.25.7.121.19.157&1003&992& & &3.5.11.17.23.83.107&2869&2858 \cr
5949&568485225&64.3.5.11.17.31.59.157&161&4&5967&572957715&4.5.17.19.107.151.1429&20853&13708 \cr
 & &512.7.17.23.31.59&12121&15104& & &32.9.7.23.149.151.331&49981&50064 \cr
\noalign{\hrule}
 & &243.5.7.11.59.103&60301&64844& & &27.13.19.23.37.101&4325&2002 \cr
5950&568533735&8.13.29.43.47.1283&1323&40&5968&573207219&4.3.25.7.11.169.173&851&332 \cr
 & &128.27.5.49.13.43&301&832& & &32.25.11.23.37.83&913&400 \cr
\noalign{\hrule}
 & &27.125.7.13.17.109&953&1772& & &11.529.29.43.79&97&570 \cr
5951&569102625&8.3.5.17.443.953&1133&1082&5969&573247147&4.3.5.19.23.79.97&957&860 \cr
 & &32.11.103.541.953&98159&95216& & &32.9.25.11.19.29.43&225&304 \cr
\noalign{\hrule}
 & &5.7.23.41.61.283&37103&44118& & &3.7.121.13.17.1021&745&2318 \cr
5952&569765315&4.27.11.19.43.3373&9553&566&5970&573353781&4.5.7.17.19.61.149&1021&594 \cr
 & &16.9.41.233.283&233&72& & &16.27.11.149.1021&149&72 \cr
\noalign{\hrule}
 & &3.7.17.43.137.271&195&76& & &5.7.11.13.19.37.163&365&2484 \cr
5953&569936577&8.9.5.13.19.43.137&781&2014&5971&573517945&8.27.25.19.23.73&17303&17372 \cr
 & &32.11.361.53.71&25631&9328& & &64.9.1331.13.43.101&12221&12384 \cr
\noalign{\hrule}
 & &5.11.13.47.71.239&297&58& & &25.41.631.887&26071&10296 \cr
5954&570243245&4.27.121.13.29.47&239&850&5972&573689425&16.9.11.13.841.31&1025&184 \cr
 & &16.3.25.17.29.239&85&696& & &256.3.25.11.23.41&759&128 \cr
\noalign{\hrule}
 & &125.7.11.13.47.97&1161&1114& & &11.19.761.3607&7235&7224 \cr
5955&570444875&4.27.5.11.43.97.557&161&5174&5973&573689743&16.3.5.7.43.1447.3607&5421&1814 \cr
 & &16.3.7.13.23.43.199&2967&1592& & &64.9.7.13.43.139.907&376551&377312 \cr
\noalign{\hrule}
 & &27.5.11.13.17.37.47&1163&58& & &9.5.11.23.29.37.47&11875&10136 \cr
5956&570713715&4.9.29.47.1163&793&370&5974&574156935&16.3.3125.7.19.181&6721&3596 \cr
 & &16.5.13.29.37.61&29&488& & &128.7.11.13.29.31.47&403&448 \cr
\noalign{\hrule}
 & &5.289.29.53.257&38243&29862& & &5.13.19.31.43.349&67&282 \cr
5957&570788005&4.27.7.79.167.229&119&110&5975&574542995&4.3.13.19.31.47.67&11&600 \cr
 & &16.3.5.49.11.17.79.167&13193&12936& & &64.9.25.11.67&99&10720 \cr
\noalign{\hrule}
 & &9.5.7.11.37.61.73&23&208& & &3.13.19.31.131.191&4205&1716 \cr
5958&570896865&32.3.13.23.61.73&2915&2122&5976&574757391&8.9.5.11.169.841&7945&8786 \cr
 & &128.5.11.53.1061&1061&3392& & &32.25.7.23.191.227&4025&3632 \cr
\noalign{\hrule}
}%
}
$$
\eject
\vglue -23 pt
\noindent\hskip 1 in\hbox to 6.5 in{\ 5977 -- 6012 \hfill\fbd 574826967 -- 583358555\frb}
\vskip -9 pt
$$
\vbox{
\nointerlineskip
\halign{\strut
    \vrule \ \ \hfil \frb #\ 
   &\vrule \hfil \ \ \fbb #\frb\ 
   &\vrule \hfil \ \ \frb #\ \hfil
   &\vrule \hfil \ \ \frb #\ 
   &\vrule \hfil \ \ \frb #\ \ \vrule \hskip 2 pt
   &\vrule \ \ \hfil \frb #\ 
   &\vrule \hfil \ \ \fbb #\frb\ 
   &\vrule \hfil \ \ \frb #\ \hfil
   &\vrule \hfil \ \ \frb #\ 
   &\vrule \hfil \ \ \frb #\ \vrule \cr%
\noalign{\hrule}
 & &9.11.169.17.43.47&2005&16& & &27.5.121.529.67&8021&86 \cr
5977&574826967&32.5.11.13.401&129&272&5995&578961405&4.9.13.43.617&2783&2770 \cr
 & &1024.3.5.17.43&5&512& & &16.5.121.23.43.277&277&344 \cr
\noalign{\hrule}
 & &3.49.19.29.47.151&925&1378& & &25.169.23.59.101&371&396 \cr
5978&574835709&4.25.13.19.29.37.53&151&1386&5996&579065825&8.9.7.11.13.23.53.101&10207&12530 \cr
 & &16.9.5.7.11.37.151&111&440& & &32.3.5.49.59.173.179&8477&8592 \cr
\noalign{\hrule}
 & &27.5.43.97.1021&24211&3356& & &9.7.121.17.41.109&9025&25702 \cr
5979&574909785&8.11.31.71.839&4579&4650&5997&579142179&4.25.361.71.181&1127&222 \cr
 & &32.3.25.19.961.241&18259&19280& & &16.3.5.49.19.23.37&851&5320 \cr
\noalign{\hrule}
 & &81.49.13.71.157&20225&31372& & &9.11.19.23.59.227&1337&1160 \cr
5980&575151759&8.25.11.23.31.809&2379&1666&5998&579421359&16.3.5.7.19.23.29.191&13&1324 \cr
 & &32.3.5.49.11.13.17.61&1037&880& & &128.5.13.29.331&1655&24128 \cr
\noalign{\hrule}
 & &3.25.13.59.73.137&77&142& & &9.5.7.11.23.29.251&13&680 \cr
5981&575307525&4.5.7.11.59.71.137&79&216&5999&580099905&16.25.13.17.251&3027&3248 \cr
 & &64.27.7.11.71.79&7821&15904& & &512.3.7.29.1009&1009&256 \cr
\noalign{\hrule}
 & &9.23.43.59.1097&6205&3668& & &5.7.11.169.37.241&837&1022 \cr
5982&576099423&8.5.7.17.23.73.131&381&1298&6000&580184605&4.27.49.31.73.241&2873&704 \cr
 & &32.3.5.11.17.59.127&935&2032& & &512.3.11.169.17.31&527&768 \cr
\noalign{\hrule}
 & &9.5.7.11.17.97.101&7387&8066& & &19.31.419.2351&2805&5156 \cr
5983&577092285&4.5.11.37.83.89.109&69&476&6001&580205641&8.3.5.11.17.31.1289&381&908 \cr
 & &32.3.7.17.23.83.89&1909&1424& & &64.9.5.11.127.227&28829&15840 \cr
\noalign{\hrule}
 & &25.23.29.53.653&351&374& & &5.17.23.47.71.89&87&2 \cr
5984&577105075&4.27.11.13.17.53.653&8299&190&6002&580621315&4.3.23.29.47.71&1335&2002 \cr
 & &16.3.5.11.19.43.193&3667&11352& & &16.9.5.7.11.13.89&63&1144 \cr
\noalign{\hrule}
 & &729.25.11.2879&1489&1390& & &27.125.7.13.31.61&3053&2992 \cr
5985&577167525&4.81.125.139.1489&51617&68992&6003&580773375&32.9.25.7.11.17.43.71&7277&248 \cr
 & &1024.49.11.71.727&35623&36352& & &512.17.19.31.383&6511&4864 \cr
\noalign{\hrule}
 & &25.13.43.109.379&11907&11528& & &5.343.11.13.23.103&351&454 \cr
5986&577321225&16.243.5.49.11.13.131&6443&12862&6004&580985405&4.27.49.11.169.227&491&1030 \cr
 & &64.9.17.59.109.379&531&544& & &16.3.5.103.227.491&1473&1816 \cr
\noalign{\hrule}
 & &27.7.17.71.2531&7975&5444& & &9.125.7.17.43.101&1159&1034 \cr
5987&577379313&8.25.11.17.29.1361&927&434&6005&581419125&4.3.7.11.19.47.61.101&5627&2150 \cr
 & &32.9.25.7.11.31.103&3193&4400& & &16.25.17.43.47.331&331&376 \cr
\noalign{\hrule}
 & &9.11.13.37.67.181&117713&115234& & &17.19.529.41.83&11633&10056 \cr
5988&577475613&4.7.53.2221.8231&3005&5226&6006&581460401&16.3.17.419.11633&2255&9378 \cr
 & &16.3.5.7.13.53.67.601&3005&2968& & &64.27.5.11.41.521&5731&4320 \cr
\noalign{\hrule}
 & &9.7.13.19.137.271&18139&15700& & &9.7.11.41.59.347&689&1730 \cr
5989&577733247&8.25.7.11.17.97.157&1083&16&6007&581699349&4.3.5.7.11.13.53.173&41&118 \cr
 & &256.3.25.17.361&323&3200& & &16.5.13.41.59.173&173&520 \cr
\noalign{\hrule}
 & &3.7.13.19.23.29.167&385&356& & &9.7.23.47.83.103&3135&766 \cr
5990&577774743&8.5.49.11.23.89.167&1233&13630&6008&582212547&4.27.5.7.11.19.383&799&664 \cr
 & &32.9.25.29.47.137&3425&2256& & &64.17.47.83.383&383&544 \cr
\noalign{\hrule}
 & &23.103.239.1021&49833&55330& & &9.49.13.17.43.139&3035&9614 \cr
5991&578081011&4.9.5.49.11.113.503&5105&3596&6009&582524397&4.5.7.11.19.23.607&3961&3354 \cr
 & &32.3.25.7.29.31.1021&2697&2800& & &16.3.13.17.23.43.233&233&184 \cr
\noalign{\hrule}
 & &11.361.41.53.67&24675&5542& & &27.25.7.31.41.97&8987&12062 \cr
5992&578141861&4.3.25.7.17.47.163&1221&1550&6010&582531075&4.9.11.19.37.43.163&497&970 \cr
 & &16.9.625.11.31.37&5625&9176& & &16.5.7.19.37.71.97&703&568 \cr
\noalign{\hrule}
 & &7.19.29.31.47.103&297&1660& & &9.25.343.7549&3787&3762 \cr
5993&578823847&8.27.5.7.11.31.83&667&418&6011&582594075&4.81.2401.11.19.541&899&44720 \cr
 & &32.9.121.19.23.29&1089&368& & &128.5.11.13.29.31.43&17329&20416 \cr
\noalign{\hrule}
 & &5.49.11.169.31.41&1615&3624& & &5.13.61.167.881&4653&5534 \cr
5994&578883305&16.3.25.11.17.19.151&903&3772&6012&583358555&4.9.5.11.13.47.2767&2623&5678 \cr
 & &128.9.7.23.41.43&387&1472& & &16.3.11.17.43.61.167&473&408 \cr
\noalign{\hrule}
}%
}
$$
\eject
\vglue -23 pt
\noindent\hskip 1 in\hbox to 6.5 in{\ 6013 -- 6048 \hfill\fbd 583509591 -- 593774775\frb}
\vskip -9 pt
$$
\vbox{
\nointerlineskip
\halign{\strut
    \vrule \ \ \hfil \frb #\ 
   &\vrule \hfil \ \ \fbb #\frb\ 
   &\vrule \hfil \ \ \frb #\ \hfil
   &\vrule \hfil \ \ \frb #\ 
   &\vrule \hfil \ \ \frb #\ \ \vrule \hskip 2 pt
   &\vrule \ \ \hfil \frb #\ 
   &\vrule \hfil \ \ \fbb #\frb\ 
   &\vrule \hfil \ \ \frb #\ \hfil
   &\vrule \hfil \ \ \frb #\ 
   &\vrule \hfil \ \ \frb #\ \vrule \cr%
\noalign{\hrule}
 & &9.49.61.109.199&1551&8200& & &9.7.23.31.103.127&2329&1608 \cr
6013&583509591&16.27.25.11.41.47&671&436&6031&587585439&16.27.17.23.67.137&671&2480 \cr
 & &128.5.121.61.109&121&320& & &512.5.11.17.31.61&5185&2816 \cr
\noalign{\hrule}
 & &5.49.11.17.47.271&10249&4656& & &9.25.13.127.1583&671&2254 \cr
6014&583545655&32.3.7.37.97.277&825&2764&6032&588044925&4.49.11.23.61.127&109&780 \cr
 & &256.9.25.11.691&3455&1152& & &32.3.5.7.13.23.109&161&1744 \cr
\noalign{\hrule}
 & &5.11.17.47.97.137&3211&3228& & &9.5.37.151.2341&3059&8646 \cr
6015&583985105&8.3.5.11.169.19.97.269&799&15594&6033&588562515&4.27.7.11.19.23.131&629&1142 \cr
 & &32.9.17.19.23.47.113&2599&2736& & &16.17.37.131.571&2227&4568 \cr
\noalign{\hrule}
 & &9.625.49.13.163&447&1672& & &5.49.961.41.61&337&624 \cr
6016&584049375&16.27.25.11.19.149&1753&1078&6034&588847945&32.3.5.7.13.61.337&803&3162 \cr
 & &64.49.121.1753&1753&3872& & &128.9.11.17.31.73&1683&4672 \cr
\noalign{\hrule}
 & &5.49.121.17.19.61&7379&7566& & &25.47.167.3001&41437&33588 \cr
6017&584095435&4.3.11.13.19.47.97.157&1785&58&6035&588871225&8.27.11.311.3767&1417&2350 \cr
 & &16.9.5.7.13.17.29.47&1363&936& & &32.9.25.11.13.47.109&1417&1584 \cr
\noalign{\hrule}
 & &27.5.17.31.43.191&3619&3046& & &27.5.7.11.19.29.103&1429&1326 \cr
6018&584313885&4.9.7.11.17.47.1523&1603&80&6036&589947435&4.81.7.11.13.17.1429&377&190 \cr
 & &128.5.49.47.229&11221&3008& & &16.5.169.19.29.1429&1429&1352 \cr
\noalign{\hrule}
 & &243.5.11.17.31.83&7969&5396& & &3.11.17.19.23.29.83&795&1114 \cr
6019&584598465&8.13.17.19.71.613&5885&4536&6037&590092899&4.9.5.17.19.53.557&203&9266 \cr
 & &128.81.5.7.11.13.107&749&832& & &16.5.7.29.41.113&287&4520 \cr
\noalign{\hrule}
 & &5.7.17.19.59.877&481&396& & &5.11.19.47.61.197&697&462 \cr
6020&584954615&8.9.7.11.13.19.37.59&941&2062&6038&590214955&4.3.7.121.17.41.197&159&38 \cr
 & &32.3.37.941.1031&38147&45168& & &16.9.7.17.19.41.53&6307&2952 \cr
\noalign{\hrule}
 & &25.17.41.59.569&8073&25498& & &11.17.23.31.43.103&793&540 \cr
6021&584974675&4.27.11.13.19.23.61&349&410&6039&590522999&8.27.5.13.17.61.103&1333&418 \cr
 & &16.9.5.13.19.41.349&2223&2792& & &32.9.11.13.19.31.43&247&144 \cr
\noalign{\hrule}
 & &125.49.11.19.457&1153&222& & &5.121.13.19.59.67&2279&1674 \cr
6022&585017125&4.3.37.457.1153&805&348&6040&590716555&4.27.13.19.31.43.53&335&1342 \cr
 & &32.9.5.7.23.29.37&851&4176& & &16.9.5.11.31.61.67&279&488 \cr
\noalign{\hrule}
 & &27.5.13.17.67.293&431&1034& & &9.5.13.17.19.53.59&545&222 \cr
6023&585690885&4.3.11.13.17.47.431&185&614&6041&590862285&4.27.25.37.53.109&253&1178 \cr
 & &16.5.37.307.431&11359&3448& & &16.11.19.23.31.109&1199&5704 \cr
\noalign{\hrule}
 & &3.49.23.29.43.139&407&260& & &9.5.13.17.19.53.59&181&142 \cr
6024&586038873&8.5.11.13.37.43.139&485&1044&6042&590862285&4.3.5.53.59.71.181&221&44 \cr
 & &64.9.25.29.37.97&2425&3552& & &32.11.13.17.71.181&1991&1136 \cr
\noalign{\hrule}
 & &13.43.47.53.421&24293&1980& & &3.5.17.31.37.43.47&1195&396 \cr
6025&586229449&8.9.5.11.17.1429&689&740&6043&591112185&8.27.25.11.31.239&3367&4042 \cr
 & &64.3.25.11.13.37.53&1221&800& & &32.7.11.13.37.43.47&143&112 \cr
\noalign{\hrule}
 & &27.49.11.191.211&703&914& & &27.125.11.89.179&133&312 \cr
6026&586500453&4.9.19.37.191.457&1477&5590&6044&591438375&16.81.25.7.11.13.19&97&178 \cr
 & &16.5.7.13.19.43.211&817&520& & &64.7.13.19.89.97&1729&3104 \cr
\noalign{\hrule}
 & &5.7.37.313.1447&1003&1188& & &27.41.359.1489&191&1298 \cr
6027&586519745&8.27.11.17.59.1447&1225&222&6045&591747957&4.11.59.191.359&871&1230 \cr
 & &32.81.25.49.11.37&567&880& & &16.3.5.13.41.59.67&871&2360 \cr
\noalign{\hrule}
 & &3.19.29.89.3989&2821&1168& & &81.7.11.269.353&1387&496 \cr
6028&586849713&32.7.13.31.73.89&513&110&6046&592246809&32.19.31.73.353&6165&4778 \cr
 & &128.27.5.11.19.73&3285&704& & &128.9.5.137.2389&11945&8768 \cr
\noalign{\hrule}
 & &9.5.121.13.43.193&361&218& & &9.5.7.11.23.43.173&1521&3424 \cr
6029&587444715&4.3.5.11.361.43.109&193&1612&6047&592851105&64.81.7.169.107&5&86 \cr
 & &32.13.31.109.193&109&496& & &256.5.13.43.107&107&1664 \cr
\noalign{\hrule}
 & &27.125.49.11.17.19&1933&8308& & &9.25.11.31.71.109&6283&4508 \cr
6030&587577375&8.9.31.67.1933&665&1268&6048&593774775&8.49.23.31.61.103&239&2130 \cr
 & &64.5.7.19.31.317&317&992& & &32.3.5.49.71.239&239&784 \cr
\noalign{\hrule}
}%
}
$$
\eject
\vglue -23 pt
\noindent\hskip 1 in\hbox to 6.5 in{\ 6049 -- 6084 \hfill\fbd 593859915 -- 603034795\frb}
\vskip -9 pt
$$
\vbox{
\nointerlineskip
\halign{\strut
    \vrule \ \ \hfil \frb #\ 
   &\vrule \hfil \ \ \fbb #\frb\ 
   &\vrule \hfil \ \ \frb #\ \hfil
   &\vrule \hfil \ \ \frb #\ 
   &\vrule \hfil \ \ \frb #\ \ \vrule \hskip 2 pt
   &\vrule \ \ \hfil \frb #\ 
   &\vrule \hfil \ \ \fbb #\frb\ 
   &\vrule \hfil \ \ \frb #\ \hfil
   &\vrule \hfil \ \ \frb #\ 
   &\vrule \hfil \ \ \frb #\ \vrule \cr%
\noalign{\hrule}
 & &9.5.11.19.233.271&387&658& & &3.25.11.61.73.163&1157&668 \cr
6049&593859915&4.81.7.43.47.233&157&76&6067&598817175&8.11.13.61.89.167&1315&522 \cr
 & &32.7.19.43.47.157&7379&4816& & &32.9.5.29.89.263&7627&4272 \cr
\noalign{\hrule}
 & &27.7.19.251.659&253&260& & &9.11.23.271.971&451&520 \cr
6050&593983719&8.5.11.13.23.251.659&6511&738&6068&599172057&16.3.5.121.13.41.271&317&46 \cr
 & &32.9.5.13.17.41.383&15703&17680& & &64.5.13.23.41.317&4121&6560 \cr
\noalign{\hrule}
 & &9.25.343.43.179&8717&6032& & &27.7.11.229.1259&38413&29600 \cr
6051&594015975&32.3.5.13.23.29.379&407&28&6069&599398569&64.25.37.107.359&6467&2508 \cr
 & &256.7.11.13.23.37&9361&1664& & &512.3.11.19.29.223&6467&4864 \cr
\noalign{\hrule}
 & &25.13.23.101.787&697&798& & &9.11.13.23.47.431&991&560 \cr
6052&594165325&4.3.5.7.17.19.41.787&1111&324&6070&599627457&32.3.5.7.13.23.991&47&944 \cr
 & &32.243.11.17.19.101&4131&3344& & &1024.5.7.47.59&295&3584 \cr
\noalign{\hrule}
 & &9.5.11.17.241.293&503&262& & &27.25.7.11.83.139&71&104 \cr
6053&594208395&4.11.131.293.503&1863&1360&6071&599635575&16.9.13.71.83.139&1615&9922 \cr
 & &128.81.5.17.23.131&1179&1472& & &64.5.121.17.19.41&3553&1312 \cr
\noalign{\hrule}
 & &3.5.11.19.37.47.109&2419&2704& & &27.5.23.151.1279&14989&16268 \cr
6054&594242385&32.11.169.37.41.59&63&470&6072&599665545&8.3.5.49.13.83.1153&137&3322 \cr
 & &128.9.5.7.13.47.59&1239&832& & &32.11.83.137.151&913&2192 \cr
\noalign{\hrule}
 & &9.7.1627.5801&845&782& & &27.19.79.113.131&1075&1414 \cr
6055&594608301&4.5.169.17.23.5801&957&4844&6073&599921181&4.9.25.7.43.79.101&3971&574 \cr
 & &32.3.5.7.11.17.29.173&9515&7888& & &16.5.49.11.361.41&2695&6232 \cr
\noalign{\hrule}
 & &27.49.19.41.577&295&484& & &5.7.11.13.19.59.107&945&1088 \cr
6056&594666009&8.5.7.121.59.577&495&82&6074&600334735&128.27.25.49.17.59&583&642 \cr
 & &32.9.25.1331.41&1331&400& & &512.81.11.17.53.107&4293&4352 \cr
\noalign{\hrule}
 & &25.7.19.23.31.251&697&108& & &3.5.23.29.47.1277&30227&29792 \cr
6057&595051975&8.27.5.17.41.251&1027&1232&6075&600490095&64.49.19.23.167.181&1287&2554 \cr
 & &256.3.7.11.13.17.79&11297&6528& & &256.9.7.11.13.19.1277&2717&2688 \cr
\noalign{\hrule}
 & &27.7.11.467.613&5707&10844& & &9.11.23.31.67.127&85&116 \cr
6058&595157409&8.7.13.439.2711&1401&1310&6076&600624783&8.3.5.11.17.23.29.127&20501&13516 \cr
 & &32.3.5.131.439.467&2195&2096& & &64.13.19.31.83.109&9047&7904 \cr
\noalign{\hrule}
 & &27.5.11.13.67.461&5033&24932& & &25.13.23.37.41.53&87&938 \cr
6059&596273535&8.7.23.271.719&495&224&6077&600997475&4.3.7.13.29.53.67&723&814 \cr
 & &512.9.5.49.11.23&1127&256& & &16.9.11.37.67.241&2651&4824 \cr
\noalign{\hrule}
 & &9.5.7.11.13.17.19.41&43&412& & &9.5.11.31.197.199&1009&14 \cr
6060&596530935&8.11.17.19.43.103&405&728&6078&601570035&4.3.7.197.1009&1415&1612 \cr
 & &128.81.5.7.13.43&387&64& & &32.5.7.13.31.283&1981&208 \cr
\noalign{\hrule}
 & &81.25.7.13.41.79&163&242& & &5.11.29.31.43.283&5635&6534 \cr
6061&596866725&4.5.7.121.13.41.163&459&74&6079&601696205&4.27.25.49.1331.23&403&928 \cr
 & &16.27.11.17.37.163&1793&5032& & &256.9.7.13.23.29.31&1449&1664 \cr
\noalign{\hrule}
 & &3.73.79.179.193&7005&7084& & &3.25.11.361.43.47&2573&1398 \cr
6062&597697647&8.9.5.7.11.23.179.467&4565&6176&6080&601904325&4.9.31.43.83.233&1645&1924 \cr
 & &512.25.7.121.83.193&21175&21248& & &32.5.7.13.37.47.233&3367&3728 \cr
\noalign{\hrule}
 & &3.5.107.401.929&137&1066& & &9.7.1331.43.167&689&2020 \cr
6063&597909045&4.5.13.41.107.137&7227&7432&6081&602148393&8.5.13.53.101.167&33&134 \cr
 & &64.9.11.13.73.929&949&1056& & &32.3.5.11.13.53.67&4355&848 \cr
\noalign{\hrule}
 & &7.11.13.19.149.211&179&30& & &3.5.7.13.41.47.229&979&3956 \cr
6064&597938341&4.3.5.7.13.179.211&745&1998&6082&602351295&8.11.23.41.43.89&2061&1588 \cr
 & &16.81.25.37.149&2025&296& & &64.9.23.229.397&1191&736 \cr
\noalign{\hrule}
 & &27.5.7.11.13.43.103&391&82& & &9.13.19.29.47.199&2849&2650 \cr
6065&598512915&4.9.5.7.13.17.23.41&43&412&6083&602959851&4.25.7.11.19.29.37.53&47&2802 \cr
 & &32.17.23.43.103&17&368& & &16.3.5.47.53.467&467&2120 \cr
\noalign{\hrule}
 & &9.5.11.29.179.233&1061&104& & &5.11.17.43.53.283&1073&342 \cr
6066&598703985&16.3.13.179.1061&2755&428&6084&603034795&4.9.11.19.29.37.53&663&344 \cr
 & &128.5.19.29.107&2033&64& & &64.27.13.17.37.43&999&416 \cr
\noalign{\hrule}
}%
}
$$
\eject
\vglue -23 pt
\noindent\hskip 1 in\hbox to 6.5 in{\ 6085 -- 6120 \hfill\fbd 603207855 -- 611127825\frb}
\vskip -9 pt
$$
\vbox{
\nointerlineskip
\halign{\strut
    \vrule \ \ \hfil \frb #\ 
   &\vrule \hfil \ \ \fbb #\frb\ 
   &\vrule \hfil \ \ \frb #\ \hfil
   &\vrule \hfil \ \ \frb #\ 
   &\vrule \hfil \ \ \frb #\ \ \vrule \hskip 2 pt
   &\vrule \ \ \hfil \frb #\ 
   &\vrule \hfil \ \ \fbb #\frb\ 
   &\vrule \hfil \ \ \frb #\ \hfil
   &\vrule \hfil \ \ \frb #\ 
   &\vrule \hfil \ \ \frb #\ \vrule \cr%
\noalign{\hrule}
 & &9.5.17.37.101.211&8393&10292& & &9.7.23.37.43.263&32461&21152 \cr
6085&603207855&8.7.11.17.31.83.109&5&114&6103&606309417&64.11.13.227.661&1145&1806 \cr
 & &32.3.5.11.19.31.83&6479&1328& & &256.3.5.7.11.43.229&1145&1408 \cr
\noalign{\hrule}
 & &27.5.13.17.73.277&7385&3784& & &121.23.359.607&3825&4432 \cr
6086&603293535&16.3.25.7.11.43.211&73&202&6104&606451879&32.9.25.121.17.277&113&718 \cr
 & &64.7.73.101.211&1477&3232& & &128.3.5.17.113.359&1695&1088 \cr
\noalign{\hrule}
 & &5.7.139.163.761&477&338& & &9.67.907.1109&2117&1210 \cr
6087&603469195&4.9.7.169.53.761&1733&550&6105&606535389&4.3.5.121.29.67.73&907&4922 \cr
 & &16.3.25.11.53.1733&8665&13992& & &16.11.23.107.907&253&856 \cr
\noalign{\hrule}
 & &9.49.97.103.137&429&292& & &9.5.7.23.31.37.73&5311&5384 \cr
6088&603626247&8.27.7.11.13.73.97&515&164&6106&606631095&16.3.7.37.47.113.673&253&4964 \cr
 & &64.5.11.41.73.103&2993&1760& & &128.11.17.23.73.113&1243&1088 \cr
\noalign{\hrule}
 & &9.5.11.13.101.929&2537&2108& & &9.25.49.17.41.79&143&568 \cr
6089&603789615&8.3.17.31.43.59.101&1661&3490&6107&607069575&16.49.11.13.41.71&969&1040 \cr
 & &32.5.11.43.151.349&6493&5584& & &512.3.5.11.169.17.19&3211&2816 \cr
\noalign{\hrule}
 & &3.25.7.17.31.37.59&29&22& & &9.5.41.283.1163&527&322 \cr
6090&603981525&4.25.11.29.31.37.59&1377&452&6108&607243005&4.3.7.17.23.31.1163&3113&376 \cr
 & &32.81.11.17.29.113&3277&4752& & &64.11.31.47.283&341&1504 \cr
\noalign{\hrule}
 & &81.11.23.41.719&23933&34306& & &3.125.11.19.23.337&83&292 \cr
6091&604113147&4.7.13.17.263.1009&1205&2214&6109&607484625&8.23.73.83.337&671&1008 \cr
 & &16.27.5.7.17.41.241&1205&952& & &256.9.7.11.61.83&5063&2688 \cr
\noalign{\hrule}
 & &5.7.29.47.53.239&155&1518& & &9.13.17.23.37.359&375&16 \cr
6092&604279235&4.3.25.11.23.31.53&6931&7644&6110&607657401&32.27.125.13.37&493&506 \cr
 & &32.9.49.13.29.239&91&144& & &128.125.11.17.23.29&1375&1856 \cr
\noalign{\hrule}
 & &169.19.29.43.151&297&254& & &19.71.331.1361&509&840 \cr
6093&604621667&4.27.11.169.127.151&1591&70&6111&607712359&16.3.5.7.509.1361&935&426 \cr
 & &16.3.5.7.37.43.127&3885&1016& & &64.9.25.7.11.17.71&3825&2464 \cr
\noalign{\hrule}
 & &81.13.361.37.43&575&242& & &3.5.41.757.1307&481&276 \cr
6094&604791603&4.9.25.121.13.19.23&511&74&6112&608480385&8.9.13.23.37.1307&757&550 \cr
 & &16.5.7.121.37.73&605&4088& & &32.25.11.13.37.757&715&592 \cr
\noalign{\hrule}
 & &9.13.19.31.67.131&3277&5500& & &81.11.19.157.229&2751&232 \cr
6095&604849401&8.125.11.29.31.113&373&402&6113&608648337&16.243.7.29.131&785&916 \cr
 & &32.3.5.11.67.113.373&6215&5968& & &128.5.29.157.229&145&64 \cr
\noalign{\hrule}
 & &28561.59.359&24871&3690& & &81.121.23.37.73&1399&3078 \cr
6096&604950541&4.9.5.7.11.17.19.41&169&118&6114&608867523&4.6561.19.1399&16571&10010 \cr
 & &16.3.5.11.169.19.59&95&264& & &16.5.7.11.13.73.227&1135&728 \cr
\noalign{\hrule}
 & &9.5.193.227.307&1285&1478& & &3.7.13.37.47.1283&435&176 \cr
6097&605248965&4.25.227.257.739&38407&19932&6115&609100401&32.9.5.11.29.1283&2077&794 \cr
 & &32.3.11.151.193.199&2189&2416& & &128.5.31.67.397&26599&9920 \cr
\noalign{\hrule}
 & &9.5.13.19.107.509&103&638& & &729.29.151.191&24991&3850 \cr
6098&605356245&4.3.11.29.103.509&5681&9080&6116&609727581&4.25.7.11.67.373&1359&986 \cr
 & &64.5.13.19.23.227&227&736& & &16.9.5.11.17.29.151&55&136 \cr
\noalign{\hrule}
 & &9.5.7.19.23.53.83&13&1232& & &9.11.13.17.61.457&125&1162 \cr
6099&605544345&32.3.49.11.13.19&15355&15368&6117&609920883&4.125.7.83.457&187&270 \cr
 & &512.5.17.37.83.113&4181&4352& & &16.27.625.7.11.17&625&168 \cr
\noalign{\hrule}
 & &3.49.17.29.61.137&3469&3518& & &27.25.11.13.71.89&24961&35114 \cr
6100&605640147&4.29.61.1759.3469&3349&103950&6118&609941475&4.97.109.181.229&169&60 \cr
 & &16.27.25.7.11.17.197&1773&2200& & &32.3.5.169.97.181&2353&1552 \cr
\noalign{\hrule}
 & &9.343.29.67.101&473&130& & &5.289.19.37.601&615&14 \cr
6101&605802141&4.5.11.13.29.43.101&497&816&6119&610516835&4.3.25.7.17.19.41&4181&3894 \cr
 & &128.3.5.7.17.43.71&6035&2752& & &16.9.11.37.59.113&1243&4248 \cr
\noalign{\hrule}
 & &3.5.11.41.101.887&6097&7208& & &3.25.7.11.23.43.107&163&912 \cr
6102&606056055&16.7.13.17.41.53.67&693&4&6120&611127825&32.9.11.19.23.163&2863&1070 \cr
 & &128.9.49.11.67&147&4288& & &128.5.7.107.409&409&64 \cr
\noalign{\hrule}
}%
}
$$
\eject
\vglue -23 pt
\noindent\hskip 1 in\hbox to 6.5 in{\ 6121 -- 6156 \hfill\fbd 611206695 -- 618599007\frb}
\vskip -9 pt
$$
\vbox{
\nointerlineskip
\halign{\strut
    \vrule \ \ \hfil \frb #\ 
   &\vrule \hfil \ \ \fbb #\frb\ 
   &\vrule \hfil \ \ \frb #\ \hfil
   &\vrule \hfil \ \ \frb #\ 
   &\vrule \hfil \ \ \frb #\ \ \vrule \hskip 2 pt
   &\vrule \ \ \hfil \frb #\ 
   &\vrule \hfil \ \ \fbb #\frb\ 
   &\vrule \hfil \ \ \frb #\ \hfil
   &\vrule \hfil \ \ \frb #\ 
   &\vrule \hfil \ \ \frb #\ \vrule \cr%
\noalign{\hrule}
 & &27.5.121.17.31.71&1033&1168& & &49.29.31.61.229&4815&1826 \cr
6121&611206695&32.121.17.73.1033&1065&992&6139&615348419&4.9.5.11.31.83.107&145&176 \cr
 & &2048.3.5.31.71.1033&1033&1024& & &128.3.25.121.29.83&6225&7744 \cr
\noalign{\hrule}
 & &3.7.13.19.191.617&47&86& & &3.25.7.13.19.47.101&8767&8908 \cr
6122&611272389&4.43.47.191.617&4797&4180&6140&615567225&8.11.13.17.19.131.797&1795&16938 \cr
 & &32.9.5.11.13.19.41.43&2255&2064& & &32.9.5.17.359.941&15997&17232 \cr
\noalign{\hrule}
 & &13.19.23.37.41.71&14873&13950& & &9.7.17.29.43.461&221&682 \cr
6123&611883467&4.9.25.23.31.107.139&11&104&6141&615682557&4.3.11.13.289.29.31&3227&530 \cr
 & &64.3.5.11.13.107.139&7645&10272& & &16.5.7.11.53.461&53&440 \cr
\noalign{\hrule}
 & &9.25.7.121.169.19&4307&82& & &27.5.11.13.19.23.73&331&34 \cr
6124&611936325&4.3.11.41.59.73&13195&13414&6142&615848805&4.13.17.19.23.331&225&212 \cr
 & &16.5.7.13.19.29.353&353&232& & &32.9.25.17.53.331&5627&4240 \cr
\noalign{\hrule}
 & &27.5.49.67.1381&951&2332& & &9.5.289.127.373&415&704 \cr
6125&612066105&8.81.5.11.53.317&229&176&6143&616059855&128.3.25.11.83.127&3811&2414 \cr
 & &256.121.229.317&38357&29312& & &512.17.37.71.103&7313&9472 \cr
\noalign{\hrule}
 & &27.7.11.79.3727&19669&21328& & &9.5.11.13.17.43.131&661&518 \cr
6126&612126207&32.9.13.17.31.43.89&3727&100&6144&616222035&4.5.7.17.37.43.661&1641&9596 \cr
 & &256.25.17.3727&425&128& & &32.3.7.547.2399&16793&8752 \cr
\noalign{\hrule}
 & &9.11.19.43.67.113&32453&20996& & &9.25.11.41.59.103&437&212 \cr
6127&612365193&8.17.23.29.83.181&105&286&6145&616663575&8.19.23.41.53.103&2065&108 \cr
 & &32.3.5.7.11.13.29.83&5395&3248& & &64.27.5.7.23.59&483&32 \cr
\noalign{\hrule}
 & &9.7.13.19.23.29.59&759&1880& & &11.13.23.277.677&19467&20144 \cr
6128&612372033&16.27.5.11.529.47&413&116&6146&616782881&32.27.7.23.103.1259&12035&14404 \cr
 & &128.5.7.29.47.59&235&64& & &256.9.5.13.29.83.277&3735&3712 \cr
\noalign{\hrule}
 & &3.11.23.29.61.457&185&272& & &9.7.11.797.1117&643&9410 \cr
6129&613600647&32.5.11.17.23.37.61&6129&5278&6147&616942557&4.5.7.643.941&2721&1780 \cr
 & &128.27.5.7.13.29.227&7945&7488& & &32.3.25.89.907&2225&14512 \cr
\noalign{\hrule}
 & &125.11.13.17.43.47&211&24& & &3.5.7.121.31.1567&799&1016 \cr
6130&614131375&16.3.25.13.43.211&531&544&6148&617170785&16.17.47.127.1567&847&720 \cr
 & &1024.27.17.59.211&12449&13824& & &512.9.5.7.121.17.47&799&768 \cr
\noalign{\hrule}
 & &3.121.19.29.37.83&41&70& & &9.71.83.103.113&2483&3410 \cr
6131&614239923&4.5.7.121.19.41.83&9657&7358&6149&617297643&4.5.11.13.31.113.191&3243&2678 \cr
 & &16.9.7.13.29.37.283&849&728& & &16.3.11.169.23.47.103&3887&4136 \cr
\noalign{\hrule}
 & &5.7.13.53.73.349&8993&9504& & &9.7.11.17.19.31.89&149&60 \cr
6132&614377855&64.27.5.11.13.17.529&9281&8176&6150&617571801&8.27.5.7.17.31.149&703&3916 \cr
 & &2048.9.7.73.9281&9281&9216& & &64.5.11.19.37.89&185&32 \cr
\noalign{\hrule}
 & &3.11.29.43.109.137&1365&4526& & &3.25.49.11.17.29.31&659&166 \cr
6133&614507883&4.9.5.7.11.13.31.73&145&548&6151&617815275&4.49.31.83.659&1089&430 \cr
 & &32.25.29.73.137&73&400& & &16.9.5.121.43.83&473&1992 \cr
\noalign{\hrule}
 & &9.5.11.19.223.293&359&136& & &9.125.11.13.23.167&4949&6574 \cr
6134&614513295&16.17.19.293.359&5901&920&6152&617920875&4.3.49.11.19.101.173&529&10 \cr
 & &256.3.5.7.23.281&6463&896& & &16.5.19.529.101&2323&152 \cr
\noalign{\hrule}
 & &243.49.11.13.361&641&290& & &25.169.23.6359&75241&71016 \cr
6135&614675061&4.9.5.11.19.29.641&5915&146&6153&617935825&16.3.11.67.269.1123&9573&8450 \cr
 & &16.25.7.169.73&73&2600& & &64.9.25.11.169.3191&3191&3168 \cr
\noalign{\hrule}
 & &5.23.41.283.461&34309&23706& & &5.29.1681.43.59&963&748 \cr
6136&615133045&4.27.11.439.3119&7093&2264&6154&618381065&8.9.11.17.1681.107&3569&1888 \cr
 & &64.9.41.173.283&173&288& & &512.3.11.43.59.83&913&768 \cr
\noalign{\hrule}
 & &121.22201.229&24955&2754& & &3.25.49.121.13.107&3127&2802 \cr
6137&615167509&4.81.5.7.17.23.31&149&242&6155&618542925&4.9.53.59.107.467&715&248 \cr
 & &16.27.5.7.121.149&189&40& & &64.5.11.13.31.53.59&1643&1888 \cr
\noalign{\hrule}
 & &27.11.71.163.179&17&196& & &9.13.23.47.67.73&1979&3520 \cr
6138&615255399&8.9.49.11.17.163&923&760&6156&618599007&128.5.11.73.1979&2997&1018 \cr
 & &128.5.49.13.19.71&3185&1216& & &512.81.37.509&4581&9472 \cr
\noalign{\hrule}
}%
}
$$
\eject
\vglue -23 pt
\noindent\hskip 1 in\hbox to 6.5 in{\ 6157 -- 6192 \hfill\fbd 618602699 -- 626242617\frb}
\vskip -9 pt
$$
\vbox{
\nointerlineskip
\halign{\strut
    \vrule \ \ \hfil \frb #\ 
   &\vrule \hfil \ \ \fbb #\frb\ 
   &\vrule \hfil \ \ \frb #\ \hfil
   &\vrule \hfil \ \ \frb #\ 
   &\vrule \hfil \ \ \frb #\ \ \vrule \hskip 2 pt
   &\vrule \ \ \hfil \frb #\ 
   &\vrule \hfil \ \ \fbb #\frb\ 
   &\vrule \hfil \ \ \frb #\ \hfil
   &\vrule \hfil \ \ \frb #\ 
   &\vrule \hfil \ \ \frb #\ \vrule \cr%
\noalign{\hrule}
 & &121.28561.179&3451&25110& & &7.17.37.353.401&1385&1422 \cr
6157&618602699&4.81.5.7.17.29.31&169&358&6175&623257859&4.9.5.17.79.277.353&19943&1940 \cr
 & &16.3.5.169.29.179&145&24& & &32.3.25.49.11.37.97&3201&2800 \cr
\noalign{\hrule}
 & &27.5.11.23.59.307&18473&18166& & &9.5.7.61.163.199&9737&2402 \cr
6158&618649515&4.5.49.11.13.29.31.293&109&2124&6176&623276955&4.49.13.107.1201&5185&10428 \cr
 & &32.9.7.59.109.293&2051&1744& & &32.3.5.11.17.61.79&869&272 \cr
\noalign{\hrule}
 & &9.25.7.19.23.29.31&1199&2924& & &5.7.169.19.31.179&2183&3366 \cr
6159&618759225&8.3.11.17.29.43.109&1271&2470&6177&623624365&4.9.5.11.17.19.37.59&31&26 \cr
 & &32.5.13.17.19.31.41&533&272& & &16.3.11.13.17.31.37.59&3009&3256 \cr
\noalign{\hrule}
 & &11.29.47.157.263&38259&45638& & &9.7.121.19.59.73&1105&16 \cr
6160&619075963&4.27.13.19.109.1201&1525&2726&6178&623812959&32.5.7.13.17.73&191&264 \cr
 & &16.9.25.19.29.47.61&1525&1368& & &512.3.11.17.191&191&4352 \cr
\noalign{\hrule}
 & &9.5.11.13.23.59.71&1123&824& & &11.17.19.29.73.83&12545&13932 \cr
6161&619992945&16.3.5.71.103.1123&23&332&6179&624301183&8.81.5.13.17.43.193&83&304 \cr
 & &128.23.83.1123&1123&5312& & &256.9.5.19.83.193&965&1152 \cr
\noalign{\hrule}
 & &9.7.11.23.167.233&323&1954& & &7.19.53.101.877&85&792 \cr
6162&620202429&4.17.19.167.977&405&572&6180&624379273&16.9.5.11.17.19.53&1247&1456 \cr
 & &32.81.5.11.13.17.19&1105&2736& & &512.3.5.7.13.29.43&5655&11008 \cr
\noalign{\hrule}
 & &9.11.17.239.1543&8563&8410& & &5.7.13.19.29.47.53&151&1386 \cr
6163&620651691&4.5.841.239.8563&7747&816&6181&624506155&4.9.49.11.47.151&925&1378 \cr
 & &128.3.5.17.29.61.127&8845&8128& & &16.3.25.11.13.37.53&111&440 \cr
\noalign{\hrule}
 & &3.7.11.41.173.379&831&380& & &27.17.41.89.373&7513&7780 \cr
6164&620985057&8.9.5.19.277.379&299&2194&6182&624734343&8.9.5.11.17.389.683&41&724 \cr
 & &32.13.19.23.1097&25231&3952& & &64.11.41.181.389&4279&5792 \cr
\noalign{\hrule}
 & &27.25.37.149.167&87&62& & &81.11.13.17.19.167&2303&1250 \cr
6165&621452925&4.81.29.31.37.167&4345&1834&6183&624798603&4.625.49.47.167&3239&936 \cr
 & &16.5.7.11.29.79.131&16037&11528& & &64.9.25.13.41.79&1025&2528 \cr
\noalign{\hrule}
 & &3.7.11.19.23.47.131&1591&2028& & &9.125.49.17.23.29&281&484 \cr
6166&621530679&8.9.169.37.43.131&1045&658&6184&625062375&8.25.7.121.23.281&1271&696 \cr
 & &32.5.7.11.13.19.37.47&185&208& & &128.3.121.29.31.41&3751&2624 \cr
\noalign{\hrule}
 & &3.25.23.557.647&1141&584& & &3.11.23.961.857&5251&5320 \cr
6167&621653775&16.7.73.163.647&405&242&6185&625094943&16.5.7.19.59.89.857&5625&374 \cr
 & &64.81.5.7.121.73&8833&6048& & &64.9.3125.11.17.19&9375&10336 \cr
\noalign{\hrule}
 & &9.13.59.113.797&15&782& & &25.7.13.19.29.499&1909&1584 \cr
6168&621691083&4.27.5.17.23.113&173&286&6186&625508975&32.9.11.19.23.29.83&697&260 \cr
 & &16.5.11.13.23.173&253&6920& & &256.3.5.13.17.41.83&4233&5248 \cr
\noalign{\hrule}
 & &5.13.19.31.37.439&187&252& & &9.5.7.11.13.17.19.43&381&854 \cr
6169&621863255&8.9.7.11.17.19.31.37&439&1142&6187&625630005&4.27.49.17.61.127&143&1180 \cr
 & &32.3.7.11.439.571&1713&1232& & &32.5.11.13.59.127&127&944 \cr
\noalign{\hrule}
 & &27.5.7.11.19.47.67&1207&1378& & &3.5.7.11.13.71.587&409&178 \cr
6170&621943245&4.3.7.13.17.53.67.71&2585&1178&6188&625780155&4.5.13.71.89.409&401&756 \cr
 & &16.5.11.13.17.19.31.47&221&248& & &32.27.7.401.409&3609&6544 \cr
\noalign{\hrule}
 & &5.179.787.883&68229&72644& & &9.49.13.23.47.101&105485&105908 \cr
6171&621954295&8.27.7.11.13.361.127&115&1766&6189&625934673&8.5.7.11.289.29.73.83&42159&254 \cr
 & &32.3.5.7.19.23.883&483&304& & &32.3.11.13.23.47.127&127&176 \cr
\noalign{\hrule}
 & &3.5.11.13.17.37.461&249&212& & &9.5.13.19.23.31.79&743&154 \cr
6172&621983505&8.9.5.11.13.17.53.83&23&2408&6190&626074605&4.3.5.7.11.79.743&3287&1058 \cr
 & &128.7.23.43.83&24983&1472& & &16.7.19.529.173&1211&184 \cr
\noalign{\hrule}
 & &3.5.17.67.83.439&25289&4124& & &7.251.593.601&2201&1950 \cr
6173&622526145&8.1331.19.1031&1665&634&6191&626182501&4.3.25.13.31.71.601&99&502 \cr
 & &32.9.5.11.37.317&951&6512& & &16.27.25.11.71.251&1775&2376 \cr
\noalign{\hrule}
 & &27.49.17.19.31.47&205&628& & &27.7.11.13.17.29.47&635&652 \cr
6174&622618353&8.3.5.19.31.41.157&3839&5110&6192&626242617&8.3.5.7.29.47.127.163&13801&27982 \cr
 & &32.25.7.11.73.349&8725&12848& & &32.5.17.37.373.823&30451&29840 \cr
\noalign{\hrule}
}%
}
$$
\eject
\vglue -23 pt
\noindent\hskip 1 in\hbox to 6.5 in{\ 6193 -- 6228 \hfill\fbd 627136059 -- 633965205\frb}
\vskip -9 pt
$$
\vbox{
\nointerlineskip
\halign{\strut
    \vrule \ \ \hfil \frb #\ 
   &\vrule \hfil \ \ \fbb #\frb\ 
   &\vrule \hfil \ \ \frb #\ \hfil
   &\vrule \hfil \ \ \frb #\ 
   &\vrule \hfil \ \ \frb #\ \ \vrule \hskip 2 pt
   &\vrule \ \ \hfil \frb #\ 
   &\vrule \hfil \ \ \fbb #\frb\ 
   &\vrule \hfil \ \ \frb #\ \hfil
   &\vrule \hfil \ \ \frb #\ 
   &\vrule \hfil \ \ \frb #\ \vrule \cr%
\noalign{\hrule}
 & &3.11.361.61.863&345&326& & &9.5.19.37.43.463&26105&25288 \cr
6193&627136059&4.9.5.19.23.163.863&1535&2398&6211&629821215&16.3.25.23.29.109.227&703&22 \cr
 & &16.25.11.109.163.307&33463&32600& & &64.11.19.23.37.109&1199&736 \cr
\noalign{\hrule}
 & &7.37.59.67.613&205&264& & &3.25.7.13.19.43.113&4249&4136 \cr
6194&627605951&16.3.5.11.37.41.613&2613&4130&6212&630090825&16.5.49.11.19.47.607&9&11524 \cr
 & &64.9.25.7.13.59.67&325&288& & &128.9.11.43.67&2211&64 \cr
\noalign{\hrule}
 & &9.49.11.109.1187&34477&23686& & &27.11.13.23.47.151&707&190 \cr
6195&627636933&4.13.23.911.1499&1205&294&6213&630234891&4.9.5.7.19.101.151&347&1012 \cr
 & &16.3.5.49.13.23.241&1205&2392& & &32.11.23.101.347&347&1616 \cr
\noalign{\hrule}
 & &9.7.17.19.109.283&341&58& & &27.25.11.13.31.211&157&184 \cr
6196&627705603&4.3.11.17.29.31.109&2429&950&6214&631370025&16.25.13.23.157.211&1617&3658 \cr
 & &16.25.7.11.19.347&347&2200& & &64.3.49.11.23.31.59&1357&1568 \cr
\noalign{\hrule}
 & &5.13.19.29.47.373&89&462& & &3.5.23.71.149.173&6067&6216 \cr
6197&627872765&4.3.5.7.11.13.47.89&461&696&6215&631407615&16.9.5.7.23.37.6067&121&86 \cr
 & &64.9.7.11.29.461&5071&2016& & &64.121.37.43.6067&192511&194144 \cr
\noalign{\hrule}
 & &3.11.23.457.1811&1285&526& & &7.961.37.43.59&8307&6974 \cr
6198&628168893&4.5.257.263.457&871&414&6216&631456763&4.9.11.13.31.71.317&43&360 \cr
 & &16.9.13.23.67.263&3419&1608& & &64.81.5.11.43.71&4455&2272 \cr
\noalign{\hrule}
 & &9.11.13.29.113.149&4025&296& & &3.7.13.17.23.61.97&165&262 \cr
6199&628407351&16.3.25.7.13.23.37&319&1124&6217&631598331&4.9.5.11.13.17.23.131&97&488 \cr
 & &128.5.11.29.281&1405&64& & &64.11.61.97.131&131&352 \cr
\noalign{\hrule}
 & &9.25.121.41.563&133&142& & &81.49.11.17.23.37&331&520 \cr
6200&628434675&4.7.11.19.41.71.563&190485&196678&6218&631614753&16.3.5.7.11.13.17.331&173&184 \cr
 & &16.27.5.17.29.83.3391&57647&57768& & &256.5.13.23.173.331&21515&22144 \cr
\noalign{\hrule}
 & &5.7.13.19.23.29.109&3043&1626& & &9.5.11.43.59.503&2291&2236 \cr
6201&628517435&4.3.5.17.19.179.271&1023&2378&6219&631674945&8.13.29.1849.59.79&69&1780 \cr
 & &16.9.11.17.29.31.41&3069&5576& & &64.3.5.13.23.79.89&7031&9568 \cr
\noalign{\hrule}
 & &27.5.23.47.59.73&36817&26962& & &5.7.17.19.29.41.47&63&22 \cr
6202&628542045&4.11.13.17.61.3347&1277&2070&6220&631757315&4.9.49.11.19.29.47&629&2050 \cr
 & &16.9.5.11.17.23.1277&1277&1496& & &16.3.25.11.17.37.41&111&440 \cr
\noalign{\hrule}
 & &5.7.17.47.113.199&117&682& & &125.121.529.79&971&846 \cr
6203&628848955&4.9.7.11.13.31.199&425&226&6221&632088875&4.9.121.23.47.971&65&2848 \cr
 & &16.3.25.11.13.17.113&65&264& & &256.3.5.13.47.89&1833&11392 \cr
\noalign{\hrule}
 & &81.125.23.37.73&6061&3064& & &3.7.23.199.6577&25885&6154 \cr
6204&628995375&16.11.19.23.29.383&7957&3150&6222&632161509&4.5.17.31.167.181&4225&1386 \cr
 & &64.9.25.7.73.109&109&224& & &16.9.125.7.11.169&5577&1000 \cr
\noalign{\hrule}
 & &9.5.17.19.73.593&781&188& & &7.19.23.31.59.113&4455&5044 \cr
6205&629205615&8.3.5.11.47.71.73&593&958&6223&632224943&8.81.5.11.13.97.113&185&1058 \cr
 & &32.71.479.593&479&1136& & &32.9.25.13.529.37&7475&5328 \cr
\noalign{\hrule}
 & &11.17.19.29.31.197&611&288& & &121.31.37.47.97&41&4518 \cr
6206&629246959&64.9.11.13.47.197&115&82&6224&632729933&4.9.31.41.251&141&110 \cr
 & &256.3.5.13.23.41.47&28905&38272& & &16.27.5.11.41.47&1107&40 \cr
\noalign{\hrule}
 & &5.19.43.223.691&1773&2464& & &11.169.61.5581&53909&59490 \cr
6207&629469905&64.9.5.7.11.43.197&493&1872&6225&632879819&4.9.5.31.37.47.661&1861&122 \cr
 & &2048.81.13.17.29&30537&17408& & &16.3.5.31.61.1861&1861&3720 \cr
\noalign{\hrule}
 & &7.11.13.19.79.419&425&444& & &5.11.169.17.19.211&8947&26712 \cr
6208&629547919&8.3.25.7.13.17.37.419&7979&7524&6226&633482135&16.9.7.23.53.389&15&38 \cr
 & &64.27.5.11.17.19.79.101&2727&2720& & &64.27.5.7.19.389&2723&864 \cr
\noalign{\hrule}
 & &25.17.43.131.263&1771&456& & &9.5.49.97.2963&93577&93092 \cr
6209&629628575&16.3.5.7.11.19.23.43&1651&3156&6227&633741255&8.7.11.17.1369.47.181&2963&114 \cr
 & &128.9.13.127.263&1651&576& & &32.3.19.37.47.2963&703&752 \cr
\noalign{\hrule}
 & &3.49.41.83.1259&249491&250750& & &3.5.83.271.1879&1617&262 \cr
6210&629803419&4.125.11.17.37.59.613&7&18&6228&633965205&4.9.49.11.83.131&6775&6194 \cr
 & &16.9.5.7.17.37.59.613&52105&52392& & &16.25.7.19.163.271&815&1064 \cr
\noalign{\hrule}
}%
}
$$
\eject
\vglue -23 pt
\noindent\hskip 1 in\hbox to 6.5 in{\ 6229 -- 6264 \hfill\fbd 634347015 -- 645174735\frb}
\vskip -9 pt
$$
\vbox{
\nointerlineskip
\halign{\strut
    \vrule \ \ \hfil \frb #\ 
   &\vrule \hfil \ \ \fbb #\frb\ 
   &\vrule \hfil \ \ \frb #\ \hfil
   &\vrule \hfil \ \ \frb #\ 
   &\vrule \hfil \ \ \frb #\ \ \vrule \hskip 2 pt
   &\vrule \ \ \hfil \frb #\ 
   &\vrule \hfil \ \ \fbb #\frb\ 
   &\vrule \hfil \ \ \frb #\ \hfil
   &\vrule \hfil \ \ \frb #\ 
   &\vrule \hfil \ \ \frb #\ \vrule \cr%
\noalign{\hrule}
 & &3.5.19.23.29.47.71&8383&7832& & &49.121.23.37.127&7391&12090 \cr
6229&634347015&16.11.71.83.89.101&7279&108&6247&640788533&4.3.5.7.13.19.31.389&1067&1656 \cr
 & &128.27.11.29.251&2761&576& & &64.27.5.11.13.23.97&2619&2080 \cr
\noalign{\hrule}
 & &3.5.19.431.5171&3663&1508& & &3.7.13.17.167.827&14707&17546 \cr
6230&635179785&8.27.11.13.19.29.37&10045&10342&6248&640963869&4.49.11.31.191.283&4509&4850 \cr
 & &32.5.49.13.41.5171&637&656& & &16.27.25.97.167.283&7075&6984 \cr
\noalign{\hrule}
 & &3.11.19.277.3659&1725&1934& & &3.125.49.11.19.167&12797&11752 \cr
6231&635491461&4.9.25.23.277.967&3659&5044&6249&641342625&16.25.13.67.113.191&171&2654 \cr
 & &32.5.13.23.97.3659&1495&1552& & &64.9.19.67.1327&3981&2144 \cr
\noalign{\hrule}
 & &3.5.17.67.167.223&83&418& & &5.7.11.733.2273&19021&21294 \cr
6232&636262485&4.11.17.19.83.223&1413&2824&6250&641451965&4.9.49.169.23.827&95&732 \cr
 & &64.9.11.157.353&3883&15072& & &32.27.5.13.19.23.61&15067&9936 \cr
\noalign{\hrule}
 & &27.7.11.19.89.181&8401&7708& & &3.5.49.31.47.599&429&1028 \cr
6233&636321609&8.3.19.31.41.47.271&493&400&6251&641466105&8.9.5.49.11.13.257&2773&568 \cr
 & &256.25.17.29.41.271&133603&131200& & &128.11.47.59.71&781&3776 \cr
\noalign{\hrule}
 & &5.19.79.137.619&819&682& & &27.7.11.23.29.463&247&710 \cr
6234&636446515&4.9.5.7.11.13.31.619&1781&76&6252&642038859&4.9.5.7.13.19.23.71&187&250 \cr
 & &32.3.7.169.19.137&507&112& & &16.625.11.13.17.71&10625&7384 \cr
\noalign{\hrule}
 & &121.29.167.1087&25865&5658& & &9.25.121.103.229&2603&5122 \cr
6235&636985261&4.3.5.7.23.41.739&1087&348&6253&642156075&4.3.11.13.19.137.197&1049&458 \cr
 & &32.9.23.29.1087&23&144& & &16.13.19.229.1049&1049&1976 \cr
\noalign{\hrule}
 & &9.125.121.31.151&3337&3458& & &11.31.1291.1459&6371&7830 \cr
6236&637201125&4.25.7.13.19.31.47.71&25469&3144&6254&642297029&4.27.5.23.29.31.277&7295&1292 \cr
 & &64.3.7.131.25469&25469&29344& & &32.3.25.17.19.1459&475&816 \cr
\noalign{\hrule}
 & &25.7.11.17.101.193&69&124& & &25.7.79.46489&23521&22968 \cr
6237&637908425&8.3.5.7.17.23.31.101&193&312&6255&642710425&16.9.25.11.29.43.547&23257&24332 \cr
 & &128.9.13.23.31.193&2691&1984& & &128.3.7.121.13.79.1789&23257&23232 \cr
\noalign{\hrule}
 & &9.7.11.59.67.233&25&674& & &3.11.29.37.67.271&503&570 \cr
6238&638286957&4.3.25.7.67.337&233&236&6256&642921213&4.9.5.11.19.271.503&185&86 \cr
 & &32.25.59.233.337&337&400& & &16.25.19.37.43.503&9557&8600 \cr
\noalign{\hrule}
 & &3.11.13.19.157.499&171&328& & &9.11.19.31.41.269&46391&52870 \cr
6239&638573793&16.27.11.13.361.41&2495&2198&6257&643112019&4.5.17.23.311.2017&31&1986 \cr
 & &64.5.7.41.157.499&205&224& & &16.3.31.311.331&311&2648 \cr
\noalign{\hrule}
 & &9.49.13.23.29.167&149&352& & &5.11.23.31.47.349&855&602 \cr
6240&638593137&64.3.7.11.13.23.149&425&334&6258&643243645&4.9.25.7.19.43.349&1703&4928 \cr
 & &256.25.17.149.167&2533&3200& & &512.3.49.11.13.131&6419&9984 \cr
\noalign{\hrule}
 & &3.11.47.239.1723&335&382& & &11.31.883.2137&1053&1084 \cr
6241&638697147&4.5.11.67.191.1723&189&1912&6259&643457111&8.81.11.13.271.883&1085&202 \cr
 & &64.27.5.7.67.239&603&1120& & &32.9.5.7.31.101.271&12195&11312 \cr
\noalign{\hrule}
 & &5.13.29.61.67.83&19&396& & &25.7.13.19.47.317&10649&4626 \cr
6242&639431585&8.9.11.19.61.67&6235&6026&6260&644009275&4.9.7.23.257.463&473&10 \cr
 & &32.3.5.23.29.43.131&2967&2096& & &16.3.5.11.43.257&257&11352 \cr
\noalign{\hrule}
 & &27.125.7.11.23.107&37&712& & &27.13.97.127.149&341&2278 \cr
6243&639552375&16.5.11.23.37.89&63&52&6261&644271381&4.11.17.31.67.127&97&30 \cr
 & &128.9.7.13.37.89&3293&832& & &16.3.5.11.17.31.97&155&1496 \cr
\noalign{\hrule}
 & &3.25.47.173.1049&589&3736& & &3.7.11.89.137.229&3495&3358 \cr
6244&639706425&16.19.31.47.467&495&962&6262&644997507&4.9.5.23.73.229.233&2473&7832 \cr
 & &64.9.5.11.13.19.37&1443&6688& & &64.11.73.89.2473&2473&2336 \cr
\noalign{\hrule}
 & &25.7.19.193.997&9559&9384& & &3.19.29.37.53.199&1067&470 \cr
6245&639799825&16.3.121.17.23.79.193&3&190&6263&645065067&4.5.11.19.37.47.97&1791&52 \cr
 & &64.9.5.11.19.23.79&1817&3168& & &32.9.5.11.13.199&33&1040 \cr
\noalign{\hrule}
 & &5.11.13.17.19.47.59&999&106& & &3.5.17.19.37.59.61&351&352 \cr
6246&640410485&4.27.11.37.53.59&9149&8500&6264&645174735&64.81.5.11.13.17.59.61&37&4978 \cr
 & &32.3.125.7.17.1307&3921&2800& & &256.11.13.19.37.131&1703&1408 \cr
\noalign{\hrule}
}%
}
$$
\eject
\vglue -23 pt
\noindent\hskip 1 in\hbox to 6.5 in{\ 6265 -- 6300 \hfill\fbd 645449805 -- 653495715\frb}
\vskip -9 pt
$$
\vbox{
\nointerlineskip
\halign{\strut
    \vrule \ \ \hfil \frb #\ 
   &\vrule \hfil \ \ \fbb #\frb\ 
   &\vrule \hfil \ \ \frb #\ \hfil
   &\vrule \hfil \ \ \frb #\ 
   &\vrule \hfil \ \ \frb #\ \ \vrule \hskip 2 pt
   &\vrule \ \ \hfil \frb #\ 
   &\vrule \hfil \ \ \fbb #\frb\ 
   &\vrule \hfil \ \ \frb #\ \hfil
   &\vrule \hfil \ \ \frb #\ 
   &\vrule \hfil \ \ \frb #\ \vrule \cr%
\noalign{\hrule}
 & &9.5.49.11.13.23.89&835&3526& & &9.5.49.11.17.19.83&169&188 \cr
6265&645449805&4.25.11.41.43.167&931&906&6283&650252295&8.3.5.7.11.169.47.83&1139&226 \cr
 & &16.3.49.19.41.43.151&6191&6536& & &32.13.17.47.67.113&7571&9776 \cr
\noalign{\hrule}
 & &17.961.39569&19769&19800& & &9.5.17.43.53.373&3103&3476 \cr
6266&646438753&16.9.25.11.17.31.53.373&27953&22&6284&650301255&8.5.11.29.53.79.107&27&292 \cr
 & &64.3.121.27953&27953&11616& & &64.27.73.79.107&5767&10272 \cr
\noalign{\hrule}
 & &5.11.53.463.479&213&266& & &25.13.361.23.241&783&422 \cr
6267&646479955&4.3.5.7.11.19.71.463&39&424&6285&650332475&4.27.5.13.23.29.211&817&5302 \cr
 & &64.9.13.19.53.71&1349&3744& & &16.9.11.19.43.241&387&88 \cr
\noalign{\hrule}
 & &23.53.59.89.101&63063&57710& & &9.7.11.53.89.199&12209&7492 \cr
6268&646497869&4.9.5.49.11.13.29.199&37&236&6286&650507319&8.7.29.421.1873&6345&6766 \cr
 & &32.3.5.7.11.29.37.59&3885&5104& & &32.27.5.17.29.47.199&2397&2320 \cr
\noalign{\hrule}
 & &9.25.1277.2251&5833&26092& & &3.25.7.13.17.71.79&223&132 \cr
6269&646768575&8.11.19.307.593&1985&1392&6287&650784225&8.9.5.11.17.79.223&4361&2354 \cr
 & &256.3.5.19.29.397&7543&3712& & &32.49.121.89.107&9523&13552 \cr
\noalign{\hrule}
 & &3.7.23.89.101.149&10115&3146& & &25.11.43.47.1171&693&478 \cr
6270&646911363&4.5.49.121.13.289&477&356&6288&650812525&4.9.5.7.121.47.239&1339&7026 \cr
 & &32.9.5.13.17.53.89&2067&1360& & &16.27.13.103.1171&351&824 \cr
\noalign{\hrule}
 & &27.7.11.23.29.467&2537&2600& & &27.25.11.127.691&8227&8918 \cr
6271&647585631&16.3.25.13.23.29.43.59&13327&12338&6289&651595725&4.5.343.11.13.19.433&5461&2766 \cr
 & &64.5.13.31.199.13327&413137&413920& & &16.3.7.13.43.127.461&3913&3688 \cr
\noalign{\hrule}
 & &529.43.71.401&3085&6138& & &3.7.13.53.73.617&1181&670 \cr
6272&647629837&4.9.5.11.23.31.617&71&688&6290&651698229&4.5.13.53.67.1181&22491&23672 \cr
 & &128.3.5.31.43.71&465&64& & &64.27.5.49.11.17.269&22865&22176 \cr
\noalign{\hrule}
 & &27.11.19.29.37.107&845&332& & &3.5.11.13.23.73.181&5919&5846 \cr
6273&647878473&8.5.169.29.37.83&963&1444&6291&651863355&4.9.11.23.37.79.1973&2675&19028 \cr
 & &64.9.5.13.361.107&247&160& & &32.25.37.67.71.107&37985&39664 \cr
\noalign{\hrule}
 & &25.7.11.13.19.29.47&237&92& & &27.343.17.41.101&2585&1556 \cr
6274&648072425&8.3.5.11.13.19.23.79&119&24&6292&651946617&8.9.5.11.17.47.389&41&194 \cr
 & &128.9.7.17.23.79&1817&9792& & &32.11.41.97.389&4279&1552 \cr
\noalign{\hrule}
 & &729.5.11.19.23.37&3791&4228& & &5.7.11.13.113.1153&4387&10602 \cr
6275&648296055&8.5.7.17.37.151.223&16929&16744&6293&652096445&4.9.7.19.31.41.107&625&646 \cr
 & &128.81.49.11.13.17.19.23&833&832& & &16.3.625.17.361.107&45125&43656 \cr
\noalign{\hrule}
 & &27.7.11.181.1723&30229&16292& & &9.7.19.41.47.283&1397&3980 \cr
6276&648363177&8.19.37.43.4073&2445&1628&6294&652773177&8.5.11.47.127.199&123&76 \cr
 & &64.3.5.11.1369.163&6845&5216& & &64.3.5.11.19.41.127&635&352 \cr
\noalign{\hrule}
 & &81.5.7.17.43.313&2761&2356& & &3.125.13.29.31.149&1331&606 \cr
6277&648656505&8.11.19.31.251.313&261&52&6295&653011125&4.9.5.1331.31.101&1363&32 \cr
 & &64.9.13.29.31.251&11687&8032& & &256.29.47.101&4747&128 \cr
\noalign{\hrule}
 & &5.7.11.13.103.1259&289&432& & &3.5.49.67.89.149&2563&2652 \cr
6278&649033385&32.27.5.289.1259&247&1012&6296&653037945&8.9.7.11.13.17.67.233&89&782 \cr
 & &256.3.11.13.17.19.23&1173&2432& & &32.289.23.89.233&5359&4624 \cr
\noalign{\hrule}
 & &25.7.17.167.1307&13959&15266& & &9.23.31.137.743&3355&3332 \cr
6279&649350275&4.27.11.289.47.449&1169&3770&6297&653192847&8.5.49.11.17.31.61.137&743&216 \cr
 & &16.3.5.7.13.29.47.167&611&696& & &128.27.5.7.11.61.743&2135&2112 \cr
\noalign{\hrule}
 & &47.137.281.359&209&72& & &3.5.49.17.61.857&14381&29326 \cr
6280&649559881&16.9.11.19.47.359&515&562&6298&653201115&4.11.31.43.73.197&221&252 \cr
 & &64.3.5.11.19.103.281&3399&3040& & &32.9.7.13.17.73.197&2847&3152 \cr
\noalign{\hrule}
 & &5.7.11.67.113.223&10757&14442& & &3.5.11.13.17.19.23.41&205&458 \cr
6281&650008205&4.3.7.29.31.83.347&11&2418&6299&653343405&4.25.19.1681.229&18837&23188 \cr
 & &16.9.11.13.961&117&7688& & &32.9.7.11.13.17.23.31&93&112 \cr
\noalign{\hrule}
 & &9.7.121.269.317&899&52& & &27.5.29.71.2351&217&2134 \cr
6282&650036079&8.3.13.29.31.269&605&202&6300&653495715&4.5.7.11.29.31.97&5187&4702 \cr
 & &32.5.121.29.101&101&2320& & &16.3.49.13.19.2351&247&392 \cr
\noalign{\hrule}
}%
}
$$
\eject
\vglue -23 pt
\noindent\hskip 1 in\hbox to 6.5 in{\ 6301 -- 6336 \hfill\fbd 653606343 -- 662265919\frb}
\vskip -9 pt
$$
\vbox{
\nointerlineskip
\halign{\strut
    \vrule \ \ \hfil \frb #\ 
   &\vrule \hfil \ \ \fbb #\frb\ 
   &\vrule \hfil \ \ \frb #\ \hfil
   &\vrule \hfil \ \ \frb #\ 
   &\vrule \hfil \ \ \frb #\ \ \vrule \hskip 2 pt
   &\vrule \ \ \hfil \frb #\ 
   &\vrule \hfil \ \ \fbb #\frb\ 
   &\vrule \hfil \ \ \frb #\ \hfil
   &\vrule \hfil \ \ \frb #\ 
   &\vrule \hfil \ \ \frb #\ \vrule \cr%
\noalign{\hrule}
 & &9.13.1993.2803&1001&992& & &27.41.719.827&3575&25904 \cr
6301&653606343&64.7.11.169.31.2803&1993&810&6319&658236591&32.25.11.13.1619&881&738 \cr
 & &256.81.5.11.31.1993&1395&1408& & &128.9.25.41.881&881&1600 \cr
\noalign{\hrule}
 & &9.5.11.13.37.41.67&25813&3502& & &7.23.29.1369.103&809&2178 \cr
6302&654046965&4.17.83.103.311&41&42&6320&658361683&4.9.7.121.23.809&481&1290 \cr
 & &16.3.7.17.41.103.311&5287&5768& & &16.27.5.11.13.37.43&1485&4472 \cr
\noalign{\hrule}
 & &3.5.19.83.139.199&91&506& & &7.31.41.43.1721&165&122 \cr
6303&654320955&4.7.11.13.19.23.139&193&54&6321&658404691&4.3.5.11.31.61.1721&85&1806 \cr
 & &16.27.7.11.23.193&14861&1656& & &16.9.25.7.11.17.43&225&1496 \cr
\noalign{\hrule}
 & &9.5.7.17.19.59.109&709&8556& & &13.23.41.223.241&261&38 \cr
6304&654322095&8.27.23.31.709&665&44&6322&658835437&4.9.19.29.41.241&715&474 \cr
 & &64.5.7.11.19.31&11&992& & &16.27.5.11.13.19.79&7505&2376 \cr
\noalign{\hrule}
 & &5.121.23.59.797&1577&5562& & &9.5.121.19.23.277&683&406 \cr
6305&654325045&4.27.19.23.83.103&22879&22132&6323&659111805&4.5.7.19.23.29.683&3483&1298 \cr
 & &32.3.11.137.167.503&22879&24144& & &16.81.11.29.43.59&1711&3096 \cr
\noalign{\hrule}
 & &9.5.49.13.289.79&11671&6536& & &9.5.107.367.373&1495&1862 \cr
6306&654450615&16.7.11.19.43.1061&3477&3950&6324&659130165&4.25.49.13.19.23.107&473&2202 \cr
 & &64.3.25.361.61.79&1805&1952& & &16.3.7.11.23.43.367&989&616 \cr
\noalign{\hrule}
 & &9.5.49.11.13.31.67&145&548& & &5.49.13.29.37.193&4969&4488 \cr
6307&654909255&8.25.7.29.67.137&1059&884&6325&659578465&16.3.5.11.17.29.4969&3717&1252 \cr
 & &64.3.13.17.137.353&6001&4384& & &128.27.7.11.59.313&17523&20032 \cr
\noalign{\hrule}
 & &11.13.19.43.71.79&1141&360& & &9.5.7.13.367.439&3569&9276 \cr
6308&655305079&16.9.5.7.13.43.163&115&158&6326&659757735&8.27.43.83.773&1507&734 \cr
 & &64.3.25.23.79.163&3749&2400& & &32.11.43.137.367&1507&688 \cr
\noalign{\hrule}
 & &25.13.17.31.43.89&829&504& & &7.23.41.199.503&351&152 \cr
6309&655469425&16.9.7.17.89.829&2145&3658&6327&660740297&16.27.7.13.19.23.41&503&440 \cr
 & &64.27.5.11.13.31.59&649&864& & &256.3.5.11.13.19.503&2145&2432 \cr
\noalign{\hrule}
 & &27.13.23.31.43.61&451&1440& & &9.25.7.17.23.29.37&1905&832 \cr
6310&656439849&64.243.5.11.13.41&145&388&6328&660780225&128.27.125.13.127&2527&902 \cr
 & &512.25.11.29.97&26675&7424& & &512.7.11.361.41&3971&10496 \cr
\noalign{\hrule}
 & &3.25.343.121.211&699&356& & &729.23.89.443&5459&4730 \cr
6311&656784975&8.9.5.121.89.233&703&98&6329&661072509&4.5.11.43.53.89.103&4347&4820 \cr
 & &32.49.19.37.233&703&3728& & &32.27.25.7.23.53.241&6025&5936 \cr
\noalign{\hrule}
 & &9.49.83.131.137&93&44& & &9.7.11.19.149.337&235&25714 \cr
6312&656914041&8.27.11.31.83.131&3151&910&6330&661154571&4.5.13.23.43.47&1011&1010 \cr
 & &32.5.7.11.13.23.137&299&880& & &16.3.25.13.23.101.337&2525&2392 \cr
\noalign{\hrule}
 & &121.289.19.23.43&7363&5064& & &11.13.17.29.83.113&205&288 \cr
6313&657102479&16.3.23.37.199.211&1615&6192&6331&661210121&64.9.5.11.13.41.113&1411&58 \cr
 & &512.27.5.17.19.43&135&256& & &256.3.5.17.29.83&15&128 \cr
\noalign{\hrule}
 & &9.125.11.13.61.67&163&38& & &3.25.11.13.17.19.191&271&356 \cr
6314&657496125&4.3.11.13.19.61.163&365&794&6332&661657425&8.5.13.89.191.271&17307&308 \cr
 & &16.5.73.163.397&11899&3176& & &64.27.7.11.641&4487&288 \cr
\noalign{\hrule}
 & &729.5.11.529.31&721&8& & &9.5.7.11.361.529&15221&12576 \cr
6315&657517905&16.5.7.11.23.103&1377&992&6333&661707585&64.27.31.131.491&8659&4598 \cr
 & &1024.81.17.31&17&512& & &256.7.121.19.1237&1237&1408 \cr
\noalign{\hrule}
 & &9.5.7.31.193.349&5207&6952& & &71.73.277.461&6993&26660 \cr
6316&657741105&16.11.31.41.79.127&349&3588&6334&661853551&8.27.5.7.31.37.43&1649&1606 \cr
 & &128.3.11.13.23.349&253&832& & &32.9.11.17.37.73.97&9603&10064 \cr
\noalign{\hrule}
 & &3.625.11.19.23.73&469&906& & &11.13.29.43.47.79&425&48 \cr
6317&657958125&4.9.5.7.67.73.151&1957&598&6335&662105873&32.3.25.17.47.79&221&174 \cr
 & &16.13.19.23.67.103&1339&536& & &128.9.5.13.289.29&1445&576 \cr
\noalign{\hrule}
 & &3.7.11.43.97.683&145&156& & &49.19.43.71.233&407&1224 \cr
6318&658071183&8.9.5.13.29.97.683&2257&6622&6336&662265919&16.9.7.11.17.37.71&261&520 \cr
 & &32.7.11.29.37.43.61&1073&976& & &256.81.5.13.17.29&11745&28288 \cr
\noalign{\hrule}
}%
}
$$
\eject
\vglue -23 pt
\noindent\hskip 1 in\hbox to 6.5 in{\ 6337 -- 6372 \hfill\fbd 662283783 -- 671552275\frb}
\vskip -9 pt
$$
\vbox{
\nointerlineskip
\halign{\strut
    \vrule \ \ \hfil \frb #\ 
   &\vrule \hfil \ \ \fbb #\frb\ 
   &\vrule \hfil \ \ \frb #\ \hfil
   &\vrule \hfil \ \ \frb #\ 
   &\vrule \hfil \ \ \frb #\ \ \vrule \hskip 2 pt
   &\vrule \ \ \hfil \frb #\ 
   &\vrule \hfil \ \ \fbb #\frb\ 
   &\vrule \hfil \ \ \frb #\ \hfil
   &\vrule \hfil \ \ \frb #\ 
   &\vrule \hfil \ \ \frb #\ \vrule \cr%
\noalign{\hrule}
 & &81.7.31.41.919&15905&12584& & &151.24649.179&25839&1190 \cr
6337&662283783&16.5.7.121.13.3181&527&3708&6355&666237821&4.81.5.7.11.17.29&157&302 \cr
 & &128.9.13.17.31.103&1751&832& & &16.3.7.11.151.157&11&168 \cr
\noalign{\hrule}
 & &3.5.59.71.83.127&99&28& & &27.29.31.97.283&3685&4522 \cr
6338&662343735&8.27.5.7.11.59.83&2921&1976&6356&666318123&4.5.7.11.17.19.67.97&1461&188 \cr
 & &128.11.13.19.23.127&2717&1472& & &32.3.5.7.11.47.487&17045&8272 \cr
\noalign{\hrule}
 & &11.31.661.2939&111231&114170& & &5.13.17.23.157.167&693&1478 \cr
6339&662453539&4.9.5.49.17.233.727&5653&9614&6357&666355885&4.9.7.11.17.23.739&3173&2000 \cr
 & &16.3.5.7.11.19.23.5653&45885&45224& & &128.3.125.11.19.167&825&1216 \cr
\noalign{\hrule}
 & &5.7.17.29.103.373&36051&25234& & &5.11.13.17.163.337&1805&3924 \cr
6340&662919845&4.3.11.31.37.61.197&1227&1030&6358&667686305&8.9.25.11.361.109&1173&898 \cr
 & &16.9.5.11.31.103.409&3069&3272& & &32.27.17.19.23.449&10327&8208 \cr
\noalign{\hrule}
 & &3.11.13.47.131.251&149&102& & &9.5.47.53.59.101&2639&2108 \cr
6341&662979603&4.9.11.13.17.131.149&1757&470&6359&667974105&8.5.7.13.17.29.31.53&243&22 \cr
 & &16.5.7.47.149.251&149&280& & &32.243.7.11.29.31&9207&3248 \cr
\noalign{\hrule}
 & &9.5.7.11.13.41.359&1513&1718& & &3.5.7.289.361.61&341&86 \cr
6342&663017355&4.7.11.13.17.89.859&675&8774&6360&668227245&4.11.17.361.31.43&3239&2898 \cr
 & &16.27.25.17.41.107&321&680& & &16.9.7.23.41.43.79&9717&7912 \cr
\noalign{\hrule}
 & &27.7.17.43.4799&9757&4640& & &27.25.7.11.13.23.43&95&256 \cr
6343&663025041&64.9.5.11.29.887&3653&782&6361&668242575&512.125.11.19.43&299&174 \cr
 & &256.13.17.23.281&6463&1664& & &2048.3.13.19.23.29&551&1024 \cr
\noalign{\hrule}
 & &81.25.11.13.29.79&2257&4142& & &5.11.13.29.103.313&1777&2292 \cr
6344&663416325&4.5.11.19.37.61.109&553&492&6362&668475665&8.3.11.29.191.1777&12543&7004 \cr
 & &32.3.7.37.41.79.109&4469&4144& & &64.9.17.37.103.113&4181&4896 \cr
\noalign{\hrule}
 & &3.25.11.23.59.593&423&170& & &3.5.11.17.23.43.241&1339&616 \cr
6345&663878325&4.27.125.17.47.59&2189&1186&6363&668568945&16.7.121.13.43.103&4275&8188 \cr
 & &16.11.47.199.593&199&376& & &128.9.25.19.23.89&1691&960 \cr
\noalign{\hrule}
 & &9.11.23.487.599&249&238& & &27.29.31.59.467&15125&12428 \cr
6346&664230501&4.27.7.17.23.83.599&11&610&6364&668793969&8.9.125.121.13.239&31&86 \cr
 & &16.5.7.11.17.61.83&2905&8296& & &32.25.11.31.43.239&5975&7568 \cr
\noalign{\hrule}
 & &25.7.121.13.19.127&129&256& & &7.113.271.3121&26235&4388 \cr
6347&664238575&512.3.5.11.13.19.43&381&854&6365&669020681&8.9.5.11.53.1097&565&532 \cr
 & &2048.9.7.61.127&549&1024& & &64.3.25.7.19.53.113&1425&1696 \cr
\noalign{\hrule}
 & &7.11.13.269.2467&2175&292& & &9.125.19.173.181&841&716 \cr
6348&664286623&8.3.25.11.13.29.73&2037&2110&6366&669315375&8.19.841.179.181&4325&924 \cr
 & &32.9.125.7.97.211&20467&18000& & &64.3.25.7.11.29.173&319&224 \cr
\noalign{\hrule}
 & &9.5.49.11.67.409&139&106& & &81.5.7.13.37.491&28303&11468 \cr
6349&664659765&4.3.53.67.139.409&2651&24752&6367&669544785&8.11.31.47.61.83&977&1890 \cr
 & &128.7.11.13.17.241&3133&1088& & &32.27.5.7.31.977&977&496 \cr
\noalign{\hrule}
 & &9.7.17.37.97.173&15029&11440& & &9.49.11.169.19.43&905&954 \cr
6350&664980687&32.5.49.11.13.19.113&519&2666&6368&669792123&4.81.5.19.43.53.181&2555&1738 \cr
 & &128.3.11.31.43.173&1333&704& & &16.25.7.11.73.79.181&14299&14600 \cr
\noalign{\hrule}
 & &3.11.13.73.79.269&25&244& & &11.13.19.31.73.109&11189&11592 \cr
6351&665517567&8.25.11.13.61.79&269&126&6369&670194239&16.9.7.23.67.73.167&589&1090 \cr
 & &32.9.5.7.61.269&915&112& & &64.3.5.7.19.31.67.109&469&480 \cr
\noalign{\hrule}
 & &81.25.7.151.311&16393&30568& & &5.7.11.13.29.31.149&225&2162 \cr
6352&665672175&16.169.97.3821&55&42&6370&670424755&4.9.125.23.29.47&1937&3938 \cr
 & &64.3.5.7.11.13.3821&3821&4576& & &16.3.11.13.149.179&179&24 \cr
\noalign{\hrule}
 & &9.25.7.169.41.61&1837&2006& & &9.25.13.17.97.139&1631&176 \cr
6353&665703675&4.25.11.17.41.59.167&2401&43524&6371&670442175&32.3.5.7.11.17.233&1009&776 \cr
 & &32.27.2401.13.31&1029&496& & &512.11.97.1009&1009&2816 \cr
\noalign{\hrule}
 & &9.5.13.59.101.191&427&528& & &25.17.529.29.103&41&534 \cr
6354&665828865&32.27.7.11.13.59.61&31&382&6372&671552275&4.3.23.41.89.103&3135&1088 \cr
 & &128.11.31.61.191&671&1984& & &512.9.5.11.17.19&171&2816 \cr
\noalign{\hrule}
}%
}
$$
\eject
\vglue -23 pt
\noindent\hskip 1 in\hbox to 6.5 in{\ 6373 -- 6408 \hfill\fbd 671553575 -- 682679745\frb}
\vskip -9 pt
$$
\vbox{
\nointerlineskip
\halign{\strut
    \vrule \ \ \hfil \frb #\ 
   &\vrule \hfil \ \ \fbb #\frb\ 
   &\vrule \hfil \ \ \frb #\ \hfil
   &\vrule \hfil \ \ \frb #\ 
   &\vrule \hfil \ \ \frb #\ \ \vrule \hskip 2 pt
   &\vrule \ \ \hfil \frb #\ 
   &\vrule \hfil \ \ \fbb #\frb\ 
   &\vrule \hfil \ \ \frb #\ \hfil
   &\vrule \hfil \ \ \frb #\ 
   &\vrule \hfil \ \ \frb #\ \vrule \cr%
\noalign{\hrule}
 & &25.49.11.19.43.61&297&1228& & &11.19.43.139.541&5209&5070 \cr
6373&671553575&8.27.121.43.307&3983&1220&6391&675813413&4.3.5.11.169.43.5209&3787&1422 \cr
 & &64.3.5.7.61.569&569&96& & &16.27.7.169.79.541&4563&4424 \cr
\noalign{\hrule}
 & &11.19.29.271.409&22389&34250& & &13.167.433.719&1445&726 \cr
6374&671795179&4.3.125.17.137.439&67&2262&6392&675890917&4.3.5.121.289.433&507&2672 \cr
 & &16.9.25.13.29.67&325&4824& & &128.9.11.169.167&143&576 \cr
\noalign{\hrule}
 & &27.29.37.139.167&11773&11440& & &3.5.7.23.107.2617&1067&1550 \cr
6375&672503823&32.3.5.11.13.29.61.193&15883&10286&6393&676245885&4.125.11.31.97.107&9&116 \cr
 & &128.5.7.37.139.2269&2269&2240& & &32.9.11.29.31.97&3007&15312 \cr
\noalign{\hrule}
 & &11.13.19.73.3391&14413&29670& & &3.5.7.11.29.61.331&37&268 \cr
6376&672574331&4.3.5.7.23.29.43.71&949&1110&6394&676297545&8.29.37.67.331&6039&3560 \cr
 & &16.9.25.13.37.43.73&1591&1800& & &128.9.5.11.61.89&267&64 \cr
\noalign{\hrule}
 & &3.25.139.173.373&1551&1924& & &25.7.11.17.23.29.31&239&36 \cr
6377&672714825&8.9.11.13.37.47.173&19685&556&6395&676654825&8.9.17.23.31.239&3509&3900 \cr
 & &64.5.31.127.139&127&992& & &64.27.25.121.13.29&297&416 \cr
\noalign{\hrule}
 & &169.67.103.577&8737&8670& & &3.5.11.13.17.67.277&249&28 \cr
6378&672937213&4.3.5.289.577.8737&9273&536&6396&676753935&8.9.5.7.11.67.83&1469&2216 \cr
 & &64.9.5.11.17.67.281&8415&8992& & &128.7.13.113.277&113&448 \cr
\noalign{\hrule}
 & &9.7.29.197.1871&43&44& & &5.13.37.43.79.83&61&18 \cr
6379&673408449&8.3.7.11.43.197.1871&23881&15410&6397&678092155&4.9.5.13.37.61.83&6083&688 \cr
 & &32.5.121.13.23.67.167&180895&179024& & &128.3.7.11.43.79&33&448 \cr
\noalign{\hrule}
 & &3.25.7.13.29.41.83&57&148& & &3.25.11.13.29.37.59&93&388 \cr
6380&673538775&8.9.5.19.29.37.83&1441&136&6398&678967575&8.9.5.11.29.31.97&697&202 \cr
 & &128.11.17.37.131&2227&26048& & &32.17.41.97.101&4141&26384 \cr
\noalign{\hrule}
 & &25.7.17.29.37.211&13&198& & &9.7.61.149.1187&437&10 \cr
6381&673548925&4.9.5.7.11.13.17.29&2627&3782&6399&679684509&4.3.5.19.23.1187&3871&2684 \cr
 & &16.3.31.37.61.71&4331&744& & &32.49.11.61.79&77&1264 \cr
\noalign{\hrule}
 & &7.11.13.337.1997&927&1070& & &7.23.29.41.53.67&605&2142 \cr
6382&673661989&4.9.5.7.103.107.337&239&1924&6400&679764379&4.9.5.49.121.17.23&3551&2378 \cr
 & &32.3.13.37.107.239&8843&5136& & &16.3.5.29.41.53.67&15&8 \cr
\noalign{\hrule}
 & &3.5.13.361.61.157&277&638& & &3.25.11.13.23.31.89&5253&8542 \cr
6383&674172915&4.11.13.29.157.277&1159&882&6401&680576325&4.9.5.17.103.4271&4453&182 \cr
 & &16.9.49.11.19.29.61&539&696& & &16.7.13.17.61.73&7259&584 \cr
\noalign{\hrule}
 & &11.17.37.41.2377&11861&28548& & &11.13.361.79.167&207&40 \cr
6384&674304983&8.9.13.29.61.409&985&1394&6402&681062239&16.9.5.11.19.23.79&833&668 \cr
 & &32.3.5.17.29.41.197&985&1392& & &128.3.49.17.23.167&2499&1472 \cr
\noalign{\hrule}
 & &121.17.23.53.269&2285&2016& & &3.5.7.23.29.71.137&137&298 \cr
6385&674512927&64.9.5.7.11.53.457&3341&6256&6403&681230445&4.71.18769.149&4095&14674 \cr
 & &2048.3.13.17.23.257&3341&3072& & &16.9.5.7.11.13.23.29&33&104 \cr
\noalign{\hrule}
 & &125.19.29.41.239&4653&2278& & &3.49.11.19.41.541&299&152 \cr
6386&674906125&4.9.11.17.41.47.67&717&20&6404&681466863&16.13.361.23.541&451&90 \cr
 & &32.27.5.47.239&27&752& & &64.9.5.11.13.23.41&195&736 \cr
\noalign{\hrule}
 & &3.13.19.61.109.137&35&22& & &3.13.47.163.2281&46375&53218 \cr
6387&674986533&4.5.7.11.61.109.137&8375&18&6405&681514899&4.125.7.11.41.53.59&2281&846 \cr
 & &16.9.625.67&67&15000& & &16.9.25.11.47.2281&75&88 \cr
\noalign{\hrule}
 & &3.13.23.43.83.211&45&1034& & &27.7.11.31.71.149&5623&26710 \cr
6388&675493923&4.27.5.11.47.211&3107&2590&6406&681805971&4.5.2671.5623&4147&1476 \cr
 & &16.25.7.13.37.239&1673&7400& & &32.9.5.11.13.29.41&2665&464 \cr
\noalign{\hrule}
 & &3.13.23.43.83.211&385&174& & &13.23.97.101.233&99&200 \cr
6389&675493923&4.9.5.7.11.23.29.83&481&2390&6407&682527599&16.9.25.11.97.233&29&262 \cr
 & &16.25.7.13.37.239&1673&7400& & &64.3.25.11.29.131&3799&26400 \cr
\noalign{\hrule}
 & &3.25.7.11.19.47.131&4833&3392& & &81.5.11.293.523&481&1004 \cr
6390&675576825&128.81.19.53.179&31&50&6408&682679745&8.3.13.37.251.293&523&230 \cr
 & &512.25.31.53.179&9487&7936& & &32.5.13.23.37.523&299&592 \cr
\noalign{\hrule}
}%
}
$$
\eject
\vglue -23 pt
\noindent\hskip 1 in\hbox to 6.5 in{\ 6409 -- 6444 \hfill\fbd 682753995 -- 693076043\frb}
\vskip -9 pt
$$
\vbox{
\nointerlineskip
\halign{\strut
    \vrule \ \ \hfil \frb #\ 
   &\vrule \hfil \ \ \fbb #\frb\ 
   &\vrule \hfil \ \ \frb #\ \hfil
   &\vrule \hfil \ \ \frb #\ 
   &\vrule \hfil \ \ \frb #\ \ \vrule \hskip 2 pt
   &\vrule \ \ \hfil \frb #\ 
   &\vrule \hfil \ \ \fbb #\frb\ 
   &\vrule \hfil \ \ \frb #\ \hfil
   &\vrule \hfil \ \ \frb #\ 
   &\vrule \hfil \ \ \frb #\ \vrule \cr%
\noalign{\hrule}
 & &27.5.49.121.853&5353&4030& & &9.13.37.271.587&39353&29326 \cr
6409&682753995&4.25.11.13.31.53.101&833&492&6427&688644333&4.11.23.29.31.43.59&793&540 \cr
 & &32.3.49.13.17.41.101&4141&3536& & &32.27.5.13.29.59.61&5307&4720 \cr
\noalign{\hrule}
 & &25.37.71.101.103&50369&43056& & &3.5.11.961.43.101&563&398 \cr
6410&683217025&32.9.11.13.19.23.241&1515&4166&6428&688647795&4.43.101.199.563&53847&58190 \cr
 & &128.27.5.101.2083&2083&1728& & &16.9.5.11.529.31.193&1587&1544 \cr
\noalign{\hrule}
 & &5.121.13.17.19.269&1377&6488& & &125.11.43.89.131&2051&1776 \cr
6411&683366255&16.81.289.811&261&550&6429&689338375&32.3.5.7.37.131.293&27&158 \cr
 & &64.729.25.11.29&3645&928& & &128.81.7.79.293&23733&35392 \cr
\noalign{\hrule}
 & &343.11.293.619&10971&11590& & &9.5.7.11.19.37.283&313&390 \cr
6412&684297691&4.9.5.49.19.23.53.61&149&586&6430&689358285&4.27.25.13.283.313&7733&92 \cr
 & &16.3.53.61.149.293&3233&3576& & &32.11.13.19.23.37&299&16 \cr
\noalign{\hrule}
 & &3.7.13.23.61.1787&95&5456& & &9.7.11.17.139.421&95&326 \cr
6413&684454953&32.5.11.19.23.31&203&234&6431&689412339&4.3.5.17.19.139.163&1177&1594 \cr
 & &128.9.5.7.11.13.29&1595&192& & &16.5.11.19.107.797&10165&6376 \cr
\noalign{\hrule}
 & &2187.121.13.199&1475&1112& & &3.5.7.31.191.1109&715&622 \cr
6414&684590049&16.729.25.59.139&17&712&6432&689470845&4.25.11.13.311.1109&36099&8374 \cr
 & &256.5.17.59.89&26255&2176& & &16.27.7.53.79.191&477&632 \cr
\noalign{\hrule}
 & &27.5.17.29.41.251&233&260& & &9.5.961.37.431&265&696 \cr
6415&684917505&8.25.13.41.233.251&2233&8058&6433&689628015&16.27.25.29.37.53&253&1178 \cr
 & &32.3.7.11.13.17.29.79&1027&1232& & &64.11.19.23.29.31&4807&928 \cr
\noalign{\hrule}
 & &3.7.37.59.67.223&1351&1128& & &9.25.11.169.17.97&227&1294 \cr
6416&684940263&16.9.49.47.59.193&17545&8474&6434&689735475&4.25.17.227.647&111&536 \cr
 & &64.5.121.19.29.223&3509&3040& & &64.3.37.67.227&15209&1184 \cr
\noalign{\hrule}
 & &3.5.11.13.29.73.151&185&34& & &3.25.11.13.131.491&1557&6958 \cr
6417&685685715&4.25.11.13.17.29.37&1611&2536&6435&689842725&4.27.5.49.71.173&559&386 \cr
 & &64.9.17.179.317&16167&5728& & &16.7.13.43.71.193&8299&3976 \cr
\noalign{\hrule}
 & &9.11.13.17.23.29.47&665&2& & &81.37.409.563&15235&30368 \cr
6418&685884771&4.3.5.7.11.19.47&3811&3944&6436&690110199&64.5.11.13.73.277&1049&1998 \cr
 & &64.17.29.37.103&103&1184& & &256.27.5.37.1049&1049&640 \cr
\noalign{\hrule}
 & &13.29.67.101.269&10027&28050& & &9.7.11.361.31.89&221&430 \cr
6419&686261771&4.3.25.11.17.37.271&101&84&6437&690227307&4.3.5.13.17.19.43.89&3131&5324 \cr
 & &32.9.5.7.11.101.271&3465&4336& & &32.1331.13.31.101&1313&1936 \cr
\noalign{\hrule}
 & &3.125.7.11.13.31.59&171&46& & &5.49.23.31.59.67&29241&12826 \cr
6420&686560875&4.27.11.13.19.23.59&149&500&6438&690529805&4.81.121.361.53&343&134 \cr
 & &32.125.19.23.149&437&2384& & &16.9.343.11.19.67&209&504 \cr
\noalign{\hrule}
 & &3.5.7.17.37.101.103&57&572& & &5.121.169.29.233&3033&1868 \cr
6421&687066135&8.9.7.11.13.19.101&1369&550&6439&690869465&8.9.121.337.467&13&350 \cr
 & &32.25.121.1369&185&1936& & &32.3.25.7.13.467&3269&240 \cr
\noalign{\hrule}
 & &3.13.31.691.823&259&950& & &3.5.11.17.83.2969&15707&16952 \cr
6422&687549837&4.25.7.19.37.823&319&504&6440&691227735&16.13.17.113.139.163&2241&122 \cr
 & &64.9.5.49.11.19.29&8265&17248& & &64.27.61.83.113&1017&1952 \cr
\noalign{\hrule}
 & &9.25.121.13.29.67&4313&5402& & &9.5.7.11.29.71.97&1937&1258 \cr
6423&687676275&4.5.13.19.37.73.227&93&1042&6441&692040195&4.11.13.17.29.37.149&485&4806 \cr
 & &16.3.19.31.37.521&16151&5624& & &16.27.5.17.89.97&89&408 \cr
\noalign{\hrule}
 & &27.11.17.29.37.127&2029&1400& & &9.11.17.173.2377&7375&9752 \cr
6424&688032279&16.25.7.11.29.2029&855&1174&6442&692084943&16.125.17.23.53.59&2541&586 \cr
 & &64.9.125.7.19.587&16625&18784& & &64.3.25.7.121.293&3223&5600 \cr
\noalign{\hrule}
 & &5.11.79.151.1049&2697&1648& & &27.19.71.83.229&73&156 \cr
6425&688243655&32.3.29.31.103.151&525&374&6443&692291961&8.81.13.19.71.73&715&634 \cr
 & &128.9.25.7.11.17.103&6489&5440& & &32.5.11.169.73.317&61685&55792 \cr
\noalign{\hrule}
 & &5.121.17.23.41.71&107&498& & &11.289.23.9479&6329&3150 \cr
6426&688611605&4.3.41.71.83.107&2967&2926&6444&693076043&4.9.25.7.23.6329&3337&2992 \cr
 & &16.9.7.11.19.23.43.107&6741&6536& & &128.3.5.7.11.17.47.71&4935&4544 \cr
\noalign{\hrule}
}%
}
$$
\eject
\vglue -23 pt
\noindent\hskip 1 in\hbox to 6.5 in{\ 6445 -- 6480 \hfill\fbd 693746235 -- 702965887\frb}
\vskip -9 pt
$$
\vbox{
\nointerlineskip
\halign{\strut
    \vrule \ \ \hfil \frb #\ 
   &\vrule \hfil \ \ \fbb #\frb\ 
   &\vrule \hfil \ \ \frb #\ \hfil
   &\vrule \hfil \ \ \frb #\ 
   &\vrule \hfil \ \ \frb #\ \ \vrule \hskip 2 pt
   &\vrule \ \ \hfil \frb #\ 
   &\vrule \hfil \ \ \fbb #\frb\ 
   &\vrule \hfil \ \ \frb #\ \hfil
   &\vrule \hfil \ \ \frb #\ 
   &\vrule \hfil \ \ \frb #\ \vrule \cr%
\noalign{\hrule}
 & &27.5.7.13.149.379&8453&8602& & &81.25.11.13.19.127&73&3502 \cr
6445&693746235&4.3.7.11.13.17.23.79.107&1895&30098&6463&698744475&4.3.17.19.73.103&1205&2956 \cr
 & &16.5.11.101.149.379&101&88& & &32.5.241.739&241&11824 \cr
\noalign{\hrule}
 & &81.23.841.443&1573&2414& & &27.5.7.29.97.263&1441&1178 \cr
6446&694084869&4.9.121.13.17.23.71&1333&1450&6464&699128955&4.5.7.11.19.29.31.131&1843&2742 \cr
 & &16.25.17.29.31.43.71&18275&17608& & &16.3.11.361.97.457&3971&3656 \cr
\noalign{\hrule}
 & &9.5.7.17.103.1259&247&1012& & &81.5.49.131.269&11&1334 \cr
6447&694420335&8.7.11.13.19.23.103&289&432&6465&699317955&4.3.11.23.29.131&2173&2150 \cr
 & &256.27.289.19.23&1173&2432& & &16.25.29.41.43.53&7685&14104 \cr
\noalign{\hrule}
 & &9.5.7.13.19.79.113&10021&1094& & &27.5.7.11.17.37.107&23&142 \cr
6448&694565235&4.7.11.547.911&417&494&6466&699614685&4.9.23.37.71.107&5369&9328 \cr
 & &16.3.13.19.139.547&547&1112& & &128.7.11.13.53.59&3127&832 \cr
\noalign{\hrule}
 & &25.7.11.19.31.613&1561&10086& & &9.5.11.289.67.73&580469&579866 \cr
6449&695034725&4.3.49.1681.223&893&1116&6467&699682005&4.49.19.61.97.137.223&410241&748 \cr
 & &32.27.19.31.41.47&1269&656& & &32.3.11.13.17.67.157&157&208 \cr
\noalign{\hrule}
 & &81.5.11.17.67.137&9347&3182& & &27.13.289.67.103&2783&974 \cr
6450&695171565&4.9.13.37.43.719&1805&2524&6468&700030539&4.121.23.103.487&323&810 \cr
 & &32.5.361.43.631&15523&10096& & &16.81.5.11.17.19.23&1311&440 \cr
\noalign{\hrule}
 & &9.5.7.17.31.59.71&551&656& & &3.13.19.43.127.173&3855&1606 \cr
6451&695394945&32.3.19.29.31.41.59&869&2698&6469&700061973&4.9.5.11.19.73.257&6055&1172 \cr
 & &128.11.361.71.79&3971&5056& & &32.25.7.173.293&2051&400 \cr
\noalign{\hrule}
 & &11.17.251.14821&30879&16058& & &3.5.7.11.13.149.313&1945&1498 \cr
6452&695653277&4.9.7.31.37.47.73&785&748&6470&700254555&4.25.49.13.107.389&7857&26918 \cr
 & &32.3.5.11.17.31.47.157&7285&7536& & &16.81.43.97.313&1161&776 \cr
\noalign{\hrule}
 & &3.11.29.37.71.277&2703&5330& & &9.5.7.121.17.23.47&29&2086 \cr
6453&696388803&4.9.5.11.13.17.41.53&575&1108&6471&700439355&4.49.23.29.149&309&358 \cr
 & &32.125.23.53.277&1219&2000& & &16.3.103.149.179&15347&1432 \cr
\noalign{\hrule}
 & &3.25.169.17.53.61&8799&5566& & &13.47.53.97.223&5289&5192 \cr
6454&696630675&4.9.5.7.121.23.419&1219&1714&6472&700476673&16.3.11.13.41.43.53.59&669&20 \cr
 & &16.11.529.53.857&5819&6856& & &128.9.5.41.43.223&1935&2624 \cr
\noalign{\hrule}
 & &3.11.13.17.19.47.107&35&22& & &9.5.11.169.289.29&33067&42362 \cr
6455&696853443&4.5.7.121.17.47.107&7657&2628&6473&701112555&4.43.59.359.769&1653&884 \cr
 & &32.9.7.13.19.31.73&651&1168& & &32.3.13.17.19.29.359&359&304 \cr
\noalign{\hrule}
 & &121.19.113.2683&1963&720& & &3.11.17.47.67.397&257&260 \cr
6456&697008521&32.9.5.11.13.19.151&2491&226&6474&701335833&8.5.13.17.67.257.397&423&1562 \cr
 & &128.3.47.53.113&2491&192& & &32.9.11.13.47.71.257&3341&3408 \cr
\noalign{\hrule}
 & &27.25.17.31.37.53&551&286& & &121.19.43.47.151&425&468 \cr
6457&697576725&4.5.11.13.17.19.29.37&339&154&6475&701588129&8.9.25.121.13.17.151&1889&74 \cr
 & &16.3.7.121.13.19.113&10283&18392& & &32.3.5.17.37.1889&28335&10064 \cr
\noalign{\hrule}
 & &3.5.73.673.947&3377&1358& & &5.11.41.61.5101&2581&2520 \cr
6458&697877445&4.7.11.73.97.307&673&1476&6476&701668055&16.9.25.7.11.29.41.89&9461&1486 \cr
 & &32.9.41.97.673&291&656& & &64.3.7.743.9461&66227&71328 \cr
\noalign{\hrule}
 & &3.25.11.13.37.1759&3829&26696& & &81.25.13.19.23.61&239&814 \cr
6459&698015175&16.7.47.71.547&297&250&6477&701745525&4.11.19.37.61.239&687&16 \cr
 & &64.27.125.7.11.71&639&1120& & &128.3.229.239&229&15296 \cr
\noalign{\hrule}
 & &7.11.13.37.47.401&48569&55290& & &9.5.11.13.29.53.71&35&178 \cr
6460&698036339&4.3.5.17.19.97.2857&507&2350&6478&702232245&4.3.25.7.29.53.89&1031&506 \cr
 & &16.9.125.169.17.47&1625&1224& & &16.11.23.89.1031&2047&8248 \cr
\noalign{\hrule}
 & &27.11.13.223.811&183&40& & &25.19.79.97.193&1963&462 \cr
6461&698273433&16.81.5.61.811&365&446&6479&702505525&4.3.7.11.13.151.193&135&58 \cr
 & &64.25.61.73.223&1525&2336& & &16.81.5.13.29.151&12231&3016 \cr
\noalign{\hrule}
 & &5.7.13.107.113.127&3033&3922& & &13.29.47.97.409&205&204 \cr
6462&698678435&4.9.37.53.113.337&143&196&6480&702965887&8.3.5.13.17.29.41.47.97&121&2934 \cr
 & &32.3.49.11.13.37.337&7077&6512& & &32.27.121.17.41.163&74817&79376 \cr
\noalign{\hrule}
}%
}
$$
\eject
\vglue -23 pt
\noindent\hskip 1 in\hbox to 6.5 in{\ 6481 -- 6516 \hfill\fbd 703016325 -- 717537645\frb}
\vskip -9 pt
$$
\vbox{
\nointerlineskip
\halign{\strut
    \vrule \ \ \hfil \frb #\ 
   &\vrule \hfil \ \ \fbb #\frb\ 
   &\vrule \hfil \ \ \frb #\ \hfil
   &\vrule \hfil \ \ \frb #\ 
   &\vrule \hfil \ \ \frb #\ \ \vrule \hskip 2 pt
   &\vrule \ \ \hfil \frb #\ 
   &\vrule \hfil \ \ \fbb #\frb\ 
   &\vrule \hfil \ \ \frb #\ \hfil
   &\vrule \hfil \ \ \frb #\ 
   &\vrule \hfil \ \ \frb #\ \vrule \cr%
\noalign{\hrule}
 & &9.25.11.479.593&159&434& & &9.11.13.283.1951&1615&1498 \cr
6481&703016325&4.27.7.31.53.479&21941&22420&6499&710595171&4.5.7.17.19.107.1951&1885&66 \cr
 & &32.5.7.19.37.59.593&2183&2128& & &16.3.25.7.11.13.19.29&551&1400 \cr
\noalign{\hrule}
 & &7.11.13.31.131.173&1425&2636& & &11.13.73.103.661&3895&4698 \cr
6482&703255553&8.3.25.11.13.19.659&1029&1688&6500&710718437&4.81.5.19.29.41.103&6059&6610 \cr
 & &128.9.25.343.211&10339&14400& & &16.27.25.73.83.661&675&664 \cr
\noalign{\hrule}
 & &7.23.37.43.2749&8675&27918& & &9.7.11.17.23.43.61&2375&362 \cr
6483&704159099&4.27.25.11.47.347&949&602&6501&710735949&4.3.125.19.43.181&1771&2314 \cr
 & &16.9.25.7.13.43.73&949&1800& & &16.25.7.11.13.23.89&325&712 \cr
\noalign{\hrule}
 & &3.11.19.41.67.409&4225&4016& & &27.13.37.127.431&6055&10754 \cr
6484&704448921&32.25.169.251.409&1107&938&6502&710869419&4.9.5.7.19.173.283&569&286 \cr
 & &128.27.5.7.41.67.251&2259&2240& & &16.7.11.13.173.569&6259&9688 \cr
\noalign{\hrule}
 & &3.5.7.13.19.29.937&6861&302& & &3.5.7.11.13.23.29.71&229&2114 \cr
6485&704731755&4.9.5.151.2287&1121&1166&6503&711065355&4.49.23.151.229&639&488 \cr
 & &16.11.19.53.59.151&8003&5192& & &64.9.61.71.229&687&1952 \cr
\noalign{\hrule}
 & &3.7.13.17.47.53.61&55&744& & &25.71.367.1093&363&730 \cr
6486&705204591&16.9.5.7.11.31.61&127&188&6504&712007525&4.3.125.121.71.73&3279&1904 \cr
 & &128.11.31.47.127&3937&704& & &128.9.7.11.17.1093&1071&704 \cr
\noalign{\hrule}
 & &7.11.79.269.431&8255&12996& & &25.7.11.13.71.401&27&428 \cr
6487&705256937&8.9.5.7.13.361.127&431&2236&6505&712486775&8.27.5.11.71.107&1121&56 \cr
 & &64.3.169.43.431&507&1376& & &128.9.7.19.59&10089&64 \cr
\noalign{\hrule}
 & &11.31.71.151.193&135&206& & &9.13.31.47.53.79&205&484 \cr
6488&705581173&4.27.5.103.151.193&667&88&6506&713753703&8.5.121.41.47.79&157&78 \cr
 & &64.9.11.23.29.103&2369&8352& & &32.3.121.13.41.157&4961&2512 \cr
\noalign{\hrule}
 & &3.125.49.19.43.47&403&528& & &9.49.17.23.41.101&2291&1850 \cr
6489&705581625&32.9.11.13.31.43.47&25&448&6507&714036771&4.25.17.23.29.37.79&2123&168 \cr
 & &4096.25.7.13.31&403&2048& & &64.3.5.7.11.37.193&2035&6176 \cr
\noalign{\hrule}
 & &11.17.31.193.631&47385&54326& & &9.25.17.31.37.163&793&4868 \cr
6490&705976051&4.729.5.13.23.1181&631&550&6508&715125825&8.13.31.61.1217&407&810 \cr
 & &16.9.125.11.13.23.631&1625&1656& & &32.81.5.11.37.61&671&144 \cr
\noalign{\hrule}
 & &3.5.7.11.13.31.37.41&881&266& & &9.5.49.29.67.167&355927&355994 \cr
6491&706110405&4.49.11.13.19.881&10897&1206&6509&715480605&4.3.5.11.13.19.23.71.109.131&270373&928 \cr
 & &16.9.17.67.641&1139&15384& & &256.13.29.167.1619&1619&1664 \cr
\noalign{\hrule}
 & &9.841.31.3011&2279&5290& & &27.5.47.257.439&133&572 \cr
6492&706498029&4.5.529.31.43.53&783&1496&6510&715861935&8.9.7.11.13.19.257&79&178 \cr
 & &64.27.5.11.17.23.29&1955&1056& & &32.7.13.19.79.89&7031&27664 \cr
\noalign{\hrule}
 & &27.11.23.31.47.71&3287&4556& & &9.49.13.29.31.139&1375&3182 \cr
6493&706646457&8.17.19.67.71.173&517&690&6511&716401413&4.3.125.11.29.37.43&973&248 \cr
 & &32.3.5.11.19.23.47.67&335&304& & &64.5.7.31.43.139&43&160 \cr
\noalign{\hrule}
 & &3.7.37.53.89.193&7205&6946& & &27.5.11.13.17.37.59&181&116 \cr
6494&707366037&4.5.11.23.131.151.193&31&162&6512&716427855&8.17.29.37.59.181&2847&3850 \cr
 & &16.81.5.11.23.31.151&34155&37448& & &32.3.25.7.11.13.29.73&1015&1168 \cr
\noalign{\hrule}
 & &61.101.313.367&8987&28080& & &27.5.37.67.2141&923&886 \cr
6495&707720231&32.27.5.11.13.19.43&5321&5794&6513&716517765&4.5.13.71.443.2141&37&2178 \cr
 & &128.3.17.313.2897&2897&3264& & &16.9.121.13.37.71&1573&568 \cr
\noalign{\hrule}
 & &9.5.11.29.31.37.43&71&114& & &11.29.31.71.1021&2699&26910 \cr
6496&708002955&4.27.11.19.29.31.71&301&598&6514&716863499&4.9.5.13.23.2699&901&1798 \cr
 & &16.7.13.19.23.43.71&5681&3976& & &16.3.5.17.29.31.53&901&120 \cr
\noalign{\hrule}
 & &9.5.49.11.23.31.41&293&76& & &9.169.47.79.127&1045&8988 \cr
6497&709046415&8.5.7.11.19.23.293&53&108&6515&717229071&8.27.5.7.11.19.107&139&158 \cr
 & &64.27.19.53.293&3021&9376& & &32.5.7.79.107.139&3745&2224 \cr
\noalign{\hrule}
 & &3.5.7.73.151.613&7657&3366& & &9.5.11.59.79.311&69&10 \cr
6498&709495395&4.27.5.11.13.17.19.31&389&1226&6516&717537645&4.27.25.11.23.311&7303&472 \cr
 & &16.11.13.389.613&389&1144& & &64.59.67.109&67&3488 \cr
\noalign{\hrule}
}%
}
$$
\eject
\vglue -23 pt
\noindent\hskip 1 in\hbox to 6.5 in{\ 6517 -- 6552 \hfill\fbd 717742095 -- 728007735\frb}
\vskip -9 pt
$$
\vbox{
\nointerlineskip
\halign{\strut
    \vrule \ \ \hfil \frb #\ 
   &\vrule \hfil \ \ \fbb #\frb\ 
   &\vrule \hfil \ \ \frb #\ \hfil
   &\vrule \hfil \ \ \frb #\ 
   &\vrule \hfil \ \ \frb #\ \ \vrule \hskip 2 pt
   &\vrule \ \ \hfil \frb #\ 
   &\vrule \hfil \ \ \fbb #\frb\ 
   &\vrule \hfil \ \ \frb #\ \hfil
   &\vrule \hfil \ \ \frb #\ 
   &\vrule \hfil \ \ \frb #\ \vrule \cr%
\noalign{\hrule}
 & &3.5.7.37.239.773&16093&12508& & &9.5.13.29.157.271&1631&1892 \cr
6517&717742095&8.49.121.19.53.59&773&234&6535&721809855&8.5.7.11.43.157.233&7453&702 \cr
 & &32.9.11.13.59.773&767&528& & &32.27.11.13.29.257&257&528 \cr
\noalign{\hrule}
 & &3.11.17.19.31.41.53&155&852& & &81.5.7.11.79.293&137&430 \cr
6518&718022217&8.9.5.11.961.71&3995&3034&6536&721839195&4.25.11.43.79.137&5549&5274 \cr
 & &32.25.17.37.41.47&1175&592& & &16.9.31.43.179.293&1333&1432 \cr
\noalign{\hrule}
 & &3.5.7.11.13.17.29.97&2421&2914& & &11.17.29.37.59.61&541&1170 \cr
6519&718032315&4.27.7.13.31.47.269&1363&1094&6537&722142949&4.9.5.11.13.61.541&1147&476 \cr
 & &16.29.31.2209.547&16957&17672& & &32.3.5.7.13.17.31.37&1209&560 \cr
\noalign{\hrule}
 & &3.7.13.29.83.1093&7&1086& & &27.25.7.13.19.619&9287&8668 \cr
6520&718222323&4.9.49.29.181&4565&4304&6538&722419425&8.5.11.13.37.197.251&119&2286 \cr
 & &128.5.11.83.269&2959&320& & &32.9.7.17.127.251&2159&4016 \cr
\noalign{\hrule}
 & &9.19.23.31.43.137&89&500& & &13.17.19.59.2917&8745&5828 \cr
6521&718248393&8.3.125.23.43.89&307&682&6539&722660497&8.3.5.11.17.31.47.53&771&686 \cr
 & &32.11.31.89.307&979&4912& & &32.9.343.11.53.257&88151&83952 \cr
\noalign{\hrule}
 & &7.11.13.19.179.211&745&1998& & &9.5.11.19.43.1787&4553&4382 \cr
6522&718328611&4.27.5.11.19.37.149&179&30&6540&722689605&4.7.11.29.43.157.313&5361&68900 \cr
 & &16.81.25.37.179&2025&296& & &32.3.25.13.53.1787&265&208 \cr
\noalign{\hrule}
 & &3.5.19.53.199.239&3003&7544& & &3.5.7.11.13.67.719&14883&15602 \cr
6523&718408905&16.9.5.7.11.13.23.41&53&152&6541&723317595&4.9.1331.29.41.269&475&11504 \cr
 & &256.7.13.19.23.53&299&896& & &128.25.19.29.719&551&320 \cr
\noalign{\hrule}
 & &25.7.11.29.61.211&1501&24& & &3.5.121.13.23.31.43&7837&7772 \cr
6524&718523575&16.3.11.19.29.79&169&150&6542&723399105&8.17.23.29.31.67.461&12943&2340 \cr
 & &64.9.25.169.79&711&5408& & &64.9.5.7.13.1849.67&1407&1376 \cr
\noalign{\hrule}
 & &3.121.13.23.37.179&2845&1272& & &7.11.37.229.1109&425&684 \cr
6525&718840551&16.9.5.37.53.569&13783&11822&6543&723534889&8.9.25.11.17.19.229&23&252 \cr
 & &64.7.11.23.179.257&257&224& & &64.81.7.17.19.23&1863&10336 \cr
\noalign{\hrule}
 & &3.5.7.19.29.31.401&1297&708& & &7.11.2789.3373&27145&3534 \cr
6526&719195505&8.9.7.29.59.1297&11825&152&6544&724361869&4.3.5.19.31.61.89&901&990 \cr
 & &128.25.11.19.43&55&2752& & &16.27.25.11.17.19.53&8075&11448 \cr
\noalign{\hrule}
 & &5.7.43.53.71.127&1317&962& & &25.13.83.97.277&4267&27258 \cr
6527&719241005&4.3.7.13.37.127.439&2343&13900&6545&724791275&4.3.7.11.17.59.251&97&90 \cr
 & &32.9.25.11.71.139&695&1584& & &16.27.5.59.97.251&1593&2008 \cr
\noalign{\hrule}
 & &5.13.361.23.31.43&847&486& & &9.5.11.13.61.1847&109727&106372 \cr
6528&719413435&4.243.5.7.121.13.23&607&152&6546&725012145&8.7.29.131.179.613&40241&40062 \cr
 & &64.81.11.19.607&6677&2592& & &32.3.7.11.29.607.40241&281687&281648 \cr
\noalign{\hrule}
 & &5.11.13.41.79.311&91&3330& & &9.11.13.23.127.193&5443&3320 \cr
6529&720240235&4.9.25.7.169.37&2057&2168&6547&725550111&16.3.5.13.83.5443&11989&15226 \cr
 & &64.3.7.121.17.271&8943&3808& & &64.19.23.331.631&11989&10592 \cr
\noalign{\hrule}
 & &5.7.11.19.29.43.79&279&274& & &5.11.13.19.529.101&363&2282 \cr
6530&720622595&4.9.11.19.29.31.43.137&1975&1766&6548&725832965&4.3.7.1331.13.163&1725&394 \cr
 & &16.3.25.31.79.137.883&21235&21192& & &16.9.25.7.23.197&197&2520 \cr
\noalign{\hrule}
 & &9.17.31.47.53.61&6875&7912& & &3.125.11.13.19.23.31&5671&2954 \cr
6531&720703593&16.625.11.23.43.47&151&366&6549&726457875&4.7.31.53.107.211&3933&7250 \cr
 & &64.3.125.23.61.151&3473&4000& & &16.9.125.7.19.23.29&87&56 \cr
\noalign{\hrule}
 & &9.5.7.13.19.73.127&35&22& & &25.11.17.29.31.173&117&202 \cr
6532&721330155&4.3.25.49.11.73.127&703&19372&6550&727088725&4.9.5.13.31.101.173&1267&748 \cr
 & &32.19.29.37.167&6179&464& & &32.3.7.11.17.101.181&2121&2896 \cr
\noalign{\hrule}
 & &625.7.23.71.101&2079&446& & &5.41.67.167.317&33337&19602 \cr
6533&721581875&4.27.25.49.11.223&391&1616&6551&727117165&4.81.121.17.37.53&317&266 \cr
 & &128.3.11.17.23.101&561&64& & &16.27.7.11.19.37.317&3591&3256 \cr
\noalign{\hrule}
 & &9.7.11.41.109.233&1415&216& & &3.5.7.13.29.53.347&27071&26724 \cr
6534&721604961&16.243.5.41.283&5123&4840&6552&728007735&8.9.11.13.17.23.107.131&1735&25382 \cr
 & &256.25.121.47.109&1175&1408& & &32.5.343.11.37.347&539&592 \cr
\noalign{\hrule}
}%
}
$$
\eject
\vglue -23 pt
\noindent\hskip 1 in\hbox to 6.5 in{\ 6553 -- 6588 \hfill\fbd 728107069 -- 734687421\frb}
\vskip -9 pt
$$
\vbox{
\nointerlineskip
\halign{\strut
    \vrule \ \ \hfil \frb #\ 
   &\vrule \hfil \ \ \fbb #\frb\ 
   &\vrule \hfil \ \ \frb #\ \hfil
   &\vrule \hfil \ \ \frb #\ 
   &\vrule \hfil \ \ \frb #\ \ \vrule \hskip 2 pt
   &\vrule \ \ \hfil \frb #\ 
   &\vrule \hfil \ \ \fbb #\frb\ 
   &\vrule \hfil \ \ \frb #\ \hfil
   &\vrule \hfil \ \ \frb #\ 
   &\vrule \hfil \ \ \frb #\ \vrule \cr%
\noalign{\hrule}
 & &41.563.31543&27313&4230& & &3.5.7.19.23.89.179&759&136 \cr
6553&728107069&4.9.5.11.13.47.191&15751&15764&6571&730993935&16.9.11.17.19.529&179&350 \cr
 & &32.3.7.19.47.563.829&15751&15792& & &64.25.7.11.17.179&85&352 \cr
\noalign{\hrule}
 & &9.25.17.191.997&13013&29962& & &81.125.23.43.73&3179&304 \cr
6554&728383275&4.7.11.169.71.211&191&1668&6572&730994625&32.11.289.19.73&725&516 \cr
 & &32.3.71.139.191&139&1136& & &256.3.25.17.29.43&493&128 \cr
\noalign{\hrule}
 & &27.5.49.11.17.19.31&149&688& & &11.23.41.107.659&27045&26 \cr
6555&728595945&32.5.17.19.43.149&451&366&6573&731431349&4.9.5.13.601&901&902 \cr
 & &128.3.11.41.61.149&6109&3904& & &16.3.5.11.13.17.41.53&901&1560 \cr
\noalign{\hrule}
 & &9.5.11.41.149.241&7031&922& & &9.5.13.31.157.257&5203&1862 \cr
6556&728773155&4.3.5.79.89.461&1271&1034&6574&731729115&4.49.121.19.31.43&243&1090 \cr
 & &16.11.31.41.47.89&1457&712& & &16.243.5.7.19.109&763&4104 \cr
\noalign{\hrule}
 & &3.5.11.17.31.83.101&13293&15010& & &81.11.13.23.41.67&3185&2242 \cr
6557&728943765&4.27.25.7.19.79.211&3113&1612&6575&731825523&4.5.49.11.169.19.59&69&1114 \cr
 & &32.11.13.31.211.283&3679&3376& & &16.3.7.23.59.557&557&3304 \cr
\noalign{\hrule}
 & &49.11.169.53.151&3683&4320& & &81.5.7.23.103.109&197&524 \cr
6558&729001273&64.27.5.11.13.29.127&4931&784&6576&732056535&8.27.5.23.131.197&803&182 \cr
 & &2048.3.49.4931&4931&3072& & &32.7.11.13.73.131&10439&2096 \cr
\noalign{\hrule}
 & &81.5.7.43.5981&3055&2926& & &27.5.37.47.3119&143&190 \cr
6559&729113805&4.27.25.49.11.13.19.47&1247&22&6577&732232035&4.3.25.11.13.19.3119&1147&1972 \cr
 & &16.121.13.19.29.43&1573&4408& & &32.13.17.19.29.31.37&7163&8432 \cr
\noalign{\hrule}
 & &11.17.29.43.53.59&16575&17822& & &3.25.37.41.47.137&31691&16484 \cr
6560&729182003&4.3.25.7.13.289.19.67&477&188&6578&732597225&8.11.13.43.67.317&137&180 \cr
 & &32.27.5.13.47.53.67&9045&9776& & &64.9.5.11.13.67.137&871&1056 \cr
\noalign{\hrule}
 & &125.49.11.79.137&599&10224& & &9.19.61.163.431&877&9066 \cr
6561&729199625&32.9.7.71.599&1045&446&6579&732809043&4.27.877.1511&12595&11084 \cr
 & &128.3.5.11.19.223&223&3648& & &32.5.11.17.163.229&2519&1360 \cr
\noalign{\hrule}
 & &81.7.11.19.47.131&31&50& & &5.13.19.37.43.373&8371&7668 \cr
6562&729622971&4.25.7.11.31.47.131&4833&3392&6580&732902105&8.27.5.11.13.71.761&3287&518 \cr
 & &512.27.31.53.179&9487&7936& & &32.9.7.11.19.37.173&1211&1584 \cr
\noalign{\hrule}
 & &3.19.41.43.53.137&85&44& & &3.11.19.29.157.257&70725&70882 \cr
6563&729665151&8.5.11.17.19.53.137&3735&3526&6581&733665867&4.9.25.7.11.23.41.61.83&3161&8224 \cr
 & &32.9.25.17.41.43.83&1275&1328& & &256.5.7.29.41.109.257&4469&4480 \cr
\noalign{\hrule}
 & &25.7.19.41.53.101&1419&754& & &3.5.13.17.19.61.191&6979&27974 \cr
6564&729747725&4.3.5.11.13.29.43.101&27&532&6582&733838235&4.7.71.197.997&191&1188 \cr
 & &32.81.7.11.19.29&2349&176& & &32.27.11.71.191&99&1136 \cr
\noalign{\hrule}
 & &25.13.29.43.1801&38709&38734& & &3.13.43.47.67.139&1135&1746 \cr
6565&729900275&4.9.11.13.17.23.29.107.181&665&2&6583&734041347&4.27.5.97.139.227&1309&2444 \cr
 & &16.3.5.7.11.19.107.181&24717&27512& & &32.7.11.13.17.47.97&1649&1232 \cr
\noalign{\hrule}
 & &5.121.13.17.43.127&421&294& & &25.11.13.17.47.257&27&248 \cr
6566&730163005&4.3.49.11.17.43.421&635&96&6584&734101225&16.27.31.47.257&1133&1180 \cr
 & &256.9.5.127.421&421&1152& & &128.3.5.11.31.59.103&5487&6592 \cr
\noalign{\hrule}
 & &25.11.23.263.439&41211&31114& & &3.5.11.13.151.2267&1139&1126 \cr
6567&730265525&4.9.19.47.241.331&15097&3770&6585&734269965&4.11.17.67.563.2267&1963&4230 \cr
 & &16.3.5.13.29.31.487&15097&9048& & &16.9.5.13.17.47.67.151&1139&1128 \cr
\noalign{\hrule}
 & &7.11.13.29.139.181&154017&156370& & &3.5.7.11.17.113.331&933&988 \cr
6568&730340611&4.9.5.19.109.157.823&257&2726&6586&734407905&8.9.7.13.19.311.331&40205&3818 \cr
 & &16.3.5.29.47.109.257&12079&13080& & &32.5.11.17.23.43.83&989&1328 \cr
\noalign{\hrule}
 & &9.5.7.11.17.79.157&229&166& & &7.13.17.19.67.373&2201&27192 \cr
6569&730598715&4.11.17.83.157.229&1569&158&6587&734560463&16.3.11.31.71.103&51&20 \cr
 & &16.3.79.229.523&523&1832& & &128.9.5.11.17.103&495&6592 \cr
\noalign{\hrule}
 & &7.47.67.71.467&2145&2612& & &3.13.19.29.179.191&365&186 \cr
6570&730879751&8.3.5.7.11.13.47.653&2465&1812&6588&734687421&4.9.5.13.31.73.191&385&1334 \cr
 & &64.9.25.11.17.29.151&74443&79200& & &16.25.7.11.23.29.31&5425&2024 \cr
\noalign{\hrule}
}%
}
$$
\eject
\vglue -23 pt
\noindent\hskip 1 in\hbox to 6.5 in{\ 6589 -- 6624 \hfill\fbd 734763575 -- 743702337\frb}
\vskip -9 pt
$$
\vbox{
\nointerlineskip
\halign{\strut
    \vrule \ \ \hfil \frb #\ 
   &\vrule \hfil \ \ \fbb #\frb\ 
   &\vrule \hfil \ \ \frb #\ \hfil
   &\vrule \hfil \ \ \frb #\ 
   &\vrule \hfil \ \ \frb #\ \ \vrule \hskip 2 pt
   &\vrule \ \ \hfil \frb #\ 
   &\vrule \hfil \ \ \fbb #\frb\ 
   &\vrule \hfil \ \ \frb #\ \hfil
   &\vrule \hfil \ \ \frb #\ 
   &\vrule \hfil \ \ \frb #\ \vrule \cr%
\noalign{\hrule}
 & &25.49.13.29.37.43&207&352& & &3.11.79.503.563&533&30 \cr
6589&734763575&64.9.5.49.11.23.37&793&58&6607&738273723&4.9.5.11.13.41.79&563&464 \cr
 & &256.3.11.13.29.61&183&1408& & &128.5.29.41.563&1189&320 \cr
\noalign{\hrule}
 & &13.1369.157.263&1705&336& & &7.11.13.17.43.1009&109839&113150 \cr
6590&734855927&32.3.5.7.11.31.263&1369&1524&6608&738316579&4.3.25.19.31.41.47.73&959&66 \cr
 & &256.9.7.1369.127&889&1152& & &16.9.7.11.31.73.137&4247&5256 \cr
\noalign{\hrule}
 & &3.5.49.17.83.709&3993&8060& & &9.25.7.11.19.2243&173&2 \cr
6591&735293265&8.9.25.1331.13.31&413&3338&6609&738339525&4.11.173.2243&12423&12250 \cr
 & &32.7.11.59.1669&1669&10384& & &16.3.125.49.41.101&1435&808 \cr
\noalign{\hrule}
 & &13.61.71.73.179&6137&6930& & &27.11.17.37.59.67&30233&13250 \cr
6592&735711301&4.9.5.7.11.17.361.71&219&716&6610&738471789&4.125.49.53.617&9&44 \cr
 & &32.27.361.73.179&361&432& & &32.9.25.7.11.617&617&2800 \cr
\noalign{\hrule}
 & &3.5.37.47.89.317&32147&17248& & &11.17.23.29.31.191&53997&46660 \cr
6593&735936105&64.49.11.17.31.61&65&54&6611&738520409&8.3.5.41.439.2333&7833&10166 \cr
 & &256.27.5.7.13.31.61&13237&14976& & &32.9.5.7.13.17.23.373&4849&5040 \cr
\noalign{\hrule}
 & &25.7.11.17.43.523&1789&1866& & &27.25.11.29.47.73&1547&4978 \cr
6594&735952525&4.3.5.311.523.1789&1039&516&6612&738780075&4.3.7.11.13.17.19.131&1205&236 \cr
 & &32.9.43.1039.1789&16101&16624& & &32.5.7.13.59.241&1687&12272 \cr
\noalign{\hrule}
 & &9.5.13.17.19.47.83&37&28& & &3.5.11.13.53.67.97&21&76 \cr
6595&737113455&8.7.17.19.37.47.83&389&1188&6613&738838815&8.9.7.13.19.53.67&33325&34144 \cr
 & &64.27.7.11.37.389&8547&12448& & &512.25.11.31.43.97&1333&1280 \cr
\noalign{\hrule}
 & &3.5.7.13.23.53.443&17&282& & &3.125.49.19.29.73&649&576 \cr
6596&737123205&4.9.7.17.47.443&253&190&6614&739097625&128.27.5.11.19.29.59&427&2138 \cr
 & &16.5.11.17.19.23.47&517&2584& & &512.7.11.61.1069&11759&15616 \cr
\noalign{\hrule}
 & &9.11.23.41.53.149&9135&1238& & &3.5.7.19.29.67.191&407&598 \cr
6597&737240229&4.81.5.7.29.619&269&298&6615&740370435&4.7.11.13.19.23.29.37&603&470 \cr
 & &16.5.149.269.619&3095&2152& & &16.9.5.11.13.23.47.67&1833&2024 \cr
\noalign{\hrule}
 & &3.5.7.19.41.71.127&809&682& & &3.125.11.23.37.211&653&722 \cr
6598&737545515&4.5.11.19.31.41.809&10863&4508&6616&740689125&4.361.37.211.653&10107&2300 \cr
 & &32.9.49.11.17.23.71&1309&1104& & &32.9.25.19.23.1123&1123&912 \cr
\noalign{\hrule}
 & &9.7.11.13.19.31.139&80975&85408& & &9.5.13.19.103.647&5863&3922 \cr
6599&737575839&64.25.17.41.79.157&741&44&6617&740714715&4.3.11.169.37.41.53&133&6386 \cr
 & &512.3.5.11.13.19.79&395&256& & &16.7.11.19.31.103&31&616 \cr
\noalign{\hrule}
 & &9.7.11.289.29.127&8303&28400& & &27.25.169.73.89&253&422 \cr
6600&737620191&32.25.361.23.71&4263&3902&6618&741145275&4.11.23.73.89.211&1425&622 \cr
 & &128.3.5.49.29.1951&1951&2240& & &16.3.25.19.211.311&4009&2488 \cr
\noalign{\hrule}
 & &9.5.7.19.59.2089&959&1130& & &9.25.169.17.31.37&1223&298 \cr
6601&737657235&4.25.49.59.113.137&33&2858&6619&741449475&4.17.31.149.1223&12705&8086 \cr
 & &16.3.11.137.1429&15719&1096& & &16.3.5.7.121.13.311&2177&968 \cr
\noalign{\hrule}
 & &9.5.13.23.109.503&491&7030& & &9.7.121.13.59.127&5513&6044 \cr
6602&737697285&4.3.25.19.37.491&1199&274&6620&742548807&8.121.37.149.1511&135&14 \cr
 & &16.11.19.109.137&2603&88& & &32.27.5.7.37.1511&4533&2960 \cr
\noalign{\hrule}
 & &27.25.19.23.41.61&3859&27434& & &9.11.19.67.71.83&1519&2300 \cr
6603&737732475&4.11.17.29.43.227&123&350&6621&742677111&8.3.25.49.23.31.83&667&418 \cr
 & &16.3.25.7.17.29.41&493&56& & &32.5.7.11.19.529.29&2645&3248 \cr
\noalign{\hrule}
 & &5.49.43.89.787&501&286& & &27.25.11.13.43.179&59&238 \cr
6604&737903005&4.3.49.11.13.89.167&225&314&6622&742952925&4.25.7.13.17.43.59&1969&1944 \cr
 & &16.27.25.13.157.167&21195&17368& & &64.243.11.17.59.179&531&544 \cr
\noalign{\hrule}
 & &7.13.19.37.83.139&231&250& & &5.11.169.167.479&2195&4032 \cr
6605&738056501&4.3.125.49.11.83.139&969&12506&6623&743534935&128.9.25.7.13.439&479&1796 \cr
 & &16.9.5.169.17.19.37&221&360& & &1024.3.449.479&1347&512 \cr
\noalign{\hrule}
 & &3.11.53.503.839&27173&514& & &27.7.661.5953&65441&59488 \cr
6606&738107733&4.29.257.937&4195&3258&6624&743702337&64.11.169.31.2111&1985&126 \cr
 & &16.9.5.181.839&543&40& & &256.9.5.7.31.397&1985&3968 \cr
\noalign{\hrule}
}%
}
$$
\eject
\vglue -23 pt
\noindent\hskip 1 in\hbox to 6.5 in{\ 6625 -- 6660 \hfill\fbd 743899845 -- 753049583\frb}
\vskip -9 pt
$$
\vbox{
\nointerlineskip
\halign{\strut
    \vrule \ \ \hfil \frb #\ 
   &\vrule \hfil \ \ \fbb #\frb\ 
   &\vrule \hfil \ \ \frb #\ \hfil
   &\vrule \hfil \ \ \frb #\ 
   &\vrule \hfil \ \ \frb #\ \ \vrule \hskip 2 pt
   &\vrule \ \ \hfil \frb #\ 
   &\vrule \hfil \ \ \fbb #\frb\ 
   &\vrule \hfil \ \ \frb #\ \hfil
   &\vrule \hfil \ \ \frb #\ 
   &\vrule \hfil \ \ \frb #\ \vrule \cr%
\noalign{\hrule}
 & &3.5.13.101.107.353&17131&18522& & &41.71.113.2273&3453&1180 \cr
6625&743899845&4.81.5.343.37.463&50479&12976&6643&747687439&8.3.5.59.71.1151&783&4972 \cr
 & &128.11.13.353.811&811&704& & &64.81.11.29.113&319&2592 \cr
\noalign{\hrule}
 & &3.5.17.19.97.1583&2233&6982& & &125.53.157.719&297&422 \cr
6626&743954595&4.7.11.17.29.3491&12125&12312&6644&747849875&4.27.11.53.157.211&213&370 \cr
 & &64.81.125.19.29.97&783&800& & &16.81.5.37.71.211&17091&21016 \cr
\noalign{\hrule}
 & &3.5.23.43.131.383&1477&488& & &3.5.17.71.109.379&64493&70052 \cr
6627&744316455&16.7.61.211.383&297&86&6645&747935655&8.121.13.41.83.211&621&1700 \cr
 & &64.27.7.11.43.61&549&2464& & &64.27.25.11.17.23.41&4059&3680 \cr
\noalign{\hrule}
 & &307.463.5237&2465&2772& & &3.5.49.11.29.31.103&2071&916 \cr
6628&744392417&8.9.5.7.11.17.29.463&3481&3464&6646&748646745&8.7.19.31.109.229&927&1144 \cr
 & &128.3.7.11.29.3481.433&804111&803648& & &128.9.11.13.103.229&687&832 \cr
\noalign{\hrule}
 & &27.125.13.361.47&1177&1502& & &9.163.467.1093&1247&154 \cr
6629&744427125&4.9.5.11.19.107.751&803&52&6647&748802277&4.3.7.11.29.43.163&2335&1088 \cr
 & &32.121.13.73.107&7811&1936& & &512.5.11.17.467&187&1280 \cr
\noalign{\hrule}
 & &3.19.29.61.83.89&711&1870& & &25.11.17.37.61.71&2267&1638 \cr
6630&744853371&4.27.5.11.17.79.83&557&854&6648&749154725&4.9.5.7.13.61.2267&1591&676 \cr
 & &16.5.7.61.79.557&2765&4456& & &32.3.7.2197.37.43&6591&4816 \cr
\noalign{\hrule}
 & &3.5.23.31.257.271&221&492& & &25.7.17.419.601&1767&1166 \cr
6631&744874665&8.9.5.13.17.41.257&2489&176&6649&749161525&4.3.25.11.17.19.31.53&63&838 \cr
 & &256.11.17.19.131&27379&2176& & &16.27.7.11.19.419&513&88 \cr
\noalign{\hrule}
 & &9.5.11.41.73.503&28379&7756& & &5.11.41.43.59.131&23&18 \cr
6632&745212105&8.7.13.37.59.277&953&1230&6650&749442485&4.9.11.23.43.59.131&893&3430 \cr
 & &32.3.5.7.13.41.953&953&1456& & &16.3.5.343.19.23.47&16121&10488 \cr
\noalign{\hrule}
 & &9.25.17.23.37.229&313&316& & &3.5.7.11.29.61.367&33&1802 \cr
6633&745412175&8.3.25.23.79.229.313&137&5588&6651&749852565&4.9.7.121.17.53&185&662 \cr
 & &64.11.127.137.313&39751&48224& & &16.5.17.37.331&12247&136 \cr
\noalign{\hrule}
 & &27.5.7.13.47.1291&9101&64& & &9.5.7.37.229.281&1379&1150 \cr
6634&745416945&128.9.19.479&325&154&6652&749987595&4.125.49.23.37.197&1749&2876 \cr
 & &512.25.7.11.13&11&1280& & &32.3.11.53.197.719&38107&34672 \cr
\noalign{\hrule}
 & &3.5.7.11.31.109.191&13&204& & &9.121.13.29.31.59&425&106 \cr
6635&745424295&8.9.5.11.13.17.109&2939&3056&6653&750901437&4.25.11.13.17.31.53&335&192 \cr
 & &256.17.191.2939&2939&2176& & &512.3.125.53.67&6625&17152 \cr
\noalign{\hrule}
 & &81.13.19.83.449&1309&230& & &23.79.479.863&477&1340 \cr
6636&745600869&4.5.7.11.17.23.449&1577&1566&6654&751106009&8.9.5.53.67.479&1817&5368 \cr
 & &16.27.5.17.19.23.29.83&667&680& & &128.3.11.23.61.79&671&192 \cr
\noalign{\hrule}
 & &7.11.19.29.43.409&98657&94800& & &25.11.17.23.29.241&51103&43128 \cr
6637&746163649&32.3.25.13.79.7589&1227&6362&6655&751492225&16.9.13.599.3931&9715&13646 \cr
 & &128.9.5.409.3181&3181&2880& & &64.3.5.29.67.6823&6823&6432 \cr
\noalign{\hrule}
 & &9.49.169.17.19.31&6809&1570& & &3.5.7.11.37.43.409&107&366 \cr
6638&746258877&4.5.11.17.157.619&403&216&6656&751580445&4.9.5.61.107.409&1573&2108 \cr
 & &64.27.5.13.31.157&157&480& & &32.121.13.17.31.61&6851&10736 \cr
\noalign{\hrule}
 & &27.5.11.23.131.167&113&140& & &9.25.11.31.97.101&697&212 \cr
6639&747208935&8.25.7.113.131.167&12351&9526&6657&751674825&8.5.11.17.31.41.53&303&148 \cr
 & &32.3.7.11.23.179.433&3031&2864& & &64.3.17.37.53.101&629&1696 \cr
\noalign{\hrule}
 & &9.5.7.11.379.569&235&334& & &729.7.11.59.227&7207&7130 \cr
6640&747230715&4.25.7.47.167.379&761&3414&6658&751789269&4.3.5.23.31.227.7207&23863&2242 \cr
 & &16.3.47.569.761&761&376& & &16.49.19.31.59.487&4123&3896 \cr
\noalign{\hrule}
 & &25.7.17.23.67.163&13959&12284& & &9.5.13.17.19.23.173&8371&5084 \cr
6641&747269425&8.27.11.17.37.47.83&335&182&6659&751851945&8.11.17.31.41.761&519&178 \cr
 & &32.3.5.7.13.37.67.83&1443&1328& & &32.3.89.173.761&761&1424 \cr
\noalign{\hrule}
 & &3.25.19.41.67.191&635&638& & &11.13.29.41.43.103&2375&612 \cr
6642&747664725&4.125.11.29.41.127.191&207&5332&6660&753049583&8.9.125.11.13.17.19&41&206 \cr
 & &32.9.11.23.31.43.127&43989&46736& & &32.3.25.17.41.103&425&48 \cr
\noalign{\hrule}
}%
}
$$
\eject
\vglue -23 pt
\noindent\hskip 1 in\hbox to 6.5 in{\ 6661 -- 6696 \hfill\fbd 753790851 -- 764421405\frb}
\vskip -9 pt
$$
\vbox{
\nointerlineskip
\halign{\strut
    \vrule \ \ \hfil \frb #\ 
   &\vrule \hfil \ \ \fbb #\frb\ 
   &\vrule \hfil \ \ \frb #\ \hfil
   &\vrule \hfil \ \ \frb #\ 
   &\vrule \hfil \ \ \frb #\ \ \vrule \hskip 2 pt
   &\vrule \ \ \hfil \frb #\ 
   &\vrule \hfil \ \ \fbb #\frb\ 
   &\vrule \hfil \ \ \frb #\ \hfil
   &\vrule \hfil \ \ \frb #\ 
   &\vrule \hfil \ \ \frb #\ \vrule \cr%
\noalign{\hrule}
 & &9.11.137.149.373&37&410& & &9.5.29.43.59.229&175&1886 \cr
6661&753790851&4.3.5.11.37.41.137&745&608&6679&758169765&4.125.7.23.41.43&3057&2068 \cr
 & &256.25.19.37.149&925&2432& & &32.3.7.11.47.1019&3619&16304 \cr
\noalign{\hrule}
 & &9.125.41.59.277&30653&29822& & &31.37.67.71.139&3695&6174 \cr
6662&753820875&4.3.5.7.13.29.31.37.151&697&5192&6680&758422781&4.9.5.343.31.739&695&44 \cr
 & &64.7.11.17.37.41.59&1309&1184& & &32.3.25.49.11.139&275&2352 \cr
\noalign{\hrule}
 & &3.25.19.29.71.257&1001&1058& & &3.13.59.487.677&627&140 \cr
6663&754057275&4.25.7.11.13.529.257&3683&342&6681&758637399&8.9.5.7.11.19.677&25979&26150 \cr
 & &16.9.11.19.23.29.127&759&1016& & &32.125.83.313.523&163699&166000 \cr
\noalign{\hrule}
 & &5.17.29.223.1373&3679&3186& & &27.11.13.19.79.131&365&346 \cr
6664&754731235&4.27.13.59.223.283&253&30&6682&759192291&4.3.5.11.13.73.131.173&1975&26714 \cr
 & &16.81.5.11.13.23.59&8437&14904& & &16.125.361.37.79&703&1000 \cr
\noalign{\hrule}
 & &3.13.31.421.1483&943&540& & &27.5.11.289.29.61&23737&24026 \cr
6665&754830687&8.81.5.23.41.421&121&1984&6683&759192885&4.5.7.11.41.293.3391&13113&125918 \cr
 & &1024.121.31.41&4961&512& & &16.9.13.29.31.47.167&5177&4888 \cr
\noalign{\hrule}
 & &27.11.23.29.37.103&2239&1572& & &289.19.1369.101&143&180 \cr
6666&754955289&8.81.11.131.2239&7009&17620&6684&759235079&8.9.5.11.13.17.37.101&5221&70 \cr
 & &64.5.43.163.881&37883&26080& & &32.3.25.7.23.227&1589&27600 \cr
\noalign{\hrule}
 & &9.125.7.11.23.379&2947&1222& & &81.7.11.23.67.79&85&76 \cr
6667&755110125&4.3.5.49.13.47.421&671&34&6685&759286143&8.9.5.11.17.19.67.79&1561&6854 \cr
 & &16.11.17.61.421&421&8296& & &32.7.19.23.149.223&2831&3568 \cr
\noalign{\hrule}
 & &27.11.37.53.1297&67&1364& & &7.11.17.271.2141&17845&32832 \cr
6668&755394849&8.121.31.37.67&1179&1300&6686&759496199&128.27.5.19.43.83&481&266 \cr
 & &64.9.25.13.31.131&10075&4192& & &512.3.7.13.361.37&14079&9472 \cr
\noalign{\hrule}
 & &27.5.7.11.23.29.109&235&7102& & &9.5.11.23.43.1553&1121&2674 \cr
6669&755747685&4.3.25.47.53.67&667&508&6687&760278915&4.3.7.19.43.59.191&109&1012 \cr
 & &32.23.29.67.127&127&1072& & &32.11.23.109.191&191&1744 \cr
\noalign{\hrule}
 & &3.5.13.37.311.337&1283&272& & &625.49.19.1307&2409&3716 \cr
6670&756182505&32.13.17.37.1283&3447&4730&6688&760510625&8.3.5.11.19.73.929&987&58 \cr
 & &128.9.5.11.43.383&4213&8256& & &32.9.7.29.47.73&1363&10512 \cr
\noalign{\hrule}
 & &11.103.331.2017&957&1060& & &5.13.73.109.1471&261&1210 \cr
6671&756421391&8.3.5.121.29.53.331&2379&4034&6689&760808555&4.9.25.121.29.109&301&26 \cr
 & &32.9.13.29.61.2017&1769&1872& & &16.3.7.11.13.29.43&301&7656 \cr
\noalign{\hrule}
 & &81.49.11.13.31.43&167&134& & &49.11.13.23.29.163&2927&1800 \cr
6672&756566811&4.27.7.13.31.67.167&2315&506&6690&761808047&16.9.25.11.13.2927&1821&1106 \cr
 & &16.5.11.23.167.463&10649&6680& & &64.27.5.7.79.607&16389&12640 \cr
\noalign{\hrule}
 & &81.7.11.29.47.89&241&560& & &13.41.587.2437&3807&27874 \cr
6673&756591759&32.9.5.49.47.241&67&2236&6691&762466627&4.81.7.11.47.181&925&1066 \cr
 & &256.5.13.43.67&871&27520& & &16.27.25.7.13.37.41&925&1512 \cr
\noalign{\hrule}
 & &9.11.23.31.71.151&12083&30046& & &9.11.17.31.47.311&56621&49430 \cr
6674&756763227&4.43.83.181.281&15553&7770&6692&762612741&4.5.41.1381.4943&1781&3162 \cr
 & &16.3.5.7.37.103.151&721&1480& & &16.3.5.13.17.31.41.137&1781&1640 \cr
\noalign{\hrule}
 & &3.7.13.37.137.547&187&224& & &25.11.29.163.587&1627&1308 \cr
6675&756958839&64.49.11.13.17.547&45&592&6693&763055975&8.3.5.109.163.1627&1221&406 \cr
 & &2048.9.5.11.17.37&2805&1024& & &32.9.7.11.29.37.109&2331&1744 \cr
\noalign{\hrule}
 & &9.5.19.47.83.227&4319&5896& & &5.13.19.29.83.257&33&62 \cr
6676&757125585&16.7.11.47.67.617&1819&4968&6694&763969765&4.3.11.13.31.83.257&1953&874 \cr
 & &256.27.7.17.23.107&12733&8832& & &16.27.7.19.23.961&4347&7688 \cr
\noalign{\hrule}
 & &3.5.7.121.109.547&29&356& & &7.11.13.17.97.463&33929&32280 \cr
6677&757510215&8.11.29.89.547&763&216&6695&764250487&16.3.5.49.37.131.269&4167&9014 \cr
 & &128.27.7.29.109&29&576& & &64.27.5.463.4507&4507&4320 \cr
\noalign{\hrule}
 & &27.11.17.31.47.103&1387&364& & &3.5.11.17.31.59.149&2337&2678 \cr
6678&757708479&8.9.7.13.19.47.73&301&310&6696&764421405&4.9.13.19.41.103.149&1441&2782 \cr
 & &32.5.49.19.31.43.73&15695&14896& & &16.11.169.19.107.131&18083&19912 \cr
\noalign{\hrule}
}%
}
$$
\eject
\vglue -23 pt
\noindent\hskip 1 in\hbox to 6.5 in{\ 6697 -- 6732 \hfill\fbd 764612805 -- 778273155\frb}
\vskip -9 pt
$$
\vbox{
\nointerlineskip
\halign{\strut
    \vrule \ \ \hfil \frb #\ 
   &\vrule \hfil \ \ \fbb #\frb\ 
   &\vrule \hfil \ \ \frb #\ \hfil
   &\vrule \hfil \ \ \frb #\ 
   &\vrule \hfil \ \ \frb #\ \ \vrule \hskip 2 pt
   &\vrule \ \ \hfil \frb #\ 
   &\vrule \hfil \ \ \fbb #\frb\ 
   &\vrule \hfil \ \ \frb #\ \hfil
   &\vrule \hfil \ \ \frb #\ 
   &\vrule \hfil \ \ \frb #\ \vrule \cr%
\noalign{\hrule}
 & &3.5.11.23.113.1783&9863&9750& & &9.5.11.17.19.61.79&3941&28036 \cr
6697&764612805&4.9.625.7.13.23.1409&847&13528&6715&770485815&8.7.43.163.563&12675&11534 \cr
 & &64.49.121.13.19.89&18601&20384& & &32.3.25.169.73.79&845&1168 \cr
\noalign{\hrule}
 & &3.5.7.11.13.127.401&363&272& & &5.49.41.59.1301&1807&4698 \cr
6698&764668905&32.9.1331.17.401&1139&2470&6716&771044155&4.81.13.29.41.139&9605&5852 \cr
 & &128.5.13.289.19.67&5491&4288& & &32.3.5.7.11.17.19.113&5763&3344 \cr
\noalign{\hrule}
 & &125.49.11.41.277&1539&4586& & &3.25.11.13.227.317&679&906 \cr
6699&765177875&4.81.19.41.2293&277&236&6717&771760275&4.9.5.7.11.13.97.151&9827&1522 \cr
 & &32.3.59.277.2293&2293&2832& & &16.7.31.317.761&761&1736 \cr
\noalign{\hrule}
 & &5.7.19.59.109.179&615&506& & &7.19.23.29.31.281&117&320 \cr
6700&765514085&4.3.25.7.11.23.41.179&177&2&6718&772761521&128.9.5.13.31.281&187&94 \cr
 & &16.9.11.23.41.59&207&3608& & &512.3.5.11.13.17.47&31161&14080 \cr
\noalign{\hrule}
 & &3.5.7.13.17.137.241&1087&2046& & &169.31.43.47.73&1609&3630 \cr
6701&766159485&4.9.5.11.17.31.1087&241&286&6719&772925387&4.3.5.121.73.1609&403&1206 \cr
 & &16.121.13.241.1087&1087&968& & &16.27.5.11.13.31.67&1485&536 \cr
\noalign{\hrule}
 & &9.59.1129.1279&39061&27550& & &5.17.29.31.67.151&6707&28602 \cr
6702&766759221&4.25.11.19.29.53.67&1129&546&6720&773090555&4.9.7.19.227.353&5863&1550 \cr
 & &16.3.7.13.19.29.1129&551&728& & &16.3.25.11.13.31.41&2255&312 \cr
\noalign{\hrule}
 & &11.13.23.431.541&1915&1374& & &19.97.307.1367&635&732 \cr
6703&766899419&4.3.5.229.383.431&407&24&6721&773449967&8.3.5.19.61.127.307&957&6790 \cr
 & &64.9.5.11.37.229&8473&1440& & &32.9.25.7.11.29.97&2233&3600 \cr
\noalign{\hrule}
 & &625.17.257.281&3641&3384& & &81.7.11.19.61.107&265&56 \cr
6704&767305625&16.9.25.11.17.47.331&5441&10116&6722&773469081&16.27.5.49.53.61&955&2278 \cr
 & &128.81.281.5441&5441&5184& & &64.25.17.67.191&28475&6112 \cr
\noalign{\hrule}
 & &27.5.7.11.169.19.23&487&696& & &5.7.11.1229.1637&44403&45632 \cr
6705&767701935&16.81.5.23.29.487&2149&286&6723&774571105&128.3.7.361.23.31.41&297&8600 \cr
 & &64.7.11.13.29.307&307&928& & &2048.81.25.11.43&3483&5120 \cr
\noalign{\hrule}
 & &125.7.43.137.149&2613&2762& & &27.13.47.107.439&31735&15238 \cr
6706&768039125&4.3.7.13.67.137.1381&915&44&6724&774913581&4.5.11.19.401.577&2171&4176 \cr
 & &32.9.5.11.61.1381&15191&8784& & &128.9.13.19.29.167&3173&1856 \cr
\noalign{\hrule}
 & &17.37.41.83.359&8343&4940& & &9.5.113.257.593&2639&326 \cr
6707&768434833&8.81.5.13.17.19.103&407&2158&6725&774959085&4.7.13.29.113.163&627&514 \cr
 & &32.3.11.169.37.83&169&528& & &16.3.11.13.19.29.257&551&1144 \cr
\noalign{\hrule}
 & &3.7.11.17.137.1429&505&454& & &5.343.11.17.41.59&817&4062 \cr
6708&768800571&4.5.11.101.227.1429&19323&3604&6726&775785395&4.3.49.19.43.677&15&34 \cr
 & &32.9.5.17.19.53.113&5989&4560& & &16.9.5.17.43.677&677&3096 \cr
\noalign{\hrule}
 & &9.25.7.19.23.1117&501&616& & &3.125.13.61.2609&5117&2508 \cr
6709&768803175&16.27.5.49.11.19.167&73&4582&6727&775851375&8.9.7.11.13.17.19.43&61&160 \cr
 & &64.11.29.73.79&2117&27808& & &512.5.7.19.43.61&817&1792 \cr
\noalign{\hrule}
 & &121.13.43.83.137&3485&2406& & &11.13.17.29.101.109&1115&738 \cr
6710&769123069&4.3.5.121.17.41.401&3483&1478&6728&776123491&4.9.5.11.41.101.223&221&890 \cr
 & &16.243.17.43.739&4131&5912& & &16.3.25.13.17.41.89&1025&2136 \cr
\noalign{\hrule}
 & &11.127.379.1453&90181&94350& & &5.7.19.31.41.919&5553&11908 \cr
6711&769309739&4.3.25.7.13.17.37.991&537&758&6729&776752585&8.9.7.13.229.617&1727&124 \cr
 & &16.9.5.179.379.991&8055&7928& & &64.3.11.13.31.157&429&5024 \cr
\noalign{\hrule}
 & &11.31.61.71.521&225&296& & &9.25.11.13.19.31.41&23&302 \cr
6712&769449791&16.9.25.11.31.37.61&103&568&6730&776994075&4.11.19.23.41.151&1209&1660 \cr
 & &256.3.5.37.71.103&3811&1920& & &32.3.5.13.23.31.83&83&368 \cr
\noalign{\hrule}
 & &81.5.49.17.2281&4541&11426& & &13.59.277.3659&145&132 \cr
6713&769529565&4.7.19.29.197.239&561&1112&6731&777387481&8.3.5.11.29.59.3659&1505&2154 \cr
 & &64.3.11.17.139.197&2167&4448& & &32.9.25.7.29.43.359&108059&104400 \cr
\noalign{\hrule}
 & &25.37.47.89.199&507&418& & &9.5.11.19.83.997&49&34 \cr
6714&769985725&4.3.11.169.19.47.199&405&206&6732&778273155&4.3.49.11.17.19.997&20303&30544 \cr
 & &16.243.5.11.13.19.103&25029&21736& & &128.23.79.83.257&5911&5056 \cr
\noalign{\hrule}
}%
}
$$
\eject
\vglue -23 pt
\noindent\hskip 1 in\hbox to 6.5 in{\ 6733 -- 6768 \hfill\fbd 778371405 -- 792545247\frb}
\vskip -9 pt
$$
\vbox{
\nointerlineskip
\halign{\strut
    \vrule \ \ \hfil \frb #\ 
   &\vrule \hfil \ \ \fbb #\frb\ 
   &\vrule \hfil \ \ \frb #\ \hfil
   &\vrule \hfil \ \ \frb #\ 
   &\vrule \hfil \ \ \frb #\ \ \vrule \hskip 2 pt
   &\vrule \ \ \hfil \frb #\ 
   &\vrule \hfil \ \ \fbb #\frb\ 
   &\vrule \hfil \ \ \frb #\ \hfil
   &\vrule \hfil \ \ \frb #\ 
   &\vrule \hfil \ \ \frb #\ \vrule \cr%
\noalign{\hrule}
 & &3.5.7.23.31.37.281&137&418& & &5.7.67.197.1699&33847&32148 \cr
6733&778371405&4.7.11.19.23.31.137&3653&594&6751&784878535&8.9.7.11.17.19.47.181&3349&6968 \cr
 & &16.27.121.13.281&117&968& & &128.3.13.289.67.197&867&832 \cr
\noalign{\hrule}
 & &13.181.373.887&1363&990& & &9.7.11.41.43.643&83&1846 \cr
6734&778492403&4.9.5.11.29.47.887&24707&20272&6752&785591037&4.3.7.11.13.71.83&3215&3176 \cr
 & &128.3.7.31.181.797&5579&5952& & &64.5.71.397.643&1985&2272 \cr
\noalign{\hrule}
 & &9.121.17.23.31.59&25777&37990& & &81.49.11.41.439&289&278 \cr
6735&778786371&4.5.29.131.149.173&2783&22302&6753&785818341&4.7.289.41.139.439&729&17270 \cr
 & &16.27.7.121.23.59&7&24& & &16.729.5.11.17.157&785&1224 \cr
\noalign{\hrule}
 & &25.13.43.139.401&341&354& & &9.25.7.19.97.271&187&478 \cr
6736&778952525&4.3.5.11.31.43.59.401&127&12558&6754&786638475&4.3.5.11.17.239.271&13547&1358 \cr
 & &16.9.7.11.13.23.127&8001&2024& & &16.7.19.23.31.97&23&248 \cr
\noalign{\hrule}
 & &3.43.67.89.1013&109&20& & &5.17.19.29.107.157&217&2886 \cr
6737&779226951&8.5.67.109.1013&473&540&6755&786781165&4.3.5.7.13.19.31.37&657&638 \cr
 & &64.27.25.11.43.109&2475&3488& & &16.27.11.13.29.31.73&9207&7592 \cr
\noalign{\hrule}
 & &25.7.19.29.59.137&21403&35178& & &3.7.23.73.137.163&117&44 \cr
6738&779403275&4.3.11.13.17.41.1259&11455&2394&6756&787368729&8.27.11.13.137.163&9805&13504 \cr
 & &16.27.5.7.19.29.79&79&216& & &1024.5.37.53.211&39035&27136 \cr
\noalign{\hrule}
 & &5.11.19.841.887&943&102& & &31.83.359.853&29835&40964 \cr
6739&779535515&4.3.17.23.41.887&95&792&6757&787922071&8.27.5.49.11.13.17.19&965&916 \cr
 & &64.27.5.11.19.23&23&864& & &64.3.25.13.17.193.229&151827&154400 \cr
\noalign{\hrule}
 & &25.49.11.13.61.73&3749&4974& & &9.25.49.19.53.71&3961&1364 \cr
6740&780054275&4.3.23.73.163.829&15483&3584&6758&788254425&8.3.11.17.19.31.233&4657&5356 \cr
 & &4096.9.7.13.397&3573&2048& & &64.11.13.103.4657&51227&42848 \cr
\noalign{\hrule}
 & &121.59.103.1061&67229&61152& & &13.37.41.71.563&1449&1462 \cr
6741&780171337&64.3.49.13.23.37.79&927&890&6759&788307533&4.9.7.17.23.37.43.563&13079&130 \cr
 & &256.27.5.49.13.89.103&31239&31360& & &16.3.5.11.13.29.41.43&1595&1032 \cr
\noalign{\hrule}
 & &27.25.11.17.41.151&637&4412& & &11.13.17.29.53.211&11151&8830 \cr
6742&781458975&8.49.13.41.1103&695&408&6760&788390317&4.27.5.7.17.59.883&703&3712 \cr
 & &128.3.5.7.13.17.139&1807&448& & &1024.9.7.19.29.37&6327&3584 \cr
\noalign{\hrule}
 & &9.5.53.157.2087&1441&646& & &3.25.11.19.71.709&527&7272 \cr
6743&781466715&4.3.11.17.19.131.157&53&104&6761&789063825&16.27.5.17.31.101&1097&1198 \cr
 & &64.11.13.19.53.131&2489&4576& & &64.31.599.1097&34007&19168 \cr
\noalign{\hrule}
 & &81.5.53.79.461&1271&1034& & &9.5.31.43.13159&26899&38896 \cr
6744&781733835&4.27.11.31.41.47.53&127&710&6762&789342615&32.11.13.17.37.727&3725&4452 \cr
 & &16.5.41.47.71.127&9017&15416& & &256.3.25.7.11.53.149&20405&19072 \cr
\noalign{\hrule}
 & &81.11.17.113.457&18247&33394& & &81.19.199.2579&16441&32560 \cr
6745&782206227&4.59.71.257.283&2465&17628&6763&789847119&32.5.11.37.41.401&221&180 \cr
 & &32.3.5.13.17.29.113&377&80& & &256.9.25.11.13.17.37&15725&18304 \cr
\noalign{\hrule}
 & &3.7.11.13.43.73.83&1739&1830& & &5.11.31.47.71.139&3447&862 \cr
6746&782392611&4.9.5.11.37.47.61.73&7&664&6764&790852315&4.9.71.383.431&611&682 \cr
 & &64.5.7.37.47.83&235&1184& & &16.3.11.13.31.47.383&383&312 \cr
\noalign{\hrule}
 & &9.5.7.11.101.2239&1081&1158& & &13.19.23.41.43.79&669&110 \cr
6747&783571635&4.27.5.23.47.101.193&9037&4598&6765&791232637&4.3.5.11.23.79.223&10105&9882 \cr
 & &16.7.121.19.47.1291&9823&10328& & &16.243.25.43.47.61&11421&12200 \cr
\noalign{\hrule}
 & &31.41.53.103.113&665&2838& & &5.19.67.277.449&1261&984 \cr
6748&784037957&4.3.5.7.11.19.43.103&53&156&6766&791634145&16.3.13.19.41.67.97&10263&6286 \cr
 & &32.9.5.7.13.43.53&3913&720& & &64.9.7.11.311.449&2799&2464 \cr
\noalign{\hrule}
 & &3.121.53.83.491&2139&2260& & &3.23.37.41.67.113&4257&3314 \cr
6749&784046967&8.9.5.23.31.113.491&385&106&6767&792479283&4.27.11.37.43.1657&23165&21574 \cr
 & &32.25.7.11.23.53.113&2825&2576& & &16.5.7.11.23.41.67.113&55&56 \cr
\noalign{\hrule}
 & &81.109.181.491&205&286& & &9.13.23.311.947&4997&2156 \cr
6750&784642059&4.5.11.13.41.109.181&4419&3002&6768&792545247&8.3.49.11.13.19.263&1403&4120 \cr
 & &16.9.5.11.19.79.491&869&760& & &128.5.7.23.61.103&6283&2240 \cr
\noalign{\hrule}
}%
}
$$
\eject
\vglue -23 pt
\noindent\hskip 1 in\hbox to 6.5 in{\ 6769 -- 6804 \hfill\fbd 792614251 -- 802798003\frb}
\vskip -9 pt
$$
\vbox{
\nointerlineskip
\halign{\strut
    \vrule \ \ \hfil \frb #\ 
   &\vrule \hfil \ \ \fbb #\frb\ 
   &\vrule \hfil \ \ \frb #\ \hfil
   &\vrule \hfil \ \ \frb #\ 
   &\vrule \hfil \ \ \frb #\ \ \vrule \hskip 2 pt
   &\vrule \ \ \hfil \frb #\ 
   &\vrule \hfil \ \ \fbb #\frb\ 
   &\vrule \hfil \ \ \frb #\ \hfil
   &\vrule \hfil \ \ \frb #\ 
   &\vrule \hfil \ \ \frb #\ \vrule \cr%
\noalign{\hrule}
 & &121.13.47.71.151&3825&1862& & &5.17.61.281.547&38203&47502 \cr
6769&792614251&4.9.25.49.17.19.71&11&60&6787&796970795&4.9.7.11.13.23.29.151&37&1094 \cr
 & &32.27.125.11.17.19&8721&2000& & &16.3.11.23.37.547&1221&184 \cr
\noalign{\hrule}
 & &3.125.7.11.13.2113&4469&2356& & &9.17.79.101.653&11977&11960 \cr
6770&793167375&8.5.11.19.31.41.109&27&182&6788&797173911&16.3.5.7.13.23.29.59.653&293&946 \cr
 & &32.27.7.13.41.109&369&1744& & &64.5.11.13.23.29.43.293&438035&438944 \cr
\noalign{\hrule}
 & &9.5.17.61.89.191&1223&2024& & &3.17.31.43.11731&285&242 \cr
6771&793258335&16.5.11.23.61.1223&459&764&6789&797508573&4.9.5.121.19.11731&16211&4480 \cr
 & &128.27.11.17.23.191&253&192& & &1024.25.7.13.29.43&5075&6656 \cr
\noalign{\hrule}
 & &9.5.11.13.43.47.61&665&622& & &25.11.13.29.43.179&1961&366 \cr
6772&793313235&4.25.7.19.47.61.311&1677&20648&6790&797986475&4.3.5.37.43.53.61&111&154 \cr
 & &64.3.7.13.29.43.89&623&928& & &16.9.7.11.1369.61&9583&4392 \cr
\noalign{\hrule}
 & &5.7.11.53.97.401&1989&1406& & &27.25.11.13.17.487&469&194 \cr
6773&793693285&4.9.13.17.19.37.401&145&256&6791&799130475&4.9.7.67.97.487&12091&18590 \cr
 & &2048.3.5.13.17.19.29&21489&17408& & &16.5.11.169.107.113&1391&904 \cr
\noalign{\hrule}
 & &5.7.11.13.17.19.491&31&354& & &9.5.19.29.167.193&11&182 \cr
6774&793757965&4.3.13.31.59.491&447&44&6792&799167645&4.5.7.11.13.29.167&519&316 \cr
 & &32.9.11.59.149&149&8496& & &32.3.11.13.79.173&1903&16432 \cr
\noalign{\hrule}
 & &3.5.49.23.151.311&9823&14488& & &3.19.31.107.4229&2409&1820 \cr
6775&793875705&16.7.11.19.47.1811&41827&43290&6793&799572801&8.9.5.7.11.13.73.107&737&226 \cr
 & &64.9.5.13.37.151.277&3601&3552& & &32.5.121.13.67.113&37855&25168 \cr
\noalign{\hrule}
 & &11.23.59.83.641&26277&26926& & &3.11.47.61.79.107&2793&5660 \cr
6776&794161181&4.3.19.23.461.13463&12033&1430&6794&799746783&8.9.5.49.11.19.283&67&32 \cr
 & &16.27.5.7.11.13.19.191&17955&19864& & &512.7.19.67.283&18961&34048 \cr
\noalign{\hrule}
 & &81.11.47.61.311&1735&1064& & &9.25.7.11.13.53.67&257&432 \cr
6777&794448567&16.9.5.7.19.47.347&2239&4354&6795&799773975&32.243.11.67.257&247&490 \cr
 & &64.49.311.2239&2239&1568& & &128.5.49.13.19.257&1799&1216 \cr
\noalign{\hrule}
 & &7.31.43.53.1607&695&912& & &3.343.17.19.29.83&28325&1516 \cr
6778&794730601&32.3.5.19.43.53.139&111&154&6796&800007369&8.25.11.103.379&1827&2342 \cr
 & &128.9.7.11.19.37.139&23769&26048& & &32.9.5.7.29.1171&1171&240 \cr
\noalign{\hrule}
 & &9.13.19.43.53.157&335&1342& & &7.31.71.223.233&5&228 \cr
6779&795396069&4.3.5.11.61.67.157&581&424&6797&800532313&8.3.5.7.19.31.71&233&264 \cr
 & &64.7.11.53.61.83&4697&2656& & &128.9.5.11.19.233&1881&320 \cr
\noalign{\hrule}
 & &3.7.11.17.841.241&3051&400& & &27.7.43.137.719&221&940 \cr
6780&795928287&32.81.25.29.113&2587&238&6798&800533881&8.5.7.13.17.47.137&473&486 \cr
 & &128.7.13.17.199&2587&64& & &32.243.5.11.17.43.47&2585&2448 \cr
\noalign{\hrule}
 & &5.11.137.149.709&44469&52664& & &27.13.23.1681.59&121&1478 \cr
6781&796004935&16.729.29.61.227&137&1906&6799&800672067&4.9.121.41.739&14605&15694 \cr
 & &64.81.137.953&953&2592& & &16.5.7.19.23.59.127&665&1016 \cr
\noalign{\hrule}
 & &3.25.7.13.67.1741&7401&1304& & &9.11.13.17.19.41.47&221&202 \cr
6782&796115775&16.9.5.163.2467&6901&5434&6800&801055827&4.11.169.289.41.101&1645&47196 \cr
 & &64.11.13.19.67.103&1133&608& & &32.27.5.7.19.23.47&115&336 \cr
\noalign{\hrule}
 & &3.5.13.17.89.2699&4601&3496& & &3.25.11.13.17.53.83&3577&4002 \cr
6783&796299465&16.19.23.43.89.107&539&450&6801&802047675&4.9.49.23.29.73.83&5035&968 \cr
 & &64.9.25.49.11.19.107&22363&23520& & &64.5.121.19.53.73&803&608 \cr
\noalign{\hrule}
 & &9.5.7.13.19.29.353&473&1292& & &5.13.37.59.5653&102663&106498 \cr
6784&796489785&8.11.17.361.29.43&10473&5050&6802&802132435&4.9.7.11.17.61.7607&2797&4810 \cr
 & &32.3.25.101.3491&3491&8080& & &16.3.5.7.13.17.37.2797&2797&2856 \cr
\noalign{\hrule}
 & &9.11.361.31.719&445&274& & &7.31.73.89.569&145&72 \cr
6785&796586571&4.5.11.19.31.89.137&719&2226&6803&802204081&16.9.5.29.89.569&151&418 \cr
 & &16.3.7.53.89.719&371&712& & &64.3.5.11.19.29.151&31559&13920 \cr
\noalign{\hrule}
 & &9.7.11.13.191.463&437&900& & &7.23.43.61.1901&12681&31042 \cr
6786&796692897&8.81.25.11.13.19.23&1289&574&6804&802798003&4.9.11.17.83.1409&665&2074 \cr
 & &32.5.7.19.41.1289&6445&12464& & &16.3.5.7.289.19.61&1445&456 \cr
\noalign{\hrule}
}%
}
$$
\eject
\vglue -23 pt
\noindent\hskip 1 in\hbox to 6.5 in{\ 6805 -- 6840 \hfill\fbd 802925123 -- 813793617\frb}
\vskip -9 pt
$$
\vbox{
\nointerlineskip
\halign{\strut
    \vrule \ \ \hfil \frb #\ 
   &\vrule \hfil \ \ \fbb #\frb\ 
   &\vrule \hfil \ \ \frb #\ \hfil
   &\vrule \hfil \ \ \frb #\ 
   &\vrule \hfil \ \ \frb #\ \ \vrule \hskip 2 pt
   &\vrule \ \ \hfil \frb #\ 
   &\vrule \hfil \ \ \fbb #\frb\ 
   &\vrule \hfil \ \ \frb #\ \hfil
   &\vrule \hfil \ \ \frb #\ 
   &\vrule \hfil \ \ \frb #\ \vrule \cr%
\noalign{\hrule}
 & &49.11.13.19.37.163&1525&594& & &27.11.13.19.41.269&245&206 \cr
6805&802925123&4.27.25.121.37.61&1007&82&6823&809076411&4.9.5.49.19.103.269&2189&232 \cr
 & &16.3.19.41.53.61&6519&488& & &64.5.49.11.29.199&9751&4640 \cr
\noalign{\hrule}
 & &9.5.121.43.47.73&441&76& & &27.11.19.43.47.71&169&40 \cr
6806&803317185&8.81.49.11.19.43&3379&6862&6824&809719713&16.9.5.169.47.71&1903&1292 \cr
 & &32.31.47.73.109&109&496& & &128.11.13.17.19.173&2249&1088 \cr
\noalign{\hrule}
 & &3.25.11.13.173.433&1007&292& & &5.7.17.29.71.661&1089&970 \cr
6807&803399025&8.5.19.53.73.173&71&936&6825&809794405&4.9.25.121.97.661&221&2204 \cr
 & &128.9.13.71.73&219&4544& & &32.3.121.13.17.19.29&1573&912 \cr
\noalign{\hrule}
 & &5.11.17.29.107.277&1413&1690& & &25.7.107.181.239&507&688 \cr
6808&803661485&4.9.25.11.169.17.157&8033&5992&6826&810024775&32.3.5.7.169.43.107&657&5258 \cr
 & &64.3.7.13.29.107.277&91&96& & &128.27.11.73.239&1971&704 \cr
\noalign{\hrule}
 & &81.25.17.19.1229&10109&10784& & &27.5.7.13.23.47.61&283&1364 \cr
6809&803858175&64.3.11.19.337.919&14659&4550&6827&810085185&8.5.7.11.13.31.283&207&1208 \cr
 & &256.25.7.13.107.137&12467&13696& & &128.9.23.31.151&151&1984 \cr
\noalign{\hrule}
 & &9.5.7.11.13.61.293&745&74& & &23.43.73.103.109&4473&3484 \cr
6810&805089285&4.25.37.149.293&2009&1716&6828&810555719&8.9.7.13.67.71.103&55&158 \cr
 & &32.3.49.11.13.37.41&287&592& & &32.3.5.7.11.13.67.79&30485&41712 \cr
\noalign{\hrule}
 & &27.5.11.17.19.23.73&725&516& & &243.5.7.11.13.23.29&30785&30694 \cr
6811&805340745&8.81.125.23.29.43&3179&304&6829&811215405&4.25.29.47.103.131.149&37&3762 \cr
 & &256.11.289.19.29&493&128& & &16.9.11.19.37.47.103&4841&5624 \cr
\noalign{\hrule}
 & &7.11.59.281.631&319&312& & &3.5.7.11.29.53.457&1103&646 \cr
6812&805523873&16.3.121.13.29.59.281&505&7644&6830&811282395&4.5.7.17.19.29.1103&2483&5238 \cr
 & &128.9.5.49.169.101&17069&20160& & &16.27.13.17.97.191&21437&13752 \cr
\noalign{\hrule}
 & &7.13.29.31.43.229&3429&3212& & &19.43.383.2593&16399&32868 \cr
6813&805573223&8.27.11.13.43.73.127&1561&7710&6831&811378223&8.9.11.529.31.83&1105&1634 \cr
 & &32.81.5.7.223.257&20817&17840& & &32.3.5.13.17.19.31.43&1105&1488 \cr
\noalign{\hrule}
 & &3.11.17.23.197.317&5053&8540& & &3.73.127.163.179&28495&682 \cr
6814&805779447&8.5.7.17.31.61.163&1089&52&6832&811499901&4.5.11.31.41.139&2413&3942 \cr
 & &64.9.5.121.13.31&429&4960& & &16.27.19.73.127&171&8 \cr
\noalign{\hrule}
 & &9.121.13.289.197&335&532& & &17.43.67.73.227&21&22 \cr
6815&806000481&8.3.5.7.121.13.19.67&6647&1460&6833&811597867&4.3.7.11.17.67.73.227&2021&14550 \cr
 & &64.25.289.23.73&1825&736& & &16.9.25.7.43.47.97&8225&6984 \cr
\noalign{\hrule}
 & &5.169.19.23.37.59&927&2284& & &27.5.49.17.31.233&185&32 \cr
6816&806105495&8.9.5.37.103.571&8671&10384&6834&812262465&64.3.25.7.37.233&2203&572 \cr
 & &256.3.11.13.23.29.59&319&384& & &512.11.13.2203&28639&2816 \cr
\noalign{\hrule}
 & &49.19.31.101.277&36531&22958& & &3.7.17.37.227.271&29815&31702 \cr
6817&807444197&4.81.11.13.41.883&1085&202&6835&812578053&4.5.7.121.67.89.131&27&104 \cr
 & &16.9.5.7.31.41.101&205&72& & &64.27.5.11.13.67.89&33165&37024 \cr
\noalign{\hrule}
 & &9.5.1331.97.139&4745&8738& & &5.23.41.97.1777&3971&4914 \cr
6818&807564285&4.3.25.13.17.73.257&1067&4408&6836&812719835&4.27.7.11.13.361.97&395&298 \cr
 & &64.11.17.19.29.97&323&928& & &16.3.5.13.361.79.149&35313&37544 \cr
\noalign{\hrule}
 & &81.25.11.19.23.83&4453&2378& & &5.11.43.47.71.103&59&414 \cr
6819&807936525&4.3.19.29.41.61.73&689&470&6837&812876515&4.9.23.47.59.103&1161&1208 \cr
 & &16.5.13.29.41.47.53&19981&15416& & &64.243.43.59.151&8909&7776 \cr
\noalign{\hrule}
 & &5.7.11.19.29.37.103&23&3834& & &17.29.1069.1543&13287&31460 \cr
6820&808446485&4.27.5.11.23.71&4099&4066&6838&813187231&8.3.5.121.13.43.103&113&102 \cr
 & &16.9.19.107.4099&4099&7704& & &32.9.11.13.17.103.113&14729&16272 \cr
\noalign{\hrule}
 & &9.7.11.47.103.241&155&2806& & &169.17.107.2647&2233&414 \cr
6821&808509933&4.5.23.31.61.103&51&52&6839&813716917&4.9.7.11.169.23.29&113&1070 \cr
 & &32.3.5.13.17.23.31.61&25415&30256& & &16.3.5.23.107.113&113&2760 \cr
\noalign{\hrule}
 & &5.83.107.131.139&42619&31746& & &9.49.13.19.31.241&495&3628 \cr
6822&808570645&4.3.11.13.17.23.37.109&525&4558&6840&813793617&8.81.5.7.11.907&1871&964 \cr
 & &16.9.25.7.11.43.53&13545&4664& & &64.11.241.1871&1871&352 \cr
\noalign{\hrule}
}%
}
$$
\eject
\vglue -23 pt
\noindent\hskip 1 in\hbox to 6.5 in{\ 6841 -- 6876 \hfill\fbd 814164255 -- 823745285\frb}
\vskip -9 pt
$$
\vbox{
\nointerlineskip
\halign{\strut
    \vrule \ \ \hfil \frb #\ 
   &\vrule \hfil \ \ \fbb #\frb\ 
   &\vrule \hfil \ \ \frb #\ \hfil
   &\vrule \hfil \ \ \frb #\ 
   &\vrule \hfil \ \ \frb #\ \ \vrule \hskip 2 pt
   &\vrule \ \ \hfil \frb #\ 
   &\vrule \hfil \ \ \fbb #\frb\ 
   &\vrule \hfil \ \ \frb #\ \hfil
   &\vrule \hfil \ \ \frb #\ 
   &\vrule \hfil \ \ \frb #\ \vrule \cr%
\noalign{\hrule}
 & &9.5.17.61.73.239&481&176& & &27.11.83.167.199&141&58 \cr
6841&814164255&32.11.13.17.37.239&1647&1460&6859&819226683&4.81.11.29.47.167&985&2822 \cr
 & &256.27.5.37.61.73&111&128& & &16.5.17.29.83.197&3349&1160 \cr
\noalign{\hrule}
 & &7.11.169.19.37.89&55977&62830& & &9.49.29.139.461&451&10 \cr
6842&814184371&4.3.5.47.61.103.397&175&222&6860&819506331&4.5.11.29.41.139&1995&2036 \cr
 & &16.9.125.7.37.61.103&7625&7416& & &32.3.25.7.11.19.509&5225&8144 \cr
\noalign{\hrule}
 & &27.5.19.23.37.373&169&682& & &5.7.13.43.163.257&2573&3132 \cr
6843&814189995&4.5.11.169.31.373&3059&1044&6861&819597415&8.27.29.31.83.257&11&268 \cr
 & &32.9.7.13.19.23.29&91&464& & &64.3.11.29.67.83&5829&29216 \cr
\noalign{\hrule}
 & &3.11.13.83.89.257&537&620& & &5.11.13.17.181.373&563&342 \cr
6844&814438911&8.9.5.11.31.179.257&79&178&6862&820620515&4.9.11.19.373.563&12833&5746 \cr
 & &32.5.31.79.89.179&5549&6320& & &16.3.169.17.41.313&1599&2504 \cr
\noalign{\hrule}
 & &11.169.17.19.23.59&2105&1782& & &9.5.13.17.19.43.101&77&94 \cr
6845&814820149&4.81.5.121.59.421&14053&10786&6863&820631565&4.5.7.11.13.43.47.101&27&532 \cr
 & &16.3.5.13.23.47.5393&5393&5640& & &32.27.49.11.19.47&539&2256 \cr
\noalign{\hrule}
 & &5.11.841.67.263&289&1026& & &5.7.13.29.37.1681&267&748 \cr
6846&815059355&4.27.289.19.841&335&506&6864&820689415&8.3.11.17.1681.89&351&1330 \cr
 & &16.3.5.11.289.23.67&289&552& & &32.81.5.7.13.17.19&1377&304 \cr
\noalign{\hrule}
 & &7.23.73.223.311&625&936& & &9.13.17.19.103.211&18179&33938 \cr
6847&815105809&16.9.625.13.23.73&1327&352&6865&821311803&4.343.53.71.239&305&66 \cr
 & &1024.3.25.11.1327&33175&16896& & &16.3.5.49.11.61.71&17395&5368 \cr
\noalign{\hrule}
 & &27.19.53.157.191&4983&5140& & &5.37.47.59.1601&1387&1386 \cr
6848&815316543&8.81.5.11.19.151.257&127&1412&6866&821321005&4.9.5.7.11.19.37.73.1601&91981&115996 \cr
 & &64.11.127.151.353&44831&53152& & &32.3.19.47.59.617.1559&29621&29616 \cr
\noalign{\hrule}
 & &27.25.121.67.149&403&202& & &11.29.67.83.463&1203&740 \cr
6849&815361525&4.9.5.13.31.101.149&703&638&6867&821343017&8.3.5.11.37.83.401&3741&670 \cr
 & &16.11.19.29.31.37.101&21793&23432& & &32.9.25.29.43.67&225&688 \cr
\noalign{\hrule}
 & &3.25.19.23.139.179&12153&12728& & &9.7.11.59.101.199&355&554 \cr
6850&815474775&16.9.19.37.43.4051&575&242&6868&821787813&4.5.7.11.59.71.277&597&184 \cr
 & &64.25.121.23.4051&4051&3872& & &64.3.5.23.199.277&1385&736 \cr
\noalign{\hrule}
 & &3.11.23.37.71.409&169&238& & &3.7.11.19.37.61.83&845&68 \cr
6851&815502237&4.7.169.17.71.409&2035&828&6869&822195759&8.5.169.17.19.61&99&1136 \cr
 & &32.9.5.11.169.23.37&169&240& & &256.9.11.13.71&71&4992 \cr
\noalign{\hrule}
 & &3.5.49.13.23.47.79&1349&954& & &9.5.11.17.239.409&7627&3128 \cr
6852&815987445&4.27.13.19.23.53.71&43099&43450&6870&822574665&16.289.23.29.263&189&478 \cr
 & &16.25.7.11.19.47.79.131&1045&1048& & &64.27.7.239.263&263&672 \cr
\noalign{\hrule}
 & &49.11.13.37.47.67&533&204& & &3.125.11.13.529.29&939&686 \cr
6853&816406591&8.3.7.169.17.37.41&45&214&6871&822661125&4.9.343.23.29.313&205&2 \cr
 & &32.27.5.17.41.107&5535&29104& & &16.5.49.41.313&12833&392 \cr
\noalign{\hrule}
 & &7.13.17.29.109.167&495&2666& & &3.11.17.43.67.509&25&42 \cr
6854&816641189&4.9.5.7.11.17.31.43&167&218&6872&822666669&4.9.25.7.11.43.509&1445&4154 \cr
 & &16.3.31.43.109.167&129&248& & &16.125.289.31.67&527&1000 \cr
\noalign{\hrule}
 & &9.5.11.17.37.43.61&551&78& & &9.5.11.13.19.53.127&5831&5246 \cr
6855&816684165&4.27.5.13.19.29.61&3487&3182&6873&822965715&4.343.17.43.61.127&899&6360 \cr
 & &16.11.29.37.43.317&317&232& & &64.3.5.49.29.31.53&1421&992 \cr
\noalign{\hrule}
 & &5.7.121.43.4489&1017&1864& & &27.25.17.179.401&5803&9412 \cr
6856&817469345&16.9.5.67.113.233&111&224&6874&823664025&8.3.5.7.13.181.829&4411&1924 \cr
 & &1024.27.7.37.233&8621&13824& & &64.11.169.37.401&1859&1184 \cr
\noalign{\hrule}
 & &9.13.17.19.59.367&4697&4330& & &3.11.169.113.1307&85&254 \cr
6857&818288523&4.5.7.11.13.19.61.433&83&5712&6875&823672707&4.5.11.17.127.1307&45&1352 \cr
 & &128.3.49.11.17.83&4067&704& & &64.9.25.169.17&1275&32 \cr
\noalign{\hrule}
 & &27.19.31.53.971&31405&12956& & &5.11.17.19.89.521&7397&2502 \cr
6858&818416089&8.5.11.41.79.571&265&186&6876&823745285&4.9.13.17.139.569&5045&4628 \cr
 & &32.3.25.31.53.571&571&400& & &32.3.5.169.89.1009&3027&2704 \cr
\noalign{\hrule}
}%
}
$$
\eject
\vglue -23 pt
\noindent\hskip 1 in\hbox to 6.5 in{\ 6877 -- 6912 \hfill\fbd 823773405 -- 835016413\frb}
\vskip -9 pt
$$
\vbox{
\nointerlineskip
\halign{\strut
    \vrule \ \ \hfil \frb #\ 
   &\vrule \hfil \ \ \fbb #\frb\ 
   &\vrule \hfil \ \ \frb #\ \hfil
   &\vrule \hfil \ \ \frb #\ 
   &\vrule \hfil \ \ \frb #\ \ \vrule \hskip 2 pt
   &\vrule \ \ \hfil \frb #\ 
   &\vrule \hfil \ \ \fbb #\frb\ 
   &\vrule \hfil \ \ \frb #\ \hfil
   &\vrule \hfil \ \ \frb #\ 
   &\vrule \hfil \ \ \frb #\ \vrule \cr%
\noalign{\hrule}
 & &3.5.7.13.19.23.1381&5699&12254& & &25.13.19.29.41.113&181&294 \cr
6877&823773405&4.7.11.41.139.557&621&908&6895&829654475&4.3.49.13.29.41.181&7095&8362 \cr
 & &32.27.23.227.557&5013&3632& & &16.9.5.7.11.37.43.113&2849&3096 \cr
\noalign{\hrule}
 & &3.11.13.17.37.43.71&279&350& & &125.17.31.43.293&39237&26638 \cr
6878&823824573&4.27.25.7.11.13.31.43&1207&2368&6896&829959125&4.3.11.19.29.41.701&125&84 \cr
 & &512.7.17.31.37.71&217&256& & &32.9.125.7.29.701&4907&4176 \cr
\noalign{\hrule}
 & &9.5.11.13.19.23.293&1645&3922& & &27.5.11.23.109.223&247&868 \cr
6879&823943835&4.25.7.13.37.47.53&6177&5848&6897&830205585&8.7.11.13.19.31.109&555&446 \cr
 & &64.3.17.29.43.53.71&51901&49184& & &32.3.5.19.31.37.223&703&496 \cr
\noalign{\hrule}
 & &5.29.1451.3917&1233&2684& & &3.25.7.19.31.2687&11651&8964 \cr
6880&824117215&8.9.5.11.29.61.137&6683&2162&6898&830887575&8.81.5.61.83.191&12617&12698 \cr
 & &32.3.23.41.47.163&11247&30832& & &32.7.11.31.37.191.907&33559&33616 \cr
\noalign{\hrule}
 & &5.11.13.17.19.43.83&1439&1356& & &3.11.17.841.41.43&137&50 \cr
6881&824242705&8.3.11.17.19.113.1439&3105&24236&6899&831785163&4.25.29.41.43.137&487&702 \cr
 & &64.81.5.23.73.83&1863&2336& & &16.27.5.13.137.487&21915&14248 \cr
\noalign{\hrule}
 & &7.121.17.19.23.131&11695&9396& & &7.11.13.529.1571&12079&5202 \cr
6882&824299553&8.81.5.17.29.2339&1049&3388&6900&831890059&4.9.7.289.47.257&6099&5980 \cr
 & &64.9.5.7.121.1049&1049&1440& & &32.27.5.13.17.19.23.107&8721&8560 \cr
\noalign{\hrule}
 & &25.19.643.2699&1397&1302& & &5.11.13.29.53.757&801&736 \cr
6883&824342075&4.3.5.7.11.31.127.643&361&3576&6901&831908935&64.9.11.23.89.757&1625&646 \cr
 & &64.9.7.11.361.149&9387&6688& & &256.3.125.13.17.19.23&8075&8832 \cr
\noalign{\hrule}
 & &3.13.17.29.53.809&5505&5012& & &9.5.7.19.43.53.61&11&806 \cr
6884&824396079&8.9.5.7.53.179.367&13753&5698&6902&832028715&4.3.7.11.13.31.61&215&212 \cr
 & &32.49.11.17.37.809&539&592& & &32.5.11.13.31.43.53&341&208 \cr
\noalign{\hrule}
 & &9.73.463.2711&18255&15544& & &7.11.13.19.67.653&111&358 \cr
6885&824661801&16.27.5.29.67.1217&29273&23188&6903&832100269&4.3.11.37.179.653&13065&11096 \cr
 & &128.11.17.31.73.401&12431&11968& & &64.9.5.13.19.67.73&365&288 \cr
\noalign{\hrule}
 & &9.5.49.13.47.613&4351&1166& & &9.11.29.239.1213&29419&5758 \cr
6886&825867315&4.11.19.47.53.229&1141&1350&6904&832322997&4.13.31.73.2879&965&1914 \cr
 & &16.27.25.7.163.229&2445&1832& & &16.3.5.11.29.31.193&965&248 \cr
\noalign{\hrule}
 & &9.5.49.19.109.181&1739&976& & &25.23.29.107.467&8109&5434 \cr
6887&826546455&32.3.7.19.37.47.61&49775&49348&6905&833233075&4.9.11.13.17.19.23.53&11861&11300 \cr
 & &256.25.11.169.73.181&9295&9344& & &32.3.25.13.29.113.409&5317&5424 \cr
\noalign{\hrule}
 & &27.47.373.1747&14819&32350& & &7.11.41.107.2467&9223&8046 \cr
6888&826919739&4.25.7.29.73.647&141&506&6906&833350133&4.27.23.41.149.401&385&16 \cr
 & &16.3.5.7.11.23.29.47&805&2552& & &128.3.5.7.11.23.149&3427&960 \cr
\noalign{\hrule}
 & &13.29.37.101.587&363&950& & &5.7.11.169.23.557&6597&2698 \cr
6889&826994363&4.3.25.121.19.29.37&1581&1174&6907&833547715&4.9.19.23.71.733&1183&450 \cr
 & &16.9.5.11.17.31.587&1705&1224& & &16.81.25.7.169.19&405&152 \cr
\noalign{\hrule}
 & &3.11.17.23.61.1051&7725&7708& & &11.13.29.53.3793&51209&58788 \cr
6890&827224233&8.9.25.41.47.103.1051&989&62&6908&833667263&8.9.23.41.71.1249&1825&3074 \cr
 & &32.25.23.31.41.43.47&31775&32336& & &32.3.25.29.41.53.73&1825&1968 \cr
\noalign{\hrule}
 & &3.5.7.19.29.79.181&40183&40002& & &9.11.13.17.53.719&7025&554 \cr
6891&827268645&4.9.11.13.19.59.113.281&4525&2378&6909&833743053&4.25.17.277.281&3081&1696 \cr
 & &16.25.11.29.41.181.281&2255&2248& & &256.3.5.13.53.79&79&640 \cr
\noalign{\hrule}
 & &5.11.47.97.3301&33327&16822& & &27.11.13.23.41.229&703&4564 \cr
6892&827709245&4.9.7.13.529.647&1117&470&6910&833771367&8.7.19.37.41.163&495&208 \cr
 & &16.3.5.7.13.47.1117&1117&2184& & &256.9.5.11.13.163&815&128 \cr
\noalign{\hrule}
 & &49.11.53.79.367&249&620& & &9.17.89.197.311&275&3074 \cr
6893&828243031&8.3.5.7.31.83.367&107&474&6911&834272739&4.25.11.29.53.89&357&622 \cr
 & &32.9.5.31.79.107&3317&720& & &16.3.5.7.17.29.311&35&232 \cr
\noalign{\hrule}
 & &3.5.7.13.17.31.1153&1081&2234& & &7.11.73.149.997&243&754 \cr
6894&829416315&4.7.23.31.47.1117&3553&4266&6912&835016413&4.243.11.13.29.149&355&1994 \cr
 & &16.27.11.17.19.47.79&8037&6952& & &16.3.5.13.71.997&923&120 \cr
\noalign{\hrule}
}%
}
$$
\eject
\vglue -23 pt
\noindent\hskip 1 in\hbox to 6.5 in{\ 6913 -- 6948 \hfill\fbd 835599765 -- 849309813\frb}
\vskip -9 pt
$$
\vbox{
\nointerlineskip
\halign{\strut
    \vrule \ \ \hfil \frb #\ 
   &\vrule \hfil \ \ \fbb #\frb\ 
   &\vrule \hfil \ \ \frb #\ \hfil
   &\vrule \hfil \ \ \frb #\ 
   &\vrule \hfil \ \ \frb #\ \ \vrule \hskip 2 pt
   &\vrule \ \ \hfil \frb #\ 
   &\vrule \hfil \ \ \fbb #\frb\ 
   &\vrule \hfil \ \ \frb #\ \hfil
   &\vrule \hfil \ \ \frb #\ 
   &\vrule \hfil \ \ \frb #\ \vrule \cr%
\noalign{\hrule}
 & &3.5.7.11.13.19.29.101&237&314& & &7.37.61.71.751&8275&7524 \cr
6913&835599765&4.9.5.13.79.101.157&343&242&6931&842418479&8.9.25.11.19.71.331&153&1502 \cr
 & &16.343.121.79.157&7693&6952& & &32.81.5.11.17.751&1377&880 \cr
\noalign{\hrule}
 & &13.17.19.41.43.113&725&1422& & &13.17.29.47.2797&14915&32634 \cr
6914&836520581&4.9.25.13.29.43.79&191&836&6932&842520731&4.9.5.49.19.37.157&11891&11188 \cr
 & &32.3.5.11.19.29.191&2101&6960& & &32.3.5.11.23.47.2797&253&240 \cr
\noalign{\hrule}
 & &5.17.41.271.887&4521&86& & &3.23.61.191.1049&429&620 \cr
6915&837713845&4.3.11.41.43.137&813&950&6933&843311031&8.9.5.11.13.23.31.61&1049&842 \cr
 & &16.9.25.11.19.271&95&792& & &32.5.11.13.421.1049&2105&2288 \cr
\noalign{\hrule}
 & &27.19.23.29.31.79&1751&250& & &121.19.59.6221&3215&3006 \cr
6916&837976779&4.9.125.17.31.103&3509&316&6934&843822661&4.9.5.11.59.167.643&3&646 \cr
 & &32.5.121.29.79&605&16& & &16.27.5.17.19.167&459&6680 \cr
\noalign{\hrule}
 & &27.5.7.11.19.31.137&241&272& & &5.13.101.173.743&753&1496 \cr
6917&838803735&32.5.7.11.17.137.241&261&124&6935&843858535&16.3.5.11.17.101.251&231&1486 \cr
 & &256.9.17.29.31.241&4097&3712& & &64.9.7.121.743&1089&224 \cr
\noalign{\hrule}
 & &3.5.17.23.239.599&301&1496& & &27.25.13.19.61.83&7511&374 \cr
6918&839639265&16.7.11.289.23.43&639&350&6936&844128675&4.3.5.7.11.17.29.37&19&16 \cr
 & &64.9.25.49.11.71&2343&7840& & &128.11.17.19.29.37&6919&1856 \cr
\noalign{\hrule}
 & &3.361.67.71.163&53&110& & &81.19.59.71.131&715&634 \cr
6919&839748453&4.5.11.19.53.67.71&1043&306&6937&844540101&4.5.11.13.59.131.317&213&82 \cr
 & &16.9.5.7.17.53.149&6307&17880& & &16.3.11.13.41.71.317&4121&3608 \cr
\noalign{\hrule}
 & &9.25.49.11.169.41&1121&1154& & &3.11.23.37.67.449&1501&1052 \cr
6920&840314475&4.3.7.13.19.41.59.577&115&11078&6938&844820889&8.11.19.67.79.263&1413&1480 \cr
 & &16.5.23.29.59.191&4393&13688& & &128.9.5.19.37.79.157&14915&15168 \cr
\noalign{\hrule}
 & &25.7.121.13.43.71&999&76& & &13.17.23.43.53.73&1289&990 \cr
6921&840414575&8.27.7.121.19.37&71&50&6939&845643461&4.9.5.11.17.73.1289&1213&1196 \cr
 & &32.9.25.19.37.71&703&144& & &32.3.5.13.23.1213.1289&19335&19408 \cr
\noalign{\hrule}
 & &9.25.11.37.67.137&893&782& & &9.7.11.13.37.43.59&25&34 \cr
6922&840566925&4.3.11.17.19.23.47.137&2077&526&6940&845665821&4.25.7.11.13.17.37.43&4189&12366 \cr
 & &16.17.23.31.67.263&4471&5704& & &16.27.5.59.71.229&1145&1704 \cr
\noalign{\hrule}
 & &9.25.13.263.1093&7097&7112& & &9.1331.31.43.53&2665&1334 \cr
6923&840817575&16.3.5.7.47.127.151.263&487&12848&6941&846304371&4.3.5.13.23.29.41.53&9331&8954 \cr
 & &512.11.73.151.487&121253&124672& & &16.7.121.31.37.41.43&287&296 \cr
\noalign{\hrule}
 & &3.125.343.13.503&21689&21186& & &5.11.23.37.79.229&23517&23288 \cr
6924&841078875&4.27.11.13.529.41.107&385&1006&6942&846749255&16.27.13.41.67.71.79&5773&2534 \cr
 & &16.5.7.121.23.41.503&943&968& & &64.3.7.23.67.181.251&45431&45024 \cr
\noalign{\hrule}
 & &11.37.739.2797&29055&1712& & &7.13.37.43.5849&20493&20450 \cr
6925&841262081&32.3.5.13.107.149&53&54&6943&846824069&4.81.25.11.13.23.37.409&1505&2356 \cr
 & &128.81.5.13.53.149&55809&47680& & &32.3.125.7.19.31.43.409&38037&38000 \cr
\noalign{\hrule}
 & &27.121.43.53.113&1843&1886& & &27.7.11.23.37.479&95303&99650 \cr
6926&841340709&4.9.11.19.23.41.53.97&655&9718&6944&847460691&4.25.13.1993.7331&2669&4662 \cr
 & &16.5.43.97.113.131&655&776& & &16.9.25.7.13.17.37.157&2669&2600 \cr
\noalign{\hrule}
 & &9.7.11.13.29.3221&1957&1264& & &11.13.19.131.2381&2435&54 \cr
6927&841521681&32.13.19.29.79.103&1815&476&6945&847462187&4.27.5.11.13.487&2549&1834 \cr
 & &256.3.5.7.121.17.19&1045&2176& & &16.3.7.131.2549&2549&168 \cr
\noalign{\hrule}
 & &9.25.13.53.61.89&323&262& & &9.125.17.23.41.47&4697&1178 \cr
6928&841630725&4.5.17.19.53.89.131&177&88&6946&847639125&4.7.11.19.31.41.61&925&966 \cr
 & &64.3.11.17.19.59.131&27379&32096& & &16.3.25.49.11.19.23.37&1813&1672 \cr
\noalign{\hrule}
 & &3.17.19.529.31.53&383&330& & &7.11.23.53.83.109&98719&94320 \cr
6929&842203443&4.9.5.11.17.19.23.383&265&58&6947&849178561&32.9.5.17.131.5807&437&6244 \cr
 & &16.25.11.29.53.383&4213&5800& & &256.3.5.7.19.23.223&3345&2432 \cr
\noalign{\hrule}
 & &27.5.11.29.31.631&1407&1748& & &27.11.43.73.911&29563&4966 \cr
6930&842394465&8.81.7.19.23.29.67&379&5048&6948&849309813&4.13.17.37.47.191&965&774 \cr
 & &128.19.379.631&379&1216& & &16.9.5.13.17.43.193&965&1768 \cr
\noalign{\hrule}
}%
}
$$
\eject
\vglue -23 pt
\noindent\hskip 1 in\hbox to 6.5 in{\ 6949 -- 6984 \hfill\fbd 849808575 -- 862097775\frb}
\vskip -9 pt
$$
\vbox{
\nointerlineskip
\halign{\strut
    \vrule \ \ \hfil \frb #\ 
   &\vrule \hfil \ \ \fbb #\frb\ 
   &\vrule \hfil \ \ \frb #\ \hfil
   &\vrule \hfil \ \ \frb #\ 
   &\vrule \hfil \ \ \frb #\ \ \vrule \hskip 2 pt
   &\vrule \ \ \hfil \frb #\ 
   &\vrule \hfil \ \ \fbb #\frb\ 
   &\vrule \hfil \ \ \frb #\ \hfil
   &\vrule \hfil \ \ \frb #\ 
   &\vrule \hfil \ \ \frb #\ \vrule \cr%
\noalign{\hrule}
 & &9.25.7.11.181.271&1777&148& & &9.5.19.529.31.61&517&5278 \cr
6949&849808575&8.37.271.1777&4125&5902&6967&855289845&4.7.11.13.29.31.47&201&202 \cr
 & &32.3.125.11.13.227&227&1040& & &16.3.7.11.29.47.67.101&74437&76328 \cr
\noalign{\hrule}
 & &9.13.29.31.59.137&2693&5390& & &25.13.17.37.53.79&243&22 \cr
6950&850194189&4.3.5.49.11.13.2693&2419&274&6968&855927475&4.243.5.11.37.79&1301&1696 \cr
 & &16.49.41.59.137&49&328& & &256.3.11.53.1301&1301&4224 \cr
\noalign{\hrule}
 & &27.5.29.83.2617&43549&32344& & &9.5.11.169.29.353&245&3638 \cr
6951&850381065&16.11.13.37.107.311&1433&2610&6969&856376235&4.25.49.13.17.107&33&58 \cr
 & &64.9.5.29.37.1433&1433&1184& & &16.3.7.11.17.29.107&107&952 \cr
\noalign{\hrule}
 & &9.5.7.19.23.37.167&1249&1082& & &25.13.103.157.163&1017&3058 \cr
6952&850570245&4.5.19.23.541.1249&65549&53106&6970&856658725&4.9.11.103.113.139&157&260 \cr
 & &16.3.11.53.59.101.167&5353&5192& & &32.3.5.11.13.113.157&339&176 \cr
\noalign{\hrule}
 & &5.13.289.45281&22649&22632& & &27.7.19.61.3911&5287&6446 \cr
6953&850603585&16.3.5.11.13.17.23.29.41.71&32033&12&6971&856708461&4.9.7.11.17.293.311&247&1930 \cr
 & &128.9.23.103.311&32033&13248& & &16.5.13.19.193.293&3809&7720 \cr
\noalign{\hrule}
 & &5.7.11.31.37.41.47&289&1746& & &3.11.37.269.2609&29279&580 \cr
6954&850953565&4.9.7.289.41.97&2585&2294&6972&856923441&8.5.19.23.29.67&185&252 \cr
 & &16.3.5.11.17.31.37.47&17&24& & &64.9.25.7.29.37&175&2784 \cr
\noalign{\hrule}
 & &9.25.11.17.113.179&2263&2212& & &9.5.11.419.4133&181&3952 \cr
6955&851051025&8.3.7.11.31.73.79.113&5063&68500&6973&857204865&32.5.11.13.19.181&12989&12894 \cr
 & &64.125.61.83.137&11371&9760& & &128.3.7.31.307.419&2149&1984 \cr
\noalign{\hrule}
 & &17.41.71.17207&147715&144804& & &9.7.53.491.523&539&1030 \cr
6956&851522809&8.3.5.11.31.953.1097&1025&72&6974&857431827&4.3.5.343.11.53.103&703&1046 \cr
 & &128.27.125.11.31.41&9207&8000& & &16.5.19.37.103.523&1957&1480 \cr
\noalign{\hrule}
 & &3.5.11.19.23.53.223&139&298& & &3.25.11.169.47.131&1593&1462 \cr
6957&852208995&4.5.11.139.149.223&207&1322&6975&858439725&4.81.5.11.13.17.43.59&2527&50042 \cr
 & &16.9.23.149.661&661&3576& & &16.7.361.131.191&1337&2888 \cr
\noalign{\hrule}
 & &9.41.47.211.233&4537&4114& & &9.13.29.409.619&335&74 \cr
6958&852633909&4.121.13.17.233.349&987&3550&6976&859009203&4.5.13.37.67.619&69&550 \cr
 & &16.3.25.7.11.17.47.71&4675&3976& & &16.3.125.11.23.67&16951&1000 \cr
\noalign{\hrule}
 & &3.25.11.17.71.857&77467&82792& & &243.25.11.13.23.43&17869&16880 \cr
6959&853379175&16.13.59.79.101.131&4285&8946&6977&859169025&32.125.107.167.211&851&21726 \cr
 & &64.9.5.7.13.71.857&91&96& & &128.9.17.23.37.71&1207&2368 \cr
\noalign{\hrule}
 & &9.11.29.53.71.79&115&754& & &5.13.139.251.379&18183&34498 \cr
6960&853482267&4.5.13.23.841.53&22137&22436&6978&859490515&4.3.11.19.29.47.367&169&198 \cr
 & &32.3.5.47.71.79.157&785&752& & &16.27.121.169.19.47&16497&18392 \cr
\noalign{\hrule}
 & &3.25.23.167.2963&40337&27812& & &25.7.41.47.2549&2313&4862 \cr
6961&853566225&8.11.17.19.193.409&1643&5310&6979&859586525&4.9.11.13.17.47.257&5759&6320 \cr
 & &32.9.5.11.31.53.59&5487&9328& & &128.3.5.169.79.443&34997&32448 \cr
\noalign{\hrule}
 & &17.19.113.149.157&81&2750& & &5.17.19.29.31.593&3663&6418 \cr
6962&853821107&4.81.125.11.113&157&182&6980&860967805&4.9.11.31.37.3209&3325&116 \cr
 & &16.27.5.7.11.13.157&455&2376& & &32.3.25.7.11.19.29&33&560 \cr
\noalign{\hrule}
 & &11.17.1697.2693&7987&10680& & &49.67.73.3593&1755&1822 \cr
6963&854593927&16.3.5.49.17.89.163&1697&1074&6981&861094787&4.27.5.13.911.3593&2303&5896 \cr
 & &64.9.5.7.179.1697&1253&1440& & &64.3.5.49.11.13.47.67&1833&1760 \cr
\noalign{\hrule}
 & &27.5.49.13.61.163&701&946& & &3.25.17.541.1249&233&308 \cr
6964&855048285&4.11.13.43.163.701&5453&3660&6982&861528975&8.7.11.17.233.1249&565&684 \cr
 & &32.3.5.7.19.41.43.61&779&688& & &64.9.5.11.19.113.233&23617&22368 \cr
\noalign{\hrule}
 & &3.5.7.121.13.31.167&603&482& & &9.7.11.19.29.37.61&335&274 \cr
6965&855059205&4.27.13.67.167.241&8365&31882&6983&861819651&4.3.5.11.19.37.67.137&1247&26 \cr
 & &16.5.7.19.239.839&4541&6712& & &16.5.13.29.43.137&559&5480 \cr
\noalign{\hrule}
 & &625.19.23.31.101&30069&41944& & &3.25.7.121.41.331&557&436 \cr
6966&855154375&16.9.49.13.107.257&11&760&6984&862097775&8.25.7.41.109.557&1881&20956 \cr
 & &256.3.5.7.11.13.19&3003&128& & &64.9.11.169.19.31&1767&5408 \cr
\noalign{\hrule}
}%
}
$$
\eject
\vglue -23 pt
\noindent\hskip 1 in\hbox to 6.5 in{\ 6985 -- 7020 \hfill\fbd 862326465 -- 871881219\frb}
\vskip -9 pt
$$
\vbox{
\nointerlineskip
\halign{\strut
    \vrule \ \ \hfil \frb #\ 
   &\vrule \hfil \ \ \fbb #\frb\ 
   &\vrule \hfil \ \ \frb #\ \hfil
   &\vrule \hfil \ \ \frb #\ 
   &\vrule \hfil \ \ \frb #\ \ \vrule \hskip 2 pt
   &\vrule \ \ \hfil \frb #\ 
   &\vrule \hfil \ \ \fbb #\frb\ 
   &\vrule \hfil \ \ \frb #\ \hfil
   &\vrule \hfil \ \ \frb #\ 
   &\vrule \hfil \ \ \frb #\ \vrule \cr%
\noalign{\hrule}
 & &3.5.7.121.13.23.227&817&2588& & &5.11.31.47.79.137&45&92 \cr
6985&862326465&8.11.13.19.43.647&395&252&7003&867301105&8.9.25.11.23.31.79&47&822 \cr
 & &64.9.5.7.19.43.79&3397&1824& & &32.27.23.47.137&27&368 \cr
\noalign{\hrule}
 & &5.7.11.13.19.47.193&73&60& & &3.5.7.11.37.53.383&117&2798 \cr
6986&862606745&8.3.25.11.47.73.193&14573&5502&7004&867477765&4.27.13.37.1399&9593&8594 \cr
 & &32.9.7.13.19.59.131&1179&944& & &16.53.181.4297&4297&1448 \cr
\noalign{\hrule}
 & &81.5.7.11.89.311&1753&29432& & &19.37.83.107.139&2613&2530 \cr
6987&863169615&16.13.283.1753&735&1018&7005&867824677&4.3.5.11.13.19.23.67.107&59&1332 \cr
 & &64.3.5.49.13.509&509&2912& & &32.27.5.11.23.37.59&6785&4752 \cr
\noalign{\hrule}
 & &27.5.13.23.73.293&1793&844& & &3.11.17.19.31.37.71&1585&4212 \cr
6988&863366985&8.3.5.11.23.163.211&3457&292&7006&868036983&8.243.5.13.19.317&4369&248 \cr
 & &64.11.73.3457&3457&352& & &128.5.17.31.257&1285&64 \cr
\noalign{\hrule}
 & &3.7.11.17.97.2267&3439&3362& & &3.5.11.361.61.239&10019&9836 \cr
6989&863543373&4.17.19.1681.97.181&1665&16&7007&868398135&8.43.233.239.2459&25025&80712 \cr
 & &128.9.5.19.37.181&6697&18240& & &128.9.25.7.11.13.19.59&2301&2240 \cr
\noalign{\hrule}
 & &9.5.13.19.23.31.109&9559&7336& & &3.5.19.61.107.467&3311&3216 \cr
6990&863824455&16.7.121.23.79.131&1457&1326&7008&868711065&32.9.7.11.43.67.467&427&40 \cr
 & &64.3.7.13.17.31.47.79&3713&3808& & &512.5.49.11.61.67&3283&2816 \cr
\noalign{\hrule}
 & &3.5.13.43.269.383&27&242& & &9.25.13.23.37.349&913&62 \cr
6991&863881395&4.81.121.13.383&5983&5600&7009&868722075&4.3.11.31.83.349&299&50 \cr
 & &256.25.7.11.31.193&14861&19840& & &16.25.11.13.23.31&341&8 \cr
\noalign{\hrule}
 & &3.13.17.19.181.379&181&142& & &9.5.37.43.61.199&161&2096 \cr
6992&864141603&4.71.32761.379&29835&2926&7010&869091705&32.7.23.131.199&165&34 \cr
 & &16.27.5.7.11.13.17.19&99&280& & &128.3.5.7.11.17.23&253&7616 \cr
\noalign{\hrule}
 & &9.7.23.31.71.271&4925&3476& & &5.13.2111.6337&65439&71776 \cr
6993&864286479&8.25.11.71.79.197&403&378&7011&869531455&64.9.11.661.2243&2113&130 \cr
 & &32.27.7.13.31.79.197&3081&3152& & &256.3.5.11.13.2113&2113&4224 \cr
\noalign{\hrule}
 & &9.25.13.19.79.197&493&98& & &49.11.53.83.367&107&474 \cr
6994&864913725&4.3.5.49.13.17.19.29&2449&4334&7012&870179387&4.3.7.11.53.79.107&249&620 \cr
 & &16.7.11.31.79.197&77&248& & &32.9.5.31.83.107&3317&720 \cr
\noalign{\hrule}
 & &5.11.29.41.101.131&1239&1690& & &25.13.23.173.673&17451&34276 \cr
6995&865241245&4.3.25.7.169.59.131&2727&5002&7013&870306775&8.9.7.11.19.41.277&125&84 \cr
 & &16.81.13.41.61.101&793&648& & &64.27.125.49.277&6615&8864 \cr
\noalign{\hrule}
 & &125.7.13.29.43.61&5893&342& & &9.7.17.61.67.199&159&40 \cr
6996&865262125&4.9.25.19.71.83&671&754&7014&871058223&16.27.5.53.61.67&935&874 \cr
 & &16.3.11.13.29.61.71&213&88& & &64.25.11.17.19.23.53&13409&15200 \cr
\noalign{\hrule}
 & &27.11.19.23.59.113&1345&2588& & &9.361.31.41.211&1177&94 \cr
6997&865303263&8.3.5.59.269.647&119&766&7015&871320069&4.3.11.47.107.211&55&266 \cr
 & &32.7.17.269.383&6511&30128& & &16.5.7.121.19.47&235&6776 \cr
\noalign{\hrule}
 & &5.7.23.41.157.167&2223&1388& & &9.25.11.17.139.149&361&86 \cr
6998&865358095&8.9.7.13.19.41.347&803&926&7016&871415325&4.3.17.361.43.139&3725&3364 \cr
 & &32.3.11.73.347.463&75993&81488& & &32.25.841.43.149&841&688 \cr
\noalign{\hrule}
 & &9.7.13.47.113.199&203&814& & &9.7.17.43.127.149&21275&36388 \cr
6999&865592091&4.49.11.29.37.199&611&810&7017&871460919&8.25.11.23.37.827&2749&3576 \cr
 & &16.81.5.11.13.37.47&495&296& & &128.3.37.149.2749&2749&2368 \cr
\noalign{\hrule}
 & &27.5.7.11.19.41.107&1139&1426& & &25.343.37.41.67&7719&4972 \cr
7000&866454435&4.11.17.23.31.67.107&10775&9234&7018&871554425&8.3.25.11.31.83.113&469&444 \cr
 & &16.243.25.19.31.431&2155&2232& & &64.9.7.31.37.67.113&1017&992 \cr
\noalign{\hrule}
 & &27.7.11.157.2657&223&29450& & &9.121.19.103.409&155&254 \cr
7001&867252771&4.25.19.31.223&567&548&7019&871649757&4.5.11.19.31.103.127&1677&280 \cr
 & &32.81.5.7.31.137&465&2192& & &64.3.25.7.13.31.43&5425&17888 \cr
\noalign{\hrule}
 & &7.17.41.149.1193&513&530& & &27.11.83.113.313&499&1742 \cr
7002&867276403&4.27.5.19.41.53.1193&25967&22946&7020&871881219&4.13.67.313.499&10735&10236 \cr
 & &16.9.5.7.11.23.149.1129&10161&10120& & &32.3.5.13.19.113.853&4265&3952 \cr
\noalign{\hrule}
}%
}
$$
\eject
\vglue -23 pt
\noindent\hskip 1 in\hbox to 6.5 in{\ 7021 -- 7056 \hfill\fbd 872243645 -- 885247069\frb}
\vskip -9 pt
$$
\vbox{
\nointerlineskip
\halign{\strut
    \vrule \ \ \hfil \frb #\ 
   &\vrule \hfil \ \ \fbb #\frb\ 
   &\vrule \hfil \ \ \frb #\ \hfil
   &\vrule \hfil \ \ \frb #\ 
   &\vrule \hfil \ \ \frb #\ \ \vrule \hskip 2 pt
   &\vrule \ \ \hfil \frb #\ 
   &\vrule \hfil \ \ \fbb #\frb\ 
   &\vrule \hfil \ \ \frb #\ \hfil
   &\vrule \hfil \ \ \frb #\ 
   &\vrule \hfil \ \ \frb #\ \vrule \cr%
\noalign{\hrule}
 & &5.7.169.239.617&1189&6& & &3.11.23.577.2003&11051&10982 \cr
7021&872243645&4.3.29.41.617&903&286&7039&877199829&4.289.19.43.257.577&621&10430 \cr
 & &16.9.7.11.13.43&11&3096& & &16.27.5.7.17.19.23.149&11305&10728 \cr
\noalign{\hrule}
 & &25.11.13.17.83.173&107&972& & &3.23.89.241.593&3367&2774 \cr
7022&872668225&8.243.5.11.17.107&10207&9802&7040&877628733&4.7.13.19.37.73.241&33583&1530 \cr
 & &32.3.169.29.59.173&1131&944& & &16.9.5.11.17.43.71&3905&17544 \cr
\noalign{\hrule}
 & &729.125.11.13.67&24559&24284& & &9.5.17.19.23.37.71&301&338 \cr
7023&873068625&8.5.169.41.467.599&98483&2748&7041&878219235&4.5.7.169.17.19.23.43&957&658 \cr
 & &64.3.7.11.229.1279&8953&7328& & &16.3.49.11.13.29.43.47&58609&56056 \cr
\noalign{\hrule}
 & &81.5.7.13.19.29.43&12319&11374& & &729.25.139.347&33229&15004 \cr
7024&873205515&4.3.121.13.47.97.127&1943&3010&7042&879046425&8.7.121.31.47.101&4617&3160 \cr
 & &16.5.7.11.29.43.47.67&517&536& & &128.243.5.11.19.79&869&1216 \cr
\noalign{\hrule}
 & &27.5.17.23.71.233&39589&40796& & &3.7.17.29.31.2741&1311&1430 \cr
7025&873222255&8.9.7.11.31.47.59.61&481&68&7043&879704763&4.9.5.11.13.19.23.29.31&1855&5482 \cr
 & &64.11.13.17.31.37.47&14911&16544& & &16.25.7.19.53.2741&475&424 \cr
\noalign{\hrule}
 & &25.13.43.101.619&549&764& & &13.17.19.53.59.67&183&506 \cr
7026&873703025&8.9.5.61.191.619&451&1406&7044&879728291&4.3.11.23.59.61.67&1275&2678 \cr
 & &32.3.11.19.37.41.61&38247&24272& & &16.9.25.11.13.17.103&927&2200 \cr
\noalign{\hrule}
 & &9.7.19.37.109.181&72215&65888& & &25.11.23.37.53.71&91&162 \cr
7027&873777681&64.5.11.13.29.71.101&109&210&7045&880636075&4.81.25.7.13.37.53&253&1178 \cr
 & &256.3.25.7.13.71.109&1775&1664& & &16.3.7.11.13.19.23.31&651&1976 \cr
\noalign{\hrule}
 & &9.5.49.121.29.113&33269&35096& & &9.17.19.29.31.337&51695&42328 \cr
7028&874319985&16.7.17.19.41.103.107&7759&5940&7046&880713441&16.5.49.11.13.37.211&29&510 \cr
 & &128.27.5.11.41.7759&7759&7872& & &64.3.25.17.29.211&211&800 \cr
\noalign{\hrule}
 & &11.19.1307.3203&2255&948& & &9.11.13.31.71.311&30989&8908 \cr
7029&874941089&8.3.5.121.19.41.79&2769&470&7047&880965657&8.7.17.19.131.233&1915&2046 \cr
 & &32.9.25.13.47.71&30033&5200& & &32.3.5.7.11.19.31.383&1915&2128 \cr
\noalign{\hrule}
 & &5.7.13.17.29.47.83&537&44& & &25.23.877.1747&4379&4356 \cr
7030&875052815&8.3.5.11.13.47.179&5559&6076&7048&880968425&8.9.5.121.29.151.877&221&656 \cr
 & &64.9.49.17.31.109&3379&2016& & &256.3.121.13.17.41.151&253011&251264 \cr
\noalign{\hrule}
 & &9.31.41.59.1297&30989&9218& & &71.32761.379&29835&2926 \cr
7031&875346597&4.7.11.19.233.419&615&848&7049&881565749&4.27.5.7.11.13.17.19&181&142 \cr
 & &128.3.5.41.53.419&2095&3392& & &16.9.5.7.11.71.181&99&280 \cr
\noalign{\hrule}
 & &5.11.19.53.97.163&49&534& & &27.5.49.11.17.23.31&1577&1118 \cr
7032&875692235&4.3.49.19.89.163&3729&632&7050&881984565&4.13.19.23.31.43.83&1781&792 \cr
 & &64.9.11.79.113&8927&288& & &64.9.11.169.19.137&3211&4384 \cr
\noalign{\hrule}
 & &5.11.13.29.157.269&85&72& & &3.5.7.103.127.643&8791&4290 \cr
7033&875701255&16.9.25.11.17.29.269&2567&392&7051&883163715&4.9.25.11.13.59.149&643&2282 \cr
 & &256.3.49.289.151&22197&36992& & &16.7.59.163.643&163&472 \cr
\noalign{\hrule}
 & &9.79.1091.1129&585&506& & &9.11.37.293.823&245&578 \cr
7034&875766429&4.81.5.11.13.23.1129&1091&38&7052&883292157&4.5.49.11.289.293&185&108 \cr
 & &16.5.11.19.23.1091&209&920& & &32.27.25.7.289.37&2023&1200 \cr
\noalign{\hrule}
 & &7.121.13.19.53.79&3483&2930& & &5.121.19.31.37.67&2533&1944 \cr
7035&875958083&4.81.5.13.19.43.293&35&22&7053&883379255&16.243.5.17.67.149&93&242 \cr
 & &16.27.25.7.11.43.293&7911&8600& & &64.729.121.17.31&729&544 \cr
\noalign{\hrule}
 & &3.121.13.19.29.337&22977&22640& & &9.343.11.53.491&703&1046 \cr
7036&876256953&32.81.5.19.23.37.283&3377&19546&7054&883663011&4.3.19.37.491.523&539&1030 \cr
 & &128.5.11.29.307.337&307&320& & &16.5.49.11.19.37.103&1957&1480 \cr
\noalign{\hrule}
 & &5.7.13.83.139.167&341&11196& & &5.11.23.47.89.167&171&4012 \cr
7037&876638945&8.9.7.11.31.311&145&166&7055&883679665&8.9.5.11.17.19.59&611&434 \cr
 & &32.3.5.11.29.31.83&319&1488& & &32.3.7.13.17.31.47&527&4368 \cr
\noalign{\hrule}
 & &343.11.37.61.103&2505&2542& & &7.17.19.23.29.587&51937&42570 \cr
7038&877113083&4.3.5.7.11.31.41.61.167&1075&5772&7056&885247069&4.9.5.11.43.167.311&29&14 \cr
 & &32.9.125.13.31.37.43&17329&18000& & &16.3.7.11.29.167.311&3421&4008 \cr
\noalign{\hrule}
}%
}
$$
\eject
\vglue -23 pt
\noindent\hskip 1 in\hbox to 6.5 in{\ 7057 -- 7092 \hfill\fbd 885581125 -- 897231335\frb}
\vskip -9 pt
$$
\vbox{
\nointerlineskip
\halign{\strut
    \vrule \ \ \hfil \frb #\ 
   &\vrule \hfil \ \ \fbb #\frb\ 
   &\vrule \hfil \ \ \frb #\ \hfil
   &\vrule \hfil \ \ \frb #\ 
   &\vrule \hfil \ \ \frb #\ \ \vrule \hskip 2 pt
   &\vrule \ \ \hfil \frb #\ 
   &\vrule \hfil \ \ \fbb #\frb\ 
   &\vrule \hfil \ \ \frb #\ \hfil
   &\vrule \hfil \ \ \frb #\ 
   &\vrule \hfil \ \ \frb #\ \vrule \cr%
\noalign{\hrule}
 & &125.11.169.37.103&1357&18& & &5.19.107.239.367&781&414 \cr
7057&885581125&4.9.13.23.37.59&275&206&7075&891602645&4.9.11.19.23.71.107&367&260 \cr
 & &16.3.25.11.59.103&177&8& & &32.3.5.13.23.71.367&897&1136 \cr
\noalign{\hrule}
 & &5.121.13.251.449&133&582& & &9.25.7.71.79.101&8107&2498 \cr
7058&886377635&4.3.7.11.19.97.251&10585&8742&7076&892251675&4.3.5.121.67.1249&403&202 \cr
 & &16.9.5.29.31.47.73&12267&18104& & &16.13.31.101.1249&1249&3224 \cr
\noalign{\hrule}
 & &9.17.19.41.43.173&649&130& & &9.5.11.13.17.41.199&7409&1546 \cr
7059&886632093&4.3.5.11.13.17.43.59&133&82&7077&892553805&4.17.31.239.773&14013&9950 \cr
 & &16.7.11.13.19.41.59&413&1144& & &16.81.25.173.199&173&360 \cr
\noalign{\hrule}
 & &3.25.7.17.23.29.149&5409&4366& & &3.5.11.13.29.113.127&23&3706 \cr
7060&886993275&4.27.29.37.59.601&39767&23540&7078&892703955&4.5.13.17.23.109&751&666 \cr
 & &32.5.7.11.13.19.23.107&2033&2288& & &16.9.23.37.751&17273&888 \cr
\noalign{\hrule}
 & &27.5.11.13.19.41.59&83&1204& & &5.121.17.29.41.73&75101&75060 \cr
7061&887277105&8.3.5.7.41.43.83&6071&4636&7079&892707145&8.27.25.13.29.53.109.139&41&3434 \cr
 & &64.13.19.61.467&467&1952& & &32.3.17.41.53.101.109&5353&5232 \cr
\noalign{\hrule}
 & &13.17.19.241.877&1661&2538& & &5.19.47.317.631&22523&37422 \cr
7062&887488043&4.27.11.47.151.241&2005&646&7080&893120555&4.243.7.11.101.223&25&52 \cr
 & &16.3.5.17.19.47.401&2005&1128& & &32.9.25.13.101.223&11817&17840 \cr
\noalign{\hrule}
 & &13.17.19.29.37.197&3223&3186& & &3.7.17.23.127.857&2809&3190 \cr
7063&887588819&4.27.11.19.59.197.293&553&3190&7081&893677029&4.5.11.17.23.29.2809&117&784 \cr
 & &16.3.5.7.121.29.59.79&21417&22120& & &128.9.5.49.11.13.53&8745&5824 \cr
\noalign{\hrule}
 & &27.17.97.127.157&4067&3940& & &9.5.49.11.29.31.41&79&166 \cr
7064&887744097&8.9.5.49.83.97.197&2563&3548&7082&894015045&4.3.11.31.41.79.83&109&232 \cr
 & &64.7.11.83.233.887&206671&204512& & &64.29.79.83.109&6557&3488 \cr
\noalign{\hrule}
 & &3.7.169.17.41.359&169&118& & &5.7.13.19.31.47.71&14283&14330 \cr
7065&888041427&4.28561.59.359&24871&3690&7083&894299315&4.27.25.7.19.529.1433&4147&31678 \cr
 & &16.9.5.7.11.17.19.41&95&264& & &16.3.11.13.23.29.47.337&7751&7656 \cr
\noalign{\hrule}
 & &125.13.289.31.61&307&18& & &9.17.53.211.523&341&182 \cr
7066&888060875&4.9.5.31.61.307&473&442&7084&894852477&4.3.7.11.13.17.31.211&53&580 \cr
 & &16.3.11.13.17.43.307&1419&2456& & &32.5.7.11.13.29.53&4147&560 \cr
\noalign{\hrule}
 & &27.25.67.73.269&127&142& & &243.5.7.127.829&3509&2294 \cr
7067&888083325&4.9.5.67.71.73.127&22327&21692&7085&895431915&4.121.29.31.37.127&249&1148 \cr
 & &32.11.17.29.71.83.269&22649&22576& & &32.3.7.11.37.41.83&3071&7216 \cr
\noalign{\hrule}
 & &27.5.49.11.29.421&793&3838& & &9.25.11.169.2141&10057&9212 \cr
7068&888387885&4.9.7.13.19.61.101&649&548&7086&895526775&8.5.49.11.47.89.113&897&1688 \cr
 & &32.11.13.59.61.137&8357&12272& & &128.3.7.13.23.89.211&14329&13504 \cr
\noalign{\hrule}
 & &71.73.173.991&42313&30030& & &27.5.49.23.71.83&1451&1454 \cr
7069&888589069&4.3.5.7.11.13.17.19.131&519&196&7087&896590485&4.9.7.23.71.727.1451&30481&330838 \cr
 & &32.9.343.131.173&3087&2096& & &16.11.17.83.163.1993&21923&22168 \cr
\noalign{\hrule}
 & &5.11.19.23.71.521&17973&19018& & &9.11.31.37.53.149&3893&4004 \cr
7070&889078685&4.9.23.37.257.1997&37229&8702&7088&896728041&8.3.7.121.13.17.31.229&625&63142 \cr
 & &16.3.19.59.229.631&13511&15144& & &32.625.131.241&31571&10000 \cr
\noalign{\hrule}
 & &3.5.101.307.1913&1417&496& & &9.625.31.37.139&16211&15064 \cr
7071&889745865&32.5.13.31.101.109&307&198&7089&896810625&16.25.7.13.29.43.269&297&28 \cr
 & &128.9.11.13.31.307&1209&704& & &128.27.49.11.29.43&13717&9408 \cr
\noalign{\hrule}
 & &9.49.1331.37.41&4825&4492& & &3.23.29.59.71.107&201&2260 \cr
7072&890435007&8.25.7.41.193.1123&7887&26&7090&896894223&8.9.5.59.67.113&209&322 \cr
 & &32.3.25.11.13.239&239&5200& & &32.5.7.11.19.23.67&6365&1232 \cr
\noalign{\hrule}
 & &5.7.17.23.151.431&549&1606& & &7.11.29.47.83.103&21441&21470 \cr
7073&890633485&4.9.11.17.23.61.73&10721&13130&7091&897226099&4.3.5.49.19.103.113.1021&105183&20 \cr
 & &16.3.5.13.71.101.151&1313&1704& & &32.27.25.13.29.31&775&5616 \cr
\noalign{\hrule}
 & &27.5.7.23.131.313&31673&10582& & &5.7.121.13.43.379&2603&708 \cr
7074&891200205&4.11.13.19.37.1667&1037&630&7092&897231335&8.3.11.13.19.59.137&1137&370 \cr
 & &16.9.5.7.13.17.19.61&1037&1976& & &32.9.5.19.37.379&703&144 \cr
\noalign{\hrule}
}%
}
$$
\eject
\vglue -23 pt
\noindent\hskip 1 in\hbox to 6.5 in{\ 7093 -- 7128 \hfill\fbd 897466325 -- 909845055\frb}
\vskip -9 pt
$$
\vbox{
\nointerlineskip
\halign{\strut
    \vrule \ \ \hfil \frb #\ 
   &\vrule \hfil \ \ \fbb #\frb\ 
   &\vrule \hfil \ \ \frb #\ \hfil
   &\vrule \hfil \ \ \frb #\ 
   &\vrule \hfil \ \ \frb #\ \ \vrule \hskip 2 pt
   &\vrule \ \ \hfil \frb #\ 
   &\vrule \hfil \ \ \fbb #\frb\ 
   &\vrule \hfil \ \ \frb #\ \hfil
   &\vrule \hfil \ \ \frb #\ 
   &\vrule \hfil \ \ \frb #\ \vrule \cr%
\noalign{\hrule}
 & &25.7.23.293.761&33099&14074& & &3.5.7.19.59.79.97&16559&12056 \cr
7093&897466325&4.3.11.17.31.59.227&293&234&7111&901973415&16.7.11.29.137.571&765&194 \cr
 & &16.27.11.13.227.293&2951&2376& & &64.9.5.11.17.29.97&957&544 \cr
\noalign{\hrule}
 & &5.7.13.17.19.41.149&1643&1188& & &3.11.13.17.337.367&3463&3830 \cr
7094&897809185&8.27.11.17.31.41.53&3419&4690&7112&901990947&4.5.337.383.3463&889&2574 \cr
 & &32.3.5.7.11.13.67.263&2893&3216& & &16.9.7.11.13.127.383&2681&3048 \cr
\noalign{\hrule}
 & &9.5.121.23.71.101&3661&8560& & &7.11.113.199.521&729&2918 \cr
7095&898060185&32.3.25.7.107.523&1111&1136&7113&902110979&4.729.113.1459&2255&796 \cr
 & &1024.11.71.101.523&523&512& & &32.27.5.11.41.199&135&656 \cr
\noalign{\hrule}
 & &3.19.47.179.1873&5143&3270& & &9.7.13.19.103.563&253&310 \cr
7096&898180293&4.9.5.19.37.109.139&187&358&7114&902366829&4.3.5.7.11.13.23.31.103&33283&2252 \cr
 & &16.11.17.37.139.179&2363&3256& & &32.83.401.563&401&1328 \cr
\noalign{\hrule}
 & &3.125.121.29.683&1717&1792& & &3.5.37.67.149.163&2915&3116 \cr
7097&898742625&512.5.7.17.101.683&89&594&7115&903112095&8.25.11.19.41.53.149&307&1332 \cr
 & &2048.27.7.11.17.89&10591&9216& & &64.9.19.37.53.307&5833&5088 \cr
\noalign{\hrule}
 & &81.47.139.1699&6479&4780& & &81.5.11.43.53.89&2231&1252 \cr
7098&899064927&8.5.11.19.31.47.239&151&1044&7116&903612105&8.5.23.53.97.313&2717&13872 \cr
 & &64.9.11.29.31.151&4379&10912& & &256.3.11.13.289.19&5491&1664 \cr
\noalign{\hrule}
 & &5.11.23.31.101.227&191&36& & &9.25.11.37.71.139&143&782 \cr
7099&899082305&8.9.11.23.101.191&4313&1990&7117&903753675&4.121.13.17.23.139&2627&570 \cr
 & &32.3.5.19.199.227&199&912& & &16.3.5.13.19.37.71&13&152 \cr
\noalign{\hrule}
 & &9.25.7.11.17.43.71&551&656& & &9.19.2213.2389&1199&1190 \cr
7100&899184825&32.3.5.11.19.29.41.43&71&544&7118&904052547&4.5.7.11.17.19.109.2213&2237&24 \cr
 & &2048.17.19.29.71&551&1024& & &64.3.5.11.109.2237&24607&17440 \cr
\noalign{\hrule}
 & &243.7.121.17.257&38425&24026& & &11.17.19.29.67.131&169&150 \cr
7101&899231949&4.25.29.41.53.293&22627&7098&7119&904355749&4.3.25.169.17.67.131&11571&2794 \cr
 & &16.3.7.1331.169.17&169&88& & &16.9.5.7.11.19.29.127&635&504 \cr
\noalign{\hrule}
 & &3.5.11.169.59.547&1541&994& & &11.13.19.43.61.127&119&1278 \cr
7102&899932605&4.7.11.23.59.67.71&1501&3042&7120&905089757&4.9.7.13.17.43.71&635&1558 \cr
 & &16.9.169.19.71.79&1349&1896& & &16.3.5.7.19.41.127&615&56 \cr
\noalign{\hrule}
 & &11.19.47.277.331&585&308& & &9.7.11.31.113.373&10209&14312 \cr
7103&900641401&8.9.5.7.121.13.331&105&226&7121&905486967&16.27.41.83.1789&2015&226 \cr
 & &32.27.25.49.13.113&36725&21168& & &64.5.13.31.41.113&533&160 \cr
\noalign{\hrule}
 & &27.5.13.47.61.179&1177&1150& & &9.25.17.41.53.109&1577&596 \cr
7104&900653715&4.125.11.23.47.61.107&6201&326&7122&905978025&8.25.17.19.83.149&1199&876 \cr
 & &16.9.11.13.23.53.163&3749&4664& & &64.3.11.73.109.149&1639&2336 \cr
\noalign{\hrule}
 & &27.25.7.11.13.31.43&1349&16& & &27.13.41.79.797&43&754 \cr
7105&900674775&32.9.5.11.19.71&527&518&7123&906100533&4.3.169.29.41.43&5335&1594 \cr
 & &128.7.17.31.37.71&1207&2368& & &16.5.11.97.797&1067&40 \cr
\noalign{\hrule}
 & &27.25.7.11.13.31.43&493&592& & &3.7.11.19.23.47.191&1465&2546 \cr
7106&900674775&32.3.5.13.17.29.37.43&497&62&7124&906201219&4.5.11.361.67.293&549&188 \cr
 & &128.7.17.31.37.71&1207&2368& & &32.9.5.47.61.293&1465&2928 \cr
\noalign{\hrule}
 & &27.125.13.19.23.47&803&2072& & &5.11.19.47.53.349&4437&598 \cr
7107&901148625&16.7.11.13.19.37.73&2209&7500&7125&908480155&4.9.13.17.23.29.47&165&212 \cr
 & &128.3.625.2209&235&64& & &32.27.5.11.17.23.53&621&272 \cr
\noalign{\hrule}
 & &7.13.23.29.31.479&297&80& & &9.121.23.131.277&5&28 \cr
7108&901289753&32.27.5.11.23.479&899&1378&7126&908880489&8.3.5.7.11.131.277&305&1136 \cr
 & &128.3.5.13.29.31.53&159&320& & &256.25.7.61.71&30317&3200 \cr
\noalign{\hrule}
 & &5.11.19.29.151.197&1507&1362& & &9.121.23.131.277&305&1136 \cr
7109&901482835&4.3.121.137.197.227&159&38&7127&908880489&32.3.5.11.23.61.71&5&28 \cr
 & &16.9.19.53.137.227&12031&9864& & &256.25.7.61.71&30317&3200 \cr
\noalign{\hrule}
 & &5.11.13.23.29.31.61&361&42& & &81.5.7.487.659&3839&544 \cr
7110&901827355&4.3.5.7.361.23.61&371&66&7128&909845055&64.9.7.11.17.349&659&650 \cr
 & &16.9.49.11.19.53&9063&392& & &256.25.13.349.659&1745&1664 \cr
\noalign{\hrule}
}%
}
$$
\eject
\vglue -23 pt
\noindent\hskip 1 in\hbox to 6.5 in{\ 7129 -- 7164 \hfill\fbd 909871677 -- 921829581\frb}
\vskip -9 pt
$$
\vbox{
\nointerlineskip
\halign{\strut
    \vrule \ \ \hfil \frb #\ 
   &\vrule \hfil \ \ \fbb #\frb\ 
   &\vrule \hfil \ \ \frb #\ \hfil
   &\vrule \hfil \ \ \frb #\ 
   &\vrule \hfil \ \ \frb #\ \ \vrule \hskip 2 pt
   &\vrule \ \ \hfil \frb #\ 
   &\vrule \hfil \ \ \fbb #\frb\ 
   &\vrule \hfil \ \ \frb #\ \hfil
   &\vrule \hfil \ \ \frb #\ 
   &\vrule \hfil \ \ \frb #\ \vrule \cr%
\noalign{\hrule}
 & &27.11.13.19.79.157&11265&11186& & &5.7.121.17.29.439&7419&6980 \cr
7129&909871677&4.81.5.7.17.19.47.751&61&26224&7147&916568345&8.3.25.29.349.2473&1599&874 \cr
 & &128.11.47.61.149&2867&9536& & &32.9.13.19.23.41.349&110331&106096 \cr
\noalign{\hrule}
 & &7.11.79.151.991&963&94& & &9.5.11.169.19.577&5729&5234 \cr
7130&910266203&4.9.47.107.991&1463&1510&7148&917109765&4.169.17.337.2617&50721&6232 \cr
 & &16.3.5.7.11.19.107.151&535&456& & &64.3.11.19.29.41.53&1189&1696 \cr
\noalign{\hrule}
 & &3.5.11.13.19.89.251&4049&5740& & &3.343.11.13.23.271&3425&5518 \cr
7131&910425945&8.25.7.11.41.4049&19809&8534&7149&917167251&4.25.49.31.89.137&247&198 \cr
 & &32.9.17.31.71.251&1581&1136& & &16.9.5.11.13.19.31.137&2945&3288 \cr
\noalign{\hrule}
 & &3.19.29.31.109.163&52745&43262& & &3.7.11.19.23.61.149&301&370 \cr
7132&910434381&4.5.7.11.97.137.223&36283&49572&7150&917507283&4.5.49.19.37.43.149&509&2322 \cr
 & &32.729.13.17.2791&47447&50544& & &16.27.5.1849.509&22905&14792 \cr
\noalign{\hrule}
 & &3.25.11.13.73.1163&42807&42092& & &3.5.11.13.31.37.373&3059&1044 \cr
7133&910541775&8.9.5.17.19.619.751&803&52&7151&917697495&8.27.7.19.23.29.37&169&682 \cr
 & &64.11.13.17.73.619&619&544& & &32.7.11.169.29.31&91&464 \cr
\noalign{\hrule}
 & &3.7.11.19.31.37.181&181&218& & &9.13.19.43.9601&3575&6026 \cr
7134&911187123&4.11.31.109.32761&34965&2204&7152&917749989&4.3.25.11.169.23.131&22909&10234 \cr
 & &32.27.5.7.19.29.37&261&80& & &16.7.17.31.43.739&3689&5912 \cr
\noalign{\hrule}
 & &3.125.7.11.37.853&519&334& & &7.11.23.587.883&213&374 \cr
7135&911323875&4.9.25.7.11.167.173&131&1706&7153&917946491&4.3.121.17.71.883&587&1470 \cr
 & &16.131.173.853&173&1048& & &16.9.5.49.71.587&315&568 \cr
\noalign{\hrule}
 & &11.23.1559.2311&198369&196058& & &27.5.11.17.41.887&161&26 \cr
7136&911520797&4.81.31.79.167.587&46483&110&7154&918084915&4.7.13.23.41.887&3025&9234 \cr
 & &16.9.5.11.23.43.47&1935&376& & &16.243.25.121.19&95&792 \cr
\noalign{\hrule}
 & &81.343.11.19.157&265&1462& & &3.29.53.139.1433&2375&2236 \cr
7137&911643579&4.9.5.49.17.43.53&247&194&7155&918451257&8.125.13.19.43.1433&99&1334 \cr
 & &16.5.13.17.19.43.97&8245&4472& & &32.9.25.11.23.29.43&2967&4400 \cr
\noalign{\hrule}
 & &5.11.17.361.37.73&1143&98& & &9.5.121.13.19.683&5899&5216 \cr
7138&911682035&4.9.49.19.37.127&425&278&7156&918576945&64.121.17.163.347&22379&19608 \cr
 & &16.3.25.17.127.139&635&3336& & &1024.3.7.19.23.43.139&22379&22016 \cr
\noalign{\hrule}
 & &25.11.13.17.43.349&363&712& & &11.13.19.523.647&170549&167832 \cr
7139&912050425&16.3.1331.13.17.89&8643&8660&7157&919381177&16.81.7.29.37.5881&5681&200 \cr
 & &128.9.5.43.67.89.433&38537&38592& & &256.3.25.13.19.23.37&2553&3200 \cr
\noalign{\hrule}
 & &5.13.19.37.41.487&1707&2188& & &7.11.17.19.103.359&2813&30456 \cr
7140&912392065&8.3.487.547.569&41&528&7158&919654967&16.81.29.47.97&2149&2410 \cr
 & &256.9.11.41.547&4923&1408& & &64.9.5.7.241.307&13815&7712 \cr
\noalign{\hrule}
 & &7.43.53.173.331&24975&39208& & &3.5.7.121.19.37.103&137&248 \cr
7141&913516639&16.27.25.169.29.37&17&22&7159&919956345&16.11.19.31.103.137&10295&3816 \cr
 & &64.9.5.11.13.17.29.37&57681&65120& & &256.9.5.29.53.71&4611&9088 \cr
\noalign{\hrule}
 & &5.7.11.19.29.59.73&333&1378& & &81.11.13.19.37.113&9947&10060 \cr
7142&913665445&4.9.7.13.37.53.73&275&236&7160&920141937&8.5.343.11.29.37.503&167&2682 \cr
 & &32.3.25.11.37.53.59&795&592& & &32.9.49.29.149.167&24883&22736 \cr
\noalign{\hrule}
 & &9.11.19.23.37.571&455&248& & &3.5.7.17.841.613&36053&36894 \cr
7143&914017401&16.5.7.11.13.31.571&115&456&7161&920226405&4.9.5.11.13.31.43.1163&5423&392 \cr
 & &256.3.25.7.13.19.23&325&896& & &64.49.121.17.29.31&847&992 \cr
\noalign{\hrule}
 & &3.19.61.241.1091&33155&33396& & &13.41.59.83.353&4743&154 \cr
7144&914211087&8.9.5.121.361.23.349&377&125612&7162&921365653&4.9.7.11.17.31.41&365&332 \cr
 & &64.13.29.31.1013&11687&32416& & &32.3.5.7.31.73.83&1085&3504 \cr
\noalign{\hrule}
 & &5.11.13.31.157.263&1369&1524& & &7.13.961.83.127&45&172 \cr
7145&915215015&8.3.13.1369.127.157&1705&336&7163&921820991&8.9.5.13.31.43.83&4697&698 \cr
 & &256.9.5.7.11.31.127&889&1152& & &32.3.7.11.61.349&11517&976 \cr
\noalign{\hrule}
 & &3.7.11.13.23.89.149&535&444& & &9.529.37.5233&90695&85462 \cr
7146&915924009&8.9.5.23.37.107.149&29&178&7164&921829581&4.5.11.13.17.19.97.173&28797&7360 \cr
 & &32.5.29.37.89.107&3959&2320& & &512.3.25.23.29.331&8275&7424 \cr
\noalign{\hrule}
}%
}
$$
\eject
\vglue -23 pt
\noindent\hskip 1 in\hbox to 6.5 in{\ 7165 -- 7200 \hfill\fbd 922489351 -- 933738549\frb}
\vskip -9 pt
$$
\vbox{
\nointerlineskip
\halign{\strut
    \vrule \ \ \hfil \frb #\ 
   &\vrule \hfil \ \ \fbb #\frb\ 
   &\vrule \hfil \ \ \frb #\ \hfil
   &\vrule \hfil \ \ \frb #\ 
   &\vrule \hfil \ \ \frb #\ \ \vrule \hskip 2 pt
   &\vrule \ \ \hfil \frb #\ 
   &\vrule \hfil \ \ \fbb #\frb\ 
   &\vrule \hfil \ \ \frb #\ \hfil
   &\vrule \hfil \ \ \frb #\ 
   &\vrule \hfil \ \ \frb #\ \vrule \cr%
\noalign{\hrule}
 & &7.31.47.151.599&61443&68540& & &3.125.7.29.89.137&1679&1946 \cr
7165&922489351&8.9.5.23.149.6827&1727&8554&7183&928192125&4.49.23.73.137.139&1595&8406 \cr
 & &32.3.5.7.11.13.47.157&2041&2640& & &16.9.5.11.23.29.467&1401&2024 \cr
\noalign{\hrule}
 & &3.5.229.379.709&341&346& & &9.7.121.13.17.19.29&229&1318 \cr
7166&923022285&4.11.31.173.379.709&88101&153668&7184&928260333&4.19.29.229.659&9581&2940 \cr
 & &32.27.13.41.251.937&133783&134928& & &32.3.5.49.11.13.67&67&560 \cr
\noalign{\hrule}
 & &5.7.11.59.97.419&83519&77796& & &49.11.19.29.53.59&291&830 \cr
7167&923205745&8.9.47.1777.2161&2353&4130&7185&928684603&4.3.5.29.53.83.97&147&118 \cr
 & &32.3.5.7.13.47.59.181&2353&2256& & &16.9.49.59.83.97&873&664 \cr
\noalign{\hrule}
 & &3.49.121.23.37.61&2319&2158& & &9.25.11.19.23.859&533&326 \cr
7168&923340957&4.9.7.13.61.83.773&298891&293480&7186&929072925&4.25.11.13.19.41.163&2247&6322 \cr
 & &64.5.11.23.29.257.1163&37265&37216& & &16.3.7.13.29.107.109&11663&21112 \cr
\noalign{\hrule}
 & &81.5.11.13.41.389&4189&868& & &3.5.17.31.41.47.61&925&966 \cr
7169&923686335&8.5.7.11.31.59.71&243&538&7187&929209035&4.9.125.7.17.23.37.47&4697&1178 \cr
 & &32.243.7.31.269&807&3472& & &16.49.11.19.31.37.61&1813&1672 \cr
\noalign{\hrule}
 & &3.5.11.23.29.37.227&1001&4220& & &49.11.19.29.31.101&2025&2924 \cr
7170&924351945&8.25.7.121.13.211&261&586&7188&929872559&8.81.25.11.17.19.43&1519&1046 \cr
 & &32.9.29.211.293&879&3376& & &32.3.5.49.17.31.523&1569&1360 \cr
\noalign{\hrule}
 & &9.25.7.17.19.23.79&341&134& & &9.5.121.19.89.101&391&410 \cr
7171&924353325&4.7.11.17.31.67.79&1879&570&7189&929956995&4.25.121.17.23.41.101&129&2654 \cr
 & &16.3.5.19.67.1879&1879&536& & &16.3.17.41.43.1327&22559&14104 \cr
\noalign{\hrule}
 & &81.5.17.19.37.191&3397&8554& & &3.5.49.43.59.499&481&422 \cr
7172&924469605&4.3.5.7.13.43.47.79&191&836&7190&930482805&4.5.7.13.37.211.499&477&22 \cr
 & &32.7.11.19.47.191&517&112& & &16.9.11.37.53.211&5883&18568 \cr
\noalign{\hrule}
 & &3.25.7.19.59.1571&2483&30508& & &3.5.19.37.103.857&377&480 \cr
7173&924572775&8.13.29.191.263&36135&35872&7191&930817695&64.9.25.13.19.29.37&3533&8492 \cr
 & &512.9.5.11.19.59.73&803&768& & &512.11.193.3533&38863&49408 \cr
\noalign{\hrule}
 & &9.7.11.19.23.43.71&65&562& & &9.13.23.53.61.107&43175&30548 \cr
7174&924573573&4.3.5.13.23.43.281&781&2186&7192&930900321&8.25.7.11.157.1091&1409&318 \cr
 & &16.11.13.71.1093&1093&104& & &32.3.25.7.53.1409&1409&2800 \cr
\noalign{\hrule}
 & &5.11.13.37.73.479&34959&17236& & &3.5.7.11.23.101.347&13&2 \cr
7175&925051985&8.3.31.43.139.271&205&66&7193&931023555&4.7.13.23.101.347&605&8586 \cr
 & &32.9.5.11.31.41.43&1763&4464& & &16.81.5.121.53&53&2376 \cr
\noalign{\hrule}
 & &9.5.7.11.13.19.23.47&3811&3944& & &1331.19.29.31.41&549&230 \cr
7176&925179255&16.3.13.17.23.29.37.103&665&2&7194&932127251&4.9.5.121.23.31.61&2337&446 \cr
 & &64.5.7.19.37.103&103&1184& & &16.27.5.19.41.223&1115&216 \cr
\noalign{\hrule}
 & &125.19.29.89.151&7887&5018& & &9.23.149.167.181&31175&6292 \cr
7177&925611125&4.3.25.11.13.193.239&1691&1416&7195&932291361&8.25.121.13.29.43&181&138 \cr
 & &64.9.19.59.89.193&1737&1888& & &32.3.25.11.13.23.181&325&176 \cr
\noalign{\hrule}
 & &3.7.13.19.137.1303&787&2090& & &27.11.19.29.41.139&851&338 \cr
7178&925936557&4.5.11.13.361.787&1953&2740&7196&932624253&4.11.169.23.37.139&615&914 \cr
 & &32.9.25.7.11.31.137&775&528& & &16.3.5.13.37.41.457&2405&3656 \cr
\noalign{\hrule}
 & &9.13.17.23.31.653&373&280& & &5.7.11.13.29.59.109&9189&31904 \cr
7179&926056521&16.3.5.7.13.17.23.373&7469&7078&7197&933427495&64.9.997.1021&1033&2030 \cr
 & &64.5.49.11.97.3539&173411&170720& & &256.3.5.7.29.1033&1033&384 \cr
\noalign{\hrule}
 & &81.5.43.127.419&781&362& & &125.289.43.601&513&88 \cr
7180&926704395&4.9.5.11.43.71.181&2581&1676&7198&933578375&16.27.5.11.17.19.43&667&752 \cr
 & &32.29.71.89.419&2059&1424& & &512.9.19.23.29.47&31349&43776 \cr
\noalign{\hrule}
 & &9.7.29.313.1621&24475&39064& & &9.5.29.359.1993&817&1176 \cr
7181&926970471&16.25.11.19.89.257&203&1488&7199&933710535&16.27.5.49.19.29.43&1993&572 \cr
 & &512.3.5.7.11.29.31&341&1280& & &128.11.13.43.1993&473&832 \cr
\noalign{\hrule}
 & &7.13.43.101.2347&201&100& & &3.17.29.37.113.151&1265&1302 \cr
7182&927564911&8.3.25.13.67.2347&16093&14418&7200&933738549&4.9.5.7.11.23.29.31.113&187&74 \cr
 & &32.243.7.121.19.89&29403&27056& & &16.5.7.121.17.23.31.37&4991&4840 \cr
\noalign{\hrule}
}%
}
$$
\eject
\vglue -23 pt
\noindent\hskip 1 in\hbox to 6.5 in{\ 7201 -- 7236 \hfill\fbd 933779691 -- 949486545\frb}
\vskip -9 pt
$$
\vbox{
\nointerlineskip
\halign{\strut
    \vrule \ \ \hfil \frb #\ 
   &\vrule \hfil \ \ \fbb #\frb\ 
   &\vrule \hfil \ \ \frb #\ \hfil
   &\vrule \hfil \ \ \frb #\ 
   &\vrule \hfil \ \ \frb #\ \ \vrule \hskip 2 pt
   &\vrule \ \ \hfil \frb #\ 
   &\vrule \hfil \ \ \fbb #\frb\ 
   &\vrule \hfil \ \ \frb #\ \hfil
   &\vrule \hfil \ \ \frb #\ 
   &\vrule \hfil \ \ \frb #\ \vrule \cr%
\noalign{\hrule}
 & &27.13.529.47.107&4895&134& & &3.7.17.23.239.479&1313&2750 \cr
7201&933779691&4.3.5.11.13.67.89&107&94&7219&940003491&4.125.7.11.13.23.101&153&958 \cr
 & &16.5.11.47.89.107&89&440& & &16.9.25.13.17.479&39&200 \cr
\noalign{\hrule}
 & &5.7.11.13.29.41.157&2599&558& & &5.11.19.43.67.313&2223&658 \cr
7202&934298365&4.9.5.23.29.31.113&1631&1066&7220&942331885&4.9.7.11.13.361.47&17201&15650 \cr
 & &16.3.7.13.23.41.233&233&552& & &16.3.25.103.167.313&1545&1336 \cr
\noalign{\hrule}
 & &9.25.11.13.71.409&16307&12732& & &9.19.37.181.823&245&578 \cr
7203&934329825&8.27.23.709.1061&665&44&7221&942488901&4.5.49.289.19.181&1441&174 \cr
 & &64.5.7.11.19.1061&1061&4256& & &16.3.7.11.17.29.131&2227&17864 \cr
\noalign{\hrule}
 & &7.11.19.37.61.283&1139&11610& & &9.11.31.41.59.127&325&2744 \cr
7204&934463453&4.27.5.7.17.43.67&2257&2860&7222&942836697&16.25.343.13.127&59&186 \cr
 & &32.3.25.11.13.37.61&75&208& & &64.3.5.7.13.31.59&35&416 \cr
\noalign{\hrule}
 & &9.7.13.37.67.461&1199&1132& & &81.25.17.61.449&2783&1258 \cr
7205&935968761&8.11.13.109.283.461&45359&4890&7223&942866325&4.9.121.289.23.37&95&194 \cr
 & &32.3.5.67.163.677&3385&2608& & &16.5.11.19.23.37.97&16169&8536 \cr
\noalign{\hrule}
 & &3.5.13.17.19.89.167&9417&10252& & &3.83.137.139.199&8327&8190 \cr
7206&936146055&8.9.11.19.43.73.233&85&4342&7224&943599693&4.27.5.7.11.13.139.757&11371&15124 \cr
 & &32.5.13.17.73.167&73&16& & &32.11.13.19.83.137.199&209&208 \cr
\noalign{\hrule}
 & &5.7.11.31.131.599&377&222& & &11.17.43.239.491&32903&11790 \cr
7207&936527515&4.3.7.11.13.29.37.131&599&2040&7225&943603309&4.9.5.13.131.2531&1069&1462 \cr
 & &64.9.5.17.37.599&629&288& & &16.3.5.13.17.43.1069&1069&1560 \cr
\noalign{\hrule}
 & &3.5.7.19.43.61.179&8319&2600& & &7.17.37.269.797&1261&30750 \cr
7208&936686415&16.9.125.13.47.59&179&946&7226&943972379&4.3.125.13.41.97&269&264 \cr
 & &64.11.43.47.179&47&352& & &64.9.25.11.97.269&2475&3104 \cr
\noalign{\hrule}
 & &9.31.883.3803&2343&1460& & &9.25.11.13.43.683&160067&162992 \cr
7209&936895671&8.27.5.11.31.71.73&115&2086&7227&944947575&32.59.61.167.2713&3737&6450 \cr
 & &32.25.7.11.23.149&3427&30800& & &128.3.25.37.43.59.101&3737&3776 \cr
\noalign{\hrule}
 & &27.7.13.23.59.281&1907&550& & &3.25.11.61.89.211&7&68 \cr
7210&936895869&4.25.11.281.1907&6313&3222&7228&945053175&8.7.11.17.89.211&1205&1116 \cr
 & &16.9.5.59.107.179&535&1432& & &64.9.5.7.17.31.241&11067&7712 \cr
\noalign{\hrule}
 & &9.5.13.19.37.43.53&1441&964& & &5.11.17.19.83.641&3807&8372 \cr
7211&937250145&8.11.19.43.131.241&527&5106&7229&945151295&8.81.7.13.17.23.47&641&158 \cr
 & &32.3.11.17.23.31.37&713&2992& & &32.27.13.79.641&351&1264 \cr
\noalign{\hrule}
 & &5.11.19.31.103.281&117&398& & &9.5.7.13.19.29.419&127&146 \cr
7212&937608485&4.9.11.13.19.31.199&89&120&7230&945408555&4.3.5.29.73.127.419&36673&55088 \cr
 & &64.27.5.13.89.199&17711&11232& & &128.7.11.169.31.313&9703&9152 \cr
\noalign{\hrule}
 & &9.11.13.31.71.331&35&178& & &9.5.11.13.47.53.59&2501&554 \cr
7213&937619397&4.3.5.7.31.89.331&781&874&7231&945745515&4.3.41.53.61.277&3589&11092 \cr
 & &16.7.11.19.23.71.89&1691&1288& & &32.37.47.59.97&97&592 \cr
\noalign{\hrule}
 & &3.13.289.19.29.151&4015&3148& & &3.11.23.41.83.367&245&122 \cr
7214&937758471&8.5.11.73.151.787&437&1224&7232&947915859&4.5.49.11.23.61.83&3303&1394 \cr
 & &128.9.5.17.19.23.73&1095&1472& & &16.9.5.7.17.41.367&85&168 \cr
\noalign{\hrule}
 & &3.17.151.193.631&605&26& & &9.25.29.31.43.109&77&3302 \cr
7215&937850883&4.5.121.13.17.151&369&386&7233&948062925&4.3.7.11.13.29.127&5735&5822 \cr
 & &16.9.121.13.41.193&1599&968& & &16.5.11.31.37.41.71&2627&3608 \cr
\noalign{\hrule}
 & &9.5.7.11.19.53.269&65&334& & &3.25.7.13.43.53.61&737&178 \cr
7216&938609595&4.3.25.11.13.53.167&779&196&7234&948804675&4.5.7.11.53.67.89&9&44 \cr
 & &32.49.19.41.167&287&2672& & &32.9.121.67.89&17889&1936 \cr
\noalign{\hrule}
 & &9.25.343.43.283&253&596& & &9.5.49.17.19.31.43&451&366 \cr
7217&939142575&8.3.25.11.23.43.149&917&158&7235&949382595&4.27.49.11.31.41.61&149&688 \cr
 & &32.7.79.131.149&11771&2096& & &128.41.43.61.149&6109&3904 \cr
\noalign{\hrule}
 & &9.17.1051.5843&8773&8756& & &9.5.7.73.157.263&957&884 \cr
7218&939571929&8.3.11.31.199.283.1051&263671&94720&7236&949486545&8.27.5.11.13.17.29.157&1&784 \cr
 & &8192.5.37.41.59.109&447515&446464& & &256.49.11.13.17&187&11648 \cr
\noalign{\hrule}
}%
}
$$
\eject
\vglue -23 pt
\noindent\hskip 1 in\hbox to 6.5 in{\ 7237 -- 7272 \hfill\fbd 949668045 -- 960637587\frb}
\vskip -9 pt
$$
\vbox{
\nointerlineskip
\halign{\strut
    \vrule \ \ \hfil \frb #\ 
   &\vrule \hfil \ \ \fbb #\frb\ 
   &\vrule \hfil \ \ \frb #\ \hfil
   &\vrule \hfil \ \ \frb #\ 
   &\vrule \hfil \ \ \frb #\ \ \vrule \hskip 2 pt
   &\vrule \ \ \hfil \frb #\ 
   &\vrule \hfil \ \ \fbb #\frb\ 
   &\vrule \hfil \ \ \frb #\ \hfil
   &\vrule \hfil \ \ \frb #\ 
   &\vrule \hfil \ \ \frb #\ \vrule \cr%
\noalign{\hrule}
 & &3.5.23.53.167.311&2387&1454& & &9.11.169.43.1327&2125&3452 \cr
7237&949668045&4.5.7.11.31.53.727&3933&7568&7255&954687591&8.3.125.17.43.863&665&1528 \cr
 & &128.9.121.19.23.43&5203&3648& & &128.625.7.19.191&25403&40000 \cr
\noalign{\hrule}
 & &3.5.13.23.37.59.97&1199&1102& & &9.7.121.17.53.139&1363&1000 \cr
7238&949703235&4.5.11.19.23.29.37.109&16461&3926&7256&954696897&16.3.125.7.29.47.53&629&484 \cr
 & &16.9.11.13.31.59.151&1023&1208& & &128.25.121.17.37.47&1739&1600 \cr
\noalign{\hrule}
 & &3.11.13.31.37.1931&13125&11978& & &25.7.13.23.71.257&3683&342 \cr
7239&950173653&4.9.625.7.11.53.113&37&62&7257&954774275&4.9.19.29.71.127&1001&1058 \cr
 & &16.25.7.31.37.53.113&2825&2968& & &16.3.7.11.13.529.127&759&1016 \cr
\noalign{\hrule}
 & &3.5.11.23.37.67.101&11537&14016& & &3.49.13.31.71.227&139&542 \cr
7240&950188305&128.9.5.73.83.139&121&536&7258&954787197&4.49.71.139.271&1705&11574 \cr
 & &2048.121.67.139&1529&1024& & &16.9.5.11.31.643&3215&264 \cr
\noalign{\hrule}
 & &27.25.13.131.827&10123&32452& & &9.5.289.23.31.103&619&2574 \cr
7241&950657175&8.7.19.53.61.191&2301&9350&7259&955074195&4.81.11.13.17.619&3335&4712 \cr
 & &32.3.25.11.13.17.59&1003&176& & &64.5.11.19.23.29.31&319&608 \cr
\noalign{\hrule}
 & &5.343.23.89.271&247&198& & &7.11.23.67.83.97&691&222 \cr
7242&951373955&4.9.7.11.13.19.23.271&3425&5518&7260&955307507&4.3.23.37.97.691&1461&770 \cr
 & &16.3.25.19.31.89.137&2945&3288& & &16.9.5.7.11.37.487&1665&3896 \cr
\noalign{\hrule}
 & &3.5.11.37.101.1543&161&1382& & &25.11.13.17.79.199&3161&8136 \cr
7243&951421515&4.5.7.23.101.691&2889&566&7261&955443775&16.9.17.29.109.113&3781&5702 \cr
 & &16.27.7.107.283&17829&856& & &64.3.19.199.2851&2851&1824 \cr
\noalign{\hrule}
 & &9.5.7.11.53.71.73&571&494& & &71.647.20807&10727&10080 \cr
7244&951832035&4.3.13.19.53.73.571&1927&920&7262&955811159&64.9.5.7.17.71.631&10439&11646 \cr
 & &64.5.23.41.47.571&26837&30176& & &256.81.11.13.73.647&10439&10368 \cr
\noalign{\hrule}
 & &27.5.13.41.101.131&11339&4774& & &3.11.17.23.53.1399&1903&2294 \cr
7245&952036605&4.9.7.11.17.23.29.31&37&26&7263&956718741&4.121.31.37.53.173&1417&3060 \cr
 & &16.13.17.23.29.31.37&12121&8584& & &32.9.5.13.17.109.173&7085&8304 \cr
\noalign{\hrule}
 & &9.5.121.179.977&41&936& & &9.37.83.89.389&3503&10890 \cr
7246&952237935&16.81.121.13.41&3007&1954&7264&956889819&4.81.5.121.31.113&559&13114 \cr
 & &64.31.97.977&97&992& & &16.13.43.79.83&3397&104 \cr
\noalign{\hrule}
 & &9.5.7.11.17.19.23.37&1457&1688& & &9.25.13.19.67.257&3949&934 \cr
7247&952434945&16.3.19.23.31.47.211&935&3074&7265&956945925&4.5.11.13.359.467&201&266 \cr
 & &64.5.11.17.29.47.53&1363&1696& & &16.3.7.11.19.67.359&359&616 \cr
\noalign{\hrule}
 & &7.19.71.79.1277&705&572& & &9.11.71.173.787&845&1058 \cr
7248&952638169&8.3.5.11.13.47.71.79&1277&4332&7266&957005379&4.3.5.169.529.787&55&842 \cr
 & &64.9.11.361.1277&209&288& & &16.25.11.13.23.421&9683&2600 \cr
\noalign{\hrule}
 & &3.5.11.19.23.73.181&4273&2282& & &9.13.289.23.1231&11809&9118 \cr
7249&952723365&4.7.73.163.4273&1881&2392&7267&957347469&4.49.17.47.97.241&989&660 \cr
 & &64.9.11.13.19.23.163&489&416& & &32.3.5.7.11.23.43.241&8435&7568 \cr
\noalign{\hrule}
 & &3.5.7.11.73.89.127&43&846& & &3.7.11.13.19.107.157&4849&4100 \cr
7250&953012445&4.27.5.43.47.89&2159&2024&7268&958500543&8.25.11.169.41.373&1413&5516 \cr
 & &64.11.17.23.43.127&989&544& & &64.9.25.7.157.197&591&800 \cr
\noalign{\hrule}
 & &3.5.67.443.2141&37&2178& & &49.529.71.521&13221&12700 \cr
7251&953205315&4.27.121.37.67&923&886&7269&958843711&8.9.25.13.71.113.127&253&14098 \cr
 & &16.121.13.71.443&1573&568& & &32.3.5.7.11.19.23.53&1007&2640 \cr
\noalign{\hrule}
 & &11.43.743.2713&1053&30896& & &121.37.167.1283&277&1560 \cr
7252&953454007&32.81.13.1931&985&946&7270&959246497&16.3.5.11.13.37.277&1597&2004 \cr
 & &128.27.5.11.43.197&985&1728& & &128.9.5.167.1597&1597&2880 \cr
\noalign{\hrule}
 & &3.25.7.31.41.1429&17537&18188& & &27.5.7.11.13.31.229&7939&9694 \cr
7253&953535975&8.13.19.41.71.4547&7337&2790&7271&959323365&4.17.31.37.131.467&20839&40338 \cr
 & &32.9.5.11.23.29.31.71&4899&5104& & &16.243.7.13.83.229&83&72 \cr
\noalign{\hrule}
 & &25.7.41.61.2179&7657&7596& & &3.7.13.19.43.59.73&5055&5566 \cr
7254&953693825&8.9.25.13.19.31.41.211&77&1348&7272&960637587&4.9.5.121.23.59.337&3139&106 \cr
 & &64.3.7.11.13.211.337&30173&32352& & &16.11.23.43.53.73&583&184 \cr
\noalign{\hrule}
}%
}
$$
\eject
\vglue -23 pt
\noindent\hskip 1 in\hbox to 6.5 in{\ 7273 -- 7308 \hfill\fbd 960706635 -- 968080113\frb}
\vskip -9 pt
$$
\vbox{
\nointerlineskip
\halign{\strut
    \vrule \ \ \hfil \frb #\ 
   &\vrule \hfil \ \ \fbb #\frb\ 
   &\vrule \hfil \ \ \frb #\ \hfil
   &\vrule \hfil \ \ \frb #\ 
   &\vrule \hfil \ \ \frb #\ \ \vrule \hskip 2 pt
   &\vrule \ \ \hfil \frb #\ 
   &\vrule \hfil \ \ \fbb #\frb\ 
   &\vrule \hfil \ \ \frb #\ \hfil
   &\vrule \hfil \ \ \frb #\ 
   &\vrule \hfil \ \ \frb #\ \vrule \cr%
\noalign{\hrule}
 & &3.5.7.17.29.67.277&1077&62& & &3.25.7.11.47.53.67&943&608 \cr
7273&960706635&4.9.31.277.359&143&134&7291&963830175&64.5.7.19.23.41.53&9&44 \cr
 & &16.11.13.31.67.359&4667&2728& & &512.9.11.19.23.41&2337&5888 \cr
\noalign{\hrule}
 & &27.25.11.13.37.269&323&158& & &3.5.13.19.433.601&517&84 \cr
7274&960713325&4.9.5.17.19.79.269&18821&6734&7292&964163265&8.9.5.7.11.13.19.47&181&2404 \cr
 & &16.7.11.13.29.37.59&203&472& & &64.7.181.601&1267&32 \cr
\noalign{\hrule}
 & &25.11.13.17.97.163&12867&10442& & &9.625.11.31.503&2379&3154 \cr
7275&960913525&4.3.17.23.227.4289&215&4074&7293&964816875&4.27.25.13.19.61.83&6919&1856 \cr
 & &16.9.5.7.23.43.97&1449&344& & &512.11.17.19.29.37&11951&7424 \cr
\noalign{\hrule}
 & &3.5.121.43.109.113&1757&1972& & &9.5.11.41.199.239&53&152 \cr
7276&961280265&8.7.11.17.29.109.251&2961&200&7294&965250495&16.19.53.199.239&3003&7544 \cr
 & &128.9.25.49.17.47&4165&9024& & &256.3.7.11.13.23.41&299&896 \cr
\noalign{\hrule}
 & &9.25.11.29.59.227&5299&3256& & &25.7.11.43.107.109&31117&1692 \cr
7277&961282575&16.5.7.121.37.757&681&76&7295&965404825&8.9.841.37.47&535&538 \cr
 & &128.3.7.19.37.227&259&1216& & &32.3.5.29.47.107.269&4089&4304 \cr
\noalign{\hrule}
 & &5.49.43.263.347&1921&186& & &27.5.13.23.71.337&5473&2278 \cr
7278&961434635&4.3.17.31.113.263&395&132&7296&965813355&4.3.169.17.67.421&805&2068 \cr
 & &32.9.5.11.79.113&1017&13904& & &32.5.7.11.23.47.67&3149&1232 \cr
\noalign{\hrule}
 & &3.5.53.89.107.127&451&184& & &17.19.61.139.353&33215&11682 \cr
7279&961489695&16.11.23.41.53.107&27&80&7297&966767101&4.9.5.7.11.13.59.73&181&38 \cr
 & &512.27.5.11.23.41&8487&2816& & &16.3.5.7.19.59.181&1239&7240 \cr
\noalign{\hrule}
 & &9.25.11.31.83.151&3935&3784& & &3.5.7.11.23.89.409&303&142 \cr
7280&961594425&16.3.125.121.43.787&41383&56992&7298&966992565&4.9.11.71.101.409&2231&1450 \cr
 & &1024.13.29.137.1427&195499&193024& & &16.25.23.29.97.101&2813&4040 \cr
\noalign{\hrule}
 & &3.5.13.19.277.937&8869&8934& & &3.7.11.29.241.599&365&1322 \cr
7281&961629045&4.9.49.181.277.1489&195833&73676&7299&967060941&4.5.73.599.661&31&630 \cr
 & &32.11.19.113.163.937&1793&1808& & &16.9.25.7.31.73&6789&200 \cr
\noalign{\hrule}
 & &243.5.49.11.13.113&551&664& & &9.121.29.113.271&3425&6406 \cr
7282&962026065&16.49.11.13.19.29.83&81&458&7300&967104963&4.3.25.11.137.3203&1189&2014 \cr
 & &64.81.19.83.229&4351&2656& & &16.19.29.41.53.137&5617&8056 \cr
\noalign{\hrule}
 & &3.11.169.17.73.139&835&406& & &3.7.97.577.823&607&1430 \cr
7283&962026923&4.5.7.13.29.139.167&675&298&7301&967312227&4.5.11.13.577.607&1139&1746 \cr
 & &16.27.125.149.167&20875&10728& & &16.9.11.13.17.67.97&2431&1608 \cr
\noalign{\hrule}
 & &3.5.31.71.151.193&667&88& & &9.17.19.43.71.109&9625&32558 \cr
7284&962156145&16.11.23.29.31.71&135&206&7302&967382739&4.125.7.11.73.223&19&54 \cr
 & &64.27.5.23.29.103&2369&8352& & &16.27.25.11.19.223&2453&600 \cr
\noalign{\hrule}
 & &7.11.13.89.101.107&11983&2460& & &5.11.13.17.19.59.71&75&134 \cr
7285&962784823&8.3.5.7.23.41.521&117&404&7303&967428605&4.3.125.13.17.67.71&209&1416 \cr
 & &64.27.5.13.23.101&115&864& & &64.9.11.19.59.67&67&288 \cr
\noalign{\hrule}
 & &49.11.13.23.43.139&485&1044& & &9.5.343.11.41.139&22711&24966 \cr
7286&963259297&8.9.5.49.23.29.97&407&260&7304&967604715&4.81.13.19.73.1747&77&4 \cr
 & &64.3.25.11.13.37.97&2425&3552& & &32.7.11.13.19.1747&1747&3952 \cr
\noalign{\hrule}
 & &13.41.79.137.167&335&198& & &27.5.11.23.29.977&109&868 \cr
7287&963366053&4.9.5.11.67.79.167&17947&21632&7305&967713615&8.9.5.7.29.31.109&1603&1558 \cr
 & &1024.3.169.131.137&1703&1536& & &32.49.19.31.41.229&62279&69616 \cr
\noalign{\hrule}
 & &9.5.7.13.23.53.193&55&634& & &27.5.49.11.53.251&4687&2090 \cr
7288&963418365&4.3.25.7.11.23.317&901&824&7306&967992795&4.25.121.19.43.109&477&2548 \cr
 & &64.17.53.103.317&5389&3296& & &32.9.49.13.43.53&43&208 \cr
\noalign{\hrule}
 & &27.25.11.13.67.149&703&638& & &27.43.47.113.157&9451&8290 \cr
7289&963609075&4.3.5.121.19.29.37.67&403&202&7307&968073147&4.5.13.47.727.829&22473&11696 \cr
 & &16.13.19.29.31.37.101&21793&23432& & &128.9.5.11.17.43.227&3859&3520 \cr
\noalign{\hrule}
 & &7.11.31.41.43.229&351&122& & &27.343.11.13.17.43&95&248 \cr
7290&963696349&4.27.7.13.31.41.61&5&1276&7308&968080113&16.3.5.11.13.19.31.43&833&586 \cr
 & &32.9.5.11.13.29&29&9360& & &64.5.49.17.31.293&1465&992 \cr
\noalign{\hrule}
}%
}
$$
\eject
\vglue -23 pt
\noindent\hskip 1 in\hbox to 6.5 in{\ 7309 -- 7344 \hfill\fbd 969786809 -- 989333895\frb}
\vskip -9 pt
$$
\vbox{
\nointerlineskip
\halign{\strut
    \vrule \ \ \hfil \frb #\ 
   &\vrule \hfil \ \ \fbb #\frb\ 
   &\vrule \hfil \ \ \frb #\ \hfil
   &\vrule \hfil \ \ \frb #\ 
   &\vrule \hfil \ \ \frb #\ \ \vrule \hskip 2 pt
   &\vrule \ \ \hfil \frb #\ 
   &\vrule \hfil \ \ \fbb #\frb\ 
   &\vrule \hfil \ \ \frb #\ \hfil
   &\vrule \hfil \ \ \frb #\ 
   &\vrule \hfil \ \ \frb #\ \vrule \cr%
\noalign{\hrule}
 & &19.89.587.977&9575&8988& & &81.25.13.17.37.59&1111&482 \cr
7309&969786809&8.3.25.7.89.107.383&913&1002&7327&976947075&4.3.25.11.13.101.241&8381&9694 \cr
 & &32.9.5.7.11.83.107.167&287595&285904& & &16.11.289.29.37.131&2227&2552 \cr
\noalign{\hrule}
 & &25.11.19.29.37.173&95427&90202& & &31.59.487.1097&1463&366 \cr
7310&969911525&4.9.7.17.23.379.461&1045&2182&7328&977123131&4.3.7.11.19.61.487&177&310 \cr
 & &16.3.5.11.17.19.23.1091&3273&3128& & &16.9.5.11.31.59.61&305&792 \cr
\noalign{\hrule}
 & &5.7.121.397.577&4203&164& & &243.7.41.107.131&10315&15686 \cr
7311&970107215&8.9.5.11.41.467&1173&1162&7329&977559597&4.5.7.11.23.31.2063&951&1112 \cr
 & &32.27.7.17.23.41.83&25461&22576& & &64.3.5.11.31.139.317&47399&50720 \cr
\noalign{\hrule}
 & &3.5.121.19.107.263&1807&492& & &25.53.67.73.151&6993&3124 \cr
7312&970442385&8.9.13.41.107.139&9835&5038&7330&978566825&8.27.25.7.11.37.71&73&2 \cr
 & &32.5.7.11.229.281&1603&4496& & &32.9.7.11.37.73&3663&112 \cr
\noalign{\hrule}
 & &25.17.19.31.3877&621&3256& & &9.7.11.17.29.47.61&137&50 \cr
7313&970510025&16.27.5.11.19.23.37&2023&2232&7331&979507683&4.3.25.7.47.61.137&341&646 \cr
 & &256.243.7.289.31&1701&2176& & &16.5.11.17.19.31.137&2603&1240 \cr
\noalign{\hrule}
 & &3.5.13.23.53.61.67&8987&9252& & &3.5.13.29.311.557&4087&3154 \cr
7314&971500335&8.27.11.19.43.67.257&115&182&7332&979598685&4.5.19.29.61.67.83&3421&666 \cr
 & &32.5.7.13.19.23.43.257&4883&4816& & &16.9.11.37.83.311&1221&664 \cr
\noalign{\hrule}
 & &9.11.43.53.59.73&2175&2132& & &27.289.29.61.71&275&14 \cr
7315&971749647&8.27.25.11.13.29.41.53&1829&344&7333&980048997&4.3.25.7.11.61.71&493&422 \cr
 & &128.5.13.29.31.43.59&1885&1984& & &16.5.7.11.17.29.211&1477&440 \cr
\noalign{\hrule}
 & &9.25.113.167.229&469&1034& & &9.5.73.107.2789&847&1942 \cr
7316&972328275&4.5.7.11.47.67.229&2147&1002&7334&980319555&4.3.7.121.107.971&325&646 \cr
 & &16.3.7.11.19.113.167&77&152& & &16.25.7.121.13.17.19&8645&16456 \cr
\noalign{\hrule}
 & &9.43.47.127.421&4625&836& & &9.5.121.19.31.307&13447&34138 \cr
7317&972511263&8.125.11.19.37.47&127&108&7335&984581235&4.7.169.17.101.113&57&44 \cr
 & &64.27.25.11.37.127&1221&800& & &32.3.7.11.13.17.19.113&1921&1456 \cr
\noalign{\hrule}
 & &9.11.19.23.113.199&65&2212& & &27.49.23.139.233&13&220 \cr
7318&972855081&8.5.7.13.79.199&703&690&7336&985504023&8.3.5.49.11.13.139&16813&17242 \cr
 & &32.3.25.19.23.37.79&925&1264& & &32.17.23.37.43.233&731&592 \cr
\noalign{\hrule}
 & &9.5.49.11.289.139&1867&338& & &27.25.7.23.43.211&1637&6490 \cr
7319&974347605&4.169.289.1867&945&2812&7337&986008275&4.125.11.59.1637&12691&5316 \cr
 & &32.27.5.7.13.19.37&247&1776& & &32.3.343.37.443&1813&7088 \cr
\noalign{\hrule}
 & &9.11.59.107.1559&80327&74014& & &9.25.11.13.23.31.43&697&722 \cr
7320&974354733&4.13.23.37.167.1609&2285&3894&7338&986453325&4.3.13.17.361.23.31.41&20525&7282 \cr
 & &16.3.5.11.13.23.59.457&2285&2392& & &16.25.11.19.331.821&6289&6568 \cr
\noalign{\hrule}
 & &9.13.29.31.73.127&4031&5240& & &9.5.11.13.23.59.113&323&436 \cr
7321&975151593&16.3.5.841.131.139&7843&10366&7339&986749335&8.3.5.13.17.19.59.109&4633&382 \cr
 & &64.5.11.23.31.71.73&1633&1760& & &32.19.41.113.191&779&3056 \cr
\noalign{\hrule}
 & &27.5.11.19.71.487&3961&1396& & &13.19.23.31.71.79&5309&3960 \cr
7322&975590055&8.17.71.233.349&429&778&7340&987806599&16.9.5.11.79.5309&1351&3958 \cr
 & &32.3.11.13.233.389&3029&6224& & &64.3.5.7.193.1979&29685&43232 \cr
\noalign{\hrule}
 & &5.17.19.29.83.251&22239&23494& & &27.7.31.227.743&7955&13156 \cr
7323&975713555&4.9.7.289.353.691&779&88&7341&988184799&8.9.5.11.13.23.37.43&797&238 \cr
 & &64.3.7.11.19.41.353&11649&9184& & &32.7.11.17.37.797&13549&6512 \cr
\noalign{\hrule}
 & &3.5.11.31.107.1783&1403&380& & &9.7.17.557.1657&43439&41782 \cr
7324&975844815&8.25.19.23.61.107&341&234&7342&988478379&4.7.121.13.359.1607&321&680 \cr
 & &32.9.11.13.19.31.61&1159&624& & &64.3.5.11.17.107.1607&17677&17120 \cr
\noalign{\hrule}
 & &27.25.7.13.23.691&2701&4774& & &13.29.43.181.337&26873&41364 \cr
7325&976227525&4.9.49.11.31.37.73&1955&1622&7343&988822367&8.27.7.11.349.383&1345&2494 \cr
 & &16.5.11.17.23.31.811&5797&6488& & &32.9.5.7.29.43.269&1345&1008 \cr
\noalign{\hrule}
 & &243.25.11.47.311&11021&3596& & &3.5.11.19.41.43.179&359&1404 \cr
7326&976781025&8.9.29.31.103.107&517&410&7344&989333895&8.81.13.179.359&15703&13376 \cr
 & &32.5.11.29.31.41.47&899&656& & &1024.11.19.41.383&383&512 \cr
\noalign{\hrule}
}%
}
$$
\eject
\vglue -23 pt
\noindent\hskip 1 in\hbox to 6.5 in{\ 7345 -- 7380 \hfill\fbd 989520945 -- 1003083939\frb}
\vskip -9 pt
$$
\vbox{
\nointerlineskip
\halign{\strut
    \vrule \ \ \hfil \frb #\ 
   &\vrule \hfil \ \ \fbb #\frb\ 
   &\vrule \hfil \ \ \frb #\ \hfil
   &\vrule \hfil \ \ \frb #\ 
   &\vrule \hfil \ \ \frb #\ \ \vrule \hskip 2 pt
   &\vrule \ \ \hfil \frb #\ 
   &\vrule \hfil \ \ \fbb #\frb\ 
   &\vrule \hfil \ \ \frb #\ \hfil
   &\vrule \hfil \ \ \frb #\ 
   &\vrule \hfil \ \ \frb #\ \vrule \cr%
\noalign{\hrule}
 & &3.5.49.43.131.239&1119&988& & &25.7.13.17.149.173&2987&738 \cr
7345&989520945&8.9.5.13.19.239.373&143&2008&7363&996925475&4.9.7.17.29.41.103&4345&7912 \cr
 & &128.11.169.19.251&35321&16064& & &64.3.5.11.23.43.79&5451&15136 \cr
\noalign{\hrule}
 & &5.529.31.47.257&187&342& & &3.125.7.47.59.137&7843&14282 \cr
7346&990417605&4.9.11.17.19.47.257&6107&7000&7364&997240125&4.7.11.23.31.37.193&1075&6066 \cr
 & &64.3.125.7.11.31.197&5775&6304& & &16.9.25.11.43.337&1011&3784 \cr
\noalign{\hrule}
 & &25.11.13.241.1151&4203&1552& & &9.11.61.197.839&7267&4750 \cr
7347&991672825&32.9.5.13.97.467&201&266&7365&998144037&4.3.125.11.169.19.43&197&362 \cr
 & &128.27.7.19.67.97&45493&32832& & &16.25.13.19.181.197&3439&2600 \cr
\noalign{\hrule}
 & &11.17.61.227.383&4817&9030& & &3.11.13.19.29.41.103&6059&6610 \cr
7348&991735987&4.3.5.7.17.43.4817&2043&2774&7366&998228517&4.5.11.13.73.83.661&3895&4698 \cr
 & &16.27.5.7.19.73.227&3591&2920& & &16.81.25.19.29.41.83&675&664 \cr
\noalign{\hrule}
 & &5.11.13.29.109.439&483&2678& & &81.7.71.137.181&41375&26714 \cr
7349&992190485&4.3.7.11.169.23.103&981&878&7367&998252829&4.125.361.37.331&273&88 \cr
 & &16.27.7.23.109.439&161&216& & &64.3.25.7.11.13.331&8275&4576 \cr
\noalign{\hrule}
 & &3.5.7.13.17.127.337&1111&3270& & &9.5.11.29.79.881&12419&13130 \cr
7350&993150795&4.9.25.7.11.101.109&403&578&7368&999093645&4.25.121.13.101.1129&2077&948 \cr
 & &16.11.13.289.31.101&1717&2728& & &32.3.13.31.67.79.101&6767&6448 \cr
\noalign{\hrule}
 & &27.25.7.13.19.23.37&2867&3608& & &7.17.31.439.617&115985&115368 \cr
7351&993180825&16.9.11.23.41.47.61&815&266&7369&999213607&16.3.5.7.11.19.23.23197&3951&19246 \cr
 & &64.5.7.11.19.41.163&1793&1312& & &64.27.11.439.9623&9623&9504 \cr
\noalign{\hrule}
 & &3.7.11.23.31.37.163&95&312& & &9.11.19.47.89.127&55775&45952 \cr
7352&993323793&16.9.5.13.19.23.163&1445&1652&7370&999264321&256.25.23.97.359&231&254 \cr
 & &128.25.7.13.289.59&17051&20800& & &1024.3.5.7.11.127.359&2513&2560 \cr
\noalign{\hrule}
 & &49.17.31.79.487&6325&32148& & &9.25.7.11.13.23.193&901&824 \cr
7353&993488279&8.9.25.11.19.23.47&707&1972&7371&999773775&16.3.13.17.53.103.193&55&634 \cr
 & &64.3.5.7.17.29.101&2929&480& & &64.5.11.17.103.317&5389&3296 \cr
\noalign{\hrule}
 & &9.7.11.17.19.23.193&841&1420& & &9.25.17.197.1327&253&338 \cr
7354&993621321&8.3.5.11.23.841.71&1409&650&7372&999927675&4.3.5.11.169.23.1327&31331&11426 \cr
 & &32.125.13.29.1409&40861&26000& & &16.17.19.29.97.197&551&776 \cr
\noalign{\hrule}
 & &17.23.359.7079&589&7668& & &9.5.7.97.137.239&8509&14674 \cr
7355&993672151&8.27.17.19.31.71&1795&1826&7373&1000462365&4.7.11.23.29.67.127&411&478 \cr
 & &32.9.5.11.19.83.359&3735&3344& & &16.3.11.23.29.137.239&253&232 \cr
\noalign{\hrule}
 & &243.13.19.59.281&21721&38300& & &3.7.29.31.37.1433&6625&3406 \cr
7356&995088159&8.25.7.29.107.383&583&1332&7374&1000983459&4.125.13.31.53.131&4959&3256 \cr
 & &64.9.5.11.29.37.53&7685&13024& & &64.9.25.11.19.29.37&825&608 \cr
\noalign{\hrule}
 & &9.5.19.59.109.181&2401&4030& & &3.5.13.361.41.347&2075&2436 \cr
7357&995229405&4.25.2401.13.19.31&5133&9592&7375&1001509665&8.9.125.7.29.41.83&28769&7106 \cr
 & &64.3.7.11.29.59.109&319&224& & &32.11.13.17.19.2213&2213&2992 \cr
\noalign{\hrule}
 & &5.49.11.37.67.149&585&1054& & &243.11.13.19.37.41&2461&2830 \cr
7358&995454845&4.9.25.7.13.17.31.37&22177&24998&7376&1001570427&4.27.5.19.23.107.283&2537&352 \cr
 & &16.3.29.67.331.431&9599&10344& & &256.11.43.59.283&12169&7552 \cr
\noalign{\hrule}
 & &3.7.11.83.127.409&1805&2386& & &5.13.31.41.67.181&4731&7396 \cr
7359&995903139&4.5.361.409.1193&3289&4482&7377&1001872105&8.3.19.31.1849.83&2211&362 \cr
 & &16.27.5.11.13.19.23.83&1495&1368& & &32.9.11.19.67.181&99&304 \cr
\noalign{\hrule}
 & &3.5.7.23.229.1801&99&1702& & &3.5.11.23.41.47.137&147&106 \cr
7360&996016035&4.27.5.11.529.37&2977&2842&7378&1001876205&4.9.5.49.47.53.137&11033&3772 \cr
 & &16.49.13.29.37.229&1073&728& & &32.7.11.17.23.41.59&413&272 \cr
\noalign{\hrule}
 & &9.47.59.107.373&26243&26350& & &5.13.289.197.271&2003&558 \cr
7361&996058827&4.3.25.7.17.23.31.59.163&19811&22256&7379&1002874795&4.9.31.271.2003&595&1408 \cr
 & &128.5.7.11.13.17.107.1801&99055&99008& & &1024.3.5.7.11.17.31&1023&3584 \cr
\noalign{\hrule}
 & &27.23.31.53.977&845&132& & &27.11.13.61.4259&26677&28690 \cr
7362&996836031&8.81.5.11.169.53&101&790&7380&1003083939&4.9.5.7.19.37.103.151&187&2144 \cr
 & &32.25.13.79.101&1313&31600& & &256.5.11.17.67.151&12835&8576 \cr
\noalign{\hrule}
}%
}
$$
\eject
\vglue -23 pt
\noindent\hskip 1 in\hbox to 6.5 in{\ 7381 -- 7416 \hfill\fbd 1003559843 -- 1016803449\frb}
\vskip -9 pt
$$
\vbox{
\nointerlineskip
\halign{\strut
    \vrule \ \ \hfil \frb #\ 
   &\vrule \hfil \ \ \fbb #\frb\ 
   &\vrule \hfil \ \ \frb #\ \hfil
   &\vrule \hfil \ \ \frb #\ 
   &\vrule \hfil \ \ \frb #\ \ \vrule \hskip 2 pt
   &\vrule \ \ \hfil \frb #\ 
   &\vrule \hfil \ \ \fbb #\frb\ 
   &\vrule \hfil \ \ \frb #\ \hfil
   &\vrule \hfil \ \ \frb #\ 
   &\vrule \hfil \ \ \frb #\ \vrule \cr%
\noalign{\hrule}
 & &121.13.37.43.401&9&1582& & &9.5.11.13.59.2657&1081&1576 \cr
7381&1003559843&4.9.7.113.401&31&370&7399&1008769905&16.13.23.47.59.197&285&482 \cr
 & &16.3.5.7.31.37&5&5208& & &64.3.5.19.23.47.241&11327&13984 \cr
\noalign{\hrule}
 & &3.7.79.151.4007&3589&418& & &27.7.11.13.163.229&551&590 \cr
7382&1003789563&4.11.19.37.79.97&1413&1510&7400&1008836829&4.9.5.11.19.29.59.229&2587&68 \cr
 & &16.9.5.11.19.151.157&1727&2280& & &32.13.17.19.29.199&9367&3184 \cr
\noalign{\hrule}
 & &9.7.11.43.59.571&95105&100244& & &3.13.289.19.53.89&145&1012 \cr
7383&1003898511&8.5.19.23.827.1319&4559&11154&7401&1010140833&8.5.11.19.23.29.53&351&86 \cr
 & &32.3.11.169.23.47.97&16393&17296& & &32.27.11.13.29.43&387&5104 \cr
\noalign{\hrule}
 & &3.25.19.529.31.43&7627&7098& & &27.5.49.11.17.19.43&11461&4366 \cr
7384&1004848725&4.9.7.169.29.43.263&6355&17578&7402&1010633085&4.9.37.59.73.157&23987&28294 \cr
 & &16.5.11.13.17.31.41.47&5863&6392& & &16.7.289.43.47.83&799&664 \cr
\noalign{\hrule}
 & &7.13.31.593.601&99&502& & &23.29.37.71.577&34353&13004 \cr
7385&1005384653&4.9.7.11.251.593&2201&1950&7403&1011024593&8.9.11.347.3251&1105&2146 \cr
 & &16.27.25.11.13.31.71&1775&2376& & &32.3.5.11.13.17.29.37&1105&528 \cr
\noalign{\hrule}
 & &3.121.19.41.3557&2947&610& & &7.121.17.23.43.71&1103&1950 \cr
7386&1005837789&4.5.7.121.61.421&629&1476&7404&1011083381&4.3.25.13.17.23.1103&1529&426 \cr
 & &32.9.17.37.41.61&629&2928& & &16.9.5.11.13.71.139&1251&520 \cr
\noalign{\hrule}
 & &3.5.13.19.37.41.179&4249&2732& & &243.7.19.113.277&517&274 \cr
7387&1006066815&8.5.7.19.607.683&4059&7474&7405&1011617019&4.11.19.47.137.277&2825&222 \cr
 & &32.9.7.11.37.41.101&707&528& & &16.3.25.37.47.113&1739&200 \cr
\noalign{\hrule}
 & &25.19.67.101.313&1419&1106& & &25.169.23.29.359&65391&57134 \cr
7388&1006083725&4.3.7.11.19.43.67.79&2191&690&7406&1011688925&4.3.49.11.53.71.307&2987&390 \cr
 & &16.9.5.49.11.23.313&1127&792& & &16.9.5.13.29.71.103&927&568 \cr
\noalign{\hrule}
 & &121.29.281.1021&1035&2056& & &67.2243.6733&78507&71774 \cr
7389&1006735609&16.9.5.11.23.29.257&339&328&7407&1011841973&4.9.11.13.17.61.2111&2245&134 \cr
 & &256.27.5.41.113.257&145205&141696& & &16.3.5.11.17.67.449&2245&4488 \cr
\noalign{\hrule}
 & &9.5.43.587.887&7337&32578& & &3.11.23.43.101.307&925&1398 \cr
7390&1007494515&4.7.11.13.23.29.179&1037&1290&7408&1011975459&4.9.25.37.233.307&437&10922 \cr
 & &16.3.5.7.17.29.43.61&1037&1624& & &16.5.19.23.43.127&95&1016 \cr
\noalign{\hrule}
 & &11.529.1681.103&2069&3750& & &25.13.17.19.31.311&511&264 \cr
7391&1007519117&4.3.625.103.2069&253&2322&7409&1012063975&16.3.7.11.17.73.311&249&62 \cr
 & &16.81.25.11.23.43&2025&344& & &64.9.7.31.73.83&6059&2016 \cr
\noalign{\hrule}
 & &9.7.169.17.19.293&15&308& & &3.49.13.431.1229&2123&894 \cr
7392&1007621433&8.27.5.49.11.169&7361&7192&7410&1012254789&4.9.7.11.13.149.193&1229&410 \cr
 & &128.5.17.29.31.433&12557&9920& & &16.5.41.193.1229&965&328 \cr
\noalign{\hrule}
 & &9.49.11.13.19.841&235&26& & &49.11.17.31.43.83&435&478 \cr
7393&1007683677&4.5.49.169.29.47&795&626&7411&1013785157&4.3.5.49.17.29.31.239&473&1992 \cr
 & &16.3.25.47.53.313&7825&19928& & &64.9.11.43.83.239&239&288 \cr
\noalign{\hrule}
 & &3.11.19.29.157.353&3627&6610& & &9.7.11.593.2467&30587&28120 \cr
7394&1007720043&4.27.5.11.13.31.661&6707&11140&7412&1013811183&16.5.7.19.37.73.419&2001&4934 \cr
 & &32.25.19.353.557&557&400& & &64.3.23.29.37.2467&851&928 \cr
\noalign{\hrule}
 & &81.625.43.463&1429&2054& & &3.625.13.17.31.79&6853&3772 \cr
7395&1007893125&4.13.79.463.1429&27577&9000&7413&1014804375&8.7.11.23.31.41.89&1125&146 \cr
 & &64.9.125.11.23.109&1199&736& & &32.9.125.7.23.73&1679&336 \cr
\noalign{\hrule}
 & &25.41.43.137.167&4103&3078& & &3.11.31.37.47.571&415363&416584 \cr
7396&1008391925&4.81.11.19.137.373&1435&2668&7414&1015807287&16.7.13.43.89.173.359&15417&20 \cr
 & &32.9.5.7.19.23.29.41&3933&3248& & &128.27.5.7.13.571&315&832 \cr
\noalign{\hrule}
 & &3.19.23.41.73.257&2551&2332& & &3.11.37.53.113.139&7771&7658 \cr
7397&1008422511&8.11.23.41.53.2551&62415&72788&7415&1016447091&4.7.11.19.53.409.547&565&18 \cr
 & &64.9.5.19.31.73.587&2935&2976& & &16.9.5.7.19.113.409&2863&2280 \cr
\noalign{\hrule}
 & &11.13.43.61.2689&31905&3052& & &81.19.61.10831&5995&4836 \cr
7398&1008614321&8.9.5.7.109.709&1009&1118&7416&1016803449&8.243.5.11.13.31.109&1525&1634 \cr
 & &32.3.5.7.13.43.1009&1009&1680& & &32.125.11.19.31.43.61&5375&5456 \cr
\noalign{\hrule}
}%
}
$$
\eject
\vglue -23 pt
\noindent\hskip 1 in\hbox to 6.5 in{\ 7417 -- 7452 \hfill\fbd 1018966485 -- 1031277265\frb}
\vskip -9 pt
$$
\vbox{
\nointerlineskip
\halign{\strut
    \vrule \ \ \hfil \frb #\ 
   &\vrule \hfil \ \ \fbb #\frb\ 
   &\vrule \hfil \ \ \frb #\ \hfil
   &\vrule \hfil \ \ \frb #\ 
   &\vrule \hfil \ \ \frb #\ \ \vrule \hskip 2 pt
   &\vrule \ \ \hfil \frb #\ 
   &\vrule \hfil \ \ \fbb #\frb\ 
   &\vrule \hfil \ \ \frb #\ \hfil
   &\vrule \hfil \ \ \frb #\ 
   &\vrule \hfil \ \ \frb #\ \vrule \cr%
\noalign{\hrule}
 & &3.5.17.19.43.67.73&9251&12390& & &7.11.19.31.97.233&2115&272 \cr
7417&1018966485&4.9.25.7.11.841.59&13433&7592&7435&1025023153&32.9.5.17.47.233&2051&1910 \cr
 & &64.49.13.19.73.101&1313&1568& & &128.3.25.7.191.293&14325&18752 \cr
\noalign{\hrule}
 & &3.7.11.13.19.53.337&4265&558& & &23.47.311.3049&27755&42372 \cr
7418&1019095077&4.27.5.19.31.853&371&466&7436&1025046359&8.9.5.7.11.13.61.107&1369&2162 \cr
 & &16.7.53.233.853&853&1864& & &32.3.5.7.23.1369.47&1369&1680 \cr
\noalign{\hrule}
 & &3.5.7.11.19.23.43.47&611&654& & &5.7.31.43.127.173&209&426 \cr
7419&1020069435&4.9.7.13.19.2209.109&65621&23650&7437&1025057005&4.3.11.19.43.71.173&11139&3700 \cr
 & &16.25.11.43.211.311&1555&1688& & &32.9.25.37.47.79&15651&6320 \cr
\noalign{\hrule}
 & &9.5.13.17.37.47.59&2173&2156& & &3.7.41.47.73.347&13&60 \cr
7420&1020366945&8.5.49.11.41.47.53.59&17&312&7438&1025069577&8.9.5.7.13.41.347&8471&3674 \cr
 & &128.3.7.11.13.17.41.53&3157&3392& & &32.11.43.167.197&32899&7568 \cr
\noalign{\hrule}
 & &25.13.19.29.41.139&891&916& & &27.7.11.23.89.241&9685&11764 \cr
7421&1020548425&8.81.11.19.29.41.229&13&13066&7439&1025626833&8.5.13.17.23.149.173&723&3256 \cr
 & &32.27.13.47.139&27&752& & &128.3.5.11.13.37.241&481&320 \cr
\noalign{\hrule}
 & &7.29.211.23827&9503&33330& & &11.13.29.31.79.101&669&1780 \cr
7422&1020581891&4.3.5.11.13.17.43.101&1189&1176&7440&1025756303&8.3.5.13.29.89.223&411&34 \cr
 & &64.9.49.17.29.41.101&12019&11808& & &32.9.17.137.223&2329&32112 \cr
\noalign{\hrule}
 & &3.25.7.79.103.239&991&682& & &9.5.11.13.31.37.139&733&1074 \cr
7423&1020990075&4.25.11.31.79.991&4223&4302&7441&1025951355&4.27.5.37.179.733&1333&2332 \cr
 & &16.9.41.103.239.991&991&984& & &32.11.31.43.53.179&2279&2864 \cr
\noalign{\hrule}
 & &243.25.7.11.37.59&5627&848& & &3.7.11.17.43.59.103&43&60 \cr
7424&1021152825&32.3.11.17.53.331&403&590&7442&1026168297&8.9.5.7.11.1849.59&8131&21074 \cr
 & &128.5.13.31.53.59&689&1984& & &32.41.47.173.257&44461&30832 \cr
\noalign{\hrule}
 & &3.7.121.37.83.131&1075&668& & &5.23.47.179.1061&4711&594 \cr
7425&1022246841&8.25.11.43.131.167&3735&1898&7443&1026512195&4.27.7.11.47.673&3193&2864 \cr
 & &32.9.125.13.73.83&2847&2000& & &128.3.11.31.103.179&3193&2112 \cr
\noalign{\hrule}
 & &13.31.107.131.181&12551&11160& & &83.113.223.491&4935&4444 \cr
7426&1022442031&16.9.5.7.11.961.163&1703&742&7444&1026934847&8.3.5.7.11.47.101.223&229&6 \cr
 & &64.3.49.11.13.53.131&1749&1568& & &32.9.7.11.101.229&7777&32976 \cr
\noalign{\hrule}
 & &25.11.17.23.37.257&457&468& & &3.125.49.11.13.17.23&305&696 \cr
7427&1022455225&8.9.13.17.23.257.457&3071&7440&7445&1027401375&16.9.625.7.29.61&2599&3026 \cr
 & &256.27.5.13.31.37.83&10881&10624& & &64.17.23.29.89.113&3277&2848 \cr
\noalign{\hrule}
 & &5.49.37.97.1163&19521&37466& & &3.5.11.13.19.151.167&2611&106 \cr
7428&1022631715&4.81.11.13.131.241&1213&490&7446&1027718835&4.7.53.151.373&715&342 \cr
 & &16.27.5.49.11.1213&1213&2376& & &16.9.5.11.13.19.53&159&8 \cr
\noalign{\hrule}
 & &5.11.19.53.59.313&18459&36926& & &9.5.11.17.19.59.109&83&26 \cr
7429&1022794795&4.9.7.37.293.499&689&190&7447&1028220435&4.3.5.11.13.17.59.83&133&428 \cr
 & &16.3.5.7.13.19.37.53&481&168& & &32.7.13.19.83.107&7553&1712 \cr
\noalign{\hrule}
 & &5.11.103.419.431&2623&3042& & &25.7.41.193.743&3113&2088 \cr
7430&1023036685&4.9.169.43.61.431&1133&6736&7448&1028887825&16.9.11.29.193.283&637&1486 \cr
 & &128.3.11.13.103.421&1263&832& & &64.3.49.13.29.743&609&416 \cr
\noalign{\hrule}
 & &243.5.11.13.43.137&151&14& & &9.7.11.13.163.701&5453&3660 \cr
7431&1023531795&4.81.7.13.43.151&505&548&7449&1029395367&8.27.5.49.19.41.61&701&946 \cr
 & &32.5.7.101.137.151&1057&1616& & &32.11.19.41.43.701&779&688 \cr
\noalign{\hrule}
 & &3.11.17.37.149.331&149&38& & &81.5.23.31.43.83&301&2210 \cr
7432&1023714483&4.19.22201.331&14245&7956&7450&1030602285&4.25.7.13.17.1849&837&1012 \cr
 & &32.9.5.7.11.13.17.37&91&240& & &32.27.11.13.17.23.31&143&272 \cr
\noalign{\hrule}
 & &3.11.29.43.149.167&181&138& & &125.121.41.1663&11571&13234 \cr
7433&1023960333&4.9.23.149.167.181&31175&6292&7451&1031267875&4.3.25.7.13.19.29.509&123&9548 \cr
 & &32.25.121.13.29.43&325&176& & &32.9.49.11.31.41&1519&144 \cr
\noalign{\hrule}
 & &243.11.23.79.211&2755&2834& & &5.29.89.157.509&7425&7336 \cr
7434&1024793451&4.5.11.13.19.29.109.211&69&4078&7452&1031277265&16.27.125.7.11.131.157&19&1394 \cr
 & &16.3.5.23.109.2039&2039&4360& & &64.3.7.17.19.41.131&42313&27552 \cr
\noalign{\hrule}
}%
}
$$
\eject
\vglue -23 pt
\noindent\hskip 1 in\hbox to 6.5 in{\ 7453 -- 7488 \hfill\fbd 1031353269 -- 1043227185\frb}
\vskip -9 pt
$$
\vbox{
\nointerlineskip
\halign{\strut
    \vrule \ \ \hfil \frb #\ 
   &\vrule \hfil \ \ \fbb #\frb\ 
   &\vrule \hfil \ \ \frb #\ \hfil
   &\vrule \hfil \ \ \frb #\ 
   &\vrule \hfil \ \ \frb #\ \ \vrule \hskip 2 pt
   &\vrule \ \ \hfil \frb #\ 
   &\vrule \hfil \ \ \fbb #\frb\ 
   &\vrule \hfil \ \ \frb #\ \hfil
   &\vrule \hfil \ \ \frb #\ 
   &\vrule \hfil \ \ \frb #\ \vrule \cr%
\noalign{\hrule}
 & &3.19.83.277.787&3025&2238& & &9.49.53.157.283&2945&2998 \cr
7453&1031353269&4.9.25.121.83.373&8771&554&7471&1038485763&4.3.5.7.19.31.157.1499&1817&1166 \cr
 & &16.49.11.179.277&1969&392& & &16.5.11.23.53.79.1499&34477&34760 \cr
\noalign{\hrule}
 & &9.5.11.23.31.37.79&2363&86& & &27.361.47.2267&687&1580 \cr
7454&1031629005&4.5.17.37.43.139&3081&2896&7472&1038533103&8.81.5.19.79.229&14077&7678 \cr
 & &128.3.13.17.79.181&2353&1088& & &32.7.11.349.2011&22121&39088 \cr
\noalign{\hrule}
 & &3.5.7.13.29.67.389&12877&12408& & &27.5.11.19.131.281&32657&5278 \cr
7455&1031703855&16.9.11.29.47.79.163&1415&52&7473&1038622365&4.7.13.289.29.113&1935&1822 \cr
 & &128.5.11.13.79.283&3113&5056& & &16.9.5.7.29.43.911&8729&7288 \cr
\noalign{\hrule}
 & &9.19.29.43.47.103&3395&4642& & &3.13.289.37.47.53&1125&836 \cr
7456&1032280317&4.5.7.11.97.103.211&10491&9976&7474&1038814257&8.27.125.11.13.19.47&17&1252 \cr
 & &64.3.7.11.13.29.43.269&2959&2912& & &64.25.11.17.313&7825&352 \cr
\noalign{\hrule}
 & &11.29.97.173.193&437&630& & &3.17.31.47.71.197&627&580 \cr
7457&1033155827&4.9.5.7.19.23.29.173&15611&3674&7475&1039332009&8.9.5.11.19.29.31.197&337&3406 \cr
 & &16.3.11.67.167.233&11189&5592& & &32.5.13.29.131.337&44147&30160 \cr
\noalign{\hrule}
 & &27.49.11.19.37.101&1207&2530& & &27.13.53.199.281&31&22 \cr
7458&1033306659&4.5.121.17.19.23.71&101&222&7476&1040261157&4.3.11.13.31.199.281&103615&99962 \cr
 & &16.3.5.23.37.71.101&355&184& & &16.5.17.23.53.151.331&28135&27784 \cr
\noalign{\hrule}
 & &3.11.19.67.73.337&2951&1940& & &27.5.7.11.101.991&3401&3536 \cr
7459&1033463409&8.5.11.13.19.97.227&511&4824&7477&1040445945&32.11.13.17.19.101.179&1843&126 \cr
 & &128.9.7.13.67.73&273&64& & &128.9.7.13.361.97&4693&6208 \cr
\noalign{\hrule}
 & &3.7.11.13.23.71.211&10509&10720& & &3.5.11.169.67.557&47435&46698 \cr
7460&1034722689&64.9.5.7.11.31.67.113&185&284&7478&1040640315&4.9.25.43.53.179.181&9581&94 \cr
 & &512.25.31.37.71.113&28675&28928& & &16.11.13.47.67.181&181&376 \cr
\noalign{\hrule}
 & &49.11.17.37.43.71&13395&15502& & &9.25.169.31.883&407&438 \cr
7461&1035061643&4.3.5.17.19.23.47.337&231&568&7479&1040858325&4.27.5.11.37.73.883&3193&1222 \cr
 & &64.9.5.7.11.19.23.71&855&736& & &16.11.13.31.37.47.103&3811&4136 \cr
\noalign{\hrule}
 & &9.49.11.37.73.79&16791&12920& & &9.5.7.11.17.23.769&1081&2764 \cr
7462&1035101529&16.27.5.17.19.29.193&18353&9632&7480&1041852735&8.7.529.47.691&10013&14850 \cr
 & &1024.7.43.18353&18353&22016& & &32.27.25.11.17.19.31&465&304 \cr
\noalign{\hrule}
 & &9.5.7.11.961.311&6063&4508& & &5.7.17.19.37.47.53&5925&6248 \cr
7463&1035588015&8.27.343.23.43.47&463&806&7481&1041947935&16.3.125.11.53.71.79&1819&1944 \cr
 & &32.13.23.31.43.463&10649&8944& & &256.729.11.17.79.107&92983&93312 \cr
\noalign{\hrule}
 & &5.19.29.313.1201&1221&1534& & &9.5.7.13.289.881&1793&2612 \cr
7464&1035640315&4.3.11.13.37.59.1201&7825&7788&7482&1042623855&8.11.289.163.653&471&182 \cr
 & &32.9.25.121.3481.313&17405&17424& & &32.3.7.11.13.157.163&1727&2608 \cr
\noalign{\hrule}
 & &17.23.71.107.349&64431&72028& & &3.5.13.23.37.61.103&593&902 \cr
7465&1036679023&8.9.11.1637.7159&6035&1124&7483&1042632435&4.11.37.41.61.593&257&24570 \cr
 & &64.3.5.11.17.71.281&1405&1056& & &16.27.5.7.13.257&2313&56 \cr
\noalign{\hrule}
 & &7.257.569.1013&73623&72610& & &9.5.19.41.71.419&1133&962 \cr
7466&1036938203&4.3.5.7.11.23.53.97.137&1013&54&7484&1042851195&4.11.13.37.41.71.103&1149&368 \cr
 & &16.81.5.23.53.1013&1863&2120& & &128.3.13.23.103.383&30797&24512 \cr
\noalign{\hrule}
 & &7.17.113.229.337&745&858& & &125.11.13.23.43.59&8051&4806 \cr
7467&1037745331&4.3.5.11.13.17.149.337&8001&10534&7485&1043024125&4.27.25.83.89.97&7429&1204 \cr
 & &16.27.7.13.23.127.229&2921&2808& & &32.9.7.17.19.23.43&1071&304 \cr
\noalign{\hrule}
 & &9.5.11.13.29.67.83&599&532& & &3.5.7.13.17.79.569&4323&5350 \cr
7468&1037766015&8.3.5.7.11.19.83.599&33031&86&7486&1043087955&4.9.125.7.11.107.131&1027&152 \cr
 & &32.17.29.43.67&731&16& & &64.11.13.19.79.107&1177&608 \cr
\noalign{\hrule}
 & &11.17.43.337.383&23&360& & &243.11.37.53.199&4595&4396 \cr
7469&1037859911&16.9.5.11.17.23.43&245&228&7487&1043108847&8.5.7.11.53.157.919&1791&64 \cr
 & &128.27.25.49.19.23&23275&39744& & &1024.9.199.919&919&512 \cr
\noalign{\hrule}
 & &11.17.43.337.383&245&228& & &3.5.7.11.13.17.61.67&699&1436 \cr
7470&1037859911&8.3.5.49.19.337.383&23&360&7488&1043227185&8.9.13.17.233.359&1159&938 \cr
 & &128.27.25.49.19.23&23275&39744& & &32.7.19.61.67.359&359&304 \cr
\noalign{\hrule}
}%
}
$$
\eject
\vglue -23 pt
\noindent\hskip 1 in\hbox to 6.5 in{\ 7489 -- 7524 \hfill\fbd 1043255367 -- 1055125357\frb}
\vskip -9 pt
$$
\vbox{
\nointerlineskip
\halign{\strut
    \vrule \ \ \hfil \frb #\ 
   &\vrule \hfil \ \ \fbb #\frb\ 
   &\vrule \hfil \ \ \frb #\ \hfil
   &\vrule \hfil \ \ \frb #\ 
   &\vrule \hfil \ \ \frb #\ \ \vrule \hskip 2 pt
   &\vrule \ \ \hfil \frb #\ 
   &\vrule \hfil \ \ \fbb #\frb\ 
   &\vrule \hfil \ \ \frb #\ \hfil
   &\vrule \hfil \ \ \frb #\ 
   &\vrule \hfil \ \ \frb #\ \vrule \cr%
\noalign{\hrule}
 & &9.7.11.23.29.37.61&65&268& & &11.29.41.191.421&39095&21834 \cr
7489&1043255367&8.5.11.13.23.61.67&1917&2170&7507&1051695469&4.9.5.7.1117.1213&1261&2378 \cr
 & &32.27.25.7.13.31.71&6603&5200& & &16.3.5.7.13.29.41.97&1261&840 \cr
\noalign{\hrule}
 & &3.25.7.31.59.1087&9119&3684& & &3.5.11.19.23.29.503&2149&1646 \cr
7490&1043764575&8.9.5.11.307.829&12173&4712&7508&1051795635&4.7.19.29.307.823&9017&14850 \cr
 & &128.7.19.31.37.47&893&2368& & &16.27.25.7.11.71.127&4445&5112 \cr
\noalign{\hrule}
 & &9.25.11.23.59.311&87&28& & &27.5.67.103.1129&25351&25454 \cr
7491&1044516825&8.27.5.7.11.29.311&247&3668&7509&1051815915&4.3.11.13.67.89.101.251&7123&1160 \cr
 & &64.49.13.19.131&32357&1568& & &64.5.13.17.29.101.419&92599&93728 \cr
\noalign{\hrule}
 & &729.11.2209.59&2623&5396& & &7.169.23.29.31.43&165&4 \cr
7492&1045124289&8.19.43.47.61.71&2077&1260&7510&1051818313&8.3.5.11.29.31.43&507&826 \cr
 & &64.9.5.7.31.61.67&14539&9760& & &32.9.5.7.169.59&45&944 \cr
\noalign{\hrule}
 & &27.11.13.421.643&959&970& & &9.5.121.73.2647&1399&4046 \cr
7493&1045184283&4.9.5.7.13.97.137.421&57211&466&7511&1052142795&4.7.289.73.1399&79&1320 \cr
 & &16.49.11.233.743&11417&5944& & &64.3.5.7.11.17.79&119&2528 \cr
\noalign{\hrule}
 & &9.125.11.289.293&2303&2678& & &9.625.19.43.229&971&154 \cr
7494&1047877875&4.3.49.11.13.17.47.103&25&586&7512&1052398125&4.5.7.11.229.971&1333&9348 \cr
 & &16.25.49.103.293&103&392& & &32.3.19.31.41.43&31&656 \cr
\noalign{\hrule}
 & &3.5.11.361.29.607&4393&4712& & &3.25.11.13.17.23.251&3027&3248 \cr
7495&1048522695&16.6859.23.31.191&469&6390&7513&1052562225&32.9.7.11.23.29.1009&13&680 \cr
 & &64.9.5.7.23.67.71&10787&6816& & &512.5.13.17.1009&1009&256 \cr
\noalign{\hrule}
 & &49.23.29.97.331&2565&248& & &3.11.41.2809.277&32485&1586 \cr
7496&1049353081&16.27.5.7.19.23.31&1&22&7514&1052759829&4.5.13.61.73.89&697&252 \cr
 & &64.9.5.11.19.31&32395&288& & &32.9.7.17.41.61&119&2928 \cr
\noalign{\hrule}
 & &27.11.13.23.53.223&119&2572& & &3.7.11.17.19.29.487&45&164 \cr
7497&1049562657&8.3.7.17.53.643&779&1150&7515&1053759399&8.27.5.29.41.487&635&148 \cr
 & &32.25.17.19.23.41&779&6800& & &64.25.37.41.127&4699&32800 \cr
\noalign{\hrule}
 & &9.7.29.37.53.293&2635&584& & &3.5.11.73.89.983&37721&49766 \cr
7498&1049744871&16.3.5.17.31.53.73&5567&638&7516&1053780915&4.67.149.167.563&365&198 \cr
 & &64.11.19.29.293&19&352& & &16.9.5.11.67.73.149&447&536 \cr
\noalign{\hrule}
 & &9.43.131.139.149&9031&8900& & &5.7.13.17.29.37.127&2537&1908 \cr
7499&1049985567&8.3.25.11.89.149.821&473&12788&7517&1054056185&8.9.13.29.43.53.59&253&124 \cr
 & &64.5.121.23.43.139&605&736& & &64.3.11.23.31.53.59&54219&43424 \cr
\noalign{\hrule}
 & &9.7.11.19.47.1697&9827&8840& & &27.631.61879&30953&30926 \cr
7500&1050186753&16.3.5.13.17.19.31.317&9419&11186&7518&1054232523&4.7.13.2209.631.2381&30305&648 \cr
 & &64.7.289.47.9419&9419&9248& & &64.81.5.7.11.19.29.47&31255&30624 \cr
\noalign{\hrule}
 & &11.241.607.653&199511&196860& & &9.13.29.167.1861&17215&36754 \cr
7501&1050779521&8.3.5.13.17.103.149.193&207&14&7519&1054500291&4.5.11.17.23.47.313&1181&2262 \cr
 & &32.27.5.7.23.103.149&76735&69552& & &16.3.5.13.17.29.1181&1181&680 \cr
\noalign{\hrule}
 & &9.25.49.19.29.173&933&2354& & &7.13.19.61.73.137&32505&1334 \cr
7502&1050936075&4.27.25.11.107.311&493&182&7520&1054795469&4.3.5.11.23.29.197&61&84 \cr
 & &16.7.11.13.17.29.107&2431&856& & &32.9.7.11.61.197&197&1584 \cr
\noalign{\hrule}
 & &9.5.49.29.41.401&5357&1748& & &9.5.7.13.529.487&297&232 \cr
7503&1051319745&8.11.19.23.41.487&401&378&7521&1054966185&16.243.7.11.29.487&3041&368 \cr
 & &32.27.7.11.401.487&487&528& & &512.23.29.3041&3041&7424 \cr
\noalign{\hrule}
 & &29.31.37.101.313&1419&2318& & &9.5.23.31.131.251&6149&6034 \cr
7504&1051543219&4.3.11.19.43.61.313&1313&2130&7522&1054986885&4.3.7.11.13.43.251.431&6157&2894 \cr
 & &16.9.5.13.61.71.101&3965&5112& & &16.11.43.47.131.1447&15917&16168 \cr
\noalign{\hrule}
 & &9.25.49.11.13.23.29&67&158& & &9.5.7.11.43.73.97&8507&1426 \cr
7505&1051575525&4.7.11.23.29.67.79&2025&208&7523&1055033595&4.3.5.23.31.47.181&817&1898 \cr
 & &128.81.25.13.67&603&64& & &16.13.19.31.43.73&589&104 \cr
\noalign{\hrule}
 & &27.25.13.37.41.79&129&2794& & &11.13.29.43.61.97&95&2718 \cr
7506&1051622325&4.81.5.11.43.127&553&338&7524&1055125357&4.9.5.11.13.19.151&1649&1220 \cr
 & &16.7.169.79.127&91&1016& & &32.3.25.17.61.97&425&48 \cr
\noalign{\hrule}
}%
}
$$
\eject
\vglue -23 pt
\noindent\hskip 1 in\hbox to 6.5 in{\ 7525 -- 7560 \hfill\fbd 1055728267 -- 1069693625\frb}
\vskip -9 pt
$$
\vbox{
\nointerlineskip
\halign{\strut
    \vrule \ \ \hfil \frb #\ 
   &\vrule \hfil \ \ \fbb #\frb\ 
   &\vrule \hfil \ \ \frb #\ \hfil
   &\vrule \hfil \ \ \frb #\ 
   &\vrule \hfil \ \ \frb #\ \ \vrule \hskip 2 pt
   &\vrule \ \ \hfil \frb #\ 
   &\vrule \hfil \ \ \fbb #\frb\ 
   &\vrule \hfil \ \ \frb #\ \hfil
   &\vrule \hfil \ \ \frb #\ 
   &\vrule \hfil \ \ \frb #\ \vrule \cr%
\noalign{\hrule}
 & &121.23.29.103.127&3233&450& & &3.5.29.37.251.263&23349&23086 \cr
7525&1055728267&4.9.25.53.61.103&1771&1462&7543&1062479235&4.9.7.17.29.43.97.181&2761&52 \cr
 & &16.3.25.7.11.17.23.43&1275&2408& & &32.11.13.17.181.251&3077&2288 \cr
\noalign{\hrule}
 & &27.25.17.19.29.167&7027&7502& & &5.11.29.61.67.163&879&716 \cr
7526&1055895075&4.9.121.17.31.7027&4355&2672&7544&1062558695&8.3.61.67.179.293&15275&4356 \cr
 & &128.5.11.13.31.67.167&4433&4288& & &64.27.25.121.13.47&6345&4576 \cr
\noalign{\hrule}
 & &9.7.13.29.79.563&145&418& & &9.11.13.29.71.401&2525&1886 \cr
7527&1056373227&4.3.5.11.19.841.79&9517&6994&7545&1062623133&4.25.13.23.29.41.101&67&600 \cr
 & &16.5.13.31.269.307&8339&12280& & &64.3.625.67.101&6767&20000 \cr
\noalign{\hrule}
 & &27.5.11.13.19.43.67&107&94& & &3.17.97.359.599&445&154 \cr
7528&1056736395&4.9.5.11.19.43.47.107&6901&2086&7546&1063807827&4.5.7.11.17.89.359&577&936 \cr
 & &16.7.47.67.103.149&7003&5768& & &64.9.5.7.11.13.577&22503&12320 \cr
\noalign{\hrule}
 & &27.5.7.83.97.139&187&104& & &3.121.361.23.353&6093&2210 \cr
7529&1057539105&16.9.5.7.11.13.17.139&779&194&7547&1063938117&4.27.5.11.13.17.677&6707&740 \cr
 & &64.11.17.19.41.97&3553&1312& & &32.25.19.37.353&37&400 \cr
\noalign{\hrule}
 & &5.11.23.29.127.227&2709&212& & &29.97.149.2539&891&3430 \cr
7530&1057591865&8.9.5.7.29.43.53&67&368&7548&1064188843&4.81.5.343.11.97&2183&2668 \cr
 & &256.3.23.53.67&201&6784& & &32.9.7.23.29.37.59&5957&8496 \cr
\noalign{\hrule}
 & &9.11.13.23.31.1153&3071&32672& & &81.23.37.15439&913&950 \cr
7531&1058028543&64.37.83.1021&1025&2046&7549&1064225709&4.25.11.19.83.15439&11993&27432 \cr
 & &256.3.25.11.31.41&1025&128& & &64.27.11.67.127.179&22733&23584 \cr
\noalign{\hrule}
 & &27.25.49.13.23.107&385&1006& & &121.41.397.541&265&276 \cr
7532&1058168475&4.125.343.11.503&21689&21186&7550&1065508697&8.3.5.11.23.41.53.397&4541&174 \cr
 & &16.9.121.529.41.107&943&968& & &32.9.19.29.53.239&40869&24592 \cr
\noalign{\hrule}
 & &3.5.11.13.19.83.313&4087&644& & &9.5.11.169.47.271&1961&2104 \cr
7533&1058774145&8.5.7.13.23.61.67&513&982&7551&1065513735&16.3.13.37.47.53.263&1595&30788 \cr
 & &32.27.19.61.491&4419&976& & &128.5.11.29.43.179&5191&2752 \cr
\noalign{\hrule}
 & &3.11.19.37.109.419&277&980& & &81.83.257.617&65&682 \cr
7534&1059521529&8.5.49.11.109.277&5447&4248&7552&1066059387&4.9.5.11.13.31.257&1079&1234 \cr
 & &128.9.7.13.59.419&1239&832& & &16.11.169.83.617&169&88 \cr
\noalign{\hrule}
 & &9.5.7.11.59.71.73&1891&1394& & &7.121.19.23.43.67&1325&1458 \cr
7535&1059586605&4.11.17.31.41.59.61&20391&20732&7553&1066370459&4.729.25.43.53.67&3653&14102 \cr
 & &32.3.7.61.71.73.971&971&976& & &16.3.5.11.13.281.641&18265&15384 \cr
\noalign{\hrule}
 & &25.49.11.13.23.263&34669&47556& & &3.7.19.529.31.163&1&22 \cr
7536&1059633575&8.9.37.937.1321&1513&2450&7554&1066541763&4.11.19.23.31.163&2169&1580 \cr
 & &32.3.25.49.17.37.89&1513&1776& & &32.9.5.11.79.241&19039&2640 \cr
\noalign{\hrule}
 & &5.289.31.59.401&4279&4680& & &9.11.13.37.43.521&12137&5364 \cr
7537&1059804905&16.9.25.11.13.59.389&527&15698&7555&1066808457&8.81.53.149.229&6149&5920 \cr
 & &64.3.17.31.47.167&501&1504& & &512.5.11.13.37.43.53&265&256 \cr
\noalign{\hrule}
 & &5.49.11.19.127.163&32591&1476& & &9.5.11.19.53.2141&149&434 \cr
7538&1059994705&8.9.13.23.41.109&987&1520&7556&1067213565&4.3.7.31.149.2141&847&1294 \cr
 & &256.27.5.7.19.47&1269&128& & &16.49.121.31.647&7117&12152 \cr
\noalign{\hrule}
 & &13.19.37.311.373&3303&3784& & &3.25.17.29.67.431&749&6578 \cr
7539&1060151417&16.9.11.43.311.367&703&230&7557&1067727075&4.25.7.11.13.23.107&149&126 \cr
 & &64.3.5.19.23.37.367&1835&2208& & &16.9.49.13.107.149&21903&11128 \cr
\noalign{\hrule}
 & &25.7.11.31.109.163&333&442& & &25.11.169.83.277&1101&1946 \cr
7540&1060245725&4.9.7.11.13.17.37.163&43&120&7558&1068506725&4.3.5.7.83.139.367&277&138 \cr
 & &64.27.5.13.17.37.43&20683&14688& & &16.9.7.23.277.367&2569&1656 \cr
\noalign{\hrule}
 & &5.7.167.419.433&2563&468& & &3.7.11.73.241.263&95&168 \cr
7541&1060440815&8.9.11.13.167.233&433&266&7559&1068827529&16.9.5.49.11.19.241&263&2432 \cr
 & &32.3.7.11.13.19.433&627&208& & &4096.361.263&361&2048 \cr
\noalign{\hrule}
 & &17.41.467.3259&5599&2340& & &125.7.11.13.83.103&179&900 \cr
7542&1060801241&8.9.5.11.13.41.509&497&1030&7560&1069693625&8.9.3125.11.179&17277&17098 \cr
 & &32.3.25.7.11.71.103&37275&18128& & &32.27.13.83.103.443&443&432 \cr
\noalign{\hrule}
}%
}
$$
\eject
\vglue -23 pt
\noindent\hskip 1 in\hbox to 6.5 in{\ 7561 -- 7596 \hfill\fbd 1070033735 -- 1084072759\frb}
\vskip -9 pt
$$
\vbox{
\nointerlineskip
\halign{\strut
    \vrule \ \ \hfil \frb #\ 
   &\vrule \hfil \ \ \fbb #\frb\ 
   &\vrule \hfil \ \ \frb #\ \hfil
   &\vrule \hfil \ \ \frb #\ 
   &\vrule \hfil \ \ \frb #\ \ \vrule \hskip 2 pt
   &\vrule \ \ \hfil \frb #\ 
   &\vrule \hfil \ \ \fbb #\frb\ 
   &\vrule \hfil \ \ \frb #\ \hfil
   &\vrule \hfil \ \ \frb #\ 
   &\vrule \hfil \ \ \frb #\ \vrule \cr%
\noalign{\hrule}
 & &5.19.841.59.227&99&128& & &25.11.169.19.23.53&8791&33966 \cr
7561&1070033735&256.9.5.11.19.29.59&333&1378&7579&1076407475&4.27.17.37.59.149&2491&3022 \cr
 & &1024.81.13.37.53&25493&41472& & &16.3.17.47.53.1511&4533&6392 \cr
\noalign{\hrule}
 & &625.7.11.13.29.59&3043&1332& & &3.11.17.29.127.521&507&14 \cr
7562&1070444375&8.9.11.13.17.37.179&9215&11542&7580&1076470923&4.9.7.11.169.127&3071&4930 \cr
 & &32.3.5.19.29.97.199&3781&4656& & &16.5.17.29.37.83&83&1480 \cr
\noalign{\hrule}
 & &243.19.47.4933&2255&2678& & &25.11.79.179.277&117&62 \cr
7563&1070456067&4.27.5.11.13.19.41.103&4933&710&7581&1077190675&4.9.5.13.31.79.277&179&1206 \cr
 & &16.25.13.71.4933&325&568& & &16.81.31.67.179&2511&536 \cr
\noalign{\hrule}
 & &3.19.23.67.73.167&725&3898& & &13.8681.9547&57715&66396 \cr
7564&1070820867&4.25.29.73.1949&1157&792&7582&1077407591&8.3.5.7.11.17.97.503&785&282 \cr
 & &64.9.5.11.13.29.89&20735&8544& & &32.9.25.7.17.47.157&51653&61200 \cr
\noalign{\hrule}
 & &81.7.11.83.2069&157&74& & &9.7.11.17.23.41.97&19015&19694 \cr
7565&1071061299&4.27.37.157.2069&1085&3154&7583&1077619851&4.5.41.43.229.3803&2793&6596 \cr
 & &16.5.7.19.31.37.83&589&1480& & &32.3.5.49.17.19.43.97&665&688 \cr
\noalign{\hrule}
 & &9.5.13.17.19.53.107&2849&2186& & &11.19.23.29.71.109&801&548 \cr
7566&1071563805&4.3.7.11.37.107.1093&1265&2014&7584&1078839817&8.9.29.89.109.137&21109&23690 \cr
 & &16.5.121.19.23.37.53&851&968& & &32.3.5.11.19.23.101.103&1545&1616 \cr
\noalign{\hrule}
 & &27.125.17.19.983&287&2662& & &125.7.17.29.41.61&71&54 \cr
7567&1071592875&4.9.7.1331.17.41&475&596&7585&1078868875&4.27.7.29.41.61.71&5665&11842 \cr
 & &32.25.11.19.41.149&1639&656& & &16.9.5.11.31.103.191&19673&24552 \cr
\noalign{\hrule}
 & &81.19.23.107.283&985&878& & &81.5.7.31.71.173&17059&11696 \cr
7568&1071856557&4.5.19.197.283.439&13223&13662&7586&1079491455&32.49.17.43.2437&165&2272 \cr
 & &16.27.7.11.23.197.1889&15169&15112& & &2048.3.5.11.17.71&187&1024 \cr
\noalign{\hrule}
 & &71.79.223.857&3233&2376& & &13.43.59.71.461&87067&93060 \cr
7569&1071941599&16.27.11.53.61.223&613&2620&7587&1079501111&8.9.5.11.47.83.1049&29323&13588 \cr
 & &128.3.5.11.131.613&21615&39232& & &64.3.7.43.59.71.79&237&224 \cr
\noalign{\hrule}
 & &7.121.41.89.347&3695&1266& & &27.5.11.23.101.313&1079&1244 \cr
7570&1072473941&4.3.5.89.211.739&9759&9020&7588&1079742015&8.9.13.83.311.313&635&3434 \cr
 & &32.9.25.11.41.3253&3253&3600& & &32.5.17.83.101.127&1411&2032 \cr
\noalign{\hrule}
 & &25.11.37.97.1087&3211&2124& & &9.5.7.169.53.383&40711&24016 \cr
7571&1072841825&8.9.5.169.19.37.59&2141&1374&7589&1080617265&32.11.19.79.3701&6405&10106 \cr
 & &32.27.13.229.2141&57807&47632& & &128.3.5.7.31.61.163&5053&3904 \cr
\noalign{\hrule}
 & &9.5.11.101.109.197&2993&3002& & &25.7.11.67.83.101&6579&188 \cr
7572&1073542635&4.19.41.73.79.101.197&7725&352&7590&1081197425&8.9.25.17.43.47&539&536 \cr
 & &256.3.25.11.19.79.103&9785&10112& & &128.3.49.11.17.47.67&987&1088 \cr
\noalign{\hrule}
 & &9.5.7.17.19.61.173&167&1448& & &27.2401.13.1283&2185&902 \cr
7573&1073714985&16.3.167.173.181&343&176&7591&1081249533&4.3.5.7.11.13.19.23.41&1283&1720 \cr
 & &512.343.11.181&1991&12544& & &64.25.41.43.1283&1025&1376 \cr
\noalign{\hrule}
 & &5.11.41.359.1327&29133&43852& & &9.25.49.17.29.199&2629&1634 \cr
7574&1074266215&8.27.13.19.83.577&2077&2654&7592&1081629675&4.3.5.11.17.19.43.239&5771&2116 \cr
 & &32.9.13.31.67.1327&2077&1872& & &32.19.529.29.199&529&304 \cr
\noalign{\hrule}
 & &3.13.23.29.109.379&3717&1210& & &25.7.19.43.67.113&25957&28782 \cr
7575&1074623043&4.27.5.7.121.29.59&1895&184&7593&1082463725&4.27.7.13.41.101.257&1661&2480 \cr
 & &64.25.11.23.379&25&352& & &128.3.5.11.31.151.257&51491&49344 \cr
\noalign{\hrule}
 & &5.7.13.53.109.409&2277&3040& & &19.29.37.173.307&3047&3354 \cr
7576&1075070815&64.9.25.11.19.23.53&197&1022&7594&1082773957&4.3.11.13.19.29.43.277&477&8510 \cr
 & &256.3.7.19.73.197&14381&7296& & &16.27.5.13.23.37.53&6095&2808 \cr
\noalign{\hrule}
 & &9.7.169.43.2351&1585&766& & &3.7.13.67.193.307&479&12452 \cr
7577&1076337171&4.5.13.43.317.383&18537&2068&7595&1083760041&8.7.11.283.479&1535&1818 \cr
 & &32.3.11.37.47.167&19129&2672& & &32.9.5.11.101.307&1111&240 \cr
\noalign{\hrule}
 & &9.5.49.37.79.167&97&88& & &7.121.19.31.41.53&169&48 \cr
7578&1076350905&16.49.11.79.97.167&7923&260&7596&1084072759&32.3.169.19.41.53&963&1210 \cr
 & &128.3.5.11.13.19.139&2717&8896& & &128.27.5.121.13.107&2889&4160 \cr
\noalign{\hrule}
}%
}
$$
\eject
\vglue -23 pt
\noindent\hskip 1 in\hbox to 6.5 in{\ 7597 -- 7632 \hfill\fbd 1084777803 -- 1094988411\frb}
\vskip -9 pt
$$
\vbox{
\nointerlineskip
\halign{\strut
    \vrule \ \ \hfil \frb #\ 
   &\vrule \hfil \ \ \fbb #\frb\ 
   &\vrule \hfil \ \ \frb #\ \hfil
   &\vrule \hfil \ \ \frb #\ 
   &\vrule \hfil \ \ \frb #\ \ \vrule \hskip 2 pt
   &\vrule \ \ \hfil \frb #\ 
   &\vrule \hfil \ \ \fbb #\frb\ 
   &\vrule \hfil \ \ \frb #\ \hfil
   &\vrule \hfil \ \ \frb #\ 
   &\vrule \hfil \ \ \frb #\ \vrule \cr%
\noalign{\hrule}
 & &9.17.1369.5179&3571&8750& & &9.7.17.19.149.359&8621&34100 \cr
7597&1084777803&4.625.7.17.3571&1573&1998&7615&1088488359&8.25.11.31.37.233&1669&894 \cr
 & &16.27.25.7.121.13.37&2275&2904& & &32.3.37.149.1669&1669&592 \cr
\noalign{\hrule}
 & &9.5.7.13.43.61.101&17&22& & &5.7.13.37.71.911&12489&646 \cr
7598&1084859685&4.3.7.11.17.43.61.101&32311&1910&7616&1088904635&4.3.7.17.19.23.181&1221&1040 \cr
 & &16.5.79.191.409&15089&3272& & &128.9.5.11.13.23.37&207&704 \cr
\noalign{\hrule}
 & &11.103.853.1123&6603&5750& & &3.25.11.169.73.107&9953&10122 \cr
7599&1085322227&4.3.125.23.31.71.103&3369&3944&7617&1089048675&4.9.7.37.107.241.269&8347&1606 \cr
 & &64.9.5.17.29.31.1123&4743&4640& & &16.11.17.73.241.491&4097&3928 \cr
\noalign{\hrule}
 & &9.5.13.19.23.31.137&2641&1606& & &81.49.11.13.19.101&33109&31228 \cr
7600&1085724315&4.11.13.361.73.139&2727&7420&7618&1089161073&8.9.37.113.211.293&30745&31078 \cr
 & &32.27.5.7.11.53.101&3333&5936& & &32.5.11.13.41.43.113.379&77695&77744 \cr
\noalign{\hrule}
 & &3.19.41.73.6367&232375&232416& & &9.5.13.19.29.31.109&581&550 \cr
7601&1086216567&64.81.125.11.169.19.269&38471&4&7619&1089169965&4.3.125.7.11.19.83.109&961&8086 \cr
 & &512.5.13.17.31.73&2015&4352& & &16.7.11.13.961.311&2177&2728 \cr
\noalign{\hrule}
 & &5.13.31.73.83.89&1197&1562& & &3.7.11.83.113.503&15&98 \cr
7602&1086590765&4.9.7.11.13.19.71.83&1333&410&7620&1089774147&4.9.5.343.11.503&377&4150 \cr
 & &16.3.5.11.19.31.41.43&2451&3608& & &16.125.13.29.83&3625&104 \cr
\noalign{\hrule}
 & &9.7.29.577.1031&136565&132526& & &9.25.49.121.19.43&1819&8716 \cr
7603&1086858549&4.5.11.13.23.43.67.191&4039&354&7621&1089898425&8.3.5.17.107.2179&3311&3226 \cr
 & &16.3.7.13.43.59.577&559&472& & &32.7.11.43.107.1613&1613&1712 \cr
\noalign{\hrule}
 & &9.5.7.11.529.593&3067&2752& & &5.53.83.179.277&35789&13794 \cr
7604&1086960105&128.43.593.3067&14283&11216&7622&1090578085&4.3.121.13.19.2753&2735&18 \cr
 & &4096.27.529.701&2103&2048& & &16.27.5.11.547&14769&88 \cr
\noalign{\hrule}
 & &3.5.11.13.19.149.179&3077&2830& & &27.25.11.13.89.127&709&434 \cr
7605&1086976605&4.25.17.149.181.283&3401&324&7623&1091022075&4.3.7.13.31.89.709&15205&12446 \cr
 & &32.81.19.179.283&283&432& & &16.5.343.127.3041&3041&2744 \cr
\noalign{\hrule}
 & &27.7.11.13.19.29.73&2881&2078& & &27.125.7.11.13.17.19&6847&1798 \cr
7606&1087107021&4.3.7.13.43.67.1039&3515&398&7624&1091215125&4.25.29.31.41.167&17451&12274 \cr
 & &16.5.19.37.67.199&2479&7960& & &16.9.7.17.361.277&277&152 \cr
\noalign{\hrule}
 & &27.5.49.13.47.269&61&208& & &3.5.7.11.17.23.41.59&3163&908 \cr
7607&1087234785&32.9.5.169.47.61&2369&5236&7625&1092432495&8.7.17.227.3163&11925&15088 \cr
 & &256.7.11.17.23.103&4301&13184& & &256.9.25.23.41.53&265&384 \cr
\noalign{\hrule}
 & &5.13.17.29.83.409&2997&3412& & &5.11.17.23.89.571&4221&5486 \cr
7608&1087831615&8.81.37.409.853&187&1040&7626&1092862595&4.9.7.13.67.89.211&40975&36544 \cr
 & &256.27.5.11.13.17.37&999&1408& & &512.3.25.11.149.571&745&768 \cr
\noalign{\hrule}
 & &3.7.11.71.113.587&3&74& & &5.121.53.103.331&2379&4034 \cr
7609&1087894731&4.9.37.113.587&2585&2698&7627&1093192045&4.3.13.61.103.2017&957&1060 \cr
 & &16.5.11.19.37.47.71&893&1480& & &32.9.5.11.13.29.53.61&1769&1872 \cr
\noalign{\hrule}
 & &3.5.13.29.199.967&2393&508& & &3.5.7.11.361.43.61&811&994 \cr
7610&1088208615&8.127.199.2393&13833&11440&7628&1093672965&4.49.11.43.71.811&36661&3078 \cr
 & &256.9.5.11.13.29.53&583&384& & &16.81.19.61.601&601&216 \cr
\noalign{\hrule}
 & &3.5.49.89.127.131&247&198& & &9.7.41.431.983&5395&12276 \cr
7611&1088309355&4.27.11.13.19.127.131&2435&994&7629&1094347359&8.81.5.11.13.31.83&983&70 \cr
 & &16.5.7.13.19.71.487&9253&7384& & &32.25.7.31.983&31&400 \cr
\noalign{\hrule}
 & &27.11.23.107.1489&8201&8178& & &5.49.13.43.61.131&33511&5346 \cr
7612&1088335413&4.81.29.47.59.107.139&1253&5060&7630&1094407405&4.243.11.23.31.47&6623&7336 \cr
 & &32.5.7.11.23.29.139.179&20155&20048& & &64.9.7.37.131.179&1611&1184 \cr
\noalign{\hrule}
 & &27.25.11.47.3119&1147&1972& & &25.7.121.17.3041&25425&26272 \cr
7613&1088453025&8.9.17.29.31.37.47&143&190&7631&1094683975&64.9.625.113.821&723&98 \cr
 & &32.5.11.13.17.19.29.31&7163&8432& & &256.27.49.113.241&27233&24192 \cr
\noalign{\hrule}
 & &49.13.41.71.587&3959&3672& & &9.1331.17.19.283&12883&35510 \cr
7614&1088478209&16.27.7.17.37.71.107&4187&3410&7632&1094988411&4.5.13.53.67.991&363&628 \cr
 & &64.9.5.11.17.31.53.79&188415&185504& & &32.3.121.13.67.157&2041&1072 \cr
\noalign{\hrule}
}%
}
$$
\eject
\vglue -23 pt
\noindent\hskip 1 in\hbox to 6.5 in{\ 7633 -- 7668 \hfill\fbd 1095007077 -- 1108778055\frb}
\vskip -9 pt
$$
\vbox{
\nointerlineskip
\halign{\strut
    \vrule \ \ \hfil \frb #\ 
   &\vrule \hfil \ \ \fbb #\frb\ 
   &\vrule \hfil \ \ \frb #\ \hfil
   &\vrule \hfil \ \ \frb #\ 
   &\vrule \hfil \ \ \frb #\ \ \vrule \hskip 2 pt
   &\vrule \ \ \hfil \frb #\ 
   &\vrule \hfil \ \ \fbb #\frb\ 
   &\vrule \hfil \ \ \frb #\ \hfil
   &\vrule \hfil \ \ \frb #\ 
   &\vrule \hfil \ \ \frb #\ \vrule \cr%
\noalign{\hrule}
 & &9.17.107.211.317&319&2& & &79.109.227.563&21717&39650 \cr
7633&1095007077&4.3.11.17.29.211&845&634&7651&1100494411&4.9.25.13.19.61.127&1199&1180 \cr
 & &16.5.11.169.317&1859&40& & &32.3.125.11.59.109.127&22125&22352 \cr
\noalign{\hrule}
 & &3.5.13.17.47.79.89&163&74& & &7.11.13.37.113.263&95&18 \cr
7634&1095464955&4.5.13.17.37.47.163&803&1602&7652&1100702603&4.9.5.13.19.37.263&635&154 \cr
 & &16.9.11.73.89.163&2409&1304& & &16.3.25.7.11.19.127&381&3800 \cr
\noalign{\hrule}
 & &3.19.23.41.89.229&605&338& & &9.5.13.59.167.191&2413&242 \cr
7635&1095499131&4.5.121.169.19.229&3105&1246&7653&1100924955&4.121.19.127.191&1111&1302 \cr
 & &16.27.25.7.11.23.89&275&504& & &16.3.7.1331.31.101&21917&10648 \cr
\noalign{\hrule}
 & &3.11.13.19.43.53.59&335&1342& & &3.13.1481.19081&721&760 \cr
7636&1095991611&4.5.121.59.61.67&2279&1674&7654&1102099479&16.5.7.19.103.19081&2691&16390 \cr
 & &16.27.31.43.53.61&279&488& & &64.9.25.11.13.23.149&10281&8800 \cr
\noalign{\hrule}
 & &169.23.29.71.137&875&4026& & &9.7.11.13.79.1549&629&398 \cr
7637&1096456621&4.3.125.7.11.61.71&221&276&7655&1102440339&4.3.17.37.199.1549&4015&632 \cr
 & &32.9.25.13.17.23.61&1525&2448& & &64.5.11.37.73.79&365&1184 \cr
\noalign{\hrule}
 & &27.37.53.139.149&2783&2730& & &3.25.11.29.31.1487&18417&18758 \cr
7638&1096585317&4.81.5.7.121.13.23.139&167&1696&7656&1102870725&4.9.7.29.83.113.877&569&308 \cr
 & &256.5.7.11.13.53.167&10855&9856& & &32.49.11.83.113.569&64297&65072 \cr
\noalign{\hrule}
 & &9.25.11.169.43.61&4447&1702& & &3.11.17.683.2879&1781&1098 \cr
7639&1097135325&4.5.13.23.37.4447&393&458&7657&1103126277&4.27.11.13.17.61.137&3415&1634 \cr
 & &16.3.131.229.4447&29999&35576& & &16.5.19.43.61.683&1159&1720 \cr
\noalign{\hrule}
 & &243.11.47.8737&203983&206656& & &5.11.41.43.59.193&13401&13208 \cr
7640&1097638047&128.169.17.71.3229&4385&7614&7658&1104140455&16.9.5.13.43.127.1489&4469&11914 \cr
 & &512.81.5.17.47.877&4385&4352& & &64.3.7.13.23.37.41.109&28231&28704 \cr
\noalign{\hrule}
 & &3.5.11.13.19.29.929&2537&2108& & &11.59.67.109.233&441&208 \cr
7641&1097980455&8.17.19.29.31.43.59&27207&1880&7659&1104337751&32.9.49.13.67.109&197&130 \cr
 & &128.9.5.47.3023&9069&3008& & &128.3.5.49.169.197&33293&47040 \cr
\noalign{\hrule}
 & &3.5.7.13.23.79.443&279&164& & &5.23.37.53.59.83&1701&482 \cr
7642&1098730815&8.27.7.13.31.41.79&3443&2416&7660&1104346955&4.243.5.7.83.241&7579&856 \cr
 & &256.11.41.151.313&47263&57728& & &64.3.11.13.53.107&321&4576 \cr
\noalign{\hrule}
 & &3.7.121.13.29.31.37&1037&36& & &3.5.11.31.149.1451&24679&26130 \cr
7643&1098776679&8.27.11.17.31.61&2075&3722&7661&1105857885&4.9.25.13.23.29.37.67&121&454 \cr
 & &32.25.83.1861&1861&33200& & &16.121.13.29.67.227&21373&23608 \cr
\noalign{\hrule}
 & &3.25.121.289.419&4747&4328& & &3.5.7.11.17.23.31.79&1879&570 \cr
7644&1098900825&16.289.47.101.541&1881&27308&7662&1105980645&4.9.25.19.23.1879&341&134 \cr
 & &128.9.11.19.6827&6827&3648& & &16.11.31.67.1879&1879&536 \cr
\noalign{\hrule}
 & &5.7.29.67.103.157&611&1332& & &9.17.71.223.457&775&1232 \cr
7645&1099708855&8.9.5.13.37.47.157&2047&308&7663&1107059193&32.25.7.11.17.31.71&1755&446 \cr
 & &64.3.7.11.13.23.89&12727&2208& & &128.27.125.13.223&375&832 \cr
\noalign{\hrule}
 & &3.25.7.11.1849.103&837&1012& & &5.7.13.19.211.607&321057&319328 \cr
7646&1099831425&8.81.121.23.31.103&7163&5300&7664&1107225665&64.27.11.17.23.47.587&5275&3514 \cr
 & &64.25.13.19.29.31.53&19981&18848& & &256.9.25.7.23.211.251&5773&5760 \cr
\noalign{\hrule}
 & &9.5.11.19.29.37.109&5747&5044& & &9.25.49.11.13.19.37&289&44 \cr
7647&1099980585&8.5.7.13.29.97.821&27&118&7665&1108332225&8.5.121.13.289.19&597&1702 \cr
 & &32.27.59.97.821&17169&13136& & &32.3.17.23.37.199&391&3184 \cr
\noalign{\hrule}
 & &3.7.11.13.17.29.743&703&40& & &9.25.169.103.283&1529&5546 \cr
7648&1099995897&16.5.7.11.19.29.37&83&468&7666&1108390725&4.3.11.13.47.59.139&159661&160178 \cr
 & &128.9.13.37.83&111&5312& & &16.67.80089.2383&18961&19064 \cr
\noalign{\hrule}
 & &9.5.7.13.29.59.157&267383&265318& & &3.25.7.29.47.1549&3379&4366 \cr
7649&1100027565&4.47.53.2503.5689&2497&6&7667&1108425675&4.5.29.31.37.59.109&5463&968 \cr
 & &16.3.11.227.5689&5689&19976& & &64.9.121.37.607&22459&11616 \cr
\noalign{\hrule}
 & &11.43.73.151.211&33195&1334& & &3.5.7.121.197.443&109&88 \cr
7650&1100128469&4.3.5.23.29.2213&1095&1118&7668&1108778055&16.5.1331.109.443&1773&442 \cr
 & &16.9.25.13.29.43.73&725&936& & &64.9.13.17.109.197&1853&1248 \cr
\noalign{\hrule}
}%
}
$$
\eject
\vglue -23 pt
\noindent\hskip 1 in\hbox to 6.5 in{\ 7669 -- 7704 \hfill\fbd 1108976211 -- 1120969239\frb}
\vskip -9 pt
$$
\vbox{
\nointerlineskip
\halign{\strut
    \vrule \ \ \hfil \frb #\ 
   &\vrule \hfil \ \ \fbb #\frb\ 
   &\vrule \hfil \ \ \frb #\ \hfil
   &\vrule \hfil \ \ \frb #\ 
   &\vrule \hfil \ \ \frb #\ \ \vrule \hskip 2 pt
   &\vrule \ \ \hfil \frb #\ 
   &\vrule \hfil \ \ \fbb #\frb\ 
   &\vrule \hfil \ \ \frb #\ \hfil
   &\vrule \hfil \ \ \frb #\ 
   &\vrule \hfil \ \ \frb #\ \vrule \cr%
\noalign{\hrule}
 & &27.7.19.23.29.463&187&250& & &3.5.11.13.19.139.197&1803&758 \cr
7669&1108976211&4.3.125.11.17.29.463&247&710&7687&1115994165&4.9.139.379.601&815&436 \cr
 & &16.625.13.17.19.71&10625&7384& & &32.5.109.163.601&17767&9616 \cr
\noalign{\hrule}
 & &7.11.97.163.911&795&116& & &5.169.19.157.443&14319&14476 \cr
7670&1109094217&8.3.5.11.29.53.163&911&882&7688&1116641305&8.9.7.11.13.19.37.43.47&53&194 \cr
 & &32.27.5.49.53.911&945&848& & &32.3.7.11.37.43.53.97&175483&172272 \cr
\noalign{\hrule}
 & &5.7.11.101.103.277&3827&4104& & &11.23.29.41.47.79&5525&1812 \cr
7671&1109427935&16.27.5.19.43.89.101&359&4186&7689&1116933521&8.3.25.13.17.41.151&29&726 \cr
 & &64.3.7.13.19.23.359&17043&11488& & &32.9.5.121.13.29&55&1872 \cr
\noalign{\hrule}
 & &11.23.29.37.61.67&1917&2170& & &9.13.23.29.103.139&8041&4010 \cr
7672&1109493803&4.27.5.7.29.31.37.71&65&268&7690&1117284363&4.5.11.17.23.43.401&327&74 \cr
 & &32.3.25.13.31.67.71&6603&5200& & &16.3.5.17.37.43.109&20165&5848 \cr
\noalign{\hrule}
 & &9.5.127.163.1193&533&660& & &7.23.31.41.43.127&773&990 \cr
7673&1111333185&8.27.25.11.13.41.163&1397&722&7691&1117489891&4.9.5.11.23.127.773&301&2620 \cr
 & &32.121.361.41.127&4961&5776& & &32.3.25.7.11.43.131&825&2096 \cr
\noalign{\hrule}
 & &7.169.29.179.181&1265&1062& & &9.5.19.43.113.269&49&220 \cr
7674&1111512493&4.9.5.11.13.23.59.181&539&358&7692&1117545705&8.25.49.11.43.113&243&3068 \cr
 & &16.3.5.49.121.59.179&2541&2360& & &64.243.7.13.59&11151&416 \cr
\noalign{\hrule}
 & &9.5.11.289.19.409&33973&6518& & &3.17.103.137.1553&2717&4270 \cr
7675&1111680405&4.53.641.3259&1309&1950&7693&1117633533&4.5.7.11.13.19.61.103&11061&1688 \cr
 & &16.3.25.7.11.13.17.53&689&280& & &64.9.5.211.1229&6145&20256 \cr
\noalign{\hrule}
 & &13.163.409.1283&1701&418& & &5.7.11.47.113.547&3177&442 \cr
7676&1111938893&4.243.7.11.19.409&1283&1580&7694&1118470045&4.9.13.17.113.353&59&280 \cr
 & &32.9.5.19.79.1283&855&1264& & &64.3.5.7.59.353&353&5664 \cr
\noalign{\hrule}
 & &81.7.41.151.317&6215&6782& & &3.17.29.61.79.157&845&924 \cr
7677&1112764149&4.5.11.113.151.3391&17009&54&7695&1118986257&8.9.5.7.11.169.17.157&1403&638 \cr
 & &16.27.11.73.233&233&6424& & &32.7.121.13.23.29.61&2093&1936 \cr
\noalign{\hrule}
 & &5.13.17.491.2053&14931&16984& & &139.29929.269&3731&33660 \cr
7678&1113865415&16.27.7.11.17.79.193&775&962&7696&1119075239&8.9.5.7.11.13.17.41&173&278 \cr
 & &64.3.25.7.13.31.37.79&17205&17696& & &32.3.13.17.139.173&221&48 \cr
\noalign{\hrule}
 & &5.7.241.269.491&1089&598& & &3.13.29.37.47.569&1159&548 \cr
7679&1114086365&4.9.5.121.13.23.269&491&854&7697&1119114321&8.19.29.37.61.137&605&468 \cr
 & &16.3.7.13.23.61.491&793&552& & &64.9.5.121.13.19.61&7381&9120 \cr
\noalign{\hrule}
 & &9.25.67.193.383&407&742& & &23.31.41.149.257&58387&47850 \cr
7680&1114328925&4.3.5.7.11.37.53.193&871&94&7698&1119419269&4.3.25.7.11.19.29.439&149&54 \cr
 & &16.11.13.47.53.67&689&4136& & &16.81.5.11.149.439&4455&3512 \cr
\noalign{\hrule}
 & &5.7.121.841.313&1007&3198& & &81.5.13.19.23.487&2149&286 \cr
7681&1114791755&4.3.121.13.19.41.53&2817&3596&7699&1120492035&4.7.11.169.19.307&487&696 \cr
 & &32.27.13.29.31.313&351&496& & &64.3.29.307.487&307&928 \cr
\noalign{\hrule}
 & &3.5.13.47.103.1181&61127&60516& & &243.13.43.73.113&34133&6674 \cr
7682&1114858095&8.27.5.11.1681.5557&11&5546&7700&1120519413&4.11.29.47.71.107&11295&3698 \cr
 & &32.121.41.47.59&4961&944& & &16.9.5.1849.251&251&1720 \cr
\noalign{\hrule}
 & &3.11.13.739.3517&3045&6562& & &27.49.17.19.43.61&143&1180 \cr
7683&1114998027&4.9.5.7.11.17.29.193&1973&1478&7701&1120883967&8.5.11.13.19.43.59&381&854 \cr
 & &16.193.739.1973&1973&1544& & &32.3.7.59.61.127&127&944 \cr
\noalign{\hrule}
 & &9.5.343.29.47.53&3587&3322& & &9.13.31.199.1553&35713&12430 \cr
7684&1115008965&4.3.7.11.17.29.151.211&4405&26&7702&1120913469&4.5.11.71.113.503&26467&30372 \cr
 & &16.5.11.13.17.881&14977&1144& & &32.3.7.19.199.2531&2531&2128 \cr
\noalign{\hrule}
 & &7.11.13.23.47.1031&12717&35740& & &9.7.11.169.17.563&375&188 \cr
7685&1115625511&8.81.5.157.1787&1129&658&7703&1120926807&8.27.125.7.169.47&2329&6554 \cr
 & &32.27.5.7.47.1129&1129&2160& & &32.5.17.29.113.137&3973&9040 \cr
\noalign{\hrule}
 & &243.125.23.1597&33553&3178& & &3.13.17.19.23.53.73&9211&8820 \cr
7686&1115704125&4.7.13.29.89.227&1275&1364&7704&1120969239&8.27.5.49.53.61.151&403&9614 \cr
 & &32.3.25.11.17.31.227&3859&5456& & &32.5.7.11.13.19.23.31&341&560 \cr
\noalign{\hrule}
}%
}
$$
\eject
\vglue -23 pt
\noindent\hskip 1 in\hbox to 6.5 in{\ 7705 -- 7740 \hfill\fbd 1121169555 -- 1136657349\frb}
\vskip -9 pt
$$
\vbox{
\nointerlineskip
\halign{\strut
    \vrule \ \ \hfil \frb #\ 
   &\vrule \hfil \ \ \fbb #\frb\ 
   &\vrule \hfil \ \ \frb #\ \hfil
   &\vrule \hfil \ \ \frb #\ 
   &\vrule \hfil \ \ \frb #\ \ \vrule \hskip 2 pt
   &\vrule \ \ \hfil \frb #\ 
   &\vrule \hfil \ \ \fbb #\frb\ 
   &\vrule \hfil \ \ \frb #\ \hfil
   &\vrule \hfil \ \ \frb #\ 
   &\vrule \hfil \ \ \frb #\ \vrule \cr%
\noalign{\hrule}
 & &9.5.11.1499.1511&34423&48022& & &3.121.13.23.101.103&7561&2842 \cr
7705&1121169555&4.13.29.1187.1847&1517&330&7723&1129110411&4.49.23.29.7561&3861&3700 \cr
 & &16.3.5.11.13.29.37.41&1517&3016& & &32.27.25.7.11.13.29.37&5075&5328 \cr
\noalign{\hrule}
 & &27.49.13.197.331&11951&4268& & &3.25.19.29.151.181&2969&7348 \cr
7706&1121495193&8.9.11.17.19.37.97&4303&3430&7724&1129453575&8.25.11.167.2969&603&3572 \cr
 & &32.5.343.13.17.331&85&112& & &64.9.11.19.47.67&737&4512 \cr
\noalign{\hrule}
 & &9.5.7.19.23.29.281&793&242& & &27.5.29.59.67.73&3773&1118 \cr
7707&1121750595&4.7.121.13.61.281&4321&1230&7725&1129747635&4.3.343.11.13.29.43&25&584 \cr
 & &16.3.5.11.29.41.149&1639&328& & &64.25.49.11.73&49&1760 \cr
\noalign{\hrule}
 & &5.49.13.31.83.137&50067&37202& & &27.25.11.13.23.509&185&68 \cr
7708&1122715685&4.9.11.19.89.5563&143&124&7726&1130018175&8.3.125.17.37.509&299&1826 \cr
 & &32.3.121.13.31.5563&5563&5808& & &32.11.13.23.37.83&83&592 \cr
\noalign{\hrule}
 & &3.23.71.419.547&6919&5662& & &9.25.121.13.31.103&1019&4646 \cr
7709&1122816507&4.11.17.19.37.71.149&145&1494&7727&1130082525&4.5.11.23.101.1019&3103&1992 \cr
 & &16.9.5.17.29.37.83&21165&8584& & &64.3.23.29.83.107&8881&21344 \cr
\noalign{\hrule}
 & &9.25.13.17.19.29.41&97&682& & &3.7.11.13.19.29.683&127&556 \cr
7710&1123337475&4.5.11.17.29.31.97&407&492&7728&1130127999&8.7.19.29.127.139&169&720 \cr
 & &32.3.121.37.41.97&3589&1936& & &256.9.5.169.139&5421&640 \cr
\noalign{\hrule}
 & &9.7.13.19.23.43.73&1411&7130& & &5.17.19.43.73.223&69&154 \cr
7711&1123457517&4.5.17.529.31.83&1551&1022&7729&1130495155&4.3.7.11.19.23.43.73&1383&10370 \cr
 & &16.3.5.7.11.17.47.73&935&376& & &16.9.5.17.61.461&549&3688 \cr
\noalign{\hrule}
 & &3.13.29.37.107.251&2023&1936& & &9.13.17.19.23.1301&3491&2190 \cr
7712&1123884879&32.7.121.13.289.251&1665&92&7730&1130820093&4.27.5.17.73.3491&2893&598 \cr
 & &256.9.5.289.23.37&4335&2944& & &16.11.13.23.73.263&803&2104 \cr
\noalign{\hrule}
 & &9.11.19.31.37.521&271&62& & &27.23.89.97.211&715&1516 \cr
7713&1124061147&4.961.271.521&741&220&7731&1131190623&8.3.5.11.13.211.379&803&1940 \cr
 & &32.3.5.11.13.19.271&1355&208& & &64.25.121.73.97&1825&3872 \cr
\noalign{\hrule}
 & &3.7.29.37.139.359&1679&1540& & &343.11.13.17.23.59&241&150 \cr
7714&1124419233&8.5.49.11.23.73.359&15921&10286&7732&1131511381&4.3.25.49.11.59.241&8073&5182 \cr
 & &32.9.11.29.37.61.139&183&176& & &16.81.5.13.23.2591&2591&3240 \cr
\noalign{\hrule}
 & &9.7.13.23.841.71&539&1462& & &27.5.41.73.2803&11687&3278 \cr
7715&1124776107&4.3.343.11.17.29.43&16259&4940&7733&1132566165&4.9.11.13.29.31.149&205&56 \cr
 & &32.5.13.19.71.229&1145&304& & &64.5.7.11.13.31.41&2387&416 \cr
\noalign{\hrule}
 & &9.5.13.53.131.277&1591&902& & &27.343.19.47.137&65801&56540 \cr
7716&1125078435&4.5.11.37.41.43.131&1209&554&7734&1133000001&8.5.11.29.257.2269&279&2548 \cr
 & &16.3.11.13.31.37.277&407&248& & &64.9.5.49.13.29.31&1885&992 \cr
\noalign{\hrule}
 & &27.25.121.37.373&241&166& & &3.49.11.13.31.37.47&551&1188 \cr
7717&1127196675&4.9.11.83.241.373&2135&1222&7735&1133221089&8.81.121.19.29.31&2405&106 \cr
 & &16.5.7.13.47.61.241&20069&25064& & &32.5.13.29.37.53&1537&80 \cr
\noalign{\hrule}
 & &11.19.139.38809&4879&33930& & &3.11.41.853.983&553&430 \cr
7718&1127440259&4.9.5.7.13.17.29.41&197&418&7736&1134489147&4.5.7.11.43.79.853&2993&6390 \cr
 & &16.3.7.11.19.29.197&29&168& & &16.9.25.7.41.71.73&5183&4200 \cr
\noalign{\hrule}
 & &25.11.23.173.1031&117&1148& & &25.29.41.73.523&1287&1328 \cr
7719&1128145975&8.9.5.7.13.41.173&1031&1564&7737&1134870775&32.9.5.11.13.29.73.83&371&6 \cr
 & &64.3.7.17.23.1031&51&224& & &128.27.7.11.53.83&15741&37184 \cr
\noalign{\hrule}
 & &81.7.11.113.1601&66215&57062& & &17.23.71.163.251&40777&22956 \cr
7720&1128354381&4.5.17.19.41.103.277&4803&6554&7738&1135785793&8.3.121.337.1913&897&2810 \cr
 & &16.3.5.19.29.113.1601&145&152& & &32.9.5.11.13.23.281&3653&7920 \cr
\noalign{\hrule}
 & &5.7.121.41.67.97&999&12736& & &9.5.121.23.43.211&9881&4678 \cr
7721&1128453865&128.27.7.37.199&1265&1066&7739&1136257155&4.3.5.41.241.2339&1477&862 \cr
 & &512.3.5.11.13.23.41&299&768& & &16.7.211.241.431&3017&1928 \cr
\noalign{\hrule}
 & &9.625.11.71.257&703&78& & &9.17.19.29.97.139&2695&2834 \cr
7722&1129033125&4.27.13.19.37.257&115&142&7740&1136657349&4.3.5.49.11.13.17.29.109&7201&194 \cr
 & &16.5.13.19.23.37.71&851&1976& & &16.19.97.109.379&379&872 \cr
\noalign{\hrule}
}%
}
$$
\eject
\vglue -23 pt
\noindent\hskip 1 in\hbox to 6.5 in{\ 7741 -- 7776 \hfill\fbd 1136691627 -- 1148201901\frb}
\vskip -9 pt
$$
\vbox{
\nointerlineskip
\halign{\strut
    \vrule \ \ \hfil \frb #\ 
   &\vrule \hfil \ \ \fbb #\frb\ 
   &\vrule \hfil \ \ \frb #\ \hfil
   &\vrule \hfil \ \ \frb #\ 
   &\vrule \hfil \ \ \frb #\ \ \vrule \hskip 2 pt
   &\vrule \ \ \hfil \frb #\ 
   &\vrule \hfil \ \ \fbb #\frb\ 
   &\vrule \hfil \ \ \frb #\ \hfil
   &\vrule \hfil \ \ \frb #\ 
   &\vrule \hfil \ \ \frb #\ \vrule \cr%
\noalign{\hrule}
 & &3.29.43.311.977&141&170& & &3.11.17.1181.1723&9177&10900 \cr
7741&1136691627&4.9.5.17.43.47.977&7337&9272&7759&1141558143&8.9.25.7.11.19.23.109&629&134 \cr
 & &64.11.19.23.29.47.61&20539&21472& & &32.5.17.19.23.37.67&7705&11248 \cr
\noalign{\hrule}
 & &81.11.17.47.1597&37975&37084& & &3.125.7.11.19.2081&7397&3008 \cr
7742&1136918673&8.25.49.17.31.73.127&5499&14186&7760&1141688625&128.25.13.47.569&7353&7922 \cr
 & &32.9.5.7.13.41.47.173&7093&7280& & &512.9.17.19.43.233&11883&11008 \cr
\noalign{\hrule}
 & &5.11.13.19.31.37.73&423&58& & &9.11.29.41.89.109&5&94 \cr
7743&1137485635&4.9.11.19.29.31.47&481&1070&7761&1141914411&4.5.29.41.47.109&649&540 \cr
 & &16.3.5.13.29.37.107&107&696& & &32.27.25.11.47.59&3525&944 \cr
\noalign{\hrule}
 & &3.47.97.193.431&145&286& & &9.5.23.61.79.229&11527&2442 \cr
7744&1137693891&4.5.11.13.29.97.193&431&1692&7762&1142175285&4.27.11.37.11527&5777&5750 \cr
 & &32.9.5.29.47.431&145&48& & &16.125.11.23.37.53.109&21571&21800 \cr
\noalign{\hrule}
 & &9.5.7.13.17.59.277&701&1364& & &343.17.31.71.89&1131&1628 \cr
7745&1137717945&8.3.11.31.277.701&65&766&7763&1142228759&8.3.49.11.13.17.29.37&213&620 \cr
 & &32.5.11.13.31.383&4213&496& & &64.9.5.13.29.31.71&585&928 \cr
\noalign{\hrule}
 & &27.125.11.19.1613&14761&15886& & &49.31.241.3121&42471&54280 \cr
7746&1137769875&4.3.11.169.29.47.509&1675&3226&7764&1142532559&16.27.5.121.13.23.59&1681&3038 \cr
 & &16.25.67.509.1613&509&536& & &64.9.5.49.31.1681&1681&1440 \cr
\noalign{\hrule}
 & &13.17.23.41.43.127&35025&35968& & &5.49.11.13.17.19.101&783&12802 \cr
7747&1138088783&256.3.25.17.281.467&1221&3556&7765&1142946805&4.27.7.29.37.173&323&286 \cr
 & &2048.9.5.7.11.37.127&11655&11264& & &16.9.11.13.17.19.173&173&72 \cr
\noalign{\hrule}
 & &3.49.13.19.79.397&3495&4048& & &9.17.43.239.727&485&246 \cr
7748&1138759167&32.9.5.7.11.13.23.233&4727&632&7766&1143120987&4.27.5.41.97.727&1673&946 \cr
 & &512.11.29.79.163&4727&2816& & &16.5.7.11.41.43.239&287&440 \cr
\noalign{\hrule}
 & &9.25.7.37.113.173&143&1068& & &243.343.11.29.43&6341&4108 \cr
7749&1139217975&8.27.11.13.89.113&1145&1258&7767&1143298233&8.49.13.17.79.373&505&132 \cr
 & &32.5.11.13.17.37.229&2519&3536& & &64.3.5.11.17.79.101&6715&3232 \cr
\noalign{\hrule}
 & &5.13.47.79.4721&833&3888& & &9.31.67.193.317&143&460 \cr
7750&1139389745&32.243.49.17.79&47&506&7768&1143656433&8.5.11.13.23.31.193&8881&13314 \cr
 & &128.9.7.11.23.47&207&4928& & &32.3.7.83.107.317&749&1328 \cr
\noalign{\hrule}
 & &9.25.7.11.17.53.73&863&12788& & &5.11.13.47.67.509&4977&622 \cr
7751&1139517225&8.7.23.139.863&55&918&7769&1146031315&4.9.7.47.79.311&1159&1018 \cr
 & &32.27.5.11.17.23&3&368& & &16.3.19.61.79.509&1159&1896 \cr
\noalign{\hrule}
 & &5.41.73.271.281&5373&5738& & &81.49.23.29.433&28529&25498 \cr
7752&1139599715&4.27.19.151.199.281&539&820&7770&1146290859&4.7.11.19.47.61.607&2165&702 \cr
 & &32.3.5.49.11.19.41.199&9751&10032& & &16.27.5.13.433.607&607&520 \cr
\noalign{\hrule}
 & &9.5.121.17.109.113&131&1112& & &27.25.7.11.13.1697&799&898 \cr
7753&1140123105&16.5.11.17.131.139&763&678&7771&1146620475&4.3.25.7.13.17.47.449&18667&22192 \cr
 & &64.3.7.109.113.139&139&224& & &128.11.17.19.73.1697&1241&1216 \cr
\noalign{\hrule}
 & &5.7.11.29.41.47.53&13271&3636& & &9.5.11.13.19.83.113&251&334 \cr
7754&1140292615&8.9.7.23.101.577&47&530&7772&1146723435&4.11.19.113.167.251&155&1992 \cr
 & &32.3.5.47.53.101&303&16& & &64.3.5.31.83.251&251&992 \cr
\noalign{\hrule}
 & &5.49.121.17.31.73&923&318& & &3.5.7.113.163.593&627&514 \cr
7755&1140472795&4.3.49.13.31.53.71&2641&1122&7773&1146859035&4.9.5.11.19.257.593&2639&326 \cr
 & &16.9.11.13.17.19.139&2641&936& & &16.7.11.13.19.29.163&551&1144 \cr
\noalign{\hrule}
 & &3.7.37.61.24071&15421&8650& & &3.5.23.841.59.67&1181&176 \cr
7756&1140893187&4.25.49.173.2203&1269&9746&7774&1146943185&32.11.841.1181&6075&6916 \cr
 & &16.27.5.11.47.443&19935&4136& & &256.243.25.7.13.19&7695&11648 \cr
\noalign{\hrule}
 & &25.11.289.83.173&3699&3526& & &5.11.13.53.107.283&4675&996 \cr
7757&1141181525&4.27.11.41.43.83.137&8483&5080&7775&1147498495&8.3.125.121.17.83&3959&6084 \cr
 & &64.3.5.17.43.127.499&16383&15968& & &64.27.169.37.107&999&416 \cr
\noalign{\hrule}
 & &7.23.59.107.1123&2595&3718& & &9.7.11.19.29.31.97&1833&4646 \cr
7758&1141409339&4.3.5.7.11.169.23.173&513&4492&7776&1148201901&4.27.7.13.23.47.101&5147&2420 \cr
 & &32.81.13.19.1123&1053&304& & &32.5.121.13.5147&5147&11440 \cr
\noalign{\hrule}
}%
}
$$
\eject
\vglue -23 pt
\noindent\hskip 1 in\hbox to 6.5 in{\ 7777 -- 7812 \hfill\fbd 1149255261 -- 1161762019\frb}
\vskip -9 pt
$$
\vbox{
\nointerlineskip
\halign{\strut
    \vrule \ \ \hfil \frb #\ 
   &\vrule \hfil \ \ \fbb #\frb\ 
   &\vrule \hfil \ \ \frb #\ \hfil
   &\vrule \hfil \ \ \frb #\ 
   &\vrule \hfil \ \ \frb #\ \ \vrule \hskip 2 pt
   &\vrule \ \ \hfil \frb #\ 
   &\vrule \hfil \ \ \fbb #\frb\ 
   &\vrule \hfil \ \ \frb #\ \hfil
   &\vrule \hfil \ \ \frb #\ 
   &\vrule \hfil \ \ \frb #\ \vrule \cr%
\noalign{\hrule}
 & &9.49.11.19.37.337&1339&1510& & &9.17.23.41.53.151&959&7150 \cr
7777&1149255261&4.5.7.13.103.151.337&20757&30130&7795&1154664837&4.25.7.11.13.23.137&193&492 \cr
 & &16.3.25.11.17.23.37.131&3275&3128& & &32.3.5.7.11.41.193&385&3088 \cr
\noalign{\hrule}
 & &9.7.11.13.29.53.83&233&350& & &13.19.61.173.443&2889&2870 \cr
7778&1149287139&4.25.49.29.83.233&309&11726&7796&1154719813&4.27.5.7.41.61.107.173&7429&3124 \cr
 & &16.3.5.11.13.41.103&103&1640& & &32.9.11.17.19.23.71.107&68373&68816 \cr
\noalign{\hrule}
 & &9.125.7.13.103.109&107&2& & &49.23.61.107.157&533&594 \cr
7779&1149364125&4.3.25.13.103.107&539&436&7797&1154880853&4.27.11.13.41.107.157&8625&8174 \cr
 & &32.49.11.107.109&107&1232& & &16.81.125.13.23.61.67&8375&8424 \cr
\noalign{\hrule}
 & &5.11.41.3721.137&117&568& & &5.7.11.13.23.79.127&3219&1126 \cr
7780&1149547135&16.9.13.3721.71&5617&5546&7798&1154948795&4.3.29.37.127.563&2123&1560 \cr
 & &64.3.13.41.47.59.137&1833&1888& & &64.9.5.11.13.37.193&1737&1184 \cr
\noalign{\hrule}
 & &81.73.433.449&18337&14440& & &27.5.11.29.139.193&71737&79382 \cr
7781&1149587721&16.9.5.11.361.1667&39493&35522&7799&1155304755&4.19.23.2089.3119&5607&53654 \cr
 & &64.73.541.17761&17761&17312& & &16.9.7.89.139.193&89&56 \cr
\noalign{\hrule}
 & &3.5.121.17.19.37.53&563&658& & &27.7.17.19.23.823&335&488 \cr
7782&1149626445&4.7.11.17.47.53.563&10525&19314&7800&1155558663&16.3.5.7.19.23.61.67&869&442 \cr
 & &16.9.25.7.29.37.421&3045&3368& & &64.5.11.13.17.67.79&9581&12640 \cr
\noalign{\hrule}
 & &9.5.11.61.113.337&3943&34138& & &3.11.13.19.83.1709&11125&11092 \cr
7783&1149855795&4.169.101.3943&1921&2022&7801&1156194897&8.125.19.47.59.83.89&2277&6178 \cr
 & &16.3.169.17.113.337&169&136& & &32.9.25.11.23.59.3089&71047&70800 \cr
\noalign{\hrule}
 & &81.5.13.37.5903&7871&13774& & &81.125.11.13.17.47&1&142 \cr
7784&1149933915&4.9.17.71.97.463&2035&2132&7802&1156852125&4.27.125.17.71&2021&104 \cr
 & &32.5.11.13.17.37.41.71&2911&2992& & &64.13.43.47&1&1376 \cr
\noalign{\hrule}
 & &27.5.19.157.2857&37739&16544& & &9.5.17.19.97.821&26599&43186 \cr
7785&1150528185&64.11.13.47.2903&1193&1710&7803&1157523795&4.11.13.67.151.397&2619&1748 \cr
 & &256.9.5.13.19.1193&1193&1664& & &32.27.19.23.97.151&453&368 \cr
\noalign{\hrule}
 & &121.19.29.41.421&2691&2270& & &3.5.7.11.17.61.967&819&148 \cr
7786&1150808131&4.9.5.13.19.23.29.227&65&616&7804&1158209745&8.27.5.49.13.17.37&10637&12932 \cr
 & &64.3.25.7.11.169.23&11661&5600& & &64.11.53.61.967&53&32 \cr
\noalign{\hrule}
 & &5.13.1849.61.157&291&2332& & &27.49.29.109.277&42535&42812 \cr
7787&1151011745&8.3.5.11.43.53.97&12831&12874&7805&1158414831&8.5.343.11.47.139.181&1853&138 \cr
 & &32.9.7.11.13.41.47.157&4653&4592& & &32.3.17.23.47.109.139&6533&6256 \cr
\noalign{\hrule}
 & &27.7.11.169.29.113&569&250& & &5.11.13.19.269.317&2523&2588 \cr
7788&1151377227&4.3.125.13.113.569&781&3626&7806&1158433705&8.3.11.841.317.647&1&318 \cr
 & &16.25.49.11.37.71&925&3976& & &32.9.29.53.647&34291&4176 \cr
\noalign{\hrule}
 & &5.11.13.19.29.37.79&1377&508& & &5.23.29.43.59.137&2439&3796 \cr
7789&1151559695&8.81.17.19.37.127&275&428&7807&1159142615&8.9.13.73.137.271&49&88 \cr
 & &64.9.25.11.107.127&5715&3424& & &128.3.49.11.73.271&39347&52032 \cr
\noalign{\hrule}
 & &3.11.31.37.41.743&52003&39386& & &3.5.11.53.71.1867&6549&2786 \cr
7790&1153055013&4.7.17.19.23.47.419&3321&2930&7808&1159210965&4.9.7.11.37.59.199&71&2260 \cr
 & &16.81.5.41.293.419&11313&11720& & &32.5.59.71.113&59&1808 \cr
\noalign{\hrule}
 & &9.11.13.23.47.829&311&518& & &27.25.11.13.61.197&2857&3182 \cr
7791&1153343763&4.7.11.13.37.47.311&6015&6158&7809&1159940925&4.3.37.43.197.2857&641&7930 \cr
 & &16.3.5.311.401.3079&124711&123160& & &16.5.13.43.61.641&641&344 \cr
\noalign{\hrule}
 & &3.5.59.89.97.151&8771&138& & &3.7.361.47.3259&4755&1496 \cr
7792&1153670955&4.9.5.49.23.179&649&604&7810&1161204513&16.9.5.11.17.19.317&2171&1316 \cr
 & &32.7.11.23.59.151&77&368& & &128.7.13.17.47.167&2171&1088 \cr
\noalign{\hrule}
 & &9.11.13.23.127.307&1435&5626& & &3.5.49.11.13.43.257&251&6 \cr
7793&1154113389&4.3.5.7.13.29.41.97&575&614&7811&1161515355&4.9.11.13.43.251&329&230 \cr
 & &16.125.7.23.97.307&875&776& & &16.5.7.23.47.251&5773&376 \cr
\noalign{\hrule}
 & &25.121.19.53.379&22797&23062& & &1331.41.61.349&33441&47750 \cr
7794&1154500325&4.9.5.13.17.19.149.887&11749&25054&7812&1161762019&4.3.125.71.157.191&91&66 \cr
 & &16.3.17.31.379.12527&12527&12648& & &16.9.5.7.11.13.71.191&22365&19864 \cr
\noalign{\hrule}
}%
}
$$
\eject
\vglue -23 pt
\noindent\hskip 1 in\hbox to 6.5 in{\ 7813 -- 7848 \hfill\fbd 1163219187 -- 1178212889\frb}
\vskip -9 pt
$$
\vbox{
\nointerlineskip
\halign{\strut
    \vrule \ \ \hfil \frb #\ 
   &\vrule \hfil \ \ \fbb #\frb\ 
   &\vrule \hfil \ \ \frb #\ \hfil
   &\vrule \hfil \ \ \frb #\ 
   &\vrule \hfil \ \ \frb #\ \ \vrule \hskip 2 pt
   &\vrule \ \ \hfil \frb #\ 
   &\vrule \hfil \ \ \fbb #\frb\ 
   &\vrule \hfil \ \ \frb #\ \hfil
   &\vrule \hfil \ \ \frb #\ 
   &\vrule \hfil \ \ \frb #\ \vrule \cr%
\noalign{\hrule}
 & &3.13.43.61.83.137&375&418& & &25.7.31.61.3539&33&28 \cr
7813&1163219187&4.9.125.11.19.83.137&1405&172&7831&1171143575&8.3.5.49.11.31.3539&5567&2028 \cr
 & &32.625.11.43.281&6875&4496& & &64.9.11.169.19.293&49517&60192 \cr
\noalign{\hrule}
 & &9.29.73.107.571&1837&1266& & &3.5.7.17.23.103.277&75433&67222 \cr
7814&1164081141&4.27.11.73.167.211&175&2146&7832&1171340205&4.19.29.61.241.313&187&126 \cr
 & &16.25.7.29.37.167&6475&1336& & &16.9.7.11.17.19.29.241&6061&5784 \cr
\noalign{\hrule}
 & &9.11.17.19.23.1583&179&8& & &9.5.29.47.71.269&487&218 \cr
7815&1164250593&16.23.179.1583&1267&2850&7833&1171437165&4.3.29.71.109.487&299&1760 \cr
 & &64.3.25.7.19.181&4525&224& & &256.5.11.13.23.109&15587&2944 \cr
\noalign{\hrule}
 & &729.5.121.19.139&667&62& & &9.13.23.37.61.193&78505&85646 \cr
7816&1164799845&4.19.23.29.31.139&1677&964&7834&1172202291&4.5.7.11.17.229.2243&2381&138 \cr
 & &32.3.13.29.43.241&3133&19952& & &16.3.5.7.17.23.2381&2381&4760 \cr
\noalign{\hrule}
 & &3.169.17.29.59.79&3955&946& & &729.169.31.307&48811&49118 \cr
7817&1165021611&4.5.7.11.43.79.113&6119&2808&7835&1172503917&4.3.7.13.19.41.367.599&373&33770 \cr
 & &64.27.5.13.29.211&1055&288& & &16.5.11.41.307.373&2255&2984 \cr
\noalign{\hrule}
 & &27.13.23.37.47.83&209&1060& & &27.5.11.19.29.1433&2105&3538 \cr
7818&1165232601&8.5.11.13.19.53.83&1133&444&7836&1172530755&4.25.841.61.421&1183&342 \cr
 & &64.3.5.121.37.103&605&3296& & &16.9.7.169.19.421&2947&1352 \cr
\noalign{\hrule}
 & &3.37.47.101.2213&114365&109148& & &9.7.11.19.29.37.83&25739&28810 \cr
7819&1166067321&8.5.13.89.257.2099&407&1692&7837&1172639853&4.5.49.43.67.3677&1221&886 \cr
 & &64.9.11.13.37.47.89&1157&1056& & &16.3.11.37.443.3677&3677&3544 \cr
\noalign{\hrule}
 & &3.5.11.13.367.1483&6093&1322& & &27.5.13.43.79.197&67&62 \cr
7820&1167439845&4.27.11.661.677&12775&5504&7838&1174461795&4.9.13.31.67.79.197&6923&61886 \cr
 & &1024.25.7.43.73&15695&3584& & &16.7.11.23.29.43.97&7469&5336 \cr
\noalign{\hrule}
 & &3.17.1681.53.257&35&88& & &3.5.67.73.83.193&9825&4264 \cr
7821&1167741951&16.5.7.11.17.41.257&477&220&7839&1175233935&16.9.125.13.41.131&2123&3248 \cr
 & &128.9.25.7.121.53&3025&1344& & &512.7.11.13.29.193&2233&3328 \cr
\noalign{\hrule}
 & &9.13.31.53.6079&205&484& & &5.13.289.19.37.89&27459&43924 \cr
7822&1168572249&8.5.121.41.6079&559&5520&7840&1175321095&8.243.79.113.139&425&286 \cr
 & &256.3.25.13.23.43&575&5504& & &32.27.25.11.13.17.113&1485&1808 \cr
\noalign{\hrule}
 & &27.5.11.31.53.479&21941&22420& & &3.49.121.169.17.23&7349&7050 \cr
7823&1168690545&8.25.11.19.37.59.593&159&434&7841&1175347173&4.9.25.7.13.47.7349&31933&34208 \cr
 & &32.3.7.19.31.37.53.59&2183&2128& & &256.11.47.1069.2903&136441&136832 \cr
\noalign{\hrule}
 & &27.5.7.361.23.149&583&222& & &9.7.19.31.79.401&255&334 \cr
7824&1169103915&4.81.11.37.53.149&7015&5054&7842&1175512653&4.27.5.7.17.167.401&4763&9272 \cr
 & &16.5.7.11.361.23.61&61&88& & &64.11.17.19.61.433&11407&13856 \cr
\noalign{\hrule}
 & &9.5.11.17.19.71.103&34207&938& & &3.11.29.31.97.409&1391&1422 \cr
7825&1169239005&4.7.67.79.433&14229&14782&7843&1176978891&4.27.11.13.79.107.409&32045&266 \cr
 & &16.27.17.19.31.389&389&744& & &16.5.7.169.17.19.29&2261&6760 \cr
\noalign{\hrule}
 & &3.17.19.23.97.541&845&1386& & &27.5.19.23.71.281&12091&3926 \cr
7826&1169554899&4.27.5.7.11.169.17.19&4393&4328&7844&1177009245&4.9.13.107.113.151&253&1216 \cr
 & &64.7.11.13.23.191.541&2483&2464& & &512.11.19.23.151&1661&256 \cr
\noalign{\hrule}
 & &7.11.17.271.3299&38475&17608& & &27.47.59.79.199&301&11440 \cr
7827&1170283961&16.81.25.19.31.71&73&98&7845&1177046991&32.9.5.7.11.13.43&3239&3196 \cr
 & &64.9.49.31.71.73&15841&20448& & &256.7.17.41.47.79&697&896 \cr
\noalign{\hrule}
 & &19.29.83.157.163&1285&1122& & &3.625.11.13.23.191&36349&11974 \cr
7828&1170353203&4.3.5.11.17.19.157.257&2905&78&7846&1177873125&4.163.223.5987&3105&2882 \cr
 & &16.9.25.7.13.17.83&3825&728& & &16.27.5.11.23.131.163&1179&1304 \cr
\noalign{\hrule}
 & &3.49.17.29.31.521&803&716& & &25.47.313.3203&2189&1014 \cr
7829&1170479121&8.11.17.73.179.521&2105&10962&7847&1177983325&4.3.11.169.199.313&3129&940 \cr
 & &32.27.5.7.11.29.421&2105&1584& & &32.9.5.7.13.47.149&819&2384 \cr
\noalign{\hrule}
 & &13.31.83.157.223&781&34230& & &343.67.167.307&38925&18356 \cr
7830&1171082939&4.3.5.7.11.71.163&157&228&7848&1178212889&8.9.25.13.173.353&17&22 \cr
 & &32.9.19.157.163&163&2736& & &32.3.5.11.17.173.353&58245&47056 \cr
\noalign{\hrule}
}%
}
$$
\eject
\vglue -23 pt
\noindent\hskip 1 in\hbox to 6.5 in{\ 7849 -- 7884 \hfill\fbd 1178336965 -- 1200498299\frb}
\vskip -9 pt
$$
\vbox{
\nointerlineskip
\halign{\strut
    \vrule \ \ \hfil \frb #\ 
   &\vrule \hfil \ \ \fbb #\frb\ 
   &\vrule \hfil \ \ \frb #\ \hfil
   &\vrule \hfil \ \ \frb #\ 
   &\vrule \hfil \ \ \frb #\ \ \vrule \hskip 2 pt
   &\vrule \ \ \hfil \frb #\ 
   &\vrule \hfil \ \ \fbb #\frb\ 
   &\vrule \hfil \ \ \frb #\ \hfil
   &\vrule \hfil \ \ \frb #\ 
   &\vrule \hfil \ \ \frb #\ \vrule \cr%
\noalign{\hrule}
 & &5.13.19.37.107.241&2453&1062& & &9.25.11.43.53.211&8149&7676 \cr
7849&1178336965&4.9.11.59.223.241&455&214&7867&1190150775&8.3.19.29.53.101.281&50165&5486 \cr
 & &16.3.5.7.11.13.59.107&231&472& & &32.5.13.79.127.211&1651&1264 \cr
\noalign{\hrule}
 & &9.25.19.29.37.257&57233&66742& & &27.11.59.101.673&96961&102920 \cr
7850&1178878275&4.1331.13.17.43.151&551&1110&7868&1191090879&16.5.31.47.83.2063&89257&81972 \cr
 & &16.3.5.121.17.19.29.37&121&136& & &128.81.7.11.23.41.311&19803&19904 \cr
\noalign{\hrule}
 & &3.5.169.19.127.193&519&326& & &9.73.463.3919&195&268 \cr
7851&1180572315&4.9.19.127.163.173&473&1940&7869&1192124529&8.27.5.13.67.3919&2563&6482 \cr
 & &32.5.11.43.97.173&7439&17072& & &32.7.11.13.233.463&1631&2288 \cr
\noalign{\hrule}
 & &9.5.17.163.9473&253&236& & &3.5.7.169.17.59.67&759&424 \cr
7852&1181235735&8.3.5.11.23.59.9473&51527&23108&7870&1192481745&16.9.11.17.23.53.59&871&812 \cr
 & &64.7.17.53.109.433&22949&24416& & &128.7.13.23.29.53.67&1537&1472 \cr
\noalign{\hrule}
 & &25.41.1061.1087&12719&13806& & &5.7.17.101.103.193&1551&200 \cr
7853&1182139675&4.9.7.13.23.41.59.79&7953&8980&7871&1194628505&16.3.125.11.47.101&493&618 \cr
 & &32.27.5.11.23.241.449&108209&109296& & &64.9.17.29.47.103&423&928 \cr
\noalign{\hrule}
 & &27.5.121.19.37.103&271&136& & &49.11.37.139.431&13325&2622 \cr
7854&1182801015&16.11.17.19.103.271&18125&15144&7872&1194765187&4.3.25.7.13.19.23.41&57&148 \cr
 & &256.3.625.29.631&18299&16000& & &32.9.5.361.23.37&1035&5776 \cr
\noalign{\hrule}
 & &3.5.11.107.113.593&551&1144& & &81.25.7.11.79.97&481&586 \cr
7855&1183043895&16.121.13.19.29.107&1845&454&7873&1194853275&4.27.5.13.37.79.293&373&22 \cr
 & &64.9.5.29.41.227&9307&2784& & &16.11.37.293.373&10841&2984 \cr
\noalign{\hrule}
 & &625.11.13.17.19.41&87&538& & &121.43.163.1409&765&644 \cr
7856&1183593125&4.3.13.17.19.29.269&123&370&7874&1194957401&8.9.5.7.17.23.43.163&121&610 \cr
 & &16.9.5.37.41.269&2421&296& & &32.3.25.7.121.23.61&4209&2800 \cr
\noalign{\hrule}
 & &81.25.11.17.53.59&1463&562& & &9.25.13.23.109.163&77&38 \cr
7857&1184116725&4.7.121.19.59.281&3429&3710&7875&1195274925&4.3.5.7.11.19.109.163&4847&20332 \cr
 & &16.27.5.49.19.53.127&931&1016& & &32.13.17.23.37.131&629&2096 \cr
\noalign{\hrule}
 & &7.11.59.463.563&513&50& & &9.5.11.13.23.59.137&19517&22668 \cr
7858&1184219267&4.27.25.7.11.19.59&1689&1556&7876&1196324415&8.27.29.673.1889&8141&10030 \cr
 & &32.81.5.389.563&1945&1296& & &32.5.7.17.29.59.1163&8141&7888 \cr
\noalign{\hrule}
 & &9.5.11.23.29.37.97&673&992& & &9.29.37.1849.67&57113&66770 \cr
7859&1184962185&64.23.31.97.673&1131&1100&7877&1196338131&4.5.7.11.41.199.607&3219&1030 \cr
 & &512.3.25.11.13.29.673&3365&3328& & &16.3.25.29.37.41.103&1025&824 \cr
\noalign{\hrule}
 & &27.5.11.13.17.23.157&3089&3034& & &25.49.61.67.239&5031&944 \cr
7860&1185076035&4.9.17.23.37.41.3089&1777&128426&7878&1196571425&32.9.49.13.43.59&1195&1342 \cr
 & &16.157.409.1777&1777&3272& & &128.3.5.11.13.61.239&143&192 \cr
\noalign{\hrule}
 & &3.25.13.17.43.1663&27467&44042& & &5.19.23.43.47.271&3111&2294 \cr
7861&1185261675&4.121.361.61.227&9&218&7879&1196704835&4.3.17.31.37.61.271&1175&3432 \cr
 & &16.9.11.19.61.109&19947&1672& & &64.9.25.11.13.31.47&3069&2080 \cr
\noalign{\hrule}
 & &289.19.41.5273&109&5382& & &9.5.7.13.17.103.167&7981&7216 \cr
7862&1187115763&4.9.13.23.41.109&3781&3740&7880&1197447615&32.11.23.41.103.347&203&4020 \cr
 & &32.3.5.11.13.17.19.199&2985&2288& & &256.3.5.7.23.29.67&1541&3712 \cr
\noalign{\hrule}
 & &81.5.7.11.13.29.101&3845&3526& & &7.121.17.19.29.151&351&200 \cr
7863&1187431245&4.25.41.43.101.769&1647&878&7881&1198011199&16.27.25.7.121.13.17&755&92 \cr
 & &16.27.41.43.61.439&18877&20008& & &128.9.125.23.151&2875&576 \cr
\noalign{\hrule}
 & &3.13.23.31.151.283&5225&5194& & &729.5.23.79.181&4355&9944 \cr
7864&1188276531&4.25.49.11.13.19.53.283&207&1208&7882&1198756665&16.3.25.11.13.67.113&553&9028 \cr
 & &64.9.5.7.19.23.53.151&1855&1824& & &128.7.37.61.79&259&3904 \cr
\noalign{\hrule}
 & &11.19.37.227.677&1587&910& & &27.25.29.101.607&17&118 \cr
7865&1188399707&4.3.5.7.13.19.529.37&15741&8864&7883&1200084525&4.5.17.29.59.607&231&376 \cr
 & &256.81.11.53.277&14681&10368& & &64.3.7.11.17.47.59&5593&20768 \cr
\noalign{\hrule}
 & &81.5.19.349.443&2665&9296& & &7.11.13.17.19.47.79&351&1150 \cr
7866&1189700865&32.3.25.7.13.41.83&443&418&7884&1200498299&4.27.25.7.11.169.23&799&722 \cr
 & &128.11.13.19.83.443&913&832& & &16.3.25.17.361.23.47&475&552 \cr
\noalign{\hrule}
}%
}
$$
\eject
\vglue -23 pt
\noindent\hskip 1 in\hbox to 6.5 in{\ 7885 -- 7920 \hfill\fbd 1201497115 -- 1217447845\frb}
\vskip -9 pt
$$
\vbox{
\nointerlineskip
\halign{\strut
    \vrule \ \ \hfil \frb #\ 
   &\vrule \hfil \ \ \fbb #\frb\ 
   &\vrule \hfil \ \ \frb #\ \hfil
   &\vrule \hfil \ \ \frb #\ 
   &\vrule \hfil \ \ \frb #\ \ \vrule \hskip 2 pt
   &\vrule \ \ \hfil \frb #\ 
   &\vrule \hfil \ \ \fbb #\frb\ 
   &\vrule \hfil \ \ \frb #\ \hfil
   &\vrule \hfil \ \ \frb #\ 
   &\vrule \hfil \ \ \frb #\ \vrule \cr%
\noalign{\hrule}
 & &5.7.13.23.29.37.107&723&1738& & &5.17.19.1369.547&1529&1206 \cr
7885&1201497115&4.3.11.13.37.79.241&717&310&7903&1209381445&4.9.11.1369.67.139&893&476 \cr
 & &16.9.5.31.239.241&7471&17208& & &32.3.7.11.17.19.47.67&3149&3696 \cr
\noalign{\hrule}
 & &5.7.11.13.17.71.199&1653&2962& & &3.19.37.83.6911&130219&125488 \cr
7886&1202165965&4.3.19.29.199.1481&2145&3626&7904&1209749817&32.11.23.31.107.1217&1197&20 \cr
 & &16.9.5.49.11.13.19.37&703&504& & &256.9.5.7.19.23.31&4991&1920 \cr
\noalign{\hrule}
 & &3.5.7.13.43.103.199&11269&6840& & &5.11.2197.17.19.31&6105&3908 \cr
7887&1203071415&16.27.25.19.59.191&559&1034&7905&1209920855&8.3.25.121.37.977&19587&16562 \cr
 & &64.11.13.43.47.191&2101&1504& & &32.9.49.169.6529&6529&7056 \cr
\noalign{\hrule}
 & &5.13.17.31.113.311&5373&4268& & &3.11.23.59.61.443&233&210 \cr
7888&1203823465&8.27.11.97.113.199&7&106&7906&1210116963&4.9.5.7.11.59.61.233&299&3544 \cr
 & &32.3.7.53.97.199&35987&9552& & &64.13.23.233.443&233&416 \cr
\noalign{\hrule}
 & &7.11.1849.8461&4381&4080& & &5.7.19.41.157.283&311&468 \cr
7889&1204617953&32.3.5.11.13.17.43.337&5939&210&7907&1211411215&8.9.5.7.13.283.311&13079&7136 \cr
 & &128.9.25.7.5939&5939&14400& & &512.3.11.29.41.223&6467&8448 \cr
\noalign{\hrule}
 & &9.25.7.23.29.31.37&4323&12352& & &9.7.11.13.17.41.193&38125&15544 \cr
7890&1204952175&128.27.11.131.193&1885&3326&7908&1211899689&16.625.29.61.67&33&28 \cr
 & &512.5.13.29.1663&1663&3328& & &128.3.125.7.11.29.67&3625&4288 \cr
\noalign{\hrule}
 & &5.41.107.137.401&7227&7432& & &9.5.49.11.107.467&715&248 \cr
7891&1205043095&16.9.11.73.401.929&137&1066&7909&1211998095&16.25.49.121.13.31&3127&2802 \cr
 & &64.3.11.13.41.73.137&949&1056& & &64.3.31.53.59.467&1643&1888 \cr
\noalign{\hrule}
 & &3.5.11.13.23.53.461&28427&44872& & &3.5.11.31.79.3001&181&214 \cr
7892&1205402055&16.7.31.71.79.131&1035&1166&7910&1212659085&4.31.107.181.3001&3159&158 \cr
 & &64.9.5.7.11.23.53.79&237&224& & &16.243.13.79.181&2353&648 \cr
\noalign{\hrule}
 & &5.7.121.41.53.131&1807&366& & &81.5.7.11.19.23.89&3757&698 \cr
7893&1205547805&4.3.5.7.11.13.61.139&123&262&7911&1212878205&4.13.289.89.349&3025&1512 \cr
 & &16.9.13.41.61.131&117&488& & &64.27.25.7.121.17&85&352 \cr
\noalign{\hrule}
 & &3.5.11.59.83.1493&5763&10660& & &3.11.289.47.2707&403&114 \cr
7894&1206351465&8.9.25.13.17.41.113&2407&418&7912&1213382973&4.9.13.19.31.2707&2465&242 \cr
 & &32.11.19.29.41.83&1189&304& & &16.5.121.17.29.31&155&2552 \cr
\noalign{\hrule}
 & &5.7.11.19.23.71.101&2559&2996& & &5.49.121.13.23.137&111&188 \cr
7895&1206484895&8.3.49.71.107.853&1313&2166&7913&1214348135&8.3.5.7.11.37.47.137&69&1576 \cr
 & &32.9.13.361.101.107&2033&1872& & &128.9.23.37.197&7289&576 \cr
\noalign{\hrule}
 & &9.11.23.29.47.389&100223&100890& & &7.17.23.577.769&43175&25488 \cr
7896&1207281339&4.81.5.19.31.53.59.61&47&4988&7914&1214442481&32.27.25.11.59.157&7501&5774 \cr
 & &32.29.31.43.47.59&1333&944& & &128.3.13.577.2887&2887&2496 \cr
\noalign{\hrule}
 & &5.11.13.43.107.367&2891&3258& & &5.13.31.47.101.127&3035&2934 \cr
7897&1207325405&4.9.5.49.59.107.181&18447&7768&7915&1214781035&4.9.25.13.31.163.607&103&2222 \cr
 & &64.27.11.13.43.971&971&864& & &16.3.11.101.103.607&3399&4856 \cr
\noalign{\hrule}
 & &25.49.19.23.37.61&89&1314& & &9.5.7.23.43.47.83&37&178 \cr
7898&1208228525&4.9.19.37.73.89&10097&9394&7916&1215298035&4.3.7.23.37.83.89&205&44 \cr
 & &16.3.7.11.23.61.439&439&264& & &32.5.11.37.41.89&1517&15664 \cr
\noalign{\hrule}
 & &3.25.7.17.191.709&13137&9592& & &9.25.13.199.2089&2507&418 \cr
7899&1208614575&16.9.5.11.29.109.151&1817&4978&7917&1215954675&4.11.19.23.109.199&699&500 \cr
 & &64.11.19.23.79.131&33143&48032& & &32.3.125.19.23.233&4427&1840 \cr
\noalign{\hrule}
 & &9.25.11.41.43.277&1387&1838& & &49.11.43.73.719&2015&1296 \cr
7900&1208668725&4.3.19.73.277.919&875&44&7918&1216491199&32.81.5.7.13.31.73&1&218 \cr
 & &32.125.7.11.19.73&365&2128& & &128.27.5.13.109&351&34880 \cr
\noalign{\hrule}
 & &7.121.31.41.1123&1919&3042& & &3.5.11.13.17.61.547&1709&2256 \cr
7901&1208951051&4.9.7.169.19.31.101&451&2370&7919&1216727655&32.9.11.17.47.1709&16853&12200 \cr
 & &16.27.5.11.13.41.79&1755&632& & &512.25.19.61.887&4435&4864 \cr
\noalign{\hrule}
 & &27.7.11.31.73.257&3299&3640& & &5.23.131.211.383&297&86 \cr
7902&1209127689&16.5.49.13.73.3299&139&3438&7920&1217447845&4.27.5.11.23.43.131&1477&488 \cr
 & &64.9.5.13.139.191&12415&4448& & &64.9.7.11.61.211&549&2464 \cr
\noalign{\hrule}
}%
}
$$
\eject
\vglue -23 pt
\noindent\hskip 1 in\hbox to 6.5 in{\ 7921 -- 7956 \hfill\fbd 1217687985 -- 1228539455\frb}
\vskip -9 pt
$$
\vbox{
\nointerlineskip
\halign{\strut
    \vrule \ \ \hfil \frb #\ 
   &\vrule \hfil \ \ \fbb #\frb\ 
   &\vrule \hfil \ \ \frb #\ \hfil
   &\vrule \hfil \ \ \frb #\ 
   &\vrule \hfil \ \ \frb #\ \ \vrule \hskip 2 pt
   &\vrule \ \ \hfil \frb #\ 
   &\vrule \hfil \ \ \fbb #\frb\ 
   &\vrule \hfil \ \ \frb #\ \hfil
   &\vrule \hfil \ \ \frb #\ 
   &\vrule \hfil \ \ \frb #\ \vrule \cr%
\noalign{\hrule}
 & &81.5.17.47.53.71&623&5128& & &3.25.11.841.41.43&7371&6346 \cr
7921&1217687985&16.7.47.89.641&1771&2412&7939&1223213475&4.243.7.13.19.29.167&113&4730 \cr
 & &128.9.49.11.23.67&12397&4288& & &16.5.7.11.13.43.113&113&728 \cr
\noalign{\hrule}
 & &11.31.109.32761&34965&2204& & &3.5.11.19.37.61.173&173&382 \cr
7922&1217693609&8.27.5.7.19.29.37&181&218&7940&1224095235&4.61.29929.191&9139&20790 \cr
 & &32.9.5.29.109.181&261&80& & &16.27.5.7.11.13.19.37&117&56 \cr
\noalign{\hrule}
 & &3.5.7.13.17.73.719&9017&330& & &9.5.11.23.41.43.61&257&2020 \cr
7923&1217960835&4.9.25.11.71.127&559&584&7941&1224377055&8.25.61.101.257&3243&3182 \cr
 & &64.11.13.43.71.73&781&1376& & &32.3.23.37.43.47.101&1739&1616 \cr
\noalign{\hrule}
 & &27.121.13.23.29.43&1097&476& & &27.5.11.13.229.277&35071&2324 \cr
7924&1218110751&8.7.17.29.43.1097&13689&4960&7942&1224574065&8.7.17.83.2063&1073&990 \cr
 & &512.81.5.169.31&2015&768& & &32.9.5.7.11.17.29.37&1073&1904 \cr
\noalign{\hrule}
 & &5.13.17.23.191.251&1773&1474& & &5.7.11.23.41.3373&8207&8658 \cr
7925&1218420515&4.9.5.11.67.197.251&43&208&7943&1224584515&4.9.7.13.23.29.37.283&8497&2540 \cr
 & &128.3.13.43.67.197&8471&12864& & &32.3.5.841.127.293&37211&40368 \cr
\noalign{\hrule}
 & &9.5.49.17.19.29.59&11&844& & &5.7.11.19.191.877&17767&22152 \cr
7926&1218599865&8.11.29.59.211&861&850&7944&1225313705&16.3.7.13.71.109.163&3547&11286 \cr
 & &32.3.25.7.17.41.211&1055&656& & &64.81.11.19.3547&3547&2592 \cr
\noalign{\hrule}
 & &27.49.31.71.419&6919&6500& & &25.11.13.353.971&4369&486 \cr
7927&1220095737&8.125.7.11.13.17.31.37&1879&10056&7945&1225377725&4.243.5.13.17.257&353&418 \cr
 & &128.3.25.419.1879&1879&1600& & &16.81.11.17.19.353&323&648 \cr
\noalign{\hrule}
 & &3.49.11.61.89.139&1845&2516& & &3.5.49.11.13.89.131&81437&88642 \cr
7928&1220238327&8.27.5.17.37.41.139&3719&1424&7946&1225419195&4.23.31.37.41.47.71&70609&37098 \cr
 & &256.41.89.3719&3719&5248& & &16.81.49.11.131.229&229&216 \cr
\noalign{\hrule}
 & &9.5.121.271.827&1077&1904& & &5.7.11.19.41.61.67&237&542 \cr
7929&1220317065&32.27.5.7.11.17.359&2309&3794&7947&1225752605&4.3.7.11.67.79.271&8125&13284 \cr
 & &128.49.271.2309&2309&3136& & &32.243.625.13.41&3159&2000 \cr
\noalign{\hrule}
 & &5.7.23.31.59.829&35541&6526& & &3.25.11.41.101.359&1891&2250 \cr
7930&1220574005&4.9.11.13.251.359&1093&2170&7948&1226460675&4.27.3125.11.31.61&1711&1414 \cr
 & &16.3.5.7.11.31.1093&1093&264& & &16.7.29.31.59.61.101&13237&13688 \cr
\noalign{\hrule}
 & &3.5.53.107.113.127&143&196& & &3.5.49.13.19.29.233&35033&782 \cr
7931&1220767815&8.5.49.11.13.107.127&3033&3922&7949&1226699565&4.17.23.53.661&319&342 \cr
 & &32.9.7.11.37.53.337&7077&6512& & &16.9.11.17.19.29.53&901&264 \cr
\noalign{\hrule}
 & &27.289.47.3329&2237&5566& & &9.5.13.719.2917&10241&36494 \cr
7932&1220880789&4.121.23.47.2237&1377&860&7950&1226933955&4.49.11.19.71.257&4181&702 \cr
 & &32.81.5.11.17.23.43&1265&2064& & &16.27.11.13.37.113&339&3256 \cr
\noalign{\hrule}
 & &7.17.41.179.1399&18029&39330& & &11.17.73.139.647&12673&13320 \cr
7933&1221804059&4.9.5.121.19.23.149&1399&1432&7951&1227675383&16.9.5.19.23.29.37.73&13033&11646 \cr
 & &64.3.5.11.23.179.1399&345&352& & &64.81.5.647.13033&13033&12960 \cr
\noalign{\hrule}
 & &3.49.29.43.59.113&15275&7664& & &3.19.83.139.1867&145&1722 \cr
7934&1222121103&32.25.7.13.47.479&583&1062&7952&1227756003&4.9.5.7.29.41.139&583&722 \cr
 & &128.9.5.11.13.53.59&2145&3392& & &16.7.11.361.41.53&8569&2968 \cr
\noalign{\hrule}
 & &27.5.49.19.71.137&5525&4202& & &81.11.19.29.41.61&7363&5716 \cr
7935&1222537995&4.125.11.13.17.19.191&2877&752&7953&1227843441&8.3.19.37.199.1429&1769&340 \cr
 & &128.3.7.11.13.47.137&517&832& & &64.5.17.29.61.199&995&544 \cr
\noalign{\hrule}
 & &3.13.17.19.37.43.61&339&820& & &25.13.17.199.1117&8997&9992 \cr
7936&1222551447&8.9.5.17.41.43.113&143&874&7954&1228113575&16.3.5.13.1249.2999&16371&22616 \cr
 & &32.5.11.13.19.23.41&2255&368& & &256.27.11.17.107.257&31779&32896 \cr
\noalign{\hrule}
 & &5.19.43.229.1307&3289&3246& & &11.13.23.29.79.163&1455&3272 \cr
7937&1222652755&4.3.11.13.19.23.229.541&27447&20414&7955&1228221137&16.3.5.11.13.97.409&133&276 \cr
 & &16.9.7.23.59.173.1307&10899&10856& & &128.9.5.7.19.23.97&4365&8512 \cr
\noalign{\hrule}
 & &25.11.13.43.73.109&29&4716& & &5.11.13.31.43.1289&2109&2324 \cr
7938&1223189825&8.9.5.11.29.131&89&56&7956&1228539455&8.3.7.19.37.83.1289&515&774 \cr
 & &128.3.7.89.131&34977&448& & &32.27.5.19.43.83.103&8549&8208 \cr
\noalign{\hrule}
}%
}
$$
\eject
\vglue -23 pt
\noindent\hskip 1 in\hbox to 6.5 in{\ 7957 -- 7992 \hfill\fbd 1228722495 -- 1241383715\frb}
\vskip -9 pt
$$
\vbox{
\nointerlineskip
\halign{\strut
    \vrule \ \ \hfil \frb #\ 
   &\vrule \hfil \ \ \fbb #\frb\ 
   &\vrule \hfil \ \ \frb #\ \hfil
   &\vrule \hfil \ \ \frb #\ 
   &\vrule \hfil \ \ \frb #\ \ \vrule \hskip 2 pt
   &\vrule \ \ \hfil \frb #\ 
   &\vrule \hfil \ \ \fbb #\frb\ 
   &\vrule \hfil \ \ \frb #\ \hfil
   &\vrule \hfil \ \ \frb #\ 
   &\vrule \hfil \ \ \frb #\ \vrule \cr%
\noalign{\hrule}
 & &3.5.7.11.13.19.59.73&191&264& & &27.25.11.31.41.131&931&94 \cr
7957&1228722495&16.9.121.19.59.191&1105&16&7975&1236269925&4.49.11.19.47.131&697&744 \cr
 & &512.5.13.17.191&191&4352& & &64.3.49.17.19.31.41&833&608 \cr
\noalign{\hrule}
 & &25.13.17.19.23.509&29429&29106& & &27.13.31.71.1601&56975&56696 \cr
7958&1228942325&4.27.5.49.11.13.29429&9937&19492&7976&1236854151&16.3.25.13.19.43.53.373&497&2948 \cr
 & &32.9.121.19.443.523&63283&63792& & &128.5.7.11.67.71.373&24991&24640 \cr
\noalign{\hrule}
 & &5.13.17.59.109.173&207&1210& & &9.25.343.11.31.47&83&358 \cr
7959&1229382115&4.9.25.121.23.173&247&272&7977&1236883725&4.7.31.47.83.179&85&132 \cr
 & &128.3.121.13.17.19.23&2783&3648& & &32.3.5.11.17.83.179&3043&1328 \cr
\noalign{\hrule}
 & &27.25.7.11.41.577&44399&56576& & &9.25.7.13.23.37.71&15367&7708 \cr
7960&1229572575&512.13.17.29.1531&577&954&7978&1237119975&8.7.121.41.47.127&115&402 \cr
 & &2048.9.17.53.577&901&1024& & &32.3.5.11.23.67.127&737&2032 \cr
\noalign{\hrule}
 & &27.7.61.107.997&445&552& & &5.13.31.41.71.211&4361&4290 \cr
7961&1229902191&16.81.5.7.23.61.89&9235&3806&7979&1237655315&4.3.25.49.11.169.31.89&387&29962 \cr
 & &64.25.11.173.1847&47575&59104& & &16.27.7.43.71.211&189&344 \cr
\noalign{\hrule}
 & &3.11.13.37.179.433&145&34& & &7.17.73.89.1601&1557&44 \cr
7962&1230268611&4.5.11.13.17.29.433&999&1432&7980&1237801943&8.9.7.11.73.173&1853&1780 \cr
 & &64.27.5.29.37.179&145&288& & &64.3.5.11.17.89.109&545&1056 \cr
\noalign{\hrule}
 & &3.343.121.41.241&533&190& & &5.7.11.17.43.53.83&657&74 \cr
7963&1230273429&4.5.121.13.19.1681&15183&16756&7981&1238032565&4.9.5.7.37.73.83&689&1054 \cr
 & &32.9.5.7.59.71.241&1065&944& & &16.3.13.17.31.37.53&1209&296 \cr
\noalign{\hrule}
 & &9.49.11.6859.37&6149&710& & &121.43.313.761&2575&35298 \cr
7964&1231101333&4.3.5.121.13.43.71&851&722&7982&1239318179&4.9.25.37.53.103&59&44 \cr
 & &16.5.361.23.37.71&355&184& & &32.3.5.11.37.53.59&10915&2544 \cr
\noalign{\hrule}
 & &27.25.7.11.19.29.43&1369&1502& & &9.5.7.11.17.53.397&457&1528 \cr
7965&1231443675&4.3.25.1369.43.751&589&1664&7983&1239420105&16.11.53.191.457&7575&2548 \cr
 & &1024.13.19.31.1369&17797&15872& & &128.3.25.49.13.101&3535&832 \cr
\noalign{\hrule}
 & &9.25.7.47.127.131&4103&4150& & &5.59.71.73.811&2057&1998 \cr
7966&1231553925&4.625.11.83.127.373&1833&8708&7984&1240006835&4.27.121.17.37.71.73&5677&494 \cr
 & &32.3.7.13.47.311.373&4849&4976& & &16.9.7.13.19.37.811&2331&1976 \cr
\noalign{\hrule}
 & &3.49.19.29.67.227&33&34& & &3.11.13.61.83.571&211&460 \cr
7967&1231883373&4.9.49.11.17.19.29.227&1465&101572&7985&1240227417&8.5.13.23.211.571&1577&1278 \cr
 & &32.5.67.293.379&1895&4688& & &32.9.19.71.83.211&4047&3376 \cr
\noalign{\hrule}
 & &9.25.7.11.17.47.89&163&262& & &37.41.661.1237&37587&13130 \cr
7968&1231998075&4.7.47.89.131.163&1633&2550&7986&1240385669&4.3.5.11.13.17.67.101&2701&984 \cr
 & &16.3.25.17.23.71.163&1633&1304& & &64.9.13.37.41.73&949&288 \cr
\noalign{\hrule}
 & &37.47.331.2141&5053&7194& & &3.7.23.31.41.43.47&36157&8164 \cr
7969&1232378869&4.3.11.31.47.109.163&2141&5520&7987&1240677753&8.11.13.19.157.173&1075&828 \cr
 & &128.9.5.11.23.2141&1035&704& & &64.9.25.23.43.157&471&800 \cr
\noalign{\hrule}
 & &5.11.13.373.4621&259&4362& & &3.17.23.31.149.229&149&242 \cr
7970&1232397595&4.3.5.7.13.37.727&4743&4708&7988&1240741923&4.121.22201.229&24955&2754 \cr
 & &32.27.11.17.31.37.107&30969&29104& & &16.81.5.7.17.23.31&189&40 \cr
\noalign{\hrule}
 & &625.11.19.9437&4031&5406& & &25.11.13.449.773&4887&5162 \cr
7971&1232708125&4.3.5.17.19.29.53.139&2553&88&7989&1240800275&4.27.29.89.181.449&625&176 \cr
 & &64.9.11.23.37.53&1219&10656& & &128.3.625.11.29.181&5249&4800 \cr
\noalign{\hrule}
 & &9.17.53.383.397&3949&2800& & &9.5.11.13.29.61.109&8251&7336 \cr
7972&1232981559&32.3.25.7.11.53.359&3613&5362&7990&1240803135&16.3.7.29.37.131.223&109&22 \cr
 & &128.49.383.3613&3613&3136& & &64.7.11.37.109.223&1561&1184 \cr
\noalign{\hrule}
 & &9.7.11.169.53.199&5921&1840& & &3.11.13.19.29.59.89&6207&956 \cr
7973&1235232999&32.3.5.13.23.31.191&2801&1592&7991&1241226129&8.9.11.239.2069&12455&10304 \cr
 & &512.5.199.2801&2801&1280& & &1024.5.7.23.47.53&42665&24064 \cr
\noalign{\hrule}
 & &3.5.7.11.13.17.47.103&2403&652& & &5.11.13.19.23.29.137&2043&3638 \cr
7974&1235689455&8.81.7.11.89.163&289&278&7992&1241383715&4.9.17.107.137.227&259&152 \cr
 & &32.289.89.139.163&22657&24208& & &64.3.7.17.19.37.227&11577&8288 \cr
\noalign{\hrule}
}%
}
$$
\eject
\vglue -23 pt
\noindent\hskip 1 in\hbox to 6.5 in{\ 7993 -- 8028 \hfill\fbd 1241815575 -- 1256176233\frb}
\vskip -9 pt
$$
\vbox{
\nointerlineskip
\halign{\strut
    \vrule \ \ \hfil \frb #\ 
   &\vrule \hfil \ \ \fbb #\frb\ 
   &\vrule \hfil \ \ \frb #\ \hfil
   &\vrule \hfil \ \ \frb #\ 
   &\vrule \hfil \ \ \frb #\ \ \vrule \hskip 2 pt
   &\vrule \ \ \hfil \frb #\ 
   &\vrule \hfil \ \ \fbb #\frb\ 
   &\vrule \hfil \ \ \frb #\ \hfil
   &\vrule \hfil \ \ \frb #\ 
   &\vrule \hfil \ \ \frb #\ \vrule \cr%
\noalign{\hrule}
 & &3.25.49.11.13.17.139&206769&206756& & &9.5.169.17.19.509&4859&4774 \cr
7993&1241815575&8.9.7.121.37.127.157.439&87635&45344&8011&1250315235&4.3.7.11.31.43.113.509&24947&135400 \cr
 & &512.5.13.17.109.439.1031&112379&112384& & &64.25.13.19.101.677&3385&3232 \cr
\noalign{\hrule}
 & &3.7.31.1237.1543&665&572& & &3.5.7.19.619.1013&559&454 \cr
7994&1242557841&8.5.49.11.13.19.1543&1237&306&8012&1250958765&4.13.19.43.227.619&99&718 \cr
 & &32.9.5.11.13.17.1237&663&880& & &16.9.11.13.227.359&11847&23608 \cr
\noalign{\hrule}
 & &5.19.37.349.1013&155&858& & &27.7.19.29.41.293&85&104 \cr
7995&1242682555&4.3.25.11.13.31.349&6579&2146&8013&1251021807&16.5.13.17.29.41.293&7239&4774 \cr
 & &16.27.17.29.37.43&783&5848& & &64.3.7.11.13.19.31.127&4433&4064 \cr
\noalign{\hrule}
 & &3.5.13.19.37.43.211&1441&2074& & &9.25.11.13.19.23.89&1241&1450 \cr
7996&1243772205&4.11.13.17.43.61.131&8055&10678&8014&1251382275&4.625.17.29.73.89&1069&444 \cr
 & &16.9.5.17.19.179.281&4777&4296& & &32.3.29.37.73.1069&39553&33872 \cr
\noalign{\hrule}
 & &81.7.41.73.733&391&5522& & &13.53.73.149.167&4389&4462 \cr
7997&1243923723&4.11.17.23.41.251&105&146&8015&1251540251&4.3.7.11.13.19.23.97.149&1837&1590 \cr
 & &16.3.5.7.11.17.23.73&187&920& & &16.9.5.7.121.53.97.167&5445&5432 \cr
\noalign{\hrule}
 & &9.13.19.53.59.179&473&294& & &9.5.17.421.3889&41929&22484 \cr
7998&1244286459&4.27.49.11.19.43.53&1475&9602&8016&1252510785&8.7.11.23.73.1823&1313&510 \cr
 & &16.25.7.59.4801&4801&1400& & &32.3.5.7.13.17.23.101&1313&2576 \cr
\noalign{\hrule}
 & &5.11.23.37.67.397&2609&1758& & &9.125.17.19.3449&7169&10076 \cr
7999&1244966195&4.3.5.67.293.2609&3293&16338&8017&1253280375&8.25.11.67.107.229&6447&722 \cr
 & &16.9.7.37.89.389&2723&6408& & &32.3.7.11.361.307&1463&4912 \cr
\noalign{\hrule}
 & &81.7.17.23.41.137&583&376& & &25.121.19.113.193&11439&11914 \cr
8000&1245272049&16.9.11.17.41.47.53&137&560&8018&1253472275&4.9.7.23.31.37.41.113&55&10564 \cr
 & &512.5.7.11.53.137&583&1280& & &32.3.5.11.19.23.139&139&1104 \cr
\noalign{\hrule}
 & &27.17.71.167.229&1739&1100& & &49.29.73.107.113&13455&1364 \cr
8001&1246301127&8.3.25.11.37.47.229&59&170&8019&1254235703&8.9.5.7.11.13.23.31&1217&876 \cr
 & &32.125.11.17.47.59&5875&10384& & &64.27.5.73.1217&1217&4320 \cr
\noalign{\hrule}
 & &13.31.107.137.211&287&3030& & &5.11.23.29.31.1103&911&6426 \cr
8002&1246498747&4.3.5.7.41.101.137&429&530&8020&1254370205&4.27.7.17.31.911&1103&1630 \cr
 & &16.9.25.11.13.41.53&4059&10600& & &16.9.5.7.163.1103&163&504 \cr
\noalign{\hrule}
 & &81.5.7.307.1433&2167&2132& & &361.47.107.691&101&792 \cr
8003&1247204385&8.27.11.13.41.197.307&8183&106&8021&1254489079&16.9.11.19.101.107&25&82 \cr
 & &32.49.11.13.53.167&8851&16016& & &64.3.25.11.41.101&3075&35552 \cr
\noalign{\hrule}
 & &27.49.13.29.41.61&701&92& & &125.11.29.73.431&153&278 \cr
8004&1247426271&8.9.7.23.41.701&503&440&8022&1254587125&4.9.11.17.29.73.139&811&3220 \cr
 & &128.5.11.503.701&27665&44864& & &32.3.5.7.17.23.811&17031&6256 \cr
\noalign{\hrule}
 & &5.7.11.23.29.43.113&75837&67108& & &27.343.23.71.83&645&988 \cr
8005&1247766905&8.3.17.19.883.1487&8249&6762&8023&1255226679&8.81.5.13.19.43.83&1309&230 \cr
 & &32.9.49.19.23.73.113&1197&1168& & &32.25.7.11.17.23.43&1075&2992 \cr
\noalign{\hrule}
 & &25.361.1369.101&35343&1118& & &121.13.43.67.277&43099&24540 \cr
8006&1247877725&4.27.7.11.13.17.43&185&202&8024&1255312201&8.3.5.7.47.131.409&201&208 \cr
 & &16.3.5.7.11.13.37.101&273&88& & &256.9.5.13.47.67.131&5895&6016 \cr
\noalign{\hrule}
 & &25.7.11.73.83.107&55233&42458& & &7.11.19.31.89.311&1535&1224 \cr
8007&1248002525&4.9.13.17.361.23.71&9125&7918&8025&1255325687&16.9.5.7.11.17.19.307&163&68 \cr
 & &16.3.125.19.37.73.107&285&296& & &128.3.289.163.307&47107&58944 \cr
\noalign{\hrule}
 & &3.5.7.17.29.89.271&7711&36166& & &5.11.13.31.181.313&129&184 \cr
8008&1248520035&4.11.169.107.701&279&422&8026&1255713745&16.3.13.23.31.43.181&743&8526 \cr
 & &16.9.13.31.107.211&22577&9672& & &64.9.49.29.743&36407&8352 \cr
\noalign{\hrule}
 & &3.7.13.137.173.193&2095&1958& & &27.5.343.47.577&63167&72428 \cr
8009&1248781989&4.5.11.13.89.173.419&42607&29880&8027&1255745295&8.13.19.43.113.953&3765&8624 \cr
 & &64.9.25.83.137.311&7775&7968& & &256.3.5.49.11.19.251&2761&2432 \cr
\noalign{\hrule}
 & &81.5.11.17.29.569&1763&3358& & &9.13.37.61.67.71&1615&1012 \cr
8010&1249703235&4.9.17.23.41.43.73&677&20&8028&1256176233&8.5.11.13.17.19.23.61&3219&18286 \cr
 & &32.5.23.43.677&677&15824& & &32.3.29.37.41.223&1189&3568 \cr
\noalign{\hrule}
}%
}
$$
\eject
\vglue -23 pt
\noindent\hskip 1 in\hbox to 6.5 in{\ 8029 -- 8064 \hfill\fbd 1256310363 -- 1268100975\frb}
\vskip -9 pt
$$
\vbox{
\nointerlineskip
\halign{\strut
    \vrule \ \ \hfil \frb #\ 
   &\vrule \hfil \ \ \fbb #\frb\ 
   &\vrule \hfil \ \ \frb #\ \hfil
   &\vrule \hfil \ \ \frb #\ 
   &\vrule \hfil \ \ \frb #\ \ \vrule \hskip 2 pt
   &\vrule \ \ \hfil \frb #\ 
   &\vrule \hfil \ \ \fbb #\frb\ 
   &\vrule \hfil \ \ \frb #\ \hfil
   &\vrule \hfil \ \ \frb #\ 
   &\vrule \hfil \ \ \frb #\ \vrule \cr%
\noalign{\hrule}
 & &3.49.11.29.73.367&285&82& & &27.5.17.53.97.107&48931&47476 \cr
8029&1256310363&4.9.5.7.11.19.41.73&367&290&8047&1262449665&8.9.11.13.83.167.293&133&34 \cr
 & &16.25.19.29.41.367&475&328& & &32.7.13.17.19.83.293&26663&25232 \cr
\noalign{\hrule}
 & &7.11.23.37.127.151&3489&35620& & &5.11.19.31.127.307&17&324 \cr
8030&1256611279&8.3.5.13.137.1163&513&650&8048&1263048655&8.81.5.17.19.127&403&232 \cr
 & &32.81.125.169.19&13689&38000& & &128.9.13.17.29.31&6409&576 \cr
\noalign{\hrule}
 & &11.29.41.307.313&75033&66130& & &25.11.29.251.631&32147&47922 \cr
8031&1256774189&4.27.5.7.17.389.397&5125&1624&8049&1263088475&4.3.49.17.31.61.163&3045&2008 \cr
 & &64.3.625.49.29.41&1875&1568& & &64.9.5.343.29.251&343&288 \cr
\noalign{\hrule}
 & &5.7.13.17.23.37.191&9801&6434& & &5.13.19.41.61.409&587&648 \cr
8032&1257254635&4.81.121.23.3217&217&3000&8050&1263292615&16.81.41.409.587&143&266 \cr
 & &64.243.125.7.31&6075&992& & &64.27.7.11.13.19.587&6457&6048 \cr
\noalign{\hrule}
 & &3.5.11.169.23.37.53&467&430& & &9.5.23.29.71.593&1919&140 \cr
8033&1257697155&4.25.11.13.43.53.467&207&482&8051&1263721545&8.3.25.7.19.23.101&6523&4402 \cr
 & &16.9.23.43.241.467&10363&11208& & &32.11.31.71.593&31&176 \cr
\noalign{\hrule}
 & &9.25.7.11.13.37.151&67&158& & &3.125.11.13.17.19.73&103&322 \cr
8034&1258332075&4.11.37.67.79.151&19599&7670&8052&1264423875&4.5.7.11.13.19.23.103&1503&268 \cr
 & &16.3.5.13.47.59.139&2773&1112& & &32.9.67.103.167&20703&2672 \cr
\noalign{\hrule}
 & &9.5.7.121.17.29.67&3739&4234& & &5.11.13.29.71.859&145071&149366 \cr
8035&1258976565&4.11.841.73.3739&19&822&8053&1264606915&4.729.7.47.199.227&325&1718 \cr
 & &16.3.19.137.3739&3739&20824& & &16.81.25.13.47.859&405&376 \cr
\noalign{\hrule}
 & &3.11.17.19.31.37.103&143&180& & &529.457.5231&7871&2640 \cr
8036&1259264919&8.27.5.121.13.31.103&37&3230&8054&1264609943&32.3.5.11.17.23.463&2351&2742 \cr
 & &32.25.13.17.19.37&325&16& & &128.9.5.457.2351&2351&2880 \cr
\noalign{\hrule}
 & &9.5.11.17.109.1373&155&1528& & &3.25.11.23.131.509&2787&2812 \cr
8037&1259363655&16.25.31.109.191&4323&1598&8055&1265234025&8.9.19.23.37.131.929&30745&3628 \cr
 & &64.3.11.17.47.131&131&1504& & &64.5.11.13.19.43.907&17233&17888 \cr
\noalign{\hrule}
 & &5.7.11.17.193.997&3&190& & &3.11.13.17.29.31.193&75&452 \cr
8038&1259395445&4.3.25.7.19.997&9559&9384&8056&1265386551&8.9.25.11.113.193&151&344 \cr
 & &64.9.121.17.23.79&1817&3168& & &128.5.43.113.151&6493&36160 \cr
\noalign{\hrule}
 & &7.23.29.43.6277&2515&3762& & &5.13.19.53.83.233&9&74 \cr
8039&1260214459&4.9.5.7.11.19.23.503&493&10&8057&1265834245&4.9.19.37.53.233&18875&18172 \cr
 & &16.3.25.11.17.19.29&323&6600& & &32.3.125.7.11.59.151&30975&26576 \cr
\noalign{\hrule}
 & &3.25.11.19.37.41.53&5017&1502& & &27.5.121.17.47.97&19&116 \cr
8040&1260285675&4.5.11.29.173.751&23161&18144&8058&1266011505&8.121.17.19.29.47&1475&582 \cr
 & &256.81.7.19.23.53&621&896& & &32.3.25.29.59.97&145&944 \cr
\noalign{\hrule}
 & &11.13.67.149.883&73&810& & &9.5.13.67.79.409&1327&2354 \cr
8041&1260543427&4.81.5.13.73.149&1005&932&8059&1266429645&4.5.11.67.107.1327&101&636 \cr
 & &32.243.25.67.233&5825&3888& & &32.3.53.101.1327&5353&21232 \cr
\noalign{\hrule}
 & &5.7.13.17.19.23.373&297&76& & &9.31.47.71.1361&10585&9224 \cr
8042&1260812735&8.27.5.7.11.361.23&1781&746&8060&1267122303&16.5.29.47.73.1153&2201&3564 \cr
 & &32.3.11.13.137.373&137&528& & &128.81.11.31.71.73&657&704 \cr
\noalign{\hrule}
 & &9.5.7.19.59.3571&5917&4796& & &25.11.17.47.73.79&5221&546 \cr
8043&1260973665&8.3.5.7.11.61.97.109&533&2822&8061&1267154075&4.3.7.13.23.47.227&7821&6232 \cr
 & &32.13.17.41.83.97&44239&26384& & &64.27.11.19.41.79&513&1312 \cr
\noalign{\hrule}
 & &3.25.7.13.19.37.263&2057&1268& & &3.5.13.289.83.271&517&296 \cr
8044&1261867425&8.121.13.17.37.317&2375&1746&8062&1267593015&16.5.11.17.37.47.83&3523&378 \cr
 & &32.9.125.121.19.97&1815&1552& & &64.27.7.11.13.271&63&352 \cr
\noalign{\hrule}
 & &5.49.13.61.73.89&4163&198& & &27.7.11.13.17.31.89&3035&2762 \cr
8045&1262269645&4.9.11.23.73.181&6227&6986&8063&1267647381&4.9.5.89.607.1381&4061&58084 \cr
 & &16.3.7.13.479.499&1437&3992& & &32.13.31.131.1117&1117&2096 \cr
\noalign{\hrule}
 & &27.7.121.289.191&1315&742& & &3.25.17.23.83.521&451&1624 \cr
8046&1262345931&4.9.5.49.17.53.263&1375&1222&8064&1268100975&16.7.11.29.41.521&855&334 \cr
 & &16.625.11.13.47.263&29375&27352& & &64.9.5.7.11.19.167&9519&2464 \cr
\noalign{\hrule}
}%
}
$$
\eject
\vglue -23 pt
\noindent\hskip 1 in\hbox to 6.5 in{\ 8065 -- 8100 \hfill\fbd 1268617141 -- 1287505711\frb}
\vskip -9 pt
$$
\vbox{
\nointerlineskip
\halign{\strut
    \vrule \ \ \hfil \frb #\ 
   &\vrule \hfil \ \ \fbb #\frb\ 
   &\vrule \hfil \ \ \frb #\ \hfil
   &\vrule \hfil \ \ \frb #\ 
   &\vrule \hfil \ \ \frb #\ \ \vrule \hskip 2 pt
   &\vrule \ \ \hfil \frb #\ 
   &\vrule \hfil \ \ \fbb #\frb\ 
   &\vrule \hfil \ \ \frb #\ \hfil
   &\vrule \hfil \ \ \frb #\ 
   &\vrule \hfil \ \ \frb #\ \vrule \cr%
\noalign{\hrule}
 & &11.23.73.149.461&267&194& & &729.5.7.23.37.59&731&2914 \cr
8065&1268617141&4.3.11.23.89.97.149&18901&5640&8083&1281082635&4.7.17.23.31.43.47&2035&702 \cr
 & &64.9.5.41.47.461&2115&1312& & &16.27.5.11.13.37.47&517&104 \cr
\noalign{\hrule}
 & &3.7.11.13.19.23.967&155&1156& & &27.5.11.17.193.263&11563&37618 \cr
8066&1269004737&8.5.289.31.967&747&220&8084&1281410955&4.7.31.373.2687&1157&1530 \cr
 & &64.9.25.11.17.83&6225&544& & &16.9.5.7.13.17.31.89&1157&1736 \cr
\noalign{\hrule}
 & &9.7.11.13.29.43.113&587&700& & &3.5.13.43.67.2281&1037&3318 \cr
8067&1269467199&8.25.49.29.43.587&9565&7458&8085&1281454395&4.9.7.17.43.61.79&253&134 \cr
 & &32.3.125.11.113.1913&1913&2000& & &16.11.23.61.67.79&4819&2024 \cr
\noalign{\hrule}
 & &11.13.41.53.61.67&1665&1082& & &121.31.563.607&221&342 \cr
8068&1269990293&4.9.5.13.37.61.541&3379&3392&8086&1281870491&4.9.13.17.19.31.607&16043&4510 \cr
 & &512.3.5.31.53.109.541&83855&83712& & &16.3.5.11.41.61.263&12505&6312 \cr
\noalign{\hrule}
 & &5.49.11.59.61.131&783&134& & &625.7.13.22543&7044733&7044642 \cr
8069&1270608955&4.27.5.7.29.61.67&473&808&8087&1282133125&4.9.11.43.47.173.757.947&31901&650 \cr
 & &64.9.11.29.43.101&11223&3232& & &16.3.25.13.19.23.47.73.173&78913&78888 \cr
\noalign{\hrule}
 & &3.25.11.17.31.37.79&2101&5246& & &3.5.169.19.79.337&121&126 \cr
8070&1270847325&4.5.121.43.61.191&2449&2754&8088&1282296795&4.27.7.121.13.79.337&26825&202 \cr
 & &16.81.17.31.79.191&191&216& & &16.25.11.29.37.101&2929&16280 \cr
\noalign{\hrule}
 & &3.5.7.13.17.157.349&157&38& & &49.11.23.31.47.71&2379&35890 \cr
8071&1271471565&4.19.24649.349&9009&15640&8089&1282432459&4.3.5.13.37.61.97&23359&23298 \cr
 & &64.9.5.7.11.13.17.23&253&96& & &16.9.5.7.11.47.71.353&353&360 \cr
\noalign{\hrule}
 & &5.7.13.1021.2741&1377&1364& & &5.7.121.41.83.89&105&16 \cr
8072&1273345255&8.81.5.7.11.17.31.1021&377&5482&8090&1282641745&32.3.25.49.41.83&1089&2314 \cr
 & &32.3.11.13.17.29.2741&493&528& & &128.27.121.13.89&351&64 \cr
\noalign{\hrule}
 & &13.31.37.127.673&40975&20112& & &11.289.191.2113&533&2646 \cr
8073&1274458081&32.3.25.11.149.419&6057&4418&8091&1282990357&4.27.49.13.41.191&935&402 \cr
 & &128.27.2209.673&2209&1728& & &16.81.5.7.11.17.67&567&2680 \cr
\noalign{\hrule}
 & &81.25.7.11.13.17.37&239&1786& & &169.289.109.241&37555&11286 \cr
8074&1274998725&4.11.19.37.47.239&9333&490&8092&1283004229&4.27.5.7.11.19.29.37&221&482 \cr
 & &16.9.5.49.17.61&7&488& & &16.3.5.7.11.13.17.241&165&56 \cr
\noalign{\hrule}
 & &9.5.11.13.17.107.109&1537&2714& & &3.5.19.23.29.43.157&1815&1796 \cr
8075&1275873885&4.3.5.17.23.29.53.59&163&3172&8093&1283331345&8.9.25.121.29.43.449&57371&43654 \cr
 & &32.13.53.61.163&9943&848& & &32.11.13.23.73.103.557&97747&98032 \cr
\noalign{\hrule}
 & &41.61.101.5059&551&5610& & &9.5.7.11.23.71.227&10811&5590 \cr
8076&1277908459&4.3.5.11.17.19.29.41&549&230&8094&1284444315&4.3.25.13.19.43.569&227&98 \cr
 & &16.27.25.17.23.61&459&4600& & &16.49.19.227.569&569&1064 \cr
\noalign{\hrule}
 & &121.13.71.11449&5719&5730& & &5.7.11.29.31.47.79&14227&4338 \cr
8077&1278658667&4.3.5.7.11.13.19.43.71.191&115&39804&8095&1285124995&4.9.7.41.241.347&323&1364 \cr
 & &32.9.25.7.23.31.107&4991&3600& & &32.3.11.17.19.31.41&323&1968 \cr
\noalign{\hrule}
 & &3.5.67.79.89.181&439&104& & &9.25.13.29.109.139&44011&56764 \cr
8078&1278974055&16.13.79.89.439&17919&16762&8096&1285183575&8.11.23.617.4001&1393&5394 \cr
 & &64.9.11.289.29.181&3179&2784& & &32.3.7.23.29.31.199&4991&3184 \cr
\noalign{\hrule}
 & &3.5.169.83.6079&3293&2786& & &5.7.19.47.79.521&1573&72 \cr
8079&1279051995&4.5.7.37.83.89.199&16533&1178&8097&1286424545&16.9.121.13.521&475&46 \cr
 & &16.9.7.11.19.31.167&22211&8184& & &64.3.25.11.19.23&33&3680 \cr
\noalign{\hrule}
 & &3.13.19.31.103.541&935&688& & &9.5.79.113.3203&341&9268 \cr
8080&1280013033&32.5.11.17.31.43.103&5113&684&8098&1286693145&8.3.5.7.11.31.331&113&218 \cr
 & &256.9.5.19.5113&5113&1920& & &32.11.31.109.113&3379&176 \cr
\noalign{\hrule}
 & &25.11.439.10607&2889&7718& & &9.7.19.23.101.463&3311&13960 \cr
8081&1280530075&4.27.25.17.107.227&109&5566&8099&1287432153&16.5.49.11.43.349&8787&8314 \cr
 & &16.9.121.23.109&10791&184& & &64.3.5.29.101.4157&4157&4640 \cr
\noalign{\hrule}
 & &9.25.11.13.53.751&427&262& & &61.83.109.2333&1125&1208 \cr
8082&1280661525&4.3.5.7.61.131.751&3419&1166&8100&1287505711&16.9.125.61.109.151&8987&2338 \cr
 & &16.11.13.53.61.263&263&488& & &64.3.5.7.11.19.43.167&66633&75680 \cr
\noalign{\hrule}
}%
}
$$
\eject
\vglue -23 pt
\noindent\hskip 1 in\hbox to 6.5 in{\ 8101 -- 8136 \hfill\fbd 1287753467 -- 1300546395\frb}
\vskip -9 pt
$$
\vbox{
\nointerlineskip
\halign{\strut
    \vrule \ \ \hfil \frb #\ 
   &\vrule \hfil \ \ \fbb #\frb\ 
   &\vrule \hfil \ \ \frb #\ \hfil
   &\vrule \hfil \ \ \frb #\ 
   &\vrule \hfil \ \ \frb #\ \ \vrule \hskip 2 pt
   &\vrule \ \ \hfil \frb #\ 
   &\vrule \hfil \ \ \fbb #\frb\ 
   &\vrule \hfil \ \ \frb #\ \hfil
   &\vrule \hfil \ \ \frb #\ 
   &\vrule \hfil \ \ \frb #\ \vrule \cr%
\noalign{\hrule}
 & &49.11.169.67.211&1641&1102& & &7.121.191.8011&15561&7550 \cr
8101&1287753467&4.3.13.19.29.67.547&825&448&8119&1295995547&4.9.25.49.13.19.151&8011&14186 \cr
 & &512.9.25.7.11.547&4923&6400& & &16.3.41.173.8011&519&328 \cr
\noalign{\hrule}
 & &7.11.29.67.79.109&43&36& & &27.25.7.13.47.449&18667&22192 \cr
8102&1288300321&8.9.11.29.43.67.109&475&2686&8120&1296251775&32.9.11.19.73.1697&799&898 \cr
 & &32.3.25.17.19.43.79&3225&5168& & &128.17.19.47.73.449&1241&1216 \cr
\noalign{\hrule}
 & &27.7.13.23.29.787&85&176& & &5.7.13.29.227.433&2817&2812 \cr
8103&1289750553&32.3.5.11.17.23.787&203&2158&8121&1296949745&8.9.7.19.29.37.227.313&5&6578 \cr
 & &128.7.11.13.29.83&83&704& & &32.3.5.11.13.19.23.37&4807&1776 \cr
\noalign{\hrule}
 & &27.5.7.11.127.977&65249&66646& & &9.5.7.11.17.361.61&2419&13826 \cr
8104&1289801205&4.7.47.71.709.919&4885&1548&8122&1297147005&4.7.31.41.59.223&4337&4560 \cr
 & &32.9.5.43.709.977&709&688& & &128.3.5.19.59.4337&4337&3776 \cr
\noalign{\hrule}
 & &9.5.43.59.89.127&159&286& & &9.5.17.101.103.163&295&194 \cr
8105&1290406995&4.27.11.13.43.53.59&23575&6052&8123&1297202085&4.3.25.17.59.97.103&3839&1414 \cr
 & &32.25.17.23.41.89&697&1840& & &16.7.11.59.101.349&4543&2792 \cr
\noalign{\hrule}
 & &3.25.11.37.67.631&10319&10504& & &3.11.17.23.193.521&3445&2924 \cr
8106&1290505425&16.5.13.17.67.101.607&67671&6364&8124&1297435359&8.5.13.289.23.43.53&4351&594 \cr
 & &128.9.37.43.73.103&7519&8256& & &32.27.11.19.53.229&4351&7632 \cr
\noalign{\hrule}
 & &27.13.23.307.521&38863&52930& & &3.5.7.121.13.29.271&7&3516 \cr
8107&1291252131&4.5.11.67.79.3533&116183&120528&8125&1298031735&8.9.5.49.293&1441&1196 \cr
 & &128.243.31.223.521&2007&1984& & &64.11.13.23.131&3013&32 \cr
\noalign{\hrule}
 & &27.11.13.53.59.107&505&184& & &3.19.67.151.2251&25&176 \cr
8108&1291848129&16.9.5.11.23.59.101&221&428&8126&1298081919&32.25.11.19.2251&603&1648 \cr
 & &128.5.13.17.101.107&505&1088& & &1024.9.5.67.103&309&2560 \cr
\noalign{\hrule}
 & &25.7.121.17.37.97&1827&1198& & &3.25.11.13.19.23.277&7565&2302 \cr
8109&1291950275&4.9.49.29.97.599&121&170&8127&1298250525&4.125.17.89.1151&513&638 \cr
 & &16.3.5.121.17.29.599&599&696& & &16.27.11.17.19.29.89&1513&2088 \cr
\noalign{\hrule}
 & &11.47.73.97.353&35991&1750& & &3.125.169.31.661&17&22 \cr
8110&1292289581&4.27.125.7.31.43&1649&1606&8128&1298617125&4.25.11.13.17.31.661&8559&34 \cr
 & &16.9.25.11.17.73.97&153&200& & &16.27.289.317&317&20808 \cr
\noalign{\hrule}
 & &25.11.19.31.79.101&387&482& & &3.11.13.19.37.59.73&1207&2156 \cr
8111&1292398525&4.9.5.31.43.101.241&3373&242&8129&1298935209&8.49.121.17.37.71&11115&19706 \cr
 & &16.3.121.43.3373&3373&11352& & &32.9.5.13.19.59.167&167&240 \cr
\noalign{\hrule}
 & &9.5.7.13.31.61.167&451&634& & &9.5.13.17.193.677&4687&1406 \cr
8112&1293188715&4.3.11.13.41.167.317&17675&35264&8130&1299423645&4.5.13.19.37.43.109&187&372 \cr
 & &512.25.7.19.29.101&9595&7424& & &32.3.11.17.19.31.109&2071&5456 \cr
\noalign{\hrule}
 & &17.41.557.3331&16895&39732& & &9.49.13.31.71.103&6677&1630 \cr
8113&1293190799&8.3.5.7.11.31.43.109&205&96&8131&1299688299&4.5.11.31.163.607&5865&812 \cr
 & &512.9.25.11.31.41&3069&6400& & &32.3.25.7.17.23.29&575&7888 \cr
\noalign{\hrule}
 & &7.13.17.23.41.887&3025&9234& & &9.121.73.83.197&4541&4292 \cr
8114&1293974227&4.243.25.121.17.19&161&26&8132&1299855447&8.3.19.29.37.197.239&415&176 \cr
 & &16.9.5.7.11.13.19.23&95&792& & &256.5.11.19.29.37.83&3515&3712 \cr
\noalign{\hrule}
 & &5.11.13.83.113.193&633&1876& & &9.121.17.23.43.71&535&1742 \cr
8115&1294255105&8.3.5.7.67.83.211&313&102&8133&1299964347&4.5.11.13.43.67.107&199&672 \cr
 & &32.9.7.17.67.313&37247&9648& & &256.3.5.7.107.199&21293&4480 \cr
\noalign{\hrule}
 & &25.29.53.67.503&309&1634& & &9.25.49.13.47.193&13717&1558 \cr
8116&1294960925&4.3.19.43.103.503&157&660&8134&1300101075&4.7.11.19.29.41.43&75&376 \cr
 & &32.9.5.11.103.157&16171&1584& & &64.3.25.19.29.47&29&608 \cr
\noalign{\hrule}
 & &17.41.563.3301&3135&6436& & &3.121.23.29.41.131&7&950 \cr
8117&1295348711&8.3.5.11.19.41.1609&415&1194&8135&1300431891&4.25.7.11.19.131&333&322 \cr
 & &32.9.25.11.83.199&16517&39600& & &16.9.5.49.19.23.37&931&4440 \cr
\noalign{\hrule}
 & &5.2197.79.1493&3219&4246& & &27.5.43.157.1427&4089&2662 \cr
8118&1295647795&4.3.11.169.29.37.193&5925&328&8136&1300546395&4.81.5.1331.29.47&43&362 \cr
 & &64.9.25.11.41.79&1845&352& & &16.121.43.47.181&5687&1448 \cr
\noalign{\hrule}
}%
}
$$
\eject
\vglue -23 pt
\noindent\hskip 1 in\hbox to 6.5 in{\ 8137 -- 8172 \hfill\fbd 1301675375 -- 1312834683\frb}
\vskip -9 pt
$$
\vbox{
\nointerlineskip
\halign{\strut
    \vrule \ \ \hfil \frb #\ 
   &\vrule \hfil \ \ \fbb #\frb\ 
   &\vrule \hfil \ \ \frb #\ \hfil
   &\vrule \hfil \ \ \frb #\ 
   &\vrule \hfil \ \ \frb #\ \ \vrule \hskip 2 pt
   &\vrule \ \ \hfil \frb #\ 
   &\vrule \hfil \ \ \fbb #\frb\ 
   &\vrule \hfil \ \ \frb #\ \hfil
   &\vrule \hfil \ \ \frb #\ 
   &\vrule \hfil \ \ \frb #\ \vrule \cr%
\noalign{\hrule}
 & &125.7.11.13.101.103&2553&2452& & &25.17.59.109.479&189&1664 \cr
8137&1301675375&8.3.25.23.37.103.613&5757&9568&8155&1309190825&256.27.7.13.479&649&170 \cr
 & &512.9.13.19.529.101&4761&4864& & &1024.3.5.11.17.59&33&512 \cr
\noalign{\hrule}
 & &5.13.23.61.109.131&9065&9174& & &9.5.11.61.131.331&377&278 \cr
8138&1302173405&4.3.25.49.11.37.131.139&4209&934&8156&1309285395&4.13.29.61.139.331&9039&560 \cr
 & &16.9.49.11.23.61.467&4203&4312& & &128.3.5.7.13.23.131&161&832 \cr
\noalign{\hrule}
 & &9.25.7.11.17.19.233&2921&1040& & &3.121.17.31.41.167&7705&7826 \cr
8139&1303862175&32.125.7.13.23.127&17&108&8157&1309837947&4.5.7.13.17.23.41.43.67&3629&8712 \cr
 & &256.27.17.23.127&381&2944& & &64.9.5.121.19.67.191&10887&10720 \cr
\noalign{\hrule}
 & &3.5.7.11.37.131.233&2119&2774& & &3.25.7.11.23.71.139&1269&9434 \cr
8140&1304400405&4.11.13.19.37.73.163&1097&4194&8158&1310849925&4.81.5.47.53.89&77&158 \cr
 & &16.9.73.233.1097&1097&1752& & &16.7.11.53.79.89&4187&712 \cr
\noalign{\hrule}
 & &9.11.13.29.73.479&9287&9766& & &27.25.11.13.289.47&2077&1102 \cr
8141&1305073341&4.11.13.19.37.251.257&3051&290&8159&1311099075&4.9.19.29.31.47.67&221&202 \cr
 & &16.27.5.19.29.37.113&3515&2712& & &16.13.17.29.31.67.101&6767&7192 \cr
\noalign{\hrule}
 & &27.25.7.347.797&14647&5278& & &9.5.7.23.37.67.73&79&286 \cr
8142&1306741275&4.49.13.29.97.151&3069&4330&8160&1311105915&4.7.11.13.37.67.79&31&438 \cr
 & &16.9.5.11.29.31.433&4763&7192& & &16.3.13.31.73.79&2449&104 \cr
\noalign{\hrule}
 & &9.13.289.29.31.43&3227&530& & &5.13.17.19.197.317&1899&2222 \cr
8143&1307109141&4.3.5.7.43.53.461&221&682&8161&1311116755&4.9.5.11.101.197.211&317&2638 \cr
 & &16.5.11.13.17.31.53&53&440& & &16.3.101.317.1319&1319&2424 \cr
\noalign{\hrule}
 & &17.961.191.419&473&54& & &9.11.23.37.79.197&2119&6650 \cr
8144&1307433773&4.27.11.31.43.191&2095&1904&8162&1311167187&4.3.25.7.11.13.19.163&197&34 \cr
 & &128.9.5.7.11.17.419&385&576& & &16.25.13.17.19.197&323&2600 \cr
\noalign{\hrule}
 & &3.5.7.31.47.83.103&667&418& & &25.11.13.31.11831&459&316 \cr
8145&1307868765&4.11.19.23.29.47.103&297&1660&8163&1311170575&8.27.17.79.11831&3901&7930 \cr
 & &32.27.5.121.23.83&1089&368& & &32.9.5.13.47.61.83&5063&6768 \cr
\noalign{\hrule}
 & &81.25.13.17.37.79&1301&1696& & &27.5.11.13.23.2953&2117&36272 \cr
8146&1308115575&64.5.13.17.53.1301&243&22&8164&1311176295&32.29.73.2267&1097&1170 \cr
 & &256.243.11.1301&1301&4224& & &128.9.5.13.29.1097&1097&1856 \cr
\noalign{\hrule}
 & &13.17.29.43.47.101&1121&900& & &9.5.11.17.19.59.139&479&2162 \cr
8147&1308211489&8.9.25.19.29.59.101&3677&748&8165&1311216885&4.5.23.47.59.479&189&2584 \cr
 & &64.3.11.17.19.3677&11031&6688& & &64.27.7.17.19.23&23&672 \cr
\noalign{\hrule}
 & &9.5.121.23.31.337&2501&7946& & &5.7.11.17.43.59.79&1257&1988 \cr
8148&1308330045&4.23.29.41.61.137&3707&4650&8166&1311768535&8.3.49.71.79.419&14811&18290 \cr
 & &16.3.25.11.29.31.337&29&40& & &32.9.5.31.59.4937&4937&4464 \cr
\noalign{\hrule}
 & &3.5.49.11.13.59.211&213&424& & &25.11.13.37.47.211&80109&87916 \cr
8149&1308452145&16.9.5.11.53.59.71&2743&3098&8167&1311771175&8.81.23.31.43.709&665&44 \cr
 & &64.13.53.211.1549&1549&1696& & &64.3.5.7.11.19.31.43&5719&2976 \cr
\noalign{\hrule}
 & &27.5.49.19.29.359&1993&572& & &3.5.13.2089.3221&1037&1052 \cr
8150&1308506535&8.11.13.359.1993&817&1176&8168&1312090455&8.13.17.61.263.3221&99&3320 \cr
 & &128.3.49.11.13.19.43&473&832& & &128.9.5.11.17.61.83&15189&11968 \cr
\noalign{\hrule}
 & &9.5.7.11.37.59.173&2291&1426& & &27.5.11.19.193.241&1547&3032 \cr
8151&1308588435&4.11.23.29.31.37.79&707&366&8169&1312364295&16.7.13.17.193.379&207&14 \cr
 & &16.3.7.23.61.79.101&7979&11224& & &64.9.49.23.379&1127&12128 \cr
\noalign{\hrule}
 & &25.121.13.107.311&1767&5188& & &9.49.11.19.29.491&2479&1940 \cr
8152&1308618025&8.3.5.11.19.31.1297&1171&126&8170&1312394391&8.5.19.29.37.67.97&1323&620 \cr
 & &32.27.7.31.1171&36301&3024& & &64.27.25.49.31.97&2425&2976 \cr
\noalign{\hrule}
 & &9.7.11.13.31.43.109&145&548& & &5.289.23.127.311&3131&36366 \cr
8153&1308980673&8.5.29.43.109.137&1365&4526&8171&1312682795&4.3.11.19.29.31.101&289&300 \cr
 & &32.3.25.7.13.31.73&73&400& & &32.9.25.289.29.101&1305&1616 \cr
\noalign{\hrule}
 & &7.11.83.239.857&1265&408& & &3.419.809.1291&195&614 \cr
8154&1309023793&16.3.5.121.17.23.83&5999&4044&8172&1312834683&4.9.5.13.307.1291&4609&1846 \cr
 & &128.9.7.337.857&337&576& & &16.11.169.71.419&1859&568 \cr
\noalign{\hrule}
}%
}
$$
\eject
\vglue -23 pt
\noindent\hskip 1 in\hbox to 6.5 in{\ 8173 -- 8208 \hfill\fbd 1312836525 -- 1329431791\frb}
\vskip -9 pt
$$
\vbox{
\nointerlineskip
\halign{\strut
    \vrule \ \ \hfil \frb #\ 
   &\vrule \hfil \ \ \fbb #\frb\ 
   &\vrule \hfil \ \ \frb #\ \hfil
   &\vrule \hfil \ \ \frb #\ 
   &\vrule \hfil \ \ \frb #\ \ \vrule \hskip 2 pt
   &\vrule \ \ \hfil \frb #\ 
   &\vrule \hfil \ \ \fbb #\frb\ 
   &\vrule \hfil \ \ \frb #\ \hfil
   &\vrule \hfil \ \ \frb #\ 
   &\vrule \hfil \ \ \frb #\ \vrule \cr%
\noalign{\hrule}
 & &27.25.7.11.13.29.67&1867&58& & &23.41.337.4159&6897&6920 \cr
8173&1312836525&4.13.841.1867&11715&12556&8191&1321692769&16.3.5.121.19.173.4159&21231&24518 \cr
 & &32.3.5.11.43.71.73&3053&1168& & &64.27.5.7.11.13.23.41.337&2079&2080 \cr
\noalign{\hrule}
 & &27.7.29.31.59.131&34765&34796& & &27.121.47.79.109&203&124 \cr
8174&1313242119&8.3.5.7.17.29.409.8699&827&9526&8192&1322210439&8.9.7.121.29.31.47&395&3356 \cr
 & &32.5.11.409.433.827&358091&359920& & &64.5.29.79.839&4195&928 \cr
\noalign{\hrule}
 & &9.25.7.31.71.379&799&1154& & &9.25.23.59.61.71&2957&1232 \cr
8175&1313831925&4.5.17.47.379.577&29667&2548&8193&1322362575&32.3.7.11.61.2957&2485&472 \cr
 & &32.3.49.11.13.29.31&319&1456& & &512.5.49.59.71&49&256 \cr
\noalign{\hrule}
 & &3.25.49.11.13.41.61&3&536& & &121.61.139.1289&4179&4300 \cr
8176&1314338025&16.9.25.61.67&287&262&8194&1322461151&8.3.25.7.43.199.1289&147&1142 \cr
 & &64.7.41.67.131&131&2144& & &32.9.5.343.43.571&73745&82224 \cr
\noalign{\hrule}
 & &3.11.269.373.397&409&782& & &27.25.343.29.197&5219&494 \cr
8177&1314515037&4.11.17.23.269.409&2285&2016&8195&1322702325&4.49.13.17.19.307&3333&2500 \cr
 & &256.9.5.7.409.457&47985&52352& & &32.3.625.11.13.101&1313&4400 \cr
\noalign{\hrule}
 & &3.13.37.79.83.139&5047&6490& & &9.121.23.101.523&1&100 \cr
8178&1315183389&4.5.49.11.59.79.103&1369&2502&8196&1323057681&8.25.11.23.523&135&388 \cr
 & &16.9.5.1369.59.139&555&472& & &64.27.125.97&12125&96 \cr
\noalign{\hrule}
 & &7.11.17.31.71.457&1755&446& & &25.11.23.311.673&3857&3546 \cr
8179&1316666813&4.27.5.13.223.457&775&1232&8197&1323841475&4.9.25.7.19.23.29.197&91&13684 \cr
 & &128.3.125.7.11.13.31&375&832& & &32.3.49.11.13.311&49&624 \cr
\noalign{\hrule}
 & &3.5.11.19.89.4723&12121&11494& & &27.25.13.29.41.127&5687&5362 \cr
8180&1317787845&4.7.17.23.31.89.821&99&722&8198&1325051325&4.9.7.121.41.47.383&145&3302 \cr
 & &16.9.11.17.361.23.31&2139&2584& & &16.5.11.13.29.47.127&47&88 \cr
\noalign{\hrule}
 & &9.25.7.11.29.43.61&67&32& & &3.49.11.169.23.211&5085&5254 \cr
8181&1317860775&64.5.29.43.61.67&5863&2982&8199&1326193869&4.27.5.11.23.37.71.113&1085&832 \cr
 & &256.3.7.11.13.41.71&2911&1664& & &512.25.7.13.31.37.113&28675&28928 \cr
\noalign{\hrule}
 & &9.7.11.13.41.43.83&3215&3176& & &9.5.7.11.17.101.223&159&346 \cr
8182&1318277961&16.3.5.41.43.397.643&83&1846&8200&1326717315&4.27.7.53.173.223&2825&3196 \cr
 & &64.5.13.71.83.397&1985&2272& & &32.25.17.47.113.173&8131&9040 \cr
\noalign{\hrule}
 & &9.5.11.13.29.37.191&145&262& & &25.13.19.59.3643&5465&9108 \cr
8183&1318808205&4.25.841.131.191&23023&1998&8201&1327235975&8.9.125.11.23.1093&10507&14632 \cr
 & &16.27.7.11.13.23.37&23&168& & &128.3.7.19.31.59.79&1659&1984 \cr
\noalign{\hrule}
 & &13.19.67.173.461&205&666& & &27.17.37.41.1907&47585&30602 \cr
8184&1319832397&4.9.5.19.37.41.173&649&130&8202&1327849821&4.5.11.13.31.107.307&17&324 \cr
 & &16.3.25.11.13.37.59&10175&1416& & &32.81.5.13.17.107&195&1712 \cr
\noalign{\hrule}
 & &11.289.107.3881&11945&8064& & &5.11.13.31.181.331&3633&1978 \cr
8185&1320133793&256.9.5.7.17.2389&3881&3286&8203&1327927315&4.3.7.11.13.23.43.173&1163&1086 \cr
 & &1024.3.31.53.3881&1643&1536& & &16.9.23.43.181.1163&8901&9304 \cr
\noalign{\hrule}
 & &49.17.29.31.41.43&689&732& & &3.7.11.23.29.37.233&6435&1076 \cr
8186&1320252521&8.3.13.17.31.41.53.61&147&550&8204&1328297817&8.27.5.121.13.269&115&3382 \cr
 & &32.9.25.49.11.53.61&11925&10736& & &32.25.19.23.89&89&7600 \cr
\noalign{\hrule}
 & &27.11.17.23.83.137&1235&674& & &31.73.157.3739&4303&564 \cr
8187&1320480117&4.9.5.13.19.137.337&17&154&8205&1328433049&8.3.13.47.73.331&1881&1550 \cr
 & &16.5.7.11.13.17.337&455&2696& & &32.27.25.11.13.19.31&3861&7600 \cr
\noalign{\hrule}
 & &5.17.19.29.163.173&447&418& & &3.25.11.101.107.149&4277&36698 \cr
8188&1320700165&4.3.11.17.361.149.163&43353&15490&8206&1328450475&4.7.13.47.59.311&2475&298 \cr
 & &16.27.5.1549.4817&41823&38536& & &16.9.25.11.13.149&3&104 \cr
\noalign{\hrule}
 & &49.11.169.89.163&1121&3240& & &9.11.169.19.37.113&833&1314 \cr
8189&1321457137&16.81.5.11.13.19.59&83&1204&8207&1329093909&4.81.49.11.13.17.73&617&760 \cr
 & &128.9.5.7.43.83&747&13760& & &64.5.49.19.73.617&17885&19744 \cr
\noalign{\hrule}
 & &3.7.11.19.43.47.149&1011&1010& & &101.52441.251&13545&38896 \cr
8190&1321655181&4.9.5.7.11.19.101.149.337&235&25714&8208&1329431791&32.9.5.7.11.13.17.43&229&502 \cr
 & &16.25.13.23.43.47.101&2525&2392& & &128.3.5.11.229.251&165&64 \cr
\noalign{\hrule}
}%
}
$$
\eject
\vglue -23 pt
\noindent\hskip 1 in\hbox to 6.5 in{\ 8209 -- 8244 \hfill\fbd 1329456375 -- 1348656359\frb}
\vskip -9 pt
$$
\vbox{
\nointerlineskip
\halign{\strut
    \vrule \ \ \hfil \frb #\ 
   &\vrule \hfil \ \ \fbb #\frb\ 
   &\vrule \hfil \ \ \frb #\ \hfil
   &\vrule \hfil \ \ \frb #\ 
   &\vrule \hfil \ \ \frb #\ \ \vrule \hskip 2 pt
   &\vrule \ \ \hfil \frb #\ 
   &\vrule \hfil \ \ \fbb #\frb\ 
   &\vrule \hfil \ \ \frb #\ \hfil
   &\vrule \hfil \ \ \frb #\ 
   &\vrule \hfil \ \ \frb #\ \vrule \cr%
\noalign{\hrule}
 & &27.125.13.157.193&427&2468& & &9.25.49.11.13.23.37&6769&4294 \cr
8209&1329456375&8.9.25.7.61.617&19&44&8227&1341665325&4.343.19.113.967&655&312 \cr
 & &64.11.19.61.617&6787&37088& & &64.3.5.13.19.113.131&2489&3616 \cr
\noalign{\hrule}
 & &3.5.17.29.1681.107&3569&1888& & &3.5.7.11.37.89.353&295&3588 \cr
8210&1330116465&64.5.29.43.59.83&963&748&8228&1342605495&8.9.25.7.13.23.59&353&178 \cr
 & &512.9.11.17.83.107&913&768& & &32.13.23.89.353&23&208 \cr
\noalign{\hrule}
 & &25.17.79.173.229&1749&2144& & &9.5.7.13.17.101.191&29047&28826 \cr
8211&1330140775&64.3.5.11.53.67.173&1343&2208&8229&1342942965&4.3.5.49.29.31.71.937&101&836 \cr
 & &4096.9.11.17.23.79&2277&2048& & &32.11.19.29.31.71.101&14839&14384 \cr
\noalign{\hrule}
 & &3.25.7.103.24623&12349&12274& & &125.19.23.73.337&671&1008 \cr
8212&1331488725&4.7.17.361.53.103.233&675&13024&8230&1343829625&32.9.125.7.11.19.61&83&292 \cr
 & &256.27.25.11.17.19.37&11951&12672& & &256.3.7.61.73.83&5063&2688 \cr
\noalign{\hrule}
 & &3.23.41.409.1151&11917&38390& & &9.5.7.19.29.61.127&517&372 \cr
8213&1331777211&4.5.11.17.349.701&117&818&8231&1344608055&8.27.11.19.31.47.61&377&1270 \cr
 & &16.9.13.349.409&1047&104& & &32.5.11.13.29.31.127&341&208 \cr
\noalign{\hrule}
 & &9.2197.31.41.53&6295&8492& & &9.5.11.13.101.2069&18163&40922 \cr
8214&1331968599&8.5.11.41.193.1259&10881&2968&8232&1344715515&4.7.37.41.79.443&365&78 \cr
 & &128.27.5.7.13.31.53&105&64& & &16.3.5.13.37.73.79&2701&632 \cr
\noalign{\hrule}
 & &49.43.59.71.151&2431&4968& & &17.29.101.113.239&2497&780 \cr
8215&1332759673&16.27.11.13.17.23.71&35&178&8233&1344759551&8.3.5.11.13.227.239&2121&986 \cr
 & &64.9.5.7.17.23.89&3519&14240& & &32.9.7.11.17.29.101&63&176 \cr
\noalign{\hrule}
 & &9.5.13.17.29.41.113&539&46& & &5.1369.349.563&403&966 \cr
8216&1336180365&4.49.11.23.41.113&87&200&8234&1344953515&4.3.5.7.13.23.31.349&8847&13838 \cr
 & &64.3.25.7.11.23.29&161&1760& & &16.27.11.17.37.983&5049&7864 \cr
\noalign{\hrule}
 & &9.7.11.59.97.337&1835&2708& & &13.19.31.43.61.67&231&1102 \cr
8217&1336555143&8.5.337.367.677&15&352&8235&1345648837&4.3.7.11.361.29.61&155&516 \cr
 & &512.3.25.11.677&677&6400& & &32.9.5.7.29.31.43&1015&144 \cr
\noalign{\hrule}
 & &5.11.17.43.79.421&381&1724& & &5.17.1933.8191&4087&4104 \cr
8218&1337178095&8.3.11.43.127.431&21&452&8236&1345822255&16.27.5.19.61.67.1933&5797&2 \cr
 & &64.9.7.113.127&791&36576& & &64.9.11.17.31.67&737&8928 \cr
\noalign{\hrule}
 & &3.5.7.23.29.97.197&1619&1716& & &3.343.19.23.41.73&643&300 \cr
8219&1338298815&8.9.7.11.13.197.1619&7975&6596&8237&1345871289&8.9.25.19.73.643&7751&4466 \cr
 & &64.25.121.13.17.29.97&2057&2080& & &32.5.7.11.23.29.337&3707&2320 \cr
\noalign{\hrule}
 & &27.7.13.23.137.173&2105&2242& & &1331.17.19.31.101&13969&39258 \cr
8220&1339367211&4.5.13.19.59.173.421&43703&29130&8238&1346057603&4.27.61.229.727&119&110 \cr
 & &16.3.25.11.29.137.971&7975&7768& & &16.3.5.7.11.17.61.727&6405&5816 \cr
\noalign{\hrule}
 & &5.7.19.23.67.1307&16269&8564& & &9.25.7.13.17.53.73&253&1072 \cr
8221&1339367855&8.3.7.11.17.29.2141&13501&10050&8239&1346702175&32.11.17.23.67.73&1391&150 \cr
 & &32.9.25.23.67.587&587&720& & &128.3.25.11.13.107&1177&64 \cr
\noalign{\hrule}
 & &5.7.11.13.41.61.107&435&358& & &5.7.11.13.17.71.223&423&1138 \cr
8222&1339373035&4.3.25.29.41.107.179&4431&44&8240&1347150805&4.9.17.47.71.569&169&9842 \cr
 & &32.9.7.11.29.211&211&4176& & &16.3.7.169.19.37&57&3848 \cr
\noalign{\hrule}
 & &243.5.7.11.103.139&611&604& & &5.49.13.23.53.347&4309&3672 \cr
8223&1339426935&8.11.13.47.103.139.151&17955&3638&8241&1347232705&16.27.5.17.31.53.139&1757&2552 \cr
 & &32.27.5.7.17.19.47.107&5029&5168& & &256.9.7.11.17.29.251&38403&40832 \cr
\noalign{\hrule}
 & &5.13.1039.19843&9889&9954& & &11.19.37.61.2857&39555&14728 \cr
8224&1340097005&4.9.7.11.29.31.79.1039&413&2704&8242&1347684041&16.27.5.7.263.293&1037&1330 \cr
 & &128.3.49.11.169.31.59&50127&49088& & &64.3.25.49.17.19.61&1225&1632 \cr
\noalign{\hrule}
 & &5.11.31.43.101.181&2569&3042& & &23.53.683.1619&43245&42562 \cr
8225&1340271515&4.9.5.7.169.101.367&1843&7348&8243&1347942163&4.9.5.13.23.961.1637&37565&86 \cr
 & &32.3.11.13.19.97.167&16199&11856& & &16.3.25.11.43.683&1419&200 \cr
\noalign{\hrule}
 & &3.5.7.2053.6221&4075&10296& & &289.23.137.1481&885&596 \cr
8226&1341029865&16.27.125.11.13.163&791&2584&8244&1348656359&8.3.5.23.59.137.149&51711&50354 \cr
 & &256.7.13.17.19.113&4199&14464& & &32.9.11.17.1481.1567&1567&1584 \cr
\noalign{\hrule}
}%
}
$$
\eject
\vglue -23 pt
\noindent\hskip 1 in\hbox to 6.5 in{\ 8245 -- 8280 \hfill\fbd 1348961705 -- 1366272875\frb}
\vskip -9 pt
$$
\vbox{
\nointerlineskip
\halign{\strut
    \vrule \ \ \hfil \frb #\ 
   &\vrule \hfil \ \ \fbb #\frb\ 
   &\vrule \hfil \ \ \frb #\ \hfil
   &\vrule \hfil \ \ \frb #\ 
   &\vrule \hfil \ \ \frb #\ \ \vrule \hskip 2 pt
   &\vrule \ \ \hfil \frb #\ 
   &\vrule \hfil \ \ \fbb #\frb\ 
   &\vrule \hfil \ \ \frb #\ \hfil
   &\vrule \hfil \ \ \frb #\ 
   &\vrule \hfil \ \ \frb #\ \vrule \cr%
\noalign{\hrule}
 & &5.7.13.41.167.433&1917&748& & &625.31.43.1627&891&736 \cr
8245&1348961705&8.27.11.17.71.433&167&266&8263&1355494375&64.81.125.11.23.43&1057&68 \cr
 & &32.3.7.17.19.71.167&1349&816& & &512.9.7.11.17.151&11781&38656 \cr
\noalign{\hrule}
 & &27.7.11.13.19.37.71&215&566& & &9.25.19.31.53.193&1595&2072 \cr
8246&1348998651&4.5.7.19.37.43.283&639&1342&8264&1355598225&16.125.7.11.29.31.37&5459&834 \cr
 & &16.9.5.11.43.61.71&215&488& & &64.3.11.53.103.139&1529&3296 \cr
\noalign{\hrule}
 & &27.5.11.17.19.29.97&9485&118& & &9.49.13.29.41.199&209&168 \cr
8247&1349269515&4.3.25.7.59.271&319&494&8265&1356490863&16.27.343.11.19.199&5945&572 \cr
 & &16.11.13.19.29.59&13&472& & &128.5.121.13.29.41&121&320 \cr
\noalign{\hrule}
 & &3.5.49.11.13.37.347&593&628& & &3.11.13.17.43.61.71&635&1558 \cr
8248&1349443095&8.7.13.157.347.593&4331&180&8266&1358197269&4.5.11.19.41.61.127&119&1278 \cr
 & &64.9.5.61.71.157&11147&5856& & &16.9.5.7.17.41.71&615&56 \cr
\noalign{\hrule}
 & &3.169.37.193.373&205&374& & &5.11.169.41.43.83&68517&77812 \cr
8249&1350441651&4.5.11.17.37.41.373&1737&2366&8267&1360128055&8.9.49.23.331.397&7975&362 \cr
 & &16.9.5.7.169.41.193&205&168& & &32.3.25.7.11.29.181&5249&1680 \cr
\noalign{\hrule}
 & &5.29.61.107.1427&171&1598& & &5.49.13.47.61.149&14877&5192 \cr
8250&1350534205&4.9.5.17.19.47.107&667&132&8268&1360577855&16.27.7.11.19.29.59&445&676 \cr
 & &32.27.11.19.23.29&6831&304& & &128.9.5.169.29.89&3393&5696 \cr
\noalign{\hrule}
 & &243.25.11.17.29.41&13&1202& & &27.47.983.1091&37829&8372 \cr
8251&1350733725&4.5.11.13.17.601&91&96&8269&1360942857&8.7.11.13.19.23.181&1145&846 \cr
 & &256.3.7.169.601&4207&21632& & &32.9.5.7.19.47.229&1145&2128 \cr
\noalign{\hrule}
 & &27.5.11.17.109.491&3627&1774& & &5.7.11.83.191.223&589&366 \cr
8252&1351087155&4.243.5.13.31.887&199&44&8270&1361059315&4.3.7.11.19.31.61.83&3725&1338 \cr
 & &32.11.13.199.887&11531&3184& & &16.9.25.19.149.223&1341&760 \cr
\noalign{\hrule}
 & &7.17.67.97.1747&3113&4860& & &81.19.47.67.281&5395&13432 \cr
8253&1351096607&8.243.5.11.97.283&431&442&8271&1361813391&16.9.5.13.23.73.83&517&562 \cr
 & &32.27.5.13.17.283.431&58185&58864& & &64.11.23.47.73.281&803&736 \cr
\noalign{\hrule}
 & &23.29.41.73.677&435&508& & &5.7.13.61.139.353&1369&396 \cr
8254&1351516087&8.3.5.841.127.677&759&82&8272&1361854585&8.9.11.13.1369.61&1007&214 \cr
 & &32.9.5.11.23.41.127&1397&720& & &32.3.19.37.53.107&37259&5136 \cr
\noalign{\hrule}
 & &3.5.121.23.139.233&29127&29822& & &3.5.11.17.19.37.691&12679&12888 \cr
8255&1351995315&4.9.7.11.13.19.31.37.73&3151&2150&8273&1362593265&16.27.5.17.31.179.409&13129&2086 \cr
 & &16.25.23.37.43.73.137&25345&25112& & &64.7.19.31.149.691&1043&992 \cr
\noalign{\hrule}
 & &9.25.19.37.83.103&9499&286& & &27.5.11.29.31.1021&221&562 \cr
8256&1352238075&4.3.5.7.11.13.23.59&19&4&8274&1363050315&4.5.13.17.281.1021&4379&726 \cr
 & &32.7.11.13.19.59&8437&112& & &16.3.121.17.29.151&187&1208 \cr
\noalign{\hrule}
 & &11.71.229.7561&93205&85644& & &121.361.529.59&7721&582 \cr
8257&1352277289&8.27.5.7.13.61.2663&2521&142&8275&1363327691&4.3.7.23.97.1103&7257&8360 \cr
 & &32.9.5.7.71.2521&2521&5040& & &64.9.5.11.19.41.59&369&160 \cr
\noalign{\hrule}
 & &9.7.47.443.1031&42887&22066& & &1681.199.4079&6119&2040 \cr
8258&1352386413&4.11.13.17.59.3299&1325&1974&8276&1364503001&16.3.5.17.29.41.211&1199&2388 \cr
 & &16.3.25.7.13.17.47.53&1325&1768& & &128.9.5.11.109.199&1199&2880 \cr
\noalign{\hrule}
 & &5.49.13.47.83.109&22847&26748& & &5.343.23.53.653&46637&28458 \cr
8259&1354290665&8.9.7.11.31.67.743&1075&332&8277&1365152005&4.27.17.31.149.313&901&3718 \cr
 & &64.3.25.11.31.43.83&2365&2976& & &16.3.11.169.289.53&5577&2312 \cr
\noalign{\hrule}
 & &3.25.11.13.289.19.23&937&508& & &121.13.31.37.757&10377&38386 \cr
8260&1354492425&8.5.19.23.127.937&43047&45968&8278&1365802867&4.9.17.1129.1153&1165&2294 \cr
 & &256.9.169.17.4783&4783&4992& & &16.3.5.17.31.37.233&1165&408 \cr
\noalign{\hrule}
 & &9.13.23.31.37.439&2989&3428& & &5.13.23.31.41.719&11165&11124 \cr
8261&1355007303&8.49.13.37.61.857&13079&18630&8279&1366204255&8.27.25.7.11.13.23.29.103&56703&38678 \cr
 & &32.81.5.7.11.23.29.41&8323&7920& & &32.81.41.83.233.461&107413&107568 \cr
\noalign{\hrule}
 & &25.19.47.101.601&38111&9864& & &125.11.157.6329&3227&3102 \cr
8262&1355149825&16.9.23.137.1657&5555&3898&8280&1366272875&4.3.7.121.47.157.461&1335&20332 \cr
 & &64.3.5.11.101.1949&1949&1056& & &32.9.5.7.13.17.23.89&18423&24752 \cr
\noalign{\hrule}
}%
}
$$
\eject
\vglue -23 pt
\noindent\hskip 1 in\hbox to 6.5 in{\ 8281 -- 8316 \hfill\fbd 1366433013 -- 1379325311\frb}
\vskip -9 pt
$$
\vbox{
\nointerlineskip
\halign{\strut
    \vrule \ \ \hfil \frb #\ 
   &\vrule \hfil \ \ \fbb #\frb\ 
   &\vrule \hfil \ \ \frb #\ \hfil
   &\vrule \hfil \ \ \frb #\ 
   &\vrule \hfil \ \ \frb #\ \ \vrule \hskip 2 pt
   &\vrule \ \ \hfil \frb #\ 
   &\vrule \hfil \ \ \fbb #\frb\ 
   &\vrule \hfil \ \ \frb #\ \hfil
   &\vrule \hfil \ \ \frb #\ 
   &\vrule \hfil \ \ \frb #\ \vrule \cr%
\noalign{\hrule}
 & &3.11.361.23.4987&3463&8450& & &5.17.101.277.577&74847&65038 \cr
8281&1366433013&4.25.169.23.3463&1881&1582&8299&1372131965&4.3.31.61.409.1049&189&220 \cr
 & &16.9.25.7.11.13.19.113&2275&2712& & &32.81.5.7.11.61.1049&80773&79056 \cr
\noalign{\hrule}
 & &3.5.11.73.233.487&13&2422& & &9.25.7.13.199.337&53&1738 \cr
8282&1366758195&4.7.13.173.233&1241&1008&8300&1373114925&4.5.7.11.13.53.79&199&186 \cr
 & &128.9.49.17.73&2499&64& & &16.3.31.53.79.199&1643&632 \cr
\noalign{\hrule}
 & &3.5.11.19.271.1609&183&88& & &3.5.49.47.83.479&11&94 \cr
8283&1366982265&16.9.121.61.1609&1349&260&8301&1373405565&4.7.11.2209.479&2781&572 \cr
 & &128.5.13.19.61.71&923&3904& & &32.27.121.13.103&1089&21424 \cr
\noalign{\hrule}
 & &3.5.13.23.59.5167&6237&11404& & &11.13.73.149.883&1005&932 \cr
8284&1367265705&8.243.5.7.11.2851&2033&818&8302&1373427913&8.3.5.11.67.233.883&73&810 \cr
 & &32.7.11.19.107.409&43763&23408& & &32.243.25.73.233&5825&3888 \cr
\noalign{\hrule}
 & &9.5.11.13.41.71.73&437&366& & &3.25.121.17.29.307&987&2522 \cr
8285&1367456805&4.27.5.13.19.23.41.61&22649&49306&8303&1373510325&4.9.5.7.13.17.47.97&989&116 \cr
 & &16.11.29.71.89.277&2581&2216& & &32.7.23.29.43.47&6923&752 \cr
\noalign{\hrule}
 & &5.13.109.137.1409&1595&186& & &125.11.13.31.37.67&369&34 \cr
8286&1367638805&4.3.25.11.29.31.109&333&442&8304&1373675875&4.9.25.11.17.37.41&1081&806 \cr
 & &16.27.11.13.17.29.37&10989&3944& & &16.3.13.23.31.41.47&1927&552 \cr
\noalign{\hrule}
 & &9.11.37.397.941&39235&53924& & &27.125.49.53.157&269&1144 \cr
8287&1368412551&8.5.7.13.17.19.59.61&37&24&8305&1376085375&16.3.7.11.13.53.269&29&400 \cr
 & &128.3.5.7.17.19.37.59&7021&6080& & &512.25.29.269&269&7424 \cr
\noalign{\hrule}
 & &25.7.13.293.2053&12331&38994& & &5.11.13.17.19.59.101&279&488 \cr
8288&1368478475&4.3.11.19.59.67.97&685&588&8306&1376201255&16.9.5.17.31.61.101&1021&494 \cr
 & &32.9.5.49.11.59.137&10549&8496& & &64.3.13.19.61.1021&3063&1952 \cr
\noalign{\hrule}
 & &11.23.61.131.677&36627&52060& & &25.11.17.271.1087&6191&5766 \cr
8289&1368706471&8.3.5.19.29.137.421&2013&5986&8307&1377147475&4.3.961.41.151.271&1105&7296 \cr
 & &32.9.5.11.41.61.73&1845&1168& & &1024.9.5.13.17.19.31&5301&6656 \cr
\noalign{\hrule}
 & &9.125.11.17.23.283&3093&1718& & &3.11.19.23.29.37.89&1337&710 \cr
8290&1369330875&4.27.23.859.1031&119&740&8308&1377162237&4.5.7.29.37.71.191&785&288 \cr
 & &32.5.7.17.37.1031&1031&4144& & &256.9.25.157.191&11775&24448 \cr
\noalign{\hrule}
 & &5.7.13.17.23.43.179&1179&14036& & &3.5.7.11.17.29.41.59&101&130 \cr
8291&1369334785&8.9.7.121.29.131&731&710&8309&1377414885&4.25.13.17.41.59.101&4149&36974 \cr
 & &32.3.5.11.17.29.43.71&957&1136& & &16.9.7.19.139.461&8759&3336 \cr
\noalign{\hrule}
 & &7.13.43.269.1301&801&500& & &5.13.109.113.1721&1595&126 \cr
8292&1369428697&8.9.125.13.89.269&297&28&8310&1377841205&4.9.25.7.11.29.109&317&208 \cr
 & &64.243.5.7.11.89&2673&14240& & &128.3.11.13.29.317&10461&1856 \cr
\noalign{\hrule}
 & &5.23.53.379.593&7009&6630& & &9.11.13.19.157.359&12245&39092 \cr
8293&1369832965&4.3.25.13.17.43.53.163&45859&7116&8311&1378244439&8.5.29.31.79.337&1161&1130 \cr
 & &32.9.121.379.593&121&144& & &32.27.25.43.113.337&43473&45200 \cr
\noalign{\hrule}
 & &3.49.11.31.151.181&195&146& & &9.11.37.41.67.137&3029&1522 \cr
8294&1370021037&4.9.5.13.73.151.181&589&770&8312&1378529757&4.3.13.67.233.761&481&280 \cr
 & &16.25.7.11.13.19.31.73&1825&1976& & &64.5.7.169.37.233&8155&5408 \cr
\noalign{\hrule}
 & &27.5.1331.29.263&2623&1292& & &5.7.11.19.29.67.97&1775&1038 \cr
8295&1370457495&8.17.19.43.61.263&497&234&8313&1378665365&4.3.125.7.19.71.173&1793&582 \cr
 & &32.9.7.13.19.61.71&8113&14768& & &16.9.11.71.97.163&639&1304 \cr
\noalign{\hrule}
 & &9.49.11.13.103.211&6005&958& & &3.89.97.139.383&4125&4508 \cr
8296&1370548179&4.3.5.13.479.1201&8525&7088&8314&1378785063&8.9.125.49.11.23.139&383&590 \cr
 & &128.125.11.31.443&13733&8000& & &32.625.7.11.59.383&6875&6608 \cr
\noalign{\hrule}
 & &9.5.11.13.23.59.157&4681&4582& & &25.13.131.139.233&10619&7590 \cr
8297&1370970315&4.5.13.23.29.31.79.151&1093&10362&8315&1378876525&4.3.125.7.11.23.37.41&1113&262 \cr
 & &16.3.11.151.157.1093&1093&1208& & &16.9.49.41.53.131&2009&3816 \cr
\noalign{\hrule}
 & &49.11.13.23.67.127&54125&39474& & &7.13.17.19.167.281&725&444 \cr
8298&1371318949&4.27.125.17.43.433&889&1276&8316&1379325311&8.3.25.13.17.19.29.37&583&4782 \cr
 & &32.3.25.7.11.17.29.127&725&816& & &32.9.5.11.53.797&8767&38160 \cr
\noalign{\hrule}
}%
}
$$
\eject
\vglue -23 pt
\noindent\hskip 1 in\hbox to 6.5 in{\ 8317 -- 8352 \hfill\fbd 1379976165 -- 1395689191\frb}
\vskip -9 pt
$$
\vbox{
\nointerlineskip
\halign{\strut
    \vrule \ \ \hfil \frb #\ 
   &\vrule \hfil \ \ \fbb #\frb\ 
   &\vrule \hfil \ \ \frb #\ \hfil
   &\vrule \hfil \ \ \frb #\ 
   &\vrule \hfil \ \ \frb #\ \ \vrule \hskip 2 pt
   &\vrule \ \ \hfil \frb #\ 
   &\vrule \hfil \ \ \fbb #\frb\ 
   &\vrule \hfil \ \ \frb #\ \hfil
   &\vrule \hfil \ \ \frb #\ 
   &\vrule \hfil \ \ \frb #\ \vrule \cr%
\noalign{\hrule}
 & &9.5.29.47.149.151&4181&2822& & &9.13.19.31.41.491&405&374 \cr
8317&1379976165&4.5.17.29.37.83.113&3599&6006&8335&1387287603&4.729.5.11.13.17.491&3511&5966 \cr
 & &16.3.7.11.13.37.59.61&29341&36344& & &16.11.17.19.157.3511&29359&28088 \cr
\noalign{\hrule}
 & &27.7.23.29.47.233&3529&14480& & &23.31.37.41.1283&1215&68 \cr
8318&1380515913&32.5.7.181.3529&1131&2398&8336&1387719743&8.243.5.17.23.41&829&1034 \cr
 & &128.3.5.11.13.29.109&1199&4160& & &32.3.11.17.47.829&9119&38352 \cr
\noalign{\hrule}
 & &9.7.11.19.23.47.97&1819&1800& & &3.25.13.17.31.37.73&199&726 \cr
8319&1380652119&16.81.25.17.23.97.107&41783&10112&8337&1387841325&4.9.121.13.73.199&7225&16058 \cr
 & &4096.5.7.47.79.127&10033&10240& & &16.25.7.289.31.37&17&56 \cr
\noalign{\hrule}
 & &9.17.19.23.107.193&19085&22752& & &3.5.7.19.23.79.383&15531&14726 \cr
8320&1380746511&64.81.5.11.79.347&14089&17906&8338&1388342445&4.9.19.31.37.167.199&86411&30800 \cr
 & &256.7.73.193.1279&8953&9344& & &128.25.7.11.13.289.23&3757&3520 \cr
\noalign{\hrule}
 & &27.7.11.13.29.41.43&125&658& & &3.125.19.23.37.229&959&7514 \cr
8321&1381809429&4.125.49.11.43.47&377&162&8339&1388512875&4.25.7.13.289.137&23&198 \cr
 & &16.81.25.13.29.47&75&376& & &16.9.11.17.23.137&137&4488 \cr
\noalign{\hrule}
 & &9.19.47.373.461&13145&4386& & &3.5.7.43.457.673&319&354 \cr
8322&1381986261&4.27.5.11.17.43.239&461&700&8340&1388637915&4.9.11.29.43.59.457&3325&788 \cr
 & &32.125.7.11.17.461&2125&1232& & &32.25.7.11.19.29.197&10835&8816 \cr
\noalign{\hrule}
 & &3.5.11.13.59.67.163&24411&23674& & &9.13.17.809.863&1969&9250 \cr
8323&1382107155&4.9.7.13.19.79.89.103&977&50&8341&1388654163&4.125.11.17.37.179&2427&2248 \cr
 & &16.25.7.19.89.977&18563&24920& & &64.3.5.37.281.809&1405&1184 \cr
\noalign{\hrule}
 & &11.169.17.67.653&9253&2070& & &7.17.43.47.53.109&136125&135394 \cr
8324&1382662853&4.9.5.17.19.23.487&201&286&8342&1389362723&4.9.125.49.121.19.509&53&878 \cr
 & &16.27.11.13.19.23.67&437&216& & &16.3.5.11.53.439.509&16797&17560 \cr
\noalign{\hrule}
 & &3.5.11.17.59.61.137&4503&3854& & &11.13.17.41.73.191&9135&1304 \cr
8325&1383041715&4.9.5.17.19.41.47.79&121&274&8343&1389712753&16.9.5.7.17.29.163&73&90 \cr
 & &16.121.19.41.47.137&1927&1672& & &64.81.25.7.29.73&5075&2592 \cr
\noalign{\hrule}
 & &5.11.29.59.61.241&46163&27342& & &3.5.49.19.29.3433&1209&2224 \cr
8326&1383437605&4.9.49.13.31.53.67&1525&118&8344&1390313505&32.9.7.13.19.31.139&8207&550 \cr
 & &16.3.25.7.13.59.61&39&280& & &128.25.11.29.283&3113&320 \cr
\noalign{\hrule}
 & &5.11.37.439.1549&50719&34476& & &31.47.337.2833&51831&35992 \cr
8327&1383822385&8.3.169.17.67.757&57&814&8345&1391028497&16.9.11.13.409.443&1585&2914 \cr
 & &32.9.11.13.17.19.37&323&1872& & &64.3.5.13.31.47.317&1585&1248 \cr
\noalign{\hrule}
 & &5.49.11.17.19.37.43&4453&2862& & &81.11.47.167.199&985&2822 \cr
8328&1384941635&4.27.7.17.53.61.73&1075&38&8346&1391698341&4.5.17.83.197.199&141&58 \cr
 & &16.9.25.19.43.73&45&584& & &16.3.5.17.29.47.197&3349&1160 \cr
\noalign{\hrule}
 & &5.11.61.293.1409&29811&47684& & &9.11.169.19.29.151&1591&70 \cr
8329&1385068135&8.3.7.13.19.131.523&3075&586&8347&1392035931&4.5.7.19.29.37.43&297&254 \cr
 & &32.9.25.13.41.293&369&1040& & &16.27.5.7.11.37.127&3885&1016 \cr
\noalign{\hrule}
 & &3.5.7.13.41.53.467&253&214& & &27.5.19.29.41.457&37759&5656 \cr
8330&1385189715&4.5.7.11.23.41.53.107&13897&16812&8348&1393751745&16.7.61.101.619&5247&914 \cr
 & &32.9.13.23.467.1069&1069&1104& & &64.9.11.53.457&583&32 \cr
\noalign{\hrule}
 & &3.5.7.11.17.19.47.79&1565&696& & &9.11.13.43.89.283&145&704 \cr
8331&1385190345&16.9.25.29.47.313&749&9826&8349&1393873767&128.3.5.121.29.89&26273&26362 \cr
 & &64.7.4913.107&289&3424& & &512.49.13.43.47.269&12643&12544 \cr
\noalign{\hrule}
 & &9.25.13.19.97.257&17&308& & &9.7.11.17.19.23.271&21463&20000 \cr
8332&1385429175&8.3.7.11.17.19.257&185&442&8350&1395188487&64.625.169.23.127&2261&1626 \cr
 & &32.5.7.13.289.37&259&4624& & &256.3.125.7.17.19.271&125&128 \cr
\noalign{\hrule}
 & &5.11.17.841.41.43&487&702& & &27.23.53.109.389&4235&4712 \cr
8333&1386308605&4.27.11.13.17.29.487&137&50&8351&1395544113&16.3.5.7.121.19.31.109&239&784 \cr
 & &16.9.25.13.137.487&21915&14248& & &512.343.11.19.239&71687&61184 \cr
\noalign{\hrule}
 & &3.5.11.17.19.53.491&291&9620& & &13.59.751.2423&1595&828 \cr
8334&1386895785&8.9.25.13.37.97&5687&5662&8352&1395689191&8.9.5.11.23.29.751&103&854 \cr
 & &32.121.19.37.47.149&7003&6512& & &32.3.5.7.23.61.103&6405&37904 \cr
\noalign{\hrule}
}%
}
$$
\eject
\vglue -23 pt
\noindent\hskip 1 in\hbox to 6.5 in{\ 8353 -- 8388 \hfill\fbd 1396392855 -- 1414490805\frb}
\vskip -9 pt
$$
\vbox{
\nointerlineskip
\halign{\strut
    \vrule \ \ \hfil \frb #\ 
   &\vrule \hfil \ \ \fbb #\frb\ 
   &\vrule \hfil \ \ \frb #\ \hfil
   &\vrule \hfil \ \ \frb #\ 
   &\vrule \hfil \ \ \frb #\ \ \vrule \hskip 2 pt
   &\vrule \ \ \hfil \frb #\ 
   &\vrule \hfil \ \ \fbb #\frb\ 
   &\vrule \hfil \ \ \frb #\ \hfil
   &\vrule \hfil \ \ \frb #\ 
   &\vrule \hfil \ \ \frb #\ \vrule \cr%
\noalign{\hrule}
 & &3.5.11.13.53.71.173&6769&4520& & &3.11.13.67.131.373&2793&2056 \cr
8353&1396392855&16.25.7.11.113.967&621&346&8371&1404469209&16.9.49.19.131.257&1055&124 \cr
 & &64.27.7.23.113.173&2599&2016& & &128.5.31.211.257&39835&13504 \cr
\noalign{\hrule}
 & &3.5.11.89.251.379&4223&1462& & &7.11.41.43.79.131&455&414 \cr
8354&1396969365&4.17.41.43.89.103&2971&1458&8372&1404887099&4.9.5.49.13.23.43.131&1501&3608 \cr
 & &16.729.41.2971&9963&23768& & &64.3.5.11.19.23.41.79&437&480 \cr
\noalign{\hrule}
 & &81.11.19.23.37.97&1885&346& & &125.149.197.383&28413&47038 \cr
8355&1397438163&4.5.11.13.29.37.173&1731&3634&8373&1405274875&4.9.7.11.29.41.811&31&10 \cr
 & &16.3.13.23.79.577&1027&4616& & &16.3.5.11.29.31.811&9889&19464 \cr
\noalign{\hrule}
 & &3.5.7.11.289.53.79&13065&16244& & &27.7.13.61.83.113&3685&3868 \cr
8356&1397599665&8.9.25.13.31.67.131&901&7876&8374&1405696383&8.9.5.11.67.113.967&61&52 \cr
 & &64.11.13.17.53.179&179&416& & &64.5.11.13.61.67.967&10637&10720 \cr
\noalign{\hrule}
 & &3.25.7.11.31.37.211&3211&2064& & &3.5.49.13.37.41.97&99&580 \cr
8357&1397648175&32.9.7.11.169.19.43&26375&24494&8375&1406008695&8.27.25.7.11.29.41&481&194 \cr
 & &128.125.37.211.331&331&320& & &32.11.13.29.37.97&29&176 \cr
\noalign{\hrule}
 & &625.7.121.19.139&35391&51484& & &121.19.41.43.347&329&450 \cr
8358&1398079375&8.3.47.61.211.251&231&20&8376&1406438539&4.9.25.7.43.47.347&1273&232 \cr
 & &64.9.5.7.11.47.61&423&1952& & &64.3.5.19.29.47.67&9715&4512 \cr
\noalign{\hrule}
 & &3.49.13.29.43.587&9565&7458& & &11.37.53.113.577&333&910 \cr
8359&1398830979&4.9.5.11.13.113.1913&587&700&8377&1406450771&4.9.5.7.13.1369.53&605&764 \cr
 & &32.125.7.587.1913&1913&2000& & &32.3.25.7.121.13.191&17381&13200 \cr
\noalign{\hrule}
 & &3.25.13.23.89.701&63&638& & &3.25.11.13.179.733&29233&47558 \cr
8360&1399073325&4.27.7.11.13.29.89&701&1702&8378&1407195075&4.7.23.31.41.43.79&773&990 \cr
 & &16.23.29.37.701&29&296& & &16.9.5.11.23.79.773&5451&6184 \cr
\noalign{\hrule}
 & &3.5.13.23.29.47.229&9717&1046& & &9.11.19.29.131.197&865&314 \cr
8361&1399889595&4.9.5.41.79.523&253&458&8379&1407746043&4.5.11.157.173.197&459&1444 \cr
 & &16.11.23.229.523&523&88& & &32.27.17.361.157&2983&816 \cr
\noalign{\hrule}
 & &3.5.7.13.17.23.43.61&1177&1238& & &9.25.7.1369.653&1159&506 \cr
8362&1399934445&4.11.13.17.43.107.619&1485&25132&8380&1407982275&4.5.7.11.19.23.37.61&1231&804 \cr
 & &32.27.5.121.61.103&1089&1648& & &32.3.19.23.67.1231&29279&19696 \cr
\noalign{\hrule}
 & &5.13.149.163.887&1503&616& & &5.11.169.31.67.73&2001&2014 \cr
8363&1400266985&16.9.5.7.11.149.167&731&1774&8381&1409317195&4.3.13.19.23.29.31.53.67&15549&32098 \cr
 & &64.3.11.17.43.887&561&1376& & &16.9.11.23.71.73.1459&13131&13064 \cr
\noalign{\hrule}
 & &49.11.13.17.19.619&403&216& & &5.7.13.31.139.719&2241&7106 \cr
8364&1400958559&16.27.49.169.19.31&6809&1570&8382&1409667805&4.27.11.17.19.31.83&26875&25298 \cr
 & &64.3.5.11.157.619&157&480& & &16.3.625.7.13.43.139&375&344 \cr
\noalign{\hrule}
 & &9.7.19.67.101.173&689&220& & &5.7.11.13.29.71.137&221&276 \cr
8365&1401317127&8.5.11.13.19.53.173&39&134&8383&1411825415&8.3.169.17.23.29.137&875&4026 \cr
 & &32.3.11.169.53.67&583&2704& & &32.9.125.7.11.17.61&1525&2448 \cr
\noalign{\hrule}
 & &9.5.7.11.19.107.199&1955&4144& & &9.5.11.19.29.31.167&6839&1996 \cr
8366&1401824655&32.3.25.49.17.23.37&361&214&8384&1412000865&8.3.7.11.499.977&20677&17746 \cr
 & &128.17.361.37.107&703&1088& & &32.19.23.29.31.467&467&368 \cr
\noalign{\hrule}
 & &3.5.7.17.29.41.661&77&128& & &81.25.7.13.79.97&373&22 \cr
8367&1402883265&256.49.11.29.661&2925&4346&8385&1412099325&4.3.5.7.11.97.373&481&586 \cr
 & &1024.9.25.13.41.53&2067&2560& & &16.13.37.293.373&10841&2984 \cr
\noalign{\hrule}
 & &25.11.37.173.797&2561&1764& & &5.11.41.773.811&161&612 \cr
8368&1402939175&8.9.49.11.13.37.197&79&1300&8386&1413666265&8.9.5.7.17.23.811&533&278 \cr
 & &64.3.25.7.169.79&1183&7584& & &32.3.7.13.23.41.139&6279&2224 \cr
\noalign{\hrule}
 & &7.121.43.59.653&1053&37474& & &27.7.13.109.5279&498157&499574 \cr
8369&1403191867&4.81.13.41.457&215&242&8387&1413784827&4.121.23.37.43.157.179&7925&13734 \cr
 & &16.3.5.121.13.41.43&533&120& & &16.9.25.7.23.43.109.317&7925&7912 \cr
\noalign{\hrule}
 & &5.11.13.41.83.577&187&228& & &9.5.7.13.29.43.277&5371&2662 \cr
8370&1403924665&8.3.121.13.17.19.577&1075&498&8388&1414490805&4.5.1331.13.41.131&667&1998 \cr
 & &32.9.25.17.19.43.83&2907&3440& & &16.27.23.29.37.131&2553&1048 \cr
\noalign{\hrule}
}%
}
$$
\eject
\vglue -23 pt
\noindent\hskip 1 in\hbox to 6.5 in{\ 8389 -- 8424 \hfill\fbd 1414665945 -- 1430037675\frb}
\vskip -9 pt
$$
\vbox{
\nointerlineskip
\halign{\strut
    \vrule \ \ \hfil \frb #\ 
   &\vrule \hfil \ \ \fbb #\frb\ 
   &\vrule \hfil \ \ \frb #\ \hfil
   &\vrule \hfil \ \ \frb #\ 
   &\vrule \hfil \ \ \frb #\ \ \vrule \hskip 2 pt
   &\vrule \ \ \hfil \frb #\ 
   &\vrule \hfil \ \ \fbb #\frb\ 
   &\vrule \hfil \ \ \frb #\ \hfil
   &\vrule \hfil \ \ \frb #\ 
   &\vrule \hfil \ \ \frb #\ \vrule \cr%
\noalign{\hrule}
 & &27.5.7.11.23.61.97&847&556& & &5.7.17.23.73.1423&11913&12278 \cr
8389&1414665945&8.9.5.49.1331.139&437&1768&8407&1421584115&4.3.49.11.361.23.877&1241&18930 \cr
 & &128.13.17.19.23.139&4199&8896& & &16.9.5.11.17.73.631&631&792 \cr
\noalign{\hrule}
 & &27.5.11.529.1801&2977&2842& & &11.31.71.151.389&1039&3240 \cr
8390&1414802565&4.49.13.29.229.1801&99&1702&8408&1422129929&16.81.5.151.1039&293&746 \cr
 & &16.9.7.11.13.23.29.37&1073&728& & &64.27.5.293.373&39555&11936 \cr
\noalign{\hrule}
 & &243.5.7.31.41.131&6059&1474& & &11.13.139.163.439&63&76 \cr
8391&1416091005&4.11.41.67.73.83&13&54&8409&1422338489&8.9.7.11.19.163.439&1555&238 \cr
 & &16.27.11.13.73.83&949&7304& & &32.3.5.49.17.19.311&45717&25840 \cr
\noalign{\hrule}
 & &5.7.13.17.19.23.419&33&128& & &11.13.17.47.59.211&5825&6624 \cr
8392&1416301705&256.3.11.13.17.419&1425&1006&8410&1422385393&64.9.25.11.13.23.233&47&8 \cr
 & &1024.9.25.19.503&2515&4608& & &1024.3.5.23.47.233&5359&7680 \cr
\noalign{\hrule}
 & &729.5.343.11.103&3277&496& & &9.5.7.11.17.19.31.41&1807&454 \cr
8393&1416516255&32.27.5.29.31.113&3709&206&8411&1422496845&4.3.5.13.31.139.227&3857&452 \cr
 & &128.103.3709&3709&64& & &32.7.13.19.29.113&3277&208 \cr
\noalign{\hrule}
 & &3.25.7.13.41.61.83&847&232& & &3.5.29.617.5303&52687&36778 \cr
8394&1416753975&16.5.49.121.29.61&117&422&8412&1423298685&4.7.19.37.47.59.71&725&396 \cr
 & &64.9.11.13.29.211&6119&1056& & &32.9.25.11.29.37.71&2343&2960 \cr
\noalign{\hrule}
 & &3.121.13.467.643&61&60& & &9.29.43.89.1427&6325&4898 \cr
8395&1417026039&8.9.5.13.61.467.643&5929&142&8413&1425354669&4.25.11.23.31.79.89&5633&6612 \cr
 & &32.5.49.121.61.71&3479&4880& & &32.3.5.19.23.29.43.131&2185&2096 \cr
\noalign{\hrule}
 & &27.125.49.11.19.41&277&236& & &3.25.49.23.101.167&2651&1524 \cr
8396&1417098375&8.125.49.11.59.277&1539&4586&8414&1425683175&8.9.11.101.127.241&575&334 \cr
 & &32.81.19.59.2293&2293&2832& & &32.25.11.23.127.167&127&176 \cr
\noalign{\hrule}
 & &81.11.13.53.2309&6601&14180& & &5.37.43.277.647&14049&13772 \cr
8397&1417492791&8.9.5.7.23.41.709&265&104&8415&1425687145&8.9.5.7.11.37.223.313&23&208 \cr
 & &128.25.13.53.709&709&1600& & &256.3.13.23.223.313&93587&85632 \cr
\noalign{\hrule}
 & &9.11.43.53.61.103&2635&7064& & &81.5.11.13.103.239&1405&1702 \cr
8398&1417576743&16.3.5.11.17.31.883&4223&1574&8416&1425693555&4.3.25.23.37.103.281&9061&1336 \cr
 & &64.5.41.103.787&3935&1312& & &64.13.17.23.41.167&6847&12512 \cr
\noalign{\hrule}
 & &27.7.17.31.43.331&8723&5510& & &17.23.1319.2767&12595&9828 \cr
8399&1417649499&4.5.11.13.19.29.31.61&567&226&8417&1427022143&8.27.5.7.11.13.23.229&1363&1156 \cr
 & &16.81.5.7.19.29.113&2147&3480& & &64.3.5.7.13.289.29.47&31161&32480 \cr
\noalign{\hrule}
 & &3.49.11.13.19.53.67&64643&52540& & &9.25.11.13.17.2609&61&160 \cr
8400&1418265849&8.5.37.71.127.509&99&28&8418&1427057775&64.125.61.2609&5117&2508 \cr
 & &64.9.5.7.11.37.509&2545&3552& & &512.3.7.11.17.19.43&817&1792 \cr
\noalign{\hrule}
 & &9.25.49.13.19.521&18349&11576& & &9.11.19.37.73.281&2291&410 \cr
8401&1418774175&16.7.59.311.1447&517&930&8419&1427643261&4.5.29.41.79.281&2455&5694 \cr
 & &64.3.5.11.31.47.311&14617&10912& & &16.3.25.13.73.491&491&2600 \cr
\noalign{\hrule}
 & &9.7.23.37.59.449&14413&2200& & &9.17.19.61.83.97&5075&842 \cr
8402&1420261983&16.25.49.11.29.71&2663&1242&8420&1427659677&4.3.25.7.19.29.421&1037&616 \cr
 & &64.27.5.23.2663&2663&480& & &64.25.49.11.17.61&1225&352 \cr
\noalign{\hrule}
 & &3.5.11.13.31.41.521&209&6564& & &3.25.11.17.23.43.103&863&3438 \cr
8403&1420399695&8.9.121.19.547&3611&1312&8421&1428684675&4.27.43.191.863&149&1012 \cr
 & &512.23.41.157&157&5888& & &32.11.23.149.191&191&2384 \cr
\noalign{\hrule}
 & &9.19.23.31.61.191&73307&62876& & &3.49.17.197.2903&1035&1868 \cr
8404&1420524873&8.11.13.1429.5639&2105&3534&8422&1429155609&8.27.5.23.197.467&803&182 \cr
 & &32.3.5.11.13.19.31.421&2105&2288& & &32.7.11.13.73.467&5137&15184 \cr
\noalign{\hrule}
 & &9.7.11.13.19.43.193&113&360& & &3.5.7.19.41.101.173&1677&242 \cr
8405&1420548129&16.81.5.7.113.193&187&380&8423&1429204035&4.9.121.13.43.173&893&1010 \cr
 & &128.25.11.17.19.113&2825&1088& & &16.5.11.19.43.47.101&517&344 \cr
\noalign{\hrule}
 & &3.11.19.23.29.43.79&221&1090& & &9.25.11.223.2591&2533&58 \cr
8406&1420655973&4.5.13.17.29.43.109&7&36&8424&1430037675&4.17.29.149.223&1155&1378 \cr
 & &32.9.5.7.13.17.109&1853&21840& & &16.3.5.7.11.13.29.53&377&2968 \cr
\noalign{\hrule}
}%
}
$$
\eject
\vglue -23 pt
\noindent\hskip 1 in\hbox to 6.5 in{\ 8425 -- 8460 \hfill\fbd 1430256971 -- 1449344325\frb}
\vskip -9 pt
$$
\vbox{
\nointerlineskip
\halign{\strut
    \vrule \ \ \hfil \frb #\ 
   &\vrule \hfil \ \ \fbb #\frb\ 
   &\vrule \hfil \ \ \frb #\ \hfil
   &\vrule \hfil \ \ \frb #\ 
   &\vrule \hfil \ \ \frb #\ \ \vrule \hskip 2 pt
   &\vrule \ \ \hfil \frb #\ 
   &\vrule \hfil \ \ \fbb #\frb\ 
   &\vrule \hfil \ \ \frb #\ \hfil
   &\vrule \hfil \ \ \frb #\ 
   &\vrule \hfil \ \ \frb #\ \vrule \cr%
\noalign{\hrule}
 & &11.169.17.167.271&2927&54& & &13.37.1217.2459&989&1470 \cr
8425&1430256971&4.27.167.2927&1547&1380&8443&1439442043&4.3.5.49.23.43.1217&689&528 \cr
 & &32.81.5.7.13.17.23&405&2576& & &128.9.5.7.11.13.43.53&21285&23744 \cr
\noalign{\hrule}
 & &3.5.11.29.337.887&141&196& & &81.149.229.521&6149&5920 \cr
8426&1430327415&8.9.49.29.47.887&779&8762&8444&1439940321&64.5.11.13.37.43.521&12137&5364 \cr
 & &32.7.13.19.41.337&533&2128& & &512.9.5.53.149.229&265&256 \cr
\noalign{\hrule}
 & &9.5.19.53.131.241&3397&10340& & &11.17.23.29.31.373&171&202 \cr
8427&1430639865&8.3.25.11.43.47.79&131&106&8445&1442241427&4.9.11.17.19.23.29.101&373&120 \cr
 & &32.11.43.47.53.131&517&688& & &64.27.5.19.101.373&2565&3232 \cr
\noalign{\hrule}
 & &9.25.19.37.83.109&4433&3452& & &3.7.11.43.53.2741&167&134 \cr
8428&1431009225&8.5.11.13.31.37.863&4883&27048&8446&1442996709&4.53.67.167.2741&3055&5796 \cr
 & &128.3.49.19.23.257&5911&3136& & &32.9.5.7.13.23.47.67&16215&13936 \cr
\noalign{\hrule}
 & &27.5.11.43.73.307&217&2& & &3.11.31.61.73.317&475&548 \cr
8429&1431054405&4.9.7.11.31.307&1387&1376&8447&1444067823&8.25.19.61.137.317&7227&36202 \cr
 & &256.7.19.31.43.73&589&896& & &32.9.11.23.73.787&787&1104 \cr
\noalign{\hrule}
 & &5.29.43.53.61.71&2169&110& & &5.7.31.43.83.373&99&56 \cr
8430&1431200605&4.9.25.11.61.241&2167&2408&8448&1444392145&16.9.49.11.83.373&4085&18 \cr
 & &64.3.7.121.43.197&2541&6304& & &64.81.5.19.43&19&2592 \cr
\noalign{\hrule}
 & &9.125.7.11.13.31.41&2051&3074& & &11.31.37.53.2161&1251&910 \cr
8431&1431304875&4.3.49.13.29.53.293&16441&17320&8449&1445062861&4.9.5.7.13.37.53.139&2651&2492 \cr
 & &64.5.29.41.401.433&12557&12832& & &32.3.5.49.11.13.89.241&56693&57840 \cr
\noalign{\hrule}
 & &11.31.47.223.401&63861&51430& & &9.121.17.163.479&3325&2846 \cr
8432&1433182421&4.3.5.7.37.139.3041&2007&1034&8450&1445439501&4.3.25.7.19.163.1423&5269&1846 \cr
 & &16.27.5.11.37.47.223&185&216& & &16.5.11.13.19.71.479&923&760 \cr
\noalign{\hrule}
 & &9.49.11.19.103.151&2531&5400& & &9.11.19.31.43.577&41561&16750 \cr
8433&1433504457&16.243.25.7.2531&2987&5518&8451&1446754221&4.125.13.23.67.139&1731&1744 \cr
 & &64.5.29.31.89.103&2759&4640& & &128.3.5.23.67.109.577&7303&7360 \cr
\noalign{\hrule}
 & &27.121.89.4931&2955&1976& & &5.7.121.23.83.179&23049&27166 \cr
8434&1433752353&16.81.5.11.13.19.197&4099&356&8452&1447146085&4.9.7.13.289.47.197&1089&290 \cr
 & &128.13.89.4099&4099&832& & &16.81.5.121.13.17.29&2349&1768 \cr
\noalign{\hrule}
 & &9.5.49.37.73.241&2873&704& & &3.11.361.29.59.71&505&144 \cr
8435&1435324905&128.5.11.169.17.37&837&1022&8453&1447203153&32.27.5.29.71.101&7847&5788 \cr
 & &512.27.7.17.31.73&527&768& & &256.7.19.59.1447&1447&896 \cr
\noalign{\hrule}
 & &27.125.11.23.1681&597&2278& & &81.5.13.17.19.23.37&5&328 \cr
8436&1435363875&4.81.11.17.67.199&15893&21320&8454&1447206345&16.9.25.13.23.41&3553&3922 \cr
 & &64.5.13.23.41.691&691&416& & &64.11.17.19.37.53&53&352 \cr
\noalign{\hrule}
 & &49.11.13.19.41.263&705&2714& & &3.13.19.29.31.41.53&35&22 \cr
8437&1435573139&4.3.5.11.19.23.47.59&819&302&8455&1447563507&4.5.7.11.29.31.41.53&10647&2432 \cr
 & &16.27.5.7.13.23.151&4077&920& & &1024.9.49.169.19&637&1536 \cr
\noalign{\hrule}
 & &27.7.31.41.43.139&4895&1082& & &5.7.13.47.79.857&29711&48276 \cr
8438&1435788963&4.9.5.7.11.89.541&5779&172&8456&1447828655&8.81.11.37.73.149&857&6370 \cr
 & &32.5.43.5779&5779&80& & &32.9.5.49.13.857&63&16 \cr
\noalign{\hrule}
 & &3.5.7.11.19.29.37.61&47&2802& & &5.43.79.269.317&1857&1540 \cr
8439&1436366085&4.9.47.61.467&2071&2132&8457&1448361905&8.3.25.7.11.269.619&5529&1196 \cr
 & &32.13.19.41.47.109&4469&9776& & &64.9.11.13.19.23.97&51129&34144 \cr
\noalign{\hrule}
 & &3.25.17.19.31.1913&693&1220& & &9.11.19.29.101.263&59&40 \cr
8440&1436615175&8.27.125.7.11.19.61&319&194&8458&1448985087&16.5.29.59.101.263&1333&1596 \cr
 & &32.7.121.29.61.97&41419&56144& & &128.3.5.7.19.31.43.59&12803&13760 \cr
\noalign{\hrule}
 & &243.13.19.43.557&2255&12838& & &25.7.11.47.83.193&1547&14472 \cr
8441&1437562971&4.9.5.49.11.41.131&43&88&8459&1449319025&16.27.49.13.17.67&1525&386 \cr
 & &64.49.121.41.43&5929&1312& & &64.9.25.61.193&549&32 \cr
\noalign{\hrule}
 & &11.13.23.41.47.227&585&358& & &3.25.11.13.337.401&817&5198 \cr
8442&1438703981&4.9.5.11.169.47.179&681&164&8460&1449344325&4.5.11.19.23.43.113&27&28 \cr
 & &32.27.41.179.227&179&432& & &32.27.7.19.23.43.113&43731&48944 \cr
\noalign{\hrule}
}%
}
$$
\eject
\vglue -23 pt
\noindent\hskip 1 in\hbox to 6.5 in{\ 8461 -- 8496 \hfill\fbd 1449366919 -- 1465544619\frb}
\vskip -9 pt
$$
\vbox{
\nointerlineskip
\halign{\strut
    \vrule \ \ \hfil \frb #\ 
   &\vrule \hfil \ \ \fbb #\frb\ 
   &\vrule \hfil \ \ \frb #\ \hfil
   &\vrule \hfil \ \ \frb #\ 
   &\vrule \hfil \ \ \frb #\ \ \vrule \hskip 2 pt
   &\vrule \ \ \hfil \frb #\ 
   &\vrule \hfil \ \ \fbb #\frb\ 
   &\vrule \hfil \ \ \frb #\ \hfil
   &\vrule \hfil \ \ \frb #\ 
   &\vrule \hfil \ \ \frb #\ \vrule \cr%
\noalign{\hrule}
 & &7.121.13.23.59.97&3735&5092& & &11.29.67.163.419&369&50 \cr
8461&1449366919&8.9.5.121.19.67.83&161&40&8479&1459711781&4.9.25.41.67.163&77&412 \cr
 & &128.3.25.7.19.23.83&2075&3648& & &32.3.5.7.11.41.103&615&11536 \cr
\noalign{\hrule}
 & &9.11.13.17.23.43.67&38719&7028& & &3.5.7.61.193.1181&249&56 \cr
8462&1449768177&8.7.31.251.1249&4515&3266&8480&1459910865&16.9.49.83.1181&3281&7348 \cr
 & &32.3.5.49.23.43.71&355&784& & &128.11.17.167.193&2839&704 \cr
\noalign{\hrule}
 & &41.127.383.727&2967&2240& & &9.49.11.13.19.23.53&15355&15368 \cr
8463&1449842287&128.3.5.7.23.43.383&211&594&8481&1460602143&16.3.5.17.23.37.53.83.113&13&1232 \cr
 & &512.81.11.43.211&38313&54016& & &512.7.11.13.17.37.113&4181&4352 \cr
\noalign{\hrule}
 & &23.37.181.9421&104993&111690& & &3.11.17.19.43.3187&1185&2002 \cr
8464&1451126051&4.9.5.7.17.53.73.283&2629&74&8482&1460720019&4.9.5.7.121.13.17.79&475&596 \cr
 & &16.3.11.37.239.283&2629&6792& & &32.125.13.19.79.149&18625&16432 \cr
\noalign{\hrule}
 & &9.5.7.11.13.103.313&5255&4942& & &9.13.17.19.29.31.43&1265&9356 \cr
8465&1452205755&4.25.49.13.353.1051&31613&83112&8483&1460886687&8.5.11.17.23.2339&1161&1178 \cr
 & &64.3.101.313.3463&3463&3232& & &32.27.5.11.19.23.31.43&253&240 \cr
\noalign{\hrule}
 & &3.11.17.19.41.3323&9269&700& & &7.11.37.43.79.151&97&54 \cr
8466&1452214137&8.25.7.13.17.23.31&1089&1004&8484&1461386003&4.27.7.11.37.79.97&731&1600 \cr
 & &64.9.5.121.31.251&8283&4960& & &512.3.25.17.43.97&7275&4352 \cr
\noalign{\hrule}
 & &5.17.31.37.47.317&232507&229362& & &9.7.13.17.37.2837&265&3102 \cr
8467&1452578005&4.3.7.11.23.43.127.919&317&156&8485&1461483387&4.27.5.11.17.47.53&889&1406 \cr
 & &32.9.13.127.317.919&14859&14704& & &16.7.19.37.53.127&2413&424 \cr
\noalign{\hrule}
 & &81.25.11.37.41.43&553&338& & &9.5.7.11.13.71.457&199&298 \cr
8468&1453020525&4.5.7.169.37.41.79&129&2794&8486&1461575115&4.5.13.149.199.457&27&29678 \cr
 & &16.3.7.11.13.43.127&91&1016& & &16.27.11.19.71&1&456 \cr
\noalign{\hrule}
 & &3.125.343.11.13.79&8633&13092& & &5.17.37.59.7879&38253&1142 \cr
8469&1453076625&8.9.5.89.97.1091&2291&3164&8487&1461987845&4.3.41.311.571&11861&11550 \cr
 & &64.7.29.79.89.113&3277&2848& & &16.9.25.7.11.29.409&22495&14616 \cr
\noalign{\hrule}
 & &3.5.7.13.17.31.43.47&2563&1952& & &27.13.67.79.787&2253&3040 \cr
8470&1453816455&64.11.17.31.61.233&9315&2092&8488&1462122441&64.81.5.13.19.751&151&902 \cr
 & &512.81.5.23.523&12029&6912& & &256.5.11.19.41.151&42845&19328 \cr
\noalign{\hrule}
 & &27.25.239.9013&4469&4544& & &5.289.361.2803&4147&1344 \cr
8471&1454022225&128.9.41.71.109.239&26125&74&8489&1462170935&128.3.5.7.11.13.19.29&499&1734 \cr
 & &512.125.11.19.37&7733&1280& & &512.9.289.499&499&2304 \cr
\noalign{\hrule}
 & &9.25.11.289.19.107&7693&7982& & &7.17.97.151.839&31707&49676 \cr
8472&1454154075&4.3.49.13.107.157.307&235&1156&8490&1462371127&8.9.11.13.271.1129&20525&23506 \cr
 & &32.5.49.289.47.157&2303&2512& & &32.3.25.7.23.73.821&41975&39408 \cr
\noalign{\hrule}
 & &9.25.11.17.181.191&8701&5624& & &3.11.37.61.73.269&201&470 \cr
8473&1454574825&16.3.7.121.19.37.113&5539&1358&8491&1462583397&4.9.5.37.47.67.73&403&2882 \cr
 & &64.49.29.97.191&1421&3104& & &16.11.13.31.47.131&1703&11656 \cr
\noalign{\hrule}
 & &3.25.13.41.83.439&1591&4994& & &9.25.7.23.31.1303&319&164 \cr
8474&1456569075&4.5.11.13.37.43.227&1297&738&8492&1463236425&8.3.5.11.29.41.1303&869&434 \cr
 & &16.9.41.227.1297&3891&1816& & &32.7.121.31.41.79&3239&1936 \cr
\noalign{\hrule}
 & &27.11.17.31.67.139&8165&10528& & &3.289.37.43.1061&38269&988 \cr
8475&1457661447&64.3.5.7.11.23.47.71&4355&1112&8493&1463540217&8.49.11.13.19.71&629&720 \cr
 & &1024.25.13.67.139&325&512& & &256.9.5.7.11.17.37&165&896 \cr
\noalign{\hrule}
 & &3.13.29.43.157.191&3685&2438& & &5.11.13.17.101.1193&25967&46248 \cr
8476&1458357771&4.5.11.23.53.67.191&725&2826&8494&1464592415&16.3.23.41.47.1129&1105&24 \cr
 & &16.9.125.23.29.157&375&184& & &256.9.5.13.17.41&41&1152 \cr
\noalign{\hrule}
 & &5.23.29.37.53.223&3417&4268& & &3.7.1331.19.31.89&3775&5542 \cr
8477&1458405505&8.3.11.17.67.97.223&4551&17080&8495&1465219371&4.25.17.89.151.163&13671&836 \cr
 & &128.9.5.7.37.41.61&3843&2624& & &32.9.5.49.11.19.31&105&16 \cr
\noalign{\hrule}
 & &27.5.11.13.19.23.173&37177&34928& & &9.7.121.17.43.263&313&50 \cr
8478&1459477305&32.9.7.37.47.59.113&209&322&8496&1465544619&4.3.25.7.17.43.313&4387&6578 \cr
 & &128.49.11.19.23.37.47&2303&2368& & &16.5.11.13.23.41.107&12305&4264 \cr
\noalign{\hrule}
}%
}
$$
\eject
\vglue -23 pt
\noindent\hskip 1 in\hbox to 6.5 in{\ 8497 -- 8532 \hfill\fbd 1466215785 -- 1480382397\frb}
\vskip -9 pt
$$
\vbox{
\nointerlineskip
\halign{\strut
    \vrule \ \ \hfil \frb #\ 
   &\vrule \hfil \ \ \fbb #\frb\ 
   &\vrule \hfil \ \ \frb #\ \hfil
   &\vrule \hfil \ \ \frb #\ 
   &\vrule \hfil \ \ \frb #\ \ \vrule \hskip 2 pt
   &\vrule \ \ \hfil \frb #\ 
   &\vrule \hfil \ \ \fbb #\frb\ 
   &\vrule \hfil \ \ \frb #\ \hfil
   &\vrule \hfil \ \ \frb #\ 
   &\vrule \hfil \ \ \frb #\ \vrule \cr%
\noalign{\hrule}
 & &9.5.29.59.137.139&10109&14140& & &9.25.11.43.109.127&709&434 \cr
8497&1466215785&8.3.25.7.11.101.919&2641&116&8515&1473241275&4.7.31.43.109.709&4825&138 \cr
 & &64.7.11.19.29.139&77&608& & &16.3.25.23.31.193&4439&248 \cr
\noalign{\hrule}
 & &3.13.37.41.137.181&1197&1156& & &3.49.17.19.41.757&39259&8222 \cr
8498&1467064911&8.27.7.289.19.37.137&6895&1826&8516&1473667797&4.11.43.83.4111&24395&20826 \cr
 & &32.5.49.11.17.83.197&106183&112880& & &16.9.5.7.13.17.41.89&445&312 \cr
\noalign{\hrule}
 & &27.5.7.17.29.47.67&3097&52& & &3.7.11.17.53.73.97&111&790 \cr
8499&1467071865&8.9.13.17.19.163&77&94&8517&1473775611&4.9.5.11.37.73.79&1649&2014 \cr
 & &32.7.11.13.47.163&2119&176& & &16.17.19.53.79.97&79&152 \cr
\noalign{\hrule}
 & &27.5.7.11.29.31.157&15517&16354& & &125.23.31.71.233&4617&742 \cr
8500&1467181485&4.5.11.13.17.37.59.263&471&3364&8518&1474394875&4.243.7.19.53.71&715&634 \cr
 & &32.3.17.841.37.157&493&592& & &16.3.5.7.11.13.53.317&16801&24024 \cr
\noalign{\hrule}
 & &11.17.53.317.467&78975&69064& & &9.25.19.41.47.179&319&460 \cr
8501&1467214529&16.243.25.13.89.97&19&116&8519&1474588575&8.3.125.11.23.29.179&533&4658 \cr
 & &128.9.5.13.19.29.89&52065&35264& & &32.13.17.23.41.137&1781&6256 \cr
\noalign{\hrule}
 & &9.13.71.173.1021&26125&46366& & &81.5.29.37.43.79&56749&41764 \cr
8502&1467290331&4.125.11.19.97.239&1721&654&8520&1476217305&8.7.121.53.67.197&4559&8640 \cr
 & &16.3.109.239.1721&26051&13768& & &1024.27.5.11.47.97&4559&5632 \cr
\noalign{\hrule}
 & &9.7.19.53.101.229&2327&2024& & &19.23.31.61.1787&203&234 \cr
8503&1467326889&16.3.7.11.13.23.53.179&1313&7220&8521&1476717829&4.9.7.13.29.61.1787&95&5456 \cr
 & &128.5.169.361.101&845&1216& & &128.3.5.11.19.29.31&1595&192 \cr
\noalign{\hrule}
 & &3.25.13.29.223.233&52283&35558& & &9.25.13.23.29.757&108009&108766 \cr
8504&1469140725&4.49.11.23.97.773&467&306&8522&1476888075&4.81.7.11.17.457.1091&377&190 \cr
 & &16.9.7.11.17.97.467&34629&41096& & &16.5.13.19.29.457.1091&8683&8728 \cr
\noalign{\hrule}
 & &3.7.11.17.157.2383&537&1846& & &9.5.343.13.37.199&6161&2794 \cr
8505&1469212437&4.9.13.71.157.179&457&1870&8523&1477422765&4.49.11.61.101.127&2193&796 \cr
 & &16.5.11.17.71.457&457&2840& & &32.3.17.43.101.199&1717&688 \cr
\noalign{\hrule}
 & &9.11.53.229.1223&66607&54470& & &25.11.31.37.43.109&7061&10440 \cr
8506&1469511549&4.5.13.43.419.1549&1949&3498&8524&1478396975&16.9.125.23.29.307&2273&1352 \cr
 & &16.3.5.11.43.53.1949&1949&1720& & &256.3.169.23.2273&52279&64896 \cr
\noalign{\hrule}
 & &23.29.41.223.241&715&474& & &3.7.11.23.31.47.191&3071&1300 \cr
8507&1469709821&4.3.5.11.13.23.79.223&261&38&8525&1478538831&8.25.13.37.83.191&1449&1034 \cr
 & &16.27.5.11.19.29.79&7505&2376& & &32.9.5.7.11.23.37.47&111&80 \cr
\noalign{\hrule}
 & &5.7.13.529.31.197&187&342& & &3.529.43.47.461&17303&42166 \cr
8508&1469924365&4.9.7.11.13.17.19.197&2645&1098&8526&1478577747&4.1331.13.29.727&423&1150 \cr
 & &16.81.5.11.529.61&891&488& & &16.9.25.11.23.29.47&725&264 \cr
\noalign{\hrule}
 & &9.25.49.11.17.23.31&1247&2318& & &9.7.29.37.131.167&323&286 \cr
8509&1469974275&4.5.7.11.19.29.43.61&39&94&8527&1478863323&4.3.11.13.17.19.131.167&185&316 \cr
 & &16.3.13.29.43.47.61&26273&14152& & &32.5.11.13.17.19.37.79&16511&17680 \cr
\noalign{\hrule}
 & &9.5.11.13.19.23.523&101&146& & &5.11.19.31.71.643&313&468 \cr
8510&1470725685&4.11.23.73.101.523&135&388&8528&1478928935&8.9.13.19.313.643&445&198 \cr
 & &32.27.5.73.97.101&7373&4656& & &32.81.5.11.89.313&7209&5008 \cr
\noalign{\hrule}
 & &5.7.13.17.29.79.83&19&396& & &125.7.11.37.4153&15553&13518 \cr
8511&1470833455&8.9.7.11.17.19.79&3067&3016&8529&1478987125&4.9.25.103.151.751&761&3014 \cr
 & &128.3.13.19.29.3067&3067&3648& & &16.3.11.103.137.761&14111&18264 \cr
\noalign{\hrule}
 & &3.17.23.29.59.733&33&700& & &3.11.13.31.173.643&25&118 \cr
8512&1471133199&8.9.25.7.11.17.59&257&274&8530&1479367461&4.25.59.173.643&1521&1694 \cr
 & &32.25.7.11.137.257&23975&45232& & &16.9.5.7.121.169.59&3003&2360 \cr
\noalign{\hrule}
 & &3.5.43.59.67.577&403&2134& & &7.11.13.31.37.1289&515&774 \cr
8513&1471168245&4.5.11.13.31.67.97&4257&2242&8531&1479961483&4.9.5.11.13.31.43.103&2109&2324 \cr
 & &16.9.121.19.43.59&121&456& & &32.27.7.19.37.83.103&8549&8208 \cr
\noalign{\hrule}
 & &5.11.13.19.29.37.101&1581&1174& & &9.13.19.593.1123&5177&4930 \cr
8514&1472247205&4.3.13.17.31.101.587&363&950&8532&1480382397&4.5.17.29.31.167.593&1123&1716 \cr
 & &16.9.25.121.17.19.31&1705&1224& & &32.3.5.11.13.29.31.1123&899&880 \cr
\noalign{\hrule}
}%
}
$$
\eject
\vglue -23 pt
\noindent\hskip 1 in\hbox to 6.5 in{\ 8533 -- 8568 \hfill\fbd 1480563315 -- 1495631797\frb}
\vskip -9 pt
$$
\vbox{
\nointerlineskip
\halign{\strut
    \vrule \ \ \hfil \frb #\ 
   &\vrule \hfil \ \ \fbb #\frb\ 
   &\vrule \hfil \ \ \frb #\ \hfil
   &\vrule \hfil \ \ \frb #\ 
   &\vrule \hfil \ \ \frb #\ \ \vrule \hskip 2 pt
   &\vrule \ \ \hfil \frb #\ 
   &\vrule \hfil \ \ \fbb #\frb\ 
   &\vrule \hfil \ \ \frb #\ \hfil
   &\vrule \hfil \ \ \frb #\ 
   &\vrule \hfil \ \ \frb #\ \vrule \cr%
\noalign{\hrule}
 & &9.5.7.11.19.43.523&2867&2886& & &25.49.79.103.149&3127&14898 \cr
8533&1480563315&4.27.5.7.13.37.43.47.61&71651&28694&8551&1485205925&4.3.7.13.53.59.191&3223&8046 \cr
 & &16.13.137.523.14347&14347&14248& & &16.81.11.149.293&3223&648 \cr
\noalign{\hrule}
 & &5.169.23.29.37.71&3267&1634& & &11.13.29.53.67.101&315&1628 \cr
8534&1480616605&4.27.5.121.19.37.43&609&1426&8552&1487325697&8.9.5.7.121.37.53&2929&3484 \cr
 & &16.81.7.11.23.29.31&891&1736& & &64.3.7.13.29.67.101&21&32 \cr
\noalign{\hrule}
 & &25.19.29.293.367&8105&2538& & &9.7.13.29.31.43.47&2147&814 \cr
8535&1481239525&4.27.125.47.1621&539&586&8553&1488023901&4.11.13.19.29.37.113&3765&2296 \cr
 & &16.3.49.11.293.1621&4863&4312& & &64.3.5.7.37.41.251&10291&5920 \cr
\noalign{\hrule}
 & &7.13.23.41.61.283&517&426& & &5.11.19.37.61.631&3119&3822 \cr
8536&1481389819&4.3.11.47.61.71.283&123&2990&8554&1488254515&4.3.5.49.13.61.3119&2627&492 \cr
 & &16.9.5.13.23.41.71&355&72& & &32.9.7.13.37.41.71&4797&7952 \cr
\noalign{\hrule}
 & &25.13.17.31.41.211&5013&1528& & &9.7.11.13.373.443&1165&164 \cr
8537&1481700025&16.9.5.13.191.557&5117&2332&8555&1488638151&8.3.5.41.233.373&443&676 \cr
 & &128.3.7.11.17.43.53&9933&3392& & &64.5.169.41.443&533&160 \cr
\noalign{\hrule}
 & &9.5.13.29.113.773&517&952& & &27.25.49.19.23.103&155&776 \cr
8538&1481875785&16.3.7.11.17.47.773&565&208&8556&1488738825&16.125.31.97.103&7659&4466 \cr
 & &512.5.11.13.47.113&517&256& & &64.9.7.11.23.29.37&1073&352 \cr
\noalign{\hrule}
 & &9.5.7.289.73.223&253&3538& & &5.49.11.13.17.41.61&477&2024 \cr
8539&1481958765&4.7.11.17.23.29.61&321&988&8557&1489583095&16.9.5.7.121.23.53&481&366 \cr
 & &32.3.13.19.61.107&1159&22256& & &64.27.13.37.53.61&999&1696 \cr
\noalign{\hrule}
 & &7.13.23.31.53.431&45375&24016& & &19.43.67.73.373&17325&42316 \cr
8540&1482122369&32.3.125.121.19.79&483&562&8558&1490488231&8.9.25.7.11.71.149&1387&1742 \cr
 & &128.9.25.7.11.23.281&7025&6336& & &32.3.5.11.13.19.67.73&165&208 \cr
\noalign{\hrule}
 & &3.13.23.41.191.211&343&554& & &5.7.11.13.37.83.97&3933&4118 \cr
8541&1482149877&4.343.41.191.277&909&1100&8559&1490924435&4.9.7.11.13.19.23.29.71&52207&4850 \cr
 & &32.9.25.7.11.101.277&48475&53328& & &16.3.25.17.37.83.97&51&40 \cr
\noalign{\hrule}
 & &5.7.11.13.43.71.97&221&276& & &27.17.19.271.631&4225&4496 \cr
8542&1482185705&8.3.169.17.23.43.97&4029&142&8560&1491299721&32.25.169.281.631&6061&9714 \cr
 & &32.9.289.71.79&711&4624& & &128.3.11.13.19.29.1619&21047&20416 \cr
\noalign{\hrule}
 & &9.25.7.11.23.3721&1291&986& & &81.5.13.841.337&4009&5764 \cr
8543&1482725475&4.5.7.17.29.61.1291&11&1026&8561&1492190505&8.3.11.19.29.131.211&14299&13342 \cr
 & &16.27.11.19.1291&3873&152& & &32.7.19.79.181.953&172493&168112 \cr
\noalign{\hrule}
 & &9.5.11.17.23.47.163&335&182& & &13.67.79.109.199&3355&1938 \cr
8544&1482748245&4.25.7.13.23.67.163&13959&12284&8562&1492536019&4.3.5.11.17.19.61.199&109&90 \cr
 & &32.27.11.13.37.47.83&1443&1328& & &16.27.25.11.17.61.109&7425&8296 \cr
\noalign{\hrule}
 & &9.13.31.227.1801&41195&14636& & &27.7.11.13.167.331&2903&1400 \cr
8545&1482815529&8.5.7.11.107.3659&1241&2418&8563&1493971479&16.3.25.49.11.2903&643&2260 \cr
 & &32.3.5.7.13.17.31.73&1241&560& & &128.125.113.643&14125&41152 \cr
\noalign{\hrule}
 & &3.11.13.37.41.43.53&161&320& & &9.11.37.43.53.179&3875&3822 \cr
8546&1483157247&128.5.7.11.23.41.43&333&118&8564&1494287883&4.27.125.49.11.13.31.37&463&1462 \cr
 & &512.9.7.23.37.59&4071&1792& & &16.5.7.13.17.31.43.463&30095&29512 \cr
\noalign{\hrule}
 & &11.31.43.137.739&11691&11218& & &3.7.11.13.43.71.163&1021&1098 \cr
8547&1484526109&4.27.71.79.137.433&54901&4420&8565&1494409917&4.27.43.61.71.1021&835&71656 \cr
 & &32.3.5.7.11.13.17.23.31&3315&2576& & &64.5.169.53.167&10855&1696 \cr
\noalign{\hrule}
 & &9.121.13.17.31.199&655&434& & &3.25.11.17.19.71.79&43&518 \cr
8548&1484687061&4.5.7.961.131.199&1177&216&8566&1494658275&4.7.37.43.71.79&209&288 \cr
 & &64.27.5.11.107.131&1605&4192& & &256.9.11.19.37.43&1591&384 \cr
\noalign{\hrule}
 & &9.5.121.23.71.167&361&1142& & &27.5.103.191.563&349&214 \cr
8549&1484911395&4.5.11.361.23.571&44239&47094&8567&1495246365&4.103.107.191.349&27819&38840 \cr
 & &16.3.13.41.47.83.167&3901&4264& & &64.9.5.11.281.971&10681&8992 \cr
\noalign{\hrule}
 & &3.7.11.169.109.349&755&428& & &11.17.19.31.37.367&38925&37778 \cr
8550&1485082599&8.5.11.107.151.349&763&414&8568&1495631797&4.9.25.13.17.173.1453&89&1364 \cr
 & &32.9.5.7.23.109.151&755&1104& & &32.3.11.13.31.89.173&3471&2768 \cr
\noalign{\hrule}
}%
}
$$
\eject
\vglue -23 pt
\noindent\hskip 1 in\hbox to 6.5 in{\ 8569 -- 8604 \hfill\fbd 1496201175 -- 1515635693\frb}
\vskip -9 pt
$$
\vbox{
\nointerlineskip
\halign{\strut
    \vrule \ \ \hfil \frb #\ 
   &\vrule \hfil \ \ \fbb #\frb\ 
   &\vrule \hfil \ \ \frb #\ \hfil
   &\vrule \hfil \ \ \frb #\ 
   &\vrule \hfil \ \ \frb #\ \ \vrule \hskip 2 pt
   &\vrule \ \ \hfil \frb #\ 
   &\vrule \hfil \ \ \fbb #\frb\ 
   &\vrule \hfil \ \ \frb #\ \hfil
   &\vrule \hfil \ \ \frb #\ 
   &\vrule \hfil \ \ \frb #\ \vrule \cr%
\noalign{\hrule}
 & &9.25.7.23.103.401&97&212& & &3.7.29.37.137.487&4239&830 \cr
8569&1496201175&8.3.5.7.53.97.401&14729&694&8587&1503379227&4.81.5.29.83.157&1001&1406 \cr
 & &32.11.13.103.347&3817&208& & &16.7.11.13.19.37.157&2717&1256 \cr
\noalign{\hrule}
 & &3.49.11.19.113.431&39713&8990& & &7.11.29.59.101.113&1251&460 \cr
8570&1496302269&4.5.29.31.151.263&7889&15516&8588&1503628511&8.9.5.11.23.101.139&1017&512 \cr
 & &32.9.343.23.431&161&48& & &8192.81.23.113&1863&4096 \cr
\noalign{\hrule}
 & &27.169.17.101.191&577&2296& & &3.5.11.41.131.1697&9231&9436 \cr
8571&1496422161&16.3.7.41.101.577&137&440&8589&1503906855&8.9.7.17.131.181.337&299&1928 \cr
 & &256.5.7.11.41.137&39319&7040& & &128.7.13.23.241.337&100763&107968 \cr
\noalign{\hrule}
 & &27.5.7.121.13.19.53&73&1162& & &9.7.11.13.103.1621&535&598 \cr
8572&1496890395&4.3.49.53.73.83&5339&5392&8590&1504169667&4.5.169.23.107.1621&1133&2754 \cr
 & &128.19.83.281.337&23323&21568& & &16.81.5.11.17.103.107&765&856 \cr
\noalign{\hrule}
 & &9.11.43.61.73.79&1205&3614& & &3.5.17.19.79.3931&44671&22156 \cr
8573&1497557259&4.3.5.13.43.139.241&4453&4582&8591&1504609905&8.11.29.31.131.191&2835&2704 \cr
 & &16.29.61.73.79.241&241&232& & &256.81.5.7.11.169.31&36673&38016 \cr
\noalign{\hrule}
 & &3.5.7.113.293.431&67&498& & &243.25.13.17.19.59&737&4868 \cr
8574&1498347795&4.9.7.67.83.293&7201&12430&8592&1505026575&8.5.11.13.67.1217&7833&1748 \cr
 & &16.5.11.19.113.379&379&1672& & &64.3.7.19.23.373&373&5152 \cr
\noalign{\hrule}
 & &81.25.11.19.3541&1343&2198& & &5.49.11.61.9161&213&458 \cr
8575&1498639725&4.9.5.7.11.17.79.157&3541&1814&8593&1506022595&4.3.71.229.9161&3549&12710 \cr
 & &16.79.907.3541&907&632& & &16.9.5.7.169.31.41&6929&2232 \cr
\noalign{\hrule}
 & &27.125.289.29.53&10153&1772& & &27.5.49.11.127.163&109&136 \cr
8576&1499151375&8.3.5.11.13.71.443&10591&6686&8594&1506308265&16.11.17.109.127.163&9751&10950 \cr
 & &32.7.17.89.3343&3343&9968& & &64.3.25.49.17.73.199&6205&6368 \cr
\noalign{\hrule}
 & &9.7.121.239.823&425&664& & &31.163.379.787&70015&58266 \cr
8577&1499421231&16.25.7.17.83.823&121&702&8595&1507173469&4.27.5.11.13.19.67.83&155&758 \cr
 & &64.27.25.121.13.17&425&1248& & &16.3.25.13.19.31.379&475&312 \cr
\noalign{\hrule}
 & &5.11.13.17.19.43.151&783&878& & &27.5.11.13.23.43.79&331&538 \cr
8578&1499525885&4.27.13.17.29.43.439&14345&23848&8596&1508318955&4.3.5.13.43.269.331&2297&1738 \cr
 & &64.9.5.11.19.151.271&271&288& & &16.11.79.331.2297&2297&2648 \cr
\noalign{\hrule}
 & &9.25.7.17.79.709&8507&9218& & &3.23.157.277.503&173&330 \cr
8579&1499694525&4.7.11.17.47.181.419&923&20616&8597&1509372723&4.9.5.11.23.173.277&2041&1006 \cr
 & &64.3.11.13.71.859&11167&24992& & &16.13.157.173.503&173&104 \cr
\noalign{\hrule}
 & &9.121.43.103.311&7397&39430& & &7.41.103.199.257&1415&2808 \cr
8580&1500009291&4.5.13.569.3943&1727&5670&8598&1511838223&16.27.5.13.257.283&1067&218 \cr
 & &16.81.25.7.11.157&3925&504& & &64.9.11.13.97.109&10791&40352 \cr
\noalign{\hrule}
 & &5.7.11.13.43.6971&5141&1830& & &25.11.13.59.67.107&14679&14746 \cr
8581&1500263765&4.3.25.13.53.61.97&11571&5654&8599&1512121325&4.9.7.13.59.73.101.233&737&30 \cr
 & &16.9.7.11.19.29.257&4883&2088& & &16.27.5.11.67.73.233&1971&1864 \cr
\noalign{\hrule}
 & &9.25.11.43.59.239&72991&73034& & &9.125.13.29.43.83&943&682 \cr
8582&1500698925&4.13.47.2809.239.1553&114169&31860&8600&1513702125&4.11.23.31.41.43.83&8775&5206 \cr
 & &32.27.5.11.53.59.97.107&5141&5136& & &16.27.25.13.19.23.137&1311&1096 \cr
\noalign{\hrule}
 & &343.11.13.113.271&22145&16608& & &3.25.11.13.17.361.23&6231&94 \cr
8583&1502027527&64.3.5.7.43.103.173&813&5242&8601&1513844475&4.9.13.31.47.67&1097&980 \cr
 & &256.9.271.2621&2621&1152& & &32.5.49.47.1097&2303&17552 \cr
\noalign{\hrule}
 & &81.7.11.23.37.283&505&494& & &17.47.359.5281&11077&5796 \cr
8584&1502075421&4.3.5.7.13.19.23.101.283&3047&16480&8602&1514807321&8.9.7.11.17.19.23.53&1077&1660 \cr
 & &256.25.11.13.103.277&33475&35456& & &64.27.5.19.83.359&2241&3040 \cr
\noalign{\hrule}
 & &27.49.31.67.547&7619&7150& & &5.17.19.43.139.157&5929&5886 \cr
8585&1503085437&4.25.7.11.13.19.31.401&67&522&8603&1515498235&4.27.49.121.19.109.157&8671&278 \cr
 & &16.9.5.11.29.67.401&2005&2552& & &16.9.7.11.13.23.29.139&7337&6552 \cr
\noalign{\hrule}
 & &3.7.11.13.29.41.421&79&530& & &11.13.83.277.461&3135&2858 \cr
8586&1503208707&4.5.13.53.79.421&539&1566&8604&1515635693&4.3.5.121.19.83.1429&849&728 \cr
 & &16.27.49.11.29.53&477&56& & &64.9.5.7.13.283.1429&64305&63392 \cr
\noalign{\hrule}
}%
}
$$
\eject
\vglue -23 pt
\noindent\hskip 1 in\hbox to 6.5 in{\ 8605 -- 8640 \hfill\fbd 1515903619 -- 1534808847\frb}
\vskip -9 pt
$$
\vbox{
\nointerlineskip
\halign{\strut
    \vrule \ \ \hfil \frb #\ 
   &\vrule \hfil \ \ \fbb #\frb\ 
   &\vrule \hfil \ \ \frb #\ \hfil
   &\vrule \hfil \ \ \frb #\ 
   &\vrule \hfil \ \ \frb #\ \ \vrule \hskip 2 pt
   &\vrule \ \ \hfil \frb #\ 
   &\vrule \hfil \ \ \fbb #\frb\ 
   &\vrule \hfil \ \ \frb #\ \hfil
   &\vrule \hfil \ \ \frb #\ 
   &\vrule \hfil \ \ \frb #\ \vrule \cr%
\noalign{\hrule}
 & &361.23.41.61.73&23409&45430& & &9.25.7.11.17.71.73&8917&9698 \cr
8605&1515903619&4.81.5.7.11.289.59&557&2622&8623&1526523075&4.3.5.7.13.37.241.373&1207&88 \cr
 & &16.243.19.23.557&557&1944& & &64.11.13.17.71.241&241&416 \cr
\noalign{\hrule}
 & &9.5.11.19.23.43.163&2863&1070& & &3.49.11.17.19.37.79&41&8 \cr
8606&1516151835&4.25.7.43.107.409&163&912&8624&1526656593&16.17.19.37.41.79&7595&4356 \cr
 & &128.3.19.163.409&409&64& & &128.9.5.49.121.31&1023&320 \cr
\noalign{\hrule}
 & &3.11.17.19.23.41.151&3409&2782& & &49.11.17.19.31.283&1&18 \cr
8607&1517766987&4.7.13.17.23.107.487&79173&66440&8625&1527352981&4.9.49.11.31.283&3835&722 \cr
 & &64.9.5.11.19.151.463&463&480& & &16.3.5.13.361.59&195&8968 \cr
\noalign{\hrule}
 & &3.5.11.29.83.3823&19823&22230& & &81.121.13.19.631&45373&44860 \cr
8608&1518323565&4.27.25.13.19.43.461&617&58&8626&1527554457&8.3.5.11.289.157.2243&281&2524 \cr
 & &16.19.29.461.617&8759&4936& & &64.17.157.281.631&4777&5024 \cr
\noalign{\hrule}
 & &25.121.13.29.31.43&981&266& & &31.211.409.571&189&220 \cr
8609&1520186525&4.9.5.7.11.19.31.109&559&464&8627&1527578599&8.27.5.7.11.211.571&925&5356 \cr
 & &128.3.7.13.29.43.109&327&448& & &64.9.125.13.37.103&34299&52000 \cr
\noalign{\hrule}
 & &27.5.7.11.19.43.179&493&1388& & &11.23.31.397.491&5427&9794 \cr
8610&1520195985&8.3.7.17.29.43.347&4009&3278&8628&1528812461&4.81.23.59.67.83&1859&50 \cr
 & &32.11.19.29.149.211&4321&3376& & &16.3.25.11.169.59&507&11800 \cr
\noalign{\hrule}
 & &5.7.11.17.31.59.127&7739&3294& & &3.7.11.29.31.53.139&2309&2278 \cr
8611&1520292235&4.27.31.61.71.109&1375&826&8629&1529897523&4.7.17.29.53.67.2309&3501&50 \cr
 & &16.3.125.7.11.59.109&327&200& & &16.9.25.389.2309&29175&18472 \cr
\noalign{\hrule}
 & &5.343.13.19.37.97&99&580& & &27.5.11.13.31.2557&29981&49286 \cr
8612&1520318345&8.9.25.49.11.19.29&481&4744&8630&1530249435&4.7.19.1297.4283&1493&2790 \cr
 & &128.3.13.37.593&1779&64& & &16.9.5.7.19.31.1493&1493&1064 \cr
\noalign{\hrule}
 & &7.13.61.79.3467&33057&29590& & &49.11.23.37.47.71&48347&28404 \cr
8613&1520380043&4.9.5.7.11.269.3673&4661&988&8631&1530645193&8.27.13.263.3719&1465&2254 \cr
 & &32.3.5.11.13.19.59.79&1947&1520& & &32.9.5.49.13.23.293&1465&1872 \cr
\noalign{\hrule}
 & &27.7.19.43.59.167&1525&1012& & &3.5.151.191.3539&40963&12122 \cr
8614&1521431289&8.25.7.11.23.61.167&9503&318&8632&1531024485&4.11.13.19.23.29.137&1017&764 \cr
 & &32.3.5.13.17.43.53&3445&272& & &32.9.19.29.113.191&2147&1392 \cr
\noalign{\hrule}
 & &5.7.23.47.131.307&25527&46618& & &3.11.19.23.89.1193&7585&5538 \cr
8615&1521610195&4.3.11.13.67.127.163&70923&71050&8633&1531178517&4.9.5.13.19.37.41.71&2047&2750 \cr
 & &16.9.25.49.11.29.47.503&20097&20120& & &16.625.11.23.71.89&625&568 \cr
\noalign{\hrule}
 & &9.11.13.23.101.509&70847&59140& & &9.5.43.179.4423&23579&45694 \cr
8616&1521757809&8.5.7.29.349.2957&2351&606&8634&1531972395&4.11.17.19.31.67.73&725&516 \cr
 & &32.3.7.29.101.2351&2351&3248& & &32.3.25.29.31.43.67&1943&2480 \cr
\noalign{\hrule}
 & &5.53.1367.4201&25723&46728& & &7.11.13.19.23.31.113&617&39150 \cr
8617&1521833255&16.9.11.29.59.887&1417&530&8635&1532341811&4.27.25.29.617&1547&1538 \cr
 & &64.3.5.13.29.53.109&1417&2784& & &16.3.5.7.13.17.29.769&7395&6152 \cr
\noalign{\hrule}
 & &9.11.29.197.2693&263825&266696& & &3.7.11.29.31.47.157&125&32 \cr
8618&1523125791&16.25.17.37.53.61.173&10179&374&8636&1532389551&64.125.7.11.29.47&157&360 \cr
 & &64.27.5.11.13.289.29&1445&1248& & &1024.9.625.157&1875&512 \cr
\noalign{\hrule}
 & &27.25.49.11.53.79&2881&1306& & &9.25.11.13.29.31.53&7177&2248 \cr
8619&1523335275&4.3.7.11.43.67.653&2635&676&8637&1533042225&16.3.11.281.7177&157&124 \cr
 & &32.5.169.17.31.67&5239&18224& & &128.31.157.7177&7177&10048 \cr
\noalign{\hrule}
 & &9.5.7.11.431.1021&3467&39202& & &9.29.41.251.571&329&242 \cr
8620&1524776715&4.17.1153.3467&1157&2310&8638&1533678021&4.3.7.121.41.47.251&1543&13340 \cr
 & &16.3.5.7.11.13.17.89&1157&136& & &32.5.7.23.29.1543&7715&2576 \cr
\noalign{\hrule}
 & &3.7.11.17.19.107.191&325&248& & &9.625.281.971&2849&2006 \cr
8621&1524865881&16.25.13.17.19.31.107&18909&26566&8639&1534786875&4.3.125.7.11.17.37.59&971&3154 \cr
 & &64.9.11.37.191.359&1077&1184& & &16.7.17.19.83.971&1577&952 \cr
\noalign{\hrule}
 & &9.5.7.11.17.19.29.47&383&416& & &3.11.13.19.29.43.151&1649&1220 \cr
8622&1525462785&64.3.5.7.13.19.29.383&289&376&8640&1534808847&8.5.17.29.43.61.97&95&2718 \cr
 & &1024.13.289.47.383&6511&6656& & &32.9.25.17.19.151&425&48 \cr
\noalign{\hrule}
}%
}
$$
\eject
\vglue -23 pt
\noindent\hskip 1 in\hbox to 6.5 in{\ 8641 -- 8676 \hfill\fbd 1535039757 -- 1553831207\frb}
\vskip -9 pt
$$
\vbox{
\nointerlineskip
\halign{\strut
    \vrule \ \ \hfil \frb #\ 
   &\vrule \hfil \ \ \fbb #\frb\ 
   &\vrule \hfil \ \ \frb #\ \hfil
   &\vrule \hfil \ \ \frb #\ 
   &\vrule \hfil \ \ \frb #\ \ \vrule \hskip 2 pt
   &\vrule \ \ \hfil \frb #\ 
   &\vrule \hfil \ \ \fbb #\frb\ 
   &\vrule \hfil \ \ \frb #\ \hfil
   &\vrule \hfil \ \ \frb #\ 
   &\vrule \hfil \ \ \frb #\ \vrule \cr%
\noalign{\hrule}
 & &9.23.37.43.59.79&323&1034& & &5.49.11.19.47.641&155&174 \cr
8641&1535039757&4.11.17.19.37.43.47&2449&4470&8659&1542653035&4.3.25.7.11.29.31.641&6063&988 \cr
 & &16.3.5.19.31.79.149&2831&1240& & &32.9.13.19.31.43.47&1333&1872 \cr
\noalign{\hrule}
 & &7.11.13.17.31.41.71&555&368& & &3.25.19.3481.311&33847&32292 \cr
8642&1535631097&32.3.5.7.23.31.37.41&2343&3614&8660&1542692175&8.81.5.11.13.17.23.181&3481&5834 \cr
 & &128.9.5.11.13.71.139&695&576& & &32.11.17.3481.2917&2917&2992 \cr
\noalign{\hrule}
 & &5.13.19.367.3391&4609&39474& & &3.13.17.23.107.947&757&1704 \cr
8643&1536953795&4.27.11.17.43.419&20375&21106&8661&1545165921&16.9.13.17.71.757&385&538 \cr
 & &16.3.125.61.163.173&29829&34600& & &64.5.7.11.269.757&58289&43040 \cr
\noalign{\hrule}
 & &3.13.29.67.103.197&6307&594& & &9.11.17.43.131.163&923&760 \cr
8644&1537591107&4.81.7.11.13.17.53&335&556&8662&1545295257&16.5.13.19.43.71.131&81&736 \cr
 & &32.5.7.53.67.139&4865&848& & &1024.81.13.23.71&14697&6656 \cr
\noalign{\hrule}
 & &9.43.53.167.449&145&304& & &9.25.11.13.43.1117&60473&45952 \cr
8645&1537976313&32.3.5.19.29.43.167&143&358&8663&1545397425&256.7.53.163.359&55&108 \cr
 & &128.11.13.19.29.179&44213&20416& & &2048.27.5.7.11.359&2513&3072 \cr
\noalign{\hrule}
 & &9.25.19.37.71.137&3451&824& & &3.25.11.37.179.283&3269&3806 \cr
8646&1538568225&16.7.17.29.103.137&1425&1562&8664&1546304925&4.7.121.37.173.467&50657&5850 \cr
 & &64.3.25.7.11.17.19.71&187&224& & &16.9.25.13.179.283&13&24 \cr
\noalign{\hrule}
 & &37.113.569.647&11685&12254& & &81.7.11.29.43.199&42869&17030 \cr
8647&1539205883&4.3.5.11.19.41.113.557&335&222&8665&1547730261&4.5.13.131.163.263&1791&1628 \cr
 & &16.9.25.11.19.37.41.67&31825&32472& & &32.9.5.11.37.131.199&655&592 \cr
\noalign{\hrule}
 & &3.7.47.61.107.239&4267&39470& & &27.7.23.193.1847&6775&6154 \cr
8648&1539673611&4.5.17.251.3947&34177&32922&8666&1549579437&4.25.17.181.193.271&187&6 \cr
 & &16.9.11.13.31.59.239&2301&2728& & &16.3.25.11.289.271&7225&23848 \cr
\noalign{\hrule}
 & &9.7.13.37.89.571&3193&4906& & &3.25.11.17.29.37.103&369&556 \cr
8649&1539968157&4.3.11.31.37.103.223&28145&40762&8667&1550028975&8.27.29.41.103.139&625&3406 \cr
 & &16.5.13.89.229.433&2165&1832& & &32.625.13.41.131&3275&8528 \cr
\noalign{\hrule}
 & &9.5.7.41.43.47.59&10013&9398& & &27.13.107.149.277&621&770 \cr
8650&1539971685&4.3.17.19.31.37.43.127&2303&110&8668&1550090061&4.729.5.7.11.23.277&299&3346 \cr
 & &16.5.49.11.31.37.47&1147&616& & &16.49.13.529.239&11711&4232 \cr
\noalign{\hrule}
 & &9.5.19.41.197.223&501&278& & &27.5.13.17.223.233&17177&14278 \cr
8651&1540001205&4.27.5.139.167.197&1691&24904&8669&1550196765&4.121.17.59.89.193&1125&932 \cr
 & &64.11.19.89.283&979&9056& & &32.9.125.59.89.233&1475&1424 \cr
\noalign{\hrule}
 & &25.7.11.13.19.41.79&657&370& & &5.11.13.17.89.1433&227&1206 \cr
8652&1540063525&4.9.125.11.19.37.73&1391&984&8670&1550212235&4.9.5.13.17.67.227&1433&1518 \cr
 & &64.27.13.41.73.107&2889&2336& & &16.27.11.23.67.1433&621&536 \cr
\noalign{\hrule}
 & &81.11.43.131.307&3277&2356& & &11.13.47.71.3253&3799&39582 \cr
8653&1540833921&8.27.11.19.29.31.113&3419&6470&8671&1552302323&4.27.29.131.733&1165&1034 \cr
 & &32.5.13.19.263.647&64961&51760& & &16.9.5.11.29.47.233&1165&2088 \cr
\noalign{\hrule}
 & &7.19.31.337.1109&11715&9356& & &9.5.11.89.131.269&1073&272 \cr
8654&1540901159&8.3.5.11.31.71.2339&69&2270&8672&1552454145&32.11.17.29.37.131&269&138 \cr
 & &32.9.25.11.23.227&5221&39600& & &128.3.17.23.29.269&493&1472 \cr
\noalign{\hrule}
 & &3.5.11.19.593.829&114871&121394& & &7.11.13.53.73.401&2515&1566 \cr
8655&1541156595&4.7.13.23.29.313.367&5985&1214&8673&1553020469&4.27.5.29.401.503&3551&8078 \cr
 & &16.9.5.49.19.29.607&4263&4856& & &16.3.5.7.53.67.577&1731&2680 \cr
\noalign{\hrule}
 & &11.19.67.283.389&2655&2722& & &3.5.13.289.43.641&1561&116 \cr
8656&1541548261&4.9.5.11.59.389.1361&18961&3990&8674&1553312865&8.7.29.223.641&10075&8514 \cr
 & &16.27.25.7.19.67.283&189&200& & &32.9.25.11.13.31.43&465&176 \cr
\noalign{\hrule}
 & &9.13.19.61.83.137&1405&172& & &27.11.43.239.509&46435&24548 \cr
8657&1541941713&8.5.13.43.61.281&375&418&8675&1553604921&8.5.17.361.37.251&979&276 \cr
 & &32.3.625.11.19.281&6875&4496& & &64.3.11.17.19.23.89&7429&2848 \cr
\noalign{\hrule}
 & &81.7.11.37.41.163&2231&1090& & &73.2663.7993&93203&101196 \cr
8658&1542229227&4.5.11.23.37.97.109&1025&174&8676&1553831207&8.27.11.37.229.937&2665&146 \cr
 & &16.3.125.29.41.97&3625&776& & &32.9.5.13.37.41.73&2665&5328 \cr
\noalign{\hrule}
}%
}
$$
\eject
\vglue -23 pt
\noindent\hskip 1 in\hbox to 6.5 in{\ 8677 -- 8712 \hfill\fbd 1554482853 -- 1575028123\frb}
\vskip -9 pt
$$
\vbox{
\nointerlineskip
\halign{\strut
    \vrule \ \ \hfil \frb #\ 
   &\vrule \hfil \ \ \fbb #\frb\ 
   &\vrule \hfil \ \ \frb #\ \hfil
   &\vrule \hfil \ \ \frb #\ 
   &\vrule \hfil \ \ \frb #\ \ \vrule \hskip 2 pt
   &\vrule \ \ \hfil \frb #\ 
   &\vrule \hfil \ \ \fbb #\frb\ 
   &\vrule \hfil \ \ \frb #\ \hfil
   &\vrule \hfil \ \ \frb #\ 
   &\vrule \hfil \ \ \frb #\ \vrule \cr%
\noalign{\hrule}
 & &27.7.11.19.23.29.59&2069&2278& & &5.7.19.857.2741&6919&6786 \cr
8677&1554482853&4.17.29.59.67.2069&1425&286&8695&1562109605&4.9.11.13.17.29.37.857&2501&70 \cr
 & &16.3.25.11.13.19.2069&2069&2600& & &16.3.5.7.29.37.41.61&6771&9512 \cr
\noalign{\hrule}
 & &3.13.29.31.101.439&385&824& & &3.7.169.17.19.29.47&2263&610 \cr
8678&1554569679&16.5.7.11.29.101.103&65&36&8696&1562443701&4.5.7.31.47.61.73&2873&7326 \cr
 & &128.9.25.7.11.13.103&7931&4800& & &16.9.5.11.169.17.37&111&440 \cr
\noalign{\hrule}
 & &3.11.23.29.31.43.53&262123&261994& & &9.5.11.17.47.59.67&3743&3448 \cr
8679&1555055139&4.289.23.101.907.1297&299367&207760&8697&1563431265&16.11.19.67.197.431&35025&6148 \cr
 & &128.9.5.49.17.29.31.37.53&5439&5440& & &128.3.25.29.53.467&13543&16960 \cr
\noalign{\hrule}
 & &9.5.7.13.41.59.157&4369&2068& & &81.49.47.83.101&13253&19976 \cr
8680&1555211385&8.3.5.7.11.17.47.257&6571&6536&8698&1563789969&16.7.11.29.227.457&635&954 \cr
 & &128.11.19.43.47.6571&422389&420544& & &64.9.5.53.127.457&24221&20320 \cr
\noalign{\hrule}
 & &27.25.121.17.19.59&277&398& & &5.169.61.97.313&1413&152 \cr
8681&1556480475&4.17.19.59.199.277&2915&1794&8699&1564957745&16.9.13.19.61.157&3641&2482 \cr
 & &16.3.5.11.13.23.53.199&4577&5512& & &64.3.11.17.73.331&40953&10592 \cr
\noalign{\hrule}
 & &3.7.37.47.89.479&209&120& & &7.11.17.31.47.821&1521&7268 \cr
8682&1556841489&16.9.5.11.19.37.479&3913&398&8700&1565821873&8.9.169.23.31.79&7055&4606 \cr
 & &64.7.11.13.43.199&2587&15136& & &32.3.5.49.17.47.83&581&240 \cr
\noalign{\hrule}
 & &243.25.7.19.41.47&1279&5896& & &9.7.29.43.127.157&16837&21390 \cr
8683&1556967825&16.11.47.67.1279&935&2214&8701&1566427779&4.27.5.23.31.113.149&3311&712 \cr
 & &64.27.5.121.17.41&121&544& & &64.5.7.11.31.43.89&1705&2848 \cr
\noalign{\hrule}
 & &81.11.41.47.907&4777&5200& & &25.11.17.31.79.137&817&1158 \cr
8684&1557279999&32.9.25.13.17.41.281&6649&376&8702&1568523275&4.3.17.19.43.137.193&93&230 \cr
 & &512.13.47.61.109&6649&3328& & &16.9.5.23.31.43.193&4439&3096 \cr
\noalign{\hrule}
 & &5.13.19.23.29.31.61&371&66& & &27.19.53.197.293&553&3190 \cr
8685&1557701795&4.3.7.11.13.29.31.53&361&42&8703&1569376269&4.3.5.7.11.29.53.79&589&1702 \cr
 & &16.9.49.361.53&9063&392& & &16.5.11.19.23.31.37&1265&9176 \cr
\noalign{\hrule}
 & &81.5.31.101.1229&7163&1018& & &29.41.743.1777&36443&36414 \cr
8686&1558439595&4.13.19.29.31.509&187&216&8704&1569849779&4.9.7.11.289.743.3313&725&37168 \cr
 & &64.27.11.17.19.509&9671&5984& & &128.3.25.7.17.23.29.101&36057&36800 \cr
\noalign{\hrule}
 & &13.31.41.157.601&19315&43956& & &27.7.11.37.137.149&1919&412 \cr
8687&1559060711&8.27.5.11.37.3863&1321&2542&8705&1570229199&8.3.19.101.103.149&4795&10552 \cr
 & &32.9.5.31.41.1321&1321&720& & &128.5.7.137.1319&1319&320 \cr
\noalign{\hrule}
 & &27.5.7.13.61.2081&8107&6460& & &9.5.13.841.31.103&77&68 \cr
8688&1559470185&8.25.121.13.17.19.67&261&14&8706&1570908105&8.7.11.13.17.29.31.103&1695&1292 \cr
 & &32.9.7.11.17.29.67&1139&5104& & &64.3.5.7.11.289.19.113&38437&39776 \cr
\noalign{\hrule}
 & &3.5.1369.139.547&893&476& & &3.11.17.23.233.523&1719&2242 \cr
8689&1561337655&8.5.7.17.19.47.547&1529&1206&8707&1572346677&4.27.11.19.23.59.191&1601&5230 \cr
 & &32.9.7.11.47.67.139&3149&3696& & &16.5.59.523.1601&1601&2360 \cr
\noalign{\hrule}
 & &9.5.49.13.19.47.61&437&2552& & &11.13.19.53.61.179&405&1564 \cr
8690&1561468545&16.11.13.361.23.29&6321&2350&8708&1572346919&8.81.5.13.17.23.53&305&358 \cr
 & &64.3.25.49.43.47&215&32& & &32.27.25.23.61.179&575&432 \cr
\noalign{\hrule}
 & &11.13.139.251.313&3807&39700& & &3.5.343.23.97.137&6077&8308 \cr
8691&1561596751&8.81.25.47.397&131&266&8709&1572553815&8.49.31.59.67.103&1485&3562 \cr
 & &32.3.5.7.19.47.131&4935&39824& & &32.27.5.11.13.59.137&1287&944 \cr
\noalign{\hrule}
 & &81.49.11.83.431&17431&17480& & &7.23.47.97.2143&14885&17028 \cr
8692&1561813407&16.5.11.19.23.83.17431&19215&1784&8710&1572959857&8.9.5.11.13.23.43.229&493&194 \cr
 & &256.9.25.7.19.61.223&28975&28544& & &32.3.5.11.17.29.43.97&8041&6960 \cr
\noalign{\hrule}
 & &3.13.23.29.97.619&1097&1716& & &3.343.11.13.289.37&473&1340 \cr
8693&1561898559&8.9.11.169.23.1097&2993&6880&8711&1573442871&8.5.7.121.13.43.67&703&144 \cr
 & &512.5.11.41.43.73&34529&52480& & &256.9.5.19.37.67&3819&640 \cr
\noalign{\hrule}
 & &27.11.23.107.2137&42385&15314& & &19.1849.107.419&715&1134 \cr
8694&1561969629&4.5.49.13.19.31.173&1089&4274&8712&1575028123&4.81.5.7.11.13.19.107&2005&3182 \cr
 & &16.9.121.19.2137&19&88& & &16.27.25.37.43.401&10025&7992 \cr
\noalign{\hrule}
}%
}
$$
\eject
\vglue -23 pt
\noindent\hskip 1 in\hbox to 6.5 in{\ 8713 -- 8748 \hfill\fbd 1575143947 -- 1593319195\frb}
\vskip -9 pt
$$
\vbox{
\nointerlineskip
\halign{\strut
    \vrule \ \ \hfil \frb #\ 
   &\vrule \hfil \ \ \fbb #\frb\ 
   &\vrule \hfil \ \ \frb #\ \hfil
   &\vrule \hfil \ \ \frb #\ 
   &\vrule \hfil \ \ \frb #\ \ \vrule \hskip 2 pt
   &\vrule \ \ \hfil \frb #\ 
   &\vrule \hfil \ \ \fbb #\frb\ 
   &\vrule \hfil \ \ \frb #\ \hfil
   &\vrule \hfil \ \ \frb #\ 
   &\vrule \hfil \ \ \frb #\ \vrule \cr%
\noalign{\hrule}
 & &13.19.47.241.563&16533&42994& & &25.11.23.29.89.97&247&732 \cr
8713&1575143947&4.9.7.11.37.83.167&1285&1786&8731&1583508025&8.3.5.13.19.23.29.61&291&146 \cr
 & &16.3.5.7.11.19.47.257&1285&1848& & &32.9.13.61.73.97&4453&1872 \cr
\noalign{\hrule}
 & &7.11.31.59.67.167&6497&16350& & &9.25.7.11.17.19.283&207&5018 \cr
8714&1575780437&4.3.25.7.73.89.109&201&310&8732&1583660925&4.81.7.13.23.193&187&380 \cr
 & &16.9.125.31.67.89&1125&712& & &32.5.11.13.17.19.23&299&16 \cr
\noalign{\hrule}
 & &27.11.17.43.53.137&259&1160& & &9.5.11.169.29.653&12257&6680 \cr
8715&1576413927&16.9.5.7.29.37.137&313&646&8733&1584174735&16.3.25.7.17.103.167&1157&1682 \cr
 & &64.5.17.19.29.313&2755&10016& & &64.13.841.89.103&2581&3296 \cr
\noalign{\hrule}
 & &3.49.11.361.37.73&425&278& & &3.7.11.43.67.2381&347&390 \cr
8716&1576673637&4.25.11.17.19.73.139&1143&98&8734&1584581691&4.9.5.7.13.347.2381&28633&11966 \cr
 & &16.9.5.49.127.139&635&3336& & &16.5.11.19.31.137.193&26441&23560 \cr
\noalign{\hrule}
 & &81.5.11.23.73.211&2899&3014& & &3.17.29.37.59.491&3637&4710 \cr
8717&1578268395&4.121.13.137.211.223&27945&962&8735&1585270587&4.9.5.59.157.3637&491&3146 \cr
 & &16.243.5.169.23.37&507&296& & &16.121.13.157.491&2041&968 \cr
\noalign{\hrule}
 & &9.7.17.23.139.461&1577&1160& & &9.5.11.17.19.47.211&2107&1480 \cr
8718&1578458007&16.3.5.19.29.83.461&373&88&8736&1585579545&16.3.25.49.37.43.47&2071&3146 \cr
 & &256.11.29.83.373&26477&47744& & &64.49.121.13.19.109&5341&4576 \cr
\noalign{\hrule}
 & &9.5.11.97.167.197&713&122& & &9.25.7.37.73.373&1207&88 \cr
8719&1579645485&4.3.11.23.31.61.97&497&2510&8737&1586769975&16.3.5.11.17.71.73&8917&9698 \cr
 & &16.5.7.23.71.251&17821&1288& & &64.13.37.241.373&241&416 \cr
\noalign{\hrule}
 & &27.25.13.43.53.79&51791&50764& & &25.1681.179.211&39897&2128 \cr
8720&1579859775&8.3.5.343.37.67.773&13&2332&8738&1587242225&32.9.7.11.13.19.31&205&422 \cr
 & &64.49.11.13.37.53&407&1568& & &128.3.5.13.41.211&13&192 \cr
\noalign{\hrule}
 & &5.13.17.19.73.1031&473&558& & &3.5.121.41.83.257&323&282 \cr
8721&1580146685&4.9.11.13.19.31.43.73&5155&994&8739&1587346365&4.9.17.19.47.83.257&12389&310 \cr
 & &16.3.5.7.31.71.1031&497&744& & &16.5.13.19.31.953&18107&3224 \cr
\noalign{\hrule}
 & &11.31.43.47.2293&46657&24426& & &9.5.13.19.251.569&4945&176 \cr
8722&1580246173&4.9.13.23.37.59.97&1405&2666&8740&1587433185&32.25.11.13.23.43&137&162 \cr
 & &16.3.5.31.37.43.281&1405&888& & &128.81.11.43.137&13563&2752 \cr
\noalign{\hrule}
 & &3.25.7.11.17.89.181&11327&4782& & &81.7.17.79.2087&767&1320 \cr
8723&1581501075&4.9.5.47.241.797&3077&908&8741&1589210847&16.243.5.11.13.17.59&79&164 \cr
 & &32.17.47.181.227&227&752& & &128.11.13.41.59.79&5863&3776 \cr
\noalign{\hrule}
 & &9.41.43.263.379&85&44& & &3.11.47.61.107.157&205&266 \cr
8724&1581574959&8.3.5.11.17.263.379&437&700&8742&1589370189&4.5.7.11.19.41.47.107&793&2826 \cr
 & &64.125.7.11.17.19.23&40375&56672& & &16.9.5.13.41.61.157&205&312 \cr
\noalign{\hrule}
 & &9.41.43.263.379&437&700& & &9.289.29.47.449&2915&10106 \cr
8725&1581574959&8.3.25.7.19.23.41.43&85&44&8743&1591778187&4.5.11.17.31.53.163&1347&446 \cr
 & &64.125.7.11.17.19.23&40375&56672& & &16.3.5.31.223.449&1115&248 \cr
\noalign{\hrule}
 & &3.5.11.13.67.101.109&2499&1186& & &27.25.49.17.19.149&79&754 \cr
8726&1582158435&4.9.49.17.109.593&3595&1742&8744&1591800525&4.13.19.29.79.149&825&676 \cr
 & &16.5.49.13.67.719&719&392& & &32.3.25.11.2197.29&2197&5104 \cr
\noalign{\hrule}
 & &27.11.19.23.89.137&11515&11002& & &5.7.11.13.17.97.193&8587&10134 \cr
8727&1582517277&4.5.49.47.137.5501&35351&3156&8745&1592876285&4.9.5.11.31.277.563&199&364 \cr
 & &32.3.7.23.29.53.263&7627&5936& & &32.3.7.13.31.199.277&8587&9552 \cr
\noalign{\hrule}
 & &3.7.121.13.17.2819&2625&194& & &9.25.49.23.61.103&4339&9386 \cr
8728&1583040459&4.9.125.49.11.97&1139&3614&8746&1593211725&4.13.361.23.4339&671&5010 \cr
 & &16.5.13.17.67.139&139&2680& & &16.3.5.11.19.61.167&209&1336 \cr
\noalign{\hrule}
 & &29.31.631.2791&1845&946& & &9.7.11.13.17.101.103&815&94 \cr
8729&1583247779&4.9.5.11.41.43.631&377&254&8747&1593250659&4.5.11.13.17.47.163&105&116 \cr
 & &16.3.5.11.13.29.43.127&8385&11176& & &32.3.25.7.29.47.163&7661&11600 \cr
\noalign{\hrule}
 & &9.5.7.13.19.47.433&473&382& & &5.13.19.773.1669&204189&208054 \cr
8730&1583409555&4.11.43.47.191.433&475&42&8748&1593319195&4.3.49.11.29.193.2347&351&1772 \cr
 & &16.3.25.7.19.43.191&955&344& & &32.81.13.443.2347&35883&37552 \cr
\noalign{\hrule}
}%
}
$$
\eject
\vglue -23 pt
\noindent\hskip 1 in\hbox to 6.5 in{\ 8749 -- 8784 \hfill\fbd 1593645207 -- 1612695395\frb}
\vskip -9 pt
$$
\vbox{
\nointerlineskip
\halign{\strut
    \vrule \ \ \hfil \frb #\ 
   &\vrule \hfil \ \ \fbb #\frb\ 
   &\vrule \hfil \ \ \frb #\ \hfil
   &\vrule \hfil \ \ \frb #\ 
   &\vrule \hfil \ \ \frb #\ \ \vrule \hskip 2 pt
   &\vrule \ \ \hfil \frb #\ 
   &\vrule \hfil \ \ \fbb #\frb\ 
   &\vrule \hfil \ \ \frb #\ \hfil
   &\vrule \hfil \ \ \frb #\ 
   &\vrule \hfil \ \ \frb #\ \vrule \cr%
\noalign{\hrule}
 & &3.7.11.29.233.1021&38695&35632& & &9.25.7.121.13.647&1073&2162 \cr
8749&1593645207&32.5.7.17.71.109.131&145&1062&8767&1602926325&4.5.7.13.23.29.37.47&459&836 \cr
 & &128.9.25.29.59.109&6431&4800& & &32.27.11.17.19.23.47&3243&5168 \cr
\noalign{\hrule}
 & &9.11.37.47.59.157&43121&58772& & &9.25.11.13.19.43.61&1&474 \cr
8750&1594727343&8.7.13.31.107.2099&1745&354&8768&1603505475&4.27.13.61.79&817&830 \cr
 & &32.3.5.7.31.59.349&1745&3472& & &16.5.19.43.79.83&83&632 \cr
\noalign{\hrule}
 & &9.11.13.19.37.41.43&231&250& & &243.5.49.11.31.79&437&778 \cr
8751&1595093643&4.27.125.7.121.41.43&5369&166&8769&1603813365&4.49.19.23.79.389&775&726 \cr
 & &16.25.49.13.59.83&4067&11800& & &16.3.25.121.23.31.389&1945&2024 \cr
\noalign{\hrule}
 & &25.59.61.113.157&116379&123046& & &81.5.7.43.59.223&6127&7030 \cr
8752&1596246475&4.9.7.11.17.47.67.193&137&2260&8770&1603904085&4.27.25.11.19.37.557&59&616 \cr
 & &32.3.5.7.67.113.137&1407&2192& & &64.7.121.19.37.59&2299&1184 \cr
\noalign{\hrule}
 & &3.7.11.23.29.43.241&1625&3044& & &25.13.17.41.73.97&40067&1158 \cr
8753&1596699951&8.125.13.241.761&3393&2632&8771&1604023525&4.3.103.193.389&143&246 \cr
 & &128.9.5.7.169.29.47&2535&3008& & &16.9.11.13.41.193&1737&88 \cr
\noalign{\hrule}
 & &9.11.13.17.19.23.167&40379&5810& & &3.17.31.379.2677&1859&9890 \cr
8754&1596707541&4.5.7.83.149.271&1023&874&8772&1604055723&4.5.11.169.17.23.43&271&288 \cr
 & &16.3.5.11.19.23.31.83&415&248& & &256.9.5.11.13.23.271&44715&38272 \cr
\noalign{\hrule}
 & &3.11.13.17.439.499&11305&5816& & &5.11.31.41.53.433&387&46 \cr
8755&1597611873&16.5.7.289.19.727&4563&928&8773&1604249845&4.9.5.23.41.43.53&1279&484 \cr
 & &1024.27.7.169.29&1827&6656& & &32.3.121.23.1279&1279&12144 \cr
\noalign{\hrule}
 & &25.13.19.83.3119&4257&1138& & &25.11.73.157.509&871&1674 \cr
8756&1598565475&4.9.5.11.19.43.569&12091&12376&8774&1604253475&4.27.5.13.31.67.157&451&334 \cr
 & &64.3.7.11.13.17.107.113&38199&39776& & &16.3.11.31.41.67.167&20541&16616 \cr
\noalign{\hrule}
 & &25.343.121.23.67&12123&10858& & &27.7.121.29.41.59&5353&1786 \cr
8757&1598903075&4.27.5.11.61.89.449&3283&5528&8775&1604283219&4.9.7.19.47.53.101&619&290 \cr
 & &64.3.49.61.67.691&2073&1952& & &16.5.19.29.53.619&3095&8056 \cr
\noalign{\hrule}
 & &3.25.11.13.17.31.283&5429&9004& & &3.5.37.73.173.229&6887&5742 \cr
8758&1599537225&8.31.61.89.2251&2071&180&8776&1605082755&4.27.11.29.37.71.97&2993&5620 \cr
 & &64.9.5.19.89.109&1691&10464& & &32.5.41.73.97.281&3977&4496 \cr
\noalign{\hrule}
 & &3.7.29.41.139.461&1995&2036& & &3.7.19.29.89.1559&1625&66 \cr
8759&1599988551&8.9.5.49.19.461.509&451&10&8777&1605487821&4.9.125.7.11.13.29&1007&1268 \cr
 & &32.25.11.19.41.509&5225&8144& & &32.5.11.19.53.317&1585&9328 \cr
\noalign{\hrule}
 & &3.5.13.17.41.61.193&33&28& & &49.11.19.233.673&279&260 \cr
8760&1600127295&8.9.7.11.13.17.41.193&38125&15544&8778&1605880969&8.9.5.13.31.233.673&13&686 \cr
 & &128.625.29.61.67&3625&4288& & &32.3.5.343.169.31&5915&1488 \cr
\noalign{\hrule}
 & &3.7.11.19.43.61.139&1405&124& & &179.37249.241&2945&40194 \cr
8761&1600216233&8.5.19.31.43.281&549&268&8779&1606884611&4.9.5.7.11.19.29.31&193&358 \cr
 & &64.9.5.31.61.67&335&2976& & &16.3.7.31.179.193&217&24 \cr
\noalign{\hrule}
 & &3.25.49.19.101.227&23449&46376& & &3.7.11.23.41.83.89&18727&44270 \cr
8762&1600877775&16.11.17.31.131.179&18767&14706&8780&1609132371&4.5.19.61.233.307&38663&32868 \cr
 & &64.9.49.19.43.383&1149&1376& & &32.9.11.23.1681.83&41&48 \cr
\noalign{\hrule}
 & &5.11.37.257.3061&513&2548& & &9.5.7.1331.23.167&3895&3728 \cr
8763&1600887695&8.27.49.13.19.257&851&80&8781&1610396865&32.25.11.19.23.41.233&67&5292 \cr
 & &256.9.5.13.23.37&117&2944& & &256.27.49.41.67&1407&5248 \cr
\noalign{\hrule}
 & &3.17.149.359.587&473&114& & &11.13.17.19.139.251&2813&450 \cr
8764&1601360067&4.9.11.17.19.43.149&4109&2470&8782&1611488021&4.9.25.11.19.29.97&71&556 \cr
 & &16.5.7.13.361.587&1805&728& & &32.3.5.29.71.139&145&3408 \cr
\noalign{\hrule}
 & &25.121.37.41.349&73&48& & &11.13.19.59.89.113&279&488 \cr
8765&1601534825&32.3.37.41.73.349&7953&4960&8783&1612167271&16.9.31.61.89.113&77&190 \cr
 & &2048.9.5.11.31.241&7471&9216& & &64.3.5.7.11.19.31.61&1891&3360 \cr
\noalign{\hrule}
 & &9.5.7.11.13.961.37&47&232& & &5.17.19.53.83.227&26609&48174 \cr
8766&1601665065&16.7.11.13.29.31.47&51&950&8784&1612695395&4.3.7.11.31.37.41.59&689&1140 \cr
 & &64.3.25.17.19.47&235&10336& & &32.9.5.7.13.19.37.53&819&592 \cr
\noalign{\hrule}
}%
}
$$
\eject
\vglue -23 pt
\noindent\hskip 1 in\hbox to 6.5 in{\ 8785 -- 8820 \hfill\fbd 1613273207 -- 1630904979\frb}
\vskip -9 pt
$$
\vbox{
\nointerlineskip
\halign{\strut
    \vrule \ \ \hfil \frb #\ 
   &\vrule \hfil \ \ \fbb #\frb\ 
   &\vrule \hfil \ \ \frb #\ \hfil
   &\vrule \hfil \ \ \frb #\ 
   &\vrule \hfil \ \ \frb #\ \ \vrule \hskip 2 pt
   &\vrule \ \ \hfil \frb #\ 
   &\vrule \hfil \ \ \fbb #\frb\ 
   &\vrule \hfil \ \ \frb #\ \hfil
   &\vrule \hfil \ \ \frb #\ 
   &\vrule \hfil \ \ \frb #\ \vrule \cr%
\noalign{\hrule}
 & &49.13.41.223.277&3663&7264& & &5.29.31.277.1303&211&66 \cr
8785&1613273207&64.9.11.37.41.227&2007&490&8803&1622384845&4.3.11.31.211.1303&667&636 \cr
 & &256.81.5.49.223&405&128& & &32.9.11.23.29.53.211&20889&19504 \cr
\noalign{\hrule}
 & &9.11.23.37.107.179&7003&380& & &9.25.11.13.29.37.47&2537&2680 \cr
8786&1613621097&8.3.5.11.19.47.149&5539&2954&8804&1622617425&16.3.125.29.43.59.67&121&5254 \cr
 & &32.7.29.191.211&40301&3248& & &64.121.37.67.71&781&2144 \cr
\noalign{\hrule}
 & &25.11.13.29.37.421&651&274& & &3.25.19.31.109.337&3861&6586 \cr
8787&1614945475&4.3.7.11.31.137.421&983&1404&8805&1622680275&4.81.11.13.19.37.89&145&64 \cr
 & &32.81.13.137.983&11097&15728& & &512.5.13.29.37.89&13949&22784 \cr
\noalign{\hrule}
 & &5.13.23.47.83.277&9207&3812& & &9.25.7.419.2459&541&716 \cr
8788&1615462615&8.27.11.23.31.953&2033&8450&8806&1622755575&8.3.179.541.2459&2041&418 \cr
 & &32.3.25.169.19.107&1235&5136& & &32.11.13.19.157.179&25597&47728 \cr
\noalign{\hrule}
 & &3.5.11.361.43.631&2119&1036& & &81.11.31.67.877&29825&28934 \cr
8789&1616177145&8.7.11.13.37.43.163&3645&2504&8807&1622982339&4.25.17.23.31.37.1193&73&702 \cr
 & &128.729.5.37.313&11581&15552& & &16.27.13.23.73.1193&15509&13432 \cr
\noalign{\hrule}
 & &27.5.7.121.67.211&1099&15236& & &49.11.43.89.787&225&314 \cr
8790&1616495265&8.49.13.157.293&3993&3700&8808&1623386611&4.9.25.43.157.787&501&286 \cr
 & &64.3.25.1331.13.37&715&1184& & &16.27.5.11.13.157.167&21195&17368 \cr
\noalign{\hrule}
 & &5.23.31.43.53.199&171&8386& & &25.11.289.107.191&2927&252 \cr
8791&1616802365&4.9.7.19.23.599&5771&5610&8809&1624230575&8.9.7.191.2927&1177&1750 \cr
 & &16.27.5.11.17.29.199&783&1496& & &32.3.125.49.11.107&49&240 \cr
\noalign{\hrule}
 & &9.5.7.13.29.53.257&2057&5396& & &9.13.31.347.1291&785&506 \cr
8792&1617561855&8.5.121.13.17.19.71&159&764&8810&1624812579&4.5.11.13.23.157.347&939&796 \cr
 & &64.3.17.19.53.191&3247&608& & &32.3.23.157.199.313&62287&57776 \cr
\noalign{\hrule}
 & &81.5.7.11.23.37.61&12149&12556& & &3.5.7.29.59.83.109&134279&128084 \cr
8793&1618844535&8.7.23.43.73.12149&11951&198&8811&1625338785&8.11.41.47.71.2857&16677&14750 \cr
 & &32.9.11.17.19.37.43&731&304& & &32.9.125.17.59.71.109&1207&1200 \cr
\noalign{\hrule}
 & &3.7.11.13.23.131.179&731&710& & &9.11.31.401.1321&1975&1634 \cr
8794&1619598981&4.5.13.17.23.43.71.179&1179&14036&8812&1625713749&4.25.19.43.79.1321&123&1198 \cr
 & &32.9.121.29.71.131&957&1136& & &16.3.19.41.79.599&11381&25912 \cr
\noalign{\hrule}
 & &5.7.79.479.1223&15863&21978& & &9.5.7.13.23.41.421&29&6 \cr
8795&1619784005&4.27.7.11.29.37.547&479&3350&8813&1625727285&4.27.13.29.41.421&6985&7406 \cr
 & &16.3.25.37.67.479&335&888& & &16.5.7.11.529.29.127&2921&2552 \cr
\noalign{\hrule}
 & &3.5.7.11.289.23.211&2529&208& & &7.29.47.227.751&5555&16224 \cr
8796&1619907135&32.27.5.13.17.281&1007&1288&8814&1626521057&64.3.5.7.11.169.101&6759&6254 \cr
 & &512.7.13.19.23.53&1007&3328& & &256.27.53.59.751&3127&3456 \cr
\noalign{\hrule}
 & &27.5.169.19.37.101&9599&9086& & &9.7.13.17.107.1093&1265&2014 \cr
8797&1619933445&4.7.11.169.29.59.331&1447&264&8815&1628308773&4.3.5.11.13.17.19.23.53&2849&2186 \cr
 & &64.3.121.331.1447&40051&46304& & &16.7.121.23.37.1093&851&968 \cr
\noalign{\hrule}
 & &9.5.23.67.103.227&1071&1298& & &3.7.11.13.17.19.23.73&583&804 \cr
8798&1621355445&4.81.5.7.11.17.59.67&5221&206&8816&1628577951&8.9.7.121.23.53.67&91561&94900 \cr
 & &16.7.11.23.103.227&11&56& & &64.25.13.19.61.73.79&1975&1952 \cr
\noalign{\hrule}
 & &9.5.11.13.19.89.149&1731&986& & &27.25.7.11.13.19.127&10609&1964 \cr
8799&1621356165&4.27.17.29.89.577&913&1490&8817&1630403775&8.3.5.10609.491&1073&1382 \cr
 & &16.5.11.17.29.83.149&493&664& & &32.29.37.103.691&25567&47792 \cr
\noalign{\hrule}
 & &27.19.29.73.1493&737&756& & &3.11.13.23.29.41.139&931&1070 \cr
8800&1621429353&8.729.7.11.29.67.73&1493&3610&8818&1630728957&4.5.49.11.13.19.41.107&9153&584 \cr
 & &32.5.11.361.67.1493&1045&1072& & &64.81.5.7.73.113&9855&25312 \cr
\noalign{\hrule}
 & &9.49.11.367.911&16571&6550& & &3.23.41.53.73.149&6467&1430 \cr
8801&1621868787&4.25.7.73.131.227&321&190&8819&1630864749&4.5.11.13.29.41.223&855&2044 \cr
 & &16.3.125.19.107.227&28375&16264& & &32.9.25.7.11.19.73&627&2800 \cr
\noalign{\hrule}
 & &19.29.37.251.317&49977&29590& & &3.11.17.109.149.179&525&674 \cr
8802&1622132429&4.81.5.11.269.617&95&712&8820&1630904979&4.9.25.7.17.179.337&763&848 \cr
 & &64.27.25.11.19.89&2225&9504& & &128.5.49.53.109.337&16513&16960 \cr
\noalign{\hrule}
}%
}
$$
\eject
\vglue -23 pt
\noindent\hskip 1 in\hbox to 6.5 in{\ 8821 -- 8856 \hfill\fbd 1632235501 -- 1654897335\frb}
\vskip -9 pt
$$
\vbox{
\nointerlineskip
\halign{\strut
    \vrule \ \ \hfil \frb #\ 
   &\vrule \hfil \ \ \fbb #\frb\ 
   &\vrule \hfil \ \ \frb #\ \hfil
   &\vrule \hfil \ \ \frb #\ 
   &\vrule \hfil \ \ \frb #\ \ \vrule \hskip 2 pt
   &\vrule \ \ \hfil \frb #\ 
   &\vrule \hfil \ \ \fbb #\frb\ 
   &\vrule \hfil \ \ \frb #\ \hfil
   &\vrule \hfil \ \ \frb #\ 
   &\vrule \hfil \ \ \frb #\ \vrule \cr%
\noalign{\hrule}
 & &13.17.37.433.461&58951&42930& & &3.5.13.23.41.79.113&739&1078 \cr
8821&1632235501&4.81.5.53.167.353&2431&2078&8839&1641541395&4.5.49.11.13.41.739&3713&18 \cr
 & &16.3.5.11.13.17.53.1039&8745&8312& & &16.9.7.11.47.79&33&2632 \cr
\noalign{\hrule}
 & &7.13.17.361.37.79&1385&4752& & &121.19.29.89.277&2867&5166 \cr
8822&1632399041&32.27.5.11.79.277&217&494&8840&1643640163&4.9.7.41.47.61.89&277&10 \cr
 & &128.3.5.7.11.13.19.31&1023&320& & &16.3.5.47.61.277&141&2440 \cr
\noalign{\hrule}
 & &169.29.43.61.127&1139&3762& & &27.5.7.89.113.173&1145&1258 \cr
8823&1632626021&4.9.11.17.19.67.127&2407&1010&8841&1644168645&4.25.7.17.37.173.229&143&1068 \cr
 & &16.3.5.19.29.83.101&7885&2424& & &32.3.11.13.17.89.229&2519&3536 \cr
\noalign{\hrule}
 & &9.25.11.13.17.29.103&1119&14& & &5.49.13.19.29.937&1333&396 \cr
8824&1633814325&4.27.5.7.29.373&409&374&8842&1644374095&8.9.5.7.11.29.31.43&1349&884 \cr
 & &16.11.17.373.409&409&2984& & &64.3.13.17.19.43.71&3053&1632 \cr
\noalign{\hrule}
 & &9.17.23.29.67.239&19&220& & &3.5.17.19.23.29.509&11&334 \cr
8825&1634142663&8.3.5.11.17.19.23.29&3571&3824&8843&1644892035&4.11.29.167.509&95&414 \cr
 & &256.19.239.3571&3571&2432& & &16.9.5.19.23.167&501&8 \cr
\noalign{\hrule}
 & &81.5.7.11.23.43.53&895&524& & &25.11.53.157.719&213&370 \cr
8826&1634624145&8.27.25.23.131.179&1927&2548&8844&1645269725&4.3.125.37.71.719&297&422 \cr
 & &64.49.13.41.47.131&25051&29344& & &16.81.11.37.71.211&17091&21016 \cr
\noalign{\hrule}
 & &9.25.23.37.83.103&19&4& & &13.59.83.103.251&1683&1580 \cr
8827&1636919775&8.3.5.19.37.83.103&9499&286&8845&1645827833&8.9.5.11.17.59.79.83&381&1030 \cr
 & &32.7.11.13.23.59&8437&112& & &32.27.25.79.103.127&10033&10800 \cr
\noalign{\hrule}
 & &5.23.43.103.3217&2793&424& & &11.19.149.227.233&6363&41080 \cr
8828&1638530695&16.3.5.49.19.43.53&103&198&8846&1647078631&16.9.5.7.13.79.101&29&50 \cr
 & &64.27.7.11.53.103&4081&864& & &64.3.125.13.29.101&10875&42016 \cr
\noalign{\hrule}
 & &343.13.19.23.841&929&492& & &11.13.53.101.2153&365&948 \cr
8829&1638758303&8.3.7.13.29.41.929&165&368&8847&1648076287&8.3.5.73.79.2153&1807&3960 \cr
 & &256.9.5.11.23.929&10219&5760& & &128.27.25.11.13.139&3753&1600 \cr
\noalign{\hrule}
 & &5.23.61.167.1399&4207&2808& & &5.11.23.29.179.251&66599&52794 \cr
8830&1638935495&16.27.7.13.167.601&4491&3322&8848&1648220365&4.9.7.13.47.109.419&1933&3190 \cr
 & &64.243.11.151.499&75349&85536& & &16.3.5.7.11.13.29.1933&1933&2184 \cr
\noalign{\hrule}
 & &25.11.13.17.149.181&3401&324& & &27.5.13.19.61.811&1023&212 \cr
8831&1639040975&8.81.11.13.19.179&3077&2830&8849&1649610495&8.81.11.31.53.61&2641&2300 \cr
 & &32.27.5.17.181.283&283&432& & &64.25.19.23.53.139&6095&4448 \cr
\noalign{\hrule}
 & &27.49.19.61.1069&70265&71912& & &9.7.11.23.29.43.83&13&680 \cr
8832&1639158633&16.5.7.13.23.47.89.101&209&120&8850&1649702439&16.5.13.17.43.83&149&66 \cr
 & &256.3.25.11.13.19.23.101&32825&32384& & &64.3.11.13.17.149&149&7072 \cr
\noalign{\hrule}
 & &3.11.13.59.67.967&1015&1886& & &11.19.73.241.449&345&104 \cr
8833&1639874379&4.5.7.11.23.29.41.59&5283&3572&8851&1650944713&16.3.5.11.13.19.23.73&449&354 \cr
 & &32.9.19.41.47.587&36613&28176& & &64.9.13.23.59.449&2691&1888 \cr
\noalign{\hrule}
 & &11.37.53.139.547&565&18& & &5.11.43.191.3659&40657&408 \cr
8834&1640107843&4.9.5.37.113.139&7771&7658&8852&1652825185&16.3.17.109.373&3225&3116 \cr
 & &16.3.5.7.19.409.547&2863&2280& & &128.9.25.19.41.43&779&2880 \cr
\noalign{\hrule}
 & &9.11.13.53.67.359&1247&1984& & &81.11.13.17.37.227&281&200 \cr
8835&1640679183&128.13.29.31.43.53&359&330&8853&1653855489&16.25.11.17.227.281&9657&5798 \cr
 & &512.3.5.11.31.43.359&1333&1280& & &64.9.5.13.29.37.223&1115&928 \cr
\noalign{\hrule}
 & &9.11.17.31.83.379&40631&3110& & &9.7.31.71.79.151&201&352 \cr
8836&1641206061&4.5.41.311.991&651&340&8854&1654110927&64.27.11.31.67.71&10465&12382 \cr
 & &32.3.25.7.17.31.41&1025&112& & &256.5.7.13.23.41.151&2665&2944 \cr
\noalign{\hrule}
 & &9.5.11.17.41.67.71&7&362& & &25.11.17.107.3307&9177&7358 \cr
8837&1641236355&4.7.11.17.67.181&1041&4118&8855&1654244075&4.3.5.7.11.13.19.23.283&3307&25098 \cr
 & &16.3.29.71.347&29&2776& & &16.9.47.89.3307&801&376 \cr
\noalign{\hrule}
 & &5.23.29.509.967&43289&15246& & &3.5.49.13.31.37.151&3933&748 \cr
8838&1641497005&4.9.7.121.73.593&2185&1966&8856&1654897335&8.27.11.17.19.23.37&581&122 \cr
 & &16.3.5.121.19.23.983&6897&7864& & &32.7.11.23.61.83&15433&1328 \cr
\noalign{\hrule}
}%
}
$$
\eject
\vglue -23 pt
\noindent\hskip 1 in\hbox to 6.5 in{\ 8857 -- 8892 \hfill\fbd 1654940735 -- 1676764859\frb}
\vskip -9 pt
$$
\vbox{
\nointerlineskip
\halign{\strut
    \vrule \ \ \hfil \frb #\ 
   &\vrule \hfil \ \ \fbb #\frb\ 
   &\vrule \hfil \ \ \frb #\ \hfil
   &\vrule \hfil \ \ \frb #\ 
   &\vrule \hfil \ \ \frb #\ \ \vrule \hskip 2 pt
   &\vrule \ \ \hfil \frb #\ 
   &\vrule \hfil \ \ \fbb #\frb\ 
   &\vrule \hfil \ \ \frb #\ \hfil
   &\vrule \hfil \ \ \frb #\ 
   &\vrule \hfil \ \ \frb #\ \vrule \cr%
\noalign{\hrule}
 & &5.7.17.23.31.47.83&19869&12584& & &27.125.343.11.131&1297&2476 \cr
8857&1654940735&16.3.7.121.13.37.179&423&830&8875&1668137625&8.3.125.619.1297&81991&80134 \cr
 & &64.27.5.11.13.47.83&297&416& & &32.7.13.17.53.103.389&85969&87344 \cr
\noalign{\hrule}
 & &27.5.7.13.19.47.151&347&1012& & &5.11.23.37.127.281&84537&93898 \cr
8858&1656546255&8.3.11.13.23.47.347&707&190&8876&1670330035&4.27.7.19.31.101.353&37&316 \cr
 & &32.5.7.19.101.347&347&1616& & &32.3.7.19.37.79.101&5757&8848 \cr
\noalign{\hrule}
 & &5.13.19.73.79.233&54281&35874& & &5.11.13.29.61.1321&8559&8614 \cr
8859&1659483085&4.9.17.31.103.1993&77&26&8877&1670847035&4.27.29.59.61.73.317&1391&17312 \cr
 & &16.3.7.11.13.31.1993&13951&8184& & &256.3.13.73.107.541&39493&41088 \cr
\noalign{\hrule}
 & &27.25.13.19.23.433&121&554& & &81.125.7.17.19.73&11269&5356 \cr
8860&1660414275&4.121.13.19.23.277&1299&1000&8878&1671161625&8.13.17.59.103.191&6615&6424 \cr
 & &64.3.125.277.433&277&160& & &128.27.5.49.11.73.103&721&704 \cr
\noalign{\hrule}
 & &9.25.121.17.37.97&54119&35606& & &27.5.7.11.13.83.149&901&1036 \cr
8861&1661078925&4.13.19.23.181.937&1613&2550&8879&1671214545&8.49.11.17.37.53.83&8667&6854 \cr
 & &16.3.25.13.17.19.1613&1613&1976& & &32.81.23.53.107.149&2461&2544 \cr
\noalign{\hrule}
 & &27.25.11.23.37.263&1037&278& & &27.7.121.107.683&523&226 \cr
8862&1661811525&4.9.5.17.37.61.139&1847&1298&8880&1671289389&4.11.113.523.683&963&280 \cr
 & &16.11.59.139.1847&8201&14776& & &64.9.5.7.107.523&523&160 \cr
\noalign{\hrule}
 & &9.5.109.137.2473&3757&11176& & &5.13.19.23.229.257&2581&396 \cr
8863&1661818905&16.3.5.11.13.289.127&959&2356&8881&1671719465&8.9.11.29.89.257&875&104 \cr
 & &128.7.17.19.31.137&3689&1216& & &128.3.125.7.13.29&5075&192 \cr
\noalign{\hrule}
 & &9.11.13.23.89.631&893&13620& & &9.5.11.23.29.37.137&20877&15808 \cr
8864&1662362559&8.27.5.19.47.227&181&46&8882&1673606385&128.27.13.19.6959&145&6814 \cr
 & &32.19.23.47.181&8507&304& & &512.5.29.3407&3407&256 \cr
\noalign{\hrule}
 & &5.71.97.109.443&83919&73346& & &9.13.43.47.73.97&9725&8464 \cr
8865&1662762845&4.3.7.11.169.31.2543&771&1772&8883&1674352017&32.25.529.73.389&645&1034 \cr
 & &32.9.13.31.257.443&3627&4112& & &128.3.125.11.23.43.47&1375&1472 \cr
\noalign{\hrule}
 & &3.25.11.83.107.227&881&1794& & &3.5.11.23.31.43.331&17&1006 \cr
8866&1663189275&4.9.13.23.227.881&581&1462&8884&1674441285&4.5.17.331.503&14319&13816 \cr
 & &16.7.13.17.23.43.83&5083&2408& & &64.9.11.37.43.157&471&1184 \cr
\noalign{\hrule}
 & &9.5.7.19.23.43.281&781&2186& & &9.5.11.23.37.41.97&1471&2506 \cr
8867&1663285365&4.3.7.11.19.71.1093&65&562&8885&1675291365&4.7.11.37.179.1471&689&2160 \cr
 & &16.5.13.281.1093&1093&104& & &128.27.5.13.53.179&6981&3392 \cr
\noalign{\hrule}
 & &9.5.11.73.191.241&5887&8056& & &27.11.19.31.61.157&7567&10550 \cr
8868&1663330185&16.5.7.11.19.841.53&27&292&8886&1675333341&4.25.7.23.31.47.211&171&884 \cr
 & &128.27.7.19.29.73&1653&448& & &32.9.5.7.13.17.19.47&1645&3536 \cr
\noalign{\hrule}
 & &9.11.17.19.23.31.73&11663&11916& & &27.7.19.439.1063&535&528 \cr
8869&1664370873&8.81.31.107.109.331&9545&716&8887&1675765287&32.81.5.11.19.107.439&163&8504 \cr
 & &64.5.23.83.107.179&14857&17120& & &512.5.11.163.1063&1793&1280 \cr
\noalign{\hrule}
 & &9.7.11.29.41.43.47&2183&2470& & &9.11.41.127.3251&3655&404 \cr
8870&1665257517&4.5.13.19.29.37.43.59&555&4&8888&1675867743&8.5.17.43.101.127&3251&1092 \cr
 & &32.3.25.1369.59&1369&23600& & &64.3.5.7.13.3251&35&416 \cr
\noalign{\hrule}
 & &3.7.13.43.139.1021&1045&2018& & &9.11.13.53.79.311&4811&1390 \cr
8871&1665987141&4.5.11.13.19.43.1009&225&784&8889&1675876059&4.5.17.79.139.283&267&128 \cr
 & &128.9.125.49.11.19&9625&3648& & &1024.3.17.89.283&25187&8704 \cr
\noalign{\hrule}
 & &3.5.19.71.281.293&2563&2776& & &3.5.19.29.43.53.89&10071&9064 \cr
8872&1666008255&16.5.11.233.293.347&549&2014&8890&1676398215&16.81.11.29.103.373&32129&6290 \cr
 & &64.9.19.53.61.347&9699&11104& & &64.5.17.361.37.89&703&544 \cr
\noalign{\hrule}
 & &81.19.23.103.457&23621&23450& & &9.49.17.19.79.149&825&676 \cr
8873&1666172187&4.9.25.7.13.529.67.79&5027&266&8891&1676696553&8.27.25.49.11.169.17&79&754 \cr
 & &16.25.49.11.13.19.457&1225&1144& & &32.11.2197.29.79&2197&5104 \cr
\noalign{\hrule}
 & &9.11.17.19.31.1681&2813&2984& & &49.11.17.31.5903&20567&20754 \cr
8874&1666353447&16.29.1681.97.373&42465&6284&8892&1676764859&4.9.7.31.131.157.1153&125599&18370 \cr
 & &128.3.5.19.149.1571&7855&9536& & &16.3.5.11.29.61.71.167&50935&49416 \cr
\noalign{\hrule}
}%
}
$$
\eject
\vglue -23 pt
\noindent\hskip 1 in\hbox to 6.5 in{\ 8893 -- 8928 \hfill\fbd 1676876773 -- 1698434793\frb}
\vskip -9 pt
$$
\vbox{
\nointerlineskip
\halign{\strut
    \vrule \ \ \hfil \frb #\ 
   &\vrule \hfil \ \ \fbb #\frb\ 
   &\vrule \hfil \ \ \frb #\ \hfil
   &\vrule \hfil \ \ \frb #\ 
   &\vrule \hfil \ \ \frb #\ \ \vrule \hskip 2 pt
   &\vrule \ \ \hfil \frb #\ 
   &\vrule \hfil \ \ \fbb #\frb\ 
   &\vrule \hfil \ \ \frb #\ \hfil
   &\vrule \hfil \ \ \frb #\ 
   &\vrule \hfil \ \ \frb #\ \vrule \cr%
\noalign{\hrule}
 & &11.13.29.281.1439&3355&4794& & &243.11.29.47.463&5377&8050 \cr
8893&1676876773&4.3.5.121.13.17.47.61&1211&2784&8911&1686847437&4.25.7.19.23.47.283&261&214 \cr
 & &256.9.7.29.61.173&10899&7808& & &16.9.7.23.29.107.283&6509&5992 \cr
\noalign{\hrule}
 & &9.11.19.43.89.233&85&4342& & &3.25.11.71.83.347&5549&344 \cr
8894&1677270771&4.5.13.17.89.167&9417&10252&8912&1687018575&16.5.11.31.43.179&639&694 \cr
 & &32.3.11.43.73.233&73&16& & &64.9.71.179.347&179&96 \cr
\noalign{\hrule}
 & &9.17.37.43.61.113&143&874& & &27.5.11.13.19.43.107&6901&2086 \cr
8895&1677914739&4.11.13.19.23.37.61&339&820&8913&1687623795&4.3.7.13.67.103.149&107&94 \cr
 & &32.3.5.11.23.41.113&2255&368& & &16.7.47.103.107.149&7003&5768 \cr
\noalign{\hrule}
 & &9.25.11.29.67.349&1495&448& & &9.25.11.29.43.547&407&668 \cr
8896&1678314825&128.3.125.7.11.13.23&1073&698&8914&1688219775&8.121.37.167.547&20223&16 \cr
 & &512.13.29.37.349&481&256& & &256.27.7.107&21&13696 \cr
\noalign{\hrule}
 & &1331.23.29.31.61&2337&446& & &25.7.13.23.41.787&1111&324 \cr
8897&1678786307&4.3.11.19.29.41.223&549&230&8915&1688370775&8.81.5.11.13.23.101&697&798 \cr
 & &16.27.5.23.61.223&1115&216& & &32.243.7.11.17.19.41&4131&3344 \cr
\noalign{\hrule}
 & &13.23.37.47.3229&40987&990& & &7.11.17.41.163.193&117&76 \cr
8898&1678954069&4.9.5.11.17.2411&1189&1222&8916&1688373071&8.9.7.11.13.17.19.163&899&410 \cr
 & &16.3.5.13.17.29.41.47&1189&2040& & &32.3.5.13.19.29.31.41&6045&8816 \cr
\noalign{\hrule}
 & &5.19.41.97.4447&235&4212& & &7.11.109.271.743&2465&2736 \cr
8899&1680143305&8.81.25.13.19.47&1639&1886&8917&1689955729&32.9.5.11.17.19.29.109&115&1084 \cr
 & &32.27.11.23.41.149&6831&2384& & &256.3.25.23.29.271&2001&3200 \cr
\noalign{\hrule}
 & &3.25.49.13.29.1213&8497&7272& & &3.25.7.11.13.47.479&113&498 \cr
8900&1680581175&16.27.841.101.293&3443&26150&8918&1690163475&4.9.5.83.113.479&2303&2782 \cr
 & &64.25.11.313.523&5753&10016& & &16.49.13.47.83.107&581&856 \cr
\noalign{\hrule}
 & &5.17.109.419.433&3289&3834& & &5.7.121.13.31.991&1593&5344 \cr
8901&1680921155&4.27.11.13.23.71.433&3953&1676&8919&1691344655&64.27.5.13.59.167&121&56 \cr
 & &32.3.59.67.71.419&3953&3408& & &1024.9.7.121.167&1503&512 \cr
\noalign{\hrule}
 & &125.19.41.61.283&41447&6072& & &9.5.7.13.17.109.223&173&282 \cr
8902&1680984625&16.3.7.11.23.31.191&1025&312&8920&1692131805&4.27.17.47.173.223&4481&1540 \cr
 & &256.9.25.11.13.41&1287&128& & &32.5.7.11.47.4481&4481&8272 \cr
\noalign{\hrule}
 & &9.5.31.251.4801&3311&11092& & &9.25.7.11.19.53.97&43&62 \cr
8903&1681046145&8.3.5.7.11.43.47.59&1357&1228&8921&1692288675&4.3.5.11.31.43.53.97&931&136 \cr
 & &64.7.23.3481.307&80063&68768& & &64.49.17.19.31.43&1333&3808 \cr
\noalign{\hrule}
 & &3.49.11.23.53.853&235&304& & &5.13.127.421.487&33&454 \cr
8904&1681367919&32.5.19.47.53.853&819&1672&8922&1692497885&4.3.5.11.13.127.227&421&294 \cr
 & &512.9.5.7.11.13.361&4693&3840& & &16.9.49.227.421&2043&392 \cr
\noalign{\hrule}
 & &3.5.37.47.59.1093&71&70& & &9.7.121.31.67.107&1709&5460 \cr
8905&1682143395&4.25.7.37.59.71.1093&65081&594&8923&1694127897&8.27.5.49.13.1709&10447&11770 \cr
 & &16.27.7.11.151.431&27153&13288& & &32.25.11.31.107.337&337&400 \cr
\noalign{\hrule}
 & &27.7.43.317.653&15631&2000& & &27.49.13.83.1187&10411&42460 \cr
8906&1682297127&32.125.343.11.29&859&516&8924&1694462679&8.5.11.29.193.359&1241&882 \cr
 & &256.3.29.43.859&859&3712& & &32.9.5.49.17.29.73&1241&2320 \cr
\noalign{\hrule}
 & &9.11.17.43.67.347&323&280& & &9.5.49.11.13.19.283&15707&2122 \cr
8907&1682506881&16.5.7.11.289.19.347&837&4654&8925&1695448755&4.7.113.139.1061&7857&7850 \cr
 & &64.27.5.7.13.31.179&18795&12896& & &16.81.25.97.157.1061&76145&76392 \cr
\noalign{\hrule}
 & &25.11.47.139.937&2907&41132& & &5.7.11.17.31.61.137&8003&12798 \cr
8908&1683390775&8.9.7.13.17.19.113&937&1210&8926&1695593515&4.81.17.53.79.151&217&2350 \cr
 & &32.3.5.121.17.937&187&48& & &16.3.25.7.31.47.53&235&1272 \cr
\noalign{\hrule}
 & &125.49.29.53.179&1381&156& & &9.7.17.409.3877&17719&44858 \cr
8909&1685128375&8.3.5.13.179.1381&3363&3542&8927&1698277203&4.11.13.29.47.2039&1325&714 \cr
 & &32.9.7.11.13.19.23.59&25783&20592& & &16.3.25.7.11.17.29.53&1325&2552 \cr
\noalign{\hrule}
 & &3.11.13.19.23.89.101&521&590& & &9.11.289.23.29.89&11425&14296 \cr
8910&1685194797&4.5.13.19.59.89.521&26361&4378&8928&1698434793&16.25.23.457.1787&951&836 \cr
 & &16.9.5.11.29.101.199&995&696& & &128.3.5.11.19.317.457&30115&29248 \cr
\noalign{\hrule}
}%
}
$$
\eject
\vglue -23 pt
\noindent\hskip 1 in\hbox to 6.5 in{\ 8929 -- 8964 \hfill\fbd 1698992801 -- 1716167115\frb}
\vskip -9 pt
$$
\vbox{
\nointerlineskip
\halign{\strut
    \vrule \ \ \hfil \frb #\ 
   &\vrule \hfil \ \ \fbb #\frb\ 
   &\vrule \hfil \ \ \frb #\ \hfil
   &\vrule \hfil \ \ \frb #\ 
   &\vrule \hfil \ \ \frb #\ \ \vrule \hskip 2 pt
   &\vrule \ \ \hfil \frb #\ 
   &\vrule \hfil \ \ \fbb #\frb\ 
   &\vrule \hfil \ \ \frb #\ \hfil
   &\vrule \hfil \ \ \frb #\ 
   &\vrule \hfil \ \ \frb #\ \vrule \cr%
\noalign{\hrule}
 & &11.17.47.61.3169&42705&11168& & &3.125.11.37.53.211&1337&412 \cr
8929&1698992801&64.9.5.13.73.349&1249&1598&8947&1706805375&8.5.7.103.191.211&14209&7524 \cr
 & &256.3.5.17.47.1249&1249&1920& & &64.9.11.13.19.1093&3279&7904 \cr
\noalign{\hrule}
 & &9.11.13.23.67.857&133&166& & &3.7.11.61.89.1361&1615&254 \cr
8930&1699659819&4.3.7.19.67.83.857&1495&4066&8948&1706828739&4.5.11.17.19.61.127&623&414 \cr
 & &16.5.7.13.361.23.107&3745&2888& & &16.9.5.7.23.89.127&345&1016 \cr
\noalign{\hrule}
 & &9.11.17.19.79.673&259&1760& & &3.5.7.67.107.2269&4719&2450 \cr
8931&1700121159&64.3.5.7.121.17.37&1027&1030&8949&1707978405&4.9.125.343.121.13&2231&856 \cr
 & &256.25.7.13.37.79.103&33475&33152& & &64.11.13.23.97.107&2231&4576 \cr
\noalign{\hrule}
 & &3.5.13.17.29.31.571&781&118& & &27.31.61.109.307&7469&820 \cr
8932&1701685635&4.5.11.59.71.571&1037&1818&8950&1708520391&8.5.7.11.31.41.97&2097&1880 \cr
 & &16.9.17.59.61.101&5959&1464& & &128.9.25.11.47.233&10951&17600 \cr
\noalign{\hrule}
 & &27.25.13.19.59.173&92579&101354& & &3.5.19.139.179.241&437&258 \cr
8933&1701762075&4.11.17.43.271.2153&6903&4750&8951&1708951485&4.9.361.23.43.241&32705&17182 \cr
 & &16.9.125.11.13.17.19.59&85&88& & &16.5.121.31.71.211&25531&17608 \cr
\noalign{\hrule}
 & &3.25.7.11.53.67.83&14711&10146& & &23.43.71.101.241&1577&1476 \cr
8934&1702083075&4.9.5.19.47.89.313&11&434&8952&1709200679&8.9.19.23.41.83.241&143&580 \cr
 & &16.7.11.19.31.313&313&4712& & &64.3.5.11.13.29.41.83&77285&87648 \cr
\noalign{\hrule}
 & &3.13.17.841.43.71&427&414& & &3.11.13.17.151.1553&41515&58598 \cr
8935&1702300899&4.27.7.17.23.43.61.71&44341&250&8953&1710230379&4.5.361.23.83.353&965&612 \cr
 & &16.125.7.11.29.139&1375&7784& & &32.9.25.17.19.23.193&10925&9264 \cr
\noalign{\hrule}
 & &81.11.13.19.71.109&4805&4024& & &243.5.7.11.101.181&2413&422 \cr
8936&1703175903&16.5.13.19.961.503&137&1098&8954&1710278955&4.3.19.101.127.211&12427&400 \cr
 & &64.9.61.137.503&8357&16096& & &128.25.289.43&215&18496 \cr
\noalign{\hrule}
 & &9.13.19.43.103.173&6655&784& & &5.7.37.101.103.127&187&702 \cr
8937&1703300391&32.3.5.49.1331.13&443&404&8955&1710929395&4.27.11.13.17.37.101&1025&692 \cr
 & &256.5.7.11.101.443&38885&56704& & &32.3.25.11.13.41.173&11245&21648 \cr
\noalign{\hrule}
 & &3.7.17.19.23.61.179&737&300& & &125.11.19.31.2113&27&182 \cr
8938&1703462271&8.9.25.7.11.67.179&157&1768&8956&1711265875&4.27.25.7.13.2113&4469&2356 \cr
 & &128.13.17.67.157&871&10048& & &32.9.19.31.41.109&369&1744 \cr
\noalign{\hrule}
 & &9.125.7.11.13.17.89&271&1954& & &9.5.11.13.23.43.269&833&1102 \cr
8939&1703827125&4.5.7.13.271.977&32701&30804&8957&1711973835&4.49.11.13.17.19.23.29&81&172 \cr
 & &32.3.17.53.151.617&8003&9872& & &32.81.7.17.19.29.43&3451&2736 \cr
\noalign{\hrule}
 & &9.5.11.17.31.47.139&91&326& & &81.53.311.1283&60203&43720 \cr
8940&1704231045&4.3.7.11.13.17.31.163&8035&5264&8958&1712962809&16.5.11.13.421.1093&1769&2862 \cr
 & &128.5.49.47.1607&1607&3136& & &64.27.5.13.29.53.61&1769&2080 \cr
\noalign{\hrule}
 & &9.5.11.41.137.613&923&310& & &3.25.7.11.13.19.1201&1673&472 \cr
8941&1704394395&4.25.11.13.31.41.71&613&162&8959&1713136425&16.5.49.19.59.239&3053&8658 \cr
 & &16.81.13.71.613&117&568& & &64.9.13.37.43.71&1591&6816 \cr
\noalign{\hrule}
 & &9.25.23.29.41.277&221&1164& & &27.13.19.37.53.131&6955&2108 \cr
8942&1704401775&8.27.5.13.17.29.97&77&2542&8960&1713206079&8.3.5.169.17.31.107&583&262 \cr
 & &32.7.11.13.31.41&77&6448& & &32.11.17.31.53.131&527&176 \cr
\noalign{\hrule}
 & &121.13.97.11173&5053&6120& & &3.25.49.169.31.89&5731&2546 \cr
8943&1704787513&16.9.5.11.13.17.31.163&6373&1940&8961&1713545925&4.5.11.13.19.67.521&279&994 \cr
 & &128.3.25.97.6373&6373&4800& & &16.9.7.31.71.521&521&1704 \cr
\noalign{\hrule}
 & &9.125.7.13.19.877&6413&4988& & &5.7.13.17.19.107.109&435&1826 \cr
8944&1705874625&8.3.5.7.121.29.43.53&247&548&8962&1714052795&4.3.25.11.29.83.109&4161&4886 \cr
 & &64.121.13.19.29.137&3509&4384& & &16.9.7.11.19.73.349&3839&5256 \cr
\noalign{\hrule}
 & &3.7.11.23.53.73.83&377&536& & &27.5.11.17.113.601&35359&32554 \cr
8945&1706147751&16.7.13.23.29.67.73&405&106&8963&1714463685&4.9.19.41.397.1861&4207&72094 \cr
 & &64.81.5.29.53.67&3915&2144& & &16.7.11.29.113.601&29&56 \cr
\noalign{\hrule}
 & &27.25.13.163.1193&1397&722& & &27.5.169.19.37.107&583&262 \cr
8946&1706377725&4.11.361.127.1193&533&660&8964&1716167115&4.9.11.19.37.53.131&6955&2108 \cr
 & &32.3.5.121.13.361.41&4961&5776& & &32.5.11.13.17.31.107&527&176 \cr
\noalign{\hrule}
}%
}
$$
\eject
\vglue -23 pt
\noindent\hskip 1 in\hbox to 6.5 in{\ 8965 -- 9000 \hfill\fbd 1716768845 -- 1732655925\frb}
\vskip -9 pt
$$
\vbox{
\nointerlineskip
\halign{\strut
    \vrule \ \ \hfil \frb #\ 
   &\vrule \hfil \ \ \fbb #\frb\ 
   &\vrule \hfil \ \ \frb #\ \hfil
   &\vrule \hfil \ \ \frb #\ 
   &\vrule \hfil \ \ \frb #\ \ \vrule \hskip 2 pt
   &\vrule \ \ \hfil \frb #\ 
   &\vrule \hfil \ \ \fbb #\frb\ 
   &\vrule \hfil \ \ \frb #\ \hfil
   &\vrule \hfil \ \ \frb #\ 
   &\vrule \hfil \ \ \frb #\ \vrule \cr%
\noalign{\hrule}
 & &5.11.19.53.139.223&207&1322& & &9.5.41.61.67.229&1753&748 \cr
8965&1716768845&4.9.19.23.53.661&139&298&8983&1726777935&8.3.11.17.229.1753&533&1220 \cr
 & &16.3.139.149.661&661&3576& & &64.5.11.13.17.41.61&187&416 \cr
\noalign{\hrule}
 & &7.11.31.37.53.367&2059&2022& & &3.25.13.61.71.409&2147&622 \cr
8966&1717892869&4.3.29.31.71.337.367&10545&98&8984&1727094525&4.19.113.311.409&869&1278 \cr
 & &16.9.5.49.19.37.71&1197&2840& & &16.9.11.71.79.311&3421&1896 \cr
\noalign{\hrule}
 & &121.31.71.6451&5101&1350& & &81.49.19.37.619&1835&22 \cr
8967&1718036771&4.27.25.71.5101&2373&2728&8985&1727138133&4.27.5.11.19.367&259&254 \cr
 & &64.81.5.7.11.31.113&2835&3616& & &16.7.11.37.127.367&1397&2936 \cr
\noalign{\hrule}
 & &27.7.169.361.149&3275&748& & &3.5.11.37.41.67.103&8773&4962 \cr
8968&1718074449&8.25.11.169.17.131&57&112&8986&1727354805&4.9.11.31.283.827&3275&5822 \cr
 & &256.3.5.7.17.19.131&655&2176& & &16.25.31.41.71.131&4061&2840 \cr
\noalign{\hrule}
 & &27.5.7.31.89.659&1199&754& & &125.13.43.79.313&2511&886 \cr
8969&1718181045&4.3.11.13.29.109.659&1085&3062&8987&1727799125&4.81.31.313.443&1837&2150 \cr
 & &16.5.7.31.109.1531&1531&872& & &16.9.25.11.31.43.167&1837&2232 \cr
\noalign{\hrule}
 & &7.11.29.31.103.241&112295&104364& & &9.5.29.41.43.751&45563&46810 \cr
8970&1718322529&8.9.5.13.37.223.607&31&638&8988&1727836965&4.3.25.7.23.31.151.283&10049&1276 \cr
 & &32.3.5.11.13.29.31.37&481&240& & &32.7.11.13.23.29.773&17779&16016 \cr
\noalign{\hrule}
 & &9.7.11.13.53.59.61&589&650& & &81.5.13.257.1277&245&3586 \cr
8971&1718439723&4.3.25.11.169.19.31.53&2183&7112&8989&1727915085&4.27.25.49.11.163&1943&2132 \cr
 & &64.5.7.19.37.59.127&3515&4064& & &32.7.11.13.29.41.67&13601&7216 \cr
\noalign{\hrule}
 & &25.11.83.127.593&11349&11476& & &13.31.71.131.461&1035&1166 \cr
8972&1718973575&8.9.13.19.97.151.593&667&74&8990&1727967683&4.9.5.11.13.23.53.461&28427&44872 \cr
 & &32.3.23.29.37.97.151&128501&135024& & &64.3.7.31.71.79.131&237&224 \cr
\noalign{\hrule}
 & &3.49.11.13.19.59.73&635&3526& & &27.11.19.31.41.241&241&538 \cr
8973&1720211493&4.5.11.13.41.43.127&589&1062&8991&1728512973&4.31.58081.269&24871&33210 \cr
 & &16.9.5.19.31.41.59&205&744& & &16.81.5.7.11.17.19.41&105&136 \cr
\noalign{\hrule}
 & &13.19.29.37.73.89&115&1272& & &27.19.23.101.1451&6479&6580 \cr
8974&1721906407&16.3.5.23.29.37.53&539&534&8992&1729155249&8.3.5.7.11.361.23.31.47&1&1082 \cr
 & &64.9.49.11.23.53.89&13409&14112& & &32.5.7.11.31.541&41657&2480 \cr
\noalign{\hrule}
 & &5.7.11.17.29.43.211&1971&494& & &9.125.7.11.13.29.53&59641&42734 \cr
8975&1722100765&4.27.11.13.19.43.73&833&844&8993&1730854125&4.19.23.43.73.929&16415&15486 \cr
 & &32.9.49.17.19.73.211&1197&1168& & &16.3.5.49.29.43.67.89&3827&3752 \cr
\noalign{\hrule}
 & &3.13.23.29.107.619&41525&1186& & &9.5.11.13.47.59.97&5983&42032 \cr
8976&1722919029&4.25.11.151.593&7337&7488&8994&1730892735&32.31.37.71.193&1217&1410 \cr
 & &512.9.121.13.23.29&363&256& & &128.3.5.31.47.1217&1217&1984 \cr
\noalign{\hrule}
 & &9.5.343.13.31.277&33&310& & &17.19.31.41.4217&273&254 \cr
8977&1723024485&4.27.25.11.13.961&11837&12188&8995&1731217661&4.3.7.13.41.127.4217&495&4712 \cr
 & &32.7.121.19.89.277&1691&1936& & &64.27.5.7.11.13.19.31&2457&1760 \cr
\noalign{\hrule}
 & &121.17.19.44087&23193&20894& & &7.11.13.19.29.43.73&2819&3330 \cr
8978&1723052221&4.27.17.31.337.859&7255&7348&8996&1731318589&4.9.5.19.29.37.2819&583&2236 \cr
 & &32.9.5.11.167.337.1451&242317&242640& & &32.3.5.11.13.37.43.53&795&592 \cr
\noalign{\hrule}
 & &9.7.11.59.149.283&731&908& & &5.11.29.823.1319&1209&386 \cr
8979&1724082129&8.3.7.17.43.227.283&36803&48380&8997&1731431515&4.3.13.31.193.1319&595&1914 \cr
 & &64.5.13.19.41.59.149&1235&1312& & &16.9.5.7.11.17.29.31&1071&248 \cr
\noalign{\hrule}
 & &27.5.59.293.739&427&1166& & &3.7.13.163.167.233&5&2 \cr
8980&1724637555&4.5.7.11.53.61.293&22909&12654&8998&1731500589&4.5.13.163.167.233&20515&18396 \cr
 & &16.9.19.31.37.739&589&296& & &32.9.25.7.11.73.373&20075&17904 \cr
\noalign{\hrule}
 & &9.5.7.19.29.61.163&4285&442& & &3.25.7.11.13.17.23.59&1623&148 \cr
8981&1725756795&4.25.13.17.19.857&5379&10904&8999&1731905175&8.9.13.17.37.541&1265&724 \cr
 & &64.3.11.29.47.163&47&352& & &64.5.11.23.37.181&181&1184 \cr
\noalign{\hrule}
 & &3.5.7.13.361.31.113&18177&14674& & &9.25.343.11.13.157&383&58 \cr
8982&1726155795&4.9.5.11.23.29.73.83&15029&6634&9000&1732655925&4.7.11.29.157.383&709&390 \cr
 & &16.7.11.19.31.107.113&107&88& & &16.3.5.13.383.709&709&3064 \cr
\noalign{\hrule}
}%
}
$$
\eject
\vglue -23 pt
\noindent\hskip 1 in\hbox to 6.5 in{\ 9001 -- 9036 \hfill\fbd 1732857027 -- 1746944599\frb}
\vskip -9 pt
$$
\vbox{
\nointerlineskip
\halign{\strut
    \vrule \ \ \hfil \frb #\ 
   &\vrule \hfil \ \ \fbb #\frb\ 
   &\vrule \hfil \ \ \frb #\ \hfil
   &\vrule \hfil \ \ \frb #\ 
   &\vrule \hfil \ \ \frb #\ \ \vrule \hskip 2 pt
   &\vrule \ \ \hfil \frb #\ 
   &\vrule \hfil \ \ \fbb #\frb\ 
   &\vrule \hfil \ \ \frb #\ \hfil
   &\vrule \hfil \ \ \frb #\ 
   &\vrule \hfil \ \ \frb #\ \vrule \cr%
\noalign{\hrule}
 & &3.11.71.757.977&745&1526& & &3.7.11.37.47.61.71&43&2300 \cr
9001&1732857027&4.5.7.109.149.977&563&414&9019&1739801679&8.25.7.23.43.47&923&1098 \cr
 & &16.9.5.7.23.109.563&38847&30520& & &32.9.13.23.61.71&39&368 \cr
\noalign{\hrule}
 & &3.5.11.13.31.67.389&1501&1112& & &3.5.11.13.47.61.283&1043&372 \cr
9002&1733059185&16.5.11.19.31.79.139&2173&468&9020&1740369345&8.9.7.13.31.47.149&7525&5588 \cr
 & &128.9.13.41.53.79&4187&7872& & &64.25.49.11.43.127&5461&7840 \cr
\noalign{\hrule}
 & &11.29.1063.5113&79985&68292& & &3.49.17.211.3301&1675&1626 \cr
9003&1733802961&8.9.5.7.17.271.941&5597&990&9021&1740580989&4.9.25.17.67.211.271&3301&286 \cr
 & &32.81.25.11.29.193&2025&3088& & &16.5.11.13.271.3301&1355&1144 \cr
\noalign{\hrule}
 & &7.1331.19.97.101&1005&326& & &7.59.1129.3733&114075&106172 \cr
9004&1734294331&4.3.5.19.67.101.163&1127&792&9022&1740612041&8.27.25.11.169.19.127&1189&46 \cr
 & &64.27.49.11.23.163&3749&6048& & &32.3.5.11.13.23.29.41&46371&20240 \cr
\noalign{\hrule}
 & &81.5.11.13.31.967&13969&43946& & &3.5.11.43.107.2293&8033&3432 \cr
9005&1736117955&4.7.43.61.73.229&7047&9670&9023&1740765345&16.9.121.13.29.277&2033&460 \cr
 & &16.243.5.7.29.967&87&56& & &128.5.19.23.29.107&437&1856 \cr
\noalign{\hrule}
 & &9.7.13.23.37.47.53&319&578& & &5.7.17.29.157.643&2047&1404 \cr
9006&1736149779&4.3.11.289.29.47.53&475&1274&9024&1741909505&8.27.5.13.23.89.157&19&176 \cr
 & &16.25.49.13.17.19.29&5075&2584& & &256.9.11.19.23.89&43263&11392 \cr
\noalign{\hrule}
 & &5.11.13.29.31.37.73&813&782& & &5.7.11.13.17.59.347&10291&1854 \cr
9007&1736162285&4.3.13.17.23.37.73.271&3249&4928&9025&1741945205&4.9.17.41.103.251&1001&750 \cr
 & &512.27.7.11.361.271&68229&69376& & &16.27.125.7.11.13.41&675&328 \cr
\noalign{\hrule}
 & &3.25.13.17.19.37.149&1431&506& & &9.11.19.37.79.317&901&980 \cr
9008&1736181525&4.81.11.17.19.23.53&185&1192&9026&1742917671&8.5.49.17.37.53.317&2607&19408 \cr
 & &64.5.11.23.37.149&253&32& & &256.3.7.11.79.1213&1213&896 \cr
\noalign{\hrule}
 & &5.169.71.103.281&173&342& & &9.121.13.23.53.101&7939&4282 \cr
9009&1736435285&4.9.19.71.173.281&9889&10062&9027&1742995683&4.3.13.17.467.2141&237281&235880 \cr
 & &16.81.11.13.19.29.31.43&72819&71896& & &64.5.121.37.53.5897&5897&5920 \cr
\noalign{\hrule}
 & &19.31.97.113.269&2925&8036& & &27.7.17.269.2017&473&1544 \cr
9010&1736671801&8.9.25.49.13.31.41&1067&452&9028&1743287049&16.3.11.43.193.269&1769&1190 \cr
 & &64.3.5.11.13.97.113&429&160& & &64.5.7.17.29.43.61&2623&4640 \cr
\noalign{\hrule}
 & &17.19.23.37.71.89&249&1100& & &11.41.79.173.283&265&186 \cr
9011&1736922487&8.3.25.11.17.83.89&483&928&9029&1744360211&4.3.5.31.53.173.283&55&228 \cr
 & &512.9.5.7.11.23.29&9135&2816& & &32.9.25.11.19.31.53&14725&7632 \cr
\noalign{\hrule}
 & &9.25.11.17.19.41.53&203&698& & &5.11.13.19.23.37.151&12609&15326 \cr
9012&1737150525&4.5.7.19.29.41.349&5777&8532&9030&1745686085&4.27.23.79.97.467&21793&15100 \cr
 & &32.27.7.53.79.109&2289&1264& & &32.9.25.19.31.37.151&155&144 \cr
\noalign{\hrule}
 & &49.11.289.19.587&123&710& & &3.11.19.41.113.601&132305&138746 \cr
9013&1737313963&4.3.5.11.17.19.41.71&15849&16172&9031&1745839491&4.5.47.173.401.563&133209&92554 \cr
 & &32.81.5.13.311.587&5265&4976& & &16.9.7.11.361.41.601&57&56 \cr
\noalign{\hrule}
 & &11.17.31.41.71.103&519&178& & &9.31.53.263.449&218075&214034 \cr
9014&1738131901&4.3.71.89.103.173&10725&1558&9032&1746152469&4.25.11.13.61.103.1039&789&68324 \cr
 & &16.9.25.11.13.19.41&225&1976& & &32.3.5.19.29.31.263&145&304 \cr
\noalign{\hrule}
 & &3.17.29.101.103.113&975&946& & &9.5.11.23.157.977&8645&13826 \cr
9015&1738622181&4.9.25.11.13.43.101.103&4487&58&9033&1746333765&4.3.25.7.13.19.31.223&737&2162 \cr
 & &16.5.7.11.13.29.641&8333&3080& & &16.7.11.23.31.47.67&3149&1736 \cr
\noalign{\hrule}
 & &9.7.23.89.97.139&383&590& & &3.13.17.37.257.277&145&132 \cr
9016&1738781163&4.5.59.89.97.383&4125&4508&9034&1746337359&8.9.5.11.17.29.37.257&3601&21842 \cr
 & &32.3.625.49.11.23.59&6875&6608& & &32.5.13.67.163.277&815&1072 \cr
\noalign{\hrule}
 & &3.5.11.13.61.97.137&3451&4906& & &9.5.13.23.31.53.79&437&358 \cr
9017&1738799205&4.7.121.13.17.29.223&305&3204&9035&1746418635&4.3.13.19.529.31.179&869&340 \cr
 & &32.9.5.7.17.61.89&357&1424& & &32.5.11.17.19.79.179&3553&2864 \cr
\noalign{\hrule}
 & &9.5.19.23.241.367&143&580& & &289.23.89.2953&2233&720 \cr
9018&1739310255&8.3.25.11.13.29.367&3473&1298&9036&1746944599&32.9.5.7.11.17.23.29&1847&890 \cr
 & &32.121.23.59.151&8909&1936& & &128.3.25.89.1847&1847&4800 \cr
\noalign{\hrule}
}%
}
$$
\eject
\vglue -23 pt
\noindent\hskip 1 in\hbox to 6.5 in{\ 9037 -- 9072 \hfill\fbd 1747056935 -- 1761547245\frb}
\vskip -9 pt
$$
\vbox{
\nointerlineskip
\halign{\strut
    \vrule \ \ \hfil \frb #\ 
   &\vrule \hfil \ \ \fbb #\frb\ 
   &\vrule \hfil \ \ \frb #\ \hfil
   &\vrule \hfil \ \ \frb #\ 
   &\vrule \hfil \ \ \frb #\ \ \vrule \hskip 2 pt
   &\vrule \ \ \hfil \frb #\ 
   &\vrule \hfil \ \ \fbb #\frb\ 
   &\vrule \hfil \ \ \frb #\ \hfil
   &\vrule \hfil \ \ \frb #\ 
   &\vrule \hfil \ \ \frb #\ \vrule \cr%
\noalign{\hrule}
 & &5.169.17.19.37.173&207&3080& & &9.47.101.181.227&1045&998 \cr
9037&1747056935&16.9.25.7.11.23.37&2941&3016&9055&1755359901&4.5.11.19.101.181.499&141&40 \cr
 & &256.3.11.13.17.29.173&319&384& & &64.3.25.11.19.47.499&9481&8800 \cr
\noalign{\hrule}
 & &3.121.23.31.43.157&287&1046& & &9.5.7.17.31.71.149&26983&25912 \cr
9038&1747287069&4.7.11.41.157.523&6095&342&9056&1756166895&16.121.31.41.79.223&1323&8236 \cr
 & &16.9.5.7.19.23.53&5035&168& & &128.27.49.29.41.71&1189&1344 \cr
\noalign{\hrule}
 & &3.5.7.19.53.61.271&1507&3642& & &121.13.59.127.149&135&14 \cr
9039&1747905285&4.9.11.53.137.607&65&542&9057&1756186861&4.27.5.7.13.59.127&5513&6044 \cr
 & &16.5.11.13.137.271&137&1144& & &32.3.5.37.149.1511&4533&2960 \cr
\noalign{\hrule}
 & &9.25.13.31.37.521&2607&3128& & &17.43.89.113.239&3945&118 \cr
9040&1747941975&16.27.5.11.13.17.23.79&24911&4924&9058&1757048413&4.3.5.59.113.263&8041&7476 \cr
 & &128.29.859.1231&35699&54976& & &32.9.7.11.17.43.89&63&176 \cr
\noalign{\hrule}
 & &3.7.13.17.19.43.461&221&178& & &3.5.7.11.59.107.241&5497&5738 \cr
9041&1747972317&4.169.289.89.461&44935&3906&9059&1757255115&4.11.19.23.59.151.239&1391&270 \cr
 & &16.9.5.7.11.19.31.43&341&120& & &16.27.5.13.23.107.239&2151&2392 \cr
\noalign{\hrule}
 & &3.13.17.23.29.59.67&11125&14134& & &3.7.11.73.127.821&533&1930 \cr
9042&1748099613&4.125.23.37.89.191&8651&1584&9060&1758254421&4.5.7.13.41.73.193&6817&7272 \cr
 & &128.9.25.11.41.211&15825&28864& & &64.9.17.41.101.401&49323&54944 \cr
\noalign{\hrule}
 & &3.13.19.61.101.383&451&698& & &3.5.11.17.773.811&533&278 \cr
9043&1748510283&4.11.41.61.101.349&1161&5000&9061&1758462915&4.11.13.41.139.773&161&612 \cr
 & &64.27.625.41.43&26875&11808& & &32.9.7.13.17.23.139&6279&2224 \cr
\noalign{\hrule}
 & &9.11.13.37.73.503&42455&7342& & &5.49.169.107.397&1305&1474 \cr
9044&1748522061&4.5.7.1213.3671&1229&2442&9062&1758842995&4.9.25.7.11.29.67.107&1961&286 \cr
 & &16.3.5.7.11.37.1229&1229&280& & &16.3.121.13.29.37.53&10527&15688 \cr
\noalign{\hrule}
 & &25.7.23.31.37.379&44847&16172& & &13.19.23.29.59.181&77035&87714 \cr
9045&1749719825&8.27.11.13.151.311&2789&3100&9063&1759354571&4.9.5.7.11.31.71.443&15067&1334 \cr
 & &64.9.25.11.31.2789&2789&3168& & &16.3.5.13.19.23.29.61&61&120 \cr
\noalign{\hrule}
 & &27.7.11.19.31.1429&1313&116& & &27.13.23.31.79.89&7&20 \cr
9046&1749854799&8.3.11.13.29.31.101&2623&6770&9064&1759599153&8.5.7.23.31.79.89&13689&3454 \cr
 & &32.5.43.61.677&3385&41968& & &32.81.11.169.157&471&2288 \cr
\noalign{\hrule}
 & &5.11.17.19.29.43.79&351&142& & &3.11.13.17.31.43.181&13205&12678 \cr
9047&1750083445&4.27.5.13.43.71.79&1837&3772&9065&1759603989&4.9.5.19.43.139.2113&8959&1606 \cr
 & &32.3.11.13.23.41.167&6847&14352& & &16.11.289.31.73.139&1241&1112 \cr
\noalign{\hrule}
 & &25.11.13.37.101.131&56079&37346& & &5.7.19.59.113.397&481&84 \cr
9048&1750130525&4.27.31.67.71.263&145&2222&9066&1760121335&8.3.49.13.19.37.59&879&242 \cr
 & &16.3.5.11.29.71.101&87&568& & &32.9.121.37.293&4477&42192 \cr
\noalign{\hrule}
 & &5.7.13.31.43.2887&15433&22098& & &121.23.29.113.193&247&3030 \cr
9049&1751008805&4.3.7.11.23.29.61.127&4005&5774&9067&1760138963&4.3.5.13.19.101.193&5313&7232 \cr
 & &16.27.5.23.89.2887&621&712& & &512.9.7.11.23.113&63&256 \cr
\noalign{\hrule}
 & &3.25.7.17.29.67.101&451&956& & &3.5.19.31.41.4861&2353&2508 \cr
9050&1751468775&8.5.11.17.29.41.239&117&202&9068&1760824335&8.9.11.13.361.41.181&7583&59206 \cr
 & &32.9.13.41.101.239&3107&1968& & &32.7.4229.7583&53081&67664 \cr
\noalign{\hrule}
 & &27.13.19.31.37.229&67&770& & &27.5.7.11.61.2777&1631&4408 \cr
9051&1751699547&4.5.7.11.13.67.229&5363&4218&9069&1760881815&16.3.5.49.19.29.233&10361&3604 \cr
 & &16.3.7.19.31.37.173&173&56& & &128.13.17.53.797&13549&44096 \cr
\noalign{\hrule}
 & &27.25.7.19.29.673&9361&8810& & &9.11.19.67.89.157&2275&3688 \cr
9052&1752138675&4.125.7.11.23.37.881&5217&14908&9070&1760975271&16.25.7.11.13.19.461&421&1884 \cr
 & &32.3.1369.47.3727&64343&59632& & &128.3.5.13.157.421&2105&832 \cr
\noalign{\hrule}
 & &125.11.23.37.1499&887&612& & &13.17.1871.4259&43587&11780 \cr
9053&1754017375&8.9.5.17.23.37.887&36421&38974&9071&1761058169&8.9.5.19.29.31.167&461&374 \cr
 & &32.3.7.121.13.43.1499&1419&1456& & &32.3.11.17.19.31.461&14291&10032 \cr
\noalign{\hrule}
 & &9.11.169.23.47.97&625&7318& & &3.5.37.43.223.331&23807&24138 \cr
9054&1754362467&4.3.625.11.3659&1417&2242&9072&1761547245&4.243.7.19.37.149.179&2453&41044 \cr
 & &16.25.13.19.59.109&2071&11800& & &32.11.19.31.223.331&341&304 \cr
\noalign{\hrule}
}%
}
$$
\eject
\vglue -23 pt
\noindent\hskip 1 in\hbox to 6.5 in{\ 9073 -- 9108 \hfill\fbd 1761752025 -- 1775512879\frb}
\vskip -9 pt
$$
\vbox{
\nointerlineskip
\halign{\strut
    \vrule \ \ \hfil \frb #\ 
   &\vrule \hfil \ \ \fbb #\frb\ 
   &\vrule \hfil \ \ \frb #\ \hfil
   &\vrule \hfil \ \ \frb #\ 
   &\vrule \hfil \ \ \frb #\ \ \vrule \hskip 2 pt
   &\vrule \ \ \hfil \frb #\ 
   &\vrule \hfil \ \ \fbb #\frb\ 
   &\vrule \hfil \ \ \frb #\ \hfil
   &\vrule \hfil \ \ \frb #\ 
   &\vrule \hfil \ \ \frb #\ \vrule \cr%
\noalign{\hrule}
 & &81.25.11.139.569&2183&1292& & &11.17.19.23.97.223&91&300 \cr
9073&1761752025&8.17.19.37.59.569&845&276&9091&1767663689&8.3.25.7.13.97.223&63&160 \cr
 & &64.3.5.169.17.23.37&6253&12512& & &512.27.125.49.13&43875&12544 \cr
\noalign{\hrule}
 & &3.5.7.11.59.103.251&61&1194& & &11.17.19.23.97.223&63&160 \cr
9074&1761752685&4.9.7.59.61.199&15499&16892&9092&1767663689&64.9.5.7.11.17.19.23&91&300 \cr
 & &32.11.41.103.1409&1409&656& & &512.27.125.49.13&43875&12544 \cr
\noalign{\hrule}
 & &25.29.31.277.283&4181&4026& & &5.11.19.53.59.541&793&252 \cr
9075&1761837725&4.3.5.11.37.61.113.277&22357&12108&9093&1767833815&8.9.7.13.53.59.61&2833&766 \cr
 & &32.9.11.79.283.1009&11099&11376& & &32.3.7.383.2833&8499&42896 \cr
\noalign{\hrule}
 & &11.17.47.193.1039&6371&42462& & &9.5.11.43.97.857&235&622 \cr
9076&1762431803&4.9.7.23.277.337&935&1004&9094&1769400765&4.25.11.47.97.311&2923&498 \cr
 & &32.3.5.11.17.251.337&5055&4016& & &16.3.37.47.79.83&3713&24568 \cr
\noalign{\hrule}
 & &9.7.41.43.59.269&337&50& & &3.25.11.19.37.43.71&609&1426 \cr
9077&1762776099&4.25.59.269.337&7093&8778&9095&1770663675&4.9.5.7.23.29.31.71&281&74 \cr
 & &16.3.5.7.11.19.41.173&865&1672& & &16.7.29.31.37.281&1967&7192 \cr
\noalign{\hrule}
 & &9.5.7.67.101.827&31211&52316& & &5.7.11.29.127.1249&4239&2006 \cr
9078&1762837335&8.11.529.29.41.59&945&1474&9096&1771025795&4.27.17.59.127.157&5203&4060 \cr
 & &32.27.5.7.121.29.67&363&464& & &32.3.5.7.121.17.29.43&561&688 \cr
\noalign{\hrule}
 & &181.38809.251&3311&42120& & &7.11.13.43.79.521&42389&5022 \cr
9079&1763131679&16.81.5.7.11.13.43&197&362&9097&1771606837&4.81.19.23.31.97&85&86 \cr
 & &64.27.7.181.197&27&224& & &16.9.5.17.23.31.43.97&20079&21080 \cr
\noalign{\hrule}
 & &25.13.23.409.577&2371&2946& & &25.7.121.13.47.137&5337&1102 \cr
9080&1764047675&4.3.491.577.2371&897&1474&9098&1772495725&4.9.5.13.19.29.593&1507&272 \cr
 & &16.9.11.13.23.67.491&5401&4824& & &128.3.11.17.29.137&87&1088 \cr
\noalign{\hrule}
 & &27.169.29.67.199&3355&1546& & &5.7.13.41.101.941&10449&43384 \cr
9081&1764315891&4.5.11.61.199.773&3027&838&9099&1772989855&16.243.11.17.29.43&2375&1882 \cr
 & &16.3.61.419.1009&25559&8072& & &64.27.125.19.941&513&800 \cr
\noalign{\hrule}
 & &27.5.49.11.127.191&3023&2134& & &5.11.29.73.97.157&7137&4324 \cr
9082&1765060605&4.5.7.121.97.3023&3629&606&9100&1773188615&8.9.5.11.13.23.47.61&103&2482 \cr
 & &16.3.19.97.101.191&1919&776& & &32.3.17.23.73.103&391&4944 \cr
\noalign{\hrule}
 & &9.13.17.29.127.241&4525&428& & &3.67.113.163.479&11275&11438 \cr
9083&1765442367&8.3.25.29.107.181&127&308&9101&1773362901&4.25.7.11.19.41.43.479&477&2 \cr
 & &64.5.7.11.107.127&3745&352& & &16.9.7.11.41.43.53&6837&25256 \cr
\noalign{\hrule}
 & &27.7.89.103.1019&1591&572& & &9.5.49.43.97.193&167&412 \cr
9084&1765481697&8.9.11.13.37.43.89&745&412&9102&1775031615&8.3.43.97.103.167&27887&23716 \cr
 & &64.5.11.43.103.149&2365&4768& & &64.49.121.79.353&9559&11296 \cr
\noalign{\hrule}
 & &9.5.11.17.37.53.107&101&434& & &9.5.19.47.163.271&5863&1798 \cr
9085&1765694205&4.7.11.17.31.53.101&5575&222&9103&1775092005&4.3.11.13.19.29.31.41&1075&524 \cr
 & &16.3.25.7.37.223&35&1784& & &32.25.11.31.43.131&20305&7568 \cr
\noalign{\hrule}
 & &9.11.41.43.67.151&161&290& & &9.5.13.73.197.211&1043&1700 \cr
9086&1765790829&4.3.5.7.23.29.67.151&25&176&9104&1775118735&8.125.7.17.149.197&2409&21034 \cr
 & &128.125.7.11.23.29&2875&12992& & &32.3.11.13.73.809&809&176 \cr
\noalign{\hrule}
 & &11.29.37.1681.89&351&1330& & &9.49.13.19.43.379&547&1364 \cr
9087&1765835027&4.27.5.7.13.19.29.37&267&748&9105&1775183319&8.3.11.31.379.547&6695&5054 \cr
 & &32.81.11.17.19.89&1377&304& & &32.5.7.11.13.361.103&1133&1520 \cr
\noalign{\hrule}
 & &3.5.7.13.17.269.283&14645&9834& & &3.11.23.47.157.317&149&8 \cr
9088&1766527035&4.9.25.11.29.101.149&1393&1132&9106&1775409537&16.11.23.149.317&4465&2826 \cr
 & &32.7.11.149.199.283&2189&2384& & &64.9.5.19.47.157&285&32 \cr
\noalign{\hrule}
 & &3.5.7.13.307.4217&1341&2876& & &121.101.131.1109&60479&73710 \cr
9089&1767154935&8.27.7.13.149.719&869&1588&9107&1775454659&4.81.5.7.13.197.307&1109&2882 \cr
 & &64.11.79.149.397&31363&52448& & &16.9.5.7.11.131.1109&45&56 \cr
\noalign{\hrule}
 & &3.25.13.61.113.263&15523&18942& & &7.31.53.317.487&185&132 \cr
9090&1767537525&4.9.5.7.11.361.41.43&23&338&9108&1775512879&8.3.5.7.11.31.37.487&1131&2278 \cr
 & &16.11.169.23.41.43&12259&3784& & &32.9.5.11.13.17.29.67&51255&66352 \cr
\noalign{\hrule}
}%
}
$$
\eject
\vglue -23 pt
\noindent\hskip 1 in\hbox to 6.5 in{\ 9109 -- 9144 \hfill\fbd 1776006855 -- 1791406773\frb}
\vskip -9 pt
$$
\vbox{
\nointerlineskip
\halign{\strut
    \vrule \ \ \hfil \frb #\ 
   &\vrule \hfil \ \ \fbb #\frb\ 
   &\vrule \hfil \ \ \frb #\ \hfil
   &\vrule \hfil \ \ \frb #\ 
   &\vrule \hfil \ \ \frb #\ \ \vrule \hskip 2 pt
   &\vrule \ \ \hfil \frb #\ 
   &\vrule \hfil \ \ \fbb #\frb\ 
   &\vrule \hfil \ \ \frb #\ \hfil
   &\vrule \hfil \ \ \frb #\ 
   &\vrule \hfil \ \ \frb #\ \vrule \cr%
\noalign{\hrule}
 & &9.5.7.19.43.67.103&961&3124& & &125.7.59.71.487&13601&20976 \cr
9109&1776006855&8.3.11.961.67.71&475&262&9127&1785037625&32.3.49.19.23.29.67&9735&11678 \cr
 & &32.25.19.961.131&4805&2096& & &128.9.5.11.59.5839&5839&6336 \cr
\noalign{\hrule}
 & &7.37.109.113.557&218277&222310& & &3.25.7.11.169.31.59&241&5246 \cr
9110&1776887371&4.9.5.11.43.47.79.307&17945&21658&9128&1785058275&4.5.13.43.61.241&289&504 \cr
 & &16.3.25.49.11.13.17.37.97&31525&31416& & &64.9.7.289.241&723&9248 \cr
\noalign{\hrule}
 & &3.5.11.17.103.6151&32955&34706& & &9.5.7.11.19.47.577&221&202 \cr
9111&1777116165&4.9.25.7.2197.37.67&77&1598&9129&1785379365&4.5.7.11.13.17.101.577&2487&1552 \cr
 & &16.49.11.13.17.37.47&2303&3848& & &128.3.13.97.101.829&80413&84032 \cr
\noalign{\hrule}
 & &9.25.11.13.59.937&109631&101194& & &3.11.29.31.59.1021&49&38 \cr
9112&1778730525&4.19.37.2663.2963&2813&150&9130&1787110413&4.49.19.31.59.1021&11639&43290 \cr
 & &16.3.25.19.29.37.97&2813&5624& & &16.9.5.13.37.103.113&20905&32136 \cr
\noalign{\hrule}
 & &81.31.67.97.109&6325&4248& & &11.13.23.29.41.457&2985&2042 \cr
9113&1778769801&16.729.25.11.23.59&1649&16576&9131&1787153797&4.3.5.13.29.199.1021&3751&9522 \cr
 & &2048.7.17.37.97&4403&1024& & &16.27.5.121.529.31&3565&2376 \cr
\noalign{\hrule}
 & &125.29.53.59.157&1539&3014& & &3.25.29.79.101.103&11457&3478 \cr
9114&1779653875&4.81.5.11.19.53.137&109&2494&9132&1787495475&4.27.5.19.37.47.67&3443&5252 \cr
 & &16.9.11.29.43.109&473&7848& & &32.11.13.19.101.313&2717&5008 \cr
\noalign{\hrule}
 & &9.5.13.17.23.43.181&121&784& & &81.5.23.31.41.151&829&1034 \cr
9115&1780244505&32.3.49.121.23.43&325&578&9133&1787744115&4.11.31.47.151.829&2755&1926 \cr
 & &128.25.7.11.13.289&1309&320& & &16.9.5.11.19.29.47.107&14993&16264 \cr
\noalign{\hrule}
 & &3.25.11.19.1369.83&50429&27604& & &243.11.17.361.109&829&370 \cr
9116&1781103225&8.67.103.211.239&85&18&9134&1788057909&4.9.5.361.37.829&2039&1210 \cr
 & &32.9.5.17.211.239&3587&11472& & &16.25.121.37.2039&10175&16312 \cr
\noalign{\hrule}
 & &9.19.529.47.419&4521&3440& & &3.5.7.11.19.227.359&493&134 \cr
9117&1781409087&32.27.5.11.23.43.137&6361&470&9135&1788363885&4.5.7.17.29.67.227&2119&4464 \cr
 & &128.25.47.6361&6361&1600& & &128.9.13.17.31.163&8313&25792 \cr
\noalign{\hrule}
 & &27.5.7.13.17.19.449&29&484& & &27.7.11.13.289.229&205&16 \cr
9118&1781656695&8.121.17.29.449&2223&2716&9136&1788673887&32.5.11.17.41.229&21&208 \cr
 & &64.9.7.11.13.19.97&97&352& & &1024.3.5.7.13.41&41&2560 \cr
\noalign{\hrule}
 & &3.5.37.113.157.181&12219&5522& & &27.13.29.43.61.67&347&3740 \cr
9119&1782172155&4.9.5.11.251.4073&3611&7684&9137&1788867639&8.3.5.11.17.43.347&2231&1586 \cr
 & &32.11.17.23.113.157&187&368& & &32.13.17.23.61.97&391&1552 \cr
\noalign{\hrule}
 & &9.25.49.17.37.257&683&242& & &9.7.29.47.67.311&5401&8200 \cr
9120&1782224325&4.121.17.257.683&35&222&9138&1789252353&16.25.11.41.47.491&10411&12666 \cr
 & &16.3.5.7.11.37.683&683&88& & &64.3.5.29.359.2111&10555&11488 \cr
\noalign{\hrule}
 & &27.7.13.29.131.191&7421&17600& & &11.31.53.83.1193&32071&31158 \cr
9121&1782821313&128.25.7.11.41.181&3393&2942&9139&1789570387&4.27.13.31.577.2467&20045&2158 \cr
 & &512.9.5.13.29.1471&1471&1280& & &16.3.5.169.19.83.211&9633&8440 \cr
\noalign{\hrule}
 & &9.25.11.13.157.353&39353&57722& & &27.5.13.61.73.229&3683&3454 \cr
9122&1783170675&4.49.19.23.29.31.59&929&492&9140&1789638435&4.3.5.11.29.73.127.157&119&2236 \cr
 & &32.3.31.41.59.929&38089&29264& & &32.7.11.13.17.43.127&8041&14224 \cr
\noalign{\hrule}
 & &27.11.17.19.29.641&469&172& & &5.13.43.61.10499&375&418 \cr
9123&1783261359&8.7.17.19.29.43.67&641&90&9141&1790027005&4.3.625.11.19.10499&11187&688 \cr
 & &32.9.5.7.67.641&67&560& & &128.27.121.43.113&3267&7232 \cr
\noalign{\hrule}
 & &9.5.11.23.961.163&151&314& & &3.5.121.19.23.37.61&4247&8502 \cr
9124&1783380555&4.3.11.23.31.151.157&239&1900&9142&1790150835&4.9.11.13.31.109.137&8471&5092 \cr
 & &32.25.19.157.239&2983&19120& & &32.13.19.43.67.197&13199&8944 \cr
\noalign{\hrule}
 & &9.5.11.13.19.29.503&10999&1442& & &9.7.11.71.89.409&2231&1450 \cr
9125&1783479555&4.3.5.7.17.103.647&283&232&9143&1791038403&4.25.7.23.29.89.97&303&142 \cr
 & &64.7.29.283.647&4529&9056& & &16.3.5.29.71.97.101&2813&4040 \cr
\noalign{\hrule}
 & &9.25.13.691.883&85747&86438& & &81.13.17.19.23.229&2075&902 \cr
9126&1784697525&4.3.5.11.19.3929.4513&71173&3478&9144&1791406773&4.27.25.11.19.41.83&233&442 \cr
 & &16.11.37.47.103.691&3811&4136& & &16.13.17.41.83.233&3403&1864 \cr
\noalign{\hrule}
}%
}
$$
\eject
\vglue -23 pt
\noindent\hskip 1 in\hbox to 6.5 in{\ 9145 -- 9180 \hfill\fbd 1791570235 -- 1813065443\frb}
\vskip -9 pt
$$
\vbox{
\nointerlineskip
\halign{\strut
    \vrule \ \ \hfil \frb #\ 
   &\vrule \hfil \ \ \fbb #\frb\ 
   &\vrule \hfil \ \ \frb #\ \hfil
   &\vrule \hfil \ \ \frb #\ 
   &\vrule \hfil \ \ \frb #\ \ \vrule \hskip 2 pt
   &\vrule \ \ \hfil \frb #\ 
   &\vrule \hfil \ \ \fbb #\frb\ 
   &\vrule \hfil \ \ \frb #\ \hfil
   &\vrule \hfil \ \ \frb #\ 
   &\vrule \hfil \ \ \frb #\ \vrule \cr%
\noalign{\hrule}
 & &5.7.13.41.137.701&3839&5274& & &27.11.13.37.12611&493&506 \cr
9145&1791570235&4.9.11.137.293.349&28741&19072&9163&1801569627&4.121.17.23.29.12611&4551&8060 \cr
 & &1024.3.41.149.701&447&512& & &32.3.5.13.17.23.31.37.41&6355&6256 \cr
\noalign{\hrule}
 & &9.5.7.17.29.83.139&1331&80& & &27.7.11.13.163.409&1283&1580 \cr
9146&1791638415&32.25.7.1331.29&1411&1614&9164&1801809009&8.5.13.79.163.1283&1701&418 \cr
 & &128.3.11.17.83.269&269&704& & &32.243.5.7.11.19.79&855&1264 \cr
\noalign{\hrule}
 & &27.5.31.37.71.163&3289&8284& & &9.121.29.43.1327&44321&12740 \cr
9147&1792021185&8.11.13.19.23.31.109&7437&8150&9165&1802043441&8.5.49.13.23.41.47&29&6 \cr
 & &32.3.25.19.37.67.163&335&304& & &32.3.7.13.29.41.47&611&4592 \cr
\noalign{\hrule}
 & &9.5.17.29.211.383&28717&26818& & &27.5.7.11.29.53.113&2831&446 \cr
9148&1792836405&4.11.13.17.23.2209.53&111&2320&9166&1805413995&4.3.7.11.19.149.223&2525&2392 \cr
 & &128.3.5.23.29.37.53&1219&2368& & &64.25.13.23.101.223&30199&35680 \cr
\noalign{\hrule}
 & &27.11.17.137.2593&1313&1016& & &25.7.11.83.89.127&5761&17064 \cr
9149&1793611809&16.13.101.127.2593&1233&1360&9167&1805936825&16.27.49.79.823&5395&6218 \cr
 & &512.9.5.13.17.101.137&1313&1280& & &64.9.5.13.83.3109&3109&3744 \cr
\noalign{\hrule}
 & &243.11.31.59.367&5&6& & &27.5.49.121.37.61&1247&400 \cr
9150&1794232539&4.729.5.31.59.367&473&22126&9168&1806536655&32.125.7.29.37.43&7381&3756 \cr
 & &16.5.11.13.23.37.43&12857&1480& & &256.3.121.61.313&313&128 \cr
\noalign{\hrule}
 & &3.49.17.61.79.149&46537&19270& & &3.1849.199.1637&1117&520 \cr
9151&1794359469&4.5.41.47.173.269&63&110&9169&1807007361&16.5.13.1849.1117&8185&6336 \cr
 & &16.9.25.7.11.41.269&11275&6456& & &2048.9.25.11.1637&825&1024 \cr
\noalign{\hrule}
 & &27.11.13.19.43.569&91&478& & &9.125.7.11.31.673&1313&640 \cr
9152&1794874653&4.3.7.11.169.19.239&4385&14018&9170&1807257375&256.625.11.13.101&243&868 \cr
 & &16.5.43.163.877&877&6520& & &2048.243.7.13.31&351&1024 \cr
\noalign{\hrule}
 & &3.5.7.11.13.173.691&271&444& & &9.25.13.53.89.131&2783&14442 \cr
9153&1794938145&8.9.7.37.271.691&4247&21320&9171&1807436475&4.27.121.23.29.83&325&296 \cr
 & &128.5.13.31.41.137&5617&1984& & &64.25.121.13.37.83&3071&3872 \cr
\noalign{\hrule}
 & &3.125.11.29.43.349&1211&164& & &9.5.11.361.67.151&85&86 \cr
9154&1795212375&8.7.29.41.43.173&25477&25650&9172&1807857315&4.25.11.17.19.43.67.151&30031&1794 \cr
 & &32.27.25.7.19.73.349&1197&1168& & &16.3.13.23.43.59.509&58351&52936 \cr
\noalign{\hrule}
 & &25.17.31.59.2311&889&114& & &3.5.11.31.571.619&13081&44486 \cr
9155&1796398075&4.3.7.19.127.2311&9295&6882&9173&1807891635&4.13.29.59.103.127&2511&1172 \cr
 & &16.9.5.11.169.31.37&3663&1352& & &32.81.31.59.293&1593&4688 \cr
\noalign{\hrule}
 & &5.7.11.19.23.59.181&671&234& & &243.13.23.41.607&1357&1964 \cr
9156&1796688355&4.9.7.121.13.59.61&529&3070&9174&1808214759&8.3.13.529.59.491&1177&410 \cr
 & &16.3.5.13.529.307&3991&552& & &32.5.11.41.107.491&5401&8560 \cr
\noalign{\hrule}
 & &3.5.7.59.61.67.71&1927&2262& & &13.29.37.53.2447&1169&792 \cr
9157&1797646515&4.9.7.13.29.41.47.61&11447&5896&9175&1809059759&16.9.7.11.167.2447&2977&530 \cr
 & &64.11.29.67.11447&11447&10208& & &64.3.5.11.13.53.229&687&1760 \cr
\noalign{\hrule}
 & &27.5.7.11.13.53.251&31&4& & &3.121.13.59.67.97&161&40 \cr
9158&1797700905&8.11.13.31.53.251&1725&1036&9176&1809458079&16.5.7.13.23.59.97&3735&5092 \cr
 & &64.3.25.7.23.31.37&4255&992& & &128.9.25.19.67.83&2075&3648 \cr
\noalign{\hrule}
 & &25.121.31.127.151&1449&1576& & &9.5.7.19.71.4261&68497&80638 \cr
9159&1798323175&16.9.7.23.31.151.197&583&130&9177&1810648035&4.11.13.23.479.1753&2237&3990 \cr
 & &64.3.5.7.11.13.53.197&10441&8736& & &16.3.5.7.11.19.23.2237&2237&2024 \cr
\noalign{\hrule}
 & &3.5.7.11.59.61.433&1189&3354& & &3.11.19.43.47.1429&9085&6634 \cr
9160&1799913885&4.9.13.29.41.43.61&1069&700&9178&1810781643&4.5.23.31.47.79.107&3515&198 \cr
 & &32.25.7.13.43.1069&5345&8944& & &16.9.25.11.19.23.37&1725&296 \cr
\noalign{\hrule}
 & &9.5.11.17.31.67.103&565&38& & &3.7.11.19.841.491&42479&60140 \cr
9161&1800229365&4.25.11.19.103.113&3207&8432&9179&1812354159&8.5.31.97.107.397&145&252 \cr
 & &128.3.17.31.1069&1069&64& & &64.9.25.7.29.31.97&2425&2976 \cr
\noalign{\hrule}
 & &5.7.29.37.79.607&4807&558& & &7.19.71.101.1901&7667&5766 \cr
9162&1800874915&4.9.11.19.23.31.79&845&1604&9180&1813065443&4.3.11.17.961.41.71&245&1026 \cr
 & &32.3.5.169.19.401&3211&19248& & &16.81.5.49.17.19.31&2635&4536 \cr
\noalign{\hrule}
}%
}
$$
\eject
\vglue -23 pt
\noindent\hskip 1 in\hbox to 6.5 in{\ 9181 -- 9216 \hfill\fbd 1813075517 -- 1830995439\frb}
\vskip -9 pt
$$
\vbox{
\nointerlineskip
\halign{\strut
    \vrule \ \ \hfil \frb #\ 
   &\vrule \hfil \ \ \fbb #\frb\ 
   &\vrule \hfil \ \ \frb #\ \hfil
   &\vrule \hfil \ \ \frb #\ 
   &\vrule \hfil \ \ \frb #\ \ \vrule \hskip 2 pt
   &\vrule \ \ \hfil \frb #\ 
   &\vrule \hfil \ \ \fbb #\frb\ 
   &\vrule \hfil \ \ \frb #\ \hfil
   &\vrule \hfil \ \ \frb #\ 
   &\vrule \hfil \ \ \frb #\ \vrule \cr%
\noalign{\hrule}
 & &11.17.31.37.79.107&873&2050& & &17.47.79.28879&14839&14040 \cr
9181&1813075517&4.9.25.17.31.41.97&7169&7866&9199&1822871359&16.27.5.11.13.19.71.79&14479&15980 \cr
 & &16.81.5.19.23.67.107&9315&10184& & &128.9.25.17.47.14479&14479&14400 \cr
\noalign{\hrule}
 & &9.5.151.173.1543&10747&15376& & &27.5.17.67.71.167&907&902 \cr
9182&1813850505&32.3.5.11.961.977&1543&1388&9200&1823191605&4.11.17.41.71.167.907&347913&196444 \cr
 & &256.11.31.347.1543&3817&3968& & &32.9.29.31.43.67.733&21257&21328 \cr
\noalign{\hrule}
 & &125.7.13.19.37.227&3233&282& & &9.7.11.89.101.293&1007&1630 \cr
9183&1815233875&4.3.25.7.47.53.61&6583&4092&9201&1825207461&4.5.11.19.53.101.163&14357&19710 \cr
 & &32.9.11.29.31.227&2871&496& & &16.27.25.49.73.293&511&600 \cr
\noalign{\hrule}
 & &3.25.43.107.5261&269&376& & &9.5.11.41.53.1697&259&218 \cr
9184&1815439575&16.5.47.269.5261&6831&19474&9202&1825352595&4.5.7.11.37.109.1697&8439&46 \cr
 & &64.27.7.11.13.23.107&2093&3168& & &16.3.23.29.37.97&24679&776 \cr
\noalign{\hrule}
 & &27.5.11.13.157.599&961&766& & &5.11.13.37.151.457&2033&372 \cr
9185&1815500115&4.9.961.383.599&3241&8632&9203&1825580185&8.3.19.31.107.457&1887&1430 \cr
 & &64.7.13.31.83.463&14353&18592& & &32.9.5.11.13.17.19.37&153&304 \cr
\noalign{\hrule}
 & &9.13.37.41.53.193&51491&28910& & &9.11.31.43.109.127&95&1238 \cr
9186&1815534981&4.5.49.11.31.59.151&1443&386&9204&1826819181&4.5.11.19.109.619&2883&8878 \cr
 & &16.3.5.7.11.13.37.193&35&88& & &16.3.23.961.193&4439&248 \cr
\noalign{\hrule}
 & &243.5.361.41.101&4379&238& & &9.7.17.29.59.997&121&1118 \cr
9187&1816304715&4.5.7.17.19.29.151&1417&1452&9205&1826983557&4.3.121.13.17.29.43&295&266 \cr
 & &32.3.121.13.17.29.109&41093&32912& & &16.5.7.11.13.19.43.59&2717&1720 \cr
\noalign{\hrule}
 & &5.7.113.607.757&53393&32148& & &3.5.49.17.31.53.89&141&230 \cr
9188&1817318545&8.9.19.47.107.499&607&286&9206&1827106365&4.9.25.7.17.23.31.47&583&4408 \cr
 & &32.3.11.13.499.607&1497&2288& & &64.11.19.29.47.53&6061&1504 \cr
\noalign{\hrule}
 & &7.11.13.37.191.257&1575&1766& & &9.5.41.47.107.197&42889&3406 \cr
9189&1818035219&4.9.25.49.11.37.883&1657&1118&9207&1827865485&4.7.11.13.131.557&213&344 \cr
 & &16.3.13.43.883.1657&37969&39768& & &64.3.7.11.13.43.71&3913&24992 \cr
\noalign{\hrule}
 & &289.19.127.2609&1441&4050& & &3.19.29.41.53.509&6991&7770 \cr
9190&1819404413&4.81.25.11.127.131&149&1292&9208&1828312221&4.9.5.7.37.53.6991&38323&24596 \cr
 & &32.9.25.17.19.149&225&2384& & &32.5.11.13.19.43.2017&30745&32272 \cr
\noalign{\hrule}
 & &5.59.101.173.353&4851&5356& & &5.49.11.13.17.37.83&8667&6854 \cr
9191&1819550855&8.9.49.11.13.103.353&391&38&9209&1829072245&4.81.5.13.23.107.149&901&1036 \cr
 & &32.3.49.17.19.23.103&33269&54096& & &32.3.7.17.23.37.53.107&2461&2544 \cr
\noalign{\hrule}
 & &3.11.17.31.227.461&42661&212& & &9.25.49.17.43.227&899&1144 \cr
9192&1819915977&8.37.53.1153&603&550&9210&1829455425&16.5.11.13.17.29.31.43&199&1134 \cr
 & &32.9.25.11.37.67&7437&400& & &64.81.7.13.29.199&3393&6368 \cr
\noalign{\hrule}
 & &5.7.17.47.151.431&913&1242& & &81.5.7.61.71.149&209&706 \cr
9193&1819990165&4.27.11.17.23.83.151&1385&26&9211&1829479365&4.27.11.19.149.353&61&88 \cr
 & &16.3.5.11.13.23.277&19113&1144& & &64.121.19.61.353&6707&3872 \cr
\noalign{\hrule}
 & &81.5.11.23.109.163&3977&14768& & &9.11.29.47.71.191&443&5096 \cr
9194&1820495655&32.9.13.41.71.97&271&368&9212&1829880657&16.49.13.71.443&27&470 \cr
 & &1024.13.23.41.271&11111&6656& & &64.27.5.7.13.47&21&2080 \cr
\noalign{\hrule}
 & &9.11.41.293.1531&331&38& & &3.11.131.173.2447&1027&1420 \cr
9195&1820798397&4.11.19.331.1531&661&870&9213&1830059913&8.5.11.13.71.79.173&493&11790 \cr
 & &16.3.5.29.331.661&9599&26440& & &32.9.25.17.29.131&2175&272 \cr
\noalign{\hrule}
 & &25.11.29.43.47.113&6417&1558& & &11.13.37.53.61.107&27&80 \cr
9196&1821274675&4.9.19.23.31.41.47&339&1118&9214&1830320921&32.27.5.11.13.37.61&1007&214 \cr
 & &16.27.13.23.43.113&621&104& & &128.9.5.19.53.107&95&576 \cr
\noalign{\hrule}
 & &11.17.29.79.4253&69055&54282& & &7.11.17.47.29753&13067&16686 \cr
9197&1822057501&4.3.5.7.83.109.1973&605&1368&9215&1830493819&4.81.17.73.103.179&7535&25972 \cr
 & &64.27.25.121.19.83&17347&21600& & &32.3.5.11.43.137.151&17673&12080 \cr
\noalign{\hrule}
 & &9.11.13.59.103.233&7015&7714& & &3.7.13.19.31.59.193&1121&1700 \cr
9198&1822316067&4.3.5.7.19.23.29.59.61&2387&6458&9216&1830995439&8.25.17.361.3481&6253&2772 \cr
 & &16.49.11.19.31.3229&28861&25832& & &64.9.7.11.169.17.37&2431&3552 \cr
\noalign{\hrule}
}%
}
$$
\eject
\vglue -23 pt
\noindent\hskip 1 in\hbox to 6.5 in{\ 9217 -- 9252 \hfill\fbd 1832044695 -- 1853648489\frb}
\vskip -9 pt
$$
\vbox{
\nointerlineskip
\halign{\strut
    \vrule \ \ \hfil \frb #\ 
   &\vrule \hfil \ \ \fbb #\frb\ 
   &\vrule \hfil \ \ \frb #\ \hfil
   &\vrule \hfil \ \ \frb #\ 
   &\vrule \hfil \ \ \frb #\ \ \vrule \hskip 2 pt
   &\vrule \ \ \hfil \frb #\ 
   &\vrule \hfil \ \ \fbb #\frb\ 
   &\vrule \hfil \ \ \frb #\ \hfil
   &\vrule \hfil \ \ \frb #\ 
   &\vrule \hfil \ \ \frb #\ \vrule \cr%
\noalign{\hrule}
 & &3.5.13.289.19.29.59&501&3256& & &9.29.47.179.839&19985&19448 \cr
9217&1832044695&16.9.11.37.59.167&1843&340&9235&1842270327&16.3.5.7.11.13.17.29.571&829&4826 \cr
 & &128.5.11.17.19.97&1067&64& & &64.11.17.19.127.829&105283&113696 \cr
\noalign{\hrule}
 & &3.7.13.17.23.41.419&221&262& & &3.7.121.43.101.167&33&134 \cr
9218&1833737997&4.169.289.131.419&51865&3024&9236&1842939021&4.9.7.1331.43.67&689&2020 \cr
 & &128.27.5.7.11.23.41&99&320& & &32.5.13.53.67.101&4355&848 \cr
\noalign{\hrule}
 & &25.11.31.107.2011&1273&738& & &625.7.23.73.251&3393&2882 \cr
9219&1834383925&4.9.5.11.19.31.41.67&2011&244&9237&1843751875&4.9.25.11.13.23.29.131&73&502 \cr
 & &32.3.61.67.2011&201&976& & &16.3.29.73.131.251&393&232 \cr
\noalign{\hrule}
 & &9.25.7.13.19.53.89&18737&21758& & &3.61.109.193.479&5723&6050 \cr
9220&1835030925&4.3.5.11.23.41.43.457&1813&442&9238&1844040309&4.25.121.59.97.479&793&27468 \cr
 & &16.49.13.17.23.37.43&5117&6808& & &32.9.7.11.13.61.109&143&336 \cr
\noalign{\hrule}
 & &27.5.19.359.1993&2279&286& & &3.5.11.127.191.461&1281&820 \cr
9221&1835224155&4.11.13.43.53.359&3969&3610&9239&1845108705&8.9.25.7.41.61.127&191&1334 \cr
 & &16.81.5.49.361.43&817&1176& & &32.7.23.29.41.191&943&3248 \cr
\noalign{\hrule}
 & &81.25.11.73.1129&152281&158194& & &3.5.17.151.191.251&1397&1850 \cr
9222&1835838675&4.19.23.181.197.773&1485&2258&9240&1845968205&4.125.11.37.127.251&18037&13338 \cr
 & &16.27.5.11.23.181.1129&181&184& & &16.27.11.13.17.19.1061&13793&15048 \cr
\noalign{\hrule}
 & &3.25.7.11.961.331&5739&12466& & &9.25.11.13.61.941&823&152 \cr
9223&1836975525&4.9.5.23.271.1913&673&682&9241&1846877175&16.3.19.823.941&1705&764 \cr
 & &16.11.23.31.673.1913&15479&15304& & &128.5.11.19.31.191&5921&1216 \cr
\noalign{\hrule}
 & &11.41.53.151.509&37&546& & &3.25.11.41.193.283&637&1486 \cr
9224&1837160677&4.3.7.13.37.41.151&795&262&9242&1847487675&4.25.49.13.41.743&3113&2088 \cr
 & &16.9.5.37.53.131&5895&296& & &64.9.7.11.13.29.283&609&416 \cr
\noalign{\hrule}
 & &5.11.13.53.179.271&777&2746& & &7.121.19.29.37.107&15&92 \cr
9225&1838248555&4.3.5.7.37.53.1373&909&946&9243&1847653423&8.3.5.11.19.23.29.37&1917&1598 \cr
 & &16.27.11.43.101.1373&37071&34744& & &32.81.17.23.47.71&31671&53392 \cr
\noalign{\hrule}
 & &3.11.289.23.83.101&7515&868& & &11.13.19.29.31.757&477&74 \cr
9226&1838819433&8.27.5.7.11.31.167&391&446&9244&1849035331&4.9.11.37.53.757&87&670 \cr
 & &32.7.17.23.167.223&1561&2672& & &16.27.5.29.37.67&335&7992 \cr
\noalign{\hrule}
 & &81.19.37.43.751&1145&394& & &5.11.41.43.61.313&2033&468 \cr
9227&1838860299&4.5.37.43.197.229&13299&23146&9245&1851352745&8.9.11.13.19.43.107&763&656 \cr
 & &16.3.11.13.31.71.163&28613&14344& & &256.3.7.13.19.41.109&9919&7296 \cr
\noalign{\hrule}
 & &3.11.83.709.947&349&598& & &11.41.47.167.523&657&1180 \cr
9228&1839027597&4.11.13.23.349.709&4735&4482&9246&1851367177&8.9.5.41.47.59.73&1913&5344 \cr
 & &16.27.5.83.349.947&349&360& & &512.3.5.167.1913&1913&3840 \cr
\noalign{\hrule}
 & &27.7.43.97.2333&57475&55142& & &3.17.29.37.43.787&14313&14806 \cr
9229&1839148227&4.25.7.121.19.79.349&9&86&9247&1851881043&4.9.11.13.43.367.673&47261&31480 \cr
 & &16.9.5.11.43.79.349&3839&3160& & &64.5.11.167.283.787&9185&9056 \cr
\noalign{\hrule}
 & &25.7.13.17.19.2503&919&1584& & &3.5.7.11.13.23.31.173&61811&61538 \cr
9230&1839266975&32.9.5.11.13.17.919&2117&1198&9248&1852085235&4.5.11.29.113.547.1061&5661&356 \cr
 & &128.3.11.29.73.599&43727&61248& & &32.9.17.29.37.89.113&95497&92208 \cr
\noalign{\hrule}
 & &27.5.47.239.1213&3667&2398& & &3.5.7.11.19.29.41.71&241&374 \cr
9231&1839459915&4.11.19.109.193.239&279&2350&9249&1852574955&4.121.17.29.71.241&801&7790 \cr
 & &16.9.25.31.47.193&965&248& & &16.9.5.17.19.41.89&89&408 \cr
\noalign{\hrule}
 & &5.121.89.127.269&16497&16052& & &81.29.61.67.193&77&106 \cr
9232&1839506735&8.27.13.47.127.4013&7&134&9250&1852870059&4.27.7.11.53.67.193&5185&26 \cr
 & &32.9.7.13.67.4013&52169&67536& & &16.5.13.17.53.61&85&5512 \cr
\noalign{\hrule}
 & &9.5.23.43.173.239&235&754& & &11.23.29.41.61.101&47439&19928 \cr
9233&1840148235&4.3.25.13.29.47.239&229&946&9251&1853333537&16.27.7.47.53.251&205&46 \cr
 & &16.11.13.29.43.229&6641&1144& & &64.9.5.7.23.41.47&329&1440 \cr
\noalign{\hrule}
 & &3.25.11.13.199.863&122543&114782& & &343.121.59.757&109&648 \cr
9234&1841879325&4.29.31.59.67.1979&1961&18&9252&1853648489&16.81.7.11.59.109&223&190 \cr
 & &16.9.31.37.53.59&3441&25016& & &64.27.5.19.109.223&55917&35680 \cr
\noalign{\hrule}
}%
}
$$
\eject
\vglue -23 pt
\noindent\hskip 1 in\hbox to 6.5 in{\ 9253 -- 9288 \hfill\fbd 1853870337 -- 1874943291\frb}
\vskip -9 pt
$$
\vbox{
\nointerlineskip
\halign{\strut
    \vrule \ \ \hfil \frb #\ 
   &\vrule \hfil \ \ \fbb #\frb\ 
   &\vrule \hfil \ \ \frb #\ \hfil
   &\vrule \hfil \ \ \frb #\ 
   &\vrule \hfil \ \ \frb #\ \ \vrule \hskip 2 pt
   &\vrule \ \ \hfil \frb #\ 
   &\vrule \hfil \ \ \fbb #\frb\ 
   &\vrule \hfil \ \ \frb #\ \hfil
   &\vrule \hfil \ \ \frb #\ 
   &\vrule \hfil \ \ \frb #\ \vrule \cr%
\noalign{\hrule}
 & &9.11.19.61.107.151&575&584& & &49.17.71.139.227&5335&1476 \cr
9253&1853870337&16.25.11.23.73.107.151&247&8058&9271&1866139079&8.9.5.11.41.71.97&323&32 \cr
 & &64.3.5.13.17.19.23.79&9085&7072& & &512.3.11.17.19.41&627&10496 \cr
\noalign{\hrule}
 & &11.13.23.47.67.179&1525&444& & &9.5.7.17.29.61.197&18623&21368 \cr
9254&1853913919&8.3.25.13.37.61.67&1049&3306&9272&1866180015&16.11.17.1693.2671&29081&300 \cr
 & &32.9.5.19.29.1049&9441&44080& & &128.3.25.13.2237&11185&832 \cr
\noalign{\hrule}
 & &9.121.71.103.233&8255&8288& & &3.25.7.11.19.73.233&1079&846 \cr
9255&1855578681&64.3.5.7.11.13.37.103.127&14413&316&9273&1866312525&4.27.13.19.47.73.83&539&1432 \cr
 & &512.5.49.29.71.79&11455&12544& & &64.49.11.13.83.179&7553&5728 \cr
\noalign{\hrule}
 & &9.7.17.23.59.1277&2497&1220& & &121.19.37.47.467&585&118 \cr
9256&1855924119&8.5.11.17.23.61.227&2303&1998&9274&1867047787&4.9.5.121.13.47.59&37&84 \cr
 & &32.27.49.37.47.227&12173&10896& & &32.27.5.7.13.37.59&1593&7280 \cr
\noalign{\hrule}
 & &27.5.7.13.19.73.109&11869&12002& & &3.25.49.121.13.17.19&597&1702 \cr
9257&1857283155&4.9.5.11.169.17.83.353&1919&61576&9275&1867190325&4.9.5.49.23.37.199&289&44 \cr
 & &64.11.19.43.101.179&18079&15136& & &32.11.289.23.199&391&3184 \cr
\noalign{\hrule}
 & &9.25.43.181.1061&2743&2562& & &121.41.53.7103&2465&4638 \cr
9258&1857996675&4.27.5.7.13.43.61.211&10679&1606&9276&1867613099&4.3.5.121.17.29.773&689&84 \cr
 & &16.11.59.61.73.181&4307&5368& & &32.9.7.13.17.29.53&1071&6032 \cr
\noalign{\hrule}
 & &5.49.11.13.29.31.59&81&458& & &361.53.89.1097&16091&16038 \cr
9259&1858291435&4.81.5.31.59.229&203&262&9277&1868012189&4.729.11.1097.16091&22855&6764 \cr
 & &16.27.7.29.131.229&3537&1832& & &32.27.5.7.11.19.89.653&10395&10448 \cr
\noalign{\hrule}
 & &3.5.49.11.47.67.73&907&2524& & &81.25.11.23.41.89&6767&442 \cr
9260&1858555545&8.5.67.631.907&621&286&9278&1869473925&4.13.17.41.67.101&1323&1424 \cr
 & &32.27.11.13.23.631&8203&3312& & &128.27.49.13.17.89&637&1088 \cr
\noalign{\hrule}
 & &3.25.13.31.89.691&10373&6902& & &7.11.13.53.101.349&4421&4770 \cr
9261&1858807275&4.7.11.17.23.29.31.41&2755&3042&9279&1870065197&4.9.5.11.2809.4421&4241&4186 \cr
 & &16.9.5.169.19.23.841&10933&10488& & &16.3.7.13.23.4241.4421&101683&101784 \cr
\noalign{\hrule}
 & &9.7.11.37.47.1543&1643&1976& & &5.121.13.17.71.197&159&38 \cr
9262&1859510961&16.13.19.31.53.1543&1593&50&9280&1870131835&4.3.5.13.17.19.53.71&517&588 \cr
 & &64.27.25.13.19.59&14573&2400& & &32.9.49.11.19.47.53&20727&16112 \cr
\noalign{\hrule}
 & &9.25.11.43.101.173&2509&2036& & &27.5.17.29.157.179&157&302 \cr
9263&1859564025&8.5.13.173.193.509&119&54&9281&1870395165&4.151.24649.179&25839&1190 \cr
 & &32.27.7.17.193.509&25959&21616& & &16.81.5.7.11.17.29&11&168 \cr
\noalign{\hrule}
 & &5.7.13.19.29.41.181&7095&8362& & &3.7.11.19.31.59.233&4669&17000 \cr
9264&1860481805&4.3.25.11.19.37.43.113&181&294&9282&1870403073&16.125.49.17.23.29&279&946 \cr
 & &16.9.49.11.37.43.181&2849&3096& & &64.9.5.11.17.31.43&255&1376 \cr
\noalign{\hrule}
 & &81.127.331.547&21329&48140& & &13.41.53.151.439&3735&4268 \cr
9265&1862533359&8.5.7.11.29.83.277&331&54&9283&1872597961&8.9.5.11.83.97.439&125&1192 \cr
 & &32.27.29.83.331&83&464& & &128.3.625.83.149&51875&28608 \cr
\noalign{\hrule}
 & &9.25.7.361.29.113&869&1886& & &9.11.47.53.71.107&1085&1406 \cr
9266&1863220275&4.5.7.11.19.23.41.79&7053&452&9284&1873488573&4.3.5.7.11.19.31.37.71&107&1598 \cr
 & &32.3.11.113.2351&2351&176& & &16.17.19.37.47.107&629&152 \cr
\noalign{\hrule}
 & &3.11.29.59.61.541&1147&476& & &49.11.19.23.73.109&1843&1734 \cr
9267&1863335463&8.7.17.29.31.37.59&541&1170&9285&1874215651&4.3.11.289.361.23.97&64435&39894 \cr
 & &32.9.5.7.13.31.541&1209&560& & &16.9.5.49.61.109.263&2367&2440 \cr
\noalign{\hrule}
 & &81.5.11.29.47.307&6355&2548& & &11.19.83.103.1049&1435&522 \cr
9268&1864154655&8.25.49.11.13.31.41&1363&912&9286&1874291309&4.9.5.7.29.41.1049&979&2168 \cr
 & &256.3.7.19.29.31.47&589&896& & &64.3.5.7.11.89.271&5691&14240 \cr
\noalign{\hrule}
 & &125.11.169.71.113&3199&4824& & &9.7.11.19.23.41.151&7871&16440 \cr
9269&1864344625&16.9.7.11.13.67.457&75&68&9287&1874888631&16.27.5.17.137.463&299&164 \cr
 & &128.27.25.17.67.457&30619&29376& & &128.13.17.23.41.137&1781&1088 \cr
\noalign{\hrule}
 & &3.5.7.19.29.167.193&519&316& & &3.169.37.127.787&881&770 \cr
9270&1864724505&8.9.19.79.173.193&11&182&9288&1874943291&4.5.7.11.13.787.881&8199&2032 \cr
 & &32.7.11.13.79.173&1903&16432& & &128.9.5.11.127.911&2733&3520 \cr
\noalign{\hrule}
}%
}
$$
\eject
\vglue -23 pt
\noindent\hskip 1 in\hbox to 6.5 in{\ 9289 -- 9324 \hfill\fbd 1875443115 -- 1893899085\frb}
\vskip -9 pt
$$
\vbox{
\nointerlineskip
\halign{\strut
    \vrule \ \ \hfil \frb #\ 
   &\vrule \hfil \ \ \fbb #\frb\ 
   &\vrule \hfil \ \ \frb #\ \hfil
   &\vrule \hfil \ \ \frb #\ 
   &\vrule \hfil \ \ \frb #\ \ \vrule \hskip 2 pt
   &\vrule \ \ \hfil \frb #\ 
   &\vrule \hfil \ \ \fbb #\frb\ 
   &\vrule \hfil \ \ \frb #\ \hfil
   &\vrule \hfil \ \ \frb #\ 
   &\vrule \hfil \ \ \frb #\ \vrule \cr%
\noalign{\hrule}
 & &3.5.7.13.23.31.41.47&177&1258& & &3.5.7.11.13.17.47.157&8989&9694 \cr
9289&1875443115&4.9.13.17.31.37.59&4715&4312&9307&1883526645&4.11.13.37.89.101.131&1377&64 \cr
 & &64.5.49.11.23.37.41&259&352& & &512.81.17.37.89&3293&6912 \cr
\noalign{\hrule}
 & &7.13.79.311.839&575&1602& & &9.3125.193.347&31923&35048 \cr
9290&1875818581&4.9.25.23.89.839&553&286&9308&1883559375&16.81.13.337.3547&11875&15422 \cr
 & &16.3.25.7.11.13.23.79&575&264& & &64.625.11.13.19.701&7711&7904 \cr
\noalign{\hrule}
 & &27.67.599.1733&1771&38& & &7.17.89.293.607&269&1782 \cr
9291&1877863203&4.7.11.19.23.599&1765&4824&9309&1883619941&4.81.11.269.607&293&6970 \cr
 & &64.9.5.67.353&1765&32& & &16.3.5.17.41.293&615&8 \cr
\noalign{\hrule}
 & &81.5.13.31.37.311&70433&55522& & &9.13.289.103.541&1265&724 \cr
9292&1878115005&4.11.17.19.23.71.337&8249&15678&9310&1884161799&8.5.11.17.23.103.181&2105&972 \cr
 & &16.9.11.13.67.73.113&7571&6424& & &64.243.25.23.421&11367&18400 \cr
\noalign{\hrule}
 & &3.5.11.101.251.449&3097&1842& & &9.11.169.61.1847&64699&47968 \cr
9293&1878128835&4.9.19.101.163.307&18017&34480&9311&1885031577&64.23.29.97.1499&1865&366 \cr
 & &128.5.43.419.431&18533&26816& & &256.3.5.29.61.373&1865&3712 \cr
\noalign{\hrule}
 & &121.29.43.59.211&1649&1860& & &9.5.11.43.263.337&11591&44954 \cr
9294&1878392263&8.3.5.17.31.43.59.97&471&6194&9312&1886510835&4.7.169.19.67.173&5805&5786 \cr
 & &32.9.17.19.157.163&24021&49552& & &16.27.5.7.11.169.43.263&169&168 \cr
\noalign{\hrule}
 & &3.125.101.49603&74657&74152& & &5.7.17.43.89.829&2145&3658 \cr
9295&1878713625&16.25.121.13.23.31.617&6401&12726&9313&1887686885&4.3.25.11.13.31.43.59&829&504 \cr
 & &64.9.7.11.13.37.101.173&24739&24864& & &64.27.7.11.59.829&649&864 \cr
\noalign{\hrule}
 & &81.5.13.361.23.43&511&74& & &3.11.17.19.29.31.197&115&82 \cr
9296&1879757685&4.9.7.19.37.43.73&575&242&9314&1887740877&4.5.17.19.23.29.31.41&611&288 \cr
 & &16.25.7.121.23.73&605&4088& & &256.9.5.13.23.41.47&28905&38272 \cr
\noalign{\hrule}
 & &27.49.37.193.199&1175&176& & &9.11.97.421.467&1943&1846 \cr
9297&1880061057&32.25.7.11.47.199&579&814&9315&1888017021&4.11.13.29.67.71.467&4947&190 \cr
 & &128.3.5.121.37.193&121&320& & &16.3.5.13.17.19.29.97&1615&3016 \cr
\noalign{\hrule}
 & &27.5.53.317.829&1351&1034& & &9.25.7.29.67.617&1031&16456 \cr
9298&1880283915&4.3.7.11.47.193.829&361&1190&9316&1888158825&16.7.121.17.1031&3515&3702 \cr
 & &16.5.49.17.361.193&17689&26248& & &64.3.5.11.19.37.617&407&608 \cr
\noalign{\hrule}
 & &9.25.7.283.4219&299&16& & &27.5.7.17.41.47.61&107&46 \cr
9299&1880513775&32.5.13.23.4219&2167&2052&9317&1888392555&4.3.5.7.23.41.47.107&12337&25498 \cr
 & &256.27.11.13.19.197&11229&18304& & &16.11.169.19.61.73&3211&6424 \cr
\noalign{\hrule}
 & &11.19.31.37.47.167&7791&58& & &3.19.31.61.89.197&3009&2420 \cr
9300&1881585827&4.3.49.29.31.53&3173&3120&9318&1889829471&8.9.5.121.17.59.197&3019&2822 \cr
 & &128.9.5.7.13.19.167&315&832& & &32.5.11.289.83.3019&250577&254320 \cr
\noalign{\hrule}
 & &729.25.2197.47&42517&60742& & &5.11.79.257.1693&819&874 \cr
9301&1881895275&4.121.17.41.61.251&13&684&9319&1890513845&4.9.7.13.19.23.79.257&18623&890 \cr
 & &32.9.11.13.19.251&209&4016& & &16.3.5.7.11.89.1693&89&168 \cr
\noalign{\hrule}
 & &9.5.11.289.59.223&481&634& & &27.7.11.17.59.907&9215&8308 \cr
9302&1882174635&4.11.13.17.37.59.317&6913&1524&9320&1891309959&8.5.7.17.19.31.67.97&2169&92 \cr
 & &32.3.31.37.127.223&1147&2032& & &64.9.5.23.97.241&5543&15520 \cr
\noalign{\hrule}
 & &9.7.11.13.1849.113&241&16400& & &3.25.19.47.61.463&2497&182 \cr
9303&1882313433&32.25.7.41.241&741&946&9321&1891574925&4.5.7.11.13.61.227&389&8334 \cr
 & &128.3.5.11.13.19.43&19&320& & &16.9.389.463&1167&8 \cr
\noalign{\hrule}
 & &9.25.47.67.2657&1919&11366& & &7.11.13.47.167.241&7857&7340 \cr
9304&1882550925&4.3.5.19.101.5683&8777&8272&9322&1893500609&8.81.5.97.241.367&167&2002 \cr
 & &128.11.19.47.67.131&1441&1216& & &32.9.7.11.13.97.167&97&144 \cr
\noalign{\hrule}
 & &9.25.121.13.17.313&1577&3142& & &3.7.11.1031.7951&10817&44840 \cr
9305&1883234925&4.3.5.17.19.83.1571&3443&4412&9323&1893618111&16.5.19.29.59.373&7821&734 \cr
 & &32.11.83.313.1103&1103&1328& & &64.9.11.79.367&367&7584 \cr
\noalign{\hrule}
 & &9.5.49.11.19.61.67&253&82& & &3.5.13.29.53.71.89&295&82 \cr
9306&1883473515&4.49.121.23.41.61&2263&3666&9324&1893899085&4.25.41.53.59.89&2871&1846 \cr
 & &16.3.13.31.41.47.73&29419&15416& & &16.9.11.13.29.59.71&177&88 \cr
\noalign{\hrule}
}%
}
$$
\eject
\vglue -23 pt
\noindent\hskip 1 in\hbox to 6.5 in{\ 9325 -- 9360 \hfill\fbd 1894483899 -- 1915478565\frb}
\vskip -9 pt
$$
\vbox{
\nointerlineskip
\halign{\strut
    \vrule \ \ \hfil \frb #\ 
   &\vrule \hfil \ \ \fbb #\frb\ 
   &\vrule \hfil \ \ \frb #\ \hfil
   &\vrule \hfil \ \ \frb #\ 
   &\vrule \hfil \ \ \frb #\ \ \vrule \hskip 2 pt
   &\vrule \ \ \hfil \frb #\ 
   &\vrule \hfil \ \ \fbb #\frb\ 
   &\vrule \hfil \ \ \frb #\ \hfil
   &\vrule \hfil \ \ \frb #\ 
   &\vrule \hfil \ \ \frb #\ \vrule \cr%
\noalign{\hrule}
 & &9.7.11.29.107.881&13983&20150& & &7.11.13.19.239.419&171&248 \cr
9325&1894483899&4.27.25.13.31.59.79&29&2&9343&1904581679&16.9.13.361.31.239&2095&1012 \cr
 & &16.25.13.29.59.79&1475&8216& & &128.3.5.11.23.31.419&713&960 \cr
\noalign{\hrule}
 & &625.43.251.281&9179&2904& & &9.25.49.13.97.137&3267&1486 \cr
9326&1895520625&16.3.25.121.67.137&843&832&9344&1904645925&4.243.25.121.743&1049&7124 \cr
 & &2048.9.11.13.137.281&13563&13312& & &32.11.13.137.1049&1049&176 \cr
\noalign{\hrule}
 & &3.7.11.19.29.53.281&281&302& & &729.11.59.4027&1649&2378 \cr
9327&1895595933&4.19.29.151.78961&81081&2120&9345&1905258267&4.11.17.29.41.59.97&6695&972 \cr
 & &64.81.5.7.11.13.53&135&416& & &32.243.5.13.29.103&2987&1040 \cr
\noalign{\hrule}
 & &169.59.397.479&2915&3312& & &81.7.19.149.1187&481&562 \cr
9328&1896115273&32.9.5.11.13.23.53.59&4787&8622&9346&1905345099&4.13.19.37.281.1187&945&242 \cr
 & &128.81.479.4787&4787&5184& & &16.27.5.7.121.13.281&3653&4840 \cr
\noalign{\hrule}
 & &3.5.361.223.1571&3493&1922& & &9.5.11.13.43.71.97&2573&1792 \cr
9329&1897053195&4.7.961.223.499&5203&1710&9347&1905667335&512.7.13.31.43.83&561&518 \cr
 & &16.9.5.121.19.31.43&1333&2904& & &2048.3.49.11.17.31.37&30821&31744 \cr
\noalign{\hrule}
 & &7.11.17.23.29.41.53&67&186& & &27.11.19.29.61.191&1937&3706 \cr
9330&1897252819&4.3.29.31.41.53.67&605&2142&9348&1906651197&4.13.17.109.149.191&25&2508 \cr
 & &16.27.5.7.121.17.31&297&1240& & &32.3.25.11.19.109&25&1744 \cr
\noalign{\hrule}
 & &9.11.101.109.1741&4071&15080& & &27.7.11.17.79.683&123157&130670 \cr
9331&1897500231&16.27.5.13.23.29.59&545&1166&9349&1907002251&4.5.73.107.179.1151&6783&1028 \cr
 & &64.25.11.13.53.109&1325&416& & &32.3.7.17.19.179.257&3401&4112 \cr
\noalign{\hrule}
 & &3.5.11.61.409.461&253&662& & &5.11.17.41.71.701&117&818 \cr
9332&1897745685&4.121.23.331.461&25327&14724&9350&1907971285&4.9.13.41.71.409&1203&4114 \cr
 & &32.9.19.31.43.409&817&1488& & &16.27.121.17.401&297&3208 \cr
\noalign{\hrule}
 & &27.7.11.13.23.43.71&7021&17900& & &9.5.29.59.137.181&14749&9500 \cr
9333&1897808913&8.25.49.17.59.179&2967&76&9351&1909245015&8.3.625.343.19.43&1991&116 \cr
 & &64.3.25.19.23.43&19&800& & &64.7.11.19.29.181&77&608 \cr
\noalign{\hrule}
 & &81.25.7.11.19.641&11293&10652& & &81.5.37.233.547&83&1082 \cr
9334&1899010575&8.27.5.23.491.2663&221&2884&9352&1909853235&4.3.83.541.547&67081&67628 \cr
 & &64.7.13.17.103.491&22763&15712& & &32.49.11.29.1369.53&21571&22736 \cr
\noalign{\hrule}
 & &27.11.23.101.2753&245&866& & &7.11.47.97.5441&6031&44118 \cr
9335&1899380043&4.5.49.433.2753&6201&8954&9353&1910024963&4.27.19.37.43.163&5441&3850 \cr
 & &16.9.7.121.13.37.53&5291&2968& & &16.9.25.7.11.5441&25&72 \cr
\noalign{\hrule}
 & &27.5.49.41.43.163&13&176& & &625.7.79.5527&3751&1776 \cr
9336&1900945935&32.5.7.11.13.41.43&527&978&9354&1910269375&32.3.25.7.121.31.37&4951&474 \cr
 & &128.3.13.17.31.163&221&1984& & &128.9.79.4951&4951&576 \cr
\noalign{\hrule}
 & &81.5.37.181.701&373&328& & &11.13.17.19.59.701&545&222 \cr
9337&1901311785&16.9.37.41.181.373&1573&56&9355&1910330851&4.3.5.11.37.109.701&249&950 \cr
 & &256.7.121.13.373&45133&11648& & &16.9.125.19.37.83&10375&2664 \cr
\noalign{\hrule}
 & &3.25.7.11.13.19.31.43&1049&26& & &27.5.7.121.17.983&475&596 \cr
9338&1901424525&4.7.169.19.1049&151&18&9356&1910819295&8.3.125.19.149.983&287&2662 \cr
 & &16.9.151.1049&151&25176& & &32.7.1331.41.149&1639&656 \cr
\noalign{\hrule}
 & &81.5.7.17.19.31.67&109&46& & &9.625.11.17.23.79&2311&494 \cr
9339&1901919285&4.9.17.19.23.67.109&9889&8750&9357&1911256875&4.3.125.13.19.2311&1343&968 \cr
 & &16.625.7.11.23.29.31&2875&2552& & &64.121.13.17.19.79&209&416 \cr
\noalign{\hrule}
 & &9.49.19.23.71.139&35&104& & &27.7.11.13.193.367&26275&31046 \cr
9340&1901924073&16.3.5.343.13.19.71&2881&1166&9358&1914349437&4.25.7.361.43.1051&193&858 \cr
 & &64.11.13.43.53.67&46163&15136& & &16.3.5.11.13.19.43.193&215&152 \cr
\noalign{\hrule}
 & &9.1331.31.47.109&2355&1024& & &9.5.7.11.17.19.29.59&73&130 \cr
9341&1902420927&2048.27.5.47.157&1027&242&9359&1914942645&4.3.25.11.13.17.59.73&8533&11542 \cr
 & &8192.121.13.79&1027&4096& & &16.7.13.23.29.53.199&4577&5512 \cr
\noalign{\hrule}
 & &3.13.23.71.167.179&185&114& & &3.5.7.11.13.29.53.83&3383&7782 \cr
9342&1903795491&4.9.5.19.37.167.179&781&2392&9360&1915478565&4.9.13.17.199.1297&1643&346 \cr
 & &64.5.11.13.23.37.71&185&352& & &16.31.53.173.199&6169&1384 \cr
\noalign{\hrule}
}%
}
$$
\eject
\vglue -23 pt
\noindent\hskip 1 in\hbox to 6.5 in{\ 9361 -- 9396 \hfill\fbd 1916305105 -- 1938804945\frb}
\vskip -9 pt
$$
\vbox{
\nointerlineskip
\halign{\strut
    \vrule \ \ \hfil \frb #\ 
   &\vrule \hfil \ \ \fbb #\frb\ 
   &\vrule \hfil \ \ \frb #\ \hfil
   &\vrule \hfil \ \ \frb #\ 
   &\vrule \hfil \ \ \frb #\ \ \vrule \hskip 2 pt
   &\vrule \ \ \hfil \frb #\ 
   &\vrule \hfil \ \ \fbb #\frb\ 
   &\vrule \hfil \ \ \frb #\ \hfil
   &\vrule \hfil \ \ \frb #\ 
   &\vrule \hfil \ \ \frb #\ \vrule \cr%
\noalign{\hrule}
 & &5.11.13.43.157.397&31093&31236& & &5.13.19.47.89.373&461&696 \cr
9361&1916305105&8.3.5.17.19.31.43.59.137&143&5748&9379&1926919865&16.3.19.29.373.461&89&462 \cr
 & &64.9.11.13.17.31.479&8143&8928& & &64.9.7.11.89.461&5071&2016 \cr
\noalign{\hrule}
 & &3.25.97.211.1249&135407&128132& & &9.5.7.13.17.19.31.47&87&134 \cr
9362&1917246225&8.43.47.67.103.311&855&1166&9380&1927152045&4.27.5.7.19.29.31.67&799&286 \cr
 & &32.9.5.11.19.53.67.103&60049&61104& & &16.11.13.17.29.47.67&319&536 \cr
\noalign{\hrule}
 & &3.49.17.29.103.257&339&4708& & &9.19.29.41.53.179&5975&1364 \cr
9363&1918379841&8.9.11.29.107.113&1043&2060&9381&1928887353&8.3.25.11.19.31.239&377&212 \cr
 & &64.5.7.11.103.149&1639&160& & &64.5.13.29.53.239&1195&416 \cr
\noalign{\hrule}
 & &9.5.13.37.263.337&2313&14782& & &9.25.11.19.89.461&923&122 \cr
9364&1918417995&4.81.19.257.389&169&88&9382&1929388725&4.5.13.61.71.461&53&408 \cr
 & &64.11.169.19.389&2717&12448& & &64.3.13.17.53.61&901&25376 \cr
\noalign{\hrule}
 & &27.25.17.29.73.79&989&836& & &9.343.11.113.503&377&4150 \cr
9365&1919113425&8.3.11.19.23.29.43.79&1241&11432&9383&1930081923&4.25.13.29.83.113&15&98 \cr
 & &128.11.17.73.1429&1429&704& & &16.3.125.49.13.29&3625&104 \cr
\noalign{\hrule}
 & &9.5.11.19.43.47.101&18745&19654& & &3.5.121.17.19.37.89&293&312 \cr
9366&1919758005&4.25.11.23.31.163.317&21573&30098&9384&1930504785&16.9.13.17.37.89.293&715&86 \cr
 & &16.27.17.23.47.101.149&1173&1192& & &64.5.11.169.43.293&7267&9376 \cr
\noalign{\hrule}
 & &5.11.13.841.31.103&477&3670& & &25.13.19.37.79.107&1683&292 \cr
9367&1919998795&4.9.25.29.53.367&3819&5356&9385&1931299175&8.9.11.17.19.37.73&17525&18328 \cr
 & &32.27.13.19.67.103&1273&432& & &128.3.25.29.79.701&2103&1856 \cr
\noalign{\hrule}
 & &9.25.11.23.89.379&1971&76& & &9.13.289.19.31.97&66527&49940 \cr
9368&1920137175&8.243.5.11.19.73&503&712&9386&1931838129&8.5.11.71.227.937&291&646 \cr
 & &128.73.89.503&503&4672& & &32.3.11.17.19.97.227&227&176 \cr
\noalign{\hrule}
 & &13.17.31.53.67.79&1665&1886& & &9.11.3559.5483&14141&46172 \cr
9369&1921904179&4.9.5.23.31.37.41.79&19&1166&9387&1931885703&8.7.17.79.97.179&5353&2310 \cr
 & &16.3.11.19.23.41.53&779&6072& & &32.3.5.49.11.53.101&2597&8080 \cr
\noalign{\hrule}
 & &9.5.11.13.59.61.83&4531&5180& & &27.7.11.13.43.1663&9139&5828 \cr
9370&1922243895&8.25.7.23.37.61.197&73&1452&9388&1932673743&8.3.169.19.31.37.47&115&3326 \cr
 & &64.3.121.23.37.73&2701&8096& & &32.5.23.47.1663&1081&80 \cr
\noalign{\hrule}
 & &25.7.13.43.103.191&2673&4012& & &5.7.11.31.67.2417&1175&1242 \cr
9371&1924511225&8.243.5.11.17.43.59&721&764&9389&1932741965&4.27.125.7.11.23.31.47&6799&8174 \cr
 & &64.9.7.17.59.103.191&531&544& & &16.9.13.47.61.67.523&37271&37656 \cr
\noalign{\hrule}
 & &9.11.43.59.79.97&361&350& & &9.25.113.139.547&46543&15268 \cr
9372&1924662069&4.25.7.361.43.59.97&12793&108&9390&1933139025&8.7.11.61.109.347&31&30 \cr
 & &32.27.5.11.19.1163&1163&4560& & &32.3.5.7.11.31.109.347&37169&38864 \cr
\noalign{\hrule}
 & &3.5.7.11.13.17.19.397&227&172& & &13.41.83.89.491&20383&20370 \cr
9373&1925388465&8.13.17.43.227.397&1899&4850&9391&1933200061&4.3.5.7.11.17.41.89.97.109&149&67758 \cr
 & &32.9.25.43.97.211&12513&16880& & &16.9.5.11.23.149.491&3427&3960 \cr
\noalign{\hrule}
 & &3.1331.127.3797&1835&1962& & &9.5.19.5041.449&5863&2668 \cr
9374&1925500467&4.27.5.1331.109.367&2033&37970&9392&1935214695&8.11.13.23.29.41.71&133&1056 \cr
 & &16.25.19.107.3797&475&856& & &512.3.7.121.19.23&847&5888 \cr
\noalign{\hrule}
 & &49.11.13.17.19.23.37&171&220& & &3.5.169.37.47.439&3479&2774 \cr
9375&1926035111&8.9.5.121.13.361.37&4369&9062&9393&1935272235&4.49.19.71.73.439&22311&13970 \cr
 & &32.3.5.17.23.197.257&3855&3152& & &16.9.5.7.11.37.67.127&4191&3752 \cr
\noalign{\hrule}
 & &9.7.17.37.61.797&5921&44290& & &9.5.11.13.23.103.127&16351&14446 \cr
9376&1926545859&4.5.31.43.103.191&493&462&9394&1936053405&4.3.11.31.83.197.233&65405&19504 \cr
 & &16.3.7.11.17.29.43.103&2987&3784& & &128.5.23.53.103.127&53&64 \cr
\noalign{\hrule}
 & &27.25.7.37.103.107&33&3778& & &3.5.121.43.59.421&7429&10674 \cr
9377&1926746325&4.81.5.11.1889&1147&742&9395&1938559755&4.27.11.17.19.23.593&91&118 \cr
 & &16.7.11.31.37.53&53&2728& & &16.7.13.17.23.59.593&13639&12376 \cr
\noalign{\hrule}
 & &11.17.29.43.8263&112463&120726& & &3.5.7.11.31.173.313&44239&4276 \cr
9378&1926840707&4.9.13.19.41.211.353&493&140&9396&1938804945&8.13.41.83.1069&5247&8650 \cr
 & &32.3.5.7.13.17.19.29.41&3895&4368& & &32.9.25.11.53.173&265&48 \cr
\noalign{\hrule}
}%
}
$$
\eject
\vglue -23 pt
\noindent\hskip 1 in\hbox to 6.5 in{\ 9397 -- 9432 \hfill\fbd 1939000121 -- 1957496749\frb}
\vskip -9 pt
$$
\vbox{
\nointerlineskip
\halign{\strut
    \vrule \ \ \hfil \frb #\ 
   &\vrule \hfil \ \ \fbb #\frb\ 
   &\vrule \hfil \ \ \frb #\ \hfil
   &\vrule \hfil \ \ \frb #\ 
   &\vrule \hfil \ \ \frb #\ \ \vrule \hskip 2 pt
   &\vrule \ \ \hfil \frb #\ 
   &\vrule \hfil \ \ \fbb #\frb\ 
   &\vrule \hfil \ \ \frb #\ \hfil
   &\vrule \hfil \ \ \frb #\ 
   &\vrule \hfil \ \ \frb #\ \vrule \cr%
\noalign{\hrule}
 & &31.101.421.1471&20525&21996& & &27.343.43.59.83&275&2812 \cr
9397&1939000121&8.9.25.13.31.47.821&8979&1694&9415&1950098031&8.3.25.11.19.37.83&49&34 \cr
 & &32.27.5.7.121.41.73&80811&67760& & &32.5.49.11.17.19.37&1615&6512 \cr
\noalign{\hrule}
 & &3.25.11.13.193.937&111&826& & &5.11.13.23.31.43.89&981&1384 \cr
9398&1939519725&4.9.5.7.37.59.193&217&748&9416&1950985465&16.9.23.89.109.173&10819&1118 \cr
 & &32.49.11.17.31.37&1147&13328& & &64.3.13.31.43.349&349&96 \cr
\noalign{\hrule}
 & &5.11.169.23.43.211&137&852& & &3.25.29.59.67.227&1085&858 \cr
9399&1939671305&8.3.13.71.137.211&6235&3492&9417&1951694925&4.9.125.7.11.13.31.59&1171&2546 \cr
 & &64.27.5.29.43.97&783&3104& & &16.13.19.31.67.1171&7657&9368 \cr
\noalign{\hrule}
 & &27.5.7.11.17.79.139&1409&2344& & &7.367.457.1663&45&412 \cr
9400&1940507415&16.7.79.293.1409&981&428&9418&1952416879&8.9.5.7.103.1663&6239&4576 \cr
 & &128.9.107.109.293&11663&18752& & &512.3.11.13.17.367&2431&768 \cr
\noalign{\hrule}
 & &27.5.41.547.641&64721&47414& & &9.5.11.13.37.59.139&6517&4988 \cr
9401&1940720445&4.61.151.157.1061&75317&84894&9419&1952617095&8.3.343.19.29.37.43&425&278 \cr
 & &16.3.11.41.167.14149&14149&14696& & &32.25.7.17.29.43.139&3655&3248 \cr
\noalign{\hrule}
 & &25.11.13.31.83.211&153&2168& & &81.5.11.263.1667&793&874 \cr
9402&1940878225&16.9.5.17.83.271&199&216&9420&1953165555&4.5.11.13.19.23.61.263&20039&3996 \cr
 & &256.243.199.271&53929&31104& & &32.27.13.29.37.691&13949&11056 \cr
\noalign{\hrule}
 & &3.25.7.11.29.67.173&4351&666& & &3.5.7.11.17.29.47.73&5453&3336 \cr
9403&1941202725&4.27.5.7.19.37.229&319&346&9421&1953662865&16.9.5.49.19.41.139&16907&6862 \cr
 & &16.11.29.37.173.229&229&296& & &64.11.29.47.53.73&53&32 \cr
\noalign{\hrule}
 & &3.5.49.11.17.71.199&1241&4226& & &3.5.7.11.23.29.43.59&409&1656 \cr
9404&1941960405&4.7.289.73.2113&45&2068&9422&1954466745&16.27.11.529.409&413&116 \cr
 & &32.9.5.11.47.73&3431&48& & &128.7.29.59.409&409&64 \cr
\noalign{\hrule}
 & &9.227.773.1231&1775&544& & &7.11.13.19.23.41.109&2241&3440 \cr
9405&1944043209&64.3.25.17.71.227&1397&2462&9423&1954905953&32.27.5.7.41.43.83&3959&556 \cr
 & &256.5.11.127.1231&635&1408& & &256.9.37.107.139&46287&13696 \cr
\noalign{\hrule}
 & &3.23.37.43.89.199&5555&3002& & &3.125.11.13.361.101&2813&5338 \cr
9406&1944295869&4.5.11.19.79.89.101&1449&470&9424&1955221125&4.5.17.19.29.97.157&529&1314 \cr
 & &16.9.25.7.23.47.79&5925&2632& & &16.9.17.529.29.73&15341&29784 \cr
\noalign{\hrule}
 & &49.11.1019.3541&44441&5490& & &9.625.17.113.181&13039&7414 \cr
9407&1944862381&4.9.5.19.61.2339&1261&1078&9425&1955818125&4.11.13.289.59.337&6335&10716 \cr
 & &16.3.5.49.11.13.19.97&1261&2280& & &32.3.5.7.11.19.47.181&893&1232 \cr
\noalign{\hrule}
 & &27.19.59.131.491&965&3454& & &5.7.11.13.17.83.277&97&90 \cr
9408&1946803707&4.3.5.11.59.157.193&491&1456&9426&1956189235&4.9.25.13.83.97.277&4267&27258 \cr
 & &128.7.13.157.491&2041&448& & &16.27.7.11.17.59.251&1593&2008 \cr
\noalign{\hrule}
 & &27.5.7.13.23.61.113&2679&1276& & &3.11.17.29.127.947&225&1172 \cr
9409&1947651615&8.81.11.13.19.29.47&395&17324&9427&1956656361&8.27.25.17.29.293&233&260 \cr
 & &64.5.61.71.79&71&2528& & &64.125.13.233.293&68269&52000 \cr
\noalign{\hrule}
 & &3.11.13.83.229.239&337&350& & &11.13.17.23.31.1129&369&760 \cr
9410&1948806717&4.25.7.11.83.239.337&6183&208&9428&1956899087&16.9.5.11.13.19.31.41&629&970 \cr
 & &128.27.13.229.337&337&576& & &64.3.25.17.19.37.97&10767&15200 \cr
\noalign{\hrule}
 & &27.5.17.283.3001&22649&22366& & &27.13.17.29.43.263&4883&4620 \cr
9411&1949104485&4.9.11.17.29.53.71.211&13555&278&9429&1956943287&8.81.5.7.11.19.29.257&275&2074 \cr
 & &16.5.139.211.2711&29329&21688& & &32.125.121.17.19.61&15125&18544 \cr
\noalign{\hrule}
 & &3.11.5711.10343&8027&2316& & &3.7.29.61.139.379&451&830 \cr
9412&1949272809&8.9.11.23.193.349&21895&18056&9430&1957046469&4.5.11.29.41.83.139&3717&314 \cr
 & &128.5.29.37.61.151&65453&48320& & &16.9.5.7.11.59.157&9735&1256 \cr
\noalign{\hrule}
 & &3.29.73.191.1607&847&760& & &5.7.11.17.37.59.137&56407&71688 \cr
9413&1949356887&16.5.7.121.19.73.191&6867&7076&9431&1957419695&16.3.13.29.103.4339&37&66 \cr
 & &128.9.5.49.11.29.61.109&58751&58560& & &64.9.11.13.37.4339&4339&3744 \cr
\noalign{\hrule}
 & &3.13.23.67.71.457&2827&1930& & &169.47.59.4177&127193&119250 \cr
9414&1950032253&4.5.11.193.257.457&871&414&9432&1957496749&4.9.125.11.31.53.373&103&1222 \cr
 & &16.9.11.13.23.67.193&193&264& & &16.3.5.11.13.31.47.103&1705&2472 \cr
\noalign{\hrule}
}%
}
$$
\eject
\vglue -23 pt
\noindent\hskip 1 in\hbox to 6.5 in{\ 9433 -- 9468 \hfill\fbd 1958111595 -- 1981776795\frb}
\vskip -9 pt
$$
\vbox{
\nointerlineskip
\halign{\strut
    \vrule \ \ \hfil \frb #\ 
   &\vrule \hfil \ \ \fbb #\frb\ 
   &\vrule \hfil \ \ \frb #\ \hfil
   &\vrule \hfil \ \ \frb #\ 
   &\vrule \hfil \ \ \frb #\ \ \vrule \hskip 2 pt
   &\vrule \ \ \hfil \frb #\ 
   &\vrule \hfil \ \ \fbb #\frb\ 
   &\vrule \hfil \ \ \frb #\ \hfil
   &\vrule \hfil \ \ \frb #\ 
   &\vrule \hfil \ \ \frb #\ \vrule \cr%
\noalign{\hrule}
 & &9.5.11.17.19.37.331&1019&26& & &243.5.7.11.107.197&583&1118 \cr
9433&1958111595&4.3.13.17.37.1019&8717&8606&9451&1972045845&4.121.13.43.53.197&4487&1926 \cr
 & &16.169.23.331.379&3887&3032& & &16.9.7.43.107.641&641&344 \cr
\noalign{\hrule}
 & &125.7.13.167.1031&56839&72036& & &81.7.37.149.631&8389&8648 \cr
9434&1958513375&8.27.23.29.113.503&1169&1430&9452&1972424601&16.3.23.47.149.8389&9335&946 \cr
 & &32.3.5.7.11.13.167.503&503&528& & &64.5.11.43.47.1867&87749&75680 \cr
\noalign{\hrule}
 & &3.121.29.37.47.107&855&884& & &9.5.121.841.431&371&470 \cr
9435&1958790471&8.27.5.121.13.17.19.107&295&2594&9453&1973654595&4.25.7.11.47.53.431&353&78 \cr
 & &32.25.13.17.59.1297&76523&88400& & &16.3.7.13.47.53.353&32383&19768 \cr
\noalign{\hrule}
 & &9.25.13.17.157.251&25773&29698& & &11.73.97.101.251&52725&28378 \cr
9436&1959513075&4.27.121.31.71.479&67643&33634&9454&1974614741&4.3.25.7.19.37.2027&29&66 \cr
 & &16.17.23.67.173.251&1541&1384& & &16.9.5.7.11.29.2027&9135&16216 \cr
\noalign{\hrule}
 & &3.5.11.29.151.2713&3737&44432& & &3.5.11.31.149.2591&923&1668 \cr
9437&1960237455&32.37.101.2777&140257&140220&9455&1974691785&8.9.11.13.31.71.139&149&490 \cr
 & &256.9.5.13.19.41.10789&204991&204672& & &32.5.49.13.139.149&1807&784 \cr
\noalign{\hrule}
 & &25.11.17.79.5309&703&6012& & &5.49.11.23.151.211&1657&3978 \cr
9438&1960746425&8.9.5.11.19.37.167&3007&3172&9456&1974904085&4.9.13.17.151.1657&3773&2116 \cr
 & &64.3.13.19.31.61.97&73749&58976& & &32.3.343.11.17.529&483&272 \cr
\noalign{\hrule}
 & &81.7.17.271.751&275&292& & &121.19.67.101.127&261&1658 \cr
9439&1961738919&8.25.11.73.271.751&387&3368&9457&1975781291&4.9.11.29.67.829&691&1520 \cr
 & &128.9.5.43.73.421&18103&23360& & &128.3.5.19.29.691&10365&1856 \cr
\noalign{\hrule}
 & &81.5.49.121.19.43&3379&6862& & &5.7.73.181.4273&1881&2392 \cr
9440&1961817165&4.5.11.31.47.73.109&441&76&9458&1976070215&16.9.5.11.13.19.23.181&4273&2282 \cr
 & &32.9.49.19.31.109&109&496& & &64.3.7.13.163.4273&489&416 \cr
\noalign{\hrule}
 & &27.7.11.17.43.1291&101&200& & &81.343.11.29.223&115&6352 \cr
9441&1961995959&16.3.25.17.101.1291&2795&1078&9459&1976399271&32.5.49.23.397&223&174 \cr
 & &64.125.49.11.13.43&875&416& & &128.3.5.23.29.223&23&320 \cr
\noalign{\hrule}
 & &81.5.7.11.13.47.103&289&278& & &27.25.13.17.29.457&649&1106 \cr
9442&1962565605&4.5.13.289.47.103.139&2403&652&9460&1977016275&4.5.7.11.17.29.59.79&457&138 \cr
 & &32.27.17.89.139.163&22657&24208& & &16.3.23.59.79.457&1357&632 \cr
\noalign{\hrule}
 & &5.29.53.59.61.71&2829&770& & &5.11.23.29.31.37.47&4119&7684 \cr
9443&1963740365&4.3.25.7.11.23.41.53&261&314&9461&1977651665&8.3.17.47.113.1373&9015&14326 \cr
 & &16.27.7.11.29.41.157&12089&8856& & &32.9.5.13.19.29.601&5409&3952 \cr
\noalign{\hrule}
 & &81.7.11.29.73.149&355&1994& & &5.11.13.19.41.53.67&1091&3264 \cr
9444&1967355621&4.5.7.71.73.997&243&754&9462&1977853735&128.3.11.17.19.1091&441&650 \cr
 & &16.243.5.13.29.71&923&120& & &512.27.25.49.13.17&4165&6912 \cr
\noalign{\hrule}
 & &3.5.11.31.409.941&3987&8692& & &9.25.11.19.23.31.59&1547&282 \cr
9445&1968604935&8.27.11.41.53.443&941&490&9463&1978200675&4.27.5.7.13.17.19.47&1181&1276 \cr
 & &32.5.49.443.941&443&784& & &32.11.17.29.47.1181&23171&18896 \cr
\noalign{\hrule}
 & &3.25.11.19.23.43.127&347&728& & &3.5.7.11.23.163.457&19&96 \cr
9446&1968827025&16.7.11.13.19.23.347&45&254&9464&1978853415&64.9.19.163.457&5075&3608 \cr
 & &64.9.5.7.127.347&1041&224& & &1024.25.7.11.29.41&1189&2560 \cr
\noalign{\hrule}
 & &3.5.17.23.37.43.211&1115&1078& & &9.25.7.121.13.17.47&10097&6478 \cr
9447&1968886365&4.25.49.11.23.211.223&73&5202&9465&1979502525&4.3.11.23.41.79.439&1781&32900 \cr
 & &16.9.49.11.289.73&2409&6664& & &32.25.7.13.47.137&137&16 \cr
\noalign{\hrule}
 & &3.7.17.41.59.2281&59&60& & &2401.11.13.73.79&7155&3284 \cr
9448&1969832823&8.9.5.41.3481.2281&122683&20038&9466&1980059081&8.27.5.49.53.821&631&190 \cr
 & &32.11.19.43.233.587&136771&143792& & &32.3.25.19.53.631&47325&16112 \cr
\noalign{\hrule}
 & &3.7.13.139.167.311&145&166& & &49.19.23.37.41.61&22005&20878 \cr
9449&1970853339&4.5.13.29.83.139.167&341&11196&9467&1981494781&4.27.5.11.13.41.73.163&25123&1814 \cr
 & &32.9.11.29.31.311&319&1488& & &16.3.5.7.37.97.907&2721&3880 \cr
\noalign{\hrule}
 & &27.25.11.19.89.157&1261&964& & &3.5.37.53.89.757&4251&466 \cr
9450&1971240975&8.13.19.97.157.241&75&3058&9468&1981776795&4.9.13.37.109.233&2075&1958 \cr
 & &32.3.25.11.97.139&139&1552& & &16.25.11.83.89.233&2563&3320 \cr
\noalign{\hrule}
}%
}
$$
\eject
\vglue -23 pt
\noindent\hskip 1 in\hbox to 6.5 in{\ 9469 -- 9504 \hfill\fbd 1982392425 -- 2001962885\frb}
\vskip -9 pt
$$
\vbox{
\nointerlineskip
\halign{\strut
    \vrule \ \ \hfil \frb #\ 
   &\vrule \hfil \ \ \fbb #\frb\ 
   &\vrule \hfil \ \ \frb #\ \hfil
   &\vrule \hfil \ \ \frb #\ 
   &\vrule \hfil \ \ \frb #\ \ \vrule \hskip 2 pt
   &\vrule \ \ \hfil \frb #\ 
   &\vrule \hfil \ \ \fbb #\frb\ 
   &\vrule \hfil \ \ \frb #\ \hfil
   &\vrule \hfil \ \ \frb #\ 
   &\vrule \hfil \ \ \frb #\ \vrule \cr%
\noalign{\hrule}
 & &9.25.13.23.79.373&13489&10132& & &27.7.11.17.23.31.79&45835&10492 \cr
9469&1982392425&8.25.7.17.41.47.149&1463&2262&9487&1990765161&8.5.43.61.89.103&23&66 \cr
 & &32.3.49.11.13.19.29.41&13079&14896& & &32.3.5.11.23.61.103&305&1648 \cr
\noalign{\hrule}
 & &9.5.49.361.47.53&187&548& & &9.7.11.19.37.61.67&1247&26 \cr
9470&1982848455&8.3.11.17.47.53.137&475&1274&9488&1991100573&4.3.7.13.29.43.61&335&274 \cr
 & &32.25.49.13.19.137&685&208& & &16.5.13.43.67.137&559&5480 \cr
\noalign{\hrule}
 & &243.5.7.11.17.29.43&6251&7466& & &23.37.127.18427&44825&63252 \cr
9471&1983272445&4.49.17.19.47.3733&865&66&9489&1991534879&8.9.25.7.11.163.251&4849&1426 \cr
 & &16.3.5.11.173.3733&3733&1384& & &32.3.11.13.23.31.373&11563&6864 \cr
\noalign{\hrule}
 & &9.5.11.13.43.67.107&1787&1228& & &9.5.43.47.61.359&833&962 \cr
9472&1983698145&8.11.107.307.1787&1635&1742&9490&1991604555&4.3.49.13.17.37.47.61&25&1012 \cr
 & &32.3.5.13.67.109.1787&1787&1744& & &32.25.7.11.13.23.37&3367&20240 \cr
\noalign{\hrule}
 & &3.25.7.11.289.29.41&507&218& & &9.11.13.17.29.43.73&83&304 \cr
9473&1984411275&4.9.7.11.169.41.109&685&316&9491&1991667249&32.11.19.29.73.83&12545&13932 \cr
 & &32.5.13.79.109.137&10823&22672& & &256.81.5.13.43.193&965&1152 \cr
\noalign{\hrule}
 & &9.5.121.37.59.167&329&34& & &27.5.7.17.19.61.107&1147&1742 \cr
9474&1985034645&4.3.7.17.37.47.167&1475&1364&9492&1992268845&4.13.19.31.37.61.67&1819&660 \cr
 & &32.25.7.11.31.47.59&1085&752& & &32.3.5.11.13.17.31.107&341&208 \cr
\noalign{\hrule}
 & &5.169.31.41.1849&3501&1738& & &25.7.19.23.71.367&4895&3546 \cr
9475&1985816755&4.9.5.11.43.79.389&14911&2074&9493&1992709075&4.9.125.7.11.89.197&1521&646 \cr
 & &16.3.13.17.31.37.61&2257&408& & &16.81.169.17.19.89&15041&11016 \cr
\noalign{\hrule}
 & &49.47.61.67.211&1305&2782& & &27.5.7.11.37.71.73&180133&174512 \cr
9476&1986008171&4.9.5.7.13.29.47.107&305&682&9494&1993459545&32.13.61.839.2953&3977&6930 \cr
 & &16.3.25.11.31.61.107&3317&6600& & &128.9.5.7.11.41.61.97&3977&3904 \cr
\noalign{\hrule}
 & &3.25.31.53.71.227&481&1162& & &5.11.29.41.43.709&909&854 \cr
9477&1986017325&4.25.7.13.37.71.83&127&198&9495&1993697365&4.9.7.29.61.101.709&1039&19522 \cr
 & &16.9.7.11.37.83.127&21497&33528& & &16.3.7.43.227.1039&4767&8312 \cr
\noalign{\hrule}
 & &5.23.31.47.71.167&179&534& & &27.7.13.31.73.359&6697&7100 \cr
9478&1986699635&4.3.47.89.167.179&3225&11638&9496&1996108569&8.25.37.71.181.359&2211&416 \cr
 & &16.9.25.11.529.43&4257&920& & &512.3.5.11.13.67.181&12127&14080 \cr
\noalign{\hrule}
 & &27.5.13.43.113.233&7511&7744& & &729.11.89.2797&197&2600 \cr
9479&1986917985&128.7.121.13.29.37.43&9&1582&9497&1996193727&16.27.25.11.13.197&2611&2314 \cr
 & &512.9.49.29.113&1421&256& & &64.7.169.89.373&2611&5408 \cr
\noalign{\hrule}
 & &3.5.11.13.29.89.359&6205&1538& & &3.7.19.23.367.593&1581&988 \cr
9480&1987511955&4.25.11.17.73.769&247&522&9498&1997199687&8.9.13.17.361.23.31&6523&1780 \cr
 & &16.9.13.17.19.29.73&1387&408& & &64.5.11.13.89.593&979&2080 \cr
\noalign{\hrule}
 & &5.11.19.71.73.367&559&486& & &3.25.11.419.5779&253&6032 \cr
9481&1987758245&4.243.13.43.71.367&125&3178&9499&1997655825&32.5.121.13.23.29&18017&18162 \cr
 & &16.27.125.7.13.227&20657&5400& & &128.9.43.419.1009&3027&2752 \cr
\noalign{\hrule}
 & &7.169.107.113.139&16895&1188& & &5.97.139.149.199&143&342 \cr
9482&1988207767&8.27.5.7.11.31.109&79&86&9500&1998922165&4.9.11.13.19.139.149&5&144 \cr
 & &32.9.31.43.79.109&30573&54064& & &128.81.5.11.13.19&20007&704 \cr
\noalign{\hrule}
 & &5.19.29.601.1201&827277&828478& & &27.5.7.11.13.361.41&757&244 \cr
9483&1988561755&4.3.7.121.17.43.53.3481&601&300&9501&2000133135&8.5.19.41.61.757&21867&25652 \cr
 & &32.9.25.121.3481.601&17405&17424& & &64.3.121.37.53.197&10441&13024 \cr
\noalign{\hrule}
 & &5.13.59.509.1019&763&1782& & &25.49.19.23.37.101&57419&36006 \cr
9484&1989103285&4.81.7.11.13.59.109&1019&928&9502&2000509525&4.3.17.67.353.857&143&210 \cr
 & &256.27.29.109.1019&3161&3456& & &16.9.5.7.11.13.17.857&11141&13464 \cr
\noalign{\hrule}
 & &9.11.19.43.73.337&7807&16794& & &11.23.67.263.449&1413&1480 \cr
9485&1989802683&4.243.37.211.311&19883&31390&9503&2001692737&16.9.5.23.37.157.449&1501&1052 \cr
 & &16.5.43.59.73.337&59&40& & &128.3.5.19.79.157.263&14915&15168 \cr
\noalign{\hrule}
 & &3.5.29.67.163.419&77&412& & &5.13.37.47.89.199&405&206 \cr
9486&1990516065&8.7.11.29.103.419&369&50&9504&2001962885&4.81.25.37.89.103&507&418 \cr
 & &32.9.25.7.41.103&615&11536& & &16.243.11.169.19.103&25029&21736 \cr
\noalign{\hrule}
}%
}
$$
\eject
\vglue -23 pt
\noindent\hskip 1 in\hbox to 6.5 in{\ 9505 -- 9540 \hfill\fbd 2002886145 -- 2025077945\frb}
\vskip -9 pt
$$
\vbox{
\nointerlineskip
\halign{\strut
    \vrule \ \ \hfil \frb #\ 
   &\vrule \hfil \ \ \fbb #\frb\ 
   &\vrule \hfil \ \ \frb #\ \hfil
   &\vrule \hfil \ \ \frb #\ 
   &\vrule \hfil \ \ \frb #\ \ \vrule \hskip 2 pt
   &\vrule \ \ \hfil \frb #\ 
   &\vrule \hfil \ \ \fbb #\frb\ 
   &\vrule \hfil \ \ \frb #\ \hfil
   &\vrule \hfil \ \ \frb #\ 
   &\vrule \hfil \ \ \frb #\ \vrule \cr%
\noalign{\hrule}
 & &9.5.13.1097.3121&2109&1012& & &3.121.59.167.563&221&342 \cr
9505&2002886145&8.27.5.11.13.19.23.37&6257&5542&9523&2013647757&4.27.13.17.19.59.167&2815&1694 \cr
 & &32.17.37.163.6257&106369&96496& & &16.5.7.121.13.17.563&221&280 \cr
\noalign{\hrule}
 & &169.289.89.461&44935&3906& & &3.11.13.17.281.983&2037&1054 \cr
9506&2003897389&4.9.5.7.11.19.31.43&221&178&9524&2014494339&4.9.7.13.289.31.97&2711&110 \cr
 & &16.3.5.11.13.17.31.89&341&120& & &16.5.11.97.2711&13555&776 \cr
\noalign{\hrule}
 & &7.11.17.73.139.151&1293&236& & &9.25.7.19.23.29.101&6523&4402 \cr
9507&2005645873&8.3.17.59.73.431&5817&1510&9525&2015957475&4.3.11.29.31.71.593&1919&140 \cr
 & &32.9.5.7.151.277&1385&144& & &32.5.7.11.19.31.101&31&176 \cr
\noalign{\hrule}
 & &3.5.41.59.167.331&1273&3692& & &9.5.7.71.109.827&25019&29154 \cr
9508&2005725945&8.13.19.67.71.167&3965&792&9526&2016048195&4.27.43.113.127.197&71&5390 \cr
 & &128.9.5.11.169.61&1859&11712& & &16.5.49.11.71.113&77&904 \cr
\noalign{\hrule}
 & &9.17.29.59.79.97&45227&6280& & &125.7.13.361.491&28341&33034 \cr
9509&2006043129&16.5.49.13.71.157&3363&6842&9527&2016230125&4.9.7.47.67.83.199&2185&1716 \cr
 & &64.3.11.19.59.311&3421&608& & &32.27.5.11.13.19.23.199&4577&4752 \cr
\noalign{\hrule}
 & &27.5.11.13.37.2809&2725&2566& & &11.13.23.31.73.271&9963&9820 \cr
9510&2006426565&4.9.125.53.109.1283&27313&95312&9528&2017054897&8.243.5.23.31.41.491&73&11366 \cr
 & &128.7.11.13.23.37.191&1337&1472& & &32.27.5.73.5683&5683&2160 \cr
\noalign{\hrule}
 & &5.7.11.13.53.67.113&18981&18914& & &9.49.11.19.43.509&1825&2756 \cr
9511&2008321315&4.27.343.19.37.113.193&559&1588&9529&2017302903&8.25.11.13.43.53.73&509&294 \cr
 & &32.9.13.37.43.193.397&132201&132784& & &32.3.5.49.13.53.509&265&208 \cr
\noalign{\hrule}
 & &27.13.23.241.1033&54175&30416& & &27.19.23.41.43.97&3955&22 \cr
9512&2009797569&32.25.11.197.1901&133&2034&9530&2017758789&4.3.5.7.11.43.113&1843&1886 \cr
 & &128.9.25.7.19.113&15029&1600& & &16.5.7.19.23.41.97&35&8 \cr
\noalign{\hrule}
 & &9.13.23.29.43.599&66061&79838& & &9.5.529.29.37.79&1859&786 \cr
9513&2010050523&4.11.19.31.191.2131&1161&970&9531&2017878435&4.27.11.169.79.131&2203&1334 \cr
 & &16.27.5.11.19.31.43.97&8835&8536& & &16.169.23.29.2203&2203&1352 \cr
\noalign{\hrule}
 & &9.25.7.43.67.443&5189&4114& & &3.125.49.17.23.281&1271&696 \cr
9514&2010145725&4.3.121.17.67.5189&4303&886&9532&2018879625&16.9.5.7.17.29.31.41&281&484 \cr
 & &16.121.13.331.443&1573&2648& & &128.121.31.41.281&3751&2624 \cr
\noalign{\hrule}
 & &9.5.17.361.29.251&497&758& & &49.19.29.37.43.47&1775&528 \cr
9515&2010205035&4.7.17.361.71.379&91113&97250&9533&2018904223&32.3.25.11.19.37.71&343&1692 \cr
 & &16.3.125.121.251.389&3025&3112& & &256.27.5.343.47&945&128 \cr
\noalign{\hrule}
 & &9.11.19.29.191.193&923&730& & &27.7.11.23.29.31.47&32825&29206 \cr
9516&2010839787&4.3.5.11.13.71.73.191&22439&8878&9534&2020411701&4.9.25.13.17.101.859&47&812 \cr
 & &16.5.19.23.193.1181&1181&920& & &32.5.7.13.29.47.101&505&208 \cr
\noalign{\hrule}
 & &5.11.13.19.23.41.157&519&266& & &5.49.19.37.59.199&1139&1044 \cr
9517&2011272835&4.3.7.13.361.41.173&471&15272&9535&2022211135&8.9.49.17.29.67.199&715&118 \cr
 & &64.9.23.83.157&747&32& & &32.3.5.11.13.29.59.67&4147&3216 \cr
\noalign{\hrule}
 & &9.11.17.41.103.283&12545&16604& & &19.31.79.157.277&8833&3570 \cr
9518&2011368447&8.5.7.13.17.193.593&849&256&9536&2023586659&4.3.5.7.121.17.31.73&597&1108 \cr
 & &4096.3.7.193.283&1351&2048& & &32.9.11.17.199.277&1683&3184 \cr
\noalign{\hrule}
 & &25.7.11.37.47.601&34119&56356& & &81.7.23.311.499&209&290 \cr
9519&2011892575&8.9.17.73.193.223&1973&1750&9537&2023819749&4.5.7.11.19.23.29.311&1663&108 \cr
 & &32.3.125.7.193.1973&9865&9264& & &32.27.19.29.1663&1663&8816 \cr
\noalign{\hrule}
 & &9.5.7.13.31.83.191&667&418& & &49.29.37.137.281&209&72 \cr
9520&2012459085&4.3.11.13.19.23.29.191&1787&2360&9538&2024056769&16.9.49.11.19.29.37&1765&344 \cr
 & &64.5.19.23.59.1787&41101&35872& & &256.3.5.11.43.353&45537&7040 \cr
\noalign{\hrule}
 & &3.5.11.47.139.1867&51319&51366& & &27.5.13.59.113.173&209&322 \cr
9521&2012523315&4.9.7.19.37.73.139.1223&10741&34510&9539&2024201205&4.3.5.7.11.13.19.23.173&37177&34928 \cr
 & &16.5.49.17.23.29.73.467&389011&389528& & &128.49.37.47.59.113&2303&2368 \cr
\noalign{\hrule}
 & &5.11.13.571.4931&2537&2394& & &5.11.31.41.59.491&1039&1380 \cr
9522&2013154715&4.9.5.7.19.43.59.571&359&8206&9540&2025077945&8.3.25.23.491.1039&533&42 \cr
 & &16.3.11.43.359.373&16039&8616& & &32.9.7.13.41.1039&7273&1872 \cr
\noalign{\hrule}
}%
}
$$
\eject
\vglue -23 pt
\noindent\hskip 1 in\hbox to 6.5 in{\ 9541 -- 9576 \hfill\fbd 2025858549 -- 2047626105\frb}
\vskip -9 pt
$$
\vbox{
\nointerlineskip
\halign{\strut
    \vrule \ \ \hfil \frb #\ 
   &\vrule \hfil \ \ \fbb #\frb\ 
   &\vrule \hfil \ \ \frb #\ \hfil
   &\vrule \hfil \ \ \frb #\ 
   &\vrule \hfil \ \ \frb #\ \ \vrule \hskip 2 pt
   &\vrule \ \ \hfil \frb #\ 
   &\vrule \hfil \ \ \fbb #\frb\ 
   &\vrule \hfil \ \ \frb #\ \hfil
   &\vrule \hfil \ \ \frb #\ 
   &\vrule \hfil \ \ \frb #\ \vrule \cr%
\noalign{\hrule}
 & &3.11.13.71.227.293&1395&1102& & &3.5.7.11.19.23.37.109&4527&1468 \cr
9541&2025858549&4.27.5.13.19.29.31.71&293&10588&9559&2035596255&8.27.37.367.503&4351&14260 \cr
 & &32.19.293.2647&2647&304& & &64.5.19.23.31.229&229&992 \cr
\noalign{\hrule}
 & &3.25.7.13.37.71.113&36703&32522& & &81.25.7.121.1187&4601&1334 \cr
9542&2026008075&4.49.289.23.101.127&6435&212&9560&2035912725&4.3.5.7.23.29.43.107&37&572 \cr
 & &32.9.5.11.13.53.101&3333&848& & &32.11.13.23.37.43&559&13616 \cr
\noalign{\hrule}
 & &81.5.7.11.181.359&14603&20938& & &5.49.11.73.79.131&903&538 \cr
9543&2026370115&4.9.17.361.29.859&1625&1624&9561&2036010515&4.3.343.43.79.269&7765&19332 \cr
 & &64.125.7.13.17.841.859&357425&357344& & &32.81.5.179.1553&14499&24848 \cr
\noalign{\hrule}
 & &125.7.17.29.37.127&75573&60698& & &9.67.131.149.173&4165&4612 \cr
9544&2027031125&4.243.11.31.89.311&29&27650&9562&2036202561&8.3.5.49.17.173.1153&2915&22516 \cr
 & &16.3.25.7.29.79&3&632& & &64.25.11.13.53.433&61919&42400 \cr
\noalign{\hrule}
 & &9.7.11.19.337.457&35445&34988& & &27.7.121.169.17.31&169&358 \cr
9545&2027836503&8.27.5.7.17.139.8747&17509&8762&9563&2036781747&4.121.28561.179&3451&25110 \cr
 & &32.5.13.17.337.17509&17509&17680& & &16.81.5.7.17.29.31&145&24 \cr
\noalign{\hrule}
 & &3.5.7.121.23.53.131&2329&2984& & &3.25.41.59.103.109&17427&12958 \cr
9546&2028848745&16.11.17.53.137.373&2925&22694&9564&2036858475&4.9.5.11.19.31.37.157&767&938 \cr
 & &64.9.25.7.13.1621&4863&2080& & &16.7.13.37.59.67.157&17353&16328 \cr
\noalign{\hrule}
 & &25.7.11.13.17.19.251&177&44& & &11.17.19.29.53.373&11997&4910 \cr
9547&2028851825&8.3.25.121.59.251&1387&1638&9565&2036938453&4.9.5.17.31.43.491&371&1102 \cr
 & &32.27.7.13.19.59.73&1593&1168& & &16.3.5.7.19.29.31.53&93&280 \cr
\noalign{\hrule}
 & &7.11.2963.8893&109629&118522& & &25.7.11.17.19.29.113&907&2370 \cr
9548&2028946843&4.9.13.19.937.3119&2965&154&9566&2037556675&4.3.125.17.79.907&609&1516 \cr
 & &16.3.5.7.11.13.19.593&2965&5928& & &32.9.7.29.79.379&3411&1264 \cr
\noalign{\hrule}
 & &125.49.13.17.1499&27071&52554& & &11.13.47.53.59.97&669&20 \cr
9549&2029083875&4.3.11.19.23.107.461&2457&2350&9567&2038606999&8.3.5.47.97.223&5289&5192 \cr
 & &16.81.25.7.13.47.461&3807&3688& & &128.9.5.11.41.43.59&1935&2624 \cr
\noalign{\hrule}
 & &9.5.17.23.29.41.97&983&15048& & &3.29.41.59.89.109&133&44 \cr
9550&2029284135&16.81.11.19.983&1261&278&9568&2041604553&8.7.11.19.29.41.109&649&540 \cr
 & &64.11.13.97.139&139&4576& & &64.27.5.7.121.19.59&5445&4256 \cr
\noalign{\hrule}
 & &27.5.17.23.79.487&4495&3784& & &5.7.19.23.31.59.73&9801&5494 \cr
9551&2030797305&16.3.25.11.23.29.31.43&917&158&9569&2042142515&4.81.121.31.41.67&2681&3952 \cr
 & &64.7.29.31.79.131&4061&6496& & &128.9.7.11.13.19.383&4979&6336 \cr
\noalign{\hrule}
 & &9.31.431.16889&231&200& & &9.7.11.19.109.1423&10679&12102 \cr
9552&2030885361&16.27.25.7.11.16889&3247&13642&9570&2042293869&4.27.7.59.181.2017&163&1430 \cr
 & &64.5.17.19.191.359&68569&51680& & &16.5.11.13.163.2017&10085&16952 \cr
\noalign{\hrule}
 & &25.23.59.139.431&4609&5304& & &9.11.19.41.71.373&7327&7966 \cr
9553&2032412825&16.3.5.11.13.17.59.419&203&4812&9571&2042395443&4.7.11.17.19.431.569&81863&24540 \cr
 & &128.9.7.13.29.401&46917&12992& & &32.3.5.71.409.1153&5765&6544 \cr
\noalign{\hrule}
 & &25.11.17.29.53.283&2789&324& & &27.5.7.11.23.83.103&3757&10148 \cr
9554&2033489425&8.81.5.53.2789&6893&7052&9572&2043937665&8.13.289.23.43.59&225&166 \cr
 & &64.27.41.43.61.113&107543&97632& & &32.9.25.13.17.43.83&731&1040 \cr
\noalign{\hrule}
 & &9.7.17.23.31.2663&1595&1068& & &5.11.211.353.499&2047&1836 \cr
9555&2033528049&8.27.5.7.11.23.29.89&1&622&9573&2044185935&8.27.5.17.23.89.499&77&422 \cr
 & &32.5.11.29.311&45095&176& & &32.9.7.11.17.89.211&1071&1424 \cr
\noalign{\hrule}
 & &3.5.11.13.361.37.71&1411&2494& & &9.5.13.17.29.47.151&19153&20258 \cr
9556&2034204315&4.13.17.29.37.43.83&963&1444&9574&2046810285&4.7.47.107.179.1447&7579&2550 \cr
 & &32.9.17.361.43.107&2193&1712& & &16.3.25.11.13.17.53.179&1969&2120 \cr
\noalign{\hrule}
 & &3.17.107.163.2287&2529&242& & &3.13.23.97.101.233&29&262 \cr
9557&2034265917&4.27.121.107.281&1585&1304&9575&2047582797&4.13.23.29.101.131&99&200 \cr
 & &64.5.121.163.317&1585&3872& & &64.9.25.11.29.131&3799&26400 \cr
\noalign{\hrule}
 & &9.5.49.13.19.37.101&7199&9118& & &3.5.7.19.31.113.293&1931&2024 \cr
9558&2035300995&4.5.13.23.47.97.313&3333&2072&9576&2047626105&16.11.19.23.293.1931&55409&18720 \cr
 & &64.3.7.11.37.101.313&313&352& & &1024.9.5.13.67.827&32253&34304 \cr
\noalign{\hrule}
}%
}
$$
\eject
\vglue -23 pt
\noindent\hskip 1 in\hbox to 6.5 in{\ 9577 -- 9612 \hfill\fbd 2048234045 -- 2071430163\frb}
\vskip -9 pt
$$
\vbox{
\nointerlineskip
\halign{\strut
    \vrule \ \ \hfil \frb #\ 
   &\vrule \hfil \ \ \fbb #\frb\ 
   &\vrule \hfil \ \ \frb #\ \hfil
   &\vrule \hfil \ \ \frb #\ 
   &\vrule \hfil \ \ \frb #\ \ \vrule \hskip 2 pt
   &\vrule \ \ \hfil \frb #\ 
   &\vrule \hfil \ \ \fbb #\frb\ 
   &\vrule \hfil \ \ \frb #\ \hfil
   &\vrule \hfil \ \ \frb #\ 
   &\vrule \hfil \ \ \frb #\ \vrule \cr%
\noalign{\hrule}
 & &5.11.13.101.113.251&70929&54614& & &3.5.7.11.13.19.31.233&17717&4418 \cr
9577&2048234045&4.27.7.37.47.71.83&6275&2938&9595&2060613555&4.49.2209.2531&52855&55386 \cr
 & &16.9.25.7.13.113.251&45&56& & &16.9.5.11.17.961.181&1581&1448 \cr
\noalign{\hrule}
 & &5.7.37.47.67.503&927&2222& & &3.5.121.13.23.29.131&387&1054 \cr
9578&2051211365&4.9.11.101.103.503&199&1310&9596&2061660315&4.27.5.11.13.17.31.43&707&1048 \cr
 & &16.3.5.103.131.199&20497&3144& & &64.7.17.43.101.131&5117&3232 \cr
\noalign{\hrule}
 & &25.49.17.29.43.79&1179&164& & &3.5.11.13.31.101.307&57&44 \cr
9579&2051533225&8.9.5.7.41.43.131&143&158&9597&2061810465&8.9.5.121.19.31.307&13447&34138 \cr
 & &32.3.11.13.41.79.131&5109&7216& & &32.7.169.17.101.113&1921&1456 \cr
\noalign{\hrule}
 & &5.7.59.263.3779&75&338& & &27.31.97.109.233&455&2552 \cr
9580&2052356005&4.3.125.169.3779&8673&12452&9598&2061957033&16.3.5.7.11.13.29.109&1165&34 \cr
 & &32.9.49.11.59.283&2547&1232& & &64.25.7.17.233&2975&32 \cr
\noalign{\hrule}
 & &9.7.11.29.41.47.53&785&202& & &9.5.7.13.41.71.173&15029&22814 \cr
9581&2052526707&4.3.5.29.41.101.157&1175&7612&9599&2062254285&4.49.11.17.19.61.113&205&1038 \cr
 & &32.125.11.47.173&173&2000& & &16.3.5.19.41.61.173&61&152 \cr
\noalign{\hrule}
 & &121.19.29.41.751&2205&2756& & &27.11.13.29.113.163&313033&316310 \cr
9582&2052866761&8.9.5.49.13.53.751&561&190&9600&2062357011&4.5.7.47.197.227.673&1131&458 \cr
 & &32.27.25.7.11.13.17.19&4725&3536& & &16.3.5.13.29.47.197.229&9259&9160 \cr
\noalign{\hrule}
 & &27.5.169.17.67.79&33583&56398& & &5.17.59.101.4073&97&198 \cr
9583&2052916515&4.11.43.71.163.173&2637&9646&9601&2063035595&4.9.11.17.97.4073&24875&19928 \cr
 & &16.9.7.11.13.53.293&2051&4664& & &64.3.125.47.53.199&31641&37600 \cr
\noalign{\hrule}
 & &3.13.17.37.97.863&5561&5658& & &27.125.7.11.17.467&247&212 \cr
9584&2053514541&4.9.17.23.37.41.67.83&16393&16060&9602&2063147625&8.25.11.13.19.53.467&423&5648 \cr
 & &32.5.11.169.41.67.73.97&47905&47888& & &256.9.47.53.353&16591&6784 \cr
\noalign{\hrule}
 & &9.5.11.23.61.2957&2485&472& & &3.25.11.17.367.401&3123&1288 \cr
9585&2053592145&16.3.25.7.23.59.71&2957&1232&9603&2064017175&16.27.5.7.17.23.347&2843&5138 \cr
 & &512.49.11.2957&49&256& & &64.49.367.2843&2843&1568 \cr
\noalign{\hrule}
 & &9.11.169.17.31.233&19413&12190& & &9.19.529.29.787&8851&6490 \cr
9586&2054416221&4.243.5.23.53.719&481&238&9604&2064545757&4.3.5.11.19.53.59.167&253&754 \cr
 & &16.5.7.13.17.23.37.53&4255&2968& & &16.5.121.13.23.29.59&1573&2360 \cr
\noalign{\hrule}
 & &9.11.13.17.29.41.79&313&5110& & &25.11.17.103.4289&50619&22294 \cr
9587&2055116349&4.5.7.73.79.313&41&354&9605&2065260725&4.3.47.71.157.359&3869&3510 \cr
 & &16.3.7.41.59.73&59&4088& & &16.81.5.13.53.71.73&50297&46008 \cr
\noalign{\hrule}
 & &81.19.31.41.1051&15521&17060& & &9.25.11.13.19.31.109&5069&5068 \cr
9588&2055828519&8.5.11.17.41.83.853&4343&5040&9606&2065667175&8.3.25.7.11.13.19.37.137.181&101&6074 \cr
 & &256.9.25.7.43.83.101&89225&90496& & &32.7.37.101.137.3037&416069&418544 \cr
\noalign{\hrule}
 & &9.49.11.29.47.311&7277&7340& & &9.5.17.19.23.37.167&23287&20114 \cr
9589&2056304943&8.5.7.11.19.29.367.383&2773&77274&9607&2065670595&4.3.5.11.29.73.89.113&1169&74 \cr
 & &32.729.5.47.53.59&4779&4240& & &16.7.29.37.89.167&623&232 \cr
\noalign{\hrule}
 & &5.7.13.31.211.691&109&108& & &3.25.31.53.97.173&15257&12250 \cr
9590&2056523105&8.27.5.13.109.211.691&2387&5842&9608&2067838725&4.3125.49.11.19.73&6417&3292 \cr
 & &32.9.7.11.23.31.109.127&22563&22352& & &32.9.7.11.23.31.823&9053&7728 \cr
\noalign{\hrule}
 & &729.43.211.311&209&520& & &9.25.13.37.97.197&5131&12814 \cr
9591&2057021487&16.5.11.13.19.43.211&119&936&9609&2068071525&4.3.5.7.43.149.733&3707&1508 \cr
 & &256.9.7.11.169.17&1859&15232& & &32.11.13.29.43.337&13717&5392 \cr
\noalign{\hrule}
 & &3.5.19.29.241.1033&59551&38584& & &7.1331.53.59.71&49933&20610 \cr
9592&2057596545&16.7.13.17.31.53.113&205&198&9610&2068532389&4.9.5.13.23.167.229&1429&742 \cr
 & &64.9.5.11.17.41.53.113&65879&66912& & &16.3.5.7.23.53.1429&1429&2760 \cr
\noalign{\hrule}
 & &25.7.11.13.19.61.71&3303&1028& & &5.7.11.19.131.2161&223&432 \cr
9593&2059282225&8.9.11.19.257.367&497&130&9611&2070810665&32.27.7.223.2161&1261&3422 \cr
 & &32.3.5.7.13.71.257&257&48& & &128.9.13.29.59.97&51507&24128 \cr
\noalign{\hrule}
 & &3.25.13.23.139.661&1407&2068& & &9.11.529.37.1069&8183&11390 \cr
9594&2060386575&8.9.7.11.13.23.47.67&451&152&9612&2071430163&4.3.5.49.11.17.67.167&20539&10690 \cr
 & &128.7.121.19.41.47&34727&57152& & &16.25.19.23.47.1069&475&376 \cr
\noalign{\hrule}
}%
}
$$
\eject
\vglue -23 pt
\noindent\hskip 1 in\hbox to 6.5 in{\ 9613 -- 9648 \hfill\fbd 2071883709 -- 2095998751\frb}
\vskip -9 pt
$$
\vbox{
\nointerlineskip
\halign{\strut
    \vrule \ \ \hfil \frb #\ 
   &\vrule \hfil \ \ \fbb #\frb\ 
   &\vrule \hfil \ \ \frb #\ \hfil
   &\vrule \hfil \ \ \frb #\ 
   &\vrule \hfil \ \ \frb #\ \ \vrule \hskip 2 pt
   &\vrule \ \ \hfil \frb #\ 
   &\vrule \hfil \ \ \fbb #\frb\ 
   &\vrule \hfil \ \ \frb #\ \hfil
   &\vrule \hfil \ \ \frb #\ 
   &\vrule \hfil \ \ \frb #\ \vrule \cr%
\noalign{\hrule}
 & &9.49.19.37.41.163&425&278& & &7.121.17.19.29.263&1693&1200 \cr
9613&2071883709&4.3.25.17.41.139.163&67081&75394&9631&2086602287&32.3.25.7.11.19.1693&115&1578 \cr
 & &16.49.11.23.1369.149&3427&3256& & &128.9.125.23.263&2875&576 \cr
\noalign{\hrule}
 & &9.11.13.19.29.37.79&8275&46274& & &3.11.13.17.137.2089&291&1798 \cr
9614&2072807451&4.25.17.331.1361&589&6216&9632&2087205549&4.9.13.17.29.31.97&2017&2420 \cr
 & &64.3.5.7.19.31.37&35&992& & &32.5.121.97.2017&10085&17072 \cr
\noalign{\hrule}
 & &5.67.89.197.353&399&46& & &9.5.11.17.29.43.199&235&496 \cr
9615&2073364915&4.3.7.19.23.67.197&65&132&9633&2088207495&32.25.11.31.47.199&507&1682 \cr
 & &32.9.5.7.11.13.19.23&3059&20592& & &128.3.169.841.31&5239&1856 \cr
\noalign{\hrule}
 & &81.19.71.83.229&715&634& & &3.5.7.13.19.23.31.113&173&264 \cr
9616&2076875883&4.5.11.13.83.229.317&73&156&9634&2089557015&16.9.5.11.31.113.173&2141&5644 \cr
 & &32.3.5.11.169.73.317&61685&55792& & &128.11.17.83.2141&36397&58432 \cr
\noalign{\hrule}
 & &3.5.7.13.17.37.41.59&277&318& & &27.13.17.23.97.157&45733&1210 \cr
9617&2076917115&4.9.13.37.53.59.277&217&550&9635&2090043189&4.5.121.19.29.83&471&442 \cr
 & &16.25.7.11.31.53.277&15235&13144& & &16.3.5.11.13.17.19.157&55&152 \cr
\noalign{\hrule}
 & &25.29.97.109.271&2769&44& & &9.5.7.23.31.41.227&6383&13420 \cr
9618&2077330175&8.3.11.13.71.271&1455&1526&9636&2090305665&8.3.25.11.13.61.491&899&574 \cr
 & &32.9.5.7.13.97.109&63&208& & &32.7.11.29.31.41.61&671&464 \cr
\noalign{\hrule}
 & &11.17.29.43.59.151&2691&1688& & &9.25.7.11.23.29.181&3&178 \cr
9619&2077480801&16.9.11.13.23.43.211&815&604&9637&2091595275&4.27.11.23.29.89&2125&4706 \cr
 & &128.3.5.13.23.151.163&6357&7360& & &16.125.13.17.181&13&680 \cr
\noalign{\hrule}
 & &27.121.23.139.199&196105&190732& & &5.7.97.383.1609&153&1762 \cr
9620&2078475201&8.5.7.13.41.431.1163&5709&106&9638&2092158565&4.9.7.17.97.881&10075&8426 \cr
 & &32.3.7.11.41.53.173&9169&4592& & &16.3.25.11.13.31.383&465&1144 \cr
\noalign{\hrule}
 & &7.11.17.67.151.157&97&90& & &243.5.13.139.953&1927&1232 \cr
9621&2079175021&4.9.5.67.97.151.157&5137&4980&9639&2092316265&32.7.11.41.47.953&251&702 \cr
 & &32.27.25.11.83.97.467&201275&201744& & &128.27.7.13.47.251&1757&3008 \cr
\noalign{\hrule}
 & &3.125.7.11.13.29.191&283&718& & &9.11.13.29.47.1193&31&1162 \cr
9622&2079202125&4.25.191.283.359&48411&53186&9640&2092737933&4.3.7.11.31.47.83&65&76 \cr
 & &16.27.7.11.29.131.163&1179&1304& & &32.5.7.13.19.31.83&12865&2128 \cr
\noalign{\hrule}
 & &7.13.29.53.107.139&135&242& & &7.121.13.67.2837&1919&918 \cr
9623&2080241891&4.27.5.7.121.53.139&89&884&9641&2092959869&4.27.11.17.19.67.101&217&520 \cr
 & &32.9.121.13.17.89&2057&12816& & &64.9.5.7.13.17.19.31&2635&5472 \cr
\noalign{\hrule}
 & &121.169.19.23.233&2997&214& & &27.107.691.1049&185&506 \cr
9624&2082137629&4.81.37.107.233&365&598&9642&2094117651&4.9.5.11.23.37.1049&691&358 \cr
 & &16.9.5.13.23.37.73&365&2664& & &16.5.11.23.179.691&1969&920 \cr
\noalign{\hrule}
 & &27.25.49.13.29.167&46321&7954& & &3.5.7.11.13.199.701&61&138 \cr
9625&2082368925&4.11.41.97.4211&2331&1880&9643&2094577485&4.9.5.13.23.61.701&5929&3184 \cr
 & &64.9.5.7.37.47.97&1739&3104& & &128.49.121.23.199&253&448 \cr
\noalign{\hrule}
 & &7.121.47.113.463&10485&11276& & &5.13.23.29.31.1559&11289&47146 \cr
9626&2082767071&8.9.5.121.233.2819&47&652&9644&2095303795&4.3.11.53.71.2143&1363&780 \cr
 & &64.3.47.163.2819&8457&5216& & &32.9.5.13.29.47.71&423&1136 \cr
\noalign{\hrule}
 & &49.17.41.139.439&729&17270& & &5.59.113.239.263&8041&7476 \cr
9627&2084050213&4.729.5.7.11.157&289&278&9645&2095338095&8.3.7.11.17.43.89.239&3945&118 \cr
 & &16.9.5.289.139.157&785&1224& & &32.9.5.7.11.59.263&63&176 \cr
\noalign{\hrule}
 & &81.5.7.11.13.53.97&5969&1888& & &3.5.11.29.131.3343&64657&45662 \cr
9628&2084187105&64.5.13.47.59.127&1067&4902&9646&2095509405&4.289.19.41.79.83&3343&9900 \cr
 & &256.3.11.19.43.97&817&128& & &32.9.25.11.17.3343&51&80 \cr
\noalign{\hrule}
 & &3.5.19.31.53.61.73&2211&2242& & &17.19.41.79.2003&28303&53820 \cr
9629&2085139515&4.9.5.11.361.53.59.67&15407&194102&9647&2095532591&8.9.5.11.13.23.31.83&505&574 \cr
 & &16.7.31.37.43.61.71&2627&2408& & &32.3.25.7.11.31.41.101&16275&17776 \cr
\noalign{\hrule}
 & &5.29.101.109.1307&4851&6158& & &7.11.29.43.83.263&17773&4056 \cr
9630&2086370635&4.9.5.49.11.29.3079&10201&964&9648&2095998751&16.3.49.169.2539&2871&5410 \cr
 & &32.3.7.10201.241&1687&4848& & &64.27.5.11.29.541&2705&864 \cr
\noalign{\hrule}
}%
}
$$
\eject
\vglue -23 pt
\noindent\hskip 1 in\hbox to 6.5 in{\ 9649 -- 9684 \hfill\fbd 2096174355 -- 2121591563\frb}
\vskip -9 pt
$$
\vbox{
\nointerlineskip
\halign{\strut
    \vrule \ \ \hfil \frb #\ 
   &\vrule \hfil \ \ \fbb #\frb\ 
   &\vrule \hfil \ \ \frb #\ \hfil
   &\vrule \hfil \ \ \frb #\ 
   &\vrule \hfil \ \ \frb #\ \ \vrule \hskip 2 pt
   &\vrule \ \ \hfil \frb #\ 
   &\vrule \hfil \ \ \fbb #\frb\ 
   &\vrule \hfil \ \ \frb #\ \hfil
   &\vrule \hfil \ \ \frb #\ 
   &\vrule \hfil \ \ \frb #\ \vrule \cr%
\noalign{\hrule}
 & &3.5.121.223.5179&3497&1682& & &27.5.121.23.71.79&5149&3016 \cr
9649&2096174355&4.13.841.223.269&1485&4982&9667&2107329345&16.121.13.19.29.271&63&184 \cr
 & &16.27.5.11.29.47.53&1363&3816& & &256.9.7.23.29.271&1897&3712 \cr
\noalign{\hrule}
 & &49.11.19.311.659&48645&16354& & &27.121.13.841.59&395&1178 \cr
9650&2098882709&4.9.5.13.17.23.37.47&1381&1862&9668&2107368549&4.5.19.29.31.59.79&1881&2780 \cr
 & &16.3.5.49.17.19.1381&1381&2040& & &32.9.25.11.361.139&3475&5776 \cr
\noalign{\hrule}
 & &81.5.49.361.293&15301&944& & &5.11.13.17.19.23.397&3103&2058 \cr
9651&2099065185&32.9.11.13.59.107&2513&3800&9669&2108758795&4.3.343.17.23.29.107&17251&15432 \cr
 & &512.25.7.19.359&1795&256& & &64.9.7.13.643.1327&40509&42464 \cr
\noalign{\hrule}
 & &3.7.23.43.271.373&17369&1330& & &9.5.13.37.103.947&8107&10948 \cr
9652&2099392827&4.5.49.11.19.1579&765&814&9670&2111274945&8.3.7.121.13.17.23.67&947&8920 \cr
 & &16.9.25.121.17.19.37&35853&24200& & &128.5.11.223.947&223&704 \cr
\noalign{\hrule}
 & &7.17.19.23.71.569&123723&124930& & &27.25.7.17.97.271&319&494 \cr
9653&2100869197&4.9.5.7.13.961.59.233&209&442&9671&2111503275&4.9.11.13.17.19.29.97&9485&118 \cr
 & &16.3.5.11.169.17.19.31.59&20119&20280& & &16.5.7.13.59.271&13&472 \cr
\noalign{\hrule}
 & &27.5.7.121.17.23.47&1871&2476& & &9.25.13.23.31.1013&16687&6612 \cr
9654&2101318065&8.17.47.619.1871&1357&30450&9672&2112636825&8.27.11.19.29.37.41&5065&3548 \cr
 & &32.3.25.7.23.29.59&145&944& & &64.5.19.887.1013&887&608 \cr
\noalign{\hrule}
 & &3.25.7.23.151.1153&13247&13178& & &3.5.11.47.2809.97&3913&1422 \cr
9655&2102293725&4.11.13.599.1019.1153&3871&2718&9673&2113028115&4.27.7.13.43.53.79&995&436 \cr
 & &16.9.49.13.79.151.1019&13247&13272& & &32.5.7.79.109.199&8611&22288 \cr
\noalign{\hrule}
 & &729.23.37.3389&4489&1100& & &3.11.13.41.43.2797&1121&7270 \cr
9656&2102464431&8.3.25.11.37.4489&893&782&9674&2115446619&4.5.19.41.59.727&1573&846 \cr
 & &32.11.17.19.23.47.67&15181&11792& & &16.9.5.121.13.19.47&517&2280 \cr
\noalign{\hrule}
 & &11.13.19.41.43.439&10797&8080& & &3.7.13.67.103.1123&2255&46052 \cr
9657&2102841169&32.3.5.41.59.61.101&957&1462&9675&2115701679&8.5.11.29.41.397&793&396 \cr
 & &128.9.11.17.29.43.61&4437&3904& & &64.9.5.121.13.61&7381&480 \cr
\noalign{\hrule}
 & &49.11.31.37.41.83&4275&8342& & &5.7.19.23.43.3217&103&198 \cr
9658&2103846899&4.9.25.19.41.43.97&1643&7172&9676&2115772645&4.9.11.23.103.3217&2793&424 \cr
 & &32.3.5.11.31.53.163&795&2608& & &64.27.49.11.19.53&4081&864 \cr
\noalign{\hrule}
 & &27.5.19.23.53.673&4121&1936& & &25.49.11.53.2963&2829&134 \cr
9659&2104292655&32.3.121.13.53.317&767&184&9677&2116100525&4.3.5.23.41.53.67&1029&1144 \cr
 & &512.11.169.23.59&9971&2816& & &64.9.343.11.13.67&819&2144 \cr
\noalign{\hrule}
 & &729.5.7.11.13.577&20519&26866& & &9.25.7.11.23.47.113&59&106 \cr
9660&2105268165&4.49.289.19.71.101&4485&1006&9678&2116300725&4.3.5.7.23.53.59.113&32549&14356 \cr
 & &16.3.5.13.23.101.503&2323&4024& & &32.121.37.97.269&9953&17072 \cr
\noalign{\hrule}
 & &3.125.7.13.17.19.191&4631&6256& & &3.19.31.61.73.269&663&724 \cr
9661&2105273625&32.7.11.289.23.421&4797&1850&9679&2116613319&8.9.13.17.31.181.269&473&4100 \cr
 & &128.9.25.11.13.37.41&1517&2112& & &64.25.11.41.43.181&44075&63712 \cr
\noalign{\hrule}
 & &9.49.29.193.853&325&528& & &9.13.19.83.11491&5831&5660 \cr
9662&2105440281&32.27.25.7.11.13.193&67&1418&9680&2120192919&8.5.343.13.17.83.283&5757&1298 \cr
 & &128.5.13.67.709&46085&4288& & &32.3.11.19.59.101.283&16697&17776 \cr
\noalign{\hrule}
 & &9.7.11.13.103.2269&13561&11398& & &3.5.13.841.67.193&517&324 \cr
9663&2105466363&4.3.13.41.71.139.191&2269&5180&9681&2120619345&8.243.5.11.13.47.67&841&3514 \cr
 & &32.5.7.37.139.2269&695&592& & &32.7.841.47.251&1757&752 \cr
\noalign{\hrule}
 & &3.11.41.61.97.263&7429&10322& & &3.5.7.17.47.131.193&46157&3058 \cr
9664&2105499363&4.13.17.19.23.41.397&8357&7920&9682&2121117285&4.11.101.139.457&2583&2444 \cr
 & &128.9.5.11.13.17.61.137&5343&5440& & &32.9.7.13.41.47.101&1313&1968 \cr
\noalign{\hrule}
 & &3.5.11.19.31.47.461&3481&3434& & &9.11.13.41.61.659&48877&16364 \cr
9665&2105707395&4.11.17.19.31.3481.101&2691&862&9683&2121180633&8.37.1321.4091&1385&2706 \cr
 & &16.9.13.23.59.101.431&137057&134472& & &32.3.5.11.37.41.277&1385&592 \cr
\noalign{\hrule}
 & &9.49.17.37.71.107&2223&404& & &49.13.947.3517&341&46062 \cr
9666&2107324233&8.81.49.13.19.101&425&506&9684&2121591563&4.27.11.31.853&1285&1274 \cr
 & &32.25.11.13.17.23.101&14443&9200& & &16.9.5.49.13.31.257&1285&2232 \cr
\noalign{\hrule}
}%
}
$$
\eject
\vglue -23 pt
\noindent\hskip 1 in\hbox to 6.5 in{\ 9685 -- 9720 \hfill\fbd 2121853125 -- 2142860895\frb}
\vskip -9 pt
$$
\vbox{
\nointerlineskip
\halign{\strut
    \vrule \ \ \hfil \frb #\ 
   &\vrule \hfil \ \ \fbb #\frb\ 
   &\vrule \hfil \ \ \frb #\ \hfil
   &\vrule \hfil \ \ \frb #\ 
   &\vrule \hfil \ \ \frb #\ \ \vrule \hskip 2 pt
   &\vrule \ \ \hfil \frb #\ 
   &\vrule \hfil \ \ \fbb #\frb\ 
   &\vrule \hfil \ \ \frb #\ \hfil
   &\vrule \hfil \ \ \frb #\ 
   &\vrule \hfil \ \ \frb #\ \vrule \cr%
\noalign{\hrule}
 & &3.3125.49.31.149&13&638& & &37.59.61.83.193&96481&84708 \cr
9685&2121853125&4.5.7.11.13.29.149&837&802&9703&2133138097&8.9.49.11.13.179.181&2545&4514 \cr
 & &16.27.13.29.31.401&5213&2088& & &32.3.5.49.37.61.509&2545&2352 \cr
\noalign{\hrule}
 & &9.25.7.17.19.43.97&539&1304& & &3.11.13.521.9551&22925&51578 \cr
9686&2121891975&16.5.343.11.43.163&12937&5928&9704&2134734459&4.25.7.17.37.41.131&521&504 \cr
 & &256.3.13.17.19.761&761&1664& & &64.9.49.37.131.521&4847&4704 \cr
\noalign{\hrule}
 & &125.7.169.113.127&867&22& & &5.7.19.37.229.379&639&506 \cr
9687&2122154125&4.3.25.11.289.113&1421&1404&9705&2135492555&4.9.11.23.37.71.379&6025&1856 \cr
 & &32.81.49.11.13.17.29&6237&7888& & &512.3.25.23.29.241&20967&29440 \cr
\noalign{\hrule}
 & &27.25.7.13.17.19.107&11891&11756& & &25.49.169.17.607&14839&14904 \cr
9688&2122909425&8.5.7.11.19.23.47.2939&12519&8054&9706&2136290975&16.81.5.11.13.17.19.23.71&21041&46 \cr
 & &32.9.11.13.23.107.4027&4027&4048& & &64.3.529.53.397&21041&50784 \cr
\noalign{\hrule}
 & &25.11.29.41.43.151&491&534& & &3.11.17.19.43.59.79&2575&1454 \cr
9689&2123048675&4.3.11.29.89.151.491&2121&460&9707&2136308757&4.25.11.43.103.727&18463&12798 \cr
 & &32.9.5.7.23.101.491&30933&37168& & &16.81.5.37.79.499&4995&3992 \cr
\noalign{\hrule}
 & &3.5.49.41.251.281&899&858& & &25.11.43.281.643&84429&96254 \cr
9690&2125451685&4.9.5.7.11.13.29.31.281&43723&29618&9708&2136576475&4.27.17.19.53.59.149&185&716 \cr
 & &16.11.23.59.251.1901&14927&15208& & &32.3.5.19.37.149.179&26671&33744 \cr
\noalign{\hrule}
 & &5.11.37.41.73.349&7953&4960& & &11.31.61.127.809&69&740 \cr
9691&2125673495&64.3.25.121.31.241&73&48&9709&2137157143&8.3.5.23.31.37.127&2427&2272 \cr
 & &2048.9.31.73.241&7471&9216& & &512.9.23.71.809&1633&2304 \cr
\noalign{\hrule}
 & &5.11.13.29.109.941&3933&772& & &9.11.17.29.193.227&775&962 \cr
9692&2126768215&8.9.11.13.19.23.193&6587&6730&9710&2138283477&4.25.13.29.31.37.227&3723&772 \cr
 & &32.3.5.7.19.673.941&2019&2128& & &32.3.5.17.37.73.193&365&592 \cr
\noalign{\hrule}
 & &81.5.7.169.23.193&4099&8446& & &9.5.11.19.31.41.179&377&212 \cr
9693&2126791485&4.3.13.41.103.4099&2849&1250&9711&2139722145&8.3.13.29.41.53.179&5975&1364 \cr
 & &16.625.7.11.37.103&4625&9064& & &64.25.11.13.31.239&1195&416 \cr
\noalign{\hrule}
 & &5.7.11.23.37.43.151&1161&104& & &3.11.13.43.311.373&703&230 \cr
9694&2127334055&16.27.13.37.1849&943&906&9712&2139907341&4.5.13.19.23.37.373&3303&3784 \cr
 & &64.81.13.23.41.151&1053&1312& & &64.9.5.11.23.43.367&1835&2208 \cr
\noalign{\hrule}
 & &25.49.29.89.673&621&17446& & &9.49.11.433.1019&61049&39832 \cr
9695&2127840925&4.27.7.11.13.23.61&355&194&9713&2140392177&16.13.41.383.1489&1745&3234 \cr
 & &16.3.5.11.13.71.97&20661&1144& & &64.3.5.49.11.41.349&1745&1312 \cr
\noalign{\hrule}
 & &3.7.11.29.31.37.277&47759&12350& & &9.5.169.17.29.571&539&46 \cr
9696&2128399581&4.25.13.19.163.293&261&554&9714&2140830315&4.49.11.13.23.571&435&136 \cr
 & &16.9.5.13.19.29.277&95&312& & &64.3.5.49.11.17.29&539&32 \cr
\noalign{\hrule}
 & &3.7.17.23.53.67.73&1391&150& & &9.7.13.31.37.43.53&55&1388 \cr
9697&2128480053&4.9.25.7.13.53.107&253&1072&9715&2140876647&8.3.5.7.11.53.347&817&1612 \cr
 & &128.11.23.67.107&1177&64& & &64.11.13.19.31.43&19&352 \cr
\noalign{\hrule}
 & &3.5.11.19.41.59.281&7439&1834& & &3.5.11.4913.19.139&13603&10962 \cr
9698&2130981765&4.7.41.43.131.173&12213&10450&9716&2140913445&4.81.7.11.29.61.223&289&278 \cr
 & &16.9.25.7.11.19.23.59&161&120& & &16.289.29.61.139.223&1769&1784 \cr
\noalign{\hrule}
 & &5.7.13.19.23.71.151&12617&11694& & &9.25.121.151.521&185&336 \cr
9699&2131710035&4.3.5.11.19.31.37.1949&1633&8112&9717&2141817975&32.27.125.7.121.37&551&3926 \cr
 & &128.9.169.23.37.71&481&576& & &128.7.13.19.29.151&3857&832 \cr
\noalign{\hrule}
 & &3.25.11.19.31.41.107&2011&244& & &81.5.7.11.23.29.103&85&806 \cr
9700&2131752975&8.5.61.107.2011&1273&738&9718&2142440685&4.25.13.17.23.29.31&4187&3288 \cr
 & &32.9.19.41.61.67&201&976& & &64.3.17.53.79.137&10823&28832 \cr
\noalign{\hrule}
 & &9.25.11.13.23.43.67&2269&12136& & &9.7.11.13.29.59.139&925&786 \cr
9701&2132012025&16.3.5.37.41.2269&3311&3496&9719&2142601461&4.27.25.7.11.13.37.131&163&2242 \cr
 & &256.7.11.19.23.41.43&779&896& & &16.5.19.59.131.163&2489&6520 \cr
\noalign{\hrule}
 & &9.121.73.139.193&26335&53162& & &9.5.49.23.29.31.47&485983&485836 \cr
9702&2132666019&4.5.19.23.229.1399&5211&26966&9720&2142860895&8.3.5.13.59.8237.9343&122507&1048 \cr
 & &16.27.97.139.193&97&24& & &128.7.11.37.43.59.131&93869&92224 \cr
\noalign{\hrule}
}%
}
$$
\eject
\vglue -23 pt
\noindent\hskip 1 in\hbox to 6.5 in{\ 9721 -- 9756 \hfill\fbd 2144138715 -- 2163550543\frb}
\vskip -9 pt
$$
\vbox{
\nointerlineskip
\halign{\strut
    \vrule \ \ \hfil \frb #\ 
   &\vrule \hfil \ \ \fbb #\frb\ 
   &\vrule \hfil \ \ \frb #\ \hfil
   &\vrule \hfil \ \ \frb #\ 
   &\vrule \hfil \ \ \frb #\ \ \vrule \hskip 2 pt
   &\vrule \ \ \hfil \frb #\ 
   &\vrule \hfil \ \ \fbb #\frb\ 
   &\vrule \hfil \ \ \frb #\ \hfil
   &\vrule \hfil \ \ \frb #\ 
   &\vrule \hfil \ \ \frb #\ \vrule \cr%
\noalign{\hrule}
 & &27.5.31.37.61.227&421&726& & &17.31.41.137.727&59013&36476 \cr
9721&2144138715&4.81.121.227.421&3317&30784&9739&2152035593&8.9.11.79.83.829&85&164 \cr
 & &512.13.31.37.107&1391&256& & &64.3.5.11.17.41.829&2487&1760 \cr
\noalign{\hrule}
 & &3.5.11.19.109.6277&141&6136& & &27.5.13.23.71.751&19807&20558 \cr
9722&2144945055&16.9.13.19.47.59&275&2498&9740&2152302165&4.19.29.71.541.683&759&1300 \cr
 & &64.25.11.1249&6245&32& & &32.3.25.11.13.19.23.683&3415&3344 \cr
\noalign{\hrule}
 & &3.7.169.17.961.37&7975&11524& & &3.5.7.11.17.19.29.199&101&1494 \cr
9723&2145260481&8.25.11.29.31.43.67&12051&1666&9741&2152958115&4.27.17.19.83.101&175&148 \cr
 & &32.9.5.49.13.17.103&721&240& & &32.25.7.37.83.101&8383&2960 \cr
\noalign{\hrule}
 & &9.5.121.47.83.101&43&58& & &27.11.41.47.53.71&1115&2222 \cr
9724&2145335445&4.3.121.29.43.47.83&8453&3250&9742&2153636397&4.5.121.53.101.223&201&12020 \cr
 & &16.125.13.29.79.107&29783&21400& & &32.3.25.67.601&601&26800 \cr
\noalign{\hrule}
 & &13.73.811.2789&11473&47730& & &3.7.11.13.17.31.1361&83&134 \cr
9725&2146523171&4.3.5.7.11.37.43.149&949&690&9743&2153892741&4.11.13.67.83.1361&6615&5254 \cr
 & &16.9.25.13.23.43.73&989&1800& & &16.27.5.49.37.67.71&22311&19880 \cr
\noalign{\hrule}
 & &27.7.11.13.19.37.113&1411&17570& & &9.25.13.29.109.233&7157&10318 \cr
9726&2146997853&4.5.49.17.83.251&2033&2034&9744&2154300525&4.3.7.11.13.17.67.421&2071&932 \cr
 & &16.9.5.17.19.107.113.251&4267&4280& & &32.19.109.233.421&421&304 \cr
\noalign{\hrule}
 & &3.5.11.19.43.89.179&1391&2436& & &49.11.19.43.59.83&18005&49806 \cr
9727&2147578455&8.9.7.13.29.107.179&775&5966&9745&2156457611&4.9.5.13.277.2767&6785&4018 \cr
 & &32.25.13.19.31.157&4867&1040& & &16.3.25.49.23.41.59&1025&552 \cr
\noalign{\hrule}
 & &9.5.121.13.127.239&2071&3644& & &9.5.31.97.107.149&433&530 \cr
9728&2148537105&8.19.109.239.911&1815&2726&9746&2157327045&4.25.31.53.149.433&3103&7722 \cr
 & &32.3.5.121.29.47.109&1363&1744& & &16.27.11.13.29.53.107&2067&2552 \cr
\noalign{\hrule}
 & &7.19.41.313.1259&18377&5544& & &27.25.7.11.13.31.103&851&634 \cr
9729&2148847351&16.9.49.11.17.23.47&465&74&9747&2157430275&4.5.13.23.37.103.317&2387&18 \cr
 & &64.27.5.31.37.47&4995&46624& & &16.9.7.11.31.317&317&8 \cr
\noalign{\hrule}
 & &3.5.7.23.53.103.163&39&1180& & &9.5.7.169.23.41.43&1639&1156 \cr
9730&2148908055&8.9.25.13.59.103&3151&2926&9748&2158626015&8.3.11.13.289.41.149&43&490 \cr
 & &32.7.11.13.19.23.137&1507&3952& & &32.5.49.11.289.43&289&1232 \cr
\noalign{\hrule}
 & &9.73.911.3593&2303&5896& & &25.7.11.17.19.23.151&929&2544 \cr
9731&2150507511&16.49.11.47.67.73&1755&1822&9749&2159424575&32.3.5.7.11.53.929&4709&1794 \cr
 & &64.27.5.11.13.47.911&1833&1760& & &128.9.13.17.23.277&2493&832 \cr
\noalign{\hrule}
 & &9.7.11.19.29.43.131&6107&5290& & &9.7.11.29.41.43.61&151&520 \cr
9732&2150921619&4.3.5.7.11.529.31.197&841&4978&9750&2161291671&16.5.7.13.29.43.151&543&9272 \cr
 & &16.5.19.841.31.131&155&232& & &256.3.19.61.181&3439&128 \cr
\noalign{\hrule}
 & &13.37.89.109.461&115721&106020& & &5.11.19.23.293.307&47&162 \cr
9733&2151109441&8.9.5.19.31.97.1193&3293&286&9751&2161972285&4.81.47.293.307&8533&5896 \cr
 & &32.3.5.11.13.19.37.89&285&176& & &64.9.7.11.23.53.67&3551&2016 \cr
\noalign{\hrule}
 & &3.7.13.23.59.5807&200725&199958& & &3.25.7.11.13.83.347&1811&2006 \cr
9734&2151267027&4.25.49.11.31.37.61.149&1593&2132&9752&2162235075&4.5.7.17.59.83.1811&3303&3718 \cr
 & &32.27.13.31.37.41.59.61&20313&20336& & &16.9.11.169.367.1811&14313&14488 \cr
\noalign{\hrule}
 & &9.37.47.101.1361&3817&46540& & &9.7.11.19.277.593&1513&1534 \cr
9735&2151402111&8.5.11.13.179.347&12625&12972&9753&2162824587&4.3.13.17.19.59.89.593&35137&23870 \cr
 & &64.3.625.23.47.101&625&736& & &16.5.7.11.31.41.59.857&50563&50840 \cr
\noalign{\hrule}
 & &5.7.11.37.131.1153&3887&46548& & &9.169.17.23.3637&28897&54754 \cr
9736&2151607535&8.27.169.23.431&1297&2590&9754&2162963907&4.7.11.37.71.3911&1565&2346 \cr
 & &32.9.5.7.37.1297&1297&144& & &16.3.5.7.17.23.37.313&1565&2072 \cr
\noalign{\hrule}
 & &3.17.19.61.89.409&2849&1690& & &81.25.11.23.41.103&365&1498 \cr
9737&2151626709&4.5.7.11.169.37.409&2897&1602&9755&2163548475&4.125.7.41.73.107&2369&6756 \cr
 & &16.9.169.89.2897&2897&4056& & &32.3.7.23.103.563&563&112 \cr
\noalign{\hrule}
 & &9.7.29.53.71.313&2257&4316& & &7.121.17.31.37.131&11975&16452 \cr
9738&2151875313&8.3.13.37.53.61.83&469&220&9756&2163550543&8.9.25.17.457.479&2687&5082 \cr
 & &64.5.7.11.37.61.67&12395&21472& & &32.27.5.7.121.2687&2687&2160 \cr
\noalign{\hrule}
}%
}
$$
\eject
\vglue -23 pt
\noindent\hskip 1 in\hbox to 6.5 in{\ 9757 -- 9792 \hfill\fbd 2164265155 -- 2188928885\frb}
\vskip -9 pt
$$
\vbox{
\nointerlineskip
\halign{\strut
    \vrule \ \ \hfil \frb #\ 
   &\vrule \hfil \ \ \fbb #\frb\ 
   &\vrule \hfil \ \ \frb #\ \hfil
   &\vrule \hfil \ \ \frb #\ 
   &\vrule \hfil \ \ \frb #\ \ \vrule \hskip 2 pt
   &\vrule \ \ \hfil \frb #\ 
   &\vrule \hfil \ \ \fbb #\frb\ 
   &\vrule \hfil \ \ \frb #\ \hfil
   &\vrule \hfil \ \ \frb #\ 
   &\vrule \hfil \ \ \frb #\ \vrule \cr%
\noalign{\hrule}
 & &5.13.17.23.31.41.67&2233&11502& & &9.5.11.13.31.67.163&12157&10690 \cr
9757&2164265155&4.81.7.11.17.29.71&161&620&9775&2178575685&4.25.13.1069.12157&697&372 \cr
 & &32.3.5.49.23.29.31&87&784& & &32.3.17.31.41.12157&12157&11152 \cr
\noalign{\hrule}
 & &9.7.13.31.269.317&605&202& & &9.5.11.13.167.2029&727&560 \cr
9758&2164996197&4.3.5.7.121.101.317&899&52&9776&2180454705&32.25.7.727.2029&1027&1002 \cr
 & &32.5.13.29.31.101&101&2320& & &128.3.7.13.79.167.727&5089&5056 \cr
\noalign{\hrule}
 & &9.25.17.29.131.149&37&3762& & &81.25.7.11.71.197&2557&2368 \cr
9759&2165145075&4.81.11.17.19.37&355&274&9777&2180922975&128.3.11.37.71.2557&107&2450 \cr
 & &16.5.11.19.71.137&781&20824& & &512.25.49.37.107&3959&1792 \cr
\noalign{\hrule}
 & &7.11.169.163.1021&419&582& & &3.25.13.19.29.31.131&373&402 \cr
9760&2165662499&4.3.13.97.419.1021&1141&120&9778&2181670725&4.9.13.19.67.131.373&3277&5500 \cr
 & &64.9.5.7.163.419&419&1440& & &32.125.11.29.113.373&6215&5968 \cr
\noalign{\hrule}
 & &9.25.7.11.19.29.227&443&2032& & &13.29.59.233.421&37925&60168 \cr
9761&2166959025&32.19.29.127.443&5217&7630&9779&2181882599&16.3.25.23.37.41.109&1947&2522 \cr
 & &128.3.5.7.37.47.109&5123&2368& & &64.9.11.13.37.59.97&3589&3168 \cr
\noalign{\hrule}
 & &3.13.41.53.107.239&55187&72200& & &25.13.19.47.73.103&9603&8428 \cr
9762&2167235031&16.25.11.361.29.173&477&9992&9780&2182201775&8.9.49.11.43.97.103&1875&2554 \cr
 & &256.9.5.53.1249&3747&640& & &32.27.625.7.11.1277&31925&33264 \cr
\noalign{\hrule}
 & &3.7.23.73.137.449&9715&286& & &5.11.13.17.29.41.151&46737&1432 \cr
9763&2168886867&4.5.11.13.23.29.67&171&148&9781&2182296545&16.81.179.577&275&302 \cr
 & &32.9.5.13.19.37.67&27417&5360& & &64.3.25.11.151.179&537&160 \cr
\noalign{\hrule}
 & &49.29.31.41.1201&555&1454& & &7.11.17.61.151.181&24525&52762 \cr
9764&2169115291&4.3.5.37.727.1201&2639&3366&9782&2182353019&4.9.25.23.31.37.109&3077&302 \cr
 & &16.27.7.11.13.17.29.37&5049&3848& & &16.3.17.23.151.181&69&8 \cr
\noalign{\hrule}
 & &5.31.37.47.83.97&2853&1706& & &27.5.7.23.31.41.79&32351&31244 \cr
9765&2170106795&4.9.5.83.317.853&37777&33022&9783&2182389615&8.11.17.31.73.107.173&207&134 \cr
 & &16.3.11.19.37.79.1021&19399&20856& & &32.9.17.23.67.107.173&18511&18224 \cr
\noalign{\hrule}
 & &13.23.47.257.601&7327&6726& & &5.19.43.67.79.101&2191&690 \cr
9766&2170584221&4.3.17.19.59.257.431&6105&1222&9784&2183812405&4.3.25.7.23.101.313&1419&1106 \cr
 & &16.9.5.11.13.37.47.59&2655&3256& & &16.9.49.11.23.43.79&1127&792 \cr
\noalign{\hrule}
 & &9.5.11.31.353.401&391&12040& & &89.78961.311&25641&53320 \cr
9767&2172130785&16.3.25.7.17.23.43&1601&1624&9785&2185561519&16.9.5.7.11.31.37.43&281&622 \cr
 & &256.49.17.29.1601&78449&63104& & &64.3.5.37.281.311&185&96 \cr
\noalign{\hrule}
 & &9.5.7.11.47.13339&6917&6422& & &81.7.11.29.107.113&905&338 \cr
9768&2172322845&4.7.169.19.47.6917&5597&1320&9786&2186935443&4.5.169.29.107.181&853&6102 \cr
 & &64.3.5.11.13.19.29.193&7163&6176& & &16.27.13.113.853&853&104 \cr
\noalign{\hrule}
 & &5.11.17.89.151.173&114145&114318& & &25.7.11.569.1997&276579&272596 \cr
9769&2173825445&4.27.25.11.29.37.73.617&1211&6764&9787&2187364025&8.9.23.79.389.2963&29&8918 \cr
 & &32.3.7.19.37.73.89.173&4161&4144& & &32.3.343.13.29.79&3871&18096 \cr
\noalign{\hrule}
 & &27.5.7.37.97.641&3553&934& & &9.5.49.53.97.193&13943&14428 \cr
9770&2174018805&4.5.11.17.19.37.467&15981&23714&9788&2187829665&8.3.53.73.191.3607&349&10472 \cr
 & &16.3.7.71.167.761&11857&6088& & &128.7.11.17.73.349&13651&22336 \cr
\noalign{\hrule}
 & &5.49.31.37.71.109&312309&312806& & &27.5.19.43.83.239&69&26 \cr
9771&2174775085&4.27.7.13.43.53.227.269&14279&22&9789&2187921915&4.81.13.23.83.239&8987&10850 \cr
 & &16.3.11.13.43.109.131&5633&3432& & &16.25.7.11.13.19.31.43&1001&1240 \cr
\noalign{\hrule}
 & &9.5.17.19.43.3481&4105&4922& & &11.43.1951.2371&949&1422 \cr
9772&2175642405&4.25.23.59.107.821&13419&7106&9790&2188013333&4.9.13.73.79.1951&2413&3440 \cr
 & &16.27.7.11.17.19.23.71&1633&1848& & &128.3.5.19.43.73.127&12065&14016 \cr
\noalign{\hrule}
 & &27.343.121.29.67&4825&5122& & &9.5.7.11.47.89.151&41&664 \cr
9773&2177288883&4.25.11.13.67.193.197&11907&1292&9791&2188608345&16.3.11.41.83.151&2491&2492 \cr
 & &32.243.5.49.13.17.19&1989&1520& & &128.7.41.47.53.83.89&3403&3392 \cr
\noalign{\hrule}
 & &5.199.557.3929&3357&572& & &5.11.13.47.53.1229&6877&6642 \cr
9774&2177510735&8.9.11.13.199.373&1075&1114&9792&2188928885&4.81.169.529.41.53&2021&8950 \cr
 & &32.3.25.43.373.557&1865&2064& & &16.9.25.23.43.47.179&7697&8280 \cr
\noalign{\hrule}
}%
}
$$
\eject
\vglue -23 pt
\noindent\hskip 1 in\hbox to 6.5 in{\ 9793 -- 9828 \hfill\fbd 2190279175 -- 2206117987\frb}
\vskip -9 pt
$$
\vbox{
\nointerlineskip
\halign{\strut
    \vrule \ \ \hfil \frb #\ 
   &\vrule \hfil \ \ \fbb #\frb\ 
   &\vrule \hfil \ \ \frb #\ \hfil
   &\vrule \hfil \ \ \frb #\ 
   &\vrule \hfil \ \ \frb #\ \ \vrule \hskip 2 pt
   &\vrule \ \ \hfil \frb #\ 
   &\vrule \hfil \ \ \fbb #\frb\ 
   &\vrule \hfil \ \ \frb #\ \hfil
   &\vrule \hfil \ \ \frb #\ 
   &\vrule \hfil \ \ \frb #\ \vrule \cr%
\noalign{\hrule}
 & &25.49.1849.967&837&1012& & &27.25.13.19.67.197&3827&2842 \cr
9793&2190279175&8.27.7.11.23.31.967&1677&710&9811&2200603275&4.5.49.29.43.67.89&2767&16368 \cr
 & &32.81.5.13.23.43.71&5751&4784& & &128.3.7.11.31.2767&30437&13888 \cr
\noalign{\hrule}
 & &9.11.23.53.67.271&1829&610& & &3.25.13.109.139.149&8993&11718 \cr
9794&2191204917&4.5.11.31.59.61.67&1431&2168&9812&2201061525&4.81.7.13.17.529.31&1705&158 \cr
 & &64.27.5.31.53.271&93&160& & &16.5.11.23.961.79&22103&6952 \cr
\noalign{\hrule}
 & &5.13.17.19.29.59.61&1435&276& & &3.25.19.23.239.281&2233&2308 \cr
9795&2191269145&8.3.25.7.13.17.23.41&6177&1298&9813&2201136225&8.7.11.23.29.281.577&2561&10710 \cr
 & &32.9.11.29.59.71&99&1136& & &32.9.5.49.11.13.17.197&32487&34672 \cr
\noalign{\hrule}
 & &11.13.17.23.197.199&15903&12554& & &25.11.19.43.97.101&24901&24084 \cr
9796&2191957339&4.27.19.23.31.6277&21473&35020&9814&2201140975&8.27.5.11.37.223.673&29&194 \cr
 & &32.3.5.17.103.109.197&1545&1744& & &32.9.29.37.97.673&9657&10768 \cr
\noalign{\hrule}
 & &5.11.43.67.101.137&2177&2166& & &2187.5.7.149.193&1769&418 \cr
9797&2192541835&4.3.5.7.361.67.137.311&129&2474&9815&2201204565&4.5.11.19.29.61.149&17177&4428 \cr
 & &16.9.19.43.311.1237&23503&22392& & &32.27.41.89.193&89&656 \cr
\noalign{\hrule}
 & &9.5.31.37.107.397&2369&12320& & &27.37.43.107.479&63745&43148 \cr
9798&2192553585&64.3.25.7.11.23.103&1369&794&9816&2201675121&8.5.7.11.19.23.61.67&333&1070 \cr
 & &256.11.1369.397&407&128& & &32.9.25.7.19.37.107&175&304 \cr
\noalign{\hrule}
 & &9.7.13.19.29.43.113&3775&3388& & &121.59.373.827&247&126 \cr
9799&2192716071&8.25.49.121.113.151&4321&3078&9817&2202174469&4.9.7.13.19.59.827&5995&4756 \cr
 & &32.81.25.11.19.29.149&2475&2384& & &32.3.5.11.19.29.41.109&22345&26448 \cr
\noalign{\hrule}
 & &49.13.73.101.467&28841&22770& & &5.11.169.17.53.263&11711&2754 \cr
9800&2193312667&4.9.5.7.11.23.151.191&611&3782&9818&2202571085&4.81.49.289.239&1325&1276 \cr
 & &16.3.5.11.13.31.47.61&7755&15128& & &32.9.25.11.29.53.239&2151&2320 \cr
\noalign{\hrule}
 & &3.13.19.29.31.37.89&4961&2202& & &3.13.41.79.107.163&4245&8632 \cr
9801&2193661587&4.9.121.37.41.367&2015&1648&9819&2203164561&16.9.5.169.83.283&24287&23540 \cr
 & &128.5.11.13.31.41.103&4223&3520& & &128.25.11.107.149.163&1639&1600 \cr
\noalign{\hrule}
 & &9.25.49.23.41.211&2059&2794& & &3.5.11.13.167.6151&3059&3092 \cr
9802&2193677325&4.3.5.11.29.41.71.127&57239&38824&9820&2203380465&8.5.7.13.19.23.167.773&7&2178 \cr
 & &64.7.13.17.23.37.211&481&544& & &32.9.49.121.773&8503&2352 \cr
\noalign{\hrule}
 & &25.7.23.601.907&1559&2466& & &5.7.11.13.37.73.163&667&282 \cr
9803&2194055675&4.9.137.601.1559&95&506&9821&2203516315&4.3.23.29.37.47.163&507&344 \cr
 & &16.3.5.11.19.23.1559&4677&1672& & &64.9.169.29.43.47&12267&17888 \cr
\noalign{\hrule}
 & &9.7.121.29.9929&46435&42926& & &3.5.7.11.41.173.269&173&278 \cr
9804&2194974243&4.5.7.169.37.127.251&1119&638&9822&2203759635&4.139.29929.269&3731&33660 \cr
 & &16.3.5.11.13.29.127.373&4849&5080& & &32.9.5.7.11.13.17.41&221&48 \cr
\noalign{\hrule}
 & &83521.41.641&54901&28620& & &11.169.17.137.509&12555&13064 \cr
9805&2195015401&8.27.5.7.11.23.31.53&289&82&9823&2203771999&16.81.5.169.23.31.71&5069&170 \cr
 & &32.3.5.11.289.31.41&465&176& & &64.27.25.17.37.137&675&1184 \cr
\noalign{\hrule}
 & &7.37.71.307.389&18095&3702& & &3.19.29.59.97.233&807&3620 \cr
9806&2196069547&4.3.5.49.11.47.617&333&284&9824&2204207727&8.9.5.59.181.269&187&718 \cr
 & &32.27.5.11.37.47.71&1269&880& & &32.11.17.269.359&6103&47344 \cr
\noalign{\hrule}
 & &11.13.19.53.101.151&327&1334& & &9.5.7.121.17.41.83&1219&6448 \cr
9807&2196159251&4.3.13.23.29.101.109&453&1870&9825&2204990865&32.5.11.13.23.31.53&2199&5644 \cr
 & &16.9.5.11.17.29.151&153&1160& & &256.3.17.83.733&733&128 \cr
\noalign{\hrule}
 & &9.5.19.529.4861&33077&57382& & &5.7.11.13.37.43.277&2117&732 \cr
9808&2198605995&4.11.13.31.97.2207&1637&570&9826&2205738535&8.3.13.29.43.61.73&605&1164 \cr
 & &16.3.5.13.19.31.1637&1637&3224& & &64.9.5.121.73.97&7081&3168 \cr
\noalign{\hrule}
 & &3.5.7.11.13.1369.107&1525&2006& & &9.5.7.13.47.73.157&17663&25042 \cr
9809&2199442245&4.125.7.17.37.59.61&45243&8132&9827&2205841365&4.7.17.19.659.1039&611&1650 \cr
 & &32.9.11.19.107.457&457&912& & &16.3.25.11.13.47.659&659&440 \cr
\noalign{\hrule}
 & &25.11.289.19.31.47&4527&4432& & &19.137.229.3701&325&4026 \cr
9810&2200106425&32.9.5.11.47.277.503&1819&2336&9828&2206117987&4.3.25.11.13.61.137&591&916 \cr
 & &2048.3.17.73.107.503&110157&109568& & &32.9.61.197.229&549&3152 \cr
\noalign{\hrule}
}%
}
$$
\eject
\vglue -23 pt
\noindent\hskip 1 in\hbox to 6.5 in{\ 9829 -- 9864 \hfill\fbd 2206917735 -- 2227942255\frb}
\vskip -9 pt
$$
\vbox{
\nointerlineskip
\halign{\strut
    \vrule \ \ \hfil \frb #\ 
   &\vrule \hfil \ \ \fbb #\frb\ 
   &\vrule \hfil \ \ \frb #\ \hfil
   &\vrule \hfil \ \ \frb #\ 
   &\vrule \hfil \ \ \frb #\ \ \vrule \hskip 2 pt
   &\vrule \ \ \hfil \frb #\ 
   &\vrule \hfil \ \ \fbb #\frb\ 
   &\vrule \hfil \ \ \frb #\ \hfil
   &\vrule \hfil \ \ \frb #\ 
   &\vrule \hfil \ \ \frb #\ \vrule \cr%
\noalign{\hrule}
 & &3.5.11.19.23.127.241&63&52& & &5.17.19.23.841.71&20229&51256 \cr
9829&2206917735&8.27.7.13.19.127.241&575&4004&9847&2217965095&16.3.11.43.149.613&3393&3350 \cr
 & &64.25.49.11.169.23&845&1568& & &64.27.25.13.29.67.149&23517&23840 \cr
\noalign{\hrule}
 & &81.5.11.17.151.193&607&148& & &25.11.13.31.37.541&261&664 \cr
9830&2207145105&8.3.11.37.193.607&151&1972&9848&2218384025&16.9.11.29.83.541&111&430 \cr
 & &64.17.29.37.151&37&928& & &64.27.5.37.43.83&3569&864 \cr
\noalign{\hrule}
 & &5.121.19.29.37.179&3009&3614& & &9.49.11.23.59.337&29621&55640 \cr
9831&2207810165&4.3.13.17.19.29.59.139&895&31614&9849&2218405959&16.5.13.19.107.1559&177&70 \cr
 & &16.9.5.11.179.479&479&72& & &64.3.25.7.59.1559&1559&800 \cr
\noalign{\hrule}
 & &5.121.19.47.61.67&5265&2116& & &27.5.13.17.23.53.61&739&2494 \cr
9832&2208063055&8.81.25.13.19.529&341&134&9850&2218500765&4.17.23.29.43.739&6039&5300 \cr
 & &32.9.11.13.23.31.67&713&1872& & &32.9.25.11.43.53.61&215&176 \cr
\noalign{\hrule}
 & &9.25.11.17.31.1693&209&1484& & &243.49.227.821&23179&17050 \cr
9833&2208222225&8.3.7.121.19.31.53&365&3386&9851&2219071869&4.9.25.11.13.31.1783&1589&194 \cr
 & &32.5.7.73.1693&511&16& & &16.5.7.11.13.97.227&1261&440 \cr
\noalign{\hrule}
 & &49.13.17.59.3457&10593&48176& & &25.7.11.17.19.43.83&247&828 \cr
9834&2208715327&32.9.11.107.3011&1345&1666&9852&2219114975&8.9.11.13.17.361.23&4105&8798 \cr
 & &128.3.5.49.11.17.269&1345&2112& & &32.3.5.53.83.821&2463&848 \cr
\noalign{\hrule}
 & &5.11.13.19.29.71.79&4619&6678& & &9.7.13.83.103.317&2167&1850 \cr
9835&2209749685&4.9.5.7.19.31.53.149&2343&488&9853&2219517027&4.3.25.7.11.37.83.197&533&3604 \cr
 & &64.27.11.31.61.71&1891&864& & &32.25.11.13.17.41.53&9911&16400 \cr
\noalign{\hrule}
 & &3.11.73.379.2423&5719&1550& & &9.5.19.373.6961&14111&7150 \cr
9836&2212225653&4.25.7.19.31.43.73&97&3042&9854&2219967315&4.3.125.11.13.103.137&9517&29108 \cr
 & &16.9.5.7.169.97&507&27160& & &32.19.31.307.383&11873&4912 \cr
\noalign{\hrule}
 & &9.7.11.13.23.59.181&823&2180& & &3.5.121.59.89.233&289&410 \cr
9837&2212763553&8.3.5.109.181.823&1687&782&9855&2220621645&4.25.289.41.59.89&19&1494 \cr
 & &32.7.17.23.109.241&4097&1744& & &16.9.17.19.41.83&969&27224 \cr
\noalign{\hrule}
 & &9.5.13.23.29.53.107&1419&1042& & &7.23.1201.11489&6345&5144 \cr
9838&2212795845&4.27.5.11.43.53.521&3317&3838&9856&2221524529&16.27.5.7.23.47.643&6677&890 \cr
 & &16.11.19.31.43.101.107&14663&15352& & &64.3.25.11.89.607&54023&26400 \cr
\noalign{\hrule}
 & &9.5.49.863.1163&47311&5024& & &9.5.7.13.17.59.541&1265&724 \cr
9839&2213090145&64.121.17.23.157&1313&1470&9857&2222041185&8.25.7.11.23.59.181&1623&148 \cr
 & &256.3.5.49.13.17.101&1313&2176& & &64.3.37.181.541&181&1184 \cr
\noalign{\hrule}
 & &3.11.19.23.29.67.79&157&394& & &3.11.17.23.31.67.83&535&876 \cr
9840&2213580237&4.11.23.67.157.197&6989&23940&9858&2224361073&8.9.5.23.67.73.107&1099&442 \cr
 & &32.9.5.7.19.29.241&1205&336& & &32.5.7.13.17.107.157&14287&8560 \cr
\noalign{\hrule}
 & &5.11.41.397.2473&59869&76146& & &25.11.13.53.59.199&329&2916 \cr
9841&2213916155&4.3.343.19.23.37.137&397&534&9859&2224625975&8.729.5.7.47.53&767&502 \cr
 & &16.9.7.23.37.89.397&5957&6408& & &32.27.7.13.59.251&1757&432 \cr
\noalign{\hrule}
 & &25.7.13.17.19.23.131&333&242& & &11.169.289.41.101&1645&47196 \cr
9842&2214027725&4.9.121.17.19.37.131&1205&236&9860&2224756391&8.27.5.7.19.23.47&221&202 \cr
 & &32.3.5.11.37.59.241&24013&11568& & &32.3.5.7.13.17.23.101&115&336 \cr
\noalign{\hrule}
 & &7.11.13.73.157.193&4181&18270& & &27.13.89.229.311&1955&1022 \cr
9843&2214184973&4.9.5.49.29.37.113&8321&12584&9861&2224810341&4.9.5.7.17.23.73.89&3091&3406 \cr
 & &64.3.121.13.53.157&159&352& & &16.11.13.17.23.131.281&36811&34408 \cr
\noalign{\hrule}
 & &25.7.13.17.23.47.53&5269&6444& & &3.11.19.23.41.53.71&1675&668 \cr
9844&2215806775&8.9.7.11.23.179.479&159&320&9862&2224915143&8.25.23.41.67.167&639&304 \cr
 & &1024.27.5.11.53.179&4833&5632& & &256.9.5.19.71.167&835&384 \cr
\noalign{\hrule}
 & &27.5.121.13.53.197&4487&1926& & &9.25.7.11.23.37.151&14981&16276 \cr
9845&2217198555&4.243.5.7.107.641&583&1118&9863&2226279825&8.5.11.13.71.211.313&1147&3468 \cr
 & &16.11.13.43.53.641&641&344& & &64.3.289.31.37.313&8959&10016 \cr
\noalign{\hrule}
 & &27.5.7.11.13.61.269&323&1324& & &5.7.11.29.31.41.157&45&1144 \cr
9846&2217430215&8.5.17.19.269.331&693&962&9864&2227942255&16.9.25.121.13.31&4553&4522 \cr
 & &32.9.7.11.13.17.19.37&323&592& & &64.3.7.13.17.19.29.157&663&608 \cr
\noalign{\hrule}
}%
}
$$
\eject
\vglue -23 pt
\noindent\hskip 1 in\hbox to 6.5 in{\ 9865 -- 9900 \hfill\fbd 2229126075 -- 2254369575\frb}
\vskip -9 pt
$$
\vbox{
\nointerlineskip
\halign{\strut
    \vrule \ \ \hfil \frb #\ 
   &\vrule \hfil \ \ \fbb #\frb\ 
   &\vrule \hfil \ \ \frb #\ \hfil
   &\vrule \hfil \ \ \frb #\ 
   &\vrule \hfil \ \ \frb #\ \ \vrule \hskip 2 pt
   &\vrule \ \ \hfil \frb #\ 
   &\vrule \hfil \ \ \fbb #\frb\ 
   &\vrule \hfil \ \ \frb #\ \hfil
   &\vrule \hfil \ \ \frb #\ 
   &\vrule \hfil \ \ \frb #\ \vrule \cr%
\noalign{\hrule}
 & &81.25.11.19.23.229&233&442& & &3.5.7.17.59.61.349&4183&3146 \cr
9865&2229126075&4.3.13.17.23.229.233&2075&902&9883&2242051035&4.5.121.13.47.59.89&171&124 \cr
 & &16.25.11.41.83.233&3403&1864& & &32.9.121.13.19.31.89&65949&60016 \cr
\noalign{\hrule}
 & &11.169.43.79.353&45&124& & &27.25.11.169.1787&2537&6398 \cr
9866&2229203119&8.9.5.11.31.43.353&5803&5846&9884&2242372275&4.5.7.13.43.59.457&2863&3078 \cr
 & &32.3.5.7.31.37.79.829&40145&39792& & &16.81.49.19.59.409&24131&22344 \cr
\noalign{\hrule}
 & &27.5.109.197.769&2809&1036& & &9.11.43.709.743&21535&52022 \cr
9867&2229219495&8.3.7.37.2809.109&817&1144&9885&2242532259&4.5.19.1369.59.73&1469&2838 \cr
 & &128.7.11.13.19.43.53&43301&64064& & &16.3.5.11.13.19.43.113&1469&760 \cr
\noalign{\hrule}
 & &9.11.13.43.59.683&49&7562& & &7.11.4091.7121&2423&47424 \cr
9868&2230076277&4.3.49.13.19.199&4505&5246&9886&2243164847&128.3.13.19.2423&1205&1218 \cr
 & &16.5.17.43.53.61&305&7208& & &512.9.5.7.19.29.241&41211&37120 \cr
\noalign{\hrule}
 & &3.11.13.31.43.47.83&18901&20358& & &121.13.29.137.359&1783&47400 \cr
9869&2230814157&4.81.169.29.41.461&13529&160&9887&2243580911&16.3.25.79.1783&1089&694 \cr
 & &256.5.41.83.163&815&5248& & &64.27.5.121.347&347&4320 \cr
\noalign{\hrule}
 & &125.11.19.23.47.79&1079&186& & &3.5.343.61.7151&14037&6886 \cr
9870&2231048875&4.3.25.13.31.79.83&1551&524&9888&2244305595&4.9.5.11.313.4679&20273&21838 \cr
 & &32.9.11.31.47.131&279&2096& & &16.121.19.61.97.179&17363&18392 \cr
\noalign{\hrule}
 & &3.5.7.17.599.2087&601&1196& & &11.13.83.379.499&729&350 \cr
9871&2231451705&8.13.23.601.2087&2863&4950&9889&2244677149&4.729.25.7.11.499&4901&7574 \cr
 & &32.9.25.7.11.23.409&6135&4048& & &16.3.49.169.29.541&18473&12984 \cr
\noalign{\hrule}
 & &9.5.7.961.47.157&415&572& & &5.7.11.37.359.439&2507&2322 \cr
9872&2233733985&8.3.25.11.13.961.83&893&68&9890&2245026245&4.27.7.23.43.109.359&561&202 \cr
 & &64.13.17.19.47.83&1577&7072& & &16.81.11.17.23.43.101&31671&34744 \cr
\noalign{\hrule}
 & &3.11.19.23.37.53.79&15&422& & &125.11.13.23.43.127&2263&612 \cr
9873&2234086899&4.9.5.53.79.211&2071&2116&9891&2245153625&8.9.11.17.31.43.73&8825&4826 \cr
 & &32.19.529.109.211&2507&3376& & &32.3.25.19.127.353&353&912 \cr
\noalign{\hrule}
 & &3.11.113.359.1669&349&10& & &243.5.11.13.67.193&841&3514 \cr
9874&2234308659&4.5.11.349.1669&8307&10052&9892&2246696595&4.7.841.193.251&517&324 \cr
 & &32.9.7.13.71.359&213&1456& & &32.81.7.11.47.251&1757&752 \cr
\noalign{\hrule}
 & &9.25.7.13.29.53.71&1031&506& & &3.25.11.71.89.431&581&2356 \cr
9875&2234375325&4.3.11.13.23.71.1031&35&178&9893&2246878425&8.7.19.31.83.431&2277&1846 \cr
 & &16.5.7.23.89.1031&2047&8248& & &32.9.11.13.23.71.83&897&1328 \cr
\noalign{\hrule}
 & &9.25.121.13.71.89&703&98& & &49.53.461.1877&58203&33770 \cr
9876&2236452075&4.5.49.13.19.37.71&1541&1086&9894&2247176309&4.9.5.11.29.223.307&265&42 \cr
 & &16.3.7.19.23.67.181&29141&10184& & &16.27.25.7.11.29.53&675&2552 \cr
\noalign{\hrule}
 & &3.11.13.23.31.71.103&513&410& & &3.19.83.661.719&2035&2696 \cr
9877&2236878501&4.81.5.11.19.23.31.41&2333&11648&9895&2248450329&16.5.11.37.337.719&2213&1494 \cr
 & &1024.7.13.19.2333&16331&9728& & &64.9.5.37.83.2213&6639&5920 \cr
\noalign{\hrule}
 & &5.7.13.17.19.31.491&447&44& & &3.5.11.29.37.97.131&519&548 \cr
9878&2236954265&8.3.5.7.11.17.19.149&31&354&9896&2249710815&8.9.5.37.131.137.173&253&6148 \cr
 & &32.9.31.59.149&149&8496& & &64.11.23.29.53.137&3151&1696 \cr
\noalign{\hrule}
 & &9.7.107.457.727&270205&274318& & &27.61.379.3607&32875&36482 \cr
9879&2239623099&4.5.11.13.37.337.4157&2477&2814&9897&2251536291&4.9.125.17.29.37.263&4697&14428 \cr
 & &16.3.5.7.67.2477.4157&165959&166280& & &32.7.11.29.61.3607&203&176 \cr
\noalign{\hrule}
 & &81.5.23.37.67.97&119&82& & &3.7.289.19.59.331&451&1572 \cr
9880&2239912845&4.27.5.7.17.23.41.97&6809&11524&9898&2251908519&8.9.11.41.131.331&2195&1864 \cr
 & &32.11.17.43.67.619&8041&9904& & &128.5.131.233.439&57509&74560 \cr
\noalign{\hrule}
 & &3.49.47.227.1429&2257&2030& & &5.121.109.127.269&7267&5922 \cr
9881&2241162147&4.5.343.29.37.47.61&8597&7524&9899&2252879035&4.9.7.169.43.47.127&3685&9146 \cr
 & &32.9.5.11.19.61.8597&163343&161040& & &16.3.5.11.13.17.67.269&663&536 \cr
\noalign{\hrule}
 & &3.7.19.67.191.439&120913&122230& & &3.25.17.53.73.457&15213&9008 \cr
9882&2241535317&4.5.7.13.17.71.131.719&209&4824&9900&2254369575&32.9.5.11.461.563&15203&10132 \cr
 & &64.9.11.17.19.67.131&1441&1632& & &256.17.23.149.661&15203&19072 \cr
\noalign{\hrule}
}%
}
$$
\eject
\vglue -23 pt
\noindent\hskip 1 in\hbox to 6.5 in{\ 9901 -- 9936 \hfill\fbd 2255598015 -- 2276334081\frb}
\vskip -9 pt
$$
\vbox{
\nointerlineskip
\halign{\strut
    \vrule \ \ \hfil \frb #\ 
   &\vrule \hfil \ \ \fbb #\frb\ 
   &\vrule \hfil \ \ \frb #\ \hfil
   &\vrule \hfil \ \ \frb #\ 
   &\vrule \hfil \ \ \frb #\ \ \vrule \hskip 2 pt
   &\vrule \ \ \hfil \frb #\ 
   &\vrule \hfil \ \ \fbb #\frb\ 
   &\vrule \hfil \ \ \frb #\ \hfil
   &\vrule \hfil \ \ \frb #\ 
   &\vrule \hfil \ \ \frb #\ \vrule \cr%
\noalign{\hrule}
 & &3.5.11.19.421.1709&18551&13920& & &25.7.11.13.31.37.79&1021&5454 \cr
9901&2255598015&64.9.25.13.29.1427&1709&11134&9919&2267590325&4.27.79.101.1021&605&8584 \cr
 & &256.19.293.1709&293&128& & &64.3.5.121.29.37&29&1056 \cr
\noalign{\hrule}
 & &3.5.11.23.31.127.151&583&130& & &9.7.17.37.89.643&613&1256 \cr
9902&2256078165&4.25.121.13.53.127&1449&1576&9920&2267734329&16.3.17.37.157.613&79&550 \cr
 & &64.9.7.13.23.53.197&10441&8736& & &64.25.11.79.613&48427&8800 \cr
\noalign{\hrule}
 & &25.23.59.71.937&209&1566& & &37.107.683.839&2321&1638 \cr
9903&2256928475&4.27.11.19.29.937&155&782&9921&2268653483&4.9.7.11.13.211.839&2399&6830 \cr
 & &16.9.5.17.23.29.31&4437&248& & &16.3.5.13.683.2399&2399&1560 \cr
\noalign{\hrule}
 & &5.43.89.179.659&403&492& & &9.5.343.19.71.109&1177&538 \cr
9904&2257183735&8.3.13.31.41.43.659&5907&22430&9922&2269577835&4.11.19.107.109.269&3479&25842 \cr
 & &32.9.5.11.179.2243&2243&1584& & &16.3.49.59.71.73&73&472 \cr
\noalign{\hrule}
 & &9.5.7.19.233.1619&715&482& & &9.7.11.17.37.41.127&115&4 \cr
9905&2257703595&4.25.11.13.241.1619&2597&978&9923&2269715679&8.3.5.11.23.41.127&2567&1624 \cr
 & &16.3.49.53.163.241&8639&13496& & &128.5.7.17.29.151&4379&320 \cr
\noalign{\hrule}
 & &23.4489.131.167&3751&738& & &27.7.11.43.67.379&5267&4966 \cr
9906&2258734619&4.9.121.31.41.167&7705&7826&9924&2270058021&4.11.13.23.67.191.229&5985&21328 \cr
 & &16.3.5.7.13.23.41.43.67&4305&4472& & &128.9.5.7.19.23.31.43&2185&1984 \cr
\noalign{\hrule}
 & &27.5.17.19.103.503&407&2158& & &11.19.31.41.83.103&1035&7514 \cr
9907&2259131445&4.11.13.37.83.503&2897&2394&9925&2270947811&4.9.5.13.289.23.41&83&206 \cr
 & &16.9.7.19.83.2897&2897&4648& & &16.3.5.13.23.83.103&69&520 \cr
\noalign{\hrule}
 & &9.25.11.13.23.43.71&227&98& & &9.11.17.19.29.31.79&205&2086 \cr
9908&2259296325&4.3.49.11.23.71.227&10811&5590&9926&2271038517&4.5.7.17.31.41.149&627&644 \cr
 & &16.5.7.13.19.43.569&569&1064& & &32.3.5.49.11.19.23.149&3427&3920 \cr
\noalign{\hrule}
 & &9.31.587.13799&2657&2626& & &9.11.13.529.47.71&33257&58120 \cr
9909&2259903627&4.13.101.2657.13799&44485&223872&9927&2271906351&16.5.7.1453.4751&1649&3102 \cr
 & &1024.3.5.7.11.31.41.53&20405&20992& & &64.3.5.7.11.17.47.97&1649&1120 \cr
\noalign{\hrule}
 & &5.169.61.163.269&27249&54796& & &9.5.7.11.13.31.1627&1689&3316 \cr
9910&2260093615&8.3.7.19.31.103.293&143&150&9928&2271934665&8.27.31.563.829&833&4 \cr
 & &32.9.25.11.13.19.31.103&21527&22320& & &64.49.17.563&119&18016 \cr
\noalign{\hrule}
 & &361.31.37.53.103&2435&3024& & &25.11.19.29.53.283&157079&164154 \cr
9911&2260391753&32.27.5.7.19.37.487&3949&434&9929&2272723475&4.3.13.43.109.251.281&3975&712 \cr
 & &128.3.49.11.31.359&3949&9408& & &64.9.25.53.89.281&2529&2848 \cr
\noalign{\hrule}
 & &9.11.19.443.2713&637&3350& & &9.17.37.281.1429&245&4532 \cr
9912&2260696779&4.25.49.11.13.19.67&1599&326&9930&2273168889&8.3.5.49.11.37.103&493&802 \cr
 & &16.3.7.169.41.163&27547&2296& & &32.7.11.17.29.401&11629&1232 \cr
\noalign{\hrule}
 & &9.125.19.29.41.89&8533&6842& & &19.101.797.1487&82665&67522 \cr
9913&2261923875&4.3.7.11.23.29.53.311&133&800&9931&2274281741&4.9.5.49.11.13.53.167&505&664 \cr
 & &256.25.49.11.19.53&2597&1408& & &64.3.25.7.11.13.83.101&14525&13728 \cr
\noalign{\hrule}
 & &25.43.563.3739&3277&462& & &9.11.289.281.283&3945&832 \cr
9914&2262936275&4.3.5.7.11.29.43.113&3739&2496&9932&2275232553&128.27.5.13.17.263&2857&562 \cr
 & &512.9.7.13.3739&819&256& & &512.281.2857&2857&256 \cr
\noalign{\hrule}
 & &27.11.29.103.2551&115&2666& & &3.25.7.11.23.37.463&4747&5902 \cr
9915&2263091589&4.5.11.23.29.31.43&303&1292&9933&2275425075&4.5.13.37.47.101.227&4167&4232 \cr
 & &32.3.17.19.31.101&589&27472& & &64.9.529.47.101.463&3243&3232 \cr
\noalign{\hrule}
 & &3.5.7.13.17.23.31.137&89&66& & &49.121.43.79.113&5833&3726 \cr
9916&2266687605&4.9.7.11.13.17.89.137&355&1426&9934&2275911869&4.81.19.23.113.307&3425&2408 \cr
 & &16.5.11.23.31.71.89&979&568& & &64.9.25.7.23.43.137&5175&4384 \cr
\noalign{\hrule}
 & &9.11.13.19.23.29.139&48055&6494& & &3.5.11.13.23.29.37.43&5353&3318 \cr
9917&2267110989&4.5.7.17.191.1373&209&1164&9935&2276267565&4.9.7.43.53.79.101&17963&12610 \cr
 & &32.3.7.11.17.19.97&97&1904& & &16.5.7.11.13.23.71.97&679&568 \cr
\noalign{\hrule}
 & &9.11.13.19.23.29.139&615&914& & &9.7.529.167.409&36003&32300 \cr
9918&2267110989&4.27.5.19.29.41.457&851&338&9936&2276334081&8.27.25.11.17.19.1091&161&26 \cr
 & &16.5.169.23.37.457&2405&3656& & &32.5.7.13.19.23.1091&5455&3952 \cr
\noalign{\hrule}
}%
}
$$
\eject
\vglue -23 pt
\noindent\hskip 1 in\hbox to 6.5 in{\ 9937 -- 9972 \hfill\fbd 2277131997 -- 2300179747\frb}
\vskip -9 pt
$$
\vbox{
\nointerlineskip
\halign{\strut
    \vrule \ \ \hfil \frb #\ 
   &\vrule \hfil \ \ \fbb #\frb\ 
   &\vrule \hfil \ \ \frb #\ \hfil
   &\vrule \hfil \ \ \frb #\ 
   &\vrule \hfil \ \ \frb #\ \ \vrule \hskip 2 pt
   &\vrule \ \ \hfil \frb #\ 
   &\vrule \hfil \ \ \fbb #\frb\ 
   &\vrule \hfil \ \ \frb #\ \hfil
   &\vrule \hfil \ \ \frb #\ 
   &\vrule \hfil \ \ \frb #\ \vrule \cr%
\noalign{\hrule}
 & &3.7.17.23.29.73.131&381&1298& & &9.5.11.19.29.37.227&109&964 \cr
9937&2277131997&4.9.11.17.29.59.127&305&19126&9955&2290785255&8.11.109.227.241&1369&1128 \cr
 & &16.5.61.73.131&5&488& & &128.3.1369.47.109&5123&2368 \cr
\noalign{\hrule}
 & &3.11.67.137.7523&55145&27608& & &9.59.83.131.397&5525&5348 \cr
9938&2278769361&16.5.7.17.29.41.269&2881&1692&9956&2292104511&8.3.25.7.13.17.191.397&27343&50072 \cr
 & &128.9.5.7.43.47.67&4515&3008& & &128.5.11.37.569.739&231583&236480 \cr
\noalign{\hrule}
 & &5.49.13.37.83.233&4883&47988& & &11.19.83.349.379&3595&3606 \cr
9939&2279004455&8.9.19.31.43.257&407&364&9957&2294505037&4.3.5.83.349.601.719&250173&758 \cr
 & &64.3.7.11.13.19.31.37&1023&608& & &16.27.7.11.361.379&133&216 \cr
\noalign{\hrule}
 & &3.5.11.59.157.1493&899&8364& & &9.7.121.13.19.23.53&1333&1450 \cr
9940&2281893735&8.9.11.17.29.31.41&1477&206&9958&2295231939&4.25.7.19.29.31.43.53&82407&73678 \cr
 & &32.7.29.103.211&1477&47792& & &16.3.5.11.13.17.197.2113&16745&16904 \cr
\noalign{\hrule}
 & &13.17.23.41.47.233&11803&2250& & &9.7.11.17.23.43.197&1655&10126 \cr
9941&2282221253&4.9.125.11.17.29.37&1081&806&9959&2295327573&4.5.23.61.83.331&985&924 \cr
 & &16.3.5.13.23.29.31.47&465&232& & &32.3.25.7.11.197.331&331&400 \cr
\noalign{\hrule}
 & &3.5.49.37.127.661&8701&5396& & &3.11.19.31.41.43.67&13515&10634 \cr
9942&2282938665&8.343.11.19.71.113&793&450&9960&2295917877&4.9.5.11.13.17.53.409&3427&3526 \cr
 & &32.9.25.13.19.61.71&13845&18544& & &16.5.13.23.41.43.53.149&9685&9752 \cr
\noalign{\hrule}
 & &81.17.47.35281&17241&18040& & &3.5.7.13.2809.599&3819&374 \cr
9943&2283351039&16.243.5.7.11.41.821&17681&15980&9961&2296736715&4.9.11.17.19.53.67&77&400 \cr
 & &128.25.11.17.47.17681&17681&17600& & &128.25.7.121.67&8107&320 \cr
\noalign{\hrule}
 & &9.11.13.17.139.751&76105&89866& & &7.23.43.257.1291&27837&27676 \cr
9944&2283926931&4.5.343.31.131.491&221&270&9962&2296961401&8.27.11.17.37.257.1031&115&9394 \cr
 & &16.27.25.7.13.17.31.131&4061&4200& & &32.3.5.7.121.17.23.61&6171&4880 \cr
\noalign{\hrule}
 & &9.25.13.17.19.41.59&869&12374& & &5.11.19.71.83.373&1377&2726 \cr
9945&2285410725&4.3.5.11.23.79.269&2753&2698&9963&2297003005&4.81.5.17.29.47.83&1969&4376 \cr
 & &16.19.71.269.2753&19099&22024& & &64.3.11.17.179.547&9299&17184 \cr
\noalign{\hrule}
 & &3.11.13.23.29.61.131&2855&5868& & &9.7.11.17.41.67.71&47959&1528 \cr
9946&2286568713&8.27.5.29.163.571&36179&40906&9964&2297730897&16.191.199.241&23115&22916 \cr
 & &32.121.13.23.113.181&1991&1808& & &128.3.5.17.23.67.337&1685&1472 \cr
\noalign{\hrule}
 & &11.41.73.127.547&28365&28912& & &27.5.11.1229.1259&183889&190034 \cr
9947&2287127887&32.3.5.13.31.61.73.139&41&114&9965&2297756835&4.13.17.29.373.7309&6825&484 \cr
 & &128.9.13.19.41.61.139&34333&35136& & &32.3.25.7.121.169.29&4901&6160 \cr
\noalign{\hrule}
 & &9.5.7.103.107.659&1409&12430& & &3.5.13.19.37.41.409&185&62 \cr
9948&2287794285&4.3.25.11.113.1409&659&584&9966&2298778365&4.25.31.1369.409&10773&23452 \cr
 & &64.73.659.1409&1409&2336& & &32.81.7.11.13.19.41&297&112 \cr
\noalign{\hrule}
 & &7.11.23.293.4409&13485&9076& & &49.11.13.59.67.83&2403&2494 \cr
9949&2287843327&8.3.5.23.29.31.2269&1491&778&9967&2298989693&4.27.7.11.29.43.67.89&1909&4790 \cr
 & &32.9.5.7.29.71.389&27619&20880& & &16.9.5.23.83.89.479&18423&19160 \cr
\noalign{\hrule}
 & &13.31.37.1849.83&10483&13554& & &3.17.47.67.103.139&805&946 \cr
9950&2288346437&4.27.11.31.251.953&44225&41366&9968&2299295883&4.5.7.11.23.43.67.139&855&118 \cr
 & &16.9.25.13.29.37.43.61&1769&1800& & &16.9.25.19.23.43.59&33925&19608 \cr
\noalign{\hrule}
 & &3.11.17.23.31.59.97&839&742& & &9.5.121.37.101.113&153&254 \cr
9951&2289159939&4.7.11.23.53.59.839&387&970&9969&2299320045&4.81.5.11.17.113.127&8251&2036 \cr
 & &16.9.5.7.43.97.839&5873&5160& & &32.17.37.223.509&8653&3568 \cr
\noalign{\hrule}
 & &3.7.121.97.9293&4465&4828& & &25.11.29.31.71.131&833&762 \cr
9952&2290510761&8.5.7.17.19.47.71.97&29&3366&9970&2299439725&4.3.5.49.17.31.127.131&4521&12116 \cr
 & &32.9.11.289.19.29&8381&912& & &32.9.11.13.17.137.233&30277&33552 \cr
\noalign{\hrule}
 & &9.25.49.121.17.101&263&2762& & &5.29.41.59.79.83&117&88 \cr
9953&2290520925&4.3.101.263.1381&539&842&9971&2299900535&16.9.11.13.59.79.83&4429&2128 \cr
 & &16.49.11.263.421&421&2104& & &512.3.7.11.19.43.103&62909&79104 \cr
\noalign{\hrule}
 & &9.5.6859.41.181&3737&3122& & &17.29.67.83.839&1125&286 \cr
9954&2290528755&4.3.7.37.101.181.223&21307&1216&9972&2300179747&4.9.125.11.13.29.67&2701&2324 \cr
 & &512.7.11.13.19.149&11473&3328& & &32.3.5.7.11.37.73.83&8103&6160 \cr
\noalign{\hrule}
}%
}
$$
\eject
\vglue -23 pt
\noindent\hskip 1 in\hbox to 6.5 in{\ 9973 -- 10008 \hfill\fbd 2301654745 -- 2318141825\frb}
\vskip -9 pt
$$
\vbox{
\nointerlineskip
\halign{\strut
    \vrule \ \ \hfil \frb #\ 
   &\vrule \hfil \ \ \fbb #\frb\ 
   &\vrule \hfil \ \ \frb #\ \hfil
   &\vrule \hfil \ \ \frb #\ 
   &\vrule \hfil \ \ \frb #\ \ \vrule \hskip 2 pt
   &\vrule \ \ \hfil \frb #\ 
   &\vrule \hfil \ \ \fbb #\frb\ 
   &\vrule \hfil \ \ \frb #\ \hfil
   &\vrule \hfil \ \ \frb #\ 
   &\vrule \hfil \ \ \frb #\ \vrule \cr%
\noalign{\hrule}
 & &5.13.29.37.61.541&5913&5372& & &7.11.31.61.83.191&3725&1338 \cr
9973&2301654745&8.81.13.17.29.73.79&41&418&9991&2308307771&4.3.25.149.191.223&589&366 \cr
 & &32.3.11.19.41.73.79&61541&38544& & &16.9.5.19.31.61.149&1341&760 \cr
\noalign{\hrule}
 & &9.11.23.61.73.227&76187&90034& & &9.5.7.19.367.1051&193&858 \cr
9974&2301662187&4.7.47.59.109.1621&2405&4026&9992&2308516245&4.27.11.13.193.367&26275&31046 \cr
 & &16.3.5.7.11.13.37.47.61&2405&2632& & &16.25.361.43.1051&215&152 \cr
\noalign{\hrule}
 & &243.5.31.53.1153&49387&11722& & &3.5.7.1331.13.31.41&731&116 \cr
9975&2301670485&4.13.29.131.5861&2079&3782&9993&2309171865&8.11.13.17.29.31.43&595&738 \cr
 & &16.27.7.11.29.31.61&1769&616& & &32.9.5.7.289.29.41&867&464 \cr
\noalign{\hrule}
 & &3.25.7.11.13.23.31.43&289787&285788& & &9.5.7.11.37.67.269&377&118 \cr
9976&2301724425&8.37.197.1471.1931&1701&230&9994&2310638715&4.13.29.59.67.269&101&168 \cr
 & &32.243.5.7.23.37.197&2997&3152& & &64.3.7.13.29.59.101&5959&12064 \cr
\noalign{\hrule}
 & &5.19.23.53.59.337&1261&924& & &11.13.23.43.59.277&399&958 \cr
9977&2302550815&8.3.7.11.13.53.59.97&6651&928&9995&2311341461&4.3.7.11.19.277.479&177&100 \cr
 & &512.27.7.29.739&19953&51968& & &32.9.25.19.59.479&4311&7600 \cr
\noalign{\hrule}
 & &5.23.29.53.83.157&1037&2574& & &125.7.11.37.43.151&13731&14204 \cr
9978&2303294405&4.9.5.11.13.17.61.83&157&92&9996&2312319625&8.3.25.7.23.53.67.199&237&3788 \cr
 & &32.3.11.17.23.61.157&561&976& & &64.9.79.199.947&74813&57312 \cr
\noalign{\hrule}
 & &5.529.73.79.151&209&186& & &27.7.11.17.29.37.61&1817&1846 \cr
9979&2303310965&4.3.11.19.23.31.73.151&1501&3180&9997&2313305379&4.3.7.13.17.23.61.71.79&39853&590 \cr
 & &32.9.5.11.361.53.79&5247&5776& & &16.5.11.23.59.3623&18115&10856 \cr
\noalign{\hrule}
 & &9.5.11.19.487.503&277&2158& & &3.7.121.31.43.683&73475&87818 \cr
9980&2303858205&4.13.83.277.503&1549&2052&9998&2313425499&4.25.19.2311.2939&8503&6192 \cr
 & &32.27.19.83.1549&1549&3984& & &128.9.5.11.19.43.773&3865&3648 \cr
\noalign{\hrule}
 & &25.11.13.19.107.317&747&430& & &5.7.11.169.43.827&14193&14752 \cr
9981&2303948075&4.9.125.13.19.43.83&2317&3058&9999&2313776465&64.9.11.13.19.83.461&827&86 \cr
 & &16.3.7.11.83.139.331&27473&23352& & &256.3.43.461.827&461&384 \cr
\noalign{\hrule}
 & &13.289.31.47.421&6831&6220& & &27.5.7.13.19.23.431&1181&1276 \cr
9982&2304532529&8.27.5.11.289.23.311&4371&11524&10000&2313842895&8.11.23.29.431.1181&1539&11452 \cr
 & &64.81.31.43.47.67&2881&2592& & &64.81.7.19.29.409&1227&928 \cr
\noalign{\hrule}
 & &5.11.13.17.61.3109&4149&7258& & &9.7.13.19.23.29.223&55801&28270 \cr
9983&2305183595&4.9.5.13.19.191.461&1309&1556&10001&2314558701&4.5.11.41.257.1361&8671&1866 \cr
 & &32.3.7.11.17.389.461&8169&7376& & &16.3.11.13.23.29.311&311&88 \cr
\noalign{\hrule}
 & &9.17.41.43.83.103&23485&53456& & &729.23.31.61.73&2255&3658 \cr
9984&2305998711&32.5.7.11.13.61.257&83&174&10002&2314566981&4.9.5.11.961.41.59&1403&442 \cr
 & &128.3.5.11.29.61.83&1769&3520& & &16.11.13.17.23.59.61&767&1496 \cr
\noalign{\hrule}
 & &9.5.7.17.37.103.113&9577&11328& & &343.37.179.1019&513&506 \cr
9985&2306093265&128.27.7.59.61.157&1133&10396&10003&2314851091&4.27.49.11.19.23.37.179&1019&4420 \cr
 & &1024.11.23.103.113&253&512& & &32.9.5.11.13.17.23.1019&6435&6256 \cr
\noalign{\hrule}
 & &9.7.11.17.29.43.157&5&124& & &9.5.49.11.307.311&18881&5066 \cr
9986&2306472399&8.3.5.11.29.31.157&395&76&10004&2315794635&4.7.17.79.149.239&1755&2308 \cr
 & &64.25.19.31.79&775&48032& & &32.27.5.13.149.577&5811&9232 \cr
\noalign{\hrule}
 & &9.7.11.17.29.43.157&395&76& & &9.11.17.313.4397&1729&2668 \cr
9987&2306472399&8.3.5.7.17.19.43.79&5&124&10005&2316247263&8.3.7.11.13.17.19.23.29&295&1252 \cr
 & &64.25.19.31.79&775&48032& & &64.5.19.23.59.313&1357&3040 \cr
\noalign{\hrule}
 & &27.7.11.29.67.571&355&382& & &3.5.7.11.13.37.43.97&141&44 \cr
9988&2306552787&4.5.7.29.71.191.571&969&3028&10006&2317219905&8.9.7.121.13.43.47&3589&14600 \cr
 & &32.3.5.17.19.191.757&64345&58064& & &128.25.37.73.97&365&64 \cr
\noalign{\hrule}
 & &11.23.37.79.3119&127881&118520& & &11.17.61.137.1483&423&1906 \cr
9989&2306559761&16.9.5.13.1093.2963&3121&158&10007&2317571597&4.9.11.47.61.953&1483&530 \cr
 & &64.3.5.13.79.3121&3121&6240& & &16.3.5.47.53.1483&235&1272 \cr
\noalign{\hrule}
 & &11.19.29.349.1091&441&650& & &25.23.29.43.53.61&27&88 \cr
9990&2307780299&4.9.25.49.13.29.349&187&1234&10008&2318141825&16.27.5.11.29.43.53&413&3328 \cr
 & &16.3.25.11.13.17.617&15425&5304& & &8192.9.7.13.59&48321&4096 \cr
\noalign{\hrule}
}%
}
$$
\eject
\vglue -23 pt
\noindent\hskip 1 in\hbox to 6.5 in{\ 10009 -- 10044 \hfill\fbd 2318284563 -- 2341307605\frb}
\vskip -9 pt
$$
\vbox{
\nointerlineskip
\halign{\strut
    \vrule \ \ \hfil \frb #\ 
   &\vrule \hfil \ \ \fbb #\frb\ 
   &\vrule \hfil \ \ \frb #\ \hfil
   &\vrule \hfil \ \ \frb #\ 
   &\vrule \hfil \ \ \frb #\ \ \vrule \hskip 2 pt
   &\vrule \ \ \hfil \frb #\ 
   &\vrule \hfil \ \ \fbb #\frb\ 
   &\vrule \hfil \ \ \frb #\ \hfil
   &\vrule \hfil \ \ \frb #\ 
   &\vrule \hfil \ \ \frb #\ \vrule \cr%
\noalign{\hrule}
 & &3.7.19.23.29.31.281&187&94& & &3.5.11.17.43.101.191&8549&10742 \cr
10009&2318284563&4.7.11.17.19.23.29.47&117&320&10027&2326783965&4.5.11.41.83.103.131&82719&65704 \cr
 & &512.9.5.11.13.17.47&31161&14080& & &64.9.7.13.43.101.191&91&96 \cr
\noalign{\hrule}
 & &9.5.49.29.101.359&559&862& & &43.191.491.577&34485&59296 \cr
10010&2318581755&4.3.5.13.43.359.431&539&1616&10028&2326800391&64.3.5.121.17.19.109&577&468 \cr
 & &128.49.11.13.43.101&473&832& & &512.27.11.13.17.577&3861&4352 \cr
\noalign{\hrule}
 & &9.25.13.31.107.239&18573&24548& & &9.5.343.257.587&2405&2992 \cr
10011&2318831775&8.27.17.361.41.151&2387&1690&10029&2328508665&32.3.25.49.11.13.17.37&4891&2116 \cr
 & &32.5.7.11.169.361.31&2527&2288& & &256.17.529.67.73&83147&67712 \cr
\noalign{\hrule}
 & &81.5.49.13.89.101&5267&1298& & &25.121.53.73.199&3761&6786 \cr
10012&2319027165&4.11.23.59.89.229&2727&17654&10030&2329041275&4.9.13.29.73.3761&1295&5056 \cr
 & &16.27.7.13.97.101&97&8& & &512.3.5.7.13.37.79&8769&23296 \cr
\noalign{\hrule}
 & &3.5.13.41.263.1103&12561&1778& & &49.11.17.71.3581&1521&2060 \cr
10013&2319261555&4.9.5.7.53.79.127&14339&19316&10031&2329701913&8.9.5.169.17.71.103&1337&414 \cr
 & &32.11.13.439.1103&439&176& & &32.81.5.7.13.23.191&15471&23920 \cr
\noalign{\hrule}
 & &9.5.11.31.37.61.67&18965&17818& & &3.11.19.29.41.53.59&20687&20600 \cr
10014&2320455555&4.25.11.59.151.3793&93&1568&10032&2331187881&16.25.11.59.103.137.151&53&702 \cr
 & &256.3.49.31.3793&3793&6272& & &64.27.5.13.53.103.137&16029&16480 \cr
\noalign{\hrule}
 & &9.13.59.227.1481&1405&638& & &9.961.107.2521&1&962 \cr
10015&2320698861&4.5.11.29.281.1481&4071&12220&10033&2333041803&4.13.37.2521&16405&16368 \cr
 & &32.3.25.13.23.47.59&1081&400& & &128.3.5.11.17.31.193&2123&5440 \cr
\noalign{\hrule}
 & &3.23.29.53.79.277&36101&60200& & &11.17.23.269.2017&63683&41496 \cr
10016&2320757799&16.25.7.13.43.2777&8899&4986&10034&2333606473&16.3.7.13.19.43.1481&6051&4570 \cr
 & &64.9.5.11.277.809&2427&1760& & &64.9.5.7.457.2017&2285&2016 \cr
\noalign{\hrule}
 & &9.169.53.83.347&6683&6514& & &25.13.17.47.89.101&621&536 \cr
10017&2321735013&4.3.41.163.347.3257&18073&151610&10035&2334218575&16.27.5.23.47.67.101&2581&18326 \cr
 & &16.5.11.31.53.15161&15161&13640& & &64.3.49.11.17.29.89&957&1568 \cr
\noalign{\hrule}
 & &3.19.149.167.1637&325&176& & &3.7.19.41.61.2341&723&436 \cr
10018&2321807847&32.25.11.13.19.1637&1341&296&10036&2336081559&8.9.109.241.2341&2255&86 \cr
 & &512.9.5.13.37.149&2405&768& & &32.5.11.41.43.109&545&7568 \cr
\noalign{\hrule}
 & &9.5.11.97.137.353&2147&1736& & &3.7.11.17.23.41.631&787&156 \cr
10019&2322053415&16.3.5.7.19.31.97.113&2549&706&10037&2336694591&8.9.7.11.13.17.787&45385&46694 \cr
 & &64.113.353.2549&2549&3616& & &32.5.29.37.313.631&5365&5008 \cr
\noalign{\hrule}
 & &3.5.7.31.43.47.353&1649&1606& & &9.47.71.223.349&63&286 \cr
10020&2322159315&4.11.17.47.73.97.353&35991&1750&10038&2337378291&4.81.7.11.13.47.71&3173&6980 \cr
 & &16.27.125.7.17.31.43&153&200& & &32.5.7.19.167.349&665&2672 \cr
\noalign{\hrule}
 & &9.25.41.271.929&443&172& & &27.5.121.61.2347&1355&992 \cr
10021&2322476775&8.3.5.43.443.929&2501&286&10039&2338632945&64.9.25.31.61.271&11063&2662 \cr
 & &32.11.13.41.43.61&793&7568& & &256.1331.13.23.37&5291&2944 \cr
\noalign{\hrule}
 & &27.25.7.11.23.29.67&331&4016& & &81.7.11.29.67.193&5185&26 \cr
10022&2322710775&32.5.29.251.331&5427&4172&10040&2338868763&4.3.5.13.17.29.61&77&106 \cr
 & &256.81.7.67.149&447&128& & &16.5.7.11.13.17.53&85&5512 \cr
\noalign{\hrule}
 & &5.37.43.61.4787&1265&3522& & &81.11.23.29.31.127&223&96 \cr
10023&2322915685&4.3.25.11.23.43.587&831&244&10041&2339747289&64.243.23.31.223&235&478 \cr
 & &32.9.11.23.61.277&6371&1584& & &256.5.47.223.239&56165&28544 \cr
\noalign{\hrule}
 & &3.5.49.13.23.71.149&703&220& & &5.7.11.29.239.877&33291&58724 \cr
10024&2324893935&8.25.7.11.19.37.149&13&162&10042&2340217495&8.243.53.137.277&377&100 \cr
 & &32.81.11.13.19.37&7733&432& & &64.27.25.13.29.137&1755&4384 \cr
\noalign{\hrule}
 & &243.13.29.41.619&20075&12028& & &3.25.11.41.67.1033&12549&1186 \cr
10025&2324995569&8.9.25.11.31.73.97&1073&1352&10043&2341062075&4.9.5.47.89.593&1271&1694 \cr
 & &128.11.169.29.37.73&5291&4672& & &16.7.121.31.41.89&2387&712 \cr
\noalign{\hrule}
 & &27.13.37.79.2267&67&2200& & &5.7.13.29.191.929&3531&3154 \cr
10026&2325880791&16.25.11.13.37.67&13667&13602&10044&2341307605&4.3.11.19.83.107.929&669&9550 \cr
 & &64.3.5.79.173.2267&173&160& & &16.9.25.19.191.223&855&1784 \cr
\noalign{\hrule}
}%
}
$$
\eject
\vglue -23 pt
\noindent\hskip 1 in\hbox to 6.5 in{\ 10045 -- 10080 \hfill\fbd 2341380475 -- 2363437989\frb}
\vskip -9 pt
$$
\vbox{
\nointerlineskip
\halign{\strut
    \vrule \ \ \hfil \frb #\ 
   &\vrule \hfil \ \ \fbb #\frb\ 
   &\vrule \hfil \ \ \frb #\ \hfil
   &\vrule \hfil \ \ \frb #\ 
   &\vrule \hfil \ \ \frb #\ \ \vrule \hskip 2 pt
   &\vrule \ \ \hfil \frb #\ 
   &\vrule \hfil \ \ \fbb #\frb\ 
   &\vrule \hfil \ \ \frb #\ \hfil
   &\vrule \hfil \ \ \frb #\ 
   &\vrule \hfil \ \ \frb #\ \vrule \cr%
\noalign{\hrule}
 & &25.49.653.2927&939&286& & &9.5.7.11.13.19.41.67&5831&5626 \cr
10045&2341380475&4.3.11.13.313.2927&3185&258&10063&2351033685&4.2401.11.13.17.29.97&135&31078 \cr
 & &16.9.5.49.169.43&169&3096& & &16.27.5.17.41.379&1137&136 \cr
\noalign{\hrule}
 & &3.25.11.17.23.53.137&21667&22058& & &125.11.289.61.97&129&146 \cr
10046&2342217075&4.41.47.137.269.461&2295&8734&10064&2351267875&4.3.5.17.43.61.73.97&6061&144 \cr
 & &16.27.5.11.17.397.461&3573&3688& & &128.27.11.19.29.43&14877&2752 \cr
\noalign{\hrule}
 & &5.11.13.67.79.619&145&882& & &5.343.13.31.41.83&3615&7018 \cr
10047&2342602405&4.9.25.49.29.619&1541&316&10065&2351966435&4.3.25.121.13.29.241&23157&16168 \cr
 & &32.3.23.29.67.79&667&48& & &64.27.31.43.47.83&1161&1504 \cr
\noalign{\hrule}
 & &25.49.19.251.401&1881&124& & &3.125.11.13.17.29.89&2819&1306 \cr
10048&2342652025&8.9.5.7.11.361.31&559&524&10066&2352904125&4.13.29.653.2819&8855&27792 \cr
 & &64.3.11.13.31.43.131&61963&38688& & &128.9.5.7.11.23.193&4053&1472 \cr
\noalign{\hrule}
 & &11.17.23.47.67.173&63&110& & &9.17.19.37.79.277&217&494 \cr
10049&2343085877&4.9.5.7.121.17.23.67&20821&10034&10067&2353713597&4.7.13.17.361.31.37&1385&4752 \cr
 & &16.3.29.47.173.443&443&696& & &128.27.5.11.31.277&1023&320 \cr
\noalign{\hrule}
 & &3.5.11.43.101.3271&2345&926& & &9.5.121.17.59.431&1543&1112 \cr
10050&2343982245&4.25.7.67.101.463&9171&2404&10068&2353835385&16.121.17.139.1543&21525&4706 \cr
 & &32.9.7.601.1019&7133&28848& & &64.3.25.7.13.41.181&7421&14560 \cr
\noalign{\hrule}
 & &9.5.121.17.19.31.43&3389&5446& & &9.5.841.37.1681&159277&151708 \cr
10051&2344393755&4.3.7.43.389.3389&257&646&10069&2353845465&8.17.19.23.83.97.101&8613&1184 \cr
 & &16.17.19.257.3389&3389&2056& & &512.27.11.29.37.83&913&768 \cr
\noalign{\hrule}
 & &9.11.13.41.157.283&1945&4492& & &9.25.41.47.61.89&407&38 \cr
10052&2344490577&8.5.11.13.389.1123&18819&6466&10070&2353878675&4.5.11.19.37.47.61&697&462 \cr
 & &32.27.17.41.53.61&2703&976& & &16.3.7.121.17.37.41&4477&952 \cr
\noalign{\hrule}
 & &3.5.19.41.181.1109&23187&22282& & &5.17.23.59.137.149&51711&50354 \cr
10053&2345518365&4.9.13.19.59.131.857&6697&7876&10071&2354537485&4.3.11.289.1481.1567&885&596 \cr
 & &32.11.37.179.181.857&31709&31504& & &32.9.5.11.59.149.1567&1567&1584 \cr
\noalign{\hrule}
 & &9.25.13.17.43.1097&29371&26576& & &9.25.11.289.37.89&91&96 \cr
10054&2345577975&32.3.5.11.23.151.1277&757&2&10072&2355400575&64.27.5.7.13.17.37.89&27317&7648 \cr
 & &128.757.1277&1277&48448& & &4096.59.239.463&110657&120832 \cr
\noalign{\hrule}
 & &3.11.13.17.31.97.107&22903&23000& & &3.11.13.17.43.73.103&2325&2104 \cr
10055&2346515457&16.125.17.23.31.37.619&531&13706&10073&2357950881&16.9.25.11.31.73.263&6901&326 \cr
 & &64.9.5.7.11.37.59.89&26255&24864& & &64.31.67.103.163&2077&5216 \cr
\noalign{\hrule}
 & &5.7.11.29.43.67.73&6399&7202& & &9.5.11.37.199.647&14603&7486 \cr
10056&2348144645&4.81.5.13.43.79.277&29&14&10074&2358111195&4.3.5.17.19.197.859&481&2096 \cr
 & &16.27.7.13.29.79.277&7479&8216& & &128.13.37.131.197&2561&8384 \cr
\noalign{\hrule}
 & &9.5.13.23.37.53.89&551&668& & &9.11.17.37.43.881&4781&2138 \cr
10057&2348287695&8.5.19.29.37.89.167&10461&13754&10075&2359012293&4.3.7.43.683.1069&30155&15812 \cr
 & &32.3.11.13.19.529.317&4807&5072& & &32.5.37.59.67.163&9617&5360 \cr
\noalign{\hrule}
 & &27.7.13.23.89.467&5375&11446& & &11.361.41.43.337&87101&72300 \cr
10058&2348766693&4.125.23.43.59.97&5227&7458&10076&2359294201&8.3.25.7.23.241.541&123&118 \cr
 & &16.3.25.11.113.5227&31075&41816& & &32.9.5.7.23.41.59.541&61065&60592 \cr
\noalign{\hrule}
 & &9.5.13.89.197.229&893&92& & &81.25.11.17.23.271&47341&67834 \cr
10059&2348808345&8.13.19.23.47.229&2481&1870&10077&2360281275&4.7.13.2609.6763&12513&5750 \cr
 & &32.3.5.11.17.23.827&4301&13232& & &16.3.125.13.23.43.97&1261&1720 \cr
\noalign{\hrule}
 & &9.7.19.53.61.607&65&542& & &81.5.13.17.23.31.37&2189&10366 \cr
10060&2349029907&4.5.7.13.19.61.271&1507&3642&10078&2361231405&4.11.23.71.73.199&1887&2690 \cr
 & &16.3.11.13.137.607&137&1144& & &16.3.5.17.37.71.269&269&568 \cr
\noalign{\hrule}
 & &3.11.29.101.109.223&221&890& & &5.7.31.41.173.307&169&138 \cr
10061&2349441699&4.5.13.17.29.89.109&1115&738&10079&2362642835&4.3.5.7.169.23.41.173&407&1842 \cr
 & &16.9.25.41.89.223&1025&2136& & &16.9.11.13.23.37.307&3663&2392 \cr
\noalign{\hrule}
 & &25.49.13.41.59.61&583&642& & &9.11.41.59.71.139&57715&31106 \cr
10062&2349877075&4.3.11.13.41.53.61.107&18821&23208&10080&2363437989&4.5.7.17.97.103.151&11&108 \cr
 & &64.9.121.29.59.967&31581&30944& & &32.27.5.11.103.151&755&4944 \cr
\noalign{\hrule}
}%
}
$$
\eject
\vglue -23 pt
\noindent\hskip 1 in\hbox to 6.5 in{\ 10081 -- 10116 \hfill\fbd 2364048597 -- 2385442875\frb}
\vskip -9 pt
$$
\vbox{
\nointerlineskip
\halign{\strut
    \vrule \ \ \hfil \frb #\ 
   &\vrule \hfil \ \ \fbb #\frb\ 
   &\vrule \hfil \ \ \frb #\ \hfil
   &\vrule \hfil \ \ \frb #\ 
   &\vrule \hfil \ \ \frb #\ \ \vrule \hskip 2 pt
   &\vrule \ \ \hfil \frb #\ 
   &\vrule \hfil \ \ \fbb #\frb\ 
   &\vrule \hfil \ \ \frb #\ \hfil
   &\vrule \hfil \ \ \frb #\ 
   &\vrule \hfil \ \ \frb #\ \vrule \cr%
\noalign{\hrule}
 & &3.227.337.10301&5491&4810& & &9.49.71.191.397&661&1058 \cr
10081&2364048597&4.5.13.289.19.37.337&1705&7434&10099&2374219197&4.49.529.71.661&10505&36426 \cr
 & &16.9.25.7.11.17.31.59&33099&43400& & &16.3.5.11.13.191.467&2335&1144 \cr
\noalign{\hrule}
 & &9.19.29.281.1697&989&1540& & &27.25.7.13.29.31.43&14927&14368 \cr
10082&2364733863&8.5.7.11.23.43.1697&975&722&10100&2374506225&64.5.11.23.29.59.449&481&186 \cr
 & &32.3.125.7.13.361.43&11375&13072& & &256.3.11.13.31.37.449&4939&4736 \cr
\noalign{\hrule}
 & &9.11.13.89.107.193&8525&998& & &9.7.157.401.599&59653&34390 \cr
10083&2365427493&4.3.25.121.31.499&1157&658&10101&2375808309&4.5.121.17.19.29.181&349&2406 \cr
 & &16.5.7.13.31.47.89&1457&280& & &16.3.181.349.401&349&1448 \cr
\noalign{\hrule}
 & &5.19.23.29.107.349&1183&1278& & &27.19.29.197.811&1315&1118 \cr
10084&2366239195&4.9.7.169.29.71.349&4807&5314&10102&2376853659&4.9.5.13.19.29.43.263&3113&8110 \cr
 & &16.3.7.11.19.23.71.2657&18599&18744& & &16.25.11.13.283.811&3113&2600 \cr
\noalign{\hrule}
 & &3.2197.31.37.313&51425&16682& & &3.5.7.11.31.197.337&2543&1164 \cr
10085&2366241501&4.25.121.17.19.439&313&126&10103&2377058145&8.9.5.31.97.2543&2167&2198 \cr
 & &16.9.25.7.11.19.313&525&1672& & &32.7.11.157.197.2543&2543&2512 \cr
\noalign{\hrule}
 & &11.31.67.71.1459&9063&6986& & &7.13.17.47.53.617&10947&21436 \cr
10086&2366697883&4.9.7.19.53.71.499&10385&904&10104&2377657009&8.3.7.23.41.89.233&1009&2640 \cr
 & &64.3.5.7.31.67.113&339&1120& & &256.9.5.11.23.1009&55495&26496 \cr
\noalign{\hrule}
 & &9.7.41.61.83.181&3151&650& & &9.5.7.19.529.751&803&52 \cr
10087&2367068949&4.3.25.13.23.83.137&61&22&10105&2377714815&8.7.11.13.529.73&2253&1450 \cr
 & &16.25.11.23.61.137&6325&1096& & &32.3.25.13.29.751&65&464 \cr
\noalign{\hrule}
 & &9.5.17.73.109.389&2519&3686& & &9.25.7.11.31.43.103&2041&2216 \cr
10088&2367883845&4.3.11.19.97.109.229&17119&14600&10106&2378705175&16.13.31.103.157.277&3397&204 \cr
 & &64.25.17.361.53.73&1805&1696& & &128.3.17.43.79.157&2669&5056 \cr
\noalign{\hrule}
 & &25.49.11.31.53.107&6663&6322& & &5.13.37.727.1361&24815&25542 \cr
10089&2368918475&4.3.5.29.107.109.2221&843&1378&10107&2379620035&4.27.25.7.11.13.43.709&2201&74 \cr
 & &16.9.13.29.53.109.281&28449&29224& & &16.9.11.31.37.43.71&7029&10664 \cr
\noalign{\hrule}
 & &27.25.13.23.59.199&18619&21206& & &9.11.289.19.29.151&2821&2670 \cr
10090&2369627325&4.529.43.433.461&495&34&10108&2380463811&4.27.5.7.11.13.29.31.89&5077&1216 \cr
 & &16.9.5.11.17.43.433&8041&3464& & &512.5.19.89.5077&25385&22784 \cr
\noalign{\hrule}
 & &27.11.13.59.101.103&3055&3022& & &3.7.19.29.73.2819&583&2236 \cr
10091&2369792997&4.9.5.169.47.101.1511&149041&106318&10109&2381161377&8.7.11.13.43.53.73&2819&3330 \cr
 & &16.5.17.53.59.103.1447&7235&7208& & &32.9.5.37.53.2819&795&592 \cr
\noalign{\hrule}
 & &27.5.7.11.169.19.71&113&932& & &9.5.11.37.83.1567&77401&91504 \cr
10092&2369862495&8.3.13.71.113.233&8555&532&10110&2382067215&32.7.17.19.29.43.157&325&168 \cr
 & &64.5.7.19.29.59&29&1888& & &512.3.25.49.13.19.43&27391&24320 \cr
\noalign{\hrule}
 & &27.19.23.37.61.89&995&1052& & &3.7.11.19.31.83.211&2501&3978 \cr
10093&2370100527&8.9.5.37.61.199.263&74137&13442&10111&2382801267&4.27.13.17.41.61.83&36119&8120 \cr
 & &32.49.11.13.17.47.89&8789&10192& & &64.5.7.19.29.1901&1901&4640 \cr
\noalign{\hrule}
 & &9.7.19.23.43.2003&5&166& & &5.7.11.23.191.1409&2637&4408 \cr
10094&2371217499&4.5.43.83.2003&1023&980&10112&2383048745&16.9.19.29.191.293&1375&2254 \cr
 & &32.3.25.49.11.31.83&6391&12400& & &64.3.125.49.11.23.29&609&800 \cr
\noalign{\hrule}
 & &9.11.13.47.197.199&15113&12946& & &125.13.997.1471&52751&71874 \cr
10095&2371350267&4.3.7.13.17.127.6473&5713&760&10113&2383203875&4.27.1331.17.29.107&221&100 \cr
 & &64.5.7.17.19.29.197&2465&4256& & &32.9.25.11.13.289.29&3179&4176 \cr
\noalign{\hrule}
 & &11.19.31.389.941&5559&6500& & &3.11.19.29.53.2473&102037&113114 \cr
10096&2371631471&8.3.125.11.13.17.19.109&1231&186&10114&2383227627&4.13.23.47.167.2459&2315&144 \cr
 & &32.9.25.17.31.1231&11079&6800& & &128.9.5.23.47.463&31947&15040 \cr
\noalign{\hrule}
 & &9.5.7.13.127.4561&4387&8948& & &3.7.11.169.227.269&1349&4300 \cr
10097&2372016465&8.3.13.41.107.2237&319&1918&10115&2383838457&8.25.11.13.19.43.71&7263&3358 \cr
 & &32.7.11.29.107.137&3973&18832& & &32.27.5.23.73.269&657&1840 \cr
\noalign{\hrule}
 & &3.5.19.107.277.281&299&22& & &3.125.361.67.263&221&154 \cr
10098&2373639315&4.5.11.13.19.23.281&2061&656&10116&2385442875&4.7.11.13.17.361.263&4221&250 \cr
 & &128.9.23.41.229&28167&1472& & &16.9.125.49.13.67&49&312 \cr
\noalign{\hrule}
}%
}
$$
\eject
\vglue -23 pt
\noindent\hskip 1 in\hbox to 6.5 in{\ 10117 -- 10152 \hfill\fbd 2387194485 -- 2407242695\frb}
\vskip -9 pt
$$
\vbox{
\nointerlineskip
\halign{\strut
    \vrule \ \ \hfil \frb #\ 
   &\vrule \hfil \ \ \fbb #\frb\ 
   &\vrule \hfil \ \ \frb #\ \hfil
   &\vrule \hfil \ \ \frb #\ 
   &\vrule \hfil \ \ \frb #\ \ \vrule \hskip 2 pt
   &\vrule \ \ \hfil \frb #\ 
   &\vrule \hfil \ \ \fbb #\frb\ 
   &\vrule \hfil \ \ \frb #\ \hfil
   &\vrule \hfil \ \ \frb #\ 
   &\vrule \hfil \ \ \frb #\ \vrule \cr%
\noalign{\hrule}
 & &3.5.13.17.19.151.251&11649&9686& & &25.7.11.19.31.2113&7783&7008 \cr
10117&2387194485&4.9.11.19.29.167.353&59&40&10135&2395772225&64.3.11.19.43.73.181&2113&11100 \cr
 & &64.5.29.59.167.353&58951&54752& & &512.9.25.37.2113&333&256 \cr
\noalign{\hrule}
 & &49.121.17.19.29.43&893&354& & &27.5.13.307.4447&1763&2684 \cr
10118&2388088549&4.3.11.17.361.47.59&2967&14000&10136&2395976895&8.9.5.11.13.41.43.61&4447&1702 \cr
 & &128.9.125.7.23.43&2875&576& & &32.23.37.41.4447&943&592 \cr
\noalign{\hrule}
 & &3.17.61.739.1039&489&550& & &9.7.11.29.43.47.59&1193&1580 \cr
10119&2388691131&4.9.25.11.17.163.739&7015&5548&10137&2396346183&8.5.7.11.29.79.1193&1457&14580 \cr
 & &32.125.11.19.23.61.73&26125&26864& & &64.729.25.31.47&2511&800 \cr
\noalign{\hrule}
 & &11.13.89.337.557&191891&198018& & &9.5.7.13.31.79.239&72913&54032 \cr
10120&2388972443&4.27.7.19.79.193.347&5785&808&10138&2396848545&32.11.17.307.4289&465&4754 \cr
 & &64.3.5.13.89.101.193&2895&3232& & &128.3.5.11.31.2377&2377&704 \cr
\noalign{\hrule}
 & &9.121.361.59.103&49187&5506& & &3.25.7.19.37.73.89&10097&9394 \cr
10121&2389044933&4.101.487.2753&1133&1620&10139&2397880275&4.25.49.11.23.61.439&89&1314 \cr
 & &32.81.5.11.101.103&505&144& & &16.9.11.73.89.439&439&264 \cr
\noalign{\hrule}
 & &3.11.23.89.113.313&157&156& & &67.109.373.881&160655&167958 \cr
10122&2389211319&8.9.11.13.23.89.113.157&226195&66526&10140&2399860739&4.9.5.7.11.23.31.43.127&113&268 \cr
 & &32.5.19.29.31.37.2381&345245&348688& & &32.3.7.11.23.43.67.113&20769&19888 \cr
\noalign{\hrule}
 & &27.7.11.529.41.53&2405&1298& & &3.5.7.29.53.107.139&89&884 \cr
10123&2389845843&4.5.121.13.37.53.59&1767&194&10141&2400279105&8.13.17.29.89.107&135&242 \cr
 & &16.3.5.19.31.59.97&28615&4712& & &32.27.5.121.17.89&2057&12816 \cr
\noalign{\hrule}
 & &3.25.7.71.97.661&221&2204& & &3.5.11.13.17.41.1607&16591&10728 \cr
10124&2389961175&8.7.13.17.19.29.71&1089&970&10142&2402569455&16.27.5.47.149.353&9881&10234 \cr
 & &32.9.5.121.13.19.97&1573&912& & &64.7.17.41.43.47.241&10363&10528 \cr
\noalign{\hrule}
 & &5.67.131.157.347&5325&5194& & &5.7.19.23.31.37.137&3653&594 \cr
10125&2390810915&4.3.125.49.53.71.347&3223&5652&10143&2403441005&4.27.5.11.13.37.281&137&418 \cr
 & &32.27.7.11.53.157.293&22561&22896& & &16.9.121.13.19.137&117&968 \cr
\noalign{\hrule}
 & &3.7.11.67.191.809&3781&5118& & &27.5.11.37.67.653&37667&61828 \cr
10126&2391490563&4.9.19.67.199.853&10505&2828&10144&2403898695&8.7.13.29.41.5381&825&4556 \cr
 & &32.5.7.11.19.101.191&505&304& & &64.3.25.11.17.29.67&493&160 \cr
\noalign{\hrule}
 & &27.7.11.13.17.41.127&25655&31622& & &9.625.49.11.13.61&2599&3026 \cr
10127&2392403013&4.5.49.97.163.733&4209&544&10145&2404276875&4.7.11.13.17.23.89.113&305&696 \cr
 & &256.3.17.23.61.163&3749&7808& & &64.3.5.29.61.89.113&3277&2848 \cr
\noalign{\hrule}
 & &9.11.13.19.53.1847&2035&188& & &3.11.17.19.29.31.251&345&554 \cr
10128&2393728623&8.5.121.37.47.53&2337&4076&10146&2405192691&4.9.5.17.23.251.277&22387&16016 \cr
 & &64.3.5.19.41.1019&1019&6560& & &128.5.7.11.13.61.367&27755&23488 \cr
\noalign{\hrule}
 & &11.13.19.23.29.1321&3675&4996& & &9.5.7.121.17.47.79&565&1118 \cr
10129&2393967719&8.3.25.49.11.19.1249&107&1356&10147&2405856915&4.25.11.13.43.47.113&7247&29478 \cr
 & &64.9.25.7.107.113&12091&50400& & &16.3.4913.7247&7247&2312 \cr
\noalign{\hrule}
 & &3.5.7.11.13.17.83.113&589&1154& & &625.7.79.6961&10393&3432 \cr
10130&2394036645&4.11.13.17.19.31.577&927&1504&10148&2405895625&16.3.25.11.13.19.547&3793&3318 \cr
 & &256.9.19.31.47.103&27683&39552& & &64.9.7.11.79.3793&3793&3168 \cr
\noalign{\hrule}
 & &169.23.41.83.181&6039&6220& & &3.7.31.73.89.569&151&418 \cr
10131&2394170441&8.9.5.11.13.61.83.311&65231&29624&10149&2406612243&4.7.11.19.31.73.151&145&72 \cr
 & &128.3.7.529.37.41.43&6923&7104& & &64.9.5.11.19.29.151&31559&13920 \cr
\noalign{\hrule}
 & &11.13.17.19.139.373&2267&4074& & &3.5.7.121.169.19.59&193&16 \cr
10132&2394761083&4.3.7.11.19.97.2267&1865&402&10150&2406949545&32.5.7.11.169.193&4307&2448 \cr
 & &16.9.5.67.97.373&873&2680& & &1024.9.17.59.73&3723&512 \cr
\noalign{\hrule}
 & &9.5.53.61.101.163&6127&49322& & &9.7.83.349.1319&9867&19100 \cr
10133&2395119555&4.7.11.13.271.557&1769&1212&10151&2407070799&8.27.25.11.13.23.191&1319&3074 \cr
 & &32.3.7.13.29.61.101&377&112& & &32.5.11.29.53.1319&1595&848 \cr
\noalign{\hrule}
 & &7.11.13.19.29.43.101&1395&4324& & &5.11.13.61.97.569&6651&7220 \cr
10134&2395385993&8.9.5.11.13.23.31.47&595&3838&10152&2407242695&8.9.25.361.61.739&943&582 \cr
 & &32.3.25.7.17.19.101&425&48& & &32.27.23.41.97.739&16997&17712 \cr
\noalign{\hrule}
}%
}
$$
\eject
\vglue -23 pt
\noindent\hskip 1 in\hbox to 6.5 in{\ 10153 -- 10188 \hfill\fbd 2407814175 -- 2436215125\frb}
\vskip -9 pt
$$
\vbox{
\nointerlineskip
\halign{\strut
    \vrule \ \ \hfil \frb #\ 
   &\vrule \hfil \ \ \fbb #\frb\ 
   &\vrule \hfil \ \ \frb #\ \hfil
   &\vrule \hfil \ \ \frb #\ 
   &\vrule \hfil \ \ \frb #\ \ \vrule \hskip 2 pt
   &\vrule \ \ \hfil \frb #\ 
   &\vrule \hfil \ \ \fbb #\frb\ 
   &\vrule \hfil \ \ \frb #\ \hfil
   &\vrule \hfil \ \ \frb #\ 
   &\vrule \hfil \ \ \frb #\ \vrule \cr%
\noalign{\hrule}
 & &3.25.13.29.31.41.67&803&468& & &81.125.11.31.701&8803&8722 \cr
10153&2407814175&8.27.5.11.169.29.73&31&814&10171&2420290125&4.5.49.11.31.89.8803&2727&32 \cr
 & &32.121.31.37.73&8833&592& & &256.27.101.8803&8803&12928 \cr
\noalign{\hrule}
 & &243.7.11.41.43.73&10925&7786& & &9.13.163.179.709&3533&2824 \cr
10154&2408086989&4.25.17.19.23.41.229&9933&9542&10172&2420319681&16.3.179.353.3533&2035&1498 \cr
 & &16.3.7.11.13.43.229.367&2977&2936& & &64.5.7.11.37.107.353&135905&126688 \cr
\noalign{\hrule}
 & &5.37.43.79.3833&72427&53262& & &5.11.23.29.149.443&63&86 \cr
10155&2408829685&4.9.11.23.47.67.269&111&158&10173&2421466795&4.9.5.7.11.29.43.443&14611&894 \cr
 & &16.27.11.23.37.67.79&1809&2024& & &16.27.19.149.769&769&4104 \cr
\noalign{\hrule}
 & &49.17.31.269.347&46255&47088& & &9.11.13.17.37.41.73&7885&218 \cr
10156&2410396289&32.27.5.11.841.31.109&46687&44982&10174&2422902339&4.3.5.13.19.83.109&1817&6068 \cr
 & &128.729.49.17.46687&46687&46656& & &32.23.37.41.79&23&1264 \cr
\noalign{\hrule}
 & &9.25.11.17.23.47.53&13813&13588& & &7.37.47.107.1861&16827&52030 \cr
10157&2410602975&8.17.19.23.43.79.727&31075&186&10175&2423972971&4.3.5.121.43.71.79&1961&1092 \cr
 & &32.3.25.11.19.31.113&589&1808& & &32.9.5.7.11.13.37.53&3445&1584 \cr
\noalign{\hrule}
 & &9.11.17.23.167.373&7409&49660& & &3.5.7.31.47.83.191&12221&15086 \cr
10158&2411222319&8.5.13.31.191.239&1361&1122&10176&2425271205&4.121.19.31.101.397&3041&3438 \cr
 & &32.3.5.11.17.31.1361&1361&2480& & &16.9.11.101.191.3041&9123&8888 \cr
\noalign{\hrule}
 & &1849.241.5413&3631&1782& & &5.43.79.373.383&24969&55226 \cr
10159&2412081517&4.81.11.241.3631&731&2900&10177&2426460115&4.3.7.29.41.53.521&737&2910 \cr
 & &32.9.25.11.17.29.43&2475&7888& & &16.9.5.11.29.67.97&6499&22968 \cr
\noalign{\hrule}
 & &625.7.13.151.281&3079&5046& & &5.13.19.41.191.251&5159&4968 \cr
10160&2413263125&4.3.841.151.3079&3729&650&10178&2427492535&16.27.5.7.11.23.67.251&779&26 \cr
 & &16.9.25.11.13.29.113&2871&904& & &64.9.11.13.19.41.67&603&352 \cr
\noalign{\hrule}
 & &25.47.53.83.467&36387&61138& & &9.25.7.11.61.2297&60671&35404 \cr
10161&2413841275&4.9.7.11.13.311.397&1717&2650&10179&2427527025&8.169.53.167.359&34455&25498 \cr
 & &16.3.25.7.13.17.53.101&1717&2184& & &32.3.5.11.19.61.2297&19&16 \cr
\noalign{\hrule}
 & &7.11.13.19.37.47.73&35&108& & &9.11.13.883.2137&1085&202 \cr
10162&2414404993&8.27.5.49.19.37.47&4121&6424&10180&2428531677&4.5.7.31.101.2137&1053&1084 \cr
 & &128.9.11.13.73.317&317&576& & &32.81.5.7.13.101.271&12195&11312 \cr
\noalign{\hrule}
 & &9.11.361.257.263&211&3038& & &5.2197.37.43.139&4191&1786 \cr
10163&2415634749&4.49.31.211.263&815&1026&10181&2429321765&4.3.11.169.19.47.127&483&1376 \cr
 & &16.27.5.7.19.31.163&1085&3912& & &256.9.7.23.43.127&2921&8064 \cr
\noalign{\hrule}
 & &81.19.47.127.263&4697&5590& & &3.13.17.443.8273&4247&4026 \cr
10164&2415994533&4.5.7.11.13.43.61.263&135&2758&10182&2429854557&4.9.11.31.61.137.443&4117&130 \cr
 & &16.27.25.49.13.197&2561&9800& & &16.5.11.13.23.61.179&15433&7160 \cr
\noalign{\hrule}
 & &3.7.11.179.211.277&897&2150& & &243.13.29.101.263&6103&2684 \cr
10165&2416725003&4.9.25.13.23.43.211&3043&1988&10183&2433462993&8.81.11.17.61.359&11395&10504 \cr
 & &32.5.7.17.23.71.179&1207&1840& & &128.5.13.17.43.53.101&3655&3392 \cr
\noalign{\hrule}
 & &121.19.59.71.251&315&334& & &3.25.49.11.29.31.67&1059&884 \cr
10166&2417258261&4.9.5.7.11.71.167.251&767&14&10184&2434919025&8.9.7.11.13.17.31.353&145&548 \cr
 & &16.3.5.49.13.59.167&1911&6680& & &64.5.17.29.137.353&6001&4384 \cr
\noalign{\hrule}
 & &9.13.29.31.83.277&823&76& & &9.5.17.19.29.53.109&3091&70 \cr
10167&2418262353&8.13.19.277.823&145&132&10185&2435092155&4.3.25.7.11.17.281&937&988 \cr
 & &64.3.5.11.19.29.823&4115&6688& & &32.13.19.281.937&12181&4496 \cr
\noalign{\hrule}
 & &5.49.29.31.79.139&1121&1170& & &27.7.19.263.2579&3085&506 \cr
10168&2418620155&4.9.25.13.19.31.59.139&4367&58&10186&2435692707&4.5.11.23.263.617&4051&2736 \cr
 & &16.3.11.13.19.29.397&2717&9528& & &128.9.19.23.4051&4051&1472 \cr
\noalign{\hrule}
 & &43.73.109.7069&148005&155962& & &3.5.121.19.23.37.83&613&2522 \cr
10169&2418665419&4.9.5.11.13.23.29.2689&3589&6278&10187&2435779005&4.11.13.37.97.613&1577&5166 \cr
 & &16.3.5.29.37.43.73.97&3589&3480& & &16.9.7.13.19.41.83&123&728 \cr
\noalign{\hrule}
 & &9.5.11.17.337.853&78913&80598& & &125.151.337.383&49979&7854 \cr
10170&2418984315&4.27.7.19.23.47.73.101&853&416&10188&2436215125&4.3.7.11.17.23.41.53&337&360 \cr
 & &256.7.13.73.101.853&9191&9344& & &64.27.5.7.11.53.337&2079&1696 \cr
\noalign{\hrule}
}%
}
$$
\eject
\vglue -23 pt
\noindent\hskip 1 in\hbox to 6.5 in{\ 10189 -- 10224 \hfill\fbd 2437300845 -- 2457162631\frb}
\vskip -9 pt
$$
\vbox{
\nointerlineskip
\halign{\strut
    \vrule \ \ \hfil \frb #\ 
   &\vrule \hfil \ \ \fbb #\frb\ 
   &\vrule \hfil \ \ \frb #\ \hfil
   &\vrule \hfil \ \ \frb #\ 
   &\vrule \hfil \ \ \frb #\ \ \vrule \hskip 2 pt
   &\vrule \ \ \hfil \frb #\ 
   &\vrule \hfil \ \ \fbb #\frb\ 
   &\vrule \hfil \ \ \frb #\ \hfil
   &\vrule \hfil \ \ \frb #\ 
   &\vrule \hfil \ \ \frb #\ \vrule \cr%
\noalign{\hrule}
 & &9.5.7.43.103.1747&709&194& & &3.13.961.79.827&33&994 \cr
10189&2437300845&4.3.97.709.1747&2975&2266&10207&2448615507&4.9.7.11.71.827&6409&620 \cr
 & &16.25.7.11.17.97.103&1067&680& & &32.5.13.17.29.31&2465&16 \cr
\noalign{\hrule}
 & &27.41.47.79.593&9625&6386& & &25.23.101.181.233&40271&63804 \cr
10190&2437402563&4.125.7.11.31.47.103&227&948&10208&2449196975&8.3.7.11.13.409.523&57&466 \cr
 & &32.3.5.11.31.79.227&1705&3632& & &32.9.7.11.13.19.233&819&3344 \cr
\noalign{\hrule}
 & &25.7.13.17.19.31.107&18909&26566& & &81.5.49.19.67.97&29&902 \cr
10191&2437414525&4.9.7.11.37.191.359&325&248&10209&2450480445&4.9.5.11.29.41.67&8303&5432 \cr
 & &64.3.25.13.31.37.359&1077&1184& & &64.7.361.23.97&437&32 \cr
\noalign{\hrule}
 & &9.5.11.17.29.97.103&57&46& & &5.11.31.43.67.499&13&486 \cr
10192&2438153685&4.27.5.17.19.23.29.97&1133&2782&10210&2451140395&4.243.5.13.31.67&5017&5368 \cr
 & &16.11.13.19.23.103.107&2461&1976& & &64.9.11.29.61.173&15921&5536 \cr
\noalign{\hrule}
 & &27.49.37.109.457&18755&1846& & &9.11.211.239.491&1069&830 \cr
10193&2438396163&4.5.7.121.13.31.71&1371&1016&10211&2451303261&4.5.11.83.491.1069&211&702 \cr
 & &64.3.11.13.127.457&1397&416& & &16.27.5.13.211.1069&1069&1560 \cr
\noalign{\hrule}
 & &9.5.361.29.31.167&16369&17204& & &3.7.11.139.241.317&2587&900 \cr
10194&2438910585&8.3.11.17.23.29.16369&2821&19190&10212&2453031273&8.27.25.13.139.199&583&3170 \cr
 & &32.5.7.13.17.19.31.101&1717&1456& & &32.125.11.53.317&125&848 \cr
\noalign{\hrule}
 & &27.17.59.113.797&173&286& & &169.841.41.421&79695&62434 \cr
10195&2438941941&4.11.13.59.173.797&15&782&10213&2453288669&4.9.5.7.11.19.23.31.53&377&842 \cr
 & &16.3.5.11.17.23.173&253&6920& & &16.3.7.11.13.19.29.421&209&168 \cr
\noalign{\hrule}
 & &27.5.7.11.19.53.233&2369&2666& & &3.5.19.29.37.71.113&50507&70742 \cr
10196&2438989245&4.7.23.31.43.103.233&1653&71786&10214&2453473515&4.7.17.31.163.2971&2871&100 \cr
 & &16.3.11.13.19.29.251&377&2008& & &32.9.25.7.11.29.31&2387&240 \cr
\noalign{\hrule}
 & &9.5.49.121.41.223&1703&142& & &9.5.7.13.359.1669&50209&8206 \cr
10197&2439398115&4.7.121.13.71.131&6021&5900&10215&2453605245&4.11.23.37.59.373&1373&2730 \cr
 & &32.27.25.59.71.223&1065&944& & &16.3.5.7.13.37.1373&1373&296 \cr
\noalign{\hrule}
 & &81.5.121.17.29.101&5759&4042& & &9.25.49.457.487&2669&1444 \cr
10198&2440105965&4.5.13.29.43.47.443&18981&45254&10216&2453712975&8.17.361.157.487&165&322 \cr
 & &16.27.1331.17.19.37&209&296& & &32.3.5.7.11.17.361.23&3971&6256 \cr
\noalign{\hrule}
 & &3.31.59.241.1847&803&1044& & &5.23.31.47.97.151&3969&11066 \cr
10199&2442411849&8.27.11.29.31.59.73&335&1928&10217&2454178085&4.81.49.11.23.503&2425&1922 \cr
 & &128.5.11.29.67.241&3685&1856& & &16.3.25.7.11.961.97&465&616 \cr
\noalign{\hrule}
 & &9.25.7.17.127.719&559&584& & &9.7.121.13.17.31.47&19475&53234 \cr
10200&2444905575&16.7.13.17.43.73.719&9017&330&10218&2454583131&4.25.19.41.43.619&16157&27918 \cr
 & &64.3.5.11.43.71.127&781&1376& & &16.27.11.47.107.151&453&856 \cr
\noalign{\hrule}
 & &3.5.11.13.17.19.3529&8927&8718& & &9.5.7.11.37.107.179&2923&3818 \cr
10201&2445014715&4.9.13.17.79.113.1453&3529&160660&10219&2455510365&4.11.23.1369.79.83&1141&228 \cr
 & &32.5.29.277.3529&277&464& & &32.3.7.19.23.79.163&12877&6992 \cr
\noalign{\hrule}
 & &3.11.13.41.43.53.61&13785&11284& & &3.7.13.17.29.71.257&495&2 \cr
10202&2445205191&8.9.5.7.169.31.919&1537&3058&10220&2455845483&4.27.5.11.13.257&2059&1802 \cr
 & &32.7.11.29.31.53.139&4031&3472& & &16.5.17.29.53.71&265&8 \cr
\noalign{\hrule}
 & &3.5.11.13.23.101.491&14363&15836& & &3.11.29.47.97.563&273&244 \cr
10203&2446571985&8.5.11.37.53.107.271&3437&522&10221&2456348169&8.9.7.13.61.97.563&6715&604 \cr
 & &32.9.7.29.271.491&1897&1392& & &64.5.17.61.79.151&59645&33184 \cr
\noalign{\hrule}
 & &81.7.83.149.349&57419&29150& & &81.7.11.289.29.47&1625&1336 \cr
10204&2447219061&4.25.11.53.67.857&261&596&10222&2456797959&16.9.125.11.13.29.167&47&272 \cr
 & &32.9.5.11.29.53.149&1595&848& & &512.5.13.17.47.167&2171&1280 \cr
\noalign{\hrule}
 & &5.7.13.289.37.503&23463&30002& & &83.881.33599&16841&16758 \cr
10205&2447253445&4.27.49.11.79.2143&851&1292&10223&2456859677&4.9.49.11.19.881.1531&17671&830 \cr
 & &32.3.11.17.19.23.37.79&4503&4048& & &16.3.5.7.19.41.83.431&16359&17240 \cr
\noalign{\hrule}
 & &9.5.49.11.19.47.113&2813&2498& & &7.11.19.43.139.281&549&268 \cr
10206&2447547795&4.7.11.19.29.97.1249&107&1356&10224&2457162631&8.9.7.11.61.67.139&1405&124 \cr
 & &32.3.29.97.107.113&3103&1552& & &64.3.5.31.67.281&335&2976 \cr
\noalign{\hrule}
}%
}
$$
\eject
\vglue -23 pt
\noindent\hskip 1 in\hbox to 6.5 in{\ 10225 -- 10260 \hfill\fbd 2459339883 -- 2484893619\frb}
\vskip -9 pt
$$
\vbox{
\nointerlineskip
\halign{\strut
    \vrule \ \ \hfil \frb #\ 
   &\vrule \hfil \ \ \fbb #\frb\ 
   &\vrule \hfil \ \ \frb #\ \hfil
   &\vrule \hfil \ \ \frb #\ 
   &\vrule \hfil \ \ \frb #\ \ \vrule \hskip 2 pt
   &\vrule \ \ \hfil \frb #\ 
   &\vrule \hfil \ \ \fbb #\frb\ 
   &\vrule \hfil \ \ \frb #\ \hfil
   &\vrule \hfil \ \ \frb #\ 
   &\vrule \hfil \ \ \frb #\ \vrule \cr%
\noalign{\hrule}
 & &9.7.121.169.23.83&445&302& & &27.5.7.121.17.31.41&7169&6812 \cr
10225&2459339883&4.5.7.11.13.23.89.151&2237&1236&10243&2470652415&8.9.5.11.13.67.107.131&317&7522 \cr
 & &32.3.5.89.103.2237&45835&35792& & &32.107.317.3761&33919&60176 \cr
\noalign{\hrule}
 & &3.5.17.71.199.683&7&206& & &3.13.23.73.97.389&645&1034 \cr
10226&2460777285&4.5.7.17.103.683&639&44&10244&2470794573&4.9.5.11.13.43.47.97&9725&8464 \cr
 & &32.9.11.71.103&3399&16& & &128.125.11.529.389&1375&1472 \cr
\noalign{\hrule}
 & &49.121.13.19.1681&15183&16756& & &9.25.13.19.173.257&1147&1166 \cr
10227&2461762303&8.9.343.59.71.241&533&190&10245&2470920075&4.25.11.13.31.37.53.173&26707&1968 \cr
 & &32.3.5.13.19.41.59.71&1065&944& & &128.3.17.41.53.1571&64411&57664 \cr
\noalign{\hrule}
 & &13.17.31.43.61.137&375&418& & &9.25.11.13.43.1787&5731&3944 \cr
10228&2461913701&4.3.125.11.17.19.31.137&6161&1914&10246&2472359175&16.121.13.17.29.521&20925&5816 \cr
 & &16.9.5.121.29.61.101&12221&10440& & &256.27.25.31.727&2181&3968 \cr
\noalign{\hrule}
 & &9.5.19.43.167.401&15457&22638& & &25.11.23.29.97.139&7593&3562 \cr
10229&2462041755&4.27.343.11.13.29.41&401&50&10247&2473119275&4.3.5.11.13.137.2531&337&348 \cr
 & &16.25.343.29.401&343&1160& & &32.9.13.29.337.2531&39429&40496 \cr
\noalign{\hrule}
 & &81.5.11.19.47.619&1189&296& & &3.5.7.23.29.89.397&264983&265012 \cr
10230&2462576985&16.3.29.37.41.619&371&248&10248&2474546655&8.7.11.19.529.41.281.317&21663&26 \cr
 & &256.7.29.31.37.53&56869&27776& & &32.9.13.19.29.83.317&18069&17264 \cr
\noalign{\hrule}
 & &27.5.121.23.79.83&191131&184574& & &27.5.11.13.41.53.59&475&1124 \cr
10231&2463497685&4.13.17.31.229.11243&4133&7110&10249&2475036135&8.9.125.19.53.281&59&1066 \cr
 & &16.9.5.17.31.79.4133&4133&4216& & &32.13.41.59.281&281&16 \cr
\noalign{\hrule}
 & &7.23.179.193.443&57159&28340& & &3.125.4489.1471&98671&85204 \cr
10232&2463995681&8.27.5.13.29.73.109&1617&1544&10250&2476244625&8.7.17.79.179.1249&7695&6446 \cr
 & &128.81.5.49.11.13.193&5005&5184& & &32.81.5.7.11.17.19.293&55377&56848 \cr
\noalign{\hrule}
 & &9.19.101.293.487&1075&1562& & &79.197.397.401&297&100 \cr
10233&2464416261&4.25.11.19.43.71.101&18187&2922&10251&2477582911&8.27.25.11.79.401&799&404 \cr
 & &16.3.5.13.487.1399&1399&520& & &64.9.5.11.17.47.101&23265&54944 \cr
\noalign{\hrule}
 & &27.5.19.47.113.181&451&404& & &5.49.121.13.47.137&69&1576 \cr
10234&2465711415&8.3.11.41.101.113.181&893&350&10252&2481494015&16.3.7.11.13.23.197&111&188 \cr
 & &32.25.7.19.41.47.101&1435&1616& & &128.9.37.47.197&7289&576 \cr
\noalign{\hrule}
 & &13.19.29.53.67.97&2135&3006& & &27.7.11.53.101.223&2825&3196 \cr
10235&2467273861&4.9.5.7.19.29.61.167&1067&4240&10253&2481741801&8.25.11.17.47.101.113&159&346 \cr
 & &128.3.25.7.11.53.97&825&448& & &32.3.5.47.53.113.173&8131&9040 \cr
\noalign{\hrule}
 & &5.49.17.29.31.659&1089&430& & &25.11.31.277.1051&9243&2318 \cr
10236&2467516765&4.9.25.121.17.29.43&659&166&10254&2481857675&4.9.13.19.31.61.79&275&314 \cr
 & &16.3.11.43.83.659&473&1992& & &16.3.25.11.61.79.157&4819&3768 \cr
\noalign{\hrule}
 & &9.11.13.71.113.239&10349&6620& & &27.5.7.43.199.307&187&1348 \cr
10237&2467818639&8.3.5.13.79.131.331&317&710&10255&2482514055&8.7.11.17.199.337&555&754 \cr
 & &32.25.71.317.331&7925&5296& & &32.3.5.13.29.37.337&13949&5392 \cr
\noalign{\hrule}
 & &3.5.7.11.13.29.53.107&251&358& & &27.121.13.17.19.181&175&148 \cr
10238&2469351885&4.5.11.13.53.179.251&3103&342&10256&2482982073&8.25.7.121.13.37.181&153&1420 \cr
 & &16.9.19.29.107.179&537&152& & &64.9.125.17.37.71&4625&2272 \cr
\noalign{\hrule}
 & &25.23.1369.3137&559&3696& & &5.13.1369.103.271&3667&144 \cr
10239&2469367975&32.3.5.7.11.13.37.43&93&92&10257&2483838305&32.9.5.19.37.193&1661&1854 \cr
 & &256.9.7.11.13.23.31.43&39897&38528& & &128.81.11.103.151&1661&5184 \cr
\noalign{\hrule}
 & &3.125.7.19.67.739&3223&1950& & &9.121.23.1681.59&3159&1802 \cr
10240&2469460875&4.9.3125.11.13.293&2881&244&10258&2484136413&4.2187.13.17.41.53&4255&32686 \cr
 & &32.11.13.43.61.67&6149&976& & &16.5.23.37.59.277&1385&296 \cr
\noalign{\hrule}
 & &9.7.11.13.31.37.239&529&290& & &81.5.7.11.29.41.67&8303&5432 \cr
10241&2469664197&4.5.11.529.29.31.37&14101&2298&10259&2484290655&16.9.49.361.23.97&29&902 \cr
 & &16.3.5.59.239.383&1915&472& & &64.11.19.23.29.41&437&32 \cr
\noalign{\hrule}
 & &3.121.31.59.3721&1127&764& & &81.13.23.37.47.59&29591&57970 \cr
10242&2470472367&8.49.23.59.61.191&6479&2880&10260&2484893619&4.5.11.17.31.127.233&2679&1282 \cr
 & &1024.9.5.11.19.23.31&2185&1536& & &16.3.5.19.31.47.641&2945&5128 \cr
\noalign{\hrule}
}%
}
$$
\eject
\vglue -23 pt
\noindent\hskip 1 in\hbox to 6.5 in{\ 10261 -- 10296 \hfill\fbd 2485092555 -- 2514542745\frb}
\vskip -9 pt
$$
\vbox{
\nointerlineskip
\halign{\strut
    \vrule \ \ \hfil \frb #\ 
   &\vrule \hfil \ \ \fbb #\frb\ 
   &\vrule \hfil \ \ \frb #\ \hfil
   &\vrule \hfil \ \ \frb #\ 
   &\vrule \hfil \ \ \frb #\ \ \vrule \hskip 2 pt
   &\vrule \ \ \hfil \frb #\ 
   &\vrule \hfil \ \ \fbb #\frb\ 
   &\vrule \hfil \ \ \frb #\ \hfil
   &\vrule \hfil \ \ \frb #\ 
   &\vrule \hfil \ \ \frb #\ \vrule \cr%
\noalign{\hrule}
 & &81.5.121.17.19.157&115669&125954& & &9.11.169.23.73.89&1425&622 \cr
10261&2485092555&4.71.103.887.1123&85547&5814&10279&2500130061&4.27.25.169.19.311&253&422 \cr
 & &16.9.7.121.17.19.101&101&56& & &16.11.19.23.211.311&4009&2488 \cr
\noalign{\hrule}
 & &3.5.7.11.13.19.31.281&2701&7134& & &5.7.11.1439.4513&23777&55368 \cr
10262&2485117635&4.9.19.29.37.41.73&23&680&10280&2500269695&16.9.13.31.59.769&1537&770 \cr
 & &64.5.17.23.29.41&667&22304& & &64.3.5.7.11.29.31.53&1537&2976 \cr
\noalign{\hrule}
 & &9.7.103.409.937&24343&17784& & &3.7.13.17.23.131.179&24149&15020 \cr
10263&2486798937&16.81.11.13.19.2213&4685&6898&10281&2503016607&8.5.7.19.31.41.751&11583&2686 \cr
 & &64.5.19.937.3449&3449&3040& & &32.81.5.11.13.17.79&2133&880 \cr
\noalign{\hrule}
 & &5.11.13.29.31.53.73&303&38& & &9.7.13.31.151.653&1133&76 \cr
10264&2486935165&4.3.13.19.29.73.101&3875&3498&10282&2503431567&8.3.11.19.103.653&343&310 \cr
 & &16.9.125.11.19.31.53&225&152& & &32.5.343.19.31.103&1957&3920 \cr
\noalign{\hrule}
 & &9.5.11.31.41.59.67&1321&1924& & &3.5.49.37.251.367&159&208 \cr
10265&2487010185&8.13.31.37.41.1321&14927&32100&10283&2505121815&32.9.5.13.37.53.251&3773&6032 \cr
 & &64.3.25.11.23.59.107&535&736& & &1024.343.11.169.29&13013&14848 \cr
\noalign{\hrule}
 & &9.5.7.17.19.23.1063&451&584& & &17.41.887.4057&95&792 \cr
10266&2487563505&16.11.17.41.73.1063&183&880&10284&2508195623&16.9.5.11.19.4057&2171&1886 \cr
 & &512.3.5.121.61.73&8833&15616& & &64.3.11.13.23.41.167&5511&9568 \cr
\noalign{\hrule}
 & &27.5.7.13.31.47.139&1121&148& & &13.19.109.151.617&2143&726 \cr
10267&2487995055&8.5.13.19.31.37.59&4653&10258&10285&2508335141&4.3.121.617.2143&763&1380 \cr
 & &32.9.11.23.47.223&2453&368& & &32.9.5.7.121.23.109&5445&2576 \cr
\noalign{\hrule}
 & &9.5.11.13.19.47.433&181&2404& & &7.13.31.47.127.149&2937&1000 \cr
10268&2488215015&8.181.433.601&517&84&10286&2508943801&16.3.125.7.11.47.89&899&276 \cr
 & &64.3.7.11.47.181&1267&32& & &128.9.5.11.23.29.31&2277&9280 \cr
\noalign{\hrule}
 & &9.11.13.31.47.1327&5243&6700& & &3.25.11.17.29.31.199&2327&138 \cr
10269&2488335993&8.25.49.11.13.67.107&1147&2538&10287&2509086525&4.9.5.13.23.31.179&649&64 \cr
 & &32.27.5.49.31.37.47&735&592& & &512.11.59.179&179&15104 \cr
\noalign{\hrule}
 & &9.7.23.47.61.599&9931&10138& & &81.7.121.157.233&10547&8450 \cr
10270&2488415517&4.37.137.599.9931&224755&142692&10288&2509712667&4.9.25.7.169.53.199&11&466 \cr
 & &32.3.5.11.23.47.79.569&6259&6320& & &16.5.11.13.199.233&995&104 \cr
\noalign{\hrule}
 & &7.121.13.331.683&765&82& & &25.7.289.131.379&319&336 \cr
10271&2489289803&4.9.5.13.17.41.331&2049&2254&10289&2510998175&32.3.5.49.11.17.29.379&4167&2 \cr
 & &16.27.49.17.23.683&621&952& & &128.27.29.463&463&50112 \cr
\noalign{\hrule}
 & &3.5.289.653.881&471&182& & &3.5.59.61.193.241&217&24 \cr
10272&2493895155&4.9.5.7.13.157.881&1793&2612&10290&2511004305&16.9.5.7.31.59.61&2651&6494 \cr
 & &32.11.157.163.653&1727&2608& & &64.11.17.191.241&3247&352 \cr
\noalign{\hrule}
 & &7.13.19.41.61.577&555&22& & &9.25.121.13.47.151&1889&74 \cr
10273&2495080133&4.3.5.7.11.19.37.61&2613&2308&10291&2511805725&4.3.5.37.47.1889&163&22 \cr
 & &32.9.11.13.67.577&603&176& & &16.11.163.1889&163&15112 \cr
\noalign{\hrule}
 & &3.5.23.37.317.617&161&156& & &3.17.19.41.53.1193&25967&22946 \cr
10274&2496693585&8.9.7.13.529.37.617&78661&126800&10292&2512024941&4.7.11.17.23.149.1129&513&530 \cr
 & &256.25.11.317.7151&7151&7040& & &16.27.5.11.19.23.53.1129&10161&10120 \cr
\noalign{\hrule}
 & &7.11.151.397.541&101575&113202& & &9.7.11.19.23.43.193&185&394 \cr
10275&2497212179&4.9.25.17.19.239.331&113&1082&10293&2513277459&4.3.5.7.23.37.43.197&27599&2186 \cr
 & &16.3.5.113.331.541&1655&2712& & &16.11.13.193.1093&1093&104 \cr
\noalign{\hrule}
 & &11.19.79.151321&67405&83916& & &9.23.41.47.6301&134035&124306 \cr
10276&2498461031&8.81.5.7.13.17.37.61&389&38&10294&2513399589&4.5.7.11.13.683.2437&3221&5658 \cr
 & &32.3.5.17.19.37.389&629&240& & &16.3.5.7.11.23.41.3221&3221&3080 \cr
\noalign{\hrule}
 & &9.25.11.13.19.61.67&365&794& & &11.289.19.107.389&45029&67392 \cr
10277&2498485275&4.3.125.67.73.397&163&38&10295&2514070823&128.81.13.37.1217&5885&7102 \cr
 & &16.19.73.163.397&11899&3176& & &512.3.5.11.53.67.107&3551&3840 \cr
\noalign{\hrule}
 & &5.11.17.31.53.1627&621&280& & &3.5.121.13.19.71.79&63&184 \cr
10278&2499405535&16.27.25.7.23.1627&901&726&10296&2514542745&16.27.5.7.23.71.79&5149&3016 \cr
 & &64.81.121.17.23.53&891&736& & &256.7.13.19.29.271&1897&3712 \cr
\noalign{\hrule}
}%
}
$$
\eject
\vglue -23 pt
\noindent\hskip 1 in\hbox to 6.5 in{\ 10297 -- 10332 \hfill\fbd 2515370325 -- 2542052673\frb}
\vskip -9 pt
$$
\vbox{
\nointerlineskip
\halign{\strut
    \vrule \ \ \hfil \frb #\ 
   &\vrule \hfil \ \ \fbb #\frb\ 
   &\vrule \hfil \ \ \frb #\ \hfil
   &\vrule \hfil \ \ \frb #\ 
   &\vrule \hfil \ \ \frb #\ \ \vrule \hskip 2 pt
   &\vrule \ \ \hfil \frb #\ 
   &\vrule \hfil \ \ \fbb #\frb\ 
   &\vrule \hfil \ \ \frb #\ \hfil
   &\vrule \hfil \ \ \frb #\ 
   &\vrule \hfil \ \ \frb #\ \vrule \cr%
\noalign{\hrule}
 & &3.25.13.1373.1879&777&1102& & &729.59.89.661&21175&21836 \cr
10297&2515370325&4.9.7.19.29.37.1373&8639&3718&10315&2530294119&8.25.7.121.53.89.103&283&162 \cr
 & &16.11.169.29.53.163&19981&14344& & &32.81.5.7.53.103.283&29149&29680 \cr
\noalign{\hrule}
 & &7.19.23.53.59.263&41&18& & &3.11.23.37.227.397&2015&3206 \cr
10298&2515724659&4.9.7.19.41.53.263&2431&590&10316&2530811877&4.5.7.11.13.31.37.229&5949&11684 \cr
 & &16.3.5.11.13.17.41.59&6765&1768& & &32.9.13.23.127.661&8593&6096 \cr
\noalign{\hrule}
 & &7.11.19.31.113.491&12015&12506& & &9.11.13.43.149.307&2783&9190 \cr
10299&2516320499&4.27.5.11.169.19.37.89&3437&226&10317&2531463363&4.3.5.1331.23.919&1537&2456 \cr
 & &16.3.5.7.89.113.491&89&120& & &64.5.23.29.53.307&3335&1696 \cr
\noalign{\hrule}
 & &3.25.7.11.17.361.71&5069&6844& & &7.121.37.173.467&50657&5850 \cr
10300&2516323425&8.7.17.29.37.59.137&37411&19170&10318&2531909149&4.9.25.13.179.283&3269&3806 \cr
 & &32.27.5.11.19.71.179&179&144& & &16.3.7.11.13.173.467&13&24 \cr
\noalign{\hrule}
 & &3.5.11.17.31.103.281&685&158& & &3.49.19.61.107.139&803&1230 \cr
10301&2516738565&4.25.11.79.103.137&3807&17918&10319&2533957629&4.9.5.7.11.41.73.139&671&302 \cr
 & &16.81.289.31.47&1269&136& & &16.5.121.61.73.151&8833&6040 \cr
\noalign{\hrule}
 & &3.625.49.11.47.53&361&626& & &5.29.53.149.2213&3267&1054 \cr
10302&2517466875&4.125.7.11.361.313&621&754&10320&2534028845&4.27.5.121.17.31.53&499&2204 \cr
 & &16.27.13.19.23.29.313&64467&57592& & &32.9.11.19.29.499&9481&1584 \cr
\noalign{\hrule}
 & &5.7.11.13.17.101.293&2385&838& & &625.7.13.17.2621&6127&8748 \cr
10303&2517920405&4.9.25.53.101.419&23147&19172&10321&2534179375&8.2187.5.11.13.557&1393&5848 \cr
 & &32.3.79.293.4793&4793&3792& & &128.27.7.17.43.199&5373&2752 \cr
\noalign{\hrule}
 & &27.5.11.17.19.59.89&13&310& & &25.11.13.181.3919&833&1158 \cr
10304&2518668405&4.25.13.31.59.89&693&782&10322&2535886925&4.3.49.17.193.3919&319&3600 \cr
 & &16.9.7.11.13.17.23.31&403&1288& & &128.27.25.49.11.29&783&3136 \cr
\noalign{\hrule}
 & &9.5.7.23.59.71.83&2451&2446& & &9.7.43.67.89.157&9509&51250 \cr
10305&2518992315&4.27.7.19.23.43.71.1223&50237&557594&10323&2536141419&4.625.37.41.257&231&26 \cr
 & &16.11.83.3359.4567&36949&36536& & &16.3.125.7.11.13.37&481&11000 \cr
\noalign{\hrule}
 & &9.5.11.19.277.967&3659&5044& & &7.17.83.491.523&253&270 \cr
10306&2519213895&8.11.13.19.97.3659&1725&1934&10324&2536344461&4.27.5.7.11.23.83.491&8891&4472 \cr
 & &32.3.25.13.23.97.967&1495&1552& & &64.3.5.11.13.17.43.523&1677&1760 \cr
\noalign{\hrule}
 & &3.29.83.409.853&187&1040& & &27.5.11.13.17.59.131&23161&15956 \cr
10307&2519240817&32.5.11.13.17.29.83&2997&3412&10325&2536541865&8.9.19.23.53.3989&695&524 \cr
 & &256.81.11.37.853&999&1408& & &64.5.131.139.3989&3989&4448 \cr
\noalign{\hrule}
 & &3.11.17.19.41.73.79&4175&4394& & &27.5.49.41.47.199&671&436 \cr
10308&2520288573&4.25.2197.17.79.167&13779&586&10326&2536673895&8.49.11.61.109.199&1551&8200 \cr
 & &16.9.5.13.293.1531&19903&35160& & &128.3.25.121.41.47&121&320 \cr
\noalign{\hrule}
 & &27.49.13.47.3119&68167&84664& & &9.5.11.19.211.1279&6929&7140 \cr
10309&2521253007&16.11.19.557.6197&35&6162&10327&2538117945&8.27.25.7.169.17.19.41&12203&31022 \cr
 & &64.3.5.7.13.19.79&79&3040& & &32.13.12203.15511&201643&195248 \cr
\noalign{\hrule}
 & &81.23.31.149.293&10153&13580& & &9.5.7.17.41.43.269&25781&32054 \cr
10310&2521322721&8.5.7.11.13.31.71.97&513&1192&10328&2539592685&4.49.11.29.31.47.127&817&3120 \cr
 & &128.27.13.19.71.149&923&1216& & &128.3.5.11.13.19.29.43&2717&1856 \cr
\noalign{\hrule}
 & &3.5.17.29.47.53.137&11&3984& & &27.121.169.43.107&25075&25118 \cr
10311&2523669465&32.9.11.53.83&2195&2204&10329&2540317923&4.25.11.17.19.59.107.661&713&2532 \cr
 & &256.5.11.19.29.439&4829&2432& & &32.3.5.19.23.31.211.661&389329&388240 \cr
\noalign{\hrule}
 & &9.5.11.19.29.47.197&5483&1018& & &9.25.29.401.971&3367&5372 \cr
10312&2525345955&4.3.29.509.5483&24883&19400&10330&2540645775&8.5.7.13.17.29.37.79&6323&1188 \cr
 & &64.25.97.149.167&14453&26720& & &64.27.11.17.6323&6323&17952 \cr
\noalign{\hrule}
 & &5.11.17.23.257.457&3071&7440& & &17.23.41.257.617&71269&87300 \cr
10313&2525740745&32.3.25.11.31.37.83&457&468&10331&2542019639&8.9.25.121.19.31.97&391&676 \cr
 & &256.27.13.31.83.457&10881&10624& & &64.3.5.11.169.17.23.31&5577&4960 \cr
\noalign{\hrule}
 & &9.29.181.199.269&3289&2482& & &27.137.211.3257&2245&1012 \cr
10314&2528857971&4.3.11.13.17.23.73.181&145&398&10332&2542052673&8.3.5.11.23.211.449&3257&3706 \cr
 & &16.5.13.17.29.73.199&949&680& & &32.5.17.23.109.3257&1955&1744 \cr
\noalign{\hrule}
}%
}
$$
\eject
\vglue -23 pt
\noindent\hskip 1 in\hbox to 6.5 in{\ 10333 -- 10368 \hfill\fbd 2542108361 -- 2560490111\frb}
\vskip -9 pt
$$
\vbox{
\nointerlineskip
\halign{\strut
    \vrule \ \ \hfil \frb #\ 
   &\vrule \hfil \ \ \fbb #\frb\ 
   &\vrule \hfil \ \ \frb #\ \hfil
   &\vrule \hfil \ \ \frb #\ 
   &\vrule \hfil \ \ \frb #\ \ \vrule \hskip 2 pt
   &\vrule \ \ \hfil \frb #\ 
   &\vrule \hfil \ \ \fbb #\frb\ 
   &\vrule \hfil \ \ \frb #\ \hfil
   &\vrule \hfil \ \ \frb #\ 
   &\vrule \hfil \ \ \frb #\ \vrule \cr%
\noalign{\hrule}
 & &13.841.157.1481&121725&110792& & &27.7.11.17.19.29.131&155&358 \cr
10333&2542108361&16.9.25.11.541.1259&6123&172&10351&2551093083&4.5.11.17.31.131.179&18767&14706 \cr
 & &128.27.5.13.43.157&1161&320& & &16.9.5.49.19.43.383&1915&2408 \cr
\noalign{\hrule}
 & &27.25.17.37.53.113&451&26& & &81.11.13.23.61.157&387797&387940 \cr
10334&2542779675&4.3.11.13.37.41.113&27853&26500&10352&2551398993&8.5.7.17.23.37.47.163.223&18603&930262 \cr
 & &32.125.7.23.53.173&1211&1840& & &32.27.13.29.43.53.373&19769&19952 \cr
\noalign{\hrule}
 & &81.343.11.53.157&247&194& & &3.125.7.13.37.43.47&57&1682 \cr
10335&2543005773&4.9.7.11.13.19.97.157&265&1462&10353&2551765125&4.9.7.19.841.43&925&6644 \cr
 & &16.5.13.17.43.53.97&8245&4472& & &32.25.11.37.151&1661&16 \cr
\noalign{\hrule}
 & &9.5.7.23.37.53.179&3053&3212& & &43.47.61.127.163&71991&51290 \cr
10336&2543132655&8.3.11.23.37.43.71.73&895&784&10354&2552039981&4.9.5.19.23.223.421&4953&4730 \cr
 & &256.5.49.11.43.71.179&5467&5504& & &16.27.25.11.13.19.43.127&3861&3800 \cr
\noalign{\hrule}
 & &9.125.49.29.37.43&7381&3756& & &121.17.73.89.191&2793&6040 \cr
10337&2543412375&8.27.7.121.61.313&1247&400&10355&2552586839&16.3.5.49.19.89.151&573&484 \cr
 & &256.25.29.43.313&313&128& & &128.9.5.7.121.19.191&665&576 \cr
\noalign{\hrule}
 & &9.11.13.19.67.1553&1413&140& & &25.49.11.169.19.59&97347&109678 \cr
10338&2544359103&8.81.5.7.11.13.157&29507&28408&10356&2552825275&4.3.29.31.37.61.877&975&98 \cr
 & &128.19.53.67.1553&53&64& & &16.9.25.49.13.31.61&279&488 \cr
\noalign{\hrule}
 & &3.17.47.67.83.191&7651&2090& & &5.11.13.17.41.47.109&6395&1272 \cr
10339&2545975947&4.5.7.11.19.47.1093&711&382&10357&2553072665&16.3.25.13.53.1279&7973&9252 \cr
 & &16.9.5.11.19.79.191&869&2280& & &128.27.7.17.67.257&17219&12096 \cr
\noalign{\hrule}
 & &3.5.7.11.19.157.739&181&104& & &27.13.53.107.1283&15221&52778 \cr
10340&2546124735&16.13.157.181.739&9405&19012&10358&2553838443&4.11.31.491.2399&1445&954 \cr
 & &128.9.5.49.11.19.97&291&448& & &16.9.5.11.289.31.53&1445&2728 \cr
\noalign{\hrule}
 & &5.7.13.23.31.47.167&27115&27828& & &3.5.7.31.3721.211&176719&169334 \cr
10341&2546333335&8.9.25.11.13.17.29.773&7567&20708&10359&2555601405&4.11.19.43.71.131.179&4517&1116 \cr
 & &64.3.7.11.23.31.47.167&33&32& & &32.9.11.31.71.4517&13551&12496 \cr
\noalign{\hrule}
 & &9.5.11.151.173.197&1609&52& & &5.13.31.43.163.181&9189&1406 \cr
10342&2547384345&8.5.13.197.1609&10951&9966&10360&2556287435&4.9.19.31.37.1021&495&526 \cr
 & &32.3.11.47.151.233&233&752& & &16.81.5.11.19.37.263&32967&39976 \cr
\noalign{\hrule}
 & &25.169.23.157.167&3893&282& & &9.5.13.19.23.73.137&66281&42854 \cr
10343&2547831325&4.3.169.17.47.229&315&484&10361&2556705645&4.7.79.839.3061&1111&1950 \cr
 & &32.27.5.7.121.229&22869&3664& & &16.3.25.7.11.13.79.101&4345&5656 \cr
\noalign{\hrule}
 & &5.7.11.19.79.4409&2639&1770& & &27.19.113.44111&65093&67240 \cr
10344&2547894965&4.3.25.49.13.19.29.59&473&948&10362&2557070559&16.9.5.7.17.1681.547&4901&22 \cr
 & &32.9.11.13.43.59.79&2537&1872& & &64.5.11.169.29.41&6929&51040 \cr
\noalign{\hrule}
 & &5.11.13.17.31.67.101&445&666& & &9.25.121.29.41.79&941&1666 \cr
10345&2549839435&4.9.25.31.37.67.89&2299&4376&10363&2557271475&4.3.49.11.17.41.941&167&2990 \cr
 & &64.3.121.19.37.547&23199&17504& & &16.5.7.13.17.23.167&2093&22712 \cr
\noalign{\hrule}
 & &3.5.11.13.29.179.229&43&186& & &3.25.7.41.157.757&1853&418 \cr
10346&2549845155&4.9.5.29.31.43.179&3563&1628&10364&2558224725&4.5.11.17.19.109.157&7749&9364 \cr
 & &32.7.11.31.37.509&15779&4144& & &32.27.7.11.41.2341&2341&1584 \cr
\noalign{\hrule}
 & &3.7.11.19.37.113.139&1339&1302& & &3.5.7.43.733.773&12441&13214 \cr
10347&2550706851&4.9.49.11.13.31.103.113&3445&41978&10365&2558239635&4.9.11.13.29.43.6607&1175&5432 \cr
 & &16.5.169.53.139.151&8003&6760& & &64.25.7.13.29.47.97&17719&15520 \cr
\noalign{\hrule}
 & &9.7.13.127.137.179&50425&71368& & &3.25.7.23.59.3593&1877&1716 \cr
10348&2550710799&16.25.11.811.2017&9129&11146&10366&2559743025&8.9.25.11.13.59.1877&1141&18034 \cr
 & &64.3.11.17.179.5573&5573&5984& & &32.7.11.71.127.163&11573&22352 \cr
\noalign{\hrule}
 & &3.125.7.23.29.31.47&1937&3938& & &3.5.47.67.83.653&2117&7678 \cr
10349&2551024875&4.7.11.13.31.149.179&225&2162&10367&2560089765&4.11.29.47.73.349&1541&1890 \cr
 & &16.9.25.23.47.179&179&24& & &16.27.5.7.11.23.29.67&2277&1624 \cr
\noalign{\hrule}
 & &3.11.13.61.71.1373&553&370& & &23.47.71.73.457&5291&5220 \cr
10350&2551032627&4.5.7.11.37.79.1373&1121&252&10368&2560490111&8.9.5.11.13.29.37.47.73&12121&1828 \cr
 & &32.9.5.49.19.37.59&32745&14896& & &64.3.5.11.17.23.31.457&1705&1632 \cr
\noalign{\hrule}
}%
}
$$
\eject
\vglue -23 pt
\noindent\hskip 1 in\hbox to 6.5 in{\ 10369 -- 10404 \hfill\fbd 2560856897 -- 2584053175\frb}
\vskip -9 pt
$$
\vbox{
\nointerlineskip
\halign{\strut
    \vrule \ \ \hfil \frb #\ 
   &\vrule \hfil \ \ \fbb #\frb\ 
   &\vrule \hfil \ \ \frb #\ \hfil
   &\vrule \hfil \ \ \frb #\ 
   &\vrule \hfil \ \ \frb #\ \ \vrule \hskip 2 pt
   &\vrule \ \ \hfil \frb #\ 
   &\vrule \hfil \ \ \fbb #\frb\ 
   &\vrule \hfil \ \ \frb #\ \hfil
   &\vrule \hfil \ \ \frb #\ 
   &\vrule \hfil \ \ \frb #\ \vrule \cr%
\noalign{\hrule}
 & &17.31.71.89.769&299&228& & &3.49.11.13.29.41.103&4573&50944 \cr
10369&2560856897&8.3.13.19.23.89.769&429&340&10387&2574378807&512.17.199.269&35&234 \cr
 & &64.9.5.11.169.17.19.23&35321&33120& & &2048.9.5.7.13.17&85&3072 \cr
\noalign{\hrule}
 & &11.43.113.191.251&1049&3810& & &3.5.121.13.29.53.71&6295&118 \cr
10370&2562398509&4.3.5.127.191.1049&429&620&10388&2574851565&4.25.13.59.1259&371&396 \cr
 & &32.9.25.11.13.31.127&10075&18288& & &32.9.7.11.53.1259&1259&336 \cr
\noalign{\hrule}
 & &27.5.121.19.23.359&413&622& & &25.7.67.157.1399&4047&2948 \cr
10371&2562683805&4.3.7.11.59.311.359&119761&113230&10389&2575314175&8.3.5.11.19.4489.71&4197&292 \cr
 & &16.5.169.23.41.67.127&21463&21976& & &64.9.19.73.1399&1387&288 \cr
\noalign{\hrule}
 & &9.11.17.29.131.401&389&790& & &3.25.7.19.23.103.109&711&1474 \cr
10372&2563880517&4.5.11.17.29.79.389&681&188&10390&2575754475&4.27.5.11.67.79.103&109&406 \cr
 & &32.3.5.47.227.389&10669&31120& & &16.7.29.67.79.109&1943&632 \cr
\noalign{\hrule}
 & &9.5.121.17.103.269&13091&628& & &3.25.13.103.113.227&17&22 \cr
10373&2564698455&8.3.5.13.19.53.157&1751&1694&10391&2576001675&4.5.11.17.103.113.227&423&12062 \cr
 & &32.7.121.17.103.157&157&112& & &16.9.17.37.47.163&29563&3912 \cr
\noalign{\hrule}
 & &3.5.11.61.227.1123&5681&6804& & &27.5.7.23.29.61.67&149&454 \cr
10374&2565779865&8.729.7.13.19.23.61&2015&11836&10392&2576097405&4.3.7.23.29.149.227&2717&2050 \cr
 & &64.5.11.169.31.269&5239&8608& & &16.25.11.13.19.41.149&31141&21320 \cr
\noalign{\hrule}
 & &3.49.11.17.31.3011&60415&32926& & &5.49.11.503.1901&60407&32742 \cr
10375&2565850749&4.5.43.101.163.281&33759&12044&10393&2576967085&4.9.17.29.107.2083&1573&4676 \cr
 & &32.9.121.31.3011&33&16& & &32.3.7.121.13.17.167&6513&2992 \cr
\noalign{\hrule}
 & &5.11.17.59.193.241&3263&18& & &25.7.121.19.43.149&31941&25534 \cr
10376&2565889645&4.9.13.241.251&5015&4774&10394&2577696275&4.27.49.169.17.751&1331&580 \cr
 & &16.3.5.7.11.17.31.59&217&24& & &32.9.5.1331.13.17.29&2871&3536 \cr
\noalign{\hrule}
 & &9.5.11.29.31.73.79&433&466& & &9.13.109.131.1543&97405&83126 \cr
10377&2566343835&4.3.5.73.79.233.433&93697&7192&10395&2577802149&4.5.7.121.23.89.467&545&78 \cr
 & &64.29.31.43.2179&2179&1376& & &16.3.25.121.13.23.109&575&968 \cr
\noalign{\hrule}
 & &7.71.79.151.433&29887&35496& & &25.7.13.361.43.73&4101&19624 \cr
10378&2567132729&16.9.7.121.13.17.19.29&71&50&10396&2577982225&16.3.7.11.223.1367&645&722 \cr
 & &64.3.25.13.17.19.29.71&16575&17632& & &64.9.5.361.43.223&223&288 \cr
\noalign{\hrule}
 & &9.5.7.13.19.61.541&2519&4514& & &9.7.13.43.179.409&16027&53246 \cr
10379&2567642805&4.3.11.37.3721.229&1517&2204&10397&2578271787&4.11.31.47.79.337&25645&15198 \cr
 & &32.11.19.29.1369.41&15059&19024& & &16.3.5.17.23.149.223&33227&15640 \cr
\noalign{\hrule}
 & &9.5.13.529.43.193&3727&1218& & &9.121.13.19.43.223&821&598 \cr
10380&2568250035&4.27.7.23.29.3727&7141&10868&10398&2579277987&4.3.11.169.19.23.821&45985&10664 \cr
 & &32.7.11.13.19.37.193&703&1232& & &64.5.17.31.43.541&9197&4960 \cr
\noalign{\hrule}
 & &27.121.17.103.449&50803&3526& & &27.31.83.107.347&9125&244 \cr
10381&2568512133&4.41.43.101.503&1133&630&10399&2579388759&8.125.31.61.73&979&1284 \cr
 & &16.9.5.7.11.101.103&505&56& & &64.3.25.11.89.107&275&2848 \cr
\noalign{\hrule}
 & &3.25.23.43.59.587&29&616& & &3.5.11.13.79.97.157&599&586 \cr
10382&2568902775&16.5.7.11.23.29.59&5283&3572&10400&2580630195&4.11.97.157.293.599&171513&3994 \cr
 & &128.9.19.47.587&141&1216& & &16.9.17.19.59.1997&37943&24072 \cr
\noalign{\hrule}
 & &81.125.13.29.673&2101&2774& & &9.5.7.13.23.79.347&1069&1622 \cr
10383&2568925125&4.27.11.19.29.73.191&84125&69248&10401&2581893405&4.5.347.811.1069&1273&462 \cr
 & &1024.125.541.673&541&512& & &16.3.7.11.19.67.1069&14003&8552 \cr
\noalign{\hrule}
 & &9.49.17.29.53.223&20185&6352& & &9.5.13.29.41.47.79&253&458 \cr
10384&2569604247&32.5.7.11.367.397&4929&8966&10402&2582632845&4.11.13.23.29.47.229&9717&1046 \cr
 & &128.3.31.53.4483&4483&1984& & &16.3.11.41.79.523&523&88 \cr
\noalign{\hrule}
 & &27.19.37.43.47.67&161&1430& & &11.13.101.113.1583&1413&170 \cr
10385&2570160267&4.5.7.11.13.19.23.67&2673&3008&10403&2583549397&4.9.5.13.17.101.157&1991&678 \cr
 & &512.243.7.121.47&1089&1792& & &16.27.5.11.113.181&135&1448 \cr
\noalign{\hrule}
 & &3.7.13.37.149.1709&5925&7634& & &25.11.37.229.1109&23&252 \cr
10386&2572128741&4.9.25.11.37.79.347&1709&1214&10404&2584053175&8.9.7.23.37.1109&425&684 \cr
 & &16.5.347.607.1709&3035&2776& & &64.81.25.17.19.23&1863&10336 \cr
\noalign{\hrule}
}%
}
$$
\eject
\vglue -23 pt
\noindent\hskip 1 in\hbox to 6.5 in{\ 10405 -- 10440 \hfill\fbd 2585473407 -- 2600649675\frb}
\vskip -9 pt
$$
\vbox{
\nointerlineskip
\halign{\strut
    \vrule \ \ \hfil \frb #\ 
   &\vrule \hfil \ \ \fbb #\frb\ 
   &\vrule \hfil \ \ \frb #\ \hfil
   &\vrule \hfil \ \ \frb #\ 
   &\vrule \hfil \ \ \frb #\ \ \vrule \hskip 2 pt
   &\vrule \ \ \hfil \frb #\ 
   &\vrule \hfil \ \ \fbb #\frb\ 
   &\vrule \hfil \ \ \frb #\ \hfil
   &\vrule \hfil \ \ \frb #\ 
   &\vrule \hfil \ \ \frb #\ \vrule \cr%
\noalign{\hrule}
 & &9.11.17.41.89.421&66263&83524& & &9.5.13.361.71.173&2623&7238 \cr
10405&2585473407&8.7.19.23.43.67.157&33&10&10423&2593985355&4.3.7.11.19.43.47.61&7&50 \cr
 & &32.3.5.7.11.19.67.157&8911&12560& & &16.25.49.11.47.61&2695&22936 \cr
\noalign{\hrule}
 & &25.11.13.17.157.271&1541&186& & &7.11.169.17.37.317&45&214 \cr
10406&2585793925&4.3.5.13.17.23.31.67&1413&542&10424&2594701109&4.9.5.11.17.107.317&3283&2106 \cr
 & &16.27.31.157.271&27&248& & &16.729.5.49.13.67&2345&5832 \cr
\noalign{\hrule}
 & &3.25.7.11.41.67.163&61&676& & &3.11.13.23.409.643&5145&9644 \cr
10407&2585819775&8.5.7.169.61.163&603&538&10425&2594892729&8.9.5.343.13.2411&1615&796 \cr
 & &32.9.13.61.67.269&2379&4304& & &64.25.49.17.19.199&64277&39200 \cr
\noalign{\hrule}
 & &3.5.17.31.41.79.101&2937&302& & &3.49.17.19.31.41.43&325&326 \cr
10408&2586033795&4.9.11.89.101.151&629&730&10426&2594979093&4.25.7.13.17.19.41.43.163&899&4554 \cr
 & &16.5.11.17.37.73.89&3293&6424& & &16.9.5.11.13.23.29.31.163&43355&43032 \cr
\noalign{\hrule}
 & &3.5.17.19.41.47.277&3&44& & &5.11.101.463.1009&51129&4366 \cr
10409&2586159255&8.9.5.11.17.19.277&1927&2782&10427&2595112685&4.9.13.19.23.37.59&505&616 \cr
 & &32.11.13.41.47.107&1177&208& & &64.3.5.7.11.13.23.101&273&736 \cr
\noalign{\hrule}
 & &81.5.29.257.857&3317&968& & &5.7.11.17.23.43.401&27&962 \cr
10410&2586824505&16.121.31.107.257&717&460&10428&2595675005&4.27.7.13.37.401&731&472 \cr
 & &128.3.5.11.23.31.239&7843&15296& & &64.9.13.17.43.59&117&1888 \cr
\noalign{\hrule}
 & &27.5.67.73.3919&2563&6482& & &81.5.11.13.107.419&44813&67858 \cr
10411&2587656915&4.7.11.73.233.463&195&268&10429&2596503195&4.7.37.41.131.1093&1877&2970 \cr
 & &32.3.5.7.11.13.67.233&1631&2288& & &16.27.5.7.11.41.1877&1877&2296 \cr
\noalign{\hrule}
 & &9.7.13.17.151.1231&19297&1630& & &3.11.19.23.31.37.157&2977&1830 \cr
10412&2588021163&4.5.7.23.163.839&13541&12702&10430&2596919259&4.9.5.13.61.157.229&16169&2200 \cr
 & &16.3.5.11.29.73.1231&803&1160& & &64.125.11.19.23.37&125&32 \cr
\noalign{\hrule}
 & &5.7.13.67.73.1163&423&88& & &3.25.11.29.961.113&1389&2114 \cr
10413&2588146015&16.9.11.13.47.1163&335&1498&10431&2598087525&4.9.7.11.31.151.463&691&970 \cr
 & &64.3.5.7.11.67.107&107&1056& & &16.5.7.97.463.691&44911&38696 \cr
\noalign{\hrule}
 & &5.7.19.23.31.43.127&11737&15822& & &9.7.41.43.149.157&24115&30866 \cr
10414&2589305845&4.27.121.23.97.293&19&272&10432&2598237117&4.5.49.11.13.23.53.61&2021&894 \cr
 & &128.9.11.17.19.293&2637&11968& & &16.3.13.43.47.61.149&611&488 \cr
\noalign{\hrule}
 & &9.25.11.169.41.151&1021&266& & &49.13.29.41.47.73&711&2720 \cr
10415&2589540525&4.5.7.13.19.41.1021&2241&4906&10433&2598615383&64.9.5.13.17.29.79&2773&308 \cr
 & &16.27.11.19.83.223&1577&5352& & &512.3.7.11.47.59&177&2816 \cr
\noalign{\hrule}
 & &3.7.11.29.193.2003&835&516& & &11.13.23.3481.227&1615&1674 \cr
10416&2589692721&8.9.5.43.167.2003&5423&12604&10434&2598925043&4.27.5.17.19.31.59.227&11703&1690 \cr
 & &64.5.11.17.23.29.137&2329&3680& & &16.81.25.169.47.83&26975&30456 \cr
\noalign{\hrule}
 & &3.11.169.41.47.241&531&80& & &9.25.17.23.31.953&689&264 \cr
10417&2589997839&32.27.5.13.59.241&10475&10234&10435&2599045425&16.27.11.13.23.31.53&23825&20536 \cr
 & &128.125.7.17.43.419&126119&136000& & &256.25.17.151.953&151&128 \cr
\noalign{\hrule}
 & &121.19.53.89.239&105&16& & &125.49.17.109.229&45741&49634 \cr
10418&2591807537&32.3.5.7.19.53.239&2403&2138&10436&2599064125&4.3.7.13.23.79.83.193&1199&618 \cr
 & &128.81.7.89.1069&7483&5184& & &16.9.11.13.103.109.193&14729&13896 \cr
\noalign{\hrule}
 & &9.5.13.2099.2111&77717&58718& & &19.31.61.72361&54145&18216 \cr
10419&2592128565&4.11.17.23.31.109.157&2117&390&10437&2599858369&16.9.5.49.11.13.17.23&269&122 \cr
 & &16.3.5.13.17.29.31.73&2117&4216& & &64.3.5.11.13.61.269&429&160 \cr
\noalign{\hrule}
 & &11.13.23.31.59.431&64665&82306& & &27.125.13.19.3119&293&3412 \cr
10420&2592715411&4.27.5.7.479.5879&2221&3658&10438&2600076375&8.9.25.293.853&539&314 \cr
 & &16.9.5.7.31.59.2221&2221&2520& & &32.49.11.157.293&7693&51568 \cr
\noalign{\hrule}
 & &81.7.13.61.73.79&6095&1276& & &5.13.17.61.173.223&3381&7172 \cr
10421&2593021977&8.5.11.23.29.53.73&3429&4966&10439&2600417495&8.3.5.49.11.13.23.163&213&928 \cr
 & &32.27.11.13.127.191&1397&3056& & &512.9.7.23.29.71&14413&52992 \cr
\noalign{\hrule}
 & &3.5.13.19.43.73.223&9537&6742& & &9.25.13.23.29.31.43&349&726 \cr
10422&2593488885&4.9.11.289.19.3371&3139&232&10440&2600649675&4.27.121.23.31.349&46625&54652 \cr
 & &64.11.17.29.43.73&187&928& & &32.125.13.373.1051&5255&5968 \cr
\noalign{\hrule}
}%
}
$$
\eject
\vglue -23 pt
\noindent\hskip 1 in\hbox to 6.5 in{\ 10441 -- 10476 \hfill\fbd 2604366765 -- 2623270221\frb}
\vskip -9 pt
$$
\vbox{
\nointerlineskip
\halign{\strut
    \vrule \ \ \hfil \frb #\ 
   &\vrule \hfil \ \ \fbb #\frb\ 
   &\vrule \hfil \ \ \frb #\ \hfil
   &\vrule \hfil \ \ \frb #\ 
   &\vrule \hfil \ \ \frb #\ \ \vrule \hskip 2 pt
   &\vrule \ \ \hfil \frb #\ 
   &\vrule \hfil \ \ \fbb #\frb\ 
   &\vrule \hfil \ \ \frb #\ \hfil
   &\vrule \hfil \ \ \frb #\ 
   &\vrule \hfil \ \ \frb #\ \vrule \cr%
\noalign{\hrule}
 & &9.5.7.11.13.17.19.179&3569&526& & &9.25.7.11.961.157&893&68 \cr
10441&2604366765&4.11.19.43.83.263&2735&2262&10459&2613944025&8.3.7.17.19.47.157&415&572 \cr
 & &16.3.5.13.29.83.547&2407&4376& & &64.5.11.13.17.19.83&1577&7072 \cr
\noalign{\hrule}
 & &9.11.67.607.647&2363&3100& & &81.5.13.31.37.433&7291&7694 \cr
10442&2604971457&8.25.17.31.139.647&111&536&10460&2614867515&4.23.317.433.3847&112871&24390 \cr
 & &128.3.31.37.67.139&4309&2368& & &16.9.5.11.31.271.331&2981&2648 \cr
\noalign{\hrule}
 & &3.11.37.47.83.547&49&1690& & &125.11.71.73.367&3279&1904 \cr
10443&2605427187&4.5.49.11.169.83&2105&1962&10461&2615471375&32.3.7.17.367.1093&363&730 \cr
 & &16.9.25.13.109.421&31575&11336& & &128.9.5.7.121.17.73&1071&704 \cr
\noalign{\hrule}
 & &27.11.17.19.101.269&8255&18914& & &9.49.11.13.19.37.59&425&278 \cr
10444&2606349339&4.9.5.49.13.127.193&1313&424&10462&2615664051&4.3.25.11.13.17.59.139&6517&4988 \cr
 & &64.5.7.169.53.101&1855&5408& & &32.5.343.17.19.29.43&3655&3248 \cr
\noalign{\hrule}
 & &5.49.13.491.1667&711&956& & &729.13.19.73.199&6161&2380 \cr
10445&2606912945&8.9.13.79.239.491&759&268&10463&2615775201&8.81.5.7.17.61.101&739&638 \cr
 & &64.27.11.23.67.239&70983&49312& & &32.5.7.11.29.61.739&97295&82768 \cr
\noalign{\hrule}
 & &3.25.13.19.23.29.211&817&5302& & &27.11.13.17.23.1733&6145&4412 \cr
10446&2607152925&4.5.11.361.43.241&783&422&10464&2616225183&8.5.11.13.1103.1229&909&194 \cr
 & &16.27.11.29.43.211&387&88& & &32.9.97.101.1229&9797&19664 \cr
\noalign{\hrule}
 & &9.5.7.11.169.61.73&181&38& & &5.13.59.131.5209&63261&37216 \cr
10447&2607610005&4.3.5.7.13.19.61.181&319&136&10465&2616923465&64.81.11.71.1163&1753&590 \cr
 & &64.11.17.19.29.181&3439&15776& & &256.27.5.59.1753&1753&3456 \cr
\noalign{\hrule}
 & &9.49.121.19.31.83&427&1340& & &25.1331.31.43.59&2989&2214 \cr
10448&2608659207&8.3.5.343.11.61.67&1963&248&10466&2616978925&4.27.49.11.41.59.61&31&1978 \cr
 & &128.13.31.61.151&9211&832& & &16.9.23.31.43.61&1403&72 \cr
\noalign{\hrule}
 & &9.17.31.61.71.127&1375&826& & &5.7.11.83.101.811&3483&4294 \cr
10449&2608825491&4.125.7.11.17.59.127&7739&3294&10467&2617466005&4.81.5.19.43.83.113&4247&6488 \cr
 & &16.27.25.61.71.109&327&200& & &64.3.31.43.137.811&3999&4384 \cr
\noalign{\hrule}
 & &5.7.11.289.47.499&81153&63058& & &9.5.11.2197.29.83&4123&6862 \cr
10450&2609498045&4.9.41.71.127.769&3757&1450&10468&2617648605&4.3.7.19.29.31.47.73&3887&5654 \cr
 & &16.3.25.13.289.29.71&2059&1560& & &16.11.169.23.73.257&1679&2056 \cr
\noalign{\hrule}
 & &27.11.19.29.37.431&10045&10342& & &243.5.11.31.71.89&97&542 \cr
10451&2609678709&4.5.49.41.431.5171&3663&1508&10469&2618056485&4.27.11.31.97.271&47515&44896 \cr
 & &32.9.49.11.13.29.37.41&637&656& & &256.5.13.17.23.43.61&44591&38272 \cr
\noalign{\hrule}
 & &27.5.7.11.23.61.179&5149&5770& & &27.13.89.191.439&181&620 \cr
10452&2610569115&4.25.7.11.19.271.577&4519&6444&10470&2619358911&8.3.5.13.31.181.191&7175&9658 \cr
 & &32.9.179.271.4519&4519&4336& & &32.125.7.11.41.439&3157&2000 \cr
\noalign{\hrule}
 & &9.5.7.11.17.23.41.47&114409&114904& & &9.343.11.19.31.131&2479&2822 \cr
10453&2610728505&16.23.53.191.271.599&197579&35250&10471&2620088163&4.11.17.37.67.83.131&7695&7826 \cr
 & &64.3.125.41.47.61.79&1975&1952& & &16.81.5.7.13.19.37.43.67&22311&22360 \cr
\noalign{\hrule}
 & &27.25.7.17.19.29.59&24779&17996& & &9.11.13.17.23.41.127&5539&10492 \cr
10454&2611285425&8.9.11.71.349.409&138125&134444&10472&2620250919&8.3.11.29.43.61.191&2921&820 \cr
 & &64.625.13.17.19.29.61&793&800& & &64.5.23.41.61.127&61&160 \cr
\noalign{\hrule}
 & &5.121.13.83.4001&2637&2758& & &11.127.577.3251&2270535&2271112 \cr
10455&2611832795&4.9.7.197.293.4001&4069&68&10473&2620530319&16.3.5.23.229.661.12343&5841&6502 \cr
 & &32.3.13.17.293.313&4981&15024& & &64.27.5.11.23.59.229.3251&36639&36640 \cr
\noalign{\hrule}
 & &9.5.11.17.19.59.277&2915&1794& & &27.17.19.31.89.109&5893&51550 \cr
10456&2613000555&4.27.25.121.13.23.53&277&398&10474&2622675051&4.25.71.83.1031&2431&3462 \cr
 & &16.13.23.53.199.277&4577&5512& & &16.3.25.11.13.17.577&7501&2200 \cr
\noalign{\hrule}
 & &9.25.7.13.109.1171&5053&802& & &5.19.103.181.1481&1431&526 \cr
10457&2613408525&4.3.5.7.31.163.401&14773&2342&10475&2622976885&4.27.53.263.1481&2353&2090 \cr
 & &16.11.17.79.1171&79&1496& & &16.9.5.11.13.19.53.181&689&792 \cr
\noalign{\hrule}
 & &243.13.757.1093&127879&137720& & &9.11.169.17.23.401&4957&11774 \cr
10458&2613759759&16.5.11.41.313.3119&3281&162&10476&2623270221&4.7.23.841.4957&2145&2812 \cr
 & &64.81.5.17.41.193&3281&6560& & &32.3.5.7.11.13.19.29.37&3857&2960 \cr
\noalign{\hrule}
}%
}
$$
\eject
\vglue -23 pt
\noindent\hskip 1 in\hbox to 6.5 in{\ 10477 -- 10512 \hfill\fbd 2623893195 -- 2643999889\frb}
\vskip -9 pt
$$
\vbox{
\nointerlineskip
\halign{\strut
    \vrule \ \ \hfil \frb #\ 
   &\vrule \hfil \ \ \fbb #\frb\ 
   &\vrule \hfil \ \ \frb #\ \hfil
   &\vrule \hfil \ \ \frb #\ 
   &\vrule \hfil \ \ \frb #\ \ \vrule \hskip 2 pt
   &\vrule \ \ \hfil \frb #\ 
   &\vrule \hfil \ \ \fbb #\frb\ 
   &\vrule \hfil \ \ \frb #\ \hfil
   &\vrule \hfil \ \ \frb #\ 
   &\vrule \hfil \ \ \frb #\ \vrule \cr%
\noalign{\hrule}
 & &3.5.7.11.41.67.827&721&106& & &9.11.17.23.149.457&11525&11272 \cr
10477&2623893195&4.49.11.53.67.103&15327&20786&10495&2635811937&16.25.457.461.1409&459&2 \cr
 & &16.9.13.19.131.547&32357&13128& & &64.27.25.17.1409&4227&800 \cr
\noalign{\hrule}
 & &9.49.11.23.29.811&19375&722& & &3.5.7.37.59.11503&25483&32032 \cr
10478&2624085387&4.625.7.361.31&277&312&10496&2636660145&64.49.11.13.17.1499&333&1166 \cr
 & &64.3.125.13.19.277&30875&8864& & &256.9.121.13.37.53&4719&6784 \cr
\noalign{\hrule}
 & &3.5.7.11.31.83.883&2037&2378& & &125.13.17.31.3079&51359&984 \cr
10479&2624112645&4.9.49.29.41.83.97&4693&8096&10497&2636778625&16.3.7.11.23.29.41&1105&1128 \cr
 & &256.11.13.361.23.97&35017&38272& & &256.9.5.13.17.41.47&1927&1152 \cr
\noalign{\hrule}
 & &9.5.11.31.41.43.97&673&394& & &3.11.13.29.127.1669&1023&646 \cr
10480&2624163795&4.5.41.43.197.673&5723&34662&10498&2637031683&4.9.121.17.19.31.127&28145&46658 \cr
 & &16.3.53.59.97.109&3127&872& & &16.5.13.41.433.569&17753&22760 \cr
\noalign{\hrule}
 & &5.7.121.169.19.193&4307&2448& & &81.7.23.31.61.107&19&88 \cr
10481&2624526905&32.9.11.17.19.59.73&193&16&10499&2638676817&16.27.7.11.19.31.61&851&1040 \cr
 & &1024.3.17.73.193&3723&512& & &512.5.11.13.19.23.37&9139&14080 \cr
\noalign{\hrule}
 & &3.11.13.59.97.1069&389&378& & &5.19.29.47.89.229&617&528 \cr
10482&2624573523&4.81.7.97.389.1069&187&7670&10500&2639033785&32.3.11.19.29.47.617&17&534 \cr
 & &16.5.11.13.17.59.389&389&680& & &128.9.17.89.617&5553&1088 \cr
\noalign{\hrule}
 & &5.19.53.61.83.103&1221&3178& & &3.5.7.17.59.71.353&143&214 \cr
10483&2625697115&4.3.5.7.11.37.61.227&1577&1272&10501&2639509845&4.5.11.13.59.107.353&1207&558 \cr
 & &64.9.19.53.83.227&227&288& & &16.9.13.17.31.71.107&1209&856 \cr
\noalign{\hrule}
 & &81.25.37.101.347&46699&46726& & &81.25.7.241.773&1693&3718 \cr
10484&2625896475&4.3.17.41.61.67.347.383&817&20350&10502&2640703275&4.11.169.241.1693&2413&720 \cr
 & &16.25.11.19.37.41.43.67&14003&14104& & &128.9.5.11.13.19.127&2717&8128 \cr
\noalign{\hrule}
 & &3.5.23.79.173.557&2301&484& & &3.11.37.283.7643&1555&1558 \cr
10485&2626319055&8.9.121.13.59.173&557&730&10503&2640985149&4.5.19.37.41.311.7643&29&7614 \cr
 & &32.5.11.59.73.557&649&1168& & &16.81.19.29.47.311&36801&47272 \cr
\noalign{\hrule}
 & &9.17.19.23.163.241&3575&7672& & &7.11.73.641.733&3159&3892 \cr
10486&2626500663&16.3.25.7.11.13.19.137&4097&3712&10504&2641043713&8.243.49.13.73.139&3665&88 \cr
 & &4096.5.13.17.29.241&1885&2048& & &128.9.5.11.13.733&65&576 \cr
\noalign{\hrule}
 & &3.11.13.17.47.79.97&35&834& & &243.7.11.17.361.23&1135&2836 \cr
10487&2626654173&4.9.5.7.13.97.139&341&244&10505&2641076361&8.5.17.23.227.709&623&1332 \cr
 & &32.7.11.31.61.139&1891&15568& & &64.9.7.37.89.227&3293&7264 \cr
\noalign{\hrule}
 & &3.5.13.23.31.41.461&1757&548& & &5.7.11.13.17.37.839&4509&1364 \cr
10488&2627900535&8.7.23.41.137.251&407&1350&10506&2641293655&8.27.121.13.31.167&4541&6712 \cr
 & &32.27.25.11.37.137&13563&2960& & &128.9.19.239.839&2151&1216 \cr
\noalign{\hrule}
 & &9.25.19.4489.137&27233&58058& & &9.11.13.31.239.277&55821&55810 \cr
10489&2629095075&4.7.11.13.29.113.241&945&2188&10507&2641301091&4.27.5.23.239.809.5581&6017&436 \cr
 & &32.27.5.49.29.547&4263&8752& & &32.5.11.23.109.547.809&298115&297712 \cr
\noalign{\hrule}
 & &3.5.121.23.29.41.53&83741&66196& & &9.5.49.59.79.257&2227&428 \cr
10490&2630644665&8.49.13.19.67.1709&37575&5104&10508&2641318785&8.7.17.79.107.131&297&1046 \cr
 & &256.9.25.11.29.167&835&384& & &32.27.11.131.523&17259&2096 \cr
\noalign{\hrule}
 & &3.7.11.169.79.853&1445&4526& & &7.121.53.83.709&2025&2938 \cr
10491&2630721093&4.5.11.13.289.31.73&1179&1252&10509&2641700677&4.81.25.11.13.53.113&26233&10492 \cr
 & &32.9.5.17.31.131.313&69037&75120& & &32.3.37.43.61.709&2623&1776 \cr
\noalign{\hrule}
 & &3.5.11.29.31.113.157&1631&1066& & &3.11.29.31.139.641&161&480 \cr
10492&2631611235&4.7.11.13.41.157.233&2599&558&10510&2643300033&64.9.5.7.23.31.139&2771&3484 \cr
 & &16.9.23.31.113.233&233&552& & &512.7.13.17.67.163&76447&56576 \cr
\noalign{\hrule}
 & &3.49.59.223.1361&33233&33456& & &11.17.23.31.79.251&2881&2892 \cr
10493&2632281519&32.9.17.41.59.167.199&5365&111034&10511&2643820399&8.3.17.31.43.67.79.241&859&1590 \cr
 & &128.5.49.11.29.37.103&19055&20416& & &32.9.5.53.67.241.859&726615&728432 \cr
\noalign{\hrule}
 & &5.17.619.50053&23479&26574& & &31.113569.751&68425&45144 \cr
10494&2633538595&4.3.17.43.53.103.443&3817&3714&10512&2643999889&16.27.25.7.11.17.19.23&337&62 \cr
 & &16.9.11.43.53.347.619&25069&24984& & &64.9.17.23.31.337&207&544 \cr
\noalign{\hrule}
}%
}
$$
\eject
\vglue -23 pt
\noindent\hskip 1 in\hbox to 6.5 in{\ 10513 -- 10548 \hfill\fbd 2644199635 -- 2668377285\frb}
\vskip -9 pt
$$
\vbox{
\nointerlineskip
\halign{\strut
    \vrule \ \ \hfil \frb #\ 
   &\vrule \hfil \ \ \fbb #\frb\ 
   &\vrule \hfil \ \ \frb #\ \hfil
   &\vrule \hfil \ \ \frb #\ 
   &\vrule \hfil \ \ \frb #\ \ \vrule \hskip 2 pt
   &\vrule \ \ \hfil \frb #\ 
   &\vrule \hfil \ \ \fbb #\frb\ 
   &\vrule \hfil \ \ \frb #\ \hfil
   &\vrule \hfil \ \ \frb #\ 
   &\vrule \hfil \ \ \frb #\ \vrule \cr%
\noalign{\hrule}
 & &5.7.11.17.37.61.179&207&1102& & &3.11.13.37.41.61.67&3379&3392 \cr
10513&2644199635&4.9.19.23.29.37.61&2275&1202&10531&2659790991&128.11.31.41.53.67.109&1665&1082 \cr
 & &16.3.25.7.13.23.601&9015&2392& & &512.9.5.31.37.109.541&83855&83712 \cr
\noalign{\hrule}
 & &9.169.29.167.359&2701&530& & &3.5.19.53.293.601&6149&5270 \cr
10514&2644466877&4.5.13.29.37.53.73&359&330&10532&2659884765&4.25.11.13.17.31.43.53&3737&2412 \cr
 & &16.3.25.11.37.73.359&2701&2200& & &32.9.17.31.37.67.101&76849&82416 \cr
\noalign{\hrule}
 & &9.5.11.13.19.43.503&379&94& & &3.25.11.13.23.41.263&2771&7486 \cr
10515&2644469685&4.3.13.47.379.503&2701&3838&10533&2659896525&4.5.11.17.19.163.197&621&424 \cr
 & &16.19.37.47.73.101&7373&13912& & &64.27.17.23.53.163&8109&5216 \cr
\noalign{\hrule}
 & &9.5.49.131.9157&4251&4906& & &5.49.79.167.823&6149&2034 \cr
10516&2645045235&4.27.49.11.13.109.223&445&94&10534&2660170555&4.9.11.13.43.79.113&1169&300 \cr
 & &16.5.47.89.109.223&24307&33464& & &32.27.25.7.43.167&135&688 \cr
\noalign{\hrule}
 & &5.49.11.19.149.347&18009&54514& & &13.41.59.191.443&2755&8514 \cr
10517&2647452115&4.27.23.29.97.281&4511&3638&10535&2660825011&4.9.5.11.19.29.41.43&443&8 \cr
 & &16.3.13.17.23.107.347&4173&3128& & &64.3.19.43.443&2451&32 \cr
\noalign{\hrule}
 & &9.5.13.23.47.53.79&2771&6958& & &9.11.13.61.109.311&7327&8260 \cr
10518&2647795995&4.5.49.13.17.71.163&223&132&10536&2661308793&8.3.5.7.17.59.61.431&56291&19996 \cr
 & &32.3.7.11.17.163.223&36349&20944& & &64.181.311.4999&4999&5792 \cr
\noalign{\hrule}
 & &25.11.23.31.59.229&190137&183038& & &27.25.23.37.41.113&1859&916 \cr
10519&2649169325&4.3.61.71.1039.1289&76199&2430&10537&2661311025&8.9.11.169.113.229&1765&296 \cr
 & &16.729.5.23.3313&3313&5832& & &128.5.11.13.37.353&3883&832 \cr
\noalign{\hrule}
 & &9.11.13.17.347.349&51509&1888& & &81.5.11.13.19.41.59&1649&4214 \cr
10520&2649612537&64.19.59.2711&1385&1326&10538&2661831315&4.3.49.17.43.59.97&12617&8446 \cr
 & &256.3.5.13.17.19.277&1385&2432& & &16.7.11.31.37.41.103&3193&2072 \cr
\noalign{\hrule}
 & &9.5.7.11.23.29.31.37&1349&884& & &25.53.1123.1789&33371&61446 \cr
10521&2650894785&8.3.13.17.19.23.37.71&2449&1598&10539&2661987275&4.3.49.11.13.17.19.151&1961&300 \cr
 & &32.13.289.31.47.79&13583&16432& & &32.9.25.7.13.37.53&333&1456 \cr
\noalign{\hrule}
 & &5.11.13.71.89.587&401&756& & &27.5.109.239.757&2809&4004 \cr
10522&2652115895&8.27.7.11.401.587&409&178&10540&2662281945&8.3.7.11.13.2809.109&3953&3626 \cr
 & &32.9.89.401.409&3609&6544& & &32.343.37.53.59.67&209509&203056 \cr
\noalign{\hrule}
 & &3.25.7.19.23.31.373&9379&1546& & &125.17.19.23.47.61&2017&858 \cr
10523&2652841275&4.31.83.113.773&1793&1710&10541&2662367875&4.3.11.13.17.47.2017&16287&5900 \cr
 & &16.9.5.11.19.163.773&5379&6184& & &32.9.25.59.61.89&531&1424 \cr
\noalign{\hrule}
 & &9.29.41.43.73.79&11225&7828& & &9.17.19.43.101.211&1315&584 \cr
10524&2653644681&8.25.19.41.103.449&989&3234&10542&2663896311&16.5.19.73.101.263&6185&1188 \cr
 & &32.3.5.49.11.19.23.43&5635&3344& & &128.27.25.11.1237&3711&17600 \cr
\noalign{\hrule}
 & &17.29.397.13567&1027&12540& & &9.25.7.23.29.43.59&481&186 \cr
10525&2655346807&8.3.5.11.13.17.19.79&3367&3348&10543&2665181925&4.27.5.7.13.31.37.43&14927&14368 \cr
 & &64.81.7.11.169.31.37&193347&200096& & &256.11.23.37.59.449&4939&4736 \cr
\noalign{\hrule}
 & &9.125.7.11.23.31.43&201&674& & &3.5.193.241.3821&103&3718 \cr
10526&2655835875&4.27.23.31.67.337&407&430&10544&2665892595&4.11.169.103.193&9765&7642 \cr
 & &16.5.11.37.43.67.337&2479&2696& & &16.9.5.7.31.3821&217&24 \cr
\noalign{\hrule}
 & &13.17.47.331.773&133125&122738& & &27.5.7.11.13.109.181&589&316 \cr
10527&2657648981&4.3.625.7.11.71.797&221&276&10545&2666078415&8.9.11.19.31.79.109&22789&4150 \cr
 & &32.9.125.13.17.23.797&18331&18000& & &32.25.13.83.1753&1753&6640 \cr
\noalign{\hrule}
 & &81.19.997.1733&15965&16962& & &3.7.23.31.41.43.101&3287&5610 \cr
10528&2659085739&4.243.5.11.31.103.257&8531&16498&10546&2666137299&4.9.5.11.17.19.43.173&17297&38582 \cr
 & &16.5.11.19.73.113.449&41245&39512& & &16.49.101.191.353&2471&1528 \cr
\noalign{\hrule}
 & &81.5.151.157.277&10879&52706& & &3.5.7.17.19.31.43.59&29&30 \cr
10529&2659569795&4.11.361.23.43.73&675&314&10547&2667313005&4.9.25.7.17.19.29.31.43&140261&2714 \cr
 & &16.27.25.11.73.157&365&88& & &16.11.23.41.59.311&3421&7544 \cr
\noalign{\hrule}
 & &3.5.19.23.47.89.97&231&254& & &9.5.7.37.283.809&145&138 \cr
10530&2659697805&4.9.7.11.19.47.89.127&55775&45952&10548&2668377285&4.27.25.23.29.37.809&21791&3184 \cr
 & &1024.25.7.23.97.359&2513&2560& & &128.7.11.29.199.283&2189&1856 \cr
\noalign{\hrule}
}%
}
$$
\eject
\vglue -23 pt
\noindent\hskip 1 in\hbox to 6.5 in{\ 10549 -- 10584 \hfill\fbd 2669291625 -- 2689504425\frb}
\vskip -9 pt
$$
\vbox{
\nointerlineskip
\halign{\strut
    \vrule \ \ \hfil \frb #\ 
   &\vrule \hfil \ \ \fbb #\frb\ 
   &\vrule \hfil \ \ \frb #\ \hfil
   &\vrule \hfil \ \ \frb #\ 
   &\vrule \hfil \ \ \frb #\ \ \vrule \hskip 2 pt
   &\vrule \ \ \hfil \frb #\ 
   &\vrule \hfil \ \ \fbb #\frb\ 
   &\vrule \hfil \ \ \frb #\ \hfil
   &\vrule \hfil \ \ \frb #\ 
   &\vrule \hfil \ \ \frb #\ \vrule \cr%
\noalign{\hrule}
 & &3.125.7.11.169.547&149&696& & &81.23.31.71.653&3535&2882 \cr
10549&2669291625&16.9.25.7.11.29.149&169&106&10567&2677602339&4.9.5.7.11.71.101.131&653&1832 \cr
 & &64.169.29.53.149&4321&1696& & &64.11.101.229.653&2519&3232 \cr
\noalign{\hrule}
 & &11.17.43.61.5443&497&540& & &27.5.169.19.37.167&341&172 \cr
10550&2669796943&8.27.5.7.11.71.5443&3893&1550&10568&2678503815&8.5.11.31.37.43.167&243&6422 \cr
 & &32.9.125.7.17.31.229&28625&31248& & &32.243.11.169.19&9&176 \cr
\noalign{\hrule}
 & &27.5.11.13.19.29.251&73&62& & &31.53.1021.1597&2303&51810 \cr
10551&2669900805&4.13.19.29.31.73.251&1815&63812&10569&2678972291&4.3.5.49.11.47.157&901&744 \cr
 & &32.3.5.7.121.43.53&4081&688& & &64.9.7.11.17.31.53&1309&288 \cr
\noalign{\hrule}
 & &25.7.29.37.59.241&87363&87362& & &5.7.11.17.19.29.743&83&468 \cr
10552&2669972725&4.9.7.121.17.361.37.59.571&212875&1033618&10570&2679477185&8.9.13.17.83.743&703&40 \cr
 & &16.3.125.13.19.29.71.131.251&263055&263048& & &128.3.5.19.37.83&111&5312 \cr
\noalign{\hrule}
 & &9.25.29.41.67.149&6529&3454& & &169.289.131.419&51865&3024 \cr
10553&2670702075&4.3.11.29.157.6529&395&76&10571&2680833649&32.27.5.7.11.23.41&221&262 \cr
 & &32.5.19.79.6529&6529&24016& & &128.9.5.11.13.17.131&99&320 \cr
\noalign{\hrule}
 & &7.19.29.31.89.251&4485&2794& & &3.5.11.17.29.61.541&10437&5252 \cr
10554&2671007213&4.3.5.7.11.13.23.31.127&13303&11652&10572&2684466345&8.9.49.11.13.71.101&1769&3082 \cr
 & &32.9.11.53.251.971&8739&9328& & &32.23.29.61.67.71&1633&1072 \cr
\noalign{\hrule}
 & &9.5.11.23.41.59.97&18421&12698& & &3.11.29.61.139.331&9039&560 \cr
10555&2671410555&4.3.5.7.169.109.907&733&902&10573&2685867393&32.9.5.7.11.23.131&377&278 \cr
 & &16.7.11.41.733.907&6349&5864& & &128.7.13.23.29.139&161&832 \cr
\noalign{\hrule}
 & &3.5.23.109.227.313&19669&14448& & &25.11.13.43.101.173&119&54 \cr
10556&2671872855&32.9.5.7.13.17.43.89&3059&6886&10574&2686036925&4.27.5.7.11.17.43.101&2509&2036 \cr
 & &128.49.11.19.23.313&931&704& & &32.3.7.13.17.193.509&25959&21616 \cr
\noalign{\hrule}
 & &9.5.13.19.37.73.89&85&6412& & &1331.13.29.53.101&6453&7990 \cr
10557&2671923735&8.25.7.13.17.229&1701&1276&10575&2686065811&4.27.5.121.17.47.239&145&944 \cr
 & &64.243.49.11.29&783&17248& & &128.3.25.29.59.239&5975&11328 \cr
\noalign{\hrule}
 & &9.121.19.29.61.73&54223&74914& & &25.13.29.61.4673&3099&1574 \cr
10558&2671973667&4.7.13.43.97.5351&2045&3306&10576&2686624525&4.3.13.29.787.1033&4697&18126 \cr
 & &16.3.5.7.19.29.43.409&2045&2408& & &16.27.7.11.19.53.61&7049&2376 \cr
\noalign{\hrule}
 & &3.49.29.43.61.239&6479&5232& & &27.7.11.13.89.1117&15181&23000 \cr
10559&2672461911&32.9.11.19.31.61.109&5975&674&10577&2686835151&16.9.125.17.19.23.47&1573&9002 \cr
 & &128.25.11.239.337&3707&1600& & &64.5.7.121.13.643&643&1760 \cr
\noalign{\hrule}
 & &9.7.11.13.47.59.107&17543&17660& & &9.5.23.541.4799&213&328 \cr
10560&2673069399&8.5.11.53.59.331.883&3293&6420&10578&2687128065&16.27.41.71.4799&1441&3358 \cr
 & &64.3.25.37.89.107.331&29459&29600& & &64.11.23.41.73.131&9563&14432 \cr
\noalign{\hrule}
 & &11.23.29.37.43.229&335&654& & &3.25.7.11.361.1289&1899&7124 \cr
10561&2673155243&4.3.5.37.67.109.229&6283&18678&10579&2687274975&8.27.13.19.137.211&1115&1628 \cr
 & &16.9.11.61.103.283&17263&7416& & &64.5.11.37.137.223&8251&4384 \cr
\noalign{\hrule}
 & &121.13.101.113.149&111&10& & &5.49.11.19.73.719&20033&19314 \cr
10562&2674944701&4.3.5.13.37.113.149&3059&1122&10580&2687596835&4.9.5.13.19.23.29.37.67&3577&242 \cr
 & &16.9.5.7.11.17.19.23&1197&15640& & &16.3.49.121.13.37.73&407&312 \cr
\noalign{\hrule}
 & &17.71.73.97.313&69&5252& & &5.7.13.19.29.71.151&3901&2844 \cr
10563&2675138071&8.3.13.23.97.101&1065&1166&10581&2687808305&8.9.13.29.47.79.83&1221&142 \cr
 & &32.9.5.11.13.53.71&6201&880& & &32.27.11.37.71.79&2133&6512 \cr
\noalign{\hrule}
 & &29.41.43.199.263&605&51732& & &9.5.11.41.139.953&251&702 \cr
10564&2675833799&8.27.5.121.479&1247&1148&10582&2688417765&4.243.5.13.139.251&1927&1232 \cr
 & &64.3.7.11.29.41.43&231&32& & &128.7.11.41.47.251&1757&3008 \cr
\noalign{\hrule}
 & &7.13.17.19.59.1543&671&450& & &5.11.19.23.53.2111&11001&12220 \cr
10565&2675850541&4.9.25.7.11.61.1543&2431&3974&10583&2689107905&8.3.25.13.361.47.193&23&9048 \cr
 & &16.3.5.121.13.17.1987&5961&4840& & &128.9.169.23.29&261&10816 \cr
\noalign{\hrule}
 & &5.11.13.23.31.59.89&171&1876& & &27.25.23.191.907&74893&54032 \cr
10566&2676933545&8.9.7.13.19.59.67&341&400&10584&2689504425&32.7.11.13.307.823&885&4876 \cr
 & &256.3.25.7.11.31.67&1407&640& & &256.3.5.11.23.53.59&3127&1408 \cr
\noalign{\hrule}
}%
}
$$
\eject
\vglue -23 pt
\noindent\hskip 1 in\hbox to 6.5 in{\ 10585 -- 10620 \hfill\fbd 2689845675 -- 2719466295\frb}
\vskip -9 pt
$$
\vbox{
\nointerlineskip
\halign{\strut
    \vrule \ \ \hfil \frb #\ 
   &\vrule \hfil \ \ \fbb #\frb\ 
   &\vrule \hfil \ \ \frb #\ \hfil
   &\vrule \hfil \ \ \frb #\ 
   &\vrule \hfil \ \ \frb #\ \ \vrule \hskip 2 pt
   &\vrule \ \ \hfil \frb #\ 
   &\vrule \hfil \ \ \fbb #\frb\ 
   &\vrule \hfil \ \ \frb #\ \hfil
   &\vrule \hfil \ \ \frb #\ 
   &\vrule \hfil \ \ \frb #\ \vrule \cr%
\noalign{\hrule}
 & &3.25.11.19.157.1093&5957&2678& & &27.5.11.13.239.587&49981&46874 \cr
10585&2689845675&4.5.7.13.19.23.37.103&1099&396&10603&2708356365&4.9.23.151.331.1019&247&3226 \cr
 & &32.9.49.11.103.157&309&784& & &16.13.19.1019.1613&19361&12904 \cr
\noalign{\hrule}
 & &9.19.53.317.937&2773&3250& & &3.5.49.11.13.19.23.59&43&204 \cr
10586&2691973827&4.125.13.47.59.937&119&1056&10604&2709922215&8.9.5.7.11.17.43.59&199&494 \cr
 & &256.3.5.7.11.13.17.59&35105&18304& & &32.13.17.19.43.199&3383&688 \cr
\noalign{\hrule}
 & &5.11.41.307.3889&10387&6498& & &9.7.23.29.251.257&36425&30652 \cr
10587&2692296365&4.9.13.17.361.41.47&185&62&10605&2710648647&8.25.7.31.47.79.97&87&242 \cr
 & &16.3.5.17.19.31.37.47&27683&15096& & &32.3.5.121.29.79.97&9559&7760 \cr
\noalign{\hrule}
 & &3.23.59.109.6073&1783&4290& & &9.139.1069.2027&3575&2506 \cr
10588&2694826947&4.9.5.11.13.59.1783&2219&436&10606&2710745613&4.3.25.7.11.13.139.179&2027&58 \cr
 & &32.7.11.13.109.317&1001&5072& & &16.5.7.13.29.2027&91&1160 \cr
\noalign{\hrule}
 & &3.5.7.23.31.137.263&901&5148& & &3.25.11.19.23.73.103&61&1326 \cr
10589&2697460815&8.27.5.7.11.13.17.53&6049&6236&10607&2710787475&4.9.5.13.17.61.103&421&506 \cr
 & &64.23.53.263.1559&1559&1696& & &16.11.13.23.61.421&793&3368 \cr
\noalign{\hrule}
 & &9.7.361.31.43.89&2875&5402& & &3.29.43.59.71.173&4873&7410 \cr
10590&2698161291&4.3.125.23.37.43.73&18601&23374&10608&2711091477&4.9.5.11.13.19.29.443&59&436 \cr
 & &16.5.11.13.19.29.31.89&377&440& & &32.19.59.109.443&2071&7088 \cr
\noalign{\hrule}
 & &9.31.37.73.3583&2365&1218& & &81.11.13.19.97.127&1943&3010 \cr
10591&2700073557&4.27.5.7.11.29.43.73&1961&1178&10609&2711128563&4.27.5.7.19.29.43.67&12319&11374 \cr
 & &16.5.7.11.19.31.37.53&1463&2120& & &16.121.47.67.97.127&517&536 \cr
\noalign{\hrule}
 & &9.7.43.71.101.139&5995&3874& & &3.13.19.29.257.491&2365&4018 \cr
10592&2700247221&4.3.5.11.13.43.109.149&22873&28684&10610&2711632443&4.5.49.11.41.43.257&1557&1270 \cr
 & &32.5.71.89.101.257&1285&1424& & &16.9.25.7.43.127.173&66675&59512 \cr
\noalign{\hrule}
 & &3.13.211.239.1373&581&792& & &27.5.7.17.23.41.179&2531&1586 \cr
10593&2700321663&16.27.7.11.13.83.239&15403&13730&10611&2711723805&4.13.17.41.61.2531&4807&37710 \cr
 & &64.5.11.73.211.1373&365&352& & &16.9.5.11.19.23.419&419&1672 \cr
\noalign{\hrule}
 & &3.11.13.29.31.47.149&1535&2898& & &27.37.71.167.229&59&170 \cr
10594&2700854013&4.27.5.7.23.149.307&2431&5858&10612&2712537747&4.9.5.17.59.71.167&1739&1100 \cr
 & &16.5.7.11.13.17.29.101&505&952& & &32.125.11.37.47.59&5875&10384 \cr
\noalign{\hrule}
 & &3.7.11.13.19.23.29.71&52207&4850& & &9.5.13.109.157.271&3053&1012 \cr
10595&2702048349&4.25.17.37.83.97&3933&4118&10613&2713009455&8.3.11.23.43.71.109&1363&1690 \cr
 & &16.9.5.17.19.23.29.71&51&40& & &32.5.11.169.23.29.47&7337&9776 \cr
\noalign{\hrule}
 & &625.49.11.71.113&703&78& & &3.5.169.19.103.547&1283&1452 \cr
10596&2702748125&4.3.49.13.19.37.113&375&262&10614&2713664265&8.9.121.19.103.1283&15863&1750 \cr
 & &16.9.125.19.37.131&4847&1368& & &32.125.7.11.29.547&725&1232 \cr
\noalign{\hrule}
 & &19.29.53.151.613&8385&9392& & &5.11.13.59.139.463&37&102 \cr
10597&2703117289&32.3.5.13.43.151.587&569&7062&10615&2714900045&4.3.11.17.37.59.463&93&556 \cr
 & &128.9.5.11.107.569&60883&31680& & &32.9.17.31.37.139&629&4464 \cr
\noalign{\hrule}
 & &27.5.169.109.1087&17&22& & &9.11.31.41.113.191&1163&52126 \cr
10598&2703189645&4.9.11.13.17.109.1087&4183&5600&10616&2715767307&4.67.389.1163&39155&38766 \cr
 & &256.25.7.11.17.47.89&61523&56960& & &16.3.5.7.13.41.71.191&455&568 \cr
\noalign{\hrule}
 & &9.25.11.169.29.223&5111&7564& & &3.13.89.397.1973&1565&408 \cr
10599&2704984425&8.3.19.29.31.61.269&119&1772&10617&2718768351&16.9.5.17.313.397&1607&1210 \cr
 & &64.7.17.269.443&52717&8608& & &64.25.121.17.1607&40175&65824 \cr
\noalign{\hrule}
 & &27.13.17.19.29.823&20251&26660& & &27.5.7.121.13.31.59&2185&356 \cr
10600&2705873391&8.9.5.7.11.31.43.263&5251&6584&10618&2718781065&8.9.25.13.19.23.89&11&236 \cr
 & &128.7.11.59.89.823&5251&4928& & &64.11.23.59.89&2047&32 \cr
\noalign{\hrule}
 & &3.5.49.19.29.41.163&683&970& & &3.7.29.71.227.277&25&4 \cr
10601&2706514755&4.25.7.97.163.683&7967&9108&10619&2718825081&8.25.71.227.277&31617&31262 \cr
 & &32.9.11.23.31.97.257&99231&94576& & &32.27.5.49.11.29.1171&8197&7920 \cr
\noalign{\hrule}
 & &9.25.13.71.13033&7979&5054& & &3.5.7.13.19.23.47.97&3333&2072 \cr
10602&2706628275&4.7.361.71.79.101&3315&4664&10620&2719466295&16.9.49.11.19.37.101&7199&9118 \cr
 & &64.3.5.7.11.13.17.19.53&7049&5984& & &64.11.23.47.97.313&313&352 \cr
\noalign{\hrule}
}%
}
$$
\eject
\vglue -23 pt
\noindent\hskip 1 in\hbox to 6.5 in{\ 10621 -- 10656 \hfill\fbd 2720158415 -- 2741682625\frb}
\vskip -9 pt
$$
\vbox{
\nointerlineskip
\halign{\strut
    \vrule \ \ \hfil \frb #\ 
   &\vrule \hfil \ \ \fbb #\frb\ 
   &\vrule \hfil \ \ \frb #\ \hfil
   &\vrule \hfil \ \ \frb #\ 
   &\vrule \hfil \ \ \frb #\ \ \vrule \hskip 2 pt
   &\vrule \ \ \hfil \frb #\ 
   &\vrule \hfil \ \ \fbb #\frb\ 
   &\vrule \hfil \ \ \frb #\ \hfil
   &\vrule \hfil \ \ \frb #\ 
   &\vrule \hfil \ \ \frb #\ \vrule \cr%
\noalign{\hrule}
 & &5.13.37.79.103.139&2103&1708& & &27.25.49.17.43.113&55141&52316 \cr
10621&2720158415&8.3.7.13.61.139.701&553&1254&10639&2732094225&8.9.11.29.41.67.823&24295&428 \cr
 & &32.9.49.11.19.61.79&10431&8624& & &64.5.11.43.107.113&107&352 \cr
\noalign{\hrule}
 & &13.19.31.47.7559&3051&4508& & &9.23.83.157.1013&3325&286 \cr
10622&2720325361&8.27.49.13.19.23.113&7985&10208&10640&2732483421&4.3.25.7.11.13.19.83&157&92 \cr
 & &512.3.5.7.11.29.1597&122969&111360& & &32.5.7.11.19.23.157&133&880 \cr
\noalign{\hrule}
 & &81.6859.59.83&9529&11770& & &27.5.121.19.23.383&413&1502 \cr
10623&2720670363&4.3.5.11.13.19.107.733&83&2116&10641&2734005285&4.3.7.19.23.59.751&4213&4964 \cr
 & &32.5.11.13.529.83&5819&1040& & &32.11.17.59.73.383&1241&944 \cr
\noalign{\hrule}
 & &27.5.13.41.59.641&8551&8756& & &29.73.79.83.197&3069&2990 \cr
10624&2721266145&8.11.13.17.59.199.503&29067&610&10642&2734590293&4.9.5.11.13.23.29.31.197&83&2084 \cr
 & &32.3.5.17.61.9689&9689&16592& & &32.3.5.13.31.83.521&6045&8336 \cr
\noalign{\hrule}
 & &27.17.107.157.353&473&490& & &9.25.7.11.17.37.251&377&118 \cr
10625&2721891573&4.3.5.49.11.43.157.353&52751&5494&10643&2735253675&4.5.13.17.29.59.251&22701&6386 \cr
 & &16.7.17.29.41.67.107&1943&2296& & &16.3.7.23.31.47.103&2369&11656 \cr
\noalign{\hrule}
 & &3.5.79.83.89.311&5611&946& & &9.49.11.17.41.809&1625&2434 \cr
10626&2722368045&4.11.31.43.89.181&711&622&10644&2735347923&4.125.49.13.17.1217&233&12 \cr
 & &16.9.11.79.181.311&181&264& & &32.3.25.233.1217&30425&3728 \cr
\noalign{\hrule}
 & &3.17.37.41.61.577&31993&32054& & &9.7.113.199.1931&70015&87394 \cr
10627&2723086299&4.11.13.17.23.31.41.47.107&5193&26800&10645&2735611011&4.5.11.19.37.67.1181&8927&56028 \cr
 & &128.9.25.11.47.67.577&12925&12864& & &32.3.7.23.29.79.113&667&1264 \cr
\noalign{\hrule}
 & &9.7.19.41.211.263&325&536& & &9.11.37.67.71.157&1855&9292 \cr
10628&2723429961&16.3.25.13.19.67.263&223&1012&10646&2735707887&8.3.5.7.11.23.53.101&13031&5254 \cr
 & &128.5.11.23.67.223&56419&21440& & &32.37.71.83.157&83&16 \cr
\noalign{\hrule}
 & &13.193.823.1319&595&1914& & &9.11.17.19.23.3721&15995&54704 \cr
10629&2723612333&4.3.5.7.11.17.29.823&1209&386&10647&2736687591&32.5.7.13.263.457&1481&1938 \cr
 & &16.9.7.13.17.31.193&1071&248& & &128.3.5.7.17.19.1481&1481&2240 \cr
\noalign{\hrule}
 & &361.23.41.53.151&52245&2266& & &49.11.13.17.127.181&15111&24890 \cr
10630&2724405269&4.243.5.11.43.103&3599&3496&10648&2738188453&4.9.5.7.19.23.73.131&381&1298 \cr
 & &64.81.19.23.59.61&3599&2592& & &16.27.5.11.19.59.127&1593&760 \cr
\noalign{\hrule}
 & &7.121.37.89.977&2727&1750& & &9.7.13.19.37.67.71&215&44 \cr
10631&2725020067&4.27.125.49.89.101&4477&472&10649&2738876049&8.5.11.13.43.67.71&5453&696 \cr
 & &64.3.25.121.37.59&177&800& & &128.3.5.7.19.29.41&145&2624 \cr
\noalign{\hrule}
 & &9.25.11.17.239.271&13547&1358& & &3.361.23.109.1009&493&516 \cr
10632&2725155675&4.3.5.7.19.23.31.97&187&478&10650&2739516729&8.9.17.361.29.43.109&3205&44 \cr
 & &16.11.17.23.31.239&23&248& & &64.5.11.17.43.641&54485&15136 \cr
\noalign{\hrule}
 & &81.289.101.1153&739&638& & &7.17.29.53.71.211&495&2 \cr
10633&2726048277&4.11.17.29.739.1153&12623&60&10651&2740069843&4.9.5.11.53.211&343&290 \cr
 & &32.3.5.13.29.971&971&30160& & &16.3.25.343.11.29&3675&88 \cr
\noalign{\hrule}
 & &9.5.49.11.23.67.73&3265&1586& & &3.5.17.127.211.401&3091&3726 \cr
10634&2728517715&4.25.13.61.67.653&511&1164&10652&2740127235&4.243.11.23.211.281&60833&7450 \cr
 & &32.3.7.13.61.73.97&793&1552& & &16.25.127.149.479&2395&1192 \cr
\noalign{\hrule}
 & &3.5.49.23.31.41.127&3127&3096& & &27.5.11.169.61.179&589&82 \cr
10635&2728754385&16.27.5.23.41.43.53.59&2431&3374&10653&2740286835&4.9.5.19.31.41.179&517&338 \cr
 & &64.7.11.13.17.53.59.241&241723&242528& & &16.11.169.31.41.47&1271&376 \cr
\noalign{\hrule}
 & &3.49.1331.13.29.37&59&422& & &3.5.11.13.19.137.491&31&460 \cr
10636&2729219493&4.49.11.29.59.211&1591&14040&10654&2741466585&8.25.19.23.31.137&169&306 \cr
 & &64.27.5.13.37.43&45&1376& & &32.9.169.17.23.31&897&8432 \cr
\noalign{\hrule}
 & &3.11.13.29.41.53.101&575&536& & &9.5.11.29.31.61.101&2419&650 \cr
10637&2730463593&16.25.23.29.41.53.67&605&2142&10655&2741675805&4.125.13.41.59.101&5929&6696 \cr
 & &64.9.125.7.121.17.23&20125&17952& & &64.27.49.121.31.41&1617&1312 \cr
\noalign{\hrule}
 & &5.7.23.41.83.997&12969&53846& & &125.11.31.131.491&5967&6308 \cr
10638&2731196755&4.9.11.13.19.109.131&1453&2870&10656&2741682625&8.27.5.13.17.19.83.131&8339&4106 \cr
 & &16.3.5.7.19.41.1453&1453&456& & &32.9.13.31.269.2053&31473&32848 \cr
\noalign{\hrule}
}%
}
$$
\eject
\vglue -23 pt
\noindent\hskip 1 in\hbox to 6.5 in{\ 10657 -- 10692 \hfill\fbd 2741925585 -- 2770372935\frb}
\vskip -9 pt
$$
\vbox{
\nointerlineskip
\halign{\strut
    \vrule \ \ \hfil \frb #\ 
   &\vrule \hfil \ \ \fbb #\frb\ 
   &\vrule \hfil \ \ \frb #\ \hfil
   &\vrule \hfil \ \ \frb #\ 
   &\vrule \hfil \ \ \frb #\ \ \vrule \hskip 2 pt
   &\vrule \ \ \hfil \frb #\ 
   &\vrule \hfil \ \ \fbb #\frb\ 
   &\vrule \hfil \ \ \frb #\ \hfil
   &\vrule \hfil \ \ \frb #\ 
   &\vrule \hfil \ \ \frb #\ \vrule \cr%
\noalign{\hrule}
 & &3.5.49.53.59.1193&112073&99088& & &9.13.23.59.97.179&1675&2442 \cr
10657&2741925585&32.11.13.37.233.563&8901&11930&10675&2756706147&4.27.25.11.37.67.97&2147&472 \cr
 & &128.9.5.11.23.43.1193&1419&1472& & &64.11.19.37.59.113&7733&3616 \cr
\noalign{\hrule}
 & &27.7.13.71.79.199&419&220& & &9.5.7.11.169.17.277&13373&4078 \cr
10658&2742481287&8.3.5.7.11.13.79.419&1349&796&10676&2757519765&4.17.43.311.2039&1385&654 \cr
 & &64.19.71.199.419&419&608& & &16.3.5.109.277.311&311&872 \cr
\noalign{\hrule}
 & &5.67.109.163.461&59861&28974& & &3.5.11.13.19.31.37.59&381&322 \cr
10659&2743846645&4.3.11.31.439.1931&17425&3816&10677&2758013115&4.9.5.7.11.13.23.31.127&2593&5662 \cr
 & &64.27.25.17.41.53&24327&6560& & &16.7.19.23.149.2593&23989&20744 \cr
\noalign{\hrule}
 & &3.5.169.29.107.349&4807&5314& & &3.25.4591.8011&22323&17732 \cr
10660&2745270645&4.5.11.19.23.107.2657&1183&1278&10678&2758387575&8.9.5.7.11.13.31.1063&4591&4976 \cr
 & &16.9.7.11.169.71.2657&18599&18744& & &256.13.31.311.4591&4043&3968 \cr
\noalign{\hrule}
 & &9.5.13.61.107.719&4697&6088& & &11.29.31.157.1777&52585&2502 \cr
10661&2745354105&16.3.7.11.3721.761&719&3002&10679&2758922221&4.9.5.13.139.809&2117&310 \cr
 & &64.7.11.19.79.719&1501&2464& & &16.3.25.29.31.73&25&1752 \cr
\noalign{\hrule}
 & &25.7.13.397.3041&1517&1524& & &7.11.13.23.313.383&30675&30988 \cr
10662&2746555175&8.3.25.13.37.41.127.397&14007&682&10680&2759974217&8.3.25.11.13.61.127.409&267&4232 \cr
 & &32.9.7.11.23.29.31.127&66033&62992& & &128.9.5.529.89.127&26289&28480 \cr
\noalign{\hrule}
 & &5.49.41.181.1511&18467&55572& & &9.25.49.11.13.17.103&31&1256 \cr
10663&2747217095&8.3.11.59.313.421&1511&1932&10681&2760582825&16.17.31.103.157&891&860 \cr
 & &64.9.7.23.59.1511&531&736& & &128.81.5.11.43.157&1413&2752 \cr
\noalign{\hrule}
 & &3.5.49.13.19.37.409&6605&1166& & &9.25.7.11.37.59.73&23&2678 \cr
10664&2747320485&4.25.11.13.53.1321&1127&2448&10682&2760894675&4.5.7.11.13.23.103&4659&4714 \cr
 & &128.9.49.17.23.53&901&4416& & &16.3.23.1553.2357&35719&18856 \cr
\noalign{\hrule}
 & &9.7.11.19.23.29.313&4615&8058& & &3.5.19.2591.3739&10643&8052 \cr
10665&2748887757&4.27.5.7.13.17.71.79&47&506&10683&2761008465&8.9.11.19.29.61.367&3197&7234 \cr
 & &16.5.11.13.23.47.71&3337&520& & &32.23.29.139.3617&83191&64496 \cr
\noalign{\hrule}
 & &25.7.17.19.127.383&54023&13002& & &9.11.37.47.61.263&3265&398 \cr
10666&2749432525&4.3.11.89.197.607&205&402&10684&2761978923&4.5.199.263.653&629&366 \cr
 & &16.9.5.11.41.67.89&24723&7832& & &16.3.17.37.61.653&653&136 \cr
\noalign{\hrule}
 & &125.11.29.101.683&89&594& & &9.5.11.13.31.61.227&449&222 \cr
10667&2750697125&4.27.25.121.29.89&1717&1792&10685&2762268795&4.27.5.13.31.37.449&1507&508 \cr
 & &2048.9.7.17.89.101&10591&9216& & &32.11.127.137.449&17399&7184 \cr
\noalign{\hrule}
 & &9.49.19.23.109.131&143&2650& & &3.11.13.19.61.67.83&6235&6026 \cr
10668&2751806043&4.3.25.11.13.53.131&3059&3884&10686&2764990371&4.5.13.23.29.43.83.131&19&396 \cr
 & &32.7.13.19.23.971&971&208& & &32.9.11.19.23.43.131&2967&2096 \cr
\noalign{\hrule}
 & &9.7.13.17.19.101.103&1369&550& & &5.7.43.101.131.139&111&806 \cr
10669&2751978411&4.25.11.17.1369.103&57&572&10687&2767859045&4.3.13.31.37.43.101&139&1452 \cr
 & &32.3.5.121.13.19.37&185&1936& & &32.9.121.31.139&121&4464 \cr
\noalign{\hrule}
 & &9.5.7.43.199.1021&403&618& & &25.19.61.83.1151&1551&26 \cr
10670&2752059555&4.27.7.13.31.103.199&23749&46100&10688&2768068675&4.3.11.13.47.1151&581&570 \cr
 & &32.25.11.17.127.461&25355&34544& & &16.9.5.7.13.19.47.83&423&728 \cr
\noalign{\hrule}
 & &5.7.19.31.53.2521&1287&1658& & &121.169.43.47.67&753&820 \cr
10671&2754431995&4.9.11.13.829.2521&30031&2300&10689&2768937743&8.3.5.13.41.43.47.251&5561&17358 \cr
 & &32.3.25.23.59.509&6785&24432& & &32.9.5.11.67.83.263&3735&4208 \cr
\noalign{\hrule}
 & &3.25.7.19.271.1019&29029&29054& & &7.11.37.41.151.157&3015&3422 \cr
10672&2754586275&4.49.11.13.29.73.199.271&27&3550&10690&2769190963&4.9.5.7.29.59.67.151&779&628 \cr
 & &16.27.25.11.29.71.199&18531&17512& & &32.3.5.19.29.41.59.157&2755&2832 \cr
\noalign{\hrule}
 & &27.5.7.361.41.197&685&398& & &9.5.7.11.19.23.31.59&3071&3484 \cr
10673&2755428165&4.9.25.137.197.199&8531&35794&10691&2769480945&8.3.11.13.31.37.67.83&295&46 \cr
 & &16.11.19.449.1627&4939&13016& & &32.5.13.23.37.59.67&871&592 \cr
\noalign{\hrule}
 & &3.19.29.79.21101&11377&9724& & &81.5.11.47.101.131&7693&488 \cr
10674&2755516287&8.11.13.17.31.79.367&1161&3610&10692&2770372935&16.49.47.61.157&243&86 \cr
 & &32.27.5.11.17.361.43&9405&11696& & &64.243.7.43.61&301&5856 \cr
\noalign{\hrule}
}%
}
$$
\eject
\vglue -23 pt
\noindent\hskip 1 in\hbox to 6.5 in{\ 10693 -- 10728 \hfill\fbd 2772820575 -- 2796989481\frb}
\vskip -9 pt
$$
\vbox{
\nointerlineskip
\halign{\strut
    \vrule \ \ \hfil \frb #\ 
   &\vrule \hfil \ \ \fbb #\frb\ 
   &\vrule \hfil \ \ \frb #\ \hfil
   &\vrule \hfil \ \ \frb #\ 
   &\vrule \hfil \ \ \frb #\ \ \vrule \hskip 2 pt
   &\vrule \ \ \hfil \frb #\ 
   &\vrule \hfil \ \ \fbb #\frb\ 
   &\vrule \hfil \ \ \frb #\ \hfil
   &\vrule \hfil \ \ \frb #\ 
   &\vrule \hfil \ \ \frb #\ \vrule \cr%
\noalign{\hrule}
 & &9.25.343.19.31.61&7811&13112& & &9.49.11.43.67.199&169&370 \cr
10693&2772820575&16.25.11.73.107.149&2451&1274&10711&2781170469&4.3.5.169.37.43.199&7123&29212 \cr
 & &64.3.49.13.19.43.73&949&1376& & &32.17.67.109.419&1853&6704 \cr
\noalign{\hrule}
 & &3.5.11.29.41.79.179&265&186& & &11.17.23.137.4721&10449&15170 \cr
10694&2774252085&4.9.25.29.31.53.179&3367&8558&10712&2781787877&4.243.5.23.37.41.43&1727&136 \cr
 & &16.7.11.13.31.37.389&14911&21784& & &64.3.5.11.17.41.157&785&3936 \cr
\noalign{\hrule}
 & &9.11.17.23.229.313&36425&35252& & &9.11.13.103.139.151&17955&3638 \cr
10695&2774544993&8.3.25.7.11.31.47.1259&529&8284&10713&2782322829&4.243.5.7.17.19.107&611&604 \cr
 & &64.5.19.529.31.109&13547&17440& & &32.13.17.19.47.107.151&5029&5168 \cr
\noalign{\hrule}
 & &9.5.11.17.29.83.137&923&1484& & &125.13.17.131.769&447&322 \cr
10696&2774921985&8.3.5.7.13.53.71.137&29&166&10714&2782914875&4.3.7.13.17.23.131.149&517&2220 \cr
 & &32.7.29.53.71.83&371&1136& & &32.9.5.11.37.47.149&14751&27824 \cr
\noalign{\hrule}
 & &9.7.17.19.31.53.83&451&4850& & &9.5.11.17.31.47.227&1421&622 \cr
10697&2774972781&4.25.7.11.17.41.97&285&166&10715&2783168685&4.5.49.11.29.31.311&10387&498 \cr
 & &16.3.125.19.83.97&125&776& & &16.3.7.13.17.47.83&1079&56 \cr
\noalign{\hrule}
 & &3.11.19.23.199.967&413&214& & &9.5.11.17.23.73.197&14231&15910 \cr
10698&2775076293&4.7.23.59.107.967&11495&10746&10716&2783370645&4.25.7.11.19.37.43.107&4913&312 \cr
 & &16.27.5.121.19.59.199&495&472& & &64.3.7.13.4913.37&3757&8288 \cr
\noalign{\hrule}
 & &729.5.7.121.29.31&1349&884& & &9.11.37.661.1151&1501&2162 \cr
10699&2775496185&8.243.11.13.17.19.71&875&1798&10717&2786850693&4.19.23.47.79.1151&175&1326 \cr
 & &32.125.7.17.19.29.31&425&304& & &16.3.25.7.13.17.23.47&10387&32200 \cr
\noalign{\hrule}
 & &27.11.31.103.2927&40081&38948& & &9.49.11.13.193.229&8507&950 \cr
10700&2775735567&8.7.13.31.107.149.269&1775&30558&10718&2787195411&4.3.25.13.19.47.181&1351&1364 \cr
 & &32.3.25.11.13.71.463&11575&14768& & &32.5.7.11.19.31.47.193&1457&1520 \cr
\noalign{\hrule}
 & &11.19.47.233.1213&1305&1258& & &3.49.11.13.23.73.79&2085&1786 \cr
10701&2776264667&4.9.5.17.19.29.37.1213&117&20504&10719&2788246461&4.9.5.11.19.47.73.139&1027&224 \cr
 & &64.81.5.11.13.233&1053&160& & &256.5.7.13.19.47.79&893&640 \cr
\noalign{\hrule}
 & &81.11.43.233.311&157&76& & &121.19.43.89.317&4925&1098 \cr
10702&2776274919&8.11.19.43.157.311&23093&9720&10720&2789052541&4.9.25.121.61.197&159&38 \cr
 & &128.243.5.7.3299&3299&6720& & &16.27.25.19.53.61&1431&12200 \cr
\noalign{\hrule}
 & &9.5.11.13.59.71.103&23&332& & &25.101.211.5237&2871&2366 \cr
10703&2776490145&8.3.11.13.23.59.83&1123&824&10721&2790142675&4.9.5.7.11.169.29.211&5237&6368 \cr
 & &128.83.103.1123&1123&5312& & &256.3.7.13.199.5237&2587&2688 \cr
\noalign{\hrule}
 & &9.7.37.83.113.127&20995&6644& & &7.11.13.19.229.641&821&180 \cr
10704&2776531023&8.5.7.11.13.17.19.151&1163&498&10722&2791779991&8.9.5.19.229.821&8333&3982 \cr
 & &32.3.13.17.83.1163&1163&3536& & &32.3.11.13.181.641&181&48 \cr
\noalign{\hrule}
 & &5.7.121.13.61.827&3277&4104& & &49.11.13.43.73.127&589&1062 \cr
10705&2777359585&16.27.5.7.13.19.29.113&17&1712&10723&2793361571&4.9.49.19.31.59.73&635&3526 \cr
 & &512.9.17.29.107&52751&2304& & &16.3.5.31.41.43.127&205&744 \cr
\noalign{\hrule}
 & &5.121.23.31.47.137&371&234& & &27.7.11.19.173.409&2021&1612 \cr
10706&2777559235&4.9.7.13.23.31.47.53&137&850&10724&2794972257&8.9.11.13.19.31.43.47&875&1348 \cr
 & &16.3.25.13.17.53.137&663&2120& & &64.125.7.31.47.337&42125&46624 \cr
\noalign{\hrule}
 & &3.2401.181.2131&135&2266& & &25.7.11.13.61.1831&7353&1802 \cr
10707&2778276333&4.81.5.11.103.181&4081&4262&10725&2795067275&4.9.5.11.17.19.43.53&427&218 \cr
 & &16.5.7.121.53.2131&605&424& & &16.3.7.17.53.61.109&2703&872 \cr
\noalign{\hrule}
 & &81.25.7.43.47.97&6061&8086& & &5.49.11.19.79.691&92903&76392 \cr
10708&2778824475&4.11.13.19.29.97.311&4361&1548&10726&2795229745&16.9.61.1061.1523&231&1292 \cr
 & &32.9.49.11.13.43.89&1001&1424& & &128.27.7.11.17.19.61&1037&1728 \cr
\noalign{\hrule}
 & &11.13.17.29.79.499&1035&1256& & &3.13.17.37.83.1373&52877&670 \cr
10709&2779141079&16.9.5.11.23.157.499&1343&3838&10727&2795528229&4.5.121.19.23.67&1073&468 \cr
 & &64.3.17.19.23.79.101&2323&1824& & &32.9.13.19.29.37&19&1392 \cr
\noalign{\hrule}
 & &9.7.17.29.37.41.59&55&314& & &27.11.13.17.43.991&2245&2804 \cr
10710&2779873677&4.5.11.17.29.59.157&429&574&10728&2796989481&8.5.449.701.991&145&846 \cr
 & &16.3.7.121.13.41.157&2041&968& & &32.9.25.29.47.449&21103&11600 \cr
\noalign{\hrule}
}%
}
$$
\eject
\vglue -23 pt
\noindent\hskip 1 in\hbox to 6.5 in{\ 10729 -- 10764 \hfill\fbd 2797180925 -- 2820633543\frb}
\vskip -9 pt
$$
\vbox{
\nointerlineskip
\halign{\strut
    \vrule \ \ \hfil \frb #\ 
   &\vrule \hfil \ \ \fbb #\frb\ 
   &\vrule \hfil \ \ \frb #\ \hfil
   &\vrule \hfil \ \ \frb #\ 
   &\vrule \hfil \ \ \frb #\ \ \vrule \hskip 2 pt
   &\vrule \ \ \hfil \frb #\ 
   &\vrule \hfil \ \ \fbb #\frb\ 
   &\vrule \hfil \ \ \frb #\ \hfil
   &\vrule \hfil \ \ \frb #\ 
   &\vrule \hfil \ \ \frb #\ \vrule \cr%
\noalign{\hrule}
 & &25.49.11.41.61.83&117&422& & &3.7.223.353.1699&1313&3784 \cr
10729&2797180925&4.9.5.13.41.83.211&847&232&10747&2808615201&16.11.13.43.101.223&945&3398 \cr
 & &64.3.7.121.29.211&6119&1056& & &64.27.5.7.13.1699&65&288 \cr
\noalign{\hrule}
 & &9.11.17.43.67.577&553&1130& & &5.11.17.79.109.349&69527&68328 \cr
10730&2797713171&4.5.7.43.67.79.113&1717&1164&10748&2809898465&16.9.13.17.73.251.277&16511&3710 \cr
 & &32.3.5.17.97.101.113&10961&8080& & &64.3.5.7.11.13.19.53.79&3021&2912 \cr
\noalign{\hrule}
 & &27.49.13.37.53.83&397&1034& & &3.5.7.11.17.257.557&165579&165836 \cr
10731&2799360837&4.11.37.47.83.397&245&162&10749&2810730615&8.9.121.97.569.3769&22829&11092 \cr
 & &16.81.5.49.47.397&1191&1880& & &64.37.47.59.569.617&1072963&1074272 \cr
\noalign{\hrule}
 & &9.25.7.13.19.23.313&869&556& & &31.151.491.1223&7535&7686 \cr
10732&2800590975&8.3.7.11.13.23.79.139&125&148&10750&2810907733&4.9.5.7.11.61.137.1223&18881&12766 \cr
 & &64.125.11.37.79.139&25715&27808& & &16.3.13.61.79.239.491&18881&19032 \cr
\noalign{\hrule}
 & &19.259081.569&124135&134946& & &3.5.7.17.43.53.691&14671&15042 \cr
10733&2800924691&4.81.5.49.11.17.37.61&509&1138&10751&2810998365&4.9.5.289.23.109.863&4967&7568 \cr
 & &16.3.5.49.11.509.569&245&264& & &128.11.43.863.4967&54637&55232 \cr
\noalign{\hrule}
 & &11.13.19.29.961.37&2727&9890& & &27.227.317.1447&16445&55514 \cr
10734&2801642701&4.27.5.23.31.43.101&1859&464&10752&2811366171&4.5.11.13.23.41.677&2043&6758 \cr
 & &128.3.11.169.29.43&129&832& & &16.9.11.31.109.227&1199&248 \cr
\noalign{\hrule}
 & &3.25.47.61.83.157&4727&2652& & &81.25.7.11.13.19.73&709&5934 \cr
10735&2801990775&8.9.13.17.29.61.163&1111&1660&10753&2811483675&4.243.23.43.709&133&110 \cr
 & &64.5.11.13.29.83.101&4147&3232& & &16.5.7.11.19.43.709&709&344 \cr
\noalign{\hrule}
 & &9.5.7.11.41.109.181&529&452& & &9.29.31.43.59.137&735&598 \cr
10736&2802800385&8.5.529.41.113.181&48441&11336&10754&2812180779&4.27.5.49.13.23.29.59&1507&86 \cr
 & &128.3.13.67.109.241&3133&4288& & &16.5.11.13.23.43.137&715&184 \cr
\noalign{\hrule}
 & &3.11.13.23.29.97.101&3053&730& & &3.5.17.19.23.43.587&29&14 \cr
10737&2803342971&4.5.11.29.43.71.73&3519&2716&10755&2812730835&4.7.17.19.23.29.587&51937&42570 \cr
 & &32.9.7.17.23.71.97&497&816& & &16.9.5.11.43.167.311&3421&4008 \cr
\noalign{\hrule}
 & &3.5.7.19.29.47.1031&559&472& & &27.5.169.31.41.97&2537&1276 \cr
10738&2803479735&16.5.7.13.19.43.47.59&27&638&10756&2812792905&8.9.5.11.13.29.43.59&49&94 \cr
 & &64.27.11.29.43.59&5841&1376& & &32.49.29.43.47.59&66787&40592 \cr
\noalign{\hrule}
 & &7.17.23.31.173.191&1537&1710& & &5.7.121.13.17.31.97&237&764 \cr
10739&2803599421&4.9.5.7.19.23.29.31.53&27313&4152&10757&2814356545&8.3.5.11.79.97.191&5301&5204 \cr
 & &64.27.11.13.173.191&297&416& & &64.27.19.31.79.1301&40527&41632 \cr
\noalign{\hrule}
 & &25.121.13.37.41.47&181&144& & &11.13.23.53.67.241&765&776 \cr
10740&2803833175&32.9.121.41.47.181&13975&8288&10758&2814696599&16.9.5.13.17.53.97.241&121&15544 \cr
 & &2048.3.25.7.13.37.43&903&1024& & &256.3.121.17.29.67&957&2176 \cr
\noalign{\hrule}
 & &49.83.331.2083&37297&64770& & &3.5.23.29.347.811&2361&1694 \cr
10741&2804086691&4.3.5.13.17.19.127.151&77&2490&10759&2815577085&4.9.7.121.347.787&1633&5450 \cr
 & &16.9.25.7.11.13.83&117&2200& & &16.25.7.11.23.71.109&3815&6248 \cr
\noalign{\hrule}
 & &5.7.11.37.101.1949&55499&16614& & &9.7.13.31.197.563&8965&2858 \cr
10742&2804114005&4.9.13.19.23.71.127&1313&1100&10760&2815919379&4.3.5.11.13.163.1429&6751&394 \cr
 & &32.3.25.11.169.23.101&845&1104& & &16.11.43.157.197&157&3784 \cr
\noalign{\hrule}
 & &27.7.13.79.97.149&737&290& & &3.5.7.121.13.289.59&97&552 \cr
10743&2805370659&4.9.5.7.11.29.67.97&395&298&10761&2816228415&16.9.11.289.23.97&767&1834 \cr
 & &16.25.29.67.79.149&725&536& & &64.7.13.23.59.131&131&736 \cr
\noalign{\hrule}
 & &25.49.17.2209.61&821&2046& & &3.5.121.13.19.61.103&3265&3018 \cr
10744&2806147925&4.3.11.17.31.47.821&7869&16900&10762&2816700315&4.9.25.121.503.653&17081&96094 \cr
 & &32.9.25.169.43.61&1521&688& & &16.19.23.29.31.2089&20677&16712 \cr
\noalign{\hrule}
 & &3.5.7.17.827.1901&37999&32296& & &11.17.31.107.4547&2537&2010 \cr
10745&2806246695&16.7.11.13.37.79.367&1719&850&10763&2820408613&4.3.5.11.43.59.67.107&23329&27282 \cr
 & &64.9.25.13.17.37.191&7067&6240& & &16.9.5.41.569.4547&2845&2952 \cr
\noalign{\hrule}
 & &13.29.53.89.1579&211&1368& & &27.7.13.19.23.37.71&365&1364 \cr
10746&2807949911&16.9.19.29.53.211&2411&2200&10764&2820633543&8.5.11.23.31.71.73&91&162 \cr
 & &256.3.25.11.19.2411&60275&80256& & &32.81.5.7.13.31.73&465&1168 \cr
\noalign{\hrule}
}%
}
$$
\eject
\vglue -23 pt
\noindent\hskip 1 in\hbox to 6.5 in{\ 10765 -- 10800 \hfill\fbd 2822122611 -- 2842289043\frb}
\vskip -9 pt
$$
\vbox{
\nointerlineskip
\halign{\strut
    \vrule \ \ \hfil \frb #\ 
   &\vrule \hfil \ \ \fbb #\frb\ 
   &\vrule \hfil \ \ \frb #\ \hfil
   &\vrule \hfil \ \ \frb #\ 
   &\vrule \hfil \ \ \frb #\ \ \vrule \hskip 2 pt
   &\vrule \ \ \hfil \frb #\ 
   &\vrule \hfil \ \ \fbb #\frb\ 
   &\vrule \hfil \ \ \frb #\ \hfil
   &\vrule \hfil \ \ \frb #\ 
   &\vrule \hfil \ \ \frb #\ \vrule \cr%
\noalign{\hrule}
 & &9.49.11.19.67.457&2665&2494& & &3.11.13.19.53.79.83&333&580 \cr
10765&2822122611&4.5.7.13.29.41.43.457&12797&456&10783&2832643671&8.27.5.29.37.53.79&2047&86 \cr
 & &64.3.5.13.19.67.191&955&416& & &32.5.23.29.43.89&3827&53360 \cr
\noalign{\hrule}
 & &9.7.11.31.331.397&1385&1394& & &3.11.41.43.113.431&2963&1670 \cr
10766&2823016581&4.5.11.17.31.41.277.331&1651&147630&10784&2833491837&4.5.11.43.167.2963&40977&38014 \cr
 & &16.3.25.7.13.19.37.127&17575&13208& & &16.9.5.29.83.157.229&99615&104248 \cr
\noalign{\hrule}
 & &3.25.7.11.13.31.1213&3049&3016& & &9.17.59.89.3527&272441&267190 \cr
10767&2823045225&16.5.7.169.29.31.3049&23769&2426&10785&2833602381&4.5.7.11.13.19.347.1103&5261&9078 \cr
 & &64.9.19.29.139.1213&2641&2784& & &16.3.5.7.17.19.89.5261&5261&5320 \cr
\noalign{\hrule}
 & &9.49.121.191.277&1915&578& & &9.5.7.11.17.73.659&1147&1658 \cr
10768&2823170427&4.5.7.121.289.383&2369&312&10786&2833742835&4.3.31.37.659.829&23725&49424 \cr
 & &64.3.5.13.17.23.103&30797&2720& & &128.25.13.73.3089&3089&4160 \cr
\noalign{\hrule}
 & &9.7.17.23.29.59.67&1425&286& & &9.5.7.13.19.83.439&373&20128 \cr
10769&2823853221&4.27.25.7.11.13.19.23&2069&2278&10787&2834980785&64.7.17.37.373&3041&3300 \cr
 & &16.25.13.17.67.2069&2069&2600& & &512.3.25.11.3041&3041&14080 \cr
\noalign{\hrule}
 & &9.5.11.97.131.449&53417&5402& & &5.13.19.83.139.199&193&54 \cr
10770&2824194285&4.7.13.37.73.587&1057&1644&10788&2835390805&4.27.5.83.193.199&91&506 \cr
 & &32.3.49.13.137.151&6713&31408& & &16.9.7.11.13.23.193&14861&1656 \cr
\noalign{\hrule}
 & &9.25.17.19.59.659&3773&2158& & &81.13.17.157.1009&1193&184 \cr
10771&2825676675&4.5.343.11.13.59.83&1377&1318&10789&2835751113&16.13.23.157.1193&12699&14740 \cr
 & &16.81.7.13.17.83.659&747&728& & &128.9.5.11.17.67.83&5561&3520 \cr
\noalign{\hrule}
 & &9.5.11.13.53.8287&6767&1520& & &3.5.11.13.19.149.467&577&5494 \cr
10772&2826322785&32.25.13.19.67.101&8287&8262&10790&2835855165&4.5.19.41.67.577&2079&806 \cr
 & &128.243.17.101.8287&1717&1728& & &16.27.7.11.13.31.41&2583&248 \cr
\noalign{\hrule}
 & &5.13.17.19.31.43.101&1169&846& & &27.5.29.347.2089&53713&6868 \cr
10773&2826619835&4.9.7.43.47.101.167&13145&8398&10791&2837916945&8.11.17.19.101.257&41&60 \cr
 & &16.3.5.7.11.13.17.19.239&717&616& & &64.3.5.11.17.41.257&4369&14432 \cr
\noalign{\hrule}
 & &5.13.43.59.61.281&319&1086& & &9.25.13.19.223.229&3493&9218 \cr
10774&2826636605&4.3.11.29.43.61.181&17051&11802&10792&2838048525&4.3.7.11.13.419.499&229&190 \cr
 & &16.9.7.289.59.281&289&504& & &16.5.7.11.19.229.499&499&616 \cr
\noalign{\hrule}
 & &9.7.19.83.149.191&68225&42746& & &3.7.361.113.3313&6233&3706 \cr
10775&2827430109&4.25.11.29.67.2729&1337&1392&10793&2838091389&4.17.23.109.113.271&16565&14058 \cr
 & &128.3.5.7.841.67.191&4205&4288& & &16.9.5.11.17.71.3313&1207&1320 \cr
\noalign{\hrule}
 & &3.7.169.29.83.331&206867&199916& & &11.29.31.239.1201&100825&114036 \cr
10776&2827548633&8.23.37.41.53.5591&747&770&10794&2838528671&8.3.25.13.17.37.43.109&93&638 \cr
 & &32.9.5.7.11.53.83.5591&27955&27984& & &32.9.5.11.13.29.31.37&481&720 \cr
\noalign{\hrule}
 & &11.13.289.223.307&3567&190& & &3.11.41.43.59.827&393&434 \cr
10777&2829287747&4.3.5.19.29.41.223&121&102&10795&2838727947&4.9.7.11.31.43.59.131&86101&16540 \cr
 & &16.9.5.121.17.29.41&4059&1160& & &32.5.29.827.2969&2969&2320 \cr
\noalign{\hrule}
 & &5.11.17.361.83.101&1569&1984& & &27.7.17.19.181.257&541&4342 \cr
10778&2829555905&128.3.19.31.101.523&3403&6534&10796&2839723299&4.9.13.17.167.541&1265&724 \cr
 & &512.81.121.41.83&3321&2816& & &32.5.11.23.167.181&835&4048 \cr
\noalign{\hrule}
 & &5.19.41.647.1123&2509&1386& & &9.11.31.41.67.337&1807&1900 \cr
10779&2830032995&4.9.7.11.13.193.647&2179&6232&10797&2841092991&8.3.25.13.19.41.67.139&337&5362 \cr
 & &64.3.11.19.41.2179&2179&1056& & &32.7.13.19.337.383&2681&3952 \cr
\noalign{\hrule}
 & &3.5.41.43.103.1039&417&622& & &9.83.97.113.347&17203&16456 \cr
10780&2830064565&4.9.43.103.139.311&67337&53020&10798&2841189849&16.121.17.113.17203&1765&15438 \cr
 & &32.5.11.289.233.241&56153&50864& & &64.3.5.17.31.83.353&6001&4960 \cr
\noalign{\hrule}
 & &9.25.49.19.59.229&2431&1286& & &9.5.841.75079&335753&339958 \cr
10781&2830216725&4.5.7.11.13.17.19.643&2773&1728&10799&2841364755&4.11.43.59.67.131.233&7085&8526 \cr
 & &512.27.13.17.47.59&2397&3328& & &16.3.5.49.13.29.43.59.109&37583&37496 \cr
\noalign{\hrule}
 & &3.11.13.41.199.809&5249&3650& & &3.11.29.41.107.677&14039&13718 \cr
10782&2831670699&4.25.29.73.181.199&7821&21034&10800&2842289043&4.11.6859.29.101.139&40401&3940 \cr
 & &16.9.5.11.13.79.809&79&120& & &32.9.5.19.4489.197&67335&59888 \cr
\noalign{\hrule}
}%
}
$$
\eject
\vglue -23 pt
\noindent\hskip 1 in\hbox to 6.5 in{\ 10801 -- 10836 \hfill\fbd 2842829505 -- 2869509951\frb}
\vskip -9 pt
$$
\vbox{
\nointerlineskip
\halign{\strut
    \vrule \ \ \hfil \frb #\ 
   &\vrule \hfil \ \ \fbb #\frb\ 
   &\vrule \hfil \ \ \frb #\ \hfil
   &\vrule \hfil \ \ \frb #\ 
   &\vrule \hfil \ \ \frb #\ \ \vrule \hskip 2 pt
   &\vrule \ \ \hfil \frb #\ 
   &\vrule \hfil \ \ \fbb #\frb\ 
   &\vrule \hfil \ \ \frb #\ \hfil
   &\vrule \hfil \ \ \frb #\ 
   &\vrule \hfil \ \ \frb #\ \vrule \cr%
\noalign{\hrule}
 & &9.5.17.41.233.389&2021&76& & &5.7.11.13.383.1489&213&1702 \cr
10801&2842829505&8.17.19.41.43.47&17941&18672&10819&2854286435&4.3.7.11.13.23.37.71&645&1126 \cr
 & &256.3.7.11.233.389&77&128& & &16.9.5.43.71.563&5067&24424 \cr
\noalign{\hrule}
 & &3.13.97.569.1321&333&236& & &729.19.31.61.109&355&374 \cr
10802&2843488167&8.27.13.37.59.1321&485&836&10820&2854954269&4.5.11.17.31.61.71.109&19&1872 \cr
 & &64.5.11.19.37.59.97&7733&9440& & &128.9.5.11.13.19.71&781&4160 \cr
\noalign{\hrule}
 & &27.7.23.31.47.449&9845&2278& & &3.25.53.67.71.151&73&2 \cr
10803&2843776971&4.5.11.17.31.67.179&3143&2406&10821&2855270325&4.53.67.73.151&6993&3124 \cr
 & &16.3.5.7.17.401.449&401&680& & &32.27.7.11.37.71&3663&112 \cr
\noalign{\hrule}
 & &3.121.13.17.59.601&477&290& & &5.13.197.271.823&1269&2254 \cr
10804&2844627357&4.27.5.11.29.53.601&8911&7316&10822&2855937565&4.27.49.23.47.823&4855&2552 \cr
 & &32.7.19.31.53.59.67&14539&16112& & &64.3.5.11.23.29.971&32043&21344 \cr
\noalign{\hrule}
 & &27.11.13.173.4261&20299&18050& & &27.11.19.43.61.193&1253&2414 \cr
10805&2846147733&4.3.25.11.361.53.383&8411&2666&10823&2856706677&4.7.11.17.61.71.179&827&3870 \cr
 & &16.5.13.19.31.43.647&27821&23560& & &16.9.5.43.71.827&827&2840 \cr
\noalign{\hrule}
 & &27.5.7.1331.31.73&63229&33934& & &9.5.11.19.37.43.191&609&1426 \cr
10806&2846390085&4.361.47.53.1193&93&1100&10824&2858000805&4.27.7.23.29.31.191&487&296 \cr
 & &32.3.25.11.19.31.47&95&752& & &64.7.23.31.37.487&11201&6944 \cr
\noalign{\hrule}
 & &3.7.157.401.2153&2725&572& & &27.11.41.43.53.103&31333&30200 \cr
10807&2846474841&8.25.11.13.109.401&157&558&10825&2858392449&16.25.41.151.31333&12571&18762 \cr
 & &32.9.5.31.109.157&327&2480& & &64.3.25.13.53.59.967&24175&24544 \cr
\noalign{\hrule}
 & &3.5.11.19.31.83.353&87&502& & &9.7.29.43.79.461&535&74 \cr
10808&2847423315&4.9.11.29.251.353&3071&812&10826&2861113059&4.3.5.37.43.79.107&14993&4802 \cr
 & &32.7.841.37.83&5887&592& & &16.2401.11.29.47&3773&376 \cr
\noalign{\hrule}
 & &5.7.13.163.193.199&757&594& & &9.125.31.79.1039&47&78 \cr
10809&2848455155&4.27.5.11.13.199.757&1141&146&10827&2862574875&4.27.13.47.79.1039&5687&7820 \cr
 & &16.3.7.73.163.757&757&1752& & &32.5.121.17.23.2209&37553&44528 \cr
\noalign{\hrule}
 & &9.5.7.131.199.347&4141&3146& & &5.7.11.19.23.29.587&5723&1614 \cr
10810&2849472045&4.3.121.13.41.101.131&16915&22778&10828&2864034635&4.3.5.19.59.97.269&15849&10244 \cr
 & &16.5.7.11.17.199.1627&1627&1496& & &32.81.13.197.587&2561&1296 \cr
\noalign{\hrule}
 & &13.17.29.37.61.197&62415&32342& & &23.37.67.191.263&118305&105508 \cr
10811&2849627261&4.9.5.19.73.103.157&493&22&10829&2864134961&8.9.5.11.13.239.2029&2929&4958 \cr
 & &16.3.11.17.19.29.73&2409&152& & &32.3.5.13.29.37.67.101&2929&3120 \cr
\noalign{\hrule}
 & &3.5.11.289.529.113&71&94& & &5.11.31.41.43.953&251&702 \cr
10812&2850466245&4.289.23.47.71.113&3267&3380&10830&2864636995&4.27.5.13.31.43.251&4981&6314 \cr
 & &32.27.5.121.169.47.71&30033&29744& & &16.3.7.11.13.17.41.293&3809&2856 \cr
\noalign{\hrule}
 & &7.17.31.67.83.139&10725&54466& & &7.11.17.31.241.293&1165&3816 \cr
10813&2851519531&4.3.25.11.13.113.241&16949&17514&10831&2865404927&16.9.5.7.31.53.233&241&24 \cr
 & &16.27.5.7.17.139.997&997&1080& & &256.27.233.241&233&3456 \cr
\noalign{\hrule}
 & &27.13.17.29.53.311&32615&27328& & &23.31.79.151.337&317&396 \cr
10814&2852267769&128.9.5.7.11.61.593&1537&4502&10832&2866312049&8.9.11.151.317.337&395&58 \cr
 & &512.7.29.53.2251&2251&1792& & &32.3.5.11.29.79.317&3487&6960 \cr
\noalign{\hrule}
 & &5.47.101.263.457&205&252& & &9.11.89.113.2879&56753&53874 \cr
10815&2852733385&8.9.25.7.41.101.263&4631&24356&10833&2866456197&4.81.19.29.41.73.103&117407&5060 \cr
 & &64.3.11.421.6089&66979&40416& & &32.5.11.23.113.1039&1039&1840 \cr
\noalign{\hrule}
 & &9.125.7.13.29.961&363&1262& & &27.5.49.13.17.37.53&2603&3058 \cr
10816&2853088875&4.27.7.121.31.631&1525&1742&10834&2866815315&4.3.7.11.19.53.137.139&2479&2108 \cr
 & &16.25.13.61.67.631&4087&5048& & &32.17.19.31.37.67.137&9179&9424 \cr
\noalign{\hrule}
 & &27.25.11.19.113.179&5239&4964& & &3.49.11.13.17.71.113&205&1038 \cr
10817&2853524025&8.9.169.17.31.73.113&595&874&10835&2867075211&4.9.5.13.41.71.173&15029&22814 \cr
 & &32.5.7.13.289.19.23.73&26299&26864& & &16.7.11.17.19.61.113&61&152 \cr
\noalign{\hrule}
 & &13.37.2081.2851&385&2466& & &9.7.11.17.29.37.227&19581&26500 \cr
10818&2853739811&4.9.5.7.11.13.37.137&2851&2440&10836&2869509951&8.27.125.53.61.107&2107&782 \cr
 & &64.3.25.7.61.2851&1281&800& & &32.5.49.17.23.43.61&6923&4880 \cr
\noalign{\hrule}
}%
}
$$
\eject
\vglue -23 pt
\noindent\hskip 1 in\hbox to 6.5 in{\ 10837 -- 10872 \hfill\fbd 2872681845 -- 2907550681\frb}
\vskip -9 pt
$$
\vbox{
\nointerlineskip
\halign{\strut
    \vrule \ \ \hfil \frb #\ 
   &\vrule \hfil \ \ \fbb #\frb\ 
   &\vrule \hfil \ \ \frb #\ \hfil
   &\vrule \hfil \ \ \frb #\ 
   &\vrule \hfil \ \ \frb #\ \ \vrule \hskip 2 pt
   &\vrule \ \ \hfil \frb #\ 
   &\vrule \hfil \ \ \fbb #\frb\ 
   &\vrule \hfil \ \ \frb #\ \hfil
   &\vrule \hfil \ \ \frb #\ 
   &\vrule \hfil \ \ \frb #\ \vrule \cr%
\noalign{\hrule}
 & &3.5.11.17.61.103.163&12099&18382& & &9.5.121.17.19.31.53&365&3386 \cr
10837&2872681845&4.9.5.7.13.37.101.109&2231&4636&10855&2889601605&4.3.25.17.73.1693&209&1484 \cr
 & &32.19.23.61.97.101&9797&6992& & &32.7.11.19.53.73&511&16 \cr
\noalign{\hrule}
 & &25.11.17.131.4691&2673&2018& & &61.121801.389&72765&49036 \cr
10838&2872885675&4.243.5.121.17.1009&8479&18764&10856&2890215929&8.27.5.49.11.13.23.41&349&778 \cr
 & &32.9.61.139.4691&1251&976& & &32.9.5.41.349.389&205&144 \cr
\noalign{\hrule}
 & &3.125.7.11.169.19.31&4531&656& & &9.11.13.37.41.1481&225&226 \cr
10839&2874246375&32.11.13.23.41.197&817&1350&10857&2891473299&4.81.25.13.37.113.1481&17119&2134 \cr
 & &128.27.25.19.23.43&387&1472& & &16.5.11.17.19.53.97.113&87397&85880 \cr
\noalign{\hrule}
 & &3.5.17.43.263.997&497&500& & &5.7.13.53.277.433&51543&74492 \cr
10840&2875143615&8.625.7.17.43.71.263&10967&342&10858&2892377215&8.27.11.23.83.1693&2017&3710 \cr
 & &32.9.7.11.19.71.997&2343&2128& & &32.9.5.7.11.53.2017&2017&1584 \cr
\noalign{\hrule}
 & &3.11.13.19.29.43.283&5763&386& & &9.25.13.41.59.409&20801&15484 \cr
10841&2876496051&4.9.17.29.113.193&1771&1510&10859&2893910175&8.3.5.49.11.31.61.79&359&194 \cr
 & &16.5.7.11.23.113.151&17365&6328& & &32.7.31.61.97.359&77903&94672 \cr
\noalign{\hrule}
 & &11.17.137.167.673&3633&3770& & &5.7.11.361.59.353&1207&558 \cr
10842&2879345029&4.3.5.7.13.17.29.167.173&1089&3928&10860&2894640595&4.9.7.17.361.31.71&473&2054 \cr
 & &64.27.5.7.121.13.491&44681&47520& & &16.3.11.13.43.71.79&16827&4472 \cr
\noalign{\hrule}
 & &27.7.11.29.163.293&425&454& & &9.125.19.23.43.137&307&682 \cr
10843&2879437869&4.9.25.7.11.17.163.227&61&754&10861&2896162875&4.3.11.19.31.137.307&89&500 \cr
 & &16.5.13.17.29.61.227&13481&9080& & &32.125.11.89.307&979&4912 \cr
\noalign{\hrule}
 & &9.121.13.23.53.167&907&1160& & &81.11.29.127.883&833&50 \cr
10844&2881982961&16.3.5.11.29.167.907&2279&442&10862&2897611299&4.3.25.49.11.17.127&725&164 \cr
 & &64.5.13.17.29.43.53&2465&1376& & &32.625.7.29.41&287&10000 \cr
\noalign{\hrule}
 & &5.11.13.29.43.53.61&1387&2058& & &13.361.41.15061&7141&7920 \cr
10845&2882558965&4.3.343.19.29.43.73&689&732&10863&2897932193&32.9.5.11.13.19.37.193&9877&738 \cr
 & &32.9.7.13.19.53.61.73&1197&1168& & &128.81.7.17.41.83&9877&5184 \cr
\noalign{\hrule}
 & &3.13.3481.67.317&1309&992& & &9.5.31.73.149.191&667&3952 \cr
10846&2883385401&64.7.11.17.31.59.67&783&3170&10864&2898122265&32.13.19.23.29.191&1023&1460 \cr
 & &256.27.5.17.29.317&1305&2176& & &256.3.5.11.29.31.73&319&128 \cr
\noalign{\hrule}
 & &5.13.53.433.1933&5177&4488& & &5.1849.53.61.97&12831&12874 \cr
10847&2883427105&16.3.11.17.31.167.433&387&46&10865&2899241245&4.3.7.13.41.43.47.61.157&291&2332 \cr
 & &64.27.17.23.43.167&65297&37152& & &32.9.7.11.41.47.53.97&4653&4592 \cr
\noalign{\hrule}
 & &9.5.7.29.31.61.167&297&130& & &125.11.17.19.47.139&16329&9796 \cr
10848&2884805595&4.243.25.11.13.29.31&8881&2806&10866&2901468625&8.3.17.31.79.5443&3393&2050 \cr
 & &16.11.23.61.83.107&2461&7304& & &32.27.25.13.29.31.41&16523&12528 \cr
\noalign{\hrule}
 & &3.11.13.53.223.569&9039&2780& & &5.11.13.19.29.53.139&59619&80354 \cr
10849&2885029719&8.9.5.13.23.131.139&223&362&10867&2902340155&4.3.7.17.167.40177&20673&19504 \cr
 & &32.23.131.181.223&4163&2096& & &128.27.17.23.53.2297&39049&39744 \cr
\noalign{\hrule}
 & &9.343.31.53.569&5555&5078& & &3.25.11.19.23.83.97&14689&9464 \cr
10850&2885934429&4.5.11.101.569.2539&985&1554&10868&2902586775&16.7.169.23.37.397&1017&166 \cr
 & &16.3.25.7.11.37.101.197&41107&39400& & &64.9.83.113.397&1191&3616 \cr
\noalign{\hrule}
 & &9.5.11.13.17.23.31.37&1081&1064& & &3.5.7.13.19.23.31.157&11561&2726 \cr
10851&2885949495&16.3.7.19.529.31.37.47&24475&388&10869&2903189835&4.11.23.29.47.1051&3645&3692 \cr
 & &128.25.11.19.89.97&9215&5696& & &32.729.5.13.71.1051&17253&16816 \cr
\noalign{\hrule}
 & &3.5.11.19.59.67.233&6879&6868& & &3.121.169.19.47.53&245&366 \cr
10852&2887488615&8.9.5.17.19.67.101.2293&11461&4&10870&2903492163&4.9.5.49.13.19.53.61&913&94 \cr
 & &64.17.73.101.157&11461&54944& & &16.5.7.11.47.61.83&581&2440 \cr
\noalign{\hrule}
 & &3.7.121.13.19.43.107&1679&862& & &29.61.599.2741&2255&486 \cr
10853&2887711827&4.13.23.73.107.431&45815&55728&10871&2904448571&4.243.5.11.41.599&2741&3848 \cr
 & &128.81.5.49.11.17.43&945&1088& & &64.9.5.13.37.2741&585&1184 \cr
\noalign{\hrule}
 & &3.5.13.17.961.907&88913&87952& & &49.289.23.79.113&167315&158388 \cr
10854&2889443505&32.11.17.23.59.137.239&949&558&10872&2907550681&8.3.5.67.109.197.307&23331&10132 \cr
 & &128.9.13.31.59.73.239&14101&14016& & &64.9.5.7.11.17.101.149&15049&15840 \cr
\noalign{\hrule}
}%
}
$$
\eject
\vglue -23 pt
\noindent\hskip 1 in\hbox to 6.5 in{\ 10873 -- 10908 \hfill\fbd 2908252457 -- 2935038015\frb}
\vskip -9 pt
$$
\vbox{
\nointerlineskip
\halign{\strut
    \vrule \ \ \hfil \frb #\ 
   &\vrule \hfil \ \ \fbb #\frb\ 
   &\vrule \hfil \ \ \frb #\ \hfil
   &\vrule \hfil \ \ \frb #\ 
   &\vrule \hfil \ \ \frb #\ \ \vrule \hskip 2 pt
   &\vrule \ \ \hfil \frb #\ 
   &\vrule \hfil \ \ \fbb #\frb\ 
   &\vrule \hfil \ \ \frb #\ \hfil
   &\vrule \hfil \ \ \frb #\ 
   &\vrule \hfil \ \ \frb #\ \vrule \cr%
\noalign{\hrule}
 & &11.23.157.211.347&279&68& & &25.49.17.841.167&1067&1908 \cr
10873&2908252457&8.9.11.17.23.31.157&8651&3470&10891&2924808775&8.9.7.11.53.97.167&935&1102 \cr
 & &32.3.5.41.211.347&205&48& & &32.3.5.121.17.19.29.53&2299&2544 \cr
\noalign{\hrule}
 & &3.5.11.23.31.59.419&17861&21106& & &125.7.11.13.97.241&981&706 \cr
10874&2908302045&4.23.53.61.173.337&8029&66330&10892&2925047125&4.9.5.13.97.109.353&187&1448 \cr
 & &16.9.5.7.11.31.37.67&469&888& & &64.3.11.17.181.353&6001&17376 \cr
\noalign{\hrule}
 & &25.121.19.23.31.71&387&962& & &5.13.5021.8963&10229&55044 \cr
10875&2909556925&4.9.121.13.31.37.43&4715&15968&10893&2925209495&8.9.11.53.139.193&4745&2622 \cr
 & &256.3.5.23.41.499&1497&5248& & &32.27.5.13.19.23.73&1971&6992 \cr
\noalign{\hrule}
 & &7.11.19.529.53.71&221&150& & &27.25.7.13.19.23.109&4697&1972 \cr
10876&2912287301&4.3.25.11.13.17.19.529&11289&2296&10894&2925857025&8.49.11.17.23.29.61&317&810 \cr
 & &64.9.5.7.41.53.71&369&160& & &32.81.5.11.61.317&2013&5072 \cr
\noalign{\hrule}
 & &9.49.29.379.601&74855&92284& & &3.5.11.19.71.13151&1219&1124 \cr
10877&2913065631&8.5.11.1361.23071&127571&126210&10895&2927215335&8.23.53.281.13151&871&14022 \cr
 & &32.3.25.7.29.53.83.601&1325&1328& & &32.9.13.19.23.41.67&4623&8528 \cr
\noalign{\hrule}
 & &25.11.13.23.71.499&459&1174& & &3.25.11.23.37.43.97&329&1396 \cr
10878&2913149525&4.27.5.17.499.587&11521&10934&10896&2928354825&8.7.37.43.47.349&621&970 \cr
 & &16.3.7.11.17.41.71.281&5901&5576& & &32.27.5.7.23.47.97&329&144 \cr
\noalign{\hrule}
 & &27.343.19.59.281&799&1730& & &9.7.11.13.23.67.211&185&284 \cr
10879&2917224261&4.3.5.7.17.47.59.173&851&2090&10897&2929285359&8.5.13.23.37.71.211&10509&10720 \cr
 & &16.25.11.19.23.37.47&6325&13912& & &512.3.25.31.37.67.113&28675&28928 \cr
\noalign{\hrule}
 & &5.19.181.269.631&6129&5860& & &9.7.13.23.29.31.173&347&550 \cr
10880&2918662105&8.27.25.181.227.293&203&22&10898&2929662099&4.3.25.11.31.173.347&425&598 \cr
 & &32.3.7.11.29.227.293&93467&76272& & &16.625.13.17.23.347&5899&5000 \cr
\noalign{\hrule}
 & &27.7.11.19.107.691&9325&22454& & &81.5.49.11.31.433&3043&1720 \cr
10881&2920585437&4.25.7.103.109.373&1033&1578&10899&2930173785&16.3.25.17.31.43.179&4387&1162 \cr
 & &16.3.5.103.263.1033&27089&41320& & &64.7.17.41.83.107&8881&22304 \cr
\noalign{\hrule}
 & &343.13.19.841.41&165&368& & &9.7.11.17.47.67.79&8339&7540 \cr
10882&2921264801&32.3.5.49.11.19.23.29&929&492&10900&2930771151&8.3.5.7.11.13.29.31.269&901&3134 \cr
 & &256.9.5.11.41.929&10219&5760& & &32.13.17.31.53.1567&21359&25072 \cr
\noalign{\hrule}
 & &5.49.13.59.103.151&81&22& & &3.5.23.79.191.563&627&1190 \cr
10883&2922641995&4.81.5.49.11.13.151&59&696&10901&2930811915&4.9.25.7.11.17.19.191&12949&16274 \cr
 & &64.243.11.29.59&7047&352& & &16.11.23.79.103.563&103&88 \cr
\noalign{\hrule}
 & &5.11.13.17.191.1259&1161&940& & &81.5.7.11.23.61.67&65&4 \cr
10884&2922900695&8.27.25.43.47.1259&5917&4658&10902&2931421185&8.27.25.7.11.13.67&1867&58 \cr
 & &32.3.17.43.61.97.137&39867&41968& & &32.13.29.1867&54143&208 \cr
\noalign{\hrule}
 & &27.5.11.169.19.613&4249&2494& & &9.5.7.121.13.61.97&839&734 \cr
10885&2922989355&4.7.13.19.29.43.607&45975&28372&10903&2931843915&4.3.61.97.367.839&9295&8456 \cr
 & &32.3.25.41.173.613&865&656& & &64.5.7.11.169.151.367&4771&4832 \cr
\noalign{\hrule}
 & &27.103.191.5503&173&5330& & &9.13.41.53.83.139&35&88 \cr
10886&2923034013&4.5.13.41.103.173&649&690&10904&2933178417&16.3.5.7.11.13.83.139&503&1582 \cr
 & &16.3.25.11.23.59.173&14927&34600& & &64.49.11.113.503&56839&17248 \cr
\noalign{\hrule}
 & &27.5.11.13.199.761&217&1972& & &81.5.29.37.43.157&16511&24466 \cr
10887&2923529895&8.7.17.29.31.761&627&134&10905&2933748315&4.9.11.13.19.79.941&157&14 \cr
 & &32.3.7.11.19.31.67&469&9424& & &16.7.79.157.941&941&4424 \cr
\noalign{\hrule}
 & &25.29.43.191.491&58681&79794& & &27.7.11.19.59.1259&2291&3550 \cr
10888&2923622675&4.9.7.11.13.31.83.101&1275&1856&10906&2934173781&4.3.25.7.19.29.71.79&5471&5036 \cr
 & &512.27.25.11.13.17.29&3861&4352& & &32.5.71.1259.5471&5471&5680 \cr
\noalign{\hrule}
 & &7.17.79.353.881&795&1676& & &3.5.43.47.151.641&2587&4510 \cr
10889&2923645193&8.3.5.17.53.79.419&881&462&10907&2934219165&4.25.11.13.41.43.199&1359&9916 \cr
 & &32.9.5.7.11.53.881&495&848& & &32.9.13.37.67.151&2613&592 \cr
\noalign{\hrule}
 & &3.23.107.601.659&901&902& & &9.5.49.13.17.19.317&1331&5452 \cr
10890&2924103597&4.11.17.23.41.53.107.659&27045&26&10908&2935038015&8.3.5.7.1331.29.47&247&82 \cr
 & &16.9.5.13.17.53.601&901&1560& & &32.121.13.19.29.41&3509&656 \cr
\noalign{\hrule}
}%
}
$$
\eject
\vglue -23 pt
\noindent\hskip 1 in\hbox to 6.5 in{\ 10909 -- 10944 \hfill\fbd 2935238845 -- 2961521277\frb}
\vskip -9 pt
$$
\vbox{
\nointerlineskip
\halign{\strut
    \vrule \ \ \hfil \frb #\ 
   &\vrule \hfil \ \ \fbb #\frb\ 
   &\vrule \hfil \ \ \frb #\ \hfil
   &\vrule \hfil \ \ \frb #\ 
   &\vrule \hfil \ \ \frb #\ \ \vrule \hskip 2 pt
   &\vrule \ \ \hfil \frb #\ 
   &\vrule \hfil \ \ \fbb #\frb\ 
   &\vrule \hfil \ \ \frb #\ \hfil
   &\vrule \hfil \ \ \frb #\ 
   &\vrule \hfil \ \ \frb #\ \vrule \cr%
\noalign{\hrule}
 & &5.7.11.19.53.67.113&168051&169294& & &1331.19.31.53.71&56181&14362 \cr
10909&2935238845&4.3.7.13.31.47.139.1801&183207&19316&10927&2950037717&4.3.43.61.167.307&5247&4940 \cr
 & &32.9.11.173.353.439&61069&63216& & &32.27.5.11.13.19.43.53&1161&1040 \cr
\noalign{\hrule}
 & &9.5.23.41.43.1609&679&2288& & &27.7.41.47.8101&4121&3980 \cr
10910&2935950345&32.3.5.7.11.13.41.97&323&938&10928&2950408503&8.9.5.7.13.41.199.317&2585&2 \cr
 & &128.49.11.17.19.67&55811&13376& & &32.25.11.47.317&7925&176 \cr
\noalign{\hrule}
 & &121.31.41.71.269&5575&5454& & &27.7.11.13.23.47.101&1435&1292 \cr
10911&2937254309&4.27.25.31.71.101.223&13181&49036&10929&2950834887&8.5.49.17.19.23.41.47&117&1010 \cr
 & &32.3.5.49.13.23.41.269&1911&1840& & &32.9.25.13.17.41.101&425&656 \cr
\noalign{\hrule}
 & &3.5.7.11.169.101.149&14031&1018& & &27.11.17.31.109.173&1739&1640 \cr
10912&2937489555&4.27.5.509.1559&10769&2974&10930&2951478783&16.3.5.17.37.41.47.173&983&1612 \cr
 & &16.121.89.1487&16357&712& & &128.13.31.41.47.983&40303&39104 \cr
\noalign{\hrule}
 & &3.5.7.121.79.2927&527&1132& & &3.5.11.13.19.23.47.67&10003&6948 \cr
10913&2937815265&8.17.31.283.2927&2923&5850&10931&2951762385&8.27.7.19.193.1429&21125&17458 \cr
 & &32.9.25.13.17.37.79&1887&1040& & &32.125.49.169.29.43&18473&17200 \cr
\noalign{\hrule}
 & &13.529.37.11549&6039&5510& & &11.83.313.10333&148051&137718 \cr
10914&2938631501&4.9.5.11.13.19.29.37.61&101&23098&10932&2952851077&4.9.7.23.41.157.1093&3445&166 \cr
 & &16.3.5.101.11549&505&24& & &16.3.5.7.13.41.53.83&3445&6888 \cr
\noalign{\hrule}
 & &9.25.7.17.19.53.109&54769&54994& & &25.7.17.41.43.563&1089&2852 \cr
10915&2938904325&4.7.11.13.17.31.383.887&151&6360&10933&2952892775&8.9.25.121.17.23.31&247&178 \cr
 & &64.3.5.11.13.31.53.151&4433&4832& & &32.3.121.13.19.31.89&52421&75504 \cr
\noalign{\hrule}
 & &3.5.7.11.13.29.43.157&2891&850& & &11.19.961.47.313&22427&22740 \cr
10916&2939621685&4.125.343.11.17.59&21113&21762&10934&2954689639&8.3.5.11.19.41.379.547&4069&6324 \cr
 & &16.27.13.17.31.43.491&4419&4216& & &64.9.13.17.31.313.379&4927&4896 \cr
\noalign{\hrule}
 & &11.13.23.47.53.359&61605&29222& & &9.121.17.19.31.271&2369&70 \cr
10917&2941250741&4.9.5.19.1369.769&1669&2438&10935&2955026547&4.5.7.17.23.31.103&183&344 \cr
 & &16.3.5.19.23.53.1669&1669&2280& & &64.3.5.43.61.103&31415&1376 \cr
\noalign{\hrule}
 & &25.11.17.31.53.383&10543&4032& & &9.5.7.1873.5009&8411&8446 \cr
10918&2941832575&128.9.7.13.31.811&607&204&10936&2955284955&4.13.41.103.647.5009&65879&762 \cr
 & &1024.27.7.17.607&4249&13824& & &16.3.11.41.53.113.127&57277&47912 \cr
\noalign{\hrule}
 & &9.5.49.61.131.167&26057&24878& & &9.13.29.43.47.431&241&370 \cr
10919&2942565885&4.343.71.367.1777&143&126024&10937&2955476043&4.3.5.29.37.241.431&1183&110 \cr
 & &64.3.11.13.59.89&767&31328& & &16.25.7.11.169.241&3133&15400 \cr
\noalign{\hrule}
 & &121.31.37.127.167&4225&474& & &3.7.19.23.29.41.271&1565&2002 \cr
10920&2943533483&4.3.25.169.79.167&8763&4588&10938&2957003763&4.5.49.11.13.271.313&405&134 \cr
 & &32.9.23.31.37.127&23&144& & &16.81.25.13.67.313&45225&32552 \cr
\noalign{\hrule}
 & &27.13.29.41.7057&125&658& & &3.25.7.11.17.47.641&2031&2456 \cr
10921&2945161323&4.125.7.47.7057&4351&2706&10939&2957718225&16.9.11.47.307.677&1477&1900 \cr
 & &16.3.25.11.19.41.229&5225&1832& & &128.25.7.19.211.677&12863&13504 \cr
\noalign{\hrule}
 & &25.17.31.47.67.71&1311&1838& & &3.5.7.13.961.37.61&1653&14146 \cr
10922&2945653325&4.3.25.19.23.71.919&27001&4026&10940&2960653605&4.9.5.11.19.29.643&427&428 \cr
 & &16.9.11.13.31.61.67&1287&488& & &32.7.11.29.61.107.643&18647&18832 \cr
\noalign{\hrule}
 & &27.125.7.13.53.181&4127&2498& & &11.53.1223.4153&55251&9568 \cr
10923&2946250125&4.3.7.13.1249.4127&9911&38800&10941&2961126377&64.9.7.13.23.877&1465&1166 \cr
 & &128.25.11.17.53.97&1649&704& & &256.3.5.7.11.53.293&1465&2688 \cr
\noalign{\hrule}
 & &7.121.19.367.499&1769&2268& & &41.59.613.1997&46345&71478 \cr
10924&2947159369&8.81.49.11.19.29.61&355&184&10942&2961245459&4.9.5.11.13.361.23.31&211&4904 \cr
 & &128.9.5.23.29.61.71&99613&83520& & &64.3.23.211.613&211&2208 \cr
\noalign{\hrule}
 & &3.5.7.13.17.107.1187&1137&682& & &27.7.11.13.19.73.79&10117&5140 \cr
10925&2947243845&4.9.11.31.379.1187&69961&59278&10943&2961429471&8.3.5.13.67.151.257&209&662 \cr
 & &16.43.107.277.1627&11911&13016& & &32.5.11.19.257.331&1655&4112 \cr
\noalign{\hrule}
 & &81.5.11.29.37.617&14231&2428& & &3.11.13.47.191.769&289&322 \cr
10926&2949392655&8.3.5.7.19.107.607&2783&962&10944&2961521277&4.7.289.23.191.769&495&4888 \cr
 & &32.121.13.19.23.37&2717&368& & &64.9.5.11.13.289.47&289&480 \cr
\noalign{\hrule}
}%
}
$$
\eject
\vglue -23 pt
\noindent\hskip 1 in\hbox to 6.5 in{\ 10945 -- 10980 \hfill\fbd 2963533739 -- 2986807747\frb}
\vskip -9 pt
$$
\vbox{
\nointerlineskip
\halign{\strut
    \vrule \ \ \hfil \frb #\ 
   &\vrule \hfil \ \ \fbb #\frb\ 
   &\vrule \hfil \ \ \frb #\ \hfil
   &\vrule \hfil \ \ \frb #\ 
   &\vrule \hfil \ \ \frb #\ \ \vrule \hskip 2 pt
   &\vrule \ \ \hfil \frb #\ 
   &\vrule \hfil \ \ \fbb #\frb\ 
   &\vrule \hfil \ \ \frb #\ \hfil
   &\vrule \hfil \ \ \frb #\ 
   &\vrule \hfil \ \ \frb #\ \vrule \cr%
\noalign{\hrule}
 & &23.107.577.2087&1097&990& & &3.49.23.41.89.241&833&110 \cr
10945&2963533739&4.9.5.11.23.577.1097&233&1498&10963&2973281829&4.5.2401.11.17.89&7121&9522 \cr
 & &16.3.7.107.233.1097&7679&5592& & &16.9.5.529.7121&7121&2760 \cr
\noalign{\hrule}
 & &5.49.23.53.9923&4829&5094& & &3.13.37.43.173.277&1405&8844 \cr
10946&2963553565&4.9.49.11.23.283.439&95&1222&10964&2973450129&8.9.5.11.13.67.281&6055&3526 \cr
 & &16.3.5.11.13.19.47.283&34827&24904& & &32.25.7.41.43.173&175&656 \cr
\noalign{\hrule}
 & &5.7.841.47.2143&17739&57266& & &3.11.29.71.107.409&9729&2132 \cr
10947&2964722635&4.243.11.19.73.137&325&188&10965&2973564561&8.27.11.13.23.41.47&409&860 \cr
 & &32.9.25.11.13.47.73&4015&1872& & &64.5.13.23.43.409&2795&736 \cr
\noalign{\hrule}
 & &9.17.29.47.59.241&67&20& & &3.11.13.19.23.29.547&7157&5516 \cr
10948&2965216041&8.3.5.17.59.67.241&2491&1606&10966&2973884199&8.7.11.13.17.197.421&1035&3596 \cr
 & &32.11.47.53.67.73&4891&9328& & &64.9.5.7.17.23.29.31&1085&1632 \cr
\noalign{\hrule}
 & &9.49.17.19.109.191&543&220& & &25.49.19.29.4409&473&948 \cr
10949&2965520817&8.27.5.7.11.181.191&5293&5102&10967&2975964775&8.3.11.43.79.4409&2639&1770 \cr
 & &32.67.79.181.2551&201529&194032& & &32.9.5.7.13.29.43.59&2537&1872 \cr
\noalign{\hrule}
 & &9.7.11.169.19.31.43&26375&24494& & &27.5.11.43.149.313&41&14 \cr
10950&2966222259&4.125.31.37.211.331&3211&2064&10968&2978005635&4.7.41.43.149.313&22437&24200 \cr
 & &128.3.5.169.19.43.331&331&320& & &64.81.25.7.121.277&3047&3360 \cr
\noalign{\hrule}
 & &7.11.17.67.149.227&2061&472& & &5.289.19.23.53.89&351&86 \cr
10951&2966378569&16.9.11.59.67.229&8645&6698&10969&2978620405&4.27.13.289.43.89&145&1012 \cr
 & &64.3.5.7.13.17.19.197&3743&6240& & &32.9.5.11.23.29.43&387&5104 \cr
\noalign{\hrule}
 & &27.5.43.47.83.131&4199&1958& & &7.37.1451.7927&855&596 \cr
10952&2966534955&4.5.11.13.17.19.43.89&747&188&10970&2979037943&8.9.5.19.149.7927&5379&2548 \cr
 & &32.9.19.47.83.89&89&304& & &64.27.5.49.11.13.163&30807&22880 \cr
\noalign{\hrule}
 & &9.5.11.23.29.89.101&97&1208& & &9.25.37.113.3167&2579&5746 \cr
10953&2967853185&16.23.89.97.151&975&1072&10971&2979276075&4.169.17.113.2579&555&2024 \cr
 & &512.3.25.13.67.151&9815&17152& & &64.3.5.11.13.17.23.37&2431&736 \cr
\noalign{\hrule}
 & &81.19.53.59.617&43&574& & &9.61.503.10789&6149&4640 \cr
10954&2969283501&4.9.7.19.41.43.53&8437&7430&10972&2979349983&64.3.5.11.13.29.43.61&5669&8048 \cr
 & &16.5.7.11.13.59.743&3715&8008& & &2048.5.503.5669&5669&5120 \cr
\noalign{\hrule}
 & &5.7.11.43.83.2161&57447&18188& & &5.13.17.19.347.409&6835&936 \cr
10955&2969354465&8.9.13.491.4547&1537&3010&10973&2979673385&16.9.25.169.1367&7021&5654 \cr
 & &32.3.5.7.13.29.43.53&1537&624& & &64.3.7.11.17.59.257&15163&7392 \cr
\noalign{\hrule}
 & &27.5.11.29.53.1301&33041&35912& & &9.13.37.431.1597&4331&54758 \cr
10956&2969460945&16.3.5.19.37.47.4489&193&8&10974&2979681003&4.11.19.61.71.131&1445&2886 \cr
 & &256.19.47.67.193&59831&24704& & &16.3.5.13.289.19.37&1445&152 \cr
\noalign{\hrule}
 & &27.5.121.13.71.197&1745&172& & &11.23.61.113.1709&88669&69870 \cr
10957&2970209385&8.25.43.197.349&5041&9966&10975&2980374661&4.3.5.7.17.53.137.239&1045&1284 \cr
 & &32.3.11.5041.151&151&1136& & &32.9.25.7.11.19.53.107&50825&53424 \cr
\noalign{\hrule}
 & &17.29.43.269.521&163293&172150& & &5.7.11.23.41.43.191&4739&3474 \cr
10958&2971018651&4.3.25.11.13.53.79.313&269&126&10976&2981770715&4.9.49.41.193.677&18607&9150 \cr
 & &16.27.5.7.53.269.313&10955&11448& & &16.27.25.23.61.809&8235&6472 \cr
\noalign{\hrule}
 & &9.25.29.31.37.397&1045&146& & &5.11.13.17.19.37.349&1207983&1206748 \cr
10959&2971217475&4.3.125.11.19.37.73&4739&3364&10977&2982192785&8.3.7.23.29.41.61.101.103&72675&146 \cr
 & &32.7.19.841.677&3857&10832& & &32.27.25.17.19.23.73&1679&2160 \cr
\noalign{\hrule}
 & &5.7.289.23.53.241&63&328& & &27.5.49.47.53.181&8671&3784 \cr
10960&2971574585&16.9.49.17.41.241&5141&7150&10978&2982511665&16.49.11.13.23.29.43&519&470 \cr
 & &64.3.25.11.13.53.97&3783&1760& & &64.3.5.11.13.29.47.173&5017&4576 \cr
\noalign{\hrule}
 & &9.11.23.31.71.593&115769&102130& & &3.5.11.961.67.281&11283&11564 \cr
10961&2971924461&4.5.7.1459.115769&57155&58614&10979&2985303255&8.9.5.49.31.59.3761&829&1124 \cr
 & &16.3.25.49.23.71.9769&9769&9800& & &64.7.281.829.3761&26327&26528 \cr
\noalign{\hrule}
 & &9.5.121.19.59.487&23791&59486& & &7.11.31.47.79.337&44577&70526 \cr
10962&2972572515&4.49.37.607.643&26983&4524&10980&2986807747&4.27.13.127.179.197&385&206 \cr
 & &32.3.121.13.29.223&2899&464& & &16.9.5.7.11.13.103.127&6695&9144 \cr
\noalign{\hrule}
}%
}
$$
\eject
\vglue -23 pt
\noindent\hskip 1 in\hbox to 6.5 in{\ 10981 -- 11016 \hfill\fbd 2987321155 -- 3014816387\frb}
\vskip -9 pt
$$
\vbox{
\nointerlineskip
\halign{\strut
    \vrule \ \ \hfil \frb #\ 
   &\vrule \hfil \ \ \fbb #\frb\ 
   &\vrule \hfil \ \ \frb #\ \hfil
   &\vrule \hfil \ \ \frb #\ 
   &\vrule \hfil \ \ \frb #\ \ \vrule \hskip 2 pt
   &\vrule \ \ \hfil \frb #\ 
   &\vrule \hfil \ \ \fbb #\frb\ 
   &\vrule \hfil \ \ \frb #\ \hfil
   &\vrule \hfil \ \ \frb #\ 
   &\vrule \hfil \ \ \frb #\ \vrule \cr%
\noalign{\hrule}
 & &5.7.13.43.179.853&519&734& & &7.17.37.43.47.337&231&568 \cr
10981&2987321155&4.3.13.173.367.853&26201&37290&10999&2998782031&16.3.49.11.37.43.71&13395&15502 \cr
 & &16.9.5.7.11.19.113.197&19503&17176& & &64.9.5.19.23.47.337&855&736 \cr
\noalign{\hrule}
 & &27.5.49.23.41.479&10183&13288& & &9.125.7.121.47.67&21733&6608 \cr
10982&2987975655&16.11.17.41.151.599&525&74&11000&3000603375&32.49.59.103.211&1551&1340 \cr
 & &64.3.25.7.17.37.151&2567&5920& & &256.3.5.11.47.67.103&103&128 \cr
\noalign{\hrule}
 & &11.13.17.59.83.251&381&1030& & &9.7.11.13.19.47.373&61175&27232 \cr
10983&2988056357&4.3.5.13.103.127.251&1683&1580&11001&3000798801&64.25.23.37.2447&1649&798 \cr
 & &32.27.25.11.17.79.127&10033&10800& & &256.3.25.7.17.19.97&1649&3200 \cr
\noalign{\hrule}
 & &27.5.11.13.19.29.281&763&6824& & &729.11.521.719&2501&3230 \cr
10984&2989012455&16.5.7.13.109.853&3969&3116&11002&3003909381&4.5.17.19.41.61.719&891&1610 \cr
 & &128.81.343.19.41&1029&2624& & &16.81.25.7.11.17.19.23&2975&3496 \cr
\noalign{\hrule}
 & &5.7.13.139.151.313&10659&9686& & &3.25.343.11.13.19.43&10081&8366 \cr
10985&2989148435&4.3.11.17.19.29.151.167&2353&516&11003&3005477475&4.5.17.19.47.89.593&16091&4824 \cr
 & &32.9.13.17.29.43.181&21199&26064& & &64.9.17.67.16091&48273&36448 \cr
\noalign{\hrule}
 & &27.5.11.37.41.1327&13&1340& & &3.5.11.19.61.79.199&223&14 \cr
10986&2989392615&8.9.25.13.37.67&1327&3002&11004&3006405435&4.5.7.61.199.223&711&284 \cr
 & &32.19.79.1327&1501&16& & &32.9.71.79.223&213&3568 \cr
\noalign{\hrule}
 & &9.5.11.19.29.97.113&1091&752& & &3.13.23.59.139.409&3337&1980 \cr
10987&2989557945&32.3.5.11.29.47.1091&8023&8342&11005&3008725473&8.27.5.11.47.71.139&2045&11914 \cr
 & &128.43.47.71.97.113&3053&3008& & &32.25.7.23.37.409&175&592 \cr
\noalign{\hrule}
 & &9.49.11.17.19.23.83&5777&4366& & &11.19.961.71.211&3281&7290 \cr
10988&2991160557&4.11.19.37.53.59.109&7525&51888&11006&3008918869&4.729.5.17.71.193&2599&682 \cr
 & &128.3.25.7.23.43.47&1175&2752& & &16.27.5.11.23.31.113&3105&904 \cr
\noalign{\hrule}
 & &9.25.71.271.691&3997&2222& & &7.11.43.59.73.211&185&228 \cr
10989&2991494475&4.7.11.101.271.571&691&420&11007&3008960647&8.3.5.11.19.37.73.211&2307&4628 \cr
 & &32.3.5.49.571.691&571&784& & &64.9.13.37.89.769&89973&105376 \cr
\noalign{\hrule}
 & &9.5.11.17.19.97.193&53&1790& & &27.25.11.23.67.263&1387&18338 \cr
10990&2993207085&4.25.11.17.53.179&7081&2406&11008&3009226275&4.9.19.53.73.173&1315&1972 \cr
 & &16.3.73.97.401&73&3208& & &32.5.17.29.53.263&1537&272 \cr
\noalign{\hrule}
 & &27.25.11.13.19.23.71&157&1192& & &3.11.13.19.251.1471&39865&67814 \cr
10991&2994881175&16.3.5.11.13.149.157&1151&994&11009&3009520371&4.5.7.17.41.67.827&1521&4268 \cr
 & &64.7.71.149.1151&8057&4768& & &32.9.5.11.169.17.97&1649&3120 \cr
\noalign{\hrule}
 & &9.5.13.19.29.9293&2167&7126& & &27.7.11.757.1913&5959&11258 \cr
10992&2995459155&4.5.7.11.13.197.509&1253&1308&11010&3010685139&4.3.11.13.59.101.433&757&1190 \cr
 & &32.3.49.109.179.509&91111&85456& & &16.5.7.13.17.101.757&1105&808 \cr
\noalign{\hrule}
 & &9.25.17.23.79.431&319&112& & &125.41.137.4289&4707&418 \cr
10993&2995460775&32.25.7.11.17.29.79&3699&1724&11011&3011414125&4.9.11.19.137.523&1015&492 \cr
 & &256.27.7.137.431&959&384& & &32.27.5.7.19.29.41&3857&432 \cr
\noalign{\hrule}
 & &27.49.11.29.31.229&4129&2512& & &49.19.53.67.911&7047&55330 \cr
10994&2996040663&32.9.31.157.4129&1925&2204&11012&3011748691&4.243.5.11.29.503&265&238 \cr
 & &256.25.7.11.19.29.157&2983&3200& & &16.9.25.7.11.17.29.53&3825&2552 \cr
\noalign{\hrule}
 & &9.17.23.73.107.109&3115&608& & &5.11.43.83.103.149&3969&8398 \cr
10995&2996073081&64.3.5.7.19.89.107&253&146&11013&3012539365&4.81.5.49.11.13.17.19&149&786 \cr
 & &256.5.11.23.73.89&445&1408& & &16.243.19.131.149&4617&1048 \cr
\noalign{\hrule}
 & &9.121.13.17.59.211&3139&9310& & &3.17.107.283.1951&1885&66 \cr
10996&2996088381&4.3.5.49.13.19.43.73&1199&1940&11014&3012989781&4.9.5.11.13.29.283&1615&1498 \cr
 & &32.25.49.11.97.109&10573&19600& & &16.25.7.17.19.29.107&551&1400 \cr
\noalign{\hrule}
 & &27.7.19.41.47.433&29755&23504& & &27.7.11.17.19.4489&811&998 \cr
10997&2996298081&32.9.5.11.13.113.541&433&5518&11015&3014439813&4.7.19.67.499.811&4861&4050 \cr
 & &128.13.31.89.433&2759&832& & &16.81.25.499.4861&37425&38888 \cr
\noalign{\hrule}
 & &9.5.11.19.31.43.239&5401&6596& & &11.41.43.83.1873&18693&20566 \cr
10998&2996310735&8.121.17.19.97.491&325&1974&11016&3014816387&4.9.7.13.31.41.67.113&3577&830 \cr
 & &32.3.25.7.13.47.491&17185&9776& & &16.3.5.343.31.73.83&10633&8760 \cr
\noalign{\hrule}
}%
}
$$
\eject
\vglue -23 pt
\noindent\hskip 1 in\hbox to 6.5 in{\ 11017 -- 11052 \hfill\fbd 3015538981 -- 3044599445\frb}
\vskip -9 pt
$$
\vbox{
\nointerlineskip
\halign{\strut
    \vrule \ \ \hfil \frb #\ 
   &\vrule \hfil \ \ \fbb #\frb\ 
   &\vrule \hfil \ \ \frb #\ \hfil
   &\vrule \hfil \ \ \frb #\ 
   &\vrule \hfil \ \ \frb #\ \ \vrule \hskip 2 pt
   &\vrule \ \ \hfil \frb #\ 
   &\vrule \hfil \ \ \fbb #\frb\ 
   &\vrule \hfil \ \ \frb #\ \hfil
   &\vrule \hfil \ \ \frb #\ 
   &\vrule \hfil \ \ \frb #\ \vrule \cr%
\noalign{\hrule}
 & &7.13.31.653.1637&921875&920238& & &27.383.433.677&1507&1940 \cr
11017&3015538981&4.3.15625.11.59.73.191&4403&15672&11035&3031371081&8.3.5.11.97.137.677&5629&7660 \cr
 & &64.9.625.7.17.37.653&10625&10656& & &64.25.11.13.383.433&325&352 \cr
\noalign{\hrule}
 & &9.5.11.13.19.23.29.37&2881&3180& & &3.11.13.73.179.541&181&38 \cr
11018&3017377935&8.27.25.37.43.53.67&253&1178&11036&3032706963&4.19.179.181.541&1971&1430 \cr
 & &32.11.19.23.31.43.67&1333&1072& & &16.27.5.11.13.73.181&181&360 \cr
\noalign{\hrule}
 & &9.23.67.97.2243&64645&84832& & &3.11.13.23.31.47.211&2147&596 \cr
11019&3017492199&64.5.7.11.241.1847&2249&402&11037&3033382209&8.19.23.31.113.149&1645&15192 \cr
 & &256.3.5.7.13.67.173&2249&4480& & &128.9.5.7.47.211&21&320 \cr
\noalign{\hrule}
 & &27.25.7.47.107.127&23&212& & &9.5.11.17.29.31.401&463&494 \cr
11020&3017777175&8.5.23.53.107.127&913&1548&11038&3033599085&4.3.5.13.17.19.401.463&26071&8014 \cr
 & &64.9.11.43.53.83&3569&18656& & &16.19.841.31.4007&4007&4408 \cr
\noalign{\hrule}
 & &25.11.13.41.59.349&563&912& & &27.5.121.13.31.461&139&4 \cr
11021&3018125825&32.3.11.13.19.41.563&3213&2650&11039&3034765305&8.11.31.139.461&381&4690 \cr
 & &128.81.25.7.17.19.53&26163&23744& & &32.3.5.7.67.127&67&14224 \cr
\noalign{\hrule}
 & &3.5.7.29.53.59.317&8927&6708& & &9.7.17.19.29.37.139&199&60 \cr
11022&3018383655&8.9.13.29.43.79.113&36757&6974&11040&3034992303&8.27.5.17.19.29.199&1529&1036 \cr
 & &32.7.11.59.89.317&89&176& & &64.7.11.37.139.199&199&352 \cr
\noalign{\hrule}
 & &11.29.43.197.1117&525&722& & &27.343.11.17.1753&26009&31840 \cr
11023&3018412133&4.3.25.7.11.361.1117&173&1290&11041&3035857671&64.9.5.31.199.839&1315&476 \cr
 & &16.9.125.19.43.173&3287&9000& & &512.25.7.17.31.263&8153&6400 \cr
\noalign{\hrule}
 & &121.13.19.139.727&2107&534& & &9.5.7.11.47.103.181&109939&95134 \cr
11024&3020171011&4.3.49.43.89.727&2277&1550&11042&3036105765&4.13.17.29.223.3659&3725&66 \cr
 & &16.27.25.49.11.23.31&19251&9800& & &16.3.25.11.13.29.149&4321&520 \cr
\noalign{\hrule}
 & &9.31.37.227.1289&53&280& & &3.5.11.41.43.53.197&57323&57528 \cr
11025&3020540769&16.5.7.31.53.1289&629&660&11043&3037234695&16.27.7.17.19.43.47.431&407&106 \cr
 & &128.3.25.7.11.17.37.53&9911&11200& & &64.11.17.37.47.53.431&20257&20128 \cr
\noalign{\hrule}
 & &7.29.31.59.79.103&1881&2780& & &81.11.17.41.67.73&20711&15820 \cr
11026&3021162319&8.9.5.7.11.19.103.139&2449&3422&11044&3037443057&8.5.7.17.113.139.149&3159&796 \cr
 & &32.3.5.11.29.31.59.79&55&48& & &64.243.13.149.199&7761&4768 \cr
\noalign{\hrule}
 & &9.5.121.361.29.53&2219&3226& & &9.169.41.83.587&12407&11660 \cr
11027&3021196365&4.7.19.29.317.1613&20377&26400&11045&3038290281&8.5.11.169.19.53.653&189&20 \cr
 & &256.3.25.49.11.41.71&10045&9088& & &64.27.25.7.53.653&27825&20896 \cr
\noalign{\hrule}
 & &27.19.47.67.1873&837199&835390& & &9.25.13.401.2591&14297&1342 \cr
11028&3025713501&4.5.1331.17.37.139.601&603&2&11046&3039048675&4.3.5.11.17.841.61&5083&4168 \cr
 & &16.9.11.17.37.67.139&2363&3256& & &64.13.289.23.521&6647&16672 \cr
\noalign{\hrule}
 & &5.13.19.23.37.43.67&2673&3008& & &9.11.13.17.19.71.103&989&350 \cr
11029&3027887785&128.243.11.37.43.47&161&1430&11047&3040021413&4.25.7.11.17.19.23.43&9483&1442 \cr
 & &512.9.5.7.121.13.23&1089&1792& & &16.3.49.29.103.109&3161&392 \cr
\noalign{\hrule}
 & &11.29.3721.2551&1701&2020& & &3.25.19.529.37.109&3663&6388 \cr
11030&3028034449&8.243.5.7.101.2551&2639&88&11048&3040176225&8.27.11.1369.1597&29089&14030 \cr
 & &128.9.5.49.11.13.29&585&3136& & &32.5.19.23.61.1531&1531&976 \cr
\noalign{\hrule}
 & &3.5.13.17.463.1973&671&2644& & &5.11.13.19.151.1483&2349&5066 \cr
11031&3028249185&8.11.61.463.661&17757&10486&11049&3042129805&4.81.17.29.149.151&299&148 \cr
 & &32.9.49.107.1973&321&784& & &32.27.13.17.23.29.37&18241&9936 \cr
\noalign{\hrule}
 & &9.5.83.229.3541&8509&9196& & &81.5.7.43.109.229&76639&53204 \cr
11032&3028670415&8.3.121.19.67.83.127&229&20&11050&3042870705&8.47.173.283.443&55&228 \cr
 & &64.5.11.67.127.229&1397&2144& & &64.3.5.11.19.47.443&8417&16544 \cr
\noalign{\hrule}
 & &9.7.13.79.139.337&671&580& & &9.25.19.47.109.139&9959&56066 \cr
11033&3030788943&8.5.11.29.61.79.337&10081&10476&11051&3044214675&4.289.23.97.433&2565&4796 \cr
 & &64.27.11.17.29.97.593&172563&173536& & &32.27.5.11.17.19.109&187&48 \cr
\noalign{\hrule}
 & &27.49.13.23.79.97&8635&6016& & &5.17.109.373.881&20371&20286 \cr
11034&3031306551&256.5.11.47.79.157&2461&6174&11052&3044599445&4.9.49.13.23.881.1567&935&19436 \cr
 & &1024.9.343.23.107&749&512& & &32.3.5.7.11.17.23.43.113&20769&19888 \cr
\noalign{\hrule}
}%
}
$$
\eject
\vglue -23 pt
\noindent\hskip 1 in\hbox to 6.5 in{\ 11053 -- 11088 \hfill\fbd 3045106575 -- 3067050357\frb}
\vskip -9 pt
$$
\vbox{
\nointerlineskip
\halign{\strut
    \vrule \ \ \hfil \frb #\ 
   &\vrule \hfil \ \ \fbb #\frb\ 
   &\vrule \hfil \ \ \frb #\ \hfil
   &\vrule \hfil \ \ \frb #\ 
   &\vrule \hfil \ \ \frb #\ \ \vrule \hskip 2 pt
   &\vrule \ \ \hfil \frb #\ 
   &\vrule \hfil \ \ \fbb #\frb\ 
   &\vrule \hfil \ \ \frb #\ \hfil
   &\vrule \hfil \ \ \frb #\ 
   &\vrule \hfil \ \ \frb #\ \vrule \cr%
\noalign{\hrule}
 & &27.25.7.29.71.313&681&884& & &7.13.17.43.71.647&1661&2868 \cr
11053&3045106575&8.81.5.13.17.71.227&3443&2308&11071&3055775177&8.3.11.13.43.151.239&399&160 \cr
 & &64.11.13.17.313.577&6347&7072& & &512.9.5.7.11.19.151&14345&25344 \cr
\noalign{\hrule}
 & &3.25.7.11.37.53.269&3601&4946& & &27.25.11.17.43.563&703&28 \cr
11054&3046364475&4.5.13.53.277.2473&2959&486&11072&3055781025&8.7.11.19.37.563&11147&9684 \cr
 & &16.243.11.269.277&277&648& & &64.9.71.157.269&19099&5024 \cr
\noalign{\hrule}
 & &7.11.37.61.89.197&1395&772& & &9.11.17.29.31.43.47&133&394 \cr
11055&3047042537&8.9.5.31.37.61.193&1157&8298&11073&3057807357&4.7.11.19.43.47.197&6045&2426 \cr
 & &32.81.13.89.461&5993&1296& & &16.3.5.13.19.31.1213&6065&1976 \cr
\noalign{\hrule}
 & &9.11.13.23.113.911&65645&64628& & &27.5.13.19.293.313&749&4818 \cr
11056&3047215743&8.5.19.23.107.151.691&237473&284232&11074&3058036605&4.81.5.7.11.73.107&2689&3224 \cr
 & &128.3.13.17.61.229.911&3893&3904& & &64.7.11.13.31.2689&18823&10912 \cr
\noalign{\hrule}
 & &27.11.13.17.59.787&3845&3058& & &5.49.361.71.487&4293&21338 \cr
11057&3047722821&4.3.5.121.17.139.769&3989&6296&11075&3058162765&4.81.7.47.53.227&1207&836 \cr
 & &64.139.787.3989&3989&4448& & &32.9.11.17.19.47.71&1683&752 \cr
\noalign{\hrule}
 & &9.11.13.19.31.4021&18389&17800& & &9.5.19.43.193.431&193&22 \cr
11058&3048090903&16.25.7.11.13.37.71.89&439&1596&11076&3058222995&4.11.37249.431&20995&16254 \cr
 & &128.3.5.49.19.71.439&21511&22720& & &16.27.5.7.13.17.19.43&119&312 \cr
\noalign{\hrule}
 & &9.5.61.463.2399&5967&6028& & &27.5.73.263.1181&1117&854 \cr
11059&3048973065&8.243.11.13.17.137.463&14099&19192&11077&3060992565&4.5.7.61.1117.1181&957&6862 \cr
 & &128.13.17.23.613.2399&14099&14144& & &16.3.11.29.47.61.73&1769&4136 \cr
\noalign{\hrule}
 & &243.7.13.19.53.137&316609&315920& & &31.43.73.163.193&7029&7060 \cr
11060&3050687367&32.5.7.11.37.43.199.359&461&2052&11078&3061243831&8.9.5.11.43.71.163.353&20221&806 \cr
 & &256.27.5.11.19.199.461&25355&25472& & &32.3.13.31.71.73.277&3601&3408 \cr
\noalign{\hrule}
 & &11.13.17.37.107.317&3283&2106& & &5.17.31.41.43.659&2013&1282 \cr
11061&3050912293&4.81.49.169.37.67&45&214&11079&3061387795&4.3.11.31.41.61.641&4161&2890 \cr
 & &16.729.5.7.67.107&2345&5832& & &16.9.5.289.19.61.73&9333&11096 \cr
\noalign{\hrule}
 & &9.19.23.37.67.313&2545&4654& & &11.13.23.41.73.311&141417&151856 \cr
11062&3051720891&4.3.5.13.67.179.509&11&190&11080&3061476847&32.9.19.827.9491&3505&5986 \cr
 & &16.25.11.13.19.509&12725&1144& & &128.3.5.19.41.73.701&3505&3648 \cr
\noalign{\hrule}
 & &3.25.7.19.29.61.173&647&1122& & &5.11.23.59.89.461&2373&55388 \cr
11063&3052719075&4.9.7.11.17.173.647&4205&36556&11081&3062199415&8.3.7.61.113.227&3333&3560 \cr
 & &32.5.13.19.841.37&1073&208& & &128.9.5.7.11.89.101&909&448 \cr
\noalign{\hrule}
 & &5.7.11.13.19.163.197&11753&11556& & &3.11.19.43.97.1171&9039&13210 \cr
11064&3053595545&8.27.5.49.19.23.73.107&2119&86&11082&3062419107&4.9.5.11.23.131.1321&6665&5224 \cr
 & &32.3.13.23.43.73.163&1679&2064& & &64.25.23.31.43.653&20243&18400 \cr
\noalign{\hrule}
 & &27.7.31.37.73.193&8789&7052& & &9.13.67.103.3793&28105&21204 \cr
11065&3054255687&8.3.11.17.37.41.43.47&1825&62&11083&3062532681&8.81.5.7.11.19.31.73&16385&544 \cr
 & &32.25.11.31.47.73&1175&176& & &512.25.17.29.113&55709&6400 \cr
\noalign{\hrule}
 & &9.25.61.131.1699&1177&522& & &27.5.169.17.53.149&197&548 \cr
11066&3054759525&4.81.5.11.29.61.107&29237&3398&11084&3062890935&8.13.17.53.137.197&1341&440 \cr
 & &16.169.173.1699&173&1352& & &128.9.5.11.149.197&197&704 \cr
\noalign{\hrule}
 & &3.11.13.29.41.53.113&55453&4490& & &3.11.13.19.67.71.79&775&574 \cr
11067&3054875109&4.5.23.449.2411&981&1430&11085&3063170253&4.25.7.11.13.31.41.79&269&126 \cr
 & &16.9.25.11.13.23.109&2725&552& & &16.9.5.49.31.41.269&39543&50840 \cr
\noalign{\hrule}
 & &27.5.11.529.3889&47&3842& & &3.125.7.13.89.1009&407&602 \cr
11068&3055062285&4.9.17.23.47.113&2759&2552&11086&3064459125&4.25.49.11.37.43.89&2259&1034 \cr
 & &64.11.17.29.31.89&2581&16864& & &16.9.121.43.47.251&32379&45496 \cr
\noalign{\hrule}
 & &3.25.169.23.47.223&473&642& & &9.5.11.29.31.71.97&5187&4702 \cr
11069&3055473525&4.9.5.11.23.43.47.107&427&1508&11087&3064749435&4.27.7.13.19.71.2351&217&2134 \cr
 & &32.7.11.13.29.61.107&19459&11984& & &16.49.11.13.19.31.97&247&392 \cr
\noalign{\hrule}
 & &3.13.59.257.5167&10841&56330& & &9.7.19.67.167.229&9425&5918 \cr
11070&3055541619&4.5.37.43.131.293&2277&2570&11088&3067050357&4.3.25.11.13.19.29.269&2077&3034 \cr
 & &16.9.25.11.23.43.257&2967&2200& & &16.25.13.31.37.41.67&10075&12136 \cr
\noalign{\hrule}
}%
}
$$
\eject
\vglue -23 pt
\noindent\hskip 1 in\hbox to 6.5 in{\ 11089 -- 11124 \hfill\fbd 3067135175 -- 3088552467\frb}
\vskip -9 pt
$$
\vbox{
\nointerlineskip
\halign{\strut
    \vrule \ \ \hfil \frb #\ 
   &\vrule \hfil \ \ \fbb #\frb\ 
   &\vrule \hfil \ \ \frb #\ \hfil
   &\vrule \hfil \ \ \frb #\ 
   &\vrule \hfil \ \ \frb #\ \ \vrule \hskip 2 pt
   &\vrule \ \ \hfil \frb #\ 
   &\vrule \hfil \ \ \fbb #\frb\ 
   &\vrule \hfil \ \ \frb #\ \hfil
   &\vrule \hfil \ \ \frb #\ 
   &\vrule \hfil \ \ \frb #\ \vrule \cr%
\noalign{\hrule}
 & &25.13.41.43.53.101&27&532& & &3.5.49.13.19.31.547&67&522 \cr
11089&3067135175&8.27.5.7.19.41.53&1419&754&11107&3078458565&4.27.7.29.67.547&7619&7150 \cr
 & &32.81.11.13.29.43&2349&176& & &16.25.11.13.19.29.401&2005&2552 \cr
\noalign{\hrule}
 & &27.49.43.199.271&38099&43472& & &125.7.11.29.41.269&6803&7992 \cr
11090&3067966881&32.7.11.13.19.31.1229&473&8130&11108&3078469625&16.27.25.7.37.6803&4189&2614 \cr
 & &128.3.5.121.43.271&121&320& & &64.3.37.59.71.1307&145077&134048 \cr
\noalign{\hrule}
 & &3.5.11.13.17.19.43.103&4669&584& & &11.361.67.71.163&1043&306 \cr
11091&3068566215&16.7.11.13.23.29.73&515&486&11109&3079077661&4.9.7.17.19.149.163&53&110 \cr
 & &64.243.5.23.73.103&1863&2336& & &16.3.5.7.11.17.53.149&6307&17880 \cr
\noalign{\hrule}
 & &9.13.53.67.83.89&33441&27880& & &3.7.11.13.47.139.157&2635&10014 \cr
11092&3069054729&16.27.5.17.41.71.157&3575&664&11110&3080120043&4.9.5.11.17.31.1669&8507&6838 \cr
 & &256.125.11.13.17.83&1375&2176& & &16.13.17.47.181.263&3077&2104 \cr
\noalign{\hrule}
 & &7.11.13.23.151.883&1691&1782& & &9.7.13.19.37.53.101&56939&12650 \cr
11093&3069725659&4.81.121.19.89.883&8257&70330&11111&3082027221&4.25.11.23.97.587&2349&8806 \cr
 & &16.3.5.13.23.359.541&5385&4328& & &16.81.5.7.17.29.37&153&1160 \cr
\noalign{\hrule}
 & &25.361.29.37.317&117&434& & &125.11.71.131.241&36093&66218 \cr
11094&3069772525&4.9.25.7.13.19.31.37&5389&4686&11112&3082118875&4.3.53.113.227.293&5061&20590 \cr
 & &16.27.7.11.17.71.317&5049&3976& & &16.9.5.7.29.71.241&203&72 \cr
\noalign{\hrule}
 & &9.7.23.59.149.241&817&4244& & &81.37.911.1129&1691351&1692262 \cr
11095&3069896319&8.3.19.43.59.1061&67&110&11113&3082471443&4.11.13.23.61.97.151.487&1235&168 \cr
 & &32.5.11.19.67.1061&20159&58960& & &64.3.5.7.169.19.151.487&514759&513760 \cr
\noalign{\hrule}
 & &27.25.11.19.47.463&261&214& & &9.1331.361.23.31&25&146 \cr
11096&3069933075&4.243.11.29.107.463&5377&8050&11114&3083310747&4.25.11.19.23.31.73&7437&958 \cr
 & &16.25.7.19.23.107.283&6509&5992& & &16.3.5.37.67.479&32093&1480 \cr
\noalign{\hrule}
 & &3.11.13.17.31.37.367&2033&2400& & &3.625.23.43.1663&519&1144 \cr
11097&3069981057&64.9.25.17.19.37.107&31&734&11115&3083825625&16.9.11.13.23.43.173&5575&3326 \cr
 & &256.5.31.107.367&535&128& & &64.25.11.223.1663&223&352 \cr
\noalign{\hrule}
 & &9.5.59.61.67.283&5423&6838& & &9.5.7.17.23.79.317&307379&312356 \cr
11098&3070828755&4.3.11.13.17.29.59.263&1351&3782&11116&3084421095&8.11.31.61.229.5039&1275&154934 \cr
 & &16.7.31.61.193.263&8153&10808& & &32.3.25.13.17.59.101&3835&1616 \cr
\noalign{\hrule}
 & &9.11.83.101.3701&2305&1396& & &5.11.29.347.5573&55809&5494 \cr
11099&3071522817&8.5.11.83.349.461&2457&2108&11117&3084460445&4.81.13.41.53.67&1073&1100 \cr
 & &64.27.7.13.17.31.461&47957&44256& & &32.3.25.11.13.29.37.67&2613&2960 \cr
\noalign{\hrule}
 & &9.5.11.13.31.73.211&19&136& & &27.5.67.389.877&3943&442 \cr
11100&3072667455&16.11.17.19.73.211&507&296&11118&3085728885&4.3.13.17.67.3943&665&3278 \cr
 & &256.3.169.17.19.37&9139&2176& & &16.5.7.11.17.19.149&24871&1192 \cr
\noalign{\hrule}
 & &5.19.23.59.113.211&961&396& & &9.121.37.191.401&995&1106 \cr
11101&3073720345&8.9.11.19.961.211&989&778&11119&3086081163&4.3.5.7.11.79.199.401&167&1036 \cr
 & &32.3.11.23.31.43.389&14663&18672& & &32.5.49.37.167.199&9751&13360 \cr
\noalign{\hrule}
 & &9.5.11.37.113.1487&1799&3286& & &9.13.17.29.73.733&26015&27494 \cr
11102&3077487765&4.7.11.31.37.53.257&3561&5948&11120&3086452629&4.3.5.121.13.43.59.233&733&34 \cr
 & &32.3.53.1187.1487&1187&848& & &16.5.121.17.43.733&605&344 \cr
\noalign{\hrule}
 & &81.13.19.101.1523&4081&4100& & &7.11.17.29.31.43.61&9215&22932 \cr
11103&3077536761&8.25.7.11.13.41.53.1523&7597&18&11121&3086722793&8.9.5.343.13.19.97&9145&10406 \cr
 & &32.9.5.7.41.71.107&20377&8560& & &32.3.25.121.31.43.59&825&944 \cr
\noalign{\hrule}
 & &5.17.47.61.73.173&5159&5394& & &5.13.17.1849.1511&8241&15796 \cr
11104&3077624155&4.3.7.11.17.29.31.67.73&519&722&11122&3087192095&8.3.11.17.41.67.359&273&86 \cr
 & &16.9.11.361.31.67.173&24187&24552& & &32.9.7.13.41.43.67&2583&1072 \cr
\noalign{\hrule}
 & &3.5.19.151.233.307&6565&6716& & &49.11.17.31.83.131&7695&7826 \cr
11105&3078336585&8.25.13.23.73.101.307&151&7524&11123&3088508269&4.81.5.343.13.19.31.43&2479&2822 \cr
 & &64.9.11.13.19.23.151&429&736& & &16.9.5.13.17.37.43.67.83&22311&22360 \cr
\noalign{\hrule}
 & &27.7.11.89.127.131&2435&994& & &3.49.121.13.361.37&4369&9062 \cr
11106&3078360747&4.5.49.71.89.487&247&198&11124&3088552467&4.49.17.23.197.257&171&220 \cr
 & &16.9.11.13.19.71.487&9253&7384& & &32.9.5.11.19.197.257&3855&3152 \cr
\noalign{\hrule}
}%
}
$$
\eject
\vglue -23 pt
\noindent\hskip 1 in\hbox to 6.5 in{\ 11125 -- 11160 \hfill\fbd 3088555659 -- 3119310117\frb}
\vskip -9 pt
$$
\vbox{
\nointerlineskip
\halign{\strut
    \vrule \ \ \hfil \frb #\ 
   &\vrule \hfil \ \ \fbb #\frb\ 
   &\vrule \hfil \ \ \frb #\ \hfil
   &\vrule \hfil \ \ \frb #\ 
   &\vrule \hfil \ \ \frb #\ \ \vrule \hskip 2 pt
   &\vrule \ \ \hfil \frb #\ 
   &\vrule \hfil \ \ \fbb #\frb\ 
   &\vrule \hfil \ \ \frb #\ \hfil
   &\vrule \hfil \ \ \frb #\ 
   &\vrule \hfil \ \ \frb #\ \vrule \cr%
\noalign{\hrule}
 & &9.7.19.59.101.433&3095&2662& & &9.7.13.23.37.61.73&2915&2122 \cr
11125&3088555659&4.3.5.7.1331.59.619&19&866&11143&3103602957&4.3.5.7.11.37.53.1061&23&208 \cr
 & &16.11.19.433.619&619&88& & &128.13.23.53.1061&1061&3392 \cr
\noalign{\hrule}
 & &9.25.7.121.13.29.43&811&6014& & &3.17.37.47.79.443&205&1534 \cr
11126&3089411325&4.3.29.31.97.811&8351&16790&11144&3103848933&4.5.13.17.41.59.79&1881&538 \cr
 & &16.5.7.23.73.1193&1679&9544& & &16.9.5.11.13.19.269&2717&32280 \cr
\noalign{\hrule}
 & &11.19.31.41.103.113&53&156& & &3.11.167.311.1811&95&406 \cr
11127&3091772321&8.3.13.31.41.53.113&665&2838&11145&3103910931&4.5.7.11.19.29.1811&713&1098 \cr
 & &32.9.5.7.11.13.19.43&3913&720& & &16.9.19.23.29.31.61&35929&16008 \cr
\noalign{\hrule}
 & &81.7.23.127.1867&57475&14534& & &25.49.11.169.29.47&1179&4 \cr
11128&3092138469&4.25.121.169.19.43&2839&2364&11146&3103925825&8.9.7.11.29.131&2015&1784 \cr
 & &32.3.169.17.167.197&32899&45968& & &128.3.5.13.31.223&6913&192 \cr
\noalign{\hrule}
 & &5.7.13.17.173.2311&341&1206& & &27.5.17.223.6067&341&118 \cr
11129&3092476205&4.9.11.31.67.2311&2261&50&11147&3104999595&4.5.11.31.59.6067&1411&4656 \cr
 & &16.3.25.7.17.19.31&31&2280& & &128.3.17.31.83.97&8051&1984 \cr
\noalign{\hrule}
 & &3.5.11.17.23.191.251&43&208& & &9.7.11.529.43.197&841&4978 \cr
11130&3092913615&32.13.17.23.43.191&1773&1474&11148&3105443187&4.3.19.841.43.131&6107&5290 \cr
 & &128.9.11.43.67.197&8471&12864& & &16.5.529.29.31.197&155&232 \cr
\noalign{\hrule}
 & &9.49.11.17.157.239&77221&53560& & &25.11.127.193.461&87&548 \cr
11131&3094409241&16.5.13.31.47.53.103&471&986&11149&3107382025&8.3.5.11.29.137.193&969&6566 \cr
 & &64.3.13.17.29.53.157&689&928& & &32.9.49.17.19.67&55811&2736 \cr
\noalign{\hrule}
 & &27.5.11.13.107.1499&58489&18016& & &11.169.61.79.347&4665&154 \cr
11132&3096386865&64.23.563.2543&1553&990&11150&3108606787&4.3.5.7.121.13.311&1961&2082 \cr
 & &256.9.5.11.23.1553&1553&2944& & &16.9.5.7.37.53.347&3339&1480 \cr
\noalign{\hrule}
 & &5.61.1823.5569&41679&69524& & &5.11.19.37.257.313&261&52 \cr
11133&3096447535&8.9.7.11.13.191.421&1873&610&11151&3110251265&8.9.5.13.29.37.257&1179&106 \cr
 & &32.3.5.7.11.61.1873&1873&3696& & &32.81.13.53.131&55809&2096 \cr
\noalign{\hrule}
 & &3.7.11.13.29.961.37&51&950& & &3.625.17.31.47.67&959&3036 \cr
11134&3096552459&4.9.25.17.19.31.37&47&232&11152&3111605625&8.9.125.7.11.23.137&31&94 \cr
 & &64.5.17.19.29.47&235&10336& & &32.11.23.31.47.137&1507&368 \cr
\noalign{\hrule}
 & &3.25.11.13.19.23.661&325&302& & &3.5.29.37.41.53.89&12551&11034 \cr
11135&3097991325&4.625.169.151.661&2907&102718&11153&3112724715&4.27.7.11.29.163.613&83927&31930 \cr
 & &16.9.7.11.17.19.23.29&493&168& & &16.5.23.31.41.89.103&713&824 \cr
\noalign{\hrule}
 & &27.5.7.11.19.29.541&443&3148& & &7.11.17.23.157.659&3075&4174 \cr
11136&3098655945&8.11.29.443.787&553&234&11154&3114960541&4.3.25.17.23.41.2087&2021&66 \cr
 & &32.9.7.13.79.443&1027&7088& & &16.9.5.11.41.43.47&1763&16920 \cr
\noalign{\hrule}
 & &13.37.59.313.349&74701&93168& & &7.11.17.961.2477&98523&81184 \cr
11137&3100036823&32.9.11.647.6791&2425&4366&11155&3115939673&64.27.41.43.59.89&19&5270 \cr
 & &128.3.25.11.37.59.97&2425&2112& & &256.9.5.17.19.31&45&2432 \cr
\noalign{\hrule}
 & &81.13.19.241.643&65&578& & &9.5.49.13.19.59.97&367&1628 \cr
11138&3100344741&4.3.5.169.289.241&6061&6230&11156&3116946105&8.3.7.11.37.59.367&25&13604 \cr
 & &16.25.7.11.17.19.29.89&24475&27608& & &64.25.19.179&5&5728 \cr
\noalign{\hrule}
 & &9.25.11.19.23.47.61&341&134& & &27.11.13.41.47.419&12245&7448 \cr
11139&3100875525&4.121.31.47.61.67&5265&2116&11157&3117421593&16.3.5.49.11.19.31.79&313&82 \cr
 & &32.81.5.13.529.31&713&1872& & &64.7.19.31.41.313&9703&4256 \cr
\noalign{\hrule}
 & &9.25.13.53.20011&129833&130310& & &7.13.17.53.109.349&2585&3348 \cr
11140&3102205275&4.125.121.29.37.83.157&4717&1092&11158&3119019631&8.27.5.11.13.31.47.53&9047&6980 \cr
 & &32.3.7.121.13.53.83.89&10043&9968& & &64.9.25.83.109.349&747&800 \cr
\noalign{\hrule}
 & &27.121.19.23.41.53&655&9718& & &25.17.71.167.619&297&322 \cr
11141&3102346467&4.3.5.11.43.113.131&1843&1886&11159&3119280275&4.27.7.11.17.23.71.167&8047&57250 \cr
 & &16.5.19.23.41.97.131&655&776& & &16.3.125.13.229.619&687&520 \cr
\noalign{\hrule}
 & &11.13.17.41.163.191&73&90& & &9.7.11.23.31.59.107&102979&86518 \cr
11142&3103057243&4.9.5.11.13.41.73.191&9135&1304&11160&3119310117&4.29.53.67.181.239&12803&3210 \cr
 & &64.81.25.7.29.163&5075&2592& & &16.3.5.7.29.31.59.107&29&40 \cr
\noalign{\hrule}
}%
}
$$
\eject
\vglue -23 pt
\noindent\hskip 1 in\hbox to 6.5 in{\ 11161 -- 11196 \hfill\fbd 3119551215 -- 3149041665\frb}
\vskip -9 pt
$$
\vbox{
\nointerlineskip
\halign{\strut
    \vrule \ \ \hfil \frb #\ 
   &\vrule \hfil \ \ \fbb #\frb\ 
   &\vrule \hfil \ \ \frb #\ \hfil
   &\vrule \hfil \ \ \frb #\ 
   &\vrule \hfil \ \ \frb #\ \ \vrule \hskip 2 pt
   &\vrule \ \ \hfil \frb #\ 
   &\vrule \hfil \ \ \fbb #\frb\ 
   &\vrule \hfil \ \ \frb #\ \hfil
   &\vrule \hfil \ \ \frb #\ 
   &\vrule \hfil \ \ \frb #\ \vrule \cr%
\noalign{\hrule}
 & &3.5.1331.37.41.103&7049&56296& & &49.13.19.37.43.163&65535&86218 \cr
11161&3119551215&16.7.19.31.53.227&1089&500&11179&3138707299&4.3.5.11.17.257.3919&35&222 \cr
 & &128.9.125.121.53&1325&192& & &16.9.25.7.37.3919&3919&1800 \cr
\noalign{\hrule}
 & &9.19.31.41.83.173&649&130& & &3.125.121.13.17.313&319&6 \cr
11162&3120799419&4.3.5.11.13.31.59.83&1189&1384&11180&3138724875&4.9.5.1331.17.29&5947&6032 \cr
 & &64.11.29.41.59.173&649&928& & &128.13.19.841.313&841&1216 \cr
\noalign{\hrule}
 & &5.7.11.13.59.97.109&57&706& & &9.25.11.31.163.251&1127&1942 \cr
11163&3122154035&4.3.5.13.19.97.353&583&678&11181&3139049925&4.5.49.23.251.971&2771&4026 \cr
 & &16.9.11.53.113.353&18709&8136& & &16.3.7.11.17.23.61.163&1037&1288 \cr
\noalign{\hrule}
 & &9.7.11.19.59.4019&3757&8300& & &25.11.13.23.181.211&4689&164 \cr
11164&3122172207&8.3.25.13.289.19.83&649&326&11182&3140254975&8.9.11.13.41.521&475&46 \cr
 & &32.11.17.59.83.163&1411&2608& & &32.3.25.19.23.41&2337&16 \cr
\noalign{\hrule}
 & &27.5.343.37.1823&6413&6278& & &3.125.11.13.157.373&24707&28632 \cr
11165&3123318555&4.121.43.53.73.1823&2481&658&11183&3140333625&16.9.5.31.797.1193&559&36424 \cr
 & &16.3.7.121.47.53.827&43831&45496& & &256.13.29.43.157&1247&128 \cr
\noalign{\hrule}
 & &5.13.37.41.79.401&959&558& & &3.25.11.13.443.661&1241&9834 \cr
11166&3123707795&4.9.5.7.13.31.79.137&407&802&11184&3140526675&4.9.121.17.73.149&15815&6982 \cr
 & &16.3.7.11.37.137.401&411&616& & &16.5.3163.3491&3491&25304 \cr
\noalign{\hrule}
 & &13.71.641.5281&2179&3102& & &81.5.41.43.53.83&229&176 \cr
11167&3124466683&4.3.11.47.641.2179&4615&2436&11185&3140951985&32.11.41.43.83.229&6625&3222 \cr
 & &32.9.5.7.13.29.47.71&2961&2320& & &128.9.125.11.53.179&1969&1600 \cr
\noalign{\hrule}
 & &729.13.73.4519&839&110& & &7.11.19.59.151.241&1391&270 \cr
11168&3126339099&4.5.11.839.4519&6683&15912&11186&3141161947&4.27.5.7.13.107.241&5497&5738 \cr
 & &64.9.13.17.41.163&697&5216& & &16.9.13.19.23.151.239&2151&2392 \cr
\noalign{\hrule}
 & &9.5.49.11.83.1553&61009&24406& & &27.25.121.17.31.73&481&294 \cr
11169&3126445245&4.169.361.12203&6225&5978&11187&3142118925&4.81.49.11.13.37.73&15709&5270 \cr
 & &16.3.25.49.13.19.61.83&793&760& & &16.5.7.17.23.31.683&683&1288 \cr
\noalign{\hrule}
 & &3.25.7.11.17.19.23.73&29713&63362& & &19.23.41.109.1609&333&1276 \cr
11170&3131880675&4.13.43.691.2437&873&1564&11188&3142301377&8.9.11.19.29.37.109&2645&574 \cr
 & &32.9.13.17.23.43.97&1261&2064& & &32.3.5.7.11.529.41&231&1840 \cr
\noalign{\hrule}
 & &5.11.19.23.2209.59&6201&4844& & &27.5.7.11.41.73.101&1039&8816 \cr
11171&3132505585&8.9.7.11.13.19.53.173&31&658&11189&3142335735&32.19.29.41.1039&909&130 \cr
 & &32.3.49.31.47.173&1519&8304& & &128.9.5.13.29.101&29&832 \cr
\noalign{\hrule}
 & &13.29.83.239.419&3675&3256& & &9.5.13.23.37.59.107&19549&17366 \cr
11172&3133512031&16.3.25.49.11.13.37.83&1557&478&11190&3142832355&4.3.13.19.113.173.457&8825&8998 \cr
 & &64.27.5.49.173.239&6615&5536& & &16.25.11.19.113.353.409&368885&369736 \cr
\noalign{\hrule}
 & &3.7.121.23.29.1849&43127&37580& & &9.7.11.37.149.823&245&578 \cr
11173&3133772103&8.5.49.61.101.1879&817&1062&11191&3144278907&4.5.343.11.289.149&10699&8166 \cr
 & &32.9.19.43.59.61.101&17877&18544& & &16.3.13.17.823.1361&1361&1768 \cr
\noalign{\hrule}
 & &3.25.23.29.31.43.47&741&506& & &81.1331.289.101&39311&68500 \cr
11174&3134116275&4.9.5.11.13.19.529.31&14093&4042&11192&3146895279&8.125.19.137.2069&153&2222 \cr
 & &16.11.17.43.47.829&829&1496& & &32.9.11.17.101.137&137&16 \cr
\noalign{\hrule}
 & &5.7.17.19.37.59.127&2859&2156& & &27.29.43.211.443&845&402 \cr
11175&3134209505&8.3.343.11.127.953&16539&27022&11193&3147142437&4.81.5.169.67.211&4085&1342 \cr
 & &32.9.37.59.149.229&2061&2384& & &16.25.11.13.19.43.61&3575&9272 \cr
\noalign{\hrule}
 & &27.5.131.233.761&2765&3526& & &9.5.49.11.13.67.149&89&782 \cr
11176&3135780405&4.25.7.41.43.79.131&13527&10252&11194&3147789645&4.5.7.17.23.89.149&2563&2652 \cr
 & &32.81.11.41.167.233&1837&1968& & &32.3.11.13.289.23.233&5359&4624 \cr
\noalign{\hrule}
 & &9.5.11.17.59.71.89&6913&1898& & &3.49.23.29.97.331&1&22 \cr
11177&3137288715&4.13.31.71.73.223&42343&47526&11195&3148059243&4.7.11.29.97.331&2565&248 \cr
 & &16.3.7.23.7921.263&6049&4984& & &64.27.5.11.19.31&32395&288 \cr
\noalign{\hrule}
 & &9.49.17.41.59.173&1651&440& & &3.5.7.11.17.19.23.367&1233&3068 \cr
11178&3137397039&16.3.5.7.11.13.59.127&289&124&11196&3149041665&8.27.7.13.19.59.137&31625&24956 \cr
 & &128.13.289.31.127&3937&14144& & &64.125.11.17.23.367&25&32 \cr
\noalign{\hrule}
}%
}
$$
\eject
\vglue -23 pt
\noindent\hskip 1 in\hbox to 6.5 in{\ 11197 -- 11232 \hfill\fbd 3150756455 -- 3182787677\frb}
\vskip -9 pt
$$
\vbox{
\nointerlineskip
\halign{\strut
    \vrule \ \ \hfil \frb #\ 
   &\vrule \hfil \ \ \fbb #\frb\ 
   &\vrule \hfil \ \ \frb #\ \hfil
   &\vrule \hfil \ \ \frb #\ 
   &\vrule \hfil \ \ \frb #\ \ \vrule \hskip 2 pt
   &\vrule \ \ \hfil \frb #\ 
   &\vrule \hfil \ \ \fbb #\frb\ 
   &\vrule \hfil \ \ \frb #\ \hfil
   &\vrule \hfil \ \ \frb #\ 
   &\vrule \hfil \ \ \frb #\ \vrule \cr%
\noalign{\hrule}
 & &5.7.11.17.31.53.293&241&24& & &9.7.19.47.181.311&24127&32164 \cr
11197&3150756455&16.3.11.17.241.293&1165&3816&11215&3166875369&8.7.11.17.23.43.1049&361&1410 \cr
 & &256.27.5.53.233&233&3456& & &32.3.5.17.361.43.47&817&1360 \cr
\noalign{\hrule}
 & &3.13.181.601.743&281&462& & &9.5.7.13.31.61.409&60559&64186 \cr
11198&3152147037&4.9.7.11.13.281.601&3715&2896&11216&3167150805&4.7.23.67.479.2633&9455&8976 \cr
 & &128.5.181.281.743&281&320& & &128.3.5.11.17.23.31.61.67&4301&4288 \cr
\noalign{\hrule}
 & &9.5.7.11.13.17.23.179&71251&43864& & &25.11.10201.1129&121303&133722 \cr
11199&3152654505&16.43.1657.5483&1913&3570&11217&3167155475&4.9.7.13.17.19.23.31.43&505&22 \cr
 & &64.3.5.7.17.43.1913&1913&1376& & &16.3.5.11.13.19.43.101&817&312 \cr
\noalign{\hrule}
 & &9.49.13.37.89.167&11033&12536& & &81.5.11.17.19.31.71&331&1870 \cr
11200&3152754423&16.11.17.59.89.1567&27&1540&11218&3167161965&4.25.121.289.331&13347&21622 \cr
 & &128.27.5.7.121.59&7139&960& & &16.9.19.569.1483&1483&4552 \cr
\noalign{\hrule}
 & &27.5.37.61.79.131&4597&56342& & &3.5.19.73.197.773&1485&2258 \cr
11201&3153288555&4.11.13.197.4597&2227&2370&11219&3168206205&4.81.25.11.73.1129&152281&158194 \cr
 & &16.3.5.17.79.131.197&197&136& & &16.19.23.181.197.773&181&184 \cr
\noalign{\hrule}
 & &27.11.17.709.881&67613&64970& & &81.5.19.43.61.157&8677&4438 \cr
11202&3153751821&4.9.5.7.13.73.89.743&2771&2836&11220&3168885645&4.3.7.19.317.8677&16027&10004 \cr
 & &32.17.73.163.709.743&11899&11888& & &32.7.11.31.41.47.61&10199&7216 \cr
\noalign{\hrule}
 & &7.23.43.73.6241&1489&1650& & &9.3125.11.10243&159809&149566 \cr
11203&3154070339&4.3.25.11.6241.1489&6843&602&11221&3168928125&4.13.17.19.53.83.647&12045&1046 \cr
 & &16.9.5.7.11.43.2281&2281&3960& & &16.3.5.11.73.83.523&6059&4184 \cr
\noalign{\hrule}
 & &7.11.13.41.151.509&795&262& & &3.31.43.47.101.167&13145&8398 \cr
11204&3154370219&4.3.5.11.53.131.509&37&546&11222&3170203251&4.5.11.13.17.19.31.239&1169&846 \cr
 & &16.9.5.7.13.37.131&5895&296& & &16.9.7.11.47.167.239&717&616 \cr
\noalign{\hrule}
 & &9.11.59.127.4253&2125&2066& & &3.11.17.23.71.3461&15611&12150 \cr
11205&3154905171&4.3.125.17.1033.4253&2587&74888&11223&3170667093&4.729.25.11.67.233&1097&6922 \cr
 & &64.5.11.13.23.37.199&22885&15392& & &16.67.1097.3461&1097&536 \cr
\noalign{\hrule}
 & &9.25.19.67.103.107&759&1274& & &9.25.11.73.97.181&9917&11908 \cr
11206&3156689925&4.27.5.49.11.13.23.67&107&1648&11224&3172110975&8.13.47.73.211.229&13083&2320 \cr
 & &128.49.11.103.107&539&64& & &256.3.5.49.13.29.89&18473&11392 \cr
\noalign{\hrule}
 & &9.49.11.13.113.443&1577&334& & &5.121.169.19.23.71&1783&1428 \cr
11207&3156870717&4.3.19.83.167.443&7595&6266&11225&3172355615&8.3.7.121.17.23.1783&5415&24896 \cr
 & &16.5.49.13.19.31.241&2945&1928& & &1024.9.5.361.389&7391&4608 \cr
\noalign{\hrule}
 & &5.11.19.47.53.1213&391&126& & &5.19.29.47.127.193&265881&265834 \cr
11208&3157554235&4.9.7.17.19.23.1213&10537&11750&11226&3173806835&4.3.7.11.23.127.1151.5779&13113&7334 \cr
 & &16.3.125.7.41.47.257&7175&6168& & &16.27.11.19.31.47.193.1151&9207&9208 \cr
\noalign{\hrule}
 & &27.11.79.311.433&221&212& & &3.5.11.47.229.1789&43301&40782 \cr
11209&3159597969&8.3.11.13.17.53.79.311&3251&170&11227&3177076155&4.9.5.7.19.43.53.971&2327&6412 \cr
 & &32.5.289.53.3251&76585&52016& & &32.49.13.53.179.229&9487&10192 \cr
\noalign{\hrule}
 & &27.29.61.127.521&19175&4066& & &7.17.29.79.107.109&5&114 \cr
11210&3160334421&4.9.25.13.19.59.107&781&610&11228&3179672027&4.3.5.19.29.79.107&1199&1092 \cr
 & &16.125.11.59.61.71&8875&5192& & &32.9.5.7.11.13.19.109&2223&880 \cr
\noalign{\hrule}
 & &9.5.11.19.41.59.139&653&598& & &81.5.19.29.53.269&14267&12922 \cr
11211&3162346605&4.13.19.23.41.59.653&74643&120604&11229&3181520835&4.3.7.11.13.29.71.1297&4045&10222 \cr
 & &32.3.11.139.179.2741&2741&2864& & &16.5.7.13.19.269.809&809&728 \cr
\noalign{\hrule}
 & &243.25.11.19.47.53&541&784& & &5.11.31.67.89.313&7569&9646 \cr
11212&3162760425&32.49.11.19.47.541&8343&18166&11230&3182244395&4.9.7.13.841.53.89&1921&1550 \cr
 & &128.81.31.103.293&9083&6592& & &16.3.25.17.841.31.113&14297&13560 \cr
\noalign{\hrule}
 & &11.169.53.97.331&65317&33210& & &5.11.59.293.3347&1349&1998 \cr
11213&3163406389&4.81.5.49.31.41.43&1649&1606&11231&3182277395&4.27.5.19.37.71.293&413&1052 \cr
 & &16.27.7.11.17.41.73.97&8687&8856& & &32.3.7.19.37.59.263&9731&6384 \cr
\noalign{\hrule}
 & &9.5.49.71.113.179&1951&14660& & &557.2017.2833&1287&730 \cr
11214&3166637985&8.3.25.733.1951&25487&23288&11232&3182787677&4.9.5.11.13.73.2833&9469&970 \cr
 & &128.7.11.41.71.331&3641&2624& & &16.3.25.17.97.557&2425&408 \cr
\noalign{\hrule}
}%
}
$$
\eject
\vglue -23 pt
\noindent\hskip 1 in\hbox to 6.5 in{\ 11233 -- 11268 \hfill\fbd 3182800257 -- 3207231027\frb}
\vskip -9 pt
$$
\vbox{
\nointerlineskip
\halign{\strut
    \vrule \ \ \hfil \frb #\ 
   &\vrule \hfil \ \ \fbb #\frb\ 
   &\vrule \hfil \ \ \frb #\ \hfil
   &\vrule \hfil \ \ \frb #\ 
   &\vrule \hfil \ \ \frb #\ \ \vrule \hskip 2 pt
   &\vrule \ \ \hfil \frb #\ 
   &\vrule \hfil \ \ \fbb #\frb\ 
   &\vrule \hfil \ \ \frb #\ \hfil
   &\vrule \hfil \ \ \frb #\ 
   &\vrule \hfil \ \ \frb #\ \vrule \cr%
\noalign{\hrule}
 & &27.7.13.19.29.2351&9925&6532& & &9.25.23.53.103.113&21057&9418 \cr
11233&3182800257&8.3.25.19.23.71.397&341&56&11251&3192286725&4.27.17.277.7019&7249&230 \cr
 & &128.5.7.11.23.31.71&7843&22720& & &16.5.11.17.23.659&187&5272 \cr
\noalign{\hrule}
 & &9.5.29.41.157.379&1969&2584& & &3.25.49.11.13.59.103&5947&22172 \cr
11234&3183708015&16.3.11.17.19.179.379&8323&1880&11252&3193615425&8.7.19.23.241.313&625&1062 \cr
 & &256.5.7.11.29.41.47&517&896& & &32.9.625.59.313&939&400 \cr
\noalign{\hrule}
 & &27.5.121.67.2909&215&148& & &7.11.157.349.757&195939&187612 \cr
11235&3183740505&8.9.25.37.43.2909&3383&6292&11253&3193829177&8.81.17.31.41.59.89&3479&1060 \cr
 & &64.121.13.17.37.199&8177&6368& & &64.27.5.49.31.53.71&57505&61344 \cr
\noalign{\hrule}
 & &9.11.31.461.2251&2059&2090& & &25.19.29.37.6269&3497&2772 \cr
11236&3184735059&4.5.121.19.29.71.2251&93&42862&11254&3195152575&8.9.7.11.13.19.37.269&425&2534 \cr
 & &16.3.841.31.739&841&5912& & &32.3.25.49.13.17.181&7059&13328 \cr
\noalign{\hrule}
 & &27.137.919.937&12875&11938& & &3.5.13.73.193.1163&107&472 \cr
11237&3185219997&4.125.47.103.127.137&4697&1272&11255&3195173865&16.13.59.107.1163&965&198 \cr
 & &64.3.5.7.11.53.61.103&60049&68320& & &64.9.5.11.107.193&107&1056 \cr
\noalign{\hrule}
 & &9.5.7.11.19.97.499&39997&45332& & &81.11.19.23.29.283&455&172 \cr
11238&3186611505&8.49.23.37.47.1619&97&1716&11256&3195534969&8.27.5.7.13.23.29.43&209&412 \cr
 & &64.3.11.13.23.47.97&1081&416& & &64.5.11.13.19.43.103&4429&2080 \cr
\noalign{\hrule}
 & &49.11.19.23.83.163&65377&86376& & &27.5.7.11.169.17.107&34621&29044 \cr
11239&3186661247&16.3.13.47.59.61.107&437&330&11257&3195537345&8.9.53.89.137.389&661&572 \cr
 & &64.9.5.11.19.23.47.61&2115&1952& & &64.11.13.53.389.661&20617&21152 \cr
\noalign{\hrule}
 & &27.5.7.121.29.961&1525&1742& & &25.7.11.19.23.29.131&2743&5232 \cr
11240&3186680805&4.125.13.29.31.61.67&363&1262&11258&3195813775&32.3.7.13.23.109.211&393&370 \cr
 & &16.3.121.61.67.631&4087&5048& & &128.9.5.13.37.131.211&7807&7488 \cr
\noalign{\hrule}
 & &5.961.47.103.137&4823&18& & &3.5.7.23.71.103.181&1997&2902 \cr
11241&3186757685&4.9.7.13.53.137&2343&2480&11259&3196621995&4.7.103.1451.1997&1359&638 \cr
 & &128.27.5.11.31.71&297&4544& & &16.9.11.29.151.1451&47883&35032 \cr
\noalign{\hrule}
 & &3.5.19.37.41.73.101&1397&5538& & &9.11.13.17.313.467&84025&62146 \cr
11242&3187679685&4.9.11.13.37.71.127&2669&2030&11260&3198075309&4.25.7.23.193.3361&39&154 \cr
 & &16.5.7.11.13.17.29.157&31871&19448& & &16.3.5.49.11.13.3361&3361&1960 \cr
\noalign{\hrule}
 & &11.41.103.163.421&83709&100970& & &27.5.7.11.37.53.157&1791&64 \cr
11243&3187744219&4.9.5.23.71.131.439&285&154&11261&3200381415&128.243.37.199&4595&4396 \cr
 & &16.27.25.7.11.19.23.71&36423&32200& & &1024.5.7.157.919&919&512 \cr
\noalign{\hrule}
 & &27.5.7.83.97.419&14471&25784& & &81.19.71.83.353&51337&76400 \cr
11244&3187833705&16.7.11.29.293.499&1275&776&11262&3201472431&32.25.11.13.191.359&4233&8182 \cr
 & &256.3.25.11.17.29.97&1595&2176& & &128.3.5.17.83.4091&4091&5440 \cr
\noalign{\hrule}
 & &3.7.19.23.37.41.229&495&208& & &11.59.107.193.239&67&126 \cr
11245&3188025561&32.27.5.11.13.23.229&703&4564&11263&3203197261&4.9.7.11.67.107.239&1475&1154 \cr
 & &256.5.7.19.37.163&815&128& & &16.3.25.7.59.67.577&11725&13848 \cr
\noalign{\hrule}
 & &27.25.7.13.23.37.61&33&292& & &9.13.17.19.29.37.79&275&428 \cr
11246&3188633175&8.81.11.23.61.73&1295&3158&11264&3203429697&8.25.11.13.29.79.107&1377&508 \cr
 & &32.5.7.11.37.1579&1579&176& & &64.81.5.17.107.127&5715&3424 \cr
\noalign{\hrule}
 & &3.7.23.137.157.307&561&538& & &25.7.43.541.787&291181&290394 \cr
11247&3189376029&4.9.11.17.137.269.307&2135&628&11265&3203896675&4.9.7.11.13.17.73.103.257&265&8 \cr
 & &32.5.7.17.61.157.269&4573&4880& & &64.3.5.11.17.53.73.103&82709&86496 \cr
\noalign{\hrule}
 & &9.25.17.19.529.83&349&824& & &31.251.277.1487&202059&209840 \cr
11248&3190941225&16.3.23.83.103.349&2915&2812&11266&3204986119&32.9.5.11.13.43.61.157&95&34 \cr
 & &128.5.11.19.37.53.349&21571&22336& & &128.3.25.11.13.17.19.157&182325&190912 \cr
\noalign{\hrule}
 & &3.5.7.23.47.61.461&153&176& & &5.49.11.13.239.383&23999&10854 \cr
11249&3191874105&32.27.5.11.17.61.461&329&1976&11267&3206998795&4.81.7.67.103.233&575&2206 \cr
 & &512.7.11.13.17.19.47&4199&2816& & &16.3.25.23.67.1103&7705&26472 \cr
\noalign{\hrule}
 & &3.13.41.67.83.359&979&620& & &9.7.11.13.19.41.457&1417&46 \cr
11250&3192242001&8.5.11.31.67.83.89&2655&4732&11268&3207231027&4.3.169.23.41.109&11647&9140 \cr
 & &64.9.25.7.11.169.59&10325&13728& & &32.5.19.457.613&613&80 \cr
\noalign{\hrule}
}%
}
$$
\eject
\vglue -23 pt
\noindent\hskip 1 in\hbox to 6.5 in{\ 11269 -- 11304 \hfill\fbd 3208473495 -- 3234973115\frb}
\vskip -9 pt
$$
\vbox{
\nointerlineskip
\halign{\strut
    \vrule \ \ \hfil \frb #\ 
   &\vrule \hfil \ \ \fbb #\frb\ 
   &\vrule \hfil \ \ \frb #\ \hfil
   &\vrule \hfil \ \ \frb #\ 
   &\vrule \hfil \ \ \frb #\ \ \vrule \hskip 2 pt
   &\vrule \ \ \hfil \frb #\ 
   &\vrule \hfil \ \ \fbb #\frb\ 
   &\vrule \hfil \ \ \frb #\ \hfil
   &\vrule \hfil \ \ \frb #\ 
   &\vrule \hfil \ \ \frb #\ \vrule \cr%
\noalign{\hrule}
 & &9.5.17.31.193.701&2915&3068& & &3.25.7.11.13.59.727&851&124 \cr
11269&3208473495&8.25.11.13.53.59.701&513&188&11287&3220191975&8.7.11.23.31.37.59&1395&1454 \cr
 & &64.27.11.19.47.53.59&52687&55968& & &32.9.5.23.961.727&961&1104 \cr
\noalign{\hrule}
 & &41.47.149.11177&137973&149150& & &11.13.17.29.43.1063&3243&2180 \cr
11270&3209173771&4.3.25.11.19.37.113.157&1937&210&11288&3222438791&8.3.5.13.23.43.47.109&30827&35772 \cr
 & &16.9.125.7.13.37.149&4625&6552& & &64.9.11.29.271.1063&271&288 \cr
\noalign{\hrule}
 & &9.25.361.71.557&67&1738& & &7.11.17.19.101.1283&351&21460 \cr
11271&3212205075&4.3.5.11.67.71.79&557&628&11289&3222858793&8.27.5.7.13.29.37&57&202 \cr
 & &32.11.67.157.557&737&2512& & &32.81.13.19.101&13&1296 \cr
\noalign{\hrule}
 & &3.5.29.41.233.773&5203&6392& & &29.41.71.181.211&220935&213514 \cr
11272&3212244015&16.121.17.43.47.233&31275&3082&11290&3224046029&4.3.5.7.11.13.101.103.151&1537&426 \cr
 & &64.9.25.23.67.139&15985&6432& & &16.9.5.7.29.53.71.103&3605&3816 \cr
\noalign{\hrule}
 & &9.25.7.29.37.1901&7097&8998& & &3.19.229.337.733&123625&123396 \cr
11273&3212642475&4.3.5.7.11.47.151.409&871&116&11291&3224365113&8.9.125.7.13.19.23.43.113&61&7414 \cr
 & &32.11.13.29.67.409&5317&11792& & &32.5.7.11.61.113.337&6893&6160 \cr
\noalign{\hrule}
 & &9.7.11.17.23.71.167&8047&57250& & &9.5.11.13.29.37.467&4747&544 \cr
11274&3212808291&4.125.13.229.619&297&322&11292&3224520585&64.5.17.29.47.101&429&934 \cr
 & &16.27.5.7.11.13.23.229&687&520& & &256.3.11.13.17.467&17&128 \cr
\noalign{\hrule}
 & &3.11.289.23.97.151&2347&2600& & &5.11.47.977.1277&22321&23598 \cr
11275&3212834097&16.25.13.17.151.2347&2457&110&11293&3225120965&4.27.5.11.13.17.19.23.101&1241&2554 \cr
 & &64.27.125.7.11.169&10647&4000& & &16.9.289.19.73.1277&5491&5256 \cr
\noalign{\hrule}
 & &3.19.151.181.2063&805&1258& & &9.5.47.59.103.251&759&2014 \cr
11276&3213879621&4.5.7.17.19.23.37.181&3839&324&11294&3226066605&4.27.11.19.23.53.103&767&664 \cr
 & &32.81.7.11.17.349&5933&33264& & &64.11.13.19.23.59.83&11869&13984 \cr
\noalign{\hrule}
 & &3.25.11.13.23.83.157&833&2742& & &3.5.7.41.43.107.163&1469&5214 \cr
11277&3214421925&4.9.49.17.157.457&13061&15730&11295&3228590715&4.9.11.13.43.79.113&815&428 \cr
 & &16.5.7.121.13.37.353&2471&3256& & &32.5.13.79.107.163&79&208 \cr
\noalign{\hrule}
 & &13.47.59.257.347&99&158& & &5.7.13.19.23.109.149&1111&306 \cr
11278&3214813771&4.9.11.13.47.79.347&215&826&11296&3229279235&4.9.11.17.19.101.149&1183&736 \cr
 & &16.3.5.7.11.43.59.79&8295&3784& & &256.3.7.11.169.17.23&663&1408 \cr
\noalign{\hrule}
 & &3.49.11.41.71.683&1725&1186& & &27.7.17.53.61.311&119&430 \cr
11279&3214940421&4.9.25.23.593.683&21197&5488&11297&3230552619&4.3.5.49.289.43.53&10109&2318 \cr
 & &128.5.343.11.41.47&235&448& & &16.5.11.19.61.919&1045&7352 \cr
\noalign{\hrule}
 & &5.11.13.19.23.41.251&47&162& & &3.25.7.19.23.73.193&5687&1248 \cr
11280&3215474405&4.81.13.41.47.251&2185&1078&11298&3232368825&64.9.5.7.121.13.47&23&68 \cr
 & &16.3.5.49.11.19.23.47&141&392& & &512.121.17.23.47&2057&12032 \cr
\noalign{\hrule}
 & &3.625.11.13.67.179&2641&3266& & &9.25.11.169.59.131&5687&5818 \cr
11281&3215623125&4.13.19.23.67.71.139&175&1098&11299&3232847475&4.3.5.1331.13.47.2909&707&188378 \cr
 & &16.9.25.7.23.61.139&4209&7784& & &16.7.101.131.719&5033&808 \cr
\noalign{\hrule}
 & &243.5.7.13.17.29.59&41&418& & &25.11.53.59.3761&118237&103662 \cr
11282&3216004155&4.9.5.7.11.19.41.59&3317&2668&11300&3234177925&4.9.49.13.19.127.443&1343&1070 \cr
 & &32.23.29.31.41.107&4387&11408& & &16.3.5.7.17.79.107.443&104991&101864 \cr
\noalign{\hrule}
 & &9.7.17.101.131.227&869&2990& & &3.625.7.11.43.521&201&674 \cr
11283&3216681027&4.3.5.11.13.23.79.131&9899&11716&11301&3234433125&4.9.5.67.337.521&1247&1768 \cr
 & &32.13.19.29.101.521&6773&8816& & &64.13.17.29.43.337&6409&10784 \cr
\noalign{\hrule}
 & &27.25.343.13.1069&7303&1958& & &27.13.961.43.223&16313&25010 \cr
11284&3217502925&4.5.11.13.67.89.109&931&486&11302&3234475179&4.9.5.11.41.61.1483&6727&688 \cr
 & &16.243.49.11.19.67&603&1672& & &128.7.961.41.43&287&64 \cr
\noalign{\hrule}
 & &9.7.13.29.313.433&28955&34584& & &81.19.47.97.461&451&8308 \cr
11285&3218949279&16.27.5.11.131.5791&5947&11738&11303&3234514761&8.11.31.41.47.67&17955&16684 \cr
 & &64.11.19.313.5869&5869&6688& & &64.27.5.7.19.43.97&301&160 \cr
\noalign{\hrule}
 & &13.67.73.197.257&1165&3726& & &5.121.23.383.607&6023&7938 \cr
11286&3219143707&4.81.5.23.233.257&6149&790&11304&3234973115&4.81.49.121.19.317&265&582 \cr
 & &16.3.25.11.13.43.79&1419&15800& & &16.243.5.7.19.53.97&32319&41128 \cr
\noalign{\hrule}
}%
}
$$
\eject
\vglue -23 pt
\noindent\hskip 1 in\hbox to 6.5 in{\ 11305 -- 11340 \hfill\fbd 3236495207 -- 3257221275\frb}
\vskip -9 pt
$$
\vbox{
\nointerlineskip
\halign{\strut
    \vrule \ \ \hfil \frb #\ 
   &\vrule \hfil \ \ \fbb #\frb\ 
   &\vrule \hfil \ \ \frb #\ \hfil
   &\vrule \hfil \ \ \frb #\ 
   &\vrule \hfil \ \ \frb #\ \ \vrule \hskip 2 pt
   &\vrule \ \ \hfil \frb #\ 
   &\vrule \hfil \ \ \fbb #\frb\ 
   &\vrule \hfil \ \ \frb #\ \hfil
   &\vrule \hfil \ \ \frb #\ 
   &\vrule \hfil \ \ \frb #\ \vrule \cr%
\noalign{\hrule}
 & &11.17.19.29.101.311&2867&2556& & &27.11.19.173.3329&1769&1560 \cr
11305&3236495207&8.9.19.47.61.71.101&1741&4420&11323&3249899631&16.81.5.13.29.61.173&2299&50 \cr
 & &64.3.5.13.17.71.1741&26115&29536& & &64.125.121.19.61&671&4000 \cr
\noalign{\hrule}
 & &5.11.19.953.3251&2373&2392& & &729.25.11.13.29.43&8953&9272 \cr
11306&3237622135&16.3.7.11.13.23.113.3251&19&3270&11324&3249900225&16.7.13.19.43.61.1279&76183&1836 \cr
 & &64.9.5.7.19.109.113&6867&3616& & &128.27.17.29.37.71&1207&2368 \cr
\noalign{\hrule}
 & &27.5.7.11.169.19.97&24127&14912& & &49.121.13.149.283&1811&5490 \cr
11307&3237699465&128.9.23.233.1049&1573&524&11325&3250105859&4.9.5.121.61.1811&11977&10166 \cr
 & &1024.121.13.23.131&3013&5632& & &16.3.5.7.13.17.23.29.59&10005&8024 \cr
\noalign{\hrule}
 & &27.11.19.31.107.173&6467&11830& & &3.7.11.13.71.79.193&135&58 \cr
11308&3238184763&4.3.5.7.11.169.29.223&113&1070&11326&3250858611&4.81.5.13.29.71.79&18073&12464 \cr
 & &16.25.107.113.223&5575&904& & &128.5.11.19.31.41.53&24149&16960 \cr
\noalign{\hrule}
 & &243.5.7.139.2741&25&164& & &9.5.11.23.347.823&9589&17570 \cr
11309&3240396495&8.9.125.41.2741&12397&12272&11327&3251339685&4.3.25.7.43.223.251&347&322 \cr
 & &256.49.11.13.23.41.59&85813&83072& & &16.49.23.43.251.347&2107&2008 \cr
\noalign{\hrule}
 & &9.5.7.19.61.83.107&4355&4526& & &9.5.7.11.13.19.29.131&2941&452 \cr
11310&3242319885&4.25.7.13.31.61.67.73&2643&17182&11328&3251393145&8.5.7.11.17.113.173&2929&1026 \cr
 & &16.3.121.71.73.881&62551&70664& & &32.27.17.19.29.101&101&816 \cr
\noalign{\hrule}
 & &9.13.41.43.79.199&16775&25332& & &27.5.11.13.29.37.157&11725&10726 \cr
11311&3242786391&8.27.25.11.61.2111&7493&15728&11329&3252139605&4.125.7.29.31.67.173&909&1034 \cr
 & &256.5.59.127.983&57997&81280& & &16.9.7.11.31.47.101.173&37541&37976 \cr
\noalign{\hrule}
 & &3.5.13.37.101.4451&2207&2244& & &121.13.17.19.37.173&1075&828 \cr
11312&3243510465&8.9.5.11.13.17.101.2207&391&11426&11330&3252213679&8.9.25.11.17.23.37.43&1&186 \cr
 & &32.11.289.23.29.197&92191&72496& & &32.27.5.23.31.43&1161&57040 \cr
\noalign{\hrule}
 & &27.25.11.13.19.29.61&553&492& & &9.31.41.59.61.79&935&1484 \cr
11313&3244301775&8.81.5.7.13.29.41.79&2257&4142&11331&3252347919&8.5.7.11.17.31.53.79&3477&1028 \cr
 & &32.7.19.37.41.61.109&4469&4144& & &64.3.7.11.19.61.257&4883&2464 \cr
\noalign{\hrule}
 & &27.5.49.11.17.43.61&117&422& & &27.5.49.37.97.137&187&2 \cr
11314&3244664115&4.243.13.17.43.211&487&244&11332&3252549195&4.7.11.17.97.137&8103&8200 \cr
 & &32.13.61.211.487&2743&7792& & &64.3.25.11.37.41.73&2993&1760 \cr
\noalign{\hrule}
 & &27.5.13.19.31.43.73&1749&514& & &17.47.113.137.263&2167&2304 \cr
11315&3244768605&4.81.11.43.53.257&9269&11548&11333&3253130897&512.9.11.47.113.197&595&5906 \cr
 & &32.11.13.23.31.2887&2887&4048& & &2048.3.5.7.17.2953&20671&15360 \cr
\noalign{\hrule}
 & &9.29.67.151.1229&715&514& & &3.125.11.13.19.31.103&1379&246 \cr
11316&3245219973&4.3.5.11.13.29.151.257&67&386&11334&3253267875&4.9.7.19.31.41.197&721&550 \cr
 & &16.5.13.67.193.257&3341&7720& & &16.25.49.11.103.197&197&392 \cr
\noalign{\hrule}
 & &125.11.31.271.281&447&3428& & &3.19.29.31.173.367&459&92 \cr
11317&3245936375&8.3.149.281.857&215&66&11335&3253469313&8.81.17.23.31.173&3245&734 \cr
 & &32.9.5.11.43.857&7713&688& & &32.5.11.17.59.367&295&2992 \cr
\noalign{\hrule}
 & &25.11.13.769.1181&247&522& & &5.7.11.19.31.113.127&9483&5546 \cr
11318&3246775675&4.9.169.19.29.1181&675&506&11336&3254304515&4.3.5.11.29.47.59.109&589&2184 \cr
 & &16.243.25.11.19.23.29&5589&4408& & &64.9.7.13.19.31.109&981&416 \cr
\noalign{\hrule}
 & &3.11.13.37.43.67.71&5453&696& & &7.11.13.17.961.199&1177&216 \cr
11319&3246838023&16.9.7.19.29.37.41&215&44&11337&3254314063&16.27.121.13.17.107&655&434 \cr
 & &128.5.11.29.41.43&145&2624& & &64.3.5.7.31.107.131&1605&4192 \cr
\noalign{\hrule}
 & &125.11.13.19.73.131&5899&12834& & &9.7.11.13.31.43.271&175&166 \cr
11320&3247833875&4.9.25.17.23.31.347&313&262&11338&3254438187&4.25.49.13.43.83.271&2697&55568 \cr
 & &16.3.31.131.313.347&9703&8328& & &128.3.5.23.29.31.151&4379&7360 \cr
\noalign{\hrule}
 & &7.31.37.41.71.139&88205&88464& & &5.7.37.47.73.733&64493&30042 \cr
11321&3248766241&32.3.5.13.19.23.59.71.97&1353&4&11339&3256825285&4.9.121.13.41.1669&1717&3290 \cr
 & &256.9.5.11.13.41.97&11349&7040& & &16.3.5.7.17.41.47.101&1717&984 \cr
\noalign{\hrule}
 & &3.11.37.67.151.263&9641&7980& & &27.25.7.11.29.2161&73&102 \cr
11322&3248801391&8.9.5.7.19.31.37.311&7571&6424&11340&3257221275&4.81.11.17.73.2161&635&1526 \cr
 & &128.7.11.19.67.73.113&8249&8512& & &16.5.7.17.73.109.127&13843&9928 \cr
\noalign{\hrule}
}%
}
$$
\eject
\vglue -23 pt
\noindent\hskip 1 in\hbox to 6.5 in{\ 11341 -- 11376 \hfill\fbd 3257393139 -- 3288920739\frb}
\vskip -9 pt
$$
\vbox{
\nointerlineskip
\halign{\strut
    \vrule \ \ \hfil \frb #\ 
   &\vrule \hfil \ \ \fbb #\frb\ 
   &\vrule \hfil \ \ \frb #\ \hfil
   &\vrule \hfil \ \ \frb #\ 
   &\vrule \hfil \ \ \frb #\ \ \vrule \hskip 2 pt
   &\vrule \ \ \hfil \frb #\ 
   &\vrule \hfil \ \ \fbb #\frb\ 
   &\vrule \hfil \ \ \frb #\ \hfil
   &\vrule \hfil \ \ \frb #\ 
   &\vrule \hfil \ \ \frb #\ \vrule \cr%
\noalign{\hrule}
 & &9.343.11.13.47.157&1535&506& & &5.7.13.19.29.31.421&6649&6402 \cr
11341&3257393139&4.3.5.121.23.47.307&931&884&11359&3271950955&4.3.5.7.11.29.61.97.109&6777&26468 \cr
 & &32.49.13.17.19.23.307&5833&6256& & &32.81.11.13.251.509&41229&44176 \cr
\noalign{\hrule}
 & &25.7.11.23.29.43.59&5283&3572& & &49.11.361.67.251&2619&2150 \cr
11342&3257444575&8.9.5.19.43.47.587&29&616&11360&3272235043&4.27.25.7.11.19.43.97&3263&62 \cr
 & &128.3.7.11.19.29.47&141&1216& & &16.9.13.31.43.251&3627&344 \cr
\noalign{\hrule}
 & &49.11.13.17.23.29.41&37&414& & &25.7.13.19.23.37.89&451&474 \cr
11343&3257547293&4.9.49.17.529.37&377&1210&11361&3273818275&4.3.7.11.13.19.41.79.89&11567&270 \cr
 & &16.3.5.121.13.29.37&55&888& & &16.81.5.41.43.269&11029&27864 \cr
\noalign{\hrule}
 & &11.19.29.67.71.113&4815&2756& & &9.7.121.13.19.37.47&1569&2050 \cr
11344&3258036001&8.9.5.11.13.19.53.107&29&1206&11362&3274330059&4.27.25.11.19.41.523&57421&51886 \cr
 & &32.81.29.53.67&53&1296& & &16.5.7.13.631.25943&25943&25240 \cr
\noalign{\hrule}
 & &27.49.19.103.1259&77407&52270& & &13.53.1151.4129&3663&57340 \cr
11345&3259690749&4.5.11.31.227.5227&905&6132&11363&3274458031&8.9.5.11.37.47.61&1489&1378 \cr
 & &32.3.25.7.11.73.181&13213&4400& & &32.3.5.11.13.53.1489&1489&2640 \cr
\noalign{\hrule}
 & &25.13.17.37.41.389&15423&10366& & &3.5.13.361.89.523&65191&31196 \cr
11346&3260374325&4.3.25.53.71.73.97&1683&3458&11364&3276676065&8.7.11.67.139.709&2225&2934 \cr
 & &16.27.7.11.13.17.19.73&5621&4104& & &32.9.25.89.139.163&2085&2608 \cr
\noalign{\hrule}
 & &27.19.31.421.487&7535&7562& & &27.5.7.121.23.29.43&5371&5156 \cr
11347&3260544381&4.5.11.361.137.199.421&481461&62566&11365&3279528945&8.9.7.23.41.131.1289&54725&114134 \cr
 & &16.3.7.961.41.109.167&36239&35752& & &32.25.11.149.199.383&29651&30640 \cr
\noalign{\hrule}
 & &27.7.11.31.197.257&181805&172598& & &289.19.29.43.479&641&90 \cr
11348&3262988421&4.5.13.211.409.2797&99&310&11366&3279845683&4.9.5.17.479.641&6167&4730 \cr
 & &16.9.25.11.13.31.2797&2797&2600& & &16.3.25.7.11.43.881&9691&4200 \cr
\noalign{\hrule}
 & &7.11.17.31.257.313&1593&1850& & &7.11.13.47.113.617&22579&6420 \cr
11349&3264215339&4.27.25.7.17.31.37.59&517&3172&11367&3280163887&8.3.5.7.67.107.337&403&66 \cr
 & &32.3.5.11.13.37.47.61&22607&14640& & &32.9.5.11.13.31.107&4815&496 \cr
\noalign{\hrule}
 & &3.29.151.257.967&19887&18920& & &25.11.23.61.67.127&5339&13086 \cr
11350&3264794103&16.9.5.7.11.29.43.947&19447&21274&11368&3282984925&4.9.19.23.281.727&5691&11030 \cr
 & &64.5.121.967.19447&19447&19360& & &16.27.5.7.271.1103&29781&15176 \cr
\noalign{\hrule}
 & &9.7.11.41.137.839&2197&5354& & &27.5.11.41.199.271&1717&722 \cr
11351&3265875459&4.2197.137.2677&891&890&11369&3283467165&4.3.11.17.361.41.101&71395&65258 \cr
 & &16.81.5.11.169.89.2677&120465&120328& & &16.5.67.109.131.487&63797&58424 \cr
\noalign{\hrule}
 & &13.361.29.103.233&3707&720& & &7.11.13.29.31.41.89&1125&146 \cr
11352&3266191903&32.9.5.11.13.19.337&853&1864&11370&3283731451&4.9.125.7.13.29.73&451&74 \cr
 & &512.3.5.233.853&853&3840& & &16.3.5.11.37.41.73&2701&120 \cr
\noalign{\hrule}
 & &5.7.13.19.23.41.401&127203&139462& & &11.13.31.67.11059&5731&5328 \cr
11353&3269046235&4.3.103.109.389.677&143&246&11371&3284644649&32.9.121.37.67.521&5585&13692 \cr
 & &16.9.11.13.41.109.677&7447&7848& & &256.27.5.7.163.1117&154035&142976 \cr
\noalign{\hrule}
 & &3.49.41.359.1511&4525&3014& & &27.7.13.29.193.239&817&1010 \cr
11354&3269340123&4.25.7.11.41.137.181&359&5976&11372&3286687131&4.3.5.13.19.43.101.239&79&638 \cr
 & &64.9.5.11.83.359&249&1760& & &16.5.11.19.29.79.101&5555&12008 \cr
\noalign{\hrule}
 & &9.5.23.31.97.1051&16159&10904& & &9.5.11.13.23.53.419&12449&9758 \cr
11355&3270969495&16.11.13.23.29.47.113&1275&194&11373&3286747035&4.5.7.11.17.41.59.211&1677&742 \cr
 & &64.3.25.11.17.29.97&1595&544& & &16.3.49.13.43.53.211&2107&1688 \cr
\noalign{\hrule}
 & &25.7.23.41.43.461&21699&1876& & &5.7.31.41.263.281&221&66 \cr
11356&3271290575&8.9.49.67.2411&1975&5258&11374&3287574955&4.3.11.13.17.263.281&3255&164 \cr
 & &32.3.25.11.79.239&869&11472& & &32.9.5.7.17.31.41&9&272 \cr
\noalign{\hrule}
 & &27.121.13.17.23.197&3305&3526& & &3.49.11.17.37.53.61&1075&38 \cr
11357&3271413717&4.5.11.41.43.197.661&127&324&11375&3288261669&4.25.7.11.19.37.43&4453&2862 \cr
 & &32.81.5.43.127.661&28423&30480& & &16.27.5.53.61.73&45&584 \cr
\noalign{\hrule}
 & &3.11.17.683.8539&15539&7000& & &3.13.41.47.107.409&1577&350 \cr
11358&3271828857&16.125.7.17.41.379&3069&416&11376&3288920739&4.25.7.13.19.83.107&23067&21338 \cr
 & &1024.9.25.11.13.31&10075&1536& & &16.9.5.11.47.227.233&7689&9080 \cr
\noalign{\hrule}
}%
}
$$
\eject
\vglue -23 pt
\noindent\hskip 1 in\hbox to 6.5 in{\ 11377 -- 11412 \hfill\fbd 3289235103 -- 3310284315\frb}
\vskip -9 pt
$$
\vbox{
\nointerlineskip
\halign{\strut
    \vrule \ \ \hfil \frb #\ 
   &\vrule \hfil \ \ \fbb #\frb\ 
   &\vrule \hfil \ \ \frb #\ \hfil
   &\vrule \hfil \ \ \frb #\ 
   &\vrule \hfil \ \ \frb #\ \ \vrule \hskip 2 pt
   &\vrule \ \ \hfil \frb #\ 
   &\vrule \hfil \ \ \fbb #\frb\ 
   &\vrule \hfil \ \ \frb #\ \hfil
   &\vrule \hfil \ \ \frb #\ 
   &\vrule \hfil \ \ \frb #\ \vrule \cr%
\noalign{\hrule}
 & &9.49.11.19.127.281&6665&12004& & &3.11.53.59.113.283&1705&1422 \cr
11377&3289235103&8.3.5.11.31.43.3001&791&2210&11395&3299944989&4.27.5.121.31.79.113&12169&1504 \cr
 & &32.25.7.13.17.31.113&24973&12400& & &256.31.43.47.283&2021&3968 \cr
\noalign{\hrule}
 & &3.25.121.841.431&353&78& & &37.4943.18047&2453&2490 \cr
11378&3289424325&4.9.11.13.841.353&371&470&11396&3300633877&4.3.5.11.83.223.18047&231&18278 \cr
 & &16.5.7.13.47.53.353&32383&19768& & &16.9.5.7.121.13.19.37&11495&6552 \cr
\noalign{\hrule}
 & &25.7.121.13.17.19.37&2461&2214& & &9.25.7.13.23.43.163&103&58 \cr
11379&3289811525&4.27.7.11.23.37.41.107&13&10606&11397&3300713325&4.5.13.29.43.103.163&1893&36938 \cr
 & &16.3.13.23.5303&15909&184& & &16.3.11.23.73.631&6941&584 \cr
\noalign{\hrule}
 & &25.11.17.23.37.827&13797&6878& & &3.5.19.89.157.829&807&22 \cr
11380&3290157475&4.27.7.19.23.73.181&25&1654&11398&3301330845&4.9.11.19.89.269&7475&7744 \cr
 & &16.3.25.7.19.827&399&8& & &512.25.1331.13.23&30613&16640 \cr
\noalign{\hrule}
 & &5.7.11.13.37.109.163&83709&101476& & &9.7.1331.53.743&22237&61616 \cr
11381&3290181895&8.9.23.71.131.1103&3161&148&11399&3302047287&32.37.601.3851&1625&2226 \cr
 & &64.3.29.37.71.109&213&928& & &128.3.125.7.13.37.53&1625&2368 \cr
\noalign{\hrule}
 & &9.5.7.13.17.361.131&3529&1826& & &9.25.11.19.23.43.71&281&74 \cr
11382&3292162965&4.11.361.83.3529&221&3750&11400&3302048475&4.5.11.19.37.43.281&609&1426 \cr
 & &16.3.625.13.17.83&83&1000& & &16.3.7.23.29.31.281&1967&7192 \cr
\noalign{\hrule}
 & &3.11.17.89.233.283&3913&5426& & &49.13.17.23.89.149&3701&7128 \cr
11383&3292268331&4.7.13.43.233.2713&12625&22644&11401&3302877487&16.81.11.89.3701&2555&6256 \cr
 & &32.9.125.7.17.37.101&13875&11312& & &512.9.5.7.17.23.73&657&1280 \cr
\noalign{\hrule}
 & &3.5.11.13.17.73.1237&11&62& & &27.25.11.31.113.127&91&36 \cr
11384&3292825965&4.5.121.13.31.1237&8343&7738&11402&3303241425&8.243.5.7.13.31.113&2501&6004 \cr
 & &16.81.31.53.73.103&5459&6696& & &64.13.19.41.61.79&62647&24928 \cr
\noalign{\hrule}
 & &31.107.433.2293&58707&12376& & &25.7.11.19.37.2441&611&1314 \cr
11385&3293346473&16.9.7.11.13.17.593&433&160&11403&3303344275&4.9.13.47.73.2441&1111&1330 \cr
 & &1024.3.5.11.17.433&2805&512& & &16.3.5.7.11.13.19.47.101&1313&1128 \cr
\noalign{\hrule}
 & &7.17.19.23.61.1039&715&324& & &9.7.13.29.61.2281&57695&81446 \cr
11386&3295898137&8.81.5.7.11.13.19.61&391&524&11404&3304737891&4.5.11.193.211.1049&377&588 \cr
 & &64.27.11.13.17.23.131&3537&4576& & &32.3.49.11.13.29.1049&1049&1232 \cr
\noalign{\hrule}
 & &9.13.29.43.59.383&7531&8680& & &9.13.29.37.113.233&2893&9650 \cr
11387&3296879703&16.3.5.7.17.31.59.443&8371&5362&11405&3305368989&4.3.25.11.13.193.263&1247&1262 \cr
 & &64.5.49.11.383.761&8371&7840& & &16.5.11.29.43.263.631&56545&55528 \cr
\noalign{\hrule}
 & &25.7.13.23.29.41.53&723&814& & &9.49.11.29.71.331&4693&4906 \cr
11388&3297364525&4.3.25.11.23.37.41.241&87&938&11406&3306097179&4.3.49.121.13.361.223&49319&31630 \cr
 & &16.9.7.11.29.67.241&2651&4824& & &16.5.13.149.331.3163&15815&15496 \cr
\noalign{\hrule}
 & &9.5.7.11.289.37.89&11609&24614& & &5.13.79.433.1487&95183&75852 \cr
11389&3297560805&4.7.13.19.31.47.397&1513&1266&11407&3306277585&8.9.49.11.17.43.509&433&76 \cr
 & &16.3.17.31.47.89.211&1457&1688& & &64.3.7.11.19.43.433&3311&1824 \cr
\noalign{\hrule}
 & &243.17.23.61.569&7367&5720& & &13.53.59.163.499&3325&3162 \cr
11390&3297806217&16.9.5.11.13.17.53.139&13087&3176&11408&3306430387&4.3.25.7.17.19.31.53.59&91&3036 \cr
 & &256.5.23.397.569&397&640& & &32.9.5.49.11.13.17.23&8415&18032 \cr
\noalign{\hrule}
 & &5.13.19.23.67.1733&699&1034& & &3.5.7.11.13.19.67.173&2209&2950 \cr
11391&3298132955&4.3.11.13.19.23.47.233&1733&1296&11409&3306738435&4.125.2209.59.173&54353&75978 \cr
 & &128.243.11.47.1733&2673&3008& & &16.81.7.13.37.67.113&999&904 \cr
\noalign{\hrule}
 & &9.11.47.59.61.197&173521&182780& & &9.25.11.17.31.43.59&1111&2114 \cr
11392&3298990959&8.5.13.19.37.73.2377&5025&40138&11410&3309072525&4.3.7.121.31.101.151&34715&1948 \cr
 & &32.3.125.7.47.61.67&875&1072& & &32.5.53.131.487&25811&2096 \cr
\noalign{\hrule}
 & &27.5.121.37.53.103&329&276& & &7.11.31.47.163.181&28425&27158 \cr
11393&3299392305&8.81.7.23.37.47.103&1193&7150&11411&3309912067&4.3.25.37.47.367.379&301&66 \cr
 & &32.25.11.13.47.1193&5965&9776& & &16.9.5.7.11.37.43.379&14023&15480 \cr
\noalign{\hrule}
 & &11.31.37.103.2539&19691&22230& & &3.5.13.23.29.31.821&4267&6406 \cr
11394&3299559989&4.9.5.7.13.19.29.31.97&11&206&11412&3310284315&4.5.17.29.251.3203&1529&1674 \cr
 & &16.3.11.19.29.97.103&1843&696& & &16.27.11.17.31.139.251&21267&22088 \cr
\noalign{\hrule}
}%
}
$$
\eject
\vglue -23 pt
\noindent\hskip 1 in\hbox to 6.5 in{\ 11413 -- 11448 \hfill\fbd 3310552245 -- 3337223175\frb}
\vskip -9 pt
$$
\vbox{
\nointerlineskip
\halign{\strut
    \vrule \ \ \hfil \frb #\ 
   &\vrule \hfil \ \ \fbb #\frb\ 
   &\vrule \hfil \ \ \frb #\ \hfil
   &\vrule \hfil \ \ \frb #\ 
   &\vrule \hfil \ \ \frb #\ \ \vrule \hskip 2 pt
   &\vrule \ \ \hfil \frb #\ 
   &\vrule \hfil \ \ \fbb #\frb\ 
   &\vrule \hfil \ \ \frb #\ \hfil
   &\vrule \hfil \ \ \frb #\ 
   &\vrule \hfil \ \ \frb #\ \vrule \cr%
\noalign{\hrule}
 & &3.5.7.11.13.37.59.101&617&2918& & &27.7.59.457.653&215&242 \cr
11413&3310552245&4.11.37.617.1459&30385&23598&11431&3327692571&4.5.7.121.43.59.653&1053&37474 \cr
 & &16.27.5.19.23.59.103&2369&1368& & &16.81.5.13.41.457&533&120 \cr
\noalign{\hrule}
 & &5.7.11.13.19.97.359&277&3672& & &3.5.7.31.83.97.127&1349&57706 \cr
11414&3311493185&16.27.13.17.19.277&2231&2478&11432&3328162635&4.11.19.43.61.71&601&558 \cr
 & &64.81.7.23.59.97&1357&2592& & &16.9.11.31.71.601&1803&6248 \cr
\noalign{\hrule}
 & &3.11.59.67.109.233&197&130& & &9.25.11.37.163.223&3869&206 \cr
11415&3313013253&4.5.11.13.59.197.233&441&208&11433&3328659675&4.53.73.103.223&163&60 \cr
 & &128.9.5.49.169.197&33293&47040& & &32.3.5.53.73.163&53&1168 \cr
\noalign{\hrule}
 & &7.47.199.223.227&149985&160654& & &37.223.313.1289&6435&5146 \cr
11416&3314207491&4.27.5.11.13.37.101.167&47&454&11434&3328923707&4.9.5.11.13.31.83.223&1481&1258 \cr
 & &16.9.5.13.47.101.227&585&808& & &16.3.5.13.17.31.37.1481&34255&35544 \cr
\noalign{\hrule}
 & &3.13.19.23.439.443&9265&924& & &3.5.13.19.37.149.163&5921&110 \cr
11417&3314471511&8.9.5.7.11.13.17.109&10601&9782&11435&3329383395&4.25.11.19.31.191&2093&2682 \cr
 & &32.5.67.73.10601&53005&78256& & &16.9.7.11.13.23.149&759&56 \cr
\noalign{\hrule}
 & &3.5.7.19.23.29.47.53&265&218& & &11.17.529.41.821&4711&4320 \cr
11418&3314686515&4.25.19.29.2809.109&5083&65142&11436&3329847103&64.27.5.7.23.41.673&6329&272 \cr
 & &16.9.7.11.13.17.23.47&143&408& & &2048.3.5.17.6329&6329&15360 \cr
\noalign{\hrule}
 & &13.19.23.53.101.109&453&1870& & &5.13.151.311.1091&11977&58938 \cr
11419&3314732837&4.3.5.11.17.19.53.151&327&1334&11437&3330239315&4.3.7.11.19.29.47.59&819&302 \cr
 & &16.9.5.17.23.29.109&153&1160& & &16.27.49.13.29.151&1323&232 \cr
\noalign{\hrule}
 & &3.25.7.19.29.73.157&1199&1356& & &9.5.13.53.233.461&1943&4972 \cr
11420&3315380775&8.9.5.11.19.29.109.113&1127&17512&11438&3330340065&8.3.11.29.53.67.113&6965&3688 \cr
 & &128.49.121.23.199&24079&10304& & &128.5.7.11.199.461&1393&704 \cr
\noalign{\hrule}
 & &11.13.89.463.563&4149&3170& & &243.25.11.19.43.61&1643&1582 \cr
11421&3317534363&4.9.5.317.461.463&923&460&11439&3330357525&4.81.7.11.19.31.53.113&17335&10286 \cr
 & &32.3.25.13.23.71.317&22507&27600& & &16.5.37.113.139.3467&128279&125656 \cr
\noalign{\hrule}
 & &9.5.49.41.73.503&403&38& & &9.7.169.269.1163&1193&2356 \cr
11422&3319581195&4.13.19.31.41.503&231&272&11440&3330882009&8.3.19.31.269.1193&4345&766 \cr
 & &128.3.7.11.13.17.19.31&6479&14144& & &32.5.11.31.79.383&59365&13904 \cr
\noalign{\hrule}
 & &9.25.49.13.17.29.47&1643&1672& & &3.125.7.11.23.29.173&2147&2522 \cr
11423&3320983575&16.3.5.49.11.19.31.47.53&31291&18836&11441&3331915125&4.11.13.19.97.113.173&1075&828 \cr
 & &128.13.17.19.29.83.277&5263&5312& & &32.9.25.23.43.97.113&4859&4656 \cr
\noalign{\hrule}
 & &37.41.43.127.401&7425&9016& & &9.25.13.17.113.593&2407&418 \cr
11424&3322019137&16.27.25.49.11.23.127&53&328&11442&3332022525&4.11.19.29.83.593&255&338 \cr
 & &256.9.49.23.41.53&10971&6272& & &16.3.5.11.169.17.19.29&551&1144 \cr
\noalign{\hrule}
 & &27.11.13.17.23.31.71&43757&72370& & &5.7.11.13.23.59.491&783&4618 \cr
11425&3322741851&4.5.49.19.47.7237&493&6744&11443&3334766435&4.27.7.23.29.2309&1129&5798 \cr
 & &64.3.5.7.17.29.281&8149&1120& & &16.9.13.223.1129&1129&16056 \cr
\noalign{\hrule}
 & &11.289.41.43.593&317&276& & &27.5.7.11.13.23.29.37&7171&340 \cr
11426&3323514161&8.3.11.289.23.43.317&7957&4470&11444&3334996665&8.25.13.17.71.101&6583&9108 \cr
 & &32.9.5.23.73.109.149&146169&134320& & &64.9.11.23.29.227&227&32 \cr
\noalign{\hrule}
 & &3.49.11.19.241.449&241&298& & &9.7.11.13.479.773&257&1030 \cr
11427&3324505107&4.149.58081.449&62491&4410&11445&3335735403&4.5.7.103.257.479&2319&1034 \cr
 & &16.9.5.49.11.13.19.23&345&104& & &16.3.11.47.103.773&103&376 \cr
\noalign{\hrule}
 & &9.11.13.757.3413&15365&15352& & &9.5.19.29.83.1621&7663&442 \cr
11428&3325145967&16.5.7.11.19.101.439.757&42987&1352&11446&3335993685&4.3.13.17.19.79.97&1621&2882 \cr
 & &256.3.49.169.19.23.89&56693&55936& & &16.11.17.131.1621&1441&136 \cr
\noalign{\hrule}
 & &5.17.23.157.10837&148049&158886& & &9.5.49.11.13.19.557&101&146 \cr
11429&3326254595&4.9.7.11.13.43.97.313&3613&170&11447&3336978645&4.49.11.73.101.557&589&5538 \cr
 & &16.3.5.7.17.43.3613&3613&7224& & &16.3.13.19.31.71.73&2201&584 \cr
\noalign{\hrule}
 & &5.961.67.10333&26871&24794& & &9.25.11.13.19.53.103&77&248 \cr
11430&3326554355&4.3.49.11.169.23.31.53&225&1868&11448&3337223175&16.7.121.31.53.103&4605&854 \cr
 & &32.27.25.7.11.13.467&27027&37360& & &64.3.5.49.61.307&15043&1952 \cr
\noalign{\hrule}
}%
}
$$
\eject
\vglue -23 pt
\noindent\hskip 1 in\hbox to 6.5 in{\ 11449 -- 11484 \hfill\fbd 3337370127 -- 3368210109\frb}
\vskip -9 pt
$$
\vbox{
\nointerlineskip
\halign{\strut
    \vrule \ \ \hfil \frb #\ 
   &\vrule \hfil \ \ \fbb #\frb\ 
   &\vrule \hfil \ \ \frb #\ \hfil
   &\vrule \hfil \ \ \frb #\ 
   &\vrule \hfil \ \ \frb #\ \ \vrule \hskip 2 pt
   &\vrule \ \ \hfil \frb #\ 
   &\vrule \hfil \ \ \fbb #\frb\ 
   &\vrule \hfil \ \ \frb #\ \hfil
   &\vrule \hfil \ \ \frb #\ 
   &\vrule \hfil \ \ \frb #\ \vrule \cr%
\noalign{\hrule}
 & &27.7.13.23.73.809&2503&3160& & &81.125.17.101.193&170753&160628 \cr
11449&3337370127&16.3.5.13.23.79.2503&803&1700&11467&3355232625&8.11.13.361.43.3089&92225&40602 \cr
 & &128.125.11.17.73.79&9875&11968& & &32.3.25.7.17.31.67.101&469&496 \cr
\noalign{\hrule}
 & &3.5.7.13.17.19.67.113&965&174& & &3.5.7.11.17.29.71.83&16263&6386 \cr
11450&3338016045&4.9.25.13.19.29.193&22099&28274&11468&3355562595&4.27.5.13.31.103.139&1189&2996 \cr
 & &16.49.11.41.67.211&3157&1688& & &32.7.29.41.103.107&4387&1648 \cr
\noalign{\hrule}
 & &27.11.31.37.41.239&9025&182& & &27.11.37.47.67.97&689&2900 \cr
11451&3338117541&4.25.7.13.361.41&407&372&11469&3356623017&8.9.25.13.29.47.53&2479&14746 \cr
 & &32.3.5.11.13.19.31.37&65&304& & &32.37.67.73.101&101&1168 \cr
\noalign{\hrule}
 & &3.25.11.13.37.47.179&2033&294& & &9.11.19.23.31.2503&895&1608 \cr
11452&3338488725&4.9.25.49.11.19.107&863&1612&11470&3356905959&16.27.5.11.19.67.179&32591&37424 \cr
 & &32.7.13.19.31.863&16397&3472& & &512.13.23.109.2339&30407&27904 \cr
\noalign{\hrule}
 & &81.5.7.13.31.37.79&50411&53966& & &9.97.137.28081&14089&13992 \cr
11453&3339542115&4.9.121.223.50411&14167&36244&11471&3358515681&16.27.11.53.73.137.193&42541&37330 \cr
 & &32.11.13.17.31.41.457&7769&7216& & &64.5.19.73.2239.3733&1362545&1361312 \cr
\noalign{\hrule}
 & &81.125.41.83.97&1591&484& & &3.19.29.103.109.181&39&142 \cr
11454&3342171375&8.3.5.121.37.43.97&2993&1178&11472&3359039811&4.9.13.19.29.71.109&905&1166 \cr
 & &32.19.31.37.41.73&2701&9424& & &16.5.11.13.53.71.181&3763&5720 \cr
\noalign{\hrule}
 & &9.25.11.13.37.2809&9301&874& & &5.11.169.31.89.131&411597&415658 \cr
11455&3344044275&4.3.13.19.23.71.131&745&958&11473&3359482555&4.9.19.29.37.41.83.137&11895&8492 \cr
 & &16.5.19.23.149.479&9101&27416& & &32.27.5.11.13.61.137.193&26441&26352 \cr
\noalign{\hrule}
 & &9.5.23.109.149.199&55941&52514& & &25.11.169.151.479&4747&522 \cr
11456&3345077565&4.27.7.121.29.31.643&2353&6854&11474&3361490275&4.9.29.47.101.151&2147&2600 \cr
 & &16.11.13.23.29.149.181&2353&2552& & &64.3.25.13.19.29.113&1653&3616 \cr
\noalign{\hrule}
 & &25.49.11.19.73.179&731&522& & &27.137.181.5021&55385&80182 \cr
11457&3345478675&4.9.25.7.17.29.43.73&2587&2888&11475&3361654899&4.5.11.19.47.53.853&685&168 \cr
 & &64.3.13.17.361.29.199&49153&47328& & &64.3.25.7.19.53.137&3325&1696 \cr
\noalign{\hrule}
 & &9.5.11.53.199.641&139&338& & &27.5.11.23.173.569&14873&24388 \cr
11458&3346510365&4.5.11.169.139.641&37&678&11476&3362115735&8.9.7.13.67.107.139&3337&3404 \cr
 & &16.3.13.37.113.139&1807&33448& & &64.13.23.37.47.71.139&84929&84064 \cr
\noalign{\hrule}
 & &3.5.7.23.41.73.463&781&2460& & &81.25.11.17.83.107&10207&9802 \cr
11459&3346607985&8.9.25.11.1681.71&7147&8828&11477&3363012675&4.5.169.29.59.83.173&107&972 \cr
 & &64.7.11.1021.2207&24277&32672& & &32.243.13.29.59.107&1131&944 \cr
\noalign{\hrule}
 & &3.5.53.67.227.277&4807&58072& & &11.29.61.83.2083&61957&84870 \cr
11460&3349249935&16.7.11.17.19.23.61&25&36&11478&3364247051&4.9.5.7.23.41.53.167&1331&1498 \cr
 & &128.9.25.7.17.19.23&15295&3264& & &16.3.5.49.1331.53.107&64735&62328 \cr
\noalign{\hrule}
 & &3.5.11.19.3481.307&2305&8138& & &7.11.13.23.101.1447&98841&80030 \cr
11461&3350271045&4.25.11.13.313.461&1377&6448&11479&3364742381&4.3.5.47.53.151.701&327&374 \cr
 & &128.81.169.17.31&4563&33728& & &16.9.5.11.17.53.109.151&72027&74120 \cr
\noalign{\hrule}
 & &7.11.841.59.877&391&450& & &27.25.17.23.41.311&3197&2090 \cr
11462&3350721451&4.9.25.7.11.17.23.877&4849&464&11480&3365307675&4.125.11.19.529.139&4097&1722 \cr
 & &128.3.5.13.17.29.373&6341&12480& & &16.3.7.17.41.139.241&1687&1112 \cr
\noalign{\hrule}
 & &125.7.11.19.73.251&377&4392& & &3.5.49.121.289.131&845&15006 \cr
11463&3350818625&16.9.25.7.13.29.61&251&176&11481&3366990165&4.9.25.169.41.61&143&82 \cr
 & &512.3.11.13.29.251&1131&256& & &16.11.2197.1681&1681&17576 \cr
\noalign{\hrule}
 & &81.17.23.29.41.89&21791&18142& & &27.5.13.79.107.227&287&394 \cr
11464&3351456891&4.7.11.23.47.193.283&4029&410&11482&3367548405&4.9.5.7.13.41.79.197&3859&1276 \cr
 & &16.3.5.17.41.79.283&1415&632& & &32.11.17.29.197.227&3349&5104 \cr
\noalign{\hrule}
 & &3.25.11.17.353.677&249&3634& & &3.25.11.361.43.263&5789&9734 \cr
11465&3351708525&4.9.5.17.23.79.83&3883&3172&11483&3368102925&4.5.7.11.31.157.827&931&774 \cr
 & &32.11.13.23.61.353&299&976& & &16.9.343.19.43.827&2481&2744 \cr
\noalign{\hrule}
 & &9.23.479.33827&17017&16810& & &3.11.31.71.79.587&1535&666 \cr
11466&3354048531&4.5.7.11.13.17.1681.479&1311&18328&11484&3368210109&4.27.5.37.307.587&6715&15004 \cr
 & &64.3.5.19.23.29.41.79&16195&17632& & &32.25.121.17.31.79&187&400 \cr
\noalign{\hrule}
}%
}
$$
\eject
\vglue -23 pt
\noindent\hskip 1 in\hbox to 6.5 in{\ 11485 -- 11520 \hfill\fbd 3370919013 -- 3397778265\frb}
\vskip -9 pt
$$
\vbox{
\nointerlineskip
\halign{\strut
    \vrule \ \ \hfil \frb #\ 
   &\vrule \hfil \ \ \fbb #\frb\ 
   &\vrule \hfil \ \ \frb #\ \hfil
   &\vrule \hfil \ \ \frb #\ 
   &\vrule \hfil \ \ \frb #\ \ \vrule \hskip 2 pt
   &\vrule \ \ \hfil \frb #\ 
   &\vrule \hfil \ \ \fbb #\frb\ 
   &\vrule \hfil \ \ \frb #\ \hfil
   &\vrule \hfil \ \ \frb #\ 
   &\vrule \hfil \ \ \frb #\ \vrule \cr%
\noalign{\hrule}
 & &9.7.11.31.173.907&965&58& & &3.7.121.17.29.37.73&1431&802 \cr
11485&3370919013&4.3.5.7.29.173.193&4535&482&11503&3383577813&4.81.11.53.73.401&21895&43148 \cr
 & &16.25.241.907&25&1928& & &32.5.7.23.29.67.151&3473&5360 \cr
\noalign{\hrule}
 & &81.5.7.11.13.53.157&11&466& & &5.289.691.3389&91377&108322 \cr
11486&3373375005&4.9.121.157.233&10547&8450&11504&3383899555&4.9.11.13.41.71.1321&4913&950 \cr
 & &16.25.169.53.199&995&104& & &16.3.25.4913.19.71&1349&2040 \cr
\noalign{\hrule}
 & &27.5.7.11.23.103.137&77507&77714& & &9.5.37.113.17987&8901&9086 \cr
11487&3373728435&4.3.5.343.13.61.179.433&1507&208&11505&3384164115&4.81.7.11.23.43.59.113&4819&40 \cr
 & &128.11.169.61.137.179&10919&10816& & &64.5.7.11.23.61.79&9821&27808 \cr
\noalign{\hrule}
 & &3.5.17.19.43.97.167&5709&2536& & &25.7.11.169.101.103&783&328 \cr
11488&3374818665&16.9.11.43.173.317&193&194&11506&3384355975&16.27.5.13.29.41.103&19&84 \cr
 & &64.11.97.173.193.317&61181&60896& & &128.81.7.19.29.41&22591&5184 \cr
\noalign{\hrule}
 & &3.25.11.53.79.977&26483&25298& & &25.343.13.97.313&1683&6142 \cr
11489&3374826675&4.5.7.11.13.71.139.373&43&738&11507&3384492475&4.9.11.17.37.83.97&545&2194 \cr
 & &16.9.7.13.41.43.373&45879&31304& & &16.3.5.37.109.1097&4033&26328 \cr
\noalign{\hrule}
 & &3.25.7.11.17.31.1109&111&230& & &9.11.19.529.41.83&845&98 \cr
11490&3375158325&4.9.125.23.37.1109&7579&33454&11508&3386151747&4.5.49.11.169.19.23&719&282 \cr
 & &16.11.13.43.53.389&5057&18232& & &16.3.5.7.13.47.719&33793&3640 \cr
\noalign{\hrule}
 & &11.169.23.89.887&20935&21822& & &5.11.17.43.109.773&4075&3966 \cr
11491&3375365851&4.3.5.53.79.89.3637&2041&1596&11509&3387552685&4.3.125.163.661.773&5559&102184 \cr
 & &32.9.7.13.19.53.79.157&111627&112784& & &64.9.17.53.109.241&2169&1696 \cr
\noalign{\hrule}
 & &9.49.13.19.29.1069&449&620& & &27.25.7.11.13.29.173&1601&2546 \cr
11492&3376845927&8.5.49.13.29.31.449&297&1718&11510&3389861475&4.5.19.67.173.1601&44187&13768 \cr
 & &32.27.11.449.859&9449&21552& & &64.3.11.13.103.1721&1721&3296 \cr
\noalign{\hrule}
 & &11.13.19.23.191.283&36001&33900& & &9.13.149.163.1193&2849&730 \cr
11493&3377826023&8.3.25.7.23.37.113.139&543&382&11511&3390003747&4.3.5.7.11.37.73.149&3523&1990 \cr
 & &32.9.113.139.181.191&20453&20016& & &16.25.11.13.199.271&6775&17512 \cr
\noalign{\hrule}
 & &11.37.67.229.541&989&1530& & &27.7.11.17.23.43.97&1885&1426 \cr
11494&3378329141&4.9.5.17.23.37.43.67&2717&164&11512&3390560019&4.5.13.529.29.31.97&5047&11352 \cr
 & &32.3.5.11.13.17.19.41&2337&17680& & &64.3.49.11.29.43.103&721&928 \cr
\noalign{\hrule}
 & &9.25.7.11.41.67.71&53153&52372& & &11.23.37.281.1289&785&504 \cr
11495&3379016025&8.23.41.2311.13093&197945&103194&11513&3390638449&16.9.5.7.11.23.37.157&1289&438 \cr
 & &32.81.5.49.11.13.59.61&6903&6832& & &64.27.5.7.73.1289&1971&1120 \cr
\noalign{\hrule}
 & &27.25.11.169.2693&703&1990& & &3.113.151.173.383&19&132 \cr
11496&3379243725&4.3.125.13.19.37.199&2693&5068&11514&3391731951&8.9.11.19.173.383&3367&80 \cr
 & &32.7.37.181.2693&1267&592& & &256.5.7.11.13.37&407&58240 \cr
\noalign{\hrule}
 & &11.13.19.31.53.757&87&670& & &9.47.73.181.607&1625&1806 \cr
11497&3379271467&4.3.5.13.19.29.31.67&477&74&11515&3392583093&4.27.125.7.13.43.607&1991&1384 \cr
 & &16.27.5.37.53.67&335&7992& & &64.7.11.13.43.173.181&15743&15136 \cr
\noalign{\hrule}
 & &11.37.41.167.1213&58851&9118& & &9.25.7.37.139.419&9503&9922 \cr
11498&3380302277&4.9.13.47.97.503&1525&3034&11516&3393994275&4.3.121.13.17.41.43.139&35&1564 \cr
 & &16.3.25.13.37.41.61&1525&312& & &32.5.7.11.289.23.43&10879&4624 \cr
\noalign{\hrule}
 & &5.49.41.487.691&94631&74664& & &11.13.31.347.2207&4275&6482 \cr
11499&3380313265&16.9.17.61.173.547&627&9926&11517&3394919957&4.9.25.7.11.13.19.463&4399&694 \cr
 & &64.27.7.11.19.709&5643&22688& & &16.3.5.7.53.83.347&1113&3320 \cr
\noalign{\hrule}
 & &27.49.11.31.59.127&725&598& & &27.5.13.17.19.53.113&511&506 \cr
11500&3380414499&4.25.11.13.23.29.31.59&3&1708&11518&3394954485&4.3.7.11.13.17.19.23.53.73&745&34394 \cr
 & &32.3.5.7.13.23.61&115&12688& & &16.5.29.73.149.593&43289&34568 \cr
\noalign{\hrule}
 & &5.11.23.29.37.47.53&14193&7378& & &9.13.361.29.47.59&1075&9394 \cr
11501&3381146395&4.9.7.17.19.23.31.83&1361&50&11519&3396572829&4.3.25.7.11.13.43.61&3337&2812 \cr
 & &16.3.25.7.31.1361&1361&26040& & &32.19.37.47.61.71&2257&1136 \cr
\noalign{\hrule}
 & &5.43.79.263.757&1629&58174& & &3.5.7.19.37.191.241&847&22048 \cr
11502&3381560635&4.9.17.29.59.181&1595&1414&11520&3397778265&64.49.121.13.53&243&296 \cr
 & &16.3.5.7.11.841.101&5887&26664& & &1024.243.11.13.37&1053&5632 \cr
\noalign{\hrule}
}%
}
$$
\eject
\vglue -23 pt
\noindent\hskip 1 in\hbox to 6.5 in{\ 11521 -- 11556 \hfill\fbd 3398052735 -- 3434237695\frb}
\vskip -9 pt
$$
\vbox{
\nointerlineskip
\halign{\strut
    \vrule \ \ \hfil \frb #\ 
   &\vrule \hfil \ \ \fbb #\frb\ 
   &\vrule \hfil \ \ \frb #\ \hfil
   &\vrule \hfil \ \ \frb #\ 
   &\vrule \hfil \ \ \frb #\ \ \vrule \hskip 2 pt
   &\vrule \ \ \hfil \frb #\ 
   &\vrule \hfil \ \ \fbb #\frb\ 
   &\vrule \hfil \ \ \frb #\ \hfil
   &\vrule \hfil \ \ \frb #\ 
   &\vrule \hfil \ \ \frb #\ \vrule \cr%
\noalign{\hrule}
 & &27.5.49.11.17.41.67&925&884& & &9.29.53.107.2311&75205&47278 \cr
11521&3398052735&8.125.49.11.13.289.37&8091&28034&11539&3420582741&4.5.7.11.169.89.307&27&116 \cr
 & &32.9.13.29.31.107.131&49387&53072& & &32.27.5.7.13.29.307&3991&1680 \cr
\noalign{\hrule}
 & &5.19.73.113.4339&11921&9774& & &3.49.13.17.31.43.79&583&1610 \cr
11522&3400279045&4.27.7.13.73.131.181&13017&3454&11540&3421108509&4.5.343.11.23.31.53&2243&1530 \cr
 & &16.81.11.157.4339&891&1256& & &16.9.25.17.53.2243&6729&10600 \cr
\noalign{\hrule}
 & &3.11.23.67.211.317&7063&7074& & &3.11.13.53.151.997&65&518 \cr
11523&3401404611&4.81.7.23.131.317.1009&1235&24442&11541&3422987139&4.5.7.169.37.997&583&414 \cr
 & &16.5.7.121.13.19.101.131&101101&99560& & &16.9.5.7.11.23.37.53&483&1480 \cr
\noalign{\hrule}
 & &9.13.29.59.73.233&733&34& & &5.11.19.29.37.43.71&343&1692 \cr
11524&3404980683&4.3.17.29.73.733&26015&27494&11542&3423283105&8.9.343.29.43.47&1775&528 \cr
 & &16.5.121.43.59.233&605&344& & &256.27.25.7.11.71&945&128 \cr
\noalign{\hrule}
 & &3.125.11.13.17.37.101&263&4888& & &343.19.23.73.313&613&774 \cr
11525&3406742625&16.11.169.47.263&2525&5418&11543&3424859459&4.9.49.43.313.613&35207&5170 \cr
 & &64.9.25.7.43.101&43&672& & &16.3.5.11.17.19.47.109&15369&7480 \cr
\noalign{\hrule}
 & &5.7.11.13.19.73.491&18103&28542& & &9.5.11.43.227.709&24427&63422 \cr
11526&3408490085&4.3.7.43.67.71.421&33&2914&11544&3425671755&4.13.19.1669.1879&18685&17016 \cr
 & &16.9.11.31.47.71&423&17608& & &64.3.5.13.37.101.709&1313&1184 \cr
\noalign{\hrule}
 & &11.13.43.2209.251&29547&65440& & &3.5.13.19.53.73.239&561&634 \cr
11527&3409368391&64.9.5.49.67.409&559&1786&11545&3425980155&4.9.11.13.17.19.53.317&84875&14818 \cr
 & &256.3.7.13.19.43.47&133&384& & &16.125.7.31.97.239&2425&1736 \cr
\noalign{\hrule}
 & &5.13.41.73.89.197&47&42& & &25.11.29.43.97.103&453&4882 \cr
11528&3410957485&4.3.7.13.41.47.73.197&5335&19716&11546&3426163675&4.3.5.29.151.2441&3913&8292 \cr
 & &32.9.5.7.11.31.53.97&46269&38192& & &32.9.7.13.43.691&8983&1008 \cr
\noalign{\hrule}
 & &27.5.7.11.19.23.751&4213&4964& & &3.289.101.109.359&775&1078 \cr
11529&3411503865&8.9.5.121.17.73.383&413&1502&11547&3426584277&4.25.49.11.17.31.359&61&2574 \cr
 & &32.7.17.59.73.751&1241&944& & &16.9.5.7.121.13.61&2379&33880 \cr
\noalign{\hrule}
 & &5.7.13.23.31.67.157&21663&30932& & &9.11.17.23.29.43.71&14705&3826 \cr
11530&3412521385&8.9.7.11.19.29.37.83&3875&5338&11548&3427178733&4.5.289.173.1913&231&58 \cr
 & &32.3.125.17.29.31.157&725&816& & &16.3.5.7.11.29.1913&1913&280 \cr
\noalign{\hrule}
 & &9.13.23.53.71.337&123445&108748& & &27.11.361.113.283&21629&62422 \cr
11531&3412540521&8.5.7.31.877.3527&429&3956&11549&3428692443&4.529.43.59.503&1017&1520 \cr
 & &64.3.7.11.13.23.31.43&2387&1376& & &128.9.5.19.529.113&529&320 \cr
\noalign{\hrule}
 & &3.5.11.29.373.1913&1435&478& & &9.49.11.19.29.1283&965&916 \cr
11532&3414331965&4.25.7.41.239.373&5499&3826&11550&3429331983&8.5.29.193.229.1283&409&6006 \cr
 & &16.9.13.41.47.1913&611&984& & &32.3.7.11.13.229.409&5317&3664 \cr
\noalign{\hrule}
 & &3.7.11.961.15383&54321&53360& & &13.17.19.59.109.127&2115&44 \cr
11533&3414887553&32.9.5.11.19.23.29.953&53857&55738&11551&3429478663&8.9.5.11.13.47.59&1241&592 \cr
 & &128.841.961.53857&53857&53824& & &256.3.5.17.37.73&365&14208 \cr
\noalign{\hrule}
 & &27.361.29.43.281&3293&2046& & &191.223.239.337&14751&60400 \cr
11534&3415417029&4.81.11.19.31.37.89&365&6844&11552&3430567999&32.9.25.11.149.151&299&148 \cr
 & &32.5.29.37.59.73&10915&1168& & &256.3.25.11.13.23.37&18975&61568 \cr
\noalign{\hrule}
 & &11.31.191.52441&6345&58786& & &25.49.17.37.61.73&101&528 \cr
11535&3415534771&4.27.5.7.13.17.19.47&229&382&11553&3431147825&32.3.25.7.11.73.101&3599&3774 \cr
 & &16.3.5.7.19.191.229&285&56& & &128.9.11.17.37.59.61&649&576 \cr
\noalign{\hrule}
 & &3.11.29.53.89.757&1625&646& & &243.7.13.29.53.101&13561&16490 \cr
11536&3417225933&4.125.13.17.19.29.53&801&736&11554&3432755781&4.3.5.13.17.71.97.191&22781&69854 \cr
 & &256.9.25.17.19.23.89&8075&8832& & &16.11.19.53.109.659&12521&9592 \cr
\noalign{\hrule}
 & &81.49.11.13.19.317&347&3140& & &3.5.11.169.23.53.101&10207&12530 \cr
11537&3418456041&8.27.5.13.157.347&217&568&11555&3433173315&4.25.7.13.59.173.179&371&396 \cr
 & &128.7.31.71.347&2201&22208& & &32.9.49.11.53.173.179&8477&8592 \cr
\noalign{\hrule}
 & &191.66049.271&58905&7144& & &5.7.37.41.71.911&2001&484 \cr
11538&3418762289&16.9.5.7.11.17.19.47&257&542&11556&3434237695&8.3.121.23.29.911&615&296 \cr
 & &64.3.7.11.257.271&33&224& & &128.9.5.11.23.37.41&207&704 \cr
\noalign{\hrule}
}%
}
$$
\eject
\vglue -23 pt
\noindent\hskip 1 in\hbox to 6.5 in{\ 11557 -- 11592 \hfill\fbd 3435039993 -- 3468568641\frb}
\vskip -9 pt
$$
\vbox{
\nointerlineskip
\halign{\strut
    \vrule \ \ \hfil \frb #\ 
   &\vrule \hfil \ \ \fbb #\frb\ 
   &\vrule \hfil \ \ \frb #\ \hfil
   &\vrule \hfil \ \ \frb #\ 
   &\vrule \hfil \ \ \frb #\ \ \vrule \hskip 2 pt
   &\vrule \ \ \hfil \frb #\ 
   &\vrule \hfil \ \ \fbb #\frb\ 
   &\vrule \hfil \ \ \frb #\ \hfil
   &\vrule \hfil \ \ \frb #\ 
   &\vrule \hfil \ \ \frb #\ \vrule \cr%
\noalign{\hrule}
 & &3.49.11.43.127.389&12695&25532& & &3.5.43.59.269.337&7093&8778 \cr
11557&3435039993&8.5.7.13.491.2539&5695&12078&11575&3449799915&4.9.7.11.19.41.43.173&337&50 \cr
 & &32.9.25.11.17.61.67&12261&6800& & &16.25.11.19.173.337&865&1672 \cr
\noalign{\hrule}
 & &7.11.19.31.239.317&4205&1818& & &121.17.19.529.167&17621&21462 \cr
11558&3436079339&4.9.5.841.101.239&779&62&11576&3452709469&4.3.49.23.67.73.263&1691&150 \cr
 & &16.3.5.19.31.41.101&205&2424& & &16.9.25.7.19.73.89&5607&14600 \cr
\noalign{\hrule}
 & &5.17.19.23.37.41.61&66703&31812& & &3.121.13.17.193.223&1125&932 \cr
11559&3437286865&8.3.7.11.13.241.733&53253&51566&11577&3452717697&8.27.125.13.223.233&17177&14278 \cr
 & &32.27.19.23.59.61.97&1593&1552& & &32.25.121.59.89.193&1475&1424 \cr
\noalign{\hrule}
 & &27.5.157.241.673&70499&91694& & &25.7.13.17.19.37.127&253&6 \cr
11560&3437680635&4.11.13.17.361.29.127&2205&3856&11578&3452942675&4.3.25.11.17.23.127&4189&5586 \cr
 & &128.9.5.49.17.19.241&931&1088& & &16.9.49.19.59.71&413&5112 \cr
\noalign{\hrule}
 & &5.11.127.677.727&471067&474702& & &9.25.7.11.13.529.29&205&2 \cr
11561&3437870315&4.3.19.61.1297.24793&75&24718&11579&3455176725&4.125.11.13.23.41&939&686 \cr
 & &16.9.25.17.61.727&1037&360& & &16.3.343.41.313&12833&392 \cr
\noalign{\hrule}
 & &5.7.13.23.37.83.107&261&154& & &9.5.11.13.19.59.479&7047&7526 \cr
11562&3438767605&4.9.49.11.13.23.29.37&11451&12118&11580&3455331165&4.2187.5.11.29.53.71&1333&114578 \cr
 & &16.27.121.73.83.347&25331&26136& & &16.31.43.59.971&1333&7768 \cr
\noalign{\hrule}
 & &25.13.289.19.41.47&11803&13248& & &9.7.11.113.131.337&379&412 \cr
11563&3438876025&128.9.5.11.19.23.29.37&12763&9248&11581&3457107423&8.3.103.131.337.379&365&28 \cr
 & &8192.3.289.12763&12763&12288& & &64.5.7.73.103.379&27667&16480 \cr
\noalign{\hrule}
 & &81.7.121.13.17.227&745&844& & &3.25.49.887.1061&793&268 \cr
11564&3441807369&8.9.5.11.13.17.149.211&545&1886&11582&3458568225&8.7.13.61.67.887&1061&5148 \cr
 & &32.25.23.41.109.211&102787&84400& & &64.9.11.169.1061&169&1056 \cr
\noalign{\hrule}
 & &11.13.29.43.97.199&113&2700& & &9.5.17.529.83.103&2915&2812 \cr
11565&3442130263&8.27.25.11.43.113&1843&1886&11583&3459652065&8.3.25.11.17.19.23.37.53&349&824 \cr
 & &32.9.25.19.23.41.97&7011&9200& & &128.11.37.53.103.349&21571&22336 \cr
\noalign{\hrule}
 & &3.5.11.169.37.47.71&655&1084& & &27.11.31.67.71.79&10465&12382 \cr
11566&3442933065&8.25.13.71.131.271&1739&36&11584&3460018221&4.5.7.13.23.41.79.151&201&352 \cr
 & &64.9.37.47.271&813&32& & &256.3.5.11.13.23.41.67&2665&2944 \cr
\noalign{\hrule}
 & &9.25.11.17.19.59.73&8533&11542& & &9.5.7.31.127.2791&441551&437614 \cr
11567&3443123475&4.3.7.19.23.29.53.199&73&130&11585&3461272605&4.11.17.61.137.211.293&310909&7068 \cr
 & &16.5.13.23.53.73.199&4577&5512& & &32.3.19.29.31.71.151&10721&8816 \cr
\noalign{\hrule}
 & &27.5.7.31.41.47.61&853&8602& & &49.19.29.41.53.59&89&2508 \cr
11568&3443539365&4.11.17.23.47.853&7791&6710&11586&3461460793&8.3.11.361.29.89&695&3276 \cr
 & &16.3.5.49.121.53.61&847&424& & &64.27.5.7.13.139&351&22240 \cr
\noalign{\hrule}
 & &81.5.49.197.881&36173&35188& & &11.101.197.15817&101525&117342 \cr
11569&3444238665&8.49.19.61.463.593&407&28650&11587&3461819339&4.27.25.31.41.53.131&9191&2248 \cr
 & &32.3.25.11.19.37.191&7067&16720& & &64.3.25.7.13.101.281&6825&8992 \cr
\noalign{\hrule}
 & &361.29.73.4507&2947&1560& & &81.5.17.313.1607&2639&4246 \cr
11570&3444416159&16.3.5.7.13.19.29.421&2981&876&11588&3463093035&4.7.11.13.29.193.313&3151&918 \cr
 & &128.9.11.13.73.271&2981&7488& & &16.27.17.23.137.193&3151&1544 \cr
\noalign{\hrule}
 & &9.25.11.13.43.47.53&1295&1196& & &5.13.67.131.6073&201761&192984 \cr
11571&3446360775&8.125.7.169.23.37.43&4631&744&11589&3464676865&16.3.7.11.17.19.37.41.43&195&12146 \cr
 & &128.3.7.11.31.37.421&13051&16576& & &64.9.5.11.13.6073&99&32 \cr
\noalign{\hrule}
 & &25.37.41.211.431&297&728& & &81.49.13.529.127&1705&8582 \cr
11572&3448937425&16.27.7.11.13.37.211&3325&1004&11590&3466441251&4.5.343.11.31.613&1361&1704 \cr
 & &128.3.25.49.19.251&12299&3648& & &64.3.11.31.71.1361&42191&24992 \cr
\noalign{\hrule}
 & &11.13.31.151.5153&7749&59240& & &9.5.7.13.41.73.283&19657&1002 \cr
11573&3449330599&16.27.5.7.41.1481&2201&2242&11591&3468542805&4.27.11.167.1787&25043&23206 \cr
 & &64.9.5.7.19.31.59.71&29323&27360& & &16.41.79.283.317&317&632 \cr
\noalign{\hrule}
 & &27.11.19.29.107.197&25241&15932& & &3.31.43.59.61.241&3175&11044 \cr
11574&3449515113&8.9.7.43.569.587&1639&1070&11592&3468568641&8.25.11.31.127.251&153&122 \cr
 & &32.5.11.107.149.587&2935&2384& & &32.9.17.61.127.251&4267&6096 \cr
\noalign{\hrule}
}%
}
$$
\eject
\vglue -23 pt
\noindent\hskip 1 in\hbox to 6.5 in{\ 11593 -- 11628 \hfill\fbd 3469578281 -- 3500997325\frb}
\vskip -9 pt
$$
\vbox{
\nointerlineskip
\halign{\strut
    \vrule \ \ \hfil \frb #\ 
   &\vrule \hfil \ \ \fbb #\frb\ 
   &\vrule \hfil \ \ \frb #\ \hfil
   &\vrule \hfil \ \ \frb #\ 
   &\vrule \hfil \ \ \frb #\ \ \vrule \hskip 2 pt
   &\vrule \ \ \hfil \frb #\ 
   &\vrule \hfil \ \ \fbb #\frb\ 
   &\vrule \hfil \ \ \frb #\ \hfil
   &\vrule \hfil \ \ \frb #\ 
   &\vrule \hfil \ \ \frb #\ \vrule \cr%
\noalign{\hrule}
 & &169.43.223.2141&379&2520& & &27.5.29.73.89.137&28237&28922 \cr
11593&3469578281&16.9.5.7.13.43.379&2563&3122&11611&3484698435&4.11.17.89.151.14461&511&13950 \cr
 & &64.3.49.11.223.233&2563&4704& & &16.9.25.7.11.17.31.73&1085&1496 \cr
\noalign{\hrule}
 & &9.7.17.23.29.43.113&6911&4312& & &9.5.11.13.59.67.137&5549&2534 \cr
11594&3471060663&16.343.11.17.6911&6371&540&11612&3484945035&4.7.11.13.31.179.181&411&590 \cr
 & &128.27.5.11.23.277&4155&704& & &16.3.5.31.59.137.181&181&248 \cr
\noalign{\hrule}
 & &11.41.53.337.431&195&236& & &9.25.11.19.37.2003&2143&2132 \cr
11595&3471839041&8.3.5.11.13.53.59.337&123&460&11613&3485069775&8.13.37.41.2003.2143&80707&1416 \cr
 & &64.9.25.13.23.41.59&13275&9568& & &128.3.121.13.23.29.59&17641&20416 \cr
\noalign{\hrule}
 & &5.49.11.31.43.967&1677&710& & &3.13.29.283.10891&9549&1342 \cr
11596&3473884645&4.3.25.7.13.1849.71&837&1012&11614&3485915043&4.27.11.13.61.1061&6011&5660 \cr
 & &32.81.11.13.23.31.71&5751&4784& & &32.5.61.283.6011&6011&4880 \cr
\noalign{\hrule}
 & &9.11.13.41.199.331&7443&21014& & &5.49.47.97.3121&1683&1438 \cr
11597&3475709523&4.81.7.19.79.827&697&130&11615&3486016555&4.9.11.17.47.97.719&49&1600 \cr
 & &16.5.13.17.19.41.79&395&2584& & &512.3.25.49.719&719&3840 \cr
\noalign{\hrule}
 & &11.17.53.97.3617&170469&180380& & &3.5.37.61.239.431&5499&3344 \cr
11598&3477264439&8.9.5.13.29.31.47.311&867&8152&11616&3487369695&32.27.11.13.19.47.61&5735&4088 \cr
 & &128.27.13.289.1019&17323&22464& & &512.5.7.13.31.37.73&6643&7936 \cr
\noalign{\hrule}
 & &81.5.31.47.71.83&159&76& & &169.19.53.103.199&9207&8200 \cr
11599&3477370905&8.243.19.31.53.71&979&1222&11617&3488240951&16.27.25.11.31.41.199&6461&106 \cr
 & &32.11.13.19.47.53.89&13091&15664& & &64.9.5.7.13.53.71&315&2272 \cr
\noalign{\hrule}
 & &9.5.19.23.149.1187&3871&2684& & &9.5.19.31.353.373&5159&16102 \cr
11600&3478010895&8.3.49.11.61.79.149&437&10&11618&3489886845&4.3.5.7.11.67.83.97&4051&2984 \cr
 & &32.5.7.11.19.23.79&77&1264& & &64.83.373.4051&4051&2656 \cr
\noalign{\hrule}
 & &9.5.7.37.457.653&22363&1798& & &81.25.11.19.73.113&9877&7052 \cr
11601&3478096755&4.7.11.19.29.31.107&5319&4570&11619&3491183025&8.7.17.41.43.73.83&261&1502 \cr
 & &16.27.5.19.197.457&197&456& & &32.9.7.29.83.751&21779&9296 \cr
\noalign{\hrule}
 & &3.5.7.11.13.19.89.137&10371&1466& & &3.5.7.11.169.31.577&69&646 \cr
11602&3478480005&4.9.11.733.3457&11525&19588&11620&3491452965&4.9.7.13.17.19.23.31&14425&4412 \cr
 & &32.25.59.83.461&27199&6640& & &32.25.577.1103&1103&80 \cr
\noalign{\hrule}
 & &9.5.11.53.97.1367&59033&2482& & &27.7.17.61.71.251&14365&946 \cr
11603&3478734765&4.13.17.19.73.239&1517&1590&11621&3492791253&4.5.11.169.289.43&147&142 \cr
 & &16.3.5.17.19.37.41.53&1517&2584& & &16.3.49.11.169.43.71&1859&2408 \cr
\noalign{\hrule}
 & &25.11.29.43.73.139&109&36& & &3.5.11.41.397.1301&1533&2834 \cr
11604&3479659975&8.9.5.11.43.109.139&111&584&11622&3494102205&4.9.5.7.13.41.73.109&3553&4988 \cr
 & &128.27.37.73.109&999&6976& & &32.11.17.19.29.43.109&23693&29648 \cr
\noalign{\hrule}
 & &3.5.313.397.1867&151&1716& & &9.5.29.43.199.313&1477&4294 \cr
11605&3479929305&8.9.11.13.151.397&20633&18670&11623&3495235005&4.5.7.19.43.113.211&7593&12452 \cr
 & &32.5.47.439.1867&439&752& & &32.3.7.11.283.2531&21791&40496 \cr
\noalign{\hrule}
 & &7.11.29.37.103.409&2925&62& & &27.11.31.47.59.137&25&1482 \cr
11606&3480574867&4.9.25.11.13.31.37&1939&2494&11624&3497748507&4.81.25.13.19.59&2603&2662 \cr
 & &16.3.5.7.29.43.277&4155&344& & &16.5.1331.361.137&1805&968 \cr
\noalign{\hrule}
 & &5.7.41.43.73.773&201909&193094& & &25.13.29.37.79.127&12587&14238 \cr
11607&3481951445&4.3.11.17.37.67.107.131&773&366&11625&3498757925&4.9.7.41.79.113.307&63373&40700 \cr
 & &16.9.61.107.131.773&7991&7704& & &32.3.25.11.37.127.499&499&528 \cr
\noalign{\hrule}
 & &11.13.29.67.83.151&707&4854& & &27.25.13.43.9277&41209&42284 \cr
11608&3482281517&4.3.7.101.151.809&7221&8030&11626&3500444025&8.3.49.11.13.841.961&469&430 \cr
 & &16.9.5.7.11.29.73.83&511&360& & &32.5.343.11.29.31.43.67&60233&60368 \cr
\noalign{\hrule}
 & &3.5.17.19.41.89.197&1661&676& & &81.5.7.13.17.37.151&7777&5210 \cr
11609&3482842785&8.11.169.17.89.151&1773&794&11627&3500451045&4.3.25.49.11.101.521&169&1394 \cr
 & &32.9.169.197.397&1191&2704& & &16.11.169.17.41.101&4141&1144 \cr
\noalign{\hrule}
 & &3.5.343.11.19.41.79&8427&8084& & &25.49.23.137.907&1559&2466 \cr
11610&3482912895&8.9.5.41.43.47.2809&119&2054&11628&3500997325&4.9.7.18769.1559&6985&25754 \cr
 & &32.7.13.17.47.53.79&2491&3536& & &16.3.5.11.79.127.163&30099&14344 \cr
\noalign{\hrule}
}%
}
$$
\eject
\vglue -23 pt
\noindent\hskip 1 in\hbox to 6.5 in{\ 11629 -- 11664 \hfill\fbd 3501885165 -- 3537170091\frb}
\vskip -9 pt
$$
\vbox{
\nointerlineskip
\halign{\strut
    \vrule \ \ \hfil \frb #\ 
   &\vrule \hfil \ \ \fbb #\frb\ 
   &\vrule \hfil \ \ \frb #\ \hfil
   &\vrule \hfil \ \ \frb #\ 
   &\vrule \hfil \ \ \frb #\ \ \vrule \hskip 2 pt
   &\vrule \ \ \hfil \frb #\ 
   &\vrule \hfil \ \ \fbb #\frb\ 
   &\vrule \hfil \ \ \frb #\ \hfil
   &\vrule \hfil \ \ \frb #\ 
   &\vrule \hfil \ \ \frb #\ \vrule \cr%
\noalign{\hrule}
 & &3.5.17.37.47.53.149&17193&17822& & &3.121.19.457.1117&4007&2890 \cr
11629&3501885165&4.9.7.11.19.53.67.521&17123&10490&11647&3520704693&4.5.289.457.4007&1881&5888 \cr
 & &16.5.7.19.1049.17123&139517&136984& & &2048.9.5.11.17.19.23&1955&3072 \cr
\noalign{\hrule}
 & &9.13.19.23.61.1123&7955&6644& & &19.89.1103.1889&99&1790 \cr
11630&3502489887&8.3.5.11.37.43.61.151&8611&10868&11648&3523311797&4.9.5.11.179.1103&2489&3026 \cr
 & &64.5.121.13.19.79.109&13189&12640& & &16.3.11.17.19.89.131&393&1496 \cr
\noalign{\hrule}
 & &19.37.53.109.863&61183&15444& & &9.13.19.29.31.41.43&583&196 \cr
11631&3504842353&8.27.11.13.17.59.61&1585&2014&11649&3523314951&8.49.11.13.29.31.53&5289&4600 \cr
 & &32.9.5.17.19.53.317&1585&2448& & &128.3.25.49.23.41.43&1225&1472 \cr
\noalign{\hrule}
 & &361.23.197.2143&973&1170& & &9.25.11.17.83.1009&111289&94136 \cr
11632&3505285813&4.9.5.7.13.361.23.139&3641&3940&11650&3523655025&16.7.1681.109.1021&1339&5808 \cr
 & &32.3.25.11.139.197.331&24825&24464& & &512.3.121.13.41.103&14729&10496 \cr
\noalign{\hrule}
 & &5.7.19.23.29.41.193&3069&5254& & &25.343.19.43.503&2721&3796 \cr
11633&3509850715&4.9.11.31.37.71.193&41&152&11651&3523904825&8.3.13.73.503.907&24255&12464 \cr
 & &64.3.11.19.31.41.71&2343&992& & &256.27.5.49.11.19.41&1107&1408 \cr
\noalign{\hrule}
 & &3.5.17.19.691.1049&1243&2212& & &25.11.13.19.23.37.61&6663&2698 \cr
11634&3511941855&8.7.11.79.113.1049&959&90&11652&3526054675&4.3.5.361.71.2221&2013&208 \cr
 & &32.9.5.49.113.137&15481&2352& & &128.9.11.13.61.71&639&64 \cr
\noalign{\hrule}
 & &27.5.11.13.19.61.157&1321&406& & &9.7.13.19.31.71.103&3175&5132 \cr
11635&3512795715&4.9.7.13.19.29.1321&451&1772&11653&3527725383&8.25.7.31.127.1283&209&426 \cr
 & &32.7.11.29.41.443&18163&3248& & &32.3.5.11.19.71.1283&1283&880 \cr
\noalign{\hrule}
 & &9.5.7.11.29.73.479&14291&15886& & &3.5.7.13.17.19.53.151&183&506 \cr
11636&3513658995&4.169.31.47.73.461&25543&3876&11654&3528482685&4.9.5.7.11.23.61.151&347&1102 \cr
 & &32.3.7.13.17.19.41.89&13243&18512& & &16.11.19.29.61.347&3817&14152 \cr
\noalign{\hrule}
 & &27.5.19.29.97.487&6149&59596& & &5.49.11.13.23.29.151&639&488 \cr
11637&3513873015&8.11.13.43.47.317&3071&1050&11655&3528620095&16.9.5.11.13.29.61.71&229&2114 \cr
 & &32.3.25.7.11.37.83&15355&1232& & &64.3.7.61.151.229&687&1952 \cr
\noalign{\hrule}
 & &9.5.7.31.41.67.131&13&54& & &11.19.54289.311&5355&59644 \cr
11638&3514003605&4.243.5.7.13.31.131&6059&1474&11656&3528730711&8.9.5.7.13.17.31.37&233&418 \cr
 & &16.11.13.67.73.83&949&7304& & &32.3.11.13.17.19.233&39&272 \cr
\noalign{\hrule}
 & &27.125.7.23.29.223&4927&202& & &9.5.31.53.59.809&66671&30266 \cr
11639&3514006125&4.5.13.29.101.379&759&1136&11657&3528991485&4.121.19.29.37.409&4045&3726 \cr
 & &128.3.11.23.71.101&7171&704& & &16.81.5.11.23.37.809&851&792 \cr
\noalign{\hrule}
 & &9.7.11.29.179.977&167&370& & &3.11.13.17.23.53.397&31&190 \cr
11640&3514623651&4.3.5.11.37.167.977&571&406&11658&3529396299&4.5.11.19.23.31.397&8957&4590 \cr
 & &16.7.29.37.167.571&6179&4568& & &16.27.25.169.17.53&325&72 \cr
\noalign{\hrule}
 & &5.13.29.47.97.409&121&2934& & &49.11.13.47.71.151&11&60 \cr
11641&3514829435&4.9.121.163.409&205&204&11659&3530736209&8.3.5.121.13.47.151&3825&1862 \cr
 & &32.27.5.121.17.41.163&74817&79376& & &32.27.125.49.17.19&8721&2000 \cr
\noalign{\hrule}
 & &5.7.11.31.271.1087&3657&4744& & &5.11.13.31.37.59.73&63&4370 \cr
11642&3515776495&16.3.5.7.11.23.53.593&3523&558&11660&3532192235&4.9.25.7.19.23.37&413&438 \cr
 & &64.27.13.23.31.271&351&736& & &16.27.49.19.59.73&931&216 \cr
\noalign{\hrule}
 & &25.23.43.71.2003&1539&464& & &5.289.19.179.719&50623&78078 \cr
11643&3516216425&32.81.19.23.29.71&8255&6622&11661&3533485955&4.3.7.11.169.23.31.71&3401&3060 \cr
 & &128.3.5.7.11.13.43.127&4191&5824& & &32.27.5.13.17.19.23.179&351&368 \cr
\noalign{\hrule}
 & &25.53.67.173.229&1343&2208& & &81.5.13.43.67.233&1199&4228 \cr
11644&3516999175&64.3.5.17.23.79.229&1749&2144&11662&3534252345&8.5.7.11.43.109.151&1809&6496 \cr
 & &4096.9.11.23.53.67&2277&2048& & &512.27.49.29.67&1421&256 \cr
\noalign{\hrule}
 & &81.5.7.71.101.173&22763&27434& & &13.29.101.131.709&35739&35870 \cr
11645&3517052805&4.3.5.11.13.17.29.43.103&1211&2806&11663&3536553683&4.9.5.11.13.17.361.29.211&48601&287944 \cr
 & &16.7.17.23.43.61.173&2623&3128& & &64.3.7.53.131.35993&35993&35616 \cr
\noalign{\hrule}
 & &3.5.11.19.1849.607&621&1228& & &9.7.13.31.127.1097&116885&131146 \cr
11646&3518545305&8.81.5.11.19.23.307&2159&5536&11664&3537170091&4.5.23.97.241.2851&2541&310 \cr
 & &512.17.23.127.173&67643&32512& & &16.3.25.7.121.31.241&3025&1928 \cr
\noalign{\hrule}
}%
}
$$
\eject
\vglue -23 pt
\noindent\hskip 1 in\hbox to 6.5 in{\ 11665 -- 11700 \hfill\fbd 3537553841 -- 3567692403\frb}
\vskip -9 pt
$$
\vbox{
\nointerlineskip
\halign{\strut
    \vrule \ \ \hfil \frb #\ 
   &\vrule \hfil \ \ \fbb #\frb\ 
   &\vrule \hfil \ \ \frb #\ \hfil
   &\vrule \hfil \ \ \frb #\ 
   &\vrule \hfil \ \ \frb #\ \ \vrule \hskip 2 pt
   &\vrule \ \ \hfil \frb #\ 
   &\vrule \hfil \ \ \fbb #\frb\ 
   &\vrule \hfil \ \ \frb #\ \hfil
   &\vrule \hfil \ \ \frb #\ 
   &\vrule \hfil \ \ \frb #\ \vrule \cr%
\noalign{\hrule}
 & &41.47.977.1879&61479&15560& & &3.11.23.29.31.41.127&2353&1330 \cr
11665&3537553841&16.243.5.11.23.389&29&52&11683&3552949587&4.5.7.13.19.23.41.181&3937&11358 \cr
 & &128.3.5.11.13.29.389&64185&24128& & &16.9.13.31.127.631&631&312 \cr
\noalign{\hrule}
 & &5.121.19.37.53.157&1755&4054& & &9.5.7.17.19.47.743&14069&20852 \cr
11666&3539046115&4.27.25.13.53.2027&11&64&11684&3553037145&8.3.5.11.13.401.1279&629&574 \cr
 & &512.9.11.13.2027&18243&3328& & &32.7.13.17.37.41.1279&19721&20464 \cr
\noalign{\hrule}
 & &9.13.23.29.67.677&19833&200& & &27.5.11.19.29.43.101&1519&1046 \cr
11667&3539771001&16.27.25.11.601&1651&1354&11685&3553594605&4.49.29.31.101.523&2025&2924 \cr
 & &64.5.13.127.677&635&32& & &32.81.25.17.43.523&1569&1360 \cr
\noalign{\hrule}
 & &5.7.13.53.181.811&2871&2806& & &3.11.17.23.61.4517&10933&60620 \cr
11668&3539864965&4.9.11.23.29.53.61.181&563&20&11686&3555253911&8.5.7.13.841.433&427&414 \cr
 & &32.3.5.23.29.61.563&38847&28304& & &32.9.5.49.23.61.433&2165&2352 \cr
\noalign{\hrule}
 & &25.31.67.79.863&60051&8126& & &9.5.7.11.59.127.137&1909&5626 \cr
11669&3540090725&4.3.17.37.239.541&177&9020&11687&3556964565&4.23.29.83.97.127&49&78 \cr
 & &32.9.5.11.41.59&26609&144& & &16.3.49.13.23.83.97&15617&8632 \cr
\noalign{\hrule}
 & &49.19.1949.1951&66285&29216& & &3.49.11.13.31.53.103&2773&1660 \cr
11670&3540126569&64.27.5.11.83.491&377&868&11688&3557362809&8.5.7.47.59.83.103&2809&2088 \cr
 & &512.9.7.11.13.29.31&11687&25344& & &128.9.5.29.47.2809&7685&9024 \cr
\noalign{\hrule}
 & &3.49.11.23.31.37.83&255&596& & &25.7.11.37.47.1063&3723&22852 \cr
11671&3540620391&8.9.5.49.17.83.149&43&790&11689&3558472225&8.3.7.17.29.73.197&19375&20616 \cr
 & &32.25.43.79.149&3397&59600& & &128.9.625.31.859&21475&17856 \cr
\noalign{\hrule}
 & &7.41.1487.8297&26335&34632& & &13.17.53.313.971&8855&60318 \cr
11672&3540902393&16.9.5.7.13.23.37.229&517&16318&11690&3559850099&4.27.5.7.11.23.1117&1549&1802 \cr
 & &64.3.11.41.47.199&9353&1056& & &16.9.5.7.17.53.1549&1549&2520 \cr
\noalign{\hrule}
 & &3.5.47.71.193.367&2717&2788& & &3.19.107.163.3581&3339&242 \cr
11673&3545445705&8.11.13.17.19.41.47.193&355&162&11691&3560004597&4.27.7.121.53.107&1259&1630 \cr
 & &32.81.5.13.17.19.41.71&8721&8528& & &16.5.121.163.1259&1259&4840 \cr
\noalign{\hrule}
 & &17.19.23.29.109.151&5295&2134& & &5.121.17.47.53.139&217&2274 \cr
11674&3545943419&4.3.5.11.97.151.353&9265&5382&11692&3561170965&4.3.5.7.31.139.379&1207&3102 \cr
 & &16.27.25.13.17.23.109&351&200& & &16.9.7.11.17.47.71&639&56 \cr
\noalign{\hrule}
 & &3.5.7.19.67.139.191&1711&374& & &9.11.13.37.239.313&15613&15374 \cr
11675&3548672085&4.11.17.19.29.59.67&3753&200&11693&3562234533&4.169.37.1201.7687&717&6970 \cr
 & &64.27.25.29.139&1305&32& & &16.3.5.17.41.239.1201&6005&5576 \cr
\noalign{\hrule}
 & &3.25.13.17.19.59.191&6615&6424& & &81.169.31.37.227&5017&2020 \cr
11676&3548889825&16.81.125.49.11.19.73&11269&5356&11694&3564191241&8.5.169.29.101.173&99&70 \cr
 & &128.7.11.13.59.103.191&721&704& & &32.9.25.7.11.101.173&27775&19376 \cr
\noalign{\hrule}
 & &7.11.19.37.173.379&1833&1454& & &9.7.11.17.41.47.157&8675&60328 \cr
11677&3549207277&4.3.7.11.13.37.47.727&101&5190&11695&3564211959&16.25.347.7541&3597&3944 \cr
 & &16.9.5.47.101.173&4545&376& & &256.3.25.11.17.29.109&2725&3712 \cr
\noalign{\hrule}
 & &3.31.41.281.3313&5605&4334& & &27.5.19.443.3137&4927&4484 \cr
11678&3549723789&4.5.11.19.59.197.281&1169&2574&11696&3564557415&8.9.5.13.361.59.379&803&2608 \cr
 & &16.9.7.121.13.59.167&45591&57112& & &256.11.13.59.73.163&105787&121472 \cr
\noalign{\hrule}
 & &9.25.49.11.19.23.67&865&262& & &11.13.31.47.71.241&4461&2260 \cr
11679&3550810725&4.125.11.19.131.173&557&1932&11697&3565093961&8.3.5.113.241.1487&177&64 \cr
 & &32.3.7.23.173.557&557&2768& & &1024.9.5.59.1487&66915&30208 \cr
\noalign{\hrule}
 & &27.31.37.73.1571&935&2506& & &3.5.11.31.181.3853&931&2922 \cr
11680&3551617827&4.9.5.7.11.17.73.179&185&1426&11698&3567165195&4.9.5.49.19.31.487&7391&7706 \cr
 & &16.25.7.11.23.31.37&1771&200& & &16.7.361.389.3853&2723&2888 \cr
\noalign{\hrule}
 & &11.436921.739&222525&214396& & &27.7.13.17.223.383&667&3124 \cr
11681&3551730809&8.9.25.7.13.19.23.31.43&661&1478&11699&3567448521&8.11.23.29.71.383&4795&4014 \cr
 & &32.3.25.7.13.661.739&325&336& & &32.9.5.7.29.137.223&685&464 \cr
\noalign{\hrule}
 & &121.17.59.73.401&3379&3438& & &9.11.17.23.37.47.53&475&1274 \cr
11682&3552659099&4.9.121.31.73.109.191&23659&2840&11700&3567692403&4.3.25.49.13.19.23.37&319&578 \cr
 & &64.3.5.31.59.71.401&1065&992& & &16.25.7.11.289.19.29&5075&2584 \cr
\noalign{\hrule}
}%
}
$$
\eject
\vglue -23 pt
\noindent\hskip 1 in\hbox to 6.5 in{\ 11701 -- 11736 \hfill\fbd 3568514719 -- 3602991645\frb}
\vskip -9 pt
$$
\vbox{
\nointerlineskip
\halign{\strut
    \vrule \ \ \hfil \frb #\ 
   &\vrule \hfil \ \ \fbb #\frb\ 
   &\vrule \hfil \ \ \frb #\ \hfil
   &\vrule \hfil \ \ \frb #\ 
   &\vrule \hfil \ \ \frb #\ \ \vrule \hskip 2 pt
   &\vrule \ \ \hfil \frb #\ 
   &\vrule \hfil \ \ \fbb #\frb\ 
   &\vrule \hfil \ \ \frb #\ \hfil
   &\vrule \hfil \ \ \frb #\ 
   &\vrule \hfil \ \ \frb #\ \vrule \cr%
\noalign{\hrule}
 & &343.11.89.10627&1887&8740& & &5.13.47.79.89.167&123&44 \cr
11701&3568514719&8.3.5.49.17.19.23.37&745&1068&11719&3587110735&8.3.5.11.13.41.47.89&2449&1734 \cr
 & &64.9.25.23.89.149&3427&7200& & &32.9.289.31.41.79&8959&5904 \cr
\noalign{\hrule}
 & &9.25.7.11.289.23.31&2249&226& & &9.5.11.13.53.67.157&6031&4304 \cr
11702&3569937525&4.13.23.31.113.173&19575&14212&11720&3587557545&32.3.37.67.163.269&1643&836 \cr
 & &32.27.25.11.17.19.29&87&304& & &256.11.19.31.53.163&3097&3968 \cr
\noalign{\hrule}
 & &25.11.23.509.1109&152437&140238& & &3.5.11.37.73.83.97&91001&108946 \cr
11703&3570342325&4.27.49.19.53.71.113&2277&1486&11721&3588048915&4.17.19.47.53.61.101&1825&666 \cr
 & &16.243.7.11.19.23.743&14117&13608& & &16.9.25.17.37.73.101&505&408 \cr
\noalign{\hrule}
 & &11.13.29.53.71.229&272375&265734& & &3.11.29.67.191.293&10647&2150 \cr
11704&3573581869&4.27.125.7.19.37.2179&377&1802&11722&3588291597&4.27.25.7.11.169.43&1465&992 \cr
 & &16.9.5.7.13.17.29.37.53&1295&1224& & &256.125.13.31.293&3875&1664 \cr
\noalign{\hrule}
 & &729.11.19.29.809&10975&10166& & &7.269.1009.1889&437&1446 \cr
11705&3574541421&4.25.11.13.17.19.23.439&261&2456&11723&3588999883&4.3.19.23.241.1889&1345&3234 \cr
 & &64.9.5.17.23.29.307&5219&3680& & &16.9.5.49.11.23.269&1449&440 \cr
\noalign{\hrule}
 & &5.11.17.41.179.521&759&8098& & &7.13.17.89.131.199&3&4 \cr
11706&3575083765&4.3.5.121.23.4049&2327&1722&11724&3589258127&8.3.13.17.89.131.199&23331&5620 \cr
 & &16.9.7.13.23.41.179&207&728& & &64.9.5.7.11.101.281&27819&16160 \cr
\noalign{\hrule}
 & &9.5.11.13.191.2909&527&428& & &25.11.17.67.73.157&12085&1566 \cr
11707&3575408265&8.13.17.31.107.2909&20567&17250&11725&3589871725&4.27.125.29.2417&1771&646 \cr
 & &32.3.125.17.23.131.157&55675&57776& & &16.3.7.11.17.19.23.29&483&4408 \cr
\noalign{\hrule}
 & &9.5.11.29.41.59.103&793&752& & &27.25.19.281.997&4011&3014 \cr
11708&3576648735&32.3.11.13.29.47.59.61&89&1858&11726&3593013525&4.81.7.11.19.137.191&281&1820 \cr
 & &128.13.47.89.929&43663&74048& & &32.5.49.13.137.281&1781&784 \cr
\noalign{\hrule}
 & &9.49.13.37.101.167&751&418& & &3.25.7.11.503.1237&631&606 \cr
11709&3577844907&4.7.11.13.19.101.751&309&1420&11727&3593268525&4.9.7.11.101.503.631&193693&123700 \cr
 & &32.3.5.71.103.751&36565&12016& & &32.25.109.1237.1777&1777&1744 \cr
\noalign{\hrule}
 & &7.59.1667.5197&30987&67366& & &5.7.61.71.151.157&44109&74426 \cr
11710&3577983787&4.9.11.13.313.2591&1765&826&11728&3593625595&4.9.11.169.17.29.199&11&210 \cr
 & &16.3.5.7.11.13.59.353&1765&3432& & &16.27.5.7.121.13.29&10179&968 \cr
\noalign{\hrule}
 & &9.7.13.97.107.421&3365&418& & &5.23.149.349.601&24707&64842 \cr
11711&3578668821&4.3.5.11.19.107.673&12121&10088&11729&3594049115&4.3.31.101.107.797&2057&1260 \cr
 & &64.5.13.17.23.31.97&1955&992& & &32.27.5.7.121.17.101&19089&32912 \cr
\noalign{\hrule}
 & &9.25.7.121.89.211&703&98& & &3.5.7.11.13.19.43.293&1545&1678 \cr
11712&3578807925&4.5.343.19.37.211&699&356&11730&3594305715&4.9.25.13.43.103.839&10291&616 \cr
 & &32.3.19.37.89.233&703&3728& & &64.7.11.41.103.251&4223&8032 \cr
\noalign{\hrule}
 & &3.5.11.13.43.151.257&551&1110& & &27.25.23.29.61.131&101783&101652 \cr
11713&3579363645&4.9.25.19.29.37.257&57233&66742&11731&3597747975&8.81.5.11.19.43.197.487&6623&14318 \cr
 & &16.1331.13.17.43.151&121&136& & &32.11.37.179.197.7159&1281461&1282864 \cr
\noalign{\hrule}
 & &9.5.13.17.47.79.97&341&244& & &9.11.29.1087.1153&533&620 \cr
11714&3581801145&8.11.17.31.47.61.79&35&834&11732&3598255881&8.3.5.11.13.31.41.1087&133&1220 \cr
 & &32.3.5.7.31.61.139&1891&15568& & &64.25.7.13.19.31.61&70525&37088 \cr
\noalign{\hrule}
 & &9.5.13.23.41.43.151&27401&33894& & &5.49.121.13.9337&8601&736 \cr
11715&3581895915&4.81.7.11.47.53.269&205&124&11733&3598339745&64.3.49.23.47.61&533&594 \cr
 & &32.5.11.31.41.53.269&8339&9328& & &256.81.11.13.41.47&3321&6016 \cr
\noalign{\hrule}
 & &9.11.17.29.97.757&6175&15778& & &3.23.29.61.163.181&77&106 \cr
11716&3583849203&4.25.343.13.17.19.23&319&1866&11734&3601165683&4.7.11.23.53.163.181&10485&19018 \cr
 & &16.3.5.49.11.29.311&311&1960& & &16.9.5.11.37.233.257&59881&48840 \cr
\noalign{\hrule}
 & &125.11.13.19.61.173&2139&236& & &5.121.19.37.43.197&82253&56238 \cr
11717&3584062625&8.3.13.23.31.59.61&1075&282&11735&3602843365&4.3.7.13.83.103.991&2167&10716 \cr
 & &32.9.25.31.43.47&11997&752& & &32.9.7.11.19.47.197&329&144 \cr
\noalign{\hrule}
 & &3.5.7.11.19.29.43.131&7501&6254& & &27.5.11.17.41.3481&20203&2798 \cr
11718&3584869365&4.11.13.19.53.59.577&9271&1692&11736&3602991645&4.9.89.227.1399&6409&6182 \cr
 & &32.9.47.59.73.127&27813&44368& & &16.11.13.17.29.89.281&8149&9256 \cr
\noalign{\hrule}
}%
}
$$
\eject
\vglue -23 pt
\noindent\hskip 1 in\hbox to 6.5 in{\ 11737 -- 11772 \hfill\fbd 3603321183 -- 3628938599\frb}
\vskip -9 pt
$$
\vbox{
\nointerlineskip
\halign{\strut
    \vrule \ \ \hfil \frb #\ 
   &\vrule \hfil \ \ \fbb #\frb\ 
   &\vrule \hfil \ \ \frb #\ \hfil
   &\vrule \hfil \ \ \frb #\ 
   &\vrule \hfil \ \ \frb #\ \ \vrule \hskip 2 pt
   &\vrule \ \ \hfil \frb #\ 
   &\vrule \hfil \ \ \fbb #\frb\ 
   &\vrule \hfil \ \ \frb #\ \hfil
   &\vrule \hfil \ \ \frb #\ 
   &\vrule \hfil \ \ \frb #\ \vrule \cr%
\noalign{\hrule}
 & &3.49.11.37.229.263&745&10476& & &9.13.89.449.773&625&176 \cr
11737&3603321183&8.81.5.11.97.149&3367&8702&11755&3614112801&32.625.11.13.773&4887&5162 \cr
 & &32.7.13.19.37.229&247&16& & &128.27.25.29.89.181&5249&4800 \cr
\noalign{\hrule}
 & &7.11.13.17.191.1109&661&1770& & &27.5.7.11.17.113.181&659&14596 \cr
11738&3604523923&4.3.5.7.59.191.661&999&338&11756&3614351895&8.17.41.89.659&427&1086 \cr
 & &16.81.5.169.37.59&5265&17464& & &32.3.7.41.61.181&41&976 \cr
\noalign{\hrule}
 & &9.7.11.19.251.1091&1675&584& & &9.5.19.47.293.307&8533&5896 \cr
11739&3605664447&16.25.7.11.19.67.73&6697&3012&11757&3614680935&16.5.7.11.19.23.53.67&47&162 \cr
 & &128.3.5.37.181.251&905&2368& & &64.81.7.47.53.67&3551&2016 \cr
\noalign{\hrule}
 & &625.11.17.59.523&5051&5574& & &11.19.37.43.83.131&3029&540 \cr
11740&3606411875&4.3.11.59.929.5051&375&55186&11758&3615479087&8.27.5.11.13.37.233&5561&3526 \cr
 & &16.9.125.41.673&369&5384& & &32.9.41.43.67.83&369&1072 \cr
\noalign{\hrule}
 & &7.11.13.23.383.409&1269&1412& & &9.5.169.19.29.863&2891&10054 \cr
11741&3606483881&8.27.23.47.353.409&3115&6292&11759&3616276365&4.3.49.11.13.59.457&5&44 \cr
 & &64.3.5.7.121.13.47.89&4895&4512& & &32.5.121.59.457&26963&1936 \cr
\noalign{\hrule}
 & &9.5.7.13.23.149.257&6061&23494& & &9.7.19.677.4463&14605&16636 \cr
11742&3606626205&4.7.11.17.19.29.691&771&1462&11760&3616676847&8.3.5.19.23.127.4159&3443&58942 \cr
 & &16.3.289.19.43.257&817&2312& & &32.11.13.313.2267&24937&65104 \cr
\noalign{\hrule}
 & &27.13.19.557.971&535&22& & &3.5.49.11.47.89.107&2533&1828 \cr
11743&3606908643&4.5.11.13.107.971&557&414&11761&3618692385&8.11.17.107.149.457&1729&90 \cr
 & &16.9.5.23.107.557&115&856& & &32.9.5.7.13.19.457&5941&912 \cr
\noalign{\hrule}
 & &9.5.7.19.37.43.379&691&506& & &3.7.11.17.19.139.349&789&3050 \cr
11744&3608889165&4.11.23.43.379.691&16095&202&11762&3619551243&4.9.25.61.139.263&5423&3056 \cr
 & &16.3.5.11.29.37.101&101&2552& & &128.25.11.17.29.191&4775&1856 \cr
\noalign{\hrule}
 & &13.29.41.43.5431&119441&103230& & &5.11.113.349.1669&8307&10052 \cr
11745&3609719581&4.9.5.7.31.37.113.151&5431&17974&11763&3620119415&8.9.7.13.71.113.359&349&10 \cr
 & &16.3.7.11.19.43.5431&209&168& & &32.3.5.7.13.71.349&213&1456 \cr
\noalign{\hrule}
 & &49.307.389.617&333&284& & &5.7.11.23.29.59.239&5283&3572 \cr
11746&3610515559&8.9.37.71.307.389&18095&3702&11764&3621066295&8.9.19.47.239.587&413&174 \cr
 & &32.27.5.7.11.47.617&1269&880& & &32.27.7.19.29.47.59&513&752 \cr
\noalign{\hrule}
 & &5.11.13.23.29.67.113&361&306& & &25.23.71.107.829&231&124 \cr
11747&3610647755&4.9.13.17.361.67.113&1175&2644&11765&3621299975&8.3.5.7.11.23.31.829&817&12 \cr
 & &32.3.25.17.19.47.661&45543&52880& & &64.9.11.19.31.43&1333&60192 \cr
\noalign{\hrule}
 & &27.5.7.11.37.41.229&158813&169802& & &9.7.13.67.149.443&463&580 \cr
11748&3611150235&4.31.47.59.109.1439&667&2106&11766&3622002111&8.5.29.67.443.463&2079&136 \cr
 & &16.81.13.23.29.31.109&27807&25288& & &128.27.7.11.17.463&7871&2112 \cr
\noalign{\hrule}
 & &5.7.17.29.47.61.73&519&722& & &5.7.11.13.41.127.139&8601&9052 \cr
11749&3611316205&4.3.5.361.47.61.173&5159&5394&11767&3622483865&8.3.5.7.13.31.47.61.73&9279&4826 \cr
 & &16.9.7.11.361.29.31.67&24187&24552& & &32.27.19.47.127.1031&19589&20304 \cr
\noalign{\hrule}
 & &5.61.157.241.313&759&446& & &27.5.169.19.61.137&341&172 \cr
11750&3612109205&4.3.11.23.61.157.223&17841&17170&11768&3622634145&8.5.11.31.43.61.137&221&84 \cr
 & &16.9.5.17.19.23.101.313&7429&7272& & &64.3.7.11.13.17.31.43&5117&10912 \cr
\noalign{\hrule}
 & &3.5.343.11.29.31.71&5039&5256& & &9.5.7.97.139.853&4609&8874 \cr
11751&3612402255&16.27.49.11.73.5039&4757&9796&11769&3622814685&4.81.7.11.17.29.419&2641&292 \cr
 & &128.31.67.71.73.79&4891&5056& & &32.11.17.19.73.139&803&5168 \cr
\noalign{\hrule}
 & &9.11.17.6859.313&8233&2096& & &3.25.7.11.17.19.29.67&43&518 \cr
11752&3613177161&32.3.19.131.8233&383&7850&11770&3624326475&4.49.29.37.43.67&1105&2178 \cr
 & &128.25.157.383&60131&1600& & &16.9.5.121.13.17.43&129&1144 \cr
\noalign{\hrule}
 & &3.13.31.83.36013&19625&16388& & &5.11.17.31.41.43.71&891&2162 \cr
11753&3613796511&8.125.17.31.157.241&12331&25506&11771&3628139405&4.81.5.121.17.23.47&4601&7384 \cr
 & &32.9.5.11.13.19.59.109&19293&16720& & &64.27.13.43.71.107&1391&864 \cr
\noalign{\hrule}
 & &9.31.41.193.1637&11323&3410& & &11.13.89.389.733&123&856 \cr
11754&3614049099&4.5.11.169.961.67&123&838&11772&3628938599&16.3.13.41.107.389&501&890 \cr
 & &16.3.13.41.67.419&5447&536& & &64.9.5.41.89.167&1503&6560 \cr
\noalign{\hrule}
}%
}
$$
\eject
\vglue -23 pt
\noindent\hskip 1 in\hbox to 6.5 in{\ 11773 -- 11808 \hfill\fbd 3631894497 -- 3660242517\frb}
\vskip -9 pt
$$
\vbox{
\nointerlineskip
\halign{\strut
    \vrule \ \ \hfil \frb #\ 
   &\vrule \hfil \ \ \fbb #\frb\ 
   &\vrule \hfil \ \ \frb #\ \hfil
   &\vrule \hfil \ \ \frb #\ 
   &\vrule \hfil \ \ \frb #\ \ \vrule \hskip 2 pt
   &\vrule \ \ \hfil \frb #\ 
   &\vrule \hfil \ \ \fbb #\frb\ 
   &\vrule \hfil \ \ \frb #\ \hfil
   &\vrule \hfil \ \ \frb #\ 
   &\vrule \hfil \ \ \frb #\ \vrule \cr%
\noalign{\hrule}
 & &27.7.121.31.47.109&395&3356& & &7.11.17.53.131.401&167&750 \cr
11773&3631894497&8.3.5.79.109.839&203&124&11791&3644443187&4.3.125.17.167.401&1219&1620 \cr
 & &64.5.7.29.31.839&4195&928& & &32.243.625.23.53&5589&10000 \cr
\noalign{\hrule}
 & &17.47.73.199.313&344045&338724& & &3.7.11.23.29.41.577&47&530 \cr
11774&3633013849&8.9.5.13.67.79.9409&313&22&11792&3644999589&4.5.11.29.41.47.53&13271&3636 \cr
 & &32.3.11.13.79.97.313&7663&6864& & &32.9.23.101.577&303&16 \cr
\noalign{\hrule}
 & &9.29.53.59.61.73&34775&69596& & &9.5.23.31.41.47.59&3515&4972 \cr
11775&3634302591&8.25.13.107.127.137&371&264&11793&3647839905&8.25.11.19.37.59.113&12121&5454 \cr
 & &128.3.5.7.11.13.53.137&7535&5824& & &32.27.11.17.23.31.101&1111&816 \cr
\noalign{\hrule}
 & &7.23.97.211.1103&17537&16434& & &81.25.11.13.43.293&12901&16124 \cr
11776&3634591261&4.9.11.13.19.71.83.97&1055&4838&11794&3648355425&8.3.7.13.19.29.97.139&829&2090 \cr
 & &16.3.5.11.19.41.59.211&7257&8360& & &32.5.11.361.29.829&10469&13264 \cr
\noalign{\hrule}
 & &81.11.13.31.53.191&911&1190& & &3.7.13.97.211.653&28475&29128 \cr
11777&3634895979&4.9.5.7.13.17.53.911&191&880&11795&3648631623&16.25.11.17.67.97.331&3717&2782 \cr
 & &128.25.11.191.911&911&1600& & &64.9.5.7.13.59.107.331&31565&31776 \cr
\noalign{\hrule}
 & &3.7.11.13.17.19.23.163&3277&6014& & &11.19.43.47.53.163&3309&5330 \cr
11778&3636413781&4.11.13.29.31.97.113&801&460&11796&3649018571&4.3.5.11.13.19.41.1103&21253&23970 \cr
 & &32.9.5.23.29.89.113&10057&6960& & &16.9.25.17.47.53.401&3609&3400 \cr
\noalign{\hrule}
 & &3.25.7.121.19.23.131&853&1898& & &729.11.23.47.421&845&424 \cr
11779&3636615675&4.5.11.13.23.73.853&11367&7102&11797&3649454919&16.27.5.11.169.23.53&47&668 \cr
 & &16.27.13.53.67.421&46163&30312& & &128.13.47.53.167&8851&832 \cr
\noalign{\hrule}
 & &121.13.59.113.347&1005&238& & &5.7.11.13.17.59.727&1363&1068 \cr
11780&3639055277&4.3.5.7.11.17.67.347&195&542&11798&3649550905&8.3.7.29.47.89.727&3835&1254 \cr
 & &16.9.25.7.13.17.271&17073&3400& & &32.9.5.11.13.19.47.59&423&304 \cr
\noalign{\hrule}
 & &9.13.23.47.107.269&42097&73570& & &9.7.11.23.29.53.149&25051&10300 \cr
11781&3640387491&4.5.7.11.43.89.1051&23&66&11799&3650238207&8.25.7.13.41.47.103&10281&23606 \cr
 & &16.3.5.7.121.23.1051&4235&8408& & &32.3.11.23.29.37.149&37&16 \cr
\noalign{\hrule}
 & &7.11.13.17.31.67.103&14513&20610& & &3.5.7.37.61.73.211&13&198 \cr
11782&3640463827&4.9.5.17.23.229.631&1417&13096&11800&3650279955&4.27.7.11.13.61.73&2321&2132 \cr
 & &64.3.5.13.109.1637&4911&17440& & &32.121.169.41.211&4961&2704 \cr
\noalign{\hrule}
 & &49.13.53.269.401&3047&450& & &625.11.41.12959&9917&3042 \cr
11783&3641765309&4.9.25.11.277.401&4531&1484&11801&3652818125&4.9.169.41.47.211&5225&4692 \cr
 & &32.3.5.7.23.53.197&345&3152& & &32.27.25.11.13.17.19.23&5967&6992 \cr
\noalign{\hrule}
 & &17.41.43.137.887&813&950& & &3.5.11.47.67.79.89&6077&6472 \cr
11784&3642045949&4.3.25.17.19.271.887&4521&86&11802&3653202135&16.11.59.67.103.809&60525&6322 \cr
 & &16.9.5.11.19.43.137&95&792& & &64.9.25.29.109.269&23403&17440 \cr
\noalign{\hrule}
 & &9.5.7.37.47.61.109&1023&7672& & &3.11.17.23.179.1583&1267&2850 \cr
11785&3642222465&16.27.49.11.31.137&2609&1090&11803&3656155371&4.9.25.7.11.17.19.181&179&8 \cr
 & &64.5.11.109.2609&2609&352& & &64.25.7.179.181&4525&224 \cr
\noalign{\hrule}
 & &3.5.17.19.31.79.307&2627&14872& & &729.25.11.13.23.61&8159&7366 \cr
11786&3642679335&16.11.169.17.37.71&279&350&11804&3656463525&4.27.11.29.41.127.199&105053&45950 \cr
 & &64.9.25.7.11.169.31&1859&3360& & &16.25.13.919.8081&8081&7352 \cr
\noalign{\hrule}
 & &81.125.7.11.4673&3349&1324& & &11.19.23.47.97.167&477&40 \cr
11787&3643187625&8.5.7.11.17.197.331&3107&3438&11805&3659823871&16.9.5.53.97.167&4183&4668 \cr
 & &32.9.13.191.197.239&45649&40976& & &128.27.47.89.389&10503&5696 \cr
\noalign{\hrule}
 & &7.23.41.73.7561&3861&3700& & &9.5.7.11.19.23.2417&1343&538 \cr
11788&3643441753&8.27.25.11.13.37.41.73&4837&10128&11806&3659833485&4.17.79.269.2417&3495&1078 \cr
 & &256.81.5.7.211.691&85455&88448& & &16.3.5.49.11.79.233&553&1864 \cr
\noalign{\hrule}
 & &27.5.7.11.13.59.457&76183&46478& & &9.5.37.97.131.173&253&6148 \cr
11789&3643645005&4.17.29.37.71.1367&1997&630&11807&3660187815&8.11.23.29.53.97&519&548 \cr
 & &16.9.5.7.17.29.1997&1997&3944& & &64.3.23.53.137.173&3151&1696 \cr
\noalign{\hrule}
 & &7.11.29.41.53.751&561&190& & &9.19.529.43.941&4147&3206 \cr
11790&3644084059&4.3.5.121.17.19.29.41&2205&2756&11808&3660242517&4.7.11.13.529.29.229&4555&2322 \cr
 & &32.27.25.49.13.17.53&4725&3536& & &16.27.5.43.229.911&4555&5496 \cr
\noalign{\hrule}
}%
}
$$
\eject
\vglue -23 pt
\noindent\hskip 1 in\hbox to 6.5 in{\ 11809 -- 11844 \hfill\fbd 3661452795 -- 3685716441\frb}
\vskip -9 pt
$$
\vbox{
\nointerlineskip
\halign{\strut
    \vrule \ \ \hfil \frb #\ 
   &\vrule \hfil \ \ \fbb #\frb\ 
   &\vrule \hfil \ \ \frb #\ \hfil
   &\vrule \hfil \ \ \frb #\ 
   &\vrule \hfil \ \ \frb #\ \ \vrule \hskip 2 pt
   &\vrule \ \ \hfil \frb #\ 
   &\vrule \hfil \ \ \fbb #\frb\ 
   &\vrule \hfil \ \ \frb #\ \hfil
   &\vrule \hfil \ \ \frb #\ 
   &\vrule \hfil \ \ \frb #\ \vrule \cr%
\noalign{\hrule}
 & &3.5.7.11.13.43.53.107&313&222& & &121.19.67.113.211&87&2234 \cr
11809&3661452795&4.9.11.37.43.53.313&485&98&11827&3672608819&4.3.11.29.67.1117&927&190 \cr
 & &16.5.49.37.97.313&11581&5432& & &16.27.5.19.29.103&13905&232 \cr
\noalign{\hrule}
 & &243.5.11.1681.163&4081&4324& & &3.25.29.53.151.211&23067&16948 \cr
11810&3662050095&8.7.121.23.47.53.163&7163&1476&11828&3672776775&8.27.5.11.19.223.233&179&44 \cr
 & &64.9.7.13.19.23.29.41&5681&6496& & &64.121.19.179.233&41707&73568 \cr
\noalign{\hrule}
 & &5.343.11.29.37.181&5969&3978& & &11.19.23.29.41.643&747&196 \cr
11811&3663828245&4.9.5.13.17.37.47.127&1267&638&11829&3675081289&8.9.49.11.83.643&2665&4408 \cr
 & &16.3.7.11.13.29.47.181&141&104& & &128.3.5.7.13.19.29.41&195&448 \cr
\noalign{\hrule}
 & &9.11.13.71.101.397&2261&2900& & &17.23.47.227.881&4797&15466 \cr
11812&3663943569&8.25.7.11.17.19.29.101&2193&5122&11830&3675161099&4.9.11.13.17.19.37.41&161&2270 \cr
 & &32.3.5.13.289.43.197&12427&15760& & &16.3.5.7.23.41.227&41&840 \cr
\noalign{\hrule}
 & &9.5.13.31.37.43.127&16247&4436& & &9.5.19.31.331.419&184253&173992 \cr
11813&3664303695&8.3.5.7.11.211.1109&17447&5842&11831&3675951945&16.7.13.23.239.8011&1257&6754 \cr
 & &32.23.73.127.239&5497&1168& & &64.3.7.11.13.307.419&3377&2912 \cr
\noalign{\hrule}
 & &9.5.19.47.53.1721&20119&20066& & &3.5.7.11.17.37.61.83&901&320 \cr
11814&3665394405&4.11.31.59.79.127.1721&52311&1040&11832&3678244185&128.25.289.53.61&5229&1996 \cr
 & &128.3.5.7.13.47.53.127&889&832& & &1024.9.7.83.499&499&1536 \cr
\noalign{\hrule}
 & &3.25.11.13.23.89.167&931&906& & &3.25.13.19.23.89.97&149&1306 \cr
11815&3666330525&4.9.49.13.19.23.89.151&835&3526&11833&3678305475&4.5.19.23.149.653&8577&8558 \cr
 & &16.5.19.41.43.151.167&6191&6536& & &16.9.11.389.653.953&254017&251592 \cr
\noalign{\hrule}
 & &81.5.49.181.1021&26767&17578& & &7.11.19.29.127.683&169&720 \cr
11816&3667375845&4.9.11.13.17.29.47.71&85&14&11834&3680160407&32.9.5.11.169.683&127&556 \cr
 & &16.5.7.13.289.29.47&3757&10904& & &256.3.5.13.127.139&5421&640 \cr
\noalign{\hrule}
 & &121.41.47.15731&10709&5022& & &3.7.13.109.337.367&3335&1046 \cr
11817&3667950077&4.81.31.41.10709&5539&5170&11835&3680316003&4.5.23.29.367.523&693&11336 \cr
 & &16.9.5.11.29.31.47.191&8091&7640& & &64.9.5.7.11.13.109&15&352 \cr
\noalign{\hrule}
 & &81.7.11.31.61.311&6367&12604& & &27.25.11.47.53.199&767&502 \cr
11818&3667985937&8.23.31.137.6367&5307&1060&11836&3680639325&4.5.11.13.59.199.251&329&2916 \cr
 & &64.3.5.23.29.53.61&1219&4640& & &32.729.7.47.251&1757&432 \cr
\noalign{\hrule}
 & &5.7.13.19.43.71.139&26373&66062& & &27.5.13.19.59.1871&39287&33682 \cr
11819&3668652715&4.3.17.29.59.67.149&781&930&11837&3680921205&4.9.11.17.1531.2311&23413&2614 \cr
 & &16.9.5.11.17.31.67.71&5797&4824& & &16.11.13.1307.1801&19811&10456 \cr
\noalign{\hrule}
 & &3.5.7.157.457.487&165&322& & &79.83.269.2087&12207&10120 \cr
11820&3668885115&4.9.25.49.11.23.457&2669&1444&11838&3681119471&16.3.5.11.13.23.79.313&813&2630 \cr
 & &32.11.17.361.23.157&3971&6256& & &64.9.25.13.263.271&71273&93600 \cr
\noalign{\hrule}
 & &81.25.11.13.19.23.29&28453&66928& & &9.7.11.13.17.29.829&39155&36284 \cr
11821&3669783975&32.37.47.89.769&429&340&11839&3681951273&8.5.17.41.47.191.193&34811&2052 \cr
 & &256.3.5.11.13.17.37.47&1739&2176& & &64.27.5.7.19.4973&4973&9120 \cr
\noalign{\hrule}
 & &9.7.71.131.6263&35685&29422& & &3.13.17.23.41.43.137&267&5350 \cr
11822&3669886269&4.81.5.13.47.61.313&131&3938&11840&3683106219&4.9.25.43.89.107&923&1012 \cr
 & &16.5.11.61.131.179&895&5368& & &32.5.11.13.23.71.107&3905&1712 \cr
\noalign{\hrule}
 & &9.11.289.19.43.157&17755&63128& & &5.7.11.23.47.53.167&28073&16182 \cr
11823&3669904359&16.5.13.53.67.607&471&136&11841&3683653435&4.9.7.29.31.67.419&495&2438 \cr
 & &256.3.13.17.53.157&689&128& & &16.81.5.11.23.31.53&81&248 \cr
\noalign{\hrule}
 & &3.19.31.47.193.229&8085&986& & &9.5.17.19.113.2243&637&1606 \cr
11824&3670516653&4.9.5.49.11.17.19.29&965&916&11842&3684026565&4.3.5.49.11.13.73.113&1517&1178 \cr
 & &32.25.17.29.193.229&493&400& & &16.13.19.31.37.41.73&16523&21608 \cr
\noalign{\hrule}
 & &9.5.11.37.43.59.79&8227&8758& & &9.17.19.29.109.401&3035&8594 \cr
11825&3670747245&4.11.19.29.37.151.433&6075&6482&11843&3684799827&4.3.5.19.607.4297&15151&19448 \cr
 & &16.243.25.7.19.151.463&100415&100008& & &64.5.11.13.17.109.139&1529&2080 \cr
\noalign{\hrule}
 & &3.5.11.19.61.73.263&3225&15974& & &9.11.61.293.2083&91885&68972 \cr
11826&3671520765&4.9.125.49.43.163&8833&11542&11844&3685716441&8.5.17.23.43.47.401&5011&4212 \cr
 & &16.7.121.29.73.199&1393&2552& & &64.81.5.13.43.5011&65143&61920 \cr
\noalign{\hrule}
}%
}
$$
\eject
\vglue -23 pt
\noindent\hskip 1 in\hbox to 6.5 in{\ 11845 -- 11880 \hfill\fbd 3686304705 -- 3721387635\frb}
\vskip -9 pt
$$
\vbox{
\nointerlineskip
\halign{\strut
    \vrule \ \ \hfil \frb #\ 
   &\vrule \hfil \ \ \fbb #\frb\ 
   &\vrule \hfil \ \ \frb #\ \hfil
   &\vrule \hfil \ \ \frb #\ 
   &\vrule \hfil \ \ \frb #\ \ \vrule \hskip 2 pt
   &\vrule \ \ \hfil \frb #\ 
   &\vrule \hfil \ \ \fbb #\frb\ 
   &\vrule \hfil \ \ \frb #\ \hfil
   &\vrule \hfil \ \ \frb #\ 
   &\vrule \hfil \ \ \frb #\ \vrule \cr%
\noalign{\hrule}
 & &3.5.961.47.5441&8207&13648& & &81.7.17.281.1367&6049&3520 \cr
11845&3686304705&32.29.31.283.853&17325&9118&11863&3702600153&128.9.5.11.17.23.263&5473&8992 \cr
 & &128.9.25.7.11.47.97&3201&2240& & &8192.13.281.421&5473&4096 \cr
\noalign{\hrule}
 & &7.13.167.431.563&469821&465880& & &3.25.7.11.17.97.389&2921&1754 \cr
11846&3687597641&16.3.5.11.19.23.613.619&1169&12930&11864&3704437275&4.7.23.97.127.877&99&778 \cr
 & &64.9.25.7.11.167.431&275&288& & &16.9.11.23.127.389&127&552 \cr
\noalign{\hrule}
 & &3.5.11.13.61.71.397&259&656& & &49.17.31.43.47.71&915&418 \cr
11847&3688128015&32.7.11.13.37.41.71&7543&12834&11865&3705368093&4.3.5.7.11.17.19.47.61&71&258 \cr
 & &128.9.19.23.31.397&2139&1216& & &16.9.5.19.43.61.71&1159&360 \cr
\noalign{\hrule}
 & &27.5.11.19.239.547&341&854& & &3.5.13.19.529.31.61&201&202 \cr
11848&3688631595&4.7.121.31.61.547&5605&1776&11866&3706255995&4.9.5.19.529.61.67.101&517&5278 \cr
 & &128.3.5.19.31.37.59&1147&3776& & &16.7.11.13.29.47.67.101&74437&76328 \cr
\noalign{\hrule}
 & &27.49.17.31.67.79&6149&3700& & &5.19.23.31.229.239&2627&1914 \cr
11849&3690390753&8.9.25.11.13.17.37.43&7049&7270&11867&3707204285&4.3.5.11.29.37.71.229&741&1886 \cr
 & &32.125.7.11.19.53.727&125875&127952& & &16.9.11.13.19.23.29.41&4147&2952 \cr
\noalign{\hrule}
 & &343.13.41.61.331&12325&7866& & &25.7.11.13.19.29.269&2561&2550 \cr
11850&3691298429&4.9.25.17.19.23.29.41&10819&9394&11868&3709180475&4.3.625.7.169.17.29.197&993&19118 \cr
 & &16.3.7.11.23.31.61.349&7843&8376& & &16.9.121.79.197.331&65207&62568 \cr
\noalign{\hrule}
 & &3.5.11.19.61.97.199&4067&16206& & &3.25.7.11.131.4903&1469&3434 \cr
11851&3691409205&4.9.5.49.37.73.83&121&194&11869&3709242075&4.5.7.11.13.17.101.113&4903&3798 \cr
 & &16.7.121.37.83.97&913&2072& & &16.9.101.211.4903&633&808 \cr
\noalign{\hrule}
 & &3.25.121.13.131.239&41&3066& & &9.25.7.11.19.59.191&23467&4558 \cr
11852&3693679275&4.9.7.41.73.131&33475&33466&11870&3709473075&4.7.31.43.53.757&375&382 \cr
 & &16.25.13.29.41.103.577&23657&23896& & &16.3.125.31.43.53.191&1643&1720 \cr
\noalign{\hrule}
 & &5.13.29.31.191.331&4653&5608& & &7.13.17.23.127.821&3663&15220 \cr
11853&3694319135&16.9.11.13.29.47.701&83&20246&11871&3709924127&8.9.5.11.17.37.761&253&508 \cr
 & &64.3.53.83.191&4399&96& & &64.3.121.23.37.127&363&1184 \cr
\noalign{\hrule}
 & &125.7.11.17.19.29.41&26559&13316& & &81.5.7.121.29.373&26619&27466 \cr
11854&3696452375&8.9.7.13.227.3329&2459&870&11872&3710609595&4.243.19.31.443.467&7975&442 \cr
 & &32.27.5.13.29.2459&2459&5616& & &16.25.11.13.17.29.467&2335&1768 \cr
\noalign{\hrule}
 & &27.121.17.19.31.113&1345&5096& & &5.11.19.41.79.1097&5513&10998 \cr
11855&3696509223&16.9.5.49.13.17.269&3293&872&11873&3713076235&4.9.13.37.41.47.149&47795&43244 \cr
 & &256.13.37.89.109&9701&61568& & &32.3.5.121.19.79.569&569&528 \cr
\noalign{\hrule}
 & &81.5.11.13.29.31.71&5123&628& & &9.25.13.41.173.179&241&2086 \cr
11856&3696656535&8.11.13.47.109.157&71&72&11874&3713717475&4.5.7.149.173.241&459&286 \cr
 & &128.9.47.71.109.157&7379&6976& & &16.27.7.11.13.17.241&1309&5784 \cr
\noalign{\hrule}
 & &9.25.13.17.23.53.61&1551&2954& & &3.11.13.19.31.61.241&4403&1270 \cr
11857&3697501275&4.27.5.7.11.13.47.211&8441&8056&11875&3714663381&4.5.7.11.17.19.37.127&739&1674 \cr
 & &64.19.23.53.211.367&6973&6752& & &16.27.7.31.37.739&6651&2072 \cr
\noalign{\hrule}
 & &25.49.121.13.19.101&333&382& & &27.25.7.13.197.307&211189&212164 \cr
11858&3697769075&4.9.5.11.19.37.101.191&5083&6128&11876&3714922575&8.9.11.29.31.59.73.263&197&460 \cr
 & &128.3.13.17.23.191.383&73153&75072& & &64.5.11.23.29.31.59.197&20677&20768 \cr
\noalign{\hrule}
 & &27.25.11.23.59.367&287&362& & &3.25.49.11.19.47.103&901&274 \cr
11859&3697791075&4.9.7.23.41.181.367&3733&430&11877&3718251075&4.49.17.53.103.137&39269&46530 \cr
 & &16.5.7.41.43.3733&12341&29864& & &16.9.5.11.47.107.367&1101&856 \cr
\noalign{\hrule}
 & &5.7.13.23.37.41.233&125&162& & &9.11.17.31.37.41.47&525&932 \cr
11860&3698969365&4.81.625.13.23.233&92389&94486&11878&3719882727&8.27.25.7.17.41.233&1403&4888 \cr
 & &16.9.7.11.17.37.227.397&34731&34936& & &128.5.7.13.23.47.61&9821&4160 \cr
\noalign{\hrule}
 & &5.11.13.23.157.1433&1653&388& & &27.5.7.19.43.61.79&71&2694 \cr
11861&3699812545&8.3.19.29.97.1433&2123&690&11879&3720581235&4.81.19.71.449&715&634 \cr
 & &32.9.5.11.19.23.193&1737&304& & &16.5.11.13.317.449&5837&27896 \cr
\noalign{\hrule}
 & &25.23.41.53.2963&1029&1144& & &9.5.7.17.439.1583&7141&7106 \cr
11862&3702194425&16.3.5.343.11.13.2963&2829&134&11880&3721387635&4.11.289.19.37.193.439&9517&1176 \cr
 & &64.9.7.13.23.41.67&819&2144& & &64.3.49.11.31.193.307&66619&67936 \cr
\noalign{\hrule}
}%
}
$$
\eject
\vglue -23 pt
\noindent\hskip 1 in\hbox to 6.5 in{\ 11881 -- 11916 \hfill\fbd 3721660267 -- 3747658959\frb}
\vskip -9 pt
$$
\vbox{
\nointerlineskip
\halign{\strut
    \vrule \ \ \hfil \frb #\ 
   &\vrule \hfil \ \ \fbb #\frb\ 
   &\vrule \hfil \ \ \frb #\ \hfil
   &\vrule \hfil \ \ \frb #\ 
   &\vrule \hfil \ \ \frb #\ \ \vrule \hskip 2 pt
   &\vrule \ \ \hfil \frb #\ 
   &\vrule \hfil \ \ \fbb #\frb\ 
   &\vrule \hfil \ \ \frb #\ \hfil
   &\vrule \hfil \ \ \frb #\ 
   &\vrule \hfil \ \ \frb #\ \vrule \cr%
\noalign{\hrule}
 & &13.31.47.349.563&97893&115346& & &27.25.7.19.43.967&1859&1042 \cr
11881&3721660267&4.9.49.11.73.107.149&16403&460&11899&3732934275&4.9.25.7.11.169.521&1829&646 \cr
 & &32.3.5.7.23.47.349&483&80& & &16.17.19.31.59.521&16151&8024 \cr
\noalign{\hrule}
 & &3.5.13.17.41.139.197&197&418& & &27.13.89.127.941&64163&19586 \cr
11882&3721760445&4.11.19.139.38809&4879&33930&11900&3733279173&4.7.11.19.307.1399&1471&25110 \cr
 & &16.9.5.7.13.17.29.41&29&168& & &16.81.5.31.1471&465&11768 \cr
\noalign{\hrule}
 & &3.7.11.37.269.1619&10643&690& & &9.25.19.751.1163&803&52 \cr
11883&3722312517&4.9.5.11.23.29.367&7361&6994&11901&3733840575&8.5.11.13.73.1163&42807&42092 \cr
 & &16.13.17.23.269.433&5083&3464& & &64.3.17.19.619.751&619&544 \cr
\noalign{\hrule}
 & &3.25.11.19.31.47.163&219&56& & &5.11.37.53.89.389&97&486 \cr
11884&3722671425&16.9.7.19.31.47.73&935&4366&11902&3734047955&4.243.5.37.89.97&389&56 \cr
 & &64.5.7.11.17.37.59&1003&8288& & &64.27.7.97.389&189&3104 \cr
\noalign{\hrule}
 & &27.11.17.29.47.541&10369&5320& & &11.263.373.3461&68085&30014 \cr
11885&3723046767&16.5.7.19.47.10369&2059&8310&11903&3734727029&4.9.5.17.43.89.349&13703&17358 \cr
 & &64.3.25.29.71.277&6925&2272& & &16.27.11.71.193.263&1917&1544 \cr
\noalign{\hrule}
 & &3.841.647.2281&1391&550& & &9.11.13.71.103.397&1265&74 \cr
11886&3723461061&4.25.11.13.107.2281&36163&20862&11904&3736496907&4.3.5.121.23.37.71&6749&6386 \cr
 & &16.9.19.841.43.61&817&1464& & &16.17.23.31.103.397&391&248 \cr
\noalign{\hrule}
 & &3.1369.47.101.191&5907&1160& & &9.13.29.97.11353&170237&158884 \cr
11887&3723722439&16.9.5.11.29.37.179&4427&764&11905&3736510713&8.11.23.37.43.107.157&3785&174 \cr
 & &128.5.19.191.233&4427&320& & &32.3.5.11.29.43.757&3785&7568 \cr
\noalign{\hrule}
 & &3.625.11.13.17.19.43&8161&10286& & &27.5.343.17.47.101&41&194 \cr
11888&3723988125&4.5.19.37.139.8161&90387&64672&11906&3736767195&4.3.343.41.97.101&2585&1556 \cr
 & &256.9.121.43.47.83&3901&4224& & &32.5.11.47.97.389&4279&1552 \cr
\noalign{\hrule}
 & &5.11.41.83.101.197&69&14& & &9.49.11.19.23.41.43&821&950 \cr
11889&3724022005&4.3.7.23.41.101.197&2877&5200&11907&3737360781&4.3.25.7.361.41.821&45623&28382 \cr
 & &128.9.25.49.13.137&28665&8768& & &16.5.23.43.617.1061&5305&4936 \cr
\noalign{\hrule}
 & &9.11.13.4913.19.31&15995&15982& & &17.19.29.31.79.163&345&182 \cr
11890&3724265259&4.5.7.289.31.61.131.457&4251&4708&11908&3739184629&4.3.5.7.13.19.23.29.79&6921&3586 \cr
 & &32.3.5.7.11.13.61.107.109.131&228445&228464& & &16.27.11.13.163.769&3861&6152 \cr
\noalign{\hrule}
 & &9.7.19.199.15643&9193&24836& & &121.281.317.347&29575&67932 \cr
11891&3726209529&8.49.29.317.887&21573&21890&11909&3740075999&8.27.25.7.169.17.37&79&40 \cr
 & &32.27.5.11.17.29.47.199&7755&7888& & &128.9.125.13.37.79&37999&72000 \cr
\noalign{\hrule}
 & &81.25.7.19.109.127&19619&18856& & &9.7.11.13.17.53.461&61535&6676 \cr
11892&3728265975&16.23.127.853.2357&81271&27060&11910&3741987249&8.5.31.397.1669&1033&636 \cr
 & &128.3.5.11.41.67.1213&49733&47168& & &64.3.5.31.53.1033&1033&4960 \cr
\noalign{\hrule}
 & &5.11.13.17.29.71.149&943&162& & &7.23.43.59.89.103&957&1580 \cr
11893&3729044605&4.81.23.29.41.149&209&238&11911&3744325319&8.3.5.11.23.29.79.103&651&1166 \cr
 & &16.27.7.11.17.19.23.41&4347&6232& & &32.9.7.121.29.31.53&33759&24592 \cr
\noalign{\hrule}
 & &81.31.37.137.293&4165&6676& & &3.7.103.1217.1423&71459&75110 \cr
11894&3729379887&8.5.49.17.137.1669&30969&2596&11912&3745863933&4.5.49.19.29.37.3761&9441&62018 \cr
 & &64.27.11.31.37.59&59&352& & &16.9.5.11.1049.2819&57695&67656 \cr
\noalign{\hrule}
 & &243.13.109.10831&1525&1634& & &9.5.13.19.31.83.131&121&10994 \cr
11895&3729449061&4.25.19.43.61.10831&5995&4836&11913&3746455245&4.121.23.31.239&351&362 \cr
 & &32.3.125.11.13.31.43.109&5375&5456& & &16.27.11.13.181.239&5973&1912 \cr
\noalign{\hrule}
 & &3.25.11.19.37.59.109&4601&1004& & &27.13.23.1849.251&421&2270 \cr
11896&3729819225&8.5.37.43.107.251&6859&16146&11914&3746671227&4.3.5.227.251.421&2093&1672 \cr
 & &32.27.13.6859.23&8303&1872& & &64.7.11.13.19.23.227&4313&2464 \cr
\noalign{\hrule}
 & &7.11.13.59.83.761&381&380& & &49.19.41.127.773&20145&17732 \cr
11897&3730343617&8.3.5.7.11.13.19.59.83.127&4969&426&11915&3747285241&8.3.5.11.13.17.31.41.79&423&28 \cr
 & &32.9.19.71.127.4969&352799&347472& & &64.27.7.13.17.31.47&10881&25568 \cr
\noalign{\hrule}
 & &25.7.11.13.29.37.139&14177&11538& & &9.11.17.71.79.397&7285&536 \cr
11898&3732403675&4.9.5.11.641.14177&403&238&11916&3747658959&16.5.31.47.67.71&1501&1836 \cr
 & &16.3.7.13.17.31.14177&14177&12648& & &128.27.17.19.31.79&589&192 \cr
\noalign{\hrule}
}%
}
$$
\vfill\eject
\vglue -23 pt
\noindent\hskip 1 in\hbox to 6.5 in{\ 11917 -- 11952 \hfill\fbd 3747979125 -- 3783512029\frb}
\vskip -9 pt
$$
\vbox{
\nointerlineskip
\halign{\strut
    \vrule \ \ \hfil \frb #\ 
   &\vrule \hfil \ \ \fbb #\frb\ 
   &\vrule \hfil \ \ \frb #\ \hfil
   &\vrule \hfil \ \ \frb #\ 
   &\vrule \hfil \ \ \frb #\ \ \vrule \hskip 2 pt
   &\vrule \ \ \hfil \frb #\ 
   &\vrule \hfil \ \ \fbb #\frb\ 
   &\vrule \hfil \ \ \frb #\ \hfil
   &\vrule \hfil \ \ \frb #\ 
   &\vrule \hfil \ \ \frb #\ \vrule \cr%
\noalign{\hrule}
 & &9.125.11.41.83.89&353&394& & &27.5.31.37.109.223&5161&4976 \cr
11917&3747979125&4.125.11.89.197.353&3621&7504&11935&3763817415&32.9.13.223.311.397&2849&50 \cr
 & &128.3.7.17.67.71.197&80869&88256& & &128.25.7.11.37.397&1985&4928 \cr
\noalign{\hrule}
 & &9.5.11.19.43.73.127&527&1924& & &11.13.19.23.29.31.67&3285&4558 \cr
11918&3749331465&8.3.5.13.17.31.37.73&4373&1672&11936&3764020403&4.9.5.13.29.43.53.73&253&124 \cr
 & &128.11.17.19.4373&4373&1088& & &32.3.5.11.23.31.53.73&1095&848 \cr
\noalign{\hrule}
 & &3.7.29.59.241.433&2875&4114& & &243.25.7.11.83.97&283&962 \cr
11919&3749507643&4.125.121.17.23.433&2169&2594&11937&3766056525&4.81.5.11.13.37.283&4067&388 \cr
 & &16.9.5.11.23.241.1297&6485&6072& & &32.49.37.83.97&259&16 \cr
\noalign{\hrule}
 & &9.125.121.169.163&1357&268& & &9.25.23.67.73.149&7579&2404 \cr
11920&3749835375&8.13.23.59.67.163&46519&44400&11938&3771327825&8.11.13.53.73.601&2831&3780 \cr
 & &256.3.25.11.37.4229&4229&4736& & &64.27.5.7.19.53.149&1007&672 \cr
\noalign{\hrule}
 & &5.49.121.13.37.263&1219&1674& & &13.961.37.41.199&28159&67560 \cr
11921&3750181435&4.27.7.11.23.31.37.53&1391&6638&11939&3771424319&16.3.5.29.563.971&8649&7678 \cr
 & &16.3.13.23.107.3319&9957&19688& & &64.27.5.11.961.349&3839&4320 \cr
\noalign{\hrule}
 & &3.25.49.11.31.41.73&2109&884& & &13.6859.101.419&349103&343656 \cr
11922&3750752775&8.9.11.13.17.19.31.37&42413&47704&11940&3773458273&16.27.37.43.59.61.97&455&418 \cr
 & &128.7.67.73.83.89&5963&5312& & &64.3.5.7.11.13.19.43.59.61&75579&75680 \cr
\noalign{\hrule}
 & &9.7.17.361.31.313&4477&1660& & &27.25.7.11.13.37.151&277&4052 \cr
11923&3751480593&8.5.7.121.31.37.83&19&12&11941&3774996225&8.3.7.11.277.1013&31487&32500 \cr
 & &64.3.5.121.19.37.83&10043&5920& & &64.625.13.23.1369&925&736 \cr
\noalign{\hrule}
 & &27.31.101.199.223&1445&4576& & &9.13.41.53.83.179&469&220 \cr
11924&3751498449&64.5.11.13.289.199&591&404&11942&3777258537&8.3.5.7.11.41.67.179&371&166 \cr
 & &512.3.13.17.101.197&2561&4352& & &32.49.11.53.67.83&539&1072 \cr
\noalign{\hrule}
 & &5.11.17.67.101.593&1657&8424& & &9.11.17.19.31.37.103&755&2438 \cr
11925&3751996985&16.81.5.11.13.1657&1571&86&11943&3777794757&4.5.19.23.37.53.151&961&258 \cr
 & &64.3.13.43.1571&61269&1376& & &16.3.5.961.43.151&6493&1240 \cr
\noalign{\hrule}
 & &37.43.1129.2089&14373&62920& & &81.11.13.41.73.109&25565&23594 \cr
11926&3752343271&16.9.5.121.13.1597&1609&3182&11944&3778803171&4.3.5.13.47.251.5113&71&682 \cr
 & &64.3.5.37.43.1609&1609&480& & &16.5.11.31.71.5113&25565&17608 \cr
\noalign{\hrule}
 & &9.5.11.31.43.5689&8843&3154& & &7.11.29.443.3821&57285&84032 \cr
11927&3753801315&4.5.11.19.37.83.239&141&1054&11945&3779805799&128.9.5.13.19.67.101&1193&322 \cr
 & &16.3.17.19.31.37.47&703&6392& & &512.3.7.19.23.1193&27439&14592 \cr
\noalign{\hrule}
 & &5.49.11.23.71.853&1313&2166& & &27.25.11.13.43.911&5457&6368 \cr
11928&3753997555&4.3.5.11.13.361.23.101&2559&2996&11946&3781173825&64.81.13.17.107.199&605&1982 \cr
 & &32.9.7.13.19.107.853&2033&1872& & &256.5.121.107.991&10901&13696 \cr
\noalign{\hrule}
 & &81.25.49.13.41.71&5293&9262& & &9.5.7.23.31.113.149&2467&2152 \cr
11929&3754971675&4.5.11.13.67.79.421&539&1566&11947&3781506015&16.23.113.269.2467&2533&66 \cr
 & &16.27.49.121.29.67&1943&968& & &64.3.11.17.149.269&2959&544 \cr
\noalign{\hrule}
 & &27.7.11.31.101.577&7955&11086& & &5.11.13.29.79.2309&541&486 \cr
11930&3755894373&4.9.5.7.23.37.43.241&6347&4016&11948&3782292085&4.243.29.541.2309&2329&20 \cr
 & &128.5.11.23.251.577&1255&1472& & &32.3.5.17.137.541&27591&2192 \cr
\noalign{\hrule}
 & &27.125.11.13.43.181&49&424& & &3.5.343.17.61.709&18709&24540 \cr
11931&3756270375&16.9.49.13.53.181&769&860&11949&3782773785&8.9.25.53.353.409&46189&45836 \cr
 & &128.5.7.43.53.769&5383&3392& & &64.7.11.13.17.19.53.1637&86761&86944 \cr
\noalign{\hrule}
 & &3.5.121.13.19.83.101&147&268& & &7.17.103.257.1201&160457&148200 \cr
11932&3758140815&8.9.49.13.19.67.101&26477&28396&11950&3783208849&16.3.25.11.13.19.29.503&15141&20674 \cr
 & &64.7.11.29.31.83.229&6641&6944& & &64.9.5.49.103.10337&10337&10080 \cr
\noalign{\hrule}
 & &11.43.101.223.353&945&3398& & &3.5.13.17.53.61.353&140415&139514 \cr
11933&3760642787&4.27.5.7.353.1699&1313&3784&11951&3783240435&4.9.25.11.23.37.79.883&398551&109174 \cr
 & &64.9.5.11.13.43.101&65&288& & &16.169.17.19.113.3527&27911&28216 \cr
\noalign{\hrule}
 & &3.361.47.107.691&25&82& & &11.23.139.271.397&205&66 \cr
11934&3763467237&4.25.19.41.47.691&101&792&11952&3783512029&4.3.5.121.23.41.397&2085&7046 \cr
 & &64.9.25.11.41.101&3075&35552& & &16.9.25.13.139.271&325&72 \cr
\noalign{\hrule}
}%
}
$$
\eject
\vglue -23 pt
\noindent\hskip 1 in\hbox to 6.5 in{\ 11953 -- 11988 \hfill\fbd 3784344993 -- 3810898091\frb}
\vskip -9 pt
$$
\vbox{
\nointerlineskip
\halign{\strut
    \vrule \ \ \hfil \frb #\ 
   &\vrule \hfil \ \ \fbb #\frb\ 
   &\vrule \hfil \ \ \frb #\ \hfil
   &\vrule \hfil \ \ \frb #\ 
   &\vrule \hfil \ \ \frb #\ \ \vrule \hskip 2 pt
   &\vrule \ \ \hfil \frb #\ 
   &\vrule \hfil \ \ \fbb #\frb\ 
   &\vrule \hfil \ \ \frb #\ \hfil
   &\vrule \hfil \ \ \frb #\ 
   &\vrule \hfil \ \ \frb #\ \vrule \cr%
\noalign{\hrule}
 & &9.11.13.17.269.643&39095&36136& & &27.25.11.17.67.449&161&26 \cr
11953&3784344993&16.5.7.17.1117.4517&20787&1798&11971&3797226675&4.5.7.13.23.67.449&1089&1156 \cr
 & &64.3.7.169.29.31.41&8323&12896& & &32.9.7.121.13.289.23&2737&2288 \cr
\noalign{\hrule}
 & &5.13.17.29.31.37.103&2355&1456& & &3.11.169.37.41.449&15&466 \cr
11954&3785828345&32.3.25.7.169.17.157&3399&526&11972&3798678741&4.9.5.13.233.449&3535&506 \cr
 & &128.9.7.11.103.263&1841&6336& & &16.25.7.11.23.101&175&18584 \cr
\noalign{\hrule}
 & &3.5.49.13.307.1291&21919&23210& & &27.7.11.841.41.53&271&590 \cr
11955&3787000035&4.25.11.13.23.211.953&639&4214&11973&3799357947&4.9.5.29.53.59.271&451&1988 \cr
 & &16.9.49.43.71.953&9159&7624& & &32.5.7.11.41.59.71&295&1136 \cr
\noalign{\hrule}
 & &3.5.7.121.19.29.541&689&326& & &9.37.41.53.59.89&86875&106522 \cr
11956&3787246155&4.13.19.53.163.541&233&774&11974&3799670859&4.625.13.17.139.241&2937&30562 \cr
 & &16.9.13.43.163.233&37979&13416& & &16.3.5.7.11.37.59.89&35&88 \cr
\noalign{\hrule}
 & &3.5.19.107.283.439&13223&13662& & &27.25.7.121.289.23&511&1546 \cr
11957&3788607315&4.81.7.11.23.107.1889&985&878&11975&3800256075&4.3.5.49.17.73.773&253&988 \cr
 & &16.5.7.11.197.439.1889&15169&15112& & &32.11.13.19.23.773&773&3952 \cr
\noalign{\hrule}
 & &5.11.19.29.31.37.109&7911&8984& & &25.13.79.179.827&961&66 \cr
11958&3788822015&16.27.11.19.293.1123&3835&25172&11976&3800747275&4.3.5.11.961.827&1989&2816 \cr
 & &128.3.5.7.13.29.31.59&1239&832& & &2048.27.121.13.17&3267&17408 \cr
\noalign{\hrule}
 & &9.11.13.19.37.59.71&89723&65270& & &27.5.11.19.23.5861&46361&41554 \cr
11959&3790043829&4.5.23.47.61.83.107&2309&4218&11977&3803466645&4.9.7.37.79.179.263&823&5800 \cr
 & &16.3.5.19.37.47.2309&2309&1880& & &64.25.29.263.823&38135&26336 \cr
\noalign{\hrule}
 & &9.25.11.31.127.389&14651&14524& & &3.5.13.23.41.137.151&2023&242 \cr
11960&3790445175&8.3.49.11.13.23.31.3631&88303&4790&11978&3804028995&4.7.121.289.23.41&20687&14040 \cr
 & &32.5.7.227.389.479&3353&3632& & &64.27.5.13.137.151&9&32 \cr
\noalign{\hrule}
 & &3.11.23.47.59.1801&19933&64714& & &9.19.29.31.53.467&1525&1496 \cr
11961&3790577307&4.13.19.31.131.643&229&360&11979&3804946479&16.3.25.11.17.31.61.467&6401&604 \cr
 & &64.9.5.13.229.643&41795&21984& & &128.5.37.61.151.173&130615&144448 \cr
\noalign{\hrule}
 & &5.49.1741.8887&20437&64872& & &3.625.7.13.29.769&1749&3634 \cr
11962&3790705415&16.9.17.53.107.191&1463&1784&11980&3805108125&4.9.125.11.23.53.79&3799&6076 \cr
 & &256.3.7.11.19.53.223&35457&26752& & &32.49.29.31.53.131&4061&5936 \cr
\noalign{\hrule}
 & &9.121.23.167.907&2279&442& & &3.5.13.17.19.23.37.71&957&658 \cr
11963&3793844043&4.3.11.13.17.23.43.53&907&1160&11981&3805616685&4.9.7.11.29.37.47.71&301&338 \cr
 & &64.5.17.29.43.907&2465&1376& & &16.49.11.169.29.43.47&58609&56056 \cr
\noalign{\hrule}
 & &27.5.17.529.53.59&7163057&7162792& & &27.121.13.19.53.89&817&340 \cr
11964&3796349985&16.7.11.13.19.9839.34273&46817&81090&11982&3806378433&8.3.5.121.17.361.43&27323&16358 \cr
 & &64.9.5.7.11.17.19.53.46817&46817&46816& & &32.89.307.8179&8179&4912 \cr
\noalign{\hrule}
 & &9.7.97.101.6151&22781&32578& & &3.7.11.19.37.131.179&4043&2580 \cr
11965&3796464861&4.49.11.13.19.109.179&4371&970&11983&3807953457&8.9.5.13.43.131.311&1367&188 \cr
 & &16.3.5.11.13.31.47.97&2015&4136& & &64.13.43.47.1367&64249&17888 \cr
\noalign{\hrule}
 & &5.11.17.23.53.3331&92053&91152& & &9.5.11.29.37.71.101&97&310 \cr
11966&3796557215&32.27.13.23.73.97.211&21147&680&11984&3808769085&4.3.25.29.31.97.101&2627&302 \cr
 & &512.81.5.7.17.19.53&1539&1792& & &16.37.71.97.151&97&1208 \cr
\noalign{\hrule}
 & &3.49.11.23.31.37.89&24289&33650& & &243.5.11.29.31.317&887&1786 \cr
11967&3796568853&4.25.7.107.227.673&711&38&11985&3808797795&4.5.19.47.317.887&2059&2376 \cr
 & &16.9.25.19.79.227&5675&36024& & &64.27.11.19.29.47.71&1349&1504 \cr
\noalign{\hrule}
 & &81.25.19.29.41.83&259&254& & &9.7.19.41.149.521&12719&12760 \cr
11968&3796982325&4.3.5.7.29.37.41.83.127&7579&1634&11986&3809798433&16.5.49.11.23.29.79.521&25365&164 \cr
 & &16.7.11.13.19.43.53.127&47117&49192& & &128.3.25.19.23.41.89&2047&1600 \cr
\noalign{\hrule}
 & &9.17.29.37.101.229&73485&70556& & &9.5.17.67.149.499&4543&3940 \cr
11969&3797064801&8.81.5.23.31.71.569&2519&3232&11987&3810860505&8.25.7.11.59.149.197&13&162 \cr
 & &512.5.11.101.229.569&2845&2816& & &32.81.11.13.59.197&23049&10384 \cr
\noalign{\hrule}
 & &361.37.67.4243&90589&66402& & &7.11.13.29.1849.71&1809&3658 \cr
11970&3797141317&4.9.7.17.31.157.577&469&2200&11988&3810898091&4.27.13.29.31.59.67&1775&20468 \cr
 & &64.3.25.49.11.31.67&4557&8800& & &32.3.25.7.17.43.71&425&48 \cr
\noalign{\hrule}
}%
}
$$
\eject
\vglue -23 pt
\noindent\hskip 1 in\hbox to 6.5 in{\ 11989 -- 12024 \hfill\fbd 3812002875 -- 3851238611\frb}
\vskip -9 pt
$$
\vbox{
\nointerlineskip
\halign{\strut
    \vrule \ \ \hfil \frb #\ 
   &\vrule \hfil \ \ \fbb #\frb\ 
   &\vrule \hfil \ \ \frb #\ \hfil
   &\vrule \hfil \ \ \frb #\ 
   &\vrule \hfil \ \ \frb #\ \ \vrule \hskip 2 pt
   &\vrule \ \ \hfil \frb #\ 
   &\vrule \hfil \ \ \fbb #\frb\ 
   &\vrule \hfil \ \ \frb #\ \hfil
   &\vrule \hfil \ \ \frb #\ 
   &\vrule \hfil \ \ \frb #\ \vrule \cr%
\noalign{\hrule}
 & &9.125.29.331.353&42109&74734& & &9.7.11.17.359.907&22475&32452 \cr
11989&3812002875&4.11.17.43.79.2477&21489&20620&12007&3836046753&8.25.49.19.29.31.61&907&318 \cr
 & &32.3.5.13.19.29.43.1031&13403&13072& & &32.3.29.53.61.907&1537&976 \cr
\noalign{\hrule}
 & &11.19.31.41.113.127&4355&852& & &3.7.11.31.383.1399&23025&20344 \cr
11990&3812185289&8.3.5.11.13.19.67.71&26781&25546&12008&3836985537&16.9.25.11.307.2543&68551&71314 \cr
 & &32.9.53.79.113.241&12773&11376& & &64.5.49.181.197.1399&6335&6304 \cr
\noalign{\hrule}
 & &19.47.61.239.293&7777&62250& & &9.125.11.13.17.23.61&5201&13924 \cr
11991&3814580771&4.3.125.7.11.83.101&63&38&12009&3837029625&8.7.23.3481.743&10285&6804 \cr
 & &16.27.5.49.11.19.83&14553&3320& & &64.243.5.49.121.17&539&864 \cr
\noalign{\hrule}
 & &9.25.7.23.31.43.79&2057&82& & &9.25.121.23.6131&2653&3478 \cr
11992&3814746075&4.3.7.121.17.41.43&5&124&12010&3839078925&4.3.7.11.23.37.47.379&14885&11266 \cr
 & &32.5.121.31.41&4961&16& & &16.5.13.37.43.131.229&73229&67784 \cr
\noalign{\hrule}
 & &3.5.49.11.29.53.307&6141&4604& & &9.125.7.11.23.41.47&3349&5276 \cr
11993&3814980015&8.9.7.11.23.89.1151&47515&24998&12011&3839306625&8.3.7.11.17.197.1319&1975&2162 \cr
 & &32.5.13.17.29.43.431&7327&8944& & &32.25.23.47.79.1319&1319&1264 \cr
\noalign{\hrule}
 & &5.29.107.331.743&115209&130724& & &3.5.11.17.41.173.193&107&148 \cr
11994&3815650495&8.27.11.17.251.2971&10939&13910&12012&3839901945&8.11.37.107.173.193&855&1048 \cr
 & &32.3.5.13.17.107.10939&10939&10608& & &128.9.5.19.37.107.131&42051&44992 \cr
\noalign{\hrule}
 & &9.5.7.17.29.67.367&2035&2402& & &25.19.89.97.937&10791&12634 \cr
11995&3818548755&4.25.7.11.37.67.1201&277&2202&12013&3842332475&4.9.11.89.109.6317&4625&14326 \cr
 & &16.3.277.367.1201&1201&2216& & &16.3.125.11.13.19.29.37&4785&3848 \cr
\noalign{\hrule}
 & &49.11.13.23.131.181&1763&2400& & &25.11.13.61.67.263&3231&12812 \cr
11996&3821288471&64.3.25.11.41.43.131&203&72&12014&3842699575&8.9.25.359.3203&551&526 \cr
 & &1024.27.7.29.41.43&32103&22016& & &32.3.19.29.263.3203&9609&8816 \cr
\noalign{\hrule}
 & &7.11.23.2209.977&23845&22074& & &9.49.11.19.173.241&289&338 \cr
11997&3822159803&4.3.5.13.19.47.251.283&64251&2254&12015&3842802117&4.3.169.289.173.241&52415&94108 \cr
 & &16.27.49.121.23.59&1593&616& & &32.5.7.11.953.3361&16805&15248 \cr
\noalign{\hrule}
 & &9.125.11.19.71.229&15967&292& & &49.11.13.17.41.787&25769&34830 \cr
11998&3822897375&8.3.5.7.73.2281&181&184&12016&3843612773&4.81.5.7.43.73.353&327&26 \cr
 & &128.7.23.181.2281&52463&81088& & &16.243.5.13.73.109&17739&4360 \cr
\noalign{\hrule}
 & &11.83.103.109.373&5427&5800& & &27.5.7.13.19.43.383&115&158 \cr
11999&3823343623&16.81.25.11.29.67.83&7087&1526&12017&3844111635&4.9.25.19.23.79.383&2821&6754 \cr
 & &64.3.25.7.19.109.373&525&608& & &16.7.11.13.31.79.307&9517&6952 \cr
\noalign{\hrule}
 & &9.5.49.17.31.37.89&53443&81748& & &25.13.29.43.53.179&111&154 \cr
12000&3826581255&8.13.107.191.4111&3297&814&12018&3844843925&4.3.5.7.11.13.29.37.179&1961&366 \cr
 & &32.3.7.11.37.107.157&1177&2512& & &16.9.7.1369.53.61&9583&4392 \cr
\noalign{\hrule}
 & &7.121.13.43.59.137&1137&370& & &9.7.13.23.43.47.101&69&22 \cr
12001&3827082259&4.3.5.7.11.37.43.379&2603&708&12019&3845027277&4.27.11.529.43.101&12737&10010 \cr
 & &32.9.19.37.59.137&703&144& & &16.5.7.121.13.47.271&1355&968 \cr
\noalign{\hrule}
 & &9.125.11.73.4241&4067&8308& & &7.11.169.23.71.181&969&12968 \cr
12002&3831213375&8.49.31.67.73.83&783&1294&12020&3846291449&16.3.17.19.23.1621&4525&2904 \cr
 & &32.27.7.29.83.647&18763&27888& & &256.9.25.121.181&225&1408 \cr
\noalign{\hrule}
 & &27.5.11.13.19.31.337&371&466& & &81.5.7.13.29.59.61&1327&442 \cr
12003&3831907365&4.7.11.13.53.233.337&4265&558&12021&3846593205&4.27.7.169.17.1327&1531&1342 \cr
 & &16.9.5.31.233.853&853&1864& & &16.11.61.1327.1531&14597&12248 \cr
\noalign{\hrule}
 & &9.125.19.1849.97&3289&8836& & &5.7.121.53.61.281&25899&11006 \cr
12004&3833670375&8.3.11.13.19.23.2209&463&430&12022&3847383155&4.3.7.89.97.5503&3091&2412 \cr
 & &32.5.13.23.43.47.463&10649&9776& & &32.27.11.67.89.281&1809&1424 \cr
\noalign{\hrule}
 & &5.13.29.31.211.311&6371&5316& & &9.25.11.13.17.31.227&259&166 \cr
12005&3834563135&8.3.23.277.311.443&377&66&12023&3849063075&4.3.7.11.13.37.83.227&527&1970 \cr
 & &32.9.11.13.23.29.277&2493&4048& & &16.5.7.17.31.83.197&1379&664 \cr
\noalign{\hrule}
 & &13.17.19.43.67.317&10509&10730& & &11.361.23.149.283&8905&62694 \cr
12006&3834850123&4.3.5.19.29.31.37.43.113&6471&55792&12024&3851238611&4.729.5.13.43.137&283&68 \cr
 & &128.27.5.11.317.719&7909&8640& & &32.27.17.137.283&3699&272 \cr
\noalign{\hrule}
}%
}
$$
\eject
\vglue -23 pt
\noindent\hskip 1 in\hbox to 6.5 in{\ 12025 -- 12060 \hfill\fbd 3853794035 -- 3880693817\frb}
\vskip -9 pt
$$
\vbox{
\nointerlineskip
\halign{\strut
    \vrule \ \ \hfil \frb #\ 
   &\vrule \hfil \ \ \fbb #\frb\ 
   &\vrule \hfil \ \ \frb #\ \hfil
   &\vrule \hfil \ \ \frb #\ 
   &\vrule \hfil \ \ \frb #\ \ \vrule \hskip 2 pt
   &\vrule \ \ \hfil \frb #\ 
   &\vrule \hfil \ \ \fbb #\frb\ 
   &\vrule \hfil \ \ \frb #\ \hfil
   &\vrule \hfil \ \ \frb #\ 
   &\vrule \hfil \ \ \frb #\ \vrule \cr%
\noalign{\hrule}
 & &5.7.169.19.53.647&1647&2882& & &27.49.11.19.71.197&62305&3796 \cr
12025&3853794035&4.27.11.13.53.61.131&895&7096&12043&3867503409&8.5.13.17.73.733&4801&4728 \cr
 & &64.3.5.11.179.887&5907&28384& & &128.3.5.17.197.4801&4801&5440 \cr
\noalign{\hrule}
 & &27.11.137.211.449&3257&3706& & &7.11.71.257.2753&14893&15390 \cr
12026&3854835171&4.9.17.109.137.3257&2245&1012&12044&3868017307&4.81.5.19.53.257.281&2431&1424 \cr
 & &32.5.11.17.23.109.449&1955&1744& & &128.27.11.13.17.89.281&98631&96832 \cr
\noalign{\hrule}
 & &121.13.529.41.113&62135&2358& & &5.7.13.19.47.89.107&2129&2148 \cr
12027&3855198061&4.9.5.289.43.131&1631&2024&12045&3869337745&8.3.5.89.107.179.2129&153729&74074 \cr
 & &64.3.7.11.17.23.233&3961&672& & &32.27.7.11.13.19.29.31.37&12617&12528 \cr
\noalign{\hrule}
 & &9.5.7.31.67.71.83&667&418& & &3.5.121.41.149.349&16169&1860 \cr
12028&3855524715&4.3.11.19.23.29.67.71&1519&2300&12046&3869654415&8.9.25.19.23.31.37&3839&3136 \cr
 & &32.25.49.529.29.31&2645&3248& & &1024.49.11.23.349&1127&512 \cr
\noalign{\hrule}
 & &3.5.11.59.601.659&5&654& & &11.43.61.113.1187&2853&4040 \cr
12029&3855634365&4.9.25.109.601&2717&2692&12047&3870081743&16.9.5.11.43.101.317&1679&1808 \cr
 & &32.11.13.19.109.673&12787&22672& & &512.3.5.23.73.101.113&25185&25856 \cr
\noalign{\hrule}
 & &5.11.361.29.37.181&1449&2522& & &9.11.17.43.149.359&4109&2470 \cr
12030&3856099115&4.9.5.7.13.23.97.181&6919&24476&12048&3871090179&4.5.7.13.19.359.587&473&114 \cr
 & &32.3.11.17.29.37.211&633&272& & &16.3.5.7.11.13.361.43&1805&728 \cr
\noalign{\hrule}
 & &27.7.11.29.109.587&11077&11050& & &27.11.17.29.53.499&1927&2510 \cr
12031&3857599053&4.25.121.13.17.19.53.587&38619&464&12049&3872396187&4.3.5.41.47.251.499&27163&3710 \cr
 & &128.9.5.7.29.53.613&3065&3392& & &16.25.7.23.53.1181&4025&9448 \cr
\noalign{\hrule}
 & &7.11.19.29.211.431&3513&496& & &9.7.13.19.23.79.137&2171&2332 \cr
12032&3858353807&32.3.11.29.31.1171&16809&17150&12050&3873584169&8.3.11.169.53.137.167&14945&83582 \cr
 & &128.9.25.343.13.431&2925&3136& & &32.5.49.529.61.79&805&976 \cr
\noalign{\hrule}
 & &7.13.47.257.3511&383&3894& & &11.17.529.53.739&169&222 \cr
12033&3859252579&4.3.11.59.257.383&1605&1222&12051&3874517141&4.3.11.169.23.37.739&1275&2014 \cr
 & &16.9.5.13.47.59.107&2655&856& & &16.9.25.13.17.19.37.53&4329&3800 \cr
\noalign{\hrule}
 & &25.11.19.29.73.349&10291&8904& & &41.97.863.1129&18711&65000 \cr
12034&3860402425&16.3.5.7.29.41.53.251&871&666&12052&3874898479&16.243.625.7.11.13&1649&1726 \cr
 & &64.27.7.13.37.67.251&225589&216864& & &64.9.5.13.17.97.863&585&544 \cr
\noalign{\hrule}
 & &25.7.13.19.31.43.67&5141&9264& & &9.7.13.23.233.883&14575&6628 \cr
12035&3860467975&32.3.5.13.53.97.193&1771&1674&12053&3875505543&8.25.11.23.53.1657&27417&14008 \cr
 & &128.81.7.11.23.31.193&15633&16192& & &128.3.13.17.19.37.103&11951&6592 \cr
\noalign{\hrule}
 & &3.61.79.97.2753&8605&5852& & &3.7.29.53.113.1063&16195&16082 \cr
12036&3860611737&8.5.7.11.19.97.1721&6557&5490&12054&3877080963&4.5.11.17.41.43.79.1063&1701&56764 \cr
 & &32.9.25.19.61.79.83&1425&1328& & &32.243.7.23.43.617&26531&29808 \cr
\noalign{\hrule}
 & &9.13.17.1033.1879&575&458& & &27.7.11.13.43.47.71&4505&5428 \cr
12037&3860662923&4.25.17.23.229.1879&29139&2804&12055&3878131257&8.9.5.17.23.47.53.59&373&3146 \cr
 & &32.3.5.11.701.883&38555&14128& & &32.5.121.13.53.373&2915&5968 \cr
\noalign{\hrule}
 & &9.5.11.17.361.31.41&38297&38824& & &7.13.23.53.73.479&288629&295272 \cr
12038&3861062865&16.5.7.19.23.211.5471&27709&354&12056&3878854343&16.27.11.19.1367.1381&5545&9646 \cr
 & &64.3.121.23.59.229&13511&8096& & &64.9.5.7.13.19.53.1109&5545&5472 \cr
\noalign{\hrule}
 & &3.5.7.17.19.293.389&4039&2574& & &9.7.11.107.113.463&65237&19400 \cr
12039&3865529955&4.27.49.11.13.19.577&575&62&12057&3879506169&16.25.89.97.733&36663&34438 \cr
 & &16.25.11.23.31.577&6347&28520& & &64.3.121.67.101.257&25957&23584 \cr
\noalign{\hrule}
 & &3.25.11.103.173.263&2911&4814& & &9.25.11.73.109.197&1403&422 \cr
12040&3866277525&4.29.41.71.83.263&173&90&12058&3879634275&4.11.23.61.197.211&13067&28500 \cr
 & &16.9.5.29.41.71.173&1189&1704& & &32.3.125.19.73.179&895&304 \cr
\noalign{\hrule}
 & &9.5.11.19.227.1811&35&1846& & &9.125.11.23.43.317&49&424 \cr
12041&3866367285&4.25.7.13.71.227&807&782&12059&3879723375&16.3.49.23.53.317&1075&1144 \cr
 & &16.3.13.17.23.71.269&19099&40664& & &256.25.7.11.13.43.53&689&896 \cr
\noalign{\hrule}
 & &3.7.11.19.23.29.1321&13865&13876& & &49.11.13.19.103.283&2647&2400 \cr
12042&3867178623&8.5.19.23.29.47.59.3469&41301&39116&12060&3880693817&64.3.25.11.283.2647&5061&8174 \cr
 & &64.9.7.11.13.127.353.3469&440563&440544& & &256.9.5.7.61.67.241&73505&77184 \cr
\noalign{\hrule}
}%
}
$$
\eject
\vglue -23 pt
\noindent\hskip 1 in\hbox to 6.5 in{\ 12061 -- 12096 \hfill\fbd 3880888605 -- 3911592797\frb}
\vskip -9 pt
$$
\vbox{
\nointerlineskip
\halign{\strut
    \vrule \ \ \hfil \frb #\ 
   &\vrule \hfil \ \ \fbb #\frb\ 
   &\vrule \hfil \ \ \frb #\ \hfil
   &\vrule \hfil \ \ \frb #\ 
   &\vrule \hfil \ \ \frb #\ \ \vrule \hskip 2 pt
   &\vrule \ \ \hfil \frb #\ 
   &\vrule \hfil \ \ \fbb #\frb\ 
   &\vrule \hfil \ \ \frb #\ \hfil
   &\vrule \hfil \ \ \frb #\ 
   &\vrule \hfil \ \ \frb #\ \vrule \cr%
\noalign{\hrule}
 & &243.5.11.17.19.29.31&51379&51658& & &25.163.587.1627&1221&406 \cr
12061&3880888605&4.27.5.23.191.269.1123&307&30628&12079&3891824675&4.3.5.7.11.29.37.587&1627&1308 \cr
 & &32.13.19.31.191.307&3991&3056& & &32.9.7.37.109.1627&2331&1744 \cr
\noalign{\hrule}
 & &9.5.97.491.1811&979&494& & &9.11.13.67.163.277&23681&21470 \cr
12062&3881362365&4.3.11.13.19.89.1811&4075&1358&12080&3893325579&4.3.5.7.13.17.19.113.199&937&1210 \cr
 & &16.25.7.89.97.163&3115&1304& & &16.25.121.17.199.937&84575&82456 \cr
\noalign{\hrule}
 & &11.17.23.31.37.787&160167&147550& & &5.11.19.53.167.421&8143&7722 \cr
12063&3882465389&4.3.25.7.13.29.227.263&341&114&12081&3893953195&4.27.121.13.17.53.479&3437&21950 \cr
 & &16.9.5.11.19.29.31.263&7627&6840& & &16.3.25.7.13.439.491&44681&52680 \cr
\noalign{\hrule}
 & &7.11.169.19.113.139&1525&282& & &5.7.11.23.37.73.163&507&344 \cr
12064&3883508629&4.3.25.7.13.19.47.61&27&638&12082&3898528865&16.3.5.7.11.169.43.73&667&282 \cr
 & &16.81.5.11.29.61&305&18792& & &64.9.13.23.29.43.47&12267&17888 \cr
\noalign{\hrule}
 & &9.125.11.17.59.313&11219&10906& & &7.17.23.37.137.281&325&66 \cr
12065&3884995125&4.3.7.11.13.17.19.41.863&6667&626&12083&3898552693&4.3.25.11.13.137.281&4403&2622 \cr
 & &16.19.41.59.113.313&779&904& & &16.9.7.11.17.19.23.37&209&72 \cr
\noalign{\hrule}
 & &27.5.13.19.37.47.67&43981&38764& & &729.17.499.631&3877&4606 \cr
12066&3885125985&8.9.7.11.61.103.881&1885&2812&12084&3902171517&4.49.47.631.3877&4147&270 \cr
 & &64.5.13.19.29.37.881&881&928& & &16.27.5.7.11.13.29.47&4277&12760 \cr
\noalign{\hrule}
 & &5.11.13.17.29.73.151&1715&2664& & &25.13.361.29.31.37&5389&4686 \cr
12067&3885552385&16.9.25.343.11.17.37&11023&9802&12085&3902581475&4.3.11.17.19.29.71.317&117&434 \cr
 & &64.3.7.169.29.73.151&91&96& & &16.27.7.11.13.17.31.71&5049&3976 \cr
\noalign{\hrule}
 & &149.58081.449&62491&4410& & &3.5.7.41.61.89.167&33263&70778 \cr
12068&3885676981&4.9.5.49.11.13.19.23&241&298&12086&3903098115&4.29.31.37.43.823&12727&36594 \cr
 & &16.3.5.13.23.149.241&345&104& & &16.9.11.13.19.89.107&4173&1672 \cr
\noalign{\hrule}
 & &5.11.23.719.4273&29367&68912& & &5.343.19.37.41.79&2593&10098 \cr
12069&3886443055&32.9.13.59.73.251&1777&2530&12087&3905084155&4.27.11.17.41.2593&1645&948 \cr
 & &128.3.5.11.13.23.1777&1777&2496& & &32.81.5.7.11.47.79&517&1296 \cr
\noalign{\hrule}
 & &3.11.19.29.37.53.109&911&2110& & &3.11.23.37.241.577&305&546 \cr
12070&3886598067&4.5.29.37.211.911&17113&9306&12088&3905137731&4.9.5.7.11.13.61.577&4097&3404 \cr
 & &16.9.5.11.47.109.157&705&1256& & &32.5.17.23.37.61.241&305&272 \cr
\noalign{\hrule}
 & &9.5.11.19.281.1471&2737&4618& & &27.125.7.37.41.109&60203&64672 \cr
12071&3887566155&4.7.17.23.281.2309&16395&22858&12089&3906464625&64.7.11.13.43.47.421&30271&8040 \cr
 & &16.3.5.7.11.1039.1093&7273&8744& & &1024.3.5.67.30271&30271&34304 \cr
\noalign{\hrule}
 & &3.5.11.169.17.59.139&735&794& & &243.25.11.53.1103&5917&21658 \cr
12072&3887643045&4.9.25.49.169.17.397&649&2224&12090&3906522675&4.9.49.13.17.61.97&95&1166 \cr
 & &128.7.11.59.139.397&397&448& & &16.5.7.11.19.53.61&1159&56 \cr
\noalign{\hrule}
 & &625.11.23.67.367&4589&3852& & &3.121.13.17.53.919&23325&25382 \cr
12073&3888135625&8.9.625.13.107.353&217&842&12091&3907421661&4.9.25.343.13.37.311&919&3124 \cr
 & &32.3.7.13.31.107.421&91357&66768& & &32.5.7.11.37.71.919&1295&1136 \cr
\noalign{\hrule}
 & &9.7.11.19.47.61.103&62015&57362& & &3.5.13.17.31.47.809&6799&6954 \cr
12074&3888228267&4.5.7.23.29.43.79.157&1695&596&12092&3907433595&4.9.169.19.47.61.523&6325&1618 \cr
 & &32.3.25.23.43.113.149&85675&77744& & &16.25.11.19.23.61.809&3355&3496 \cr
\noalign{\hrule}
 & &5.13.19.59.197.271&603&1958& & &81.13.83.197.227&275&472 \cr
12075&3890043755&4.9.11.19.59.67.89&6775&5084&12093&3908395881&16.9.25.11.13.59.227&1379&664 \cr
 & &32.3.25.11.31.41.271&2255&1488& & &256.5.7.59.83.197&413&640 \cr
\noalign{\hrule}
 & &5.11.29.53.191.241&27&292& & &3.5.7.43.79.97.113&5809&1854 \cr
12076&3891230585&8.27.73.191.241&5887&8056&12094&3909624285&4.27.37.43.103.157&1325&2486 \cr
 & &128.3.7.19.841.53&1653&448& & &16.25.11.53.113.157&2915&1256 \cr
\noalign{\hrule}
 & &11.17.31.43.67.233&783&550& & &25.11.19.29.131.197&1219&4494 \cr
12077&3891369581&4.27.25.121.17.29.67&3029&2666&12095&3910405675&4.3.7.11.19.23.53.107&655&522 \cr
 & &16.9.5.13.29.31.43.233&377&360& & &16.27.5.23.29.53.131&621&424 \cr
\noalign{\hrule}
 & &5.121.17.47.83.97&7119&11020& & &7.59.131.197.367&5353&6270 \cr
12078&3891813145&8.9.25.7.11.19.29.113&14711&16364&12096&3911592797&4.3.5.11.19.53.101.367&10281&9170 \cr
 & &64.3.7.47.313.4091&28637&30048& & &16.9.25.7.19.23.131.149&10925&10728 \cr
\noalign{\hrule}
}%
}
$$
\eject
\vglue -23 pt
\noindent\hskip 1 in\hbox to 6.5 in{\ 12097 -- 12132 \hfill\fbd 3911797197 -- 3945662721\frb}
\vskip -9 pt
$$
\vbox{
\nointerlineskip
\halign{\strut
    \vrule \ \ \hfil \frb #\ 
   &\vrule \hfil \ \ \fbb #\frb\ 
   &\vrule \hfil \ \ \frb #\ \hfil
   &\vrule \hfil \ \ \frb #\ 
   &\vrule \hfil \ \ \frb #\ \ \vrule \hskip 2 pt
   &\vrule \ \ \hfil \frb #\ 
   &\vrule \hfil \ \ \fbb #\frb\ 
   &\vrule \hfil \ \ \frb #\ \hfil
   &\vrule \hfil \ \ \frb #\ 
   &\vrule \hfil \ \ \frb #\ \vrule \cr%
\noalign{\hrule}
 & &9.7.11.19.23.12917&6355&6562& & &3.25.7.11.13.43.1217&11247&2728 \cr
12097&3911797197&4.5.7.11.17.19.31.41.193&6617&138&12115&3928749825&16.9.121.23.31.163&79&200 \cr
 & &16.3.13.17.23.41.509&8653&4264& & &256.25.23.79.163&12877&2944 \cr
\noalign{\hrule}
 & &7.13.47.83.103.107&323&426& & &9.5.19.41.71.1579&825&754 \cr
12098&3912355811&4.3.13.17.19.47.71.83&2247&8140&12116&3929980995&4.27.125.11.13.19.29.41&1579&19454 \cr
 & &32.9.5.7.11.19.37.107&1881&2960& & &16.29.71.137.1579&137&232 \cr
\noalign{\hrule}
 & &3.11.43.73.179.211&585&2554& & &9.121.61.67.883&28037&25826 \cr
12099&3912377403&4.27.5.13.211.1277&14881&13604&12117&3930006069&4.3.11.529.37.53.349&1&1220 \cr
 & &32.13.19.23.179.647&8411&6992& & &32.5.23.61.349&8027&80 \cr
\noalign{\hrule}
 & &81.5.11.17.19.2719&2129&590& & &7.11.13.29.43.47.67&1135&1746 \cr
12100&3912545835&4.25.11.17.59.2129&32129&21096&12118&3930729803&4.9.5.7.11.29.97.227&43&188 \cr
 & &64.9.361.89.293&5567&2848& & &32.3.43.47.97.227&681&1552 \cr
\noalign{\hrule}
 & &11.31.61.83.2267&9937&15000& & &29.37.47.137.569&605&468 \cr
12101&3913936961&16.3.625.19.31.523&47&78&12119&3931247743&8.9.5.121.13.47.569&1159&548 \cr
 & &64.9.5.13.19.47.523&104481&83680& & &64.3.5.121.19.61.137&7381&9120 \cr
\noalign{\hrule}
 & &27.5.17.29.103.571&943&808& & &27.11.23.37.47.331&113&140 \cr
12102&3914299215&16.23.29.41.101.571&309&880&12120&3931985079&8.5.7.37.47.113.331&9869&5688 \cr
 & &512.3.5.11.23.101.103&2323&2816& & &128.9.5.7.71.79.139&39263&44480 \cr
\noalign{\hrule}
 & &13.23.61.157.1367&2489&1122& & &3.5.49.31.107.1613&2465&852 \cr
12103&3914435941&4.3.11.13.17.19.61.131&4715&19782&12121&3932485935&8.9.25.49.17.29.71&5507&7282 \cr
 & &16.27.5.7.23.41.157&287&1080& & &32.11.17.331.5507&61897&88112 \cr
\noalign{\hrule}
 & &81.13.19.389.503&3095&10486& & &125.29.37.139.211&13&198 \cr
12104&3914709669&4.3.5.49.13.107.619&13079&17252&12122&3933752125&4.9.25.11.13.29.139&541&1266 \cr
 & &32.5.11.19.29.41.227&12485&19024& & &16.27.11.211.541&297&4328 \cr
\noalign{\hrule}
 & &3.25.49.11.113.857&9129&18556& & &9.361.23.61.863&1535&2398 \cr
12105&3914797425&8.9.5.17.179.4639&1937&2702&12123&3933853461&4.5.11.19.61.109.307&345&326 \cr
 & &32.7.13.149.179.193&28757&37232& & &16.3.25.23.109.163.307&33463&32600 \cr
\noalign{\hrule}
 & &9.11.197.199.1009&1415&1612& & &29.59.479.4801&83745&55484 \cr
12106&3916026873&8.3.5.11.13.31.199.283&1009&14&12124&3934750769&8.9.5.11.13.97.1861&2161&3422 \cr
 & &32.7.13.283.1009&1981&208& & &32.3.5.11.29.59.2161&2161&2640 \cr
\noalign{\hrule}
 & &27.25.11.17.61.509&1349&1196& & &5.7.121.47.53.373&7767&12002 \cr
12107&3919160025&8.3.5.11.13.19.23.61.71&163&508&12125&3934920605&4.9.17.47.353.863&895&1694 \cr
 & &64.13.19.71.127.163&117221&99104& & &16.3.5.7.121.179.353&1059&1432 \cr
\noalign{\hrule}
 & &27.5.361.137.587&9541&8954& & &125.7.97.199.233&98363&75762 \cr
12108&3919219965&4.7.121.361.29.37.47&3425&7044&12126&3935399125&4.27.19.23.31.61.167&187&20 \cr
 & &32.3.25.11.37.137.587&185&176& & &32.3.5.11.17.19.31.61&30039&10736 \cr
\noalign{\hrule}
 & &3.49.19.103.13633&6765&6868& & &5.19.31.37.71.509&9401&9432 \cr
12109&3921927807&8.9.5.49.11.17.19.41.101&13633&5254&12127&3937886135&16.9.5.7.17.19.71.79.131&11167&4422 \cr
 & &32.5.37.41.71.13633&2911&2960& & &64.27.11.13.67.79.859&748189&750816 \cr
\noalign{\hrule}
 & &11.13.23.593.2011&639&2650& & &3.13.31.41.107.743&97063&74030 \cr
12110&3922208147&4.9.25.53.71.593&18377&13052&12128&3940785069&4.5.11.29.673.3347&1391&18126 \cr
 & &32.3.13.17.23.47.251&4267&2256& & &16.9.11.13.19.53.107&1007&264 \cr
\noalign{\hrule}
 & &3.7.19.79.89.1399&40183&39560& & &27.17.29.43.71.97&7657&63250 \cr
12111&3924711231&16.5.11.13.23.43.79.281&19&1008&12129&3941932851&4.125.11.13.19.23.31&261&974 \cr
 & &512.9.5.7.11.19.281&4215&2816& & &16.9.25.11.29.487&275&3896 \cr
\noalign{\hrule}
 & &27.5.47.61.73.139&5291&9044& & &3.7.11.13.17.31.47.53&9047&6980 \cr
12112&3927345615&8.7.11.13.17.19.37.73&417&94&12130&3942209271&8.5.7.17.83.109.349&2585&3348 \cr
 & &32.3.11.13.37.47.139&481&176& & &64.27.25.11.31.47.83&747&800 \cr
\noalign{\hrule}
 & &3.11.251.359.1321&1199&122& & &3.5.7.11.289.53.223&403&180 \cr
12113&3928121637&4.121.61.109.251&3565&3816&12131&3945123105&8.27.25.7.13.289.31&19847&30728 \cr
 & &64.9.5.23.31.53.109&37789&52320& & &128.23.89.167.223&3841&5696 \cr
\noalign{\hrule}
 & &3.5.7.11.17.59.3391&1577&4968& & &9.49.11.13.19.37.89&95489&71920 \cr
12114&3928354815&16.81.19.23.59.83&3391&3332&12132&3945662721&32.5.17.29.31.41.137&1351&2622 \cr
 & &128.49.17.19.23.3391&437&448& & &128.3.5.7.17.19.23.193&4439&5440 \cr
\noalign{\hrule}
}%
}
$$
\eject
\vglue -23 pt
\noindent\hskip 1 in\hbox to 6.5 in{\ 12133 -- 12168 \hfill\fbd 3948466665 -- 3970810877\frb}
\vskip -9 pt
$$
\vbox{
\nointerlineskip
\halign{\strut
    \vrule \ \ \hfil \frb #\ 
   &\vrule \hfil \ \ \fbb #\frb\ 
   &\vrule \hfil \ \ \frb #\ \hfil
   &\vrule \hfil \ \ \frb #\ 
   &\vrule \hfil \ \ \frb #\ \ \vrule \hskip 2 pt
   &\vrule \ \ \hfil \frb #\ 
   &\vrule \hfil \ \ \fbb #\frb\ 
   &\vrule \hfil \ \ \frb #\ \hfil
   &\vrule \hfil \ \ \frb #\ 
   &\vrule \hfil \ \ \frb #\ \vrule \cr%
\noalign{\hrule}
 & &3.5.11.13.17.19.41.139&15839&13476& & &27.625.11.37.577&169&206 \cr
12133&3948466665&8.9.19.47.337.1123&4607&5500&12151&3962908125&4.9.5.11.169.103.577&629&6976 \cr
 & &64.125.11.17.271.337&8425&8672& & &512.17.37.103.109&11227&4352 \cr
\noalign{\hrule}
 & &49.121.37.41.439&3775&702& & &5.49.11.17.19.29.157&1103&318 \cr
12134&3948494627&4.27.25.7.13.41.151&4477&5338&12152&3963318205&4.3.11.17.19.53.1103&1225&2328 \cr
 & &16.9.5.121.17.37.157&785&1224& & &64.9.25.49.53.97&4365&1696 \cr
\noalign{\hrule}
 & &7.11.13.113.181.193&5&6& & &27.343.11.13.41.73&589&9850 \cr
12135&3951376429&4.3.5.7.13.113.181.193&2669&19140&12153&3963698739&4.25.19.31.41.197&143&1128 \cr
 & &32.9.25.11.17.29.157&12325&22608& & &64.3.5.11.13.19.47&47&3040 \cr
\noalign{\hrule}
 & &9.5.11.31.43.53.113&745&272& & &9.25.121.41.53.67&5911&3164 \cr
12136&3951751815&32.25.17.31.53.149&609&1034&12154&3963714975&8.3.7.23.53.113.257&779&440 \cr
 & &128.3.7.11.29.47.149&7003&12992& & &128.5.7.11.19.41.257&1799&1216 \cr
\noalign{\hrule}
 & &9.25.23.31.41.601&657&56& & &125.2401.47.281&51129&45254 \cr
12137&3953032425&16.81.25.7.41.73&1573&1748&12155&3963750875&4.9.7.1331.13.17.19.23&199&1530 \cr
 & &128.121.13.19.23.73&8833&15808& & &16.81.5.289.23.199&23409&36616 \cr
\noalign{\hrule}
 & &125.7.223.20261&3807&24068& & &5.7.11.37.53.59.89&1341&1786 \cr
12138&3953427625&8.81.7.11.47.547&51433&51950&12156&3964426235&4.9.7.11.19.37.47.149&5671&2822 \cr
 & &32.3.25.19.1039.2707&51433&49872& & &16.3.17.47.53.83.107&11703&14552 \cr
\noalign{\hrule}
 & &25.7.19.31.89.431&2277&1846& & &9.125.11.13.19.1297&5149&6524 \cr
12139&3953853925&4.9.25.11.13.23.71.89&581&2356&12157&3964442625&8.7.13.361.233.271&3025&498 \cr
 & &32.3.7.13.19.23.31.83&897&1328& & &32.3.25.121.83.233&2563&1328 \cr
\noalign{\hrule}
 & &3.17.841.137.673&589&1430& & &25.19.67.347.359&877&858 \cr
12140&3954593091&4.5.11.13.17.19.31.137&2853&1346&12158&3964535725&4.3.5.11.13.67.359.877&2573&122838 \cr
 & &16.9.5.31.317.673&1585&744& & &16.9.31.59.83.347&2573&4248 \cr
\noalign{\hrule}
 & &11.17.89.251.947&6669&9430& & &3.5.11.19.47.71.379&37937&41274 \cr
12141&3955991171&4.27.5.13.19.23.41.89&4735&5678&12159&3964906605&4.27.5.59.643.2293&14413&2948 \cr
 & &16.3.25.17.19.167.947&1425&1336& & &32.7.11.29.59.67.71&3953&3248 \cr
\noalign{\hrule}
 & &9.59.107.179.389&62345&62524& & &3.5.19.41.137.2477&20723&26340 \cr
12142&3956224527&8.3.5.49.11.29.37.59.337&35263&84802&12160&3965293065&8.9.25.17.23.53.439&451&26 \cr
 & &32.29.109.179.197.389&3161&3152& & &32.11.13.23.41.439&10097&2288 \cr
\noalign{\hrule}
 & &81.7.11.61.101.103&497&394& & &27.5.7.61.89.773&71179&33176 \cr
12143&3957894171&4.49.61.71.101.197&309&4640&12161&3965803065&16.11.13.17.29.53.79&125&96 \cr
 & &256.3.5.29.103.197&985&3712& & &1024.3.125.11.53.79&21725&27136 \cr
\noalign{\hrule}
 & &11.19.29.37.127.139&64139&12450& & &3.25.7.31.251.971&2771&4026 \cr
12144&3958808821&4.3.25.31.83.2069&957&1112&12162&3966559275&4.9.5.11.17.31.61.163&1127&1942 \cr
 & &64.9.5.11.29.83.139&415&288& & &16.49.17.23.61.971&1037&1288 \cr
\noalign{\hrule}
 & &11.13.17.97.103.163&155&1494& & &3.4913.19.31.457&4251&4708 \cr
12145&3958963723&4.9.5.11.31.83.163&2183&390&12163&3967340847&8.9.11.13.17.19.107.109&15995&15982 \cr
 & &16.27.25.13.37.59&675&17464& & &32.5.7.61.107.109.131.457&228445&228464 \cr
\noalign{\hrule}
 & &5.7.31.59.127.487&35497&26352& & &81.5.11.17.61.859&3703&592 \cr
12146&3959263735&32.27.49.11.61.461&10915&15986&12164&3968438265&32.27.7.11.529.37&413&116 \cr
 & &128.3.5.37.59.7993&7993&7104& & &256.49.29.37.59&63307&6272 \cr
\noalign{\hrule}
 & &29.53.71.131.277&891&1168& & &27.121.13.41.43.53&4343&3236 \cr
12147&3959891449&32.81.11.53.73.131&7085&142&12165&3968447769&8.11.1849.101.809&369&1480 \cr
 & &128.9.5.13.71.109&545&7488& & &128.9.5.37.41.809&4045&2368 \cr
\noalign{\hrule}
 & &3.11.17.19.23.29.557&2845&3282& & &49.17.37.131.983&65283&17116 \cr
12148&3960021021&4.9.5.17.29.547.569&851&3586&12166&3968912633&8.3.11.47.389.463&255&208 \cr
 & &16.11.23.37.163.569&6031&4552& & &256.9.5.11.13.17.389&21395&14976 \cr
\noalign{\hrule}
 & &121.29.31.79.461&1065&1384& & &9.11.13.23.71.1889&225&21004 \cr
12149&3961622401&16.3.5.11.71.173.461&2449&2622&12167&3970056519&8.81.25.59.89&3575&3634 \cr
 & &64.9.5.19.23.31.71.79&6745&6624& & &32.625.11.13.23.79&625&1264 \cr
\noalign{\hrule}
 & &11.19.31.743.823&783&40& & &11.29.43.137.2113&303&1810 \cr
12150&3961837231&16.27.5.11.19.29.31&4859&5030&12168&3970810877&4.3.5.29.43.101.181&2099&2244 \cr
 & &64.3.25.43.113.503&64887&90400& & &32.9.11.17.181.2099&27693&33584 \cr
\noalign{\hrule}
}%
}
$$
\eject
\vglue -23 pt
\noindent\hskip 1 in\hbox to 6.5 in{\ 12169 -- 12204 \hfill\fbd 3972079125 -- 4005175265\frb}
\vskip -9 pt
$$
\vbox{
\nointerlineskip
\halign{\strut
    \vrule \ \ \hfil \frb #\ 
   &\vrule \hfil \ \ \fbb #\frb\ 
   &\vrule \hfil \ \ \frb #\ \hfil
   &\vrule \hfil \ \ \frb #\ 
   &\vrule \hfil \ \ \frb #\ \ \vrule \hskip 2 pt
   &\vrule \ \ \hfil \frb #\ 
   &\vrule \hfil \ \ \fbb #\frb\ 
   &\vrule \hfil \ \ \frb #\ \hfil
   &\vrule \hfil \ \ \frb #\ 
   &\vrule \hfil \ \ \frb #\ \vrule \cr%
\noalign{\hrule}
 & &9.125.7.59.83.103&571&3146& & &3.5.49.59.251.367&803&436 \cr
12169&3972079125&4.5.121.13.83.571&399&316&12187&3994653705&8.5.7.11.73.109.251&2183&3438 \cr
 & &32.3.7.11.19.79.571&16511&9136& & &32.9.37.59.109.191&7067&5232 \cr
\noalign{\hrule}
 & &81.7.11.13.19.2579&65041&16040& & &19.31.67.137.739&2431&11610 \cr
12170&3973050081&16.5.193.337.401&265&72&12188&3995352509&4.27.5.11.13.17.31.43&10541&7906 \cr
 & &256.9.25.53.401&21253&3200& & &16.9.59.67.83.127&4897&9144 \cr
\noalign{\hrule}
 & &3.5.7.13.31.37.2539&11&206& & &11.13.47.733.811&625&108 \cr
12171&3975198045&4.11.37.103.2539&19691&22230&12189&3995385823&8.27.625.13.811&893&82 \cr
 & &16.9.5.7.13.19.29.97&1843&696& & &32.9.25.19.41.47&9225&304 \cr
\noalign{\hrule}
 & &3.25.7.89.149.571&17893&21890& & &9.49.11.37.113.197&23951&1690 \cr
12172&3975316275&4.125.11.29.199.617&81329&41454&12190&3995560107&4.5.7.169.43.557&591&592 \cr
 & &16.9.49.47.167.487&23547&27272& & &128.3.5.37.43.197.557&2785&2752 \cr
\noalign{\hrule}
 & &27.25.11.13.17.2423&10741&11066& & &3.67.97.223.919&11275&10356 \cr
12173&3975961275&4.3.121.17.23.467.503&101&8450&12191&3995656689&8.9.25.11.41.67.863&1159&5474 \cr
 & &16.25.169.101.467&1313&3736& & &32.5.7.17.19.23.41.61&112217&92720 \cr
\noalign{\hrule}
 & &17.157.431.3457&1807&1650& & &25.19.23.43.67.127&9143&63882 \cr
12174&3976721923&4.3.25.11.13.17.139.431&1011&27004&12192&3997315475&4.27.7.169.41.223&1265&1634 \cr
 & &32.9.5.43.157.337&1935&5392& & &16.3.5.7.11.13.19.23.43&231&104 \cr
\noalign{\hrule}
 & &9.5.121.29.89.283&26273&26362& & &27.7.11.101.137.139&1825&296 \cr
12175&3977153235&4.3.49.13.43.47.269.283&145&704&12193&3998630097&16.9.25.37.73.137&4231&4094 \cr
 & &512.5.49.11.29.47.269&12643&12544& & &64.23.73.89.4231&149431&135392 \cr
\noalign{\hrule}
 & &3.7.11.19.23.961.41&375&214& & &3.25.11.19.47.61.89&4379&196 \cr
12176&3977412747&4.9.125.11.31.41.107&5131&6308&12194&3999680025&8.49.11.19.29.151&2759&1098 \cr
 & &32.125.7.19.83.733&10375&11728& & &32.9.7.31.61.89&21&496 \cr
\noalign{\hrule}
 & &27.5.7.11.29.67.197&101&884& & &3.11.43.79.127.281&6335&2938 \cr
12177&3978904545&8.7.11.13.17.67.101&229&240&12195&4000548387&4.5.7.13.113.127.181&5869&5688 \cr
 & &256.3.5.13.17.101.229&23129&28288& & &64.9.5.79.113.5869&17607&18080 \cr
\noalign{\hrule}
 & &17.361.199.3259&1439&4698& & &27.5.11.17.361.439&49849&80534 \cr
12178&3980096117&4.81.29.199.1439&1615&176&12196&4000802355&4.67.79.601.631&2601&44878 \cr
 & &128.9.5.11.17.19.29&261&3520& & &16.9.289.19.1181&1181&136 \cr
\noalign{\hrule}
 & &3.7.11.19.379.2393&13395&15788& & &5.71.89.197.643&171&814 \cr
12179&3980590383&8.9.5.361.47.3947&97291&80324&12197&4002170245&4.9.11.19.37.71.89&2501&2572 \cr
 & &64.17.43.59.97.467&246089&254048& & &32.3.11.37.41.61.643&6771&7216 \cr
\noalign{\hrule}
 & &5.11.17.19.37.73.83&237&1150& & &5.7.13.17.19.113.241&2783&1314 \cr
12180&3982610995&4.3.125.17.23.37.79&399&2524&12198&4002297845&4.9.5.7.121.19.23.73&23377&11882 \cr
 & &32.9.7.19.23.631&14513&1008& & &16.3.13.97.241.457&1371&776 \cr
\noalign{\hrule}
 & &9.11.89.103.4391&6601&2210& & &9.5.13.17.37.73.149&3865&4312 \cr
12181&3984977403&4.5.7.13.17.23.41.103&1255&3624&12199&4002355305&16.3.25.49.11.73.773&247&2072 \cr
 & &64.3.25.13.151.251&37901&10400& & &256.343.11.13.19.37&3773&2432 \cr
\noalign{\hrule}
 & &169.17.41.97.349&6431&498& & &9.5.11.19.149.2857&1401&1456 \cr
12182&3987646429&4.3.59.83.97.109&7735&2838&12200&4003642665&32.27.7.13.19.149.467&57205&4906 \cr
 & &16.9.5.7.11.13.17.43&4257&280& & &128.5.11.17.223.673&11441&14272 \cr
\noalign{\hrule}
 & &27.11.19.37.97.197&4303&3430& & &49.137.349.1709&9405&7696 \cr
12183&3989787219&4.3.5.343.13.197.331&11951&4268&12201&4003908433&32.9.5.11.13.19.37.137&1195&6538 \cr
 & &32.5.7.11.17.19.37.97&85&112& & &128.3.25.7.239.467&17925&29888 \cr
\noalign{\hrule}
 & &3.11.29.47.89.997&2089&2000& & &27.131.563.2011&1595&416 \cr
12184&3991121607&32.125.11.997.2089&261&10706&12202&4004566641&64.3.5.11.13.29.563&197&760 \cr
 & &128.9.25.29.53.101&7575&3392& & &1024.25.13.19.197&64025&9728 \cr
\noalign{\hrule}
 & &27.5.29.37.59.467&7843&9436& & &27.11.13.289.37.97&763&304 \cr
12185&3991189815&8.5.7.11.23.29.31.337&153&8&12203&4004710281&32.7.13.17.19.37.109&757&2610 \cr
 & &128.9.11.17.31.337&5729&21824& & &128.9.5.19.29.757&21953&6080 \cr
\noalign{\hrule}
 & &3.121.23.37.67.193&325&526& & &5.7.13.17.23.47.479&159&320 \cr
12186&3994554003&4.25.121.13.193.263&1269&304&12204&4005175265&128.3.25.13.17.47.53&5269&6444 \cr
 & &128.27.5.19.47.263&40185&16832& & &1024.27.11.179.479&4833&5632 \cr
\noalign{\hrule}
}%
}
$$
\eject
\vglue -23 pt
\noindent\hskip 1 in\hbox to 6.5 in{\ 12205 -- 12240 \hfill\fbd 4006490839 -- 4035369987\frb}
\vskip -9 pt
$$
\vbox{
\nointerlineskip
\halign{\strut
    \vrule \ \ \hfil \frb #\ 
   &\vrule \hfil \ \ \fbb #\frb\ 
   &\vrule \hfil \ \ \frb #\ \hfil
   &\vrule \hfil \ \ \frb #\ 
   &\vrule \hfil \ \ \frb #\ \ \vrule \hskip 2 pt
   &\vrule \ \ \hfil \frb #\ 
   &\vrule \hfil \ \ \fbb #\frb\ 
   &\vrule \hfil \ \ \frb #\ \hfil
   &\vrule \hfil \ \ \frb #\ 
   &\vrule \hfil \ \ \frb #\ \vrule \cr%
\noalign{\hrule}
 & &13.89.967.3581&22793&23760& & &9.61.71.281.367&6919&10222 \cr
12205&4006490839&32.27.5.11.23.89.991&23321&21274&12223&4019787333&4.11.17.19.37.71.269&805&2154 \cr
 & &128.3.121.967.23321&23321&23232& & &16.3.5.7.17.23.37.359&30515&47656 \cr
\noalign{\hrule}
 & &5.41.53.191.1931&9889&234& & &27.25.7.13.29.37.61&79&506 \cr
12206&4007240165&4.9.11.13.29.31.41&1777&1790&12224&4020450525&4.3.5.11.23.29.37.79&479&2074 \cr
 & &16.3.5.11.31.179.1777&55087&47256& & &16.17.61.79.479&8143&632 \cr
\noalign{\hrule}
 & &7.11.29.79.22721&1107&1184& & &7.13.239.401.461&1593&4400 \cr
12207&4008143447&64.27.37.41.22721&9085&13636&12225&4020541889&32.27.25.11.59.239&273&922 \cr
 & &512.9.5.7.23.79.487&11201&11520& & &128.81.5.7.13.461&81&320 \cr
\noalign{\hrule}
 & &3.5.11.61.101.3943&1921&2022& & &25.7.13.361.59.83&6627&1258 \cr
12208&4008315795&4.9.5.11.17.61.113.337&3943&34138&12226&4021783675&4.3.5.17.19.37.2209&539&354 \cr
 & &16.169.17.101.3943&169&136& & &16.9.49.11.17.47.59&1309&3384 \cr
\noalign{\hrule}
 & &5.11.13.61.73.1259&63567&5678& & &9.13.19.37.41.1193&2047&2750 \cr
12209&4008523805&4.9.7.17.167.1009&1597&1430&12227&4023143163&4.125.11.23.89.1193&7585&5538 \cr
 & &16.3.5.7.11.13.17.1597&1597&2856& & &16.3.625.13.37.41.71&625&568 \cr
\noalign{\hrule}
 & &9.7.29.101.103.211&817&110& & &5.11.13.73.127.607&621&14 \cr
12210&4010325291&4.5.11.19.29.43.211&91&120&12228&4023660355&4.27.7.11.13.23.73&607&680 \cr
 & &64.3.25.7.11.13.19.43&6175&15136& & &64.3.5.7.17.23.607&1173&224 \cr
\noalign{\hrule}
 & &9.169.23.29.59.67&623&4510& & &3.25.11.13.17.71.311&13377&16798 \cr
12211&4010346171&4.3.5.7.11.41.67.89&767&1102&12229&4025918325&4.9.343.169.37.227&1555&34 \cr
 & &16.11.13.19.29.41.59&451&152& & &16.5.49.17.37.311&37&392 \cr
\noalign{\hrule}
 & &27.7.23.47.67.293&8525&63902& & &121.17.29.181.373&25191&19942 \cr
12212&4010789979&4.25.11.31.89.359&567&412&12230&4027352989&4.81.169.17.59.311&4925&362 \cr
 & &32.81.5.7.103.359&1795&4944& & &16.3.25.59.181.197&4925&1416 \cr
\noalign{\hrule}
 & &3.5.7.11.289.61.197&1327&15418& & &11.13.23.29.157.269&3285&326 \cr
12213&4011214515&4.13.17.593.1327&553&774&12231&4028225773&4.9.5.13.29.73.163&4235&2116 \cr
 & &16.9.7.43.79.593&3397&14232& & &32.3.25.7.121.529&4025&528 \cr
\noalign{\hrule}
 & &9.7.11.37.97.1613&7055&7462& & &5.49.121.29.43.109&2067&1748 \cr
12214&4011816501&4.5.49.13.17.41.83.97&4071&94&12232&4029437335&8.3.7.11.13.19.23.43.53&6649&930 \cr
 & &16.3.13.23.47.59.83&14053&39176& & &32.9.5.23.31.61.109&1891&3312 \cr
\noalign{\hrule}
 & &3.7.29.31.41.71.73&1275&988& & &3.5.7.13.19.23.29.233&307&858 \cr
12215&4011844137&8.9.25.13.17.19.29.71&43&682&12233&4030584285&4.9.7.11.169.23.307&1793&8854 \cr
 & &32.11.13.17.19.31.43&2717&11696& & &16.121.19.163.233&163&968 \cr
\noalign{\hrule}
 & &5.11.17.43.73.1367&1085&282& & &9.5.19.37.103.1237&21279&2224 \cr
12216&4012097155&4.3.25.7.17.31.43.47&1971&3146&12234&4030646985&32.27.41.139.173&2405&2266 \cr
 & &16.81.121.13.31.73&2511&1144& & &128.5.11.13.37.41.103&533&704 \cr
\noalign{\hrule}
 & &3.7.11.19.23.127.313&991&1930& & &9.125.13.529.521&1337&1268 \cr
12217&4012744197&4.5.7.11.19.193.991&635&828&12235&4030781625&8.3.25.7.13.23.191.317&7279&4796 \cr
 & &32.9.25.23.127.991&991&1200& & &64.11.29.109.251.317&300949&294176 \cr
\noalign{\hrule}
 & &5.11.17.37.43.2699&2003&15498& & &5.49.13.43.59.499&477&22 \cr
12218&4014991915&4.27.7.17.41.2003&925&1078&12236&4032092155&4.9.7.11.43.53.59&481&422 \cr
 & &16.3.25.343.11.37.41&1715&984& & &16.3.11.13.37.53.211&5883&18568 \cr
\noalign{\hrule}
 & &27.7.17.29.71.607&26779&16460& & &5.343.29.89.911&12987&13432 \cr
12219&4015649169&8.9.5.61.439.823&2387&1564&12237&4032464065&16.27.343.13.23.37.73&6281&1822 \cr
 & &64.5.7.11.17.23.31.61&7843&9760& & &64.9.11.23.571.911&6281&6624 \cr
\noalign{\hrule}
 & &9.5.11.13.19.107.307&235&1156& & &9.23.29.31.53.409&3355&6052 \cr
12220&4016282985&8.3.25.11.289.19.47&7693&7982&12238&4033937961&8.3.5.11.17.53.61.89&6923&7228 \cr
 & &32.49.13.47.157.307&2303&2512& & &64.7.11.13.17.23.43.139&56287&57824 \cr
\noalign{\hrule}
 & &3.7.11.13.17.19.41.101&2315&6254& & &3.25.23.31.37.2039&5671&4524 \cr
12221&4016641629&4.5.7.17.53.59.463&2071&1170&12239&4034314425&8.9.5.13.23.29.53.107&31&176 \cr
 & &16.9.25.13.19.59.109&2725&1416& & &256.11.13.31.53.107&7579&13696 \cr
\noalign{\hrule}
 & &25.11.29.37.53.257&3053&3372& & &9.11.23.71.109.229&77663&100334 \cr
12222&4019216575&8.3.37.43.53.71.281&8995&10956&12240&4035369987&4.13.17.37.227.2099&2525&426 \cr
 & &64.9.5.7.11.43.83.257&2709&2656& & &16.3.25.17.37.71.101&2525&5032 \cr
\noalign{\hrule}
}%
}
$$
\eject
\vglue -23 pt
\noindent\hskip 1 in\hbox to 6.5 in{\ 12241 -- 12276 \hfill\fbd 4037041547 -- 4063557775\frb}
\vskip -9 pt
$$
\vbox{
\nointerlineskip
\halign{\strut
    \vrule \ \ \hfil \frb #\ 
   &\vrule \hfil \ \ \fbb #\frb\ 
   &\vrule \hfil \ \ \frb #\ \hfil
   &\vrule \hfil \ \ \frb #\ 
   &\vrule \hfil \ \ \frb #\ \ \vrule \hskip 2 pt
   &\vrule \ \ \hfil \frb #\ 
   &\vrule \hfil \ \ \fbb #\frb\ 
   &\vrule \hfil \ \ \frb #\ \hfil
   &\vrule \hfil \ \ \frb #\ 
   &\vrule \hfil \ \ \frb #\ \vrule \cr%
\noalign{\hrule}
 & &49.11.37.47.59.73&1395&1454& & &27.2401.11.13.19.23&1283&1720 \cr
12241&4037041547&4.9.5.7.31.47.73.727&16211&944&12259&4051104057&16.9.5.343.43.1283&2185&902 \cr
 & &128.3.13.29.31.43.59&11687&8256& & &64.25.11.19.23.41.43&1025&1376 \cr
\noalign{\hrule}
 & &3.23.73.193.4153&4595&442& & &3.49.11.29.103.839&2091&7138 \cr
12242&4037301573&4.5.13.17.193.919&657&3938&12260&4052352381&4.9.17.29.41.43.83&517&3920 \cr
 & &16.9.11.13.73.179&143&4296& & &128.5.49.11.43.47&235&2752 \cr
\noalign{\hrule}
 & &5.11.13.53.61.1747&2331&584& & &81.5.11.13.19.29.127&815&1598 \cr
12243&4038356465&16.9.7.13.37.61.73&1747&954&12261&4052717955&4.3.25.11.13.17.47.163&20447&18328 \cr
 & &64.81.7.53.1747&81&224& & &64.7.17.23.29.79.127&2737&2528 \cr
\noalign{\hrule}
 & &9.41.43.277.919&875&44& & &5.169.19.41.47.131&4411&1746 \cr
12244&4039151121&8.3.125.7.11.41.43&1387&1838&12262&4052876035&4.9.11.13.19.97.401&5371&158 \cr
 & &32.5.7.19.73.919&365&2128& & &16.3.11.41.79.131&869&24 \cr
\noalign{\hrule}
 & &3.49.43.71.9001&2761&6240& & &3.59.139.269.613&4235&3966 \cr
12245&4039567791&64.9.5.11.13.43.251&2527&268&12263&4056960891&4.9.5.7.121.613.661&877&6394 \cr
 & &512.7.11.361.67&24187&2816& & &16.5.7.11.23.139.877&8855&7016 \cr
\noalign{\hrule}
 & &5.7.11.17.23.47.571&71&258& & &7.13.17.19.31.61.73&957&430 \cr
12246&4039907795&4.3.5.23.43.71.571&611&2244&12264&4057497899&4.3.5.7.11.13.29.43.61&589&204 \cr
 & &32.9.11.13.17.43.47&117&688& & &32.9.17.19.29.31.43&261&688 \cr
\noalign{\hrule}
 & &11.17.79.491.557&917&426& & &3.5.17.151.167.631&1599&968 \cr
12247&4040223451&4.3.7.11.71.131.557&24817&14730&12265&4057541385&16.9.5.121.13.41.167&12863&11018 \cr
 & &16.9.5.13.23.83.491&3735&2392& & &64.7.11.19.677.787&164483&151648 \cr
\noalign{\hrule}
 & &9.25.11.13.17.83.89&3243&4322& & &9.5.7.13.23.67.643&52063&51460 \cr
12248&4040504325&4.27.5.11.23.47.2161&1577&12382&12266&4057583985&8.25.11.13.31.83.4733&579&4154 \cr
 & &16.19.23.41.83.151&2869&7544& & &32.3.961.67.83.193&16019&15376 \cr
\noalign{\hrule}
 & &3.11.23.47.193.587&197&390& & &9.11.13.41.131.587&4799&63880 \cr
12249&4041429843&4.9.5.11.13.23.47.197&587&2702&12267&4057623999&16.5.1597.4799&1601&3198 \cr
 & &16.7.193.197.587&197&56& & &64.3.5.13.41.1601&1601&160 \cr
\noalign{\hrule}
 & &31.61.419.5101&91845&66286& & &9.5.49.17.19.41.139&16907&6862 \cr
12250&4041670229&4.9.5.11.13.23.131.157&2341&3782&12268&4058913285&4.11.17.29.47.53.73&5453&3336 \cr
 & &16.3.5.23.31.61.2341&2341&2760& & &64.3.7.19.41.53.139&53&32 \cr
\noalign{\hrule}
 & &3.7.121.19.31.37.73&25&234& & &23.59.61.181.271&13605&2926 \cr
12251&4042448949&4.27.25.11.13.31.73&1691&2494&12269&4060294627&4.3.5.7.11.19.23.907&197987&198372 \cr
 & &16.5.13.19.29.43.89&3827&15080& & &32.9.37.61.271.5351&5351&5328 \cr
\noalign{\hrule}
 & &3.169.23.31.53.211&2765&1122& & &3.361.73.191.269&65219&13840 \cr
12252&4042553853&4.9.5.7.11.17.79.211&14167&2080&12270&4061972361&32.5.49.1331.173&3933&2722 \cr
 & &256.25.13.31.457&457&3200& & &128.9.7.19.23.1361&9527&4416 \cr
\noalign{\hrule}
 & &9.11.13.17.199.929&4345&3416& & &81.11.17.31.41.211&95&2416 \cr
12253&4044792609&16.3.5.7.121.17.61.79&199&1858&12271&4062137607&32.5.17.19.41.151&3231&664 \cr
 & &64.5.61.199.929&61&160& & &512.9.83.359&359&21248 \cr
\noalign{\hrule}
 & &3.5.17.83.131.1459&4257&2798& & &3.25.11.79.157.397&7553&2372 \cr
12254&4045245285&4.27.11.43.131.1399&10625&49532&12272&4062292575&8.7.13.79.83.593&40545&44696 \cr
 & &32.625.7.17.29.61&7625&3248& & &128.9.5.17.37.53.151&33337&28992 \cr
\noalign{\hrule}
 & &9.3125.19.67.113&1591&1534& & &9.5.49.17.29.37.101&817&2112 \cr
12255&4045753125&4.3.13.37.43.59.67.113&385&7052&12273&4062361905&128.27.7.11.17.19.43&1235&74 \cr
 & &32.5.7.11.13.41.1849&41041&29584& & &512.5.13.361.37&4693&256 \cr
\noalign{\hrule}
 & &7.19.73.281.1483&115&396& & &25.13.841.89.167&427&414 \cr
12256&4045963607&8.9.5.11.19.23.1483&2747&1702&12274&4062429475&4.9.25.7.23.61.89.167&16393&14168 \cr
 & &32.3.529.37.41.67&35443&72816& & &64.3.49.11.169.529.97&185367&186208 \cr
\noalign{\hrule}
 & &5.7.13.83.101.1061&67009&21054& & &9.25.11.13.289.19.23&173&11098 \cr
12257&4046935165&4.3.121.29.113.593&1623&4900&12275&4063477275&4.3.11.31.173.179&425&598 \cr
 & &32.9.25.49.11.541&2705&11088& & &16.25.13.17.23.179&179&8 \cr
\noalign{\hrule}
 & &9.49.11.13.149.431&1229&410& & &25.23.73.131.739&10689&64636 \cr
12258&4049842797&4.5.7.41.431.1229&2123&894&12276&4063557775&8.3.7.11.13.113.509&2687&3930 \cr
 & &16.3.5.11.41.149.193&965&328& & &32.9.5.7.131.2687&2687&1008 \cr
\noalign{\hrule}
}%
}
$$
\eject
\vglue -23 pt
\noindent\hskip 1 in\hbox to 6.5 in{\ 12277 -- 12312 \hfill\fbd 4065072363 -- 4097660705\frb}
\vskip -9 pt
$$
\vbox{
\nointerlineskip
\halign{\strut
    \vrule \ \ \hfil \frb #\ 
   &\vrule \hfil \ \ \fbb #\frb\ 
   &\vrule \hfil \ \ \frb #\ \hfil
   &\vrule \hfil \ \ \frb #\ 
   &\vrule \hfil \ \ \frb #\ \ \vrule \hskip 2 pt
   &\vrule \ \ \hfil \frb #\ 
   &\vrule \hfil \ \ \fbb #\frb\ 
   &\vrule \hfil \ \ \frb #\ \hfil
   &\vrule \hfil \ \ \frb #\ 
   &\vrule \hfil \ \ \frb #\ \vrule \cr%
\noalign{\hrule}
 & &9.11.19.643.3361&101&63758& & &11.13.289.23.4297&933&3364 \cr
12277&4065072363&4.71.101.449&3361&3810&12295&4084388737&8.3.17.23.841.311&10725&3572 \cr
 & &16.3.5.127.3361&635&8& & &64.9.25.11.13.19.47&1175&5472 \cr
\noalign{\hrule}
 & &27.5.71.211.2011&47&1964& & &3.5.17.31.43.61.197&1011&26 \cr
12278&4067116785&8.5.47.211.491&11011&12066&12296&4084758555&4.9.13.31.43.337&7385&7106 \cr
 & &32.3.7.121.13.2011&847&208& & &16.5.7.11.13.17.19.211&4009&8008 \cr
\noalign{\hrule}
 & &3.5.11.13.107.113.157&4727&12072& & &27.49.19.23.37.191&335&286 \cr
12279&4071825615&16.9.11.29.163.503&2033&3500&12297&4085793117&4.5.11.13.19.37.67.191&10535&3468 \cr
 & &128.125.7.19.29.107&3325&1856& & &32.3.25.49.13.289.43&7225&8944 \cr
\noalign{\hrule}
 & &9.7.169.47.79.103&3961&880& & &81.7.11.47.53.263&21985&5416 \cr
12280&4071828033&32.3.5.7.11.13.17.233&2903&1738&12298&4086064521&16.9.5.677.4397&17431&13034 \cr
 & &128.121.79.2903&2903&7744& & &64.343.19.17431&17431&29792 \cr
\noalign{\hrule}
 & &11.19.23.29.131.223&387&164& & &3.5.7.11.47.83.907&989&754 \cr
12281&4072379839&8.9.11.23.41.43.131&8003&2370&12299&4086629085&4.11.13.23.29.43.907&17289&3572 \cr
 & &32.27.5.53.79.151&20385&66992& & &32.9.13.17.19.47.113&6441&3536 \cr
\noalign{\hrule}
 & &3.5.7.11.37.167.571&36283&27098& & &9.13.19.31.127.467&373&94 \cr
12282&4075081395&4.7.13.17.797.2791&14949&4588&12300&4087161117&4.13.19.47.127.373&6919&24450 \cr
 & &32.9.11.17.31.37.151&1581&2416& & &16.3.25.11.17.37.163&6031&37400 \cr
\noalign{\hrule}
 & &3.5.13.17.59.67.311&45&3998& & &27.625.11.361.61&121871&103754 \cr
12283&4075404645&4.27.25.17.1999&27199&26774&12301&4087648125&4.7.47.2593.7411&2409&5002 \cr
 & &16.11.59.461.1217&5071&9736& & &16.3.7.11.41.47.61.73&2993&2632 \cr
\noalign{\hrule}
 & &3.7.11.17.89.107.109&166433&153700& & &9.7.361.37.43.113&84487&55220 \cr
12284&4076253489&8.25.29.53.149.1117&6741&1156&12302&4088804769&8.5.11.13.67.97.251&105&38 \cr
 & &64.9.5.7.289.29.107&435&544& & &32.3.25.7.19.97.251&2425&4016 \cr
\noalign{\hrule}
 & &9.5.11.169.17.47.61&14507&9322& & &11.13.17.19.29.43.71&395&954 \cr
12285&4077261045&4.3.11.59.79.89.163&42347&63440&12303&4089435493&4.9.5.11.17.29.53.79&3569&8992 \cr
 & &128.5.13.17.47.53.61&53&64& & &256.3.5.43.83.281&4215&10624 \cr
\noalign{\hrule}
 & &9.7.13.37.239.563&1003&670& & &5.49.11.19.67.1193&139141&153144 \cr
12286&4077480771&4.5.13.17.59.67.563&717&154&12304&4092866855&16.27.61.709.2281&77&2204 \cr
 & &16.3.5.7.11.17.59.239&1003&440& & &128.9.7.11.19.29.61&1769&576 \cr
\noalign{\hrule}
 & &27.5.13.37.181.347&66953&20108& & &3.25.19.43.67.997&9493&9450 \cr
12287&4078372545&8.11.23.41.71.457&1697&3330&12305&4093108725&4.81.625.7.11.67.863&54223&3598 \cr
 & &32.9.5.37.41.1697&1697&656& & &16.49.11.13.43.97.257&36751&38024 \cr
\noalign{\hrule}
 & &25.109.821.1823&4147&4968& & &27.5.7.17.37.71.97&5941&946 \cr
12288&4078461175&16.27.5.11.13.23.29.109&3439&3646&12306&4093667235&4.7.11.13.17.43.457&883&426 \cr
 & &64.3.11.19.29.181.1823&10317&10208& & &16.3.13.43.71.883&883&4472 \cr
\noalign{\hrule}
 & &11.13.157.389.467&669&1058& & &27.13.17.29.41.577&275&302 \cr
12289&4078516013&4.3.13.529.223.467&471&5600&12307&4093678251&4.25.11.13.17.29.41.151&46737&1432 \cr
 & &256.9.25.7.23.157&5175&896& & &64.81.5.179.577&537&160 \cr
\noalign{\hrule}
 & &7.11.169.19.29.569&115373&131874& & &27.13.37.367.859&48517&80300 \cr
12290&4079822747&4.3.31.113.709.1021&865&156&12308&4094190711&8.25.7.11.29.73.239&1351&1278 \cr
 & &32.9.5.13.31.113.173&17515&24912& & &32.9.25.49.29.71.193&86975&89552 \cr
\noalign{\hrule}
 & &11.29.149.85849&19159&66690& & &5.13.31.43.167.283&1107&1064 \cr
12291&4080488819&4.27.5.49.13.17.19.23&293&638&12309&4094929345&16.27.5.7.19.31.41.283&6149&206 \cr
 & &16.9.11.13.17.29.293&221&72& & &64.9.11.13.19.43.103&1881&3296 \cr
\noalign{\hrule}
 & &9.25.47.233.1657&4213&4072& & &9.5.29.47.179.373&2977&1112 \cr
12292&4082806575&16.3.5.11.233.383.509&1157&8&12310&4095153945&16.3.13.139.179.229&7975&16906 \cr
 & &256.11.13.89.509&12727&65152& & &64.25.11.29.79.107&5885&2528 \cr
\noalign{\hrule}
 & &25.17.19.37.79.173&77597&44022& & &121.19.673.2647&65863&15570 \cr
12293&4083357925&4.3.11.13.23.29.47.127&173&126&12311&4095509869&4.9.5.7.9409.173&253&426 \cr
 & &16.27.7.11.29.127.173&8613&7112& & &16.27.5.11.23.71.97&8165&20952 \cr
\noalign{\hrule}
 & &9.125.19.73.2617&301&12784& & &5.17.29.43.67.577&66297&17252 \cr
12294&4083501375&32.25.7.17.43.47&1971&3146&12312&4097660705&8.3.49.11.19.41.227&1105&1392 \cr
 & &128.27.121.13.73&121&2496& & &256.9.5.7.13.17.19.29&1729&1152 \cr
\noalign{\hrule}
}%
}
$$
\eject
\vglue -23 pt
\noindent\hskip 1 in\hbox to 6.5 in{\ 12313 -- 12348 \hfill\fbd 4099256343 -- 4132674077\frb}
\vskip -9 pt
$$
\vbox{
\nointerlineskip
\halign{\strut
    \vrule \ \ \hfil \frb #\ 
   &\vrule \hfil \ \ \fbb #\frb\ 
   &\vrule \hfil \ \ \frb #\ \hfil
   &\vrule \hfil \ \ \frb #\ 
   &\vrule \hfil \ \ \frb #\ \ \vrule \hskip 2 pt
   &\vrule \ \ \hfil \frb #\ 
   &\vrule \hfil \ \ \fbb #\frb\ 
   &\vrule \hfil \ \ \frb #\ \hfil
   &\vrule \hfil \ \ \frb #\ 
   &\vrule \hfil \ \ \frb #\ \vrule \cr%
\noalign{\hrule}
 & &81.7.13.29.127.151&24857&5680& & &25.7.169.23.73.83&23947&23778 \cr
12313&4099256343&32.5.49.53.67.71&11049&7766&12331&4121483275&4.9.49.11.73.311.1321&66419&1690 \cr
 & &128.3.11.29.127.353&353&704& & &16.3.5.11.169.17.3907&3907&4488 \cr
\noalign{\hrule}
 & &9.13.17.19.23.53.89&7073&1990& & &3.5.7.11.23.311.499&1663&108 \cr
12314&4099983381&4.5.11.89.199.643&23141&34086&12332&4122595785&8.81.499.1663&209&290 \cr
 & &16.3.13.19.23.73.317&317&584& & &32.5.11.19.29.1663&1663&8816 \cr
\noalign{\hrule}
 & &125.11.13.61.3761&2033&1728& & &3.13.2809.61.617&725&1342 \cr
12315&4100900375&128.27.25.11.13.19.107&2051&874&12333&4123170987&4.25.11.29.53.3721&1993&1728 \cr
 & &512.3.7.361.23.293&105773&123648& & &512.27.5.11.29.1993&89685&81664 \cr
\noalign{\hrule}
 & &27.5.13.17.23.43.139&7&146& & &3.5.121.41.157.353&323&282 \cr
12316&4101447285&4.3.5.7.13.23.43.73&341&556&12334&4124153715&4.9.17.19.47.157.353&2101&5278 \cr
 & &32.7.11.31.73.139&511&5456& & &16.7.11.13.17.19.29.191&38773&33592 \cr
\noalign{\hrule}
 & &3.5.11.19.193.6779&18009&15886& & &3.49.11.13.29.67.101&497&816 \cr
12317&4101667845&4.81.169.19.23.29.47&3425&17114&12335&4125224103&32.9.343.17.67.71&473&130 \cr
 & &16.25.29.43.137.199&42785&31784& & &128.5.11.13.17.43.71&6035&2752 \cr
\noalign{\hrule}
 & &23.97.101.109.167&4203&14000& & &27.13.19.41.79.191&1639&1990 \cr
12318&4101700193&32.9.125.7.23.467&291&176&12336&4125770181&4.5.11.41.79.149.199&1719&1520 \cr
 & &1024.27.25.7.11.97&7425&3584& & &128.9.25.11.19.149.191&1639&1600 \cr
\noalign{\hrule}
 & &3.5.7.17.19.29.43.97&205&86& & &7.121.43.277.409&201&208 \cr
12319&4102324485&4.25.19.29.41.1849&11817&34408&12337&4126244353&32.3.121.13.43.67.277&43099&24540 \cr
 & &64.9.11.13.17.23.101&3333&9568& & &256.9.5.7.47.131.409&5895&6016 \cr
\noalign{\hrule}
 & &25.169.17.19.31.97&13203&3190& & &3.7.13.17.89.97.103&1025&1012 \cr
12320&4103577725&4.81.125.11.29.163&791&2584&12338&4126772559&8.25.11.17.23.41.89.103&7137&9506 \cr
 & &64.3.7.17.19.29.113&609&3616& & &32.9.25.49.13.41.61.97&4575&4592 \cr
\noalign{\hrule}
 & &5.7.121.31.43.727&27323&28656& & &169.19.53.79.307&2403&3430 \cr
12321&4104100385&32.9.5.11.89.199.307&427&1762&12339&4127448299&4.27.5.343.13.53.89&271&418 \cr
 & &128.3.7.61.307.881&53741&58944& & &16.9.5.7.11.19.89.271&28035&23848 \cr
\noalign{\hrule}
 & &243.5.7.11.13.31.109&353&3026& & &9.49.11.13.29.37.61&901&230 \cr
12322&4109590485&4.5.7.13.17.89.353&109&1656&12340&4127662539&4.3.5.49.17.23.37.53&793&58 \cr
 & &64.9.23.89.109&2047&32& & &16.13.17.29.53.61&17&424 \cr
\noalign{\hrule}
 & &25.19.73.103.1151&13409&15366& & &3.5.7.11.19.79.2381&13113&13078 \cr
12323&4110825275&4.3.11.13.23.53.73.197&135&62&12341&4127832555&4.27.13.19.31.47.79.503&669061&248050 \cr
 & &16.81.5.11.13.23.31.53&55809&62744& & &16.25.121.41.281.2381&2255&2248 \cr
\noalign{\hrule}
 & &27.5.7.109.167.239&1133&370& & &49.19.31.313.457&523&66 \cr
12324&4111238565&4.3.25.11.37.103.239&559&8284&12342&4128306301&4.3.49.11.313.523&14535&11092 \cr
 & &32.11.13.19.43.109&6149&304& & &32.27.5.17.19.47.59&6345&16048 \cr
\noalign{\hrule}
 & &729.5.11.37.47.59&1843&104& & &49.23.37.181.547&10053&16750 \cr
12325&4113787095&16.243.5.13.19.97&12173&10912&12343&4128492893&4.9.125.23.67.1117&481&5104 \cr
 & &1024.7.11.31.37.47&217&512& & &128.3.25.11.13.29.37&10725&1856 \cr
\noalign{\hrule}
 & &9.5.7.31.37.59.193&22825&23018& & &27.5.7.19.23.73.137&1543&1562 \cr
12326&4114180035&4.3.125.11.17.31.83.677&84767&142&12344&4130062965&4.7.11.71.73.137.1543&22275&90364 \cr
 & &16.17.29.37.71.79&2059&10744& & &32.81.25.121.19.29.41&5945&5808 \cr
\noalign{\hrule}
 & &17.19.43.53.5591&869&138& & &25.7.169.73.1913&3069&1156 \cr
12327&4115630147&4.3.11.23.79.5591&2361&3230&12345&4130119175&8.9.7.11.289.31.73&25&536 \cr
 & &16.9.5.17.19.23.787&3935&1656& & &128.3.25.17.31.67&527&12864 \cr
\noalign{\hrule}
 & &9.13.17.257.8053&54529&50160& & &13.29.61.293.613&495&118 \cr
12328&4116476169&32.27.5.11.19.31.1759&731&1028&12346&4130468173&4.9.5.11.59.61.293&2117&1238 \cr
 & &256.5.17.19.31.43.257&4085&3968& & &16.3.29.59.73.619&12921&4952 \cr
\noalign{\hrule}
 & &27.49.17.29.59.107&16027&22340& & &3.5.17.29.31.37.487&6747&7376 \cr
12329&4117584807&8.5.11.17.31.47.1117&13137&5852&12347&4130765655&32.9.5.13.31.173.461&1067&3082 \cr
 & &64.3.7.121.19.29.151&2869&3872& & &128.11.23.67.97.173&127501&142784 \cr
\noalign{\hrule}
 & &11.289.19.79.863&7231&7440& & &7.31.197.277.349&34515&34238 \cr
12330&4117958977&32.3.5.7.17.31.79.1033&11311&20712&12348&4132674077&4.9.5.7.13.17.19.31.53.59&3047&13060 \cr
 & &512.9.5.863.11311&11311&11520& & &32.3.25.11.53.277.653&21549&21200 \cr
\noalign{\hrule}
}%
}
$$
\eject
\vglue -23 pt
\noindent\hskip 1 in\hbox to 6.5 in{\ 12349 -- 12384 \hfill\fbd 4134528125 -- 4163064997\frb}
\vskip -9 pt
$$
\vbox{
\nointerlineskip
\halign{\strut
    \vrule \ \ \hfil \frb #\ 
   &\vrule \hfil \ \ \fbb #\frb\ 
   &\vrule \hfil \ \ \frb #\ \hfil
   &\vrule \hfil \ \ \frb #\ 
   &\vrule \hfil \ \ \frb #\ \ \vrule \hskip 2 pt
   &\vrule \ \ \hfil \frb #\ 
   &\vrule \hfil \ \ \fbb #\frb\ 
   &\vrule \hfil \ \ \frb #\ \hfil
   &\vrule \hfil \ \ \frb #\ 
   &\vrule \hfil \ \ \frb #\ \vrule \cr%
\noalign{\hrule}
 & &3125.49.13.31.67&803&2322& & &27.7.19.23.61.823&869&442 \cr
12349&4134528125&4.27.11.13.43.67.73&427&310&12367&4146416379&4.9.11.13.17.79.823&335&488 \cr
 & &16.3.5.7.31.43.61.73&2623&1752& & &64.5.11.13.61.67.79&9581&12640 \cr
\noalign{\hrule}
 & &3.49.27889.1009&24539&3350& & &3.7.1331.19.73.107&41839&34028 \cr
12350&4136580147&4.25.7.53.67.463&6327&9878&12368&4148179959&8.49.43.47.139.181&5445&3424 \cr
 & &16.9.5.11.19.37.449&38665&10776& & &512.9.5.121.107.139&695&768 \cr
\noalign{\hrule}
 & &9.11.17.23.89.1201&1865&182& & &81.5.37.251.1103&145919&130934 \cr
12351&4137566301&4.5.7.13.373.1201&607&594&12369&4148642205&4.17.41.3559.3851&109197&48694 \cr
 & &16.27.5.7.11.373.607&13055&14568& & &16.9.11.97.251.1103&97&88 \cr
\noalign{\hrule}
 & &3.25.23.53.167.271&23465&19624& & &81.7.169.19.43.53&2555&1738 \cr
12352&4137621225&16.125.11.13.361.223&1431&2806&12370&4149231723&4.5.49.11.169.73.79&905&954 \cr
 & &64.27.13.19.23.53.61&2223&1952& & &16.9.25.53.73.79.181&14299&14600 \cr
\noalign{\hrule}
 & &49.47.151.11903&5423&6480& & &27.5.11.13.17.47.269&2603&1258 \cr
12353&4139303959&32.81.5.7.11.17.29.47&7103&1510&12371&4149242955&4.289.19.37.47.137&725&1014 \cr
 & &128.3.25.151.7103&7103&4800& & &16.3.25.169.19.29.137&7163&5480 \cr
\noalign{\hrule}
 & &29.37.1847.2089&3879&64460& & &27.7.17.23.89.631&1007&3410 \cr
12354&4140044959&8.9.5.11.293.431&1517&638&12372&4150093941&4.5.11.17.19.23.31.53&195&518 \cr
 & &32.3.121.29.37.41&121&1968& & &16.3.25.7.11.13.37.53&14575&3848 \cr
\noalign{\hrule}
 & &27.11.19.47.67.233&7891&7720& & &9.25.13.17.19.23.191&1709&2684 \cr
12355&4140365031&16.3.5.11.13.47.193.607&31&548&12373&4150396575&8.3.11.17.19.61.1709&12415&7288 \cr
 & &128.5.13.31.137.607&94085&113984& & &128.5.11.13.191.911&911&704 \cr
\noalign{\hrule}
 & &9.7.163.191.2111&6655&8122& & &81.5.11.43.53.409&103&4396 \cr
12356&4140471069&4.5.1331.31.131.191&35&156&12374&4152554505&8.5.7.43.103.157&1573&5178 \cr
 & &32.3.25.7.11.13.31.131&18733&12400& & &32.3.121.13.863&11219&176 \cr
\noalign{\hrule}
 & &243.3125.7.19.41&1279&1846& & &3.125.11.13.73.1061&157&1218 \cr
12357&4140871875&4.3.13.19.41.71.1279&125&1474&12375&4153417125&4.9.7.13.29.73.157&7427&4034 \cr
 & &16.125.11.67.1279&1279&5896& & &16.49.1061.2017&2017&392 \cr
\noalign{\hrule}
 & &27.11.13.19.47.1201&11455&12656& & &3.25.11.19.31.83.103&5617&10348 \cr
12358&4140895473&32.5.7.11.13.29.79.113&14429&21618&12376&4154172825&8.5.11.13.41.137.199&30193&10602 \cr
 & &128.9.5.47.307.1201&307&320& & &32.9.19.31.109.277&831&1744 \cr
\noalign{\hrule}
 & &7.17.19.31.37.1597&34551&24538& & &9.5.7.11.857.1399&3699&2300 \cr
12359&4141607099&4.9.7.11.349.12269&9799&2470&12377&4154337495&8.243.125.11.23.137&1399&1274 \cr
 & &16.3.5.11.13.19.41.239&17589&9560& & &32.49.13.23.137.1399&2093&2192 \cr
\noalign{\hrule}
 & &5.7.59.877.2287&1471&594& & &25.49.13.19.31.443&3071&5346 \cr
12360&4141768435&4.27.11.1471.2287&5959&19198&12378&4155262475&4.243.7.11.31.37.83&8417&7270 \cr
 & &16.3.29.59.101.331&8787&2648& & &16.9.5.11.19.443.727&727&792 \cr
\noalign{\hrule}
 & &9.11.29.31.173.269&1381&1040& & &3.5.7.19.43.193.251&1287&470 \cr
12361&4141839537&32.5.13.29.173.1381&621&244&12379&4155682755&4.27.25.11.13.47.193&3857&9068 \cr
 & &256.27.23.61.1381&31763&23424& & &32.7.13.19.29.2267&2267&6032 \cr
\noalign{\hrule}
 & &9.5.7.11.17.31.2269&13319&25254& & &3.5.11.169.29.53.97&23821&25358 \cr
12362&4143318795&4.81.19.23.61.701&739&800&12380&4157346765&4.5.7.11.31.41.83.409&171&2216 \cr
 & &256.25.23.701.739&84985&89728& & &64.9.19.41.83.277&34071&50464 \cr
\noalign{\hrule}
 & &3.49.17.59.157.179&2325&5368& & &3.7.11.43.641.653&1335&1976 \cr
12363&4143534423&16.9.25.11.31.59.61&7&52&12381&4157685609&16.9.5.13.19.89.653&241&412 \cr
 & &128.5.7.11.13.31.61&24583&3520& & &128.5.13.89.103.241&107245&85696 \cr
\noalign{\hrule}
 & &11.13.101.281.1021&1495&1596& & &9.7.13.23.37.47.127&25&116 \cr
12364&4143711143&8.3.5.7.169.19.23.1021&609&4496&12382&4160207961&8.3.25.23.29.37.127&235&616 \cr
 & &256.9.49.19.29.281&4959&6272& & &128.125.7.11.29.47&3625&704 \cr
\noalign{\hrule}
 & &9.5.13.29.41.59.101&3287&3278& & &7.11.23.41.223.257&64961&7650 \cr
12365&4144871835&4.11.19.29.41.59.149.173&166481&88458&12383&4161409021&4.9.25.13.17.19.263&223&1012 \cr
 & &16.3.7.17.19.23.641.1399&611363&610232& & &32.3.5.11.17.23.223&15&272 \cr
\noalign{\hrule}
 & &3.61.97.131.1783&893&890& & &7.13.47.349.2789&3663&874 \cr
12366&4146154323&4.5.19.47.61.89.97.131&2871&8788&12384&4163064997&4.9.7.11.19.23.37.47&845&142 \cr
 & &32.9.5.11.2197.19.29.47&388455&386672& & &16.3.5.11.169.23.71&17963&1560 \cr
\noalign{\hrule}
}%
}
$$
\eject
\vglue -23 pt
\noindent\hskip 1 in\hbox to 6.5 in{\ 12385 -- 12420 \hfill\fbd 4164450165 -- 4200522755\frb}
\vskip -9 pt
$$
\vbox{
\nointerlineskip
\halign{\strut
    \vrule \ \ \hfil \frb #\ 
   &\vrule \hfil \ \ \fbb #\frb\ 
   &\vrule \hfil \ \ \frb #\ \hfil
   &\vrule \hfil \ \ \frb #\ 
   &\vrule \hfil \ \ \frb #\ \ \vrule \hskip 2 pt
   &\vrule \ \ \hfil \frb #\ 
   &\vrule \hfil \ \ \fbb #\frb\ 
   &\vrule \hfil \ \ \frb #\ \hfil
   &\vrule \hfil \ \ \frb #\ 
   &\vrule \hfil \ \ \frb #\ \vrule \cr%
\noalign{\hrule}
 & &243.5.103.107.311&517&410& & &9.49.23.43.53.181&519&470 \cr
12385&4164450165&4.27.25.11.41.47.311&11021&3596&12403&4183977357&4.27.5.47.53.173.181&8671&3784 \cr
 & &32.29.31.41.103.107&899&656& & &64.11.13.23.29.43.173&5017&4576 \cr
\noalign{\hrule}
 & &9.5.7.13.17.19.47.67&547&6644& & &3.11.257.619.797&2209&4600 \cr
12386&4165135065&8.5.11.19.151.547&1407&1462&12404&4184041983&16.25.23.2209.257&619&666 \cr
 & &32.3.7.17.43.67.547&547&688& & &64.9.5.23.37.47.619&4255&4512 \cr
\noalign{\hrule}
 & &11.13.43.53.67.191&725&2826& & &9.5.29.37.79.1097&241&154 \cr
12387&4170503909&4.9.25.13.29.43.157&3685&2438&12405&4184522955&4.3.7.11.37.241.1097&37525&60562 \cr
 & &16.3.125.11.23.53.67&375&184& & &16.25.19.79.107.283&5377&4280 \cr
\noalign{\hrule}
 & &3.11.13.37.71.3701&22743&25370& & &9.25.11.53.61.523&1599&1016 \cr
12388&4170964083&4.27.5.7.11.361.43.59&16523&1000&12406&4184876025&16.27.5.13.41.61.127&9823&7322 \cr
 & &64.625.7.13.31.41&8897&20000& & &64.7.11.13.19.47.523&1729&1504 \cr
\noalign{\hrule}
 & &5.11.19.47.173.491&101931&78854& & &81.5.109.113.839&583&256 \cr
12389&4171975445&4.3.61.89.443.557&11275&38298&12407&4185255015&512.27.5.11.53.113&431&134 \cr
 & &16.9.25.11.13.41.491&533&360& & &2048.53.67.431&22843&68608 \cr
\noalign{\hrule}
 & &9.11.19.29.113.677&85&28& & &37.89.107.11881&170235&182116 \cr
12390&4173053049&8.3.5.7.11.17.29.677&1175&856&12408&4186282231&8.27.5.11.13.97.4139&3961&178 \cr
 & &128.125.7.17.47.107&85493&56000& & &32.9.5.11.17.89.233&3961&7920 \cr
\noalign{\hrule}
 & &9.11.19.29.113.677&1175&856& & &9.5.7.11.29.47.887&779&8762 \cr
12391&4173053049&16.3.25.19.47.107.113&85&28&12409&4189119165&4.5.11.13.19.41.337&141&196 \cr
 & &128.125.7.17.47.107&85493&56000& & &32.3.49.13.19.41.47&533&2128 \cr
\noalign{\hrule}
 & &81.25.13.23.61.113&33&28& & &3.5.7.11.13.23.61.199&867&1322 \cr
12392&4173539175&8.243.5.7.11.13.23.113&5893&22052&12410&4192142955&4.9.289.23.61.661&9751&50072 \cr
 & &64.7.37.71.83.149&74053&98272& & &64.49.11.199.569&569&224 \cr
\noalign{\hrule}
 & &7.17.239.257.571&33&4030& & &5.19.23.29.127.521&38979&36566 \cr
12393&4173631427&4.3.5.11.13.31.257&571&714&12411&4192671955&4.9.23.47.61.71.389&77231&22382 \cr
 & &16.9.7.17.31.571&9&248& & &16.3.7.11.17.361.31.59&34751&31416 \cr
\noalign{\hrule}
 & &9.49.11.311.2767&191&502& & &7.13.289.41.3889&3823&66 \cr
12394&4174464987&4.7.191.251.2767&505&2262&12412&4193349251&4.3.7.11.41.3823&2055&1768 \cr
 & &16.3.5.13.29.101.191&38077&7640& & &64.9.5.11.13.17.137&495&4384 \cr
\noalign{\hrule}
 & &5.49.43.61.73.89&8073&11062& & &9.49.11.17.211.241&2575&76 \cr
12395&4175199595&4.27.13.23.73.5531&4189&1342&12413&4193529417&8.3.25.19.103.211&2483&2792 \cr
 & &16.9.11.23.59.61.71&7029&10856& & &128.13.19.191.349&66659&15808 \cr
\noalign{\hrule}
 & &3.107.223.227.257&64961&6622& & &9.25.19.53.83.223&4337&62 \cr
12396&4176080637&4.7.11.13.19.43.263&4815&15436&12414&4193676675&4.31.223.4337&2057&2280 \cr
 & &32.9.5.17.107.227&17&240& & &64.3.5.121.17.19.31&527&3872 \cr
\noalign{\hrule}
 & &9.125.11.349.967&13949&5246& & &9.7.17.101.137.283&2611&2200 \cr
12397&4176352125&4.25.13.29.37.43.61&349&726&12415&4193897841&16.3.25.49.11.101.373&923&188 \cr
 & &16.3.121.37.61.349&671&296& & &128.5.13.47.71.373&87655&59072 \cr
\noalign{\hrule}
 & &3.25.49.13.19.43.107&671&146& & &27.25.11.23.79.311&7303&472 \cr
12398&4176442725&4.7.11.13.61.73.107&3375&3268&12416&4195770975&16.59.67.79.109&69&10 \cr
 & &32.27.125.11.19.43.61&671&720& & &64.3.5.23.67.109&67&3488 \cr
\noalign{\hrule}
 & &3.5.7.13.23.37.59.61&1061&4896& & &27.7.11.13.23.43.157&1165&606 \cr
12399&4180652385&64.27.17.61.1061&8195&9842&12417&4196563371&4.81.5.101.157.233&17365&1508 \cr
 & &256.5.7.11.19.37.149&1639&2432& & &32.25.13.23.29.151&3775&464 \cr
\noalign{\hrule}
 & &9.5.49.409.4637&6613&1976& & &3.7.11.13.89.113.139&323&10380 \cr
12400&4181855265&16.3.5.7.13.17.19.389&3421&1636&12418&4197962769&8.9.5.13.17.19.173&1529&1412 \cr
 & &128.11.19.311.409&3421&1216& & &64.5.11.19.139.353&1765&608 \cr
\noalign{\hrule}
 & &5.13.31.41.197.257&3049&3306& & &3.5.49.17.19.23.769&865&1634 \cr
12401&4182714835&4.3.13.19.29.197.3049&22675&16962&12419&4198982235&4.25.361.23.43.173&1989&6314 \cr
 & &16.9.25.11.19.257.907&8163&8360& & &16.9.7.11.13.17.41.43&1599&3784 \cr
\noalign{\hrule}
 & &9.7.11.67.251.359&1271&2678& & &5.11.13.19.241.1283&13755&10622 \cr
12402&4183850979&4.3.13.31.41.103.251&643&110&12420&4200522755&4.3.25.7.11.47.113.131&80829&82004 \cr
 & &16.5.11.31.103.643&19933&4120& & &32.27.49.13.19.83.1283&1323&1328 \cr
\noalign{\hrule}
}%
}
$$
\eject
\vglue -23 pt
\noindent\hskip 1 in\hbox to 6.5 in{\ 12421 -- 12456 \hfill\fbd 4201158685 -- 4230217285\frb}
\vskip -9 pt
$$
\vbox{
\nointerlineskip
\halign{\strut
    \vrule \ \ \hfil \frb #\ 
   &\vrule \hfil \ \ \fbb #\frb\ 
   &\vrule \hfil \ \ \frb #\ \hfil
   &\vrule \hfil \ \ \frb #\ 
   &\vrule \hfil \ \ \frb #\ \ \vrule \hskip 2 pt
   &\vrule \ \ \hfil \frb #\ 
   &\vrule \hfil \ \ \fbb #\frb\ 
   &\vrule \hfil \ \ \frb #\ \hfil
   &\vrule \hfil \ \ \frb #\ 
   &\vrule \hfil \ \ \frb #\ \vrule \cr%
\noalign{\hrule}
 & &5.19.41.47.53.433&7631&10122& & &27.121.83.103.151&4795&5248 \cr
12421&4201158685&4.3.5.7.13.19.241.587&13409&2256&12439&4217367033&256.9.5.7.41.103.137&121&806 \cr
 & &128.9.7.11.23.47.53&693&1472& & &1024.7.121.13.31.41&8897&6656 \cr
\noalign{\hrule}
 & &3.5.13.19.23.31.37.43&1243&1724& & &25.19.29.2809.109&5083&65142 \cr
12422&4202889015&8.5.11.19.31.113.431&40549&8154&12440&4217643275&4.9.7.11.13.17.23.47&265&218 \cr
 & &32.27.23.41.43.151&1359&656& & &16.3.5.11.13.17.53.109&143&408 \cr
\noalign{\hrule}
 & &125.11.23.61.2179&157&96& & &5.7.17.19.349.1069&775&6708 \cr
12423&4203563375&64.3.125.157.2179&8723&10902&12441&4217680705&8.3.125.13.19.31.43&1529&2346 \cr
 & &256.9.11.13.23.61.79&1027&1152& & &32.9.11.13.17.23.139&16263&4048 \cr
\noalign{\hrule}
 & &27.5.11.29.41.2381&137&182& & &3.25.7.13.23.97.277&168611&180686 \cr
12424&4204048365&4.3.7.13.41.137.2381&763&6380&12442&4217761275&4.11.43.103.191.1637&1385&252 \cr
 & &32.5.49.11.13.29.109&637&1744& & &32.9.5.7.43.191.277&573&688 \cr
\noalign{\hrule}
 & &27.11.29.37.79.167&1985&148& & &5.121.17.29.79.179&11949&14006 \cr
12425&4204358433&8.5.29.1369.397&29467&28098&12443&4217765365&4.3.7.47.79.149.569&561&8 \cr
 & &32.9.7.79.223.373&2611&3568& & &64.9.11.17.47.149&423&4768 \cr
\noalign{\hrule}
 & &9.23.37.41.59.227&71&298& & &27.7.13.43.89.449&5975&9118 \cr
12426&4205656467&4.23.37.59.71.149&697&660&12444&4221919611&4.25.47.89.97.239&9933&1300 \cr
 & &32.3.5.11.17.41.71.149&13277&11920& & &32.3.625.7.11.13.43&625&176 \cr
\noalign{\hrule}
 & &3.7.19.23.137.3347&78089&101518& & &3.5.11.13.59.61.547&1843&8954 \cr
12427&4208012403&4.11.31.193.229.263&22743&66940&12445&4222760685&4.5.1331.19.37.97&8307&9638 \cr
 & &32.9.5.7.361.3347&57&80& & &16.9.13.19.61.71.79&1349&1896 \cr
\noalign{\hrule}
 & &17.41.487.12401&2061&10340& & &3.11.169.47.71.227&277&230 \cr
12428&4209383039&8.9.5.11.41.47.229&3005&2776&12446&4224571923&4.5.11.23.71.227.277&4563&658 \cr
 & &128.3.25.11.347.601&95425&115392& & &16.27.7.169.47.277&277&504 \cr
\noalign{\hrule}
 & &27.25.13.47.59.173&2849&76& & &9.5.11.29.41.43.167&19&148 \cr
12429&4209621975&8.3.7.11.19.37.173&5395&1006&12447&4226413455&8.3.5.11.19.29.37.41&1469&2426 \cr
 & &32.5.13.83.503&83&8048& & &32.13.37.113.1213&15769&66896 \cr
\noalign{\hrule}
 & &125.169.17.19.617&3773&4248& & &243.89.241.811&165&76 \cr
12430&4210022375&16.9.5.343.11.13.17.59&191&4&12448&4227018777&8.729.5.11.19.811&30167&39088 \cr
 & &128.3.343.59.191&11269&65856& & &256.7.97.311.349&108539&86912 \cr
\noalign{\hrule}
 & &9.7.11.13.29.71.227&6557&9560& & &3.5.37.41.47.59.67&1235&3982 \cr
12431&4210743537&16.3.5.29.79.83.239&845&1562&12449&4227674205&4.25.11.13.19.59.181&647&828 \cr
 & &64.25.11.169.71.79&1027&800& & &32.9.11.13.19.23.647&36879&52624 \cr
\noalign{\hrule}
 & &3.5.11.17.31.59.821&2337&2678& & &9.29.31.263.1987&2849&3112 \cr
12432&4212013245&4.9.13.19.41.103.821&595&226&12450&4228202871&16.3.7.11.29.31.37.389&13&2710 \cr
 & &16.5.7.13.17.19.103.113&13699&11752& & &64.5.11.13.37.271&3523&65120 \cr
\noalign{\hrule}
 & &29.97.421.3557&71995&31158& & &3.5.11.19.29.193.241&193&358 \cr
12433&4212459061&4.27.5.7.121.17.577&1261&1624&12451&4228729395&4.179.37249.241&2945&40194 \cr
 & &64.9.49.13.17.29.97&1989&1568& & &16.9.5.7.11.19.29.31&217&24 \cr
\noalign{\hrule}
 & &11.13.43.47.61.239&78831&44732& & &11.41.67.349.401&15353&11514 \cr
12434&4213374737&8.9.19.53.211.461&1195&188&12452&4228838933&4.3.13.19.41.101.1181&201&980 \cr
 & &64.3.5.47.211.239&633&160& & &32.9.5.49.13.67.101&4949&9360 \cr
\noalign{\hrule}
 & &5.11.169.409.1109&4093&5202& & &5.7.11.13.19.79.563&1891&2454 \cr
12435&4216035395&4.9.289.409.4093&5523&1430&12453&4229540315&4.3.7.13.19.31.61.409&6715&1398 \cr
 & &16.27.5.7.11.13.17.263&3213&2104& & &16.9.5.17.31.79.233&3961&2232 \cr
\noalign{\hrule}
 & &3.5.7.17.83.149.191&2637&4048& & &3.11.13.31.47.67.101&3411&3310 \cr
12436&4216343145&32.27.11.23.149.293&1865&2158&12454&4229733651&4.27.5.31.67.331.379&10309&76 \cr
 & &128.5.11.13.23.83.373&8579&9152& & &32.169.19.61.331&4303&18544 \cr
\noalign{\hrule}
 & &27.5.11.169.53.317&6119&10682& & &5.7.11.19.331.1747&166833&179062 \cr
12437&4216462965&4.5.49.11.29.109.211&2067&1748&12455&4229949955&4.27.13.37.71.97.167&545&11312 \cr
 & &32.3.7.13.19.23.53.211&3059&3376& & &128.9.5.7.13.101.109&12753&6464 \cr
\noalign{\hrule}
 & &5.7.11.29.37.59.173&3053&3348& & &5.343.19.131.991&1353&362 \cr
12438&4216562735&8.27.7.11.29.31.43.71&2249&4450&12456&4230217285&4.3.11.19.41.131.181&12245&11466 \cr
 & &32.9.25.13.43.89.173&5785&6192& & &16.27.5.49.11.13.31.79&9207&8216 \cr
\noalign{\hrule}
}%
}
$$
\eject
\vglue -23 pt
\noindent\hskip 1 in\hbox to 6.5 in{\ 12457 -- 12492 \hfill\fbd 4234083139 -- 4270381731\frb}
\vskip -9 pt
$$
\vbox{
\nointerlineskip
\halign{\strut
    \vrule \ \ \hfil \frb #\ 
   &\vrule \hfil \ \ \fbb #\frb\ 
   &\vrule \hfil \ \ \frb #\ \hfil
   &\vrule \hfil \ \ \frb #\ 
   &\vrule \hfil \ \ \frb #\ \ \vrule \hskip 2 pt
   &\vrule \ \ \hfil \frb #\ 
   &\vrule \hfil \ \ \fbb #\frb\ 
   &\vrule \hfil \ \ \frb #\ \hfil
   &\vrule \hfil \ \ \frb #\ 
   &\vrule \hfil \ \ \frb #\ \vrule \cr%
\noalign{\hrule}
 & &11.13.19.59.61.433&4185&4042& & &9.5.13.43.47.59.61&78319&101926 \cr
12457&4234083139&4.27.5.31.43.47.59.61&323&2214&12475&4255043715&4.11.289.41.113.271&3807&826 \cr
 & &16.729.5.17.19.41.47&32759&29160& & &16.81.7.289.47.59&289&504 \cr
\noalign{\hrule}
 & &7.61.1787.5557&73953&35054& & &3.121.13.17.29.31.59&335&192 \cr
12458&4240263293&4.81.11.17.83.1031&1885&854&12476&4255108143&128.9.5.11.29.59.67&425&106 \cr
 & &16.27.5.7.13.17.29.61&1885&3672& & &512.125.17.53.67&6625&17152 \cr
\noalign{\hrule}
 & &5.11.37.541.3853&56403&86158& & &3.25.17.43.149.521&371&374 \cr
12459&4241902555&4.27.23.1873.2089&2197&4070&12477&4256009925&4.5.7.11.289.43.53.521&137381&684 \cr
 & &16.9.5.11.2197.23.37&2197&1656& & &32.9.7.19.37.47.79&20461&42864 \cr
\noalign{\hrule}
 & &3.11.13.83.283.421&9429&14060& & &5.343.121.73.281&47&558 \cr
12460&4242324801&8.9.5.7.13.19.37.449&283&166&12478&4256755195&4.9.49.31.47.281&1573&730 \cr
 & &32.5.7.19.37.83.283&703&560& & &16.3.5.121.13.31.73&31&312 \cr
\noalign{\hrule}
 & &3.11.23.1129.4951&7553&7300& & &3.49.29.37.137.197&3245&2468 \cr
12461&4242566361&8.25.7.13.73.83.1129&207&5852&12479&4257001959&8.5.7.11.59.137.617&171&788 \cr
 & &64.9.5.49.11.13.19.23&1911&3040& & &64.9.5.11.19.59.197&3245&1824 \cr
\noalign{\hrule}
 & &3.25.13.31.229.613&553&592& & &3.5.23.79.181.863&3157&2294 \cr
12462&4242894825&32.5.7.31.37.79.613&34579&56034&12480&4257312765&4.5.7.11.31.37.41.181&1211&306 \cr
 & &128.9.11.151.229.283&9339&9664& & &16.9.49.11.17.31.173&25823&45672 \cr
\noalign{\hrule}
 & &9.17.29.31.109.283&2429&950& & &27.5.7.13.41.79.107&3859&1276 \cr
12463&4242912309&4.3.25.7.19.283.347&341&58&12481&4257649305&8.3.11.17.29.107.227&287&394 \cr
 & &16.25.11.29.31.347&347&2200& & &32.7.11.17.29.41.197&3349&5104 \cr
\noalign{\hrule}
 & &9.5.7.109.191.647&20653&20108& & &9.5.11.13.23.107.269&569&822 \cr
12464&4243016295&8.11.19.191.457.1087&133&324&12482&4260027915&4.27.5.137.269.569&1177&2522 \cr
 & &64.81.7.11.361.1087&35739&34784& & &16.11.13.97.107.569&569&776 \cr
\noalign{\hrule}
 & &9.5.7.11.23.37.1439&68731&21926& & &25.41.83.89.563&3591&3796 \cr
12465&4243200885&4.13.17.19.311.577&1733&2310&12483&4262853025&8.27.5.7.13.19.73.563&11&574 \cr
 & &16.3.5.7.11.17.19.1733&1733&2584& & &32.3.49.11.19.41.73&10731&3344 \cr
\noalign{\hrule}
 & &3.7.23.961.41.223&295&418& & &27.5.11.13.19.59.197&8651&214 \cr
12466&4243842309&4.5.7.11.19.31.59.223&611&12546&12484&4263258285&4.3.19.41.107.211&2003&2006 \cr
 & &16.9.13.17.19.41.47&799&5928& & &16.17.41.59.107.2003&34051&35096 \cr
\noalign{\hrule}
 & &5.11.19.23.173.1021&1517&1770& & &9.29.43.47.59.137&587&1950 \cr
12467&4245374155&4.3.25.37.41.59.1021&1023&2&12485&4263628923&4.27.25.13.137.587&9541&8954 \cr
 & &16.9.11.31.37.59&19647&248& & &16.5.7.121.13.29.37.47&4235&3848 \cr
\noalign{\hrule}
 & &3.5.19.113.241.547&3029&5176& & &3.125.11.19.41.1327&66007&9632 \cr
12468&4245488535&16.13.233.241.647&1243&1890&12486&4264148625&64.7.43.149.443&1625&1476 \cr
 & &64.27.5.7.11.113.233&2563&2016& & &512.9.125.13.41.43&559&768 \cr
\noalign{\hrule}
 & &25.7.11.13.383.443&5729&5346& & &27.25.13.17.37.773&13699&25724 \cr
12469&4245966725&4.243.7.121.13.17.337&443&404&12487&4266554175&8.9.7.19.59.103.109&34595&20098 \cr
 & &32.81.17.101.337.443&27297&27472& & &32.5.11.13.17.37.773&11&16 \cr
\noalign{\hrule}
 & &81.11.13.23.37.431&443&850& & &9.5.7.13.19.29.31.61&307&242 \cr
12470&4248424323&4.27.25.13.17.23.443&707&9482&12488&4266748395&4.7.121.19.29.31.307&66795&4474 \cr
 & &16.7.11.17.101.431&101&952& & &16.3.5.61.73.2237&2237&584 \cr
\noalign{\hrule}
 & &9.7.13.17.529.577&2491&2270& & &3.11.13.139.163.439&1555&238 \cr
12471&4249759059&4.5.7.47.53.227.577&1083&506&12489&4267015467&4.5.7.13.17.139.311&63&76 \cr
 & &16.3.5.11.361.23.47.53&19133&20680& & &32.9.5.49.17.19.311&45717&25840 \cr
\noalign{\hrule}
 & &5.7.11.361.73.419&221&582& & &3.25.11.19.29.41.229&13&13066 \cr
12472&4251134195&4.3.5.7.13.17.97.419&1387&708&12490&4268004675&4.25.13.47.139&891&916 \cr
 & &32.9.13.17.19.59.73&1989&944& & &32.81.11.47.229&27&752 \cr
\noalign{\hrule}
 & &3.7.11.17.29.107.349&13683&13190& & &3.5.7.37.71.113.137&373&586 \cr
12473&4252732869&4.9.5.107.1319.4561&273691&214336&12491&4270201635&4.5.37.113.293.373&48177&6028 \cr
 & &512.11.17.139.179.197&35263&35584& & &32.9.11.53.101.137&3333&848 \cr
\noalign{\hrule}
 & &243.5.73.199.241&68327&20368& & &9.7.1331.127.401&1139&2470 \cr
12474&4253723505&32.7.19.43.67.227&55&12&12492&4270381731&4.5.7.13.17.19.67.127&363&272 \cr
 & &256.3.5.7.11.19.227&2497&17024& & &128.3.121.289.19.67&5491&4288 \cr
\noalign{\hrule}
}%
}
$$
\eject
\vglue -23 pt
\noindent\hskip 1 in\hbox to 6.5 in{\ 12493 -- 12517 \hfill\fbd 4270935435 -- 4293502983\frb}
\vskip -9 pt
$$
\vbox{
\nointerlineskip
\halign{\strut
    \vrule \ \ \hfil \frb #\ 
   &\vrule \hfil \ \ \fbb #\frb\ 
   &\vrule \hfil \ \ \frb #\ \hfil
   &\vrule \hfil \ \ \frb #\ 
   &\vrule \hfil \ \ \frb #\ \ \vrule \hskip 2 pt
   &\vrule \ \ \hfil \frb #\ 
   &\vrule \hfil \ \ \fbb #\frb\ 
   &\vrule \hfil \ \ \frb #\ \hfil
   &\vrule \hfil \ \ \frb #\ 
   &\vrule \hfil \ \ \frb #\ \vrule \cr%
\noalign{\hrule}
 & &3.5.13.23.61.67.233&7871&8206& & &3.5.121.23.37.47.59&5561&2788 \cr
12493&4270935435&4.11.13.17.61.373.463&9287&4194&12506&4283078745&8.5.17.37.41.67.83&207&3278 \cr
 & &16.9.37.233.251.373&9287&8952& & &32.9.11.23.67.149&201&2384 \cr
\noalign{\hrule}
 & &25.11.2609.5953&129&65354& & &243.23.199.3851&8791&4940 \cr
12494&4271128675&4.3.41.43.797&377&420&12507&4283124561&8.81.5.13.19.59.149&199&1738 \cr
 & &32.9.5.7.13.29.41&533&29232& & &32.5.11.59.79.199&869&4720 \cr
\noalign{\hrule}
 & &5.11.13.19.23.13681&150539&164124& & &81.5.59.389.461&2159&2620 \cr
12495&4274696855&8.9.841.47.97.179&689&152&12508&4285066455&8.25.17.127.131.389&407&18 \cr
 & &128.3.13.19.47.53.97&7473&6208& & &32.9.11.37.127.131&4847&22352 \cr
\noalign{\hrule}
 & &3.25.7.17.29.83.199&209&806& & &11.13.17.53.107.311&3175&65556 \cr
12496&4275012525&4.5.11.13.17.19.31.83&199&216&12509&4287508511&8.27.25.127.607&311&296 \cr
 & &64.27.11.13.19.31.199&5301&4576& & &128.9.5.37.127.311&4699&2880 \cr
\noalign{\hrule}
 & &3.5.121.13.37.59.83&147&268& & &11.29.67.83.2417&237&2180 \cr
12497&4275154455&8.9.49.13.37.59.67&217&550&12510&4287658903&8.3.5.11.79.83.109&4089&4958 \cr
 & &32.25.343.11.31.67&10385&5488& & &32.9.5.29.37.47.67&1665&752 \cr
\noalign{\hrule}
 & &3.25.7.11.53.61.229&2347&11622& & &5.7.17.43.359.467&957&838 \cr
12498&4275561675&4.9.11.13.149.2347&22631&7880&12511&4289402005&4.3.11.29.43.419.467&3965&16116 \cr
 & &64.5.7.53.61.197&197&32& & &32.9.5.11.13.17.61.79&11297&8784 \cr
\noalign{\hrule}
 & &3.121.31.293.1297&21715&18492& & &9.25.11.47.79.467&24909&11984 \cr
12499&4276376313&8.9.5.11.23.43.67.101&149&104&12512&4291578225&32.27.7.361.23.107&1625&1264 \cr
 & &128.13.43.67.101.149&195637&184384& & &1024.125.7.13.23.79&2093&2560 \cr
\noalign{\hrule}
 & &13.17.19.59.61.283&135&902& & &81.13.29.313.449&43765&37928 \cr
12500&4276752883&4.27.5.11.19.41.283&6757&5978&12513&4291578369&16.9.5.11.431.8753&5321&14074 \cr
 & &16.3.49.11.29.61.233&7689&11368& & &64.11.17.31.227.313&7037&5984 \cr
\noalign{\hrule}
 & &3.7.169.29.89.467&1375&1206& & &9.7.13.17.199.1549&4015&632 \cr
12501&4277705523&4.27.125.7.11.67.467&12169&494&12514&4291778673&16.3.5.7.11.13.73.79&629&398 \cr
 & &16.5.11.13.19.43.283&5377&18920& & &64.5.17.37.73.199&365&1184 \cr
\noalign{\hrule}
 & &5.49.11.17.61.1531&93&212& & &121.13.23.313.379&597&976 \cr
12502&4278708665&8.3.7.11.31.53.1531&1275&2806&12515&4291806233&32.3.23.61.199.313&9669&2470 \cr
 & &32.9.25.17.23.31.61&1035&496& & &128.9.5.11.13.19.293&2637&6080 \cr
\noalign{\hrule}
 & &3.13.29.37.293.349&6517&1980& & &27.5.49.13.19.37.71&113&22 \cr
12503&4279148679&8.27.5.343.11.19.37&3205&4528&12516&4292268435&4.7.11.19.37.71.113&1551&6472 \cr
 & &256.25.7.283.641&49525&82048& & &64.3.121.47.809&38023&3872 \cr
\noalign{\hrule}
 & &27.5.13.23.37.47.61&815&266& & &27.7.11.17.29.59.71&1075&1004 \cr
12504&4281878835&4.3.25.7.13.19.37.163&2867&3608&12517&4293502983&8.25.17.29.43.59.251&2989&1278 \cr
 & &64.11.41.47.61.163&1793&1312& & &32.9.25.49.43.61.71&2623&2800 \cr
\noalign{\hrule}
 & &343.17.19.29.31.43&4565&5382& & & & & \cr
12505&4282770373&4.9.5.11.13.17.23.31.83&473&938& & & & & \cr
 & &16.3.7.121.13.23.43.67&8107&7176& & & & & \cr
\noalign{\hrule}
}%
}
$$

\vfill
\eject
}
\bye